\documentclass[12pt]{amsbook}

\usepackage{amsmath}
\usepackage{amssymb}
\usepackage{upref}
\usepackage{mathrsfs}
\usepackage{amsmidx}
\usepackage{graphicx}
\usepackage{amscd,epsfig}
\usepackage{tikz-cd}
 \usepackage{color}

\newcommand{\red }{\color{red}}

\newcommand{\lab}{\label}
\newcommand{\ben}{\begin{enumerate}}
\newcommand{\een}{\end{enumerate}}

\newcommand{\bea}{\begin{eqnarray}}
\newcommand{\ba}{\begin{array}}
\newcommand{\bean}{\begin{eqnarray*}}
\newcommand{\ea}{\end{array}}
\newcommand{\eea}{\end{eqnarray}}
\newcommand{\eean}{\end{eqnarray*}}
\newcommand{\beq}{\begin{equation}}
\newcommand{\eeq}{\end{equation}}
\newcommand{\bthm}{\begin{thm}}
\newcommand{\ethm}{\end{thm}}
\newcommand{\blem}{\begin{lem}}
\newcommand{\elem}{\end{lem}}
\newcommand{\bprop}{\begin{prop}}
\newcommand{\eprop}{\end{prop}}
\newcommand{\bcor}{\begin{cor}}
\newcommand{\ecor}{\end{cor}}
\newcommand{\brem}{\begin{rem}}
\newcommand{\erem}{\end{rem}}
\newcommand{\bdfn}{\begin{dfn}}
\newcommand{\edfn}{\end{dfn}}
\newcommand{\edf}{\end{fact}}
\newcommand{\bdf}{\begin{fact}}
\newcommand{\eobs}{\end{obs}}
\newcommand{\bobs}{\begin{obs}}
\newcommand{\bproper}{\begin{proper}}
\newcommand{\eproper}{\end{proper}}
\newcommand{\bcla}{\begin{cla}}
\newcommand{\ecla}{\end{cla}}
\newcommand{\bex}{\begin{ex}}
\newcommand{\eex}{\end{ex}}

\newcommand{\nl}{\newline}

\alph{enumii} \roman{enumiii}
\newtheorem{thm}{Theorem}[section]
\newtheorem{prop}[thm]{Proposition}
\newtheorem{lem}[thm]{Lemma}

\newtheorem{cor}[thm]{Corollary}
\newtheorem{dfn}[thm]{Definition}
\newtheorem{rem}[thm]{Remark}
\newtheorem{fact}[thm]{Fact}
\newtheorem{obs}[thm]{Observation}
\newtheorem{proper}[thm]{Property}
\newtheorem{cla}[thm]{Claim}
\newtheorem{ex}[thm]{Example}

\def\bpf{\, {\rm {\large Proof. }}} 

\def\bal{\begin{aligned}}  \def\eal{\end{aligned}}
\def\cN{\mathbb N}             \def\D{\mathbb D}   
\def\cR{\mathbb R}
\def\N{\mathbb N}              \def\Z{\mathbb Z}     \def\R{\mathbb R}
\def\Q{\mathbb Q }             \def\C{\mathbb C}     \def\cC{\mathcal C}
\def\oc{\widehat{\C}}
\def\1{1\!\!1}                 \def\ind{{\rm ind}}   \def\cF{\mathcal F} 
\def\mT{\mathbb T}             \def\cB{\mathcal B}   \def\mD{\mathfrak D}
\def\and{\text{ and }}         \def\F{{\mathfrak F}}
\def\tif{\text{ if }}         
          \def\Con{\text{Con}}   \def\cV{\mathcal V}
           \def\cO{\mathcal O}
\def\Comp{\text{{\rm Comp}}}  \def\diam{\text{\rm {diam}}}
\def\dist{\text{{\rm dist}}}   \def\Dist{\text{{\rm Dist}}}
\def\Crit{{\rm Crit}} 
\def\Sing{{\rm Sing}}          \def\Tr{{\rm Tr}} \def\Per{{\rm Per}}

\def\h{{\rm h}}
\def\hmu{\h_\mu}           \def\htop{{\rm h}_{\rm top}}
\def\H{\text{{\rm H}}}     \def\im{\text{{\rm Im}}} \def\re{\text{{\rm re}}}
\def\HD{\text{{\rm HD}}}   \def\DD{\text{{\rm DD}}}
\def\BD{\text{{\rm BD}}}         \def\PD{\text{{\rm PD}}}  
\def\Int{\text{{\rm Int}}} \def\ep{\text{e}}
         \def\P{\text{{\rm P}}}      \def\Id{{\rm Id}}
 \def\Mod{{\rm Mod}}
              \def\const{{\rm const}}     \def\sh{\sharp}
\def\abs{\prec}            \def\supp{{\rm supp}}
\def\A{\mathcal A}                     

\def\L{{\mathcal L}}                \def\Pa{\mathcal P}
                          \def\cI{\mathcal I}
\def\cA{\mathcal A}        \def\cP{\mathcal P}        \def\cS{\mathcal S}
             \def\cY{\mathcal Y}     \def\cU{\mathcal U}

\def\a{\alpha}                \def\b{\beta}             \def\d{\delta}
\def\De{\Delta}               \def\e{\varepsilon}          \def\vep{\varepsilon}
\def\f{\phi}                  \def\vp{\varphi}
\def\g{\gamma}                \def\Ga{\Gamma}           \def\l{\lambda}
\def\La{\Lambda}              \def\om{\omega}           \def\Om{\Omega}
\def\Sg{\Sigma}               \def\sg{\sigma}
               \def\th{\theta}           
\def\ka{\kappa}

\def\bi{\bigcap}              \def\bu{\bigcup}
\def\({\bigl(}                \def\){\bigr)}
\def\lt{\left}                \def\rt{\right}

\def\ld{\ldots}               \def\bd{\partial}         \def\^{\tilde}
\def\-{\hat}
\def\endpf{\hfill{$\endpfsuit$}}
\def\epf{\endpf}

\def\es{\emptyset}            \def\sms{\setminus}
\def\sbt{\subset}             \def\spt{\supset}
\def\gek{\succeq}             \def\lek{\preceq}
\def\eqv{\Leftrightarrow}     \def\lekra{\Longleftrightarrow}
\def\Lra{\Longrightarrow}     \def\imp{\Rightarrow}
\def\lra{\longrightarrow}     \def\lmt{\longmapsto}
\def\comp{\asymp}
\def\upto{\nearrow}           \def\downto{\searrow}
\def\sp{\medskip}             \def\fr{\noindent}        \def\nl{\newline}
\def\ov{\overline}            \def\un{\underline}
               \def\PC(F){{\rm PC}(F)}
\def\ess{{\rm ess}}

\def\om{\omega}

\def\arg{\text{arg}}
\def\Arg{\text{Arg}}
\def\re{{\rm Re}}

\def\endpf{{${\mathbin{\hbox{\vrule height 6.8pt depth0pt width6.8pt  }}}$}}

\def\1{1\!\!1}
\def\PC(f){{\rm PC}(f)} 

\def\cM{\mathcal M}
\def\osc{\text{{\rm osc}}}
\def\du{\bigoplus}

\def\bit{\begin{itmize}}  \def\ben{\begin{enumerate}}
\def\eit{\end{itmize}}    \def\een{\end{enumerate}}

\numberwithin{section}{chapter}
 \numberwithin{equation}{chapter}

\makeindex{(S)}

\makeindex{(N)}

\begin{document}


\title{Ergodic Theory, Geometric Measure Theory, \\ Conformal Measures \\ and \\ the Dynamics of Elliptic Functions}


\author{Janina Kotus and Mariusz Urba\'nski}
\address{Faculty of
Mathematics and Information Science, Warsaw University of
Technology, Warsaw 00-661, Poland} \curraddr{}
\email{J.Kotus@impan.pl}

\address{Department of Mathematics,
University of North Texas, P.O. Box 311430, Denton, TX 76203-1430, USA},
\curraddr{} 
\email{urbanski@unt.edu}
\thanks{}



\maketitle



\setcounter{page}{4}

\tableofcontents


\mainmatter



\section{Introduction}

The ultimate goal of our book is to present a unified approach to the dynamics, ergodic theory, and geometry of elliptic functions from $\C$ to $\oc$. We consider elliptic functions as a most regular class of transcendental meromorphic functions. Poles form an essential feature of such functions but the set of critical values is finite and an elliptic function is ``the same'' on its of its fundamental regions. In a sense this is the class of transcendental meromorphic functions which resembles rational functions most. On the other hand, the differences are huge. We will touch on them in the course of this introduction. In order to comprehensively cover the dynamics and geometry of elliptic functions we make large preparations. This is done in the first two parts of the book: Part 1, ``Ergodic Theory and Measures'' and Part 2,''Geometry and Conformal Measures''. We intend our book to be as self contained as possible and we use essentially all major results of Part~1 and Part~2 in Part~3 and Part~4 dealing with elliptic functions.

This book can be thus treated as a fairly comprehensive account of dynamics, ergodic theory, and fractal geometry of elliptic functions but also as a reference book (with proofs) for many results of geometric measure theory, finite and infinite abstract ergodic theory, Young's towers, measure--theoretic Kolmogorov--Sinai entropy, thermodynamic formalism, geometric function theory (in particular Koebe's Distortion Theorems and Riemann--Hurwitz Formulas), various kinds of conformal measures, conformal graph Directed Markov systems and iterated function systems, classical general theory of elliptic functions, and topological dynamics of transcendental meromorphic functions.

Up to Chapter~\ref{Selected Preliminaries} all the material contained in this book, after being substantially processed, collects, with virtually all proofs, the results that are known and have been published. However, Chapter~\ref{infinite-measure} contains a material on infinite ergodic theory, that up to our best knowledge, is not included, with proofs, in any book. Also, Section~\ref{NiceSetsGeneral} treating nice sets is strongly processed and goes, in many aspects, far beyond the existing knowledge and use of nice sets and nice families of conformal dynamics. 

All the results from Chapter~\ref{first-outlook} throughout the end of the book are actually new although they borrow, use, and apply, a lot from the previous results, methods and techniques. Up to our best knowledge Section~\ref{Simple Examples of Elliptic Functions}, is purely original, providing large classes of a variety of simple examples of dynamically various kinds of elliptic functions. 

\

In Part~1, we first bring up Chapter~\ref{geometric-measure}, with proofs, some basic and fundamental concepts and theorems from abstract and geometric measure theory. These include in particular the three classical covering theorems: $4r$, Besicovitch's, and Vitali's type. We also include a short section on probability theory: conditional expectations and martingales theorems. We devote quite a big space for treating Hausdorff and packing measures. In particular we formulate and prove Converse Frostman's Lemmas which form an indispensable tool for proving that a Hausdorff or packing measure is finite , positive, or infinite. Some of them are frequently called, in particular in fractal geometry literature, mass redistribution principle but these lemmas involve no mass redistribution. We then deal with Hausdorff, packing, and box counting, dimensions of sets and measures, and provide tools to calculate and estimate them.

In the next four chapters of Part~1 we deal with classical ergodic theory, both finite (probability) and, what we want to emphasize, infinite one as well. 

\sp Already in Chapter~\ref{invariant-measure} we deal with both finite and infinite invariant measures. We start with quasi--invariant measures and early on, in the second section of this chapter, we introduce the powerful concept of the first return map. This concept along with the concept of nice sets (see Section~\ref{NiceSetsGeneral}) will form our most fundamental tool in Part~4 of our book devoted to present a refined ergodic theory of elliptic functions. We introduce in this chapter the notions of ergodicity and conservativity (always satisfied for finite invariant measures), and prove Poincare\'e's Recurrence Theorem, Birkhoff's Ergodic Theorem, and Hopf's Ergodic Theorem, the last one pertaining to infinite measures. We also provide a powerful, though perhaps somewhat neglected by ergodic community, tool of proving the existence of invariant $\sg$--finite measures absolutely continuous with respect to a given quasi--invariant measures. It stems from the work of Marco Martens. 

We then prove Bogolubov--Krylov Theorem about the existence of Borel probability invariant measures for continuous dynamical systems acting on compact metrizable topological spaces. 

\sp Chapter~\ref{finite-measure} is devoted to stochastic laws for measurable endomorphisms preserving a probability measure that are finer than mere Birkhoff's Ergodic Theorem. Under appropriate hypotheses we prove the Law of Iterated Logarithm. We then describe another powerful method of ergodic theory, namely L. S. Young towers (recently also frequently called Kakutani's towers) which she developed in \cite{lsy1} and \cite{lsy2}; see also \cite{gouezel3} for further progress. With appropriate assumptions imposed on the first return time Young's  construction yields
the exponential decay of correlations, the Central Limit Theorem, and
the Law of Iterated Logarithm follows too.

\sp In Chapter~\ref{infinite-measure} we deal with refined stochastic laws for dynamical systems preserving an infinite measure. This is primarily Darling--Kac Theorem. We make use of some recent progress on this theorem and related issues, mainly due to Zweim\"muller, Thaler, Theresiu, Melbourne, Gou\"ezel, Bruin, Aaronson, and others, but we do not go into most recent subtleties and developments of this branch of infinite ergodic theory. We do not need them for our applications to elliptic functions. 

\sp In Chapter~\ref{ch8} we provide a classical account of Kolmogorov–-Sinai metric entropy for measure preserving dynamical systems. We prove Shannon--McMillan--Breimann Theorem and, based on Abramov's formula define the concept of Krengel's entropy of a conservative system preserving a (possibly infinite) invariant measure. 

\sp In Chapter~\ref{Thermodynamic Formalism}, the last chapter of Part 1, we collect and prove basics concepts and theorems of classical thermodynamic formalism. This includes topological pressure, variational principle, and equilibrium states. We have written this chapter very carefully, taking care of some inaccuracies which have persisted in expositions of thermodynamic formalism since 1970's and 1980's. We have not dealt with Gibbs states in this chapter. 

We provide lots of examples in this Part 1 of the book. Mainly of invariant measures, ergodic invariant measures, and counter examples of infinite ergodic theory. We do emphasize once more that we treat the latter one we care and detail.  

\

Part~2, Geometry and Conformal Measures, has its first chapter, Chapter~\ref{geometric-theory}, devoted to some selected issues of geometric function theory. Its character is entirely classical, meaning that no dynamics is involved. We deal here at length with extremal lengths and moduli of topological annuli. However, the central theme is about various versions of Koebe's Distortion Theorems. These theorems, proved in early years of 19th century by Koebe and Bieberbach, form a beautiful, elegant, and powerful tool of complex analysis. We prove them carefully and provide their many versions of analytic and geometric character. These theorems form an absolutely indispensable tool for non--expanding holomorphic dynamics and their applications very frequently occur throughout all the book; most notably when dealing with holomorphic inverse branches, conformal measures, and  Hausdorff and packing measures. The Riemann--Hurwitz Formula, appropriate in the context of transcendental meromorphic functions, which we treat at length in the last section of Chapter~\ref{geometric-theory}, is an elegant and probably the best tool to control the topological structure of connected components of inverse images of open connected sets under meromorphic maps. Especially to be sure that such connected components are simply connected. Another reason why we devoted a lot of time to the Riemann--Hurwitz Formula is that in the standard monographs on Riemann surfaces this formula is usually formulated only for compact surfaces and its proofs are somewhat sketchy; we do need, as a matter of fact almost 
exclusively, the non--compact case. Our approach to the Riemann--Hurwitz Formula stems from that of A. Beardon's in \cite{Bea} designed to deal with rational functions of the Riemann sphere. We modify it to fit to our context of transcendental meromorphic functions.

\sp In Chapter~\ref{CMHDIM} we encounter for the first time the beautiful, elegant, and powerful concept of conformal measures due to S. Patterson (see \cite{Pat1} and also \cite{Pat2}) in the context of Fuchsian groups, and D. Sullivan in the context of all Kleinian groups and rational functions (see \cite{Su1}--\cite{Su6}). We however, motivated by \cite {DU1} do not restrict ourselves to conformal dynamical systems only and present a fairly complete account of general conformal measures. The last two sections of this chapter deal with conformal measures in the sense of Sullivan. 

\sp Chapter~\ref{Markov-systems} deals with conformal graph directed Markov systems, its special case of iterated function systems, and thermodynamic formalism of countable alphabet subshifts of finite type, frequently also called topological Markov 
chains. This theory started in the papers \cite {MU1}, \cite{MU5}, and the book \cite{MU2}. It was in \cite {MU1} and \cite {MU2} where the concept of conformal measures due to S. Patterson and D. Sullivan was adapted to the realm of conformal graph directed Markov systems and iterated function systems.
We present some elements of this theory in Chapter~\ref{Markov-systems}, primarily those related to conformal measures and a version of Bowen's Formula for the Hausdorff dimension of limit sets of such systems. In particular we get an almost cost free, effective, lower estimate for the Hausdorff dimension of such limit sets. More about conformal graph directed Markov systems can be found 
in many papers and books; we bring up here some of them: \cite{MU3}--\cite{MU6}, \cite{MPU}, \cite{MSzU}, \cite{U3}, \cite{U4}, \cite{CTU}, \cite{HD_Spectrum}, \cite{CF Subsystems}, and \cite{CU-Porosity}. 
 Afterwards, in Part 3 and especially in Part 4, we apply the techniques developed here to get a quite good, explicit estimate from below of Hausdorff dimensions of Julia sets of all elliptic functions and to explore stochastic properties of invariant versions of conformal measures for parabolic and subexpanding elliptic functions. Getting these stochastic properties is possible for us by combining together several powerful methods which we have already mentioned. Namely, having proved the existence of nice and pre-nice sets, it turns out that the holomorphic inverse branches of the first return map they generate form a conformal iterated function system. So, the whole theory of conformal graph directed Markov systems applies and we also enhance it by L. S. Young's towers techniques developed in \cite{lsy1} and \cite{lsy2}; see also \cite{gouezel3}. 
 
\

In Part~3 we pass to elliptic functions and stay with them till the end of the book. Its first chapter, i.e. Chapter~\ref{elliptic-theory}, interesting on its own, is devoted to present some good account of the classical theory of elliptic functions. Almost no dynamics is involved here. We will actually not need this chapter anywhere in the book except in Chapter~\ref{examples} where we provide a lot of examples of elliptic functions, mainly Weierstrass $\wp$ functions, but not only them. This chapter makes a heavy use of Section~\ref{weierstrassII} of Chapter~\ref{elliptic-theory} which in turn relies on all previous sections of this chapter. We primarily follow here the classical books \cite{Du} and \cite{JS}. We would also like to bring reader's attention to the books \cite{AE} and \cite{Lang-Elliptic}.

\sp In Chapter~\ref{topological-picture} we provide a relatively short and condensed account of topological dynamics of all meromorphic functions with emphasize on Fatou domains including Baker domains that are exclusive for transcendental functions and do not occur for rational functions. We do this for all meromorphic functions and not merely for elliptic ones. In particular, we provide a complete proof of Fatou's classification of Fatou periodic components. We analyze the structure of these components and the structure of their boundaries in greater detail. Particularly, we provide a very detailed qualitative and quantitative description of the local behavior around rationally indifferent periodic points and the structure of corresponding Leau--Fatou flower petals including the Fatou Flower Petal Theorem. 

We also bring up the definitions of Speiser class $\cS$ and Eremenko--Lyubich class $\cB$, which play a seminal role in the recent development of the theory of iteration of transcendental meromorphic functions, and we prove some structural theorems about their Fatou components. 

In the last section of this chapter, Section~\ref{NiceSetsGeneral}, Nice Sets for Analytic Maps, we introduce and thoroughly study the objects related to the powerful concept of nice sets, which will be our indispensable tool in last part of the book, Part~\ref{EFB}, Elliptic Functions B, leading, along with the techniques of Young's towers, the first return map method, and the theory of conformal graph directed Markov systems to such stochastic laws as
the exponential decay of correlations, the Central Limit Theorem, and
the Law of Iterated Logarithm follows too fo large classes of elliptic functions.  

Up to our best knowledge there is no systematic book account of topological dynamics of transcendental meromorphic functions. Some results, with and without proofs, can be found in \cite{BKL1}--\cite{BKL1} and in \cite{Be}. 
Essentially all results in Chapter~\ref{topological-picture}, the one we are discussing now, of our book are supplied with proofs. 

\sp Through this whole Chapter~\ref{geometry-and-dynamics} we deal with general non--constant elliptic functions, i. e. impose no assumptions on a given non--constant elliptic function. We first give some basic preliminary facts about such functions. The rest of this chapter is devoted to analyze in greater detail fractal properties of any non--constant elliptic function. Following the paper \cite{KU2}, by associating to a given elliptic function an infinite alphabet conformal iterated function systems, and heavily utilizing its $\th$ number, we provide a strong, somewhat surprising, lower bound for the Hausdorff dimension of the Julia sets of all non--constant elliptic functions. In particular, this estimate shows that the Hausdorff dimension of the Julia sets of any non--constant elliptic function is strictly larger than $1$. We also provide a simple closed formula for the Hausdorff dimension of the set of points escaping to infinity of an elliptic function.
In the last section of this chapter we prove that no conformal measure of an elliptic function charges the set of escaping points.

\sp In Chapter~\ref{first-outlook} we define the class of non--recurrent and, more notably, the  class of compactly non--recurrent elliptic functions. This is the class of elliptic functions which will be dealt with by us since now through the end of the book in greatest detail. Its history goes back to the papers \cite{U1} and \cite{U2} by the second named author, and \cite {KU3}. One should also mention the paper \cite{CJY}. Similarly, as in all these papers, our treatment of non--recurrent elliptic functions is based on, in fact possible at all, due to an appropriate version of the breakthrough Ma\~ne's Theorem proven in \cite{M1} in the context of rational functions. In our setting of elliptic functions this is Theorem~\ref{mnt6.3a}. The first section of the current chapter is entirely devoted to proving this theorem, its first most fundamental consequences, and some other results surrounding it. The next two sections of this chapter, also relying on Ma\~ne's Theorem, provide us with further refined technical tools to study the structure of Julia sets and holomorphic inverse branches. 

The last section of this chapter has somewhat different character. It systematically defines and describes various subclasses of, mainly compactly non--recurrent, elliptic functions we will be dealing with in Part~4 of the book. Mostly, these classes of elliptic functions will be defined in terms of how strongly expanding these
functions are. We would like to add that while in the theory of rational functions such classes pop up in a natural and fairly obvious way, and for example metric and topological definitions of expanding rational functions describe the same class of functions, in the theory of iteration of transcendental meromorphic functions such classification is by no means obvious, topological and metric analogs of rational function realm concepts do not usually coincide, and the definitions of expanding, hyperbolic, topologically hyperbolic, subhyperbolic, etc, functions vary from author to author. Our definitions seem to us pretty natural and fit well for our purpose of detailed investigation of dynamical and geometric properties of elliptic functions they define.  

\sp The purpose of Chapter~\ref{examples} is to provide examples of elliptic
functions with prescribed properties of the orbits of critical
points (an values). We are primarily focused on constructing examples of various classes of compactly non--recurrent elliptic functions discerned in Section~\ref{DDCoEF}. All these examples are either Weierstrass 
$\wp_\La$ elliptic functions or their modifications. The dynamics of such functions depends heavily on the lattice $\La$ and varies drastically from $\La$ to $\La$. 

The first three sections of this chapter have a preparatory character and respectively describe basic dynamical and geometric properties of all Weierstrass $\wp_\La$ elliptic functions, the ones generated by square lattices, and triangular lattices. 

In Section~\ref{Simple Examples of Elliptic Functions} we provide simple constructions of many classes of elliptic functions discerned in Section~\ref{DDCoEF}. We essentially cover all of them. All these examples stem from 
Weierstrass $\wp$ functions.

We then, starting with Section~\ref{connectedjulia}, provide also some different, interesting by their own, and historically first examples of various kinds of Weierstrass $\wp$ elliptic functions and their modifications. These come from the series of papers \cite{HK1}, \cite{HK2}, \cite{HK3}, \cite{HKK}, \cite{HL} by Jane Hawkins and her collaborators. 

\

Part~4 is entirely covered by getting dynamical, geometric/fractal, and stochastic properties of dynamical systems generated by compactly non--recurrent elliptic functions, primarily subexpanding and parabolic ones.

\sp In Chapter~\ref{conformal-measure} we deal systematically with one of the primary concepts of the book, namely that of (Sullivan's) $h$--conformal (as always $h=\HD(J(f))$) measures for compactly non--recurrent elliptic functions. We will prove their existence for this class of elliptic functions. 

In Section~\ref{Elliptic Regularity} we introduce an important class of regular compactly non--recurrent elliptic functions. For this class of elliptic functions we  prove uniqueness and atomlessness of $h$--conformal measures along with their first basic stochastic properties such as ergodicity and conservativity.  
 


As we have already mentioned in this introduction, conformal measures were first defined and introduced by Samuel Patterson in in his seminal paper \cite{Pat1} (see also \cite{Pat2}) in the context of Fuchsian groups. Dennis Sullivan extended this concept to all Kleinian groups in \cite{Su1}--\cite{Su3}. He then, in the papers \cite{Su4} --\cite{Su6}, defined conformal measures for all rational functions of the Riemann sphere $\oc$; he also proved their existence therein. Both Patterson
and Sullivan came up with conformal measures in order to get an understanding of geometric measures, i.e. Hausdorff and packing ones. Although already Sullivan noticed that there are conformal measures for Klenian groups that are not equal, nor even equivalent to any Hausdorff and packing (generalized) measure, the main purpose to deal with them is to understand Hausdorff and packing. Chapter~\ref{Markov-systems}, Graph Directed Markov Systems, and Part~\ref{EFA}, Elliptic Functions A, and, especially, the current Part~\ref{EFB}, Elliptic Functions B, of our book, provide a good evidence. 

Conformal measures, in the sense of Sullivan have been studied in in the context of rational functions in greater detail in \cite{DU2}, where, in particular, the structure of the set of their exponents was examined and fairly clarified.

Since then conformal measures in the context of rational functions have been studied in numerous research works. We list here only very few of them appearing in the early stages of the development of their theory: \cite{DU LMS}, \cite{DU3}, \cite{DU4}. Subsequently the concept of conformal measures, in the sense of Sullivan, has been extended to countable alphabet iterated functions systems in \cite{MU1} and to conformal graph directed Markov systems in \cite{MU2}. These were treated at length in Chapter~\ref{Markov-systems}, Graph Directed Markov Systems. It was furthermore  extended to some transcendental meromorphic dynamics in \cite{KU1}, \cite{UZ1}, and \cite{MyU1}; see also \cite{UZ2}, \cite{MyU2}, \cite{BKZ1}, and \cite{BKZ1}. Our current construction fits well in this line of development.

Last, the concept of conformal measures found its place also in random dynamics; we cite only \cite{MSU}.

The results of Chapter~\ref{Hausdorff-and-packing} can be viewed from two perspectives. The first one is that we provide a geometrical characterization of the $h$--conformal measure $m_h$, which, with the absence of parabolic points, turns out to be a normalized packing measure, and the second one is that we give a complete description of geometric, Hausdorff and packing, measures of the Julia sets $J(f)$. All of this is contained in Theorem~\ref{t1020907}, which gives a simple clear picture. Due the fact that that the Hausdorff dimension of the Julia set of an elliptic function is strictly larger than $1$, this picture is even simpler an clearer then for non--recurrent rational functions of \cite{U1}, comp. \cite{DU3}.

\sp Throughout this whole Chapter~\ref{invariant-p.s.n.r.} $f:\C\lra\oc$ is assumed to be a compactly non--recurrent regular elliptic function. This chapter is in a sense a core of our book. Taking fruits of what has been done in all previous chapters,  we prove in it the existence and uniqueness, up to a multiplicative constant, of a $\sg$--finite $f$--invariant measure $\mu_h$ equivalent to the $h$-- conformal measure $m_h$. Furthermore, still heavily relying on what has been done in all previous chapters, particularly nice sets, first return map techniques, Young's towers, we provide here a systematic account of ergodic and refined stochastic properties of the dynamical system $(f,\mu_h)$ generated by such subclasses of compactly non--recurrent regular elliptic functions as normal subexanding elliptic functions of finite character and parabolic elliptic functions. 

\sp In Chapter~{Dynamical Rigidity of Compactly Non--Recurrent Regular Elliptic Functions} we deal with the problem of dynamical rigidity of compactly non--recurrent regular elliptic functions.
The issue at stake is whether a weak conjugacy such as Lipschitz one on Julia sets can be promoted to much better one such as affine on the whole complex plane $\C$. This topic has a long history and goes back at least to the seminal paper \cite{Su3} of Dennis Sullivan treating among others the dynamical rigidity of conformal expanding repellers. Its systematical account can be found in \cite{PU2}. Many various, in many contexts, smooth and conformal dynamical rigidity theorems then followed. The literature abounds. 

Our approach in this chapter stems from the original Sullivan's 
article \cite{Su3}. It is also influenced by \cite{PU1}, where the case of tame rational functions has been actually done, and \cite{SU}, where the equivalence of the statements (4) and (1) of Theorem~\ref{t1042706} was established for all tame transcendental meromorphic functions. Being tame meaning that the closure of the postsingular set does not contain the whole Julia set; in particular each non--recurrent elliptic function is tame. We would however like to emphasize that, unlike \cite{SU}, we chose in our book the approach which does not make use of dynamical rigidity results for conformal iterated function
systems proven in \cite{MPU}. Such approach is therefore on the one hand somewhat more economical in tools but on the other hand uses the existence of conformal measures $m_h$ and their invariant versions $\mu_h$ which are in general not available in the general context of tame transcendental meromorphic functions dealt with in \cite{SU}. Our approach in particular  provides us with one more way of proving dynamical rigidity, and it can be possibly used in some other contexts. We would also like to add that doing dynamical rigidity for transcendental functions, like in \cite{SU} or in the current chapter, is substantially more involved than in the case of rational functions of \cite{PU1}; the first much bigger difficulty being caused by infinite degree of transcendental functions.

In Appendix~A, A Quick Review of Some Selected Facts from Complex Analysis of One Complex Variable, we collect for the convenience of the reader many basic and fundamental theorems of complex analysis. We provide no proofs but we give detailed references (arbitrarily chosen) where the proofs can be found. We use these theorems throughout the book without directly referring to them. 

\sp {\bf Acknowledgment} {We are very indebted to Jane Hawkins for producing for us four images of Julia sets of various Weierstrass elliptic $\wp$ functions that are included in Chapter~\ref{examples} of our book. We also thank her very much for fruitful discussions about the dynamics and Julia sets of Weierstrass elliptic $\wp$ functions.} 

\backmatter

\part{Ergodic Theory and Measures}

\chapter{Geometric Measure Theory}\label{geometric-measure}

\section{Measures, Integrals and Measure Spaces}

\sp This section has an introductory character. It collects a minimum of knowledge from abstract measure theory needed in subsequent chapters of the book. Most, commonly well known, theorems are brought up without proofs. A full account of measure theory can be found in many books, for ex. in \cite{Cohen Measure Theory}, \cite{Folland Real Analysis}, \cite{Royden Real Analysis}

\bdfn  
A family $\mathfrak F$ of subsets of a set $X$ is said to be
a $\sg$--algebra\index{(N)}{$\sg$--algebra} if the  following conditions are
satisfied:
\beq\label{FP2.1.1}
 X \in \mathfrak F,
\eeq
\beq\label{FP2.1.2} A \in \mathfrak F\, \Rightarrow A^c,
\in \mathfrak F,
\eeq
\beq\label{FP2.1.3} \{A_i\}_{i=1}^\infty  \sbt \mathfrak F \Longrightarrow
\bigcup_{i=1}^\infty A_i \in \mathfrak F.
\eeq
\edfn 

\sp\fr It  follows from this definition that $\es \in \mathfrak F$,
that the $\sg$-algebra  $\mathfrak F$  is closed under  countable
intersections and under subtractions of sets. If (\ref{FP2.1.3}) is
assumed only  for finite subfamilies of $\mathfrak F$ then $\mathcal
F$ is called an  algebra. The elements of the $\sg$-algebra
$\mathfrak F$ are frequently called measurable
sets.\index{(N)}{measurable set}

\bdfn For any family $\mathfrak F$ of subsets  of $X$, we denote by
$\sg(\mathfrak F)$   the least $\sg$-algebra that contains $\F$, and we call  it the $\sg$-algebra generated by $\F$.
\edfn

\bdfn A function on a $ \sg$-algebra $\mathfrak F$, $ \mu: \mathfrak F
\to [0, +\infty]$, is said to  be
$\sg$-additive\index{(N)}{additive function} or countably additive if for any countable
subfamily $\{A_i\}_{i=1}^\infty$ of $\mathfrak F$ consisting of
mutually  disjoint  sets we have
\beq\label{FP2.1.4}
\mu\lt(\bigcup_{i=1}^\infty A_i\rt)= \sum_{i=1}^\infty  \mu(A_i).
\eeq
We  say  then that $\mu$ is a measure.\index{(N)}{measure} 
\edfn

\fr If we  consider in (\ref{FP2.1.4}) only finite  families of sets, we
say $\mu$ is additive.  The two  notions of  additivity and of
$\sg$-additivity  make  sense for a $\sg$-algebra as well as for an
algebra, provided  that in the case of an algebra one considers only
families $\{ A_i\}_{i=1}^\infty  \sbt  \mathfrak F$ such that
$\bigcup_{i=1}^\infty A_i \in  \mathfrak F$. The simplest consequences
of the definition of measure are the following:
\beq\label{FP2.1.5}
 \mu(\es)=0;
\eeq
\beq\label{FP2.1.6}
\mbox {If} \  \,\, A, B \in \mathfrak F\,\, \mbox{and} \ \,\, A \sbt B \
\,\, \mbox{then} \,\, \mu(A)  \leq \mu(B);
\eeq
\beq\label{FP2.1.7}
\mbox {If} \,\, A_1  \sbt  A_2 \sbt \ldots \,\, \mbox{and}
\,\,\{A_i\}_{i=1}^\infty \sbt \mathfrak F \,\, \mbox{then} \,\,
 \mu(\bigcup_{i=1}^\infty A_i)= \sup_{i} \mu(A_i)=\lim_{i \to \infty}
\mu(A_i).
\eeq

\bdfn We say that the triple $(X, \mathfrak F, \mu)$\index{(S)}{$(X,
\mathfrak F,  \mu)$} with a $\sg$-algebra $\mathfrak F$  and $\mu$, a
measure on $\mathfrak F$, is a measure space. If $\mu(X)$ =1, the
triple $(X, \mathfrak F, \mu)$ is called a probability space and
$\mu$ a probability measure.\index{(N)}{probability measure} \edfn

\bdfn
We say that $\varphi: X \to \mathbb R$ is a measurable
function,\index{(N)}{measurable function} if  $\varphi^{-1}(J) \in
\mathfrak F$ for every interval $J \sbt \mathbb R$, equivalently for
every Borel set $J \sbt \mathbb R$. 
\edfn

Throughout the book, for any set $A\sbt X$, we denote by $\1_A$ the characteristic function of the set $A$:
$$
\1_A(x)=
\begin{cases}
1  &\text{ if } x\in A \\
0  &\text{ if } x\notin A.
\end{cases}
$$
A step function is a linear combination of (finitely many) characteristic functions. It is easy to see that any non-negative measurable function $\varphi: X \to \mathbb R$ can be represented as the pointwise limit of a monotone increasing sequence of non--negative step functions, say 
$$
\vp=\lim_{n\to\infty}\vp_n.
$$
The integral of $\vp$ against the measure $\mu$ is then defined as:
$$
\int_X\vp\,d\mu:=\lim_{n\to\infty}\int_X\vp_n\,d\mu.
$$
It is easy to see that this definition is independent of the choice of a sequence $(\vp)_{n=1}^\infty$ of  monotone increasing non--negative step functions. Writing any (not necessarily non-negative) measurable function $\varphi: X \to \mathbb R$ in its canonical form
$$
\vp=\vp_+-\vp_-
$$
where 
$$
\vp_+:=\max\{\vp,0\} \  \text{ and } \vp_-:=-\min\{\vp,0\}
$$
we say that the function $\varphi$ is $\mu$ integrable\index{(N)}{integrable
function} if 
$$
\int_X\vp_+\,d\mu<+\infty \  \text{ and }  \int_X\vp_-\,d\mu<+\infty.
$$
We then define the integral of $\vp$ against the measure $\mu$ to be
$$
\int_X\vp\,d\mu:=\int_X\vp_+\,d\mu - \int_X\vp_-\,d\mu.
$$
The integral of $\vp$ is also frequently denoted by 
$$
\mu(\vp).
$$
Since $|\varphi|=\vp_+-\vp_-$ we see that $\vp$ is integrable if and only if $|\vp|$ is, i.e. if $\int_X|\varphi| d\mu< \infty$. We then write $\varphi \in
L^1(\mu)$. We bring up now two fundmental properties of integrals, the onece that make integrals so powerful and convenient tools.

\

\bthm[Lebesgue Monotone Convergence Theorem]\label{monotone}\index{(N)}{Lebesgue Monotone Convergence Theorem} 
Suppose that $(\vp)_{n=1}^\infty$ is a monotone-increasing sequence of
integrable, real-valued, functions on a probability space $(X,\mathfrak{F},\mu)$. Denote its limit by $\vp$. Then
$$
\int_X\vp\,d\mu= \lim_{n\to \infty}\int_X\vp_n\,d\mu
$$
In particular, the above limit exists. As a matter of fact it is enough to assume only that the sequence $(\vp)_{n=1}^\infty$ is monotone-increasing on a measurable set whose complement is of measure zero.
\ethm

\bthm[LebesgueDominated Convergence Theorem]\label{dominated}\index{(N)}{Lebesgue Dominated Convergence Theorem} 
Suppose that $(\phi_n)_{n=1}^\infty$ is a sequence of measurable, real-valued,
functions on a probability space $(X, \mathfrak{F},\mu)$, that 
$|\phi_n|\leq g $ for an integrable function $g$, and that the sequence  $(\phi_n)_{n=1}^\infty$ converges $\mu$-a.e. to a function $\varphi: X \to \mathbb R$. Then the function $\phi$ is $\mu$-integrable and
$$
\int_X\vp\,d\mu= \lim_{n\to \infty}\int_X\vp_n\,d\mu.
$$
\ethm

\sp More generally than $L^1(\mu)$, for every $1\leq p <\infty$ we write
$$
||\varphi||_p:=\left(\int_X|\varphi|^p d\mu\right)^\frac{1}{p},
$$
and we say that  $\varphi$  belongs to $L^p(\mu)=L^p(X, \mathfrak F,
\mu)$\index{(S)}{$L^p(X, \mathfrak F, \mu)$}. If
$$
\inf_{\mu(E)=0}\{\sup_{X \sms E}| \varphi|\} < \infty,
$$
then we denote the latter
expression by $||\varphi||_\infty$\index{(S)}{$\mid \mid \varphi
\mid  \mid_\infty$}, we say that the function $\varphi$ is essentially bounded, and we write that $\varphi \in L^\infty$\index{(S)}{$L^\infty$}. The numbers $||\varphi||_p$\index{(S)}{$\mid
\mid \varphi \mid \mid_p$}, $ 1\leq p <\infty$ are called
$L^p$-norms of $\varphi$. The vector spaces $ L^p(X, \mathfrak F,
\mu)$ become Banach spaces when endowed with respective norms $|| \cdot||_p$.

\

\bdfn
A measure space $(X,\F,\mu)$ and the measure $\mu$ are called
\begin{itemize}

\sp\item finite if $\mu(X)<+\infty$,

\sp\item probablity if $\mu(X)=1$,

\sp\item infinite if $\mu(X)=+\infty$,

\sp\item $\sg$-finite if the space $X$ can be expressed as a countable union of measurable sets with finite measure $\mu$.
\end{itemize}
\edfn

\fr Given two measures $\mu$ and $\nu$ on the same measurable space $(X,\F)$, we say that $\mu$ is absolutely continuous\index{(N)}{absolutely  continuous measure} with respect to $\nu$ if for any set $A$ in $\F$, $\nu(A)=0$ entails $\mu(A)=0$. The famous Radon--Nikodym Theorem gives the following.

\sp\bthm\label{Radon_Nikodym}
Let $(X,\F)$ be a measurable space, and let $\mu$ and $\nu$ be two $\sg$-finite measures on $(X,\F)$. Then the following statements are equivalent:

\begin{itemize}
\item[(a)] $\mu$ is absolutely continuous with respect to $\nu$ ($\nu(A)=0$ entails $\mu(A)=0$).

\sp\item[(b)] $\forall_{\e>0} \, \exists_{\d>0} \forall_{A\in\F}\ [\nu(A)<\d \, \imp \mu(A)<\e]$.

\sp\item[(c)] There exists a unique (up to sets of measure zero) measurable function $\rho:X\to[0,+\infty)$ such that
$$
\mu(A)=\int_A\rho\,d\nu 
$$
for every $A\in \F$.
\end{itemize}
\ethm

\fr We then write 
$$
\mu \abs \nu
$$ 
in order to indicate that a measure $\mu$ is
absolutely continuous with respect to $\nu$. The unique function $\rho:X\to[0,+\infty)$ apearing in the item (c) is denoted by $\frac{d\mu}{d\nu}$ and is called the Radon-Nikodym derivative of $\mu$ with respect to $\nu$.

\sp\fr We say that two measures $\mu$ and $\nu$ on the same measurable space $(X,\F)$ are equivalent\index{(N)}{equivalent measures} if each one is absolutely continuous with respect to the other. To denote this fact, we frequently write 
$$
\mu\comp\nu.
$$
On the other hand, there is a concept somehow opposite to equivalency or even absolute continuity of measures. Namely, we say that two measures $\mu$ and $\nu$ on $(X,\F)$ are (mutually) singular if there exists a set $Y\in\F$ such that
$$
\mu(X\sms Y)=0 \  \text{ while } \ \nu(Y)=0.
$$
We then write that 
$$
\mu\perp\nu.
$$

\section[Measures on Metric Spaces]{Measures on Metric Spaces: \\ (Metric) Outer Measures and weak$^*$ Convergence}  

In this section we will show how to construct measures starting of with functions of sets that are required to satisfy much weaker conditions than those defining a measure. These are called outer measures. At the end of the section we also deal with weak$^*$ of measures and Riesz Representation Theorem. Again, we refer, for example, to \cite{Cohen Measure Theory}, \cite{Folland Real Analysis}, \cite{Royden Real Analysis} for compplete accounts. 

\bdfn\label{def_outer_measure}
An outer measure on a set $X$\index{(N)}{outer measure}  is a
function $\mu$ defined on all  subsets of $X$ taking values in
$[0,\infty]$ such that
\beq\label{FP8.1.1}
 \mu(\es)=0;
\eeq
\beq\label{FP8.1.2}
\mbox {If}\ \,\,
 A \sbt B \  \,\, \mbox{then} \ \,\, \mu(A)  \leq \mu(B);
\eeq
\beq\label{FP8.1.3}
 \mu(\bigcup_{n=1}^\infty  A_n)\leq \sum_{n=1}^\infty  \mu(A_n)
\eeq
for any  countable family $\{ A_n\}_{n=1}^\infty$ of
subsets of $X$. 
\edfn

\fr A subset $A$ of $X$ is called $\mu$-measurable or simply measurable
with respect to the outer $ \mu$ if and only if
\beq\label{FP8.1.4}
\mu(B) \geq \mu(A\cap B) + \mu(B\sms A)
\eeq
for all  sets $B \sbt  X$. The opposite inequality follows
immediately from (\ref{FP8.1.3}). One can immediately check that if $\mu(A)=0$
then $A$ is $\mu$-measurable.

\bthm\label{FPt8.1.2} If $\mu$ is an outer  measure on $X$, then the
family $\mathfrak F$ of all $\mu$-measurable sets is a $\sg$-algebra,
and restriction of $\mu$  to $\mathfrak F$ is a measure. \ethm

\bpf Obviously $X \in  \mathfrak F$. By symmetry
(\ref{FP8.1.4}), $A \in  \mathfrak F$ if and only if $A^c\in \mathcal
F$. So the conditions (\ref{FP2.1.1}) and (\ref{FP2.1.2}) of the
definition of $\sigma$-algebra are satisfied. To check the condition
(\ref{FP2.1.3}) that $\mathfrak F$ is closed under countable union,
suppose  that $A_1, A_2, \ldots \in \mathfrak F$ and let $B \sbt X$
be any set. Applying (\ref{FP8.1.4}) in  turn to $A_1, A_2, \ldots $
we  get for all $k\geq 1$
$$\begin{aligned}
\mu(B) & \geq  \mu(B \cap A_1)+ \mu(B \sms A_1)\\
      & \geq  \mu(B \cap A_1)+ \mu((B \sms A_1) \cap A_2)+ \mu(B \sms A_1\sms A_2) \\
& \geq \ldots\\
& \geq  \sum_{j=1}^k\mu \left( \left(B \sms \bigcup_{i=1}^{j-1}
A_i\right)\cap A_j\right) + \mu \left(B \sms \bigcup_{j=1}^{k}
 A_j\right)\\
& \geq  \sum_{j=1}^k\mu\left(\left(B \sms \bigcup_{i=1}^{j-1}
A_i\right)\cap A_j\right) + \mu\left(B \sms \bigcup_{j=1}^{\infty}
 A_j\right)
 \end{aligned}
$$
and therefore
\beq\label{FP8.1.5}
\mu(B)\geq  \sum_{j=1}^k\mu\left(\left(B \sms \bigcup_{i=1}^{j-1}
A_i\right)\cap A_j\right) + \mu\left(B \sms \bigcup_{j=1}^{\infty}
 A_j\right).
\eeq
Since
$$B \cap \bigcup_{j=1}^{\infty}
 A_j= \bigcup_{j=1}^{\infty} \left( B \sms
\bigcup_{i=1}^{j-1}
 A_i\right)\cap A_j,$$
using (\ref{FP8.1.3}) we thus get
$$  \mu(B)\geq \mu \left( \bigcup_{j=1}^\infty \left(B \sms \bigcup_{i=1}^{j-1}
A_i\right)\cap A_j\right) + \mu\left(B \sms \bigcup_{j=1}^{\infty}
 A_j\right).
$$
Hence condition (\ref{FP2.1.3}) is also satisfied, and $\mathfrak F$
is a $\sigma$-algebra. To see that $\mu$ is a measure on $\mathcal
F$, that, is that condition (\ref{FP2.1.4}) is satisfied, consider
mutually disjoint sets $A_1, A_2, \ldots \in \mathfrak F$ and apply
(\ref{FP8.1.5}) to $B= \bigcup_{j=1}^\infty A_j$. We get
$$ \mu \left( \bigcup_{j=1}^\infty A_j\right) \geq
\sum_{j=1}^\infty \mu(A_j).$$ Combining  this with (\ref{FP8.1.3})
we conclude that $\mu$ is a measure on $\mathfrak F$. \endpf

\

\bdfn Let $(X, \rho)$\index{(S)}{$(X, \rho)$} be a  metric
space\index{(N)}{measure space}. An outer measure $\mu$ on $X$  is
said to be a metric outer measure\index{(N)}{metric outer measure}
if
\beq\label{FP8.1.6}
\mu(A\cup B)= \mu(A)+ \mu(B)
\eeq
for all positively separated sets\index{(N)}{positively separated
sets}  $A, B \sbt X$, that is, satisfying the following condition
$$\rho(A, B)=\inf\{\rho(x,y): \,\, x \in A, \, \, y \in B\} >0.$$
We assume the convention that $ \rho(A, \es)= \rho(A, \es )=\infty
$.
\edfn

Recall that the Borel  $\sg$-algebra\index{(N)}{Borel $\sg$-algebra}  on $X$ is the
$\sigma$-algebra generated by open, or equivalently closed, sets. We
want to show  that if $\mu$ is a metric  outer measure then the
family of  all $\mu$-measurable  sets contains this $\sg$-algebra.
The  proof is based on the following lemma.

\blem\label{FPl8.1.3}  Let $\mu$ be a metric outer measure on $(X,
\rho)$. Let $\{A_n\}_{n=1}^\infty$ be an ascending  sequence of
subsets of $X$, and denote $A:=\bigcup_{n=1}^\infty A_n$. If
$\rho(A_n, A\sms A_{n+1})>0$ for all $n \geq 1$, then 
$$ 
\mu(A)=\lim_{n \to \infty} \mu(A_n).
$$
\elem

\bpf  By (\ref{FP8.1.2})  it is sufficient to show
\beq\label{FP8.1.7}
\mu(A)\leq \lim_{n \to \infty}\mu(A_n).
\eeq
If $\lim_{n \to \infty}\mu(A_n)=\infty$, there is nothing to  prove.
So, suppose that
\beq\label{FP8.1.8}
 \lim_{n \to \infty}\mu(A_n)=\sup_{n \to \infty}\mu(A_n)< \infty.
\eeq
Let $B_1=A_1$ and $B_n=A_n \sms A_{n-1}$  for $ n \geq 2$. If $ n
\geq m+2$, then $ B_m \sbt A_m$ and $B_n \sbt A \sms A_{n-1} \sbt A
\sms A_{m+1}$. Thus $B_m$ and $B_n$ are positively  separated, and
applying  (\ref{FP8.1.6})  we get for every $j \geq 1$
\beq\label{FP8.1.9}
\mu\left( \bigcup_{i=1}^j B_{2i-1}\right)=\sum_{i=1}^j \mu(B_{2i-1})
\quad \mbox{and} \quad \mu\left( \bigcup_{i=1}^j
B_{2i}\right)=\sum_{i=1}^j \mu(B_{2i}).
\eeq
We also  have for every  $ n \geq 1$
\beq\label{FP8.1.10}
\begin{aligned}
\mu(A)=& \mu\left(\bigcup_{k=n}^\infty A_{k}\right)= \mu\left(A_n
\cup
\bigcup_{k=n+1}^\infty B_{k}\right)\\
&\leq  \mu( A_n) + \sum_{k=n+1}^\infty  \mu( B_k) \\
&\leq \lim_{l \to
\infty} \mu(A_l) + \sum _{k=n+1}^\infty \mu(B_k).
\end{aligned}
\eeq
Since the sets $\bigcup_{i=1}^j B_{2i-1}$ and  $\bigcup_{i=1}^j
B_{2i}$ appearing in  (\ref{FP8.1.9}) are both  contained in
$A_{2j}$, it  follows from (\ref{FP8.1.8}) and (\ref{FP8.1.9}) that
the series $\sum_{k=1}^\infty \mu(B_k)$ converges. Therefor
 (\ref{FP8.1.7}) follows  immediately from (\ref{FP8.1.10}). The proof
is complete.
\endpf

\

\bthm\label{FPt8.1.4} If $\mu$ is a metric outer measure on $(X,
\rho)$, then all Borel subsets of $X$ are $\mu$-measurable. \ethm

\bpf Since the Borel sets form the least $\sigma$-algebra
containing all closed subsets of $X$, it follows from
Theorem~\ref{FPt8.1.2}  that it is enough to check (\ref{FP8.1.4})
for every  non-empty  closed set $A \sbt X$ and every $B \sbt X$.
For al $n \geq 1$, let $B_n=\{ x \in B\sms A: \rho(x, A) \geq
1/n\}$, then $ \rho(B \cap A, B_n) \geq 1/n$ and by (\ref{FP8.1.6})
\beq\label{FP8.1.11}
\mu(B \cap A) + \mu( B_n)= \mu( B \cap A) \cup B_n) \leq \mu(B).
\eeq
The sequence $\{B_n\}_{n=1}^\infty $  is asending and, since $A$
is closed, $ B \sms A=\bigcup_{n=1}^\infty B_n$. In order to apply
Lemma~\ref{FPl8.1.3} we shall now show that
$$ \rho(B_n, (B \sms A)\sms B_{n+1}) > 0$$
for all $n \geq 1$.  And indeed, if $ x \in (B\sms A)\sms B_{n+1}$,
then  there exists $z \in A$  with $\rho(x, z)< 1/(n+1)$. Thus, if
$y \in B_n$, then
$$ 
\rho(x,y)\geq  \rho(y, z)- \rho(x, z) > \frac1n -\frac1{n(n+1)}>0.$$
Applying now Lemma~\ref{FPl8.1.3} with $A_n=B$ shows that $\mu(A\sms
B)=\lim_{n \to  \infty} \mu(B_n)$. Thus (\ref{FP8.1.4})  follows
from (\ref{FP8.1.11}). The proof is complete. 
\endpf

\sp\fr This theorem, as well as many other reasons disseminated all over the mathematics, many of which we will encounter in this book, justifies the following definition.

\bdfn
Any measure on a metric space that is defined on its $\sg$-algebra of Borel sets {\rm(}or larger{\rm)} is called a \index{(N)}{Borel measure} Borel measure.
\edfn

\fr Let us list the following well known properties of finite Borel measures.

\bthm\label{t_regularity_of_Borel_measures}
Any finite Borel measure $\mu$ on a metric space $X$ is both outer and inner regular. Outer regularity means that
$$
\mu(A)=\inf\{\mu(G):G\spt A \text{ and } G \text{ is open}\,\}
$$
while inner regularity means that
$$
\mu(A)=\sup\{\mu(F):F\sbt A \text{ and } F \text{ is closed}\,\}
$$
In addition, if the space $X$ is completely metrizable, then the closed sets involved in the concept of inner regularity can be replaced by compact ones.
\ethm

\

Given a metric space $(X,\rho)$ we denote by $M(X)$ the collection of all Borel probability measures on $X$. We denote by $C(X)$ the vector space of all real-valued continuous functions on $X$ and by $C_b(X)$ its vector subspace consisting of all bounded elements of $C(X)$. Let us record the following easy theorem.

\bthm
If $(X,\rho)$ is a metric space, the two measures $\mu$ and $\nu$ in $M(X)$ are equal if and only if
$$
\nu(g)=\mu(f)
$$
for all functions $g\in C_b(X)$.
\ethm

\fr If $X$ is compact, then $C(X)$ becomes a Banach space if endowed with the supremum metric. Denote by $C^*(X)$ the dual of $C(X)$. Endow $C^*(X)$ with the weak$^*$ topology. This means that 

\sp \centerline{a net $(F_\l)_{\l\in\La}$ in $C^*(X)$ converges to an element $F\in C^*(X)$}

\sp\fr if and only if 

\centerline{the net $(F_\l(g))_{\l\in\La}$ converges to $F(g)$}

\sp\fr for every $g\in C(X)$. $M(X)$, the space of all Borel probability measures on $X$, can be then natuaraly viewed as a subset of $C^*(X)$: every measure $\mu\in M(X)$ induces the functional 
$$
C(X)\ni g\longmapsto \mu(g).
$$
We will frequently use the following.

\bthm\label{Riess}
Let $X$ be a compact metrizable space. Consider $C^*(X)$ with its weak$^*$-topology. Then

\sp\begin{itemize}
\item[(a)] $M(X)$ is a convex compact subset of $C^*(X)$. 

\sp\item[(b)] $M(X)$ is a metrizable space. In particular, proving continuity or convergence one can restrict oneself to sequences only (as opposite to nets).

\sp\item[(c)] (Riesz Representation Theorem) Every $F$, a non-negative element in $C^*(X)$, such that $F(\1)=1$ is (uniquely) represented by an element in $M(X)$, precisely meaning that there exists $\mu\in M(X)$ such that
$$
F(g)=\mu(g)
$$
for all $g\in C(X)$.
\end{itemize}
\ethm

\sp\section[Covering Theorems] {Covering Theorems: $4r$, Besicovitch's, and Vitali's Type; \\ Lebesgue Density Theorem}

In this section we prove first the $4r$--Covering Theorem. Following the
arguments of \cite{MSzU} we prove it for all metric spaces. If we do not care of $4r$ and $5r$ suffices (almost always) a shorter less involved proof is possible. It can be found for example in \cite{Heinonen LAMS}. Then, following \cite{Mat}we will prove Besicovitch's Covering Theorem, and, as its fairly straightforward consequence, Vitali's Type Covering Theorem. We finally deduce from the latter The Lebesgue Density Points Theorem. All these theorems are classical and can be found in many sources with extended discussions. More applications of covering theorems will appear in further
sections of this chapter and throughout the entire book. For every
ball $B:=B(x, r)$, we put $r(B)=r$ and $c(B)=x$. 

\bthm\label{t1hg47} {\rm($4r$-Covering
Theorem).}\index{(N)}{$4r$-Covering Theorem}  Suppose $(X, \rho)$ is
a metric space and $\mathcal B$ is a family of open balls in $X$
such that $\sup\{ r(B): B\in\mathcal B\}<+\infty.$ Then there is a
family $\mathcal B'\subseteq\mathcal B$ consisting of mutually
disjoint balls such that $\bigcup_{B\in\mathcal B}
B\subseteq\bigcup_{B\in\mathcal B'} 4B.$ In addition, if the metric
space $X$ is separable, then $\mathcal B'$ is countable. 
\ethm

\bpf Fix an arbitrary $M>0$. Suppose that there is given
a family $\mathcal B'_M\subseteq\mathcal B$ consisting of mutually
disjoint balls such that
\begin{itemize}
\item[(a)] $r(B)>M$ for all $B\in\mathcal B'_M$,

\, \item[(b)] $\bigcup_{B\in\mathcal B'_M} 4B\supseteq\bigcup\{B:
  B\in\mathcal B\,\,\text{and}\,\,r(B)>M\}.$ 
\end{itemize}
We shall show that then there exists a family $\mathcal
B''_M\subseteq\mathcal B$ with the following properties:
\begin{itemize}
\item[(c)] $\mathcal B''_M\subseteq\mathfrak F:=\{B\in\mathcal B: 3M/4<r(B)\le M\},$

\,\item[(d)] $\mathcal B'_M\cup \mathcal B''_M$ consists of mutually disjoint balls,

\,\item[(e)] $\bigcup_{B\in \mathcal B'_M\cup \mathcal B''_M}4B\supseteq\bigcup\{B: B\in\mathcal B\,\,\text{and}\,\, r(B)>3M/4\}$.
\end{itemize}
Indeed, put
\beq\label{1hg49}
\mathcal B'''_M=\{B\in\mathfrak F: B\cap\bigcup_{D\in\mathcal B'_M}
D=\emptyset\}.
\eeq
Consider $B\in\mathfrak F\setminus\mathcal B'''_M$. Then there exists
$D\in\mathcal B'_M$ such that $B\cap D\neq\emptyset$. Hence,
$r(B)\le M<r(D)$, and in consequence,
$$
\rho(c(B), c(D))<r(B)+r(D)\le M+r(D)<r(D)+r(D)=2r(D),
$$
and
$$
B\subseteq B(c(D), r(B)+2r(D))\subseteq B(c(D), 3r(D))=3D\subseteq
4D.
$$
Therefore,
\beq\label{2hg49}
\bigcup_{B\in\mathfrak F\setminus\mathcal B'''_M}
B\subseteq\bigcup_{B\in\mathcal B'_M} 4B.
\eeq
So, if $\mathcal B'''_M=\emptyset$ we are done with the proof by
setting $\mathcal B''_M=\emptyset$. Otherwise, fix an arbitrary
$B_0\in\mathcal B'''_M$, and further, proceeding by transfinite
induction, fix some
$\mathcal B_{\alpha}\in\mathcal B'''_M$ such that
$$
c(B_{\alpha})\in c(\mathcal
B'''_M)\setminus\bigcup_{\gamma<\alpha}\frac{8}{3} B_{\gamma}
$$
for some some ordinal number $\gamma\ge 0$, as long as the difference on the
right-hand side above is nonempty. 
This  procedure terminates at some ordinal number $\lambda$. First,
we claim that the balls $(B_{\alpha})_{\alpha<\lambda}$ are mutually
disjoint. Indeed, fix $0\le\alpha<\beta<\lambda$. Then
$c(B_{\beta})\notin\frac{8}{3} B_{\alpha}$. So, 
$$
\rho(c(B_{\beta}),
c(B_{\alpha}))\ge\frac{8}{3}
r(B_{\alpha})>\frac{8}{3}\cdot\frac{3}{4}M=2M
$$ 
and
$$
r(B_{\beta})+r(B_{\alpha})\le M+M=2M.
$$ 
Thus $B_{\beta}\cap
B_{\alpha}=\emptyset.$ Now, if $B\in\mathcal B'_M$ and
$0\le\alpha<\lambda$, then $B_{\alpha}\in\mathcal B'''_M$, and by
(\ref{1hg49}), $B_{\alpha}\cap B=\emptyset$. Thus we proved item (d)
with $\mathcal B''_M=\{B_{\alpha}\}_{\alpha<\lambda}$. Item (c) is
obvious since $B_{\alpha}\in\mathcal B'''_M\subseteq\mathfrak F$ for
all $0\le\alpha<\lambda$. It remains to prove item (e). By the
definition of $\lambda$, $c(\mathcal
B'''_M)\subset\bigcup_{\gamma<\lambda}\frac{8}{3}
B_{\gamma}=\bigcup_{B\in\mathcal B''_M}\frac{8}{3} B$. Hence, if
$x\in B$ and $B\in\mathcal B'''_M,$ then there exists $D\in\mathcal
B''_M$ such that $c(B)\in\frac{8}{3} D$. Therefore,
$$
\begin{aligned}
\rho(x, c(D))&\le\rho(x, c(B))+\rho(c(B), c(D))\le r(B)+\frac{8}{3}r(D)\\
&\le M+\frac{8}{3}r(D) <\frac{4}{3} r(D)+\frac{8}{3} r(D)\\
&=4r(D).
\end{aligned}
$$
Thus, $x\in 4D$, and consequently, $\bigcup \mathcal
B'''_M\subseteq\bigcup_{D\in\mathcal B''_M} 4D.$ Combining this and
(\ref{2hg49}), we get that $\bigcup_{B\in\mathfrak F}
B\subseteq\bigcup_{B\in\mathcal B'_M\cup\mathcal B''_M}4B.$ This and
(b) immediately imply (e). The properties (c), (d) and (e) are
established. Now, take $S=\sup\{r(B): B\in\mathcal B\}+1<+\infty$,
and define inductively the sequence $(\mathcal B'_{(3/4)^n
S})_{n=0}^{\infty}$ by declaring $\mathcal B'_S=\emptyset$ and
$\mathcal B'_{(3/4)^{n+1}S}=\mathcal B'_{(3/4)^n S}\cup\mathcal
B''_{(3/4)^n S}.$ Then
$$
\mathcal B'=\bigcup_{n=0}^{\infty}\mathcal B'_{(3/4)^n S}.
$$
It then follows directly from (d) and our inductive definition that
$\mathcal B'$ consists of mutually disjoint balls. It follows from
(e) that $\bigcup_{B\in\mathcal B'} 4B\supseteq\bigcup\{B\in\mathcal
B: r(B)>0\}=\bigcup\mathcal B$. The first part of our theorem is
thus proved. The last part follows immediately from the fact that
any family of mutually disjoint open subsets of a separable space is
countable.
\endpf

\sp\brem\label{r1hg51} Assume the same as in Theorem~\ref{t1hg47}
(no separability of $X$ is required) and suppose that there exists a
finite Borel measure $\mu$ on $X$ such that $\mu(B)>0$ for all
$B\in\mathcal B'$. Then $\mathcal B'$ is countable. 
\erem

\sp We now shall prove {\it Besicovitch's Covering Theorem}. We consider it as one of the most powerful geometric tools when dealing with some aspects of fractal sets. We deduce easily from it two fundamental classical theorems: {\it Vitali-Type Covering Theorem}, and {\it Lebesgue Density Points Theorem}. For the proof of  Besicovitch's Covering Theorem we introduce two concepts. First,

\bdfn\label{d1_2017_09_19}
Let $(X,\rho)$ be a metric space. A collection
$\cB=\{B(x_i,r_i)\}_{i=1}^\infty$ of open balls centered at a set $A\sbt X$, meaning that $x_i\in A$ for all $i\ge 1$,
is said to be a packing\index{(N)}{packing} of $A$ if and
only if for any pair $i\ne j$
$$
\rho(x_i,x_j) \geq  r_i + r_j.
$$
This property is not in general equivalent to the requirement that
all the balls  $B(x_i, r_i)$ be mutually disjoint. It is obviously
so if $X$ is a Euclidean space. We call the number
$$
r(\cB):=\sup\{r_i:i\ge 1\}
$$
the radius of packing $\cB$. 
\edfn

\fr This notion has a far reaching meaning. It is the key concept to define packing measures and dimensions, which will be done in Section~\ref{Hausdorff and packing measures, Hausdorff and packing dimension}. The other notion we need is this.

\sp For any $x\in \R^n$, any $0<r\le \infty$ and any $0<\a<\pi$ by $\Con(x,\a,r)$ we will denote any solid central cone with vertex $x$, radius $r$ and
angle $\a$. That is, with the above data, for an arbitrary ray $l$ emanating from $x$, we denote
$$
\Con(x,\a,r)=\Con(l,x,\a,r):=\{y\in\R^n:
0<|y-x|<r,\angle(y-x,l)\le\a\}\cup\{x\}.
$$
The proof of Theorem~\ref{th:6.5.1} makes substantial use of the following obvious geometric observation.

\sp

\bobs\label{obsx:6.5.2} Let $n\ge 1$ be an integer. Then there exists
$\a(n)>0$ so small that the following holds.
If $x\in\R^n$, $0<r<\infty$, if $z\in B(x,r)\sms B(x,r/3)$ and if $x\in
\Con(z,\a(n),\infty)$, then the set $\Con(z,\a(n),\infty)\sms B(x,r/3)$
consists
of two connected components: one of $z$ and one of "$\infty$", and the one
containing $z$ is contained in $B(x,r)$.
\eobs

\sp \index{(N)}{Besicovitch's Covering Theorem}
\bthm[ Besicovitch's Covering Theorem]\label{th:6.5.1} 
Let $n\ge 1$ be an integer. Then there exists an integer constant $b(n)\ge 1$ such that the following holds:

\, If $A$ is a bounded subset of $\R^n$ then for any function $r:A\to (0,\infty)$ there exists $\{x_k\}_{k=1}^\infty$, a countable subset of $A$, such that the collection 
$$
\mathcal{B}(A,r):=\{B(x_k,r(x_k)):k\ge 1\}
$$ 
covers $A$ and it can be decomposed into $b(n)$ packings of $A$.
\ethm

\bpf In the sequel we
consider balls in $\R^n$. We will construct the sequence $\{x_k:k= 
1,2,\ld\}$ inductively. Let
$$
a_0:=\sup\{r(x):x\in A\}.
$$
If $a_0=\infty$ then one can find $x\in A$ with $r(x)$ so large that
$B(x,r(x))\spt A$ and the proof is finished.

\sp\fr If $a_0<\infty$ choose $x_1\in A$ so that $r(x_1)>a_0/2$. Fix $k\ge 1$ and
assume that the points $x_1,x_2,\ld,x_k$ have been already chosen. If
$A\sbt B(x_1,r(x_1))\cup\ld\cup B(x_k,r(x_k))$ then the selection process
is finished. Otherwise put
$$
a_k:=\sup\big\{r(x): x\in A\sms \(B(x_1,r(x_1))\cup\ld\cup B(x_k,r(x_k))\)\big\}
$$
and take 
\beq\label{6.5.1}
x_{k+1}\in A\sms \(B(x_1,r(x_1))\cup\ld\cup B(x_k,r(x_k))\)
\eeq
such that
\beq\label{6.5.2}
r(x_{k+1})> a_k/2.
\eeq
In order to shorten notation from now on throughout this proof we will write
$r_k$ for $r(x_k)$. By \eqref{6.5.1} we have $x_l\notin
B(x_k,r_k)$ for all pairs $k,l$ with $k<l$. Hence
\beq\label{6.5.3}
\|x_k -x_l\|\ge r(x_k).
\eeq
It follows from the construction of the sequence $(x_k)$ that
\beq
r_k>a_{k-1}/2 \ge r_l/2,
\eeq\label{6.5.4}
and therefore $r_k/3 + r_l/3 < r_k/3 + 2r_k/3=r_k$. By combining this and
\eqref{6.5.3} we obtain
\beq\label{6.5.5}
B(x_k,r_k/3)\cap B(x_l,r_l/3)=\es
\eeq
for all pairs $k,l$ with $k\ne l$ since then either $k<l$ or $l<k$.

\sp\fr Now we shall show that the balls $\{B(x_k,r_k):k\ge 1\}$ cover $A$.
Indeed, if the selection process stops after finitely many steps this claim is
obvious. Otherwise it follows from \eqref{6.5.5} that $\lim_{k\to\infty}r_k=0$ and
if $x\notin \bu_{k=1}^\infty B(x_k,r_k)$ for some $x\in A$ then by
construction $r_k>a_{k-1}/2 \ge r(x)/2$ for every $k\ge 1$. The contradiction
obtained proves that $\bu_{k=1}^\infty B(x_k,r_k)\spt A$.

\sp\fr The main step of the proof is given by the following.

\sp{\bf Claim:} For every $z\in\R^n$ and any cone $\Con(z,\a(n),\infty)$
($\a(n)$ given by Observation~\ref{obsx:6.5.2}), we have that
$$
\#\{k\ge 1: z\in B(x_k,r_k)\sms B(x_k,r_k/3) \;\text{and}\; x_k\in
\Con(z,\a(n),\infty)\} \le (12)^n.
$$
Denote this $Q$. Our task is to estimate its cardinality from above. If
$Q=\es$, there is nothing to prove. Otherwise let $i=\min Q$. If $k\in Q$ and
$k\ne i$ then $k>i$ and therefore $x_k\notin B(x_i,r_i)$. Therefore, by
Observation~\ref{obsx:6.5.2} applied with $x=x_i$, $r=r_i$, and by the the definition of $Q$, we get that $\|z-x_k\|\ge 2r_i/3$. Hence
\beq\label{6.5.6}
r_k\ge \|z-x_k\|\ge 2r_i/3.
\eeq
On the other hand, by \eqref{6.5.4}, we have that $r_k<2r_i$, and therefore $B(x_k,r_k/3)
\sbt B(z,4r_k/3)\sbt B(z,8r_i/3)$. Thus, using \eqref{6.5.5}, \eqref{6.5.6} and
the fact that the $n$-dimensional volume of balls in $\R^n$ is proportional to
the $n^{\text{th}}$ power of radii, we obtain $\#Q\le (8r_i/3)^n/(2r_i/9)^n
=12^n$. The proof of the claim is finished. 

\sp\fr Clearly, there exists an integer $c(n)\ge 1$ such that for every $z\in\R^n$
the space $\R^n$ can be covered by at most $c(n)$ cones of the form
$\Con(z,\a(n),\infty)$. Therefore, it follows from the above claim that for every
$z\in\R^n$,
$$
\#\{k\ge 1:z\in B(x_k,r_k)\sms B(x_k,r_k/3)\}\le c(n)(12)^n.
$$
Thus applying \eqref{6.5.5}
\beq\label{6.5.7}
\#\{k\ge 1:z\in B(x_k,r_k)\le 1+c(n)(12)^n.
\eeq
Since the ball $\ov B(0,3/2)$ is compact, it contains a finite subset $P$ such
that 
$$
\bu_{x\in P}B(x,1/2)\spt \ov B(0,3/2).
$$
Now for every $k\ge 1$ consider the composition of the map $\R^n\ni x\mapsto
r_kx\in\R^n$ and the translation determined by the vector from $0$ to
$x_k$. Call by $P_k$ the image of $P$ under this affine map. Then
$\#P_k=\#P$, $P_k\sbt\ov B(x_k,3r_k/2)$ and
\beq\label{6.5.8}
\bu_{x\in P_k}B(x,r_k/2)\spt \ov B(0,3r_k/2).
\eeq
Consider now two integers $1\le k<l$ such that
\beq\label{6.5.9}
B(x_k,r_k)\cap B(x_l,r_l)\ne \es.
\eeq
Let $y\in\R^n$ be the only point lying on the interval joining $x_l$ and
$x_k$ at the distance $r_k-r_l/2$ from $x_k$. As $x_l\notin B(x_k,r_k)$, by
\eqref{6.5.9} we have $\|y-x_l\|\le
r_l+r_l/2=3r_l/2$ and therefore by \eqref{6.5.8} there exists $z\in P_l$ such that
$\|z-y\|<r_l/2$. Consequently $z\in B(x_k,r_l/2 + r_k -r_l/2)=B(x_k,r_k)$.
Thus, applying \eqref{6.5.7}, with $z$ being the elements of $P_l$, we obtain the
following
\beq\label{6.5.10}
\#\{1\le k\le l-1:B(x_k,r_k)\cap B(x_l,r_l)\ne\es\}\le\#P(1+c(n)(1)2^n)
\eeq
for every $l\ge 1$. Putting 
$$
b(n):=\#P(1+c(n)(12)^n)+1,
$$ 
this property allows us to decompose the set 
$\N$ of positive integers into $b(n)$ subsets $\N_1,\N_2,\ld,\N_{b(n)}$ in the
following inductive way. For every $k=1,2,\ld,b(n)$ set $\N_k(b(n))=\{k\}$ and
suppose that for every $k=1,2,\ld,b(n)$ and some $j\ge b(n)$, the  mutually disjoint
families $\N_k(j)$ have been already defined so that
$$
\N_1(j)\cup ... \cup \N_{b(n)}(j)=\{1,2,\ld,j\}.
$$
Then by \eqref{6.5.10} there exists at least one $1\le k\le b(n)$ such that
$B(x_{j+1},r_{j+1})\cap B(x_i,r_i)=
\es$ for every $i\in \N_k(j)$. We set 
$$
\N_k(j+1):=\N_k(j)\cup \{j+1\}
$$ 
and
$$
\N_l(j+1)=\N_l(j)
$$ 
for all $l\in\{1,2,\ld,b(n)\}\sms \{k\}$. Putting now for every $k=1,2,\ld,b(n)$
$$
\N_k:=\N_k(b(n))\cup \N_k(b(n)+1)\cup\ld,
$$
we see from the inductive construction that these sets are mutually disjoint,
that they cover $\N$ and that for every $k=1,2,\ld,b(n)$ the families of balls
$\{B(x_l,r_l):l\in \N_k\}$ are also mutually disjoint. The proof of 
Besicovitch's Covering Theorem is finished.
\epf

\sp\fr We would like to emphasize here that the same statement
remains true
if open balls are replaced by closed ones. It also remains true if instead of balls one
considers $n$-dimensional cubes. And in this latter case it is even better: namely, the proof based on the same
idea, is technically considerably easier. There are further, frequently useful, generalizations, especially, a theorem of Morse. The reader is advised to consult the book \cite{G} of Guzman on such topics. 

\sp As we have already mentioned, we can easily deduce from  Besicovitch's
Covering Theorem some other fundamental facts. 
 
\sp\index{(N)}{Vitali-type covering theorem}
\bthm[Vitali--Type Covering Theorem]\label{th:6.5.3}
Let $\mu$ be a probability Borel measure on $\R^n$, let $A\subset\R^n$ be
a Borel set and let 
$\mathcal{B}$ be a family of closed balls such that each point of $A$ is
the center of arbitrarily small balls of $\mathcal{B}$, that is 
$$
\inf\{r:B(x,r)\in\mathcal{B}\}=0 
$$
for all $x\in A$. Then there there exists a  countable (finite or infinite) collection $\mathcal B(A)$ of mutually disjoint balls in $\mathcal{B}$ 
such that
$$
\mu\lt(A\setminus \bigcup\{B\in \mathcal B(A)\}\rt)=0. 
$$
\ethm

\bpf We assume $A$ is bounded, leaving the unbounded case to the reader. 
We may assume $\mu(A)>0$.  The measure $\mu$ restricted to a compact
ball $B(0,R)$ such that $A\subset B(0,R/2)$ 
is Borel hence regular. Hence there exists an open set 
$U\subset\R^n$ containing $A$ and such that
$$
\mu(U)\le (1+(4b(n))^{-1})\mu(A), 
$$
where $b(n)$ is as in  Besicovitch's Covering Theorem~\ref{th:6.5.1}. By 
that theorem applied for closed balls we can decompose $\mathcal{B}$ in
packings $\mathcal{B}_1,...,\mathcal{B}_{b(n)}$ of $A$ contained in $U$,
i.e. each $\mathcal{B}_i$ consists of disjoint balls and 
$$
A\subset\bigcup_{i=1}^{b(n)}\bigcup\mathcal{B}_i\subset U. 
$$
Then, $\mu(A)\le \sum_{i=1}^{b(n)}\mu\(\bigcup\mathcal{B}_i\)$ and
consequently there exists an $i$ such that 
$$
\mu(A)\le b(n)\mu\(\bigcup\mathcal{B}_i\).
$$
Further, for some finite subfamily $\mathcal B'_i$ of $\mathcal B_i$,
$$
\mu(A)\le 2b(n)\mu\(\bigcup\mathcal{B}'_i\). 
$$
Letting $A_1=A\setminus (\bigcup\mathcal{B}'_i\)$ we get
$$
\aligned
\mu(A_1)&\le \mu\(U\setminus \bigcup B'_i\)=\mu(U)-\mu\(\bigcup B'_i\) \\
&\le\bigl(1+\frac14(b(n))^{-1}- \frac12(b(n))^{-1}\bigr)\mu(A)\\
&=u\mu(A) 
\endaligned
$$
with $u:= 1-\frac14(b(n))^{-1}<1$.
Now, consider $A_1$ in the role of $A$ before. Since
$$
A_1\subset \R^n\setminus\(\bigcup\mathcal{B}'_i\)
$$ 
which is open, we find a packing, playing the role of $B'_i$ contained in 
$$
\R^n\setminus\(\bigcup\mathcal{B}'_i\),
$$
so disjoint from
$\bu\mathcal B'_i$. We then get the measure of a non--covered remnant bounded
above by $u\mu(A_1)\le u^2\mu(A)$. 
We can continue, constructing consecutively packings that exhaust the whole set $A$ except at most a set of measure 0. The proof is complete.
\epf

\sp Now we shall prove two quite straightforward consequences of
 Besicovitch's Covering Theorem (Theorem~\ref{th:6.5.1}), the first one being
the celebrated, and to some extent counter-intuitive, Density Points
Theorem. It in fact follows from Vitali's Type
Covering Theorem (Theorem~\ref{th:6.5.3}), which itself is a consequence of
 Besicovitch's Covering Theorem.

\sp\index{(N)}{Lebesgue density points}
\bthm[Lebesgue Density Theorem]\label{th:6.5.4}
Let $\mu$ be a probability Borel measure on $\R^n$ and let $A\subset\R^n$ be a Borel set. Then the limit
$$
\lim_{r\to 0} \frac{\mu(A\cap B(x,r))}{\mu(B(x,r))} 
$$
exists and is equal to 1 for $\mu$--almost every point $x\in\R^n$.
\ethm

\bpf First of all, for every Borel set $B\sbt X$ and every $x\in X$, we have obviously that
$$
\lim_{s\upto r}\mu(B\cap B(x,s))=\mu(B\cap B(x,r))
$$
and 
$$
\lim_{s\downto r}\mu(B\cap B(x,s)\ge \mu(B\cap B(x,r))).
$$
Therefore, the function
$$
X\ni x\lmt \mu(B\cap B(x,r))\in\R
$$
is lower semi--continuous, thus Borel measurable. Hence, the function
$$
X\ni x\lmt\frac{\mu(A\cap B(x,r))}{\mu(B(x,r))}
\in\R
$$
is also Borel measurable. Furthermore, since 
$$
\lim_{\Q\ni s\upto r}\mu(B\cap B(x,s))=\mu(B\cap B(x,r)),
$$
it follows that the set of points $x\in X$ for which the limit 
\beq\label{120190902}
\lim_{r\to 0}\frac{\mu(A\cap B(x,r))}{\mu(B(x,r))}
\eeq
exists is the same as the set of points $x\in X$ for which the limit 
$$
\lim_{\Q\ni r\to 0}\frac{\mu(A\cap B(x,r))}{\mu(B(x,r))}
$$
exists. Since the set $\Q$ of rational numbers is countable, we thus conclude that the set of points in $X$ for which in \eqref{120190902} exists is Borel measurable. 
 
Seeking contradiction, suppose now the set of points in $A$ where this limit 
is either not equal to $1$ or does not exist has positive measure $\mu$. Then there exists $0\le a<1$ and Borel $A'\subset A$ of positive measure $\mu$ such that for every $x\in A'$ there exists a sequence $(r_i(x))_{i=1}^\infty$ of positive radii converging to $0$ such that
$$
\frac{\mu\(A'\cap B(x,r_i(x))\)}{\mu\(B(x,r_i(x))\)}<a
$$
for all $i\ge 1$. Given an open set $U$ containing $A$, let 
$$
\mathcal B_U:=\big\{B(x,r_i(x)):x\in A', \, B(x,r_i(x))\sbt U\big\}.
$$ 
Let then $\cB_U(A')$be the corresponding collection of balls whose existence is asserted in Vitali's Type Covering Theorem (Theorem~\ref{th:6.5.3}). Then 
$$
\mu(A')=\sum_{B\in\mathcal B_U(A')}\mu(A'\cap B)
\le a\sum_{B\in\mathcal B_U(A')}\mu(B)
\le a\mu(U). 
$$
Since measure $\mu$ is regular, this yields $\mu(A')\le a\mu(A')$.
This contradiction finishes the proof.
\epf

\sp\fr The second consequence of Besicovitch's Covering Theorem (Theorem~\ref{th:6.5.1}), \index{(N)}{Besicovitch's Covering Theorem} which we have mentioned above, is the following technical, but very useful and frequently applied, lemma, suitable to prove that one given measure is absolutely continuous with respect to the other. We follow the proof from \cite{DU LMS}.

\blem\lab{l3123006} Let $\mu$ and $\nu$ be Borel probability
measures on $X$, a bounded  subset of a Euclidean space $\R^d$, $d\ge 1$. Suppose
that there is a constant $M >0$ and for every point $x \in Y$
there is  a converging to zero sequence $\(r_j(x)\)_{i=0}^\infty$ of
positive radii such that for all $ j \geq 1$ and all $x\in X$,
$$
\mu(B(x, r_j(x))\leq M \nu(B(x,r_j(x)).
$$
Then  the measure $\mu$  is absolutely continuous with respect to
$\nu$ and the Radon-Nikodym derivative satisfies
$$
d\mu/ d\nu \leq Mb(d),
$$
where $b(d)$ the constant coming from  Besicovitch's Covering Theorem, i.e. Theorem~\ref{th:6.5.1}. 
\elem

\bpf  Consider an arbitrary Borel set $E\sbt X$, and fix
$\varepsilon
>0$. Since  $\lim_{j\to \infty} r_j(x)=0$ and since measure $\nu$ is
regular, for every $x\in E$ there exists  a radius $r(x)$  being  of
the form $r_j(x)$ such that  
$$
\nu\lt(\bu_{x\in E}B_e(x, r(x)\rt)\sms E)<\varepsilon.
$$
Now, by  Besicovitch's Covering Theorem (Theorem~\ref{th:6.5.1}) we
can choose a countable subcover $\{B(x_i,
r_i(x))\}_{i=1}^{\infty}$ from the cover $\{B(x_i, r_i(x))\}_{x\in
E}$ of $E$, of multiplicity bounded above by $b(d)$. Therefore, we obtain
$$
\begin{aligned}
\mu(E) 
& \leq \sum_{i=1}^{\infty} \mu(B(x_i, r_i(x)))
  \leq M \sum_{i=1}^{\infty} \nu(B(x_i,r_i(x)))\\
&\leq  Mb(d)  \nu \lt(\bu_{i=1}^{\infty}B(x_i,r_i(x))\rt)\\
&\leq  Mb(d)( \varepsilon +\nu(E)).
\end{aligned}
$$
Letting $\varepsilon \searrow 0 $  we thus obtain $\mu(E) \leq Mb(d)
\nu(E).$ Therefore $\mu$ is  absolutely continuous  with respect  to
$\nu$  with the Radon--Nikodym derivative bounded above by $Mb(d)$.
\endpf

\sp\section{Conditional Expectations and Martingale Theorems}\label{CEaMT}

The content of this section belongs rather to the probability theory than to the classical measure theory. Its culmination (for us), i.e. Theorem~\ref{t120191123}, is however similar to the Lebesgue Density Theorem, i.e. Theorem~\ref{th:6.5.4}, so it is natural to place it here. This chapter is about conditional expectations and martingales and it is closely modeled on Chapter form Patrick Billingsley's book \cite{Bi}. 

We start with conditional expectations. Let $(X,\F,\mu)$ be a probability space. Let $\mD$ be a sub $\sg$--algebra of $\F$. Let 
$$
\phi: X \lra \mathbb R
$$ 
be a measurable function, integrable with respect to the measure $\mu$. We denote by 
$$
E(\phi|\mD)=E_\mu(\phi|\mD)
$$
the (conditional) expected value of $\phi$ with respect to the $\sg$--algebra $\mD$. This is the only function (up to sets of measure zero) that is measurable with respect to the $\sg$--algebra $\mD$ and such that
$$
\int_DE_\mu(\phi|\mD)\,d\mu=\int_D\phi\,d\mu
$$
for every sets $D\in\mD$. Its existence for non--negative integrable functions is a straightforward consequence of Radon--Nikodym Theorem. In the general case one sets
$$
E_\mu(\phi|\mD):=E_\mu(\phi_+|\mD)-E_\mu(\phi_-|\mD)
$$
Uniqueness is obvious.

\sp Conditional expectations exhibit several natural properties. We list below some of the them, most basic ones. Their proofs are straightforward and are omitted.

\bprop\label{expvalprops}
Let $(X,\mathcal{A},\mu)$ be a probability space, let $\mathcal{B}$ and $\mathcal{C}$ denote some sub--$\sg$--algebras of $\mathcal{A}$ and let $\varphi\in L^1(X,\mathcal{A},\mu)$. Then the following hold.

\begin{itemize}
\item[(a)] If $\varphi\geq0$ $\mu$-a.e., then 
$$
E(\varphi|\mathcal{B})\geq 0 \  \  \  \mu-a.e.
$$

\, 

\item[(b)] If $\varphi_1\geq\varphi_2$ $\mu$--a.e., then 
$$
E(\varphi_1|\mathcal{B})\geq E(\varphi_2|\mathcal{B}) \  \  \mu-a.e.
$$

\, 

\item[(c)] $\bigl|E(\varphi|\mathcal{B})\bigr|\leq E\bigl(|\varphi|\,\bigl|\bigr.\,\mathcal{B}\bigr)$. 

\

\item[(d)] The functional $E(\cdot|\mathcal{B})$ is linear. In other words, 
             for any $c_1,c_2\in\R$ and $\varphi_1,\varphi_2\in L^1(X,\mathcal{A},\mu)$, we have that
           \[
					 E(c_1\varphi_1+c_2\varphi_2|\mathcal{B})=c_1\,E(\varphi_1|\mathcal{B})+c_2\,E(\varphi_2|\mathcal{B}).
					 \]
					 
\, 

\item[(e)] If $\varphi$ is already $\mathcal{B}$-measurable, then $E(\varphi|\mathcal{B})=\varphi$. In particular, we have that 
$$
E\bigl(E(\varphi|\mathcal{B})\bigl|\mathcal{B}\bigr)=E(\varphi|\mathcal{B}).$$ 
Also, if $\varphi=c\in\R$ is a constant function, then $E(\varphi|\mathcal{B})=\varphi=c$.
						 
\

\item[(f)] If $\mathcal{B}\spt\mathcal{C}$ then
$$
E\bigl(E(\varphi|\mathcal{B})\bigl|\mathcal{C}\bigr)
=E(\varphi|\mathcal{C}).
$$
\end{itemize}
\eprop

We will now determine conditional expectations of an arbitrary integrable functions $\varphi$ with respect to various sub-$\sg$--algebras that are of particular interest and are simple enough.

\bex\label{expneg}{\rm 
Let $(X,\mathcal{A},\mu)$ be a probability space. The family $\mathcal{B}$
of all measurable sets that are either of null or of full measure
constitutes a sub--$\sg$--algebra of $\mathcal{A}$. Let $\varphi\in L^1(X,\mathcal{A},\mu)$. Since $E(\varphi|\mathcal{B})$ is $\mathcal{B}$--measurable, 
$$
E(\varphi|\mathcal{B})^{-1}(\{t\})\in\mathcal{B}
$$
for each $t\in\R$, meaning that the set $E(\varphi|\mathcal{B})^{-1}(\{t\})$ 
is either of measure zero or of measure one. Also bear in mind that
\[
X=E(\varphi|\mathcal{B})^{-1}(\R)
=\bigcup_{t\in\R}E(\varphi|\mathcal{B})^{-1}(\{t\}).
\]
Since the above union consists of mutually disjoint sets of measure zero and one, it follows that only one of these sets can be of measure one. In other words, there exists a unique $t\in\R$ such that
$$
E(\varphi|\mathcal{B})^{-1}(\{t\})=A
$$ 
for some $A\in\mathcal{A}$ with $\mu(A)=1$.
Because the function $E(\varphi|\mathcal{B})$ is unique up to a set of measure zero,
we may assume without loss of generality that $A=X$.
Hence $E(\varphi|\mathcal{B})$ is a constant function. Therefore, 
\[
E(\varphi|\mathcal{B})=\int_X E(\varphi|\mathcal{B})\,d\mu=\int_X \varphi\,d\mu.
\]
}
\eex

\bex\label{measpart}{\rm 
Let $(X,\mathcal{A})$ be a measurable space and let $\a$ 
be a countable measurable partition of $X$. The sub--$\sg$--algebra $\sg(\a)$ of $\mathcal{A}$ generated by $\a$ 
is the family of all sets which can be represented as a union of elements of $\a$. 
When $\a$ is finite, so is $\sg(\a)$. When $\a$ is countably infinite, $\sg(\a)$ is uncountable, in fact of cardinality continuum.
Let $\mu$ be a probability measure on $(X,\mathcal{A})$. 
Let $\varphi\in L^1(X,\mathcal{A},\mu)$. Since $E(\varphi|\sg(\a))$ is $\mathcal{B}$--measurable, 
$$
E(\varphi|\mathcal{B})^{-1}(\{t\})\in\sg(\a)
$$
for each $t\in\R$, i.e.
$$
E\(\varphi|\sg(\a)\)^{-1}(\{t\})\in\sg(\a).
$$ 
This means that the set $E\(\varphi|\sg(\a)\)^{-1}(\{t\})$ is a union of elements of $\a$. 
This further means that the conditional expectation function
$E(\varphi|\sg(\a))$ is constant on each element of $\a$. 
Let $A\in\a$. If $\mu(A)=0$ then $E(\varphi|\sg(\a))|_{A}=0$. Otherwise, 
\beq\label{120191123}
E(\varphi|\sg(\a))|_{A}
=\frac{1}{\mu(A)}\int_{A} E(\varphi|\sg(\a))\,d\mu
=\frac{1}{\mu(A)}\int_{A_n}\varphi\,d\mu.
\eeq
In summary, the conditional expectation $E(\varphi|\mathcal{B})$ of a function $\varphi$ with respect
to a sub-$\sg$-algebra generated by a countable measurable partition is constant on each
element of that partition. More precisely, on any given element of the partition,
$E(\varphi|\mathcal{B})$ is equal to the mean value of $\varphi$ on that element. In particular, if $\a$ is a trivial partition, i.e. consisting of sets of measure zero and one only, then
\beq\lab{1_mu_2014_11_11}
E(\varphi|\sg(\a))=\int_{X}\varphi\,d\mu \  \  \  \  \mu-a.e.
\eeq
}
\eex

The next result is a special case of a theorem originally due to Doob and is commonly 
called the Martingale Convergence Theorem. In order to discuss it, we first define the martingale itself.

\bdfn\label{martin}
Let $(X,\mathcal{A},\mu)$ be a probability space. Let $(\mathcal{A}_n)_{n=1}^\infty$ be a sequence of sub--$\sg$--algebras of $\mathcal{A}$. Let also $(\varphi_n:X\lra\R)_{n=1}^\infty$ 
be a sequence of random variables, that is, a sequence of $\mathcal{A}$--measurable functions. 
The sequence 
$$
\((\varphi_n,\mathcal{A}_n)\)_{n=1}^\infty
$$ 
is called a {\em martingale} \index{(N)}{martingale} if and only if the following conditions are satisfied:
\begin{itemize}
\item[(a)] $(\mathcal{A}_n)_{n=1}^\infty$ is an ascending sequence, that is, 
           $\mathcal{A}_{n+1}\spt\mathcal{A}_n$ for all $n\in\N$.
           
\,

\item[(b)] $\varphi_n$ is $\mathcal{A}_n$--measurable for all $n\in\N$.
   
\,

\item[(c)] $\varphi_n\in L^1(\mu)$ for all $n\in\N$.
   
\,

\item[(d)] $E(\varphi_{n+1}|\mathcal{A}_n)=\varphi_n$ $\mu$--a.e. for all $n\in\N$.
\end{itemize}
\edfn

The main, and frequently referred to as the simplest, convergence theorem concerning martingales is this.

\bthm [Martingale Convergence Theorem]\label{martintheorem}
\index{(N)}{Martingale Convergence Theorem}
Let $(X,\mathcal{A},\mu)$ be a probability space. 
If $((\varphi_n,\mathcal{A}_n))_{n=1}^\infty$ is a martingale 
such that 
\[
\sup\{\|\varphi_n\|_1:n\in\N\}<+\infty,
\]
then there exists $\widehat{\varphi}\in L^1(X,\mathcal{A},\mu)$ such that 
\[
\lim_{n\to\infty}\varphi_n(x)=\widehat{\varphi}(x)\  \ \text{ for }\,\mu\text{--a.e. }x\in X
\]
and 
\[
\|\widehat{\varphi}\|_1\leq\sup\{\|\varphi_n\|_1:n\in\N\}<+\infty.
\]
\ethm

\sp\fr This is a special case of Theorems~35.5 in Billingsley's book \cite{Bi} proved therein. Its proof is just too long and too involved to be reproduced here. We omit it. 

One natural martingale is formed by the conditional expectations of a function with respect to an ascending sequence of sub--$\sg$--algebras.

\bprop\label{martexp}
Let $(X,\mathcal{A},\mu)$ be a probability space and let $(\mathcal{A}_n)_{n=1}^\infty$ be 
an ascending sequence of sub--$\sg$--algebras of $\mathcal{A}$. For any $\varphi\in L^1(X,\mathcal{A},\mu)$, 
the sequence 
$$
\((E(\varphi|\mathcal{A}_n),\mathcal{A}_n)\)_{n=1}^\infty
$$ 
is a martingale. 
\eprop

\bpf 
Indeed, set 
$$
\varphi_n:=E(\varphi|\mathcal{A}_n)
$$ 
for all $n\in\N$. Condition~(a) in Definition~\ref{martin} 
is automatically fulfilled. Conditions~(b) and~(c) follow from the very definition of the conditional expectation 
function. Regarding condition~(d), a straightforward application of Proposition~\ref{expvalprops}(f) gives
\[
E(\varphi_{n+1}|\mathcal{A}_n)
=E\bigl(E(\varphi|\mathcal{A}_{n+1})|\mathcal{A}_n\bigr)
=E(\varphi|\mathcal{A}_n)
=\varphi_n
\]
$\mu$--a.e. for all $n\in\N$.
So $((E(\varphi|\mathcal{A}_n),\mathcal{A}_n))_{n=1}^\infty$ is a martingale. 
\epf

\sp With the hypotheses of this proposition, by using Proposition~\ref{expvalprops}(c), we see that
\[
\sup_{n\in\N}\|\varphi_n\|_1
=\sup_{n\in\N}\int_X \bigl|E(\varphi|\mathcal{A}_n)\bigr|\,d\mu
\leq\sup_{n\in\N}\int_X E\bigl(|\varphi|\,\bigl|\bigr.\,\mathcal{A}_n\bigr)\,d\mu
=\int_X |\varphi|\,d\mu
<\infty.
\]
According to Theorem~\ref{martintheorem}, there thus exists $\widehat{\varphi}\in L^1(X,\mathcal{A},\mu)$ such that 
\[
\lim_{n\to\infty}E(\varphi|\mathcal{A}_n)(x)=\widehat{\varphi}(x) \text{ for } \mu\text{-a.e. }x\in X
\hspace{0.5cm}\text{ and }\hspace{0.5cm}
\|\widehat{\varphi}\|_1\leq\|\varphi\|_1. 
\]
What is $\widehat{\varphi}$? This is the question we will address now. For this we need the concept of uniform integrability and a convergence theorem it entails.

\bdfn\label{unifintseq}
Let $(X,\mathcal{A},\mu)$ be a measure space. A sequence of measurable functions $(f_n)_{n=1}^\infty$ is called uniformly integrable index{(N)}{ uniformly integrable} if and only if 
\[
\lim_{M\to\infty}\,\sup_{n\in\N}\int_{\{|f_n|\geq M\}}|f_n|\,d\mu=0.
\]
\edfn

\fr The following theorem is classical in measure theory. It is proved for example as Theorem~16.14 in Billingsley's book ~\cite{Bi}. 

\bthm\label{unifintthm}
Let $(X,\mathcal{A},\mu)$ be a finite measure space and $(f_n)_{n=1}^\infty$ a sequence 
of measurable functions that converges pointwise $\mu$-a.e. to a function $f$. 
\begin{itemize}
\item[(a)] If $(f_n)_{n=1}^\infty$ is uniformly integrable, then $f_n\in L^1(\mu)$ for all $n\in\N$
and $f\in L^1(\mu)$. Moreover, 
\[
\lim_{n\to\infty}\|f_n-f\|_1=0
\ \ \ \textrm{and}\ \ \ 
\lim_{n\to\infty}\int_X f_n\,d\mu
=\int_X f\,d\mu.
\]
\item[(b)] If $f,f_n\in L^1(\mu)$ and $f_n\geq0$ $\mu$-a.e. for all $n\in\N$, then 
$\lim_{n\to\infty}\int_X f_n\,d\mu=\int_X f\,d\mu$ implies that 
$(f_n)_{n=1}^\infty$ is uniformly integrable. 
\end{itemize} 
\ethm

We shall now prove uniform integrability of the martingale appearing in Proposition~\ref{martexp}.
 
(see Definition~\ref{unifintseq}).

\blem\label{martunifinteg}
Let $(X,\mathcal{A},\mu)$ be a probability space and let $(\mathcal{A}_n)_{n=1}^\infty$ be a sequence of sub--$\sg$--algebras of $\mathcal{A}$. Then for very $\varphi\in L^1(X,\mathcal{A},\mu)$, the sequence $(E(\varphi|\mathcal{A}_n))_{n=1}^\infty$ is uniformly integrable.
\elem

\bpf
Without loss of generality, we may assume that $\varphi\geq0$. 
Let $\e>0$. Since the measure $\nu$ on $(X,\mathcal{A}$ given by the formula
$$
\nu(A):=\int_A \varphi\,d\mu
$$ 
is absolutely continuous with respect to $\mu$, it follows from the Radon--Nikodym Theorem (Theorem~\ref{Radon_Nikodym}) that there exists $\d>0$ such that 
\begin{equation}\label{dernt}
A\in\mathcal{A},\ \mu(A)<\d\ \ \ \Longrightarrow \ \ \ \int_A \varphi\,d\mu<\e.
\end{equation}
Consider any 
$$
M>\frac1\d\int_X \varphi\,d\mu.
$$
For each $n\in\N$, let 
\[
X_n(M):=\bigl\{x\in X:E(\varphi|\mathcal{A}_n)(x)\geq M\bigr\}.
\]
Observe that $X_n(M)\in\mathcal{A}_n$ since $E(\varphi|\mathcal{A}_n)$ 
is $\mathcal{A}_n$-measurable. Therefore, by Tchebyschev's Inequality, we get that
\[
\mu(X_n(M))
\leq\frac{1}{M}\int_{X_n(M)}E(\varphi|\mathcal{A}_n)\,d\mu
=\frac{1}{M}\int_{X_n(M)}\varphi\,d\mu
\leq\frac{1}{M}\int_X \varphi\,d\mu<\d
\]
for all $n\in\N$. Consequently, by~(\ref{dernt}),
\[
\int_{X_n(M)}E(\varphi|\mathcal{A}_n)\,d\mu=\int_{X_n(M)}\varphi\,d\mu<\e
\]
for all $n\in\N$. Thus,
\[
\sup_{n\in\N}\int_{\{E(\varphi|\mathcal{A}_n)\geq M\}}E(\varphi|\mathcal{A}_n)\,d\mu\leq\e.
\]
Therefore,
\[
\lim_{M\to\infty}\sup_{n\in\N}\int_{\{E(\varphi|\mathcal{A}_n)\geq M\}}E(\varphi|\mathcal{A}_n)\,d\mu=0,
\]
that is, $(E(\varphi|\mathcal{A}_n))_{n=1}^\infty$ is uniformly integrable.
\epf

\sp\bthm
[Martingale Convergence Theorem for Conditional Expectations]\label{DMCT}
\index{(N)}{Martingale Convergence Theorem for Conditional Expectations}
Let $(X,\mathcal{A},\mu)$ be a probability space and let $\varphi\in L^1(X,\mathcal{A},\mu)$. Let $(\mathcal{A}_n)_{n=1}^\infty$ be an ascending sequence of sub--$\sg$--algebras of $\mathcal{A}$ and let
\[
\mathcal{A}_\infty:=\sg\Bigl(\bigcup_{n=1}^\infty\mathcal{A}_n\Bigr). 
\]
Then
$$
\lim_{n\to\infty}E(\varphi|\mathcal{A}_n)=E(\varphi|\mathcal{A}_\infty)\  \ \mu\text{-a.e. on }X
$$
and
$$ 
\lim_{n\to\infty}\bigl\|E(\varphi|\mathcal{A}_n)-E(\varphi|\mathcal{A}_\infty)\bigr\|_1=0.
$$
\ethm

\bpf 
Let 
$$
\varphi_n:=E(\varphi|\mathcal{A}_n).
$$
It follows from Proposition~\ref{martexp} and Lemma~\ref{martunifinteg} that $((\varphi_n,\mathcal{A}_n))_{n=1}^\infty$ is a uniformly integrable martingale such that 
\[
\lim_{n\to\infty}\varphi_n=\widehat{\varphi} \  \  \  \mu\text{--a.e. on }X
\]
for some $\widehat{\varphi}\in L^1(X,\mathcal{A},\mu)$. 
For all $n\in\N$ the function $\varphi_n$ is $\mathcal{A}_\infty$--measurable since it is $\mathcal{A}_n$--measurable and $\mathcal{A}_n\sbt\mathcal{A}_\infty$. Thus $\widehat{\varphi}$ is $\mathcal{A}_\infty$-measurable, too. Moreover, it follows from Theorem~\ref{unifintthm} that 
\[
\lim_{n\to\infty}\|\varphi_n-\widehat{\varphi}\|_1=0
\ \ \ \textrm{and}\ \ \
\lim_{n\to\infty}\int_A \varphi_n\,d\mu=\int_A \widehat{\varphi}\,d\mu 
\]
for all $A\in\mathcal{A}$. Therefore, it just remains to show that
$$
\widehat{\varphi}=E(\varphi|\mathcal{A}_\infty).
$$
Let $k\in\N$ and $A\in\mathcal{A}_k$. If $n\geq k$, then $A\in\mathcal{A}_n\sbt\mathcal{A}_\infty$ 
and thus
\[
\int_A \varphi_n\,d\mu
=\int_A E(\varphi|\mathcal{A}_n)\,d\mu
=\int_A \varphi\,d\mu
=\int_A E(\varphi|\mathcal{A}_\infty)\,d\mu.
\]
Letting $n\to\infty$, this yields
\[
\int_A \widehat{\varphi}\,d\mu
=\int_A E(\varphi|\mathcal{A}_\infty)\,d\mu._k.
\]
Since $k$ was arbitrary, this entails
\[
\int_B \widehat{\varphi}\,d\mu
=\int_B E(\varphi|\mathcal{A}_\infty)\,d\mu
\]
for all $B\in\bigcup_{k=1}^\infty\mathcal{A}_k$. Finally, since $\bigcup_{k=1}^\infty\mathcal{A}_k$ is a $\pi$--system generating $\mathcal{A}_\infty$
and since both $\widehat{\varphi}$ and $E(\varphi|\mathcal{A}_\infty)$ are $\mathcal{A}_\infty$--measurable, we conclude that 
$$
\widehat{\varphi}=E(\varphi|\mathcal{A}_\infty) \  \  \ \mu\text{{\rm--a.e. in }} X.
$$  
\epf

\sp Recall that countable measurable partitions of a measurable space are defined and systematically treated in Section~\ref{Partitions}; they form the key concept for all of Chapter~\ref{ch8}. As an immediate consequence of   
Martingale Convergence Theorem for Conditional Expectations, i.e. Theorem~\ref{DMCT} and formula \eqref{120191123}, we get the following theorem somewhat similar to the Lebesgue Density Theorem, i.e. Theorem~\ref{th:6.5.4} from the previous section.

\bthm\label{t120191123}
Let $(X,\mathcal{A},\mu)$ be a probability space and let $(\a_n)_{n=1}^\infty$ be a sequence of finer and finer (more precisely, $\a_{n+1}$ is finer than $\a_n$ for all $n\ge 1$) countable measurable partitions of $X$  which generates the $\sg$--algebra $\mathcal{A}$, i.e $\sg\(\bu_{n\ge 1}\a_n\)=\mathcal{A}$. Then for every set $A\in\mathcal{A}$ and for $\mu$--a.e. $x\in A$, we have that
$$
\lim_{n\to\infty}\frac{\mu(A\cap \a_n(x))}{\mu(\a_n(x))}=\1_A(x)=1.
$$
\ethm

\sp\section{Hausdorff and packing measures, Hausdorff and packing dimension}
\label{Hausdorff and packing measures, Hausdorff and packing dimension}

\fr In this section we introduce the basic geometric concepts on metric spaces. These are Hausdorff measures, Hausdorff dimension, packing measures and packing dimensions. \index{(N)}{Hausdorff measure} \index{(N)}{Hausdorff dimension} \index{(N)}{packing measure} \index{(N)}{packing dimension}  We prove their fundamental properties. While Hausdorff measures and Hausdorff dimensions were introduced quite early, in 1919, by Felix Hausdorff in \cite{Hausdorff}, it took several decades more for packing measures and packing dimension to have been defined. It had been done in stages in \cite{Tricot}, \cite{Taylor_Tricott} and \cite{Su5}. There are now plenty of books on these concepts; we refer the reader for example to \cite{Fal1},  \cite{Fal2}, \cite{Fal3}, \cite{Mat}, and \cite{PU2}. Interesting, not only for historical reasons is also the classical book \cite{Rogers} by C. A. Rogers, which appeared first time in 1970. The 1998 edition is particularly interesting because of Falconer's comments it contains.

Let $\varphi:[0,+\infty)\lra [0,+\infty)$ be a function with the following properties:

\sp\begin{itemize}
\item $\varphi$ is non--decreasing, meaning that $s\le t\, \imp \varphi(s)\le \varphi(t)$.

\sp\item $\varphi(0)=0$ and $\varphi$ is continuous at $0$. 

\sp\item $\varphi((0,+\infty))\sbt (0,+\infty)$.
\end{itemize}
Any function $\varphi:[0,+\infty)\to [0,+\infty)$ with such properties is refered to in the sequel as a gauge function.

\sp\fr Let $(X,\rho)$ be a metric space. For every $\delta>0$ define
\beq\label{FP8.2.1}
\H_\varphi^\delta(A)
:= \inf \Big\{ \sum_{i=1}^\infty \varphi (\diam (U_i))\Big\},
\eeq\index{(S)}{$\H_\varphi^\delta$}

\fr\fr where the infimum  is taken  over all countable covers $\{U_i\}_{i=1}^\infty$ of $A$  of diameter not exceeding $\delta$.

\sp\fr We shall check that for every $\d>0$, $\H_\varphi^\delta$ is an outer measure.
Conditions (\ref{FP8.1.1}) and  (\ref{FP8.1.2}) of Definition~\ref{def_outer_measure}, defining the concept of outer measures, are obviously
satisfied  with $\mu= \H_\varphi^\delta$. To verify (\ref{FP8.1.3})
let $\{A_n\}_{n=1}^\infty$ be a countable family of subsets
of $X$.  Given $\varepsilon >0$   for every $n \geq 1$  we can find
a countable cover $\{U^n_i\}_{i=1}^\infty$  of $A_n$ with
diameters not exceeding $\delta$ such that
$$
\sum_{i=1}^\infty\varphi(\diam (U_i))
\leq \H^\delta_\varphi(A_n)+\frac{\varepsilon }{2^n}.
$$ 
Then the family $\{U^n_i:i,n\ge 1\}$ covers $\bigcup_{n=1}^\infty A_n$, and
$$
\H^\delta_\varphi \left(\bigcup_{n=1}^\infty A_n\right) 
\leq \sum_{n=1}^\infty \sum_{i=1}^\infty \varphi( \diam (U^n_i))
\leq \sum_{n=1}^\infty \Big(\H^\delta_\varphi (A_n) +\frac{\varepsilon }{2^n}\Big)
=\sum_{n=1}^\infty \H^\delta_\varphi (A_n) + \varepsilon.
$$ 
Thus, letting $\varepsilon \downto 0$, the formula (\ref{FP8.1.3}) follows, proving that
$\H^\delta_\varphi$ is an outer measure. Define
\beq\label{FP8.2.2}
\H_\varphi(A):=\sup_{\delta>0 }\Big\{\H_\varphi^\delta(A)\Big\}
=\lim_{\delta \to 0}\H_\varphi^\delta(A).
\eeq\index{(S)}{$\H_\varphi$}

\fr\fr The limit exists since
$\H_\varphi^\delta(A)$ increases as $\delta$ decreases, though it may happen to be infinite. Since all
$\H_\varphi^\delta$  are outer measures. It is therefore immediate that $\H_\varphi$ is an outer measure too.  Moreover, $\H_\varphi$ is a metric measure, since if $A$ and $B$ are two
positively  separated sets in $X$, then no set of diameter less than
$\rho(A, B)$  can intersect both $A$ and $B$. Consequently
$$ 
\H_\varphi^\delta(A\cup B)=\H_\varphi^\delta(A)+\H_\varphi^\delta(B)
$$ 
for all $ \delta <\rho(A,B)$. Letting $
\delta \downto 0$ we get the same formula for $\H_\varphi$, which is
just (\ref{FP8.1.6}) with $\mu=\H_\varphi$. The  metric outer
measure $\H_\varphi$ is called  the  Hausdorff  outer
measure\index{(N)}{Hausdorff outer measure}  associated to the gauge function
$\varphi$. Its restriction to the
$\sg$-algebra  of $\H_\varphi$-measurable sets, which by
Theorem~\ref{FPt8.1.4} includes all the Borel sets, is called  the
\index{(N)}{Hausdorff measure} Hausdorff  measure associated to the function $\varphi$. We should add that even if $E\sbt X$ is not a Borel set, nor even $\H_\varphi$-measurable, we nevertheless commonly referto $\H_\varphi(E)$ as the Hausdorff measure of $E$ rather than Hausdorff outer measure of $E$.

\

\fr As an immediate consequence of the definition of the Hausdorff
measure  and the  properties  of the function  $\varphi$, we get the
following

\bprop\label{FPp8.2.1} 
For any gauge function $\varphi$ the Hausdorff  measure $\H_\varphi$ is atomless. 
\eprop

\fr A particularly important  role is played by the gauge functions of the form
$$
\varphi_t(r)=r^t
$$
for $t>0$. In this case  the corresponding  outer Hausdorff measure is
denoted by  $\H^t$. So, $\H^t$ can be briefly defined as:

\bdfn\label{d2_mu_2014_11_07}
Given $t\ge 0$, the $t$-dimensional outer Hausdorff
measure\index{(N)}{Hausdorff measure} $\H^t(A)$ of the set $A$ is
equal to
$$
\H^t(A)=\sup_{\delta>0}\inf\Big\{\sum_{i=1}^\infty \diam^t(A_i) \Big\},
$$\index{(S)}{$\H^t$}

\fr where infimum is taken over all countable covers
$\{A_i\}_{i=1}^\infty$ of $A$ by the sets with
diameters $\le\delta$. 
\edfn

\brem\label{r1_mu_2014_11_07}
Since $\diam(\ov A)=\diam(A)$ for every set $A\sbt X$, we may, in the Definition~\ref{d2_mu_2014_11_07}, restrict ourselves to closed sets $A_i$ only.
\erem

Having defined Hausdorff measures, we now pass to defined the dual concept, that is that of packing measures. As it was said at the beginning of this section, while Hausdorff measures were introduced quite early, in 1919 by Felix Hausdorff in \cite{Hausdorff}, it took several decades more for packing measures to have been defined. It had been done in stages in \cite{Tricot}, \cite{Taylor_Tricott} and \cite{Su5}. We do it now. We recall that in Definition~\ref{d1_2017_09_19} we have introduced the concept of packing. We will also use it now. For every $A\sbt X$ and every $\delta >0$ let
\beq\label{FP8.3.1}
\Pi_\varphi^{* \delta} (A)
:= \sup \Big\{ \sum_{i=1}^\infty \varphi(\diam (r_i))\Big\},
 \eeq
where the supremum is taken  over all  packings  $\{B(x_i, r_i)\}_{i=1}^\infty$ of $A$ with radii not exceeding $\delta $. Let
\beq\label{FP8.3.2}
\Pi_\varphi^{*} (A)
:=\inf_{\delta
> 0} \Big\{\Pi_\varphi^{* \delta}(A)\Big\}
=\lim_{\delta \to 0} \Pi_\varphi^{* \delta} (A).
\eeq
The limit exists since $ \Pi_\varphi^{* \delta} (A)$ decreases as $
\delta$ decreases. Although the function $\Pi^*_\varphi$ satisfies
condition (\ref{FP8.1.2}) of outer measures, however in contrast to the case of
Hausdorff measures, this function need not to be subadditive, i.e.
conditions (\ref{FP8.1.3}) in general fails. In  order to obtain  an
outer measure we make one step more and we put
\beq\label{FP8.3.3}
\Pi_\varphi (A):= \inf \Big\{ \sum_{i=1}^\infty \Pi_\varphi^{* \delta}(A_i)\Big\},
\eeq\index{(S)}{$\Pi_\varphi $}

\fr where the infimum  is taken  over all countable covers $\{A_i\}{i=1}^\infty $  of $A$.
Analogously as in the case of Hausdorff measure, one checks, with similar
arguments, that
$\Pi_\varphi$ is already an outer measure. Furthermore, it is a metric
outer measure on $X$. It will be called the outer packing measure\index{(N)}{outer packing measure}, associated to the gauge function $\varphi$. Its restriction to the
$\sigma$-algebra of $\Pi_\varphi$-measurable sets, which by
Theorem~\ref{FPt8.1.4} includes all Borel sets, will  be called
the packing measure associated to the gauge function
$\varphi$.\index{(N)}{packing measure} 

\sp\fr In the case of gauge functions 
$$
\varphi_t(r)=r^t, 
$$
where $t>0$, the definition of  the outer packing measure takes the following form.

\bdfn\label{packing} 
The $t$-dimensional outer packing
measure\index{(N)}{packing measure} $\Pi^t(A)$\index{(S)}{$\Pi^t$}
of a set $A\sbt X$ is given by
$$ \Pi^t(A)\index{(S)}{$\Pi^t$}=
\inf_{\cup A_i=A}\Bigl\{\sum_i \Pi^t_*(A_i)\Bigr\}
$$
($A_i$ are arbitrary subsets of $A$), where
$$
\Pi^t_*(A)= \sup_{\delta>0}\sup\Bigl\{\sum_{i=1}^\infty
r_i^t\Bigr\}.
$$
Here the second supremum is taken over all packings
$\{B(x_i,r_i)\}_{i=1}^\infty$ of the set $A$ by open balls centered
at $A$ with radii which do not exceed $\delta$. 
\edfn

\

\fr From now on through the book, in order to get more meaningful geometric consequences, we assume that for a given gauge function $\phi:[0,+\infty)\to[0,+\infty)$ there  exists a function $C_\phi: (0, \infty) \to (0, \infty)$ such
that, for every $ a  \in (0, \infty)$  and every  $t>0$ sufficiently
small (depending  on $a$)
\beq\label{FP(8.2.3)}
  C_\phi(a)^{-1} \phi(t) \leq \phi(at) \leq C_\phi(a) \phi(t).
\eeq
We frequently refer to such gauge functions as evenly varying\index{(N)}{evenly varying function}.
Since $(at)^r= a^r t^r$, all the gauge functions $ \phi$ of the form $ r \mapsto
r^t$  satisfy (\ref{FP(8.2.3)}) with $C_\phi(a)=a^t$.

\sp\fr We now shall establish a simple, but crucial for geometrical consequences, relation between Hausdorff and packing measures.

\sp

\bprop\label{FPp8.3.1} For every set $A \sbt X$ it holds that
$\H_\varphi(A)\leq C_\varphi(2) \Pi_\varphi(A)$. \eprop

\bpf  First we shall now show that, for every  set $A
\sbt X$ and every $ \delta >0$
\beq\label{FP8.3.4}
\H^{2 \delta}_\varphi( A) \leq C_\varphi(2)
\Pi^{*\delta}_\varphi(A).
\eeq
Indeed, if there is no finite maximal (in the sense of inclusion)
packing of the set $A$ of the form $\{B(x_i,\d)\}_{i=1}^\infty$, then for
every $k \geq 1$ there is a packing $ \{B(x_i,\d)\}_{i=1}^k$ of $A$, and therefore 
$$
\Pi^{*\delta}_\varphi(A) \geq
\sum_{i=1}^k\varphi(\d) =k \varphi(\d).
$$ 
Since $\varphi(\d) >0$, this
yields $\Pi^{*\delta}_\varphi(A)=\infty$, and
(\ref{FP8.3.4}) holds. Otherwise, let $ \{B(x_i,\d)\}_{i=1}^l$ be a finite maximal  packing of $A$. Then the collection $\{B(x_i, 2\d)\}$ covers $A$,  and therefore
$$
\H^{2 \delta}_\varphi( A)
\leq \sum_{i=1}^l\varphi(2\d) 
\leq C_\varphi(2) l \varphi(\d)
\leq  C_\varphi(2)\Pi^{*\delta}_\varphi(A),
$$  
Hence, (\ref{FP8.3.4}) is satisfied. Thus, letting
$\d\downto 0$, we get that
\beq\label{FP8.3.5}
\H_\varphi( A) \leq C_\varphi(2) \Pi_\varphi^*(A).
\eeq
So, if $ \{A_n\}_{n\geq 1}$ is countable cover of $A$, then
$$\H_\varphi(A) \leq \sum_{n=1}^\infty \H_\varphi(A_i) \leq
C_\varphi(2) \sum_{n=1}^\infty \Pi^*_\varphi(A_i).$$ Hence, applying
(\ref{FP8.3.3}), the lemma  follows. \endpf

\sp

\bdfn $\HD(A)\index{(S)}{$\HD$}$, the Hausdorff dimension\index{(N)}{Hausdorff dimension} of the set $A$, is defined to be
\beq\label{FP8.4.1}
\HD(A):=\inf\{t:\H^t(A)=0\}=\sup\{t:\H^t(A)=\infty\}.
\eeq
Likewise, $\PD(A)\index{(S)}{$\PD$}$, the packing dimension\index{(N)}{packing dimension} of the set $A$, is defined to be
\beq\label{FP8.4.1B}
\PD(A):=\inf\{t:\Pi^t(A)=0\}=\sup\{t:\Pi^t(A)=\infty\}.
\eeq
\edfn

\sp

\fr The following theorem is the immediate consequence of the definition
of Hausdorff and packing dimensions, and the corresponding outer measures.

\sp

\bthm\label{tFPt8.4.1} The Hausdorff and packing dimensions are monotone increasing
functions of sets; that is, if $A\sbt B$ then 
$$
\HD(A)\le \HD(B) \  \  \text{ and } \  \  \PD(A)\le \PD(B).
$$
\ethm

\sp\fr We shall prove the following theorem, commonly refered to as the $\sg$-stability Hausdorff and packing dimensions.

\sp\bthm\label{FPt8.4.2} If $\{A_n\}_{n=1}^\infty$  is a countable
family of  subsets of $X$, then
$$
\HD\Big(\bigcup_{n=1}^\infty A_n\Big)=\sup_{n\ge 1}\{\HD(A_n)\}
$$
and 
$$
\PD\Big(\bigcup_{n=1}^\infty A_n\Big)=\sup_{n\ge 1}\{\PD(A_n)\}
$$

\ethm

\bpf We shall prove only the Hausdorff dimension part. The proof for the packing dimension is analogous. Inequality 
$$
\HD\Big(\bigcup_{n=1}^\infty A_n\Big)\ge \sup_{n\ge 1}\{\HD(A_n)\}
$$
is an immediate consequence of Theorem~\ref{tFPt8.4.1}.
Thus, if $\sup_n\{\HD(A_n)\}=\infty$, there is nothing to prove.
So, suppose that 
$$
s:=\sup_{n\ge 1}\{\HD(A_n)\}
$$ 
is finite, and consider an arbitrary $t>s$. Then, in view of (\ref{FP8.4.1}), $\H^t(A_n)=0$ for every $n \geq 1$, and therefore, since $\H^t$ is an outer measure, 
$$
\H^t\Big(\bigcup_{n=1}^\infty A_n\Big)=0.
$$
Hence, by (\ref{FP8.4.1}) again,
$$
\HD\Big(\bigcup_{n=1}^\infty A_n\Big) \leq t.
$$
The proof is complete. 
\endpf

\sp As an immediate consequence of this theorem,
Proposition~\ref{FPp8.2.1} and formula (\ref{FP8.4.1}) we get
following.

\bprop\label{FPp8.4.3} The Hausdorff dimension  of any  countable
set is equal to $0$. \eprop

\fr These  are the most basic, transparent, and also probably most useful  properties
of Hausdorff and packing measures and dimensions. We will apply  them frequently in the sequel. 

\sp\section{Hausdorff and Packing Measures; Frostmann's Converse Type Theorems}

\fr In this section we derive several geometrical consequences of
Theorem~\ref{t1hg47}. Their meaning is to tell us when a Hausdorff
measure or packing measure is positive, finite, zero, or infinity. We refer to them as Frostman's Converse Theorems. Somewhat strangely, these theorems are frequently called the Mass Redistribution Principle in the fractal literature, as if to indicate that such measures would always have to appear in the process of an iterative construction. At the end of the section we formulate Frostman's Direct Theorem and compare it with Frostman's Converse Theorems. As already mentioned, the advantage of the latte theorems is that these provide tools to calculate or at least to estimate, both
Hausdorff and packing measures and dimensions. We recall that in this section, as in the entire book, we keep 
$$
\varphi:[0, +\infty)\lra [0, +\infty),
$$
an evenly varying gauge function, i.e. satisfying formula \eqref{FP(8.2.3)}. We start with the following.

\sp

\index{(N)}{Frostman's Converse Theorem}
\bthm[Frostman's Converse Theorem for Generalized Hausdorff Measures]\label{t23hg53} Let $\varphi:[0, +\infty)\to [0, +\infty)$ be a
continuous evenly varying gauge function. Let $(X,\rho)$ be an
arbitrary metric space and $\mu$ a Borel probability measure on $X$.
Fix a Borel set $A\subseteq X$. Assume that there exists a constant $c\in (0, +\infty]$
\emph{($1/+\infty=0$)} such that
\begin{itemize}
\item[(1)]
$$
\limsup_{r\to 0}\frac{\mu (B(x, r))}{\varphi(r)}\ge c 
$$
for all points $x\in A$ except for countably many perhaps. Then the
Hausdorff measure $\H_{\varphi}$, corresponding to the gauge function
$\varphi$, satisfies 
$$
\H_{\varphi}(E)\le c^{-1} C_{\varphi}(8)\mu(E)
$$ 
for every Borel set $E\subseteq A$. In particular
$$
\H_{\varphi}(A)<+\infty \ \ (\H_{\varphi}(A)=0 \ \text{ if } \ c=+\infty),
$$
\item[(2)] If, conversly,
$$
\limsup_{r\to 0}\frac{\mu (B(x, r))}{\varphi(r)}\le c<+\infty
$$
for all $x\in A$, then 
$$
\mu(E)\le \H_{\varphi}(E)
$$ 
for every Borel
set $E\subseteq A$. In particular, 
$$
\H_{\varphi}(A)>0
$$ 
whenever $\mu(E)>0$. 
\end{itemize}
\ethm

\bpf Part (1). Since $\H_{\varphi}$ of any countable set is equal to
$0$, we may assume without loss of generality that $E$ does not
intersect the exceptional countable set. Fix $\varepsilon>0$. Then fix
$\delta>0$. Since measure $\mu$ is regular, there exists an open set
$G\supseteq E$ such that 
$$
\mu (G)\le\mu (E)+\varepsilon.
$$
Further, for every $x\in E$ there exists $r(x)\in (0, \delta)$ such that
$B(x, r(x))\subseteq G$ and 
$$
(c^{-1}+\varepsilon)\mu (B(x,r(x)))\ge\varphi (r(x)))>0.
$$
By virtue of both, $4r$--Covering Theorem, i.e. Theorem~\ref{t1hg47}, and of Remark~\ref{r1hg51}, there
exists $\{x_k\}_{k=1}^{\infty}$, a sequence of points in $E$, such that
$$
B(x_i, r(x_i))\cap B(x_j, r(x_j))=\emptyset\qquad\text{for $i\neq
j$}
$$
and
$$
\bigcup_{k=1}^{\infty} B(x_k, 4r(x_k))\supseteq\bigcup_{x\in E} B(x,
r(x))\supseteq E.
$$
Hence
$$
\begin{aligned}
{\H}^{2\delta}_{\varphi} (E)&\le\sum_{k=1}^{\infty} \varphi (2\cdot
4r(x_k))\le \sum_{k=1}^{\infty} c_{\varphi}(8)\varphi(r(x_k)) 
\le c_{\varphi}(8)\sum_{k=1}^{\infty}(c^{-1}+\varepsilon)\mu(B(x_k, r(x_k)))\\
&=c_{\varphi}(8)(c^{-1}+\varepsilon)\mu\(\bigcup_{k=1}^{\infty}B(x_k,r(x_k))\) 
\le c_{\varphi}(8)(c^{-1}+\varepsilon)\mu (G)\\ 
&\le c_{\varphi}(8)(c^{-1}+\varepsilon)(\mu(E)+\varepsilon).
\end{aligned}
$$
So, letting $\delta\downto 0$, we get,
$$
{\H}_{\varphi}(E)\le
c_{\varphi}(8)(c^{-1}+\varepsilon)(\mu(E)+\varepsilon)
$$
and since $\varepsilon>0$ was arbitrary, we finally get
$$
{\H}_{\varphi}(E)\le c_{\varphi}(8)c^{-1}\mu(E)
$$
This finishes the first part of the proof.

\sp\fr Part (2). Now we deal with the second part of our theorem. Fix an arbitrary $s>c$. Note that for every $r>0$ the function 
$$
X\ni x \longmapsto\frac{\mu (B(x, r))}{\varphi(r)}\quad\text{is Borel
measurable.}
$$
For every $k\ge 1$ consider the function
$$
X\ni x\longmapsto\varphi_k(x):=\sup\left\{\frac{\mu (B(x,r)}{\varphi(r)}:\,\,r\in
\mathbb Q\cap (0, 1/k]\right\},
$$
where $\mathbb Q$ denotes the set of rational numbers. This
function is Borel--measurable as the supremum of countably many
measurable functions. Let
$$
A_k=A\cap\varphi_k^{-1} ((0, s])\quad\text{for $k\ge 1$}.
$$
All $A_k$, $k\ge 1$, are then Borel subsets of $X$. Fix an arbitrary $r\in (0, 1/k)$. Then pick $r_j\searrow r$, $r_j\in
Q$. Since the function $t\mapsto\mu(B(x, t))$ is non--decreasing and
the function $\varphi$ is continuous, we get for every $x\in A_k$
that
$$
\frac{\mu(B(x, r))}{\varphi(r)}
\le\limsup\limits_{j\to\infty}\frac{\mu(B(x, r_j))}{\varphi(r_j)}\le s.
$$
Now fix $k\ge 1$. Then fix a Borel set $F\subseteq A_k$. Our first objective is to prove the assertion of (2) for the set $F$. To do this fix $r<1/k$ and then $\{F_i\}_1^{\infty}$, a countable cover of $F$ by subsets of $F$ that are closed relative to $F$ and have diameters lesser than
$r/2$. For every $i\ge 1$ pick $x_i\in F_i$. Then $F_i\subset B(x_i,
\diam (F_i))$. Since also all the sets $F_i$, $i\ge 1$, are Borel in $X$, we therefore have,
$$
\sum_{i=1}^{\infty}\varphi(\diam(F_i))\ge
s^{-1}\sum_{i=1}^{\infty}\mu(B(x_i, \diam (F_i)))\ge
s^{-1}\sum_{i=1}^{\infty}\mu(F_i)\ge s^{-1}\mu (F).
$$
Hence, invoking Remark~\ref{r1_mu_2014_11_07}, we get that
\beq\label{1_mu_2014_11_07}
{\H}_{\varphi}(F)\ge s^{-1}\mu(F).
\eeq
Moving further, by our hypothesis we have that
$$
\bigcup_{k=1}^{\infty} A_k\cap A=A.
$$
Define inductively
$$
B_1:=A_1\cap A
$$
and
$$
B_{k+1}:=A_{k+1}\cap \(A\setminus\bigcup_{j=1}^kA_j\cap A\).
$$
Obviously, the family $\{B_k\}_{1}^{\infty}$ consists of mutually disjoint sets, with each set $B_k$ contained in $A_k$, $k\ge 1$, and
$$
\bigcup_{k=1}^{\infty} B_k=\bigcup_{k=1}^{\infty} A_k=A.
$$
Hence, if $E$ is a Borel subset of $A$, then applying \eqref{1_mu_2014_11_07} for sets $F=E\cap B_k$, $k\ge 1$, we get
$$
{\H}_{\varphi} (E)=\bigcup_{k=1}^{\infty} {\H}_{\varphi} (E\cap
B_k)\ge s^{-1}\sum_{k=1}^{\infty}\mu(E\cap B_k)=s^{-1}\mu (E).
$$
Letting $s\searrow c$ then finishes the proof.
\endpf

\

\fr Now, let us prove the corresponding theorem for packing measures.

\index{(N)}{Frostman's Converse Theorem}
\bthm[Frostman's Converse Theorem for Generalized Packing Measures]\label{t24hg59} 
Let $\varphi: [0, \infty)\to [0, \infty)$ be a continuous evenly varying gauge function. Let $(X,\rho)$ be an arbitrary metric space and let $\mu$ be a Borel probability measure on $X$. Fix a Borel set $A\subset X$ and assume that there exists $c\in (0, +\infty]$, $(1/+\infty=0)$ such that
\begin{itemize}
\item[(1)]
$$
\liminf_{r\to 0}\frac{\mu(B(x, r))}{\varphi(r)}\le c\quad\text{for
all $x\in A$}. 
$$
Then 
$$
\mu(E)\le {\Pi}_{\varphi}(E)
$$ 
for every Borel set $E\subseteq
A$, where, we recall, ${\Pi}_{\varphi}$ denotes the packing measure
corresponding to the gauge function $\varphi$. In particular, if
$\mu(E)>0$, then 
$$
{\Pi}_{\varphi}(E)>0.
$$
\item[(2)]
If conversly,
$$
\liminf_{r\to 0}\frac{\mu(B(x, r))}{\varphi(r)}\ge c\quad\text{for
all $x\in A$},
$$
then 
$$
{\Pi}_{\varphi}(E)\le c^{-1}\mu(E)
$$ 
for every Borel set
$E\subseteq A$. In particular, if $\mu(E)<+\infty$, then
$$
{\Pi}_{\varphi}(E)<+\infty. 
$$
\end{itemize}
\ethm

\bpf Part (1). Let $\varepsilon>0$. Fix an arbitrary subset $F\subseteq A$. Define a decreasing sequence
$(G_n)_{n\geq1}$ of open sets containing $F$ as follows. By our hypothesis, for every $x\in A$ there exists $0<r_1(x)<1$ such that
\[
\frac{\mu\bigl(B(x,r_1(x))\bigr)}{\phi(r_1(x))}\leq c+\varepsilon.
\]
Take the family of balls $\bigl\{B\bigl(x,\frac{1}{4}r_1(x)\bigr)\bigr\}_{x\in F}$. According
to the $4r$-Covering Theorem (Theorem~\ref{t1hg47}), there is a countable set $F_1\subseteq F$ such
that the subfamily $\bigl\{B\bigl(x,\frac{1}{4}r_1(x)\bigr)\bigr\}_{x\in F_1}$ consists of
mutually disjoint balls satisfying
\[
F\subseteq\bigcup_{x\in F}B\Bigl(x,\frac{1}{4}r_1(x)\Bigr)\subseteq\bigcup_{x\in F_1}B\bigl(x,r_1(x)\bigr).
\]
Let $G_1:=\bigcup_{x\in F_1}B(x,r_1(x))$. For the inductive step, suppose that $G_n$ has been defined
for some $n\geq1$. By our hypothesis again, for every $x\in A$ there exists some $0<r_{n+1}(x)<\frac{1}{n+1}$ such that
$B(x,r_{n+1}(x))\subseteq G_n$ and
\begin{eqnarray}\label{thm3.1starpack}
\frac{\mu\bigl(B(x,r_{n+1}(x))\bigr)}{\phi(r_{n+1}(x))}\leq c+\varepsilon.
\end{eqnarray}
Consider the family of balls $\bigl\{B\bigl(x,\frac{1}{4}r_{n+1}(x)\bigr)\bigr\}_{x\in F}$. According
to the $4r$-Covering Theorem (Theorem~\ref{t1hg47}), there exists a countable set $F_{n+1}\subseteq F$ such
that the subfamily $\bigl\{B\bigl(x,\frac{1}{4}r_{n+1}(x)\bigr)\bigr\}_{x\in F_{n+1}}$ consists of
mutually disjoint balls satisfying
\[
F\subseteq\bigcup_{x\in F}B\Bigl(x,\frac{1}{4}r_{n+1}(x)\Bigr)\subseteq\bigcup_{x\in F_{n+1}}B\bigl(x,r_{n+1}(x)\bigr).
\]
Let 
$$
G_{n+1}:=\bigcup_{x\in F_{n+1}}B(x,r_{n+1}(x)).
$$
It is clear that  $G_{n+1}$ is an open set and $F\subseteq G_{n+1}\subseteq G_n$.
Moreover, for all pairs $x,y\in F_{n+1}\subseteq F$ we know that \[d(x,y)\geq\frac{1}{4}\max\{r_{n+1}(x),r_{n+1}(y)\}\geq
\frac{1}{8}(r_{n+1}(x)+r_{n+1}(y)).\]
Therefore the collection $\{(x,\frac{1}{8}r_{n+1}(x))\}_{x\in F_{n+1}}$ forms an $(\frac{1}{n+1})$-packing of $F$.
Using~(\ref{thm3.1starpack}), it follows that
\begin{eqnarray*}
\Pi_\phi^{*\frac{1}{n+1}}(F)
&\geq&\sum_{x\in F_{n+1}}\phi\Bigl(\frac{1}{8}r_{n+1}(x)\Bigr)
\geq(C_\phi(8))^{-1}\sum_{x\in F_{n+1}}\phi(r_{n+1}(x))\\
&\geq&(C_\phi(8))^{-1}\sum_{x\in F_{n+1}}\frac{\mu\bigl(B(x,r_{n+1}(x))\bigr)}{c+\varepsilon}\\
&=&\frac{(C_\phi(8))^{-1}}{c+\varepsilon}\mu\Bigl(\bigcup_{x\in F_{n+1}}B(x,r_{n+1}(x))\Bigr)\\
&\geq&\frac{(C_\phi(8))^{-1}}{c+\varepsilon}\mu(G_{n+1}).
\end{eqnarray*}
Letting $n$ increase to infinity, we thus obtain that
\[
\Pi_\phi^*(F)
\geq\frac{(C_\phi(8))^{-1}}{c+\varepsilon}\inf_{n\geq1}\mu(G_n)
=\frac{(C_\phi(8))^{-1}}{c+\varepsilon}\lim_{n\to\infty}\mu(G_n)
=\frac{(C_\phi(8))^{-1}}{c+\varepsilon}\mu(G_F),
\]
where $G_F:=\cap_{n\geq1}G_n$ is a $G_\delta$-set and therefore, in particular, is a Borel set.
Consequently, for every Borel set $E\subseteq A$ we have
\begin{eqnarray*}
\Pi_\phi(E)
&=&\inf\Bigl\{\sum_{k=1}^\infty \Pi_\phi^*(A_k):\{A_k\}_{k=1}^\infty \mbox{ is a cover of }E\Bigr\}\\
&=&\inf\Bigl\{\sum_{k=1}^\infty \Pi_\phi^*(A_k):\{A_k\}_{k=1}^\infty \mbox{ is a partition of }E\Bigr\}\\
&\geq&\frac{(C_\phi(8))^{-1}}{c+\e}\inf\Bigl\{\sum_{k=1}^\infty\mu(G_{A_k}):\{A_k\}_{k=1}^\infty \mbox{ is a partition of }E\Bigr\}\\
&\geq&\frac{(C_\phi(8))^{-1}}{c+\e}\inf\Bigl\{\mu\Bigl(\bigcup_{k=1}^\infty G_{A_k}\Bigr):\{A_k\}_{k=1}^\infty \mbox{ is a partition of }E\Bigr\}\\
&\geq&\frac{(C_\phi(8))^{-1}}{c+\e}\inf\Bigl\{\mu(E):\{A_k\}_{k=1}^\infty \mbox{ is a partition of }E\Bigr\}\\
&=&\frac{(C_\phi(8))^{-1}}{c+\e}\mu(E).
\end{eqnarray*}
Since this holds for all $\varepsilon>0$, we deduce that $\Pi_\phi(E)\geq(C_\phi(8)c)^{-1}\mu(E)$. The proof of Part (1) is complete.

\sp\fr Part (2). The sequence of functions $(\psi_k)_{k=1}^\infty$, where
\[
X\ni x\longmapsto\psi_k(x):=\inf\left\{\frac{\mu(B(x,r))}{\phi(r)}:r\in\Q\cap\left(0,\frac1k\right]\right\}
\]
forms an increasing sequence of measurable functions.
Let $0<s<c$. For each $k\geq1$, let 
$$
A_k:=\psi_k^{-1}([s,+\infty)).
$$
As $(\psi_k)_{k=1}^\infty$ is increasing, so is the sequence $(A_k)_{k=1}^\infty$.
Moreover, since $s<c$, it follows from our hypothesis that 
$$
\bigcup_{k=1}^\infty A_k\supseteq A.
$$
Furthermore, since $[s,+\infty)$ is a Borel subset of $\R$, the measurability of $\psi_k$
ensures that $A_k$ is a Borel subset of $X$. Fix $k\geq1$. Choose some arbitrary $r\in(0,1/k]$
and pick a sequence $(r_j)_{j\geq1}\in\Q$ such that $r_j$ increases to $r$. Since $\mu$ is a measure and $\phi$ is continuous, we deduce that  for all $x\in A_k$,
\[
\frac{\mu(B(x,r))}{\phi(r)}
=\lim_{j\to\infty}\frac{\mu(B(x,r_j))}{\phi(r_j)}
\geq\psi_k(x)
\geq s.
\]
Thus, if $x\in A_k$ then 
$$
\inf\Bigl\{\frac{\mu(B(x,r))}{\phi(r)}:r\in(0,1/k]\Bigr\}\geq s.
$$
Now, fix any set $F\subseteq A_k$, any $r\in(0,\frac1k]$, and let
$\{(x_i,t_i)\}_{i\geq1}$ be an $r$--packing of $F$. Then
\[
\sum_{i=1}^\infty\phi(t_i)
\leq s^{-1}\sum_{i=1}^\infty\mu(B(x_i,t_i))
=s^{-1}\mu\Bigl(\bigcup_{i=1}^\infty B(x_i,t_i)\Bigr)
\leq s^{-1}\mu(F_r),
\]
where $F_r$ denotes the open $r$--neighborhood of $F$.
Taking the supremum over all $r$--packings yields
$$
\Pi_\phi(F)\leq \Pi_\phi^*(F)\leq \Pi_\phi^{*r}(F)\leq s^{-1}\mu(F_r).
$$
Thus, we have that  $P_\phi(F)\leq s^{-1}\mu(F_r)$ for all $r\in(0,{1}/{k}]$ and each subset $F\subseteq A_k$.
Consequently, $\Pi_\phi(F)\leq s^{-1}\mu(F_0)=s^{-1}\mu(\overline{F})$ for all $F\subseteq A_k$.
In particular, if $C$ is a closed subset of $E$, then
\[\Pi_\phi(C\cap A_k)\leq s^{-1}\mu(\overline{C\cap A_k})\leq s^{-1}\mu(C)\leq s^{-1}\mu(E).\]
As this holds for all integers $k\geq1$ and closed sets $C\subseteq E\subseteq A=\cup_{k\geq1}A_k$,
we deduce that $\Pi_\phi(C)\leq s^{-1}\mu(E)$. By the regularity of $\mu$,  on taking the supremum
over all closed sets $C$ contained in $E$, we conclude that 
$$
\Pi_\phi(E)\leq s^{-1}\mu(E).
$$
Letting $s$ increase to $c$ finishes the proof.
\endpf

\sp Replacing $\varphi(r)$  by $r^t$, $\H_\varphi$ by  $\H^t$ and
$\Pi_\varphi$ by $\Pi^t$  in Theorem~\ref{t23hg53} and
Theorem~\ref{t24hg59} respectively,  we immediately get the
following two results.

\

\index{(N)}{Frostman's Converse Theorem}
\bthm[Frostman's Converse Theorem for Hausdorff Measures]\label{tncp12.1.} Fix $t>0$ arbitrary. Let $(X,\rho)$ be a
metric space and $\mu$ a Borel probability measure on $X$.
Fix a Borel set $A\subseteq X$. Assume that there exists a constant $c\in (0, +\infty]$
\emph{($1/+\infty=0$)} such that
\begin{itemize}
\item[(1)]
$$
\limsup_{r\to 0}\frac{\mu (B(x, r))}{r^t}\ge c 
$$
for all points $x\in A$ except for countably many perhaps. Then the
Hausdorff measure $\H_t$ satisfies 
$$
\H_t(E)\le c^{-1} 8^t\mu(E)
$$ 
for every Borel set $E\subseteq A$. In particular
$$
\H_t(A)<+\infty \ \ (\H_t(A)=0 \ \text{ if } \ c=+\infty),
$$
\item[(2)] If, conversly,
$$
\limsup_{r\to 0}\frac{\mu (B(x,r))}{r^t}\le c<+\infty
$$
for all $x\in A$, then 
$$
\mu(E)\le \H_t(E)
$$ 
for every Borel
set $E\subseteq A$. In particular, 
$$
\H_t(A)>0
$$ 
whenever $\mu(E)>0$. 
\end{itemize}
\ethm

\fr and

\index{(N)}{Frostman's Converse Theorem}
\bthm[Frostman's Converse Theorem for Packing Measures]\label{tncp12.2.}
Fix $t>0$ arbitrary. Let $(X,\rho)$ be ametric space and $\mu$ a Borel probability measure on $X$. Fix a Borel set $A\subset X$ and assume that there exists $c\in (0, +\infty]$, $(1/+\infty=0)$ such that
\begin{itemize}
\item[(1)]
$$
\liminf_{r\to 0}\frac{\mu(B(x,r))}{r^t}\le c\quad\text{for
all $x\in A$}. 
$$
Then 
$$
\mu(E)\le {\Pi}_t(E)
$$ 
for every Borel set $E\subseteq A$. In particular, if $\mu(E)>0$, then 
$$
{\Pi}_t(E)>0.
$$
\item[(2)]
If conversly,
$$
\liminf_{r\to 0}\frac{\mu(B(x, r))}{r^t}\ge c\quad\text{for
all $x\in A$},
$$
then 
$$
{\Pi}_t(E)\le c^{-1}\mu(E)
$$ 
for every Borel set
$E\subseteq A$. In particular, if $\mu(E)<+\infty$, then
$$
{\Pi}_t(E)<+\infty. 
$$
\end{itemize}
\ethm

\fr In the opposite direction to Frostman Converse Theorems, there is the following well known:

\bthm[Frostman Direct Lemma]\label{Frostman_Direct}
If $X$ be either a Borel subset of a Euclidean space $\R^d$, $d\ge 1$, or an arbitrary compact metric space. If $t>0$ and $\H_t(X)>0$, then there exists a Borel probability measure $\mu$ on $X$ such that
$$
\mu(B(x,r))\le r^t
$$
for every point $x\in X$ and all radii $r>0$.
\ethm

\fr This is a very interesting theorem although Frostman Converse Theorems seem to be more suitable to estimate and calculate Hausdorff and packing measures and dimensions. 

\section{Hausdorff Dimension of Measures} 

In this section we define the concepts of dimensions, both Hausdorff and packing, of Borel measures. We then provide tools to calculate and to estimate them. We also establish some relations between them. The dimensions of measures \index{(N)}{dimension of measure} play an important role in both fractal geometry and dynamical systems. We start this section with the following simple but crucial consequence of Theorem~\ref{tncp12.1.}. 

\

\index{(N)}{Volume Lemma}
\bthm[Volume Lemma for Hausdorff Measures]\label{th:6.6.3} 
Suppose that $\mu$ is a Borel probability measure on a metric space
$(X,\rho)$ and that $A$ is a bounded Borel subset of $\R^n$. Then
\begin{itemize}
\item[(a)] If $\mu(A)>0$ and there exists $\th_1$ such that for every $x\in A$
$$
\liminf_{r\to 0}\frac{\log\mu(B(x,r))}{\log r} \ge \th_1
$$
then $\HD(A)\ge \th_1$.

\sp\item[(b)] If there exists $\th_2$ such that for every $x\in A$
$$
\liminf_{r\to 0}\frac{\log\mu(B(x,r))}{\log r} \le \th_2
$$
then $\HD(A)\le \th_2$.
\end{itemize}
\ethm

\bpf (a) Take any $0<\th<\th_1$. Then, by the assumption,
$$
\limsup_{r\to 0}\mu(B(x,r))/r^\th =0.
$$ 
It therefore follows from Theorem~\ref{tncp12.1.} (2) that $\H^\th(A)=+\infty$.
Hence $\HD(A)\ge\th$. Consequently, $\HD(A)\ge\th_1$.

(b) Take now an arbitrary $\th>\th_2$. Then by the assumption 
$$
\limsup_{r\to 0}\mu(B(x,r))/r^\th=+\infty.
$$
Therefore applying Theorem~\ref{tncp12.1.} (1), we obtain
$\H^\th(A)=0$. Thus $\HD(A)\le\th$ and consequently 
$\HD(A)\le\th_2$. The proof is finished.
\epf

\

\fr Similarly, one proves the following consequence of
Theorem~\ref{tncp12.2.}. 

\

\index{(N)}{Volume Lemma}
\bthm[Volume Lemma for Packing Measures]\label{6.6.3pack}
Suppose that $\mu$ is a Borel probability measure on
$\R^n$, $n\ge 1$, and $A$ is a bounded Borel subset of $\R^n$.
\begin{itemize}
\item[(a)] If $\mu(A)>0$ and there exists $\th_1$ such that for every $x\in A$
$$
\limsup_{r\to 0}\frac{\log\mu(B(x,r))}{\log r} \ge \th_1
$$
then $\PD(A)\ge \th_1$.

\item[(b)] If there exists $\th_2$ such that for every $x\in A$
$$
\limsup_{r\to 0}\frac{\log\mu(B(x,r))}{\log r} \le \th_2
$$
then $\PD(A)\le \th_2$.
\end{itemize}
\ethm

We will now apply Theorem~\ref{th:6.6.3} (a) to get quite a general lower bound for Hausdorff dimension, the one which is a generalization of a result due to McMullen (\cite{McM1}) and whose proof is taken from \cite{MU_Non-Dense-Orbits}. Although, this result is usually applied in a dynamical context, it really does not require any dynamics to formulate and to prove it. 

\sp As always in this section, let $(X,\rho)$ be a metric space and let $\mu$ be a Borel probability upper Ahlfors measure \index{(N)}{upper Ahlfors measure} on $X$, meaning there exist constants $h>0$ and $C \geq 1$  such that  for every $x \in X$  and $r>0$
\beq\label{2.1}
\mu(B(x, r)) \leq C r^h.
\eeq
We then call $h$ the exponent of $\mu$. For any integer $k \geq 1$ let $E_k$ be a finite collection of compact subsets of $X$ each element of which has positive measure $\mu$. We denote:
\beq\label{2.2}
K:=\bigcup_{F\in E_1}F.
\eeq
We assume that
\beq\label{2.3}
\text{If $k\ge 1$ and $F, G \in E_k$ and $F \neq G$, then $\mu(F \cap G) =0$.}
\eeq
\beq\label{2.4}
\text{Every set $F\in E_{k+1}$  is contained in a unique element $ G \in E_k$.}
\eeq
For every integer $k\ge 1$ and every set $F\in E_k$, define 
 \beq\label{2.5}
 \text{density}\lt(\bigcup_{D\in E_{k+1}}\!\!\! D, F\rt)
 :=\frac{\mu\lt(D\cap \bigcup_{D\in E_{k+1}}D \rt)}{\mu(F)}
\eeq
and assume that
\beq\label{2.6}
\Delta_k:= \inf\lt\{\text{density}\lt(\bigcup_{D\in E_{k+1}}\!\!\!D, F\rt): \, F\in E_k\rt\} >0
\eeq
for every $k \geq 1$. Put also
$$
d_k:= \sup\big\{\diam(F): \, F\in E_k\big\}.
$$
Suppose that $d_k<1$   for every $k \geq 1$ and that
\beq\label{2.7}
\lim_{k \to \infty} d_k=0.
\eeq
We then call the collection
$$
\{E_k\}_{k=1}^\infty,
$$
a McMullen's sequence of sets. Let
$$
E_\infty:=\bigcap_{k=1}^\infty \bigcup_{F\in E_k}\!\!\! F.
$$
We shall prove the following generalization of a McMullen's result from \cite{McM1}) whose proof is taken from \cite{MU_Non-Dense-Orbits}

\bprop\label{p1_mcmullen_HD_estimate}
If $\{E_k\}_{k=1}^\infty$ is a McMullen's sequence of subsets of a metric space $(X,\rho)$ endowed with a Borel probability upper Ahlfors measure $\mu$ having exponent $h$, then
$$
\HD (E_\infty) \geq h - \limsup_{k \to \infty}\frac{\sum_{j=1}^{k-1} \log \Delta_j}{\log d_k}.
$$
\eprop

\bpf  We construct inductively  a sequence of Borel probability measures $\{\nu_k\}_{k=1}^\infty$ on $K$  as follows.

Put $\nu_1:=\mu$ and define $\nu_{k+1}$ by putting for each Borel set $A \sbt K$
\beq\label{2.8}
\nu_{k+1}(A)
:=\sum_{F\in E_k} \frac{\mu\lt(A\cap F \cap \bigcup_{D\in E_{k+1}}\! D\rt)}{\mu\lt(F \cap \bigcup_{D\in E_{k+1}}\! D\rt)}\nu_k(F).
\eeq
This  definition makes sense since by (\ref{2.5}) and (\ref{2.6}) we see that $\mu\lt(F \cap \bigcup_{D\in E_{k+1}}\! D\rt)>0$.
By induction we get for every $k\geq 1$ that 
\beq\label{2.9}
\nu\lt(\bigcup_{D\in E_k}\!\!\! D\rt)=1,
\eeq
and it follows  from properties (\ref{2.2})-(\ref{2.4}) that $\nu_{k+1}$   is a Borel probability measure indeed. In view of (\ref{2.8}) and (\ref{2.3}) we have that  $\nu_{k+1}(F)=\nu_k(F)$  for each $F\in E_k$. Hence using (\ref{2.3})  and (\ref{2.4}) we  conclude by  induction  that $\nu_n(F)=\nu_k(F)$ for every $n \geq k$. Since $\lim_{k \to \infty} d_k=0$, we therefore obtain a unique  probability measure $\nu$ on $K$
(being the weak limit of measures $\nu_k$) such that
 \beq\label{2.10}
 \nu(F)=\nu_k(F)
 \eeq
for every $F \in E_k$. Looking now at (\ref{2.9}) and the definition of the set $E$, one gets
 \beq\label{2.11}
 \nu(E_\infty)=1.
 \eeq
Making use of (\ref{2.8}) and (\ref{2.10}) one easily estimates for every $F\in E_k$ that 
 \beq\label{2.12}
\nu(F) \leq \frac{\mu(F)}{\Delta_{k-1}\ldots \Delta_1}.
 \eeq
In view  of Theorem~\ref{th:6.6.3}, Volume Lemma for Hausdorff Measures, in order to prove that $\HD(E)\geq \delta$ for some $\delta \geq 0$ it is enough to show that 
  \beq\label{2.13}
  \liminf_{r \to 0} \frac{\log \nu (B(x,r))}{\log r} \geq \d
  \eeq
  for $\nu$--$a.e.$ $x \in E$. 
  
Now consider $x \in E_\infty$  and $0<r < \sup_{k \geq 1} (d_k)$ arbitrary. Then there exists an integer $k=k(r)\ge 1$ such that $d_{k+1}\leq r \leq d_k$. Let $\tilde{B}(x,r)$ be the union of all sets in $E_{k+1}$ which meet $B(x,r)$. Then $ \tilde{B}(x,r) \sbt B(x, 2r)$  and using (\ref{2.12})  and (\ref{2.1}) we get
  \begin{eqnarray*}\begin{aligned}
  \frac{\log \nu (B(x, r))}{\log r}&  \geq \frac{ \log \mu(\tilde{B}(x, r))- \sum_{j=1}^{k-1} \log \Delta_j}{
  \log r}\\
      & \geq \frac{ \log C + h \log 2 + h \log r}{\log r}- \frac{\sum_{j=1}^{k-1}\log \Delta_j}{\log d_j}.
       \end{aligned}
  \end{eqnarray*}
Since $\lim_{r \to 0} k(r) =\infty$, we therefore  obtain
$$ \liminf_{r \to 0} \frac{\log \nu (B(x, r))}{\log r}\geq  h- \limsup_{k \to \infty}\frac{\sum_{j=1}^{k-1} \log \Delta_j}{\log d_k}.
$$
In view of (\ref{2.11}) and by applying Theorem~\ref{th:6.6.3} (a), this finishes the  proof.
\qed

\sp Now we define the following main concepts of this section.

\bdfn\label{d6.4.10}
Let $\mu$ be a Borel measure on a metric space $(X,\rho)$. We write
$$
\HD_\star(\mu):=\inf\{\HD(Y): \mu( Y)>0\}\ \; \  {\text{and}}\ \;  \
\HD^\star(\mu)=\inf\{\HD(Y): \mu(X\sms Y)=0\}.
$$
Of course 
$$
\HD_\star(\mu)\le \HD^\star(\mu),
$$
and, in the case when $\HD_\star(\mu)=\HD^\star(\mu)$, we call this common
value the \index{(N)}{Hausdorff dimension!of measure} 
{\it Hausdorff dimension} of the measure $\mu$ and we denote it by $\HD(\mu)$.
\edfn

An analogous definition can be formulated for packing dimension,
with respective notation $\PD_\star(\mu)$, $\PD^\star(\mu)$, $\PD(\mu)$ and the name \index{(N)}{packing dimension of measure} {\it packing dimension} of the
measure~$\mu$. 

\

\fr The next definition introduces concepts that are effective tools to
calculate the dimensions introduced above. 

\

\bdfn\label{d220130510}
Let $\mu$ be a Borel probability measure on a metric space $(X,\rho)$.
For every point $x\in X$ we define the {\it lower and upper pointwise
dimension} \index{(N)}{dimension pointwise lower}
\index{(N)}{dimension pointwise upper} \!\!\!\!\!
of the measure $\mu$ at the point $x\in X$ respectively as
$$
\un d_\mu(x):=\liminf_{r\to 0} {\log\mu(B(x,r))\over \log r}\  \  \ {\text {and} }\  \  \
\ov d_\mu(x):=\limsup_{r\to 0} {\log\mu(B(x,r))\over \log r}.
$$
In the case when both numbers $\un d_\mu(x)$ and $\ov d_\mu(x)$ are
equal, we denote their common value by $d_\mu(x)$. We then obviously have
$$
d_\mu(x)=\lim_{r\to 0}{\log\mu(B(x,r))\over \log r},
$$
and we call $d_\mu(x)$ the pointwise dimension \index{(N)}{dimension pointwise} of the measure $\mu$ at the point $x\in X$. 
\edfn

\sp The following theorem about Hausdorff and packing dimensions of a
Borel measure $\mu$, follows easily from Theorems~\ref{th:6.6.3} and
~\ref{6.6.3pack}. 

\sp\bthm\label{t6.6.4.}
If $\mu$ is a Borel probability measure on a metric space $(X,\rho)$, then
$$
\HD_\star(\mu)=\ess\inf \un d_\mu,\ \ \ \HD^\star(\mu)=\ess\sup \un
d_\mu
$$
and
$$
\PD_\star (\mu)=\ess\inf\ov d_\mu, \ \ \ \PD^\star (\mu)=\ess\sup\ov d_\mu.
$$
\ethm
\bpf 
Recall that the $\mu$--essential infimum ess\,inf of a measurable function $\phi$
and the $\mu$-essential supremum ess\,sup of this function
are respectively defined by
$$
\ess\inf(\phi):=\sup_{\mu(N)=0}\inf_{x\in X\setminus N}\phi(x)\ \ \ {\text{and}}\ \ \
\ess\sup(\phi):=\inf_{\mu(N)=0}\sup_{x\in X\setminus N}\phi(x).
$$
Put $\phi_*::=\ess\inf \phi$. We shall prove that
\beq\label{220130510}
\mu(\phi^{-1}((0,\phi_*))=0\ \ \ {\text{ and}}\ \ \
\mu(\phi^{-1}((0,\th))>0
\eeq
for all $\th> \phi_*$. Indeed, if we had
$\mu(\phi^{-1}((0,\phi_*))>0$, then there would 
exist $\th<\phi_*$ with $\mu(\phi^{-1}((0,\th])>0$. Hence, for every
measurable set $N\sbt X$ with $\mu(N)=0$, we would have
$\inf_{X\setminus N}\phi\le\th$. Thus $\ess\,\inf \phi\le\th$, which
is a contradiction, and the first part of \eqref{220130510} is proved.

\sp\fr For the second part, proceeding also by the way of contradiction,
assume that there exists $\th> \phi_*$ with
$$
\mu(\phi^{-1}((0,\th))=0.
$$
Then for $N:=\phi^{-1}((0,\th))$, we would
have $\inf_{X\setminus N}(\phi)\ge\th$.  Hence $\ess\,\inf \phi\ge\th$,
which is a contradiction, and this finishes the proof of formula \eqref{220130510}.

\sp\fr This formula, applied to the function $\phi:=\un d_\mu$, tells us
that for every Borel set $A\sbt X$, with $\mu(A)>0,$ there exists a
Borel set $A'\sbt A$ with $\mu(A')=\mu(A)>0$ such that for every $x\in
A'$ we have  $\un d_\mu(x)\ge d_{\mu*}$. Hence,
$$
\HD(A)\ge\HD(A')\ge d_{\mu*}
$$ 
by Theorem~\ref{th:6.6.3}(a). Thus, 
\beq\label{420130510}
\HD_\star(\mu)\ge d_{\mu*}.   
\eeq
On the other hand for every $\th>\th_1$ we have $\mu\(\{x\in X:\un
d_\mu(x)<\th\}\)>0$. Hence, by
Theorem~\ref{th:6.6.3}(b), 
$$
\HD(\{x:\un d_\mu(x)<\th\})\le\th.
$$
Therefore $\HD_\star(\mu)\le\th$. So, letting $\th\downto\th_1$, we
get 
$$
\HD_\star(\mu)\le d_{\mu*}.
$$
Along with \eqref{420130510} we thus
conclude that $\HD_\star(\mu)=d_{\mu*}$.

\sp\fr One proceeds similarly to prove that $\HD^\star(\mu)=\ess\sup
\un d_\mu(x)$ and to obtain corresponding results for 
packing dimension. For the latter, one refers to
Theorem~\ref{6.6.3pack} instead of Theorem~\ref{th:6.6.3}.
\endpf

\bdfn\label{d220191130}
A Borel probability measure  $\mu$ on a metric space $(X,\rho)$ is called dimensional exact \index{(N)}{dimensional exact} if and only if for $\mu$--a.e. $x\in X$, $d_\mu(x)$, the pointwise dimension \index{(N)}{dimension pointwise} of the measure $\mu$ at $x$, exists and is $\mu$ almost everywhere constant.
\edfn

\sp As an immediate consequence of Theorem~\ref{t6.6.4.}, we get the following.

\sp\bprop\label{p320191130}
If $\mu$ is a Borel probability dimensional exact measure on a metric space $(X,\rho)$, then both $\HD(\mu)$ and $\PD(mu)$ exist, and moreover
$$
\HD(\mu)=\PD(\mu)=d_\mu,
$$
where $d_\mu$ is the $\mu$ almost everywhere constant, value of the pointwise dimension of $\mu$.
\eprop

\section{Box--Counting Dimensions}

\fr We shall now examine a slightly different type of dimension, namely, the
box--counting dimension. This dimension, as we will shortly see, is not given by means of any outer measure. It behaves worse: it is not $\sg$-- stable and a set and its closure have the same box--counting dimension. Its definition is however substantially simpler than those of Hausdorff and packing dimensions, is frequently easier to calculate or to estimate, it also frequently agrees with Hausdorff and packing dimension, and is widely used in physical literature. 

\sp\bdfn\label{d120141231}
Let $0<r<1$, and let $A\subseteq X$ be a bounded set.
Define $N(A,r)$ to be the minimum number of balls of radius at most $r$
with centers in $A$ needed to cover $A$. Then the {\em upper and lower 
box--counting}
(or, simpler, {\em box}) dimensions of $A$ are respectively defined to be
\[
\overline{\BD}(A):=\limsup_{r\to0}\frac{\log N(A, r)}{-\log r}
\]
and
\[
\underline{\BD}(A):=\liminf_{r\to0}\frac{\log N(A, r)}{-\log r}.
\]
If these two quantities are equal, their common value is called  the \index{(N)}{box dimension}{\em box--counting dimension}, or simply {\em box dimension},
 of $A$, and we denote it by $\BD(A)$.  
\edfn

\fr As said, the box--counting dimensions do not share all the congenial properties of the
Hausdorff dimension. In particular, they are not $\sigma$--stable.
To see this, observe that
$$
\BD(\mathbb Q\cap[0,1])=1\neq0=\sup\{\BD(\{q\}):q\in\mathbb Q\cap[0,1]\}.
$$
Box--counting dimension is however easily seen to be finitely stable, i. e.

\sp\bprop\label{p720131014B}
If $(X,\rho)$ is a metric space and $F_1$, $F_2$,\ld,$F_n$ is a finite collection of subsets of $X$, then
$$
\ov{\BD}\(F_1\cup F_2\cup\ld\cup F_n\)
=\max\big\{\ov{\BD}(F_1), \ov{\BD}(F_2),\ld, \ov{\BD}(F_n)\big\}
$$
and the same formula holds for the lower box--counting dimension.
\eprop

\fr The terminology ``box--counting'' comes from the fact that in
Euclidean spaces we may use boxes from a lattice  rather
than balls to cover the set under scrutiny. Indeed, let $n\geq1$, $X=\R^n$
and let $\mathcal{L}(r)$ be any lattice in $\R^n$ consisting of cubes (boxes)
with edges of length $r$. For any $A\subseteq X$, define
$$
L(A,r)=\mathrm{card}\{C\in\mathcal{L}(r):C\cap A\neq\emptyset\}.
$$
\bprop\label{sem3prop2.6}
If $A $ is a bounded subset of $\R^n$, then
\[
\overline{\BD}(A)=\limsup_{r\to0}\frac{\log L(A,r)}{-\log r}
\]
and
\[
\underline{\BD}(A)=\liminf_{r\to0}\frac{\log L(A,r)}{-\log r}.
\]
\eprop
\bpf
Without loss of generality, let $0<r<1$. Select points $x_i\in A$ so that
\[A\subseteq\bu_{i=1}^{N(A,r)}B(x_i,r).\] Fix $1\leq i\leq N(A,r)$ momentarily.
If $C\in\mathcal{L}(r)$ is such that $d(C,x_i)<r$, where $d$ denotes the standard Euclidean
metric on $\R^n$, one immediately verifies that $C\subseteq B(x_i,r+r\sqrt{n})=B(x_i,(1+\sqrt{n})r)$. Thus,
for any given $1\leq i\leq N(A,r)$, we have that
$$
\begin{aligned}
\#\{C\in\mathcal{L}(r):d(C,x_i)<r\}
&=\frac{\lambda(B(x_i,(1+\sqrt{n})r))}{\lambda(\mbox{cube of side $r$})}
\leq\frac{c_n\bigl[(1+\sqrt{n})r\bigr]^n}{r^n}\\
&=c_n(1+\sqrt{n})^n,
\end{aligned}
$$
where $\lambda$ denotes the Lebesgue measure on $\R^n$ and $c_n$ denotes the volume of the unit ball in $\R^n$.
Since every $C\in L(A,r)$ admits at least one $1\leq i\leq N(A,r)$ such that $d(C,x_i)<r$,
we deduce that $L(A,r)\leq N(A,r)c_n(1+\sqrt{n})^n$. Therefore
\[
\log L(A,r)\leq\log\bigl(c_n(1+\sqrt{n})^n\bigr)+\log N(A,r).
\]
Hence
\[
\frac{\log L(A,r)}{-\log r}
\leq\frac{\log\bigl(c_n(1+\sqrt{n})^n\bigr)}{-\log r}+\frac{\log N(A,r)}{-\log r}.
\]
So
\[
\limsup_{r\to0}\frac{\log L(A,r)}{-\log r}\leq\ov{\BD}(A)
\ \text{ and }\
\liminf_{r\to0}\frac{\log L(A,r)}{-\log r}\leq\un{\BD}(A).
\]
For the opposite inequality, again let $0<r<1$ and for each $C\in L(A,r)$
choose $x_C\in C\cap A$. Then $C\cap A\subseteq B(x_C,r\sqrt{n})$. Thus, the family
of balls 
$$
\big\{B(x_C,r\sqrt{n}):C\in L(A,r)\big\}
$$
covers $A$. Therefore
$N(A,r\sqrt{n})\leq L(A,r)$. It then follows that
\[
\ov{\BD}(A)\leq\limsup_{r\to0}\frac{\log L(A,r)}{-\log r}
\ \text{ and }\
\un{\BD}(A)(A)\leq\liminf_{r\to0}\frac{\log L(A,r)}{-\log r}.
\]
\epf

\sp Now, we return to our general setting. Let $(X,\rho)$ be again a metric space,
let $A\subseteq X$ and let $r>0$. Further, define $P(A,r)$ to be the supremum of
the cardinalities of all packings of $A$ of the form
$\{B(x_i,r)\}_{i=1}^\infty$, so
\[
P(A,r):=\sup\{\#\{B(x_i,r)\}_{i=1}^\infty\}.
\]
Such packings will be called in the sequel $r$--packings of $A$. We shall prove the following technical, though interesting in itself, fact.

\sp\blem\label{sem3lem2.7}
If $A\subseteq X$ and $r>0$ then $N(A,2r)\leq P(A,r)\leq N(A,r)$.
\elem
\bpf
The first inequality certainly holds if $P(A,r)=\infty$. So assume this is not
the case and let $\{(x_i,r)\}_{i=1}^k$ be a $r$-packing of $A$
which is maximal in the sense of inclusion. Then $\{B(x_i,2r)\}_{i=1}^k$ is a cover of $A$ and consequently $N(A,2r)\leq P(A,r)$.
For the second inequality, there is nothing to prove if $N(A,r)=\infty$. So, again, let
$\{(x_i,r)\}_{i=1}^k$ be a finite $r$--packing of $A$ and assume that 
$$
\{B(y_j,r)\}_{j=1}^\ell
$$
is a finite cover of $A$ with centers in $A$. Then, for each $1\leq i\leq k$ there exists
$1\leq j(i)\leq \ell$ such that 
$$
x_i\in B(y_{j(i)},r).
$$
We will show that $k\leq \ell$. In
order to do this, it is enough to show that the function $i\mapsto j(i)$ is injective.
But for each $1\leq j\leq \ell$, the cardinality of the set 
$$
\{\{x_i\}_{i=1}^k\cap B(y_j,r)\}
$$
is at most $1$ (otherwise, $\{(x_i,r)\}_{i=1}^k$ would not be an $r$-packing), and so
$i\mapsto j(i)$ is injective, as required. Thus, $P(A,r)\leq N(A,r)$.
\epf

\

\fr These inequalities have the following immediate implications.

\

\bcor\label{sem3cor2.8}
If $X$ is a metric space and $A\subseteq X$, then
\[
\ov{\BD}(A)=\limsup_{r\to0}\frac{\log P(A,r)}{-\log r}
\]
and
\[
\un{\BD}(A)=\liminf_{r\to0}\frac{\log P(A,r)}{-\log r}.
\]
\ecor

\sp\fr As an immediate consequence of this corollary and of the second part of Definition~\ref{d120141231}, we get the following.

\sp\bcor\label{c220131231}
If $X$ is a metric space and $A\subseteq X$, then
$$
\HD(A)\le \un{\BD}(A) \le \PD(A)\le \ov{\BD}(A)
$$
\ecor
We end this section with the following.

\sp\bprop\label{p720131014}
Let $(X,\rho)$ be a metric space endowed with a finite Borel measure $\mu$
such that 
$$
\mu(B(x,r))\ge Cr^t
$$
for some constant $C>0$, all $x\in X$, and all radii $0\le r\le 1$.
Then
$$
\ov{\BD}(X)\le t.
$$
If, on the other hand, $\mu(X)>0$ and 
$$
\mu(B(x,r))\le Cr^t
$$
for some constant $C<+\infty$, all $x\in X$, and all radii $0\le r\le 1$,
then
$$
\un{\BD}(X)\ge\HD(X)\ge t.
$$
Finally if $\mu(X)>0$ and 
$$
C^{-1}r^t\le \mu(B(x,r))\le Cr^t
$$
for some constant $C\in [1,+\infty)$, all $x\in X$, and all radii $0\le r\le 1$,
then
$$
\BD(X)=\PD(X)=\HD(X)=t.
$$
\eprop

\bpf
We start with the first inequality. Let $\{B(x,r)\}_{i=1}^k$ be an $r$-packing 
of $X$. Then
$$
kr^t
\le C^{-1}\sum_{i=1}^k\mu(B(x_i,r))
\le C^{-1}.
$$
Hence $k\le C^{-1}r^{-t}$. Therefore $P(X,r)\le C^{-1}r^{-t}$. Consequently,
$$
\log P(X,r)\le -\log C-t\log r.
$$
In conjunction with the first formula of Corollary~\ref{sem3cor2.8}, this yields 
$$
\ov{\BD}(X)\le t.
$$
The second assertion of our proposition directly follows from the first inequality of Corollary~\ref{c220131231} and from item (2) of Frostman's Converse Theorem for Hausdorff Measures 
(Theorem~\ref{tncp12.1.}).

\sp\fr The last assertion of our proposition is now an immediate consequence of the two first assertions and Corollary~\ref{c220131231}.

\epf

\chapter{Invariant Measures, Finite and
Infinite}\label{invariant-measure}

In this chapter  we begin to investigate the class of all measurable dynamical systems that posses an invariant measure. This means we are in the field of ergodic theory. We do want to emphasize that most of this chapter pertains to all the systems regardless of whether the reference invariant measure is finite or infinite. With one notable exception though, namely that of Birkhoff's Ergodic Theorem proved in \cite{BET}. This theorem, formulated in the realm of probability spaces, is important for at least 
three--folded reasons. Firstly, establishing equality of time averages and space averages, it yields profound theoretical and, one could say, philosophical, consequences. Secondly, it has countless applications throughout various types of dynamical systems; our book is an evidence of this. And thirdly, it is an indispensable tool to develop the ergodic theory with infinite invariant measures. 
Otherwise, unless explicitly stated (rarely), we do not assume
the  reference measure to be finite. 

Invariant measures, finite and infinite abound but we postpone examples to further chapters. Indeed, Section~\ref{EoIaEM}, Examples of Invariant and Ergodic Measures, contains a fairly large collection of invariant and ergodic measures. Further examples, dealt in detail in this book, are provided in Section~\ref{CGDMS}, Conformal Graph Directed Markov Systems, particularly Theorem~\ref{t2.2.4}, and especially in Chapter~\ref{invariant-p.s.n.r.}, Conformal Invariant Measures for Compactly Non--Recurrent Regular Elliptic Functions, which is almost entirely devoted to investigate (ergodic) invariant measures, both finite and infinite, for dynamical systems generated by elliptic meromorphic functions.

\section[Quasi--Invariant Transformations]{Quasi--Invariant
  Transformations, Ergodicity and Conservativity}\label{QEC}

\fr In this section $(X, \mathfrak{F}, m)$ is a measure space and a map 
$T:X \lra X$ is a measurable with respect the $\sg$--algebra $\mathfrak{F}$. 

\bdfn\label{d120190610}
If $T:(X,\mathfrak{F}) \lra (X,\mathfrak{F})$ is a measurable map, then
a measure $m$ on $(X,\mathfrak{F})$ is called quasi--invariant.\index{(N)}{quasi--invariant measure} with respect to the map $T$ if and only if
$$ 
m \circ T^{-1} \prec m.
$$
Equivalently, $m(T^{-1}(A))=0$ whenever $A$ is measurable and
$m(A)=0$.\index{(N)}{quasi-invariant measure}  We would like to
emphasize that in this section we do not yet assume the measure $m$ to
be $T$--invariant. This concept will be defined in the next section.
\edfn

\bdfn\label{mu_d1_2014_11_08}
If $m$ is a quasi--invariant measure for a measurable map 
$T:X \lra X$, then we say that the map $T:X \lra X$ (or the measure $m$) is
ergodic\index{(N)}{ergodic map} \index{(N)}{ergodic measure}
if and only if for every measurable set $A\sbt X$ the following implication holds:
$$
T^{-1}(A)=A \  \  \Lra  \  \  [m(A)=0 \  \text{ or } \
m(X\sms A)=0]. 
$$
\edfn

It is easy to prove the following.

\bprop\label{mu_d1_2014_11_08B}
If $m$ is a quasi--invariant measure for a measurable map 
$T:X \lra X$, then the map $T:X\lra X$ (or the measure $m$) is
ergodic\index{(N)}{ergodic map} \index{(N)}{ergodic measure}
if and only if for every measurable set $A\sbt X$:
$$
m(T^{-1}(A)\div A)=0 \  \  \Lra  \  \  [m(A)=0 \  \text{ or } \
m(X\sms A)=0]. 
$$
\eprop

\fr The second important concept in the theory of
quasi--invariant measures\index{(N)}{quasi--invariant measure}
\index{(N)}{quasi--invariant map} is conservativity. We introduce it
now.

\bdfn\label{d1jn1} Let $T:(X,\mathfrak{F}) \lra (X,\mathfrak{F})$ be a measurable map. A measurable set $W \subset X$  is called  a
wandering set\index{(N)}{wandering set} if and only if the sets
$$
\{T^{-k}(W)\}_{n=0}^\infty
$$ 
are mutually disjoint. 
\edfn

\fr One way of constructing wandering sets is now described.
\blem\label{wandconst}
Let $T:(X,\mathfrak{F}) \lra (X,\mathfrak{F})$ be a measurable map. If  $A\in\mathfrak{F}$, the set 
$$
W_A:=A\backslash\bigcup_{n=1}^\infty T^{-n}(A)
$$
is wandering with respect to $T$.
\elem

\bpf
Suppose for a contradiction that $W_A$ is not wandering for $T$, that is,
$$
T^{-k}(W_A)\cap T^{-l}(W_A)\neq\emptyset
$$ 
for some $0\leq k<l$. This means that 
$$
T^{-k}\(W_A\cap T^{-(l-k)}(W_A)\)\neq\emptyset,
$$ 
and thus, on the one hand,
$$
W_A\cap T^{-(l-k)}(W_A)\neq\emptyset.
$$
On the other hand, by the very definition of $W_A$,
$$
W_A\cap T^{-(l-k)}(W_A)=\emptyset.
$$
This contradiction finishes the proof.
\epf

\bdfn\label{d1jn1-20190610}  A measurable map $T:(X,\mathfrak{F}) \lra (X,\mathfrak{F})$ is called conservative \index{(N)}{conservative map} if  and only if there is no wandering set $W\sbt X$  with
$m(W)>0$. 
\edfn

\fr For every set $B\sbt X$, let
$$
\begin{aligned} 
\hat B_\infty
:=&\Big\{x\in X: T^n(x)\in B\,\,\,\mbox{ for infinitely many}\,\,\,  n \geq
1\Big\}\\
=&\lt\{x\in X: \sum_{n=1}^\infty\1_B \circ T^n(x) =+\infty \rt\}
 \\
&=\bigcap_{k=1}^\infty\bigcup_{n=k}^\infty T^{-n}(B).
\end{aligned}
$$ \index{(S)}{$B_\infty$}
and 
$$
B_\infty:=B\cap \hat B_\infty.
$$
Obviously,
\beq\label{520190610}
T^{-1}(\hat B_\infty)=\hat B_\infty=T(\hat B_\infty) 
\eeq
and
\beq\label{620190610}
T^{-1}(X\sms \hat B_\infty)=X\sms \hat B_\infty=T(X\sms \hat B_\infty).
\eeq
Notice also that if $W$ is wandering then 
\beq\label{120190611}
W \cap\bigcup_{n=1}^\infty T^{-n}(W)
=\es 
=W \cap\bigcup_{n=1}^\infty T^n(W),
\eeq
In particular 
$$
W_\infty=\es.
$$
The celebrated Poincar\'e's Recurrence Theorem \index{(N)}{Poincar\'e's Recurrence Theorem} says that if the measure $m(X)$ is finite, then $m(B\sms B_\infty)=0$. We shall prove now its generalization in two respects, assuming only that the measure merely is quasi--invariant and need not be finite. The Poincar\'e Recurrence Theorem, immediately following from this generalization, will be stated as Theorem~\ref{Poincare}.

\sp\bthm\label{t2jn1} {\rm (Halmos' Recurrence
Theorem).}\index{(N)}{Halmos' Recurrence Theorem} Let $m$ a quasi--invariant measure for a measurable map $T: X\lra X$. Fix a measurable set $A\sbt X$. Then 
$$
m(B\sms B_\infty)=0
$$  for all measurable sets $B \sbt A$ if and only if
$$
m(A\cap W)=0
$$ 
for all measurable wandering sets $W\sbt X$. 
\ethm

\bpf Note that if $m(A)=0$, then $m(B\sms B_\infty)=0$
for all measurable sets $B \sbt A$ and $m(A\cap W)=0$ for all
measurable wandering sets $W\sbt X$. Thus, our  equivalence is
trivially satisfied, and we may assume without loss of generality
that $m(A)>0$.

Proving implication $(\Rightarrow)$, suppose for a contradiction  that
$m(A\cap W) >0$ for some measurable wandering set $W \sbt X $. Then
$(A\cap W)_\infty=\es$, and therefore, 
$$
m\(A\cap W\sms (A\cap W)_\infty\)>0.
$$

Let us now prove the implication $(\Leftarrow)$. Fix a measurable set $B\sbt A$ and for all $ n\geq 0$ let,
$$ 
B_n:= B \cap T^{-n}(B)\cap \bigcup_{k=n+1}^\infty T^{-k}(X \sms B).
$$
Suppose for a contrary that 
$$
m(B \sms B_\infty)>0.
$$
Then there exists  $n\geq 0$  such that $m(B_n)>0$. But, as $B_n\sbt B \sbt
A$, it follows from our hypotheses and Lemma~\ref{wandconst} that $m(W_{B_n})=0$.  This means that $m(B_n\sms
B_n^-)=0$, or equivalently 
$$
m\lt(B_n \sms  \bigcup_{k=1}^\infty T^{-k}(B_n)\rt)=0.
$$
So, (as $m(B_n)>0)$  there exists $x\in B_n \cap
\bigcup_{k=1}^\infty T^{-k}(B_n)$. Hence, 
$$
x\in B \cap T^{-k}(B_n)\sbt B \cap T^{-(n+ k)}(B),
$$
for some $ k\geq 1$.  Thus $x \notin B_n$. This contradiction  finishes  proof that $m(B\sms B_\infty)=0$. We are done. 
\endpf

\sp \fr Taking in the Theorem~\ref{t2jn1} $A=X$, we immediately get the
following  corollary.

\bcor\label{c1j58} 
A measurable transformation $T:(X,\mathfrak{F})\lra (X,\mathfrak{F})$ of a measure space $(X,\mathfrak{F})$, possessing a quasi--invariant measure $m$, is conservative with respect to $m$ if and only if
$$
m(B\sms B_\infty)=0
$$ 
for all $B \in \mathfrak{F}$. 
\ecor

\sp As an immediate consequence of this corollary and Borel--Cantelli
Lemma, we get the following.

\bcor\label{Cor5.4a} If a quasi--invariant transformation $T:X\lra
X$of a  measure space $(X, \mathfrak{F}, m)$ is conservative, then 
$$
\sum_{n=0}^\infty m(T^{-n}(A))=+\infty
$$ 
for any set $A\in\mathfrak{F}$ with  $m(A)>0$. 
\ecor

We now  shall  prove the following  characterization of ergodicity
and conservativity. 

\bthm\label{t1j59}  
A transformation $T: X\lra X$ of a measure space $(X, \mathfrak{F})$, possessing a quasi-invariant measure $m$, is conservative and ergodic with respect to $m$ if and only if for every set $A\in\F$ with $m(A)>0$,
$$ 
\sum_{n=1}^\infty \1_A \circ T^n= + \infty
$$
$m$--\textit{a.e.} on $X$. Equivalently: if
and only if 
$$
\mu(X\backslash \hat A_\infty)=0
$$ 
for every $A\in\F$ such that $\mu(A)>0$.

\ethm

\bpf Assume first that the map $T:X\lra X$ is conservative and
ergodic. Fix $A\in\F$ such that $\mu(A)>0$. It then follows from \eqref{520190610} that either $m(\hat A_\infty)=0$ or $m\(X\sms \hat A_\infty\)=0$. If $m\(X\sms \hat A_\infty\)=0$, we are done. So, suppose that
$$
m(\hat A_\infty)=0.
$$
By Corollary~\ref{c1j58} this implies that $m(A)=0$. This contradiction
 finishes the proof of our implication. 

\sp Now, let us prove the converse. We first show that the map $T:X\lra X$ is conservative. Keep $A\in \mathfrak{F}$ with $m(A)>0$, and assume that $m(X\backslash \hat A_\infty)=0$. Since $A_\infty = A\cap\hat A_\infty$, we  thus  get that $m(A\sms A_\infty)=0$. It therefore follows from Corollary~\ref{c1j58} that $T:X\to X$ is conservative. 

Now, let us turn to ergodicity. Keep $A\in\mathfrak{F}$ and assume $T^{-1}(A)=A$. Then $\hat A_\infty=A$ and $\widehat{(X\sms A)}_\infty=X\sms A$. If $m(A)>0$, then, as  by assumption $\mu(X\sms A_\infty)=0$, we get that $\mu(X\sms A)=0$. Likewise, if $m(X\sms A)>0$, we get that $\mu(A)=0$. We are done.
\endpf

\fr With a a little stronger hypotheses, we get yet another consequence of Theorem~\ref{t1j59}. 

\bdfn
A measurable transformation $T:(X,\F,m)\to (X,\F,m)$ is called non--singular if for every set $F\in \F$,
$$
m(T^{-1}(F))=0 \  \iff  \ m(F)=0.
$$
\edfn
Of course, every non--singular map is quasi--invariant. The, just announced, consequence of Theorem~\ref{t1j59} is this.

\bcor\label{c1_mu_2014_11_10} 
If $T:X\lra X$ is a non--singular transformation of a measure space $(X,\mathfrak{F},m)$, then the following are equivalent.

\sp\begin{itemize}  
\item[(a)] The map $T:X\lra X$ is conservative and ergodic.

\sp\item[(b)] For every set $A\in\F$ with $m(A)>0$
$$ 
\sum_{n=1}^\infty \1_A \circ T^n(x)= + \infty
$$
for $m$--\textit{a.e.} $x\in X$. 

\sp\item[(c)]
$$
\mu(X\backslash \hat A_\infty)=0
$$ 
for every $A\in\F$ such that $\mu(A)>0$.

\sp\item[(d)] 
$$
m\lt(X\sms \bu_{n=0}^\infty T^{-n}(A)\rt)=0
$$
for every set $A\in\F$ with $m(A)>0$. 
\end{itemize}
\ecor

\bpf
Because of Theorem~\ref{t1j59}, our only task is to show that (d)$\,\imp\,$(b). Indeed, for every $n\ge 1$ let 
$$
A_n^*:=X\sms \bu_{k=0}^\infty T^{-k}(T^{-n}(A))
$$
By non--singularity of $T$, $m(T^{-n}(A))>0$ for every $n\ge 0$. Hence, by (d), $m(A_n^*)=0$ for every $n\ge 0$. Therefore,
$$
m\lt(\bu_{n=0}^\infty A_n^*\rt)=0.
$$
But for a point $x\in X$ to be in the complement of this set means that $T^j(x)\in A$ for infinitely many $j$s. This in turn means that (b) holds for the point $x$.
\epf

Suppose now that $(X, \mathfrak{F}, m)$ is a probability space and that
$T:X\lra X$ is a measurable  quasi--invariant  mapping. Assume further, that $T(A)\in \F$ for every set $A\in\F$. We call $T:X\lra X$ weakly metrically  exact\index{(N)}{weakly metrically exact map} if and only if
\beq\label{1jb1-20190916}
\limsup_{n \to \infty} m(B\cap T^n(A))=m(B)
\eeq
for all measurable set $A, B\sbt X$ with $m(A)>0$. 

If $m$ is a probability space, then this just means that
\beq\label{1jb1}
\limsup_{n \to \infty} m(T^n(A))=1
\eeq
for all measurable set $A\sbt X$ with $m(A) >0$.

\sp We shall prove the following theorem.

\bthm\label{t320190918}
Each weakly metrically exact mapping $T:(X,\mathfrak{F})\lra (X, \mathfrak{F})$, possessing a quasi--invariant probability measure $m$, is conservative and ergodic with respect to $m$. 
\ethm

\bpf Proving ergodicity suppose that $m(A)>0$ and $T^{-1}(A)= A$. Then,
$$
T^n(A)\sbt A
$$ 
for all $n\ge 0$. Hence 
\[
m(X\sms A)
=\limsup_{n\to\infty}\mu\((X\sms A)\cap T^n(A)\)
\leq m(A\cap (X\sms A))=0.
\]
So $m(X\sms A)=0$, and $T$ is ergodic.

In order to prove conservativity of $T$ note that if $W\sbt X$ is a wandering set then by, \eqref{120190611}, we have that
$$
\bigcup_{n=0}^\infty T^n(W)\sbt X\sms W.
$$
Therefore, 
$$
m(W)=
\limsup_{n \to \infty} m(W\cap T^n(W))
\leq m(W\cap (X\sms W))=0.
$$
Hence, $m(W)=0$. The map $T:X\to X$ is thus conservative and we are done.
\endpf

\sp Now,  let $(X, \mathfrak{F}, m)$ be a $ \sigma$--finite measure space
and keep $T:X \lra X$ an $\mathfrak F$-measurable
transformation with respect to which the measure $m$ is quasi-invariant. Let $f\in L^1(m)$ be a non--negative function. Then the measure $(fm)\circ T^{-1}$ is absolutely continuous with repect to $m$, and we  put
\beq\label{1_mu_2014_11_15}
{\mathcal L}_{m}(f) := \frac{d((fm)\circ T^{-1})}{dm}.
\eeq\index{(S)}{${\mathcal L}_{m}(f)$} 
to be the Radon--Nikodym derivative of $fm\circ T^{-1}$ with respect to $m$. For any function $f\in L^1(m)$ we then define
$$
{\mathcal L}_{m}(f) := {\mathcal L}_{m}(f_+)-{\mathcal L}_{m}(f_-),
$$
where $f=f_+-f_-$ is the canonical decomposition of $f$ into its positive and negative parts $f_+=\max\{f,0\}$ and $f_-=-\min\{f,0\}$ respectively. Obviously
$$
{\mathcal L}_{m}: L^1(m) \to  L^1(m)
$$ 
is a linear operator and it is characterized by the property that
\beq\label{1j157}
 \int_{X}{\mathcal L}_{m}(f) g dm = \int_X f \cdot g\circ T dm
\eeq
for all $f \in L^1(m)$ and $g \in L^{\infty}(m)$. Taking $g=\1$, we
see that ${\mathcal L}_m: L^1(m) \to L^1(m)$ is bounded and
$\|{\mathcal L}_m\|_{L^1(m)}=1$. The bounded linear operator
\beq\lab{1_2017_11_21}
{\mathcal L}_m: L^1(m) \lra L^1(m)
\eeq
is called the transfer, or Perron--Frobenius, operator \index{(N)}{transfer operator} \index{(N)}{Perron--Frobenius operator} associated to the  quasi--invariant measure $m$. We  shall prove now one characterization  more of
ergodicity  and conservativity, this time in terms of the transfer
operator ${\mathcal L}_m$.

\sp\bthm\label{t1j159} 
Let $(X, \mathfrak{F}, \mu)$  be a $\sigma$--finite measure space and let $T: X \lra X$ be an $\mathfrak F$--measurable transformation with respect to which the measure $m$ is quasi--invariant. Then the transformation $T:X \lra X$ is ergodic and conservative with respect to $m$ if and only if 
$$
\sum_{n=0}^\infty {\mathcal L}_m^n(f)=+\infty
$$ 
$m$--a.e. for every non--negative (a.e.) function $f \in L^1(m)$ with $\int_Xf dm>0$. 
\ethm

\bpf ($\imp$) Proving by contradiction, suppose that there
exists a function $ f \in L^1(m)$ with the following properties
\begin{itemize}

\sp\item [(a)] $ f \geq 0$\,  $m$--{\it a.e.},

\sp\item [(b)] $\int_X f dm >0$,

\sp\item [(c)] $m\(\{ x \in X:\,\, \sum_{n=0}^\infty {\mathcal
L}_m^n(f)(x)< +\infty \}\)> 0.$
\end{itemize}
By virtue of (c)  there  exists a measurable set 
$$ 
B \sbt \lt\{x \in X:
\, \sum_{n=0}^\infty {\mathcal L}^n_m(f)(x)< \infty \rt\}
$$ 
such that
$m(B)\in (0, \infty)$ and
$$ 
\sup\lt\{\sum_{n=0}^\infty {\mathcal L}_m^n(f)(x):x \in B\rt\}
< +\infty.
$$
In particular
$$ 
\int_B \left(\sum_{n=0}^\infty {\mathcal L}_m^n(f)\right) dm <+\infty.
$$ 
Using also (\ref{1j157}) we therefore get that
\beq\label{1j159}
\begin{aligned}
\int_X f \left(\sum_{n=1}^\infty \1_B \circ T^n \right) dm
&=\sum_{n=1}^\infty\int_Xf\1_B\circ T^ndm
 =\sum_{n=1}^\infty\int_X\1_B{\mathcal L}_m^n(f)dm \\
&=\sum_{n=1}^\infty\int_B{\mathcal L}_m^n(f)dm \\
&=\int_B\left(\sum_{n=1}^\infty {\mathcal L}_m^n(f)\right) dm \\
&<+\infty.
\end{aligned}
\eeq
Since (b) implies that $f$ is positive on a set of positive measure
$m$, formulas (\ref{1j159}) and (a) yield that $ \sum_{n=1}^\infty
\1_B  \circ T^n$ is finite  on a set of    positive  measure $m$.
So, $T$  is not conservative and ergodic by virtue of Theorem~\ref{t1j159}. This contradiction finishes the proof of our implication.

\sp Proceeding in the opposite direction, assume that $T$ is not ergodic
and conservative. By Theorem~\ref{t1j59} this means that there exist two sets
$C, F \in \mathfrak F$ such that $m(C)>0$, $ m(F) >0$ and
$$\sum_{n=0}^\infty  \1_F\circ T^n(x) < \infty$$
for all $x \in C$. Therefore, there exists a measurable subset $D$ of $C$  such
that $m(D) \in (0, + \infty)$ and
$$ 
M:=\sup\left\{ \sum_{n=0}^\infty \1_F\circ T^n(x) : \, \,   x \in D\right\}<
 \infty.
 $$
Since $\1_D\ge 0$ and since $\1_D\in L^1(m)$, using
(\ref{1j157}), we get that
$$
\begin{aligned}
\int_F\sum_{n=0}^\infty{\mathcal L}_m^n(f) dm 
&= \int_X\1_F\sum_{n=0}^\infty {\mathcal L}^n_m(\1_D)dm 
=  \int_X \1_D \left(\sum_{n=0}^
 \infty \1_F \circ T^n\right) dm\\
 &=  \int_D \sum_{n=0}^\infty \1 \circ T^ndm\\ 
 &\leq
 M < + \infty.
\end{aligned}
$$
Hence 
$$ 
m\lt(\lt\{x \in F: \sum_{n=0}^\infty {\mathcal L}_m^n (f)(x)<+\infty\rt\}\rt) =m(F)>0
$$ and we are done by contrapositive. 
\endpf

\

\section[Invariant Measures: First Return Map (Inducing)] {Invariant Measures: First Return Map (Inducing); \\ Poincar\'e Recurrence Theorem}

In this section we define and extensively investigate a very powerful tool of ergodic theory, the one constituted by the first return time and map, also named inducing. This tool is extremely useful in both the ergodic theory of transformations preserving probability measures as well as those preserving infinite measures. In the former case, after inducing, the obtained map is usually somewhat worse in regard to a minor aspect but it is much better in regard of a major, more wanted, aspect. In the infinite case, inducing allows us to study systems preserving infinite measure by means of systems preserving probability measures. And this is a big advantage indeed, especially since for such maps we have the tool of Birkhoff's Ergodic Theorem.

Throughout this section $(X,\mathfrak{F}, \mu)$ is a measure space
and $T:(X,\mathfrak{F})\lra (X,\mathfrak{F})$  is a measurable conservative map preserving measure\index{(N)}{preserving measure map}\index{(N)}{invariant
measure} $\mu$, which means that 
$$
\mu \circ T^{-1}=\mu.
$$ 
This is the central equation of ergodic theory. We then also say that the measure $\mu$ is $T$--invariant. This concept will be explored extensively throughout the current chapter and entire book. Since obviously every measure preserving map is non--singular, as an immediate consequence of Corollary~\ref{c1_mu_2014_11_10}, we get the following.

\bthm\label{t2_mu_2014_11_10}
If $T:X\lra X$ is a measure preserving transformation of a
measure space $(X,\mathfrak{F},\mu)$, then the following are equivalent.

\sp\begin{itemize}  
\item[(a)] The map $T:X\to X$ is conservative and ergodic.

\sp\item[(b)] For every set $A\in\F$ with $\mu(A)>0$,
$$ 
\sum_{n=1}^\infty \1_A \circ T^n= + \infty
$$
$\mu$-\textit{a.e.} on $X$. 

\sp\item[(c)]
$$
\mu(X\backslash \hat A_\infty)=0
$$ 
for every $A\in\F$ such that $\mu(A)>0$.

\sp\item[(d)] 
$$
\mu\lt(X\sms \bu_{n=0}^\infty T^{-n}(A)\rt)=0
$$
for every set $A\in\F$ with $\mu(A)>0$. 
\end{itemize}
\ethm

\fr As an immediate consequence of this theorem, we get the following.

\bcor\label{ce_mu_2014_11_10}
If $T:X\lra X$ is a conservative and ergodic measure preserving transformation of a measure space $(X,\mathfrak{F},\mu)$, then the measure space $(X,\mathfrak{F},\mu)$ is $\sg$--finite if and only if there exists at least one set $F\in\F$ with $0<\mu(F)<+\infty$.
\ecor

\sp We know from Theorem~\ref{t320190918} that a sufficient condition for ergodicity and conservativity of a quasi--invariant measure is that it is weakly metrically exact. If the quasi--invariant measure is invariant, then we can say a little bit more. Namely, we say that if $(X,\mathfrak{F},\mu)$ is a measure space, then a measurable map $T:X \lra X$ preserving measure $\mu$ such that $T(A)\in \F$ for every set $A\in\F$, is called metrically exact\index{(N)}{metrically exact map} if and only if
\beq\label{1jb1-20190918A}
\lim_{n \to \infty} \mu(B\cap T^n(A))=\mu(B)
\eeq
for all measurable set $A, B\sbt X$ with $\mu(A)>0$. 

If $\mu$ is a probability space, then this just means that
\beq\label{1jb1-20190918B}
\lim_{n \to \infty} \mu(T^n(A))=1
\eeq
for all measurable set $A\sbt X$ with $\mu(A)>0$. Noting that for invariant measures $\mu$,
$$
\mu(T(A))\ge \mu(A)
$$
for every measurable set $A\sbt X$ such that $\mu(A)$ is also measurable, and applying Theorem~\ref{t320190918}, we get the following.

\bprop\label{p120190918}
If $(X,\mathfrak{F},\mu)$ is a probability space and $T:X \lra X$ is a measurable map preserving measure $\mu$, then $T$ is metrically exact
if an only if it is weakly metrically exact. If this holds, the map $T$ is ergodic and conservative.
\eprop

\sp\fr We also have the following immediate result.

\bprop\label{p220190918} 
Let $(X,\mathfrak{F},m)$  be a probability space and let $T: X\lra X$ be an $\mathfrak F$--measurable transformation with respect to which the measure $m$ is quasi--invariant and weakly metrically exact. Then any $T$--invariant probability measure on $(X,\mathfrak{F})$ absolutely continuous with with respect to $m$ is metrically exact. 
\eprop

\sp In the context of finite measures, as an immediate consequence of the very definition of conservativity and Corollary~\ref{c1j58}, we get the following celebrated theorem.

\index{(N)}{Poincar\'e's Recurrence Theorem} 
\bthm[Poincar\'e Recurrence Theorem]\label{Poincare}
If $(X,\mathfrak{F},\mu)$ is a finite measure space, then every measurable map $T:X \lra X$ preserving measure $\mu$ is conservative. This means that
$$
\mu(F\sms F_\infty)=0
$$
for every set $F\in\F$.
\ethm

This theorem has profound philosophical and physical consequences, particularly because of its apparent contradiction with the Second Law of Thermodynamics, growth of entropy. It is also a frequently used tool in everyday dealing with ergodic theory.

\sp We now define and investigate the main concept of this section, the inducing scheme. So, fix $F\in\mathfrak{F}$ with 
$$
0<\mu(F)< +\infty.
$$ 
Consider the function 
$$
\tau_F: F_\infty\lra \N=\{1, 2, 3, \ldots\}
$$
given by the formula
$$
\tau_F(x):=\min\{n \geq 1:T^n(x)\in F\},
$$\index{(S)}{$\tau_F$}
and the map 
$$
T_F:F_\infty\to F_\infty,
$$
\index{(S)}{$T_F$} given by the formula,
$$
T_F(x):=T^{\tau_F(x)}(x).
$$
Since, Corollary~\ref{c1j58}, 
$$
\mu(F\sms F_\infty)=0,
$$
we may, somewhat informally, say that the map 
$$
T_F:F\lra F
$$
and all its iterates, are defined $m$ \textit{a.e.} on $F$. The number $\tau_F(x)\geq 1$  is
called  the first return time to the  set $F$\index{(N)}{first
return time}, and  the map $T_F:F\to F$ is called the first return
map\index{(N)}{first return map} or the induced map.

\

\fr Let $\phi: X \lra \mathbb R$ be a measurable function. For any $n
\geq 1$  we define  $n$'th Birkhoff's sum of the
function\index{(N)}{Birkhoff's sum of the function} $\varphi$ by
\beq\label{srednia}
S_n\phi=\phi+\phi \circ T +\ldots + \phi\circ
T^{n-1}.\index{(S)}{$S_n\phi$}
\eeq
Given a  function $g:X\lra \mathbb{R}$, let $g_F: F \to
\mathbb{R}$ be defined by the formula,
\beq\label{suma}
g_F(x)= \sum_{j=0}^{\tau_F(x)-1}g \circ T^j(x).\index{(S)}{$g_F(x)$}
\eeq
Let $\mu_F$ be the conditional measure   on $F$, \textit{i.e.}
$$
\mu_F(A)=\frac{\mu(A)}{\mu(F)}
$$\index{(S)}{$\mu_F$}
for every measurable  set $A\sbt F$. We shall prove the following.

\bthm\label{t1j64} 
Suppose $T: X\lra X$ is a measurable transformation of a measure space $(X, \mathfrak{F},m)$. 

\sp\begin{itemize}
\item[(a)] Let $\mu$ be a $T$--invariant measure on $X$. If $F\in\F$ satisfies 

\sp\centerline{$0<\mu(F)<+\infty$ \ and \ $\mu(F \sms F_\infty)=0$, } 

\sp then measure $\mu_F$ is $T_F$--invariant on $F$.

\sp\item[(b)] Conversely, if $\nu$ is a probability $T_F$--invariant measure on $\(F,\mathfrak{F}|_F\)$, then there exists a $T$--invariant measure $\mu$ on $(X, \mathfrak{F})$ such that 
$$
\mu_F= \nu.
$$
In fact, the formula
\beq\label{1j64b}
\begin{aligned}
\mu(B):=\sum_{k=0}^\infty \nu\left(F\cap T^{-k}(B) \sms
\bigcup_{j=1}^k T^{-j}(F)\right)
=\sum_{k=0}^\infty \nu\(\{ x\in F \cap T^{-k}(B):\tau_F(x)
>k\}\). 
\end{aligned}
\eeq
defines such a $T$--invariant measure on $X$.

\sp\item[(c)] Consequently,
$$
\mu\lt(X\sms \bigcup_{k=1}^\infty T^{-k}(F)\rt)=0.
$$
\item[(d)]
In particular, the measure $\mu$ is $\sigma$--finite.
\end{itemize}
\ethm

\bpf We first prrove (a). For every $B \in \mathfrak{F}|_F$,
we have
\beq\label{1j64}
\begin{aligned}
\mu(T^{-1}_F(B))
=& \sum_{n=1}^\infty \mu(\tau_F^{-1}(n)\cap T^{-1}_F(B))
=\sum_{n=1}^\infty \mu(\tau_F^{-1}(n)\cap T^{-n}(B))\\
=& \sum_{n=1}^\infty \mu(F\cap T^{-n}(B)\sms \bigcup_{k=1}^{n-1}T^{-j}(F))\\
= & \sum_{n=1}^\infty \mu(F\cap T^{-1}(B_{n-1})),
\end{aligned}
\eeq
where 
$$
B_0:=B 
\  \  \  {\rm and} \  \  \ 
B_n
:=T^{-n}(B)\sms \bigcup_{k=0}^{n-1} T^{-k}(F)
$$
for all $n\geq 1$. Observe that 
$$
\mu(B_n)\leq\mu(T^{-n}(B))=\mu(B)\leq\mu(A)<+\infty
$$ 
for every $n\geq0$. Since 
$$
T^{-1}(B_{n-1})=(F\cap
T^{-1}(B_{n-1}))\cup B_n
$$ 
and  
$$
(F\cap T^{-1}(B_{n-1}))\cap B_n=\es
$$
we have that
\beq\label{5.3a}
\mu(F\cap T^{-1}(B_{n-1}))=\mu(T^{-1}(B_{n-1})) - \mu(B_n)\leq
\mu(B_{n-1})-\mu(B_n).
\eeq
Therefore, by (\ref{1j64}),
\beq\label{1aj64}
\mu(T^{-1}_F(B))= \lim_{n\to \infty}(\mu(B)-\mu(B_n))\leq \mu(B).
\eeq
Hence, also
\beq\label{2j64}
\mu(T^{-1}_F(B))= \mu(F)-\mu(T^{-1}_F(F \sms B))\geq \mu(F)-\mu(F
\sms B)=\mu(B).
\eeq
These last two formulas complete the proof of item (a).

\sp Proving item (b), we shall show first that the measure $\mu$  given
by (\ref{1j64b}) is $T$--invariant. Indeed,
$$\begin{aligned}
\mu \circ T^{-1}(B))=& \sum_{k=0}^{\infty} \nu(\{x\in F\cap
T^{-(k+1)}(B):\tau_F(x)> k\})\\
 = & \sum_{k=0}^{\infty} \nu(\{x\in F\cap
T^{-(k+1)}(B):\tau_F(x)> k+1\})+\\
& \  \  \  \  \  \  \  \  \  \  \  \  +\sum_{k=0}^{\infty} \nu(\{x\in F\cap T^{-(k+1)}(B):\tau_F(x)= k+1\})\\
= & \mu(B)-\nu(B\cap F) + \sum_{k=1}^{\infty} \nu(\{x\in F \cap
T^{-k}(B):\tau_F(x)= k \}) \\
= &  \mu(B)-\nu(B\cap F) + \sum_{k=1}^{\infty} \nu(\{x\in F:\tau_F(x)= k \mbox{ and}\,\,\ x\in T^{-1}_F(B)\}) \\
 = &  \mu(B)-\nu(B\cap F) + \sum_{k=1}^{\infty} \nu(\{x\in F:\tau_F(x)= k \mbox{ and}\,\,\ x\in T^{-1}_F(F\cap B)\}) \\
= &  \mu(B)-\nu(B\cap F)+ \nu(T^{-1}_F(F\cap B))\\
=& \mu(B).
\end{aligned}
$$
If $B \sbt F$, then \eqref{1j64b} reduces to
$\mu(B)=\nu(F\cap B)=\nu(B)$. Hence, 
$$
\mu_F(B)=\frac{\mu(B)}{\mu(F)}=\frac{\nu(B)}{\nu(F)}=\nu(B).
$$
Item (b) is thus proved. 

\sp Proving (c), formula \eqref{1j64b} gives
\[ 
\mu\Bigl(X\Bigl\backslash\Bigr.\bigcup_{k=1}^\infty T^{-k}(F)\Bigr)
=\sum_{n=0}^\infty\nu\Bigl(\Bigl\{x\in F\Bigl\backslash\Bigr.\bigcup_{k=1}^\infty T^{-(n+k)}(F):\tau_F(x)>n\Bigr\}\Bigr) 
=\sum_{n=0}^\infty\nu(\emptyset)
=0,
\]
and item (c) is proved.

Item (d) is now an immediate consequence of Theorem~\ref{t2_mu_2014_11_10} and Corollary~\ref{ce_mu_2014_11_10}. 
\endpf

\sp\brem\label{r1j64a.1}
 It follows from (\ref{1j64}) and (\ref{2j64}) that
 $\lim_{n \to \infty }\mu(B_n)=0$. In particular, taking $B=F$, we
 get that 
$$
\lim_{n \to \infty }\mu\lt(T^{-n}(F)\sms\bigcup_{k=0}^{n-1}T^{-k}(F)\rt)=0.
 $$
\erem

\fr In the  second part of Theorem~\ref{t1j64} we  have defined the
invariant  measure $\mu$  by the formula (\ref{1j64b}). We   now
shall show  that this  choice of the invariant measure $\mu$ is
actually unique.

\bprop\label{p1j64c}
 If $T$ is a measure preserving transformation of a
 measure space $(X,\mathfrak{F}, \mu)$, $F \in \mathfrak{F}$ with
 $ 0< \mu(F) < +\infty$ and 
$$
\mu\lt(X\sms \bigcup_{n=1}^{\infty}T^{-n}(F)\rt)=0,
$$
then for every set $B \in\mathfrak{F}$ we have that,
$$
\mu(B)
=\sum_{k=0}^\infty \mu\left( F \cap T^{-k}(B)\sms \bigcup_{j=1}^k
 T^{-j}(F)\right)
=\sum_{k=0}^\infty\mu\(\{x\in F\cap T^{-k}(B):\tau_F(x)> k\}\).
$$
\eprop \bpf Let $B \in\mathfrak{F}$. Assume first that $\mu(B)<+\infty$. For every $k \geq 0$ put,
$$
B_k:=T^{-k}(B)\sms \bigcup_{j=0}^k T^{-j}(F).
$$
Observe that $\mu(B_k)\leq\mu(T^{-j}(B))=\mu(B)<+\infty$ for every $j\geq0$.
As $\mu$ is $T$--invariant and $T^{-1}(B_{k-1})\spt B_j$ for all $k\ge 1$, we get for every $n\ge 0$ that
$$
\begin{aligned} \mu(B \sms F)-m(B_n)
=&  \sum_{j=1}^{n} \mu (B_{j-1}) -\mu (B_j)\\
=& \sum_{j=1}^{n} \mu (T^{-1}(B_{j-1}))- \mu(B_j)
= \sum_{j=1}^{n} \mu (T^{-1}(B_{j-1})\sms B_j) \\
=&\sum_{j=1}^n \mu\left(( T^{-j}(B)\sms \bigcup_{i=1}^{j}
T^{-i}(F))\sms (T^{-j}(B) \sms \bigcup_{i=0}^{j}
T^{-i}(F))\right) \\
=& \sum_{j=1}^n\mu\left(F\cap T^{-j}(B)\sms
\bigcup_{i=1}^{j} T^{-i}(F)\right).
\end{aligned}
$$
Then
\begin{eqnarray*} 
\mu(B)-\mu(B_n)
&=&\mu(B\cap F)+\mu(B\backslash F)-\mu(B_n) \\ 
&=&\sum_{j=0}^n\mu\Bigl(F\cap T^{-j}(B)\Bigl\backslash\Bigr.\bigcup_{k=1}^{j} T^{-k}(F)\Bigr).
\end{eqnarray*}
Letting $n\to \infty$, we see that  in order to complete the proof,
it is enough to show that
$$ 
\lim_{n \to \infty}\mu(B_n)=0.
$$
Since $\mu\(X \sms \bigcup_{k=0}^{\infty} T^{-k}(F)\)=0$, we have that
\beq\label{1j64d}
\mu\lt(X\sms \bigcup_{k=0}^{\infty} F^{(k)}\rt)=0,
\eeq
where 
$$
F^{(k)}:=T^{-k}(F)\sms \bigcup_{j=0}^{k-1} T^{-j}(F).
$$
The usefulness of the $F^{(k)}$s lies in their mutual disjointness.
Indeed, relation~(\ref{1j64d}) implies that 
$\mu\(B\backslash\bigcup_{k=0}^\infty F^{(k)}\)=0$. 
Then
\beq\label{bba}
\mu(B)
=\mu\Bigl(B\cap\bigcup_{k=0}^\infty F^{(k)}\Bigr) 
=\mu\Bigl(\bigcup_{k=0}^\infty B\cap F^{(k)}\Bigr) 
=\sum_{k=0}^\infty\mu(B\cap F^{(k)}).
\eeq
Fix $\varepsilon >0$. Since $\mu(B)<\infty$, there exists $k_\varepsilon\in\N$ so large that
\beq\label{elleps}
\sum_{k=k_\varepsilon+1}^\infty\mu\(B\cap F^{(k)}\)<\frac{\varepsilon}{2}.
\eeq
Relation~(\ref{1j64d}) also ensures that 
$\mu\(B_n\backslash\bigcup_{k=0}^\infty F^{(k)}\)=0$ for every $n\geq0$. So,
like for $B$,
\[
\mu(B_n)
=\sum_{k=0}^\infty\mu\(B_n\cap F^{(k)}\).
\]
But
\begin{eqnarray*}
B_n\cap F^{(k)}
&=&\left[T^{-n}(B)\Bigl\backslash\Bigr.\bigcup_{i=0}^n T^{-i}(F)\right]
   \cap
   \left[T^{-k}(F)\Bigl\backslash\Bigr.\bigcup_{j=0}^{k-1}T^{-j}(F)\right] \\
&=&\bigl[T^{-n}(B)\cap T^{-k}(A)\bigr] 
   \Bigl\backslash\Bigr.
	 \bigcup_{i=0}^{\max\{n,k-1\}}T^{-i}(F) \\
&=&\left\{
   \begin{array}{lcl}
	 \emptyset & \mbox{ if } & k\leq n \\
	 \displaystyle
	 T^{-n}\bigl(B\cap T^{-(k-n)}(F)\bigr)
   \Bigl\backslash\Bigr.
	 \bigcup_{i=0}^{k-1}T^{-i}(F)
	           & \mbox{ if } & k>n \\   
	 \end{array} 
   \right. \\
&=&\left\{
   \begin{array}{lcl}
	 \emptyset & \mbox{ if } & k\leq n \\
	 \displaystyle
	 \left[
	 T^{-n}\bigl(B\cap T^{-(k-n)}(F)\bigr)
   \Bigl\backslash\Bigr.
	 \bigcup_{i=n}^{k-1}T^{-i}(F)
   \right]
	 \Bigl\backslash\Bigr.
	 \bigcup_{i=0}^{n-1}T^{-i}(F)
	           & \mbox{ if } & k>n \\   
	 \end{array} 
   \right. \\
&=&\left\{
   \begin{array}{lcl}
	 \emptyset & \mbox{ if } & k\leq n \\
	 \displaystyle
	 T^{-n}\bigl(B\cap F^{(k-n)}\bigr)
	 \Bigl\backslash\Bigr.
	 \bigcup_{i=0}^{n-1}T^{-i}(F)
	           & \mbox{ if } & k>n. \\   
	 \end{array} 
   \right.	
\end{eqnarray*}
Consequently, 
\[
\mu(B_n)
=\sum_{k>n}\mu(B_n\cap F^{(k)})
=\sum_{k=n+1}^\infty\mu\Bigl(T^{-n}\bigl(B\cap F^{(k-n)}\bigr)
	 \Bigl\backslash\Bigr.
	 \bigcup_{i=0}^{n-1}T^{-i}(F)\Bigr)
=\sum_{k=1}^\infty\mu\bigl(B^{(k)}_n\bigr),
\]
where for each $n,k\geq 0$,
$$
B^{(k)}_n:= T^{-n}(B\cap F^{(k)})\sms \bigcup_{j=0}^{n-1} T^{-j}(F) 
\sbt T^{-(n+k)}(F)\sms \bigcup_{j=0}^{n+k+1} T^{-j}(F)
=B^{(n+k)}.
$$ 
Thus, by Remark~\ref{r1j64a.1}, we have  for every  $k\geq 0$ that
\beq\label{2j64d}
\lim_{n \to \infty}\mu(B^{(k)}_n)=0.
\eeq
Since for every $n \geq 0$, we have that
$$
\mu(B^{(k)}_n)\leq
\mu(T^{-n}(B^{(k)}))=\mu(B^{(k)}),
$$
we therefore get from \eqref{elleps} for every $n \geq 0$, that
$$ 
\sum_{k=k_\varepsilon+1}^\infty  \mu(B^{(k)}_n)
 < \frac{\varepsilon}{2}.$$
By virtue of (\ref{2j64d})  there exists $n_\varepsilon \geq 1$ so
large that for all $1\leq k \leq k_\varepsilon$ and  all $n \geq
n_\varepsilon$,
$$
\mu(B^{(k)}_n)\leq \frac{\varepsilon}{2k_\varepsilon}.
$$
So, far  all $n\geq n_\varepsilon$, we  have
$$ 
\begin{aligned}
\mu(B_n)
&=\mu\lt(\bigcup_{k=1}^\infty B^{(k)}_n\rt)
=\sum_{k=1}^\infty\mu\(B^{(k)}_n\)
=\sum_{k=1}^{k_\varepsilon} \mu( B^{(k)}_n)
+\sum_{k=k_\varepsilon+1}^\infty\mu( B^{(k)}_n)\\
&\leq \sum_{k=0}^{k_\varepsilon}
\frac{\varepsilon}{2k_\varepsilon}+ \frac{\varepsilon}{2} \\
&= \frac{\varepsilon}{2}+\frac{\varepsilon}{2} 
=\varepsilon.
\end{aligned}
$$
The proof is thus complete for all $B\in\F$ with $\mu(B)<+\infty$. 

Now, let $B\in\F$ be any set. Since 
$$
\mu(F^{(k)})\leq\mu(T^{-k}(F))=\mu(F)<+\infty
$$ 
for all $k\geq 0$, the sets 
$$
\(B\cap F^{(k)}\)_{k=0}^\infty
$$
are all of finite measure. Thus, the first part of this proof shows that 
$$
\mu\(B\cap F^{(k)}\)
=\sum_{n=0}^\infty \mu\left(F\cap T^{-n}\(B\cap F^{(k)}\)\sms \bigcup_{j=1}^n T^{-j}(F)\right).
$$   
By~(\ref{bba}), the mutual disjointness of the sets $F^{(k)}$s and Theorem~\ref{t1j64} (c) along with its formula \eqref{1j64b} (to have the last equality), we then conclude that
\[
\begin{aligned}
\mu(B)
&=\sum_{k=0}^\infty\mu\(B\cap F^{(k)}\)
=\sum_{k=0}^\infty\sum_{n=0}^\infty \mu\left(F\cap T^{-n}\(B\cap F^{(k)}\)\sms \bigcup_{j=1}^n T^{-j}(F)\right) \\
&=\sum_{n=0}^\infty\sum_{k=0}^\infty\mu\left(F\cap T^{-n}\(B\cap F^{(k)}\)\sms \bigcup_{j=1}^n T^{-j}(F)\right) \\
&=\sum_{n=0}^\infty\mu\left(F\cap T^{-n}(B)\sms \bigcup_{j=1}^n T^{-j}(F)\right).
\end{aligned}
\]
The proof is complete.
\endpf

\bprop\label{p1j65} 
Let $T:(X,\mathfrak{F}) \lra (X,\mathfrak{F})$ be a measurable map preserving a measure $\mu$ on $(X,\mathfrak{F})$. Let $\varphi_F: F \lra \mathbb{R}$ be a measurable function. Assume that $F\in\F$ is such that $0<\mu(F)<+\infty$ and 
$$
\mu\lt(X\sms \bigcup_{k=1}^\infty T^{-k}(F)\rt)=0.
$$
Then,
\begin{itemize}
\item[(a)] $\varphi_F\in L^1(\mu_F)$ whenever $\varphi\in L^1(\mu)$.

\,

\item[(b)] If $\varphi\geq0$ or $\varphi\in L^1(\mu)$, then 
$$ 
\int_F \varphi_F d\mu_F= \frac{1}{\mu(F)}\int_X \varphi d\mu.
$$
\end{itemize}
If, in addition, $T$ is conservative and ergodic, then the above two statements apply to all sets $F\in\F$ such that $0<\mu(F)<+\infty$.
\eprop

\bpf Suppose first that $\varphi=\1_B$ for some $B\in\F$ such that $0< \mu(B)< +\infty$. In view of Proposition~\ref{p1j64c}, we  have
$$\begin{aligned}
\int_X \1_Bd\mu=\mu(B)=& \sum_{k=0}^{\infty} \mu( \{x \in
T^{-k}(B)\cap F: \,\, \tau_F(x)>k \})\\
=& \sum_{n=1}^{\infty} \sum_{j=0}^{n-1} \mu( \{x \in F\cap
T^{-j}(B): \,\, \tau_F(x)=n
\})\\
=&  \sum_{n=1}^{\infty} \sum_{j=0}^{n-1} \int_{\tau_F^{-1}(n)}
\1_{T^{-j}}(B)d \mu\\
 =& \sum_{n=1}^{\infty}
\sum_{j=0}^{n-1}\int_{\tau_F^{-1}(n)} \1_B
\circ T^j d\mu \\
=& \sum_{n=1}^{\infty}\int_{\tau_F^{-1}(n)} S_n\1_B  d\mu\\
=& \sum_{n=1}^{\infty}\int_{\tau_F^{-1}(n)} \varphi_F  d\mu\\
=& \int_F \varphi_F  d\mu= \mu(F) \int_F \varphi_F  d\mu_F,
\end{aligned}
$$
and we are done in this case. 

If $\varphi:X\lra \mathbb{R}$ is a simple measurable
function, \textit{i.e.} $\varphi= \sum_{i=1}^n
a_i\varphi^{(i)}$, where all $a_i \in \mathbb{R}$, $1 \leq i \leq
n$, and all $\varphi^{(i)}$, $1 \leq i \leq n$, are characteristic
functions of some measurable sets with positive and finite measures,
then
$$\begin{aligned}
\int_X  \varphi d\mu =&\sum_{i=1}^n a_i\int_X  \varphi^{(i)} d\mu
= \mu(F) \sum_{i=1}^n a_i \int_F  \varphi^{(i)}_F d\mu_F\\
=& \mu(F) \int_F \sum_{i=1}^n  a_i \varphi^{(i)}_F d\mu_F\\
=& \mu(F) \int_F \(\sum_{i=1}^n \varphi^{(i)}\)_F d\mu_F\\
=& \mu(F) \int_F \varphi_F d\mu_F
\end{aligned}
$$
We thus are done in this case as well.

The next case is to consider an
arbitrary  non-negative measurable  function  $\varphi: X\to [0,
+\infty)$. Then  $\varphi$ is a point-wise monotone increasing limit
of non-negative step functions, say
$(\varphi^{(n)})^\infty_{n=1}$. Then also, the sequence $(\varphi_F^{(n)})^\infty_{n=1}$ converges point-wise in a monotone increasing way to $\varphi_F$. Hence, applying  Lebesgue Monotone Convergence Theorem twice, we
then get that
$$
\int_X  \varphi d\mu 
=\lim_{n \to \infty}\int_X  \varphi^{(n)} d\mu\\
=\lim_{n \to \infty}  \mu(F) \int_F  \varphi^{(n)}_F d\mu_F
=\mu(F)\int_F  \varphi_F d\mu_F.
$$
Since $|\varphi_F|\leq  |\varphi|_F$, we have in particular
shown that if $\varphi:X\to \mathbb{R}$ is $\mu$--integrable, then
$\varphi_F: F\to \mathbb{R}$ is $\mu_F$--integrable.
 Moreover, writing  $\varphi=\varphi^+  - \varphi^-$, where
$\varphi^+=\max\{\varphi,  0\}$ and $\varphi^-=\max\{-\varphi, 0\}$,
we have that both  functions  $\varphi^-$ and   $\varphi^-$ are
$\mu$-integrable and
$$
\begin{aligned}
\int_X  \varphi d\mu
= &  \int_X  \varphi^+ d\mu- \int_X  \varphi^-d\mu
=\mu(F) \int_F \varphi^+ d\mu_F - \mu(F) \int_F\varphi^-d\mu_F\\
=& \mu(F) \int_F  (\varphi_F^+ - \varphi^-_F) d\mu_F\\
=& \mu(F) \int_F  \varphi_F d\mu_F.   
\end{aligned}
$$
The proof is complete. \endpf

\sp\fr Observe that if $\varphi \equiv 1$, then $\varphi_F \equiv \tau_F$,
and  therefore,  and as an immediate consequence of this
proposition, we get the following celebrated result.

\bthm\label{t1j69} {\rm (Kac's Lemma).}\index{(N)}{Kac's  Lemma} Let
$T:(X, \mathfrak{F}, \mu) \to (X, \mathfrak{F}, \mu)$ be a measurable map
preserving measure $\mu$. If $F\in \mathfrak{F}$ is such that $0< \mu(F) < + \infty$ and 
$$
\mu\lt(X\sms \bigcup_{k=1}^\infty T^{-k}(F)\rt)=0,
$$ 
then
$$ 
\int_F \tau_F d \mu_F=\frac{\mu(X)}{\mu(F)}.
$$
In particular, 

\sp\begin{itemize}
\item[(a)] $\mu(X)<+\infty \ \iff \ \int_F \tau_F d \mu<+\infty$.

\sp\item[(b)] If $ \mu$ is a probability measure, then
$$ 
\int_F \tau_F d \mu_F=\frac{1}{\mu(F)}.
$$
\end{itemize}
If, in addition, $T$ is conservative and ergodic, then the above two statements apply to all sets $F\in\F$ such that $0<\mu(F)<+\infty$.
\ethm

\fr Item (a) of this theorem tells us that the measure $\mu$ is finite if and only if the first return time is integrable. This is a very powerful tool to check whether an invariant measure is finite or infiite. It will be heavily explored in Chapter~\ref{invariant-p.s.n.r.} in the context of the dynamics of elliptic functions.

\bprop\label{p1j70} Suppose that $T:(X, \mathfrak{F}, \mu)  \lra (X,
\mathfrak{F}, \mu)$ is a measure preserving transformation. 

\,

\begin{itemize}
\item[(a)] If $T$ is ergodic and  conservative with respect to $\mu$, then for every set $F \in
\mathfrak{F}$ with $0< m(F)< +\infty$, we have that $T_F: F\lra F$ is
ergodic with respect to $\mu_F$. 

\,

\item[(b)] Conversely, if for some such set $F$, the map $T_F: F \lra F$  is ergodic with respect to $\mu_F$ and 
$$
\mu\lt(X\sms
\bigcup_{n=1}^\infty T^{-n}(F)\rt)=0,
$$
then the map $T:X\lra X$ is ergodic with respect to $\mu$.
\end{itemize}
\eprop

\bpf Suppose that $T:X\to X$ is ergodic and conservative. If
$A\sbt F$  is $T_F$-invariant and $\mu(A)>0$, then for $\mu$ almost all $x \in F
\sms A$, we  have by Theorem~\ref{t1j59}, that
$$ 
0=\sum_{n=1}^\infty \1_A \circ T^n_F(x)
= \sum_{n=1}^\infty \1_A \circ T^n(x)=+\infty.
$$
This contradiction shows  that the map $T_F: F\to
F$ is ergodic. 

\sp In order to prove the converse, suppose that
$T^{-1}(B)=B$  and $m(B)>0$. Since $\mu(X\sms \bigcup_{n=1}^\infty
T^{-n}(F))=0$, there exists $k \geq 1$ such that $\mu(B\cap
T^{-k}(F))>0$. But then
$$
\mu(B\cap F)
=\mu(T^{-k}(B\cap F))
=\mu(T^{-k}(B)
\cap T^{-k}(F))
=\mu(B \cap T^{-k}(F))>0.
$$ 
Since $m_F$ is a probability measure and $T_F:F\to F$ is  ergodic, this implies that
$$
1=m_F(\bigcup_{n=0}^\infty T^{-n}(B \cap F))\leq
m_F(\bigcup_{n=1}^\infty T^{-n}(B\cap F)).$$ This means that $m(F
\sms \bigcup_{n=0}^\infty T^{-n}(B\cap F))=0$. As $m(X \sms
\bigcup_{n=0}^\infty T^{-n}( F))=0$, this implies that 
$$
m(X \sms
\bigcup_{n=0}^\infty T^{-n}(B\cap F))=0.
$$
Consequently, 
$$
m(X \sms B)
=m(X\sms  \bigcup_{n=1}^\infty T^{-n}(B))\leq  m(X\sms
\bigcup_{n=1}^\infty T^{-n}(B\cap F))=0.
$$
We are done.
\endpf

\

\section{Ergodic Theorems; Birkhoff's, von Neumann's and Hopf's}

In this section we deal with Birkhoff's, von Neumann's, and Hopf's Ergodic
Theorems. These are central pillars of abstract ergodic theory. Birkhoff's Ergodic Theorem, proved first time in \cite{BET}, establishing equalities of time and space averages, has profound theoretical and philosophical consequences. It abounds in applications. The present book is an evidence of this. The original proof provided by G. D. Birkhoff in \cite{BET} was very long and complicated. Since then some simplifications have been made. We provide a short simple proof of Birkhoff's Ergodic Theorem, taken from \cite{KH}. This theorem concerns measure--preserving dynamical
systems acting on probability spaces. As its first, fairly
straightforward, consequence we prove $L^p$--von Neumann's Ergodic
Theorem. Secondly, as a more involved  consequence of Birkhoff's
Ergodic Theorem, by utilizing
the induced procedure described  in the previous section, we  prove
Hopf's  Ergodic Theorem, which holds for general measure-preserving
transformations, whose $\sg$--finite invariant measure can be
infinite. We would like  to remark that Hopf's Ergodic Theorem is
not proved in  Aaronson's book~\cite{Aa}. 

\sp The formulation and the proof of Birkhoff's Ergodic Theorem we provide below utilizes the concept, widely used in the probability theory, of an expected value with respect to a sub $\sg$--algebra which we introduced and studied in Section~\ref{CEaMT}.  

Let $(X,\F,\mu)$ be a probability space and $T:X \lra X$ be a measurable map preserving the probability measure $\mu$. Let 
$$
\mathcal I\index{(S)}{$\mathcal I$}
:=\{A \in \mathfrak{F}:\mu(T^{-1}(A)\triangle A)=0\}. 
$$  
It is easy to check that $\cI$ is a $\sigma$-algebra and
we call it the $\mathcal I$ the $\sigma$-algebra of $T$-invariant (mod 0) sets. Let us record the following obvious.

\bthm\label{t1_mu_2014_11_11}
A measurable map $T:X \lra X$ preserving a probability measure $\mu$ is ergodic if and only if $\mathcal I$, the $\sigma$--algebra of $T$--invariant (mod 0) subsets of $X$, is trivial, i.e. it consists of sets of measure one and zero only.
\ethm

Recall that for any integer $n \geq 1$, the $n$'th Birkhoff sum of the function $\phi$ was defined as
$$
S_n\phi=\phi+\phi \circ T +\ldots + \phi\circ T^{n-1}.
$$

\sp\bthm[Birkhoff's Ergodic Theorem]\label{Birkhoff}\index{(N)}{Birkhoff's  Ergodic Theorem} 
If $T:X \lra X$ is a measure preserving endomorphism of a probability space $(X,\mathfrak{F}, \mu)$  and if $\phi :X \lra \mathbb R$ is an integrable
function, then
$$
\lim_{n \to \infty}\frac{1}{n}S_n\phi (x)
=E_\mu(\phi|\mathcal I)\quad \mbox{for}\quad  \mu-a.e.\, \, x \in X.
$$ 
If, in addition, $T$ is
ergodic, then
$$
\lim_{n \to \infty}\frac{1}{n}S_n\phi (x)
=\int_X\phi\, d\mu\quad
\mbox{for}\quad  \mu-a.e.\, \, x \in X.
$$
\ethm
\fr{\sl Proof.}
Let $f \in L^1(\mu)$. For every integer $n\ge 1$ let
$$
F_n:=\max\Big\{
\sum_{i=0}^{k-1}f \circ T^i:  \,\, 1 \leq k \leq n \Big\}.
$$
Of course the sequence $(F_n)_{n=1}^\infty$ is monotone increasing. Then, for every $x \in X$, we have 
$$
F_{n+1}(x)-F_n(T(x))=f(x)-\min\{0,F_n(T(x))\}\geq f(x).
$$
The sequence $(F_{n+1}(x)-F_n\circ T)_{n=1}^\infty$ is monotone decreasing, since $(F_n)_{n=1}^\infty$ is monotone increasing. Define
$$
A:=\lt\{x\in X: \,\,  \sup_{n\ge 1}\Big\{\sum_{i=0}^n f(T^i(x))\Big\}=+\infty\rt\}.
$$
Note that  $A \in \mathcal I$. If $x \in A$, then
$F_{n+1}-F_n(T(x))$ monotonously  decreases to  $f(x)$ as $n \to
\infty$. The Lebesgue Monotone Convergence Theorem implies then that
\beq\label{PU2.2.4}
0 \leq \int_A (F_{n+1}-F_n)d \mu= \int_A( F_{n+1}-F_n\circ T)d \mu
\longrightarrow \int_Af\,d\mu.
\eeq
Notice that $\frac{1}{n}\sum_{k=0}^{n-1}\phi \circ T^k \leq F_n/n$;
so outside $A$, we have
\beq\label{PU2.2.5}
\limsup_{n \to \infty} \frac{1}{n}\sum_{k=0}^{n-1}f\circ T^k \leq 0.
\eeq
Therefore, if the conditional expectation  value $f_{\mathcal I}$
of $f$ is negative $a.e.$,  that is if 
$$
\int_C f d \mu=\int_C f_{\mathcal I} d \mu <0
$$ 
for all $C \in \mathcal I$ with $\mu(C)>0$,
then, as $A \in \mathcal I$, (\ref{PU2.2.4}) implies that
$\mu(A)=0$. Hence (\ref{PU2.2.5}) holds $a.e.$  Now, fix $\varepsilon>0$ and let
$$
f:=\phi - \phi_{\mathcal I}-\varepsilon,
$$
Then obviously
$$
f_{\mathcal I}= -\varepsilon <0.
$$
The equality $\phi_{\mathcal I} \circ T=\phi_{\mathcal I}$ entails
$$ 
\frac{1}{n}\sum_{k=0}^{n-1} f \circ T^k
=\left( \frac{1}{n}\sum_{k=0}^{n-1}\phi \circ T^k\right) -\phi_{\mathcal I}-\varepsilon.
$$ 
So (\ref{PU2.2.5}) yields
$$ 
\limsup_{n \to \infty}\frac{1}{n}\sum_{k=0}^{n-1}\phi \circ
T^k\leq \phi_{\mathcal I}+ \varepsilon\quad a.e.
$$ 
Call the sets of points where this inequality holds by $X_+(\varepsilon)$. Replacing $ \phi$ by $-\phi$ gives
$$ 
\liminf_{n \to \infty} \frac{1}{n}\sum_{k=0}^{n-1}\phi \circ T^k \geq
\phi_{\mathcal I}- \varepsilon\quad a.e.
$$ 
Call the sets of points where this inequality holds by $X_-(\varepsilon)$.
Then denoting 
$$
X_*:=\bi_{n=1}^\infty \(X_+(1/n)\cap X_-(1/n)\),
$$
we have $\mu(X_*)=1$ and 
$$
\lim_{n\to
\infty}\frac{1}{n}\sum_{k=0}^{n-1}\phi \circ T^k=\phi_{\mathcal I}
$$
on $X_*$. The proof of the first part of our theorem is thus complete. The second part, i.e. the one concerning ergodic maps, directly follows from the first part along with Theorem~\ref{t1_mu_2014_11_11} and formula \eqref{1_mu_2014_11_11}. 
\endpf

\sp As the first consequence of Birkhoff's Ergodic
Theorem we prove von Neumann's Ergodic Theorem which asserts that the
converges in Birkhoff's Ergodic Theorem is not only almost everywhere
but also in $L^p(\mu)$ whenever the input function $\phi$ is in $L^p(\mu)$. It was first proved by John von Neumann in \cite{vN} and had many generalizations and extensions since then. 

\sp\bthm\label{von Neumann}{\rm ($L^p$--von Neumann's  Ergodic
Theorem).}\index{(N)}{von Neumann's Ergodic Theorem} If $T:X \lra X$ is
a measurable endomorphism of a probability space $(X,
\mathfrak{F},\mu)$ preserving measure $\mu$ and if $\phi :X \lra \mathbb R$ is a measurable function that belongs to $L^p(\mu)$
with some $p\ge 1$, then $E(\phi|\mathcal I)\in L^p(\mu)$ and
$$
\lim_{n \to \infty}\frac{1}{n}S_n\phi=E(\phi|\mathcal I),
$$ 
where the convergence and equality are understood in the Banach space
$L^p(\mu)$. If, in addition, $T$ is ergodic, then
$$
\lim_{n \to \infty}\frac{1}{n}S_n\phi (x)=\int \phi\, d\mu
$$
and, as above, both convergence and equality are taken in $L^p(\mu)$.
\ethm

\bpf As an auxiliary step consider first an essentially bounded function
$\psi\in L^p(\mu)$. This means that
$$
\|\psi\|_\infty:=\ess\sup\{|\psi(x)|:x\in X\}<+\infty.
$$
Then also $||n^{-1}S_n\psi\le||\psi||_\infty<+\infty$ for all $n\ge 1$. Hence
\beq\label{2_mu_2-14_11_11}
\lt\|\frac{1}{n}S_n\psi-\psi^*\rt\|_\infty\le||\psi||_\infty<+\infty.
\eeq
But by Birkhoff's Ergodic Theorem
(Theorem~\ref{Birkhoff}), the sequence $\(n^{-1}S_n\psi\)_1^\infty$
converges almost everywhere to $\psi^*:=E(\psi|\mathcal I)$. Hence
$\psi^*\in L^\infty(\mu)\sbt L^p(\mu)$ and 
$$
\lim_{n\to\infty}\lt|\frac{1}{n}S_n\psi(x)-\psi^*(x)\rt|=0
$$
$\mu$--almost everywhere. Therefore, looking up at \eqref{2_mu_2-14_11_11},  Lebesgue's Dominated Convergence Theorem yields,
$$
\lim_{n\to\infty}\lt\|\frac{1}{n}S_n\psi-\psi^*\rt\|_p=0,
$$
where by $||\cdot||_p$, we denoted the $L^p$ norm in the Banach space
$L^p(\mu)$. Consequently, we got the following.

\sp\fr{\bf Claim 1:} $\(n^{-1}S_n\psi\)_1^\infty$ is a Cauchy sequnce in $L^p(\mu)$.

\sp Passing to the general case, consider an arbitrary function
$\phi\in L^p(\mu)$. We aim to prove 

\sp\fr{\bf Claim 2:} $\(n^{-1}S_n\phi\)_1^\infty$
is a Cauchy sequnce in $L^p(\mu)$. 

\sp\fr To do this fix $\vep>0$
arbitrary. There then exists a function $\phi_\e\in L^\infty(\mu)$
such that 
\beq\label{620121208}
\|\phi-\phi_\vep\|_p<\vep/4.
\eeq
In consequence
\beq\label{120121210}
0\le \lt\|\frac{1}{n}S_n\phi-\frac{1}{n}S_n\phi_\e\rt\|_p\le \vep/4
\eeq
for all $n\ge 1$.
By virtue of Claim~1 there exists $N\ge 1$ such that if $k,l\ge N$,
then
$$
\lt\|\frac{1}{l}S_l\phi_\e-\frac{1}{k}S_k\phi_\e\rt\|_p<\vep/2.
$$
From this and \eqref{120121210} we get for such $k$ and $l$ that
$$
\aligned
\lt\|\frac{1}{l}S_l\phi-\frac{1}{k}S_k\phi\rt\|_p
&\le \lt\|\frac{1}{l}S_l\phi-\frac{1}{l}S_l\phi_\vep\rt\|_p
    +\lt\|\frac{1}{l}S_l\phi_\vep-\frac{1}{k}S_k\phi_\vep\rt\|_p
    +\lt\|\frac{1}{k}S_k\phi_\vep-\frac{1}{k}S_k\phi\rt\|_p \\
&\le \frac{\vep}{4}+\frac{\vep}{2}+\frac{\vep}{4} \\
&=\vep.
\endaligned
$$
Claim~2 is proved.

\sp\fr Let $\hat\phi$ be the limit of the sequence in
$$
\(n^{-1}S_n\phi\)_1^\infty
$$ 
in $L^p(\mu)$; the latter is a Banach space, so complete. Since convergence in the $L^p$ norm entails convergence in measure and since any sequence
convergent in measure contains a subsequence converging almost
everywhere, invoking Birkhoff's Ergodic Theorem
(Theorem\ref{Birkhoff}), we conclude that $\hat\phi=E(\phi|\mathcal
I)$. The last assertion of our theorem is now also an immediate
consequence of Theorem\ref{Birkhoff}.
\endpf

\sp Of particular significance is the case of $p=2$. Then
$L^p(\mu)$ is a Hilbert space with the inner product
$$
(\phi,\psi)=\int_X\phi\bar\psi\,d\mu.
$$
For every $\phi\in L^2(\mu)$ define
\beq\label{220121210}
U_T(\phi):=\phi\circ T.
\eeq
Then
$$
\int_X|U_T(\phi)|^2\,d\mu
=\int_X|\phi|^2\circ T\,d\mu
=\int_X|\phi|\,d\mu
<+\infty.
$$
So, $U_T(\phi)\in L^2(\mu)$ and \eqref{120121210} thus defines a bounded,
with norm $\le 1$, linear oprator from $L^2(\mu)$ into itself. It is
called the Koopman operator associated to the measure preserving map
$T:X\to X$. In fact, the same calculation gives the following.

\sp

\bthm\label{unitary}
$U_T(\phi):L^2(\mu)\to L^2(\mu)$ is a unitary operator, meaning that
$$
(U_Tg,U_Th)=(g,h)
$$
for all $g,h\in L^2(\mu)$.
\ethm

\

\sp\fr Denote by $U_T^*:L^2(\mu)\to L^2(\mu)$ the
operator conjugate to $U_T$. Let
$$
L_{inv}^2(\mu):=\big\{\phi\in L^2(\mu):U_T(\phi)=\phi\big\} 
\  \text{ and } \
L_{*inv}^2(\mu):=\big\{\phi\in L^2(\mu):U_T^*(\phi)=\phi\big\}
$$
be the respective spaces of all $U_T$-invariant and  $U_T^*$-invariant
elements of $L^2(\mu)$. Obviously both $L_{inv}^2(\mu)$ and
$L_{*inv}^2(\mu)$ are closed vector subspace of $L^2(\mu)$. We shall 
prove the following. 

\

\blem\label{l120121210}
We have
$$
L_{inv}^2(\mu)=L_{*inv}^2(\mu).
$$
\elem

\fr{\sl Proof.}
Let $g\in L_{inv}^2(\mu)$. Then $U_Tg=g$ and since, by
Theorem~\ref{unitary}, $U_T^*U_T=\Id$, we therefore get that
$U_T^*g=U_T^*U_Tg=g$, meaning that $g\in L_{*inv}^2(\mu)$. The inclusion 
$$
L_{inv}^2(\mu)\sbt L_{*inv}^2(\mu)
$$
is proved. Conversly, if $g\in L_{*inv}^2(\mu)$, then
\beq\label{220121210B}
\aligned
\|U_Tg-g\|_2^2
&=(U_Tg-g,U_Tg-g) 
=\|U_Tg\|_2^2-(g,U_Tg)-(U_Tg,g)+\|g\|_2^2 \\
&=\|g\|_2^2-(U_T^*g,g)-(g,U_T^*g)+\|g\|_2^2 \\
&=2\|g\|_2^2-2(g,g) 
=2\|g||_2^2-2\|g\|_2^2 \\
&=0.
\endaligned
\eeq
Hence $U_Tg=g$ meaning that $g\in L_{inv}^2(\mu)$ and the inclusion 
$$
L_{*inv}^2(\mu)\sbt L_{inv}^2(\mu)
$$
is also proved. We are done.
\endpf

\sp We will need
to know what is the orthogonal complement of $L_{inv}^2(\mu)$ in
$L^2(\mu)$ is. In fact, let 
$G_T$ be the closed vector subspace of $L^2(\mu)$ generated by all
coboundaries, i.e. by the vectors of the form
$$
u-u\circ T, \ \ u\in L^2(\mu)
$$
Of course all coboundaries form a vector space, so $G_T$ is just the
closure of all coboundaries.

\

\fr The orthogonal complement of $L_{inv}^2(\mu)$ is described by the following. 

\

\blem\label{l220121210} We have that
$$
L^2(\mu)=L_{inv}^2(\mu)\du G_T,
$$
i.e. the Hilbert space $L^2(\mu)$ can be represented as the orthogoal
sum of its closed vector subspaces $L_{inv}^2(\mu)$ and $G_T$.
\elem

\fr{\sl Proof.}
Fix a coboundary $\psi=u-u\circ T$. Then for every $\phi\in
L_{inv}^2(\mu)$ we have,
\beq\label{320121210}
\aligned
(\psi,\phi)
&=(u-U_Tu,\phi)
=(u,\phi)-(U_Tu,\phi) \\
&=(u,\phi)-(U_Tu,\phi) \\
&=(u,\phi)-(U_Tu,U_T\phi)
=(u,\phi)-(u,\phi) \\
&=0.
\endaligned
\eeq
This means that $\phi\perp\psi$, whence
\beq\label{420121210}
L_{inv}^2(\mu)\sbt G_U^{\perp}.
\eeq
Now fix $\phi\in G_U^{\perp}$. Then for every $\psi\in\in L^2(\mu)$ we
have $(\phi,\psi-U_T\psi)=0$, or equivalently,
$(\phi,\psi)=(\phi,U_T\psi)=(U_T^*\phi,\psi)$. This means that
$(U_T^*\phi-\phi,\psi)=0$. Hence, $U_T^*\phi-\phi=0$, or equivalently,
$\phi\in L_{inv}^2(\mu)$. Therefore,
$$
G_U^{\perp}\sbt L_{inv}^2(\mu).
$$
Along with \eqref{420121210} this finishes the proof.
\endpf

\sp\fr In subsequent chapters we will need the following characterization
of coboundaries.

\sp\blem\label{l320121210}
If $T:X\lra X$ is a measurable map preserving a probability measure
$\mu$ and if $g\in L^2(\mu)$, then the following two statements are
equivalent.
\begin{itemize}
\item[(a)] The function $g$ is a coboundary, i.e. $g=u-u\circ T$ for
some $u\in L^2(\mu)$.

\,

\item[(b)] The sequence $(S_ng)_1^\infty$ is bounded in the Hilbert
  space $L^2(\mu)$.
\end{itemize}
\elem

\bpf (a)$\imp$(b). For every $n\ge 1$ we have $S_ng=u-u\circ
T^n$. Therefore,
$$
\|S_ng\|_2
=\|u-u\circ T^n\|_2
\le \|u\|_2+\|u\circ T^n\|_2
=2\|u\|_2
$$
and thus the implication (a)$\imp$(b) is established.

\sp\fr (b)$\imp$(a). By our hypothesis there exists $M\ge 1$ such that
$\|S_ng\|_2\le M$ for all $n\ge 1$. Consequently,
$$
\lt\|\frac1n\sum_{j=1}^nS_jg\rt\|_2\le M
$$
for all $n\ge 1$. Since $L^2(\mu)$, as a Hilbert space, is reflexive,
its closed ball $\ov B(0,M)$ is weakly compact in $L^2(\mu)$. There
thus exist $u\in \ov B(0,M)$ and $(n_k)_1^\infty$, an increasing
sequence of positive integers, such that the sequence
$\(n_k^{-1}\sum_{j=1}^{n_k}S_jg\)_1^\infty$ converges weakly to
$u$. But then
\beq\label{1j289}
\lim_{k\to\infty}\frac1{n_k}\sum_{j=1}^{n_k}S_j(g\circ T)
=u\circ T \  \  \text{ weakly in } \  \ L^2(\mu)
\eeq
and
$$
\frac1{n_k}\sum_{j=1}^{n_k}S_j(g\circ T)
=\frac1{n_k}\lt(\sum_{j=1}^{n_k}S_{j+1}g-n_kg\rt)
=\frac1{n_k}\sum_{j=1}^{n_k}S_jg+\frac1{n_k}\(S_{n_k+1}g-g\)-g.
$$
Taking the weak limits of both sides of this equation and invoking
\eqref{1j289} we get $u\circ T=u-g$, or equivalently, $g=u-u\circ
T$. The proof of the implication (b)$\imp$(a) and the whole lemma is complete.
\endpf

\sp\bcor\label{c1fpn25}
If $T:X\lra X$ is a measurable map preserving a probability measure
$\mu$ and if $g\in L^2(\mu)$, then the following two statements are
equivalent.
\begin{itemize}
\item [(a)] The function $g$ is a  coboundary, i.e.  $g=u - u\circ T$ for some $ u \in L^2(\mu)$.

\,

\item [(b)]  There exists $l \geq 1$ such that $S_lg$ is  a  coboundary with respect to the dynamical system $(T^l, \mu)$, i.e. $S_lg=u-u\circ T^l$ for some $ u \in L^2(\mu)$.

\,

\item [(c)]  For every $l \geq 1$, $S_l$ is  a  coboundary with respect to the dynamical system $(T^l, \mu)$, i.e. $S_lg=u-u\circ T^l$ for some $ u \in L^2(\mu)$.
\end{itemize}
\ecor

\fr {\sl  Proof.} (a)$\Rightarrow$(c). We have 
$$ 
S_lg = \sum_{j=0}^{l-1} g \circ T^j= \sum_{j=0}^{l-1}(u \circ T^j- u \circ T^{j+1}) = u - u \circ T^l
$$ 
and the implication
(a)$\Rightarrow$(c) is established. 

The implication (c)$\Rightarrow$(b) is obvious.  

We shall prove   (b)$\Rightarrow$(a).  \nl Indeed, given an integer $n \geq 0$  write uniquely
$ n= lk_n+r_n, \, k=k_n \geq 0, \, \,  0 \leq r_n \leq l-1$.Then
\beq\label{1fpn26}
S_ng=S^{(l)}_{k_n}(S_lg) + \sum_{j-n- r_n}^{n-1} g \circ f^j= S^{(l)}_{k}(S_lg) +  S_r (g \circ f^{n-r_n}),
\eeq
where $S^{(l)}_n (h) $ is the  $n$-th  Birkhoff's sum of the function $h : X \to \mathbb R$ with respect to  to the dynamical  system $T^l: X \to X$.
It  follows  from Lemma~\ref{l320121210} that  the norms $|S^{(l)}_{k_n}(S_lg)\|_2$ are  uniformly bounded from above, say  by $M$.  It then  follows from
(\ref{1fpn26}),  that
$$ \|S_ng \|_2 \leq M + \| S_r(g\circ f^{n-r_n})\|_2\leq M + r \|g\|_2 \leq M +(l-1)\|g_2\|.$$
So, we  conclude from Lemma~\ref{l320121210} again that $g$ is a  coboundary in $L^2(\mu)$. The  proof is complete. \qed

\sp Coming back to Theorem~\ref{Birkhoff}, as we have already said, it implies that the {\it time average} exists 
for $\mu$-almost every $x \in X$. If additionally  $T$ is ergodic,
Theorem~\ref{Birkhoff} implies  that the {\it time average is equal to
the space average}.

\sp

\fr As a consequence of Theorem~\ref{Birkhoff} we get the following statement about visiting a measurable set of positive measure.

\sp

\bthm\label{t1_mu_2014_11-18}
Suppose that $(X, \mathfrak{F}, \mu)$ is  a
probability  space and $T:X \lra X $ is a measurable ergodic map preserving
measure $\mu$. Fix $F\in \mathfrak{F}$ with $\mu(F)>0$. Then for $\mu$-a.e. $x\in X$,
$$
\lim_{n\to\infty}\frac1n\#\big\{0\le j\le n-1:T^j(x)\in F\big\}=\mu(F).
$$
In particular, the set
$$
\big\{n\ge 0:T^n(x)\in F\big\}
$$
is infinite.
\ethm

\bpf
The first assertion is an immediate consequence of Theorem~\ref{Birkhoff} applied to the function $\vp:=\1_F$. The second assertion is an immediate consequence of the first one.
\epf

\sp We shall now prove two further little, technical, slightly surprising, but useful consequences of  Theorem~\ref{Birkhoff}

\sp\bprop\label{p1j71} Suppose that $(X, \mathfrak{F}, \mu)$ is  a
probability  space and $T:X \lra X $ is a measurable ergodic map preserving
measure $\mu$. Fix  $F \in \mathfrak{F}$ with $\mu(F)
>0$. For every $x\in X$ let $(k_n)_{n=1}^\infty$ be the sequence of
consecutive  visits of $x$ to $F$ under the action $T$. Then for $\mu$-{\it a.e.} $x
\in  X$
$$ \lim_{n \to \infty}\frac{k_{n+1}}{k_n}=1.$$
 \eprop
\bpf  Note that $S_{k_n} \1_F(x)=n$. If therefore follows
from Theorem~\ref{Birkhoff} and ergodicity of $T$ that
$$\begin{aligned}
\lim_{n \to \infty}\frac{k_{n+1}}{k_n}
=& \lim_{n \to
\infty}\left(\frac{n}{k_n}\frac{k_{n+1}}{n+1}\right)
= \lim_{n \to \infty}\frac{1}{k_n}S_{k_n} \1_F(x)
\frac{1}{\frac{1}{k_{n+1}}S_{k_{n+1}} \1_F(x)} \\
=& \frac{\lim_{n \to \infty} \left(\frac{1}{k_n}S_{k_n}
\1_F(x)\right)}{\lim_{n \to \infty}\lt(\frac{1}{k_{n+1}}S_{k_{n+1}}\1_F(x)\rt)}\\
= &  \frac{\mu(F)}{\mu(F)}=1.
\end{aligned}$$
We are done. \endpf

\sp\bprop\label{p1j73} 
Suppose that $(X, \mathfrak{F}, \mu)$ is  a
probability  space and $T:X \lra X $ is  a measurable  map preserving
measure $\mu$. If $f \in L^1(\mu)$, then for $\mu$-{\it a.e.} $ x\in
X$,
$$ \lim_{n\to \infty} \frac{1}{n}f(T^n(x))=0.$$
\eprop

\bpf It  follows from Theorem~\ref{Birkhoff} that for
$\mu$-{\it a.e.} $x\in X$,
$$\begin{aligned}
\lim_{n\to \infty} \frac{f(T^n(x))}{n}= & \lim_{n\to \infty} \frac{f(T^n(x))}{n+1}
= \lim_{n\to \infty} \frac{S_{n+1}f(x)- S_{n}f(x))}{n+1}\\
=& \lim_{n\to \infty} \frac{S_{n+1}f(x)}{n+1}- \lim_{n\to \infty} \frac{S_n f(x)}{n+1} \\
=& \lim_{n\to \infty} \frac{S_{n+1}f(x)}{n+1}- \lim_{n\to \infty} \frac{S_n f(x)}{n}\\
&=0. 
 \end{aligned}
$$
We are done. \endpf

\sp As an application of  Theorem~\ref{Birkhoff} to the ergodic theory
of measurable maps preserving  an infinite measure, we shall prove
the following, remarkable and useful in applications, theorem.

\sp\bthm\label{t2j73} {\rm (Hopf's Ergodic Theorem).}\index{(N)}{Hopf's
Ergodic Theorem} Suppose that $(X, \mathfrak{F}, \mu)$ is a $\sg$--finite measure space and that $T:X\lra X$ is a measurable ergodic and
conservative map preserving measure $\mu$.

Consider two measurable functions  $f,g \in L^1(\mu)$,
where $g\ge 0$ $\mu$--{\it a.e.} is such that
$\mu(g^{-1}((0,+\infty)))>0$ so that in particular $\int g d\mu>0$. Then 
$$
\lim_{n \to \infty} \frac{S_n f(x)}{S_n g(x)}=\frac{\int f
d\mu}{\int g d\mu}.
$$ 
\ethm

\bpf Since the measure  $\mu$ is  $\sg$--finite, there are
countably many mutually  disjoint sets $\{X_j\}_{j=1}^\infty$ such
that  for all $j \geq 1$, $0< \mu (X_j) < +\infty$ and $ \mu(X \sms
\bigcup_{j=1}^\infty X_j)=0$. Let 
$$
T_j=T_{X_j},
$$
where, we recall, $T_{X_j}:X_j\to X_j$ is the first return map of $T$ from $X_j$ to $X_j$. For  all
$\phi:\, X\to \mathbb{R}$, let $\phi_j:=\phi_{X_j}:X_j\to
\mathbb{R}$. Given $x\in X$, let
$$
S_n^{(j)}\phi_j(x):=\sum_{i=0}^{n-1}\phi_j(T^i_j(x)).
$$ 
If $x\in X_j$  and $n \geq 1$, let $j_n\geq 1 $  be the largest integer $k
\geq 0$ such that 
$$
\sum_{i=0}^{k-1}\tau_{X_j}(T_j^i(x))\leq n.
$$
Then for all $x\in X_j$ we have that
$$S_n \phi(x)= S^{(j)}_{T^{j_n}_i}\phi_j(x)+S_{\Delta  n}  \phi_j(x),$$
where  $$\Delta_n:= n - \sum_{i=0}^{k-1}\tau_{X_j}(T_j^i(x))\geq
0.$$ Then
\beq\label{1j75}
\frac{S_n \phi(x)}{j_n}= \frac{S^{(j)}_{T^{j_n}_i}\phi(x)}{j_n}+
\frac{S_{\Delta n}\phi(x)}{j_n}.
\eeq
Now
$$ \left| \frac{S_{\Delta n}\phi(x)}{j_n}\right|\leq \frac{1}{j_n}S_{\Delta
n}|\phi|(x)\leq  \frac{1}{j_n}|\phi|(T^{j_n}_j(x)).$$ It therefore
follows from Proposition~\ref{p1j73} and  Proposition~\ref{p1j65}
that if $\phi \in L^1(\mu)$, then for  $\mu$-{\it a.e.}, $ x\in
X_j$,
$$
\lim_{n \to \infty} \left| \frac{S_{\Delta n}\phi(x)}{j_n}\right|=0.
$$
Therefore, applying Theorem~\ref{Birkhoff}, Proposition~\ref{p1j65},
Proposition~\ref{p1j70}, and (\ref{1j75}), we get  for $\mu$-{\it
a.e.} $ x\in X_j$, and all $n \geq 1$, that
$$
\frac{S_nf(x)}{S_ng(x)}
=\frac{\frac{S_nf(x)}{j_n}}{ \frac{S_n g(x)}{j_n}}
=\frac{\frac{S^{(j)}_{j_n}f_j(x)}{j_n}+ \frac{S_{\Delta
 n}f(x)}{j_n}}{\frac{S^{(j)}_{j_n}g_j(x)}{j_n}+ \frac{S_{\Delta n}g(x)}{j_n}}\longrightarrow \frac{  \int_{X_j} f_j d\mu_{X_j}}
  {\int_{X_j} g_j d\mu_{X_j}} 
= \frac{  \int_{X} fd\mu}{\int_{X} g d\mu}.
$$
We are done. \endpf 

\sp As we shall seee in Section~\ref{CountET} this theorem being powerful and
interesting in itself, somewhat surprisingly, rules out any hope for a more 
direct version  of Birkhoff's Ergodic Theorem in the case  of
infinite measures. The world of infinite invariant measures is indeed
very different from the one of finite measures.

\sp Now we want to derive one additional  consequence  of
Theorem~\ref{Birkhoff}. We shall prove the following.

\sp\bthm\label{t1j74.1} Suppose that $\mathfrak{F}$ is a $\sg$--algebra
on a set $X$, $T:X\lra X$ is a measurable map, and $\mu_1$  and
$\mu_2$ are  two $\sg$--finite  measures  invariant under $T$. 

If $T:X\lra X$ is ergodic and conservative with respect to both measures
$\mu_1$ and $\mu_2$, then either $\mu_1$ and $\mu_2$ are mutually
singular or else they coincide up to a  positive multiplicative
constant (if both are probabilistic, then they are equal). 
\ethm

\bpf Suppose first that $\mu_1$ and $\mu_2$ are both
probability measures. If $\mu_1\neq \mu_2$, then there exists a set
$ F \in  \mathfrak{F}$ such that $\mu_1(F) \neq \mu_2(F)$. Let, for
$i=1,2$,
$$
X_i=\lt\{x\in  X:  \frac{1}{n} S_n\1_F(x)\lra  \int_X \1 _F
 d\mu_i=\mu_i(F)\rt\}.
$$
Since $\mu_1(F)\neq \mu_2(F)$ we  have that
  $$ X_1\cap X_2= \es.$$
But, by Theorem~\ref{Birkhoff}, $\mu_1(X_1)=1$ and $\mu_2(X_2)=1$.
So $\mu_1(X_2)=0$ and $\mu_2(X\sms X_2)=0$. Thus $\mu_1$ and $\mu_2$
are mutually  singular, and we are done in this case. Now, consider
the general case. Since both measures  $\mu_1$ and $\mu_2$ are  both
$\sg$--finite, there are  countably many measurable disjoint sets
$(Y_n)_{n=0}^\infty$  such that
$$ X=\bigcup_{n=1}^{\infty} Y_n$$
and $$ \mu_1(Y_n),\,\,\, \mu_2(Y_n)< +\infty$$ for all  $n \geq 1$.
Assume that   without loss of generality that neither measure
$\mu_1$ nor $\mu_2$ vanish.  Suppose first that  for some $n \geq
1$, $ \mu_1(Y_n)>0$, and $\mu_1$, $\mu_2$ coincide on $Y_n$ up to a
positive multiplicative constant. We may assume without  loss of
generality that $n=1$ and that  
$$
{\mu_1}_{|Y_1}={\mu_2}_{|Y_1}.
$$
It
then immediately  follows  from Proposition~\ref{p1j64c} that
$\mu_1=\mu_2$, and we are done in this case. 

Now, assume that
$\mu_1$  and $\mu_2$  do not coincide on $X$ up to a positive
multiplicative constant. If $\mu_1(Y_n)=0$, set  $Z_n=Y_n$, and if
$\mu_2(Y_n)=0$, set $Z_n=\es$. 

Now  consider the case  when
$$
\mu_1(Y_n), \mu_2(Y_n) \neq 0.
$$
From what we  have already proved, it  follows  that $ {\mu_1}_{|Y_1}$ and  $ {\mu_2}_{|Y_1}$ do not do coincide up to any positive multiplicative  factor. Hence
$$
{\mu_1}_{|Y_1}\neq  {\mu_2}_{|Y_1}.
$$ 
Combining Proposition~\ref{p1j70} and, the already proved case of probability
measures, we thus  conclude that the measures ${\mu_1}_{|Y_1}$ and $
{\mu_2}_{|Y_1}$ are mutually singular. This  means that there exists
a set $Z_n \sbt Y_n$ such that $\mu_1(Z_n)=0$ and $\mu_2(Y_n\sms
Z_n)=0$. Setting  
$$
Z:=\bigcup_{n=1}^\infty Z_n,
$$
we thus  have
$$
\mu_1(Z)=\sum_{n=1}^\infty \mu_1(Z_n)=0
$$ 
and  
$$
\mu_2(X\sms Z)=
\sum_{n=1}^\infty \mu_2(Y_n\sms Z)
=\sum_{n=1}^\infty \mu_2(Y_n\sms Z_n)=0.
$$ 
So,the measures $\mu_1$ and $\mu_2$ are mutually singular, and we are
done. 
\endpf

\bthm\label{5.2}
Suppose that $\mathfrak{F}$ is a $\sg$--algebra on a set $X$ and $T:X\lra X$ is a measurable map. Let $\mu$ be $\sg$--finite measure invariant under $T$. 

If $T:X\lra X$ is conservative with respect to $\mu$,
then $\mu$ is ergodic if and only if there is no $\sg$--finite measure $\nu$ on $(X,\mathfrak{F})$ invariant under $T$ such that  
$$
\nu\prec\hspace{-0.1cm}\prec\mu 
\  \  \  {\rm and} \  \  \  
\nu\neq t\mu.
$$
for every real $t>0$. 
\ethm

\bpf
First, suppose that $\mu$ is ergodic. Let $\nu$ be a $\sg$--finite measure on $(X,\mathfrak{F})$ invariant under $T$ and such that 
$$
\nu\prec\hspace{-0.1cm}\prec \mu
$$
We claim that $\nu$ is ergodic too.
Indeed, suppose by way of contradiction that there exists $A\in\mathfrak{F}$ such that 
$$
T^{-1}(A)=A
$$ 
with 
$$
\nu(A)>0\  \  \  {\rm and} \  \  \ 
\nu(X\backslash A)>0.
$$
Since $\nu\prec\hspace{-0.1cm}\prec\mu$, it follows that $\mu(A)>0$ and $\mu(X\backslash A)>0$. 
This contradicts the ergodicity of $\mu$. So $\nu$ is ergodic. 

Since $\mu$ is conservative and $\nu$ is absolutely continuous with respect to $\mu$, the $T$--invariant measure $\nu$ is conservative too. So, if 
$\nu\neq t\mu$ for every real $t>0$, then Theorem~\ref{t1j74.1}
asserts that $\nu\bot\mu$. This contradicts the hypothesis that $\nu\prec\hspace{-0.1cm}\prec\mu$. Hence,
$$
\nu=t\mu
$$
for some real $t>0$.

For the converse implication, suppose that $\mu$ is not ergodic (but still $T$--invariant and conservative by hypothesis). Then there exists  $A\in\mathfrak{F}$ such that 
$$
T^{-1}(A)=A
$$ 
with 
$$
\mu(A)>0\  \  \  {\rm and} \  \  \ 
\mu(X\backslash A)>0.
$$
Let $\mu^*_A$ be the $\sg$--finite measure on $(X,\mathfrak{F})$ defined by the formula
$$
\mu^*_A(F)=\mu(F\cap A)
$$
for all $F\in \mathfrak{F}$. One immediately
verifies that $\mu^*_A$ is a $T$--invariant probability measure such that
$$
\mu^*_A\prec\hspace{-0.1cm}\prec\mu 
\  \  \  {\rm and} \  \  \ \mu^*_A\neq t\mu
$$
for any real $t>0$.
\epf

\  

\section[Absolutely Continuous  $\sg$-Finite  Invariant Measures] {Absolutely Continuous  $\sg$-Finite  Invariant Measures; Marco Martens's Approach}

In this section we establish  a very useful, relatively easy to
verify, sufficient condition  for a quasi-invariant measure  to admit
an  absolutely  continuous $\sg$--finite  invariant measure. This
condition goes back to the work \cite{Mar}  of  Marco Martens. It
has been used many times, notably in \cite{KU1}, and obtained its
final form in \cite{SU}. In contrast to Martens's paper \cite{Mar}, where $\sg$--compact metric  spaces form the  setting, the sufficient condition  in
\cite{SU} is stated  for abstract measure spaces, and the proof
utilizes  the concept of  Banach limits rather than weak
convergence. We  start with  the following  uniqueness result.

\bthm\label{t1j77} Suppose that $T:(X, \mathfrak{F}, m)\lra (X,
\mathfrak{F}, m)$ is a measurable map  of a $\sg$--finite measure
space $X$. Suppose also that the measure $m$ is quasi--invariant. Then up
to a positive multiplicative constant there exists at most one
non-zero $\sg$-finite $T$--invariant measure $\mu$  absolutely
continuous with respect  to the measure $m$. 
\ethm

\bpf Suppose that $\mu_1$ and $\mu_2$  are   $\sg$-finite
non-zero $T$-invariant  measures absolutely continuous with respect
to the measure $m$. Since $m$ is ergodic and conservative, so are
the measures $\mu_1$ and $\mu_2$. It now follows from
Theorem~\ref{t1j74.1} that if $\mu_1$ and $\mu_2$ do not coincide up
to a positive multiplicative constant, then these two measures are
mutually singular. But this means that there exists a measurable set
$Y \sbt X$  such that $\mu_1(Y)=0$ and $\mu_2(X\sms Y)=0$. So
\beq\label{1j79}
\mu_1\lt(\bigcup_{n=0}^\infty T^{-n}(Y)\rt)
\leq  \sum_{n=0}^\infty
\mu_1(T^{-n}(Y))
\leq \sum _{n=0}^\infty 0=0.
\eeq
On the other hand $\mu_2(Y)>0$, so $m(Y)>0$, and $m(X \sms
\bigcup_{n=0}^\infty T^{-n}(Y))=0$ by virtue of Theorem~\ref{t1j59}.
Since $\mu_1$  is absolutely  continuous with respect to $m$, this
implies that $\mu_1(X\sms \bigcup_{n=0}^\infty T^{-n}(Y))=0$. Along
with (\ref{1j79}) this  gives  that $\mu_1(X)=0$. This contradiction
finishes the proof. 
\endpf

\sp We now introduce the concept of Marco Martens map.

\bdfn\label{d:mmmap} 
Let $T:(X,\F)\lra (X,\F)$ be a measurable transformation. 
Let also $m$ be a quasi--invariant probability measure on $(X,\F)$ with respect to $T$. 
The transformation $T$ is called a  Marco Martens map
\index{(N)}{Marco Martens map} 
if it admits a countable family $\{X_n\}_{n=0}^\infty$ of subsets of $X$
with the following properties:
\begin{itemize}
\item[(a)]
$X_n\in\F$ for all $n\geq0$.
\vspace{1.5mm}
\item[(b)]
$m\(X\backslash\bigcup_{n=0}^\infty X_n\)=0$.
\vspace{1.5mm}
\item[(c)]
For all $m,n\geq0$, there exists $j\geq0$ such that $m(X_m\cap T^{-j}(X_n))>0$.
\vspace{1.5mm}
\item[(d)]
For all $j\geq0$ there exists $K_j\geq1$ such that for all $A,B\in\F$ 
with $A\cup B\sbt X_j$ and for all $n\geq0$,
\[
m(T^{-n}(A))\,m(B)\leq K_j\,m(A)\,m(T^{-n}(B)).
\]
\item[(e)] $\displaystyle\sum_{n=0}^\infty m(T^{-n}(X_0))=+\infty$.
\item[(f)] $\displaystyle T\Bigl(\bigcup_{j=l}^\infty Y_j\Bigr)\in\F$ 
for all $j\geq0$, where $\displaystyle Y_j:=X_j\Bigl\backslash\Bigr.\bigcup_{i<j}X_i$.
\item[(g)] $\displaystyle\lim_{l\rightarrow\infty}m\lt(T\Bigl(\bigcup_{j=l}^\infty Y_j\Bigr)\rt)=0$.
\end{itemize}
The family $\{X_n\}_{n=0}^\infty$ is called a Marco Martens cover.
\edfn

\brem\label{r:mmmapf1}

\

\begin{itemize}
\item[(1)] Of course, condition~(b) follows from the stronger hypothesis that $\bigcup_{n=0}^\infty X_n=X$. 

\,

\item[(2)] Condition~(c) implies that $m(X_n)>0$ for all $n\geq0$.  

\,

\item[(3)] In light of Corollary~\ref{Cor5.4a}, if $T$ is conservative with respect to $\mu$ then condition~(e) is fulfilled. 
           
\,
           
\item[(4)] In conditions~(f-g), note that 
$$
\bigcup_{j=l}^\infty Y_j
=\bigcup_{j=0}^\infty X_j\backslash\bigcup_{i<l}X_i
\sbt X\backslash\bigcup_{i<l}X_i
$$.
\item[(5)] If the map $T:X\lra X$ is finite--to--one, then condition~(g) is satisfied. For then, 
$$
\bigcap_{l=1}^\infty T\lt(\bigcup_{j=l}^\infty Y_j\rt)=\emptyset.
$$
\end{itemize}
\erem

Let $l^\infty$ denote the Banach space of all bounded real--valued sequences 
$x=(x_n)_{n=1}^\infty$ with norm 
$$
\|x\|_\infty:=\sup_{n\in\N}|x_n|.
$$
Recall that a Banach limit is a shift--invariant positive bounded/continuous linear functional 
\beq\lab{220200106}
l_B:l^\infty\lra\R
\eeq
which extends the usual limits. More precisely, for all sequences 
$$
x=(x_n)_{n=1}^\infty, \  
y=(y_n)_{n=1}^\infty\in l^\infty
$$ 
and $\alpha,\beta\in\R$, the following properties hold:

\,

\begin{itemize}
\item[(a)] $l_B(\alpha x+\beta y)=\alpha\,l_B(x)+\beta\,l_B(y)$ (linearity).
\vspace{2mm}
\item[(b)] $\|l_B\|:=\sup\bigl\{|l_B(x)|:\,\|x\|_\infty\leq1\bigr\}<\infty$ (continuity/boundedness).
\vspace{2mm}
\item[(c)] If $x\geq0$, i.e. $x_n\geq0$ for all $n\in\N$, then $l_B(x)\geq0$ (positivity).
\vspace{2mm}
\item[(d)] $l_B(\sg(x))=l_B(x)$, where $\sg:l^\infty\to l^\infty$ is the (left) shift map 
defined by $(\sg(x))_n=x_{n+1}$ for all $n\in\N$ (shift-invariance).
\vspace{2mm}
\item[(e)] If $x$ is a convergent sequence, then $l_B(x)=\displaystyle\lim_{n\to\infty}x_n$.
\end{itemize}

\sp\fr It follows from properties~(a),~(c) and~(e) that a Banach limit also satisfies:

\,

\begin{itemize}
\item[(f)] $\displaystyle\liminf_{n\to\infty}x_n\leq l_B(x)\leq\limsup_{n\to\infty}x_n$.

\,

\item[(g)] If $x\leq y$, i.e. $x_n\leq y_n$ for all $n\in\N$, then $l_B(x)\leq l_B(y)$.
\end{itemize}

\sp\bthm\label{t1h75} 
Let $(X,\F,m)$ be a probability space and $T:X\lra X$ a Marco Martens map with Marco Martens cover $\{X_j\}_{j=0}^\infty$. Then 

\begin{itemize}
\item There exists a $\sigma$--finite $T$--invariant measure $\mu$ on $X$ which is equivalent to $m$. 

\,

\item In addition, $0<\mu(X_j)<\infty$ for each $j\geq0$. 
\end{itemize}

\,

\fr A measure $\mu$ with the above properties can be constructed as follows.

\,

\begin{itemize} 
\item Let $l_B:l^\infty\lra\R$ be a Banach limit and let $Y_j:=X_j\backslash\bigcup_{i<j}X_i$ for every $j\geq0$. 

\,

\item For each $A\in\F$, set
\beq\label{5.9a}
m_n(A):=\frac{\displaystyle\sum_{k=0}^n m(T^{-k}(A))}{\displaystyle\sum_{k=0}^n m(T^{-k}(X_0))}.
\eeq

\item If $A\in\mathcal{A}$ and $A\sbt Y_j$ for some $j\geq0$, then
$(m_n(A))_{n=1}^\infty\in l^\infty$ and set
\beq\label{5.9b}
\mu(A):=l_B\bigl((m_n(A))_{n=1}^\infty\bigr).
\eeq

\,

\item For a general $A\in\F$, set
\[
\mu(A):=\sum_{j=0}^\infty\mu(A\cap Y_j).
\]
\item In addition, if $(m_n(A))_{n=1}^\infty\in l^\infty$ for some $A\in\F$, 
then 
\beq\label{eq:muequ} 
\mu(A)
=l_B\bigl((m_n(A))_{n=1}^\infty\bigr)
-\lim_{l\rightarrow\infty}l_B\Bigl(\Bigl(m_n\Bigl(A\cap\bigcup_{j=l}^\infty Y_j\Bigr)\Bigr)_{n=0}^\infty\Bigr).
\eeq
\item In particular, if $A\in\F$ is contained in a finite union 
of sets $X_j$, $j\geq0$, then
\beq\label{120190618}
\mu(A)=l_B\bigl((m_n(A))_{n=1}^\infty\bigr). 
\eeq
\end{itemize}

\,

Finally, if $T$ is ergodic and conservative with respect to $m$, 
then $\mu$ is also ergodic and conservative, and unique up to a positive multiplicative constant.
\ethm

In order to prove Theorem~\ref{t1h75}, we need several lemmas.

\blem\label{l:dumeasure} 
Let $(Z,\mathcal{F})$ be a measurable space such that:
\begin{itemize}
\item[(a)] $\displaystyle Z=\bigcup_{j=0}^\infty Z_j$ for some mutually disjoint sets $Z_j\in\mathcal{F}$; and 
\item[(b)] $\nu_j$ is a finite measure on $Z_j$ for each $j\geq0$.
\end{itemize}
Then the set function $\nu:\mathcal{F}\lra[0,\infty]$ defined by
\[
\nu(F):=\sum_{j=0}^\infty\nu_j(F\cap Z_j)
\] 
is a $\sigma$--finite measure on $Z$.
\elem

\bpf Clearly, $\nu(\emptyset)=0$. 
Let $F\in\mathcal{F}$ and let $\{F_n\}_{n=1}^\infty$ be a partition 
of $F$ into sets in $\mathcal{F}$. Then
\begin{align*}
\nu(F)
=&\sum_{j=0}^\infty\nu_j(F\cap Z_j)
=\sum_{j=0}^\infty\nu_j\Bigl(\bigcup_{n=1}^\infty(F_n\cap Z_j)\Bigr) \\
=&\sum_{j=0}^\infty\sum_{n=1}^\infty\nu_j(F_n\cap Z_j)
=\sum_{n=1}^\infty\sum_{j=0}^\infty\nu_j(F_n\cap Z_j) \\
&=\sum_{n=1}^\infty\nu(F_n),
\end{align*}
where the order of summation could be changed since all terms
involved are non-negative. Thus, $\nu$ is a measure. Moreover, by definition, 
$$
Z=\bigcup_{j=0}^\infty Z_j
$$ 
and $\nu(Z_j)=\nu_j(Z_j)<\infty$ for all $j\geq0$. Therefore $\nu$ is $\sg$--finite.
\epf

\sp From this point on, all lemmas rely on the same main hypotheses as Theorem~\ref{t1h75}.

\blem\label{mn}
For all $n,j\geq0$ and all $A,B\in\F$ 
with $A\cup B\sbt X_j$, we have
\[
m_n(A)\,m(B)\leq K_j\,m(A)\,m_n(B).
\]
\elem

\bpf
This follows directly from the definition of $m_n$ and condition~(d) of Definition~\ref{d:mmmap}.
\epf

\blem\label{l:muxjfin} 
For every $j\geq0$, we have
$(m_n(X_j))_{n=1}^\infty\in l^\infty$ and $\mu(Y_j)\leq\mu(X_j)<\infty$.
\elem

\bpf Fix $j\geq0$. In virtue of condition~(c) of Definition~\ref{d:mmmap}, there
exists $q\geq0$ such that $m(X_j\cap T^{-q}(X_0))>0$. By Lemma~\ref{mn} and the definition of $m_n$, 
for all $n\geq0$ we have that
\begin{eqnarray}
m_n(Y_j)
\leq m_n(X_j)
&\leq&K_j\frac{m(X_j)}{m\bigl(X_j\cap T^{-q}(X_0)\bigr)}m_n\bigl(X_j\cap T^{-q}(X_0)\bigr) \nonumber \\
&\leq&K_j\frac{m(X_j)}{m\bigl(X_j\cap T^{-q}(X_0)\bigr)}m_n(T^{-q}(X_0)) \nonumber \\
&=&K_j\frac{m(X_j)}{m\bigl(X_j\cap T^{-q}(X_0)\bigr)}
   \frac{\sum_{k=0}^{n+q}m(T^{-k}(X_0))}{\sum_{k=0}^nm(T^{-k}(X_0))} \nonumber \\
&=&K_j\frac{m(X_j)}{m\bigl(X_j\cap T^{-q}(X_0)\bigr)}
     \left[1+\frac{\sum_{k=n+1}^{n+q}m(T^{-k}(X_0))}{\sum_{k=0}^nm(T^{-k}(X_0))}\right] \nonumber \\
&\leq&K_j\frac{m(X_j)}{m\bigl(X_j\cap T^{-q}(X_0)\bigr)} 
     \left[1+\frac{q}{m(X_0)}\right]. \label{NCPIIp110}
\end{eqnarray}
Consequently, $(m_n(X_j))_{n=1}^\infty\in l^\infty$ and properties~(g) and~(e) 
of a Banach limit yield that
\[
\mu(Y_j)
\leq K_j\frac{m(X_j)}{m\bigl(X_j\cap T^{-q}(X_0)\bigr)}\left[1+\frac{q}{m(X_0)}\right]
<+\infty.
\]
Since $X_j=\bigcup_{i=0}^j Y_i$ and the $Y$s are mutually disjoint, we deduce that 
\[
\mu(Y_j)
\leq\sum_{i=0}^j\mu(X_j\cap Y_i)
=\sum_{i=0}^\infty\mu(X_j\cap Y_i) 
=:\mu(X_j)
\leq\sum_{i=0}^j\mu(Y_i)
<+\infty.
\]
\epf

Now, for every $j\geq0$, set $\mu_j:=\mu|_{Y_j}$.

\blem\label{l:mumcomp} 
For every $j\geq0$ such that $\mu(Y_j)>0$ and
for every measurable set $A\sbt Y_j$, we have
\[
K_j^{-1}\frac{\mu(Y_j)}{m(Y_j)}m(A)
\leq\mu_j(A)
\leq K_j\frac{\mu(Y_j)}{m(Y_j)}m(A).
\]
\elem

\bpf 
This follows from the definition of the measure $\mu$ 
and by setting $B=Y_j$ in Lemma~\ref{mn} and using  
properties~(a) and~(g) of a Banach limit.
\epf

\blem\label{l:mujcam} 
For each $j\geq0$, $\mu_j$ is a finite measure on $Y_j$.
\elem

\bpf Let $j\geq0$. Assume without loss of
generality that $\mu_j(Y_j)>0.$ Let $A\sbt Y_j$ be a
measurable set and $(A_k)_{k=1}^\infty$ a countable
measurable partition of $A$. Using termwise operations
on sequences, for every $l\in\N$ we have
\begin{eqnarray*}
\left(\sum_{k=1}^\infty m_n(A_k)\right)_{n=1}^\infty - \sum_{k=1}^l\left(m_n(A_k)\right)_{n=1}^\infty 
&=&\left(\sum_{k=1}^\infty m_n(A_k)\right)_{n=1}^\infty - \left(\sum_{k=1}^l m_n(A_k)\right)_{n=1}^\infty \\ 
&=&\left(\sum_{k=l+1}^\infty m_n(A_k)\right)_{n=1}^\infty.
\end{eqnarray*}
It therefore follows from Lemma~\ref{mn} (with $A=A_k$ and $B=Y_j$) that
\begin{eqnarray*} 
\left\|\left(\sum_{k=1}^\infty m_n(A_k)\right)_{n=1}^\infty - \sum_{k=1}^l\left(m_n(A_k)\right)_{n=1}^\infty\right\|_\infty 
&=&\left\|\left(\sum_{k=l+1}^\infty m_n(A_k)\right)_{n=1}^\infty\right\|_\infty \\
&\leq&\left\|\frac{K_j}{m(Y_j)}\left(m_n(Y_j)\sum_{k=l+1}^\infty m(A_k)\right)_{n=1}^\infty\right\|_\infty \\
&=&\frac{K_j}{m(Y_j)}\left\|\left(m_n(Y_j)\sum_{k=l+1}^\infty m(A_k)\right)_{n=1}^\infty\right\|_\infty.
\end{eqnarray*}
Since $(m_n(Y_j))_{n=1}^\infty\in l^\infty$ by Lemma~\ref{l:muxjfin} 
and since $\lim_{l\rightarrow\infty}\sum_{k=l+1}^\infty m(A_k)=0$, we 
conclude that 
\[
\lim_{l\rightarrow\infty}\left\|\left(\sum_{k=1}^\infty m_n(A_k)\right)_{n=1}^\infty - \sum_{k=1}^l(m_n(A_k))_{n=1}^\infty\right\|_\infty=0.
\] 
This means that 
\[
\left(\sum_{k=1}^\infty m_n(A_k)\right)_{n=1}^\infty
=\sum_{k=1}^\infty\(m_n(A_k)\)_{n=1}^\infty
\] 
in $l^\infty$. Hence, using the continuity of the Banach limit $l_B:l^\infty\rightarrow\R$, we get
\begin{eqnarray*}
\mu(A) 
&=&l_B\bigl((m_n(A))_{n=1}^\infty\bigr)
 = l_B\lt(\Bigl(m_n\Bigl(\bigcup_{k=1}^\infty A_k\Bigr)\Bigr)_{n=1}^\infty\rt)
 = l_B\lt(\Bigl(\sum_{k=1}^\infty m_n(A_k)\Bigr)_{n=1}^\infty\rt) \\
&=&\sum_{k=1}^\infty l_B\Big((m_n(A_k))_{n=1}^\infty\Big) \\
&=& \sum_{k=1}^\infty\mu(A_k).
\end{eqnarray*}
So $\mu_j$ is countably additive. Also, $\mu_j(\emptyset)=0$. Thus $\mu_j$ is a measure. 
By Lemma~\ref{l:muxjfin}, $\mu_j$ is finite.
\epf

\sp Combining Lemmas~\ref{l:dumeasure},~\ref{l:muxjfin},~\ref{l:mumcomp} 
and~\ref{l:mujcam}, and condition~(b) of Definition~\ref{d:mmmap}, 
we get the following.

\blem\label{l:musfmeas} 
We have that $\mu$ is a $\sigma$--finite measure on $X$
equivalent to $m$. In addition, 
$$
\mu(Y_j)\leq\mu(X_j)<\infty
\  \  \  {\rm and} \  \  \
\mu(X_j)>0
$$ 
for all integers $j\geq0$.
\elem

\blem\label{l:mueqpf} 
Formula~(\ref{eq:muequ}) holds.
\elem

\bpf 
Fix $A\in\F$ such that $(m_n(A))_{n=1}^\infty\in l^\infty$. 
Then for every $l\in\N$ we have that
\begin{align*}
l_B\bigl((m_n(A))_{n=1}^\infty\bigr)
=&l_B\lt(\sum_{j=0}^l(m_n(A\cap Y_j))_{n=1}^\infty\rt)
 +l_B\lt(\Bigl(m_n\Bigl(\bigcup _{j=l+1}^\infty A\cap Y_j\Bigr)\Bigr)_{n=1}^\infty\rt) \\
=&\sum_{j=0}^l l_B\bigl((m_n(A\cap Y_j))_{n=1}^\infty\bigr)
 +l_B\lt(\Bigl(m_n\Bigl(A\cap\bigcup_{j=l+1}^\infty Y_j\Bigr)\Bigr)_{n=1}^\infty\rt).
\end{align*}
Letting $l\rightarrow\infty$, we thus obtain that
\begin{align*}
l_B\bigl((m_n(A))_{n=1}^\infty\bigr)
=&\sum_{j=0}^\infty l_B\bigl((m_n(A\cap Y_j))_{n=1}^\infty\bigr)
 +\lim_{l\rightarrow\infty}l_B\lt(\Bigl(m_n\Bigl(A\cap\bigcup_{j=l+1}^\infty Y_j\Bigr)\Bigr)_{n=1}^\infty\rt) \\
=&\sum_{j=0}^\infty\mu(A\cap Y_j)
 +\lim_{l\rightarrow\infty}l_B\lt(\Bigl(m_n\Bigl(A\cap\bigcup_{j=l}^\infty Y_j\Bigr)\Bigr)_{n=1}^\infty\rt) \\
=&\mu(A)+\lim_{l\rightarrow\infty}l_B\lt(\Bigl(m_n\Bigl(A\cap\bigcup_{j=l}^\infty Y_j\Bigr)\Bigr)_{n=1}^\infty\rt).
\end{align*}
This establishes formula~(\ref{eq:muequ}). In particular, if $A\sbt\bigcup_{j=0}^k X_j$ for some $k\in\N$, then 
$$
A\cap\bigcup_{j=l}^\infty Y_j
\sbt\bigl(\bigcup_{j=0}^k X_j\bigr)\cap\bigl(X\backslash\bigcup_{i<l}X_i\bigr)
=\emptyset
$$ 
for all $l>k$. In that case, the equation above reduces to
\[
l_B\bigl((m_n(A))_{n=1}^\infty\bigr)=\mu(A). 
\]
\epf

\blem\label{l:mutinv}
The $\sigma$--finite measure $\mu $ is $T$--invariant.
\elem

\bpf
Let $i\geq0$ be such that $m(Y_i)>0$. Fix a measurable set $A\subset Y_i$.
By definition, 
$$
\mu(A)=l_B\bigl((m_n(A))_{n=1}^\infty\bigr).
$$
Furthermore, for all $n\geq0$ notice that
\[
\bigl|m_n(T^{-1}(A))-m_n(A)\bigr|
=\frac{\bigl|m(T^{-(n+1)}(A))-m(A)\bigr|}{\sum_{k=0}^nm(T^{-k}(X_0))}
\leq\frac{1}{\sum_{k=0}^nm(T^{-k}(X_0))}.
\]
Thus, $(m_n(T^{-1}(A)))_{n=1}^\infty\in l^\infty$ because $(m_n(A))_{n=1}^\infty\in l^\infty$. Moreover,
by condition~(e) of Definition~\ref{d:mmmap}, it follows from the above and properties~(a),~(e) and~(g) of a Banach 
limit that 
$$
l_B\bigl((m_n(T^{-1}(A)))_{n=1}^\infty\bigr)=l_B\bigl((m_n(A))_{n=1}^\infty\bigr)=\mu(A).
$$  
 
Keep $A$ a measurable subset of $Y_i$. Fix $l\in\N$. We then have
\begin{align*}
m_n\lt(T^{-1}(A)\cap\bigcup_{j=l}^\infty Y_j\rt)
&=\frac{\sum_{k=0}^n m\bigl(T^{-k}\bigl(T^{-1}(A)\cap\bigcup_{j=l}^\infty Y_j\bigr)\bigr)}{\sum_{k=0}^n m(T^{-k}(X_0))} \\
&\leq\frac{\sum_{k=0}^n m\bigl(T^{-(k+1)}\bigl(A\cap T(\bigcup_{j=l}^\infty Y_j)\bigr)\bigr)}{\sum_{k=0}^n m(T^{-k}(X_0))} \\
&\leq m_{n+1}\Bigl(A\cap T\Bigl(\bigcup_{j=l}^\infty Y_j\Bigr)\Bigr)
      \cdot\frac{\sum_{k=0}^{n+1}m(T^{-k}(X_0))}{\sum_{k=0}^n m(T^{-k}(X_0))} \\
&\leq K_i\frac{m_{n+1}(Y_i)}{m(Y_i)}
      \cdot m\Bigl(A\cap T\Bigl(\bigcup_{j=l}^\infty Y_j\Bigr)\Bigr)
			\cdot\frac{\sum_{k=0}^{n+1}m(T^{-k}(X_0))}{\sum_{k=0}^n m(T^{-k}(X_0))},
\end{align*}
where the last inequality sign holds by Lemma~\ref{mn} since $A\sbt Y_i$. 
When $n\rightarrow\infty$, the last quotient on the right-hand side approaches $1$.
Therefore
\[
0
\leq l_B\Bigl(\Bigl(m_n\Bigl(T^{-1}(A)\cap\bigcup_{j=l}^\infty Y_j\Bigr)\Bigr)_{n=1}^\infty\Bigr)
\leq K_i\frac{\mu(Y_i)}{m(Y_i)}m\Bigl(T\Bigl(\bigcup_{j=l}^\infty Y_j\Bigr)\Bigr).
\]
Hence, by virtue of condition~(g) of Definition~\ref{d:mmmap},
\[
0
\leq\lim_{l\rightarrow\infty}l_B\lt(\Bigl(m_n\Bigl(T^{-1}(A)\cap\bigcup_{j=l}^\infty Y_j\Bigr)\Bigr)_{n=1}^\infty\rt)
\leq K_i\frac{\mu(Y_i)}{m(Y_i)}\lim_{l\rightarrow\infty}m\lt(T\Bigl(\bigcup_{j=l}^\infty Y_j\Bigr)\rt)=0.
\]
So
\[
\lim_{l\rightarrow\infty}l_B\lt(\Bigl(m_n\Bigl(T^{-1}(A)\cap\bigcup_{j=l}^\infty Y_j\Bigr)\Bigr)_{n=1}^\infty\rt)
=0.
\]
It thus follows from Lemma~\ref{l:mueqpf} that
\[
\mu(T^{-1}(A))
=l_B\bigl((m_n(T^{-1}(A)))_{n=1}^\infty\bigr)
=l_B\bigl((m_n(A))_{n=1}^\infty\bigr)
=\mu(A).
\]
For an arbitrary $A\in\F$, write $A=\bigcup_{j=0}^\infty A\cap Y_j$ and
observe that
\[
\mu(T^{-1}(A))
=\mu\lt(\bigcup_{j=0}^\infty T^{-1}(A\cap Y_j)\rt)
=\sum_{j=0}^\infty\mu\bigl(T^{-1}(A\cap Y_j)\bigr)
=\sum_{j=0}^\infty\mu(A\cap Y_j)
=\mu(A).
\]
We are done.
\epf

\vspace{2mm}

\noindent{\bf Proof of Theorem~\ref{t1h75}:} Combining 
Lemmas~\ref{l:muxjfin},~\ref{l:musfmeas},~\ref{l:mueqpf},~\ref{l:mutinv},
and~\ref{t1j77}, we obtain the statement of Theorem~\ref{t1h75}. 
\qed

\sp\brem\label{r5.21a} In the course of the proof  of
Theorem~\ref{t1h75} we have shown that $$ 0< \inf\{m_n(A): n\geq
1\}\leq \sup\{m_n(A): n\geq 1\} < +\infty $$ for all $j\geq 0$
and all measurable sets $A\sbt X_j$, such that $m(A)>0$. \erem

\sp\chapter{Probability (Finite) Invariant Measures; Basic Properties and Existence}\label{finite-measure 1} 

\section{Basic Properties of Probability Invariant Measures}\label{BPoPIM}

Invoking Theorem~\ref{Poincare}, as an immediate consequence of Theorem~\ref{t1j74.1}, we get the following. 

\sp\bthm\label{t1j74.1-2019-11-18} 
Suppose that $\mathfrak{F}$ is a $\sg$--algebra
on a set $X$, $T:X\lra X$ is a measurable map, and $\mu_1$  and
$\mu_2$ are two $T$--invariant probabbility measures on $\mathfrak{F}$. 

If $T:X\lra X$ is ergodic with respect to both measures
$\mu_1$ and $\mu_2$, then either $\mu_1$ and $\mu_2$ are mutually
singular or are equal. 
\ethm

Invoking Theorem~\ref{Poincare} again, as an immediate consequence of Theorem~\ref{5.2}, we get the following. 

\sp\bthm\label{5.2-2019-11-18}
Suppose that $\mathfrak{F}$ is a $\sg$--algebra on a set $X$ and $T:X\lra X$ is a measurable map. If $\mu$ is a $T$--invariant probability measure on $\mathfrak{F}$,
then $\mu$ is ergodic if and only if there is no $T$--invariant probability measure $\nu$ on $\mathfrak{F}$ such that  
$$
\nu\prec\hspace{-0.1cm}\prec\mu 
\  \  \  {\rm and} \  \  \  
\nu\neq \mu.
$$
\ethm

Recall that in a vector space the extreme points of a convex set are those points which cannot be represented as a non--trivial convex combination of two distinct points of the set. 
More precisely, let $V$ be a vector space and $C$ be a convex subset of $V$. 
A vector $v\in C$ is called an extreme point of $C$ if and only if the only combination of distinct vectors $v_1,v_2\in C$ such that 
$$
v=\a v_1+(1-\a)v_2
$$ 
for some $\a\in[0,1]$ is a combination with $\a=0$ or $\a=1$.

If $T:(X,\mathfrak{F})\lra(X,\mathfrak{F})$ is a measurable transformation, then $M(X,\mathfrak{F})$ is a convex subset of the vector space $SM(X,\mathfrak{F})$ of all signed measures on $(X,\mathfrak{F})$, and $M(T,\mathfrak{F})$ is a convex subset of $M(X,\mathfrak{F})$, thus a convex subset of the vector space $SM(X,\mathfrak{F})$.

\sp\bthm\label{ext}
Let $T:(X,\mathfrak{F})\lra(X,\mathfrak{F})$ be a measurable transformation.
The set $E(T,\mathfrak{F})$ of ergodic $T$--invariant measures of $T$ coincides with the set of extreme points of the set of $T$--invariant probability measures $M(T,\mathfrak{F})$.
\ethm

\bpf
Suppose for a contradiction that some measure $\mu\in E(T,\mathfrak{F})$ is not an extreme point of $M(T,\mathfrak{F})$. Then there exist measures $\mu_1\neq\mu_2$ in $M(T,\mathcal{A})$ and $0<\a<1$
such that 
$$
\mu=\a\mu_1+(1-\a)\mu_2.
$$
It follows immediately that
$$
\mu_1\prec\hspace{-0.1cm}\prec\mu
\  \  \  {\rm and} \  \  \
\mu_2\prec\hspace{-0.1cm}\prec\mu.
$$
By Theorem~\ref{5.2-2019-11-18}, we deduce from the ergodicity of $\mu$ that $$
\mu_1=\mu=\mu_2.
$$
This contradicts the fact that $\mu_1\neq\mu_2$. Thus $\mu$ is an extreme point of $M(T,\mathfrak{F})$.

In order to prove the converse implication, let $\mu\in M(T,)\backslash E(T,\mathfrak{F})$.
We want to show that $\mu$ is not an extreme point of $M(T,\mathfrak{F})$. Since $\mu$ is
not ergodic, there exists a set $A\in\mathfrak{F}$ such that 
$$
T^{-1}(A)=A
$$
with 
$$
\mu(A)>0
\  \  \  {\rm and} \  \  \
\mu(X\backslash A)>0.
$$ 
Observe that now $\mu$ can be written as the following non--trivial convex combination
of the $T$--invariant conditional measures $\mu_A$ and $\mu_{X\backslash A}$. Indeed, for every $B\in \mathfrak{F}$, we have that
\begin{eqnarray*}
\mu(B)
&=&\mu(A\cap B)+\mu((X\backslash A)\cap B)
=\mu(A)\,\mu_A(B)+\mu(X\backslash A)\,\mu_{X\backslash A}(B)\\
&=&\mu(A)\,\mu_A(B)+(1-\mu(A))\,\mu_{X\backslash A}(B).
\end{eqnarray*}
Thus
$$
\mu=\mu(A)\,\mu_A+(1-\mu(A))\,\mu_{X\backslash A}
$$
So, $\mu$ is a non--trivial convex combination of two
distinct $T$--invariant probability measures and thus is 
not an extreme point of $M(T,\mathfrak{F})$.
\epf

\sp

\section[Existence of Borel Probability Invariant Measures] {Existence of Borel Probability Invariant Measures; Bogolubov--Krylov Theorem}\label{existinv} 

\sp In general, a given measurable transformation $T:(X,\F)\lra (X,\F)$
may not admit any invariant measure. However, in one situation, quite common, although in general failing in this book, namely, for $X$ a compact metrizable space and the Borel $\sigma$--algebra $\mathcal{B}=\mathcal{B}_X$, on $X$, every continuous transformation does have a probability invariant measure. This is the celebrated Bogolubov--Krylov Theorem. We call all continuous self--maps of compact metrizable spaces topological dynamical systems. Before proving the Bogolubov--Krylov Theorem, we briefly examine the map 
$$
M(X)\ni\mu\longmapsto\mu\circ T^{-1}.
$$
By the Riesz Representation Theorem, the set $M(X)$ can be identified with a subset of the normed
space $C^*(X)$, the space of all bounded linear functionals on $X$. The set $M(X)$ is obviously convex. If $C^*(X)$ is endowed with the weak* topology, then by the Banach--Alaoglu Theorem, the set $M(X)$ becomes compact and metrizable.
The following lemma will be helpful in proving the existence of invariant measures.
\blem\label{cont}
Let $T:X\lra X$ be a continuous transformation of a compact metrizable space $X$. Then the map
$S:M(X)\lra M(X)$, given by the formula
$$
S(\mu):=\mu\circ T^{-1},
$$
is a continuous affine map.
\elem

\bpf The proof of affinity is left to the reader; here we concentrate on the
continuity of $S$. Let $(\mu_n)_{n=1}^\infty$ be a sequence in $M(X)$ which
weak* converges to $\mu$. Then, for any $f\in C(X)$, we have that
\begin{eqnarray*}
\lim_{n\to\infty}\int_X f\,d(S(\mu_n))
&=&\lim_{n\to\infty}\int_X f\,d(\mu_n\circ T^{-1})
=\lim_{n\to\infty}\int_X f\circ T\,d\mu_n\\
&=&\int_X f\circ T\,d\mu
=\int_X f\,d(\mu\circ T^{-1})\\
&=&\int_X f\,d(S(\mu)).
\end{eqnarray*}
Since $f$ was chosen arbitrarily in $C(X)$, the sequence $(S(\mu_n))_{n=1}^\infty$
weak* converges to $S(\mu)$. Thus, $S$ is continuous.
\epf

We come now to the main result of this section, namely, showing that continuous map on a compact metric space admits at least one invariant measure. This theorem is not very difficult to prove, but it is obviously important. For this reason, we  provide two different proofs. The first one depends on some basic functional analysis, whereas the second one is rather more constructive.

\sp\bthm[Bogolubov--Krylov Theorem]\label{KBT}
Any topological dynamical system, i.e. any continuous transformation $T:X\lra X$ of a compact metrizable space $X$ admits a $T$--invariant Borel probability measure.
\ethm

\bpf
We treat $M(X)$ as a subspace of the vector topological space $C^*(X)$ endowed with the weak* topology By Lemma~\ref{cont} we know that the map given by the formula
$$
S(\mu)=\mu\circ T^{-1}
$$ 
is a continuous affine map
of the convex, compact set $M(X)$. Thus, by Schauder's Fixed Point Theorem
(cf. Theorem V.10.5 in \cite{DS}) the map $S$ has a fixed point. In other words,  there exists
$\mu\in M(X)$ such that $\mu\circ T^{-1}=\mu$. 
\epf

\sp \noindent{\sl An alternative Proof.}
Suppose that $\mu_0$ is an arbitrary Borel probability measure on $X$ (for example, a Dirac point mass supported at a point of $X$, defined in Example~\ref{Dirac} below). Construct the sequence of Borel probability measures $(\mu_n)_{n=1}^\infty$,
where
\[
\mu_n=\frac{1}{n}\sum_{j=0}^{n-1}\mu_0\circ T^{-j}.
\]
The set of Borel probability measures on a compact metric space $X$ being weak*
compact, the sequence $(\mu_n)_{n=1}^\infty$ has at least one weak* limit point. Let $\mu_\infty$
be such a limit point. We claim that the measure $\mu_\infty$ is $T$-invariant. To show this,
let $(\mu_{n_k})_{k=1}^\infty$ be a subsequence
of the sequence $(\mu_n)_{n=1}^\infty$ which
weak* converges to $\mu_\infty$.
The weak* convergence of the subsequence means that for any function $f\in C(X)$ we have that
\[
\int_X f\,d\mu_\infty=\lim_{k\to\infty}\int_X f\,d\mu_{n_k}.
\]
On the other hand,
$$
\begin{aligned}
\left|\int_X f\circ T\,d\mu_{n_{k}}-\int_X f\,d\mu_{n_k}\right|
&=\left|\frac{1}{n_k}\sum_{j=0}^{n_k-1}\int_X f\circ T\,d\left(\mu_0\circ T^{-j}\right)-
         \frac{1}{n_k}\sum_{j=0}^{n_k-1}\int_X f\,d\left(\mu_0\circ T^{-j}\right)\right|\\
&=\frac{1}{n_k}\left|\sum_{j=0}^{n_k-1}\left(\int_X f\circ T^{j+1}\,d\mu_0-
                                              \int_X f\circ T^j\,d\mu_0
                                   \right)\right|\\
&=\frac{1}{n_k}\left|\int_X f\circ T^{n_k}\,d\mu_0-
                      \int_X f\,d\mu_0\right|\\
&\leq\frac{2}{n_k}\|f\|_\infty.
\end{aligned}
$$
Recall that, since $X$ is a compact metrizable space, the function $f\in C(X)$ is necessarily bounded. Therefore, passing to the limit as $k$ tends to infinity, we obtain that
\begin{eqnarray*}
\left|\int_X f\circ T\,d\mu_\infty-\int_X f\,d\mu_\infty\right|
&=&\left|\lim_{k\to\infty}\int_X f\circ T\,d\mu_{n_{k}}-\lim_{k\to\infty}\int_X f\,d\mu_{n_k}\right|\\
&=&\lim_{k\to\infty}\left|\int_X f\circ T\,d\mu_{n_{k}}-\int_X f\,d\mu_{n_k}\right|\\
&\leq&2\|f\|_\infty\lim_{k\to\infty}\frac{1}{n_k}\\
&=&0.
\end{eqnarray*}
Thus, 
$$
\int_X f\,d(\mu_\infty\circ T^{-1})=\int_X f\circ T\,d\mu_\infty
=\int_X f\,d\mu_\infty
$$ 
for all $f\in C(X)$ and so, we obtain that 
$$
\mu_\infty\circ T^{-1}=\mu_\infty.
$$
\endpf

\sp We now invoke Krein--Milman's, a classical theorem of functional analysis, Theorem to deduce from Bogolubov--Krylov Theorem somewhat more; namely that every topological dynamical system admits a Borel probability invariant ergodic measure. 

Recall that the convex hull of a subset $F$, denoted by $co(F)$, of a vector space $V$ is the set convex subset of $V$ containing $F$. It is equal to the intersection of all convex subset of $V$ containing $F$ and coincides with the set of all convex combinations of vectors from $F$. 

The proof of the following theorem can be found as Theorem~V.8.4 in Dunford and Schwartz's book \cite{DS}.

\sp\bthm[Krein--Milman Theorem]\label{KMT}
\index{(N)}{Krein-Milman Theorem}
If $K$ is a compact subset of a locally convex topological vector 
space $V$ and $E(K)$ is the set of its all extreme points, then 
$$
\overline{\mbox{co}}(E(K))\spt K.
$$
Consequently, 
$$
\overline{\mbox{co}}(E(K))=\overline{\mbox{co}}(K).
$$
In particular, if $K$ is convex then 
$$
\overline{\mbox{co}}(E(K))=K.
$$ 
\ethm

We can now prove the following.

\bthm\label{KBT-E}
Any topological dynamical system, i.e. any continuous transformation $T:X\lra X$ of a compact metrizable space $X$ admits at least one $T$--invariant Borel probability ergodic measure on $X$.
\ethm

\bpf
It follows from Theorem~\ref{KBT} and Theorem~\ref{KMT} that the set $E(T,\cB)$ of all ergodic $T$--invariant Borel probability measures on $X$ is non--empty.
\epf

\sp

\section{Examples of Invariant and Ergodic Measures}\label{EoIaEM} 

In order to prove that a given measure $\mu$ is invariant for some measurable dynamical system $T$, it is frequently useful to be able to reduce checking the equation $mu(T^{_1}(A))=\mu(A)$ for some smaller and convenient class of sets. One, especially important class of such sets is given by the concept of $\pi$-systems. We define it now and prove its fundamental, and extremely useful, property below.

\bdfn
Let $X$ be a set. A non-empty family $\mathcal{P}\sbt\mathcal{P}(X)$ is a
\index{(N)}{$\pi$-system}{\em $\pi$-system} on $X$ 
if and only if 
$$
P_1\cap P_2\in\mathcal{P}
$$ 
for all $P_1, P_2\in\mathcal{P}$. 
In other words, a $\pi$--system is a collection that is closed under finite intersections.
\edfn

The following proposition is classical, see for ex. Theorem~3.3 in \cite{Bi} for its proof.

\bprop\label{equalpi}
Let $X$ be a set and $\mathcal{P}$ be a $\pi$--system on $X$. 
Let $\mu$ and $\nu$ be probability measures on
$(X,\sg(\mathcal{P}))$. Then
\[
\mu=\nu
\ \ \ \Longleftrightarrow\ \ \ 
\mu(P)=\nu(P),\ \forall\, P\in\mathcal{P}.
\]
\eprop

\fr As its immediate consequence we get the following. 

\bprop\label{pressemi}
Let $T:(X,\F,\mu)\lra(X,\F,\mu)$ be a measurable endomorphism of a probability space $(X,\F,\mu)$. If a $\pi$--system $\mathcal{P}$ generates the $\sg$--algebra $\F$, then 
\[
\mu\circ T^{-1}=\mu 
\ \ \ \Longleftrightarrow\ \ \ 
\mu\circ T^{-1}(P)=\mu(P),\ \forall\, P\in\mathcal{P}.
\]
\eprop

We now may deal with actual examples.

\bex\label{Dirac}{\rm
Let $X$ be a set and let $T:X\lra X$ be a map. Fix an arbitrary point $\xi\in X$. The Dirac point mass supported at $\xi$, called also $\d$--Dirac measure \index{(N)}{$\d$--Dirac measure} supported at $\xi$. It is formally defined by the following formula:
$$
\delta_{\xi}(A):=
\begin{cases}
1 \  \  &{\rm if } \ \ \xi\in A \\
0 \  \  &{\rm if } \ \ \xi\notin A
\end{cases} 
$$
for every $A\sbt X$.

Now assume that $\xi$ is a fixed point point of $T$
, i.e. $T(\xi)=\xi$. Then the measure $\delta_{\xi}$ is $T$--invariant, that is, 
$$
\delta_{\xi}(T^{-1}(A))=\delta_{\xi}(A)
$$ 
for each $A\in\mathcal{A}$, since $\xi\in T^{-1}(A)$ if and only if $\xi\in A$. 
This example easily generalizes to invariant measures supported on periodic orbits. Namely, if $T^k(\xi)=\xi$ with some $k\ge 1$, then the the measure 
$$
\frac1k\sum_{j=0}^{k-1}\d_{T^j}(\xi)
$$
is $T$--invariant.} It is also obvious that this measure is ergodic.
\eex

\bex\label{Talphameaspres}{\rm 
Let $S^1:=[0,1]/(0=1)$.
Let $\a\in\R$ and define
the map $T_\a:S^1\lra S^1$ by 
\[
T_\a(x)=x+\a\hspace{0.1cm}(\mathrm{mod}1).
\]
Thus $T_\a$ is the rotation of the unit circle by the angle $2\a$.
The topological dynamics of $T_\a$ are radically different depending on whether the number $\a$
is rational or irrational. So will be the ergodicity of $T_\a$ with respect 
to the Lebesgue measure $\l$. However, it is fairly easy 
to see that $T_\a$ preserves Lebesgue measure $\l$, irrespective of the nature of $\a$. Indeed, 
$$
T_\a^{-1}(x)=x-\a\hspace{0.1cm}(\mathrm{mod} 1),
$$
and so $|\mathrm{det}\,DT_\a^{-1}(x)|=1$ for all $x\in S^1$. Therefore
\[
\l(T_\a^{-1}(B))=\int_B|\mathrm{det}\,DT_\a^{-1}(x)|\,d\l(x)
=\int_B\1 d\l=\l(B)
\]
for all $B\in\mathcal{B}(S^1)$, 
i.e. $T_\a$  preserves $\l$. 
\sp We shall prove the following.

\bthm\label{irr_rot_erg}
If $T_\a:S^1\lra S^1$ is the map of the circle defined defined as above, i.e.
$$
T_\a(x:=x+\a (\mathrm{mod}\ 1),
$$
then $T_\a$ is ergodic with respect to the Lebesgue measure $\l$ if and
only if $\a\in\R\backslash\Q$. 
\ethm

\bpf
First, assume that
$\a\notin\Q$. We want to show that $\l$ is ergodic with respect to $T_\a$.
For this, it suffice show that if 
$$
f\circ T_\a=f
$$
and $f\in L^2(\l)$, then $f$ is $\l$--a.e. constant. Consider the Fourier
series representation of $f$, which is given by
\[
f(x)
=\sum_{k\in\Z}a_k e^{2\pi ikx}.
\]
Then
\[
f\circ T_\a(x)
=\sum_{k\in\Z}a_k e^{2\pi ik(x+\a)}
=\sum_{k\in\Z}a_k e^{2\pi ik\a}e^{2\pi ikx}.
\]
Since we assumed that $f\circ T_\a=f$, we deduce from the uniqueness
of the Fourier series representation that 
$$
a_k e^{2\pi ik\a}=a_k
$$
for all $k\in\Z$. Hence, for each $k$ we have either 
$$
a_k=0
\  \  \  {\rm or} \  \  \ 
e^{2\pi ik\a}=1.
$$
The latter equality holds if and only if $k\a\in\Z$. Since $\a\notin\Q$, this 
occurs only when $k=0$. Thus $f(x)=a_0$ for $\lambda$--a.e. $x\in S^1$, that is, $f$ is $\l$--a.e. constant. This implies that $\l$ is ergodic.

The converse implication is obvious.
\epf
}
\eex

\bex\label{xing}{\rm
Fix $q\in\N$ and consider the map $T_q:S^1\lra S^1$ defined by the formula
$$
T_q(x):=qx\hspace{0.1cm}(\mathrm{mod}1).
$$
We claim that $\l$ is $T$--invariant. Let $I$ be a proper subinterval of 
$\S^1$. Then $T_q^{-1}(I)$ consists of $q$ mutually disjoint intervals (arcs) of length $\frac{1}{q}\l(I)$. Consequently,
\[
\l(T_q^{-1}(I))=q\cdot\frac{1}{q}\l(I)=\l(I).
\]
Since the family of all proper subintervals of $S^1$ forms a $\pi$--system which generates the Borel $sg$--algebra $\mathcal{B}(\S^1)$ and since $T_q$ preserves the Lebesgue measure of all
proper subintervals, Lemma~\ref{pressemi} asserts that $T_q$ preserves $\l$.

\sp We shall prove the following.

\bthm\label{t120191122} For every integer $q\ge 2$ the map $T_q:S^1\lra S^1$ is ergodic with respect to the Lebesgue measure $\l$.
\ethm

\bpf
It is possible to demonstrate the ergodicity of $T_q:S^1\lra S^1$ with respect to $\l$ in a similar way that we did for $T_\a$ in Theorem~\ref{irr_rot_erg}. However, in this example we will provide a different proof  which is more geometric and whose idea is applicable to many more classes of examples, particularly to full Markov maps dealt with in the next example.

Let $A\in\mathcal{B}(S^1)$ be a Borel set such that
$$
T_q^{-1}(A)=A
\  \  \  {\rm and} \  \  \ \lambda(A)>0.
$$
The surjectivity of $T_q$ then implies that $T_n(A)=A$. 
In order to establish ergodicity, we need to show that $\lambda(A)=1$.
By Lebesgue's Density Theorem, Theorem~\ref{th:6.5.4}, we have for $\lambda$--almost every point of $x\in A$ that
\[
\lim_{r\to 0}\frac{\lambda\bigl(A\cap B(x,r)\bigr)}{\lambda(B(x,r))}=1.
\]
Since $\lambda(A)>0$, there is at least one which we denote by $x$. Set 
\[
r_k:=1/(2q^k).
\] 
Then $T_q^k$ is injective on each open arc of length 
$2r_k$. So, the map $T_q^k|_{B(x,r_k)}$ is injective. Also,
\[
T_q^k(B(x,r_k))=S^1\backslash\{T_q^k(x+r_k)\}.
\]
Thus,
\[
\l\bigl(T_n^k(B(x,r_k))\bigr)=1.
\]
Therefore
\begin{eqnarray*}
\lambda(A)
&=&\lambda(T_q^k(A))
\geq\frac{\lambda\bigl(T_q^k(A\cap B(x,r_k))\bigr)}{\lambda\bigl(T_q^k(B(x,r_k))\bigr)}
=\frac{q^k\lambda\bigl(A\cap B(x,r_k)\bigr)}{q^k\lambda\bigl(B(x,r_k)\bigr)} \\ 
&=&\frac{\lambda\bigl(A\cap B(x,r_k)\bigr)}{q\lambda\bigl(B(x,r_k)\bigr)} 
\xrightarrow[\ \, k\to\infty \ ]{} 
1
\end{eqnarray*}
Consequently, $\lambda(A)=1$. This proves the ergodicity of $\l$.
\epf
}
\eex

\bex\label{Markovmapsvar}
{\rm The tent map $T:[0,1]\to[0,1]$ is defined by the formula:
\[
  T(x):=\left\{
    \begin{array}{lcl}
      2x    \    &{\rm  if }  \  \  x\in[0,1/2] \\ 
      2-2x  \   &{\rm  if } & x\in[1/2,1].
    \end{array}
  \right.
 \]
The family of all intervals $\{[a,b),(a,b):0<a<b<1\}$ forms a $\pi$-system that generates the Borel $\sg$--algebra $\mathcal{B}([0,1])$. Since the preimage of any such interval consists of two disjoint subintervals (one on each side of the tent) of half the length of the original interval, one readily sees by applying Proposition~\ref{pressemi} that Lebesgue measure 
on the interval $[0,1]$ is invariant under the tent map.
}
\eex

\bex{\rm 
In fact, the last two examples naturally generalize to a much larger family of maps. Let $T:[0,1]\lra[0,1]$ be a piecewise linear map of the unit interval 
that admits a ``partition'' 
$\mathcal{P}=\{p_j\}_{j=0}^q$, where $1\leq q<+\infty$ and 
$0=p_0<p_1<\ldots<p_{q-1}<p_q=1$, with the following properties:
\begin{itemize}
\item[(1)] 
$$
[0,1]=I_1\cup\ldots\cup I_q,
$$
where $I_j=[p_{j-1},p_j]$'s are the successive intervals of monotonicity of $T$.

\,

\item[(2)] 
$$
T(I_j)=[0,1]
$$ 
for all $1\leq j\leq q$.

\,

\item[(3)] The map $T$ is linear on $I_j$ for all $1\leq j\leq q$.
\end{itemize}

\fr Such a map $T$ will be called a {\em full Markov map} \index{(N)}{full Markov map}. 
We claim such a $T$ preserves the Lebesgue measure $\l$, i.e.
$$
\l\circ T^{-1}=\l
$$
Indeed, it is easy to see that the absolute value of the slope of 
the restriction of $T$ to the interval $I_j$ is $1/(p_j-p_{j-1})$. 
Therefore the absolute value of the slope of the 
corresponding inverse branch of $T$ is $p_j-p_{j-1}$. 
Let $J\sbt(0,1)$ be any interval. Then 
$$
T^{-1}(J)=\bigcup_{j=1}^q T|_{I_j}^{-1}(J),
$$
where
$T|_{I_j}^{-1}(J)$ is a subinterval of $\mathrm{Int}(I_j)$ 
of length $(p_j-p_{j-1})\cdot\l(I)$.
Since 
$$
\mathrm{Int}(I_j)\cap\mathrm{Int}(I_k)=\emptyset
$$
for all $1\leq j<k\leq q$, it follows that
\[
\l(T^{-1}(J))=\sum_{j=1}^q\l\bigl(T|_{I_j}^{-1}(J)\bigr)=\sum_{j=1}^q(p_j-p_{j-1})\cdot\l(I)=(p_q-p_0)\l(J)=\l(J).
\]
In addition, 
$$
0\leq\l(T^{-1}(\{0\}))\leq\l(\{p_j:0\leq j\leq k\})=0.
$$ 
So, 
$$
\l(T^{-1}(\{0\}))=0=\l(\{0\}).
$$
Similarly, 
$$
\l(T^{-1}(\{1\}))=0=\l(\{1\}).
$$
It follows that $\l(T^{-1}(J))=\l(J)$ for every interval $J\sbt[0,1]$.
Since the family of all intervals forms a $\pi$--system that generates 
$\mathcal{B}([0,1])$, Lemma~\ref{pressemi} asserts that the Lebesgue measure 
is invariant under any full Markov map. 

\sp We shall prove the following.

\bthm\label{t220191122} 
Every full Markov map $T:[0,1]\lra[0,1]$ is ergodic with respect to  Lebesgue measure $\l$ on $[0,1]$. 
\ethm

\bpf
As we did in the proof of Theorem~\ref{t120191122}, we would like to use Lebesgue's Density Theorem. However, now 
for all $r>0$ and all $k\in\N$, 
the restriction of $T^k$ to the ball $B(p_j,r)$ is not one--to--one when 
$p_j\ne 0,1$ is a point of continuity for a full Markov map $T$; for example, the point $1/2$ for the tent map. Nevertheless, despite the potential lack of injectivity, the proof presented below is motivated by the one above using Lebesgue's Density Theorem. 

For each $n\in\N$, let 
$$
\mathcal{P}_n:=\big\{I_j^{(n)}\,|\,1\leq j\leq q^n\big\}
$$
be the partition of the interval $[0,1]$ into 
the successive intervals of monotonicity of $T^n$. In particular, $I_j^{(1)}=I_j$ for all $1\leq j\leq q$.
For each $1\leq j<q^n$, let $p_j^{(n)}$ be the unique point in $I_j^{(n)}\cap I_{j+1}^{(n)}$.
For all $x\in[0,1]\backslash\{p_j^{(n)}:1\leq j<q^n\}$, let $I^{(n)}(x)$ be the 
unique element of $\mathcal{P}_n$ containing $x$. 
We will need two claims. Firstly:

\vspace{2mm}

{\em Claim $1^0$.} For every $n\in\N$, the map $T^n:[0,1]\lra [0,1]$ is a full Markov map under 
the partition $\mathcal{P}_n$. In addition, $\mathcal{P}_{n+1}$ is finer than $\mathcal{P}_n$. 

\vspace{2mm}

\bpf We will proceed by induction. Suppose that $T^n:[0,1]\to[0,1]$ is a full Markov map under the partition 
$\mathcal{P}_n$. It is obvious that $T^{n+1}$ is piecewise linear. 
Fix $I_j^{(n)}\in\mathcal{P}_n$. For all $1\leq i\leq q$, consider
\[
I_{j,i}^{(n+1)}:=T|_{I_i}^{-1}(I_j^{(n)}).
\] 
Define
\[
\widetilde{\mathcal{P}}_{n+1}
:=\bigl\{I_{j,i}^{(n+1)}:1\leq j\leq q^n \ {\rm and} \ 1\leq i\leq q\bigr\}.
\]
Then
\begin{eqnarray*}
\bigcup_{j=1}^{q^n}\bigcup_{i=1}^q I_{j,i}^{(n+1)} 
&=&\bigcup_{j=1}^{q^n}\bigcup_{i=1}^q T|_{I_i}^{-1}(I_j^{(n)})
 =\bigcup_{i=1}^q\bigcup_{j=1}^{q^n} T|_{I_i}^{-1}(I_j^{(n)}) \\
&=&\bigcup_{i=1}^q T|_{I_i}^{-1}\Bigl(\bigcup_{j=1}^{q^n} I_j^{(n)}\Bigr)
 =\bigcup_{i=1}^q T|_{I_i}^{-1}([0,1]) \\ 
&=&\bigcup_{i=1}^q I_i
 =[0,1].
\end{eqnarray*}
Thus $\widetilde{\mathcal{P}}_{n+1}$ is a cover of $[0,1]$.
Obviously, the interiors of the intervals in $\widetilde{\mathcal{P}}_{n+1}$
are mutually disjoint and $\widetilde{\mathcal{P}}_{n+1}$ is finer than $\mathcal{P}_n$. 
For all $1\leq j\leq q^n$ and all $1\leq i\leq q$ we further have
\[
T^{n+1}(I_{j,i}^{(n+1)})
=T^n\bigl(T(I_{j,i}^{(n+1)})\bigr)
=T^n\bigl(T(T|_{I_i}^{-1}(I_j^{(n)}))\bigr)
=T^n(I_j^{(n)})
=[0,1].
\]
So $\mathcal{P}_{n+1}=\widetilde{\mathcal{P}}_{n+1}$ and $T^{n+1}$ is a full Markov map under the partition $\mathcal{P}_{n+1}$, which is finer than 
$\mathcal{P}_n$. This completes the inductive step. Since the claim clearly holds when  $n=1$, Claim~$1^0$ has been established for all $n\in\N$.
\epf

\vspace{2mm}

{\em Claim $2^0$.} If $A\in\mathcal{B}([0,1])$, then 
\[
\lim_{n\to\infty}\frac{\l\bigl(A\cap I^{(n)}(x)\bigr)}{\l\bigl(I^{(n)}(x)\bigr)}=1 \  \mbox{ for } \  \l\mbox{-a.e. } x\in A. 
\]  

\vspace{2mm}

\bpf 
Let
\[
m:=\min\bigl\{|\text{slope}(T|_{I_i})|:1\leq i\leq q\bigr\}>1.
\]
Then,
\[
\mathrm{diam}(\mathcal{P}_n)
:=\sup\bigl\{\mathrm{diam}(I_j^{(n)}):1\leq j\leq q^n\bigr\}
\leq m^{-n}\xrightarrow[\ \, n\to\infty \ ]{}  0.
\]
Therefore, the $\sg$--algebra $\sg\bigl(\bigcup_{n=1}^\infty\mathcal{P}_n\bigr)$ contains all open, relative to $[0,1]$ subintervals of $[0,1]$. Hence 
$$
\sg\bigl(\bigcup_{n=1}^\infty\mathcal{P}_n\bigr)=\mathcal{B}([0,1]).
$$
Claim~$2^0$ now directly follows from this, the last assertion of Claim~$1^0$, and Theorem~\ref{t120191123}.
\epf 

\sp Concluding the proof of Theorem~\ref{t220191122}, let $A$ be a Borel subset of $[0,1]$ such that $T^{-1}(A)=A$ and $\l(A)>0$. 
By the surjectivity of $T$, we know that $T(A)=A$. Fix any $x\in A$ satisfying  Claim~$2^0$. For each $n\in\N$, let $m_n$ be the slope of $T^n|_{I^{(n)}(x)}$. Using both claims, 
we obtain that
\begin{eqnarray*}
\lambda(A)
=\lambda(T^n(A))
\geq\frac{\lambda\bigl(T^n(A\cap I^{(n)}(x))\bigr)}{\lambda\bigl(T^n(I^{(n)}(x))\bigr)}
=\frac{m_n\,\lambda\bigl(A\cap I^{(n)}(x)\bigr)}{m_n\,\lambda\bigl(I^{(n)}(x)\bigr)} 
=\frac{\lambda\bigl(A\cap I^{(n)}(x)\bigr)}{\lambda\bigl(I^{(n)}(x)\bigr)} 
\xrightarrow[\ \, n\to\infty \ ]{}1.
\end{eqnarray*}
Consequently, $\lambda(A)=1$. This proves the ergodicity of $\l$.
\epf
}
\eex

\sp 

\bex{\rm
Let $T:(X,\mathcal{A})\to(X,\mathcal{A})$ and $S:(Y,\mathcal{B})\to(Y,\mathcal{B})$
be measurable transformations for which there exists a measurable
transformation $h:(X,\mathcal{A})\to(Y,\mathcal{B})$ such that $h\circ T=S\circ h$.
We will show that every $T$-invariant measure generates an $S$-invariant push down under $h$.
Let $\mu$ be a $T$-invariant measure on $(X,\mathcal{A})$. 
Recall that the push down of $\mu$ under $h$ is the measure
$\mu\circ h^{-1}$ on $(Y,\mathcal{B})$. It follows from the $T$-invariance of $\mu$ that
\[
(\mu\circ h^{-1})\circ S^{-1}
=\mu\circ(S\circ h)^{-1}
=\mu\circ(h\circ T)^{-1}
=(\mu\circ T^{-1})\circ h^{-1}
=\mu\circ h^{-1}.
\]
That is, $\mu\circ h^{-1}$ is $S$-invariant.}
\eex

\bex\label{exprodmeas}{\rm 
Let $(X,\mathcal{A},\mu)$ and $(Y,\mathcal{B},\nu)$ be two probability spaces, and let $T:X\lra X$ and $S:Y\lra Y$ be two measure--preserving dynamical systems. The Cartesian product of $T$ and $S$ is the map $T\times S:X\times Y\lra X\times Y$ defined by the formula:
\[
(T\times S)(x,y)=(T(x),S(y)).
\] 

The direct product $\sg$-algebra $\sg(\mathcal{A}\times\mathcal{B})$ on $X\times Y$
is the $\sg$--algebra generated by the $\pi$--system of measurable rectangles
\[
\mathcal{A}\times\mathcal{B}:=\{A\times B:A\in\mathcal{A},\, B\in\mathcal{B}\}.
\]
The direct \index{(N)}{product measure} product measure, commonly defined in measure theory,
$\mu\otimes\nu$ on $(X\times Y,\sg(\mathcal{A}\times\mathcal{B}))$
is the measure uniquely determined by its values on its generating $\pi$--system:
\[
(\mu\otimes\nu)(A\times B):=\mu(A)\nu(B).
\]
We claim that the product map 
$$
T\times S:(X\times Y,\sg(\mathcal{A}\times\mathcal{B}),\mu\otimes\nu)
\to(X\times Y,\sg(\mathcal{A}\times\mathcal{B}),\mu\times\nu)
$$ 
is measure--preserving, i.e.
$$
(\mu\otimes\nu)\circ(T\times S)^{-1}(F)=(\mu\otimes\nu)(F)
$$
for every set $F\in \sg(\mathcal{A}\times\mathcal{B})$.
 
Thanks to Proposition~\ref{pressemi}
it suffices to show that
$$
(\mu\otimes\nu)\circ(T\times S)^{-1}(A\times B
)=(\mu\otimes\nu)(A\times B)
$$ 
for all $A\times B\in\mathcal{A}\times\mathcal{B}$. And indeed,
\begin{eqnarray*}
(\mu\otimes\nu)\circ(T\times S)^{-1}(A\times B)
&=&(\mu\otimes\nu)\bigl(T^{-1}(A)\times S^{-1}(B)\bigr)
=\mu(T^{-1}(A))\,\nu(S^{-1}(B))\\
&=&\mu(A)\nu(B)
=(\mu\otimes\nu)(A\times B).
\end{eqnarray*}
}
\eex  

\sp Our final example pertains to the shift map introduced by formula \eqref{120191121} below. In the special case when the set $E$ is countable, it will be treated at length in Section~\ref{SOFToIA;TP}. 
In fact, we now will look at this map in a more general context.

\bex\label{onesided}{\rm 
Let $(E,\F,P)$ be a probability space. Consider the one--sided product set
$$
E^\N:=\prod_{k=1}^\infty E.
$$
The product $\sg$--algebra $\F^\infty$ on $E^\infty$ is the $\sg$--algebra generated by the $\pi$--system $\mathcal{S}$ of all (finite) cylinders (also called rectangles), i.e. the sets of the form
\[
\prod_{k=1}^n E_k\times\prod_{l=n+1}^\infty E
=\Bigl\{\om=(\om_j)_{j=1}^\infty\in E^N:\om_k\in E_k\,\, \forall\, 1\leq k\leq n\Bigr\}, 
\]
where $n\in\N$ and $E_k\in\F$ for all $1\leq k\leq n$. 
Commonly in measure theory, the product measure $\mu_P$ on $\mathcal{F}_\infty$, frequently referred to as the Bernoulli measure generated by $P$, is the unique probability measure that gives to a cylinder the value
\begin{equation}\label{prodmeassemi}
\mu_P\Bigl(\prod_{k=1}^n E_k\times\prod_{l=n+1}^\infty E\Bigr)
=\prod_{k=1}^n P(E_k).
\end{equation}

Let $\sg:E^\N\lra E^\N$ be the left shift map, which is defined by the formula
\beq\label{120191121}
\sg\((\omega_n)_{n=1}^\infty\):=(\omega_{n+1})_{n=1}^\infty.
\eeq
As we have said above, in the special case when the set $E$ is countable, it will be treated at length in Section~\ref{SOFToIA;TP}. We shall prove that 

\centerline{The product measure $\mu_P$ is $\sg$--invariant.}

\sp Indeed, since the cylinder sets form a $\pi$--system
which generates the product $\sg$--algebra $\F^\infty$, in light of Lemma~\ref{pressemi} 
it is sufficient to show that 
$$
\mu_P(\sg^{-1}(S))=\mu_P(S)
$$ 
for all cylinder sets $S\in\mathcal{S}$. And we have
\begin{eqnarray*}
\mu_P\circ\sg^{-1}\Bigl(\prod_{k=1}^n E_k\times\prod_{l=n+1}^\infty E\Bigr)
&=&\mu_P\Bigl(E\times\prod_{k=1}^n E_k\times\prod_{l=n+2}^\infty E\Bigr)
=P(E)\prod_{k=1}^n P(E_k)\\
&=&\prod_{k=1}^n P(E_k)
=\mu_P\Bigl(\prod_{k=1}^n E_k\times\prod_{l=n+1}^\infty E\Bigr).
\end{eqnarray*}
This completes the proof that the product measure is shift--invariant.

\sp We shall prove the following.

\bthm\label{t120191125}
If $(E,\F,P)$ is a probability space, then the shift map $\sg:E^\N\lra E^\N$ is ergodic with respect to the Bernoulli probability measure shift--invariant measure $\mu_P$ generated by $P$.
\ethm

\bpf
Given an integer $k\ge 1$ let $\F_k'$ be the $\sg$--algebra on $E^k$ generated by all cylinders of length $k$, i.e. by all sets of the form
$$
\prod_{j=1}^k E_j,
$$
where $E_j\in\F$ for all $1\leq j\leq k$. Furthermore, let $\F_k$ be the sub--$\sg$--algebra of of $\F^\infty$ generated by all sets of the form
\[
F\times\prod_{l=k+1}^\infty E, 
\]
where $F\in\F'$. 
Let $A\in\F^\infty$ be a set such that
$$
\sg^{-1}(A)=A
\  \  \  {\rm and} \  \  \ \mu_P(A)>0.
$$
Our task is to show that $\mu_P(A)=1$. By the Martingale Convergence Theorem for Conditional Expectations, i.e. by Theorem~\ref{DMCT}, we have that
$$
E_{\mu_P}(\1_{A}|{\mathfrak{F}_k})(x)\xrightarrow[\ \, k\to\infty \ ]{}  
 \1_A(x)
$$ 
for $\mu_P$--a.e. $x\in E^{\mathbb N}$. For every $k \geq 1$, let  
$$
A_k:=\lt\{x \in E^{\mathbb N}:\big|E_{\mu_P}(\1_A|{\mathfrak{F}_k})(x)-1\big|\leq \frac{1}{4}\rt\} \in \mathfrak{F}_k.
$$
So,
$$ 
\bal
\forall \, \varepsilon\in(0, 1/8)\, \,\,  \exists  N_{\varepsilon}\, \,\,   &\forall k \geq N_{\varepsilon}\, \, \,  \\
&\mu_P\lt(\lt\{x \in E^{\mathbb N}:\, \, \big|E_{\mu_P}(\1_A|{\mathfrak{F}_k})(x)-\1_{A}(x)\big|\geq \frac{1}{8}\rt\}\rt) \leq \frac{\varepsilon}{2}.
\eal
$$
Denote the latter set by $B_k$. If $x \in A\sms  A_k$, then
$$
\big|E_{\mu_P}(\1_A|{\mathfrak{F}_k})(x)-\1_{A}(x)\big|>\frac{1}{4}.
$$
Hence, $ x \in B_k$. This means that $ A\sms A_k \sbt B_k$. Therefore,
\begin{equation}\label{1v3}
\mu_P(A\sms A_k)\leq\mu_P(B_k)\leq \frac{\varepsilon}{2}.
\end{equation}
If, on the other hand $x \in A_k\sms A$, then
$$ 
\big|E_{\mu_P}(\1_A|{\mathfrak{F}_k})(x)-\1_{A}(x)\big|
=\big|E_{\mu_P}(\1_A|{\mathfrak{F}_k}(x))\big|\geq \frac{3}{4} > \frac{1}{8}.
$$
Hence, $x \in B_k$. This mean that $A_k\sms A\sbt B_k$. Therefore, 
$$
\mu_P(A_k\sms A) \leq \mu_P(B_k)\leq \frac{\varepsilon}{2}.
$$
Together with (\ref{1v3}) this gives that
\begin{equation}\label{2v3}
\mu_P(A_k \bigtriangleup A)\leq \varepsilon.
\end{equation}
Now, if $ n \geq k+1$, then, abusing slightly notation, we  have $ \sigma^{-n}(A_k)=E^{n}\times A_k \in \mathfrak{F}_n$. So,
\begin{equation}\label{3v3}
\mu_P( A_k\cap \sigma^{-n}(A_k))=\mu_P(A_k)\cdot \mu_P(\sigma^{-n}(A_k))=\mu_P^2(A_k).
\end{equation}
Hence,
$$\begin{aligned}
|\mu_{P}^2(A)&-\mu_{P}(A)| \leq |\mu_{P}^2(A)-\mu_{P}^2(A_k)|+|\mu_{P}^2(A_k)-\mu_{P}(A)|
\\
 & \leq (\mu_{P}(A)+\mu_{P}(A_k))|\mu_{P}(A)-\mu_{P}(A_k)|+\big|\mu_{P}\(A_k \cap \sigma^{-n}(A)\)-\mu(A)\big|
 \\
& \leq 2 \mu_{P}(A_k \vartriangle A)+\mu_{P}(A\vartriangle (A_k\cap \sigma^{-n}(A_k)))\\
& \leq 2 \mu_{P}(A_k \vartriangle A)+\mu_{P}(A\sms  A_k)+ \mu_{P}(A\sms \sigma^{-n}(A)) +\mu_{P}\((A_k\cap \sigma^{-n}(A_k))\sms A\) 
\\
&=2 \mu_{P}(A_k \vartriangle A)+\mu_{P}(A\sms  A_k)+\mu_{P}(\sigma^{-n}(A)\sms \sigma^{-n}(A_k)) +\mu_{P}\((A_k\cap \sigma^{-n}(A_k)\sms A\)\\
&= 2 \mu_{P}(A_k\vartriangle A)+\mu_{P}(A\sms  A_k)+\mu_{P}\(\sigma^{-n}(A)\sms A_k)\) +\mu_{P}\((A_k\cap \sigma^{-n}(A_k))\sms A\)\\
&= 2 \mu_{P}(A_k \vartriangle A)+\mu_{P}(A\sms  A_k)+\mu_{P}(A\sms A_k) +\mu_{P}\((A_k\cap \sigma^{-n}(A_k)\)\sms A)\\
& \leq 4\mu_{P}(A_k \vartriangle A)+\mu_{P}(A_k\sms  A)\\
& \leq 5 \mu_{P}(A_k \vartriangle A)=5 \varepsilon.
\end{aligned}
$$
So, letting $ \varepsilon\searrow 0$, we get $\mu_P^2(A)=\mu_P(A)$. As $\mu_{\underline{p}}(A)>0$, this yields $ \mu_P(A)=1$. The proof is complete. 
\epf
}
\eex

Out of Example~\ref{xing} and/or Example~\ref{onesided} we easily get the classical result of Borel that almost all Lebesgue numbers in the interval $[0,1]$ are normal. This is done in the following example.

\bex\label{e120191125}{\rm
Fix an integer $q\ge 2$, set
$$
E_q:=\{0,1,2,\ld,q-1\}
$$
and consider the map $\pi:E_q^\N\lra [0,1]$, called the coding map, defined by the formula
$$
\pi(\om):=\sum_{k=1}^\infty\frac{\om_k}{q^k}.
$$
Since for each integer $k\ge 1$, the function, $p_k:E_q^\N\lra [0,1]$ defined by the formula
$$
p_k(\om):=\frac{\om_k}{q^k}
$$
is continuous, since
$$
\pi=\sum_{k=1}^\infty p_k,
$$
and since this series is uniformly convergent, we conclude that the function 
$\pi:E_q^\N\lra [0,1]$ is continuous, thus Borel measurable. Let $Z$ be the (countable) set of the endpoints of all elements of all partitions $\mathcal P_n$, $n\ge 1$, coming from Example~\ref{xing}. Then $E^\N\sms \pi^{-1}(Z)$ consists of all sequences in $E^\N$ that have no tail consisting either of $0$s only or of $(q-1)$s only, and for each $x\in [0,1]\sms Z$ the set $\pi^{-1}(x)$ is a singleton. In particular, the map
$$
\pi|_{E^\N\sms \pi^{-1}(Z)}:E^\N\sms \pi^{-1}(Z)\lra [0,1]
$$
is one--to--one. Let now $F\in \bigcup_{n=1}^\infty\mathcal{P}_n$ be arbitrary. Then there exists an integer $n\ge 1$ such that $F\in \mathcal{P}_n$, meaning that
$$
F=\lt[\frac{i}{q^n},\frac{i+1}{q^n}\rt]
$$
for some integer $i\in \{0,1,2,\ld,q^n-1\}$. Then, up to a countable set, 
$$
\pi^{-1}(F)=[\om],
$$
where $\om\in E^n$ is such that $\pi(\om0^\infty)=i/q^n$. Hence
$$
\l(F)=\frac{1}{q^n}=\mu_{P_q}([\om])=\mu_{P_q}\circ\pi^{-1}(F),
$$
where $P_q$ is the probability measure on $E_q$ given by the formula
$$
P_q(B):=\frac{\#B}{\#E_q}=\frac{\#B}{q}.
$$
Since $\bigcup_{n=1}^\infty\mathcal{P}_n$ is a $\pi$--system generating $\mathcal{B}([0,1])$, we conclude from Proposition~\ref{equalpi} that 
\beq\label{120191125}
\mu_{P_q}\circ\pi^{-1}=\l.
\eeq
Now fix $i\in \{0,1,2,\ld,q-1\}$. Then, Theorem~\ref{t1_mu_2014_11-18} 
affirms that
$$
\lim_{n\to\infty}\frac{1}{n}\#\bigl\{1\leq j\leq n:\om_j=i\bigr\}
=\lim_{n\to\infty}\frac{1}{n}\#\bigl\{0\leq j\le n-1:\sg^j(\om)\in[i]\bigr\} =\mu_{P_q}([i]) 
= \frac{1}{q} 
$$
for $\mu_{P_q}$--a.e. $\om\in E_q^\N$. Since the coding map is measure--preserving, we deduce that $\l$--almost every number 
between $0$ and $1$ has a $q$--adic expansion whose digits are equal to $i$ with a asymptotic frequency of $1/q$. Such numbers are normal to base $q$. Since the intersection of any countable family of sets of measure $1$ on a probability space is a set of measure one, we thus have the following main result of this example.

\bthm\label{t120191125B}
Lebesgue almost every number in $[0,1]$ is normal with respect to every base $q\ge 2$. 
\ethm

\sp\fr Of course this theorem has obvious refinements such as the appropriate value of asymptotic frequencies of fixed digits placed at fixed (mod $\!\!q$) positions. The ultimate generality is Birkhoff's Ergodic Theorem, i.e. Theorem~\ref{Birkhoff}.  
}
\eex

\bex\label{ex7_gauss}{\rm 
The map $G:[0,1]\to[0,1]$ defined by  
  \[
  G(x):=\left\{
    \begin{array}{ccl}
      0                                                & \mbox{ if } & x=0 \\ 
      \displaystyle\lt\{\frac{1}{x}\rt\} & \mbox{ if } & x>0,
    \end{array}
  \right.
  \]
where $\{x\}$ denotes the fractional part of a non--negative number $x$, 
is called the {\em Gauss map}. \index{(N)}{Gauss map} This map does not preserve Lebesgue measure $\l$ on $[0,1]$. However, let $\mu_G$ be the Borel probability measure on $[0,1]$ defined by the formula
\[
\mu_G(B):=\frac{1}{\log 2}\int_B\frac{1}{1+x}\,dx
\]
for every Borel set $B\sbt [0,1]$. This means that $\mu_G(B)$ is uniquely determined by the property that
$$
\frac{\text{d}\mu_G}{\text{d}\l}(x)=\frac{1}{\log 2}\frac{1}{1+x}
$$
for all $x\in[0,1]$. we have the following. The measure $\mu_G$ is known as the {\em Gauss measure}. \index{(N)}{Gauss measure}

\sp\bthm\label{t520191125}
The Borel probability measure $\mu_G$ (equivalent to Lebesgue measure $\l$) on $[0,1]$ is $G$--invariant. 
\ethm

\bpf 
Let $[a,b]\sbt(0,1)$ be arbitrary. Then 
\[
G^{-1}([a,b])=\bigcup_{n=1}^\infty\Bigl[\frac{1}{b+n},\frac{1}{a+n}\Bigr].
\]
Hence,
\begin{eqnarray*}
\mu_G(G^{-1}([a,b]))
&=&\frac{1}{\log 2}\sum_{n=1}^\infty\int_{\frac{1}{b+n}}^{\frac{1}{a+n}}\frac{1}{1+x}\,dx \\
&=&\frac{1}{\log 2}\sum_{n=1}^\infty\Bigl[\log\Bigl(1+\frac{1}{a+n}\Bigr)-\log\Bigl(1+\frac{1}{b+n}\Bigr)\Bigr] \\
&=&\frac{1}{\log 2}\sum_{n=1}^\infty\Bigl[\log(a+n+1)-\log(a+n)-\log(b+n+1)+\log(b+n)\Bigr] \\
&=&\frac{1}{\log 2}\bigl[\log(b+1)-\log(a+1)\bigr] \\
&=&\frac{1}{\log 2}\int_{a}^{b}\frac{1}{1+x}\,dx \\
&=&\mu_G([a,b]).
\end{eqnarray*}
So, by applying Proposition~\ref{pressemi}, we conclude that the Gauss measure $\mu_G$ is $G$--invariant
under the Gauss map. 
\epf

\sp\fr The Gauss measure $\mu_G$ is $G$ also ergodic with respect to $G$. This can be proved directly by improving the reasoning from Theorem~\ref{t220191122}. The new difficulty is however now two--folded: although the map is ``full'' and ``Markov'', it is not linear and there are infinitely many pieces of monotonicity. Ergodicity also follows from Theorem ~\ref{t2.2.4}. In fact the Gauss map $G$, more precisely its inverse branches 
$$
[0,1]\ni x\lmt g_n(x):=\frac1{x+n}, \  \  n\in\N,
$$
form a conformal iterated function system in the sense of Chapter~\ref{Markov-systems}, Graph Directed Markov Systems, and all proved therein applies to the system $\{g_n\}_{n\in\N}$.
}
\eex

\sp\chapter{Probability (Finite) Invariant Measures; Finer Properties}\label{finite-measure}

\sp In this chapter we deal with  stochastic laws for measurable endomorphisms preserving a finite measure that are finer than mere Birkhoff's Ergodic
Theorem. Of course, in order to obtain them we need more hypotheses. Right up front we say that having such a measure, we can
always normalize it, i.e. multiply by a constant factor, to produce a
probability measure. So, from now on throughout the entire chapter any
invariant measure is understood to be a probability one. In the first
section of this chapter we deal with the Law of Iterated Logarithm
while, the second section is devoted to describe the method of L. S. Young
towers which she developed in \cite{lsy1} and \cite{lsy2}; see also \cite{gouezel3} for further progress. With appropriate assumptions imposed on the first return time her construction yields
the exponential decay of correlations, the Central Limit Theorem, and
the Law of Iterated Logarithm follows too.

\sp\section{The Law of Iterated Logarithm}

\sp We shall show that in this section that under mild conditions, if a
first return map satisfies the Law
of Iterated Logarithm, then so does the original map. More precisely, let 
$$
T:X\lra X
$$ 
be a measurable dynamical system preserving a probability measure $\mu$ on $X$. We say that a $\mu$--integrable function 
$$
g:X\lra\R
$$ 
satisfies the Law of
Iterated Logarithm if there exists a positive number $A_g$ such that
\begin{displaymath}
  \limsup_{n\to\infty}\frac{S_{n}g(x)-n\int gd\mu}{\sqrt{n\log\log n}}
  = A_g
\end{displaymath}
for $\mu$-a.e. $x\in X$.
From now on we assume without loss of generality that 
\begin{displaymath}
  \mu(g)=\int g d\mu=0.  
\end{displaymath}
Keep a  measurable set $A\sbt X$ with $\mu(A)>0$. Given a point $x\in
A$, the sequence 
$(\tau_n(x))_{n=1}^\infty$ is then defined as follows.
\begin{displaymath}
  \tau_1(x):=\tau_A(x)\quad\textrm{ and }\quad
  \tau_n(x)=\tau_{n-1}(x)+\tau(T^{\tau_{n-1}(x)}(x)),
\end{displaymath}
where, we recall, $\tau_A(x)\ge 1$ is the first return time of $x$ to $A$. What we are up to is the following theorem. Its proof is taken from
\cite{SkorulUPUZ} and \cite{SUZ}.

\sp

\begin{thm}\label{LILfrominducing}
Let $T:X\lra X$ be a measurable dynamical system preserving a
probability measure $\mu$ on $X$. Assume that the dynamical system
($T,\mu$) is ergodic. Fix $A$, a measurable subset of $X$ having a
positive measure $\mu$. Let $g:X\lra\R$ be a measurable function such
that the function $g_A:A\lra\R$ satisfies the Law of Iterated Logarithm
with respect to the dynamical system $(T_A,\mu_A)$. If in addition,
  \beq
    \label{eq:6}
    \int_A|g_A|^{2+\g} d\mu<+\infty
  \eeq
with  some $\g>0$, then the function $g:X\lra\R$ satisfies the Law
of Iterated Logarithm with respect to the original dynamical system
$(T,\mu)$ and $A_g=A_{g_A}$. 
\end{thm}

\bpf 
Since the Law of Iterated Logarithm holds for a point $x\in X$ if and
only if it holds for $T(x)$, in virtue of ergodicity of $T$, it
suffices to prove our theorem for almost all points in $A$.
By our assumptions there exists a positive constant $A_{\hat g}$ such that
\begin{displaymath}
  \limsup_{n\to\infty}\frac{S_{\tau_n}g(x)}{\sqrt{n\log\log n}}=
 A_{\hat g}.
\end{displaymath}
for $\mu_A$-a.e. $x\in A$.
Since, by Kac's Lemma,
\beq
  \label{eq:2}
  \lim_{n\to\infty}\frac{\tau_n}{n}=\int_X\tau dm=\int_A\tau dm=1,
\eeq
$\mu_A$-a.e. on $A$, we thus have
\beq
  \label{eq:1}
\limsup_{n\to\infty}\frac{S_{\tau_n}g}{\sqrt{n\log\log n}}=
\limsup_{n\to\infty}\frac{S_{\tau_n}g}{\sqrt{\tau_n\log\log\tau_n}}
  =A_{\hat g}.
\eeq
$\mu_A$-a.e. on $A$.
Now, for every $n\in\N$ and (almost) every $x\in A$ let $k=k(x,n)$ be the
positive integer uniquely determined by the condition that
  \begin{displaymath}
    \tau_k(x)\leq n < \tau_{k+1}(x).
  \end{displaymath}
Since
  \begin{displaymath}
    S_ng(x)=S_{\tau_k(x)}g(x)+S_{n-\tau_k(x)}g(T^{\tau_k(x)}(x)),
  \end{displaymath}
we have that
  \beq
    \label{eq:4}
    \frac{S_{n}g(x)}{\sqrt{n\log\log n}}
    =\frac{S_{\tau_k(x)}g(x)}{\sqrt{n\log\log n}}+
    \frac{S_{n-\tau_k(x)}g(x)}{\sqrt{n\log\log n}}
  \eeq
  Since by \eqref{eq:2}
  \begin{displaymath}
    \lim_{n\to\infty}\frac{\tau_{k+1}(x)}{\tau_k(x)}=1,
  \end{displaymath}
  we get from \eqref{eq:1} that,
  \begin{displaymath}
    \limsup_{n\to\infty}\frac{S_{\tau_k}g(x)}{\sqrt{n\log\log n}}=
    \limsup_{n\to\infty}\frac{S_{\tau_k}g(x)}{\sqrt{k\log\log k}}=
   A_{\hat g}.
  \end{displaymath}
  Because of this and because of \eqref{eq:4}, we are only left to
  show that 
  \beq
    \label{eq:5}
    \lim_{n\to\infty}\frac{S_{n-\tau_{k(n)}}g(x)}{\sqrt{n\log\log n}}=0.
  \eeq
$\mu_A$-a.e. on $A$. To do this, note first that
  \begin{displaymath}
    \frac{S_{\tau_{k+1}-\tau_k}|g|(T^{\tau_k}(x))}
    {\sqrt{k\log\log k}}=
    \frac{|\hat g|(T_A^k(x))}{\sqrt{k\log\log k}}.
  \end{displaymath}
Take an arbitrary $\varepsilon\in (0,\g)$. Since
\beq
\aligned
  \mu\Big(\Big\{x\in A:|& g_A|(T^k_A(x)) \geq \varepsilon\sqrt{k\log\log k}\Big\}\Big)=\\
&= \mu\Big(\lt\{x\in A:|g_A|(x)\geq \varepsilon\sqrt{k\log\log k}\rt\}\Big)\\
&= \mu\Big(\lt\{x\in A:|g_A|^{2+\varepsilon}(x)
      \geq \varepsilon^{2+\varepsilon}(k\log\log
      k)^{1+\varepsilon/2}\rt\}\Big) \\
&\leq \frac{\int |g_A|^{2+\varepsilon} d\mu}{
      \varepsilon^{2+\varepsilon}(k\log\log k)^{1+\varepsilon/2}},
\endaligned
\eeq
using \eqref{eq:6} we conclude that
  \begin{displaymath}
    \sum_{k=1}^\infty\mu\Big(\lt\{x\in A:|g_A|(x)\geq 
    \varepsilon\sqrt{k\log\log k}\rt\}\Big)<\infty.
  \end{displaymath}
So, applying Borel--Cantelli Lemma, formula \eqref{eq:5} follows. We are done.
\epf

\sp\section[Decay of Correlations and the Central Limit
  Theorems] {Decay of Correlations and the Central Limit
  Theorems; \nl Lai-Sang Young Towers}\label{Young-Abstract}

\sp Let $T:X\lra X$ be a measurable dynamical system preserving a
probability measure $\mu$ on $X$. We say that a $\mu$--integrable function $g:X\lra\R$ satisfies the Central Limit Theorem if there exists $\sigma>0$ such that
$$
\frac{\sum_{j=0}^{n-1}g\circ T^j-n\int gd\mu}{\sqrt{n}} 
\xrightarrow[\ n\to\infty \ ]{}\mathcal N(0,\sigma)
$$ 
in distribution. $\mathcal N(0,\sigma)$ is here the normal (Gaussian) distribution
with $0$ mean and variance $\sigma$. This precisely means that
for every $t\in \R$,
$$
\mu\biggl(\Big\{x\in X:\frac{\sum_{j=0}^{n-1}g\circ T^j(x)-n\int gd\mu}{\sqrt{n}}
\le t\Big\}\biggr)\xrightarrow[\ n\to\infty \ ]{} {1\over \sigma
  \sqrt{2\pi}}\int_{-\infty}^t\exp\(-u^2/2\sigma^2\)\, du. 
$$
Equivalently: for every Borel set $B\sbt\R$ with $S(\bd B)=0$ (we recall that $S$ denotes here Lebesgue measure on $\R$)
$$
\mu\biggl(\Big\{x\in X:\frac{\sum_{j=0}^{n-1}g\circ T^j(x)-n\int gd\mu}{\sqrt{n}}
\in B\Big\}\biggr)\xrightarrow[\ n\to\infty \ ]{} {1\over \sigma
  \sqrt{2\pi}}\int_B\exp\(-u^2/2\sigma^2\)\, du. 
$$
Another important stochastic feature of a dynamical system is the rate of decay
of correlations it yields. Let $\psi_1$ and $\psi_2$ be real square $\mu$-integrable
functions on $X$. For every positive integer $n$ the \it $n$-th
correlation \rm of the pair $\psi_1,\psi_2$, is the number
\beq\label{correlations}
C_n(\psi_1,\psi_2)
:=\int \psi_1\cdot(\psi_2\circ f^n)\,d\mu - \int
   \psi_1\,d\mu_\phi \int \psi_2\,d\mu
\eeq
 provided the above integrals exist. Notice that, due to the
$T$-invariance of $\mu$, we can also write
$$
C_n(\psi_1,\psi_2)=\int (\psi_1-E\psi_1)\((\psi_2-E\psi_2)\circ
T^n\)\, d\mu,
$$
where we put $E\psi=\int \psi\, d\mu$.

\

We shall now describe a pwerful tool, commonly refered to as Young's
tower technique, which provides a way to prove the Central Limit
Theorem and estimate from above the decay of correlations in many
sufficiently ``regular'' dynamical systems 
exhibiting some sufficient expanding or hyperbolic features. Let
\begin{itemize}
\item $(\Delta_0,\mathcal{B}_0,m)$ be a 
measure space with a finite measure $m$,

\, \item $\mathcal{P}_0$ be a countable measurable partition of $\De_0$,

\, \item $T_0:\Delta_0\lra\Delta_0$ be a measurable map such that for every 
  $\Ga\in\mathcal{P}_0$ the map 
$$
T_0:\Ga\to \Delta_0
$$
is injective and $T_0(\Ga)$ can be reperesented as a union of elements of $\mathcal{P}_0$. 

\, \item (Big Images Property)
$$
\inf\{m(T_0(\Ga)):\Ga\in \mathcal{P}_0\}>0.
$$
\end{itemize}
Furthermore,
\begin{itemize}
\item We assume that the partition $\mathcal{P}_0$ is generating, i.e. for every $x,y\in\Delta_0$ there exists $s\ge 0$ such that the points $T_0^s(x)$ and $T_0^s(y)$ are in different elements of the partition $\mathcal{P}_0$. 
  
\, \item We denote by $s=s(x,y)$ the smallest integer with this property and we call it a separation time for the pair $x,y$.

\, \item We assume also that for each $\Ga\in\mathcal{P}_0$ the map
$({T_0}\vert_{\Ga})^{-1}$ is measurable and that $Jac_{m}(T_0)$,
the Jacobian of $T_0$ with respect to the measure $m$, is well--defined and
positive a.e. in $\Ga$. 

\, \item (Bounded Distortion Property) With some constants $0<\b<1$ and $C>0$,
\beq\label{youngcond0}
\left |\frac{Jac_{m} T_0(x)}{Jac_{m} T_0(y)}-1\right |
\le C\beta^{s(x,y)}
\eeq
for all $\Ga\in\mathcal{P}_0$ and all $x,y\in \Ga$.

\, \item We have also a function $R:\Delta_0\lra\{1, 2,\ld\}$ (the first return time to $\De_0$) which is constant on each element of the partition  
$\mathcal{P}_0$.

\, \item  Finally, let 
$$
\Delta:=\{(z,n)\in\Delta_0\times\mathbb{N}\cup\{0\}:0\le n<R(z)\}
$$
where each point $z\in\Delta_0$ is identified with $(z,0)$. 

\sp The map $T:\De\lra\De$ is defined as follows.
$$
T(z,n):=
\begin{cases}
(z, n+1)~~{\rm if}~~ n+1<R(z)\\
(T_0(z),0)~~{\rm if}~~ n+1=R(z)
\end{cases}
$$
We assume that the map $T:\De\lra\De$ is topologically mixing, meaning that for all $A\in \mathcal{P}_0$, all $k\in\{0,1,\ld,R(A)-1\}$, all $B\in\mathcal{P}_0$, and all $l\in\{0,1,\ld,R(B)-1\}$ there exists an integer $N\ge 0$ such that
$$
T^n(A\times\{k\})\spt B\times \{l\}
$$
for all $n\ge N$.
\end{itemize}

\sp The measure $m$ is spread over the whole space $\Delta$ by putting
\beq\label{520180515}
\tilde m_{|\Delta_0}:=m 
\  \  \text{ and }  \  \
\tilde m_{|\Ga\times\{j\}}
:=m_{|\Ga}\circ\pi^{-1}_j
\eeq
for all $\Ga\in \mathcal{P}_0$, where $\pi_j(z,0)=(z,j)$. Thus, we have the following.

\bobs\lab{120200124}
Themeasure $\^m$ is finite if and only if 
$$
\int_{\Delta_0}Rdm<+\infty. 
$$
\eobs
The separation time 
$$
s((x,k),(y,l))
$$ 
between any two points $(x,k)$ and $(y,l)$ in $\De$ is defined to be $s(x,y)$ if $k=l$ and $x,y$ are in the same set of the partition $\mathcal{P}$. Otherwise we set $s(x,y)=0$. We define
the space
$$
C_\beta(\De):=\big\{\varphi:\Delta\to \mathbb{R}:\exists \ C_\varphi>0~~{\rm such~
  that}~~|\varphi(x)-\varphi(y)|
\le C_\varphi\beta^{s(x,y)}~~\forall x,y\in\Delta\big\}.
$$
Note that, in particular,
$$
|\varphi(x)-\varphi(y)|
\le C_\varphi
$$
for all $x,y\in\Delta$. In particular, all functions in $C_\beta(\De)$ are uniformly bounded.

\sp We refer to the pentuple $\mathcal{Y}=(\De_0,m,T_0,\mathcal{P}_0,R)$ as a Young tower. The following basic result has been essentially proved in \cite{lsy2}; see \cite{gouezel3} for a refinement which we incorporated in our hypotheses. 

\sp
\begin{thm}\label{lsyoung0}
If $\mathcal{Y}=(\De_0,m,T_0,\mathcal{P}_0,R))$ is a Young tower and
$\int_{\De_0}R dm<+\infty$, then 
\begin{itemize}
\item There exists a unique probability $T$--invariant measure $\nu$, absolutely continuous with respect to $\tilde m$. 

\, \item The Radon--Nikodym derivative $d\nu/d\tilde m$ is bounded
from below by a positive constant. 

\, \item The dynamical system $(T,\nu)$ is metrically exact, thus ergodic. 
\end{itemize}
\end{thm}

\sp In order to formulate the further results, we need one notion more. This is the concept of the first entrance time. In our context this is the function $E:\De \lra \{0, 1, 2, \ldots, \infty\}$, given by  the formula:
\beq\label{2_2018_08_04}
E_T(x):= \min\big\{ n \in \{0, 1, 2, \ldots, \infty\}: \,  \, T^n(z) \in \De_0\big\}.
\eeq
Its distribution is determined by the distribution of $R:\De_0\lra \{0, 1, 2, \ldots\}$; for every integer $n\ge 1$:
\beq\label{120180406}
\^m\(\{x\in\De:E_T(x)>n\}\)
=\sum_{k=n+1}^\infty m\(\{x\in\De_0:R(x)>k\}\)
\eeq
The following finer stochastic properties of the dynamical system
$(T,\nu)$ have been also essentially proved in \cite{lsy2}; see \cite{gouezel3} for a refinement which we incorporated in our hypotheses.

\sp\begin{thm}\label{lsyoung}
Let $\mathcal{Y}=(\De_0,m,T_0,\mathcal{P}_0,R))$ be a Young tower. Then
the following hold.
\begin{enumerate}

\sp\item If $\^m\(\{x\in\De:E_T(x)>n\}\)\lek n^{-\alpha}$ with some $\alpha>0$, then
\beq\label{cov_poly}
|C_n(\psi,g)|
=\lt|\int (\psi\circ T^n) g d\nu-\int\psi d\nu\int g d\nu\rt|
=\cO(n^{-\alpha})
\eeq
for all functions $\psi\in
L^\infty$ and $g\in C_\beta(\De)$.

\sp\item If $\^m\(\{x\in\De:E_T(x)>n\}\)\lek \theta^n$ (implied by $m\(\{x\in\De_0:R(x)>n\}\)\lek \theta^n)$) for some $0<\theta<1$, then there
exists $0<\tilde\theta<1$ such that for all functions $\psi\in
L^\infty$ and $g\in C_\beta$ we have,  
\beq\label{cov}
|C_n(\psi,g)|
=\lt|\int (\psi\circ T^n) g d\nu-\int\psi d\nu\int g d\nu\rt|
=\cO(\tilde\theta^n).
\eeq
\item If $\^m\(\{x\in\De:E_T(x)>n\}\)\lek n^{-\alpha})$ with some $\alpha>1$ (in particular, if $m(R>n)\lek \theta^n)$), then the Central Limit
  Theorem is satisfied for all functions $ g\in C_\beta(\De)$, that
  are not cohomologous to a constant in $L^2(\nu)$. This means, we recall, that there exists $\sigma>0$ such that
$$
\frac{\sum_{j=0}^{n-1}g\circ T^j-n\int gd\nu}{\sqrt{n}}\xrightarrow[\ n\to\infty \ ]{}
\mathcal N(0,\sigma)
$$ 
in distribution, where, as noted above, $\mathcal N(0,\sigma)$ is here the normal (Gaussian) distribution with $0$ mean and variance $\sigma$. 
\end{enumerate}
\end{thm}

\sp\fr We have also the following.

\sp\begin{thm}\label{LILII}
Let $\mathcal{Y}=(\De_0,m,T_0,\mathcal{P}_0,R))$ be a Young tower.
If 
$$
m\(\{x\in\De_0:R(x)>n\}\)\lek n^{-\alpha}
$$ 
with some $\alpha>2$ (in
particular, if $m(R>n)\lek \theta^n$) with some $\th\in(0,1)$, then the Law of Iterated Logarithm holds for all functions $g\in C_\beta(\De)$, that
are not cohomologous to a constant in $L^2(\nu)$. This means, we recall, that
there exists a constant $A_g\in (0,+\infty)$ such that
\begin{displaymath}
\limsup_{n\to\infty}\frac{S_{n}g(x)-n\int gd\nu}{\sqrt{n\log\log n}}= A_g
\end{displaymath}
for $\nu$-a.e. $x\in \De$.
\end{thm}

\bpf Take $\g>0$ so small that
$$
2+\g<\a.
$$
Recalling that each function in $g\in C_\beta(\De)$ is uniformly bounded, we then  get
$$
\begin{aligned}
\int_{\De_0}|g_{\De_0}|^{2+\g}\,d\nu
&\lek\int_{\De_0}R^{2+\g}\,d\nu 
\lek \int_{\De_0}R^{2+\g}\,dm
=\int_0^{+\infty}m\(R^{2+\g}>t\)\,dt \\
&=\int_0^{+\infty}m\(R>t^{\frac1{2+\g}}\,dt
\lek \int_0^{+\infty}t^{-\frac{\a}{2+\g}} \\
&<+\infty.
\end{aligned}
$$
A direct application of Theorem~\ref{LILfrominducing} and stochastic law results in \cite{MU2} (for countable alphabet subshifts of finite type) then completes the proof. 
\epf

\sp\fr The status of this theorem is however somewhat different than that
of Theorem~\ref{lsyoung}. It was not proved in \cite{lsy1} or \cite{lsy12} but was concluded in \cite{SUZ} from Theorem~\ref{LILfrominducing} and the
stochastic law results in \cite{MU2} exactly as we did above. It also follows from \cite{MN} and \cite{gouezel-IPAS}, where the Invariance Principle Almost Surely was established for Young's Towers from which the Law of Iterated Logarithm follows in a well known standard way.

\

\section{Rokhlin's Natural Extension}

\sp The main theorem of this section, Theorem~\ref{t1ms33}, permits
one to replace  sometimes an endomorphism 
$T:X\to X$, preserving some probability measure, which is not necessarily invertible, with the (invertible) automorphism
$\widetilde{\vp}:\widetilde{X}\to \widetilde{X}$ . This may turn out to
be of great advantage in some proofs, since dealing with invertible
maps is in many cases easier than dealing with non--invertible
ones. Natural extensions share many properties with the original
maps. These can be defined in the context of Lebesgue spaces. Lebesgue spaces form the core of ergodic theory and are central in the descriptive set theory. The following theorem charactrizes them.

\bthm\label{t_Lebesgue_Spaces}
If $(X_1,\rho_1)$ and $(X_2,\rho_2)$ are two Polish (separable, completely metrizable) topological spaces and $\mu_1$ and $\mu_2$ are two atomless Borel probability measures respectively on $X_1$ and $X_2$, then the measure spaces $(X_1,\cB_{X_1},\mu_1)$ and $(X_2,\cB_{X_2},\mu_2)$ are isomorphic. Also their completions  $(X_1,\^\cB_{X_1},\^\mu_1)$ and $(X_2,\^\cB_{X_2},\^\mu_2)$ are isomorphic.
\ethm

\fr This theorem leads us to the following.

\bdfn
Any complete probability space which is isomorphic to a convex combination of a (complete) measure space appearing in Theorem~\ref{t_Lebesgue_Spaces} and a collection of (countably many) atoms is called a Lebesgue space. 

\sp\fr In other words: A Lebesgue space is any probability space which is isomorphic to a convex combination of Lebesgue measure on $[0,1]$ and a collection of (countably many) atoms in $[0,1]$.
\edfn

\sp Let $(X,\mathcal{B}, \mu)$ be a Lebesgue space. In addition, let
$T:X\lra X$ be a surjective measurable map which preserves the measure
$\mu$. We will shortly define a new map related to $T$, but first we
form the space 
\[
\widetilde{X}
:=\big\{(x_n)_{n=0}^\infty\in X^{\N_0}:T(x_{n+1})=x_n\ \text
{ for all } n\geq0\big\}\subset X^{\N_0}. 
\]
Then, for every $k\geq0$ let $\pi_k:\widetilde{X}\to X$ denote the
projection onto the $k$-th coordinate of $\widetilde{X}$, that is 
\[
\pi_k((x_n)_{n=0}^\infty):=x_k.
\]
Finally, endow the space $\widetilde{X}$ with the smallest
$\sigma$-algebra $\widetilde{\mathcal{B}}$ for which every projection
$\pi_k:\widetilde{X}\to X$ is measurable. 

\

\bdfn
Define the measurable map $\widetilde{T}:\widetilde{X}\lra
\widetilde{X}$ by setting 
\[
\widetilde{T}((x_n)_{n=0}^\infty):=(T(x_0), x_0, x_1, x_2, \ldots).
\]
We refer to the map $\widetilde{T}$ as \index{(N)}{Rokhlin's natural 
extension}{\em Rokhlin's natural extension of $T$}. 
\edfn

\sp\fr Rokhlin's natural extension has the following properties.

\sp\bthm\label{t1ms31}
If $(X,\mathcal{B}, \mu)$ is a Lebesgue space and 
$T:X\lra X$ is a surjective measurable map which preserves the measure $\mu$, then the map $\widetilde{T}:\widetilde{X}\lra \widetilde{X}$ is invertible
and its inverse $\widetilde{T}^{-1}:\widetilde{X}\lra \widetilde{X}$ is
given by the formula 
\[ 
\widetilde{T}^{-1}((x_n)_{n=0}^\infty):=(x_{n+1})_{n=0}^\infty.
\]
Furthermore, for each $n\geq0$, the following diagram commutes:
\[\begin{tikzcd}
\widetilde{X} \arrow{r}{\widetilde{T}} \arrow[swap]{d}{\pi_n} &\widetilde{X}\arrow{d}{\pi_n} \\
{X} \arrow{r}{T} & {X}
\end{tikzcd}
\]
\ethm

\sp\fr This theorem obvious. After passing to completed $\sigma$--algebras,
we deduce the following result. 

\sp\bthm\label{t2ms31}
If $(X,\mathcal{B}, \mu)$ is a Lebesgue space and 
$T:X\lra X$ is a surjective measurable map which preserves the measure $\mu$, then there exists a unique probability measure $\tilde{\mu}$ on the space
$(\widetilde{X}, \widetilde{\mathcal{B}})$ such that 
\[
\tilde{\mu}\circ \pi_n^{-1}=\mu, \ \text{ for every } n\geq0.
\]
\ethm

\bpf
This follows directly from the Daniel--Kolmogorov Consistency Theorem
(see Theorem~3.6.4 in \cite{Parthasarathy}). 
\epf

\sp

In light of Theorem \ref{t2ms31}, for every set $A\in
\mathcal{B}(X)$ and every $n\geq0$ we have that 
$$
\aligned
\tilde{\mu}\circ \widetilde{T}^{-1}(\pi_n^{-1}(A))
&= \tilde{\mu}\circ (\pi_n\circ \widetilde{T})^{-1}(A) 
= \tilde{\mu}\circ (T\circ \pi_n)^{-1}\\
&=\tilde{\mu}\circ \pi_n^{-1}\circ T^{-1}(A)\\
&=\mu \circ T^{-1}(A)\\
&=\mu(A)\\
&=\tilde{\mu}(\pi_n^{-1}(A)).
\endaligned
$$
It therefore follows from the definition of
$\widetilde{\mathcal{B}}$ that $\tilde{\mu}\circ
\widetilde{\vp}^{-1}=\tilde{\mu}$. Considering this along with Theorem
\ref{t1ms31} and Theorem \ref{t2ms31}, we obtain the following
result. 

\sp\bthm\label{t1ms33}
Let $(X,\cB, \mu)$ be a Lebesgue space and let $T:X\lra X$ be a surjective measurable map which  preserves the measure $\mu$. If
$\widetilde{\vp}:\widetilde{X}\lra\widetilde{X}$ is Rokhlin's natural
extension of $T:X\to X$, then we obtain the following: 
\begin{itemize}
\item[(a)] For every $n\geq0$ the following diagram commutes.
\[\begin{tikzcd}
\widetilde{X} \arrow{r}{\widetilde{T}} \arrow[swap]{d}{\pi_n} &\widetilde{X}\arrow{d}{\pi_n} \\
{X} \arrow{r}{T} & {X}
\end{tikzcd}
\]
\item[(b)] There exists a unique measure $\tilde{\mu}$ on the space
    $(\widetilde{X}, \widetilde{F})$ such that 
  \[
  \tilde{\mu}\circ \pi_n^{-1}={\mu},\ \text{ for all }n\geq0.
  \]
  \item[(c)]  We have that $\tilde{\mu}\circ
    \widetilde{\vp}^{-1}=\tilde{\mu}$, that is, the probability measure
    $ \tilde{\mu}$ is $\widetilde{\vp}$-invariant. 
\end{itemize}
\ethm

\sp Note that the surjectivity assumption is in fact not essential,
since the sets $T^n(X)$, for each $n\geq0$, are measurable (as $X$ is
a Lebesgue space), $\mu\left(\bigcap_{n=0}^\infty T^n(X)\right)=1$ and
the map 
$$
T:\bigcap_{n=0}^\infty T^n(X)\lra \bigcap_{n=0}^\infty T^n(X)
$$
is clearly surjective. 

\sp\fr This theorem permits one to sometimes replace the  endomorphism
$T$, which is not necessarily invertible, with the automorphism
$\widetilde{\vp}:\widetilde{X}\to \widetilde{X}$. This may turn out to
be of great advantage in some proofs, since dealing with invertible
maps is frequently easier than dealing with non-invertible
ones. Natural extensions share many properties with the original
maps, for example ergodicity as it is shown in the next theorem. 

\sp\bthm\label{t2ms33}
If $(X,\cB,\mu)$ is a Lebesgue space and $T:X\lra X$ is a
measurable map preserving the measure $\mu$, then the natural
extension measure $\tilde{\mu}$ on $\widetilde{X}$ is ergodic with
respect to $\widetilde{\vp}$ if and only if the measure $\mu$ is ergodic
with respect to the map $T:X\lra X$. 
\ethm
\bpf
Suppose first that $T:X\to X$ is not ergodic. Then there exists a set
$A\in\mathfrak{ F}$ such that $T^{-1}(A)=A$ and $0<\mu(A)<1$. It then
follows from Theorem \ref{t1ms33} (b) that
$\tilde{\mu}(\pi_0^{-1}(A))=\mu(A)\in(0,1)$. Further, from Theorem
\ref{t1ms33} (a), applied with $n=0$, it follows that 
\[
\widetilde{\vp}^{-1}(\pi_0^{-1}(A))
=(\pi_0\circ \widetilde{\vp})^{-1}(A)=(T\circ \pi_0)^{-1}(A)
=\pi_0^{-1}(T^{-1}(A)) 
= \pi_0^{-1}(A).
\]
Therefore $\widetilde{\vp}$ is not ergodic either.

\sp\fr Now assume that $T:X\to X$ is ergodic. Birkhoff's Ergodic
Theorem then yields that 
\[
\lim_{n\to\infty}\frac1n \sum_{j=0}^{n-1} g\circ T^j= \int_X g\ d\mu,
\]
for every function $g\in L^1(X, \cB, \mu)$ and the convergence
is understood in the Banach space $L^1({X}, {\cB},
{\mu})$. Invoking Theorem \ref{t1ms33}, this implies that 
\begin{eqnarray}\label{1ms35}
\lim_{n\to\infty}\frac1n \sum_{j=0}^{n-1} G\circ \widetilde{\vp}^j=
\int_X G\ d\tilde{\mu}, 
\end{eqnarray}
where $G((x_n)_{n=0}^\infty):=g(x_0)=g\circ \pi_0((x_n)_{n=0}^\infty)$
and the convergence is understood in the Banach space
$L^1(\widetilde{X}, \widetilde{\cB}, \tilde{\mu})$. Now for
every $n\geq 0$ let
$\widetilde{\cB}_n:=\pi_n^{-1}(\cB)$. Since $T\circ
\pi_{n+1}=\pi_n$, we obtain that 
\[
\widetilde{\cB}_{n+1}
=\pi_{n+1}^{-1}(\cB)\supseteq \pi_{n+1}^{-1}(T^{-1}(\cB))
=(T\circ \pi_{n+1})^{-1}(\cB)=\widetilde{\cB}_n.
\] 
Moreover, by definition $\widetilde{\cB}$ is the
$\sigma$-algebra generated by the sequence
$(\widetilde{\mathfrak{ F}}_n)_{n=0}^\infty$. It therefore follows from
the Martingale Convergence Theorem that for any $F\in
L^1(\widetilde{X}, \widetilde{\cB}, \tilde{\mu})$, 
\beq\label{2ms35}
F=\lim_{n\to\infty} E_{\tilde{\mu}}(F|\widetilde{\cB}_n),
\eeq
where the convergence is again understood in the Banach space
$L^1(\widetilde{X}, \widetilde{\cB}, \tilde{\mu})$. However,
for each $n\geq0$, the function
$E_{\tilde{\mu}}(F|\widetilde{\cB}_n)$ depends only upon the
$n$-th coordinate of a point in $\widetilde{X}$ and can therefore be
represented as $f_n\circ \pi_n$, where $f_n\in L^1(X, \cB,
\mu)$. Now, for any $j\geq0$, we deduce that 
\[
f_n\circ \pi_n\circ \widetilde{\vp}^{j+n}=f_n\circ(\pi_n\circ
\widetilde{\vp}^n)\circ \widetilde{\vp}^j=f_n\circ \pi_0\circ
\widetilde{\vp}^j, 
\]
and it therefore follows from (\ref{1ms35}) that
\[
\lim_{k\to\infty} \frac1k \sum_{j=0}^{k-1}
E_{\tilde{\mu}}(F|\widetilde{\cB}_n)\circ \widetilde{\vp}^j
=\int_X, E_{\tilde{\mu}}(F|\widetilde{\cB}_n)\ d\tilde{\mu}
= \int_X F\ d\tilde{\mu}
\]
 where the convergence is understood in the Banach space
 $L^1(\widetilde{X}, \widetilde{\cB}, \tilde{\mu})$. Combining
 this with (\ref{2ms35}), we conclude that 
\[
\lim_{k\to\infty} \frac1k \sum_{j=0}^{k-1} F\circ
\widetilde{\vp}^j=\int_X F\ d\tilde{\mu}, 
\]
where once more the convergence is understood in the Banach space
$L^1(\widetilde{X}, \widetilde{\cB}, \tilde{\mu})$. Therefore
the transformation $\widetilde{\vp}:\widetilde{X}\to \widetilde{X}$ is
ergodic and we are done. 
\epf

\chapter{Infinite Invariant Measures; Finer Properties}\label{infinite-measure}

\sp In this chapter we deal with invariant measures that are infinite. The
concepts introduced and explored here are meaningful just for such
measures. Even if they make sense for finite measures, they are trivial
then. 

\section{Counterexamples to Ergodic Theorems}\label{CountET}

\sp Hopf's Ergodic Theorem (Theorem~\ref{t2j73}) being powerful and
interesting in itself, somewhat surprisingly, rules out any hope for a more 
direct version  of Birkhoff's Ergodic Theorem in the case  of
infinite measures. The world of infinite invariant measures is indeed
very different from the one of finite measures. The following
proposition provides the first indication of this.

\sp\bprop\label{p1j153} 
Suppose that $(X, \mathfrak{F}, \mu)$ is
$\sigma$--finite  measure space and that  $T:X\lra X $ is a
measurable map  preserving measure measure $\mu$. If
$\mu(X)=+\infty$  then for  every function $ f \in L^1(\mu)$ we  have
that  
$$
\lim_{n \to \infty} \frac{1}{n}S_n f =0  \quad \mu-a.e.
$$
\eprop

\bpf  Since the  measure  $\mu$ is $\sigma$-finite  and
$\mu(X)=\infty$, there exists  a sequence $(A_n)_{n=1}^\infty$  of
 measurable sets such that $0<\mu(A_n)<+\infty$ for all $n\ge 1$ and such that
   \beq\label{1j153}
\lim_{n \to \infty}\mu(A_n)=\infty.
\eeq
Since $|f| \in L^1(\mu)$, the  theorem conclude from
Theorem~\ref{t2j73} that  for every $k\geq 1$ we have
$$
\varlimsup_{n \to \infty}\frac{1}{n} S_n|f|
\leq \varlimsup_{n\to\infty}\frac{S_n|f|}{S_n\1_{A_k}}
=\frac{\int |f|d \mu}{\int\1_{A_k}\, d \mu}
=\frac{ \int|f| \,d \mu }{\mu(A_k)}\quad  \mu-a.e.
$$ So
by \ref{1j153}, $\limsup_{n \to \infty} S_n |f|=0$ almost
everywhere. Since $|\frac{1}{n}S_nf|\leq \frac{1}{n}S_n|f|$ we are thus
done.
\endpf

\sp As a matter of fact the last two results of this section show that
no matter  how desired they could be, Hopf's Ergodic Theorem, i.e. Theorem~\ref{t2j73}, rules out any hope for even much weaker forms of Theorem~\ref{Birkhoff} in  the case of infinite invariant measures. Indeed, for all $f\in L^1(\mu)$ Proposition~\ref{p1j153} (a corollary of Hopf's Ergodic Theorem) implies that 
\[
\lim_{n\to\infty}\frac{1}{a_n}S_nf(x)=0\, \mbox{ for $\mu$--a.e. } x\in X
\]
if there exists $C>0$ such that $a_n\geq Cn$ for all $n$. We will now show that 
there are no constants $a_n>0$ such that
\[
\lim_{n\to\infty}\frac{1}{a_n}S_nf(x)=\int_X f\,d\mu
\, \mbox{ for $\mu$--a.e. } x\in X,\ \ \ \forall f\in L^1(\mu).
\]
We will accomplish this in two steps. Their statements and proofs are 
inspired by Jon Aaronson's book \cite{Aa}. The first step will require 
the following proposition.

\sp\fr Since we will simultaneously deal with Birkhoff's sum with respect
to the maps $ T:X\lra X$ and the first return maps $\tau_F: F\lra F$, in order do discern between them we put
\beq\label{suma-2}
S_{F,n}(g):= \sum_{j=0}^{n-1} g \circ T^j_F,
\eeq\index{(S)}{$S_{F,n}$}
where $g:F\lra \mathbb R$ is a measurable function.

\sp\bprop\label{p2.3.1in{A}} Suppose that $(X, \mathfrak{F}, \mu)$ is
a probability space, $T:X\lra X $ is an ergodic map preserving
measure $\mu$, and  $f:X\lra \mathbb R$ is measurable. Let $a:[0,+
\infty)\to [0,+\infty)$ be continuous, strictly increasing function
satisfying $ \frac{a(x)}{x}\searrow 0$ as $x \to+\infty$. If
$$ 
\int_Xa(|f|)d\mu <+\infty,
$$
then   
$$\lim_{n \to \infty}\frac{a(|S_nf|)}{n}=0\quad \mu-a.e.\,\, \mbox{in}\,\, X.
$$
\eprop

\bpf We first shall show that 
\beq\label{1_mu_2014_11_12}
a(x+y)\leq a(x)+a(y)
\eeq
for al $x, y\ge 0$.
Indeed, assume without loss of generality that $x\le y$. Then
$$
a(x+y)
\le \frac{a(y)}{y}(x+y)
=a(y)+\frac{a(y)}{y}x
\le a(y)+\frac{a(x)}{x}x
=a(x)+a(y),
$$
and \eqref{1_mu_2014_11_12} is proved. Dealing directly with the proposition,
we first establish i under the
additional assumption that $a(0)=0$. Given $ \varepsilon >0$, we will prove that
$$
\limsup_{n \to \infty}\frac{a(|S_nf|)}{n}\leq \varepsilon \,\,\, \  a.e.
$$
In order to do this, we firs claim that there are two measurable functions $g,h: X \to \mathbb R$ such that $|f|=g+h$, $\sup_X(g)<+\infty$ and $\int_X a(h)\,d \mu
<\vep$. This is because for the sets $A_M:=\{x \in X: \,\, |f(x)|\geq M
\}$, we have
$$
\lim_{M\to+\infty} a( |f| \1_{A_M})=0 \  \   \  a.e.,
$$ 
whence, by the Lebesgue Dominated Convergence Theorem, which is applicable since  $a(|f|)$ is integrable, we have that 
$$
\lim_{M\to+\infty}\int_X a( |f| \1_{A_M})d \mu =0,
$$  
and so, for a sufficiently large $M>0$, setting
$$ 
g:=|f|\1_{A_M^c} \  \text{ and } \  \, h:=|f|\1_{A_M},
$$
gives the required representation of $|f|$.
Using \eqref{1_mu_2014_11_12} and Birkhoff's Ergodic Theorem (Theorem~\ref{Birkhoff}), we have that
$$
\aligned
\frac{a(|S_n f|)}{n}
&\leq  \frac{a(|S_n g|)}{n} + \frac{a(|S_nh|)}{n}\\
& \leq \frac{a(| M n|)}{n} + \frac{S_n a(h)}{n}\xrightarrow[\ n\to\infty \ ]{} 
  \int_X  a(h)
d\mu < \varepsilon.
\endaligned
$$
Hence, letting $\vep\downto 0$, our assertion follows (under the assumption that $a(0)=0$). 
We now shall establish the proposition without assuming $a(0)=0$. To do this we fix any $m >0$ such that $a(m)/m=\a\in(0,1)$. We then define
$$   
\^a(x):=
\begin{cases}
\a x  &\text{ if } \ 0\le x \le m \\
a(x)  &\text{ if } \  x\ge m.
\end{cases}
$$
It is straightforward to verify that $\^{a}$ is continuous,
increasing, $a \equiv \^{a}$ on $[m, \infty)$, that $\frac{\^{a}}{n}
\searrow 0$ as $ x \to \infty$ and $\^a(0)=0$. Therefore, by the previous step ($a(0)=0$), we have  
$$
\lim_{n\to \infty}\frac{\^{a}(| S_n f|)}{n}= 0
$$ 
$a.e.$ in $X$. Since $a(x)\le \max\{a(m),\^a(x)\}\le a(m)+\^a(x)$, we thus get
$$
0\le\limsup_{n\to\infty}\frac{a(| S_n f|)}{n} 
\leq \limsup_{n\to\infty}\frac{a(m)+\^{a}(|S_nf|)}{n}
=\lim_{n\to \infty}\frac{\^{a}(| S_n f|)}{n}=0.
$$
The proof is complete.
\endpf

\sp Let $L_+^1(\mu)$ denote the set of functions $f\in L^1(\mu)$ such that $\int_X f d \mu>0$.
The next two theorems rule out any hope for any pointwise version of Birkhoff's Ergodic Theorem in the context of infinite invariant measures.

\sp\bthm\label{t2.4.1in{A}} Suppose that $(X, \mathfrak{F}, \mu)$ is a
measurable space with a $\sigma$--finite infinite measure $\mu$, and that
$T:X\lra X$ is an ergodic conservative transformation preserving $\mu$.
If $a(n)\upto  \infty$ and $ \frac{a(n)}{n}\searrow 0$ as $n \to \infty$,
then the following hold.
\begin{itemize}
\item[(1)] If there exists a set $F \in \mathfrak{F}$ such that $0<
\mu(F)< \infty$ and $\int_F a( \tau_F)\, d \mu< +\infty$  then
$$
\lim_{n\to\infty}\frac{S_nf}{a(n)}=+ \infty \,\,\,\, \  \mu-a.e.\,\, \forall f \in
L^1(\mu)_+.
$$
\item [(2)] Otherwise, 
$$\liminf_{n\to\infty}\frac{S_nf}{a(n)}=0
\,\,\,\, \  \mu-a.e.\,\, \forall f \in L^1(\mu)_+.
$$
\end{itemize}
\end{thm}

\bpf 
We begin with proving (1). Suppose that  $F \in \mathfrak{F}$, $ 0<
\mu(F)< +\infty$ and that $\int_F a(\tau_F) d \mu <+\infty $. Then by
Proposition~\ref{p2.3.1in{A}},
\beq\label{2_mu_2014_11_12}
\lim_{n \to \infty} \frac{a( S_{F,n}((\1_F)_F)}{n}
=0\,\, \, \  \  a.e. \,\, \mbox{on} \,\, F.
\eeq
For every $x\in F_\infty$ and every $n\ge 1$ fix the only integer $k_n(x)\ge 1$ such that
$$ 
\sum_{k=0}^{k_n(x)-1}\tau_F \circ T_F^k(x) 
\leq n 
< \sum_{k=0}^{k_n(x)}\tau_F \circ T_F^k(x),
$$ 
Since also $S_{F,k_n(x)}((\1_F)_F)\equiv k_n(x)$, we get that 
$$
\frac{S_n \1_F}{a(n)(x)}
\geq\frac{S_{F,k_n(x)}((\1_F)_F)(x)}{a(
 S_{F,k_n(x)+1}((\1_F)_F))(x)}=\frac{k_n(x)}{a(
  S_{F,k_n(x)+1}((\1_F)_F))(x)} \xrightarrow[\ n\to\infty \ ]{} 
 +\infty
$$ 
where the divergence to $+\infty$ is due to \eqref{2_mu_2014_11_12}.
Clearly, the set of all points $x\in X$ on which this divergence occurs is $T$
invariant, and, as containing $F_\infty$, it must be equal to $X$ modulo $\mu$. Thus
$$ 
\lim_{n \to \infty}\frac{S_n \1_F}{a(n)}= +\infty 
$$ 
$\mu$-a.e. on $X$, and (1) follows from Hopf's Ergodic Theorem
(Theorem~\ref{t2j73}). 

\sp\fr We now prove (2). To this end, we note first, that for every $f\in L_+^1(\mu)$,
$$ 
\frac{ S_n\circ T}{a(n)}\leq \left( 1+ \frac{1}{n}\right)
\frac{S_{n+1}f }{ a(n+1)}
$$ 
Hence 
$$
\liminf_{n \to \infty}\frac{S_nf}{a(n)}\circ T
\le \liminf_{n \to \infty}\frac{S_nf}{a(n)}.
$$ 
So, this lower limit is $\mu$--a.e constant on $X$ by ergodicity of $T$
with respect to the measure $\mu$. Denote it by $2\varepsilon\ge
0$. Now, seeking a contradiction, suppose that the hypothesis of condition (2) is satisfied but its assretion fails for the function $f$. Then $\varepsilon>0$. Fix arbitrarily a set $F\in \mathfrak{F}$ such that $0<\mu(F)<+\infty$. It then follows
from Egorov's Theorem that there exists a measurable set $G\sbt F$
such that $0< \mu(G)<\infty$, and 
$$ 
S_n \1_F(x) \geq \varepsilon a(n)
$$
forall $x \in G$ and all $n \geq 1$ large enough . Hence
$$ 
S_{\tau_G(x)}\1_F(x) \geq \varepsilon a ( \tau_G(x))
$$
forall $x \in G$. Using Hopf's Ergodic Theorem
(Theorem~\ref{t2j73}) this implies, that for $\mu$-a.e. $x \in G$, we have
$$
\aligned 
\frac1n
S_{n,G}(a\circ\tau_G)(x)
&= \frac{1}{n} \sum_{k=0}^{n-1}a( \tau_G(T^k(x)) \\
&\leq \frac{1}{n \varepsilon} \sum_{k=0}^{n-1}S_{ \tau_G(T^k(x))}\1_F (T^k_G(x))\\
&=\frac{1}{n \varepsilon} S_{S_n\tau_G(x)}\1_F (x) \\
&=\frac{1}{\varepsilon}\frac{S_{S_n\tau_G(x)}\1_F(x)}{S_{S_n\tau_G(x)}\1_F(x)}
\xrightarrow[\ n\to\infty \ ]{} 
  \frac{\mu(F)}{\varepsilon \mu(G)}.
\endaligned
$$ 
It thus follows that $\int_G a( \tau_G) d\mu < \infty$
contracting the hypothesis (2). \endpf

\sp\bthm\label{t2.4.2in{A}} Suppose that $(X,\mathfrak{F},\mu)$ is a $\sigma$--finite measure space and $T:X\lra X $ is an ergodic conservative
transformation preserving measure $\mu$. Let $(a_n)_{n=1}^\infty$
be a sequence such that $a_n>0$ for all $n\in\N$. Then,

\begin{itemize}
\item[(1)] either
$$\liminf_{n\to\infty}\frac{S_nf}{a_n}=0\,\,\,\,\, \forall f \in
L^1(\mu)_+.$$
\item [(2)] or \ $\exists\, n_k\to \infty$ such that
$$\liminf_{n\to\infty}\frac{S_{n_k}f}{a_{n_k}}=+\infty
\,\,\,\,a.e.\,\,\,\, \forall f \in L^1(\mu)_+.$$
\end{itemize}
\end{thm}

\bpf
If $(a_n)_{n=1}^\infty$ is bounded, then~(2) holds. Indeed, let $f\in L_+^1(\mu)$.
Since $\int_X f\,d\mu>0$, there 
exists $\e>0$ and $B\in\F$ such that 
$$
\mu(B)>0 
\  \  \  {\rm and} \  \  \ 
f_B\geq\e.
$$ 
As $T$ is conservative and ergodic, Theorem~\ref{t1j59} affirms that 
$\mu(X\backslash B_\infty)=0$. Therefore for $\mu$--a.e. $x\in X$ there exists 
a sequence $(n_k(x))_{k=1}^\infty$ such that $n_k(x)\nearrow\infty$ and 
$$
T^{n_k(x)}(x)\in B.
$$
Hence $S_{n_k(x)}f(x)\geq k\e$. In fact, $S_nf(x)\geq k\e$ 
for any $n\geq n_k(x)$ since $f\in L_+^1(\mu)$. Therefore,
$$
\lim_{n\to\infty}S_nf(x)=+\infty.
$$
As $(a_n)_{n=1}^\infty$ is bounded, it follows that 
$$\lim_{n\to\infty}\frac{S_nf(x)}{a_n}=\infty
$$ 
for $\mu$--a.e. $x\in X$. So~(2) holds for this function $f$ with $(n_k)_{k=1}^\infty=(n)_{n=1}^\infty$. By Hopf's Ergodic Theorem, this actually holds for every $f\in L_+^1(\mu)$,
with the same sequence. 

We can thereby restrict our attention to the case 
$$
\limsup_{n\to\infty}a_n=+\infty.
$$ 
Suppose that~(1) does not hold. That is, there exists a set $A\in\mathcal{A}$ with $\mu(A)>0$ 
and a function $f\in L_+^1(\mu)$ such that 
\beq\label{convAliminf2}
F(x):=\liminf_{n\to\infty}\frac{S_nf(x)}{a_n}>0\  \mbox{ for $\mu$--a.e. } \ x\in A.
\eeq
By Hopf's Ergodic Theorem, this actually holds for every $f\in L_+^1(\mu)$,
with the same set $A$. Then for $\mu$--a.e. $x\in A$, we get that
\begin{eqnarray*}
0\leq
\limsup_{n\to\infty}\frac{a_n}{n}
&=&\limsup_{n\to\infty}\Bigl[\frac{a_n}{S_nf(x)}\frac{S_nf(x)}{n}\Bigr] \\
&\leq&\limsup_{n\to\infty}\frac{a_n}{S_nf(x)}\limsup_{n\to\infty}\frac{S_nf(x)}{n} \\
&=&\Bigl[\liminf_{n\to\infty}\frac{S_nf(x)}{a_n}\Bigr]^{-1}\lim_{n\to\infty}\frac{S_nf(x)}{n}
=0
\end{eqnarray*}
by~(\ref{convAliminf2}) and Proposition~\ref{p1j153}. 
Thus, $a_n=o(n)$ as $n\to\infty$. For every $n\in\N$, set
\[
\overline{a}_n=\max_{1\leq k\leq n}a_k.
\]
Clearly, $a_n\leq\overline{a}_n$ for all $n\in\N$ and $\overline{a}_n\nearrow\infty$ 
as $n\nearrow\infty$. Moreover, for each $n\in\N$ there is $1\leq k(n)\leq n$ such that 
$\overline{a}_n=a_{k(n)}$. Note that $k(n)\rightarrow\infty$ as $n\to\infty$. Then
\[
\liminf_{n\to\infty}\frac{S_nf}{a_n}
\geq\liminf_{n\to\infty}\frac{S_nf}{\overline{a}_n}
=\liminf_{n\to\infty}\frac{S_nf}{a_{k(n)}}
\geq\liminf_{n\to\infty}\frac{S_{k(n)}f}{a_{k(n)}}
\geq\liminf_{n\to\infty}\frac{S_nf}{a_n},
\]
where the last inequality follows from the fact that the lim inf of a subsequence of
a sequence is greater than or equal to the lim inf of the full sequence. Hence 
\beq\label{ineqan}
\liminf_{n\to\infty}\frac{S_nf}{\overline{a}_n}
=\liminf_{n\to\infty}\frac{S_nf}{a_n}
>0\ \ \ \mu\mbox{--a.e. on }A,\ \ \ \forall\,f\in L_+^1(\mu).
\eeq
Next, set 
$$
f_n:=\frac{\overline{a}_n}{n},
$$
and let $1=n_0<n_1<\ldots$ be defined by
\[
\{n_k\}_{k\in\N}=\bigl\{j\ge2:f_i>f_j,\ \forall\,1\leq i\leq j-1\bigr\}.
\]
For every $k\geq0$,
\[
f_{n_k}>f_{n_{k+1}},\ \ \ n_k f_{n_k}\leq n_{k+1}f_{n_{k+1}},
\]
whence
\[
0<\frac{n_k}{n_{k+1}}\leq\frac{f_{n_{k+1}}}{f_{n_k}}<1
\]
and there thus exists $\alpha_k\in(0,1]$ such that
\[
\left(\frac{n_k}{n_{k+1}}\right)^{\alpha_k}
=\frac{f_{n_{k+1}}}{f_{n_k}}.
\]
Define
\[
f(x)=\frac{f_{n_k}n_k^{\alpha_k}}{x^{\alpha_k}},\ \ \ x\in[n_k,n_{k+1}],\ k\in\N,
\]
and
\[
a(x)=xf(x).
\]
Evidently,
\[
a(n_k)=\overline{a}_{n_k}, \ \ \ \forall\,k\in\N.
\]
By definition of the $n_k$, we have that for $k\in\N$, $n\in[n_k,n_{k+1})$,
\[
f_n\geq f_{n_k}, \mbox{ hence }f_n\geq f(n),
\]
whereby 
\[
a(n)\leq\overline{a}_n,\ \ \ \forall\,n\in\N.
\]
Hence, following~(\ref{ineqan}),
\[
\liminf_{n\to\infty}\frac{S_nf}{a(n)}>0,\ \ \ \mu\mbox{-a.e. on }A,\ \ \ \forall\,f\in L_+^1(\mu).
\]
It is evident that 
\[
a(n)\nearrow\infty,\ \ \ \frac{a(n)}{n}\searrow0\ \ \ \mbox{ as }n\nearrow\infty.
\]
So by Theorem~\ref{t2.4.1in{A}},
\[
\lim_{n\to\infty}\frac{S_nf}{a(n)}=+\infty\  \  \ \mu\mbox{--a.e. on } \ X,\ \ \ \forall\,f\in L_+^1(\mu).
\]
Then~(2) follows since $a_{n_k}\leq\overline{a}_{n_k}=a(n_k)$.
\qed \epf 

\

\section{Weak Ergodic Theorems}\label{WET}

\sp Although the previous sections has ruled out any hope for a direct counterpart of Birkhoff's Ergodic Theorem in the realm of infinite invariant measures, there are nevertheless some weak version of this theorem. By saying weak we mean the convergence in the sense of distribution and along a subsequence. We provide below some such theorems taken from various papers of Jon Aaronson, primarily in Jon Aaronson's book \cite{Aa} and the papers cited therein.

If $(X,\mathfrak{F}, \nu)$ is a probability space and
$$
\(Y_n:X\lra[-\infty,+\infty]\)_{n=1}^\infty
$$ 
is a sequence of random variables (measurable functions) with respect the $\sigma$--algebra $\mathfrak{F}$, then $(Y_n)_{n=1}^\infty$ is said to converge in distribution to some
probability distribution $P$ on $\mathbb R$ if and only if the sequence
$(\nu \circ Y^{-1}_n )_1^\infty$ of distributions of $Y_n$, $n\ge 1$, converges weakly to $P$. We then write:
$$
Y_n\xrightarrow{\  \  \,  {w}\hspace*{0.4cm}}P.
$$
If $P$ is the probability distribution of some random variable $Y:(\Om,\cA,Q)\lra[-\infty,+\infty]$ defined on some probability space $(\Om,\cA,Q)$, then we also say that the sequence $(Y_n)_{n=1}^\infty$ converges in distribution to $Y$, and we write
$$
Y_n\xrightarrow{\  \  \,  {w}\hspace*{0.4cm}}Y.
$$
Explicitly this means that
$$
\lim_{n\to\infty}\int_Xg\circ Y_n\,dm=E(g\circ Y):=\int_{\Om}g\circ Y\,dQ
$$
for every continuous function $g: [-\infty, \infty]\lra \mathbb R$.

If $(X,\mathfrak{F}, m)$  is a measure space ($m$
can be infinite), then a  sequence $(Y_n)_{n=1}^\infty$
 of real--valued measurable functions on $X$  is said to converge
 strongly in distribution to some probability distribution $P$ on
 $\mathbb R$ if and only if the sequence $(\nu \circ Y^{-1}_n)_{n=1}^\infty$
 converges weakly to $P$ for every probability measure $\nu$
 absolutely  continuous with respect to $m$. We then write:
$$
Y_n\xrightarrow{\  \  \,  {{\mathcal L}}\hspace*{0.4cm}}P.
$$

As above, if $P$ is the probability distribution of some random variable 
$$
Y:(\Om,\cA,Q)\lra\R
$$ 
defined on some probability space $(\Om,\cA,Q)$, then we also say that the sequence $(Y_n)_{n=1}^\infty$ converges strongly in distribution to $Y$, and we write
$$
Y_n\xrightarrow{\  \  \,  {{\mathcal L}}\hspace*{0.4cm}}Y.
$$
Explicitly this means that
$$
\lim_{n\to\infty}\int_Xhg\circ Y_n\,dm=E(g\circ Y):=\int_{\Om}g\circ Y\,dQ
$$
for every continuous function $g: [-\infty, \infty]\lra \mathbb R$ and every integrable function $h:X\lra[0,+\infty)$ such that $\int_Xh\,dm=1$. Furthermore, by the standard approximation argument, this means that 
\beq\label{120190614}
\lim_{n\to\infty}\int_Xhg\circ Y_n\,dm=E(g\circ Y):=\int_{\Om}g\circ Y\,dQ
\eeq
for every continuous function $g: [-\infty, \infty]\lra \mathbb R$ and every $m$--essentially bounded measurable function $h:X\lra[0,+\infty)$ such that $\int_Xh\,dm=1$. This in turn is equivalent to having that
\beq\label{220190614}
\lim_{n\to\infty}\int_Xhg\circ Y_n\,dm=E(g\circ Y):=\int_{\Om}g\circ Y\,dQ\int_Xh\,dm
\eeq
for every continuous function $g: [-\infty, \infty]\lra \mathbb R$ and every
function $h\in L^{\infty}(m)$.

\bthm\label{Proposition3.6.1}
Let $T:(X,\mathfrak{F},m)\lra (X,\mathfrak{F},\mu)$ be an ergodic transformation of a $\sigma$--finite measure space $(X,\mathfrak{F},\mu)$, and let $f:X\lra\mathbb R$ be a bounded measurable function. 

Suppose that $(n_k)_1^\infty$ is a sequence of positive integers diverging to $+\infty$ and that $(d_k)_1^\infty$ is a sequence of positive real numbers. 

Then there exist a subsequence $\(n_{k_l}\)_{l=1}^\infty$ diverging to $+\infty$, and a random variable $Y$ taking values in $[-\infty,\infty]$ such that
$$ 
\frac{1}{d_{k_l}}\sum_{j=0}^{n_{k_l}-1} f \circ T^j\xrightarrow{\  \  \,  {{\mathcal L}}\hspace*{0.4cm}} Y. 
$$
\ethm

\bpf As usually
$$
S_nf=\sum_{j=0}^{n-1}f\circ T^j.
$$
Because of \eqref{220190614} and since for every $h\in L^{\infty}(m)$, the function
$$
L^{\infty}(\mu)\ni k\longmapsto \int_Xhk\,d\mu\in\R
$$
belongs to $L^{\infty}(\mu)^*$, in order to prove our theorem it is sufficient to obtain a random variable $Y$ taking values in $[-\infty,+\infty]$, and  $m_l:= n_{k_l} \to \infty$  such that for each continuous function $g: [-\infty,+\infty]\to \mathbb R$
$$
g\circ\left(\frac{S_{n_{k_l}}f}{d_{k_l}}\right)\xrightarrow[\ l\to\infty \ ]{} 
 E(g\circ Y) \quad \text{weak * in}\,\,  L^{\infty}(m).
$$ 
For any $g \in C([-\infty,+\infty])$, the Banach space of real--valued continuous functions defined on $[-\infty,+\infty]$ endowed with the supremum norm $\|\cdot\|_\infty$, we have for all integers 
$k \geq 1$ that
$$ 
\left\|g\circ\left(\frac{S_{n_k}}{d_{k}}\right) \right\|_{L^{\infty}(\mu)} 
\leq \|g\|_{\infty}.
$$
Hence, because of Banach--Aaloglou Theorem there exists an increasing  subsequence $k'\to \infty$ and a function $\lambda(g)\in L^{\infty}(\mu)$ such that
$$ 
g\circ\left(\frac{S_{n_{k'}}}{d_{k'}}\right)\xrightarrow[\ k'\to\infty \ ]{} 
 \lambda(g) \quad \text{weak * \ in} \  L^{\infty}(\mu).
$$
We claim that the limit function $\lambda(g)$ is $T$--invariant. Indeed, putting for every $\varepsilon$:
$$
\om_g(\varepsilon)
:= \sup\big\{|g(y+h)-g(y)|:y,h \in \mathbb R, \, |h| <\varepsilon\big\},
$$
we have by the uniform continuity of $g$ that
$$
\lim_{\varepsilon\to 0} \om_g(\varepsilon)= 0.
$$
Hence
$$
\begin{aligned} 
\left\|g\circ\left(\frac{S_{n_k}}{d_{k}}\right)\circ T -g \circ \left(\frac{S_{n_k}}{d_{k}}\right)\right\|_\infty 
& =\left\|g\circ\left(\frac{S_{n_k}}{d_{k}}+ \frac{f \circ T^{n_k}- f}{d_k}\right) - g\circ\left(\frac{S_{n_k}}{d_{k}}\right)\right\|_\infty\\
&\leq \om_g\left(\frac{2||f||_\infty}{d_k}\right)\xrightarrow[\ k\to\infty \ ]{}.
0
\end{aligned}
$$
Thus, $\l(g)\circ T=\l(g)$. By ergodicity of $T$ this yields that the function $\lambda(g)$ is constant.

Since the Banach space $C([-\infty,+\infty])$ is separable, there exists $\mathcal G \sbt C([-\infty, \infty])$ be a countable, dense set in $C([-\infty,+\infty])$. A standard diagonalization argument yields an unbounded increasing susequence $(n_{k_l})_{l=1}6\infty$ and constants $\{\lambda(g)\in \mathbb R:g \in \mathcal G\}$ such that
\beq\label{120190615}
g\circ\left(\frac{S_{n_{k_l}}f}{d_{k_l}}\right)\xrightarrow[\ l\to\infty \ ]{} 
 \l(g) \quad \text{weak * \ in}\,\,  L^{\infty}(m)
\eeq
for all $g \in \mathcal G$. We claim that that this property can be extended to all elements of $C([-\infty,+\infty])$. More precisely, we claim that for all $h\in C([-\infty,+\infty])$ there exists  $\lambda(h) \in \mathbb R$ such that
$$
h\circ\left(\frac{S_{n_{k_l}}f}{d_{k_l}}\right)\xrightarrow[\ l\to\infty \ ]{} 
 \l(h) \quad \text{weak * \ in}\,\,  L^{\infty}(m)
$$
for all $h \in C([-\infty,+\infty])$.

In order to prove this, let $h\in C([-\infty,+\infty])$. Since the set $\mathcal G$ is dense in $C([-\infty,+\infty])$, there is a sequence  $(g_j)_{j=1}^\infty \in \mathcal G$ such that
$$
\lim_{j\to+\infty} \|g_j-h|_\infty=0.
$$
Because of \eqref{120190615}, for any probability measure $\nu$ on $X$ absolutely continuous with respect to $\mu$ with the Radon--Nikodym derivative $d\nu/d\mu\in L^{\infty}(\mu)$ and all $i, j\ge 1$, we have that
$$ 
|\lambda(g_i)-\lambda (g_j)|
\xleftarrow{\  \  \,  {{l\to\infty}}\hspace*{0.4cm}} 
\left|\int_X \left(g_i\left(\frac{S_{n_{k_l}}}{d_{k_l}}\right) - g_j\left(\frac{S_{n_{k_l}}}{d_{k_l}}\right) \right) d\nu \right| 
\leq \|g_i-g_j\|_\infty.
$$
Since, as convergent, the sequence $(g_j)_{j=1}^\infty$ is fundamental (Cauchy), we therefore conclude that the sequence $\(\l(g_j)\)_{j=1}^\infty$ is fundamental in $\R$. Thus there exists a limit
$$
\lambda(h)=:\lim_{j\to \infty}\lambda(g_j)\in \mathbb R,
$$ and
$$
|\lambda(g_j)-\lambda(h)|\leq \|g_j-h\|_\infty
$$
for all $j\geq 1$. For any probability measure $\nu$ on $X$ absolutely continuous with respect to $\mu$ with the Radon--Nikodym derivative $d\nu/d\mu\in L^{\infty}(\mu)$ and all integers $j\geq 1$, by applying the triangle inequality, we get that
$$
\begin{aligned} 
\lt|\int_X\!\!h\circ\left(\frac{S_{n_{k_l}}}{d_{k_l}}\right) d\nu-\lambda(h)\rt|
&\leq \left|\int_X \!\! h\circ\left(\frac{S_{n_{k_l}}}{d_{k_l}}\right) d\nu - \!\!\int_X g_j\circ\left(\frac{S_{n_{k_l}}}{d_{k_l}}\right) d\nu \right|+ \\
&\  \  \  \  \  \  \  \  \  \  \  \  \  \  \  \ +\left|\int_X g_j\left(\frac{S_{n_{k_l}}}{d_{k_l}}\right) d\nu- \lambda(g_j)\right|+|\lambda(g_j)-\lambda(h)|\\
& \leq  2\|h- g_j\|_\infty +  \left|\int_X g_j\circ\left(\frac{S_{n_{k_l}}}{d_{k_l}}\right) d\nu- \lambda(g_j)\right|
\xrightarrow{\  \  \,{{l\to\infty}}\hspace*{0.1cm}} 
2\|h- g_j\|_\infty.
\end{aligned}
$$
By letting $j \to \infty$ we thus conclude that
$$ 
h\circ\left(\frac{S_{n_{k_l}}}{d_{k_l}}\right)\xrightarrow[\ n\to\infty \ ]{} 
 \lambda(h) \quad \text{weak * \ in}\,\,  L^{\infty}(\mu).
$$
Since it also follows from this that$ \lambda: C([-\infty, +\infty])\lra \mathbb R$ is linear and positive, whence $ \exists $ a random variable $Y$ taking values in $[-\infty, \infty]$ such that
$$
\lambda(h)=E(h(Y)).
$$
for all $h\in C([-\infty, +\infty])$. The proof is complete.
\qed

\bcor\label{Corollary3.6.2}
Let $T:(X,\mathfrak{F},\mu)\lra (X,\mathfrak{F},\mu)$ be a conservative ergodic  measure preserving transformation of o a $\sigma$--finite measure space $(X,\mathfrak{F},m)$.
Suppose that $(n_k)_1^\infty$ is a sequence of positive integers diverging to $+\infty$ and that $(d_k)_1^\infty$ is a sequence of positive real numbers. 

Then there exist a subsequence $\(n_{k_l}\)_{l=1}^\infty$ diverging to $+\infty$, and a random variable $Y$ taking values in $[-\infty,\infty]$ such that
$$ 
\frac{1}{d_{k_l}}\sum_{j=0}^{n_{k_l}-1} f \circ T^j\xrightarrow{\  \  \,  {{\mathcal L}}\hspace*{0.4cm}}\lt(\int_Xf d\mu\rt)Y. 
$$
\ecor

\fr {\sl Proof.} Fix a set $A \in \mathfrak{F}$ such that $\mu(A)=1$. By Proposition~\ref{Proposition3.6.1}, there exist a subsequence $\(n_{k_l}\)_{l=1}^\infty$ diverging to $+\infty$, and a random variable $Y$ taking values in $[-\infty,\infty]$ such that
$$ 
\frac{1}{d_{k_l}}\sum_{j=0}^{n_{k_l}-1}\1_A \circ T^j\xrightarrow{\  \  \,  {{\mathcal L}}\hspace*{0.4cm}} Y. 
$$
Since also, by Hopf's Ergodic Theorem (Theorem~\ref{t2j73}), 
$$
\lim_{l\to+\infty}
\frac{\frac{1}{d_{k_l}}\sum_{j=0}^{n_{k_l}-1}f \circ T^j}
{\frac{1}{d_{k_l}}\sum_{j=0}^{n_{k_l}-1}\1_A \circ T^j}
=\frac{\int_Xf d\mu}{\int_X\1_A d\mu}
=\int_Xf d\mu
$$
$\mu$--a.e. on $X$, the statement of our corollary follows.
\qed

\sp 
\section{Darling--Kac Theorem; Abstract Version}\label{D-K AV}

In the Section~\ref{CountET}, Counterexamples to Ergodic Theorems, we have provided an evidence for any reasonable version of Birkhoff's Ergodic Theorem (i.e. pointwise) not to hold in the realm of infinite invariant measures. In the previous section, Section~\ref{WET}, Weak Ergodic Theorems, we provided some very weak counterparts of Birkhoff's Ergodic Theorem, namely, convergence in distribution along some subsequences of ``times''. In this section, we will strengthen the results of the previous one by stating, see Darling--Kac Theorems, Theorem~\ref{t1j163} and Theorem~\ref{t1j165}, due to \cite\cite{TZ}, the distributional convergence of weighted averages along the whole sequence of ``times'' and also some pointwise limsup results for such weighted averages. There is however a cost: we need to assume the existence of some special sets, namely: the sets $Y$ of Theorems~\ref{t1j163} and ~\ref{t1j165}) for distributional convergence and so called Darling--Kac sets for limsup results. The latter,
are particularly difficult for proving to exist, even in the simplest case of parabolic Markov maps of an interval. The former are easier provable to exist and we will do this for parabolic elliptic functions in Section~\ref{Darling_Kac-elliptic}. We would like to mention that in early papers 
of Denker and Aaronson (see \cite{Aa} and the references therein), the Darling--Kac sets were also used to prove the so called Pointwise Uniform Dual Ergodic Theorem. More recently, Ian Melbourne and Thresiu (see \cite{MT}) have provided fairly simple sufficient conditions, not requiring the existence of Darling--Kac sets for this theorem to hold. Because our book is anyway quite large, we will not (unfortunately) deal with Dual Ergodic Theorems any more in this book. We would like however to add that 
Darling--Kac Theorems have been also the objects of intensive refine research throughout several last decades, has attracted the attention of research (see for example \cite{Aa}, \cite{TZ}, and more) and working on then is  is still the process of progress. 

\sp A function $u : [a, +\infty) \lra  (0, +\infty)$ is called  regularly
varying of index\index{(N)}{regularly varying function of index $s$}
$s \in \mathbb R$ if $$ \lim_{t \to \infty}\frac{u(ct)}{u(t)}= c^s$$
for all  $c >0$. We  then write that $u \in \mathcal{R}_s$. Given a
sequence $(a_n)_1^\infty$ of positive numbers we  form the function
$\hat{a}:[1, +\infty) \lra (0, +\infty)$ as 
$$
\hat{a}(t)=a_{[t]}
$$ 
and we say that the sequence $(a_n)_1^\infty$ is regularly varying  of
index $s \in \mathbb R$ if the function $\hat{a}$ enjoys this
property. We then  also write that 
$$
(a_n)_1^\infty\in\mathcal{R}_s.
$$
The regularly  varying functions  and sequences of
index  $0$ are also referred to as slowly varying functions.

If $(X,\mathfrak{F}, \nu)$ is a probability  space and $f:X \lra
\mathbb R$ is an $\mathfrak F$-measurable function, following
probabilistic custom, we call $f$ a random variable and the
probability measure $\nu\circ f^{-1}$ on $\mathbb R$, the
distribution  of the random variable $f$. Frequently any Borel
probability measure on $\mathbb R$ is referred to as a distribution.
Notice that every distribution $\nu$ on $\mathbb  R$ is uniquely
determined  by its moments, i.e. the integrals
$$
\int_{-\infty}^{+\infty} x^k d\nu(x), \quad  k=1,2,3, \ldots. 
$$
Given $\a \in [0,1]$ the normalized Mittag--Leffler distribution
$\mathcal{M}_\a$\index{(S)}{$\mathcal{M}_\a$} of order $\a$ is
characterized by its moments
$$ \int_{-\infty}^{+\infty} x^k  d \mathcal{M}_\a(x)=
k!\frac{(\Gamma(1+\a))^k}{\Gamma(1+k \a)}, \quad  k=1,2, 3 \ldots
.$$ Given a Borel probability measure $\nu$ on $ \mathbb R$, we
denote by $a \nu +b$, $a>0$, the probability  distribution of the
random variable $aX_\nu+b$, where $X_\nu$ is any  random variable
whose distribution is equal to $\nu$, $i.e.$
   $$ ( a \nu +b)(F)=\nu\left(\frac{F-b}{a}\right)$$
for any real-valued $\nu$-integrable function $F: X \to \mathbb R$.

Every distribution of the form $a \mathfrak{M}_\a+ b$, $ a >0$, is
called a  Mittag--Leffler distribution.  Every  Mittag--Leffler
distribution \index{(N)}{Mittag--Leffler distribution} is
stable,\index{(N)}{stable Mittag--Leffler distribution} meaning that
the convolution of two Mittag--Leffler distribution of same order $
\a \in [0,1]$ is again a Mittag--Leffler distribution of order $\a$.
In other words, the distribution of the sum of two independent
random variables with Mittag--Leffler distribution is again  a
Mittag--Leffler distribution.

Now let $(X,\mathfrak{F}, \mu)$ be a $\sigma$--finite  measure space.
Assume that 
$$
\mu(X)=+\infty.
$$
Let $T:X\to X$ be an $\mathfrak {F}$--measurable map preserving measure $\mu$. For every  measurable
set $F$, the wandering rates\index{(N)}{wandering rates}  of $F$
are defined  to be
\beq\label{1201407}
w_n (F):= \mu \left( \bigcup_{k=0}^n T^{-k}(F)\right), \quad  n \geq 0.
\eeq\index{(S)}{$w_n (F)$}

\fr We shall  prove right now the following  lemma which establishes some
formula  for wandering rates which is useful in many calculations.

\sp\blem\label{l1j160} If  $(X,\mathfrak{F}, \mu)$ is a $\sigma$--finite
measure space and  $T:X\lra X$  is  an ergodic  conservative
measurable map preserving measure $\mu$, then  
$$ 
w_n(F)=\mu(F)\sum_{k=0}^n \mu_F( \tau_F  > k)
$$ 
for every  $ F \in \mathfrak{F}$ with  $ 0 <\mu(F)<  +\infty$, and all $n \geq 0$.
\elem

\bpf Since 
$$\bigcup_{k=0}^n T^{-k}(F)
= \bigcup_{k=0}^nT^{-k}(F)\cap \big\{x \in X:  \tau_F(T^k(x))> n-k\big\}
$$ 
and the constituents of the  latter union are mutually disjoint, we get that
$$
\aligned  
w_n(F)
= &\mu\left(\bigcup_{k=0}^n T^{-k}(F)\right)
= \mu\left( \bigcup_{k=0}^n T^{-k}(F)\cap\{x\in X:\tau_F(T^k(x))> n-k\}\right)\\
=& \sum_{k=0}^n \mu (T^{-k}(F)\cap \{ x \in X: \tau_F(T^k(x))>n-k\}) \\
=& \sum_{k=0}^n \mu (T^{-k}(F)\cap T^{-k}(\{ x \in X:\tau_F(x)>n-k\}))\\
= & \sum_{k=0}^n \mu (T^{-k}(F \cap \{ x \in X: \tau_F(x)> n-k\}))
\\= & \sum_{k=0}^n \mu (T^{-k}(\{ x \in F: \tau_F(x)>
n-k\}))\\
=& \sum_{k=0}^n \mu (\{ x \in F: \tau_F(x)>
n-k\})\\
=& \sum_{j=0}^n \mu (\{ x \in F: \tau_F(x)>
j\})\\
=& \sum_{j=0}^n \mu (\tau_F> j)\\
=& \mu(F) \sum_{j=0}^n \mu_F( \tau_F>j).
\endaligned$$
We are done. \endpf

\sp \fr We say that two number sequnces $(a_n)_1^\infty$ and
$(b_n)_1^\infty$ are asymtotically equivalent if 
$$
\lim_{n\to\infty}\frac{a_n}{b_n}=1.
$$
We then write that
$$
a_n\sim b_n.
$$
A measurable function $f : X \lra  [0, +\infty)$ is said to be supported on $F$ if $$
f^{-1}(0)\supset X \sms F.
$$ 
Furthermore, it is called uniformly sweeping\index{(N)}{uniformly sweeping function} on $F$ if there exists an integer $N\geq 0$ such that
$$
\inf_F\lt\{\sum_{n=0}^N \mathcal{L}_\mu^n(f)\rt\}>0,  
$$ 
where, we recall, 
$$
\mathcal{L}_\mu:L_\mu^1\lra L_\mu^1
$$ 
is the transfere operator of the map $T$ with respect to the measure $\mu$ defined in Section~\ref{QEC}. We shall now formulate the theorem  that is commonly refereed as  Darling--Kac Theorem. As we have already said this theorem is  in a sense a weak version  of
Theorem~\ref{Birkhoff} for (special) systems with infinite invariant measure, and it has a long history. It took on particularly 
simple, clear, and elegant form in \cite{TZ}. Thaler and Zweim\"uler
weakened the assumptions so much that these became  in  a sense local,
in particular getting rid of  dual  ergodicity, an annoying hypothesis present in Jon Aaronson's book \cite{Aa}. We now quote Theorem~3.1 from \cite{TZ} 

\bthm\label{t1j163} \rm{({Darling--Kac Theorem I})}
\index{(N)}{Darling-Kac Theorem I} Let $(X,\mathfrak{F}, \mu)$ be a
$\sigma$--finite  measure space and let $T:X \lra X$ be a measurable
map preserving measure $\mu$. Assume that there exists some set $Y
\in \mathfrak{F}$ with the following properties:
\begin{itemize}
\item [(a)]   $0< \mu(Y) < +\infty$,

\, 

\item  [(b)]  $(w_n(Y))_{n=1}^\infty \in \mathcal{ R}_{1-\a}$ with
some $\a \in [0,1]$,

\, 

\item [(c)] There exists a uniformly sweeping function $H:X \lra
 [0, +\infty)$ on $Y$ such that
$$
\frac{1}{w_n(Y)}\sum_{k=0}^{n-1} \mathcal{L}^k_\mu (
   \1_{Y^c_k}) \xrightarrow[\ n\to\infty \ ]{} 
 H,  
$$
uniformly on $Y$, where \index{(S)}{$Y^c_k$} 
$$
Y^c_k=Y^c_k(T):=T^{-k}(Y)\cap
\bigcap_{j=0}^{k-1}T^{-j}(X \sms Y)
$$ 
for $k \geq 1$ and $Y^c_0:=Y$.
\end{itemize}
Then for every function  $f\in L^1(\mu)$  with $ \int f
d\mu\neq 0$, the sequence $\lt(\frac{1}{a_n}S_nf \rt)_1^\infty$
converges strongly, with respect to the measure $\mu$, to the
Mittag--Leffler  distribution $(\int f d\mu)\mathcal{ M}_\a$; in symbols:
$$ 
\frac{1}{a_n}S_nf \xrightarrow{\  \  \,  {{\mathcal L}}\hspace*{0.4cm}}\lt(\int f d\mu\rt)\mathcal{ M}_\a
$$
where
$$
a_n:=\frac{1}{\mu(Y)}\int_Y S_n(\1_Y)\, d\mu_Y \sim
 \frac{1}{\Gamma(1+\a)\Gamma(2 - \a)}\cdot \frac{n}{w_n(Y)}\quad
 \mbox {as}  \quad n \lra \infty.
$$
\ethm

\

\fr A proof of this theorem can be found in \cite{TZ}. This paper contains also some other stochastic asymptotic laws (like arcsin law), relevant in the context of infinite measures.

\sp\fr For every $k \geq 1$, let
\beq\label{5j164}
Y_k=Y_k(T):=\tau_Y^{-1}(k)=\{x\in Y:\, \tau_Y(x)=k\}.
\eeq\index{(S)}{$Y_k$}
We shall prove  the following, quite often relatively easily
verifiable, sufficient condition for the requirement (c) of
Theorem~\ref{t1j163} to hold. 

\sp\blem\label{l1j164} 
Let $(X,\mathfrak{F},\mu)$ be a $\sigma$--finite measure space and let $T:X \lra X$ be a measurable map preserving measure $\mu$. If there exists a (measurable) function $\hat{H}:X\to[0, +\infty)$, supported on some set $Y\in\mathfrak{F}$, such that the sequence 
$$
\left(\frac{1}{\mu_Y(Y_k)}\mathcal{L}_\mu^k(\1_{Y_k})\right)_{k=1}^\infty
$$
converges uniformly on $Y$ to $\hat{H}|_Y$, then condition (c) of Theorem~\ref{t1j163} holds with the function 
$$
H:=(\mu(Y))^{-1}\hat{H}:X\lra[0,+\infty),
$$
although we  do not claim that $H:X\to[0,+\infty)$ is uniformly sweeping. Obviously however, if $\hat{H}$ is uniformly sweeping, then so is $H$. 
\elem

\bpf By virtue of (\ref{1j157}) we  have for every  set
$ F \in \mathfrak {F}$  and any  $g \in L^{1}(\mu)$, that
$$  \int g \1_F d \mu= \int g \circ  T \cdot \1_F\circ T d \mu=
\int g\circ T\cdot\1_{T^{-1}(F)} d \mu= \int g \mathcal {L}_\mu(
\1_{T^{-1}(F)})d \mu.$$ This implies  that
\beq\label{1j164}
 \1_F=\mathcal{L}_\mu
(\1_{T^{-1}(F)})\quad \mu \, a.e.
\eeq
Hence, for every  $n \geq 0$, we have $\mu$--$a.e.$, that 
\beq\label{3j164}
\1_{Y^c_n}
=\mathcal{L}_\mu(\1_{T^{-1}(Y^c_n)})
=\mathcal{L}_\mu(\1_{Y^c_{n+1}\cup Y_{n+1} })
= \mathcal{L}_\mu(\1_{Y_{n+1}}+ \1_{Y_{n+1}^c }) 
= \mathcal{L}_\mu(\1_{Y_{n+1}})+ \mathcal {L}_\mu(\1_{Y_{n+1}^c }).
\eeq
By Remark~\ref{r1j64a.1}, we have that
\beq\label{2j164}
\lim_{n\to \infty} \mu(Y^c_n)=0.
\eeq
Since $Y^c_{n+1}\sbt T^{-1}(Y^c_n)$ and since $\mathcal{L}_\mu$ is a positive operator, we get.
\beq\lab{1_2017_11_27}
\mathcal{L}_\mu( \1_{Y^c_{n+1}}) 
\leq\mathcal{L}_\mu( \1_{T^{-1}(Y^c_n)}).
\eeq 
It follows, by immediate induction, from (\ref{3j164}) that
\beq\label{1j164.a}
\1_{Y^c_n}= \sum_{j=1}^k \mathcal{L}_\mu^j( \1_{Y_{n+j}})+
\mathcal{L}_\mu^k(\1_{Y_{n+k}^c})
\eeq
for all integers $k\ge 1$. By  \eqref{1_2017_11_27} and (\ref{1j164}),  we have that
$$ 
\mathcal{L}_\mu^{k+1}(\1_{Y^c_{n+k+1}}) \leq
  \mathcal{L}_\mu^{k+1}(\1_{T^{-1}(Y^c_{n+k})})
  =\mathcal{L}_\mu^k(\mathcal{L}_\mu(\1_{T^{-1}(Y^c_{n+k})}))
  =\mathcal{L}_\mu^{k}(\1_{T^{-1}(Y^c_{n+k+1})}),
$$
meaning that the sequence $
(\mathcal{L}_\mu^k(\1_{Y^c_{n+k}}))_{k=1}^\infty$ is decreasing. Let
$$
g:=\lim_{k \to \infty}\mathcal{L}_{\mu}^k( \1_{Y^c_{n+k}}) \geq 0.
$$
By Lebesgue's Monotone Convergence Theorem and \eqref{2j164}, we thus have
$$
\int g \,d\mu
= \lim_{k\to \infty}
 \int \mathcal{L}_\mu^k(\1_{Y^c_{n+k}})\, d\mu 
= \lim_{k \to \infty} \int \1_{Y^c_{n+k}} d\mu
= \lim_{k \to \infty} \mu(\1_{Y^c_{n+k}})=0.
$$ Hence $g=0$ $\mu
\, a.e.$, and (\ref{1j164.a}) yields
$$\1_{Y^c_n}=\sum_{j=1}^\infty \mathcal{L}_\mu^j(\1_{Y^c_{n+j}})$$
$\mu$--$a.e.$ on $X$. Therefore,
\beq\label{2j164.a}
\frac{1}{\mu_Y(\tau^{-1}_Y([n+1,
\infty)))}\mathcal{L}_\mu^n(\1_{Y^c_n})= \frac{1}{\mu_Y(\tau_Y>k)}
\sum_{j=1}^\infty \mathcal{L}_\mu^{n+j}(\1_{Y_{n+j}}).
\eeq
Put  
$$\hat{H}
:=\lim_{n \to \infty}\frac{1}{\mu (Y_n)}
\mathcal{L}_\mu^n (\1_{Y_n}):X\lra [0,+\infty).
$$ 
Then, for every $\varepsilon >0$ there exists $n_\varepsilon \geq 1$ such that for all $ n \geq n_\varepsilon$
$$\hat{H}(x)-\varepsilon \leq \frac{1}{
\mu_Y(Y_{n+j})}\mathcal{L}_\mu^{n+j}(\1_{Y_{n+j}})(x)\leq
\hat{H}(x)+ \varepsilon$$ for all $j \geq 1$ and all $x \in Y$.
Equivalently,
$$(\hat{H}(x)-\varepsilon)\mu_Y(Y_{n+j})\leq
\mathcal{L}_\mu^{n+j}(\1_{Y_{n+j}})(x)\leq
(\hat{H}(x)+\varepsilon)\mu_Y(Y_{n+j}).
$$
Summing up over all $j=1,2, \ldots,$ we thus obtain,
$$(\hat{H}(x)-\varepsilon)\mu_Y(\tau_Y>n )\leq \sum_{j=1}^\infty
\mathcal{L}_\mu^{n+j}(\1_{Y_{n+j}})(x)\leq
(\hat{H}(x)+\varepsilon)\mu_Y(\tau_Y>n).
$$
In  virtue of (\ref{2j164.a}) this means  that
$$({H}(x)-\varepsilon)\mu_Y(\tau_Y>n )\leq
\mathcal{L}_\mu^{n}(\1_{Y_n^c})(x)\leq
({H}(x)+\varepsilon)\mu_Y(\tau_Y>n).
$$
for all $n \geq n_\varepsilon$ and all $x \in Y$. Hence, fixing $ l
\geq n_\varepsilon+1$ and summing up over all $n=n_\varepsilon,
n_\varepsilon+1, \ldots,l-1$, we  get for all $x\in Y$ that 
$$
(\hat{H}(x)-\varepsilon)\sum_{n=n_\varepsilon}^{l-1}\mu_Y(\tau_Y>n )\leq
\sum_{n=n_\varepsilon}^{l-1}\mathcal{L}_\mu^{n}(\1_{Y_{n}^c})(x)\leq
(\hat{H}(x)+\varepsilon) \sum_{n=n_\varepsilon}^{l-1}
\mu_Y(\tau_Y>n).
$$
Dividing by $w_l(Y)$, and using Lemma~\ref{l1j160}, we  then get
$$
\aligned \frac{(\hat{H}(x)-\varepsilon)}{\mu(Y)}\cdot \frac{\sum_{n=n_\varepsilon}^{l-1}\mu_Y(\tau_Y>n
)}{\sum_{n=0}^{l-1}\mu_Y(\tau_Y>n )}
 &\leq
\frac{1}{w_l(Y)}\sum_{n=n_\varepsilon}^{l-1}\mathcal{L}_\mu^{n}(\1_{Y_n^c})(x)
\cdot
\frac{\sum_{n=n_\varepsilon}^{l-1}\mathcal{L}_\mu^{n}(\1_{Y_n^c})(x)}{\sum_{n=0}^{l-1}\mathcal{L}_\mu(\1_{Y_n^c})(x)}\\
&\leq \frac{(\hat{H}(x)+\varepsilon)}{\mu(Y)}\cdot
\frac{\sum_{n=n_\varepsilon}^{l-1}\mu_Y(\tau_Y>n
)}{\sum_{n=0}^{l-1}\mu_Y(\tau_Y>n )}.
\endaligned$$
Since $\lim_{l \to \infty}w_l(Y)=+\infty$, taking $ l \geq
l_\varepsilon\geq n_\varepsilon$ large  enough, we  will have
$$1-\varepsilon  \leq \frac{\sum_{n=n_\varepsilon}^{l-1}\mu_Y(\tau_Y>n
)}{\sum_{n=0}^{l-1}\mu_Y(\tau_Y>n )}\leq 1
$$
and
$$ 1 -\varepsilon   \leq \frac{\sum_{n=n_\varepsilon}^{l-1}\mathcal{L}_\mu^{n}(\1_{Y_n^c})(x)}{\sum_{n=0}^{l-1}\mathcal{L}_\mu^{n}(\1_{Y_n^c})(x)}
\leq 1$$ for all $x \in Y$. Thus,
$$(1-\varepsilon)\frac{(\hat{H}(x)-\varepsilon)}{\mu(Y)}\leq
\frac{1}{w_l(Y)}\sum_{n=0}^{l-1} \mathcal{L}_\mu^{n}(\1_{Y_n^c})(x)
\leq\frac{1}{(1-\varepsilon)}\frac{(\hat{H}(x)+\varepsilon)}{\mu(Y)}
$$
for all $x \in X$. Hence  the sequence
$$ \left( \frac{1}{w_l(Y)}\sum_{n=0}^{l-1}\mathcal{L}_\mu^{n}(\1_{Y_n^c})(x)
                        \right)_{l=1}^\infty$$
converges to $(\mu(Y))^{-1}\hat{H}|_Y $ uniformly on $Y$. We are done.
\endpf

\sp The Darling--Kac Theorem also holds  if the hypotheses of
Theorem~\ref{t1j163} are verified for some iterate  of $T$ only.
Indeed, we have the following. 

\bthm\label{t1j165} \rm{(Darling-Kac Theorem
II)}\index{(N)}{Darling-Kac Theorem II} Let $(X, \mathfrak{F},\mu)$
be a  $\sigma$-finite  measure space and let $T:X \to X$ be a
measurable  map  preserving the measure  $\mu$. Assume that there exists some integer $q \geq 1$ and some set $Y \in \mathfrak{F}$ and with the following  properties:
\begin{itemize}
\item [(a)]   $0< \mu(Y) < +\infty$,

\, 

\item  [(b)]  $(w_n(T^q,Y))_{n=1}^\infty \in \mathcal{ R}_{1-\a}$ with
some $\a \in [0,1]$,

\, 

\item [(c)] There exists a function $H:X \lra [0, +\infty)$, uniformly  sweeping on $Y$ with respect to $T^q$, such that
$$  
\frac{1}{w_n(T^q,Y)}\sum_{k=0}^{n-1} \mathcal{L}^{qk}_\mu (
   \1_{Y^c_k}(T^q))\xrightarrow[\ n\to\infty \ ]{} 
 H|_Y
$$
uniformly on $Y$, where 
$$
Y^c_k(T^q):=T^{-kq}(Y)\cap
\bigcap_{j=0}^{k-1}T^{-qj}(X \sms Y)
$$ 
for  $k \geq 1$ and $Y^c_0(T^q):=Y$.
\end{itemize}
Then  for every function  $ f \in L^1(\mu)$  with $ \int f
d\mu\neq 0$, the sequence $\lt(\frac{1}{b_n}S_nf \rt)_1^\infty$
converges strongly, with respect to the measure $\mu$, to the
Mittag--Leffler  distribution $(\int f d\mu)\mathcal{ M}_\a$; in symbols:
$$ 
\frac{1}{a_n}S_nf \xrightarrow{\  \  \,  {{\mathcal L}}\hspace*{0.4cm}}\lt(\int f d\mu\rt)\mathcal{ M}_\a, 
$$
where
\beq\lab{2_2017_11_28}
b_n:=\frac{q}{\mu(Y)}\int_Y S_{E(n/q)}(\1_Y) d \mu_Y \sim
 \frac{1}{\Gamma(1+\a)\Gamma(2 - \a)}\cdot \frac{n}{w_{E(n/q)}(T^q,Y)}.
\eeq
\ethm

 \bpf Let
$$a_n=\frac{1}{\mu(Y)}\int_YS_n^{T^q}(\1_Y)\sim\frac{1}{\Gamma(1+\a)
\Gamma(2-\a)}\cdot \frac{n}{w_n(T^q, A)}.$$ Theorem~\ref{t1j163}
tells us that the sequence $\frac{1}{a_n}S^{T^q}_n( S_qf)$ converges
strongly with respect to the  measure $\mu$, to the Mittag-Leffler
distribution $(\int S_q (f)\, d\mu) \mathcal{M}_\a$. This means  that
for  every probability  measure $\nu$ absolutely   continuous with
respect to  $\mu$, and every bounded  uniformly continuous  function
$\phi: \mathbb R \to  \mathbb R$, we have that
$$  
\lim_{n \to \infty}\int_X \phi \circ \left( \frac{1}{a_n}S^{T^q}_n(
S_q(f))\right) d \nu
=  \int_X \phi\, d\(\mu (S_q(f)) {\mathcal M}_\a\).
$$ Equivalently,
$$  
\lim_{n \to \infty}\int_X \phi \circ \left( \frac{1}{a_n}
S_{qn}(f)\right) d \nu
=  \int_X \phi \, d\(\mu (S_q(f)) {\mathcal M}_\a\).
$$ 
Writing this equality for the function $\phi\circ\frac{1}{q}$, which is also bounded and uniformly continuous, rather than $\phi$, we get
\beq\label{1j166}
\begin{aligned} 
\lim_{n \to \infty}\int_X \phi\lt(\frac{1}{q a_n}S_{qn}(f)\rt)d \nu
&=\lim_{n \to \infty}\int_X \lt(\phi\circ\frac{1}{q}\rt)\circ \lt(\frac{1}{qa_n}S_{qn}(f)\rt)d\nu \\
&= \int_X\lt(\phi\circ\frac{1}{q}\rt)d\(\mu(S_q(f){\mathcal M}_\a\) \\
&=\int_X \phi(x)\, d\(q\mu(f){\mathcal M}_\a(qx)\)
=\int_X \phi\, d\lt(q\mu(f)\frac{1}{q} {\mathcal M}_\a\rt) \\
&= \int_X \phi\, d\(\mu( f) {\mathcal M}_\a\). 
\end{aligned}.
\eeq
Now, for every integer $n\geq q$ write uniquely $n=k_nq +r_n$, where
with integers $k \geq 1$ and $0 \leq r_n \leq q-1$. Note that
$k_n=E(n/q)$ and
$$
\left\|\frac{1}{q a_{E(n/q)}}S_n(f) -
\frac{1}{qa_{k_n}}S_{qk_n}(f)\right\|_\infty 
=\left\|\frac{1}{qa_{k_n}}S_{r_n}(f\circ T^{qk_n})\right\|_\infty
\leq  \frac{r_n ||f||_\infty}{q a_{k_n}}
\leq \frac{\|f\|_\infty}{a_{k_n}}
= \frac{\|f\|_\infty}{a_{E(n/q)}}.
$$
As $\lim_{j \to \infty} a_j =\infty$, it follows from this formula, applied in
(\ref{1j166}), that
 $$ \lim_{n \to \infty} \int_X \phi\circ \left(\frac{1}{q
 a_{E(n/q)}}S_n(f)\right)d \nu
  =\int_X \phi\,  d ( \mu (f)){\mathcal
 M}_\a).$$
This precisely means that the sequence $\frac{1}{b_n}S_n(f)$, $b_n=q
a_{E(n/q)}$, converges strongly to  the Mittag-Leffler distribution
$\mu(f) {\mathcal M }_\a$. We are done.
\endpf

\sp\fr This is precisely the form of Darling--Kac Theorem which we will
prove in section~\ref{Darling_Kac-elliptic} for some classes of parabolic elliptic functions. Lemma~\ref{l1j164} now takes on the following form.

\sp\blem~\label{l1j168} 
Let $(X,\mathfrak{F},\mu)$ be a $\sigma$--finite  measure space and let $T:X \lra X$ be a measurable map preserving measure $\mu$. 
If there exists a (measurable) function $\hat{H}:X\to[0, +\infty)$, supported on some set $Y\in \mathfrak{F}$, such that the sequence
$$ 
\left(\frac{1}{\mu_Y(Y_k(f^q))}\mathcal{L}^{qk}_\mu(\1_{Y_k(f^q)
})\right)^\infty_{k=1}
$$ 
converges uniformly on $Y$ to $\hat H|_Y$, then condition (c) in Theorem~\ref{t1j165} holds with the function $H:=(\mu(Y))^{-1}\hat{H}:X\to[0,+\infty)$, although we do not claim that $H:X\to[0,+\infty)$ is uniformly sweeping. Obviously however, if $\hat{H}$ is uniformly sweeping, then so is $H$. 
\elem

\sp We end this section by stating without proof a weak pointwise version of ergodic theorem in the realm of infinite invariant measures. The existence of Darling--Kac sets is needed now, and also some kind of mixing. We say that $A\in \mathfrak{F}$ is called a Darling--Kac set if and only if $0<\mu(A)<+\infty$ and there  exists a sequence $(a_n)_{n=1}$ of positive reals such that 
$$
\lim_{n\to\infty}\frac{1}{a_n}\sum_{j=0}^{n-1}\mathcal{L}^{j}_\mu(\1_A)
=\mu(A)
$$
uniformly on $A$.

\bthm\label{Theorem 4.} 
Let $(X,\mathfrak{F},\mu)$ be a $\sigma$--finite measure space and let $T:X \lra X$ be a measurable map preserving measure $\mu$. Assume that a Darling--Kac set exists. 
Suppose that $(\varphi(n))_{n=1}^\infty$ is an increasing sequence of positive real numbers such that the sequence $\lt(\frac{\varphi(n)}{n}\rt)_{n=1}^\infty$ is decreasing. Then:
\begin{itemize}
\item [(a)] If \, $ \sum_{n=1}^\infty \frac{1}{n}\exp(-\beta \varphi(n))< +\infty$ for all $\beta >1$, then
$$
\limsup_{n\to \infty} \frac{1}{ n^{\alpha}h\lt(\frac{n}{\varphi(n)}\rt)\varphi(n)^{1-\alpha}}\sum_{k=1}^n f\circ T^k\leq  K_\alpha\int_Xf d \mu
$$  
$\mu$--a.e. for every  $ f \in L^1_{+}$.

\

\item [(b)] If \, $\sum_{n=1}^\infty \frac{1}{n}\exp(-\b\varphi(n))=+\infty$ for all $\b<1$, then
$$
\limsup_{n\to \infty} \frac{1}{ n^{\alpha}h\lt(\frac{n}{\varphi(n)}\rt)\varphi(n)^{1-\alpha}}\sum_{k=1}^n f\circ T^k\geq  K_\alpha\int_Xf d \mu
$$   
$\mu$--a.e. for every $ f \in L^1_{+}$.
  
\

\item [(c)]
$$
\limsup_{n\to \infty} \frac{1}{ n^{\alpha}h\lt(\frac{n}{L_2(n)}\rt)L_2(n)^{1-\alpha}}\sum_{k=1}^n f\circ T^k= K_\alpha\int_Xf d \mu
$$ 
$\mu$--a.e. for every $ f \in L^1_{+}$.
\end{itemize}
\ethm

\sp\section{Points of Infinite Condensation; Abstract
Setting}\label{abstract}

\sp In this very short section we provide a framework for dealing with
infinite invariant measures in a topological setting. We pay
special attention to the concept  of points  of infinite
condensation introduced in \cite{U2}. So, assume  that $X$ is a
separable locally compact metrizable space and that 
$$
T:X \lra \-X
$$  
is a continuous map,  where 
$$
\-X:=X\cup\{ \infty\}
$$ 
is the Alexander
(one--point) compactification of $X$. We want to emphasize that we by no means  assume that $T:X\to \-X$ to be extendable to a
continuous  map  from $\-X$ to $\-X$. We  however set $T(\infty)=
\infty$. 

\sp Given a Borel $\sg$--finite measure $\mu$ on $X$ let $\hat
X_\mu(\infty)$ be the set of all points $x\in X$ such that
$\mu(U)=+\infty$ for every open set $U$ containing $x$. Following
\cite{U2} we call ${X}_\mu(\infty)$\index{(S)}{$X_\mu(\infty)$} the
set of points of infinite condensation\index{(N)}{points of infinite
condensation} of $\mu$. Notice that $X_\mu(\infty)$ is a closed
subset of $\hat X$ and ${X}_\mu(\infty)=\es$ if and only if the
measure $\mu$ is finite. 

Denote by $\cM_T^\infty$\index{(S)}{$\cM_T^\infty$} the family of all Borel
ergodic conservative $T$--invariant measures $\mu$ for which
\beq\label{120190925}
\mu(X\sms X_\mu(\infty))>0.
\eeq
Note that if $\mu\in\cM_T^\infty$,
then $X_\mu(\infty)$ is forward invariant, i.e.
\beq\label{220190925}
T(X_\mu(\infty))\sbt X_\mu(\infty).
\eeq
The following simple proposition shows that $X_\mu(\infty)$ is measurably negligible and, in the topologically transitive case, it is also topologically negligible.

\sp\bprop\label{p1w3} 
If $X$ is a locally compact separable metric space, if $T:X\lra \hat{X}$ is a continuous map, and if
$\mu\in\cM_T^\infty$, then $\mu(X_\mu(\infty))=0$. If, in addition,
$T:X\lra \hat{X}$ is topologically transitive, then $X_\mu(\infty)$
is a nowhere dense subset of $X$. 
\eprop

\bpf Proving the first assertion of this proposition, suppose for the contrary that
$$
\mu(X_\mu(\infty))\not=0.
$$
Then, by ergodicity and conservativity of $\mu$, we have that
$$
\mu\Big(X\sms\bu_{n=0}^\infty T^n(X_\mu(\infty))\Big)=0.
$$
Along with \eqref{220190925}, this implies that 
$$
\mu\(X\sms X_\mu(\infty))\)=0.
$$
This contradicts formula \eqref{120190925} and finishes the proof of the first assertion of the proposition.

Suppose now that the map $T:X\to X$ is topologically transitive and seeking contradiction assume that the set $X_\mu(\infty)$ is
not nowhere dense. This means that $\Int\(X_\mu(\infty)\)\ne\es$. Hence, $\Int\(X_\mu(\infty)\)$ contains a transitive point and we conclude that  
$$
\bu_{n=0}^\infty T^n(X_\mu(\infty))
$$ 
is a dense subset of $X$. Since however
$$
\ov{\bu_{n=0}^\infty T^n(X_\mu(\infty))}\sbt X_\mu(\infty)
$$ 
and since $X\sms X_\mu(\infty)\ne\es$ (by \eqref{120190925}),
we get a contradiction which finishes the proof. 
\endpf

\sp\bprop\label{p1w28} 
If $T:X\lra \hat{X}$ is a continuous map, $X$ is
 a locally compact separable metric space, and if $\mu\in\cM_T^\infty$,
then $\mu$ is $\sg$-finite. 
\eprop

\bpf By the definition of $X_\mu(\infty)$, for every
point $x\in X\sms X_\mu(\infty)$ there exists an open ball $B_x\sbt
X\sms X_\mu(\infty)$ such that $x\in B_x$ and 
$$
\mu(B_x)<+\infty.
$$
Since $X\sms X_\mu(\infty)$ is a separable metric space,
Lindel\"of's theorem implies that there exists a countable set
$Y\sbt X$ such that 
$$
\bu_{y\in Y}B_y=X\sms X_\mu(\infty).
$$
Thus 
$$
\{X_\mu(\infty)\}\cup \{B_y\}_{y\in Y}
$$ 
forms a countable cover of $X$ by Borel sets of finite measure as $\mu(X_\mu(\infty))=0$ by Proposition~\ref{p1w3}. 
\endpf

\chapter{Measure--Theoretic Entropy}\label{ch8}

We shall deal in this chapter with measure--theoretic entropy of a measurable transformation preserving a probability space. Measure--theoretic
entropy is also sometimes known as metric entropy or  Kolmogorov–-Sinai metric entropy. It was introduced by A.~ Kolmogorov and Ya.~Sinai in late 1950's; see \cite{Si}. Since then its account has been presented in virtually every textbook on ergodic theory. Its introduction to dynamical systems was motivated by the concept of Ludwig Boltzmann entropy of statistical mechanics and Claude Shannon work on information theory; see \cite{Sh1} and \cite{Sh2}. 

We will encounter
three stages in the definition of metric entropy. It is defined by
partitioning the underlying measurable space with measurable sets. Indeed, whereas one cannot
generally partition a topological space into open sets, it is no problem to partition a measurable space into measurable sets.

\sp\section{Partitions}\label{Partitions}

\fr Let $(X,\mathfrak{ F})$ be a measurable space, and let $\mathcal{A}$ be a countable measurable partition \index{(N)}{measurable partition} of $X$, that is,
$\mathcal{A}=\{A_k\}_{k\geq1}$ with each $A_k\in\mathfrak{ F}$ such that

\sp\begin{itemize} 
\item $A_i\cap A_j=\emptyset$ 
for all $i\neq j$ 

\sp\fr and 

\sp\item $\bigcup_{k\ge1}A_k=X$. 
\end{itemize}

\sp\fr For each $x\in X$ denote by
$\mathcal{A}(x)$ the unique element (atom) of the partition $\mathcal{A}$ that
contains the point $x$. In the remainder of this chapter, we shall use
the calligraphic letters  $\mathcal{A},\mathcal{B},\mathcal{C},\ldots$
to denote partitions, with the notable exception of $\mathfrak{ F}$ which
will denote a 
$\sigma$--algebra on the space $X$. 

Also, if $X$ is a metrizable topological space, then, as usually, $\mathcal{B}(X)$
denotes the Borel $\sigma$--algebra on $X$. 

We shall denote the set of all countable
measurable partitions on the space $(X,\mathfrak{ F})$ by
$\mathrm{Part}(X,\mathfrak{ F})$. 
Moreover, it  will always be implicitly understood that partitions
are countable (finite or infinite) and measurable.

\sp\bdfn
Let $(X,\mathfrak{ F})$ be a measurable space and
$\mathcal{A},\mathcal{B}\in\mathrm{Part}(X,\mathfrak{ F})$.
We say that partition $\mathcal{B}$ is {\em finer} than
partition $\mathcal{A}$, or that $\mathcal{A}$ is {\em coarser} \index{(N)}{finer partition} \index{(N)}{coarser partition} than
$\mathcal{B}$, which will be denoted by $\mathcal{A}\leq\mathcal{B}$, if for
every atom $B\in\mathcal{B}$ there exists some atom $A\in\mathcal{A}$
such that $B\subseteq A$.
\edfn

\

\fr Equivalently, $\mathcal{B}$ is a refinement of $\mathcal{A}$ if
$\mathcal{B}(x)\subseteq\mathcal{A}(x)$ for all $x\in X$.
We now introduce for partitions the analogue of the join of two
covers.

\

\bdfn
Given $\mathcal{A},\mathcal{B}\in\mathrm{Part}(X,\mathfrak{ F})$, the partition
\[
\mathcal{A}\vee\mathcal{B}:=\{A\cap B:A\in\mathcal{A}, B\in\mathcal{B}\}
\]
is called the {\em join} of $\mathcal{A}$ and $\mathcal{B}$. \index{(N)}{join of partitions}
\edfn

\

\fr The basic properties of the join are given in  the following
lemma. Their proofs are left to the 
reader as an exercise.

\

\blem\label{lempartord}
Let 
$\mathcal{A},\mathcal{B},\mathcal{C},
\mathcal{D}\in\mathrm{Part}(X,\mathfrak{ F})$.    
Then
\begin{itemize}
\item[(a)] $\mathcal{A}\vee\{X\}=\mathcal{A}$;

\,

\item[(b)] $\mathcal{A}\vee\mathcal{B}=\mathcal{B}\vee\mathcal{A}$;

\,

\item[(c)] $\mathcal{A}\leq\mathcal{A}\vee\mathcal{B}$
            and $\mathcal{B}\leq\mathcal{A}\vee\mathcal{B}$;
            
\,

\item[(d)] $\mathcal{A}\leq\mathcal{B}$ if and only if
            $\mathcal{A}\vee\mathcal{B}=\mathcal{B}$;
            
\,

\item[(e)] If $\mathcal{A}\leq\mathcal{C}$ and $\mathcal{B}\leq\mathcal{D}$,
             then $\mathcal{A}\vee\mathcal{B}\leq\mathcal{C}\vee\mathcal{D}$.
\end{itemize}
\elem

\sp\section{Information and Conditional Information Functions}

\sp\ Let $(X,\mathfrak{ F},\mu)$ be a probability space. As such, the set $X$ may be
construed as the set of all possible states (or outcomes) of an experiment, while
the \mbox{$\sigma$-algebra} $\mathfrak{ F}$ consists of the set of all
possible events, and 
$\mu(E)$ is the probability that event $E\in\mathfrak{ F}$ take place.
Imagine that this experiment is conducted using an instrument which due to some limitation
can only provide measurements accurate up to the atoms 
of a partition $\mathcal{A}=\{A_k\}_{k\geq1}\in\mathrm{Part}(X,\mathfrak{ F})$.
In other words, this instrument can
only tell us which atom of $\mathcal{A}$ the outcome of the experiment
falls into. Any observation made through this instrument will
therefore be of the form 
$A_k$ for a unique $k$. If the experiment were
conducted today, the probability that its outcome belongs to $A_k$, i.e., the
probability that the experiment result in observation $A_k$ being made with our
instrument, would be given by $\mu(A_k)$.

\sp We would like to introduce a function that describes the
information that our instrument 
would give us about the outcome of the experiment. So, let $x\in
X$. Intuitively, the smaller the atom 
of the partition to which $x$ belongs, the more information our
instrument provides us about $x$. 
In particular, if $x$ lies in an atom of full measure, then our
instrument gives us essentially no 
information about $x$. Moreover, because our instrument does not
distinguish points which belong to 
a common atom of the partition, the sought information function must
be constant on every atom. 

\

\bdfn
Let $(X,\mathfrak{ F},\mu)$ be a probability space and
$\mathcal{A}\in\mathrm{Part}(X,\mathfrak{ F})$. 
The function $I_\mu(\mathcal{A}):X\lra[0,\infty]$ defined by
\[
I_\mu(\mathcal{A})(x):=-\log\mu(\mathcal{A}(x))
\]
is called the {\em information function} of the partition $\mathcal{A}$. \index{(N)}{information function of a partition}\index{(N)}{$I_\mu(\mathcal{A})$}
\edfn

\sp  As the function $t\mapsto-\log t$ is a decreasing function, for any $x\in X$ the smaller $\mu(\mathcal{A}(x))$ is, the larger
$I_\mu(\mathcal{A})(x)$ is, that is,  the smaller 
the measure of the atom $\mathcal{A}(x)$ is, the more information the
partition $\mathcal{A}$ gives us about $x$.
In particular, the finer the partition, the more information it gives
us about every point in the space.

\sp We now enumerate some of the basic properties of the
information function. Their proofs are straightforward and are again
left to the reader. 

\sp\blem
Let $(X,\mathfrak{ F},\mu)$ be a probability space and
$\mbox{Meas}(X,\mathfrak{ F})$ be the set 
of all measurable functions on $(X,\mathfrak{ F})$.
\begin{itemize}
\item[(a)] The map
\[
\begin{array}{cccc}
I_\mu:&\mathrm{Part}(X,\mathfrak{ F})&\lra&\mbox{Meas}(X,\mathfrak{ F})\\
      &   \mathcal{A}&\longmapsto&I_\mu(\mathcal{A})
\end{array}
\]
is an increasing function. In other words, if
$\mathcal{A}\leq\mathcal{B}$ then $I_\mu(\mathcal{A})\leq I_\mu(\mathcal{B})$.

\,

\item[(b)] $I_\mu(\mathcal{A})(x)=0$ if and only if $\mu(\mathcal{A}(x))=1$.

\,

\item[(c)] $I_\mu(\mathcal{A})(x)=\infty$ if and only if
  $\mu(\mathcal{A}(x))=0$. 
  
\,

\item[(d)] $I_\mu(\mathcal{A})(x)=I_\mu(\mathcal{A})(y)$ if
  $\mathcal{A}(x)=\mathcal{A}(y)$, that is, $I_\mu(\mathcal{A})$ is
  constant over each atom of $\mathcal{A}$. 
\end{itemize}
\elem

\sp More advanced properties of the information function will be
presented below. Meanwhile, we introduce a function which describes
the information given by a 
partition  $\mathcal{A}$ given that a partition  $\mathcal{B}$  
has already been applied. 

\sp\bdfn
The {\em conditional information function} of a partition $\mathcal{A}$ given a partition $\mathcal{B}$ is defined by \index{(N)}{conditional information function of a partition}\index{(S)}{$I_\mu(\mathcal{A}|\mathcal{B})$}
\[
I_\mu(\mathcal{A}|\mathcal{B})(x)
:=-\log\mu_{\mathcal{B}(x)}(\mathcal{A}(x)).
\]
Note that
\begin{eqnarray*}
I_\mu(\mathcal{A}|\mathcal{B})(x)
&=&-\log\frac{\mu\bigl(\mathcal{A}(x)\cap\mathcal{B}(x)\bigr)}{\mu(\mathcal{B}(x))}  
=-\log\frac{\mu\bigl((\mathcal{A}\vee\mathcal{B})(x)\bigr)}{\mu(\mathcal{B}(x))}
\\&=&I_\mu(\mathcal{A}\vee\mathcal{B})(x)-I_\mu(\mathcal{B})(x).
\end{eqnarray*}
It is implicitly understood that $\frac{0}{0}=0$ and
$\infty-\infty=\infty$ (why?). 
\edfn

\sp For any partition $\mathcal{A}$, observe that
$$
I_\mu(\mathcal{A}|\{X\})=I_\mu(\mathcal{A}),
$$
i.e., the information function coincides with the conditional information function with respect to the trivial partition 
$\{X\}$. Note further that $I_\mu(\mathcal{A}|\mathcal{B})$ is
constant over each atom of 
$\mathcal{A}\vee\mathcal{B}$.

\sp We shall now provide some advanced properties of the conditional
information function. 
Note that some of these properties hold pointwise, while others hold
atomwise only, that is, 
after  integrating over atoms. In particular, the reader should
compare statements $(viii)$ and $(ix)$. 

\sp\bthm\label{sem2thm7.0}
Let $(X,\mathfrak{ F},\mu)$ be a probability space and
$\mathcal{A},\mathcal{B},\mathcal{C} 
\in\mathrm{Part}(X,\mathfrak{ F})$. The following statements hold:
\begin{itemize}
  \item[(1)]
    $I_\mu(\mathcal{A}\vee\mathcal{B}|\mathcal{C})
=I_\mu(\mathcal{A}|\mathcal{C})  
      +I_\mu(\mathcal{B}|\mathcal{A}\vee\mathcal{C})$. 
      
\,
  
  \item[(2)]
    $I_\mu(\mathcal{A}\vee\mathcal{B})
=I_\mu(\mathcal{A})+I_\mu(\mathcal{B}|\mathcal{A})$. 

\,

  \item[(3)] If $\mathcal{A}\leq\mathcal{B}$, then
    $I_\mu(\mathcal{A}|\mathcal{C})\leq
    I_\mu(\mathcal{B}|\mathcal{C})$. 
    
\,

  \item[(4)] If $\mathcal{A}\leq\mathcal{B}$, then
    $I_\mu(\mathcal{A})\leq I_\mu(\mathcal{B})$. 

\,

  \item[(5)] If $\mathcal{B}\leq\mathcal{C}$, then
             $\displaystyle\int_{A\cap B}I_\mu(\mathcal{A}|\mathcal{B})\,d\mu
             \geq\int_{A\cap B}I_\mu(\mathcal{A}|\mathcal{C})\,d\mu$,
             $\forall\,A\in\mathcal{A},\forall\,B\in\mathcal{B}$. 
            
             Note: In general,
             $\mathcal{B}\leq\mathcal{C}\not\Rightarrow
             I_\mu(\mathcal{A}|\mathcal{B})\geq
             I_\mu(\mathcal{A}|\mathcal{C})$. 
             
\,

  \item[(6)] $\displaystyle\int_C
    I_\mu(\mathcal{A}\vee\mathcal{B}|\mathcal{C})\,d\mu 
              \leq\int_C I_\mu(\mathcal{A}|\mathcal{C})\,d\mu+\int_C I_\mu(\mathcal{B}|\mathcal{C})\,d\mu$, $\forall\,C\in\mathcal{C}$.

              Note: In general, $I_\mu(\mathcal{A}\vee\mathcal{B}|\mathcal{C})\not\leq I_\mu(\mathcal{A}|\mathcal{C})+I_\mu(\mathcal{B}|\mathcal{C})$.
              
\,

  \item[(7)] $\displaystyle\int_X I_\mu(\mathcal{A}\vee\mathcal{B})\,d\mu
               \leq\int_X I_\mu(\mathcal{A})\,d\mu+\int_X I_\mu(\mathcal{B})\,d\mu$.

               Note: In general, $I_\mu(\mathcal{A}\vee\mathcal{B})\not\leq I_\mu(\mathcal{A})+I_\mu(\mathcal{B})$.

\,

  \item[(8)] $\displaystyle\int_{A\cap B}\hspace{-0.2cm}I_\mu(\mathcal{A}|\mathcal{C})\,d\mu
                 \leq\int_{A\cap B}\hspace{-0.2cm}I_\mu(\mathcal{A}|\mathcal{B})\,d\mu+\int_{A\cap B}\hspace{-0.2cm}I_\mu(\mathcal{B}|\mathcal{C})\,d\mu,
                 \forall\,A\in\mathcal{A}, \forall\,B\in\mathcal{B}$.

                Note: In general, $I_\mu(\mathcal{A}|\mathcal{C})\not\leq I_\mu(\mathcal{A}|\mathcal{B})+I_\mu(\mathcal{B}|\mathcal{C})$.
                
\,

  \item[(9)] $I_\mu(\mathcal{A})\leq I_\mu(\mathcal{A}|\mathcal{B})+I_\mu(\mathcal{B})$.
\end{itemize}
\ethm

\bpf
First, notice that $(2)$ follows directly from $(1)$ by setting $\mathcal{C}=\{X\}$.
Similarly, $(4)$ follows directly from $(3)$, and $(7)$ from $(6)$.
It is also easy to see that $(6)$ follows upon combining $(1)$ and $(5)$, since
$\mathcal{C}\leq\mathcal{A}\vee\mathcal{C}$. It therefore only remains to prove parts
$(1)$, $(3)$, $(5)$, $(8)$ and $(9)$.

$(1)$ Let $x\in X$. Then
\begin{eqnarray*}
I_\mu(\mathcal{A}\vee\mathcal{B}|\mathcal{C})(x)
&=&-\log\frac{\mu\bigl((\mathcal{A}\vee\mathcal{B}\vee\mathcal{C})(x)\bigr)}{\mu(\mathcal{C}(x))}\\
&=&-\log\frac{\mu\bigl(\mathcal{B}(x)\cap(\mathcal{A}\vee\mathcal{C})(x)\bigr)}{\mu(\mathcal{C}(x))}\\
&=&-\log\left(\frac{\mu\bigl(\mathcal{B}(x)\cap(\mathcal{A}\vee\mathcal{C})(x)\bigr)}{\mu\bigl((\mathcal{A}\vee\mathcal{C})(x)\bigr)}
    \cdot\frac{{\mu\bigl((\mathcal{A}\vee\mathcal{C})(x)\bigr)}}{\mu(\mathcal{C}(x))}\right)\\
&=&-\log\frac{\mu\bigl(\mathcal{B}(x)\cap(\mathcal{A}\vee\mathcal{C})(x)\bigr)}{\mu\bigl((\mathcal{A}\vee\mathcal{C})(x)\bigr)}
   -\log\frac{{\mu\bigl((\mathcal{A}\vee\mathcal{C})(x)\bigr)}}{\mu(\mathcal{C}(x))}\\
&=&I_\mu(\mathcal{B}|\mathcal{A}\vee\mathcal{C})(x)+I_\mu(\mathcal{A}|\mathcal{C})(x).
\end{eqnarray*}

$(3)$ If $\mathcal{A}\leq\mathcal{B}$, then $\mathcal{A}\vee\mathcal{B}=\mathcal{B}$.
It follows from $(1)$ that
\[
I_\mu(\mathcal{B}|\mathcal{C})
=I_\mu(\mathcal{A}\vee\mathcal{B}|\mathcal{C})
=I_\mu(\mathcal{A}|\mathcal{C})+I_\mu(\mathcal{B}|\mathcal{A}\vee\mathcal{C})
\geq I_\mu(\mathcal{A}|\mathcal{C}).
\]

$(5)$ Suppose that $\mathcal{B}\leq\mathcal{C}$. Let $A\in\mathcal{A}$ and $B\in\mathcal{B}$.
The function $k:[0,1]\to[0,\infty)$ defined by $k(t)=-t\log t$ when $t\in(0,1]$ and $k(0)=0$
is concave, i.e. 
$$
k\bigl(tx+(1-t)y\bigr)\geq tk(x)+(1-t)k(y)
$$ 
for all $t\in[0,1]$ and all $x,y\in[0,1]$.
Consequently, 
$$
k\Big(\sum_{n=1}^\infty a_n b_n\Big)\geq\sum_{n=1}^\infty a_n k(b_n)
$$
whenever $\sum_{n=1}^\infty a_n=1$ and $a_n\geq0$ for all $n\geq1$. Then
\[
k\Bigl(\sum_{C\in\mathcal{C}}\mu_B(C)\frac{\mu(A\cap C)}{\mu(C)}\Bigr)
\geq\sum_{C\in\mathcal{C}}\mu_B(C) k\Bigl(\frac{\mu(A\cap C)}{\mu(C)}\Bigr)
\]
Since $\mathcal{B}\leq\mathcal{C}$, either $C\cap B=C$ or $C\cap B=\emptyset$. Thus, either $\mu_B(C)=\frac{\mu(C)}{\mu(B)}$
or $\mu_B(C)=0$. So the left--hand side
of the previous inequality simplifies to
\begin{eqnarray*}
k\Bigl(\sum_{C\in\mathcal{C}}\mu_B(C)\frac{\mu(A\cap C)}{\mu(C)}\Bigr)
&=&k\Bigl(\sum_{C\subseteq B}\frac{\mu(A\cap C)}{\mu(B)}\Bigr)
=k\Bigl(\frac{\mu(A\cap B)}{\mu(B)}\Bigr)\\
&=&-\frac{\mu(A\cap B)}{\mu(B)}\log\frac{\mu(A\cap B)}{\mu(B)},
\end{eqnarray*}
whereas the right--hand side reduces to
\begin{eqnarray*}
\sum_{C\in\mathcal{C}}\mu_B(C)k\Bigl(\frac{\mu(A\cap C)}{\mu(C)}\Bigr)
=\sum_{C\subseteq B}\frac{\mu(C)}{\mu(B)}k\Bigl(\frac{\mu(A\cap C)}{\mu(C)}\Bigr)
=\sum_{C\subseteq B}-\frac{\mu(A\cap C)}{\mu(B)}\log\frac{\mu(A\cap C)}{\mu(C)}.
\end{eqnarray*}
Hence the inequality becomes
\[
-\frac{\mu(A\cap B)}{\mu(B)}\log\frac{\mu(A\cap B)}{\mu(B)}
\geq\sum_{C\subseteq B}-\frac{\mu(A\cap C)}{\mu(B)}\log\frac{\mu(A\cap C)}{\mu(C)}.
\]
Multiplying both sides by $\mu(B)$ yields
\[
-\mu(A\cap B)\log\frac{\mu(A\cap B)}{\mu(B)}
\geq\sum_{C\subseteq B}-\mu(A\cap C)\log\frac{\mu(A\cap C)}{\mu(C)}.
\]
Therefore
\begin{eqnarray*}
\int_{A\cap B}I_\mu(\mathcal{A}|\mathcal{B})\,d\mu
&=&
-\mu(A\cap B)\log\frac{\mu(A\cap B)}{\mu(B)}
\geq\sum_{C\subseteq B}
-\mu(A\cap C)\log\frac{\mu(A\cap C)}{\mu(C)}\\
&=&\sum_{C\subseteq B}
\int_{A\cap C}I_\mu(\mathcal{A}|\mathcal{C})\,d\mu\\
&=&\int_{A\cap B}I_\mu(\mathcal{A}|\mathcal{C})\,d\mu.
\end{eqnarray*}
$(8)$ Using parts $(3)$ and $(1)$ in turn, we have that
\[
I_\mu(\mathcal{A}|\mathcal{C})
\leq I_\mu(\mathcal{A}\vee\mathcal{B}|\mathcal{C})
=I_\mu(\mathcal{B}|\mathcal{C})+I_\mu(\mathcal{A}|\mathcal{B}\vee\mathcal{C}).
\]
Since $\mathcal{B}\leq\mathcal{B}\vee\mathcal{C}$, part $(v)$ ensures that
\[
\int_{A\cap B}I_\mu(\mathcal{A}|\mathcal{B})\,d\mu
\geq\int_{A\cap B}I_\mu(\mathcal{A}|\mathcal{B}\vee\mathcal{C})\,d\mu,
\forall\,A\in\mathcal{A},\forall\,B\in\mathcal{B}.
\]
Therefore
\[
\int_{A\cap B}I_\mu(\mathcal{A}|\mathcal{C})\,d\mu
\leq\int_{A\cap B}I_\mu(\mathcal{A}|\mathcal{B})\,d\mu+\int_{A\cap B}I_\mu(\mathcal{B}|\mathcal{C})\,d\mu,
\forall\,A\in\mathcal{A},\forall\,B\in\mathcal{B}.
\]
$(ix)$ Using parts $(iv)$ and $(ii)$ in succession, we get
\[
I_\mu(\mathcal{A})\leq I_\mu(\mathcal{A}\vee\mathcal{B})
=I_\mu(\mathcal{A}|\mathcal{B})+I_\mu(\mathcal{B}).
\]
\epf

\

\section{Entropy and Conditional Entropy for Partitions}

\sp The information function associated with a partition gives us the amount of information
that can be gathered from the partition about each and every outcome of the experiment. It it obviously useful
to encompass the information given by a partition within a single number rather than a function.
A natural way to achieve this is to calculate the average information given by the partition.
This means integrating the information function over the entire space. The resulting integral
is called the entropy of the partition. This is the first stage in the definition of the entropy
of an endomorphism.

\sp\bdfn
Let $(X,\mathfrak{ F},\mu)$ be a probability space and
$\mathcal{A}\in\mathrm{Part}(X,\mathfrak{ F})$. 
The  {\em entropy of $\mathcal{A}$} with \index{(N)}{entropy of a partition}
respect to the measure $\mu$  is defined to be 
\[
\H_\mu(\mathcal{A})
:=\int_X I_\mu(\mathcal{A})\,d\mu
=\sum_{A\in\mathcal{A}}-\mu(A)\log\mu(A),
\]
\index{(S)}{$\H_\mu(\mathcal{A})$} where it is implicitly understood that $0\cdot\infty=0$, since null sets do not contribute to the integral.
\edfn

\sp\fr The entropy of a partition is equal to zero if and only if the
partition has an atom of full measure (which implies that all other
atoms are null).  
In particular, $\H_\mu(\{X\})=0$.
Moreover, the entropy of a partition is small if the partition contains
one atom with nearly full measure (so all other atoms have small measure).
Using calculus, it is also possible to show that if the partition $\mathcal{A}$ is finite, then
$$0
\leq \H_\mu(\mathcal{A})\leq\log\#\mathcal{A}
$$ 
and that 
$$
\H_\mu(\mathcal{A})=\log\#\mathcal{A}
$$ 
if and only if
$$
\mu(A)=1/\#\mathcal{A}
$$ 
for all $A\in\mathcal{A}$. In other words, on average we gain the most information from carrying out
an experiment when the potential events have an equal probability of occurring.
Similarly, the conditional entropy of a partition $\mathcal{A}$ given a partition $\mathcal{B}$ is the average conditional information of $\mathcal{A}$ given
$\mathcal{B}$.

\sp\bdfn\label{condentpartition}
Let $(X,\mathfrak{ F},\mu)$ be a probability space and
$\mathcal{A},\mathcal{B}\in\mathrm{Part}(X,\mathfrak{ F})$. 
The {\em conditional entropy of $\mathcal{A}$ given $\mathcal{B}$} is
defined to be 
\[
\H_\mu(\mathcal{A}|\mathcal{B})
:=\int_X I_\mu(\mathcal{A}|\mathcal{B})d\mu
=\sum_{A\in\mathcal{A}}\sum_{B\in\mathcal{B}}-\mu(A\cap
B)\log\frac{\mu(A\cap B)}{\mu(B)}. 
\]
\edfn

\sp\fr Note that $\H_\mu(\mathcal{A})=\H_\mu(\mathcal{A}|\{X\})$. Also, if
for any set $B\in\mathfrak{ F}$ with $\mu(B)>0$, we define a new probability measure on
the space $(B, \mathfrak{ F}|_B)$, called the {\em conditional measure}
of $\mu$ on $B$, by setting 
$$
\mu_B(C):=\frac{\mu(C)}{\mu(B)}
$$ 
and define a new partition $\mathcal{A}|_B$  of $B$ by setting 
$$
\mathcal{A}|_B:=\big\{A\cap
B:A\in\mathcal{A}\big\},
$$
it follows that  
$$
\H_\mu(\mathcal{A}|\mathcal{B}) 
=\sum_{B\in\mathcal{B}}\H_{\mu_B}(\mathcal{A}|_{B})\mu(B).
$$
Indeed,
\begin{eqnarray*}
\H_\mu(\mathcal{A}|\mathcal{B})
&=&\sum_{A\in\mathcal{A}}\sum_{B\in\mathcal{B}}-\mu(A\cap B)\log\frac{\mu(A\cap B)}{\mu(B)}
=\sum_{B\in\mathcal{B}}\sum_{A\in\mathcal{A}}-\frac{\mu(A\cap B)}{\mu(B)}\log\frac{\mu(A\cap B)}{\mu(B)}\cdot\mu(B)\\
&=&\sum_{B\in\mathcal{B}}\sum_{A\in\mathcal{A}}-\mu_B(A)\log\mu_B(A)\cdot\mu(B)\\
&=&\sum_{B\in\mathcal{B}}\H_{\mu_B}(\mathcal{A}|_B)\mu(B).
\end{eqnarray*}
Hence the conditional entropy of $\mathcal{A}$ given $\mathcal{B}$ is
equal to the weighted average 
of the entropies of the partitions of each atom $B\in\mathcal{B}$ into
the sets $A\cap B$, 
$A\in\mathcal{A}$.

\sp Of course, the properties of entropy
(resp. conditional entropy) are inherited from the properties of the
information function (resp. 
the conditional information function) via integration.

\bthm\label{sem2thm7.1}
Let $(X,\mathfrak{ F},\mu)$ be a probability space and $\mathcal{A},\mathcal{B},\mathcal{C}
\in\mathrm{Part}(X,\mathfrak{ F})$. The following statements hold:
\begin{itemize}
  \item[(i)] $\H_\mu(\mathcal{A}\vee\mathcal{B}|\mathcal{C})=\H_\mu(\mathcal{A}|\mathcal{C})+\H_\mu(\mathcal{B}|\mathcal{A}\vee\mathcal{C})$.

\,
  
  \item[(2)] $\H_\mu(\mathcal{A}\vee\mathcal{B})=\H_\mu(\mathcal{A})+\H_\mu(\mathcal{B}|\mathcal{A})$.
  
\,
  
  \item[(3)] If $\mathcal{A}\leq\mathcal{B}$, then $\H_\mu(\mathcal{A}|\mathcal{C})\leq \H_\mu(\mathcal{B}|\mathcal{C})$.
  
\,
  
  \item[(4)] If $\mathcal{A}\leq\mathcal{B}$, then $\H_\mu(\mathcal{A})\leq \H_\mu(\mathcal{B})$.
  
\,
  
  \item[(5)] If $\mathcal{B}\leq\mathcal{C}$, then $\H_\mu(\mathcal{A}|\mathcal{B})\geq \H_\mu(\mathcal{A}|\mathcal{C})$.
  
\,
  
  \item[(6)] $\H_\mu(\mathcal{A}\vee\mathcal{B}|\mathcal{C})\leq \H_\mu(\mathcal{A}|\mathcal{C})+\H_\mu(\mathcal{B}|\mathcal{C})$.
  
\,
  
  \item[(7)] $\H_\mu(\mathcal{A}\vee\mathcal{B})\leq \H_\mu(\mathcal{A})+\H_\mu(\mathcal{B})$.
  
\,
  
  \item[(8)] $\H_\mu(\mathcal{A}|\mathcal{C})\leq
    \H_\mu(\mathcal{A}|\mathcal{B})+\H_\mu(\mathcal{B}|\mathcal{C})$. 
    
\,
  
  \item[(9)] $\H_\mu(\mathcal{A})\leq
    \H_\mu(\mathcal{A}|\mathcal{B})+\H_\mu(\mathcal{B})$. 
\end{itemize}
\ethm

\bpf
All the statements follow from their counterparts in
Theorem~\ref{sem2thm7.0} after integration or summation over atoms. 
For instance, let us prove $(5)$. If $\mathcal{B}\leq\mathcal{C}$,
then it follows from Theorem~\ref{sem2thm7.0}$(5)$ that 
\begin{eqnarray*}
\H_\mu(\mathcal{A}|\mathcal{B})
&=&\int_{X}I_\mu(\mathcal{A}|\mathcal{B})\,d\mu
=\sum_{A\in\mathcal{A}}\sum_{B\in\mathcal{B}}\int_{A\cap
  B}I_\mu(\mathcal{A}|\mathcal{B})\,d\mu \\
&\geq&\sum_{A\in\mathcal{A}}\sum_{B\in\mathcal{B}}\int_{A\cap
  B}I_\mu(\mathcal{A}|\mathcal{C})\,d\mu \\
&=&\int_{X}I_\mu(\mathcal{A}|\mathcal{C})\,d\mu \\
&=&\H_\mu(\mathcal{A}|\mathcal{C}).
\end{eqnarray*}
\epf

\sp\section{Entropy of a (Probability) Measure--Preserving Endomorphism}

\sp So far in this chapter we have studied partitions of a space, but of course, this is not our ultimate objective. Our ultimate aim is however
to study measure--theoretical dynamical systems. So let 
$$
T:X\lra X
$$ 
be a
measure--preserving endomorphism of a probability space
$(X,\mathfrak{ F},\mu)$, and let $\mathcal{A}=\{A_k\}_{k\geq1}$ be a
countable measurable partition of $X$. Observe that
$$
T^{-1}\mathcal{A}:=\{T^{-1}(A):A\in\mathcal{A}\}
$$ 
is also a
countable measurable partition of $X$.

\sp Recall that the set $X$ can be thought of as representing the set of all possible
outcomes (or states) of an experiment, 
while the \mbox{$\sigma$-algebra} $\mathfrak{ F}$ consists of the set of
all possible events, and 
$\mu(E)$ is the probability that event $E$ takes place.
Recall also that a partition $\mathcal{A}=\{A_k\}_{k\geq1}$ can be
thought of as the set of all 
observations that can be made with a given instrument.
The action of $T$ on $(X,\mathfrak{ F},\mu)$ 
may be conceived as the passage of one unit
of time (for instance, a day). Today would naturally taken as reference point
for time $0$. Suppose that we conduct the experiment with our instrument
tomorrow. The resulting observation would be one of the $A_k$s,
say $A_{k_1}$, on day $1$. Due to the passage of time (in other words,
one iteration of $T$), 
in order to make observation $A_{k_1}$ at time $1$, our
measure-theoretic system would have 
to be in one of the states of $T^{-1}(A_{k_1})$ today. The probability
of making observation $A_{k_1}$ on day $1$ is thus $\mu(T^{-1}(A_{k_1}))$.
Assume now that we conduct the same experiment for $n$ consecutive
days, starting today. 
What is the probability that we make the sequence of observations
$A_{k_0},A_{k_1},\ldots,A_{k_{n-1}}$ on those successive days?
We would make those observations precisely if our system is in one of
the states of 
$$
\bigcap_{m=0}^{n-1}T^{-m}(A_{k_m})
$$ 
today. Therefore, the probability that
our observations are respectively $A_{k_0},A_{k_1},\ldots,A_{k_{n-1}}$
on $n$ successive 
days starting today is 
$$
\mu\Big(\bigcap_{m=0}^{n-1}T^{-m}(A_{k_m})\Big).
$$
Given
this discussion, we claim it is thus natural to introduce for all
$0\leq m\leq n-1$ the partitions 
\[
\mathcal{A}_m^n:=\bigvee_{i=m}^{n-1}T^{-i}\mathcal{A}
                =T^{-m}\mathcal{A}\vee\cdots\vee T^{-(n-1)}\mathcal{A}.
\]
 \index{(S)}{\mathcal{A}_m^n} If $m\geq n$, we define $\mathcal{A}_m^n$ to be the trivial partition
$\{X\}$. To 
shorten the notation, we shall write $\mathcal{A}^n$ in lieu of
$\mathcal{A}_0^n$. 
Let us give basic properties of the operator $T^{-1}$ on partitions.

\

\blem\label{lemT}
Let $T:X\lra X$ be a measurable transformation of a measurable space
$(X,\mathfrak{ F})$, and 
$\mathcal{A},\mathcal{B}\in\mathrm{Part}(X,\mathfrak{ F})$.
Then the following statements hold:
\begin{itemize}\label{lempartt}
\item[(1)] The operator $T^{-1}$ commutes with the operator $\vee$, that is, $T^{-1}(\mathcal{A}\vee\mathcal{B})=T^{-1}\mathcal{A}\vee T^{-1}\mathcal{B}$.

\,

\item[(2)] $T^{-1}(\mathcal{A}_m^n)=(T^{-1}\mathcal{A})_m^n$ for all $m,n\geq0$.

\,

\item[(3)] $(\mathcal{A}\vee\mathcal{B})_m^n=\mathcal{A}_m^n\vee\mathcal{B}_m^n$
            for all $m,n\geq0$.
            
\,

\item[(4)] $(\mathcal{A}_k^l)_m^n=\mathcal{A}_{k+m}^{l+n-1}$.

\,

\item[(5)] $T^{-1}$ preserves the order $\leq$, that is, if $\mathcal{A}\leq\mathcal{B}$, then $T^{-1}\mathcal{A}\leq T^{-1}\mathcal{B}$.

\,

\item[(6)] More generally, if $\mathcal{A}\leq\mathcal{B}$ then $\mathcal{A}_m^n\leq\mathcal{B}_m^n$
            for all $m,n\geq0$.
            
\,

\item[(7)] $(T^{-1}\mathcal{A})(x)=T^{-1}\bigl(\mathcal{A}(T(x))\bigr)$ for all $x\in X$.
\end{itemize}
\elem

\bpf The proof of assertions $(1)$ and $(5)$ are left to the reader. In order to prove
$(2)$ by using $(1)$ repeatedly, we get that
\begin{eqnarray*}
T^{-1}(\mathcal{A}_m^n)
=T^{-1}\Bigl(\bigvee_{i=m}^{n-1}T^{-i}\mathcal{A}\Bigr)
=\bigvee_{i=m}^{n-1}T^{-1}(T^{-i}\mathcal{A})
=\bigvee_{i=m}^{n-1}T^{-i}(T^{-1}\mathcal{A})
=(T^{-1}\mathcal{A})_m^n.
\end{eqnarray*}

$(3)$ Again by using $(1)$ repeatedly, we obtain that
\begin{eqnarray*}
(\mathcal{A}\vee\mathcal{B})_m^n
&=&\bigvee_{i=m}^{n-1}T^{-i}(\mathcal{A}\vee\mathcal{B})
=\bigvee_{i=m}^{n-1}(T^{-i}\mathcal{A}\vee T^{-i}\mathcal{B}) \\
&=&\Bigl(\bigvee_{i=m}^{n-1}T^{-i}\mathcal{A}\Bigr)\bigvee\Bigl(\bigvee_{i=m}^{n-1}T^{-i}\mathcal{B}\Bigr) \\
&=&\mathcal{A}_m^n\vee\mathcal{B}_m^n.
\end{eqnarray*}

Dealing with $(4)$, by using $(1)$, it follows that
\begin{eqnarray*}
(\mathcal{A}_k^l)_m^n
=\bigvee_{j=m}^{n-1}T^{-j}(\mathcal{A}_k^l)
&=&\bigvee_{j=m}^{n-1}T^{-j}\Bigl(\bigvee_{i=k}^{l-1}T^{-i}\mathcal{A}\Bigr)\\
&=&\bigvee_{j=m}^{n-1}\bigvee_{i=k}^{l-1}T^{-(i+j)}\mathcal{A}
=\bigvee_{s=k+m}^{l+n-2}T^{-s}\mathcal{A}
=\mathcal{A}_{k+m}^{l+n-1}.
\end{eqnarray*}

$(6)$ Suppose $\mathcal{A}\leq\mathcal{B}$. Using $(5)$ repeatedly and Lemma~\ref{lempartord}$(v)$, we obtain
\[
\mathcal{A}_m^n
=\bigvee_{i=m}^{n-1}T^{-i}\mathcal{A}
\leq\bigvee_{i=m}^{n-1}T^{-i}\mathcal{B}
=\mathcal{B}_m^n.
\]

$(7)$ Let $x\in X$. Let $T^{-1}(A)=(T^{-1}\mathcal{A})(x)$, that is,
$A\in\mathcal{A}$ is such that $x\in T^{-1}(A)$.
Then $T(x)\in A$, i.e. $A=\mathcal{A}(T(x))$. Hence,
$$
(T^{-1}\mathcal{A})(x)=T^{-1}(A)=T^{-1}\bigl(\mathcal{A}(T(x))\bigr).
$$
\epf

\sp We now describe the behavior of the operator $T^{-1}$ with respect
to the information function. 

\sp\blem\label{lemparti}
Let $T:X\lra X$ be a measure--preserving endomorphism of a probability
space $(X,\mathfrak{ F},\mu)$, and let 
$\mathcal{A},\mathcal{B}\in\mathrm{Part}(X,\mathfrak{ F})$.
Then the following statements hold:
\begin{itemize}
\item[(1)]
  $I_\mu(T^{-1}\mathcal{A}|T^{-1}\mathcal{B})=I_\mu(\mathcal{A}|\mathcal{B})\circ
  T$. 
  
\,

\item[(2)] $I_\mu(T^{-1}\mathcal{A})=I_\mu(\mathcal{A})\circ T$.
\end{itemize}
\elem

\bpf It is clear that $(2)$ follows from $(1)$ by setting
  $\mathcal{B}=\{X\}$. To prove $(1)$, let $x\in X$. By
  Lemma~\ref{lempartt}$(1)$ and $(7)$ and the assumption that $\mu$
  is $T$-invariant, 
we have
\begin{eqnarray*}
I_\mu(T^{-1}\mathcal{A}|T^{-1}\mathcal{B})(x)
&=&-\log\frac{\mu\bigl((T^{-1}\mathcal{A}\vee T^{-1}\mathcal{B})(x)\bigr)}{\mu\bigl((T^{-1}\mathcal{B})(x)\bigr)}
=-\log\frac{\mu\bigl(\bigl(T^{-1}(\mathcal{A}\vee\mathcal{B})\bigr)(x)\bigr)}{\mu\bigl((T^{-1}\mathcal{B})(x)\bigr)}\\
&=&-\log\frac{\mu\bigl(T^{-1}\bigl((\mathcal{A}\vee\mathcal{B})(T(x))\bigr)\bigr)}{\mu\bigl(T^{-1}\bigl(\mathcal{B}(T(x))\bigr)\bigr)}
=-\log\frac{\mu\bigl((\mathcal{A}\vee\mathcal{B})(T(x))\bigr)}{\mu\bigl(\mathcal{B}(T(x))\bigr)}\\
&=&I_\mu(\mathcal{A}|\mathcal{B})(T(x)) \\
&=&I_\mu(\mathcal{A}|\mathcal{B})\circ T(x).
\end{eqnarray*}
\epf

\sp A more intricate property of the information function is given in
the following lemma. 

\sp\blem\label{sem2lem7.2i}
Let $T:X\lra X$ be a measure--preserving endomorphism of a probability
space $(X,\mathfrak{ F},\mu)$,
If $\mathcal{A}\in\mathrm{Part}(X,\mathfrak{ F})$, 
then 
$$
I_\mu(\mathcal{A}^n)=\sum_{j=1}^n
I_\mu(\mathcal{A}|\mathcal{A}^j_1)\circ T^{n-j}
$$ 
for all integers $n\geq1$. 
\elem

\bpf
We will prove this lemma by induction. For $n=1$, since
$\mathcal{A}_1^1$ is by definition equal to the trivial partition
$\{X\}$, we have 
\[
I_\mu(\mathcal{A}^1)
=I_\mu(\mathcal{A})
=I_\mu(\mathcal{A}|\{X\})
=I_\mu(\mathcal{A}|\mathcal{A}_1^1)
=I_\mu(\mathcal{A}|\mathcal{A}_1^1)\circ T^{1-1}.
\]
Now suppose that the lemma holds for some $n\geq1$. Then, in light of
Theorem~\ref{sem2thm7.0}$(ii)$ 
and Lemma~\ref{lemparti}$(ii)$,
we obtain that
\begin{eqnarray*}
I_\mu(\mathcal{A}^{n+1})
&=&I_\mu(\mathcal{A}^{n+1}_1\vee\mathcal{A})
=I_\mu(\mathcal{A}^{n+1}_1)+I_\mu(\mathcal{A}|\mathcal{A}^{n+1}_1)
=I_\mu(T^{-1}(\mathcal{A}^n))+I_\mu(\mathcal{A}|\mathcal{A}^{n+1}_1)\\
&=&I_\mu(\mathcal{A}^n)\circ T+I_\mu(\mathcal{A}|\mathcal{A}^{n+1}_1)
=\sum_{j=1}^n I_\mu(\mathcal{A}|\mathcal{A}^{j}_1)\circ T^{n-j}\circ T+I_\mu(\mathcal{A}|\mathcal{A}^{n+1}_1)\\
&=&\sum_{j=1}^n I_\mu(\mathcal{A}|\mathcal{A}^{j}_1)\circ T^{n+1-j}+I_\mu(\mathcal{A}|\mathcal{A}^{n+1}_1)
   \circ T^{n+1-(n+1)}\\
&=&\sum_{j=1}^{n+1}I_\mu(\mathcal{A}|\mathcal{A}^{j}_1)\circ T^{n+1-j}.
\end{eqnarray*}
\epf

\sp We now turn our attention to the effect that our measure--theoretic dynamical system $T$ has
on entropy. In particular, observe that because the system is measure-preserving, conducting the experiment
today or tomorrow (or at any time in the future) gives us the same amount of average information about the
outcome. This is the meaning of the second of the following properties
of entropy. 

\sp\blem\label{lempartH}
If $T:X\lra X$ is a measure--preserving endomorphism of a probability
space $(X,\mathfrak{ F},\mu)$, and 
$\mathcal{A},\mathcal{B}\in\mathrm{Part}(X,\mathfrak{ F})$,
then
\begin{itemize}
\item[(1)]
  $\H_\mu(T^{-1}\mathcal{A}|T^{-1}\mathcal{B})=\H_\mu(\mathcal{A}|\mathcal{B})$. 

\,

\item[(2)] $\H_\mu(T^{-1}\mathcal{A})=\H_\mu(\mathcal{A})$.

\,

\item[(3)] $\H_\mu(\mathcal{A}^n|\mathcal{B}^n)\leq n
  \H_\mu(\mathcal{A}|\mathcal{B})$. 
\end{itemize}
\elem

\bpf Part $(2)$ follows from $(1)$ by taking $\mathcal{B}=\{X\}$.
$(i)$ Using Lemma~\ref{lemparti}$(1)$ and the $T$-invariance of $\mu$, we obtain
\begin{eqnarray*}
\H_\mu(T^{-1}\mathcal{A}|T^{-1}\mathcal{B})
&=&\int_X I_\mu(T^{-1}\mathcal{A}|T^{-1}\mathcal{B})\,d\mu
=\int_X I_\mu(\mathcal{A}|\mathcal{B})\circ T\,d\mu\\
&=&\int_X I_\mu(\mathcal{A}|\mathcal{B})\,d\mu \\
&=&\H_\mu(\mathcal{A}|\mathcal{B}).
\end{eqnarray*}
$(3)$ We first prove that 
$$
\H_\mu(\mathcal{A}^n|\mathcal{B}^n)
\leq\sum_{j=0}^{n-1}\H_\mu(T^{-j}\mathcal{A}|T^{-j}\mathcal{B}).
$$
This statement clearly holds when $n=1$. Suppose that it holds for some
$n\geq1$. Then, using Theorem~\ref{sem2thm7.1}$(i)$ and $(v)$, we have that
\begin{eqnarray*}
\H_\mu(\mathcal{A}^{n+1}|\mathcal{B}^{n+1})
&=&\H_\mu\bigl(\mathcal{A}^n\vee T^{-n}\mathcal{A}|\mathcal{B}^n\vee T^{-n}\mathcal{B}\bigr)\\
&=&\H_\mu\bigl(\mathcal{A}^n|\mathcal{B}^n\vee T^{-n}\mathcal{B}\bigr)
   +\H_\mu\bigl(T^{-n}\mathcal{A}|\mathcal{A}^n\vee\mathcal{B}^n\vee T^{-n}\mathcal{B}\bigr)\\
&\leq&\H_\mu(\mathcal{A}^n|\mathcal{B}^n)+\H_\mu(T^{-n}\mathcal{A}|T^{-n}\mathcal{B})\\
&\leq&\sum_{j=0}^{n-1}\H_\mu(T^{-j}\mathcal{A}|T^{-j}\mathcal{B})+\H_\mu(T^{-n}\mathcal{A}|T^{-n}\mathcal{B})\\
&=&\sum_{j=0}^n \H_\mu(T^{-j}\mathcal{A}|T^{-j}\mathcal{B}).
\end{eqnarray*}
By induction, the above statement holds for all $n\geq1$. By $(1)$, we get
\[
\H_\mu(\mathcal{A}^n|\mathcal{B}^n)
\leq\sum_{j=0}^{n-1}\H_\mu(T^{-j}\mathcal{A}|T^{-j}\mathcal{B})
=\sum_{j=0}^{n-1}\H_\mu(\mathcal{A}|\mathcal{B})
=n \H_\mu(\mathcal{A}|\mathcal{B}).
\]
 \epf

\sp The average information gained by conducting
an experiment on $n$ consecutive days, using the partition  $\mathcal{A}$,
is given by the entropy $\H_\mu(\mathcal{A}^n)$ since $\mathcal{A}^n$
has for atoms the sets 
$$
\bigcap_{m=0}^{n-1}T^{-m}(A_{k_m}),
$$ 
where $A_{k_m}\in\mathcal{A}$ for all $m$.
Not surprisingly, the average information gained by conducting
the experiment on $n$ consecutive days using the partition $\mathcal{A}$
is equal to the sum of the average conditional information gained by
performing $\mathcal{A}$ 
on day $j+1$ given that the outcome of performing $\mathcal{A}$ over
the previous $j$ days is known, summing 
from the first day to the last day. This is summarised in the next lemma.

\sp\blem\label{sem2lem7.2}
If $T:X\lra X$ is a measure--preserving endomorphism of a probability
space $(X,\mathfrak{ F},\mu)$, $\mathcal{A}\in\mathrm{Part}(X,\mathfrak{ F})$, 
then 
$$
\H_\mu(\mathcal{A}^n)=\sum_{j=1}^n
\H_\mu(\mathcal{A}|\mathcal{A}^j_1)
$$
for all integers $n\geq1$. 
\elem

\bpf We deduce from Lemma~\ref{sem2lem7.2i} and the $T$--invariance of $\mu$ that 
\begin{eqnarray*}
\H_\mu(\mathcal{A}^n)
=\int_X I_\mu(\mathcal{A}^n)\,d\mu
&=&\sum_{j=1}^n \int_X I_\mu(\mathcal{A}|\mathcal{A}^j_1)\circ T^{n-j}\,d\mu\\
&=&\sum_{j=1}^n \int_X I_\mu(\mathcal{A}|\mathcal{A}^j_1)\,d\mu
=\sum_{j=1}^n \H_\mu(\mathcal{A}|\mathcal{A}^j_1).
\end{eqnarray*}
\epf

\sp In view of Theorem~\ref{sem2thm7.1}$(5)$, observe that
$$
\H_\mu(\mathcal{A}|\mathcal{A}^{j+1}_1)\leq \H_\mu(\mathcal{A}|\mathcal{A}^j_1),
$$
since $\mathcal{A}^{j+1}_1\geq\mathcal{A}^j_1$. 
So the sequence 
$$
\(\H_\mu(\mathcal{A}|\mathcal{A}^j_1)\)_{j\geq1}
$$ 
decreases to some limit
which we shall denote by 
$$
\H_\mu(T,\mathcal{A}).
$$
Consequently, the corresponding sequence of
Ces\`aro averages 
$$
\Bigl(\frac{1}{n}\sum_{j=1}^n \H_\mu(\mathcal{A}|\mathcal{A}^j_1)\Bigr)_{n\geq1}
=\Bigl(\frac{1}{n}\H_\mu(\mathcal{A}^n)\Bigr)_{n\geq1}
$$
decreases to the same limit. Thus, the following definition makes sense.
This is the second stage in the definition of the entropy of an
endomorphism.

\sp\bdfn\label{entwrtpartition}
If $T:X\lra X$ is a measure--preserving endomorphism of a probability
space $(X,\mathfrak{ F},\mu)$, and $\mathcal{A}\in\mathrm{Part}(X,\mathfrak{ F})$, then the quantity $\H_\mu(T,\mathcal{A})$ defined by
\[
\h_\mu(T,\mathcal{A})
:=\lim_{n\to\infty}\H_\mu(\mathcal{A}|\mathcal{A}^n_1)
=\lim_{n\to\infty}\frac{1}{n}\H_\mu(\mathcal{A}^n)
=\inf_{n\ge 1}\big\{\H_\mu(\mathcal{A}|\mathcal{A}^n_1)\big\}
=\inf_{n\ge 1}\Big\{\frac{1}{n}\H_\mu(\mathcal{A}^n)\Big\},
\]
is called the {\em entropy of $T$ with respect to $\mathcal{A}$}. \index{(N)}{entropy of $T$} \index{(N)}{entropy of an endomorphism}\index{(S)}{$\h_\mu(T,\mathcal{A})$}
\edfn

\sp The following theorem lists some of the basic properties of
$\H_\mu(T,\cdot)$. 

\sp\bthm\label{sem2thm7.3}
Let $T:X\lra X$ be a measure--preserving endomorphism of a probability
space $(X,\mathfrak{ F},\mu)$.
If $\mathcal{A},\mathcal{B}\in\mathrm{Part}(X,\mathfrak{ F})$,
then the following statements hold:
\begin{itemize}
  \item[(1)] $\h_\mu(T,\mathcal{A})\leq \H_\mu(\mathcal{A})$.
  
  \,
  
  \item[(2)] $\h_\mu(T,\mathcal{A}\vee\mathcal{B})\leq \h_\mu(T,\mathcal{A})+\h_\mu(T,\mathcal{B})$.
   
  \,
  
  \item[(3)] If $\mathcal{A}\leq\mathcal{B}$, then $\H_\mu(T,\mathcal{A})\leq \h_\mu(T,\mathcal{B})$.
   
  \,
  
  \item[(4)] $\h_\mu(T,\mathcal{A})\leq \h_\mu(T,\mathcal{B})+\H_\mu(\mathcal{A}|\mathcal{B})$.
   
  \,
  
  \item[(5)] $\h_\mu(T,T^{-1}\mathcal{A})=\h_\mu(T,\mathcal{A})$.
   
  \,
  
  \item[(6)] For all $k\geq1$, $\h_\mu(T,\mathcal{A}^k)=\h_\mu(T,\mathcal{A})$.
   
  \,
  
  \item[(7)] For all $k\geq1$, $\H_\mu(T^k,\mathcal{A}^k)=k\H_\mu(T,\mathcal{A})$.
   
  \,
  
  \item[(8)] If $T$ is invertible and $k\geq1$, then $\h_\mu(T,\mathcal{A})=\h_\mu\bigl(T,\bigvee_{i=-k}^k T^i\mathcal{A}\bigr)$.
   
  \,
  
  \item[(9)] If $\mathcal{B}_1,\mathcal{B}_2,\ldots\in\mathrm{Part}(X,\mathfrak{ F})$ are such that
              $\displaystyle\lim_{n\to\infty}\H_\mu(\mathcal{A}|\mathcal{B}_n)=0$, then
\[
\h_\mu(T,\mathcal{A})\leq\liminf_{n\to\infty}\h_\mu(T,\mathcal{B}_n).
\]
  \item[(10)] If $\displaystyle\lim_{n\to\infty}\H_\mu(\mathcal{A}|\mathcal{B}^n)=0$,
then $\h_\mu(T,\mathcal{A})\leq \h_\mu(T,\mathcal{B})$.
\end{itemize}
\ethm

\bpf
$(1)$ This follows from the fact that
$\h_\mu(T,\mathcal{A})=\lim_{n\to\infty}\frac{1}{n}\H_\mu(\mathcal{A}^n)$ and that
 $\bigl(\frac{1}{n}\H_\mu(\mathcal{A}^n)\bigr)_{n\geq1}$ is a decreasing sequence
with first term given by $\H_\mu(\mathcal{A})$.

\sp $(2)$ Using Theorem~\ref{sem2thm7.1}$(vii)$, we get
\begin{eqnarray*}
\h_\mu(T,\mathcal{A}\vee\mathcal{B})
&=&\lim_{n\to\infty}\frac{1}{n}\H_\mu((\mathcal{A}\vee\mathcal{B})^n)
=\lim_{n\to\infty}\frac{1}{n}\H_\mu(\mathcal{A}^n\vee\mathcal{B}^n)\\
&\leq&\lim_{n\to\infty}\frac{1}{n}\Bigl[\H_\mu(\mathcal{A}^n)+\H_\mu(\mathcal{B}^n)\Bigr]\\
&=&\lim_{n\to\infty}\frac{1}{n}\H_\mu(\mathcal{A}^n)+\lim_{n\to\infty}\frac{1}{n}\H_\mu(\mathcal{B}^n)\\
&=&\h_\mu(T,\mathcal{A})+\h_\mu(T,\mathcal{B}).
\end{eqnarray*}

$(3)$ If $\mathcal{A}\leq\mathcal{B}$, then $\mathcal{A}^n\leq\mathcal{B}^n$ for all $n\geq1$.
Consequently, $\H_\mu(\mathcal{A}^n)\leq \H_\mu(\mathcal{B}^n)$ for all $n\geq1$. Therefore,
\[
\h_\mu(T,\mathcal{A})
=\lim_{n\to\infty}\frac{1}{n}\H_\mu(\mathcal{A}^n)\\
\leq\lim_{n\to\infty}\frac{1}{n}\H_\mu(\mathcal{B}^n)\\
=\h_\mu(T,\mathcal{B}).
\]

$(4)$ Calling upon Theorem~\ref{sem2thm7.1}$(ix)$ and Lemma~\ref{lempartH}$(iii)$, we obtain that
\begin{eqnarray*}
\h_\mu(T,\mathcal{A})
&=&\lim_{n\to\infty}\frac{1}{n}\H_\mu(\mathcal{A}^n)
\leq\liminf_{n\to\infty}\frac{1}{n}\Bigl(\H_\mu(\mathcal{A}^n|\mathcal{B}^n)+\H_\mu(\mathcal{B}^n)\Bigr)\\
&=&\liminf_{n\to\infty}\frac{1}{n}\H_\mu(\mathcal{A}^n|\mathcal{B}^n)
      +\lim_{n\to\infty}\frac{1}{n}\H_\mu(\mathcal{B}^n)\\
&\leq&\H_\mu(\mathcal{A}|\mathcal{B})+\h_\mu(T,\mathcal{B}).
\end{eqnarray*}

$(5)$ By Lemma~\ref{lempartt}$(2)$, we know that $(T^{-1}\mathcal{A})^n=T^{-1}(\mathcal{A}^n)$ for
all $n\geq1$. Then, using Lemma~\ref{lempartH}$(ii)$, we deduce that
\begin{eqnarray*}
\h_\mu(T,T^{-1}\mathcal{A})
&=&\lim_{n\to\infty}\frac{1}{n}\H_\mu\bigl((T^{-1}\mathcal{A})^n\bigr)\\
&=&\lim_{n\to\infty}\frac{1}{n}\H_\mu\bigl(T^{-1}(\mathcal{A}^n)\bigr)
=\lim_{n\to\infty}\frac{1}{n}\H_\mu(\mathcal{A}^n)
=\h_\mu(T,\mathcal{A}).
\end{eqnarray*}

$(6)$ By Lemma~\ref{lempartt}$(iv)$, we know that $(\mathcal{A}^k)^n=\mathcal{A}^{n+k-1}$ and hence
\begin{eqnarray*}
\h_\mu(T,\mathcal{A}^k)
&=&\lim_{n\to\infty}\frac{1}{n}\H_\mu\bigl((\mathcal{A}^k)^n\bigr)
=\lim_{n\to\infty}\frac{1}{n}\H_\mu(\mathcal{A}^{n+k-1})\\
&=&\lim_{n\to\infty}\frac{n+k-1}{n}\cdot\frac{1}{n+k-1}\H_\mu(\mathcal{A}^{n+k-1})\\
&=&\lim_{n\to\infty}\frac{n+k-1}{n}\cdot\lim_{n\to\infty}\frac{1}{n+k-1}\H_\mu(\mathcal{A}^{n+k-1})\\
&=&\lim_{m\to\infty}\frac{1}{m}\H_\mu(\mathcal{A}^m)
=\h_\mu(T,\mathcal{A}).
\end{eqnarray*}

$(7)$ Let $k\geq1$. Using part~$(vi)$, we have that
\begin{eqnarray*}
\h_\mu(T^k,\mathcal{A}^k)
&=&\lim_{n\to\infty}\frac{1}{n}\H_\mu\Bigl(\bigvee_{j=0}^{n-1}T^{-kj}(\mathcal{A}^k)\Bigr)
=\lim_{n\to\infty}\frac{1}{n}\H_\mu\Bigl(\bigvee_{j=0}^{n-1}T^{-kj}\Bigl(\bigvee_{i=0}^{k-1}T^{-i}\mathcal{A}\Bigr)\Bigr)\\
&=&\lim_{n\to\infty}\frac{1}{n}\H_\mu\Bigl(\bigvee_{l=0}^{kn-1}T^{-l}\mathcal{A}\Bigr)\\
&=&k\lim_{n\to\infty}\frac{1}{kn}\H_\mu(\mathcal{A}^{kn})
=k\,\h_\mu(T,\mathcal{A}).
\end{eqnarray*}

$(8)$ The proof is similar to that of part $(6)$ and is thus left to the reader.

\sp $(9)$
By part $(4)$, for each $n\geq1$ we have,
\[
\h_\mu(T,\mathcal{A})\leq \h_\mu(T,\mathcal{B}_n)+\H_\mu(\mathcal{A}|\mathcal{B}_n).
\]
So if $(\mathcal{B}_n)_{n\geq1}$ are partitions such that
$\displaystyle\lim_{n\to\infty}\H_\mu(\mathcal{A}|\mathcal{B}_n)=0$, then
\begin{eqnarray*}
\h_\mu(T,\mathcal{A})
&\leq&\liminf_{n\to\infty}\bigl[\h_\mu(T,\mathcal{B}_n)+\H_\mu(\mathcal{A}|\mathcal{B}_n)\bigr]\\
&=&\liminf_{n\to\infty}\h_\mu(T,\mathcal{B}_n)+\lim_{n\to\infty}\H_\mu(\mathcal{A}|\mathcal{B}_n)
=\liminf_{n\to\infty}\h_\mu(T,\mathcal{B}_n).
\end{eqnarray*}

$(10)$ Suppose $\displaystyle\lim_{n\to\infty}\H_\mu(\mathcal{A}|\mathcal{B}^n)=0$. By parts $(9)$ and $(6)$, we have
\[
\h_\mu(T,\mathcal{A})
\leq\liminf_{n\to\infty}\h_\mu(T,\mathcal{B}^n)
=\lim_{n\to\infty}\h_\mu(T,\mathcal{B})
=\h_\mu(T,\mathcal{B}).
\]
\epf

\sp The entropy of $T$ is defined in a similar way to topological
entropy. The third and last 
stage in the definition of the entropy of an endomorphism consists of
passing to a supremum. 

\sp\bdfn\label{mte}
If $T:X\lra X$ is a measure--preserving endomorphism of a probability
space $(X,\mathfrak{ F},\mu)$, then
the {\em measure--theoretic entropy}
of $T$, \index{(N)}{measure--theoretic entropy
of $T$} \index{(N)}{measure--theoretic entropy
of an endomorphism} denoted $\h_\mu(T)$, \index{(S)}{ $\h_\mu(T)$} is defined by 
\[\h_\mu(T):=\sup\Big\{\h_\mu(T,\mathcal{A}): \mathcal{A} \ \text{ is a
  finite partition of }X\Big\}. 
\]
\edfn

\sp The following theorem is a useful tool for calculating the entropy
of an endomorphism. 

\sp\bthm\label{sem2thm7.4}
If $T:X\to X$ is a measure-preserving endomorphism of a probability
space $(X,\mathfrak{ F},\mu)$, then
\begin{itemize}
  \item[(1)] For all $k\geq1$, $\h_\mu(T^k)=k \h_\mu(T)$.
  
  \,
  
  \item[(2)] If $T$ is invertible, then $\H_\mu(T^{-1})=\H_\mu(T)$.
\end{itemize}
\ethm

\bpf
$(1)$ Let $k\geq1$. Then, by Theorem~\ref{sem2thm7.3}$(vii)$,
\begin{eqnarray*}
k\h_\mu(T)
&=&\sup\{k\h_\mu(T,\mathcal{A}):\mathcal{A}\text{ a finite partition}\}\\
&=&\sup\{\h_\mu(T^k,\mathcal{A}^k):\mathcal{A} \text{ a finite partition}\}\\
&\leq&\sup\{\h_\mu(T^k,\mathcal{B}):\mathcal{B}\text{ a finite
  partition}\}=\h_\mu(T^k). 
\end{eqnarray*}
On the other hand, by Theorem~\ref{sem2thm7.3}$(3)$ and $(7)$,
\[
\h_\mu(T^k,\mathcal{A})
\leq \h_\mu(T^k,\mathcal{A}^k)
=k \h_\mu(T,\mathcal{A}).
\]
Passing to the supremum over all finite partitions $\mathcal{A}$ of
$X$ on both sides, 
we obtain the desired inequality, namely,  $\h_\mu(T^k)\leq k\h_\mu(T)$.

\sp\fr $(2)$ To distinguish the action of $T$ from the action of
$T^{-1}$ on a partition, we shall 
use the respective notations $\mathcal{A}_T^n$ and $\mathcal{A}_{T^{-1}}^n$. Using
Lemmas~\ref{lempartH}$(ii)$ and~\ref{lemT}$(ii)$ in turn, we deduce that
\begin{eqnarray*}
\h_\mu(\mathcal{A}_{T^{-1}}^n)
&=&\H_\mu\Bigl(\bigvee_{i=0}^{n-1}(T^{-1})^{-i}\mathcal{A}\Bigr)
=\H_\mu\Bigl(\bigvee_{i=0}^{n-1}T^i\mathcal{A}\Bigr)
=\H_\mu\Bigl(T^{-(n-1)}\Bigl(\bigvee_{i=0}^{n-1}T^{i}\mathcal{A}\Bigr)\Bigr)\\
&=&\H_\mu\Bigl(\bigvee_{i=0}^{n-1}T^{-(n-1-i)}\mathcal{A}\Bigr)\\
&=&\H_\mu\Bigl(\bigvee_{j=0}^{n-1}T^{-j}\mathcal{A}\Bigr)
=\H_\mu(\mathcal{A}_T^n).
\end{eqnarray*}
It follows that $\h_\mu(T^{-1},\mathcal{A})=\h_\mu(T,\mathcal{A})$ for every
partition $\mathcal{A}$ and thus, passing to the supremum on both sides,
we conclude that $\h_\mu(T^{-1})=\h_\mu(T)$.
\epf

\sp Our goal now is to provide tools for actual calculating the entropy of an endomorphism. Its very definition requires to take the supremum over a huge set of all finite partitions. Our task is to reduce this to some sequences of partitions or even to a single partition. The following result, towards this end, is purely measure-theoretical. It says that given a
finite partition $\mathcal{A}$ of a compact metric space $X$ and given
any partition 
$\mathcal{C}$ of $X$ of sufficiently small diameter, we can group the atoms of
$\mathcal{C}$ together in such a way that we nearly construct
partition $\mathcal{A}$. 
It is worth noticing that $\mathcal{C}$ may be countably infinite.

\sp\blem\label{sem2lem7.5}
Suppose that $\mu$ is a Borel probability measure on a compact metric space $X$.
Suppose further that $\mathcal{A}=\{A_1,A_2,\ldots,A_n\}$ is a finite
partition of $X$ 
into Borel sets. Then for all $\varepsilon >0$ there exists $\delta
>0$ so that for every 
Borel partition $\mathcal{C}$ with $\mathrm{diam}(\mathcal{C})<\delta$
there is a Borel partition 
$\mathcal{B}=\{B_1,B_2,\ldots,B_n\}\leq\mathcal{C}$ such that
$$
\mu(B_i\Delta A_i)<\varepsilon
$$ 
for all $1\leq i\leq n$.
\elem

\bpf
Fix $\varepsilon>0$. Since $\mu$ is regular, there exists for each
$1\leq i\leq n$ a 
compact set $K_i\subseteq A_i$ such that $\mu(A_i\backslash K_i)<\varepsilon/n$.
As usual, let $d$ denote the metric on $X$ and let
$$
\theta:=\min\{d(K_i,K_j):i\neq j\}.
$$ 
Then $\theta>0$, as the sets $K_i$ are compact and disjoint. Let
$\delta=\theta/2$ and let $\mathcal{C}$ be a partition with
$\mathrm{diam}(\mathcal{C})<\delta$. 
For each $1\leq i\leq n$, define
\[
B_i:=\bigcup_{{C\in\mathcal{C}}\atop {C\cap K_i\neq\emptyset}}C.
\]
Clearly, the $B_i$s are Borel sets and $B_i\supseteq K_i$ for each $i$.
Moreover, due to the choice of $\delta$,  $B_i\cap B_j=\emptyset$ for
all $i\neq j$. 
However, the family of pairwise disjoint Borel sets $\{B_i\}_{i=1}^n$
may not cover $X$ completely. Indeed, there may be
some sets $C\in\mathcal{C}$ such that
$C\cap\bigcup_{i=1}^nK_i=\emptyset$. Simply take 
all those sets and put them into one of the $B_i$'s, say $B_1$. Then
the resulting family 
$\mathcal{B}:=\{B_i\}_{i=1}^n$ is a Borel partition of $X$. Clearly,
$\mathcal{B}\leq\mathcal{C}$. It remains to show that 
$$
\mu(B_i\Delta A_i)<\varepsilon
$$ 
for all $1\leq i\leq n$. But
\begin{eqnarray*}
\mu(B_i\Delta A_i)
&=&\mu(B_i\backslash A_i)+\mu(A_i\backslash B_i)
=\mu\lt(\lt(X\backslash\bigcup_{j\neq i}B_j\rt)\backslash A_i\rt)+\mu(A_i\backslash K_i)\\
&\leq&\mu\lt(\lt(X\backslash\bigcup_{j\neq i}K_j\rt)\backslash A_i\rt)+\mu(A_i\backslash K_i) \\
&=&\mu\lt(\lt(\bigcup_{k=1}^n A_k\backslash\bigcup_{j\neq i}K_j\rt)\backslash A_i\rt)+\mu(A_i\backslash K_i)\\
&=&\mu\lt(\bigcup_{k\neq i}A_k\backslash\cup_{j\neq i}K_j\rt)+\mu(A_i\backslash K_i)\\
&\leq&\mu\lt(\bigcup_{j\neq i}A_j\backslash K_j\rt) +\mu(A_i\backslash K_i)\\
&=&\sum_{j=1}^n\mu(A_j\backslash
K_j)<n\cdot\frac{\varepsilon}{n}=\varepsilon.
\end{eqnarray*}
\epf

\sp From the above result we will show that the conditional entropy of a partition $\mathcal{A}$
given a partition $\mathcal{C}$ is as small as desired provided $\mathcal{C}$ has a small
enough diameter. Indeed, from Theorem \ref{sem2thm7.1} $(v)$, given partitions $\mathcal{A}, \mathcal{B}$ and $\mathcal{C}$ as in the above lemma, we have that
$\H_\mu(\mathcal{A}|\mathcal{C})\leq \H_\mu(\mathcal{A}|\mathcal{B})$, where the partition
$\mathcal{B}$ is designed to resemble the partition $\mathcal{A}$. In order to estimate the conditional
entropy $\H_\mu(\mathcal{A}|\mathcal{B})$, we must estimate the contribution of all atoms of
the partition $\mathcal{A}\vee\mathcal{B}$. There are essentially two kinds of atoms to be
taken into account, namely,  atoms of the form $A_i\cap B_i$ and atoms of  the form $A_i\cap B_j$, for $i\neq j$.
Intuitively, because $A_i$ looks like $B_i$ (after all, $\mu(A_i\Delta B_i)$ is small),
the information provided by $A_i$ given that measurement $\mathcal{B}$ resulted in $B_i$
is small. On the other hand, since $A_i$ is nearly disjoint from $B_j$ when $i\neq j$ (after
all, $A_i$ is close to $B_i$ and $B_i\cap B_j=\emptyset$), the information obtained from
getting $A_i$ given that observation $B_j$ occurred is small. This is what
we now prove rigorously. First, let us make one definition which will prove useful
here and also in the sequel. Recall that a function $\psi:(a,b)\to\R$,
where $-\infty\leq a<b\leq\infty$, is concave if and only if
$$
\psi\bigl(tx+(1-t)y\bigr)\geq t\psi(x)+(1-t)\psi(y)
$$ 
for all $t\in[0,1]$ and all $x,y\in(a,b)$.

\sp\bdfn\label{functionk}
Let the function $k:[0,1]\to[0,1]$ be defined by
\[
k(t):=\begin{cases}
0  &{\rm if} \ x=0, \\
          -t\log t &{\rm if} \ t\in(0, 1].
        \end{cases}
\]
\edfn

\sp\fr Note that the function $k$ is continuous, concave, is increasing
on the interval $[0,e^{-1}]$ and 
decreasing on the interval $[e^{-1},1]$.

\sp\blem\label{sem2lem7.6}
Let $\mu$ be a Borel probability measure on a compact metric space $X$, and
let $\mathcal{A}$ be a finite Borel partition of $X$. Then
for every $\varepsilon>0$, there exists $\delta>0$ such that
$$
\H_\mu(\mathcal{A}|\mathcal{C})<\varepsilon
$$ 
for every Borel partition
$\mathcal{C}$ with $\mathrm{diam}(\mathcal{C})<\delta$.
\elem

\bpf Let $\mathcal{A}=\{A_1,A_2,\ldots,A_n\}$ be a finite
  partition of $X$. 
 Fix $\varepsilon>0$ and let
$0<\overline{\varepsilon}<\min\{e^{-1},1-e^{-1}\}$ be so
small that 
$$
\max\big\{k(\overline{\varepsilon}),k(1-\overline{\varepsilon})\big\}<\varepsilon/n.
$$
Then there exists $\tilde{\varepsilon}>0$ such that
\[
0<\frac{\tilde{\varepsilon}}{\mu(A_i)-\tilde{\varepsilon}}<\overline{\varepsilon}
\ \ \ \text{ and }\ \ \
\frac{\mu(A_i)-\tilde{\varepsilon}}{\mu(A_i)+\tilde{\varepsilon}}>1-\overline{\varepsilon}.
\]
for all $1\leq i\leq n$ such that $\mu(A_i)>0$. Let $\delta>0$ be the number ascribed to
$\tilde{\varepsilon}$ in Lemma \ref{sem2lem7.5}. Let $\mathcal{C}$ be a partition with
$\mathrm{diam}(\mathcal{C})<\delta$, and let 
$$
\mathcal{B}=\{B_1,B_2,\ldots,B_n\}\leq\mathcal{C}
$$
be such that 
$$
\mu(A_i\Delta B_i)\leq\tilde{\varepsilon}
$$ 
for all $1\leq i\leq n$, also as prescribed
in Lemma \ref{sem2lem7.5}. Then, for all $1\leq i\leq n$ we have that
\[
|\mu(A_i)-\mu(B_i)|\leq\mu(A_i\Delta B_i)\leq\tilde{\varepsilon}.
\]
Therefore,
\begin{eqnarray*}
0<
\mu(A_i)-\tilde{\varepsilon}
&\leq&\mu(A_i)-\mu(A_i\Delta B_i)
\leq\mu(B_i)\\&\leq&\mu(A_i)+\mu(A_i\Delta B_i)
\leq\mu(A_i)+\tilde{\varepsilon}
\end{eqnarray*}
for all $i$. Moreover,
\[
\mu(A_i\cap B_i)
=\mu(A_i)-\mu(A_i\backslash B_i)
\geq\mu(A_i)-\mu(A_i\Delta B_i)
\geq\mu(A_i)-\tilde{\varepsilon}
\]
for all $i$. Hence
\[
\frac{\mu(A_i\cap B_i)}{\mu(B_i)}
\geq\frac{\mu(A_i)-\tilde{\varepsilon}}{\mu(A_i)
  +\tilde{\varepsilon}}>1-\overline{\varepsilon}.    
\]
for all $i$ such that $\mu(A_i)>0$. By our choice of $\overline{\varepsilon}$,
the function $k$ is decreasing on the interval
$[1-\overline{\varepsilon},1]$ and thus 
\beq\label{k1}
k\lt(\frac{\mu(A_i\cap B_i)}{\mu(B_i)}\rt)\leq
k(1-\overline{\varepsilon})<\frac{\varepsilon}{n} 
\eeq
for all $i$ such that $\mu(A_i)>0$.
Suppose now that $i\neq j$. Since $\mathcal{A}=\{A_k\}_{k=1}^n$ is a
partition of $X$, we know that 
$A_i\cap B_j\subseteq B_j\backslash A_j\subseteq A_j\Delta B_j$. Hence, 
\[
\frac{\mu(A_i\cap B_j)}{\mu(B_j)}
\leq\frac{\mu(A_j\Delta B_j)}{\mu(A_j)-\mu(A_j\Delta B_j)}
\leq\frac{\tilde{\varepsilon}}{\mu(A_j)-\tilde{\varepsilon}}
<\overline{\varepsilon}
\]
for all $i$ such that $\mu(A_i)>0$. By our choice of $\overline{\varepsilon}$,
the function $k$ is increasing on the interval
$[0,\overline{\varepsilon}]$ and hence 
\beq\label{k2}
k\lt(\frac{\mu(A_i\cap B_j)}{\mu(B_j)}\rt)\leq
k(\overline{\varepsilon})<\frac{\varepsilon}{n} 
\eeq
for all $i$ such that $\mu(A_i)>0$. Furthermore, note that
$k\bigl(\frac{\mu(A_i\cap B_j)}{\mu(B_j)}\bigr)=0$ 
if $\mu(A_i)=0$. Then, by 
Theorem~\ref{sem2thm7.1}$(v)$ and $(\ref{k1})$ and $(\ref{k2})$, we have
\begin{eqnarray*}
\H_\mu(\mathcal{A}|\mathcal{C})
&\leq&\H_\mu(\mathcal{A}|\mathcal{B})
=\sum_{A\in\mathcal{A}}\sum_{B\in\mathcal{B}}-\mu(A\cap B)\log\frac{\mu(A\cap B)}{\mu(B)}\\
&=&\sum_{i,j=1}^n\mu(B_j)k\lt(\frac{\mu(A_i\cap B_j)}{\mu(B_j)}\rt)\\
&=&\sum_{i=1}^n\mu(B_i)k\lt(\frac{\mu(A_i\cap B_i)}{\mu(B_i)}\rt)
+\sum_{{\begin{array}{c}i,j=1\\i\neq j\end{array}}}^n\mu(B_j)k\lt(\frac{\mu(A_i\cap B_j)}{\mu(B_j)}\rt)\\
&<&\sum_{i=1}^n\mu(B_i)\frac{\varepsilon}{n}
+\sum_{i=1}^n\sum_{j=1}^n\mu(B_j)\frac{\varepsilon}{n}\\
&=&\frac{\varepsilon}{n}+n\cdot\frac{\varepsilon}{n}\\
&<&2\varepsilon.
\end{eqnarray*}
 \epf

\sp From the above lemma, we can infer that any sequence of partitions
whose diameters 
tend to $0$ provide asymptotically as much information as any given
partition can.

\sp\bcor\label{sem2cor}
Let $\mu$ be a Borel probability measure on a compact metric space $X$.
If $(\mathcal{A}_n)_{n\geq1}$ is a sequence of Borel partitions of $X$ such that
$$
\lim_{n\to\infty}\mathrm{diam}(\mathcal{A}_n)=0,
$$
then
$$
\lim_{n\to\infty}\H_\mu(\mathcal{A}|\mathcal{A}_n)=0
$$ for every finite Borel partition~$\mathcal{A}$ of $X$.
\ecor

\bpf
Let $\mathcal A$ be a finite Borel partition of $X$. Then, by
Lemma~\ref{sem2lem7.6}, 
for every $\varepsilon>0$ there exists a $\delta>0$ such that if
$\mathrm{diam}(\mathcal{C})<\delta$ then
$\H_\mu(\mathcal{A}|\mathcal{C})<\varepsilon$. 
Since $\mathrm{diam}(\mathcal{A}_n)\to 0$, it follows that
$\H_\mu(\mathcal{A}|\mathcal{A}_n)\to0$.
\epf

\sp This result about conditional entropy of partitions allows us to
deduce the following fact about entropy of endomorphisms.

\sp\bthm\label{sem2thm7.9}
Let $\mu$ be a Borel probability measure on a compact metric space $X$.
If $T:X\to X$ is a measure--preserving transformation of $(X,\mathcal{B}(X),\mu)$
and $(\mathcal{A}_n)_{n\geq1}$ is a sequence of finite Borel partitions of $X$
such that 
$$
\lim_{n\to\infty}\mathrm{diam}(\mathcal{A}_n)=0,
$$
then
$$
\h_\mu(T)=\lim_{n\to\infty}\h_\mu(T,\mathcal{A}_n).
$$
\ethm

\bpf
Let $\mathcal A$ be a finite partition consisting of Borel sets. By
Corollary~\ref{sem2cor}, we know that 
$\lim_{n\to\infty}\H_\mu(\mathcal{A}|\mathcal{A}_n)=0$. So, by
Theorem~\ref{sem2thm7.3}$(9)$, it follows that 
\[
\h_\mu(T,\mathcal{A})
\leq\liminf_{n\to\infty}\h_\mu(T,\mathcal{A}_n)
\leq\limsup_{n\to\infty}\h_\mu(T,\mathcal{A}_n)
\leq \H_\mu(T).
\]
Since this is true for any finite Borel partition $\mathcal{A}$ of $X$,
we deduce from a passage to the supremum that
\[
\h_\mu(T)
\leq\liminf_{n\to\infty}\h_\mu(T,\mathcal{A}_n)
\leq\limsup_{n\to\infty}\h_\mu(T,\mathcal{A}_n)
\leq \h_\mu(T).
\]
\epf

\bcor\label{sem2cor7.110}
Let $\mu$ be a Borel probability measure on a compact metric space $X$.
If $T:X\lra X$ is a measure--preserving transformation of $(X,\mathcal{B}(X),\mu)$ and $\mathcal{A}$ is a finite Borel
partition of $X$ such that
$$
\lim_{n\to\infty}\mathrm{diam}(\mathcal{A}^n)=0,
$$ 
then 
$$
\h_\mu(T)=\h_\mu(T,\mathcal{A}).
$$
\ecor

\bpf
By Theorems~\ref{sem2thm7.9} and~\ref{sem2thm7.3}$(6)$, we have that
\[
\h_\mu(T)
=\lim_{n\to\infty}\h_\mu(T,\mathcal{A}^n)
=\lim_{n\to\infty}\h_\mu(T,\mathcal{A})
=\h_\mu(T,\mathcal{A}).
\]
\epf

\sp We now will introduce a notion, very classical and
interesting in itself, which guarantees the hypotheses of Corollary~\ref{sem2cor7.110} to be satisfied and which will play an important
role later on, notably when dealing with variational principle and
equilibrium states.

\sp\bdfn
Let $(X,d)$ be a compact metric space. A continuous dynamical system
$T:(X,d)\lra (X,d)$ is said to be \index{(N)}{positively expansive}positively expansive provided that there exists $\d>0$ such that for every $x,y\in X$, $x\neq y$ there exists an integer $n=n(x,y)\ge 0$ with
\[
d\bigl(T^n(x),T^n(y)\bigr)>\d.
\]
The constant $\d$ is called an expansive constant for $T$ and
$T$ is then also said to be $\d$--expansive. Equivalently, $T$ is $\d$--expansive if
\[
\sup\Big\{d\bigl(T^n(x),T^n(y)\bigr):n\geq 0\Big\}\leq\d \ \Longrightarrow \ x=y.
\]
In other words, $\d$--expansiveness means that two forward $T$-orbits
that remain forever within 
a distance $\d$ from each other originate from the same point (and are therefore
only one orbit).
\edfn

\sp

\brem
Let us record the following.
\begin{itemize} 
\item[(1)] If $T$ is $\d$--expansive, then $T$ is ${\d'}$--expansive for any $0<{\d'}<\d$.

\,

\item[(2)] The expansiveness of $T$ is independent of topologically
  equivalent metrics, although particular expansive 
constants generally depend on the metric chosen. That is, if two metrics $d$ and
$d'$  generate the same topology on $X$,
then $T$ is expansive when $X$ is equipped with the metric $d$ if and
only if $T$ is expansive when $X$ is equipped with the metric $d'$.
\end{itemize}
\erem

\sp We record now the following fact which will follow from a somewhat stronger fact, namely Corollary~\ref{c120190930} proven in the next chapter, asserting the same, but for covers and not merely partitions. 

\bprop\label{expansive_generators}
If $(X,d)$ is a compact metric space and $T:X\lra X$ is a positively
expansive continuous map with an expansive constant $\d>0$, then 
$$
\lim_{n\to\infty}\mathrm{diam}(\mathcal{A}^n)=0
$$
for every partition $\mathcal{A}$ with $\mathrm{diam}(\mathcal{A})\leq\delta$.
\eprop

\sp\fr As an immediate consequence of this proposition and
Corollary~\ref{sem2cor7.110}, we get the following.

\bcor\label{sem2thm7.11}
Let $T:X\lra X$ be an expansive dynamical system preserving a Borel
probability measure $\mu$. 
If $\mathcal{A}$ is a finite partition with
$\mathrm{diam}(\mathcal{A})\leq\delta$, where 
$\delta>0$ is an expansive constant for $T$, then 
$$
\h_\mu(T)=\h_\mu(T,\mathcal{A}).
$$
\ecor

\

\section{Shannon--McMillan--Breiman Theorem}

\sp\fr The sole goal of this section is to prove 
Shannon--McMillanBreiman Theorem, i.e. Theorem~\ref{sem2thm8.4}. It sheds lots of light on what entropy really is and provides both  a very useful tool for further theoretical investigations of entropy and for its actual calculations.

We begin this section with the following  purely measure--theoretical result.

\sp\blem\label{sem2lem8.1}
Let $T:X\lra X$ be a measure--preserving endomorphism of a probability
space $(X,\mathfrak{ F},\mu)$. 
Let $\mathcal{A}\in\mathrm{Part}(X,\mathfrak{ F})$. Let
$$
f_n:=I_\mu(\mathcal{A}|\mathcal{A}_{1}^{n})
$$ 
for each $n\geq1$ and let 
$$
f^*:=\sup_{n\geq1}f_n.
$$
Then, for all $\lambda\in\R$ and all $A\in\mathcal{A}$, we have that 
$$
\mu\big(\big\{x\in A:f^*(x)>\lambda\big\}\)\leq\min\{\mu(A),e^{-\lambda}\}. 
$$
\elem

\bpf
Let $A\in\mathcal{A}$ and fix $n\geq1$. Let also 
$$
f^A_n:=-\log
E(\chi_A|\sigma(\mathcal{A}^n_1)),
$$ 
where $\sigma(\mathcal{A}^n_1)$ is the sub-$\sigma$-algebra generated
by the countable 
partition $\mathcal{A}^n_1$. Let $x\in A$. Then, 
\begin{eqnarray*}
\begin{aligned}
f^A_n(x)
&=-\log E(\chi_A|\sigma(\mathcal{A}^n_1))(x)
 =-\log\frac{\int_{\mathcal{A}^n_1(x)}\chi_A\,d\mu}{\mu(\mathcal{A}^n_1(x))}\\
&=-\log\frac{\mu(A\cap\mathcal{A}^n_1(x))}{\mu(\mathcal{A}^n_1(x))}\\
&=-\log\frac{\mu\bigl(\mathcal{A}(x)\cap\mathcal{A}^n_1(x)\bigr)} 
   {\mu(\mathcal{A}^n_1(x))}\\
&= I_\mu(\mathcal{A}|\mathcal{A}^n_1)(x)\\
&=f_n(x).
\end{aligned}
\end{eqnarray*}
Hence, $f_n=\displaystyle\sum_{A\in\mathcal{A}}\chi_Af_n^A$.
Fix $A\in\mathcal{A}$ and for $n\geq1$ and $\lambda\in\R$ consider the set
\[
B_n^{A,\lambda}:=\bigl\{x\in X:\max_{1\leq i<n} f_i^A(x)\leq\lambda,f_n^A(x)>\lambda\bigr\}.
\]
The family $\{B_n^{A,\lambda}\}_{n\geq1}$ consists of mutually disjoint sets. Also, recall that
$\mathcal{A}_1^n\leq\mathcal{A}_1^{n+1}$, and thus
$\sigma(\mathcal{A}_1^n)\subseteq\sigma(\mathcal{A}_1^{n+1})$, for each $n\geq1$.
By definition, each function $f_n^A$ is measurable with respect to
$\sigma(\mathcal{A}^n_1)$. Consequently, $B_n^{A,\lambda}\in\sigma(\mathcal{A}^n_1)$. Then
\begin{eqnarray*}
\begin{aligned}
\mu(B_n^{A,\lambda}\cap A)
&=\int_{B_n^{A,\lambda}}\chi_A\,d\mu=\int_{B_n^{A,\lambda}}E(\chi_A|\sg(\mathcal{A}^n_1))\,d\mu\\
&=\int_{B_n^{A,\lambda}}\exp(-f_n^A)\,d\mu\leq\int_{B_n^{A,\lambda}}e^{-\lambda}\,d\mu\\
&=e^{-\lambda}\mu(B_n^{A,\lambda}).
\end{aligned}
\end{eqnarray*}
Since 
$$
\begin{aligned}
\mu(\{x\in A:f^*(x)>\lambda\})
&=\mu(\{x\in A:\exists n\geq1 \mbox{ such that } f_n(x)>\lambda\}) \\
&=\mu(\{x\in A:\exists n\geq1 \mbox{ such that } f^A_n(x)>\lambda\}),
\end{aligned}
$$
we have that
\begin{eqnarray*}
\begin{aligned}
\mu(\{x\in A:f^*(x)>\lambda\})
&=\mu\Bigl(\bigcup_{n=1}^\infty B_n^{A,\lambda}\cap A\Bigr)
=\sum_{n=1}^\infty\mu(B_n^{A,\lambda}\cap A)\\
&\leq \sum_{n=1}^\infty e^{-\lambda}\mu(B_n^{A,\lambda})
=e^{-\lambda}\sum_{n=1}^\infty\mu(B_n^{A,\lambda})\\
&=e^{-\lambda}\mu\Bigl(\bigcup_{n=1}^\infty B_n^{A,\lambda}\Bigr) \\
&\leq e^{-\lambda}.
\end{aligned}
\end{eqnarray*}
\epf

\sp\bcor\label{sem2cor8.2}
Let $T:X\to X$ be a measure--preserving endomorphism of a probability
space $(X,\mathfrak{ F},\mu)$. Let $\mathcal{A}$ be a partition of $X$ with finite entropy. Let 
$$
f_n:=I_\mu(\mathcal{A}|\mathcal{A}_{1}^{n})
$$ 
for all
$n\geq1$ and let 
$$
f^*:=\sup_{n\geq1}f_n.
$$
Then the function $f^*$ belongs to $L^1(X,\mathfrak{ F},\mu)$ and
$$
\int_{X}f^*\,d\mu\leq \H_\mu(\mathcal{A})+1.
$$
\ecor

\bpf
Since $f^*\geq0$, we have $\int_{X}|f^*|\,d\mu=\int_{X}f^*\,d\mu$. Thus,
\begin{eqnarray*}
\int_X f^*\,d\mu
&=&\sum_{A\in\mathcal{A}}\int_A f^*\,d\mu
=\sum_{A\in\mathcal{A}}\int_0^\infty\mu(\{x\in A:f^*(x)>\lambda\})\,d\lambda\\
&\leq&\sum_{A\in\mathcal{A}}\int_0^\infty\min\{\mu(A),e^{-\lambda}\}\,d\lambda\\
&=&\sum_{A\in\mathcal{A}}\left[\int_0^{-\log\mu(A)}\mu(A)\,d\lambda+\int_{-\log\mu(A)}^\infty e^{-\lambda}\,d\lambda\right]\\
&=&\sum_{A\in\mathcal{A}}\Bigl[-\mu(A)\log\mu(A)+\left[-e^{-\lambda}\right]_{-\log\mu(A)}^\infty\Bigr]\\
&=&\sum_{A\in\mathcal{A}}-\mu(A)\log\mu(A)+\sum_{A\in\mathcal{A}}\mu(A)\\
&=&\H_\mu(\mathcal{A})+1<+\infty.
\end{eqnarray*}
\epf

\sp\bcor\label{sem2cor8.3}
The sequence $(f_n)_{n\geq1}$, defined in the previous corollary, converges $\mu$-a.e. and also in
$L^1(X,\mathfrak{ F},\mu)$. 
\ecor

\bpf
Recall that $\mathcal{A}_1^n\leq\mathcal{A}_1^{n+1}$, and thus
$\sigma(\mathcal{A}_1^n)\subseteq\sigma(\mathcal{A}_1^{n+1})$, for each $n\geq1$.
For any $x\in A\in\mathcal{A}$, we have 
$$
f_n(x)=f_n^A(x)=-\log E_\mu(\chi_A|\sg(\mathcal{A}^n_1))(x)
$$
and a version of Doob's Martingale Convergence Theorem for expected values (see for example Theorem~35.6 in \cite{Bi}) guarantees that the 
limit 
$$
\lim_{n\to\infty}E_\mu(\chi_A|\mathcal{A}^n_1)
$$ 
exists $\mu$-a.e.. Hence the sequence of non--negative functions
$(f_n)_{n\geq1}$ converges to some limit function $\hat{f}\geq 0$ $\mu$--a.e.. 

Since $|f_n|=f_n\leq f^*$ for all
$n$, we have $|\hat{f}|=\hat{f}\leq f^*$ and thus 
$$
|f_n-\hat{f}|\leq2 f^*
$$
$\mu$--a.e.. So, by applying Lebesgue's Dominated Convergence Theorem to the sequence $(|f_n-\hat{f}|)_{n\geq1}$, we obtain
\[
\lim_{n\to\infty}\|f_n-\hat{f}\|_1
=\lim_{n\to\infty}\int_X|f_n-\hat{f}|\,d\mu
=\int_X\lim_{n\to\infty}|f_n-\hat{f}|\,d\mu
=0
\]
i.e. $f_n\to\hat{f}$ in $L^1(X,\mathfrak{ F},\mu)$.
\epf

\sp We can now prove the main result of this section, the famous
Shannon--McMillan--Breiman Theorem. \index{(N)}{Shannon--McMillan--Breiman Theorem}

\sp\bthm[Shannon--McMillan--Breiman Theorem]\label{sem2thm8.4}

Let $T:X\lra X$ be a measure-preserving endomorphism of a probability
space $(X,\mathfrak{ F},\mu)$. 
Let $\mathcal{A}$ be a partition of $X$ with finite entropy. Then the following limits exist
\[
f:=\lim_{n\to\infty}I_\mu(\mathcal{A}|\mathcal{A}_{1}^{n})
\  \ {\rm  and } \  \
\lim_{n\to\infty}\frac{1}{n}\sum_{j=0}^{n-1}f\circ T^j=E_\mu(f|\mathcal{I}_\mu)\ \text{ $\mu$-a.e.}.
\]
Moreover, the following hold:
\begin{itemize}
   \item[(1)] $\displaystyle\lim_{n\to\infty}\frac{1}{n}I_\mu(\mathcal{A}^n)=E_\mu(f|\mathcal{I}_\mu)$ in $L^1(\mu)$ and $\mu$--a.e.
   
   \,
   
   \item[(2)] $\displaystyle \h_\mu(T,\mathcal{A})=\lim_{n\to\infty}\frac{1}{n}\H_\mu(\mathcal{A}^n)
                                    =\int_X E_\mu(f|\mathcal{I}_\mu)\,d\mu=\int_X f\,d\mu$.
\end{itemize}
\ethm

\bpf
The first sequence of functions 
$$
(f_n)_{n\geq1}=(I_\mu(\mathcal{A}|\mathcal{A}_{1}^{n}))_{n\geq1}
$$
converges to an integrable function $f$ by Corollary~\ref{sem2cor8.3}. The second limit exists by virtue of
Birkhoff's Ergodic Theorem. 

In order to prove the remaining two statements, let us first assume that $(1)$ holds and derive $(2)$ from it. Then we will prove $(1)$.
In fact, a.e. convergence in $(i)$ is not necessary to deduce
$(2)$. So, suppose that
$$
\lim_{n\to\infty}\frac{1}{n}I_\mu(\mathcal{A}^n)=E_\mu(f|\mathcal{I}_\mu)
$$ 
in $L^1(\mu)$.
The convergence in $L^1$ entails the convergence of the corresponding integrals, that is,
$$
\lim_{n\to\infty}\int_X\frac{1}{n}I_\mu(\mathcal{A}^n)\,d\mu=\int_X E_\mu(f|\mathcal{I}_\mu)\,d\mu.
$$
Then
\begin{eqnarray*}
\h_\mu(T,\mathcal{A})
=\lim_{n\to\infty}\frac{1}{n}\H_\mu(\mathcal{A}^n)
=\lim_{n\to\infty}\int_X\frac{1}{n}I_\mu(\mathcal{A}^n)\,d\mu
=\int_X E_\mu(f|\mathcal{I}_\mu)\,d\mu=\int_X f\,d\mu.
\end{eqnarray*}
This establishes $(2)$. 

In order to prove $(1)$, first notice that by
Lemma~\ref{sem2lem7.2i} we obtain that 
\[
I_\mu(\mathcal{A}^n)
=\sum_{k=1}^n I_\mu(\mathcal{A}|\mathcal{A}^k_1)\circ T^{n-k}\\
=\sum_{j=0}^{n-1}I_\mu(\mathcal{A}|\mathcal{A}_1^{n-j})\circ T^j\\
=\sum_{j=0}^{n-1}f_{n-j}\circ T^j.
\]
Then, by the triangle inequality,
\begin{eqnarray*}
\Bigl|\frac{1}{n}I_\mu(\mathcal{A}^n)-E_\mu(f|\mathcal{I}_\mu)\Bigr|
&=&\left|\frac{1}{n}\sum_{j=0}^{n-1}(f_{n-j}\circ T^{j}-f\circ T^{j})\right.\\
&&\hspace{4cm}+\left.\frac{1}{n}\sum_{j=0}^{n-1}f\circ T^{j}-E(f|\mathcal{I}_\mu)\right|\\
&\leq&\Bigl|\frac{1}{n}\sum_{j=0}^{n-1}(f_{n-j}-f)\circ T^{j}\Bigr|
   +\Bigl|\frac{1}{n}S_n f-E(f|\mathcal{I}_\mu)\Bigr|\\
&\leq&\frac{1}{n}\sum_{j=0}^{n-1}|f_{n-j}-f|\circ T^{j}
   +\Bigl|\frac{1}{n}S_n f-E_\mu(f|\mathcal{I}_\mu)\Bigr|.
\end{eqnarray*}
The second term on the right-hand side tends to $0$ $\mu$-a.e. by
Birkhoff's Ergodic Theorem. Furthermore, observe that
\beq
\aligned
\int_X\frac{1}{n}S_nf\,d\mu
&=\frac{1}{n}\sum_{j=0}^{n-1}\int_X f\circ T^{j}\,d\mu
 =\frac{1}{n}\sum_{j=0}^{n-1}\int_X f\,d\mu \\
&= \int_X f\,d\mu\\
&=\int_X \hspace{-0.1cm}E_\mu(f|\mathcal{I}_\mu)d\mu.
\endaligned
\eeq
Thus, the second term converges to $0$ in $L^1(\mu)$.
Let us now investigate the first term on the right-hand side.
Set $g_n:=|f_n-f|$. Then
\begin{eqnarray*}
\lim_{n\to\infty}\Bigl\|\frac{1}{n}\sum_{j=0}^{n-1}g_{n-j}\circ T^j-0\Bigr\|_1
&=&\lim_{n\to\infty}\int_X\Bigl|\frac{1}{n}\sum_{j=0}^{n-1}g_{n-j}\circ T^j-0\Bigr|\,d\mu \\
&=&\lim_{n\to\infty}\frac{1}{n}\sum_{j=0}^{n-1}\int_X g_{n-j}\circ T^j\,d\mu
=\lim_{n\to\infty}\frac{1}{n}\sum_{j=0}^{n-1}\int_X g_{n-j}\,d\mu\\
&=&\lim_{n\to\infty}\frac{1}{n}\sum_{i=1}^{n-1}\int_X g_{i}\,d\mu\\
&=&0,
\end{eqnarray*}
where the last equality sign was written since $f_i\to f$ in $L^1(X,\mathfrak{ F},\mu)$ according to Corollary~\ref{sem2cor8.3},
whence $g_i\to0$ in $L^1(X,\mathfrak{ F},\mu)$ and so do the Ces\`aro averages of the $g_i$s.
This ensures the convergence of the functions 
$$
\frac{1}{n}\sum_{j=0}^{n-1}g_{n-j}\circ T^j
$$
to zero in $L^1$ and thus the convergence of the first term on the right--hand side to $0$ in
$L^1$. It only  remains to show convergence $\mu$--a.e. of the same term. To this end, for each $N\geq1$ let
$$
G_N:=\sup_{n\geq N}g_n.
$$
 The sequence of functions $(G_N)_{N\geq1}$ is decreasing and bounded
below by 0, so it converges to some function. As $f_n\to f$ $\mu$-a.e., we know that
$g_n=|f_n-f|\to0$ $\mu$-a.e.. It follows that 
$$
G_N\xrightarrow[\ N\to\infty \ ]{} 0
$$ 
$\mu$--a.e.. Also, the functions $G_N$ are bounded above by an integrable function since
\[
0\leq G_N\leq G_1=\sup_{n\geq1}g_n\leq\sup_{n\geq1}\(|f_n|+|f|\)\leq2f^* \in L^1(\mu).
\]
Fix momentarily $N\geq1$. Then for any $n>N$ we have
\begin{eqnarray*}
\frac{1}{n}\sum_{j=0}^{n-1}g_{n-j}\circ T^j
&=&\frac{1}{n}\sum_{j=0}^{n-N}g_{n-j}\circ T^j
+\frac{1}{n}\sum_{j=n-N+1}^{n-1}g_{n-j}\circ T^j\\
&\leq&\frac{n-N}{n}\cdot\frac{1}{n-N}\sum_{j=0}^{n-N}G_N\circ T^j
+\frac{1}{n}\sum_{j=n-N+1}^{n-1}G_1\circ T^j.
\end{eqnarray*}
Let $F_N=\sum_{j=0}^{N-2}G_1\circ T^j$. Using Birkhoff's Ergodic Theorem,
we deduce that
\begin{eqnarray*}
\limsup_{n\to\infty}\frac{1}{n}\sum_{j=0}^{n-1}g_{n-j}\circ T^j
&\leq&\lim_{n\to\infty}\frac{1}{n-N}\sum_{j=0}^{n-N}G_N\circ T^j
    +\limsup_{n\to\infty}\frac{1}{n}F_N\circ T^{n-N+1}\\
&=&E(G_N|\mathcal{I}_\mu)+\limsup_{n\to\infty}\frac{1}{n}F_N\circ T^{n-N+1}\ \ \ \mu\mbox{-a.e.}\\
&=&E(G_N|\mathcal{I}_\mu)\ \ \ \mu\mbox{-a.e.}.
\end{eqnarray*}
But since each $G_N$ is a  non-negative function
uniformly bounded by $G_0\in L^1(\mu)$ and since the sequence $(G_N)_{N\geq 1}$ converges to zero $\mu$-a.e.,
Lebesgue's Dominated Convergence Theorem implies that
\[
\lim_{N\to\infty}\int_X E(G_N|\mathcal{I}_\mu)\,d\mu=\lim_{N\to\infty}\int_X G_N\,d\mu
=\int_X\lim_{N\to\infty}G_N\,d\mu=0.
\]
However, the sequence of expected values $(E(G_N|\mathcal{I}_\mu))_{N\geq1}$ is decreasing since
the sequence $(G_N)_{N\geq1}$ is decreasing. Therefore, 
$$
\int_X E_\mu(G_N|\mathcal{I}_\mu)\,d\mu\downto 0,
$$
as $n\to\infty$, and hence 
$$
E_\mu(G_N|\mathcal{I}_\mu)\xrightarrow[\ N\to\infty \ ]{} 0
$$ 
$\mu$-a.e.. Thus,
\[
\limsup_{n\to\infty}\frac{1}{n}\sum_{j=0}^{n-1}g_{n-j}\circ T^j=0\, \ \ \mu\mbox{--a.e.},
\]
thereby establishing the a.e. convergence of the first term on the
right--hand side. 
\epf
As an immediate consequence of this theorem and Birkhoff's  Ergodic
Theorem (Theorem~\ref{Birkhoff}), we get the following.

\sp\bcor[Ergodic case of Shannon-McMillan-Breiman Theorem]\label{SMBTEC}
Let $T:X\to X$ be an ergodic endomorphism of a probability space
$(X,\mathfrak{ F},\mu)$. 
Let $\mathcal{A}$ be a partition of $X$ with finite entropy. Then
\[ 
\hmu(T,\mathcal{A})=\lim_{n\to\infty}\frac{1}{n}I_\mu(\mathcal{A}^n)(x)\
\ \ \mbox{ for $\mu$--a.e. $x\in X$}. 
\]
\ecor  

\sp The right-hand side in the above equality can be viewed as a local
entropy at a point $x$. The theorem 
then states that at $\mu$-a.e. $x\in X$ the local entropy exists and is equal to
the (global) entropy of the endomorphism . 
Moreover, the theorem affirms that if $\H_\mu(T,\mathcal{A})>0$, then
$\mu(\mathcal{A}^n(x))\to0$ with exponential rate
$e^{-\H_\mu(T,\mathcal{A})}$ for $\mu$-a.e. $x\in X$. 

\

\section{Abramov's Formula and Krengel's Entropy \\ (Infinite Measures Allowed)}

\sp\fr In previous chapters we have devoted a good amount of time to study
induced maps. There is a celebrated Abramov's Formula which relates
the entropy of an induced system and the original one. It was originally proved by L. M. Abramov in \cite{Ab}. We quote if here without a proof.

\sp\bthm[Abramov's Formula]\label{abramovabstract}
 If $T:X\to X$ is an ergodic measure preserving transformation of a
 probability space $(X,\mathfrak{F},\mu)$, then for every set $F\in
 \mathfrak{F}$ with $ 0< \mu(F) < +\infty$, we have that
$$
\h_{\mu_F}(T_F)=\frac1{\mu(F)}\hmu(T).
$$
\ethm

\sp As an immediate consequence (take $X:=F\cup G$) of this theorem, we get
the following. 

\sp \bcor[Krengel's Entropy]\label{KrengelEntropy}
If $T:X\to X$ is a conservative ergodic measure preserving transformation of a
 measure space $(X,\mathfrak{F},\mu)$, then for all sets $F$ and $G$ in
$\mathfrak{F}$ with $0<\mu(F),\mu(G) < +\infty$, we have that
$$
\mu(F)\h_{\mu_F}(T_F)=\mu(G)\h_{\mu_G}(T_G).
$$
This common value is called the Krengel's entropy \index{(N)}{Krengel's entropy} of the map $T:X\lra
X$ and is denoted simply by $\hmu(T)$. In the case when the measure
$\mu$ is a probability one, it 
coincides with the standard entropy of $T$ with respect to $\mu$.
\ecor

\chapter{Thermodynamic Formalism}\label{Thermodynamic Formalism}

In this chapter we introduce the fundamental concepts of thermodynamic formalism such as topological pressure and topological entropy and we establish their basic properties. We then, in the last section of this chapter, relate them with Kolmogorov--Sinai metric entropies by proving the Variational Principle which is the corner stone of thermodynamic formalism. 
This principle naturally leads to the concept of equilibrium states and measures of maximal entropy. We deal with them at length, particularly with the problem of existence of equilibrium states. We do not touch the issue of its uniqueness as this requires a more involved and lengthy apparatus and holds only for some special systems such as open transitive distance expanding maps in the sense of the book \cite{PU2}. However, some considerations of Chapter touch on the issue of uniqueness though in somewhat different setting. 

Thermodynamic formalism originated in the late years 1960s with the works of David Ruelle. Its foundations, classical concepts and theorem were obtained throughout 1970s in the works of David Ruelle \cite {Ru}, Rufus Bowen \cite{Bow}, Peter Walters \cite{W1}, and Yakov Sinai \cite{Si2}. The more recent and modern expositions can be found for example in \cite{W}, \cite{PU2}, or \cite{Ru}. We should also mention the paper \cite{Mis} by Michal Misiurewicz where an elegant, short, and simple proof of the Variational Principle was provided. This is the proof we reproduce in the last section of this chapter.

\section{Topological Pressure}\label{ch10}

\subsection{Covers of a Set}\label{ch6topent}

Let $X$ be a nonempty set. Recall that a family ${\mathcal U}$ of subsets
of $X$ is said to form a \index{(N)}{cover of a set} of $X$ if and only if 
$$
X\subseteq\bigcup_{U\in{\mathcal U}}U.
$$ 
Recall further that ${\mathcal V}$ is
a said to be {\em subcover} of ${\mathcal U}$ if ${\mathcal V}$ is
itself a cover and 
$$
{\mathcal V}\subseteq{\mathcal U}.
$$ 
We will always denote covers by calligraphic letters, $\mathcal U, V,
W$ and so on. 

\sp Let us begin by introducing a useful way of obtaining a new
cover from two existing covers. 

\bdfn\label{ch6join}
If ${\mathcal U}$ and ${\mathcal V}$ are
covers of $X$, then their {\em join}, \index{(N)}{join of covers} denoted
${\mathcal U}\vee{\mathcal V}$, is the cover
\[
{\mathcal U}\vee{\mathcal V}:=\bigl\{U\cap V:{U\in{\mathcal
    U}},{V\in{\mathcal V}}\bigr\}. 
\]
\edfn

\brem\label{ch6assoc}
The join operation is commutative (that is, ${\mathcal U}\vee{\mathcal
  V}={\mathcal V}\vee{\mathcal U}$) 
and associative; in other words, 
$$
({\mathcal U}\vee{\mathcal
  V})\vee{\mathcal W}={\mathcal U}\vee({\mathcal V}\vee{\mathcal
  W}).
$$ 
\erem  
Thanks to the associativity of the join, this operation  extends
naturally to any 
finite collection $\{{\mathcal U}_j\}_{j=0}^{n-1}$ of covers of
$X$. That is,  we have 
\[
\bigvee_{j=0}^{n-1}{\mathcal U}_j
:={\mathcal U}_0\vee\cdots\vee{\mathcal U}_{n-1}
=\left\{\bigcap_{j=0}^{n-1} U_j:U_j\in{\mathcal U}_j,0
\leq j<1\right\}.
\]
It is also useful to be able to compare covers. For this purpose, we
introduce the following relation on the collection of all covers of a
set $X$.

\

\bdfn\label{ch6finercover}
Let $\mathcal U$ and $\mathcal V$ be covers of $X$.
We say that $\mathcal V$ is {\em finer} than, or a
\index{(N)}{cover refinement}{\em refinement} of, the cover 
$\mathcal U$, and denote this by  
$$
{\mathcal U}\prec{\mathcal V},
$$ 
if and only if every element of ${\mathcal V}$ 
is a subset of an element of ${\mathcal U}$. That is, for every set
$V\in{\mathcal V}$ 
there exists a set $U\in{\mathcal U}$ such that $V\subseteq U$.
It is also sometimes said that ${\mathcal V}$ is {\em inscribed} in
${\mathcal U}$.
\edfn

\

\blem\label{ch6refine}
Let ${\mathcal U}$, ${\mathcal V}$, ${\mathcal W}$ and ${\mathcal X}$
be covers of $X$. Then: 
\begin{itemize}
\item[$(a)$] The refinement relation $\prec$ is reflexive (i.e.,
  ${\mathcal U}\prec{\mathcal U}$) 
           and transitive (i.e., if ${\mathcal U}\prec{\mathcal V}$
           and ${\mathcal V}\prec{\mathcal W}$, then 
           ${\mathcal U}\prec{\mathcal W}$).
           
\,

\item[$(b)$] $ {\mathcal U}\prec{\mathcal U}\vee{\mathcal V}$.
          
\,

\item[$(c)$] If ${\mathcal V}$ is a subcover of ${\mathcal U}$,
          
\,
             then ${\mathcal U}\prec{\mathcal V}$.
\item[$(d)$] ${\mathcal U}$ is a subcover of ${\mathcal
    U}\vee{\mathcal U}$. Hence, from~$(c)$ and~$(b)$, 
            we deduce that
            \[
            {\mathcal U}\prec{\mathcal U}\vee{\mathcal U}\prec{\mathcal U}.
            \]
Nevertheless,
${\mathcal U}$ is not equal to ${\mathcal U}\vee{\mathcal U}$ in general.
          
\,

\item[$(e)$] If ${\mathcal U}\prec{\mathcal V}$ \underline{or}
  ${\mathcal U}\prec{\mathcal W}$, 
then ${\mathcal U}\prec{\mathcal V}\vee{\mathcal W}$.
        
\,

\item[$(f)$]  If ${\mathcal U}\prec{\mathcal W}$ \underline{and}
  ${\mathcal V}\prec{\mathcal W}$, 
then ${\mathcal U}\vee{\mathcal V}\prec{\mathcal W}$.
       
\,

\item[$(g)$]  If ${\mathcal U}\prec{\mathcal W}$ and ${\mathcal
    V}\prec{\mathcal X}$, then 
${\mathcal U}\vee{\mathcal V}\prec{\mathcal W}\vee{\mathcal X}$.
\end{itemize}
\elem

\bpf
All of these properties can be proved directly and are left to the
reader. As a hint, observe that property $(e)$ is a consequence
of~$(b)$ and~$(a)$ (transitivity), while 
property $(g)$ follows from~$(e)$ and~$(f)$.
\epf

\sp\brem
The relation $\prec$ does not constitute a
           partial order relation on the collection of all covers of
           $X$. This is because although it is reflexive and
           transitive,  it is not antisymmetric, that is, 
           ${\mathcal U}\prec{\mathcal V}\prec{\mathcal U}$ does not necessarily imply that ${\mathcal U}={\mathcal V}$; see Lemma \ref{ch6refine} $(d)$.
\erem

If $X$ is a metric space, then the maximum size of the elements
of a cover is encompassed by the notion of diameter of the cover.

\sp\bdfn\label{ch6diam}
If $(X,d)$ is a metric space, then the
\index{(N)}{cover diameter}{\em diameter} of a cover ${\mathcal U}$ of $X$
is defined by 
\[
\mathrm{diam}({\mathcal U}):=\sup\{\mathrm{diam}(U):U\in{\mathcal U}\},
\]
where
\[
\mathrm{diam}(U)
:=\sup \{d(x,y):x,y\in U\}.
\]
\edfn

It is also often of interest to know that all balls of some
specified radius are each contained in at least one element of a given
cover. Such a radius is known as a Lebesgue number for the cover. 

\sp\bdfn\label{ch6Lebnum}
A number $\delta>0$ is said to be a \index{(N)}{Lebesgue number}
Lebesgue number for a cover ${\mathcal U}$ of a metric space $(X,d)$ if every subset of $X$ of diameter not exceeding
$2\delta$ is contained in an element of ${\mathcal U}$.
\edfn

\sp\fr It is clear that if $\delta_0$ is a Lebesgue number for a cover
${\mathcal U}$, then 
so is any $0<{\delta}<\delta_0$. One can easily prove by contradiction
that every open cover of a compact metric space admits such a number. Recall
that an open cover is simply a cover whose elements are all open
subsets of the space. 

\

\subsection{Dynamical Covers}
In this section, we now add  a dynamical aspect to the above
discussion. Let $X$ be a nonempty set and let $T:X\lra X$ be a 
map. We will define covers that are induced by the dynamics of the map
$T$. First,  let us define the preimage of a cover under a map. 

\bdfn\label{ch6preim}
Let $X$ and $Y$ be non-empty sets. Let $h:X\lra Y$ be a map and let
${\mathcal V}$ be a cover of $Y$. The \index{(N)}{cover  preimage}{\em
  preimage} of 
${\mathcal V}$ under the map $h$ is the cover consisting of all the preimages of
the elements of ${\mathcal V}$ under $h$, that is, 
\[
h^{-1}{\mathcal V}:=\{h^{-1}(V):V\in{\mathcal V}\}.
\]
\edfn

\sp\fr We will now show that, as far as set operations go, the operator $h^{-1}$ behaves well with respect to cover operations. 

\blem\label{ch6tinv}
Let $h:X\lra Y$ be a map, and ${\mathcal U}$ and ${\mathcal V}$ be
covers of $Y$. The following assertions hold. 
\begin{itemize}
\item[$(a)$] The operation $h^{-1}$ preserves the refinement relation, that is,
\[
{\mathcal U}\prec{\mathcal V}\ \Longrightarrow\ h^{-1}{\mathcal
  U}\prec h^{-1}{\mathcal V}. 
\]
Moreover, if ${\mathcal V}$ is a subcover of ${\mathcal U}$ then
$h^{-1}{\mathcal V}$ is a subcover of $h^{-1}{\mathcal U}$.

\,

\item[$(b)$]$\;$The map $h^{-1}$ respects the join operation, that is,
\[
h^{-1}({\mathcal U}\vee{\mathcal V})=h^{-1}{\mathcal U}\vee h^{-1}{\mathcal V}.
\]
By induction, operation $h^{-n}$ for each $n\in\N$ also enjoys these properties.
\end{itemize}
\elem

\bpf
These assertions are straightforward to prove and are left to the
reader. 
\epf

\sp\fr We now introduce covers that follow the orbits of a given map by
indicating to 
which elements of a given cover the successive iterates of the map belong.

\bdfn\label{ch6iter}
Let $T:X\lra X$ be a map and ${\mathcal U}$ be a cover of $X$. For
every $n\in\N$, 
define the \index{(N)}{cover!dynamical}{\em dynamical cover}
\[
{\mathcal U}^n:=\bigvee_{j=0}^{n-1}T^{-j}{\mathcal U}
={\mathcal U}\vee T^{-1}{\mathcal U}\vee\cdots\vee T^{-(n-1)}{\mathcal U}.
\]
\edfn

\sp\fr A typical element of ${\mathcal U}^n$ is of the form
$$
U_0\cap T^{-1}(U_1)\cap T^{-2}(U_2)\cap\ldots\cap T^{-(n-1)}(U_{n-1})
$$ 
for some $U_0,U_1,U_2,\ldots, U_{n-1}\in{\mathcal U}$. This element
is the set of all points whose iterates fall successively into the elements
$U_0,U_1,U_2,\ldots$, and $U_{n-1}$.

\blem\label{ch6dyncovprop}
Let $T:X\lra X$ be a map and let $\mathcal U$ and $\mathcal V$ be some covers of $X$. Fix $n\in\N$. Then:
\begin{itemize}
\item[$(a)$] If ${\mathcal U}\prec{\mathcal V}$, then ${\mathcal
    U}^n\prec{\mathcal V}^n$. 
    
\,

\item[$(b)$] $({\mathcal U}\vee{\mathcal V})^n={\mathcal
    U}^n\vee{\mathcal V}^n$. 
\end{itemize}
\elem

\bpf
The first property follows directly from Lemmas~\ref{ch6tinv}$(a)$
and~\ref{ch6refine}$(g)$. The second is a consequence of
Lemma~\ref{ch6tinv}$(b)$. 
\epf


%

\subsection{Definition of Topological Pressure via Open Covers}\label{ch9deftopprescov}

We are now closer to the definition of \index{(N)}{topological pressure}topological pressure. It will involve a potential.

 Recall that a topological dynamical system $T:X\lra X$ is a self--transformation $T$ of 
a compact metrizable space $X$. Let $\varphi:X\lra\R$ be a real--valued
continuous function. In the context of topological pressure (for historical and physical reasons), such a function is usually referred to as a \index{(N)}{potential} potential.

The topological pressure of a potential $\varphi$ with respect to the transformation $T$ is defined in two stages. The first stage is to define topological pressure relative to an open cover, and then to take appropriate supremum over such covers.

\subsection{First Stage: Pressure of a Potential Relative to an Open Cover}

Recall that the $n$-th Birkhoff sum of a potential $\varphi$ at a point $x\in X$ is given by
\[
S_n\varphi(x)=\sum_{j=0}^{n-1}\varphi(T^j(x)).
\]
This is the sum of the values of the potential $\varphi$ at the first $n$
iterates of $x$ under $T$.
\bdfn\label{ch9snfsupset}
For every set $Y\sbt X$ and every $n\in\N$, define
\[
\overline{S}_n\varphi(Y)
:=\sup\big\{S_n\varphi(y):y\in Y\big\}
\hspace{0.3cm} \text{ and } \hspace{0.3cm}
\underline{S}_n\varphi(Y):=\inf\big\{S_n\varphi(y):y\in Y\big\}.
\]
\edfn

Now, let ${\mathcal U}$ be an open cover of $X$.
We start thermodynamic formalism with the following definition.

\bdfn\label{ch9znf}
Let $T:X\lra X$ be a topological dynamical system and let $\varphi:X\lra\R$ be a potential.
Let ${\mathcal U}$ be an open cover of $X$. For each $n\in\N$, define the $n$-th level functions,
frequently called \index{(N)}{partition function} partition functions, of ${\mathcal U}$ with respect 
to the potential $\varphi$ by
\[
Z_n(\varphi,{\mathcal U})
:=\inf\Bigl\{\sum_{V\in{\mathcal V}}e^{\overline{S}_n\varphi(V)}
  :{\mathcal V} \text{ is a subcover of } {\mathcal U}^n\Bigr\}
\]
and
\[
z_n(\varphi,{\mathcal U})
:=\inf\Bigl\{\sum_{V\in{\mathcal V}}e^{\underline{S}_n\varphi(V)}
  :{\mathcal V} \text{ is a subcover of } {\mathcal U}^n\Bigr\}.
\]
\edfn

\,

Note that if $\varphi\equiv 0$, then both numbers $Z_n(0,{\mathcal U})$ and $z_n(0,{\mathcal U})$ are equal to the minimum number of elements of ${\mathcal U}^n$ required to cover $X$. We then frequently write simply  $Z_n({\mathcal U})$ for $Z_n(0,{\mathcal U})$.

\sp \brem\label{ch9entpress}
\
\begin{itemize}
\item[(a)] It is sufficient to take the infimum over all finite subcovers
           since the exponential function takes only positive values and every subcover
					 has itself a finite subcover. However, this infimum may not be achieved if ${\mathcal U}$ is infinite.

\,

\item[(b)] In general, $Z_n(\varphi,{\mathcal U})\neq Z_1(\varphi,{\mathcal U}^n)$ and $z_n(\varphi,{\mathcal U})\neq z_1(\varphi,{\mathcal U}^n)$.

\,

\item[(c)] If $\varphi\equiv0$, then $Z_n(0,{\mathcal U})=z_n(0,{\mathcal U})=Z_n({\mathcal U})$ 

\,
           for all $n\in\N$ and any open cover $\mathcal{U}$ of $X$.
\item[(d)] If $\varphi\equiv c$ for some $c\in\R$, then $Z_n(c,{\mathcal U})=z_n(c,{\mathcal U})=e^{nc}Z_n({\mathcal U})$ 
           for all $n\in\N$ and every open cover $\mathcal{U}$ of $X$.

\,

\item[(e)] 
            For all open covers $\mathcal{U}$ of $X$ and all $n\in\N$, we have 
						\[e^{n\inf(\varphi)}Z_n({\mathcal U})\leq Z_n(\varphi,{\mathcal U})\leq e^{n\sup(\varphi)}Z_n({\mathcal U})\]
						and
						\[e^{n\inf(\varphi)}Z_n({\mathcal U})\leq z_n(\varphi,{\mathcal U})\leq e^{n\sup(\varphi)}Z_n({\mathcal U}).\]
\end{itemize}
\erem

We need the following.

\bdfn\label{ch9osc}
The \index{(N)}{oscillation} oscillation of $\varphi$ with respect 
to an open cover ${\mathcal U}$ is defined to be
\[
\mathrm{osc}(\varphi,{\mathcal U})
:=\sup\big\{|\varphi(y)-\varphi(x)|:U\in{\mathcal U},\, x,y\in U\big\}.
\]
\edfn

\sp\fr Note that $\mathrm{osc}(\varphi,\cdot)\leq2\|\varphi\|_\infty$. Also, 
$\mathrm{osc}(c,\cdot)=0$ for all $c\in\R$.

\blem\label{oscsnun}
Let $T:X\lra X$ be a topological dynamical system and let $\varphi:X\lra\R$ be a potential. Then, for every $n\in\N$ and every open cover ${\mathcal U}$ of $X$, we have that
\[
\mathrm{osc}(S_n\varphi,{\mathcal U}^n)\leq n\,\mathrm{osc}(\varphi,{\mathcal U}).
\]
\elem

\bpf
Let 
$$
V:=U_0\cap\ldots\cap T^{-(n-1)}(U_{n-1})\in{\mathcal U}^n
$$ 
and let $x,y\in V$.  
For each $0\leq j\le n-1$, we have that $T^j(x),T^j(y)\in U_j\in{\mathcal U}$. Hence, for all $0\leq j\le n-1$,
\[
\bigl|\varphi(T^j(x))-\varphi(T^j(y))\bigr|\leq\mathrm{osc}(\varphi,{\mathcal U}).
\]
Therefore
\[
\bigl|S_n\varphi(x)-S_n\varphi(y)\bigr|
\leq\sum_{j=0}^{n-1}\bigl|\varphi(T^j(x))-\varphi(T^j(y))\bigr| 
\leq n\,\mathrm{osc}(\varphi,{\mathcal U}).
\]
Since this is true for all $x,y\in V$ and all $V\in{\mathcal U}^n$, the result follows. 
\qed \epf

We now look at the relationship between the $Z_n$'s and the $z_n$'s.

\blem\label{ch9sem2lem9.4var2}
Let $T:X\lra X$ be a topological dynamical system and let $\varphi:X\lra\R$ be a potential. Then, for all $n\in\N$ and all open covers ${\mathcal U}$ of $X$, the following inequalities hold:
\[
z_n(\varphi,{\mathcal U})
\leq 
Z_n(\varphi,{\mathcal U})
\leq
e^{n\,\mathrm{osc}(\varphi,{\mathcal U})}z_n(\varphi,{\mathcal U}).
\]
\elem

\bpf
The left inequality is obvious. To ascertain the right one, let 
${\mathcal W}$ be a subcover of ${\mathcal U}^n$. Then
$$
\begin{aligned}
\sum_{W\in{\mathcal W}}e^{\overline{S}_n\varphi(W)}
&\leq\exp\Bigl(\sup\bigl\{\,\overline{S}_n\varphi(W)-\underline{S}_n\varphi(W):W\in{\mathcal W}\bigr\}\Bigr)
      \sum_{W\in{\mathcal W}}e^{\underline{S}_n\varphi(W)} \\
&\leq e^{\mathrm{osc}(S_n\varphi,{\mathcal U}^n)}\sum_{W\in{\mathcal W}}e^{\underline{S}_n\varphi(W)} \\
 &\leq e^{n\,\mathrm{osc}(\varphi,{\mathcal U})}\sum_{W\in{\mathcal W}}e^{\underline{S}_n\varphi(W)}.
\end{aligned}
$$
Taking on both sides the infimum over all subcovers of ${\mathcal U}^n$ results in the right inequality. 
\epf

In the next few results, we will see that the $Z_n$'s and the $z_n$'s have distinct properties.
\blem\label{ch9sem2lem9.4var}
Let $T:X\lra X$ be a topological dynamical system and let $\varphi:X\lra\R$ be a potential.
If 
${\mathcal U}\prec{\mathcal V}$, then for all $n\in\N$ we have that
\[
Z_n(\varphi,{\mathcal U})e^{-n\,\mathrm{osc}(\varphi,{\mathcal U})}\leq Z_n(\varphi,{\mathcal V})
\hspace{0.5cm} \text{ while } \hspace{0.5cm}
z_n(\varphi,{\mathcal U})\leq z_n(\varphi,{\mathcal V}).
\]
\elem

\bpf
Fix $n\in\N$. Let $i:{\mathcal V}\to{\mathcal U}$ be a map such that $V\sbt i(V)$
for all $V\in{\mathcal V}$. The map $i$ induces a map $i_n:{\mathcal V}^n\to{\mathcal U}^n$ 
in the following way. For every 
$$
W:=V_0\cap\ldots\cap T^{-(n-1)}(V_{n-1})\in{\mathcal V}^n,
$$
define 
\[
i_n(W):=i(V_0)\cap\ldots\cap T^{-(n-1)}(i(V_{n-1})).
\]
Observe that 
$$
W\sbt i_n(W)\in{\mathcal U}^n
$$ 
for all $W\in{\mathcal V}^n$.
Moreover, if $x\in W$ and $y\in i_n(W)$, then for each $0\leq j\le n-1$
we have that $T^j(x)\in V_{j}\sbt i(V_j)\ni T^j(y)$. 
So $T^j(x),T^j(y)\in i(V_j)$ for all $0\leq j<n$. Hence, $x,y\in i_n(W)\in{\mathcal U}^n$,and thus
\[
S_n\varphi(x)\geq S_n\varphi(y)-\mathrm{osc}(S_n\varphi,{\mathcal U}^n).
\]
Taking the supremum over all $x\in W$ on the left--hand side and
over all $y\in i_n(W)$ on the right-hand side yields
\[
\overline{S}_n\varphi(W)\geq \overline{S}_n\varphi(i_n(W))-\mathrm{osc}(S_n\varphi,{\mathcal U}^n).
\]
Now, let ${\mathcal W}$ be a subcover of ${\mathcal V}^n$.
Then $i_n({\mathcal W}):=\{i_n(W):W\in{\mathcal W}\}$ is a subcover of ${\mathcal U}^n$. Therefore
\begin{eqnarray*}
\sum_{W\in{\mathcal W}}e^{\overline{S}_n\varphi(W)}
&\geq&e^{-\mathrm{osc}(S_n\varphi,{\mathcal U}^n)}\sum_{W\in{\mathcal W}}e^{\overline{S}_n\varphi(i_n(W))} \\
&=&e^{-\mathrm{osc}(S_n\varphi,{\mathcal U}^n)}\sum_{Y\in i_n({\mathcal W})}e^{\overline{S}_n\varphi(Y)} \\
&\geq&e^{-\mathrm{osc}(S_n\varphi,{\mathcal U}^n)}Z_n(\varphi,{\mathcal U}).
\end{eqnarray*}
Taking the infimum over all subcovers ${\mathcal W}$ of ${\mathcal V}^n$ 
on the left--hand side and using Lemma~\ref{oscsnun}, we conclude that
\[
Z_n(\varphi,{\mathcal V})
\geq e^{-\mathrm{osc}(S_n\varphi,{\mathcal U}^n)}Z_n(\varphi,{\mathcal U})
\geq e^{-n\,\mathrm{osc}(\varphi,{\mathcal U})}Z_n(\varphi,{\mathcal U}).
\]
The proof of the inequality for the $z_n$'s is left to the reader.
\epf

\fr The proof of the following lemma  is left to the reader as an exercise.

\blem\label{ch9sem2lem9.3join}
Let $T:X\lra X$ be a topological dynamical system and let $\varphi:X\lra\R$ be a potential. Let ${\mathcal U}$ and ${\mathcal V}$ be open covers of $X$ and let $n\in\N$. Then
\[
Z_n(\varphi,{\mathcal U}\vee{\mathcal V})
\leq\min\bigl\{Z_n(\varphi,{\mathcal U})\cdot Z_n({\mathcal V}),Z_n({\mathcal U})\cdot Z_n(\varphi,{\mathcal V})\bigr\}.
\]
and
\[
z_n(\varphi,{\mathcal U}\vee{\mathcal V})
\leq\min\bigl\{e^{n\,\mathrm{osc}(\varphi,{\mathcal U})}z_n(\varphi,{\mathcal U})\cdot Z_n({\mathcal V}),
               Z_n({\mathcal U})\cdot e^{n\,\mathrm{osc}(\varphi,{\mathcal V})}z_n(\varphi,{\mathcal V})\bigr\}.
\]
\elem


\,

\fr In order to define topological pressure we need the following well known concept. Recall from 

\bdfn\label{1d20190927} A sequence $\{a_n\}_{n=1}^\infty$ consisting of real numbers is said to be subadditive\index{(N)}{subadditive sequence} if and only if 
$$
a_{n+m}\le a_n+a_n
$$
for all integers $m,n\ge 1$. Likewise, a sequence $\{b_n\}_{n=1}^\infty$ consisting of positive real numbers is said to be submultiplicative\index{(N)}{submultiplicative sequence} if and only if 
$$
b_{n+m}\le b_nb_n
$$
for all integers $m,n\ge 1$. 
\edfn

\fr We immediately get the following.

\bobs\label{o220190927}
If $\{a_n\}_{n=1}^\infty$ is a submultiplicative sequence of positive real numbers, then $\{\log(a_n)\}_{n=1}^\infty$ is a subadditive sequence of real numbers. 
\eobs

\sp\fr Subadditive sequences of real numbers possess the following incredibly helpful property. The reader should take note that this benign--looking lemma is one of the foundation stones of the theory of topological
topological pressure. 

\blem\label{ch6subadditive}
If $(a_n)_{n=1}^\infty$ is a subadditive sequence of real numbers,
then the sequence
$\bigl(\frac{1}{n}a_n\bigr)_{n=1}^\infty$ converges and
\[
\lim_{n\to\infty}\frac{1}{n}a_n=\inf_{n\in\N}\lt\{\frac{1}{n}a_n\rt\}.
\]
If, moreover, $(a_n)_{n=1}^\infty$ is bounded from below, then
$\inf_{n\in\N}\frac{1}{n}a_n\geq0$.
\elem

\bpf
Fix $m\in\N$. By the division algorithm, every $n\in\N$ can be
uniquely written in the form $n=km+r$, 
where $0\leq r<m$. The subadditivity of the sequence implies that
\[
\frac{a_n}{n}=\frac{a_{km+r}}{km+r}\leq\frac{a_{km}+a_r}{km+r}
\leq\frac{ka_m+a_r}{km}
=\frac{a_m}{m}+\frac{a_r}{km}.
\]
Notice that
\[
-\infty<\min\big\{a_s:0\leq s<m\big\}
\leq a_r\leq\max\big\{a_s:0\leq s<m\big\}<+\infty
\]
for all $n\in\N$. Therefore, as $n$ tends to infinity, also $k$ tends
to infinity and therefore 
$a_r/k$ approaches 0 by the Sandwich Theorem. Hence,
\[
\limsup_{n\to\infty}\frac{a_n}{n}\leq\frac{a_m}{m}.
\]
Since $m\in\N$ was chosen arbitrarily, taking the infimum over $m$ yields that
\[
\limsup_{n\to\infty}\frac{a_n}{n}\leq\inf_{m\in\N}\frac{a_m}{m}.
\]
Thus,
\[
\limsup_{n\to\infty}\frac{a_n}{n}
\leq\inf_{m\in\N}\lt\{\frac{a_m}{m}\rt\}
\leq\liminf_{n\to\infty}\frac{a_n}{n}
\leq\limsup_{n\to\infty}\frac{a_n}{n}.
\]
Consequently,
\[
\lim_{n\to\infty}\frac{a_n}{n}=\inf_{m\in\N}\lt\{\frac{a_m}{m}\rt\}.
\]
\epf

The application of this lemma in the current section is due to the following.

\blem\label{ch9sem2lem9.3}
Let $T:X\lra X$ be a topological dynamical system and let $\varphi:X\lra\R$ be a potential. If ${\mathcal U}$ is an open cover of $X$, the sequence
$(Z_n(\varphi,{\mathcal U}))_{n=1}^\infty$ is submultiplicative.
\elem

\bpf
Fix $m,n\in\N$, let ${\mathcal V}$ be a subcover of ${\mathcal U}^m$ and
${\mathcal W}$ a subcover of ${\mathcal U}^n$.
Note that ${\mathcal V}\vee T^{-m}({\mathcal W})$ is a subcover of
${\mathcal U}^{m+n}$ since it is a cover and
\[
{\mathcal V}\vee T^{-m}({\mathcal W})
\sbt{\mathcal U}^m\vee T^{-m}({\mathcal U}^n)
={\mathcal U}^{m+n}.
\]
Take arbitrary $V\in{\mathcal V}$ and
$W\in{\mathcal W}$. Then for every $x\in V\cap T^{-m}(W)$,
we have $x\in V$ and $T^m(x)\in W$ and hence
\[
S_{m+n}\varphi(x)
=S_m\varphi(x)+S_n\varphi(T^m(x))
\leq \overline{S}_m\varphi(V)+\overline{S}_n\varphi(W).
\]
Taking the supremum over all $x\in V\cap T^{-m}(W)$, we deduce that
\[
\overline{S}_{m+n}\varphi(V\cap T^{-m}(W))
\leq\overline{S}_m\varphi(V)+\overline{S}_n\varphi(W).
\]
Therefore
\begin{eqnarray*}
Z_{m+n}(\varphi,{\mathcal U})
&\leq&\sum_{E\in{\mathcal V}\vee T^{-m}({\mathcal W})}e^{\overline{S}_{m+n}\varphi(E)} 
\leq\sum_{V\in{\mathcal V}}\sum_{W\in{\mathcal W}}e^{\overline{S}_{m+n}\varphi(V\cap T^{-m}(W))} \\
&\leq&\sum_{V\in{\mathcal V}}\sum_{W\in{\mathcal W}}e^{\overline{S}_m\varphi(V)}e^{\overline{S}_n\varphi(W)} \\
&=&\sum_{V\in{\mathcal V}}e^{\overline{S}_m\varphi(V)}\sum_{W\in{\mathcal W}}e^{\overline{S}_n\varphi(W)}.
\end{eqnarray*}
Taking the infimum of the right-hand side over all subcovers ${\mathcal V}$ of ${\mathcal U}^m$ and
over all subcovers ${\mathcal W}$ of ${\mathcal U}^n$ gives
\[
Z_{m+n}(\varphi,{\mathcal U})
\leq Z_m(\varphi,{\mathcal U})Z_n(\varphi,{\mathcal U}).
\]
\epf

\fr We immediately get from this lemma and Observation~\ref{o220190927} the following fact.

\bcor\label{ch9sem2lem9.3log}
If $T:X\lra X$ is a topological dynamical system and $\varphi:X\lra\R$ is a potential, then the sequence 
$$
\(\log Z_n(\varphi,{\mathcal U})\)_{n=1}^\infty
$$
is subadditive for every open cover ${\mathcal U}$ of $X$.
\ecor

Thanks to this fact, we can define the topological pressure of a potential
with respect to an open cover. This constitutes the first step in the
definition of the topological pressure of a potential.

\bdfn\label{ch9presswrtcov}
Let $T:X\lra X$ be a topological dynamical system and let $\varphi:X\lra\R$ be a potential.
The \index{(N)}{topological pressure with respect to a cover} topological pressure of the potential $\varphi$ with respect to an open cover ${\mathcal U}$ of $X$, 
denoted by $\P(T,\varphi,{\mathcal U})$, is defined to be 
\[
\P(T,\varphi,{\mathcal U})
:=\lim_{n\to\infty}\frac{1}{n}\log Z_n(\varphi,{\mathcal U})
=\inf_{n\in\N}\lt\{\frac{1}{n}\log Z_n(\varphi,{\mathcal U})\rt\},
\]
where the existence of the limit and its equality with the infimum follow immediately from Lemma~\ref{ch6subadditive} and Corollary~\ref{ch9sem2lem9.3log}. 

If $\varphi\equiv 0$, we simply write
$$
\htop(T,{\mathcal U})
$$
for $P(T,0,{\mathcal U})$ and we call this quantity the \index{(N)}{topological entropy of a cover} topological entropy of $T$ with respect to the cover ${\mathcal U}$. 
\edfn

\,

\fr It is also possible to define similar quantities using the $z_n(\varphi,{\mathcal U})$'s
rather than the $Z_n(\varphi,{\mathcal U})$'s.
\bdfn\label{ch9presswrtcov2}
Let $T:X\lra X$ be a topological dynamical system. Given a potential $\varphi:X\lra\R$ and an open cover ${\mathcal U}$ of $X$,
let
\[
\underline{p}(T,\varphi,{\mathcal U})
:=\liminf_{n\to\infty}\frac{1}{n}\log z_n(\varphi,{\mathcal U})
\hspace{0.25cm} and \hspace{0.25cm}  
\overline{p}(T,\varphi,{\mathcal U})
:=\limsup_{n\to\infty}\frac{1}{n}\log z_n(\varphi,{\mathcal U}).
\]
\edfn

\brem\label{ch9presrelrmk}
Let ${\mathcal U}$ be an open cover of $X$.
\begin{itemize}
\item[(a)]  $\P(T,0,{\mathcal U})=\underline{p}(T,0,{\mathcal U})=\overline{p}(T,0,{\mathcal U})
             =\htop(T,{\mathcal U})$ by Remark~\ref{ch9entpress}(c).
\vspace{0.5mm}

\,

\item[(b)]  By Remark~\ref{ch9entpress}(e),
\[
-\infty
<\htop(T,{\mathcal U})+\inf\varphi
\leq \P(T,\varphi,{\mathcal U})
\leq\htop(T,{\mathcal U})+\sup\varphi
<\infty.
\]
These inequalities also hold with $\P(T,\varphi,{\mathcal U})$ replaced by 
$\underline{p}(T,\varphi,{\mathcal U})$ and $\overline{p}(T,\varphi,{\mathcal U})$, respectively. 
\vspace{0.5mm}

\,

\item[(c)]  Using Lemma~\ref{ch9sem2lem9.4var2},
\[ 
\underline{p}(T,\varphi,{\mathcal U})
\leq\overline{p}(T,\varphi,{\mathcal U})
\leq \P(T,\varphi,{\mathcal U})
\leq\underline{p}(T,\varphi,{\mathcal U})+\mathrm{osc}(\varphi,{\mathcal U}).
\]
\end{itemize}
\erem

The topological pressure respects the refinement relation and is subadditive
with respect to the join operation. It has the following properties.

\bprop\label{ch9sem2lem9.4}
Let $T:X\lra X$ be a topological dynamical system and let $\varphi:X\lra\R$ 
be a potential. Let ${\mathcal U}$ and ${\mathcal V}$ be open covers of $X$. If ${\mathcal U}\prec{\mathcal V}$, then
\begin{itemize}
\item[(a)] 
           $$
					 \P(T,\varphi,{\mathcal U})-\mathrm{osc}(\varphi,{\mathcal U})\leq \P(T,\varphi,{\mathcal V})
					 $$
					 
	\,
	
	\fr 			 while
					 \[
					 \underline{p}(T,\varphi,{\mathcal U})\leq\underline{p}(T,\varphi,{\mathcal V}) 
					 \hspace{0.5cm} \text{ and } \hspace{0.5cm}  
					 \overline{p}(T,\varphi,{\mathcal U})\leq\overline{p}(T,\varphi,{\mathcal V}). 
					 \]
\item[(b)] 
$$
\P(T,\varphi,{\mathcal U}\vee{\mathcal V})
            \leq\min\bigl\{\P(T,\varphi,{\mathcal U})+\htop(T,{\mathcal V}),
                           \P(T,\varphi,{\mathcal V})+\htop(T,{\mathcal U})\bigr\}
                           $$ 
                           
	\,
	
		 \fr       whereas
$$
\begin{aligned}
					 \overline{p}(T,\varphi,{\mathcal U}\vee{\mathcal V})
           \leq\min\Bigl\{\overline{p}(T,\varphi,{\mathcal U})+&\mathrm{osc}(\varphi,{\mathcal U})+\htop(T,{\mathcal V}), \\
           &\overline{p}(T,\varphi,{\mathcal V})+\mathrm{osc}(\varphi,{\mathcal V})+\htop(T,{\mathcal U})\Bigr\}.
			     \end{aligned}
			     $$
\end{itemize}
\eprop

\bpf
Part~(a) is an immediate consequence of Lemma~\ref{ch9sem2lem9.4var} 
while (b) follows directly from Lemma~\ref{ch9sem2lem9.3join}.
%
\epf

\sp We shall prove the following.

\blem\label{ch9sem2cor9.10}
Let $T:X\lra X$ be a topological dynamical system and let $\varphi:X\lra\R$ 
be a potential. If ${\mathcal U}$ is an open cover of $X$, then 
\[
\underline{p}(T,\varphi,{\mathcal U}^n)=\underline{p}(T,\varphi,{\mathcal U})
\hspace{0.5cm} \text{ and } \hspace{0.5cm}
\overline{p}(T,\varphi,{\mathcal U}^n)=\overline{p}(T,\varphi,{\mathcal U})
\]
whereas 
$$
\P(T,\varphi,{\mathcal U}^n)\leq \P(T,\varphi,{\mathcal U})
$$ 
for all $n\in\N$. 

In addition, if ${\mathcal U}$ is an open partition of $X$, then 
$$
\P(T,\varphi,{\mathcal U}^n)=\P(T,\varphi,{\mathcal U})
$$ 
for all $n\in\N$.
\elem

\bpf
Fix $n\in\N$. For all $k\in\N$ and all $x\in X$, we already know that
\[
S_{k+n-1}\varphi(x)=S_k\varphi(x)+S_{n-1}\varphi(T^k(x)).
\]
Therefore
\[
S_k\varphi(x)-\|S_{n-1}\varphi\|_\infty
\leq S_{k+n-1}\varphi(x)\leq S_k\varphi(x)+\|S_{n-1}\varphi\|_\infty.
\]
Hence, for any subset $Y$ of $X$,
\beq\label{ch9ergsumny}
\overline{S}_k\varphi(Y)-\|S_{n-1}\varphi\|_\infty
\leq\overline{S}_{k+n-1}\varphi(Y)
\leq\overline{S}_k\varphi(Y)+\|S_{n-1}\varphi\|_\infty
\eeq
and
\beq\label{ch9ergsumny2}
\underline{S}_k\varphi(Y)-\|S_{n-1}\varphi\|_\infty
\leq\underline{S}_{k+n-1}\varphi(Y)
\leq\underline{S}_k\varphi(Y)+\|S_{n-1}\varphi\|_\infty.
\eeq
We claim that
\beq\label{ch9zphin}
e^{-\|S_{n-1}\varphi\|_\infty}Z_k(\varphi,{\mathcal U}^n)
\leq Z_{k+n-1}(\varphi,{\mathcal U})
\eeq
and
\beq\label{ch9zphin2}
e^{-\|S_{n-1}\varphi\|_\infty}z_k(\varphi,{\mathcal U}^n)
\leq z_{k+n-1}(\varphi,{\mathcal U})
\leq e^{\|S_{n-1}\varphi\|_\infty}z_k(\varphi,{\mathcal U}^n).
\eeq
Let us first prove~(\ref{ch9zphin}). 
Recall that $({\mathcal U}^n)^k\prec{\mathcal U}^{k+n-1}\prec({\mathcal U}^n)^k$
for all $k\in\N$ (cf.~Lemma~\ref{ch6dyncovprop}(d)). However, this is 
insufficient to declare that a subcover ${\mathcal U}^{k+n-1}$
is also a subcover of $({\mathcal U}^n)^k$, or vice versa. We need to remember that 
${\mathcal U}\vee{\mathcal U}\spt{\mathcal U}$ and thus $({\mathcal U}^n)^k\spt{\mathcal U}^{k+n-1}$,
that is, ${\mathcal U}^{k+n-1}$ is a subcover of $({\mathcal U}^n)^k$.
Let ${\mathcal V}$ be a subcover of
${\mathcal U}^{k+n-1}$. Then ${\mathcal V}$ is a subcover
of $({\mathcal U}^n)^k$. Using the left inequality in~(\ref{ch9ergsumny}) with $Y$ replaced by each
$V\in{\mathcal V}$ successively, we obtain
\[
e^{-\|S_{n-1}\varphi\|_\infty}Z_k(\varphi,{\mathcal U}^n)
\leq e^{-\|S_{n-1}\varphi\|_\infty}\sum_{V\in{\mathcal V}}e^{\overline{S}_k\varphi(V)}
\leq\sum_{V\in{\mathcal V}}e^{\overline{S}_{k+n-1}\varphi(V)}.
\]
Taking the infimum over all subcovers ${\mathcal V}$ of ${\mathcal U}^{k+n-1}$ yields~(\ref{ch9zphin}).
Similarly, using the left inequality in~(\ref{ch9ergsumny2}),
we get that
\[
e^{-\|S_{n-1}\varphi\|_\infty}z_k(\varphi,{\mathcal U}^n)
\leq e^{-\|S_{n-1}\varphi\|_\infty}\sum_{V\in{\mathcal V}}e^{\underline{S}_k\varphi(V)}
\leq\sum_{V\in{\mathcal V}}e^{\underline{S}_{k+n-1}\varphi(V)}.
\]
Taking the infimum over all subcovers ${\mathcal V}$ of ${\mathcal U}^{k+n-1}$ yields
the left inequality in~(\ref{ch9zphin2}). Regarding the right inequality, 
since ${\mathcal U}^{k+n-1}\prec({\mathcal U}^n)^k$, there exists a map 
$i:({\mathcal U}^n)^k\to{\mathcal U}^{k+n-1}$ such that $W\sbt i(W)$
for all $W\in({\mathcal U}^n)^k$. Let ${\mathcal W}$ be a subcover of $({\mathcal U}^n)^k$. 
Then $i({\mathcal W})$ is a subcover of ${\mathcal U}^{k+n-1}$ and, 
using the right inequality in~(\ref{ch9ergsumny2}), we deduce that
\begin{eqnarray*}
\sum_{W\in{\mathcal W}}e^{\underline{S}_k\varphi(W)}
&\geq&\sum_{W\in{\mathcal W}}e^{\underline{S}_k\varphi(i(W))} 
   =  \sum_{Z\in i({\mathcal W})}e^{\underline{S}_k\varphi(Z)} \\
&\geq&\sum_{Z\in i({\mathcal W})}e^{\underline{S}_{k+n-1}\varphi(Z)-\|S_{n-1}\varphi\|_\infty} \\
&\geq&e^{-\|S_{n-1}\varphi\|_\infty}z_{k+n-1}(\varphi,{\mathcal U}).
\end{eqnarray*}
Taking the infimum over all subcovers of $({\mathcal U}^n)^k$ on the left-hand side 
gives the right inequality in~(\ref{ch9zphin2}). So~(\ref{ch9zphin}) and~(\ref{ch9zphin2})
always hold. 

Moreover, if ${\mathcal U}$ is a partition then ${\mathcal U}\vee{\mathcal U}={\mathcal U}$
and thus $({\mathcal U}^n)^k={\mathcal U}^{k+n-1}$ for all $k\in\N$. Let ${\mathcal W}$ be a 
subcover of $({\mathcal U}^n)^k$. Using the right inequality in~(\ref{ch9ergsumny}), 
we conclude that
\[
\sum_{W\in{\mathcal W}}e^{\overline{S}_k\varphi(W)}
\geq\sum_{W\in{\mathcal W}}e^{\overline{S}_{k+n-1}\varphi(W)-\|S_{n-1}\varphi\|_\infty} 
\geq e^{-\|S_{n-1}\varphi\|_\infty}Z_{k+n-1}(\varphi,{\mathcal U}).
\]
Taking the infimum over all subcovers of $({\mathcal U}^n)^k$ on the left-hand side 
gives  
\beq\label{ch9zphinr}
Z_k(\varphi,{\mathcal U}^n)
\geq e^{-\|S_{n-1}\varphi\|_\infty}
Z_{k+n-1}(\varphi,{\mathcal U}).
\eeq
Finally, for the passage from the $z_n$'s to $\overline{p}$, it follows from~(\ref{ch9zphin2}) that
\[
\frac{k}{k+n-1}\cdot\frac{1}{k}\log z_k(\varphi,{\mathcal U}^n)-\frac{\|S_{n-1}\varphi\|_\infty}{k+n-1}
\leq\frac{1}{k+n-1}\log z_{k+n-1}(\varphi,{\mathcal U})
\]
and
\[
\frac{1}{k+n-1}\log z_{k+n-1}(\varphi,{\mathcal U})
\leq\frac{k}{k+n-1}\cdot\frac{1}{k}\log z_k(\varphi,{\mathcal U}^n)+\frac{\|S_{n-1}\varphi\|_\infty}{k+n-1}.
\]
Taking the $\limsup$ as $k\to\infty$ in these two relations yields
\[
\overline{p}(T,\varphi,{\mathcal U}^n)\leq\overline{p}(T,\varphi,{\mathcal U})\leq\overline{p}(T,\varphi,{\mathcal U}^n).
\]
Similarly, one deduces from~(\ref{ch9zphin}) that $\P(T,\varphi,{\mathcal U}^n)\leq \P(T,\varphi,{\mathcal U})$ and,
when ${\mathcal U}$ is a partition, it ensues from~(\ref{ch9zphinr}) that $\P(T,\varphi,{\mathcal U}^n)\geq \P(T,\varphi,{\mathcal U})$.
\qed \epf

\subsection{Second Stage: The Pressure of a Potential}

We now give the definition of topological pressure of the potential:

\bdfn\label{ch9pressdfn}
Let $T:X\to X$ be a topological dynamical system and $\varphi:X\to\R$ a potential. The
\index{(N)}{topological pressure of a potential} topological pressure of the potential $\varphi$, 
denoted $\P(T,\varphi)$, is defined by
\[
\P(T,\varphi)
:=\sup\bigl\{\P(T,\varphi,{\mathcal U})-\mathrm{osc}(\varphi,{\mathcal U})
  :{\mathcal U} \text{ is an open cover of } X\bigr\}.
\]
\edfn

This definition may look a little bit awkward because of the term $\mathrm{osc}(\varphi,{\mathcal U})$ which appears in it. This is because of
Proposition~\ref{ch9sem2lem9.4}(a), due to which taking the supremum of the pressure relative to all covers does not always lead to a 
quantity that has natural desired properties. Our definition is nevertheless purely topological and if $\varphi\equiv 0$, then this term vanishes, thus disappears. We then write
$$
\htop(T)
:=\P(T,0)=\sup\bigl\{\P(T,0,{\mathcal U})\bigr\}
=\sup\bigl\{\htop(T,{\mathcal U})\bigr\},
$$
where, as above, the supremum is taken over all open covers of $X$. The quantity $\htop(T)$ is called the \index{(N)}{topological entropy} topological entropy of $T$.

In light of Proposition~\ref{ch9sem2lem9.4}(a), we may define the counterparts 
$\underline{p}(T,\varphi)$ and $\overline{p}(T,\varphi)$ of $\P(T,\varphi)$ 
by simply taking the supremum over all covers.
\bdfn\label{ch9pressdfn2}
Let $T:X\lra X$ be a topological dynamical system and $\varphi:X\lra\R$ a potential. Define
\[
\underline{p}(T,\varphi)
:=\sup\bigl\{\underline{p}(T,\varphi,{\mathcal U})
  :{\mathcal U} \text{ is an open cover of } X\bigr\}.
\]
and
\[
\overline{p}(T,\varphi)
:=\sup\bigl\{\overline{p}(T,\varphi,{\mathcal U})
  :{\mathcal U} \text{ is an open cover of } X\bigr\}.
\]
\edfn

Clearly, $\underline{p}(T,\varphi)\leq\overline{p}(T,\varphi)$. 
In fact, $\underline{p}(T,\varphi)$ and $\overline{p}(T,\varphi)$ are just 
other expressions of the topological pressure.

\bthm\label{pequalsP}
For any topological dynamical system $T:X\lra X$ and potential $\varphi:X\lra\R$, it turns out that 
$$
\underline{p}(T,\varphi)=\overline{p}(T,\varphi)=\P(T,\varphi).
$$
\ethm

\bpf
From a rearrangement of the right inequality in Remark~\ref{ch9presrelrmk}(c), it follows that 
$$
\P(T,\varphi)\leq\underline{p}(T,\varphi)\leq\overline{p}(T,\varphi).
$$ 
In order oo prove that 
$$
\overline{p}(T,\varphi)\leq \P(T,\varphi),
$$
let $({\mathcal U}_n)_{n=1}^\infty$ be a sequence of 
open covers such that 
$$
\lim_{n\to\infty}\overline{p}(T,\varphi,{\mathcal U}_n)=\overline{p}(T,\varphi).
$$
Each open cover ${\mathcal U}_n$ has a Lebesgue number $\delta_n>0$. The compactness 
of $X$ guarantees that there are finitely many open balls of radius $\min\{\delta_n,1/(2n)\}$ 
that cover $X$. These balls thereby constitute a refinement of ${\mathcal U}_n$ of diameter at most $1/n$. 
Thanks to Proposition~\ref{ch9sem2lem9.4}(a), this means that we may assume without loss of 
generality that the sequence $({\mathcal U}_n)_{n=1}^\infty$ is such that 
$$
\lim_{n\to\infty}\mathrm{diam}({\mathcal U}_n)=0.
$$ 
Since $\varphi$ is uniformly continuous, it ensues that 
$$
\lim_{n\to\infty}\mathrm{osc}(\varphi,{\mathcal U}_n)=0.
$$ 
Consequently, using the left inequality in Remark~\ref{ch9presrelrmk}(c), we conclude that
$$
\begin{aligned}
\P(T,\varphi)
&\geq\sup_{n\in\N}\bigl\{\P(T,\varphi,{\mathcal U}_n)-\mathrm{osc}(\varphi,{\mathcal U}_n)\bigr\} 
\geq\sup_{n\in\N}\bigl\{\overline{p}(T,\varphi,{\mathcal U}_n)-\mathrm{osc}(\varphi,{\mathcal U}_n)\bigr\} \\
&\geq\lim_{n\in\N}\bigl\{\overline{p}(T,\varphi,{\mathcal U}_n)-\mathrm{osc}(\varphi,{\mathcal U}_n)\bigr\} \\
&=\overline{p}(T,\varphi). 
\end{aligned}
$$
\epf

\brem\label{ch9presprop} 
We want to record the following straightforward properties of topological pressure

\begin{itemize}
\item[(a)] $\P(T,0)=\htop(T)$. 

\,

\item[(b)] By Remark~\ref{ch9presrelrmk}(b), 
\[\htop(T)+\inf(\varphi)-\mathrm{osc}(\varphi,X)
             \leq \P(T,\varphi)\leq\htop(T)+\sup(\varphi).
\]
\item[(c)] $\P(T,\varphi)=+\infty$ if and only if $\htop(T)=+\infty$, according to part~(b).
\end{itemize}
\erem

We now shall show that topological pressure does not increase under ``factorization''.

\bprop\label{ch9sem2prop9.8}
Suppose that $S:Y\lra Y$ is a topological factor of $T:X\lra X$ via the factor continuous surjection $h:X\lra Y$.
Then for every potential $\varphi:Y\lra\R$, we have that 
$$
\P(S,\varphi)\leq \P(T,\varphi\circ h).
$$
\eprop

\bpf
Let ${\mathcal V}$ be an open cover of $Y$. Observe
that 
$$
h^{-1}({\mathcal V}^n_S)=(h^{-1}({\mathcal V}))^n_T
$$ 
for all $n\in\N$. Moreover, letting $C$ be the collection of all subcovers of ${\mathcal V}_S^n$, the map 
$$
C\ni {\mathcal C}\longmapsto h^{-1}({\mathcal C})
$$
defines a bijection between subcovers of ${\mathcal V}_S^n$ and subcovers of $h^{-1}({\mathcal V}_S^n)=(h^{-1}({\mathcal V}))_T^n$ which preserves cardinalities,
that is, 
$$
\#h^{-1}({\mathcal C})=\#{\mathcal C}
$$ 
for all ${\mathcal C}\in C$, since $h$ is a
surjection. So 
$$
Z_n(S,\varphi,{\mathcal V})=Z_n(T,\varphi\circ h,h^{-1}({\mathcal V}))
$$
and therefore
\[
\P(S,\varphi,{\mathcal V})=\P(T,\varphi\circ h,h^{-1}({\mathcal V})).
\]
Also, observe that $\mathrm{osc}(\varphi\circ h,h^{-1}({\mathcal V}))=\mathrm{osc}(\varphi,{\mathcal V})$.
Then
\begin{eqnarray*}
\P(T,\varphi\circ h)
&\geq&\P(T,\varphi\circ h,h^{-1}({\mathcal V}))-\mathrm{osc}(\varphi\circ h,h^{-1}({\mathcal V})) \\
&=&\P(S,\varphi,{\mathcal V})-\mathrm{osc}(\varphi,{\mathcal V}).
\end{eqnarray*}
Passing to the supremum over all open covers ${\mathcal V}$ of $Y$ yields that
$\P(T,\varphi\circ h)\geq \P(S,\varphi)$.
\epf

\sp\fr An immediate but important consequence of this lemma is the following result.

\bcor\label{ch9sem2cor9.9}
If $T:X\lra X$ and $S:Y\lra Y$ are topologically conjugate dynamical systems via
a conjugacy $h:X\lra Y$, then 
$$
\P(S,\varphi)=\P(T,\varphi\circ h)
$$
for all potentials $\varphi:Y\lra\R$.
\ecor

We will now prove a result showing that topological pressure is determined by any sequence of covers whose diameters tend to zero. 
\bprop\label{ch9sem2lem9.5}
If $T:X\lra X$ is a topological dynamical system and $\varphi:X\lra\R$ is a continuous potential, then all the following quantities are all equal:
\begin{itemize}
\item[(a)] \hspace{0.1cm} $\P(T,\varphi)$.
\vspace{1.5mm}
\item[(b)] \hspace{0.1cm} $\overline{p}(T,\varphi)$.
\vspace{1.5mm}
\item[(c)] \hspace{0.1cm} $\displaystyle\lim_{\varepsilon\to0}\bigl\{\sup\bigl\{\P(T,\varphi,{\mathcal U}):
            {\mathcal U} \text{ open cover with }\mathrm{diam}({\mathcal U})\leq\varepsilon\bigr\}\bigr\}$.
\vspace{1.5mm}
\item[(d)] \hspace{0.1cm} $\sup\bigl\{\overline{p}(T,\varphi,{\mathcal U}):
            {\mathcal U} \text{ open cover with }\mathrm{diam}({\mathcal U})\leq\delta\bigr\}$
					  for any $\delta>0$.
\vspace{1.5mm}
\item[(e)] \hspace{0.1cm} $\displaystyle\lim_{\varepsilon\to0}\P(T,\varphi,{\mathcal U}_\varepsilon)$ for any family 
of open covers $({\mathcal U}_\varepsilon)_{\varepsilon\in(0,\infty)}$ such that
$\displaystyle\lim_{\varepsilon\to0}\mathrm{diam}({\mathcal U}_\varepsilon)=0$.
\vspace{1.5mm}
\item[(f)] \hspace{0.1cm} $\displaystyle\lim_{\varepsilon\to0}\overline{p}(T,\varphi,{\mathcal U}_\varepsilon)$ for any family 
of open covers $({\mathcal U}_\varepsilon)_{\varepsilon\in(0,\infty)}$ such that
$\displaystyle\lim_{\varepsilon\to0}\mathrm{diam}({\mathcal U}_\varepsilon)=0$.
\vspace{1.5mm}
\item[(g)] \hspace{0.1cm} $\displaystyle\lim_{n\to\infty}\P(T,\varphi,{\mathcal U}_{n})$ for any sequence of open covers 
$({\mathcal U}_n)_{n=1}^\infty$ such that $\displaystyle\lim_{n\to\infty}\mathrm{diam}({\mathcal U}_{n})=0$.
\item[(h)] \hspace{0.1cm} $\displaystyle\lim_{n\to\infty}\overline{p}(T,\varphi,{\mathcal U}_{n})$ for any sequence of open covers 
$({\mathcal U}_n)_{n=1}^\infty$ such that $\displaystyle\lim_{n\to\infty}\mathrm{diam}({\mathcal U}_{n})=0$.
\end{itemize}
Note that $\overline{p}$ can be replaced by $\underline{p}$ in the statements above.
\eprop

\bpf
We already know that (a)$=$(b) by Lemma~\ref{pequalsP}. It is clear that (b)$\geq$(d).  
It is also obvious that (d)$\geq$(f) and (c)$\geq$(e) for any family 
$({\mathcal U}_\varepsilon)_{\varepsilon\in(0,\infty)}$ as described, 
and that (d)$\geq$(h) and (c)$\geq$(g) for any sequence $({\mathcal U}_n)_{n=1}^\infty$ as specified. 
It thus suffices to prove that (f)$\geq$(b), that (h)$\geq$(b), that (e)$\geq$(a), that (g)$\geq$(a)
and that (b)$\geq$(c).    

We will prove that (g)$\geq$(a). The proofs of the other inequalities are similar.
Let ${\mathcal V}$ be any open cover of $X$. Since
$\lim_{n\to\infty}\mathrm{diam}({\mathcal U}_{n})=0$, there exists $N\in\N$ such that
$$
{\mathcal V}\prec{\mathcal U}_{n}
$$ 
for all $n\geq N$.
By Proposition~\ref{ch9sem2lem9.4}(a), we obtain that for all sufficiently large $n$,
\[
\P(T,\varphi,{\mathcal U}_n)\geq \P(T,\varphi,{\mathcal V})-\mathrm{osc}(\varphi,{\mathcal V}).
\]
We immediately deduce that 
\[
\liminf_{n\to\infty}\P(T,\varphi,{\mathcal U}_n)\geq \P(T,\varphi,{\mathcal V})-\mathrm{osc}(\varphi,{\mathcal V}).
\]
As the open cover ${\mathcal V}$ was chosen arbitrarily, passing to the supremum over all open covers allows us to conclude that
\begin{eqnarray*}
\liminf_{n\to\infty}\P(T,\varphi,{\mathcal U}_n)\geq \P(T,\varphi).
\end{eqnarray*}
But $\lim_{n\to\infty}\mathrm{osc}(\varphi,{\mathcal U}_n)=0$ since
$\lim_{n\to\infty}\mathrm{diam}({\mathcal U}_{n})=0$ and $\varphi$ 
is uniformly continuous. Therefore,
\begin{eqnarray*}
\P(T,\varphi)
&=&\sup_{{\mathcal V}}\bigl[\P(T,\varphi,{\mathcal V})-\mathrm{osc}(\varphi,{\mathcal V})\bigr] 
\geq\limsup_{n\to\infty}\bigl\{\P(T,\varphi,{\mathcal U}_n)-\mathrm{osc}(\varphi,{\mathcal U}_n)\bigr\} \\
&=&\limsup_{n\to\infty}\P(T,\varphi,{\mathcal U}_n)-\lim_{n\to\infty}\mathrm{osc}(\varphi,{\mathcal U}_n) \\
&=&\limsup_{n\to\infty}\P(T,\varphi,{\mathcal U}_n) \\
&\geq&\liminf_{n\to\infty}\P(T,\varphi,{\mathcal U}_n)\\
&\geq& \P(T,\varphi).
\end{eqnarray*}
Hence $\P(T,\varphi)=\lim_{n\to\infty}\P(T,\varphi,{\mathcal U}_n)$.
\epf


\sp We can now obtain a slightly stronger estimate than Remark~\ref{ch9presprop}(b) for the difference between topological entropy and topological pressure when the underlying space is metrizable.

\bcor\label{ch9cor2.1.19}
If $T:X\lra X$ is a topological dynamical system and $\varphi:X\lra\R$ is a continuous potential, then 
$$
\htop(T)+\inf(\varphi)
\leq \P(T,\varphi)\leq\htop(T)+\sup(\varphi).
$$
\ecor

\bpf
The upper bound was already mentioned in Remark~\ref{ch9presprop}(b). In order to derive the lower
bound, we  return to Remark~\ref{ch9presrelrmk}(b). Let $({\mathcal U}_n)_{n=1}^\infty$ be
a sequence of open covers of $X$ such that $\lim_{n\to\infty}\mathrm{diam}({\mathcal U}_{n})=0$.
According to Remark~\ref{ch9presrelrmk}(b), for each $n\in\N$ we have
\[
\htop(T,{\mathcal U}_n)+\inf(\varphi)\leq \P(T,\varphi,{\mathcal U}_n).
\]
Passing to the limit $n\to\infty$ and using Proposition~\ref{ch9sem2lem9.5},
we conclude that
\[
\htop(T)+\inf(\varphi)\leq \P(T,\varphi).
\]
\epf

\sp The preceding lemma characterized the topological pressure of a potential as the limit
of the topological pressure of the potential relative to a sequence of covers. An
even better result would be the characterization of the topological pressure of
a potential as the topological pressure of that potential with respect to a single cover. As might by now be expected, such a characterization exists when the system is expansive. We will therefore now devote a little bit of time to look closer at expansive maps. 

Firstly, the expansiveness of a system can  be expressed in terms of the following \index{(N)}{dynamical (or Bowen) metric} ``dynamical'', called also Bowen, metrics.

\bdfn\label{dnmetric}
Let $T:(X,d)\lra (X,d)$ be a topological dynamical system. For every $n\in\N$,
let $d_n:X\times X\lra[0,\infty)$ be the metric define by the following formula:
\[
d_n(x,y):=\max\bigl\{d\bigl(T^j(x),T^j(y)\bigr):0\leq j\le n-1\bigr\}.
\]
\edfn

\fr Although this notation does not make explicit the dependence on $T$, it is crucial to remember that the metrics $d_n$ arise from the dynamics of the system $T$. It is in this sense that they are dynamical metrics.

Observe that $d_1=d$ and that for each $x,y\in X$ we have 
$$
d_n(x,y)\geq d_m(x,y)
$$ 
whenever $n\geq m$. Moreover, it is worth noticing that 

\sp \centerline{ all the metrics $d_n$, $n\in\N$, are topologically
equivalent.}

\sp\fr Indeed, given a sequence $(x_k)_{k=1}^\infty$ in $X$, the
continuity of $T$ ensures that
\begin{eqnarray*}
\lim_{k\to\infty}d(x_k,y)=0
&\Longleftrightarrow&
\lim_{k\to\infty}d\bigl(T^j(x_k),T^j(y)\bigr)=0,
\ \ \ \forall\, 0\leq j\le n-1, \ \ \ \forall\, n\in\N \\ 
&\Longleftrightarrow&\lim_{k\to\infty}d_n(x_k,y)=0, \ \ \ \forall\, n\in\N.
\end{eqnarray*}

Furthermore, it is easy to see that a dynamical system $T:(X,d)\lra(X,d)$ is $\delta$--expansive if and only if 
\[
\sup\big\{d_n(x,y):n\in\N\big\}\leq\delta \ \Longrightarrow \ x=y.
\]

\bprop\label{sem1prop5.0}
A topological dynamical system $T:(X,d)\lra(X,d)$
is $\delta$--expansive if and only if for every 
$\eta>0$ there exists $N=N(\eta)\in\N$ such that 
\[
d(x,y)>\eta \ \Longrightarrow \ d_N(x,y)>\delta.
\]
\eprop

\bpf
The implication $(\imp)$ is obvious. For the converse one, suppose by way of contradiction that $T:(X,d)\lra(X,d)$ is a
$\delta$--expansive system and the assertion of our proposition fails. 
Then there exist $\eta>0$ and two sequences 
$$
(x_n)_{n=0}^\infty 
\  \  {\rm and} \  \ 
(y_n)_{n=0}^\infty
$$ 
in $X$ such that 
$$
d(x_n,y_n)>\eta
\  \  {\rm but} \ \ 
d_n(x_n,y_n)\leq\delta.
$$ 
Since $X$ is compact, we may assume, by passing to subsequences if necessary, that the sequences
$(x_n)_{n=0}^\infty$ and $(y_n)_{n=0}^\infty$ converge to, say, $x$
and $y$, respectively. On one hand, this implies that
\[
d(x,y)=\lim_{n\to\infty}d(x_n,y_n)\geq\eta>0,
\]
and hence 
\beq\label{120190930}
x\neq y.
\eeq 
On the other hand, if we fix momentarily any $N\in\N$, then for all $n\geq N$, we have that
\[
d_N(x_n,y_n)\leq d_n(x_n,y_n)\leq\delta.
\]
Therefore,
\[
d_N(x,y)=\lim_{n\to\infty}d_N(x_n,y_n)\leq\delta.
\]
Since $d_N(x,y)\leq\delta$ for every $N\in\N$, the
$\delta$--expansiveness of the system implies that $x=y$, contrary to \eqref{120190930}. We are done.
\epf

\sp As an immediate consequence of this proposition, we get the following.

\bcor\label{c120190930}
If $(X,d)$ is a compact metric space and $T:X\lra X$ is an expansive topological dynamical system with an expansive constant $\delta$,
then 
$$
\lim_{n\to\infty}\mathrm{diam}(\mathcal{U}^n)=0
$$
for every cover $\mathcal{U}$ with $\mathrm{diam}(\mathcal{U})\leq\delta$. 
\ecor

\fr We can now easily prove the following.

\bthm\label{ch9sem2cor9.11}
If $(X,d)$ is a compact metric space and $T:X\lra X$ is an expansive topological dynamical system with an expansive constant $\delta$, then 
\[
\P(T,\varphi)
=\underline{p}(T,\varphi,{\mathcal U})
=\overline{p}(T,\varphi,{\mathcal U})
\]
for any finite open cover ${\mathcal U}$ of $X$ with
$\mathrm{diam}({\mathcal U})\leq\delta$. 
Moreover,
\[
\P(T,\varphi)=\P(T,\varphi,{\mathcal U})
\]
for any finite open partition ${\mathcal U}$ of $X$ with
$\mathrm{diam}({\mathcal U})\leq\delta$.
\ethm

\bpf 
It follows from Corollary~\ref{c120190930}, Proposition~\ref{ch9sem2lem9.5} and~\ref{ch9sem2cor9.10}, that
\[
\P(T,\varphi)
=\lim_{n\to\infty}\overline{p}(T,\varphi,{\mathcal U}^n)
=\lim_{n\to\infty}\overline{p}(T,\varphi,{\mathcal U})
=\overline{p}(T,\varphi,{\mathcal U}).
\]
A similar argument leads to the statements for $\underline{p}(T,\varphi,{\mathcal U})$ and ${\mathcal U})$ being a partition.
 \epf

For the final result of this section, we study the behavior of
topological pressure with respect to the iterates of the system.
\bthm\label{ch9sem2prop9.7}
If $(X,d)$ is a compact metric space and $T:X\lra X$ is a topological dynamical system, then 
$$
\P(T^n,S_n\varphi)=n\,\P(T,\varphi)
$$
for every $n\in\N$. 
\ethm

\bpf
Fix $n\in\N$. Let ${\mathcal U}$ be an open cover of $X$.
The action of the map $T^n$ on ${\mathcal U}$ until time $j-1$ will be denoted by ${\mathcal U}_{T^n}^j$. Note that ${\mathcal U}^{mn}=({\mathcal U}^n)_{T^n}^m$ for all $m\in\N$. Furthermore, for all $x\in X$,
\[
S_{mn}\varphi(x)=\sum_{k=0}^{mn-1}\varphi\circ T^k(x)
=\sum_{j=0}^{m-1}(S_n\varphi)\circ T^{jn}(x)
=\sum_{j=0}^{m-1}(S_n\varphi)\circ (T^n)^j(x).
\]
Hence 
$$
\overline{S}_{mn}\varphi(Y)=S^{T^n}_m(\overline{S}_n\varphi)(Y)
$$ 
for all subsets $Y$ of $X$, and in particular for all $Y\in{\mathcal U}^{mn}=({\mathcal U}^n)_{T^n}^m$, where
$$
S^{T^n}_m\psi(x)=\sum_{j=0}^{m-1}\psi((T^n)^j(x)).
$$
Thus,
\[
Z_{mn}(T,\varphi,{\mathcal U})=Z_m(T^n,S_n\varphi,{\mathcal U}^n.
\]
for all $m\in\N$. Therefore,
\begin{eqnarray*}
\P(T,\varphi,{\mathcal U})
&=&\lim_{m\to\infty}\frac{1}{mn}\log Z_{mn}(T,\varphi,{\mathcal U})
=\frac{1}{n}\lim_{m\to\infty}\frac{1}{m}\log Z_m(T^n,S_n\varphi,{\mathcal U}^n)\\
&=&\frac{1}{n}\P(T^n,S_n\varphi,{\mathcal U}^n).
\end{eqnarray*}
Let $({\mathcal U}_k)_{k=1}^\infty$ be a sequence of open covers such that
$$
\lim_{k\to\infty}\mathrm{diam}({\mathcal U}_k)=0.
$$
Then
$$
\lim_{k\to\infty}\mathrm{diam}(({\mathcal U}_k)^n)=0
$$ 
since $\mathrm{diam}(({\mathcal U}_k)^n)\leq\mathrm{diam}({\mathcal U}_k)$
for all $k\in\N$. It thus follows from Proposition~\ref{ch9sem2lem9.5} that
\[
\P(T^n,S_n\varphi)
=\lim_{k\to\infty}\P(T^n,S_n\varphi,({\mathcal U}_k)^n)
=n\lim_{k\to\infty}\P(T,\varphi,{\mathcal U}_k)
=n\,\P(T,\varphi).
\]
\epf

\section{Bowen's Definition of Topological Pressure}\label{ch9defpresbowen}

In this section we provide the characterization of topological pressure, due to Rufus Bowen, by means of separated or spanning sets. Although this characterization does depend on a metric and is not immediately seen to be topologically invariant, it nevertheless has both, theoretical and practical (calculation of pressures) advantages. We will define and discuss separated and  spanning sets now and will then characterize topological pressure in terms of them. 

\bdfn\label{ch6nesep}
A subset $F$ of $X$ is said to be \index{(N)}{separated set, $(n,\e)$-separated}{\em $(n,\e)$--separated} if $F$ is 
\mbox{$\e$-separated} with respect to the metric $d_n$, which is to say that
$$
d_n(x,y)\geq\e
$$ 
for all $x,y\in F$ with $x\neq y$.
\edfn

\brem\label{ch6nesepprop}
\
\begin{itemize}
\item[(1)] If $F$ is an $(n,\varepsilon)$--separated set and $m\geq n$,
then $F$ is also $(m,\varepsilon)$--separated.

\sp\item[(2)] If $F$ is an $(n,\varepsilon)$--separated set and
$0<\varepsilon'<\varepsilon$, then $F$ is also $(n,\varepsilon')$--separated.

\sp\item[(3)] Given that $X$ is a compact metric space, any
$(n,\e)$--separated set is finite. Indeed, 
let $F$ be an $(n,\e)$--separated set, and consider the family of balls
$$
\{B_n(x,\e/2):x\in F\}.
$$
If the intersection of  $B_n(x,\e/2)$ and $
B_n(y,\e/2)$ is non-empty for 
some $x,y\in F$, then there exists 
$$
z\in B_n(x,\e/2)\cap B_n(y,\e/2),
$$ 
and it follows that
\[
d_n(x,y)\leq d_n(x,z)+d_n(z,y)<\e/2+\e/2=\e.
\]
As $F$ is an $(n,\e)$--separated set, this inequality implies that $x=y$. This
means that the balls 
$$
\{B_n(x,\e/2): x\in F\}
$$ 
are mutually
disjoint. Hence, as we are in a compact metric space, there can 
only be finitely many of them.
\end{itemize}
\erem

\sp\fr The largest separated sets will be especially useful in describing the
complexity of the dynamics that the system exhibits.

\bdfn\label{ch6maxnesep}
A subset $F$ of $X$ is called a \index{(N)}{separated set!maximal
  $(n,\e)$-separated}{\em maximal $(n,\e)$--separated set} if 
for any $(n,\e)$-separated set $F'$ with $F\subseteq F'$, we have
$F=F'$. In other words, 
no strict superset of $F$ is $(n,\e)$-separated.
\edfn

\

\fr The counterpart of the notion of a separated set is the concept of the
spanning set. 
\bdfn\label{ch6nespan}
A subset $E$ of $X$ is said to be an \index{(N)}{spanning
set!$(n,\e)$-spanning}{\em $(n,\e)$-spanning set} if
\[
\bigcup_{x\in E}B_n(x,\e)=X.
\]
That is, the orbit of every point in the space
is $\e$-shadowed by the orbit of a point of $E$ until at least time $n-1$.
\edfn

\

\fr The smallest spanning sets play a special role in describing the
complexity of the dynamics that the system possesses.

\

\bdfn\label{ch6maxnespan}
A subset $E$ of $X$ is called a \index{(N)}{spanning set!minimal
$(n,\e)$--spanning}{\em minimal $(n,\e)$--spanning set} if 
for any $(n,\e)$--spanning set $E'$ with $E\supseteq E'$, we have
$E=E'$. In other words, no strict subset of $E$ is $(n,\e)$-spanning.
\edfn

\

\brem\label{ch6nespanprop}
\
\begin{itemize}

\item[(1)] If $E$ is an $(n,\varepsilon)$--spanning set and $m\leq n$, then
           $E$ is also $(m,\varepsilon)$--spanning.

\sp\item[(2)] If $E$ is an $(n,\varepsilon)$--spanning set and
$\varepsilon'>\varepsilon$, then $E$ is also $(n,\varepsilon')$--spanning.

\sp\item[(3)] Given that we are in a compact metric space, any minimal
$(n,\e)$--spanning set $E\sbt X$ is finite since the open cover $\{B_n(x,\varepsilon):x\in E\}$ of the compact metric space $X$ admits a finite subcover.
\end{itemize}
\erem

\fr The next lemma describes two useful relations between separated
and spanning sets. 

\blem\label{ch6sepspan}
The following statements hold:
\begin{itemize}
\item[$(a)$]
Every maximal $(n,\e)$--separated set is an $(n,\e)$--spanning set.
\item[$(b)$] Every $(n,2\varepsilon)$--separated set can be embedded into any $(n,\varepsilon)$--spanning set.
\end{itemize}
\elem

\bpf
(a). Let $F$ be a maximal $(n,\e)$--separated set. By way of
contradiction, suppose that $F$ is not an $(n,\e)$--spanning set. Then
there would exist a point 
$$y
\in X\backslash\bigcup_{x\in F}B_n(x,\e).
$$
It is then easy to verify that the set 
$$
F\cup\{y\}
$$ 
is $(n,\e)$-separated, hence contradicting the maximality of $F$. Therefore, if $F$ is a maximal $(n,\e)$--separated set, then  $F$ is also
$(n,\varepsilon)$--spanning. 

\sp\fr (b). Let $F$ be an $(n,2\varepsilon)$--separated
set and $E$ an $(n,\varepsilon)$-spanning set. For each $x\in F$,
choose $i(x)\in E$ such that 
$$
x\in B_n(i(x),\varepsilon).
$$
We claim that the map 
$$
i:F\lra E
$$ 
is injective. In order to show this, let $x,y\in F$ be such that 
$$
i(x)=i(y)=:z.
$$
Then
$x,y\in B_n(z,\varepsilon)$. Therefore $d_n(x,y)<2\varepsilon$. Since
$F$ is a $(n,2\varepsilon)$-separated, we deduce that $x=y$, that is,
the map $i$ is injective and the proof is complete. 
\epf

\sp  We are now ready to provide the announced characterization of
topological pressure that is based on the concepts of separated and
spanning sets. To allege notation, for any $n\in\N$ and $Y\sbt X$, let 
\[
\Sigma_n(Y)=\sum_{x\in Y}e^{S_n\varphi(x)}.
\] 
%
%

\bthm\label{ch9sem2thm9.12}
Let $(X,d)$ be a compact metric space and let $T:X\lra X$ be a topological dynamical system. For all $n\in\N$ and for all $\varepsilon>0$, let $F_n(\varepsilon)$ be a maximal
$(n,\varepsilon)$--separated set and $E_n(\varepsilon)$ be a minimal $(n,\varepsilon)$--spanning set. Then
\begin{eqnarray*}
\P(T,\varphi)
&=&\lim_{\varepsilon\to0}\liminf_{n\to\infty}\frac{1}{n}\log \Sigma_n(E_n(\varepsilon))  
 = \lim_{\varepsilon\to0}\limsup_{n\to\infty}\frac{1}{n}\log \Sigma_n(E_n(\varepsilon)) \\
&=&\lim_{\varepsilon\to0}\limsup_{n\to\infty}\frac{1}{n}\log \Sigma_n(F_n(\varepsilon)) 
 = \lim_{\varepsilon\to0}\liminf_{n\to\infty}\frac{1}{n}\log \Sigma_n(F_n(\varepsilon)).
\end{eqnarray*}
\ethm

\bpf
Fix $\varepsilon>0$. Let ${\mathcal U}_\varepsilon$ be an open cover of $X$
consisting of balls of radius $\varepsilon/2$. Fix $n\in\N$. Let ${\mathcal U}$ be a subcover of ${\mathcal U}_\varepsilon^n$ such that 
$$
Z_n(\varphi,{\mathcal U}_\varepsilon)
\geq e^{-1}\sum_{U\in{\mathcal U}}\exp(\overline{S}_n\varphi(U)).
$$
For each $x\in F_n(\varepsilon)$, let $U(x)$ be an element of the cover ${\mathcal U}$
which contains $x$ and define the function $i:F_n(\varepsilon)\to{\mathcal U}$ by setting
$$
i(x)=U(x).
$$
Since $E_n(\varepsilon)$ and ${\mathcal U}$ is a subcover of ${\mathcal U}_\varepsilon^n$, it follows that this function is injective. Therefore,
\begin{eqnarray*}
Z_n(\varphi,{\mathcal U}_\varepsilon)
\geq e^{-1}\sum_{U\in{\mathcal U}}e^{\overline{S}_n\varphi(U)}
\geq e^{-1}\sum_{x\in F_n(\varepsilon)}e^{\overline{S}_n\varphi(U(x))}
\geq e^{-1}\sum_{x\in F_n(\varepsilon)}e^{S_n\varphi(x)}.
\end{eqnarray*}
Since this is true for all $n\in\N$, we deduce that
\[
\P(T,\varphi,{\mathcal U}_\varepsilon)
=\lim_{n\to\infty}\frac{1}{n}\log Z_n(\varphi,{\mathcal U}_\varepsilon)
\geq\limsup_{n\to\infty}\frac{1}{n}\log\sum_{x\in F_n(\varepsilon)}e^{S_n\varphi(x)}.
\]
Letting $\varepsilon\to0$ and using Proposition~\ref{ch9sem2lem9.5} yields that
\beq\label{ch9presineq1}
\P(T,\varphi)
\geq\limsup_{\varepsilon\to0}\limsup_{n\to\infty}\frac{1}{n}\log\sum_{x\in E_n(\varepsilon)}e^{S_n\varphi(x)}.
\eeq

On the other hand, if ${\mathcal V}$ is an arbitrary open cover of $X$, if $\delta({\mathcal V})>0$ is a
Lebesgue number for ${\mathcal V}$, if $0<\varepsilon<\delta({\mathcal V})/2$, and if $n\in\N$, then for all integers
$0\leq k\le n-1$ and for all $x\in F_n(\varepsilon)$, we have that
$$
T^k(B_n(x,\varepsilon))
\sbt B(T^k(x),\varepsilon).
$$
Consequently,
$$
\mathrm{diam}\left(T^k(B_n(x,\varepsilon))\right)\leq2\varepsilon<\delta({\mathcal V}).
$$
Hence, for all integers $0\leq k\le n-1$, the set $T^k(B_n(x,\varepsilon))$ is
contained in at least one element of the cover ${\mathcal V}$. Denote one of these elements by
$V_k(x)$. It follows that 
$$
B_n(x,\varepsilon)\sbt T^{-k}(V_k(x))
$$ 
for each integer $0\leq k\le n-1$ or,
in other words, 
$$
B_n(x,\varepsilon)\sbt\bigcap_{k=0}^{n-1}T^{-k}(V_k(x)).
$$
But this latter intersection is an element of ${\mathcal V}^n$. Let us denote it by 
$$
V(x).
$$
Since $F_n(\varepsilon)$ is a maximal $(n,\varepsilon)$--separated set, by
Lemma~\ref{ch6sepspan} it is also
$(n,\varepsilon)$--spanning. So, the family 
$$
\{B_n(x, \varepsilon)\}_{x\in F_n(\varepsilon)}
$$ 
is an open cover of $X$. Each one of these balls is contained in the corresponding set
$V(x)$. Hence the family $\{V(x)\}_{x\in F_n(\varepsilon)}$ is also an open cover of
$X$. Therefore it is a subcover of ${\mathcal V}^n$. 
Consequently,
\[
Z_n(\varphi,{\mathcal V})
\leq\sum_{x\in F_n(\varepsilon)}e^{\overline{S}_n\varphi(V(x))}
\leq e^{n\,\mathrm{osc}(\varphi,{\mathcal V})}\sum_{x\in F_n(\varepsilon)}e^{S_n\varphi(x)},
\]
where the last inequality is due to Lemma~\ref{oscsnun}. It follows that
\[
\P(T,\varphi,{\mathcal V})
\leq\mathrm{osc}(\varphi,{\mathcal V})
    +\liminf_{n\to\infty}\frac{1}{n}\log\sum_{x\in E_n(\varepsilon)}e^{S_n\varphi(x)}.
\]
Since ${\mathcal V}$ is independent of $\varepsilon>0$, we deduce that
\[
\P(T,\varphi,{\mathcal V})-\mathrm{osc}(\varphi,{\mathcal V})
\leq\liminf_{\varepsilon\to0}\liminf_{n\to\infty}\frac{1}{n}\log\sum_{x\in F_n(\varepsilon)}e^{S_n\varphi(x)}.
\]
Then, as ${\mathcal V}$ was chosen to be an arbitrary open cover of $X$, we conclude that
\beq\label{ch9presineq2}
\P(T,\varphi)
\leq\liminf_{\varepsilon\to0}\liminf_{n\to\infty}\frac{1}{n}\log\sum_{x\in E_n(\varepsilon)}e^{S_n\varphi(x)}.
\eeq
The inequalities~(\ref{ch9presineq1}) and~(\ref{ch9presineq2}) taken together complete the proof of our theorem.
\epf

\sp In Theorem~\ref{ch9sem2thm9.12}, the topological pressure of the system is expressed in terms of a specific family of maximal separated (resp. minimal spanning) sets. However, to derive theoretical results, it is sometimes simpler to use the following quantities.

\bdfn\label{ch6sem1defn7.5gen}
Let $(X,d)$ be a compact metric space and let $T:X\lra X$ be a topological dynamical system. For all $n\in\N$ and $\varepsilon>0$, let 
\[
\P_n(T,\varphi,\varepsilon):=\sup\Bigl\{\Sigma_n(E_n(\varepsilon))\,:\,E_n(\varepsilon) \text{ maximal
$(n,\varepsilon)$--separated set}\Bigr\} 
\]
and
\[
Q_n(T,\varphi,\varepsilon):=\inf\Bigl\{\Sigma_n(F_n(\varepsilon))\,:\,F_n(\varepsilon) \text{ minimal
$(n,\varepsilon)$--spanning set}\Bigr\}. 
\]
Thereafter, let
\[
\underline{\P}(T,\varphi,\varepsilon):=\liminf_{n\to\infty}\frac{1}{n}\log \P_n(T,\varphi,\varepsilon),
\hspace{0.4cm}
\overline{\P}(T,\varphi,\varepsilon):=\limsup_{n\to\infty}\frac{1}{n}\log \P_n(T,\varphi,\varepsilon)
\]
and
\[
\underline{Q}(T,\varphi,\varepsilon):=\liminf_{n\to\infty}\frac{1}{n}\log Q_n(T,\varphi,\varepsilon),
\hspace{0.4cm}
\overline{Q}(T,\varphi,\varepsilon):=\limsup_{n\to\infty}\frac{1}{n}\log Q_n(T,\varphi,\varepsilon).
\]
\edfn

The following are simple but key observations about these quantities. 

\brem\label{ch9sem2cor9.12r}
Let $(X,d)$ be a compact metric space and let $T:X\lra X$ be a topological dynamical system. Let $m\leq n\in\N$ and $0<\varepsilon<\varepsilon'$. The following relations hold:
\begin{itemize}
\item[(a)] By Remark~\ref{ch6nesepprop}(a):
$$
\P_m(T,\varphi,\varepsilon)\leq \P_n(T,\varphi,\varepsilon)e^{(n-m)\|\varphi\|_\infty}
$$ 

\item[(b)] 
$$
e^{-n\|\varphi\|_\infty}\leq \P_n(T,\varphi,\varepsilon)\leq r_n(\varepsilon)e^{n\|\varphi\|_\infty} 
$$
and
$$
\P_n(T,0,\varepsilon)=r_n(\varepsilon).
$$ 
\item[(c)] By Remark~\ref{ch6nespanprop}(a):
$$
Q_m(T,\varphi,\varepsilon)
\leq Q_n(T,\varphi,\varepsilon)e^{(n-m)\|\varphi\|_\infty}.
$$ 

\item[(d)] 
$$
e^{-n\|\varphi\|_\infty}\leq Q_n(T,\varphi,\varepsilon)\leq s_n(\varepsilon)e^{n\|\varphi\|_\infty}
$$
and 
$$
Q_n(T,0,\varepsilon)=s_n(\varepsilon).
$$ 
\item[(e)] By Remarks~\ref{ch6nesepprop} and~\ref{ch6nespanprop}(b):
$$
\P_n(T,\varphi,\varepsilon)\geq \P_n(T,\varphi,\varepsilon')
$$ 
and 
$$
Q_n(T,\varphi,\varepsilon)\geq Q_n(T,\varphi,\varepsilon')
$$ 

\item[(f)] By Lemma~\ref{ch6sepspan}:
$$
0<Q_n(T,\varphi,\varepsilon)\leq P_n(T,\varphi,\varepsilon)<\infty.
$$ 

\item[(g)] 
$$
\underline{\P}(T,\varphi,\varepsilon)
\leq\overline{\P}(T,\varphi,\varepsilon)
$$ 
and 
$$
\underline{Q}(T,\varphi,\varepsilon)
\leq\overline{Q}(T,\varphi,\varepsilon).
$$
\item[(h)] By~(b):
$$
-\|\varphi\|_\infty<\underline{\P}(T,\varphi,\varepsilon)
\leq\underline{r}(\varepsilon)+\|\varphi\|_\infty
$$ 
and 
$$
-\|\varphi\|_\infty<\overline{\P}(T,\varphi,\varepsilon)\leq\overline{r}(\varepsilon)+\|\varphi\|_\infty.
$$  
\item[(i)]  By~(d):
$$
-\|\varphi\|_\infty<\underline{Q}(T,\varphi,\varepsilon)\leq\underline{s}(\varepsilon)+\|\varphi\|_\infty
$$ 
           and 
$$
-\|\varphi\|_\infty<\overline{Q}(T,\varphi,\varepsilon)\leq\overline{s}(\varepsilon)+\|\varphi\|_\infty.
$$ 
\item[(j)] By~(e):
$$
\underline{\P}(T,\varphi,\varepsilon)\geq\underline{\P}(T,\varphi,\varepsilon')
$$ 
and 
$$
\overline{\P}(T,\varphi,\varepsilon)\geq\overline{\P}(T,\varphi,\varepsilon').
$$
\item[(k)] By~(e):
$$
\underline{Q}(T,\varphi,\varepsilon)\geq\underline{Q}(T,\varphi,\varepsilon')
$$ 
and 
$$
\overline{Q}(T,\varphi,\varepsilon)\geq\overline{Q}(T,\varphi,\varepsilon').
$$
\item[(l)] By~(f):
$$
-\|\varphi\|_\infty\leq\overline{Q}(T,\varphi,\varepsilon)\leq\overline{\P}(T,\varphi,\varepsilon)\leq\infty
$$ 
and 
$$
-\|\varphi\|_\infty\leq\underline{Q}(T,\varphi,\varepsilon)\leq\underline{\P}(T,\varphi,\varepsilon)\leq\infty
$$ 
\end{itemize}
\erem

\sp We will now describe a relationship between $\P_n$'s, $Q_n$'s, and the cover--related quantities $Z_n$'s and $z_n$'s.
\blem\label{zPQ}
If $(X,d)$ is a compact metric space and $T:X\lra X$ is a topological dynamical system, then the following relations hold.
\begin{itemize}
\item[(a)] If ${\mathcal U}$ is an open cover of $X$ with Lebesgue number $2\delta$, then
\[
z_n(T,\varphi,{\mathcal U})\leq Q_n(T,\varphi,\delta)\leq \P_n(T,\varphi,\delta).
\]
\item[(b)] If $\varepsilon>0$ and ${\mathcal V}$ is an open cover of $X$ with 
$\mathrm{diam}({\mathcal V})\leq\varepsilon$, then
\[
Q_n(T,\varphi,\varepsilon)\leq \P_n(T,\varphi,\varepsilon)\leq Z_n(T,\varphi,{\mathcal V}).
\]
\end{itemize}
\elem

\bpf We already know that 
$$
Q_n(T,\varphi,\delta)\leq P_n(T,\varphi,\delta).
$$

(a) Let ${\mathcal U}$ be an open cover with Lebesgue number $2\delta$ 
and let $F$ be an $(n,\delta)$--spanning set. Then the dynamic balls 
$\{B_n(x,\delta):x\in F\}$ form a cover of $X$. For every $0\leq i\le n-1$ the ball 
$B(T^i(x),\delta)$, which has diameter at most $2\delta$, is contained
in an element of ${\mathcal U}$. Therefore 
$$
B_n(x,\delta)=\bigcap_{i=0}^{n-1}T^{-i}(B(T^i(x),\delta))
$$
is contained in an element of 
$$
{\mathcal U}^n=\bigvee_{i=0}^{n-1}T^{-i}({\mathcal U}). 
$$
That is, 
$$
{\mathcal U}^n\prec\{B_n(x,\delta):x\in F\}.
$$ 
Then there exists a map 
$$
i:\{B_n(x,\delta):x\in F\}\to{\mathcal U}^n
$$ 
such that
$B_n(x,\delta)\sbt i(B_n(x,\delta))$ for every $x\in F$. Let 
${\mathcal W}$ be a subcover of $\{B_n(x,\delta):x\in F\}$. 
Then
\begin{eqnarray*}
\Sigma_n(F)=\sum_{x\in F}e^{S_n\varphi(x)}
&\geq&\sum_{x\in F}e^{\underline{S}_n\varphi(B_n(x,\delta))} 
\geq\sum_{W\in{\mathcal W}}e^{\underline{S}_n\varphi(W)} 
\geq\sum_{W\in{\mathcal W}}e^{\underline{S}_n\varphi(i(W))} \\
&\geq&\sum_{Z\in i({\mathcal W})}e^{\underline{S}_n\varphi(Z)} 
\geq z_n(T,\varphi,{\mathcal U}).
\end{eqnarray*}
Since $F$ is an arbitrary $(n,\delta)$--spanning set, 
it ensues that 
$$
Q_n(T,\varphi,\delta)\geq z_n(T,\varphi,{\mathcal U}).
$$

(b) Let ${\mathcal V}$ be an open cover with $\mathrm{diam}({\mathcal V})\leq\varepsilon$
and let $E$ be an $(n,\varepsilon)$--separated set. Then no element of 
the cover ${\mathcal V}^n$ can contain more than one element of $E$. 
Let ${\mathcal W}$ be a subcover of ${\mathcal V}^n$ and let a map $i:E\to{\mathcal W}$ be such that
the fact that $x\in i(x)$ for all $x\in E$. Then 
\[
\Sigma_n(E)=\sum_{x\in E}e^{S_n\varphi(x)}
\leq\sum_{x\in E}e^{\overline{S}_n\varphi(i(x))} 
=\sum_{W\in i(E)}e^{\overline{S}_n\varphi(W)} 
\leq\sum_{W\in{\mathcal W}}e^{\overline{S}_n\varphi(W)}.
\]
As ${\mathcal W}$ is an arbitrary subcover of ${\mathcal V}^n$, it follows that 
$\Sigma_n(E)\leq Z_n(T,\varphi,{\mathcal V})$.
Since $E$ is an arbitrary $(n,\varepsilon)$-separated set, 
we deduce that 
$$
\P_n(T,\varphi,\varepsilon)\leq Z_n(T,\varphi,{\mathcal V}).
$$
\epf

\fr These inequalities have the following immediate consequences.

\bcor\label{zPQc}
If $(X,d)$ is a compact metric space and $T:X\lra X$ is a topological dynamical system, then the following relations hold.
\begin{itemize}
\item[(a)] If ${\mathcal U}$ is an open cover of $X$ with Lebesgue number $2\delta$, then
\[
\underline{p}(T,\varphi,{\mathcal U})
\leq\underline{Q}(T,\varphi,\delta)
\leq\underline{\P}(T,\varphi,\delta).
\]
\item[(b)] If $\varepsilon>0$ and ${\mathcal V}$ is an open cover of $X$ with 
$\mathrm{diam}({\mathcal V})\leq\varepsilon$, then
\[
\overline{Q}(T,\varphi,\varepsilon)
\leq\overline{\P}(T,\varphi,\varepsilon)\leq \P(T,\varphi,{\mathcal V}).
\]
\end{itemize}
\ecor

\fr We can then surmise new expressions for the topological pressure.

\bcor\label{ch9sem2cor9.12}
If $(X,d)$ is a compact metric space and $T:X\lra X$ is a topological dynamical system, then the following equalities hold.  
\begin{eqnarray*}
\P(T,\varphi)
=\lim_{\varepsilon\to0}\underline{\P}(T,\varphi,\varepsilon)  
=\lim_{\varepsilon\to0}\overline{\P}(T,\varphi,\varepsilon)  
=\lim_{\varepsilon\to0}\underline{Q}(T,\varphi,\varepsilon)  
=\lim_{\varepsilon\to0}\overline{Q}(T,\varphi,\varepsilon).  
\end{eqnarray*}
\ecor

\bpf
Let $({\mathcal U}_\varepsilon)_{\varepsilon\in(0,\infty)}$ be a family of open covers such that 
$$
\lim_{\varepsilon\to\infty}\mathrm{diam}({\mathcal U}_\varepsilon)=0.
$$
Let $\delta_\varepsilon>0$ be a Lebesgue number for ${\mathcal U}_\varepsilon$. Then 
$$
\lim_{\varepsilon\to\infty}\delta_\varepsilon=0,
$$ 
as $\delta_\varepsilon\leq\mathrm{diam}({\mathcal U}_\varepsilon)$.
Using Proposition~\ref{ch9sem2lem9.5} and Corollary~\ref{zPQc}(a), we deduce that
\beq\label{PpQP}
\P(T,\varphi)
=\lim_{\varepsilon\to\infty}\underline{p}(T,\varphi,{\mathcal U}_\varepsilon)
\leq\lim_{\varepsilon\to0}\underline{Q}(T,\varphi,\varepsilon)
\leq\lim_{\varepsilon\to0}\underline{\P}(T,\varphi,\varepsilon).
\eeq
On the other hand, using Proposition~\ref{ch9sem2lem9.5} and Corollary~\ref{zPQc}(b), we obtain
\beq\label{QPPP}
\lim_{\varepsilon\to0}\overline{Q}(T,\varphi,\varepsilon)
\leq\lim_{\varepsilon\to0}\overline{\P}(T,\varphi,\varepsilon)
\leq\lim_{\varepsilon\to0}\sup\Big\{\P(T,\varphi,{\mathcal V}):\mathrm{diam}({\mathcal V})\leq\varepsilon\Big\}
=\P(T,\varphi).
\eeq
Combining~(\ref{PpQP}) and~(\ref{QPPP}) allows us to conclude. 
\epf

\sp Corollary~\ref{ch9sem2cor9.12} is useful to derive theoretical results. Nevertheless, 
in practice, Theorem~\ref{ch9sem2thm9.12} is simpler to use, as only one family
of sets is needed. Sometimes a single sequence of sets is enough. 

\bthm\label{ch9sem2cor9.11bowenseqpre}
Let $(X,d)$ be a compact metric space and let $T:X\lra X$ be an expansive topological dynamical system with an expansive constant $\delta$. If ${\mathcal U}$ is an open cover of $X$ whose Lebesgue number is $2\eta$ with some $\eta\in (0,\d/2)$ (in particular
$\mathrm{diam}({\mathcal U})\leq \delta$), then the following statements hold
for all $0<\varepsilon\leq\eta$:
\begin{itemize}
\item[(a)] If $(E_n(\varepsilon))_{n=1}^\infty$ is a sequence of maximal 
$(n,\varepsilon)$-separated sets in $X$, then 
\[
\P(T,\varphi)
=\lim_{n\to\infty}\frac{1}{n}\log \Sigma_n(E_n(\varepsilon)).
\]
\item[(b)] If $(F_n(\varepsilon))_{n=1}^\infty$ is a sequence of minimal 
$(n,\varepsilon)$-spanning sets in $X$, then 
\[
\P(T,\varphi)
\leq\liminf_{n\to\infty}\frac{1}{n}\log \Sigma_n(F_n(\varepsilon)).
\]
\item[(c)] 
$\displaystyle \P(T,\varphi)=\lim_{n\to\infty}\frac{1}{n}\log \P_n(T,\varphi,\varepsilon)$.
\vspace{1.5mm}
\item[(d)] 
$\displaystyle \P(T,\varphi)=\lim_{n\to\infty}\frac{1}{n}\log Q_n(T,\varphi,\varepsilon)$.
\end{itemize}
\ethm

\bpf
We will prove~(a) and leave it to the reader to show the other parts using similar arguments.
  
It follows from Theorem~\ref{ch9sem2cor9.11} that 
$$
\P(T,\varphi)=\underline{p}(T,\varphi,{\mathcal U})
$$ 
Fix $0<\varepsilon\leq\eta$. Observe that $2\varepsilon$ is also a 
Lebesgue number for ${\mathcal U}$. Choose any sequence $(E_n(\varepsilon))_{n=1}^\infty$
of maximal $(n,\varepsilon)$--separated sets. Since maximal $(n,\varepsilon)$--separated sets 
are $(n,\varepsilon)$--spanning sets, it follows from Lemma~\ref{zPQ}(a) that 
$$
z_n(T,\varphi,{\mathcal U})\leq Q_n(T,\varphi,\varepsilon)\leq \Sigma_n(E_n(\varepsilon)).
$$
Therefore,
\beq\label{htopsequnderpres}
\P(T,\varphi)
=\underline{p}(T,\varphi,{\mathcal U})
=\liminf_{n\to\infty}\frac{1}{n}\log z_n(T,\varphi,{\mathcal U})
\leq\liminf_{n\to\infty}\frac{1}{n}\log \Sigma_n(E_n(\varepsilon)).
\eeq

On the other hand, since $\diam({\mathcal U})\le \d$, it follows from Corollary~\ref{c120190930} that there exists $N\in\N$ such that 
$$
\mathrm{diam}({\mathcal U}^k)\leq\varepsilon
$$
for all integers $k\geq N$. It ensues from Lemma~\ref{zPQ}(b) that 
$$
\Sigma_n(E_n(\varepsilon))\leq P_n(T,\varphi,\varepsilon)\leq Z_n(T,\varphi,{\mathcal U}^k)
$$
for all $k\geq N$. Consequently, 
\[
\limsup_{n\to\infty}\frac{1}{n}\log \Sigma_n(E_n(\varepsilon))
\leq\lim_{n\to\infty}\frac{1}{n}\log Z_n(T,\varphi,{\mathcal U}^k)
=\P(T,\varphi,{\mathcal U}^k)
\]
for all $k\geq N$. It then follows from Proposition~\ref{ch9sem2lem9.5}(g) that
\beq\label{htopseqoverpres}
\limsup_{n\to\infty}\frac{1}{n}\log \Sigma_n(E_n(\varepsilon))
\leq\lim_{k\to\infty}\P(T,\varphi,{\mathcal U}^k)
=\P(T,\varphi).
\eeq 
Combining~(\ref{htopsequnderpres}) and~(\ref{htopseqoverpres}) gives (a).\epf

\sp The ultimate result of this section is the following.

\bthm\label{ch6sem1thm7.2bowennewpres}
Let $(X,d)$ be a compact metric space. If $T:X\to X$ is a $\delta$--expansive topological dynamical system on a compact metric space $(X,d)$,
then items (a)--(d) of Theorem~\ref{ch9sem2cor9.11bowenseqpre} for every $0<\e\le \delta/4$.
\ethm

\bpf
In view of Theorem~\ref{ch9sem2cor9.11bowenseqpre}, it suffices to show that there exists ${\mathcal U}$, an open cover of $X$, whose Lebesgue number is $\d/2$. So, since a Lebesgue number of the open cover  
$$
{\mathcal U}:=\big\{B(x,\d):x\in X\big\}
$$
is $\d/2$, we are done.
\epf

\section{Basic Properties of Topological Pressure}\label{ch9proppres}

In this section we give some of the most basic properties of topological pressure.
First, we show that the addition or subtraction of a constant to the potential increases or decreases
the pressure of the potential by that same constant.

\bprop\label{ch9prop1tp}
If $T:X\lra X$ is a topological dynamical system and $\varphi:X\lra\R$ is a continuous function, then for every constant $c\in\R$, we have that 
$$
\P(T,\varphi+c)=\P(T,\varphi)+c.
$$
\eprop

\bpf
For each $n\in\N$ and each $\varepsilon>0$,
let $E_n(\varepsilon)$ be a maximal $(n,\varepsilon)$--separated set. Then
\begin{eqnarray*}
\P(T,\varphi+c)
&=&\lim_{\varepsilon\to0}\limsup_{n\to\infty}\frac{1}{n}\log\sum_{x\in E_n(\varepsilon)}e^{S_n(\varphi+c)(x)} \\
&=&\lim_{\varepsilon\to0}\limsup_{n\to\infty}\frac{1}{n}\log\sum_{x\in E_n(\varepsilon)}e^{S_n\varphi(x)}e^{nc} \\
&=&\lim_{\varepsilon\to0}\limsup_{n\to\infty}\frac{1}{n}\Bigl[\log\Bigl(\sum_{x\in E_n(\varepsilon)}e^{S_n\varphi(x)}\Bigr)+nc\Bigr] \\
&=&\P(T,\varphi)+c.
\end{eqnarray*}
\epf

\sp Next, we shall show that topological pressure, as a function of the potential, is increasing.

\bprop\label{ch9pressentineq}
If $T:X\lra X$ is a topological dynamical system and $\varphi,\psi:X\lra\R$ are continuous functions, such that  $\varphi\leq\psi$, then 
$$
\P(T,\varphi)\leq \P(T,\psi).
$$
In particular,
\[
\htop(T)+\inf\varphi\leq \P(T,\varphi)\leq \htop(T)+\sup\varphi.
\]
\eprop

\bpf
That $\P(T,\varphi)\leq \P(T,\psi)$ whenever $\varphi\leq\psi$ is clear from the characterization of pressure given in Theorem~\ref{ch9sem2thm9.12}.
The second statement was proved in Corollary~\ref{ch9cor2.1.19}.
\epf

\sp In general, it is not the case that $\P(T,c\,\varphi)=c \P(T,\varphi)$. For example, suppose that
$\P(T,0)\neq0$. Then the equation $\P(T,c\,\cdot0)=c\,\P(T,0)$ only holds when $c=1$.

\section{Examples}\label{ch9examples}

\bex\label{ch9exmp1}
{\rm Let $E$ be a finite alphabet and let $\sg:E^\N\lra E^\N$ be the corresponding shift map; see formula \eqref{120191121} for its definition and Section~\ref{SOFToIA;TP} its treatement at length. Let $\widetilde{\varphi}:E\lra\R$ be an arbitrary function. Then the function $\varphi:E^\N\lra \R$ defined by 
$$
\varphi(\omega):=\widetilde{\varphi}(\omega_1)
$$
is a continuous function on $E^\infty$ which depends only upon the first coordinate $\om_1$ of the word $\omega\in E^\infty$. We will show that
\[
\P(\sigma,\varphi)=\log\sum_{e\in E}\exp(\widetilde{\varphi}(e)).
\]
It is immediate to see that the
shift map $\sg:E^\N\lra E^\N$ is expansive and any number $\delta\in(0,1)$ is an expansive constant when $E^\N$ is endowed with the metric given by formula \eqref{120191105}. Note that 
$$
{\mathcal U}=\big\{[e]:e\in E\big\}
$$ 
is a (finite) open cover of $E^\N$, and furthermore, it
 is the partition of $E^\N$ into initial
$1$--cylinders. Since $\mathrm{diam}({\mathcal U})=s<1$,
in light of Theorem~\ref{ch9sem2cor9.11}, we have that
$$
\P(\sigma,\varphi)=\P(\sigma,\varphi,{\mathcal U}).
$$
In order to compute $\P(\sigma,\varphi,{\mathcal U})$, observe that
${\mathcal U}^n=\{[\om]:\om\in E^n\}$ is the partition of $E^\infty$
into initial cylinders of length $n$. Then
\begin{eqnarray*}
\P(\sigma,\varphi)=\P(\sigma,\varphi,{\mathcal U})
&=&\lim_{n\to\infty}\frac{1}{n}\log Z_n(\varphi,{\mathcal U})
 = \lim_{n\to\infty}\frac{1}{n}\log\sum_{U\in{\mathcal U}^n}e^{\overline{S}_n\varphi(U)} \\
&=&\lim_{n\to\infty}\frac{1}{n}\log\sum_{\om\in E^n}e^{\overline{S}_n\varphi([\om])} \\
&=&\lim_{n\to\infty}\frac{1}{n}\log\sum_{\om_1\ldots\om_n\in E^n}
   \exp\bigl(\widetilde{\varphi}(\om_1)+\ldots+\widetilde{\varphi}(\om_n)\bigr) \\
&=&\lim_{n\to\infty}\frac{1}{n}\log
   \sum_{\om_1\in E}\exp(\widetilde{\varphi}(\om_1))
   \cdots\sum_{\om_n\in E}\exp(\widetilde{\varphi}(\om_n)) \\
&=&\lim_{n\to\infty}\frac{1}{n}\log\Bigl(\sum_{e\in E}\exp(\widetilde{\varphi}(e))\Bigr)^n \\
&=&\log\sum_{e\in E}\exp(\widetilde{\varphi}(e)).
\end{eqnarray*}}
\eex

\sp\bex\label{ch9exmp2}{
Let $E$ be a finite alphabet and let $\sg:E^\N\lra E^\N$ be the corresponding shift map; see, as in the previous example, formula \eqref{120191121} for its definition and Section~\ref{SOFToIA;TP} its treatement at length. Let $\widetilde{\varphi}:E^2\lra\R$ be an arbitrary function. Then the function
$\varphi:E^\N\lra\R$ defined by 
$$
\varphi(\omega):=\widetilde{\varphi}(\omega_1,\omega_2)
$$
is a continuous function on $E^\N$ which depends only upon the first two coordinates of the word $\omega\in E^\N$.

As in the previous example, 
$$
\P(\sigma,\varphi)=\P(\sigma,\varphi,{\mathcal U}),
$$
where ${\mathcal U}=\{[e]:e\in E\}$ is the (finite) open partition of $E^\infty$ into initial $1$--cylinders and
\[
\P(\sigma,\varphi)=\P(\sigma,\varphi,{\mathcal U})
=\lim_{n\to\infty}\frac{1}{n}\log\sum_{\om\in E^n}e^{\overline{S}_n\varphi([\om])}.
\]
But in this case
\begin{eqnarray*}
\sum_{\om\in E^n}e^{\overline{S}_n\varphi([\om])}\!\!\!
&=&\!\!\!\sum_{\om\in E^n}
   \exp\Big(\widetilde{\varphi}(\om_1,\om_2)+\widetilde{\varphi}(\om_2,\om_3)+\ldots +\widetilde{\varphi}(\om_{n-1},\om_n)+\max_{e\in E}\big\{\exp(\widetilde{\varphi}(\om_n,e)\big\}\Big) \\           
&=&\!\!\!\sum_{\om_1\in E}\sum_{\om_2\in E}e^{\widetilde{\varphi}(\om_1,\om_2))}
   \sum_{\om_3\in E}e^{\widetilde{\varphi}(\om_2,\om_3)}
	 \cdots\sum_{\om_n\in E}e^{\widetilde{\varphi}(\om_{n-1},\om_n)} 
\,\cdot\,\max_{e\in E}\big\{\exp(\widetilde{\varphi}(\om_n,e))\big\}.
%
\end{eqnarray*}
Since
\[
m:=\min_{e,f\in E}\big\{\exp(\widetilde{\varphi}(f,e))\big\}
\leq\max_{e\in E}\big\{\exp(\widetilde{\varphi}(\om_n,e))\big\}
\leq\max_{e,f\in E}\big\{\exp(\widetilde{\varphi}(f,e))\big\}=:M
\]
for all $n\in\N$ and all $\om_n\in E$, we have that
\[
\sum_{\om\in E^n}e^{\overline{S}_n\varphi([\om])}
\asymp\sum_{\om_1\in E}\sum_{\om_2\in E}e^{\widetilde{\varphi}(\om_1,\om_2)}
   \sum_{\om_3\in E}e^{\widetilde{\varphi}(\om_2,\om_3)}
	 \cdots\sum_{\om_n\in E}e^{\widetilde{\varphi}(\om_{n-1},\om_n)} \\
%
\]
for all $n$, with uniform constant of comparability $C:=\max\{m^{-1},M\}$.

Let $A:E^2\lra(0,+\infty)$ be the positive matrix whose entries are 
$$
A_{ef}:=\exp(\widetilde{\varphi}(e,f)).
$$
Equip this matrix with the $L^1$ norm 
$$
\|A\|:=\sum_{e\in E}\sum_{f\in E}A_{ef}.
$$
It is then easy to prove by induction that
\[
\|A^{n-1}\|
=\sum_{\om_1\in E}\sum_{\om_2\in E}e^{\widetilde{\varphi}(\om_1,\om_2)}
   \sum_{\om_3\in E}e^{\widetilde{\varphi}(\om_2,\om_3)}\cdots 
	 \sum_{\om_n\in E}e^{\widetilde{\varphi}(\om_{n-1},\om_n)} 
\]
for all $n\geq2$, and hence
\[
\sum_{\om\in E^n}e^{\overline{S}_n\varphi([\om])}\asymp\|A^{n-1}\|.
\]
Therefore
\begin{eqnarray*}
\P(T,\varphi)
&=&\lim_{n\to\infty}\frac{1}{n}\log\sum_{\om\in E^n}e^{\overline{S}_n\varphi([\om])} \\
&=&\lim_{n\to\infty}\frac{1}{n}\log\|A^{n-1}\|
=\log\Big(\lim_{n\to\infty}\|A^n\|^{1/n}\Big)
=\log r(A),
\end{eqnarray*}
where $r(A)$ is the spectral radius of $A$, that is, the largest eigenvalue of $A$ in modulus.}
\eex

\sp\section{The Variational Principle and Equilibrium States}\label{ch10B}

In this chapter, we will state and prove a fundamental result of thermodynamic formalism known as the
{\em Variational Principle}. \index{(N)}{Variational Principle} This deep result establishes a crucial relationship
between topological dynamics and ergodic theory, by way of a formula
linking topological pressure and measure theoretic entropy. The Variational Principle in its classical form and full generality was proved in \cite{W1} and \cite{Bow}.
The proof we shall present follows that of Michal Misiurewicz~\cite{Mis} which is particularly elegant, short and simple. We further introduce the concept of equilibrium states, give some sufficient conditions for their existence such as upper semicontinuity of metric entropy function and entailing it expansivness. We single out the special class of equilibrium states, those of potentials that are identically equal to zero, and following tradition, we call them measures of maximal entropy. We do not deal in this chapter with the issue of uniqueness of equilibrium states. We however provide an example of a topological dynamical system with positive and finite topological entropy which does not have any measure of maximal entropy. 

\sp\subsection{The Variational Principle}

For any topological dynamical system $T:X\lra X$, subject to a potential $\varphi:X\lra\R$ and equipped with a $T$--invariant Borel probability measure $\mu$, the quantity
$$
\hmu(T)+\int_X\varphi\,d\mu
$$ 
is called the \index{(N)}{free energy} free energy of the system
$T$ with respect to $\mu$. The variational principle 
states that the topological pressure of a system is the supremum of
the free energies generated by that system.

\sp\bthm[Variational Principle]\label{sem2thm9.13}
If $T:X\lra X$ is a continuous map of a compact metrizable space 
$X$ and $\varphi:X\lra\R$ is a continuous function, then
\[
\P(T,\varphi)=\sup\Bigl\{\hmu(T)+\int_X\varphi\,d\mu:\mu\in M(T)\Bigr\},
\]
where $M(T)=M(T,\cB)$ is the set of all $T$--invariant Borel probability measures on $X$.
\ethm

\sp\fr The proof will be given in two parts. In the first part, we will show that
$$
\P(T,\varphi)\geq \H_\mu(T)+\int\varphi\,d\mu
$$ 
for every measure $\mu\in M(T)$. The second part  shall consist of the proof of the inequality 
$$
\sup\Big\{\H_\mu(T)+\int_X\varphi\,d\mu:\mu\in M(T)\Big\}\geq \P(T,\varphi).
$$
The first part is relatively easier to prove than the second.
For the proof of Part I, we will need Jensen's Inequality.
Recall that a function $\psi:(a,b)\to\R$, where $-\infty\leq a<b\leq\infty$,
is said to be {\em convex} if 
$$
\psi\bigl(tx+(1-t)y\bigr)\leq t\psi(x)+(1-t)\psi(y)
$$ 
for all $t\in[0,1]$ and all $x,y\in(a,b)$.

\sp\bthm[Jensen's Inequality]\label{chapter11thm2}
Let $\mu$ be a probability measure on a measurable space $(X,{\mathfrak F})$
and let $\psi:(a,b)\to\R$ be a convex function. If $f\in
L^1(X,{\mathfrak F},\mu)$ and $f(X)\subseteq(a,b)$, then
\[
\psi\Bigl(\int_X f\,d\mu\Bigr)\leq\int_X\psi\circ f\,d\mu.
\]
\ethm

\sp We shall also need the following lemma, which states that any
finite Borel partition 
${\mathcal A}$ of $X$ can be, from a metric entropy viewpoint,
approximated as closely as desired 
by a finite Borel partition ${\mathcal B}$ whose elements are all (but one)
compact and contained in those of ${\mathcal A}$.

\sp\blem\label{chapter11lem1}
Let $\mu\in M(X)$, let ${\mathcal A}:=\{A_1, \ldots, A_s\}$ be a
finite partition of 
$X$ into Borel sets and let $\varepsilon>0$. Then there exist compact sets
$B_i\subseteq A_i$, 
$1\leq i\leq s$, such that the partition
${\mathcal B}:=\{B_1, \ldots, B_s, X\backslash(B_1\cup\ldots\cup
B_s)\}$ satisfies 
\[
\H_\mu({\mathcal A}|{\mathcal B})\leq \varepsilon>0.
\]
\elem

\bpf
Let the measure $\mu$ and the partition ${\mathcal A}$ be as stated and
let $\varepsilon>0$. Recall from Definition~\ref{functionk} the continuous
non-negative function 
$k:[0,1]\to[0,1]$ defined by
\[
k(t):=
\begin{cases}
          -t\log t & \text{ if } t\in(0,1] \\
          0 &  \text{ if } t=0.
\end{cases}
\] 
The continuity of $k$ at $0$ implies that there
exists $\delta>0$ such that $k(t)<\varepsilon/s$ when $0\leq t<\delta$. Since
$\mu$ is regular and $X$ 
is compact, for each $1\leq i\leq s$ there exists a compact set
$B_i\subseteq A_i$ such that 
$\mu(A_i\backslash B_i)<\delta$. Thus, $k(\mu(A_i\backslash
B_i))<\varepsilon/s$ for all 
$1\leq i\leq s$. By Definition~\ref{condentpartition} of conditional entropy,
it follows that
\begin{eqnarray*}
\H_\mu({\mathcal A}|{\mathcal B})
&=&\sum_{j=1}^s\sum_{i=1}^s-\mu(A_i\cap B_j)\log\frac{\mu(A_i\cap
  B_j)}{\mu(B_j)} \\ 
&&+\sum_{i=1}^s-\mu\bigl(A_i\cap (X\backslash\cup_{j=1}^s B_j)\bigr)
\log\frac{\mu\bigl(A_i\cap (X\backslash\cup_{j=1}^s
  B_j)\bigr)}{\mu(X\backslash\cup_{j=1}^s B_j)} \\ 
&=&\sum_{j=1}^s-\mu(B_j)\log\frac{\mu(B_j)}{\mu(B_j)} \\
&&+\sum_{i=1}^s-\mu\bigl(A_i\cap(\cup_{j=1}^s A_j\backslash B_j)\bigr)
\log\frac{\mu\bigl(A_i\cap(\cup_{j=1}^s A_j\backslash
  B_j)\bigr)}{\mu\bigl(\cup_{j=1}^s A_j\backslash B_j\bigr)} \\ 
&=&0+
\sum_{i=1}^s-\mu(A_i\backslash B_i)
\log\frac{\mu(A_i\backslash B_i)}{\mu\bigl(\cup_{j=1}^s A_j\backslash
  B_j\bigr)} \\ 
&=&\sum_{i=1}^s-\mu(A_i\backslash B_i)
\Bigl[\log\mu(A_i\backslash B_i)-\log\mu\bigl(\cup_{j=1}^s
A_j\backslash B_j\bigr)\Bigr] \\ 
&=&\sum_{i=1}^s k(\mu(A_i\backslash B_i))+\sum_{i=1}^s
\mu(A_i\backslash B_i)\log\mu\bigl(\cup_{j=1}^s A_j\backslash
B_j\bigr) \\ 
&\leq&\sum_{i=1}^s k(\mu(A_i\backslash B_i))\leq s\cdot\frac{\varepsilon}{s}=\varepsilon.
\end{eqnarray*}
\epf

\sp\bpf[Proof of Part I]
Recall that our aim is to establish  the inequality
\beq\label{chapter11eq6C=0}
\P(T,\varphi)\geq \hmu(T)+\int_X\varphi\,d\mu,\ \text{ for all } \ \mu\in M(T).
\eeq
To that end, let $\mu\in M(T)$ be arbitrary. We claim that it is
sufficient to prove that  
\beq\label{chapter11eq6C}
\P(T,\varphi)\geq \hmu(T)+\int_X\varphi\,d\mu+\log 2.
\eeq
 Indeed, suppose that (\ref{chapter11eq6C}) holds. Then, rather than
 considering directly the system $(X,T)$ under the potential 
$\varphi$, we may alternatively consider the higher--iterate system
$(X,T^n)$ under the potential 
$$
S_n\varphi=\sum_{k=0}^{n-1}\varphi\circ T^k.
$$
Since the measure $\mu$ is $T$--invariant, 
it is also $T^n$--invariant for every integer $n\geq1$. Using successively
Theorem~\ref{ch9sem2prop9.7}, 
inequality~(\ref{chapter11eq6C}) with the quadruple
$(X,T^n,S_n\varphi,\mu)$ instead of 
$(X,T,\varphi,\mu)$, and Theorems~\ref{sem2thm7.4},
we would then obtain that 
\begin{eqnarray*}
n\P(T,\varphi)
=\P(T^n,S_n\varphi)
\geq
\hmu(T^n)+\hspace{-0.1cm}\int_X\hspace{-0.1cm}
S_n\varphi\,d\mu+\hspace{-0.1cm}\log 2   
=n\H_\mu(T)+n\hspace{-0.1cm}\int_X\hspace{-0.1cm}\varphi\,d\mu+\hspace{-0.1cm}\log 2.
\end{eqnarray*}
Dividing by $n$ and letting $n$ tend to infinity would then
yield~(\ref{chapter11eq6C=0}). 
 Of course, to obtain~(\ref{chapter11eq6C}) it suffices to show that
\beq\label{chapter11eq6CA}
\P(T,\varphi)\geq \hmu(T,{\mathcal A})+\int_X\varphi\,d\mu+\log 2
\eeq
for all finite Borel partitions ${\mathcal A}$ of $X$ (see
Definition~\ref{mte}). So let ${\mathcal A}$ 
be such a partition and let $\varepsilon>0$. To get~(\ref{chapter11eq6CA}),
it is enough to prove that 
\beq\label{chapter11eq6CAe}
\P(T,\varphi)\geq \hmu(T,{\mathcal A})+\int_X\varphi\,d\mu+\log 2-2\e.
\eeq
Because of Corollary~\ref{ch9sem2cor9.12}, it suffices to demonstrate that
\beq\label{chapter11eq5CAedB}
\underline{\P}(T,\varphi,\delta)
\geq \hmu(T,{\mathcal A})+\int_X\varphi\,d\mu+\log 2-2\e
\eeq
for all sufficiently small $\delta>0$. We will in fact show more, namely that
\beq\label{chapter11eq5CAed}
\ov{\P}(T,\varphi,\delta)
\geq \hmu(T,{\mathcal A})+\int_X\varphi\,d\mu+\log 2-2\e
\eeq
for all sufficiently small $\delta>0$. In light of Definition~\ref{ch6sem1defn7.5gen} and in view of Definition~\ref{entwrtpartition} 
of the relative entropy $\hmu(T,{\mathcal A})$,
it is sufficient to prove that
\beq\label{chapter11eq5CAedn}
\frac1n\log\sum_{y\in F_n(\delta)}\exp(S_n\varphi(y))
\geq\frac1n \H_\mu({\mathcal A}^n)+\frac1n\int_X S_n\varphi\,d\mu+\log 2-2\e
\eeq
for all sufficiently small $\delta>0$, all large enough $n\in\N$ and
every $(n,\delta)$--separated 
set $F_n(\delta)$. So given a finite partition ${\mathcal A}:=\{A_1, \ldots, A_s\}$ and $\varepsilon>0$, let 
$$
{\mathcal B}:=\{B_1, \ldots, B_s, X\backslash(B_1\cup\ldots\cup B_s)\}
$$ 
be the partition given by Lemma~\ref{chapter11lem1}. 
Thus, $\H_\mu({\mathcal A}|{\mathcal B})\leq \varepsilon$. Fix
$n\geq1$. By Theorem~\ref{sem2thm7.1}(9) 
and Lemma~\ref{lempartH}(3), we know that
\beq\label{chapter11eq1a}
\H_\mu({\mathcal A}^n)\hspace{-0.1cm}
\leq \H_\mu({\mathcal B}^n)+\H_{\mu}({\mathcal A}^n|{\mathcal B}^n)
\leq \H_\mu({\mathcal B}^n)+n\H_{\mu}({\mathcal A}|{\mathcal B})
\leq \H_\mu({\mathcal B}^n)+n\varepsilon.
\eeq
From~(\ref{chapter11eq5CAedn}) and~(\ref{chapter11eq1a}), it thus
suffices to establish that 
\beq\label{chapter11eq5CBedn}
\log\sum_{y\in F_n(\delta)}\exp(S_n\varphi(y))
\geq \H_\mu({\mathcal B}^n)+\int_X S_n\varphi\,d\mu+(\log 2-\e)n
\eeq
for all sufficiently small $\delta>0$, all large enough $n\in\N$ and
every $(n,\delta)$--separated 
set $F_n(\delta)$. In order to prove this inequality, we will estimate
the term 
$$
\H_\mu({\mathcal B}^n)+\int S_n\varphi\,d\mu
$$
from above.
Since the logarithm function is concave (so its negative is convex),
Jensen's Inequality (Theorem~\ref{chapter11thm2}) implies that
\beq\label{chapter11eq2a}
\aligned
\H_\mu({\mathcal B}^n)+\int_X S_n\varphi\,d\mu
&\leq\sum_{B\in {\mathcal B}^n}\mu(B)\Bigl[-\log\mu(B)+S_n\varphi(B)\Bigr] 
=\sum_{B\in{\mathcal B}^n}\mu(B)\log\frac{\exp(S_n\varphi(B))}{\mu(B)} \\
&=\int_X\log\frac{\exp\bigl(S_n\varphi({\mathcal B}^n(x))\bigr)}{\mu({\mathcal
    B}^n(x))}\,d\mu(x) \\ 
&\leq\log\int_X\frac{\exp\bigl(S_n\varphi({\mathcal
    B}^n(x))\bigr)}{\mu({\mathcal B}^n(x))}\,d\mu(x) \\ 
&=\log\sum_{B\in {\mathcal B}^n}\exp(S_n\varphi(B)).
\endaligned
\eeq
Since each set $B_i\in{\mathcal B}$ is compact, it follows that
$d(B_i,B_j)>0$ for all 
$1\leq i\neq j\leq s$. As $\varphi$ is continuous, let
$0<\delta<\frac12\min\{d(B_i,B_j):1\leq i\neq j\leq s\}$ be such that
\beq\label{chapter11eq3a}
d(x,y)<\delta\hspace{0.5cm}\Longrightarrow\hspace{0.5cm}|\varphi(x)-\varphi(y)|
      <\varepsilon.  
\eeq
Now consider an arbitrary maximal $(n,\delta)$-separated set
$F_n(\delta)$ and fix temporarily $B\in{\mathcal B}$. 
According to Lemma~\ref{ch6sepspan}, each maximal
$(n,\delta)$-separated set is an $(n,\delta)$-spanning set. 
So for every $x\in B$ there exists $y\in F_n(\delta)$ such that $x\in
B_n(y,\delta)$ and therefore 
$|S_n\varphi(x)-S_n\varphi(y)|<n\varepsilon$ by~(\ref{chapter11eq3a}). As the
set $F_n(\delta)$ is finite, 
it follows that there exists $y_B\in F_n(\delta)$ such that
\beq\label{chapter11eq4a}
S_n\varphi(B)
\leq S_n\varphi(y_B)+n\varepsilon\hspace{0.5cm}\mbox{ and
    }\hspace{0.5cm}B\cap B_n(y_B,\delta) 
\neq\emptyset.
\eeq
Moreover, since $d(B_i,B_j)>2\delta$ for each $1\leq i\neq j\leq s$,
any dynamic ball $B_k(z,\delta)$, 
$k\geq1$, $z\in X$, intersects at most one $B_i$ and perhaps
$X\backslash\cup_{j=1}^s B_j$. Hence 
\[
\#\bigl\{B\in{\mathcal B}:B\cap B_k(z,\delta)\neq\emptyset\bigr\}\leq2
\]
for all $k\geq1$ and $z\in X$. Thus,
\[
\#\bigl\{B\in{\mathcal B}^n:B\cap B_n(z,\delta)\neq\emptyset\bigr\}\leq2^n
\]
for all $z\in X$. So, the function $f:{\mathcal B}^n\lra F_n(\delta)$
defined by $f(B)=y_B$ 
is at most $2^n$--to--one. Consequently, by~(\ref{chapter11eq4a}) we obtain that
\begin{eqnarray*}
\begin{aligned}
2^n\sum_{y\in F_n(\delta)}\exp(S_n\varphi(y))
&\geq\sum_{B\in{\mathcal B}^n}\exp(S_n\varphi(y_B))
\geq\sum_{B\in{\mathcal B}^n}\exp(S_n\varphi(B)-n\varepsilon) \\
&=\sum_{B\in{\mathcal B}^n}\exp(S_n\varphi(B))\cdot e^{-n\varepsilon}.
\end{aligned}
\end{eqnarray*}
Multiplying both sides by $2^{-n}$, then taking the logarithm of both sides and
applying~(\ref{chapter11eq2a}) yields
\begin{eqnarray*}
\log\sum_{y\in F_n(\delta)}\exp(S_n\varphi(y))
&\geq&\log\sum_{B\in{\mathcal B}^n}\exp(S_n\varphi(B))-n\varepsilon-n\log2 \\
&\geq&\H_\mu({\mathcal B}^n)+\int_X S_n\varphi\,d\mu+n(-\log2-\varepsilon).
\end{eqnarray*}
This inequality, which is nothing other than the sought
inequality~(\ref{chapter11eq5CBedn}) holds for all
$$
0<\delta<\frac12\min\big\{d(B_i,B_j):1\leq i\neq j\leq s\big\},
$$
all $n\in\N$ and all maximal $(n,\delta)$--separated sets $F_n(\delta)$. This concludes the proof.
\epf

\sp  Let us now move on to the proof of Part II of the variational principle.
For this, we shall need the following four lemmas.
The first lemma states that given any Borel probability measure $\mu$
there exist 
finite Borel partitions of arbitrarily small diameters whose atoms
have boundaries with zero $\mu$--measure.

\sp\blem\label{chapter11lem3}
Let $\mu\in M(X)$. For every $\varepsilon>0$,
there exists a finite Borel partition ${\mathcal A}$ of $X$ such that
$$
\mathrm{diam}({\mathcal A})<\varepsilon
\ \ \ {\rm and}  \  \  \
\mu(\partial A)=0
$$
for each $A\in\mathcal {A}$. 
\elem

\bpf
Let $\varepsilon>0$ and let $E:=\{x_1, \ldots, x_s\}$ be an
$(\varepsilon/4)$-spanning set of $X$. 
Since for each fixed $1\leq i\leq s$ the sets 
$$
\{x\in X:d(x,x_i)=r\},
$$
where $\varepsilon/4<r<\varepsilon/2$, are mutually disjoint, only countably many of them may have
positive $\mu$--measure. Hence, there exists $\varepsilon/4<t<\varepsilon/2$ such that for every
$1\leq i\leq s$ we have
\begin{eqnarray}\label{ch112.4.5}
\mu\bigl(\{x\in X:d(x,x_i)=t\}\bigr)=0.
\end{eqnarray}
Now let the set $A_1$ be defined by 
$$
A_1:=\{x\in X:d(x,x_1)\leq t\},
$$
and for each $2\leq i\leq s$ define the sets $A_2, \ldots, A_s$ inductively by setting
\[
A_i:=\{x\in X:d(x,x_i)\leq t\}\backslash(A_1\cup\cdots\cup A_{i-1}).
\]
Since $t<\varepsilon/2$, the family ${\mathcal A}:=\{A_1, \ldots, A_s\}$ is a
Borel partition 
of $X$ with diameter smaller than $\varepsilon$. Noting that
$\partial(A\backslash B)\subseteq\partial A\cup\partial B$ and
$\partial(A\cup B)\subseteq\partial A\cup\partial B$, 
it follows from~(\ref{ch112.4.5}) that 
$$
\mu(\partial A_i)=0
$$
for each $1\leq i\leq s$. 
\epf

\sp The second lemma states that given any finite Borel partition
${\mathcal A}$ whose atoms have boundaries with 
zero $\mu$-measure, the entropy of ${\mathcal A}$ as a function of the
underlying Borel probability measure is continuous at $\mu$.

\sp\blem\label{chapter11lem4}
Let $\mu\in M(X)$. If ${\mathcal A}$ is a finite Borel partition of $X$ such that $\mu(\partial A)=0$ for every $A\in{\mathcal A}$, then the function
$\H_{(\cdot)}({\mathcal A}):M(X)\lra [0, \infty]$ defined by 
$$
\H_{(\cdot)}({\mathcal A})(\nu)=\H_\nu(\mathcal{A})
$$ 
is continuous at $\mu$.
\elem

\bpf
This follows directly from the fact that, according to the Portmanteau
Theorem, a sequence of Borel probability measures $(\mu_n)_{n\geq1}$ converges
weakly to a measure $\mu$ if and only if 
$\lim_{n\to\infty}\mu_n(A)=\mu(A)$ for every Borel set $A$ with
$\mu(\partial A)=0$. 
\epf

\sp In the third announced lemma, we show that the entropy of ${\mathcal A}$ as a function of the underlying Borel probability measure is 
a concave function.

\sp \blem\label{chapter11lem45}
For any finite Borel partition ${\mathcal A}$ of $X$, the function
$H_{(\cdot)}({\mathcal A})$ is concave. 
\elem

\bpf
Let ${\mathcal A}$ be a finite Borel partition of $X$, and $\mu$ and $\nu$
Borel probability measures on $X$. 
Let also $t\in(0,1)$. Since the function $k(x)=-x\log x$ is concave,
for any $A\in{\mathcal A}$ we have 
\begin{eqnarray*}
k\bigl(t\mu(A)+(1-t)\nu(A)\bigr)\geq tk(\mu(A))+(1-t)k(\nu(A)),
\end{eqnarray*}
and hence
\begin{eqnarray*}
\H_{t\mu+(1-t)\nu}({\mathcal A})
&=&\sum_{A\in{\mathcal A}}k\bigl(t\mu(A)+(1-t)\nu(A)\bigr) 
\geq t\sum_{A\in{\mathcal A}}k(\mu(A))+(1-t)\sum_{A\in{\mathcal A}}k(\nu(A)) \\
&=&t\H_\mu({\mathcal A})+(1-t)\H_\nu({\mathcal A}).
\end{eqnarray*}
\epf

\sp Finally, the fourth lemma is a generalization of, already proven, Krylov--Bogolyubov
Theorem asserting that any continuous self-map of a compact metrizable
spaces has at least one Borel probability invariant measure. 

\sp\blem\label{chapter11lem5}
Let $T:X\lra X$ be a topological dynamical system. If $(\mu_n)_{n\geq1}$ is a sequence of measures in $M(X)$, then any weak* limit point of the
sequence $(m_n)_{n\geq1}$, where  
$$
m_n:=\frac1n\sum_{i=0}^{n-1}\mu_n\circ T^{-i}, 
$$
is a $T$--invariant measure.
\elem

\bpf
By the compactness of $M(X)$ in the weak$^*$ topology, we know that the sequence
$(m_n)_{n\geq1}$ has accumulation points. 
So let $(m_{n_j})_{j\geq1}$ be a subsequence which converges
weakly$^*$ to, say, $m\in M(X)$. Let $f\in C(X)$ be arbitrary. We obtain:
\begin{eqnarray*}
\left|\int_X f\circ T\,dm-\int_X f\,dm\right|
&=&\lim_{j\to\infty}\left|\int_X f\circ T\,dm_{n_j}-\int_X f\,dm_{n_j}\right| \\
&=&\lim_{j\to\infty}\left|\frac{1}{n_j}\int_X\sum_{i=0}^{n_j-1}(f\circ
  T^{i+1}-f\circ T^{i})\,d\mu_{n_j}\right| \\ 
&=&\lim_{j\to\infty}\left|\frac{1}{n_j}\int_X(f\circ
  T^{n_j}-f)\,d\mu_{n_j}\right| \\ 
&\leq&\lim_{j\to\infty}\frac{2\|f\|_\infty}{n_j}=0.
\end{eqnarray*}
Thus, 
$$
\int_X f\circ T\,dm=\int_X f\,dm,
$$
and so, the Borel probability measure $m$ is $T$--invariant.
\epf

\sp \bpf[Proof of Part II]
Fix $\varepsilon>0$. Let $(F_n(\varepsilon))_{n\geq1}$ be a sequence of maximal $(n,\varepsilon)$--separated sets 
in $X$. For every $n\geq1$, define the measures $\mu_n$ and $m_n$ by
\[
\mu_n=\frac{\sum_{x\in F_n(\varepsilon)}\delta_x
  \exp(S_n\varphi(x))}{\sum_{x\in F_n(\varepsilon)}\exp(S_n\varphi(x))} 
\ \text{ and }\
m_n=\frac1n\sum_{k=0}^{n-1}\mu_n\circ T^{-k},
\]
where $\delta_x$ denotes the Dirac measure concentrated at the point
$x$. For ease of exposition, 
define 
$$
s_n:=\sum_{x\in F_n(\varepsilon)}\exp(S_n\varphi(x))
$$ and
$$
\mu(x):=\mu(\{x\}).
$$
Let $(n_i)_{i\geq1}$ be an increasing 
sequence of natural numbers such that $(m_{n_i})_{i\geq1}$ converges
weakly to, say, $m$, 
and such that
\begin{eqnarray}
\  \  \  \  \  \  \lim_{i\to\infty}\frac1{n_i}\log\hspace{-0.25cm}\sum_{x\in
  F_{n_i}(\varepsilon)}\hspace{-0.25cm}\exp(S_n\varphi(x)) 
=\limsup_{n\to\infty}\hspace{-0.2cm}\sum_{x\in
  F_n(\varepsilon)}\hspace{-0.2cm}\exp(S_n\varphi(x))
\end{eqnarray}
By Lemma~\ref{chapter11lem5}, the limit measure $m$ belongs to
$M(T)$. Also, in view of 
Lemma~\ref{chapter11lem3}, there exists a finite Borel partition
${\mathcal A}$ such that 
$$
\mathrm{diam}({\mathcal A})<\varepsilon
\ \ \ {\rm and} \  \  \ 
m(\partial A)=0
$$ 
for all $A\in{\mathcal A}$. Since 
$\#(A\cap F_n(\varepsilon))\leq1$ for all $A\in{\mathcal A}^n$, we obtain that
\beq\label{chapter11eq7}
\aligned
\H_{\mu_n}({\mathcal A}^n)+\int_X S_n\varphi\,d\mu_n
&=\sum_{x\in F_n(\varepsilon)}\mu_n(x)\Bigl[-\log \mu_n(x)+S_n\varphi(x)\Bigr] \\
&=\sum_{x\in F_n(\varepsilon)}\frac{\exp(S_n\varphi(x))}{s_n}
  \Bigl[-\log\frac{\exp(S_n\varphi(x))}{s_n}+S_n\varphi(x)\Bigr] \\
&=\frac{1}{s_n}\sum_{x\in
  F_n(\varepsilon)}\hspace{-0.2cm}\exp(S_n\varphi(x))\Bigl[-S_n\varphi(x)+\log
s_n+S_n\varphi(x)\Bigr] \\ 
&=\log s_n=\log\sum_{x\in F_n(\varepsilon)}\exp(S_n\varphi(x)).
\endaligned
\eeq
Now, fix $M\in\N$ and let $n\geq 2M$. For $j=0, 1, \ldots, M-1$, define
$s(j):=\bigl\lfloor\frac{n-j}{M}\bigr\rfloor-1$, where $\lfloor
r\rfloor$ denotes the integer part of $r$. Note that
\begin{eqnarray*}
\bigvee_{k=0}^{s(j)}T^{-(kM+j)}({\mathcal A}^M)
&=&T^{-j}({\mathcal A}^M)\vee T^{-(M+j)}({\mathcal A}^M)\vee\cdots\vee
T^{-(s(j)M+j)}({\mathcal A}^M) \\ 
&=&T^{-j}({\mathcal A})\vee T^{-(j+1)}({\mathcal A})\vee\cdots \vee
T^{-(s(j)M+j+M-1)}({\mathcal A}) \\ 
&=&T^{-j}({\mathcal A})\vee\ldots \vee T^{-((s(j)+1)M+j-1)}({\mathcal A})
\end{eqnarray*}
and
\[
(s(j)+1)M+j-1=\Bigl\lfloor\frac{n-j}{M}\Bigr\rfloor M+j-1\leq n-j+j-1=n-1.
\]
Observe also that 
$$
(n-1)-\((s(j)+1)M+j\)
\leq n-1-\lt(\bigl\lfloor\frac{n-j}{M}\bigr\rfloor M+j\rt)
\leq n-j-\lt(\frac{n-j}{M}-1\rt)M-1=M-1.
$$
Setting 
$$
R_j:=\{0, 1, \ldots, j-1\}\cup\{(s(j)+1)M+j, \ldots, n-1\},
$$
we have $\# R_j\leq 2M$ and
\[
{\mathcal A}^n=\bigvee_{k=0}^{s(j)}T^{-(kM+j)}({\mathcal A}^M)\vee
\bigvee_{i\in R_j}T^{-i}({\mathcal A}). 
\]
Hence, using Theorem~\ref{sem2thm7.1}(7), we get that
\begin{eqnarray*}
\H_{\mu_n}({\mathcal A}^n)
&\leq&\sum_{k=0}^{s(j)}\H_{\mu_n}\bigl(T^{-(kM+j)}({\mathcal A}^M)\bigr)
      +\H_{\mu_n}\Bigl(\bigvee_{i\in R_j}T^{-i}({\mathcal A})\Bigr) \\
&\leq&\sum_{k=0}^{s(j)}\H_{\mu_n\circ T^{-(kM+j)}}({\mathcal A}^M)
      +\log\#\Bigl(\bigvee_{i\in R_j}T^{-i}({\mathcal A})\Bigr) \\
&\leq&\sum_{k=0}^{s(j)}\H_{\mu_n\circ T^{-(kM+j)}}({\mathcal
  A}^M)+\log(\#{\mathcal A})^{\# R_j} \\ 
&\leq&\sum_{k=0}^{s(j)}\H_{\mu_n\circ T^{-(kM+j)}}({\mathcal A}^M)+2M\log\#{\mathcal A}.
\end{eqnarray*}
Summing over all $j=0, 1, \ldots, M-1$ and using
Lemma~\ref{chapter11lem45}, we obtain 
\begin{eqnarray*}
M\H_{\mu_n}({\mathcal A}^n)
&\leq&\sum_{j=0}^{M-1}\sum_{k=0}^{s(j)}\H_{\mu_n\circ T^{-(kM+j)}}({\mathcal A}^M)
      +2M^2\log\#{\mathcal A} \\
&\leq&\sum_{l=0}^{n-1}\H_{\mu_n\circ T^{-l}}({\mathcal A}^M)+2M^2\log\#{\mathcal A} \\
&\leq&n\H_{\frac1n\sum_{l=0}^{n-1}\mu_n\circ T^{-l}}({\mathcal A}^M)+2M^2\log\#{\mathcal A} \\
&=&n\H_{m_n}({\mathcal A}^M)+2M^2\log\#{\mathcal A}.
\end{eqnarray*}
Adding $M\int S_n\varphi\,d\mu_n$ to both sides and
applying~(\ref{chapter11eq7}) yields 
\[
M\log\sum_{x\in F_n(\varepsilon)}\exp(S_n\varphi(x))
\leq n\H_{m_n}({\mathcal A}^M)+M\int_X S_n\varphi\,d\mu_n+2M^2\log\#{\mathcal A}.
\]
Dividing both sides by $Mn$ and since $\frac1n\int
S_n\varphi\,d\mu_n=\int\varphi\,dm_n$, it follows that
\[
\frac1n\log\sum_{x\in F_n(\varepsilon)}\exp(S_n\varphi(x))
\leq \frac1M\H_{m_n}({\mathcal
  A}^M)+\int_X\varphi\,dm_n+\frac{2M}{n}\log\#{\mathcal A}. 
\]
Since $\partial T^{-1}(A)\subseteq T^{-1}(\partial A)$ for every set
$A\subseteq X$, 
the $m$--measure of the boundaries of the partition ${\mathcal A}^M$ is
equal to zero. 
Therefore, remembering that each maximal $(n,\varepsilon)$--separated set is $(n,\varepsilon)$--spanning, on letting $n$ tend to infinity along the subsequence $(n_i)_{i\geq1}$, 
we conclude from the above inequality and from Lemma~\ref{chapter11lem4} that
\[
\varliminf_{n\to\infty}\frac1n\log\sum_{x\in F_n(\varepsilon)}\exp(S_n\varphi(x))\leq\frac1M\H_m({\mathcal A}^M)+\int_X\varphi\,dm.
\]
Letting $M$ tend to infinity, we obtain that
\[
\varliminf_{n\to\infty}\frac1n\log\sum_{x\in F_n(\varepsilon)}\exp(S_n\varphi(x))
\leq \h_m(T,{\mathcal A})+\int_X\varphi\,dm
\leq\sup\Bigl\{\hmu(T)+\int_X\varphi\,d\mu:\mu\in M(T)\Bigr\}.
\]
Finally, letting $\varepsilon>0$ tend to zero and applying
Theorem~\ref{ch9sem2thm9.12} yields 
the desired inequality. This finishes the proof of Part II and
consequently completes 
the proof of the variational principle.
\epf

\sp Let us now state some immediate consequences of the variational principle. The first consequence concerns topological entropy of a topological dynamical system. Namely: the topological 
entropy of  such a system is the supremum of its all measure--theoretic entropies. 

\sp\bcor\label{sem2cor9.140}
If $T:X\lra X$ is a continuous map of a compact metrizable space, then 
$$
\h_{\mathrm{top}}(T)=\sup\{\hmu(T):\mu\in M(T)\}.
$$
\ecor

\bpf
This follows directly from Theorem~\ref{sem2thm9.13} on letting $\varphi\equiv0$.
\epf

\sp Furthermore, the topological pressure of  a topological dynamical system is determined by the supremum of the free energy 
of the system with respect to its ergodic measures.

\sp\bcor\label{sem2cor9.14}
If $T:X\lra X$ is a continuous map of a compact metrizable space, then 
$$
\P(T,\varphi)=\sup\Big\{\hmu(T)+\int_X\varphi\,d\mu:\mu\in E(T)\Big\},
$$
where $E(T)$ denotes the subset of all ergodic measures in $M(T)$.
\ecor

\bpf
Let $\mu\in M(T)$. Since $M(T)$ is a compact convex metrizable space,
the Chocquet representation theorem 
(see \cite{Phelps}) implies that $\mu$ has a decomposition in terms of
the extreme points of $M(T)$, which, according to Theorem~\ref{ext}, are
exactly the ergodic measures $E(T)$ in M(T). This means that there exists a unique probability measure $m$ on the Borel $\sg$--algebra of $M(T)$ such that $m(E(T))=1$ and
\[
\int_X f\,d\mu = \int_{E(T)}\Bigl(\int_X f\,d\nu\Bigr)dm(\nu)
\]
for all $f\in C(X)$. In particular ,
$$
\int_X\varphi\,d\mu=\int_{E(T)}\Big(\int_X\varphi\,d\nu\Big)\,dm(\nu).
$$ 
Moreover, we have that
\[\h_\mu(T)=\h_{\int_{E(T)}\nu\,dm(\nu)}(T)=\int_{E(T)}\h_{\nu}(T)\,dm(\nu).\]
It therefore follows that 
\[
\h_\mu(T)+\int_X\varphi\,d\mu
=\int_{E(T)}\Bigl[\h_\nu(T)+\int_X\varphi\,d\nu\Bigr]\,dm(\nu).
\]
Suppose by a way of contradiction that
$$
\h_\nu(T)+\int_X\varphi\,d\nu<\h_\mu(T)+\int_X\varphi\,d\mu
$$ 
for every $\nu\in E(T)$. Then we would have that
\begin{eqnarray*}
\hmu(T)+\int_X\varphi\,d\mu
&=&\int_{E(T)}\Bigl[\h_\nu(T)+\int_X\varphi\,d\nu\Bigr]\,dm(\nu) 
<\int_{E(T)}\Bigl[\h_{\mu}(T)+\int_X\varphi\,d\mu\Bigr]\,dm(\nu) \\
&=&\h_{\mu}(T)+\int_X\varphi\,d\mu.
\end{eqnarray*}
This is a contradiction and consequently there exists $\nu\in E(T)$ such that
$$
\h_\nu(T)+\int_X\varphi\,d\nu\geq \h_\mu(T)+\int_X\varphi\,d\mu.
$$
So, an application of the Variational Principle, i.e. Theorem~\ref{sem2thm9.13}, finishes the proof.
\epf

\sp Finally, we will show that the pressure of any subsystem is at most the
pressure of the entire system. 

\sp\bcor\label{sem2cor9.15}
If $T:X\lra X$ is a topological dynamical system, $\varphi:X\lra\R$ is a
continuous potential, and $Y$ is a closed subset of $X$ such that $T(Y)\subseteq Y$, then 
$$
\P\(T|_Y,\varphi|_Y\)\leq \P(T,\varphi).
$$
\ecor

\bpf
Each Borel probability $T|_Y$--invariant measure $\nu$ on $Y$ can be extended to the Borel probability $T$--invariant measure on $X$ by the following formula:
$$
{\nu^*}(B)=\nu(B\cap Y).
$$ 
Then
$$
\h_{\nu^*}(T)=\h_{\nu}(T)
\  \  \  {\rm and} \  \  \
\int_X\varphi\,d{\nu^*}=\int_Y\varphi\,d\nu. 
$$
Therefore, by virtue of the Variational Principle, i.e. Theorem~\ref{sem2thm9.13}, we get that 
$$
\P\(T|_Y,\varphi|_Y\)\leq \P(T,\varphi).
$$
\epf

\sp\subsection{Equilibrium States}
In light of the Variational Principle, the measures that maximize the
free energy of the system, 
that is, the measures which respect to which the free energy of the
system coincides with its topological pressure, are given a special name.

\sp\bdfn
If $T:X\lra X$ is a topological dynamical system and $\varphi:X\lra\R$ is a
continuous potential, then
a measure $\mu \in M(T)$ is called an \index{(N)}{equilibrium state}{\em
equilibrium state} for the potential $\varphi$ if and only if 
\[
\P(T,\varphi)=\h_\mu(T)+\int_X\varphi\,d\mu.
\]
\edfn

\sp Notice that if a potential $\varphi$ has an equilibrium state, then, after invoking Theorem~\ref{KBT-E}, the same reasoning as that of Corollary~\ref{sem2cor9.14} shows that $\varphi$ has an ergodic equilibrium state. 

When $\varphi=0$, the equilibrium states are also called
\index{(N)}{measure of maximal entropy}{\em measures of maximal entropy},
that is, measures for which 
$$
\h_\mu(T)=\h_{\mathrm{top}}(T).
$$
In particular, if $\h_{\mathrm{top}}(T)=0$, then every invariant measure 
is a measure of maximal entropy for $T$.
Note  that this is the case for homeomorphisms of the unit circle 
among other examples.

\sp It is natural at this point to wonder whether equilibrium
states exist for all topological dynamical systems. The answer, as the
following example demonstrates, is that they do not.  

\sp \bex
{\rm We will describe a topological 
dynamical system with finite topological entropy and with no measure of maximal entropy. Let 
$$
(T_n:X_n\to X_n)_{n=1}^\infty
$$ be a sequence of topological dynamical systems  with the property that
\[
\h_{\mathrm{top}}(T_n)<\h_{\mathrm{top}}(T_{n+1})
\hspace{0.5cm}\text{ and }\hspace{0.5cm}
\sup_{n\geq1}\lt\{\h_{\mathrm{top}}(T_n)\rt\}<+\infty.
\]
Let $\bigoplus_{n=1}^{\infty}X_n$ denote the topological disjoint union of the spaces $X_n$, and let 
$$
X:=\{\omega\}\cup\bigoplus_{n=1}^{\infty}X_n
$$ 
be the one--point compactification of $\bigoplus_{n=1}^{\infty}X_n$.
Define the map $T:X\to X$ by
\[
T(x)=\left\{
       \begin{array}{cl}
         T_n(x) & \hbox{if $x\in X_n$;} \\
         \omega & \hbox{if $x=\omega$.}
       \end{array}
     \right.
\]
Then $T:X\lra X$ is continuous. Suppose that $\mu$ is an ergodic measure of
maximal entropy for $T$. 
Then 
$$
\mu(\{\omega\})\in\{0,1\}
$$ 
since $T^{-1}(\{\om\})=\{\om\}$. But if $\mu(\{\omega\})=1$, then 
we would have 
$$
\mu\lt(\bigoplus_{n=1}^{\infty}X_n\rt)=0.
$$
Hence, on the one hand, we would have $\h_\mu(T)=0$, while, on 
the other hand,
$$
\h_\mu(T)=h_{\mathrm{top}}(T)\geq\sup_{n\geq1}\lt\{\h_{\mathrm{top}}(T_n)\rt\}>0.
$$
This contradiction implies that 
$$
\mu(\{\omega\})=0.
$$ 
Similarly, 
$$
\mu(X_n)\in\{0,1\}
$$ 
for all $n\geq1$ 
since $T^{-1}(X_n)=X_n$. Therefore, there exists a unique $N\geq1$
such that $\mu(X_N)=1$. It then follows that 
\[
\h_{\mathrm{top}}(T)
=\h_\mu(T)
=\h_\mu(T_N)
\leq \h_{\mathrm{top}}(T_N)
<\sup_{n\geq 1}\lt\{\h_{\mathrm{top}}(T_n)\rt\}
\leq \h_{\mathrm{top}}(T).
\]
This contradiction shows that there is no measure of maximal entropy
for the system $T$. }
\eex

Given that equilibrium states do not always exist, we would like to
find conditions under which they 
do exist. But since the function 
$$
\mu\longmapsto\int_X\varphi\,d\mu
$$ 
is continuous in the weak* topology on  the compact space $M(T)$, the function 
$$
\mu\longmapsto \h_\mu(T)
$$ 
cannot be continuous in general. Otherwise, 
the sum of these last two functions would be continuous and would
hence attain a maximum on the compact space $M(T)$, 
that is, equilibrium states would always exist. Nevertheless, the
function 
$$
\mu\longmapsto \h_\mu(T)
$$
is sometimes upper semicontinuous and this is sufficient
to ensure the existence of an equilibrium state. Let us first recall
the notion of upper (and lower) semicontinuity.

\sp\bdfn\label{defnupplowsemicont}
Let $X$ be a topological space. A function $f:X\lra [-\infty,+\infty]$ is
{\em upper semicontinuous} if and only if for all $x\in X$ it holds that
\[
\limsup_{y\to x} f(y)\leq f(x).
\]
Equivalently, $f$ is upper semicontinuous if and only if the set 
$$
\{x\in X:f(x)<r\}
$$ 
is open in $X$ for all $r\in\R$. 

A function $f:X\to[-\infty,+\infty]$ is called {\em lower semicontinuous} if and only if $-f$ is upper semicontinuous. 
\edfn

\sp  Evidently, a function $f:X\to[-\infty,+\infty]$ is continuous if and
only if it is both upper and lower semicontinuous. A perhaps most significant property of upper semicontinuous functions is the following well--known fact:

\sp\blem\label{lemuppsemicontattain}
If $f:X\lra[-\infty,+\infty]$ is an upper semicontinuous function on a
compact topological space $X$, then $f$ attains its upper bound on $X$.
\elem

\sp  One class of dynamical systems for which the function $\mu\mapsto
\h_\mu(T)$ is upper semicontinuous are all expansive maps.

\sp\bthm\label{sem2thm9.16}
If $T:X\lra X$ is a (positively) expansive topological dynamical system, then the function 
$$
M(T)\ni\mu\longmapsto\h_\mu(T)
$$ 
is upper semicontinuous. 
Hence each continuous potential $\varphi:X\lra\R$ has an equilibrium state.
\ethm

\bpf
Fix $\delta>0$ an expansive constant for $T$ and let $\mu\in M(T)$.
According to Lemma~\ref{chapter11lem3}, there exists a finite Borel partition ${\mathcal A}$ of $X$ with the property that 
$$
\mathrm{diam}({\mathcal A})<\delta
\  \  \  {\rm and} \  \  \ \mu(\partial A)=0
$$ 
for each $A\in{\mathcal A}$. Fix $\varepsilon>0$. Notice that since 
$$
\h_\mu(T)\geq \h_\mu(T,{\mathcal A})
=\inf_{n\to\infty}\lt\{\frac{1}{n}\H_\mu({\mathcal A}^n)\rt\}
$$ 
by Definitions~\ref{mte} and~\ref{entwrtpartition}, so there exists
$m\geq1$ such that 
\[
\frac{1}{m}\H_\mu({\mathcal A}^m)\leq \H_\mu(T)+\frac{\varepsilon}{2}.
\]
Now let $(\mu_n)_{n\geq1}$ be a sequence of measures in $M(T)$
converging weakly to $\mu$. 
Since $\mathrm{diam}({\mathcal A})<\delta$, it follows from
Theorem~\ref{sem2thm7.11} that 
\[
\h_{\mu_n}(T)=\h_{\mu_n}(T,{\mathcal A})
\]
for all $n\geq1$. Moreover, by Lemma~\ref{chapter11lem4} (with ${\mathcal
A}$ replaced by 
${\mathcal A}^m$), we have that
\[
\lim_{n\to\infty}\H_{\mu_n}({\mathcal A}^m)=\H_\mu({\mathcal A}^m).
\]
Therefore, there exists $N\geq1$ such that for all $n\geq N$, we have
\[
\frac{1}{m}\bigl|\H_{\mu_n}({\mathcal A}^m)-\H_\mu({\mathcal
  A}^m)\bigr|\leq\frac{\varepsilon}{2}. 
\]
Hence, for all $n\geq N$, we deduce that
\[
\h_{\mu_n}(T)
=\h_{\mu_n}(T,{\mathcal A})
\leq\frac{1}{m}\H_{\mu_n}({\mathcal A}^m)
\leq\frac{1}{m}\H_{\mu}({\mathcal A}^m)+\frac{\varepsilon}{2}
\leq \h_\mu(T)+\varepsilon.
\]
Consequently, 
$$
\limsup_{n\to\infty}\h_{\mu_n}(T)\leq \h_\mu(T)
$$ 
for any sequence $(\mu_n)_{n\geq1}$ in $M(T)$ converging weakly to $\mu$. Thus,
$$\limsup_{\nu\to\mu}\h_{\nu}(T)\leq \h_\mu(T),
$$
or, in other words, the function
$$
M(T)\ni\mu\longmapsto \h_\mu(T)
$$ 
is upper semicontinuous.
Since the function $\mu\mapsto\int_X\varphi\,d\mu$ is continuous in
the weak topology on the 
compact space $M(T)$, it follows that the function 
$$
\mu\longmapsto\h_\mu(T)+\int_X\varphi\,d\mu
$$ 
is upper semicontinuous. Lemma~\ref{lemuppsemicontattain} then yields that each continuous potential $\varphi$ admits an equilibrium state. 
\epf

\sp Finally, we show that the pressure function is Lipschitz
continuous.

\sp\bthm\label{sem2thm9.17}
If $T:X\lra X$ is a topological dynamical, then the pressure function $\P:C(X)\to\R$ is Lipschitz continuous with
Lipschitz constant $1$. 
\ethm

\bpf Let $\psi,\varphi\in C(X)$. Let also $\varepsilon>0$.
By the Variational Principle, there exists $\mu\in M(T)$ such that
\[
\P(T,\psi)\leq \h_\mu(T)+\int_X\psi\,d\mu+\varepsilon.
\]
Then
\begin{eqnarray*}
\P(T,\psi)
&\leq&
\h_\mu(T)+\int_X\varphi\,d\mu+\int_X\bigl(\psi-\varphi\bigr)\,d\mu+\varepsilon \\
&\leq&\P(T,\varphi)+\|\psi-\varphi\|_\infty+\varepsilon.
\end{eqnarray*}
Since this is true for all $\varepsilon>0$, we conclude that
\[
\P(T,\psi)-\P(T,\varphi)\leq\|\psi-\varphi\|_\infty.
\]
\epf

\part{Geometry and Conformal Measures}

\chapter{Geometric Function Theory}\label{geometric-theory}

\sp From now on, throughout this and the next part, essentially all our
considerations take place  in the (extended) complex plane and the maps are
assumed to be analytic. We collect in this chapter some selected theorems from geometric measures theory which will form for us indispensable tools in the following chapter, particularly when we will deal with the dynamics of elliptic functions. 

The most important for us are (various versions of) Koebe's Distortion Theorems. These will be used most frequently. Section~\ref{ELaMoTA} about extremal lengths and moduli of topological annuli is a preparation for formulating Koebe's Distortion Theorems in the context where ``Koebe's space/collar'' is an arbitrary topological annulus and and its results will be used throughout the book. 

The Riemann--Hurwitz Formula, which we treat at length in the last section of this chapter, is an elegant and probably the best tool to control the topological structure of connected components of inverse images of open connected sets under meromorphic maps. Especially to be sure that such connected components are simply connected.  

Given a set $A\sbt {\mathbb C}$ and $r>0$, the symbol
$$
B_e(A,r)
$$
\index{(S)}{$B_e(A,r)$} denotes the Euclidean open
ball ($r$-neighborhood) of the set $A$. Throughout the entire book
$f^*$\index{(S)}{$f^*$}, $\diam_s$\index{(S)}{$\diam_s$} and
$B_s(A,r)$\index{(S)}{$B_s(A,r)$} denote respectively the
derivatives, diameters and open balls 
defined by means of the spherical metric whereas
$f'$\index{(S)}{$f'$} and $\diam_e$\index{(S)}{$\diam_e$} are
considered in the Euclidean sense. The spherical
distance\index{(N)}{spherical distance} between any two points $x$
and $y$ in $\oc$ is denoted by $|x-y|_*$\index{(S)}{$\mid
\cdot \mid_*$}, where 
$$
\oc:= {\mathbb C} \cup \{\infty\},
$$
endowed with the natural conformal structure is the extended complex plane, frequently also called Riemann sphere. We emphasize that when calculating
$f^*$\index{(S)}{$f^*$}\index{(N)}{spherical derivative}, we
consider the spherical metric in the domain and in the codomain.

\sp


\bdfn\lab{dcomp} If $H:D\to {\mathbb C}$ is an analytic map, $z\in
{\mathbb C}$, and $r>0$, then by \index{(S)}{$\Comp(z,H,r)$}
$$
\Comp(z,H,r)
$$
we denote the only connected component of $H^{-1}(B_e(H(z),r))$ that contains $z$. 
\edfn

\section{Extremal Lengths and Moduli of Topological Annuli}\label{ELaMoTA}

The concepts of extremal length and moduli of topological annuli are intimately connected with the theory of quasiconformal maps. We will present a concise but self--contained account of these concepts without using quasiconformal maps. The results of this section will be used in the next one, on Koebe's Distortion Theorems and throughout the entire book. More complete and comprehensive exposition of these topics can be found in many books on, already mentioned, theory of quasiconformal maps. Among them are 
\cite{Ah}, \cite{Le}, \cite{LV}, \cite{AIM}, \cite{Gardiner}, \cite{Hubbard}, \cite{FdM}, \cite{FM}, and \cite{BF}. 

Let $D\sbt \mathbb C$ be a domain, i.e. an open connected set. We
say that a  curve $\g$ in $\mathbb C$ is piecewise
rectifiable\index{(N)}{piecewise rectifiable curve} if it can be
divided into countably many pieces with disjoint interiors each of
which  is rectifiable.  Let $\Ga$ be a family of piecewise
rectifiable curves in $D$. A Borel measurable function $\rho: D\to
[0,+\infty]$ is said to be admissible if
$$
A(\rho):=\int \int_D \rho^2dS\in (0,+\infty).
$$
\index{(S)}{$A(\rho)$} 

\fr Clearly $\rho(z)|dz|$ should be thought of as a measurable,
though conformal,  change of the standard Riemannian metric $|dz|$.
For any  $\g \in \Gamma$ we  thus define the length of $\g$ with
respect to the metric $\rho$\index{(S)}{$L_\rho(\g)$} as
$$
L_\rho(\g):=\int_\g\rho|dz|.
$$
Then 
$$
L(\rho)=L_\Ga(\rho):=\inf\{L_\rho(\g):\, \g \in \Ga\}.
$$\index{(S)}{$L_\rho$} 

\fr The inverse extremal length of $\Ga$ is defined as
\beq\lab{1_2017_12_11}
\l^{-1}(\Ga):=\inf\left\{\frac{A(\rho)}{L^2(\rho)}\right\},
\eeq
where the infimum is taken over all admissible metrics $\rho: D \lra
[0,+\infty]$. The quantity 
$$
\l(\Ga):= \frac{1}{\l^{-1}(\Ga)}
$$ 
is then called the
extremal length\index{(N)}{extremal length} of $\Ga$. Note that $\l^{-1}(\Ga)$ is a conformal invariant. More precisely,

\bobs\label{o1j215} If $D \sbt \mathbb C$ is a domain, $\Ga$ is a
family of rectifiable curves in $D$  and $\varphi: D \lra \mathbb C$
is a conformal isomorphism  onto  $\varphi(D)$, then
$$
\l^{-1}(\varphi(\Ga)):=\l^{-1}(\Ga),
$$ 
where $\varphi(\Ga)=\{
\varphi(\g):   \g \in \Ga\}$. 
\eobs

\

\fr Given $w \in \mathbb C$, $0 \leq R_1 \leq R_2 \leq +\infty$,
we define the geometric annulus as \index{(N)}{geometric annulus}

$$
A(w; R_1, R_2)=\{z \in \mathbb C: R_1 < |z-w| < R_2\}.
$$\index{(S)}{$A(w; R_1, R_2)$}

 Let $\Ga $ be the family  of all closed rectifiable  curves in $A(w;
R_1, R_2)$ separating $S_{R_1}:=\{z \in \mathbb C:   |z|=R_1 \}$
from $S_{R_2}$. We define $\Mod(A(w; R_1, R_2))$, the  modulus of
$A(w; R_1, R_2)$\index{(N)}{modulus of geometric annulus}, as
$$
\Mod(A(w; R_1, R_2)):= 2 \pi \l^{-1}(\Ga).
$$\index{(S)}{$\Mod(A(w; R_1,R_2))$}

\fr We shall prove  the following.

\sp\bthm\label{t2j215} $\Mod(A(w; R_1, R_2))=\log (R_2/R_1)$. \ethm

\bpf Assume without loss of generality, that $w=0$. Put
$A=A(0;R_1, R_2)$. Let $\rho$ be an  admissible Riemannian metric on
$A$. Working in polar coordinates $( \theta, r)$, we have for every
$r \in (R_1, R_2)$ that
$$
L(\rho)\leq \int_0^{2\pi} r \rho(re^{it})\, d \theta.
$$ 
Hence
\beq\label{1j217}
\begin{aligned} \int_{R_1}^{R_2}\int_0^{2\pi}  \rho( re^{it}) d
\theta dr & =\int_{R_1}^{R_2}\int_0^{2\pi} r  \rho( r
e^{it}) d \theta\frac{dr}{r}\\
   &\geq  \int_{R_1}^{R_2} L(\rho)\frac{dr}{r}=L(\rho) \log
(R_2/R_1).
\end{aligned}
\eeq
On the other hand, using Cauchy-Schwarz inequality, we get that
$$
\begin{aligned}  \left(\int_{R_1}^{R_2}\int_0^{2\pi} \rho(
r e^{it}) d \theta dr\right)^2 &
=\left(\int_{R_1}^{R_2}\int_0^{2\pi}
\sqrt{r}\rho( r e^{it})\frac{1}{\sqrt{r}} d \theta dr\right)^2\\
& \leq \int_{R_1}^{R_2}\int_0^{2\pi} r  \rho^2( r e^{it}) d \theta
 dr \int_{R_1}^{R_2}\int_0^{2\pi} \frac{1}{r}   d \theta
dr\\ & = 2\pi \log (R_2/R_1) A(\rho). \end{aligned}
$$
Combining   this with (\ref{1j217}), we get $A(\rho)/L^2(\rho)\geq
\log(R_2/R_1)$   whence
\beq\label{2j217}
\l^{-1}(\Ga) \geq  \frac{1}{2\pi}\log(R_2/R_1).
\eeq
To prove the opposite inequality,  consider the Riemannian metric
$\tau$ on $A$ that is defined as follows,
 \beq\label{4j217}
\tau(r e^{it})= 1/r.
 \eeq
So,
\beq\label{3j217}
A(\tau)
=\int_{R_1}^{R_2}\int_0^{2\pi} r\tau^2( r e^{it}) \,d \theta dr
=\int_{R_1}^{R_2}\int_0^{2\pi} \frac{1}{r} d \theta dr=2 \pi
\log(R_2/R_1).
 \eeq
But for every curve $\g : S^1 \to A$ belonging to $\Ga$, in fact for
every curve $\g : S^1 \to \mathbb C \sms \{0\}$ belonging to
$\Ga_{\mathbb C \sms \{0\}}$, we have
\beq\label{6j217}
\begin{aligned}
L_\tau(\g)
  =\int_\g \tau |dz|
& = \int_0^{2\pi}\frac{1}{|\tau(\theta)|}|
   \g'(\theta)|d \theta
   =\int_0^{2\pi}\left|\frac{\g'(\theta)}{\g(\theta)}\right|d \theta \\
&\geq \left|\int_0^{2\pi} \frac{\g'(\theta)}{\g(\theta)} d \theta\right|\\
& \geq 2\pi |\ind_0\g|\geq 2\pi, 
\end{aligned}
\eeq 
where $\ind_0\g$ denotes the index of the closed curve $\g$ with respect to the point $0$. Thus,
\beq\label{5j217}
L(\tau) \geq 2 \pi.
\eeq
Hence, invoking (\ref{3j217}), we  get that
$$ \l^{-1}(\Ga)\leq \frac{A(\tau)}{L^2(\tau)}\leq \frac{1}{2\pi}\log
(R_2/R_1).$$ Combining this with (\ref{2j217}), we thus get
$$ \l^{-1}(\Ga)=\frac{1}{ 2\pi} \log (R_2/R_1).$$
So, finally,
$$\Mod(A)=2 \pi \l^{-1}(\Ga)=\log (R_2/R_1).$$
The proof is complete. \endpf

\

\brem \label{r1j218} It follows from the proof of this theorem that
$L(\tau)=2\pi$ and $$ \frac{A(\tau)}{L^2(\tau)}= \Mod(A(w;R_1,
R_2))=\log (R_2/R_1).$$\erem

\sp\fr As an immediate consequence of this theorem we get the following.

\sp\bcor\label{c1j219} Two annuli $A(w;R_1, R_2)$ and  $A(z;R_1',
R_2')$ are conformally equivalent if and only if their  moduli are
the same, meaning that $R_2'/R_1'=R_2/R_1$.
\ecor

\sp\fr This corollary can be proved  without  using the concept of extremal
length; a particulary appealing short argument is given at the
beginning  of Chapter VII of Z. Nehari's book \cite{N}. We will
primarily need the concepts of  moduli and extremal length for other
purposes.

\sp We  say that an open subset of the complex plane $ \mathbb C$ is a
topological annulus\index{(N)}{topological annulus} if its
complement consists of exactly two connected components one of
which is bounded. It is well--known from topology that any two
annuli are homeomorphic.  Our nearest goal is to classify them up
to conformal  equivalency. 

Given an annulus $A$ let  $\Ga_A$  be the
family  of all closed rectifiable curves in $A$ that separate the
two connected components of its complement. We define    the modulus
of $A$\index{(N)}{modulus of topological annulus} to be
\beq\lab{2_2017_12_11}
\Mod(A):= 2 \pi \l^{-1}(\Ga_A).\index{(S)}{$\Mod( A)$}
\eeq
Note that this definition extends the one for geometric  annuli
$A(w;R_1,R_2)$. We shall  prove the following.

\sp\bthm\label{t1j221} 
Two annuli are conformally equivalent if and only if their  moduli  are equal. 
\ethm

\sp\fr Let $\mathbb H$ be the open upper half--plane, i.e.
$$
\mathbb H:=\{z\in\C:\im(z)>0\}.
$$
In order to prove Theorem~\ref{t1j221}, we will need the following.

\sp\bprop\label{p1j220} 
For every $\kappa>0 $ the holomorphic  map $
\Pi_\kappa: \mathbb H \lra  \mathbb C$  given by the formula
$$ 
\Pi_\kappa (z)= \exp\left( \frac{\kappa}{\pi} i \log z \right)
$$ 
is a holomorphic covering map  from $\mathbb H$ onto 
$A(0;e^{-\kappa}, 1)$, where 
$$
\log : \mathbb H \lra \mathbb C
$$ 
is the restriction to $\mathbb H$ of the principal branch of logarithm from $\C\sms (-\infty,0]$, to $\C$, i.e. the one sending $1$ to $0$. 

Furthermore, the group of deck  transformations  of $\Pi_\kappa$ coincides with ${<g>}$, the group generated by $g$, where $g:\mathbb H\lra\mathbb H$ is given by the formula
$$
g(z):=\exp\lt(\frac{2\pi^2}{\kappa}\rt)z.
$$

In particular, the annulus $A(0; e^{-\kappa}, 1)$ and the quotient space ${\mathbb H/<g>}$ are conformally isomorphic
\eprop

\bpf  First  observe  that $\log: \mathbb H\lra \mathbb C$
is a conformal   covering of $\mathbb H$ onto  its image $\{w \in
\mathbb C: 0< \im(w) < \pi\}$. So,
$$  
i\frac{\kappa}{\pi}\{w\in\mathbb C:  0< \im(w) < \pi\}
=\{\xi\in \mathbb C:  - \kappa <   \re(\xi) < 0\}.
$$ 
Furthermore, $\exp(\{ \xi\in \mathbb C:  - \kappa < \re \xi
< 0\})=A(0; e^{-\kappa}, 1)$. So,  $\Pi_\kappa$  as a composition of
conformal covering maps, is a conformal covering map onto its image
$A(0; e^{-\kappa}, 1)$. The first part of  the proposition  is thus
proved. In  order to prove the second part we note, that
$$
\begin{aligned}
\Pi_\kappa(g(z))
&=\exp\left( \frac{\kappa}{\pi} i \log \left(\exp
\left(\frac{2 \pi^2}{\kappa}\right)z\right)\right)
=\exp\left( i\frac{\kappa}{\pi}\left(\frac{2 \pi^2}{\kappa} +\log z\right)\right)\\
& =\exp\left(  \frac{\kappa}{\pi} i \log z \right)\exp( 2 \pi i)
=\exp \left(\frac{\kappa}{\pi} i\log z\right) \\
&=\Pi_\kappa(z).
\end{aligned}
$$
So, 
$$
\Pi_{\kappa}\circ g = \Pi_\kappa,
$$
and therefore $g$ belongs to
the deck group  of the covering map $ \Pi_\kappa$.  Now, if $
\Pi_\kappa(z)= \Pi_\kappa(w)$, then  
$$ 
i\frac{\kappa}{\pi} \log z-i \frac{\kappa}{\pi}\log w= 2 \pi i n
$$ 
for some integer $ n \in \mathbb Z$. This means  that 
$$ 
\log (w/z)= \frac{2 \pi^2}{\kappa}n.
$$
So,  
$$ 
w = \exp\lt(\frac{2\pi^2}{\kappa}n\rt)z = g^n(z).
$$
Thus, the group of deck  transformations  of $\Pi_\kappa$ coincides with
${<g>}$, the group  generated by $g$. Since $A(0; e^{-\kappa}, 1)$,
as  the range of $\Pi_\kappa$, is conformally isomorphic to the
quotient of $\mathbb H$ by  the deck group of $\Pi_\kappa$, we  are
thus done.
\endpf

\sp\bprop\label{p1j220a} 
The function $ \Pi_\infty: \mathbb H\lra
\mathbb C$, given by the formula $\Pi_\infty(z)=e^{iz}$, is a
holomorphic covering map from $\mathbb H $ onto $ A(0;0,1)$. In
particular, $ A(0;0,1)$ and ${\mathbb H /<g>}$ are  conformally
equivalent, where $g:\mathbb H\lra\mathbb H$ is the translation given by the formula 
$$
g(z)= z + 2 \pi.
$$
\eprop

\bpf Note that $ i \mathbb H=\{ w \in \mathbb C:  \re(w)
<0\}$, and so 
$$ 
\Pi_\infty(\mathbb H)=e^{i\mathbb H}=\{ \xi \in
\mathbb C:   0< |\xi| <1\}= A(0;0,1).
$$ 
Since $\Pi_\infty$ is
clearly a conformal covering map,  we are  thus done with the first
part. 

\sp In order too prove the second part, note that 
$$
\Pi_\infty(g(z))= e^{i(z+2\pi)}=e^{iz}=\Pi_\infty(z).
$$
Hence, $\Pi_\infty$ is a member of the deck
group of  $\Pi_\infty$. If 
$$
\Pi_\infty(w)=\Pi_\infty(z),
$$
then $ iz=iw+ 2 \pi i n$ for some $n \in  \mathbb Z$. This means that $w=z + 2\pi n= g^n(z)$. Thus, the group  of the deck transformation of
$\Pi_\infty$ is equal to $<g>$. So, since $A(0;0,1)= \Pi_\infty(\mathbb
H)$ and the quotient of $\mathbb H$ by the deck group of
$\Pi_\infty$ are conformally equivalent, we are thus done.
\endpf

\sp{\bf Proof of Theorem~\ref{t1j221}}. Since the  moduli of  annuli
are conformally invariant, it is enough  to show  that any two
annuli with the same moduli are conformally isomorphic. And by
virtue of Theorem~\ref{t2j215}, in order to do this, it is
sufficient to prove that  any annulus $A$ is conformally  equivalent
to a geometric  annulus $A(0;R_1,1)$.  

If $\mathbb C \sms A$ is a
singleton, say $\xi$, then the translation $z \mapsto  z + \xi $
establishes a  conformal equivalence between $A(0;0, \infty)$ and
$A$.  So, we may assume  that $\mathbb C \sms A$ contains  at least
two  points.  The Koebe--Poincar\'e's  Uniformization Theorem
asserts that there are exactly three  conformally  distinct  simply
connected Riemann surfaces,  namely  $\oc$, $ \mathbb
C$, and the upper half--plane $\mathbb H$. The universal
cover  of $A$  cannot  be $\oc$ since $A$ is not  compact
(purely  topological  obstacle), and it  cannot be $\mathbb C $
because  of the  Great Picard Theorem which asserts that every entire function from
$\mathbb C$ to $\mathbb C$  missing at  least two points is 
constant. So, the universal cover  of $A$ must be the half--plane
$\mathbb H$.  Let 
$$
\Pi: \mathbb H\lra A
$$ 
be a conformal  universal
cover of $A$. The deck group  $\Ga$ of $\Pi$ being isomorphic to the
fundamental group  of $A$  must be  isomorphic to $\mathbb Z$, the
group of integers. Let $g$   be a generator of $\Ga$. Being one
element of a deck group, $g$ cannot have fixed points, but $g$
extends  continuously  to $\ov{\mathbb H}$ (even biholomorphically
to $\oc$ because of the Schwartz Reflection Principle)
mapping it univalently onto itself. Keep  for this extension the
same symbol  $g$. Then $g :\ov{\mathbb  H} \to \ov{\mathbb  H}$ has
either exactly one fixed point or exactly two distinct fixed points.
In the former  case, applying  conformal conjugacy, we may assume without
loss of generality  that $g$  fixes $\infty$  and is of the form 
$$g(z)= z + 2\pi.$$ In the latter, also because of conformal
conjugacy, we may assume that $g(0)=0$ and $g(\infty)=\infty$, and
then $g:\ov{\mathbb H} \lra \ov{\mathbb H} $ must be of the form
$$ g(z)=\l z$$
with some $\l >0$. Let us deal first with the latter case. Writing
$$
\l= \exp\(2 \pi^2/\kappa\),
$$
we conclude from
Proposition~\ref{p1j220} that the annulus $A(0; e^{-\kappa},1)$ is conformally
isomorphic to $\mathbb H /<g>$ which is conformally isomorphic to
$A$ as $<g>$ is the deck group  of $\Pi: \mathbb H \to A$. The same
argument, but based on Proposition~\ref{p1j220a}  works also for the
case when $g(z)=z + 2 \pi$. The proof is  complete. 
\endpf

\sp Now, we shall easily prove two monotonicity   relation between
moduli of annuli. We say that an annulus $A$ is  essentially
contained in an annulus $B$ if $A$ separates the (two) boundary
components of $\mathbb C  \sms B$.  Equivalently, the inclusion $i:A
\to  B$ induces an isomorphism of  fundamental groups.

\

\bthm\label{t1j223} If $\{A_n\}_{n=1}^\infty$ are mutually disjoint
annuli, all essentially  contained in an annulus $A\sbt\C$, then 
$$
\Mod (A)\geq \sum_{n=1}^\infty \Mod(A_n).
$$
\ethm

\bpf  Let $\rho$ be an admissible metric for $A$. Since
all annuli $A_n$ are  mutually disjoint, we   have $A(\rho)\geq
\sum_{n=1}^\infty A(\rho|A_n)$. Also, $L(\rho) \leq L( \rho|_{A_n})$
for  all $n \geq 1$. Hence,
$$\frac{A(\rho)}{L^2(\rho)} \geq
 \sum_{n=1}^\infty\frac{A(\rho|_{A_n})}{L^2(\rho|_{A_n})} \geq
 \sum_{n=1}^\infty \Mod(A_n).$$
 Thus, 
 $$
 \Mod(A) \geq \sum_{n=1}^\infty \Mod(A_n).
 $$
The prof is complete. 
\endpf

\sp

\fr As an immediate consequence of this theorem we get the following.

\bcor\label{c1j225} 
If an annulus $A$ is essentially contained in an annulus $B\sbt\C$, then 
$$ 
\Mod(A) \leq \Mod(B).
$$
\ecor

\sp\fr Now, given $r \in (0,1)$ let
$$
B(r):=\mathbb D \sms [0,r].
$$
Clearly $B(r)$  is an annulus and, following tradition, we denote
its modulus by $\mu(r)$\index{(S)}{$\mu(r)$}. We call the map
$$
(0,1)\ni r\longmapsto \mu(r)\in[0,+\infty],
$$
the Gr\"otzsch modulus function. We have immediately from Corollary~\ref{c1j225} the following.

\bobs\label{o1_2017_10_35}
The Gr\"otzsch modulus function $(0,1)\ni r\longmapsto \mu(r)\in[0,+\infty]$ is monotone decreasing.
\eobs

We shall now show  that $B(r)$  is in a sense extremal amongst of all annuli
 that separate $0$ and $r$  from $S^1= \partial{\mathbb D}$ (see
 \cite{LV}, comp. \cite{Grotzsch}).

\bthm\label{t2j225} {\rm(Gr\"otzsch Module
Theorem)}\index{(N)}{Gr\"otzsch Module Theorem} If a ring $A\sbt
\mathbb C$ separates $0$ and $r$ from $\partial{\mathbb D}$, then
$$
\Mod (A) \leq \mu(r).
$$
\ethm

\bpf By virtue of Theorem~\ref{t1j221} there exists a
conformal homemorphism $\varphi: B(r) \lra A(0;1,e^{\mu(r)})$. Let
$$
B_{+}(r):= \{ z \in B(r):\, \im z > 0\}.
$$ 
Since $\partial{B}_{+}(r)$ is a Jordan curve, by virtue of Caratheodory's  Theorem, the map $\varphi_{+}:= \varphi_{ |B_{+}(r)}$ extends homeomorphically to a
map from $\ov{B_{+}(r)}$ to $\mathbb C$. By the Schwarz Reflection
Principle we can then  extend the map $ \varphi_{+}$ to a
holomorphic map $\varphi^*: \mathbb D \to  \mathbb C$ by setting
$$
\varphi^*(z):= \ov{ \varphi(\ov{z})}
$$ 
for all $z\in\{w \in \mathbb
D:\im(w) < 0\}$.  But $\varphi$ and $\varphi^*$  are holomorphic
and $ \varphi^*|_{B_{+}(r)}= \varphi$, whence  
$$
\varphi^*|_{B(r)}=\varphi.
$$
So, we can now get rid of the notation $ \varphi^*$ and
can consider  only a holomorphic map 
$$
\varphi: \mathbb D \lra \ov{A(0;1, e^{\mu(r)})}
$$
such that $\varphi|_{B(r)}$ is a conformal
homeomorphism onto $A( 0;1, e^{\mu(r)})$ and $\varphi$ is symmetric
with respect to the $x$-axis, meaning that
$$
\varphi(\ov{z})=\ov{\varphi(z)}.
$$
Now, consider the metric $\tau$ on $A( 0;1, e^{\mu(r)})$ given by the formula~(\ref{4j217}). Define the pull--back metric on $B(r)$ by the formula
$$
\varphi^*(\tau|dz|):=
\tau(\varphi(z))|\varphi'(z)||dz|=\frac{|\varphi'(z)|}{|\varphi(z)|}|dz|,
$$
i.e. put
$$
\rho(z):=\frac{|\varphi'(z)|}{|\varphi(z)|}.
$$
Since the map $\C\ni z\longmapsto \ov{z}$ is an Euclidean isometry, we conclude
from the Chain Rule that $|\varphi'(\ov{z})|=|\varphi'(z)|$, and
therefore,
\beq\label{1j227}
\rho(\ov{z})= \rho({z}),
\eeq
which just means that the metric $\rho$ is symmetric with respect to
the $x$--axis. Consider $\rho$ as a metric defined on $B$. Since the
area of the segment $[0,r]$ with respect to the metric $\rho$ is
equal to  zero, we have
$$ A_\rho(B)\leq A_\rho(B(r))=A_\tau(A(0;1,e^{\mu(r)}))=: A(\tau).$$
It therefore  follows from (\ref{3j217}) that
\beq\label{1j229}
 \Mod( A)\leq 2 \pi
\frac{A(\rho)}{L^2(\rho)} \leq 4 \pi^2 \frac{\mu(r)}{L^2(\rho)}.
\eeq
Now consider an arbitrary curve $\g \in \Ga_A$. Since $\g$ separates
the points $0$ and $r$  from the unit circle $S^1$, it can be
divided into two subarcs $\g_1$ and $ \g_2$ with disjoint interiors
both having one common  endpoint on the segment  $(-1, 0)$ and one
on the segment $(r, 1)$. Let $\hat{\g}_1$  be  the curve  resulting
from $\g_1$ by  replacing each point $z$   on $\g_1$ with $\im z <0$
by its  conjugate $\ov{z}$ and let $ \hat{\g_2}$ be the curve
resulting  from $\g_2$ by replacing  each point $z$ on $\g_2$ with
$\im z>0$ by its conjugate $\ov{z}$. Since  the metric $\rho$ is
symmetric with respect  to the $x$-axis  (see(\ref{1j227})), we have
$L_\rho(\hat{\g_i})=L_\rho(\g_i)$, $i=1,2$. Obviously
$\varphi(\hat{\g_1} \cup \hat{\g_2})\in \Ga_{\mathbb C\sms \{0\}}$.
It then follows from (\ref{6j217}) that
$$\begin{aligned}
 L(\rho)&\geq L_\rho(\g)=L_\rho(\g_1 \cup \g_2)=
L_\rho(\g_1) + L_\rho(\g_2)=L_\rho(\hat{\g_1})+L_\rho(\hat{\g_2})\\
 & = L_\rho(\hat{\g_1} \cup \hat{\g_2})= L_\tau( \varphi(\hat{\g_1} \cup
 \hat{\g_2})) \\
 &\geq 2 \pi.
 \end{aligned} 
 $$
Inserting this to (\ref{1j229}), we finally  get that $\Mod(A)\leq
\mu(r)$. The proof  is complete. 
\endpf

\sp Now we shall  prove the following estimates for $\mu(r)$. 

\sp\bprop\label{p1j231} 
For all $r\in(0,1)$, we have $ \mu(r)< \log 4 -
\log r$. 
\eprop

\bpf Let $R>0$. Let $ M: \oc \lra \oc$  be the unique linear affine  map (M\"obius transformation)
sending the points $0$, $r$ and $1$ respectively to $-1, 1$ and $R$.
Since all M\"obius transformations  preserve cross-ratios, we get for
all $ z \in \oc$ that
$$ 
\frac{1}{r} \cdot \frac{z-r}{z-1}
= \frac{ 1+R}{2} \cdot \frac{M(z)-1}{M(z)-R}.
$$ 
We now require $R$ to be such that $M(-1)=-R$. Then, we get the
following  quadratic equation for $R$
\beq\label{1j231}
\frac{1+R}{2 }\cdot
\frac{1+R}{2R}=\frac{1}{r}\cdot\frac{1+r}{2}\  \  \Leftrightarrow \  \  r
(1+R)^2=2(1+r)R  \  \  \Leftrightarrow \  \  rR^2-2R+r =0.
\eeq
With $R \in \mathbb R$ being a solution to this equation the map $M$
maps the unit circle  $\{z\in\C:|z|=1\}$ onto the circle $\{z\in\C:|z|=R\}$ as
$M(S^1)$ is a  circle, $ M(\ov{\mathbb R})= \ov{\mathbb R}$,
$M(S^1)$ intersects $\ov{ \mathbb R}$ at right angles (as  M\"obius
map preserve angles), $M(-1)=-R$ and $M(1)=R$. The
equation~(\ref{1j231}) has two solutions 
$$
R=r^{-1}(1-\sqrt{1-r})  
\  \ {\rm and }  \  \
R=r^{-1}(1+\sqrt{1-r})\in (0,1).
$$
We fix the later one. Then
$M(B(0,1))=B(0, R)$ as it can be easily checked by evaluating
$M(\infty)$.  As also  $M(0)=-1$ and $M(r)=1$, we get that $M([0,
r])=[1,-1])$. Therefore $\hat{M}$, the  restriction  of $M$ to
$B(r)$, establishes  a conformal isomorphism between $B(r)$ and
$$
G(r):=M(B(r))=B(0, R)\sms [-1, 1].
$$
Now  consider the rational function $H: \oc\lra
\oc$ given by the formula 
$$
H(z)=\frac{1}{2}\left(z+\frac{1}{z}\right).
$$
This is a rational function of degree 2 and,
as $H(z)=H(1/z)$, it is injective on $\{z \in \oc: \,
|z|>1\}$. If  $|z|= \rho$, then writing $ z =  \rho e^{i\theta}$,
$0\leq \theta < 2 \pi$, we get 
$$
H(z)=\rho e^{it} +
\frac{1}{\rho}e^{- i \theta}=\left( \rho+ \frac{1}{\rho}\right)\cos
\theta +i\left( \rho- \frac{1}{\rho}\right)\sin \theta.
$$  
Hence if $\rho=1$, then 
$$
H(\{z\in\C:|z|=1\})=[-1, 1],
$$
and if $\rho>1$, then 
$$
H(\{z\in\C:|z|=\rho\})= E_\rho,
$$
the ellipse  with axes points  $\pm
\left(\rho+ 
\frac{1}{\rho}\right)$ and $\pm \left(\rho- \frac{1}{\rho}\right)$.
This is so as
$$ 
\frac{(( \rho+\frac{1}{\rho})\cos\theta)^2}{(\rho+\frac{1}{\rho})^2} +
 \frac{(( \rho-\frac{1}{\rho})\sin
\theta)^2}{(\rho-\frac{1}{\rho})^2} =1.
$$ Thus, $\hat{H}$, the map
$H$ restricted to the annulus $A(0;1,\rho)$ establishes a conformal
isomorphism  between the annulus $A(0;1,\rho)$ and $E_\rho\sms
[-1,1]$. Now, if $\rho=\frac{4}{r}$, then
$$\
\rho-\frac{1}{\rho}=\frac{4}{r}-\frac{r}{4}= \frac{1}{r}\left(4 -
\frac{r^2}{4}\right).
$$
We check that $4 - \frac{r^2}{4}\geq
2(1+\sqrt{1-r^2})$. This inequality means that
$$
2-\frac{r^2}{4}\geq \sqrt{1-r^2} \  \  \Longleftrightarrow \  \  4 -
r^2+\frac{r^4}{16}\geq 1- r^2,$$ which is true. Hence
$\rho-\frac{1}{\rho}\geq 2 R$, and therefore, $G(r)\varsubsetneq
E_{4/r}\sms [-1,1]$. Applying Corollary~\ref{c1j225}, we thus obtain
$$
\begin{aligned}
\mu(r)
&=\Mod(B(r))
=\Mod(G(r))
\leq \Mod(E_{4/r}\sms[-1,1]) \\
&=\Mod(A(0; 1, 4/r))
=\log4 - \log r.
\end{aligned} 
$$ 
We are done. 
\epf

\sp For every set $B\sbt \mathbb C$ and every point $ z \in  \mathbb C$,
let
$$ 
r_z(B)= \sup \{ |w-z|:\,  w \in B\}.
$$\index{(S)}{$r_z(B)$}

\fr We are now in position to prove easily  the following. 

\sp\bthm\label{t1j233} If $B$ is compact connected subset of the
complex plane $\mathbb C$, $ x \in B$ and $R > r_x(B)$, then
$$ 
r_x(B)\leq 4 \exp\(-\Mod ( ( B(x,R)\sms B)_{*})\),
$$
where $( B(x,R)\sms B)_{*}$  is the connected component (an annulus)
separating  $x$  from $ \partial{B}(x, R)$. 
\ethm

\bpf By Theorem~\ref{t2j225} and Proposition~\ref{p1j231}, we
get  that
$$ 
\Mod((B(x,R)\sms B)_{*})\leq \mu ( r_x(B))< \log 4 - \log r_x(B).$$
Exponentiating  this inequality,  we  get that
$$ 
r_x(B)\leq 4 \exp\( -\Mod ((B(x,R)\sms B)_{*})\).
$$
We are done. 
\epf

\sp We know, by Observation~\ref{o1_2017_10_35} that the limit $\lim_{r\upto 1}\mu(r)$ exists. We need to know that the value of this limit is equal to zero. We will prove it now. In order to do this, we shall first provide an uper estimate of the modulus of an annulus which is also interesting on its own.

\sp\bthm\label{l_Lehto_6.2}
Let $A\sbt\C$ be an annulus, let $\d$ be the minimum of spherical diameters of  boundary components of $A$, and let $\varepsilon$ be the spherical distance between these (two) components. If  $0<\varepsilon < \delta$, then
\beq\label{(L_6.6)}
\Mod(A) \leq \frac{ \pi^2}{\log\left(\frac{ \tan(\delta/2)}{\tan(\varepsilon/2)}\right)}.
\eeq
\ethm  

\bpf Let $A_0^c$ be the bounded connected component of $\C\sms A$ and let $A_*^c$ be the unbounded connected component of $\C\sms A$. Since $\varepsilon < \delta$, there  exist two points $x\in A_0^c$ and $y\in A_*^c$ such that $|y-x|_s<\d$. Since 
$$
\big|\tan(|y-x|_s)-(-\tan(|y-x|_s)\big|_s=|y-x|_s,
$$
by an affine change of coordinates, which corresponds to a rotation of the Riemann sphere $\oc$, we may assume without loss of generality that 
$$
x=\tan(|y-x|_s) 
\  \  {\rm and} \  \
y=-\tan(|y-x|_s).
$$
Both components of the complement of the annulus
$$
\tilde A:=A\(0;\tan(|y-x|_s/2), \tan( \delta/2)\)
$$
then contain points of $A_0^c$ and $A_*^c$. It follows that every rectifiable (in fact every (continuous)) curve $\g$ in $A$ which separates $A_0^c$ and $A_*^c$, i.e. which belongs to $\Ga_A$, intersects both connected components of $\C\sms\tilde A$. Then denote by $\g_1$ and $\g_2$ some two connected components of $\g\cap\tilde A$ that are separated by $\bd B(0,\tan(|y-x|_s/2)$. Denote the endpoints of the first of these two components by $\xi_1$ and $\xi_2$. Set
$$
\tau(z):=
\begin{cases}
1/|z| & {\rm if} \  \ z\in \tilde A \\
0     & {\rm if} \  \ z\in A\sms \tilde A,
\end{cases},
$$
i.e. $\tau:A\lra[0,+\infty)$ is the same measurable Riemannian metric as the one defined in \eqref{4j217}. Calculating in the same way as in \eqref{6j217}, we get
$$
\begin{aligned}
\int_{\g_1} \tau |dz|
&\geq \left|\int_0^{2\pi} \frac{\g'(\theta)}{\g(\theta)} d\theta\right|
=\big|\log(\xi_2)-\log(\xi_1)\big|
\ge \big|\log|\xi_2|-\log|\xi_1|\big| \\
&=\big|\log\(\tan(|y-x|_s/2)\)-\log\(\tan(\delta/2)\)\big| \\
&=\log\lt(\frac{\tan(|y-x|_s/2)}{\log\(\tan(\delta/2)}\)\rt),
\end{aligned}
$$
where $\log(\xi_2)$ and $\log(\xi_1)$ are some appropriate choices of logarithms respectively of $\xi_2$ and $\xi-1$. 
Hence,
$$
L_\tau(\g)
=\int_\g \tau |dz|
\ge\int_{\g_1} \tau |dz| + \int_{\g_2} \tau |dz|
\ge 2\log\lt(\frac{\tan(|y-x|_s/2)}{\log\(\tan(\delta/2)}\rt).
$$
Thus
\beq\lab{3_2017_12_11}
L(\tau)\ge 2\log\lt(\frac{\tan(|y-x|_s/2)}{\log\(\tan(\delta/2)}\rt).
\eeq
On the other hand, similarly as in \eqref{3j217}, with $R_1:=\tan(|y-x|_s/2)$ and $R_2:=\tan(\delta/2)$, we get
$$
\begin{aligned}
A(\tau)
&\le \int_A|z|^{-1}dS(z)
\le\int_{R_1}^{R_2}\int_0^{2\pi} r\tau^2( r e^{it}) d \theta dr
=\int_{R_1}^{R_2}\int_0^{2\pi} \frac{1}{r} d \theta dr
=2 \pi\log(R_2/R_1) \\
&=2 \pi\log\lt(\frac{\tan(|y-x|_s/2)}{\log\(\tan(\delta/2)}\rt).
\end{aligned}
$$
It directly follows from this formula, \eqref{3_2017_12_11}, along with formulas \eqref{2_2017_12_11} and \eqref{1_2017_12_11}, that
$$ 
\Mod(A) 
\leq 2\pi \frac{A(\tau)}{L^2(\tau)}
\leq \frac{\pi^2}{\log\lt(\frac{\tan(|y-x|_s/2)}{\log\(\tan(\delta/2)}\rt)}
$$
So, letting $|y-x|_s$ converge to $\varepsilon$, we get that
$$ 
\Mod(A) 
\leq 
\frac{\pi^2}{\log \left(\frac{\tan(\delta/2)}{\tan( \varepsilon/2)}\right)}, 
$$
and the proof is complete. 
\qed

\sp Now as an immediate consequence of this theorem, we get the result we were after.

\sp\bprop\label{p2_2017_10_27}
If $(0,1)\ni r\longmapsto \mu(r)\in[0,+\infty]$ is the Gr\"otzsch modulus function, then 
$$
\lim_{r\upto 1}\mu(r)=0.
$$
\eprop

\fr We end this section with the following.

\sp\bthm\label{t1j235} 
If $A$ and $B$ are two annuli contained in the complex plane $\C$, and $f: A \to B$ is a conformal covering map, then
$$
\Mod(B)= \deg(f) \Mod(A).
$$
\ethm

\bpf  Because of Theorem~\ref{t1j221} and
Theorem~\ref{t2j215}, we may assume  without loss of generality that
$B=A(0;1,R)$. Let  $d=\deg(f)$. Consider the map $E_d:
A(0;1,R^{1/d}) \lra B$, given by the formula
$$
E_d(z):=z^d.
$$
This is a covering map of degree $d$. Fix
a point  $w \in B$ and $w_1 \in A$, $w_2\in A(0;1,R^{1/d})$ such
that 
$$
f(w_1)=w \  \  {\rm and} \ \  E_d(w_2)=w
$$
Consider the fundamental groups
$$
\pi_1(A, w_1), \  \pi_1( A(0;1,R^{1/d}), w_1), \  \  {\rm and} \ \  \pi_1(B,w).
$$
Since both maps $f$ and $E_d$ are coverings with degree $d$, we get that  
$$
f_*(\pi_1(A, w_1))=d \pi_1(B,w)= (E_d)_*(\pi_1(A(0;1,R^{1/d}), w_1)).
$$
So, there exists a continuous map $H:A \lra A(0;1,R^{1/d})$ such that $E_d\circ H = f$ and $H(w_1)=w_2$. Since both $f$ and $E_d$ are covering map, so is also $H$. In particular the degree, $\deg (H)$ is well--defined and
$$
d=\deg(f)=\deg(E_d\circ H)=\deg(E_d) \deg(H)=d \deg(H).
$$ So,
$\deg(H)=1$, and as $H: A \lra A(0;1,R^{1/d})$ is covering, it must
be a homeomorphism. Thus, $H: A \to A(0;1,R^{1/d})$  is a conformal
isomorphism  and so,
$$
\Mod(A)=\Mod( A(0;1,R^{1/d}))= \frac{1}{d} \log R
=\frac{1}{\deg(f)}{\Mod(B)}.
$$ Equivalently, $\Mod(B)= \deg(f)
\Mod(A)$ and we are done. \endpf

\sp

\section{Koebe's Distortion Theorems}\label{KDT}

\, This section is entirely devoted to formulate and to prove various
versions of Koebe's Distortion Theorems. They are truly amazing
features of univalent holomorphic functions and form one of the main
indispensable tools when dealing with dynamical and, especially geometric, aspects, of meromorphic maps in the complex plane. Koebe's Distortion Theorems will be one of the most frequently invoked theorems in Part~\ref{EFA}, Elliptic Functions A, and Part~\ref{EFB}, Elliptic Functions B, of our book. 

Koebe's $\frac{1}{4}$--Theorem, i.e. Theorem~\ref{one-quater} was conjectured by Paul Koebe in 1907 and was proved by Ludwig Bieberbach in 1916 in \cite{Bie}. The proof is a fairly easy consequence of Bieberbach's Coefficient Inequality, i.e. Theorem~\ref{t1.2inGC}, obtained in \cite{Bie}. Likewise, analytic Koebe's Distortion Theorems, such as Theorem~\ref{th17.46in{H}A} and Theorem~\ref{th17.46in{H}} follow from Bieberbach's Coefficient Inequality. Our proofs are standard and closely follow the exposition of \cite{Hi}.

\sp All other distortion theorems proved in this section are relatively simple consequence of the ones mentioned above and some results of the previous section, Extremal Lengths and Moduli of Topological Annuli.

\sp We start with proving Theorem~\ref{t1.1inGC} and Theorem~\ref{t1.2inGC}.
The latter will form a crucial ingredient in the proof of the full
version of the first Koebe's Distortion
Theorem, i.e. Theorem~\ref{one-quater}, following it. Let $\mathcal S$ 
denote the class of all univalent holomorphic functions
$f:B(0,1)\lra\C$ such that 
$$
f(0)=0 \  \text{ and }  \  f'(0)=1.
$$

\sp\bthm\label{t1.1inGC}{\rm(Area Theorem)}\index{(N)}{Area Theorem}
Let $g:B(0,1)\to\C$ be a univalent meromorphic function with a simple
pole at $0$. Assume that the residue of $g$ at $0$ is equal to
$1$, so that the function $g$ can  be represented in a form
$$
g(z)=1/z+ b_0+b_1z+\ldots
$$ 
Then 
$$
\sum_{n=1}^\infty n|b_n|^2\leq 1.
$$ 
\ethm

\bpf For every $0<r <1$ put $D_r=\mathbb C \sms
g(B_e(0,r))$. By  Green's theorem we get, 
\beq\label{120130513}
S(D_r)
=\int\int_{D_r} dxdy
=\frac{1}{2i}\int_{\partial D_r}\ov{z}dz 
=- \frac{1}{2i}\int_{\partial B_e(0,r)}\ov{g}dg. 
\eeq
Recall that 
$$
\frac1{2\pi i}\int_{\partial B(0,r)}z^k\bar z^ldz=\d_{kl}.
$$
So, substituting the power series expansions for $g$ and $g'$ to
\eqref{120130513}, and performing the integration, we  obtain
$$
S(D_r)
=\pi\left(\frac{1}{r^2}-\sum_{n=1}^\infty n|b_n|^2r^{2n}\right).
$$ 
Since $S( D_r)\geq 0$, taking the limit as $r \upto 1$, yields
the desired result. 
\endpf

\sp \bthm[Bieberbach's Coefficient Inequality]\label{t1.2inGC} If $f(z)=z + \sum_{n=2}^\infty a_n z^n \in
\mathcal S$, then $|a_2|\leq 2$. 
\ethm

\bpf Note that the formula
$$
h(z):=\frac{f(z^2)}{z^2}
     = 1 + \sum_{n=2}^\infty a_n z^{2n-2}
     = 1 + \sum_{n=2}^\infty a_n z^{2(n-1)}
$$
defines a holomorphic function from $B_e(0,1)$ to $\C\sms\{0\}$ such
that $h(0)=1$. Let 
$$
\sqrt h:B_e(0,1)\lra\C\sms\{0\})
$$ 
be the square root of $h$ uniquely determined by the requirement that $0\mapsto 1$.
Let
$$
g(z)=\frac1{z\sqrt{h(z)}}=\frac1z-\frac12a_2z+ \ldots,
$$
and this series contains only odd powers of $z$ so that the function
$g$ is odd. If now
$g(z_1)=g(z_2)$, then $f(z_1^2)=f(z_2^2)$, so $ z_1^2=z_2^2$, and
$z_1=\pm z_2$. But since $g$ is odd, we get that $z_1=z_2$. Thus $g$
is univalent, whence Theorem~\ref{t1.1inGC}  gives  that $|a_2| \leq
2$. The proof is complete.
\endpf

\sp

\bthm\label{one-quater}{\rm(Koebe's
$\frac{1}{4}$--Theorem).}\index{(N)}{Koebe's ${1\over 4}$-Theorem}
If $w\in {\mathbb C}$, $r>0$, and $H:B_e(w,r)\lra {\mathbb C}$  is an
arbitrary univalent analytic function, then 
$$
H(B_e(w,r))\spt B_e(H(w),4^{-1}|H'(w)|r).
$$ 
\ethm  

\bpf Precomposing $H$ with the scaled (by factor $1/r$)
translation moving $0$ to $w$ and $B_e(0,r)$ onto $B_e(w,r)$, and the
postcomposing with a scaled (by factor $r$) translation moving $H(w)$
to $0$, we may assuming without loss of generality that $H\in\cS$. 
Fix now a point $c\in \C\sms H(B(0,1))$. Then an immediate inspection
shows that the function
$$
B_e(0,r)\ni z\longmapsto\frac{c H(z)}{c-H(z)}\in\C
$$
belongs to $\cS$ and takes on the form
\beq\label{320130513}
B_e(0,r)\ni z\longmapsto z+\left(a_2+ \frac{1}{c}\right) z^2+\ld,
\eeq
where, as usually, $a_2$ is the coefficient of $z^2$ in the Taylor
series expansion of $H$ about $0$.
Applying now 
Theorem~\ref{t1.2inGC} twice, to the function in \eqref{320130513} and
to $H$, we obtain
$$ \frac{1}{|c|}\
\leq |a_2|+\left|a_2+\frac{1}{c}\right| 
\leq 2 +2=4.
$$
The proof is complete.
\endpf

\sp The Koebe's function 
$$
B_e(0,1)\ni z\lmt f(z):=\frac{z}{1-z^2}=\sum_{n=1}^\infty nz^n\in \C
$$ 
maps univalently the ball unit $B_e(0,1)$ to the slit plane $\mathbb  C \sms (-\infty, 1/4]$. This shows that the number $1/4$ is optimal in the above theorem.

\sp We  will now prove a series of Koebe's Distortion Theorems. First we
establish Lemma~\ref{l17.4.1}, which will form a crucial 
ingredient in the  proof of full version of Theorem~\ref
{th17.46in{H}}, following it.

\sp\blem\label{l17.4.1} If $f \in \mathcal S$ and $|z| <1$, then
\beq\label{17.4.19}
\left|\frac{1}{2}(1-|z|^2)\frac{f''(z)}{f'(z)}-\ov{z}\right|\leq 2
\eeq
for all $z\in B_e(0,1)$.
\elem

\bpf We fix $z\in B_e(0,1)$. The parameter $w\in B_e(0,1)$ will be a
variable throughout the proof.   
The M\"obius transformation 
$$
w\longmapsto \frac{w+z}{1+\ov{z}w}
$$
maps the ball $B_e(0,1)$ onto itself, sending $0$ to $z$. It follows  that
for any choice of constants $C_1$ and $C_2$ the function 
$$
g(w)=C_1+ C_2 f\left(\frac{w+z}{1+\ov{z}w}\right)
$$ 
is univalent in $B_e(0,1)$.
Now we specify $C_1$ and $C_2$ so that $g(0)=0$ and $g'(0)=1$, that
is, so that $g\in \mathcal S$. These conditions give
$$ 
C_1
=-\frac{f(z_0)}{f'(z)(1-|z|^2)}, \quad C_2
=\frac{1}{f'(z)(1-|z_0|^2)}. 
$$
Direct calculations give
$$
g'(z)=\frac{1}{f'(z)}\frac{1}{(1+\ov{z}w)^2}
       f'\left(\frac{w+z}{(1+\ov{z}w}\right)
$$
and
$$
g''(w)=
\frac{1}{f'(z)}\frac{2 \ov{z}}{(1+\ov{z}w)^3}
f'\left(\frac{w+z}{1+\ov{z}w}\right)
+ \frac{1}{f'(z)}\frac{1-|z|^2}
{(1+\ov{z}w)^4}f''\left(\frac{w+z}{1+\ov{z}w}\right).
$$
It thus follows that 
$$
g(w)=w+\left(\frac{1}{2}(1-|z|^2)
\frac{f''(z)}{f'(z)}-\ov{z}\right)w^2+\ldots.
$$ 
Since $g \in \mathcal S$, the coefficient of $z^2$ is in modulus
$\leq 2$ because of Theorem~\ref{t1.2inGC}. And this just
is the assertion of our lemma.
\endpf

\sp Now we formulate and prove the following theorem which is the
central one among Koebe's distortion theorems.

\sp\bthm\label{th17.46in{H}A}{\rm(Koebe's Distortion Theorem, Analytic
Version I)}\index{(N)}{Koebe's Distortion Theorem (analytic version)}
If $f \in \cS$, then
 \beq\label{17.4.18A}
 \frac{1-r}{(1+r)^3} 
\leq |f'(z)|
\leq  \frac{1+r}{(1-r)^3}
\leq r,
\eeq
and there is a choice of argument of $f'(z)$ such that
\beq\label{17.4.18bA}
 |{\rm arg} f'(z)| \leq 2 \log \frac{1+r}{1-r}
\eeq
for all $z\in\ov B(0,r)$.
 \ethm

\bpf By Lemma~\ref{l17.4.1} for every $t\in B_e(0,1)$ there
exists $\eta(t)\in \ov B_e(0,1)$ such that
$$
\frac{f''(t)}{f'(t)}-\frac{2\ov{t}}{1-|t|^2}
=4\frac{\eta(t)}{1-|t|^2}. 
$$
We integrate this expression along the line segment from $0$ to $z$ 
in $B(0,1)$ to obtaine
\beq\label{17.4.20}
\log f'(z)+ \log(1-|z|^2)= 4 \int_0^z\frac{\eta(t)}{1-|t|^2}dt,
\eeq
where the modulus of the right-hand side expression does not exceed
\beq\label{120130514}
4\int_0^r\frac{ds}{1-|s|^2}=2 \log \frac{1+r}{1-r}
\eeq
for all $z\in \ov B(0,r)$. Taking then the real parts of both sides of
\eqref{17.4.20} we get the inequality
$$
2 \log \frac{1-r}{1+r}
\leq \log\((1-|z|^2)|f'(z)|\)
\leq 2\log\frac{1+r}{1-r}.
$$ 
Since the function $\ov B(0,r)\ni z\mapsto\log |f'(z)|$ is harmonic,
it assumes its maximum and minimum values on $\bd B(0,r)$. This yields
(\ref{17.4.18A}). Formula (\ref{17.4.18bA}) is what we get by equating
imaginary parts in (\ref{17.4.20}) and using \eqref{120130514}. 
\endpf

\

\fr Part (\ref{17.4.18bA}) is called the
 Rotation Theorem.\index{(N)}{Rotation Theorem} It was
discovered by Bieberbach in 1919. The estimates in (\ref{17.4.18A})
are the best possible; equality is reached for the functions
\beq\label{17.4.12}
f_\b(z)=\frac{z}{(1+e^{i\b}z)^2},
\eeq
where $\b$ is an arbitrary real number.

\sp Let $\mathcal U$ denote the class of all univalent and holomorphic
functions from the unit disk $B_e(0,1)$ into $\C$. With obvious
translations and rescalings, as an immediate consequence of
Theorem~\ref{th17.46in{H}A}, we get the following.

\sp\bthm [Koebe's Distortion Theorem, Analytic
Version II]\label{th17.46in{H}}\index{(N)}{Koebe's Distortion Theorem (Analytic version)}
If $f \in \mathcal U$, then
 \beq\label{17.4.18}
 \frac{1-r}{(1+r)^3} 
\leq \frac{|f'(z)|}{|f'(0)|}
\leq  \frac{1+r}{(1-r)^3}
\leq r,
\eeq
and there is a holomorphic branch of argument of $f'(z)$ such that
\beq\label{17.4.18b}
\lt|{\rm arg}\lt(\frac{f'(z)}{f'(0)}\rt)\rt| 
\leq 2 \log \frac{1+r}{1-r}
\eeq
for all $z\in\ov B_e(0,r)$.
 \ethm

\sp\fr  As an immediate consequence of Theorem~\ref{th17.46in{H}}, we get
the following two facts.

\sp\bthm[Koebe's Distortion Theorem, Analytic
Version III]\label{tKoebe-1}
There exists a monotone increasing continuous function 
$$
K:[0,1)\lra [1,+\infty)
$$ 
such that $K(0)=1$ and with the following property. If $w\in {\mathbb C}$,
$R>0$, and $H:B_e(w,R)\lra {\mathbb C}$ is an arbitrary univalent
analytic function, then  
\beq\label{Koebe-1}
\lt|\frac{|H'(z)|}{|H'(w)|}-1\rt|\le K(r/R)|z-w|.
\eeq
for every $r\in[0,R]$ and for all $z\in \ov B_e(w,r)$
\ethm

\sp

\bthm [Koebe's Distortion Theorem I, Euclidean
Version]\label{Euclid-I}\index{(N)}{Koebe's Distortion Theorem I, Euclidean
Version}
There exists a monotone increasing continuous function 
$$
k:[0,1)\to [1,\infty)
$$ 
such that $k(0)=1$ and for any $w\in {\mathbb C}$, any $r >0$, all $t\in [0,1)$, and any univalent analytic function $H:B_e(w,r)\lra {\mathbb C}$, we have that
$$
\sup\big\{|H'(z)|:z\in B_e(w,tr)\big\} \le k(t) \inf\big\{|H'(z)|:z\in
B_e(w,tr)\big\}.
$$
We put $K:=k(1/2)$. 
\ethm

\sp

\bthm\label{Spher-I}{\rm (Koebe's Distortion Theorem I, Spherical
Version)}\index{(N)}{Koebe's Distortion Theorem I (spherical
version)}  
Given a number $s>0$ there exists a monotone increasing continuous function $k_s:[0,1)\to[1,\infty)$ such that for any $w\in \ov{\mathbb C}$, any $r >0$, all $t\in [0,1)$, and any univalent analytic function $H:B_n(w,r)\lra
\ov{\mathbb C}$ such that the complement $\ov{\mathbb C}\sms
H(B_n(z,r))$ contains a spherical ball of radius $s$, we have that 
$$
\sup\big\{|H^*(z)|:z\in B_n(w,tr)\big\} 
\le k_s(t) \inf\big\{|H^*(z)|:z\in B_n(w,tr)\big\}.
$$
where $B_n$ stands in here for either Euclidean 
or spherical ball (necessarily spherical  if $z=\infty$).
\ethm 
\bpf By rescaling, we may assume without loss of generality
that $r=1$. Let 
$$
M_1:\oc\to\oc
$$ 
be a M\"obius transformation (with spherical
derivative bounded above by $2$ and below by $1/2$)
sending $0$ to $w$ and let 
$$
B:=M_1^{-1}(B_n(w,1)),
$$
which is a
Euclidean ball centered at $0$ with finite radius bounded uniformly
above and below. Let 
$$
B_s(\xi,s)
$$ 
be a ball of radius $s$ disjoint from $H(B_n(w,1))$. Let
$$
M_2:\oc\to\oc
$$ 
be a M\"obius transformation sending $\xi$ to $\infty$
and $\oc\sms B_e(\xi,s)$ onto the unit ball $B_e(0,1)$. Then 
$$
M_2\circ H\circ M_1:B\lra B_e(0,1)\sbt\C
$$ 
is a univalent analytic function. Noting now that the Euclidean and
spherical derivatives $(M_2\circ H\circ M_1)'$ and $(M_2\circ H\circ
M_1)^*$ are universally uniformly comparable on $B$, and the spherical
derivatives of $M_1$ and $M_2$ are universally uniformly bounded above and
separated from zero, our theorem follows from Theorem~\ref{tKoebe-1}.
\epf

\sp

\fr Employing the Mean Value Inequality, the following two lemmas are
straightforward respective consequences of the latter two distortion theorems and Koebe's $\frac{1}{4}$--Theorem (Theorem~\ref{one-quater}).

\sp

\blem\lab{lncp12.9.} 
Suppose that $D\sbt {\mathbb C}$ is an open set, $z\in D$ and $H:D\lra {\mathbb C}$ is an analytic map which has an analytic inverse $H_z^{-1}$ defined on $B_e(H(z),2R)$ for some $R>0$. Then for every $0\le r\le R$
$$
B_e\lt(z,\frac14r|H'(z)|^{-1}\rt)\sbt H_z^{-1}(B_e(H(z),r))\sbt
B_e(z,Kr|H'(z)|^{-1}).
$$
\elem

\sp 
\fr and
\sp

\blem\lab{lncp12.9s.} Suppose that $D\sbt \ov{\mathbb C}$ is an open
set, $z\in D$ and $H:D\lra \ov{\mathbb C}$ is an analytic map which
has an analytic inverse $H_z^{-1}$ defined on $B_s(H(z),2R)$ for
some $R>0$ avoiding a spherical ball of some radius $s$. Then for
every $0\le r\le R$
$$
B_s\(z,k_s^{-1}(1/2)r|H^*(z)|^{-1}\) \sbt H_z^{-1}\(B_s(H(z),r)\) \sbt
B_s\(z,k_s(1/2)r|H^*(z)|^{-1}\).
$$
\elem

\fr Since conformal homeomorphisms preserve moduli of annuli, employing
the Riemann Mapping Theorem, we directly obtain from Theorem~\ref{Euclid-I} and
Theorem~\ref{Spher-I} the respective two more geometric versions of Koebe's
Distortion Theorems that involve moduli of annuli.

\bthm\label{Euclid-II}{\rm (Koebe's Distortion Theorem II, Euclidean
version)}\index{(N)}{Koebe's Distortion Theorem II (Euclidean
version)} There exists a function $w:(0,+\infty)\to [1,\infty)$ such
that for any two open topological disks $Q_1\sbt Q_2$ with
$\text{Mod}(Q_2\sms Q_1)\ge t$ and any univalent analytic function
$H:Q_2\to \ov{\mathbb C}$ we have
$$
\sup\{|H'(\xi)|:\xi\in Q_1\} \le w(t) \inf\{|H'(\xi)|:\xi\in Q_1\}.
$$
\ethm

{\sl Proof.} Let $R:B(0,1)\lra Q_2$ be a Riemann Mapping (conformal homeomorphism) such that $R(0)\in Q_1$. Then, by Theorem~\ref{t1j235}, $\Mod(R^{-1}(Q_2\sms \ov Q_1))\ge t$. Let
$$
r:=\sup\big\{|R^{-1}(z)|:z\in Q_1\big\}<1.
$$
Then by Theorem~\ref{t2j225} (Gr\"otzsch Module Theorem), $\Mod(R^{-1}(Q_2\sms \ov Q_1))\le \mu(r)$. Hence, $t\le \mu(r)$. Therefore, by virtue of Proposition~\ref{p2_2017_10_27}, there exists $s(t)<1$ such that
$$
r<s(t).
$$
Thus, applying Theorem~\ref{Euclid-I}, Koebe's Distortion Theorem I, Euclidean
Version, (in fact Theorem~\ref{th17.46in{H}} would directly apply), we get for all $\xi, z\in Q_1$ that
$$
\frac{|H'(\xi)|}{|H'(z)|}
=\frac{\big|\((H\circ\R)\circ R^{-1}\)'(\xi)\big|} {\big|\((H\circ\R)\circ R^{-1}\)'(z)\big|}
=\frac{\big|((H\circ\R)'(R^{-1}(\xi))\big|}{\big|((H\circ\R)'(R^{-1}(z))\big|}\cdot
 \frac{|R'(R^{-1}(z))|}{|R'(R^{-1}(z))|}
 \le k^2(s(t)).
$$
So, setting $w(t):=k^2(s(t))$ finishes the proof.
\endpf

\bthm\label{Spher-II}{\rm (Koebe's Distortion Theorem II, spherical
version).}\index{(N)}{Koebe's Distortion Theorem II (spherical
version)}
Given a number $s>0$ there exists a function 
$$
w_s:(0,+\infty)\lra
[1,\infty)
$$ 
such that for any two open topological disks $Q_1\sbt
Q_2$ with $\text{Mod}(Q_2\sms Q_1)\ge t$ and any univalent analytic
function $H:Q_2\lra \ov{\mathbb C} $ such that the complement
$\ov{\mathbb C}\sms H(Q_2)$ contains a ball of radius $s$, we have that
$$
\sup\big\{|H^*(\xi)|:\xi\in Q_1\big\} 
\le w_s(t) \inf\big\{|H^*(\xi)|:\xi\in Q_1\big\}.
$$
\ethm

\

\section{Local Properties of Critical Points of Holomorphic Functions}\label{Local Properties of Critical Points of Holomorphic Functions}

This short technical section provides a quite good, non commonly known, description of what is going on near critical points of holomorphic maps.

\sp Given an  analytic function $H$ defined  throughout an open set
$D\sbt {\mathbb C}$, we  put
$$\index{(N)}{critical points}\index{(S)}{$\Crit(f)$}
\Crit(H):=\{z\in D:H'(z)=0\}.
$$
Its image, $H(\Crit(H))$, is called the set of critical
values\index{(N)}{set of critical values} of $H$. Suppose now that
$c\in\Crit(H)$.
Then there exists $R=R(H,c)>0$ and $A=A(H,c)\index{(S)}{$A(f,c)$}\ge
1$ such that
$$
A^{-1}|z-c|^{p_c} \le |H(z)-H(c)| \le A|z-c|^{p_c}
$$
and
$$
A^{-1}|z-c|^{p_c-1} \le |H'(z)| \le A|z-c|^{p_c-1}
$$
for every $z\in \Comp(c,H,R)$, and that
$$
H(\Comp(c,H,R))=B_e(H(c),R), 
$$ 
where
$$
p_c=p(H,c)\ge 2
$$\index{(S)}{$p_c$} is the order, or multiplicity, of 
$H$ at the critical point $c$\index{(N)}{order of a critical
point}\index{(N)}{multiplicity of critical point}, called also the order, or multiplicity, of the critical point $c$. In particular
$$\Comp(c,H,R)\sbt B_e(c,(R/A)^{1/p_c}).$$  Moreover, by taking
$R>0$ sufficiently small, we can ensure that the two above
inequalities hold for every $z\in B_e(c,(R/A)^{1/p_c})$ and the ball
$B_e(c,(R/A)^{1/p_c})$ can be expressed as a union of $p_c$ closed
topological disks with smooth boundaries and mutually disjoint
interiors such that the map $H$ restricted to each of these
interiors, is injective.

\sp\fr In the sequel we will need  the following technical lemma proven in
\cite{U1} as Lemma~2.11.

\

\blem\lab{lncp12.11.} Suppose that an analytic map $Q\circ H:D\lra
\mathbb C$, a radius $R>0$ and a point $z\in D$ are such that
\begin{enumerate}
\item [(a)]
$$\Comp(H(z),Q,2R)\cap \Crit(Q)=\es$$
 \and
$$\Comp(z,Q\circ H,R)\cap \Crit(H)\ne\es.$$
\item [(b)] If $c$ belongs to the last intersection and
$$
\diam_e\(\Comp(z,Q\circ H,R)\)\le (AR(H,c))^{1/p_c}
$$
then
$$
|z-c|\le KA^2|(Q\circ H)'(z)|^{-1}R.
$$
\end{enumerate}
\elem

\bpf In view of Lemma~\ref{lncp12.9.}
$$
\Comp(H(z),Q,R)\sbt B_e(H(z),KR|Q'(H(z))|^{-1}).
$$
So, since $H(c)\in \Comp(H(z),Q,R)$, we get $$H(c)\in
B_e(H(z),KR|Q'(H(z))|^{-1}).$$ Thus, using this  and  (b) we obtain
$$\begin{aligned}
A^{-1}|z-c|^{p_c} &\le |H(z)-H(c)| 
              \le KR|Q'(H(z))|^{-1}\\
              & =KR|(Q\circ H)'(z)|^{- 1}|H'(z)| \\
&\le KR |(Q\circ H)'(z)|^{-1} A|z-c|^{p_c-1}. \end{aligned}
$$
So, $|z-c|\le KA^2|(Q\circ H)'(z)|^{-1}R$. 
\endpf

\section{Proper Analytic Maps and Degree}

In this short section we bring up, with proofs, the concept of (holomorphic) proper maps and of the degree of such maps. We follow closely Section~4.2 of Otto Forster's book \cite{Fo}. 

Let $X$ and $Y$ be locally compact Hausdorff topological spaces. A continuous  map $f : X\lra Y$  is said to be proper if and only if $f^{-1}(K)$ is compact for every  compact subset $K$ of $Y$. We shall  prove the following well--known result.

\bprop\label{p1ms.2.1} 
Every proper map between locally compact Hausdorff topological spaces is closed,  meaning that images of all closed sets are closed.
\eprop

{\sl Proof.} Let $X$ and $Y$ be locally compact Hausdorff topological spaces and let $f:X \to Y$ be a proper map. Let $ F \sbt X$ be a closed set. Let $w \in \ov{f(F)}$.   Since $Y$ is a $T_2$ locally compact  space, there exists $W \sbt Y$, an open  neighbourhood of $w$, such that $\ov{W}$ is compact. Now, if $V$
 is an arbitrary open neighbourhood of $w$, then so is $W \cap V$, whence
 $$ (f(F)\cap W)\cap V=f(F)\cap (W\cap V)\neq \es.$$
 Thus
 $$ w \in\ov{f(F)\cap W}\sbt \ov{f(F)\cap \ov{W}}.$$
 Since also $f(F)\cap \ov{W}=f(F\cap f^{-1}(\ov{W}))$ and $F \cap f^{-1}(W)$ is a compact set, we therefore get that
 $$ w \in\ov{f(F\cap f^{-1}(\ov{W}))}=  f(F\cap f^{-1}(\ov{W})) \sbt f(F).$$
 The proof is complete. 
 \qed 

\blem\label{lemma4.21Forster}
Suppose $X$ and $Y$ are locally compact Hausdorff topological spaces and $p:X \lra Y$ is a proper, discrete map. Then the following hold:

\begin{itemize}
\item [(a)]  For every point $y \in Y$ the set $p^{-1}(y)$ is finite.

\, \item [(b)]  If $y \in Y$ and $V$ is a neighbourhood  of $p^{-1}(y)$, then there exists a neighbourhood $U$  of $y$ with $p^{-1}(U) \sbt V$.

\, \item [(c)]  If $D\sbt X$ is a closed discrete set, then $p(D)$ is also discrete.
\end{itemize}
\elem

\fr{\sl Proof.}
The item (a) follows from the  fact that $p^{-1}(y)$ is  a compact discrete subset of $Y$. 

For item (b) note that $p(X\sms V)$ is closed by applying Proposition~\ref {p1ms.2.1} and as the set $X\sms V$ is closed. Also $y \notin p(Y\sms V)$. Thus $U:=Y\sms p(X\sms V)$ is open neighbourhood of $x$ such that $p^{-1}(U) \sbt V$.

Proving (c), fix $z\in D$. Since $D$ is closed and discrete and since, by (a), the set $p^{-1}(p(z))$ is finite, there exists $U$, an open neighborhood of $p^{-1}(p(z))$ such that
$$
U\cap D=D\cap p^{-1}(p(z)).
$$
Then, by Proposition~\ref{p1ms.2.1}, $p(X\sms U)$ is a closed set and $p(z)\notin 
p(X\sms U)$. So, $Y\sms p(X\sms U)$ is an open neighborhood of $p(z)$ and 
$$
p(z)\in \(Y\sms p(X\sms U)\)\cap p(D)
\sbt p(U\cap D)
\sbt p(D)\cap\{p(z)\}
=\{p(z)\}.
$$
Thus, 
$$
\(Y\sms p(X\sms U)\)\cap p(D)=\{p(z)\},
$$
proving that $p(D)$ is discrete.
\qed

\sp \bthm\label{theorem4.24Forster}
Suppose $X$ and $Y$ are Riemann surfaces and $f:X \lra Y$ is a proper non--constant holomorphic map. Then there exists a natural number $n\ge 1$ such that $f$ takes every value $w \in Y$, counting multiplicities, $n$ times, i.e 
$$
\sum_{z\in f^{-1}(w)}q_z(f)=n,
$$
where $q_z(f)$ is the order of $f$ at $z$; it is equal to $1$ if $z$ is not a critical point of $f$. 
\ethm

\bpf The set $\Crit(f)$ of branch (critical) points of the map $f:X\to Y$ is closed and discrete. Since the map $f$ is proper and discrete, $B:=f(A)$ is also closed and discrete respectively by Proposition~\ref{p1ms.2.1} and Lemma~\ref{lemma4.21Forster} (c). Let  
$$
Y':=Y\sms B,
$$
$$
X':=X\sms f^{-1}(B)\sbt X\sms A
$$  
Since, 
$$
f|_{X'}:X'\to Y'
$$
is a proper discrete (as holomorphic) covering (as a proper local homeomorphism between locally compact Hausdorff spaces) map and the space $X'$ is pathwise connected, it follows from Lemma~\ref{lemma4.21Forster} (a) that there exists a natural number $n\ge 1$ such that
\beq\lab{1_2018_01_9}
\#\(f|_{X'}\)^{-1}(y)=n
\eeq
for all $y\in Y'$. 

We are thus only left to deal with points $b \in B$. 
Then there exist mutually disjoint neighborhoods $U_x$ of respective points $x\in f^{-1}(b)$, and neighborhoods $V_x$ of $b$, $x\in f^{-1}(b)$, such that $V_x\sms \{b\}\sbt Y'$ and for every $w\in V_x\sms \{b\}$ the set 
$$
\#\(f^{-1}(w)\cap U_x\)=q_x(f).
$$ 
By Lemma~\ref{lemma4.21Forster} (b) we can find a neighborhood 
$$
V \sbt \bi_{x\in f^{-1}(b)}V_x
$$
of $b$ such that 
$$
f^{-1}(V) \sbt \bu_{x\in f^{-1}(b)}U_x.
$$
Then, for every point $w \in V\sms \{b\}\sbt Y'$, we have
$$
\#\(f^{-1}(w)\)=\sum_{x\in f^{-1}(b)}q_x(f).
$$ 
Combining this with \eqref{1_2018_01_9}, we thus get
$$
\sum_{x\in f^{-1}(b)}q_x(f)=n,
$$
and the proof is complete.
\qed

\sp\section{Riemann--Hurwitz Formula}\label{RHF}

In this section we present a relation between holomorphic maps, their degrees, critical points and topological structure of images and preimages. It is known as Riemann--Hurwitz Formula. This Formula has a long history and is treated in many textbooks on Riemann surfaces and algebraic geometry. It is usually formulated for compact Riemann surfaces. We do need bigger generality. We therefore provide here a complete self--contained exposition which suffices for our needs throughout the book. 

In particular, the topological structure, mentioned above, enters our considerations in terms of Euler characteristic. We make a brief introduction to it in a restricted scope and, probably, in the most elementary way; in particular we speak of triangulations only rather than to deal with a more general and flexible concept of CW--complexes. 

The results of this section, i.e. various versions of the Riemann--Hurwitz Formula and their corollaries, will be instrumental when dealing with elliptic functions, particularly the non--recurrent ones. These formulas provided an elegant and probably the best tool to control the topological structure of connected components of inverse images of open connected sets under meromorphic maps. Especially to be sure that such connected components are simply connected.

Our exposition stems in this section from that of Allan Beardon's book \cite{Bea}. 

\subsection{The Euler Characteristic of Plane Bordered Surfaces}\label{EC 1}

We devote this subsection to a brief and of restricted scope introduction to the Euler characteristic. Let $S\sbt\oc$ be a bordered connected surface (though its boundary (border) can be empty and then $S=\oc$), i.e. a closure of an open connected set in $\oc$ whose boundary $\bd S$ (in $\oc$) consists of a finite number of simple closed (Jordan) curves. A triangulation \index{(N)}{triangulation} $T$ of $S$ is a partition  of $S$ into a finite number of mutually disjoint subsets called vertices, edges, and faces, respectively denoted by $V$, $E$, and $F$, with the following properties:
\begin{itemize}
\item [(1)] each vertex, element of $V$, is a point of $S$;

\, \item [(2)] for each edge $e\in E$, there are a compact interval $[a,b]\sbt\R$ and a homeomorphism $\varphi:[a,b]\lra \ov e$ such that
$$
\varphi\((a,b)\)=e 
\  \  {\rm and} \ \ 
\varphi(\{a,b\})\sbt V;
$$
\item [(3)] for each face $f\in F$, there are a closed triangle $\De\sbt\C$ and a homeomorphism $\varphi:\De\lra \ov f$ which repspectively maps the edges and vertices of $\De$  (in the usual sense) to $E$ and $V$, and such that 
$$
\varphi(\Int\De)=f.
$$
\end{itemize}

\fr Of course if (2) or (3) holds then for any other closed interval or respectively a triangle, there would exist the required homeomorphism $\varphi$. We stress that $T$ partitions $S$ into mutually disjoint subsets of $S$. Each such subset is either a vertex, an edge or a face and we call each of theses a simplex of $T$ of dimension 0, 1 and 2 respectively. For any simplex $s$ of dimension $m$, the Euler characteristic $ \chi(s)$ of $s$ is defined to be
$$
(-1)^m.
$$
More generally, if $S'$ is any subset of $S$ comprising  a union of simplices, say $s_1, \ldots, s_k$, where $s_j$ has dimension $n_j$, we then call $S'$ a subcomplex of $S$ relative to $T$, and we define
\beq\lab{1_2017_12_12}
\chi(S',T):=\sum_{j=1}^k\chi(s_j)=\sum_{j=1}^k(-1)^{n_j}.
\eeq
In particular, if the triangulation $T$ of $S$ consists of faces $F$, edges $E$, and  vertices $V$, then the Euler characteristic $\chi(S)$ is, by definition
$$
\chi(S)=\chi(S,T):=\#F-\#E+\#V.
$$ 
The crucial fact well--known in algebraic topology (which we accept here without proof) is that $\chi(S,T)$ is independent of the particular triangulation $T$ used. In particular it is thus a topological invariant, meaning that two homeomorphic bordered surfaces have the same Euler characteristics. So, we can compute $\chi(S)$ by using any convenient triangulation we choose. Of course, we need to know for this, and we indeed do by obvious geometry, that every such surface admits a triangulation. 

For each edge $e\in E$, the closeure $\ov e$ is a subcomplex of $S$ relative to it and 
\beq\lab{2_2017_12_12}
\chi(\ov e,T)=1.
\eeq
The boundary $\bd S$ of $S$ is its subcomplex relative to any traingulation of $S$ and it is again evident that
$$
\chi(\bd S,T)=0,
$$
as every connected component of $\bd S$ is a Jordan curve which is a subcomplex of any triangulation of $S$ and whose Euler characteristic is evidently equal to zero. The latter formula above permits us to get rid of the triangulation $T$ and to express it in the following form:
\beq\lab{3_2017_12_12}
\chi(\bd S)=0.
\eeq
We say that the Euler characteristic of $\bd S$ is zero.

Calculations of $\chi$ can often be simplified by making use of the above and 
the following simple idea extending it. If $T$ is a triangulation of $S$ and $S_1$ and $S_2$ are some two subcomplexes of $S$ relative to $T$, then also $S_1\cup S_2$ and $S_1\cap S_2$ are subcomplexes of $S$ relative to $T$, and from the formula \eqref{1_2017_12_12}, we obtain
\beq\label{(5.3.1)}
\chi(S_1\cup S_2)+\chi( S_1 \cap S_2)= \chi(S_1) +\chi(S_2).
\eeq

Let us illustrate this with some simple examples. Let $\Int S=S\sms\bd S$ denote the  interior of $S$. The interior $\Int S$ of $S$ is its subcomplex relative to any triangulation of $S$. Then, by \eqref{(5.3.1)} and \eqref{3_2017_12_12}:
$$ 
\chi(S)
=\chi(\Int S,T)+\chi(\partial{S})
=\chi(\Int S,T). 
$$
In particular $\chi(\Int S,T)$ is again independent of $T$ and we can speak of $\chi(\Int S)$, calling it the Euler characteristic of $\Int S$. We thus have
\beq\lab{4_2017_12_12}
\chi(\Int S)=\chi(S).
\eeq
However, a warning. We have not established yet that $\chi(\Int S)$ is a topological invariant, i.e. we do not know that if $\hat S$ is another bordered surface in $\oc$, and $\Int\hat S$ is homeoemorphic with $\Int S$, then $\chi(\Int\hat S)=\chi(\Int S)$ (we know this though if $\hat S$ and $S$ are homeomorphic). We will address this issue later in this section.
 
By constructing  explicit triangulations it is immediate  that $\chi(\hat{\mathbb C})=2$, and that $\chi(D)=1$  for any closed topological disk $D$; these  computations are indeed simple but important. 

Now suppose  that $S$ is the complement in  ($\hat{\mathbb C}$) of $k$ mutually disjoint open topological disks $D_1, \ldots, D_k$, whose boundaries are Jordan curves, so  $D$ is  of connectivity $k$, meaning that its complemet has $k$ connected components. We can triangulate the sphere  such that each of the sets $S$, $D_1, \ldots, D_k$ is its subcomplex. Then (\ref{(5.3.1)}) yields
$$ 
2 =\chi(\hat{\mathbb C})
=\chi(S)+\sum_{j=1}^k \chi(D_j)
=\chi(S)+k.
$$
So,
\beq\label{1_2017_12_19}
\chi(S)=2-k.
\eeq
Note that for any such subdomain $S$ of the sphere $\oc$, we have that
\begin{itemize}
  \item [(a)] $\chi(S)=2$ if and only if $S$ is the sphere $\hat{\mathbb C}$ ($k=0$);
  
\, \item[(b)] $\chi(S)=1$ if and only if $S$ is simply connected, but not being $\hat{\mathbb C}$, i.e. $S$ is a closed topological disk, equivalently if $k=1$;
   
\, \item [(c)] $\chi(S)=0$ if and only if $S$ is doubly connected, i.e. if $k=1$;
\end{itemize}
while in all other cases, i.e. if $k\ge 3$, 
$$
\chi(S)<0.
$$ 

\subsection{The Riemann--Hurwitz Formula for Borded Surfaces in $\oc$} 

Our first main goal is to obtain and to prove the formula, commonly called the Riemann--Hurwitz formula, for the situation we start to describe now. Let $S_1$ and $S_2$ be two domains in $\oc$, i.e. open connected subsets of $\oc$. Let 
$$ 
R:S_1\lra S_2, 
$$
be an analtytic which has a (unique) continuous extension from $\ov S_1$ to $\ov S_2$. Keep the same symbol $R$ for this extension. Let $V$ be an open connected subset of $S_2$ such that $\ov V\sbt S_2$. Assume also that
\beq\lab{1_2017_12_14B}
\bd V\cap R(\Crit(R))=\es.
\eeq
Let $U$ be a connected component of $R^{-1}(V)$. 
Toward the end of later dealing with Euler characteristics and Riemann--Hurwitz Formula, we will establish some properties of the map $R$ in relation to $V$ and $U$ that will turn out to be also useful in further stages of the proof of the Riemann--Hurwitz formula. 

\sp\blem\label{l1_2017_12_14}
The map $R|_U:U\lra V$ is a proper analytic surjection with
\beq\lab{2_2017_12_14}
R(\bd U)=\bd V 
\  \  {\rm and} \  \
\ov U\cap R^{-1}(\bd V)=\bd U.
\eeq
\elem

\bpf  Since $R(U)\sbt V$ and since $V \cap \partial{V}=\es$, we have that

 \beq\label{1ms3}
 \ov{U}\cap R^{-1}(\partial{V})\sbt \partial{U}.
 \eeq
Seeking contradiction,  suppose  that
$$V \cap R(\partial{U}) \neq\es.$$
Then, there exists $z \in \partial{U}$ such that $R(z)\in V$. Since $R$ is continuous and $V$ is open, there thus exists $\delta>0$  such  that
$$B(z, \delta) \sbt S_1 \quad\text{and} \quad  R(B(z, \delta)) \sbt V.  $$
But then the set $U\cup B(z, \delta)$ is  connected (as $U\cap B(z, \delta)\neq \es$) and  $R(U \cup B(z, \delta))\sbt V$. It therefore  follows from the definition of $U$  that
 $$ U \cup B(z, \delta)\sbt U.$$
In particular, $ z \in U$ contrary to the facts that $U$ is open and $ z \in \partial{U}$. Hence
$$  V \cap R(\partial{U})=\es,$$
and, as $R(\partial{U})\sbt \ov{V}$, we must have that
 \beq\label{2ms3}
 R(\partial{U})\sbt \partial{U}.
 \eeq
Therefore, $\partial{U} \sbt R^{-1}(\partial{V})$, and as $\partial{U} \sbt \ov{U}$, we conclude that
$$\partial{U} \sbt \ov{U}\cap R(\partial{V}).$$
Along with (\ref{1ms3}), this entails the right--hand side  part of (\ref{2_2017_12_14}). Now,
\beq\label{1ms4}
\partial{R(U)} \sbt \ov{ R({U})}= R(\ov U),
 \eeq
the equality holding because $R(\ov{U})$ is closed as compact ($\ov{U}$ is compact as a closed subset of the  compact set $\ov{S}_1$).  Since $R(U)$ is open, we have that $\partial{R}(U)\cap R(U)=\es$, and invoking (\ref{1ms4}), we conclude that
$$\partial{R}(U)\sbt R(\partial{U}).$$
Along with (\ref{2ms3}) this implies that
 \beq\label{2ms4}
 \partial{R}(U)\sbt \partial{V}.
 \eeq
Now, seeking contradiction, suppose that
$$ V \sms R(U)\neq  \es.$$
Fix two points: $x \in  R(U)$ and $y\in  V \sms R(U).$ Since $V$ is arcwise connected there exists a homeomorphism
$$ 
h: [0,1] \lra V
$$
such that $h(0)=x$ and $h(1)=y$. Let
$$ 
s:= \inf\big\{ t\geq 0: \, h(t)\in V \sms R(U)\big\}. 
$$
Then,  since $\ov{\mathbb C}\sms R(U)$ is  closed and $h([0,1]) \sbt V$, and since $h(0)\in R(U)$, we have that
$$ R(s) \in (V \sms R(U)) \cap \ov{R(U)}= V\cap(\ov{\mathbb C} \sms R(U)) \cap \ov{R(U)}= V \cap \partial{R}(U),$$
contrary to (\ref{2ms4}) since $V \cap \partial{V}=\es$. Thus, $R: U \lra V$  is a surjection. This in turn implies that $R(\ov{U})= \ov{R(U)}=\ov{V}$.   As $R(U)$ is open, this yields $\partial{V} \sbt R(\partial{U})$. Along with  (\ref{2ms3}), this establishes the left-hand side of (\ref{2_2017_12_14}). We are thus  left to show only that $R:U\to V$ is a proper map. Indeed, let $K \sbt V$ be a  compact  set. Then $\ov{U}\cap R^{-1}(K)$ is  compact. We are to show that $U\cap R^{-1}(K)$ is compact. For this it thus  suffices to prove that  $\ov{U}\cap R^{-1}(K)=U\cap R^{-1}(K)$, and for this in turn, it suffices to show that
\beq\lab{1ms5}       
\partial{U} \cap R^{-1}(K)=\es.
\eeq
But, utilizing  (\ref{2ms3}), we get
$$ 
R(\partial{U}\cap R^{-1}(K))
\sbt  R(\partial{U}) \cap R(R^{-1}(K))
\sbt R(\partial{U})\cap K 
\sbt \partial{V} \cap  V 
=\es.
$$
Hence, (\ref{1ms5}) holds, and the proof is complete. \qed

\sp\fr Having this lemma proved, we may, by virtue of Theorem~\ref{theorem4.24Forster}, speak about $\deg\(R|_U\)\geq 1$, the degree of the map $R|_U:U\to V$. Now, wee shall  prove the following.

\sp\blem\label{l1ms5} 
If in addition $\ov V$ is a bordered surface in $\oc$, then $\ov U$ is a bordered surface too, i.e. $\partial{U}$ is a disjoint union of Jordan cures. In addition, $R|_{\partial{U}}:\partial{U}\lra \partial{V}$ is a covering map of finite degree.
\elem

\bpf  Because of Lemma~\ref{l1_2017_12_14} and formula (\ref{1_2017_12_14}) the map
$$
R|_{\partial{U}}:\partial{U}\lra \partial{V}
$$
is a locally homeomorphic surjection. Since $\partial{U}$ is a  compact set, $R|_{\partial{U}}:\partial{U}\lra \partial{V}$ is thus  a covering map of finite degree whose value we denote by $q$. If $\Gamma$ is a connected  component of $\partial{V}$ (we know that then $\Gamma$ is a Jordan curve), then $ R:R^{-1}(\Gamma) \to \Gamma$ is also a covering  map of degree $\leq q$. Let $L$ be a  connected  component  of $R^{-1}(\Gamma)$. We claim that $R(L)=\Gamma$. Otherwise there would exist a point $x \in L$ such that  $R(x) \in \partial_\Ga{R}(L)$. But, because of (\ref{1_2017_12_14B}) there exists an open topological arc $\alpha$ containing $x$ such that $R(\alpha)\sbt \Gamma$ and $R(\a)$ is an open topological arc containining $R(x)$. But  then $L\cup \alpha$ is connected  and $R(L\cup \alpha)\sbt \Gamma$. Hence $L\cup \alpha=L$, and therefore $R(x)\in R(L) \cup R (\alpha)$, yielding $R(x)\in \Int_\Ga R(L)$. This contradiction gives that
$$
R(L)=\Gamma.
$$
Thus, the map $R: L\to \Gamma$ is a covering surjection. In consequence it is at most  $q$--to--$1$ and  $R^{-1}(\Gamma)$ has at most $q$ connected components. We are thus left to show that $L$ is a Jordan curve. But $L$ is compact and $R: L\to \Gamma$ is locally homeomorphic, so $L$ is a compact connected $1$-dimensional topological manifold. It is well known that then $L$ is a Jordan curve. Another argument, less topological but more analytic, would be this. Fix a point $\xi \in  \Gamma$. Then the set  $L\cap R^{-1}(\xi)$  has at most $q$ elements  and each element $z \in R^{-1}(\xi)$ gives  rise to a local  inverse map $R^{-1}_z$ from  a sufficiently small neighborhood  of $\xi$ onto  a neighborhood  of $z$; in particular, $R^{-1}_z(\xi)=z$. Continuing  $R^{-1}_z$ analytically from $\xi$ to $\xi$ along $\Gamma$, we traverse a path $\Gamma_z$ in $L$ whose  other endpoint $\hat{z}$ belongs to $R^{-1}(\xi)$. The path $\Gamma_z$ is either a closed topological arc or a Jordan curve depending on whether $\hat{z}\neq z$ or $\hat{z}=z$. Furthermore, the map 
$$
R^{-1}(\xi) \ni z \longmapsto \hat{z} \in R^{-1}(\xi)
$$ 
is a bijection. We then continue analytically getting points $z, \hat{z}, \hat{\hat{z}}, \ldots $, until  with at most $q$ iterates we reach $z$ again. Then the consecutive closed topological arcs $\Gamma_z ,\Gamma_{\hat{z}}, \Gamma_{\hat{\hat{z}}}, \ldots$, have exactly one common endpoint, $ \hat{z}, \hat{\hat{z}},\ldots$, respectively. So, their union, up to  the second to the last element is a closed  topological arc again. Its union with the last closed  topological arc in the sequence
$\Gamma_z, \Gamma_{\hat{z}}, \Gamma_{\hat{\hat{z}}}, \ldots$, is a Jordan curve since these two arc have two common endpoints  and no other endpoints. The proof is complete \qed

\sp\blem\label{l1ms7}
With the hypotheses of Lemma~\ref{l1ms5} (in particular
$\ov V$ is a bordered surface in $\oc$), we have that
$$
\deg\(R|_{\partial{U}}\)=\deg(R|_U).
$$
\elem

\bpf Fix $\xi \in \partial{V}$. Since $\partial{V}$ contains no critical values of $R$,  there exists $\varepsilon >0$ such that all analytic branches of $R^{-1}$ are defined on $B(\xi, \varepsilon)$.  Since, by Lemma~\ref{l1ms5}, $\partial{U}$ is a finite disjoint union of Jordan curves, and since for every $z \in \partial{U}\cap R^{-1}(\xi)$, the map $R:R^{-1}_z(B(z,\varepsilon))\lra B(z,\varepsilon)$ is a homeomorphism, there exists $\Gamma_z \sbt R^{-1}_z(B(z, \varepsilon))\cap \ov{U}$, an open  relative to $\ov{U}$,  neighbourhood of  $z $ in $\ov{U}$  such that $R(\Gamma_z)$ is an open, relative to $\ov{V}$, neighbourhood of $\xi$ in $\ov{V}$. Hence
$$ V_\xi:= \bigcap \{ R(\Gamma_z):\, z \in \partial{U}\cap R^{-1}(\xi)\}$$
is also an open, relative to $\ov{V}$, neighbourhood of $\xi$ in $\ov{V}$. Fix $y \in V_\xi\cap V$. Since both families
$$ {\mathcal F}_\xi := \{ R^{-1}_z(B( \xi, \varepsilon)): \,  z \in \partial{U} \cap R^{-1}(\xi)\}$$
 and
$$ 
{\mathcal F}_y := \{ R^{-1}_x(B( \xi, \varepsilon)):\,  x \in U\cap R^{-1}(y)\}
$$
consist of mutually disjoint sets, we have that
$$
\deg(R|_{\partial U})=\#{\mathcal F}_\xi  \quad \text{and} \quad  \deg(R|_U)=\#{\mathcal F}_y.
$$
Thus, in order to complete the proof, it suffices to show that
\beq\label{1ms7}
{\mathcal F}_\xi={\mathcal F}_y.
\eeq
Indeed, if $z \in \partial{U} \cap R^{-1}(\xi),$
then there exists a (unique) point $x \in \Gamma_z$ such that $R(x)=y$. Furthermore, $ x \in \ov{U}\cap R^{-1}(V) =U$, so $x\in U\cap R^{-1}(y)$, and $
R^{-1}_x(y)=x=R^{-1}_z(y)$.  Thus $R^{-1}_z(B(\xi,\varepsilon))=R^{-1}_x(B(\xi, \varepsilon)) \in  {\mathcal F}_y$,  and the inclusion
\beq\label{2ms7}
 {\mathcal F}_\xi \sbt {\mathcal F}_y
\eeq
is proved.
For the proof of the opposite inclusion, suppose that $x \in U \cap R^{-1}(y)$. Then
$R^{-1}_x(V_\xi\cap V)\sbt U$, whence
$$
z:= R^{-1}_x(\xi)\in R^{-1}_x(\ov{V_\xi\cap V})
= \ov{R^{-1}_x(V_\xi\cap V)}\sbt \ov{U}.
$$
Therefore, $z \in \partial{U}\cap R^{-1}(\xi)$ and $ R^{-1}_z(\xi)=z=R^{-1}_x(\xi)$. Thus,  $R^{-1}_z(B(\xi, \varepsilon))=R^{-1}_x(B(\xi, \varepsilon)).$
So,  $R^{-1}_x(B(\xi, \varepsilon))\in {\mathcal F}_\xi$, and the inclusion $
{\mathcal F}_y \sbt {\mathcal F}_\xi$ is proved. Along with (\ref{2ms7}), this completes the proof of formula (\ref{1ms7}), and, simultaneously, the proof of Lemma~\ref{l1ms7}.
\endpf

\sp Now we turn to prepare the appropriate data from $R$. It will involve a contribution from the critical points of $R$ and in order to quantify this, we introduce the {\it deficiency} of $R$ at a point $z$ belonging to $S_1$ as
\beq\label{(5.4.3)}
\delta_R(z):=p_z-1,
\eeq
where, we recall, $p_z$ is the order of $z$ with respect to $R$; if $z$ is not a critical point of $R$, the $p_z=1$ and $\d_R(z)=0$.
For any set $A\sbt S_1$ we  define the {\it total  deficiency} of $R$  {\it over} $A$ as
$$ 
\delta_R(A):=\sum_{z \in A} \delta_R(z).
$$
The function $A\longmapsto\delta_R(A)$ is additive, meaning that for disjoint sets $A$ and $B$ in $S_1$ we have
$$ 
\delta_R (A\cup B)=\delta_R(A)+\delta_R(B).
$$
Frequently, but only when $R$ is  understood, we omit the suffix $R$ and use $\delta$ instead of $\delta_R$. 

 We  are now ready to relate the quantities $\chi(U)$, $\chi(V)$, $\mbox{deg}\(R|_U\)$ and $\delta_R(U)$. We shall prove the following.

\bthm[Riemann--Hurwitz Formula for Bordered Surfaces in $\oc$]\label{Riemann--Hurwitz Formula}
Let $S_1$ and $S_2$ be two domains in $\oc$, i.e. open connected subsets of $\oc$. Let 
$$ 
R:S_1\lra S_2, 
$$
be an analtytic which has a (unique) continuous extension from $\ov S_1$ to $\ov S_2$. Let $V$ be an open connected subset of $S_2$ such that $\ov V\sbt S_2$ and $\ov V$ is a bordered surface in $\oc$. Assume also that
\beq\lab{1_2017_12_14}
\bd V\cap R(\Crit(R))=\es.
\eeq
If $U$ is a connected component of $R^{-1}(V)$, then
\beq\label{(5.4.4)}
\chi(U)+\delta_R(U)=\deg\(R|_U\)\chi(V).
\eeq
\ethm

{\sl Proof.} We triangulate  the closure of $V$ ensuring (as we may) that all   critical values of $R$ in $V$  are vertices of the triangulation. Indeed, since each connected component on $\bd V$ is a Jordan curve, given any triangulation, we can always connect with a closed topological arc each critical value of $R$ with at least two distinct vertices of $T$, bounding a $1$--dimensional simplex on one connected componnent of $\bd V$, and to form in this way a larger triangulation with the required property. Denote such triangulation by $T$, and, as always, it vertices, edges, and faces respectively by $V$, $E$ and $F$.

Now we construct a triangulation $T_U$ of $\bd U$ in the following way. Denote $\deg\(R|_U\)$ be $m$. Let $\De\in F$. Since $\De$ is a topological closed disk, there exists an open topological disk $\De'$ (whose  closure is a topological closed disk) containg $\ov\De\sms V$ and disjoint from the union of critical values of $R$ and vertices of $\De$. Then, using the Monodromy Theorem, we see
that there are  exactly $m$ distinct analytic branches $R_j^{-1}:\De'\to S_1$ (so, $j=1, 2, \ld, m$) of $R^{-1}$ such that $R_j^{-1}(\De)\sbt U$. We define 
$$
F_U:=\big\{R_j^{-1}(\De):\De\in F,\, j=1, 2, \ld, m\big\}.
$$
These are to be the faces of the ultimate triangulation $T_U$. The edges of every face $R_j^{-1}(\De)$ are to be the sets $R_j^{-1}(e)$, where $e$ is, one of the three, edges of $\De$. So, we define
$$
E_U:=\big\{R_j^{-1}(e):\De\in F,\, e\in \bd\De,\, j=1, 2, \ld, m\big\}.
$$
Finally, the set $V_U$ consists of the endpoints of the closures of all elements of $E_U$. Note that the endpoints of the closure of each element $e$ of $E_U$ are distinct since these are mapped by $R$ onto two different points, namely the endpoints of $R(e)in E$. It follows immediately from this construction that the sets $F_U$ labeled as faces, $E_U$ labeled as edges, and $V_U$ labeled as vertices, form a triangulation of $\ov U$. We denote it by $T_U$.

It is immediate from this construction that 
$$
\#F_U=m\#F, \  \  \  \#E_U=m\#E
$$
and 
$$
\#V_U= m\#V -\sum_{c\in V_U\cap\Crit(R)}(p_c-1)=m\#V-\delta_R(U).
$$
Therefore,
$$\begin{aligned}
\chi(U)
&=\chi(U,T_U)
=\#F_U-\#E_U+\#V_U
= m\#F - m\#F+(m\#V-\delta_R(U))\\
& =m(\#F-\#E+\#V)-\delta_R(U)\\
&=m\chi(V)-\delta_R(U).
\end{aligned}
$$
The proof is complete.
\qed

\sp\fr Now we shall derive a bunch of fairly straightforward consequences of this theorem.

\bcor\label{c1_2017_12_22}
With the hypotheses of Theorem~\ref{Riemann--Hurwitz Formula}, if 
$V$ is simply connected, i.e. conformally equivalent to the unit disk $\D$, then
$$
\delta_U(R))\geq \deg(R|_U)-1,
$$
with equality holding if and only if $U$ is simply connected.
\ecor

\bcor\label{Riemann Hurwitz Simply Connected_1}
With the hypotheses of Theorem~\ref{Riemann--Hurwitz Formula}, if both $U$ and $V$ are simply connected, i.e. conformally equivalent to the unit disk $\D$, then
\beq\label{1_2017_10_12}
\d_U(R)=\deg(R|_U)-1.
\eeq
In particular,
\beq\label{2_2017_10_12B}
\#\(\Crit(R|_U)\)\le \deg(R|_U)-1. 
\eeq
\ecor

\bcor\label{Riemann Hurwitz Simply Connected_2}
With the hypotheses of Theorem~\ref{Riemann--Hurwitz Formula}, if
$V$ is simply connected , i.e. conformally equivalent to the unit disk $\D$, and $U$ contains no critical points of $R$, then  $R|_U:U\to V$ is a conformal homeomorphism. In particular, $U$ is conformally equivalent to the unit disk $\D$.
\ecor

{\sl Proof.}
Since $U$ is not compact (as otherwise $V$ would be compact), $\chi(U)\le 1$. We therefore conclude from Theorem~\ref{Riemann--Hurwitz Formula} that $\deg(R|_U)\le 1$. Thus $\deg(R|_U)=1$, whence $R|_U:U\to V$ is a conformal homeomorphism. This in turn immediately yields the second assertion of our corollary.
\endpf

\bcor\label{Riemann Hurwitz Simply Connected_3}
With the hypotheses of Theorem~\ref{Riemann--Hurwitz Formula}, if $V$ is simply connected ,i.e. conformally equivalent to the unit disk $\D$, and the map $R|_U:U\lra V$ has only one critical point, denote it by $c$, then

\begin{itemize}
\item[(a)]
the Riemann surface $U$ is simply connected, i.e. conformally equivalent to the unit disk $\D$, and 

\, \item[(b)]
$$
\deg\(R|_U\)=p_c.
$$
\end{itemize}
\ecor

{\sl Proof.}
Obviously $\deg\(R|_U\)\ge p_c$. Since $U$ is orientable and not compact (as otherwise $V$ would be compact), $\chi(U)\le 1$. We therefore conclude from Theorem~\ref{Riemann--Hurwitz Formula} that $\deg\(R|_U\)\le p_c$. Thus $\deg\(R|_U\)=p_c$, proving item (b). Applying Theorem~\ref{Riemann--Hurwitz Formula} once more, it now follows that $\chi(U)=1$, and item (a) is also proved.
\endpf

\subsection{Euler Characteristic: The General $\oc$ Case}\label{EC 2}

The Euler characteristic of a domain $D$ in the complex sphere has been defined in Subsection~\ref{EC 1} whenever its boundary $\partial{D}$ consists of a  finite number of Jordan curves. In general, however, the boundary of a domain $D$ is much more complicated than this; we may not be able triangulate  the closure of $D$ and in these circumstances,
$\chi(D)$ is as  yet undefined. As this is likely to be so in the case of major interest to us (when $D$ is a component of a  Fatou set), this present us with a problem which we must now tackle. We  propose to show that  given  any domain $D$, we can define $\chi(D)$ as  the limiting value of the Euler characteristic of smooth subdomains  which exhaust $D$ and, once this has been done, we can then use the Euler characteristic as a tool to study the way in which $R$  maps one component of the Fatou set onto another. We shall restrict our discussion  to subdomains of $\ov{\mathbb C}$; nevertheless. The following development is closely related to  the construction of the ideal boundary components of a Riemann  surface.

    Let $D\sbt\oc$ be a domain, i.e. an open connected set. A subdomain $\De$ of $D$ is said to be a {\it regular subdomain} of $D$ if:
\begin{itemize}
    \item [(1)]  $\ov\De$ is a bordered surface in $D$, i.e. $\ov\De\sbt D$ and $\bd\De$ is a finite  union of  mutually disjoint Jordan curves, say $ \gamma_1, \ldots, \gamma_n$, and
        
\, \item [(2)] $\oc\sms\De$ consists of $n$ closed topological disks, say $\Ga_1, \ldots, \Ga_n$ respectively bounded  by $ \gamma_1, \ldots, \gamma_n$, and $\Ga_j\cap(\oc\sms D)\ne\es$ for every$j=1,2,\ld,n$. 
\end{itemize}

\, For example, the unit opemn disk $\D=\{z\in\C:|z| <1\}$ is  a regular  subdomain of $\mathbb C$, whereas the annulus $\{z\in\C: 1< |z| <2\}$ is not. 

A crucial observation is that, according to the previous section, the Euler characteristic $\chi(\De)$ is defined  for each  regular subdomain  $\De$  of $D$.
        
Of course, if a  subdomain $\De$ of $D$ satisfies (1) but not (2), we can  adjoint to $\De$ those sets $\Ga_j$ which do not intersect the complement of  $D$ to form a regular subdomain of $D$, which we denote by $\De_*$. Obviously,  $\chi(\De_*)\geq \chi(\De)$; in fact $\chi(\De_*)=\chi(\De)+k$, where $k\ge 0$ is the number of ajoined sets disks $\Ga_j$.

We want to consider $D$ as  the limit  of regular  subdomains,  and as   no  canonical sequence  of subdomains of $D$ presents itself, it is best to reject the idea of a sequential limit  and to consider instead convergence with respect  to the directed set (or net) of regular subdomains. There is no need for great generality here and the  details  are quite simple and explicit. First, we shall prove the following.

\blem\label{lemma5.5.1}
If $D$ is a proper subdomain of $\oc$, then
\begin{itemize}
\item [(a)] any compact subset of  $D$ lies in some regular subdomain od $D$;  and

\,\item [(b)] if $\De_1$ and $\De_2$ are regular subdomains of $D$, then there is a regular subdomain $\De$ of $D$ which contains $\De_1\cup\De_2$.
\end{itemize}
\elem

\bpf If $\xi$ is an, arbirarily chosen point in $D$, then the map $\oc\ni z\mapsto (z-\xi)^{-1}\in\oc$ sends $\xi$ to $\infty$, so we may assume without loss of generality that $\infty\in D$. Let $n$  be a positive  integer, and cover the plane with a square grid (including the axes), each square having diameter $1/n$.  

Now let $K_n$  be the union of those closed  squares in the grid that intersect $\oc\sms D$. If $n\ge 1$ is large enough, say $n\ge q$, then $\infty\notin K_n$, and let $D_n$ be the connected component of  $\oc\sms K_n$ that contains $\infty$. Then it is easy to see that $\(D_n\)_{n=q}^\infty$ is an ascending sequence of regular subdomains of $D$  whose union $\bu_{n=q}^\infty D_n$ is $D$.

Having this, the rest of the proof is straightforward. Given any compact subset $K$ of $D$, the family $\{D_n\}_{n=q}^\infty$ is an open cover of $K$ and so is covered  by a finite collections of the sets $D_n$. As the sequence $\(D_n\)_{n=q}^\infty$ is ascending, this  finite collections contains a largest domain, say $D_m$, and then $D_m$ contains $K$. This proves (a). 

Finally, (b) follows from (a) for if $\De_1$ and $\De_2$ are regular subdomains of $D$, then $\ov\De_1\cup\ov\De_2$ is a  compact subset of $D$ and so by (a), it lies in some regular subdomain of $D$. 
\qed

\sp Lemma~\ref{lemma5.5.1} says that the class $\mathcal R (D)$ of all regular subdomains of $D$ is a net with the order given by the (direct) inclusion relation.

\sp Our next task is to show that the function
$$
\mathcal R (D)\ni\De\longmapsto\chi(\De)
$$ 
is monotone decreasing, and so, tends to a limit, which may be $-\infty$. We shall prove the following.

\blem\label{lemma5.5.2} 
The Euler characteristic function
$$
\mathcal R (D)\ni\De\longmapsto\chi(\De)
$$
is decreasing. Explicitly, if $\De_1$ and $\De_2$ are regular  subdomains of $D$ such that $\De_1\sbt \De_2$, then 
$$
\chi(\De_2) \leq \chi(\De_1).
$$
\elem

\fr{\sl Proof.}
Let $V_1, \ldots, V_m$ be the connected components of $\oc\sms\De_1$, and let $W_1, \ldots, W_n$ be the connected components of of $\oc\sms\De_2$. Since $\De_1 \sbt \De_2$, we have
$$ 
W_1\cup\ldots \cup W_n \sbt V_1\cup \ldots \cup V_m.
$$
For each $j\in \{1, \ldots, m\}$ fix a point $z_j\in V_j\sms D$. As $z_j$ is not in $D$, it lies in some set $W_k$, $1\le k\le n$, and so $W_k$ (being connected) lies in $V_j$. It follows that each set $V_j$, $j\in \{1, \ldots, m\}$, contains some set $W_k$ with $k\in \{1, \ldots, n\}$. Hence $m \leq n$. The given inequality now follows since, using \eqref{1_2017_12_19}, we get that
$$ 
\chi(\De_2)= 2- n \leq 2-m= \chi(\De_1).
$$
The proof is complete.
\qed
 
\sp Lemma~\ref{lemma5.5.2} says that the  function 
$$
\chi:\mathcal R(D) \lra \big\{2, 1, 0, -1, -2, \ldots, -\infty\big\},
$$
defined on the net $\mathcal R(D)$, is monotone decreasing. Thus, we have the following.

\bprop\label{defn5.5.3} 
For any subdomain $D$ of $\oc$, the limit 
$$
\chi(D):=\lim_{\De\in \mathcal R(D)}\chi(\De)
$$
exists and
$$
\chi(D)=\inf\big\{\chi(\De): \De\in \mathcal R(D)\big\}.
$$
\eprop

\sp\fr Quite explicitly, either:
\begin{itemize}
\item [(1)] $ \chi(D)=-\infty$, and there are regular subdomains $\De_n$ of $D$ with $\chi(\Omega_n) \lra -\infty$ or, equivalently, with the connectivity of $\De_n$ diverging to $+\infty$; or else
    
\item [(2)] there is some regular subdomain $D_*$ of $D$ such that
$$
\chi(D_*)=\chi(D) > -\infty, 
$$
and then (from Lemma~\ref{lemma5.5.2}) $\chi(\De)=\chi(D)$ whenever $\De$ is a regular subdomain which contains $D_*$. 
\end{itemize}

If $D$ is simply connected  domain, then $\partial{D}$
is connected  and each regular subdomain  $\De$ can only  have  one  complementary component: thus  $\chi(D)=1$  for a simply connected  domain $D$, regardless  of the nature  of $\partial{D}$. More generally, if  $D$ has  connectivity $k$, then  $\chi(D)=2-k$  for  al sufficiently  large  regular subdomains $\De$, and so 
\beq\lab{2_2017_12_19}
\chi(D)=2-k,
\eeq
again  irrespective  of the  complexity  of $\partial{D}$. In particular, the current definition of the Euler characteristic of an open domain in $\oc$ coincides with the one of formula \eqref{4_2017_12_12} in the case $D$ is a bordered surface. Also w now know the the Euler characteristic of open domains is a topological invariant. 

\subsection{Riemann--Hurwitz Formula: The General $\oc$ Case} 
The main theorem of this subsection and the entire Section~\ref{RHF} is the following.

\bthm[General Riemann--Hurwitz Formula]\label{General Riemann--Hurwitz Formula}
Let $S_1$ and $S_2$ be two domains in $\oc$, i.e. open connected subsets of $\oc$. Let 
$$ 
R:S_1\lra S_2, 
$$
be an analytic function which has a (unique) continuous extension from $\ov S_1$ to $\ov S_2$. Let $V$ be an open connected subset of $S_2$. If $U$ is a connected component of $R^{-1}(V)$, then
\beq\label{(5.4.4)B}
\chi(U)+\delta_R(U)=\deg\(R|_U\)\chi(V).
\eeq
\ethm

{\sl Proof.} 
Fix  a point $w \in V \sms R(\Crit(R))$. Because of Lemma~\ref{lemma5.5.1}(a) there exists $\Delta_0$, a regular subdomain of $U$ such that
\begin{itemize}
\item [(1)] $\Crit (R) \cap U \sbt \Delta_0$ and

\,\item [(2)] $\ov{U} \cap R^{-1}(w) \sbt \Delta_0$.
\end{itemize}
Again because of Lemma~\ref{lemma5.5.1} (a) and since $R(\ov{\Delta}_0)$ is  a  compact  subset of $V$, there exists a regular subdomain $\Delta_1$ of $V$ such that
\beq\label{1ms8}
R(\ov{\Delta}_0) \sbt  \Delta_1.
\eeq
We shall prove the following

\

\fr{\bf Claim 1.}$$ \Delta_2:=U \cap R^{-1}(\Delta_1)$$
{\it is a regular subdomain of $U$ and $ \Delta_0\sbt \Delta_2$.}

\

\bpf  Since $\Delta_0$ is  connected and since, by (\ref{1ms8}), $R(\Delta_0) \sbt \Delta_1$, there exists $\Delta_2^{'}$, a unique connected component of $R^{-1}(\Delta_1)$ such that
\beq\label{2ms8}
{\Delta}_0 \sbt  \Delta_2'.
\eeq
 Since every connected component of $R^{-1}(\Delta_1)$ is either disjoint  from $U$ or is contained in $U$ and since any such component $\Delta_2^{''}$ contained in $U$, intersects $R^{-1}(w)$ (as  by Lemma~\ref{l1_2017_12_14}, $R(\Delta_2^{'})=\Delta_1$) it follows from (2) that $\Delta_2^{''}\cap \Delta_0\neq \es$. But then, by (\ref{2ms8}),
 $\Delta_2^{''}\cap \Delta_2^{'}\neq \es$. Consequently $\Delta_2^{''}= \Delta_2^{'}$. We have thus proved that
 $$    \Delta_2^{'}= U \cap R^{-1}( \Delta_1).         $$
 So, $\Delta_2$ is a subdomain of $U$ and, invoking (\ref{2ms8}), $\Delta_0\sbt \Delta_2$. Furthermore, since, by Lemma~\ref{l1ms5}, $\Delta_2^{'}$ is a bordered surface, and since $\Delta_2 =\Delta_2^{'}$, we get that $\Delta_2$ is a bordered surface. We are thus left to  show   that the subdomain $\Delta_2$  of $U$ is regular. Indeed, let $W$  be a connected component of $ \ov{\mathbb C}\sms \Delta_2$. Seeking contradiction suppose that $W \cap (\ov{\mathbb C}\sms U)=\es$. This means that
 \beq\label{1ms9}
W \sbt U.
\eeq
We know that $R(W)$ is a connected set contained in $ V\sms  \Delta_1$. Let $W_1$ be the connected component of $V \sms \Delta_1$ containing $R(W)$ and let $W_2$ be the connected component of $\ov{\mathbb C}\sms \Delta_1$ containing $W_1$. In particular $W_2$ is a  connected  component of
 $\ov{\mathbb C}\sms \Delta_1$, and therefore, since $\Delta_1$ is a regular subdomain of $V$, we  have that
 \beq\label{2ms9}
W_2 \cap (\ov{\mathbb C}\sms V)\neq \es.
\eeq
 Since $W_1\sbt V \cap W_2$, we have that
 $$   \partial{W}_1 \sbt \partial{V}\cup \partial{W}_2.   $$
 Seeking contradiction, suppose that
 $$ \partial{W}_1\cap \partial{V}=\es.$$
 Thus then
 $$ \partial{W}_1\sbt  \partial{W}_2.$$
 Hence $W_2=W_1\cup (W_2\sms \ov{W_1})$. Since $W_2$ is connected and $W_1\neq \es$, this implies that $W_2\sms\ov{W}_1=\es$. Equivalently  $W_2\sbt \ov{W}_1$. Since both $W_2$ and $W_1$ are open, this yields $W_2 \sbt W_1$. Thus, $ W_2=W_1$, whence $W_2\sbt V$, contrary to (\ref{2ms9}). We have thus proved that
 \beq\label{1ms10}
\partial{W}_1\cap \partial{V}\neq \es.
\eeq

\

Now, we  shall prove the following.

\

\fr{\bf  Claim~2.} {\it $W$ is  a connected  component of $R^{-1}(W_1)$ contained in $U$.}

\

\fr {\sl Proof.} $W \sbt U$ by (\ref{1ms9}). Furthermore,
$$W\sbt U\cap R^{-1}(W_1)\sbt U\cap R^{-1}(V \sms \Delta_1)= U\sms \Delta_2\sbt \ov{\mathbb C}\sms \Delta_2.$$
 Therefore, since $W$  is a connected   component of $\ov{\mathbb C}\sms \Delta_2$, it is also a connected component of $U\cap R^{-1}(W_1)$. Since each connected component
of $R^{-1}(W_1) $ is either contained in $U$ or is disjoint from $U$, $W$ is thus a connected component of $R^{-1}(W_1)$ contained in $U$. Claim~2 is  proved. \qed

\sp\fr It follows from this claim and Lemma~\ref{l1_2017_12_14} that $R(\partial{W})=\partial{W}_1$. Thus, invoking  (\ref{1ms10}), we get that $R(\partial{W})\cap \partial{V}\neq \es$.
Hence
 \beq\label{2ms10}
 R(\partial{W})\cap (\ov{\mathbb C} \sms {V})\neq \es.
 \eeq
But $\partial{W}\sbt\partial{\Delta}_2 \sbt \ov{U}\cap \ov{\Delta_2}\sbt  \ov{U}\cap R^{-1}(\ov{\Delta_1})\sbt \ov{U}\cap R^{-1}(V) \sbt U$. So,  $R(\partial{W})\sbt R(U)=V$, contrary to (\ref{2ms10}). The proof of Claim~1 is complete. \qed

\sp\fr Now, since, by Claim~1, $ \Delta_0\sbt \Delta_2\sbt U$, it follows  from Lemma~\ref{lemma5.5.2} that
  \beq\label{1ms13}
 \chi(\Delta_0)\geq \chi(\Delta_2)\geq \chi(U),
  \eeq
 and from (1) of  this proof that
  \beq\label{2ms13}
 \delta_R(\Delta_0)=\delta_R(\Delta_2)=\delta_R(U).
  \eeq
Because of Claim~1, Theorem~\ref{Riemann--Hurwitz Formula} is applicable to give, with the use of items (2) and (1) of this proof, along with \eqref{1ms13}, that
\beq\label{3ms13}
\deg(R|_U\)\chi(\Delta_1)
= \deg(R|_{\Delta_2})\chi(\Delta_1)
= \chi(\Delta_2)+\delta_R(\Delta_2)
= \chi(\Delta_2)+\delta_R(U)
\geq \chi(U)+\delta_R(U).
\eeq
Now,  the only requirement on $\Delta_1$ is (\ref{1ms5}), so  because of Lemma~\ref{lemma5.5.1} and Lemma~\ref{lemma5.5.2}, the infinimum  of $\chi(\Delta_1)$ over all regular subdomains $\Delta_1$ of $V$ is $\chi(V)$. Therefore,
  \beq\label{4ms13}
  \deg\(R|_U\)\chi(V) \geq \chi(U) +\delta_R(U).
    \eeq
On the other hand,  (\ref{3ms13}) and  (\ref{1ms13}) along with Lemma~\ref{lemma5.5.2}, give
\beq\label{5ms13}
\deg\(R|_U\)\chi(V) 
\leq \deg\(R|_U\)\chi(\Delta_1)
=\chi(\Delta_2)+\delta_R(U)
\leq \chi(\Delta_0)+\delta_R(U).
\eeq
Since the only requirements on $\Delta_0$ are (1) and (2), because of Lemma~\ref{lemma5.5.1} and Lemma~\ref{lemma5.5.2}, we see that the supremum of $\chi(\Delta_0)$ over all such domains $\Delta_0$ is $\chi(U)$. Therefore (\ref{5ms13}) yields
$$
\deg\(R|_U\) \chi(V) \leq \chi(U)+\delta_R(U).
$$
Along with (\ref{4ms13}) this completes the proof of Theorem~\ref{General Riemann--Hurwitz Formula}. \qed

\sp
\brem 
An important straightforward observation is that all the corollaries, Corollary~\ref{c1_2017_12_22} through Corollary~\ref{Riemann Hurwitz Simply Connected_3} hold true under the, weaker, hypotheses of 
Theorem~\ref{General Riemann--Hurwitz Formula}, rather than merely those of 
Theorem~\ref{Riemann--Hurwitz Formula}. In particular no Jordan curves needed nor critical values off the boundaries. The proofs are the same. We list all these corollaries below for the convenience of the reader and ease of reference.
\erem

\bcor\label{c1_2017_12_22_G}
With the hypotheses of Theorem~\ref{General Riemann--Hurwitz Formula}, if 
$V$ is simply connected, then
$$
\delta_U(R))\geq \deg(R|_U)-1,
$$
with equality holding if and only if $U$ is simply connected.
\ecor

\bcor\label{General Riemann Hurwitz Simply Connected_1}
With the hypotheses of Theorem~\ref{General Riemann--Hurwitz Formula}, if both $U$ and $V$ are simply connected and non of them is conformally equivalent to $\oc$, then
\beq\label{1_2017_10_12B}
\d_U(R)=\deg(R|_U)-1.
\eeq
In particular,
\beq\label{2_2017_10_12}
\#\(\Crit(R|_U)\)\le \deg(R|_U)-1. 
\eeq
\ecor

\bcor\label{General Riemann Hurwitz Simply Connected_2}
With the hypotheses of Theorem~\ref{General Riemann--Hurwitz Formula}, if
$V$ is conformally equivalent to the unit disk $\D$ and $U$ contains no critical points of $R$, then  $R|_U:U\to V$ is a conformal homeomorphism. In particular, $U$ is conformally equivalent to the unit disk $\D$.
\ecor

\bcor\label{General Riemann Hurwitz Simply Connected_3}
With the hypotheses of Theorem~\ref{General Riemann--Hurwitz Formula},  if $V$ is conformally equivalent to the unit disk $\D$ and if
$R|_U:U\lra V$ has only one critical point, denote it by $c$, then

\begin{itemize}
\item[(a)]
the Riemann surface $U$ is is conformally equivalent to the unit disk $\D$, and 

\, \item[(b)]
$$
\deg\(R|_U\)=p_c.
$$
\end{itemize}
\ecor

\chapter{Conformal Measures and Hausdorff Dimension \\ of \\ Invariant Measures}\label{CMHDIM}

In this chapter we encounter for the first time conformal measures, a beautiful, elegant and powerful tool of conformal dynamics due to S. Patterson (\cite{Pat1}, \cite{Pat2})and D. Sullivan (\cite{Su1}--\cite{Su6}). We however, motivated by \cite {DU1} do not restrict ourselves to conformal dynamical systems only. The results of this chapter are needed to analyze the dynamics and geometry of elliptic functions in Part\ref{EFA} and especially in Part\ref{EFB} but this chapter is also interesting on its own. It is very much stimulated by the development of the theory of iteration of rational functions of the Riemann sphere $\oc$, and to some extent (see Section~\ref{RI}, Ruelle's Inequality and Section~\ref{PT}, Pesin Theory), by the theory of non--hyperbolic smooth dynamical systems. 

\section[Mane's Partition, Pesin Theory, Volume Lemmas] {Mane's Partition, Pesin Theory, \\ Volume Lemmas and Hausdorff Dimension of Invariant Measures}\label{Volume_Lemma}

\sp This is the first section in the book where we intimately connect dynamical/ergodic notions such as Kolmogorov--Sinai metric entropy and Lyapunov exponents with geometric ones such as Hausdorff dimension. We discuss the nature and some history of these relations with appropriate subsections. 

\sp \subsection{Settings} 
Let $Y$ be a Riemann (complex analytic) surface, either compact or open. In our applications it will be always $\C$, $\oc$, or an Euclidean torus $2$--dimensional $\mT$. Let $X$ be a compact subset of $Y$. We say that $f  
\in {\mathcal A}(X)$ if and only if 

\sp\begin{itemize}
\item[(a)] $f: X\lra X$ is a continuous map. 

\sp\item[(b)] $f: X\lra X$ can be 
holomorphically extended to $U(f)=U(f,X)$, an open neighborhood of
of $X$ in $Y$. 

\sp\item[(c)] The map $f:U(f)\lra Y$ has no critical points. 

\sp\item[(d)] If $V\sbt Y$ is an open set intersecting $X$ and all
iterates $f^n:V\lra Y$ are well--defined, then they do not form a normal family.
\end{itemize}

\sp\fr The first technical result of this section is the following.

\blem\label{l10.3.1} 
If $X\sbt Y$ is compact and $f\in {\mathcal A}(X)$, then the series 
$$
\sum_{n=1}^\infty |(f^n)'(z)|^{\frac 13}
$$ 
diverges for all $z\in X$. 
\elem

\bpf By our hypothesis there exists $\e>0$ such that for every $w\in X$ we have that $B(w,\e)\sbt U(f)$ and the map $f$ restricted to the ball
$B(w,\e)$ is 1--to--1. Since, by decreasing the set $U(f)$ if necessary, the map $f:U(f)\lra Y$ is uniformly continuous, there exists $0<\a<1$ such that for every $x\in X$, we have that
\beq\label{10.3.1}
f(B(x,\a\e))\sbt B(f(x),\e).
\eeq
Suppose for a contrary that the series $\sum_{n=1}^\infty |(f^n)'(z)|^{\frac
13}$ converges for some $z\in X$. Then there exists $q\ge 1$ such
that 
$$
\sup_{n\ge q}(2|(f^n)'(z)|)^{\frac 13} < 1.
$$
Choose $0<\e_1=\e_2=\ld=\e_q<\a\e$ so small that for every $n=1,2,\ld,q$,
\beq\label{10.3.2}
\text{the map $f^n$ restricted to the ball $B(z,\e_n)$ is 1--to--1.}
\eeq
and
\beq\label{10.3.3}
f^n(B(z,\e_n))\sbt B(f^n(z),\e).
\eeq
For every $n\ge q$ define $\e_{n+1}$ inductively by
\beq\label{10.3.4}
\e_{n+1}=(1-(2|(f^n)'(z)|)^{\frac 13})\e_n.
\eeq
Then $0< \e_n <\a\e$ for every $n\ge 1$.
Assume that \eqref{10.3.2}
and \eqref{10.3.3} are satisfied for some $n\ge q$. Then by 
Theorem~\ref{th17.46in{H}} (Koebe's Distortion Theorem, analytic
version) and \eqref{10.3.4} the set $f^n(B(z,\e_{n+1}))$ is
contained in the ball centered at $f^n(z)$ and of radius
$$
\e_{n+1}|(f^n)'(z)|\frac 2{(1-\e_{n+1}/\e_n)^3}
=\frac {2\e_{n+1}|(f^n)'(z)|}{2|(f^n)'(z)|}
=\e_{n+1}
<\a \e.
$$
Therefore, since $f$ is injective on $B(f^n(z),\e)$, formula \eqref{10.3.2} is satisfied for $n+1$ and using also \eqref{10.3.1} we get
$$
f^{n+1}(B(z,\e_{n+1}))=f\(f^n(B(z,\e_{n+1}))\)\sbt
f(B(f^n(z),\a \e))\sbt B(f^{n+1}(z),\e).
$$
Thus \eqref{10.3.3} is satisfied for $n+1$. Since the sequence $(\e_n)_{n=1}^\infty$ is monotone decreasing, it has a limit. Denote it by $\e_\infty$.
Since the series
$\sum_{k=1}^\infty|(f^k)'(z)|^{\frac 13}$ converges, it follows from
\eqref{10.3.4} that $\e_\infty>0$. Clearly \eqref{10.3.2} and \eqref{10.3.3}
remain true with $\e_n$ replaced by $\e_\infty$. It thus follows that the family
$$
\lt\{f^n|_{B(z,\e_\infty/2)}\rt\}_{n=1}^\infty
$$ 
is normal, contrary to the assumption (d). We are done. 
\epf

\sp\fr Let $\mu$ be a Borel probability $f$-invariant ergodic
measure on $X$. The characteristic Lyapunov exponent $\chi_\mu(f)$ of
the map $f$ with respect to the measure $\mu$ is this.
$$
\chi_\mu(f):=\int_X\log|f'(x)|\,d\mu(x)<+\infty.
$$
As an immediate consequence of Lemma~\ref{l10.3.1} and of Birkhoff's
Ergodic Theorem, we get the following.

\bcor\label{c10.3.2}
If $Y$ is a Riemann surface, $X\sbt Y$ is compact and $f\in {\mathcal A}(X)$, then 
$$
\chi_\mu(f)\ge 0
$$
for every Borel probability $f$--invariant ergodic measure $\mu$ on $X$.
\ecor

\sp\fr We would like to remark that such inequality was proved by Feliks Przytycki in \cite{P3} for all rational functions $:\oc\lra\oc$ and all $f$--invariant measures supported on the Julia set of $f$. 

\subsection{Ruelle's Inequality}\label{RI}

\fr The first, famous, relation between entropy and Lyapunov exponent
is the following version of Ruelle's Inequality proven by David Ruelle in \cite{Ru2} in the context of smooth diffeomorphisms of multi--dimensional smooth manifolds. 

\bthm[Ruelle's Inequality]\label{t9.1.1}
If $Y$ is a Riemann surface, $X\sbt Y$ is compact, and $f\in\A(X)$, and $\mu$ is a Borel probability $f$-invariant ergodic measure on $X$, then  
$$
\hmu(f)\le 2\max\{0,\chi_\mu(f)\}.
$$
\ethm

\bpf Since $X\sbt Y$ is compact, we may assume without loss of generality that $Y=\C$. Let $S\sbt\C$\, be a square so large that
$X\sbt(1/2)S$. Consider a decreasing to zero sequence 
$(a_k)_{k=1}^\infty$ of positive numbers and $(\Pa_k)_{k=1}^\infty$,
an increasing sequence of partitions of $S$ consisting of squares with
with edges of length equal to $a_k$. For every $g\in\A(X)$ preserving
measure $\mu$, for every $x\in X$, and every integer $k\ge 1$ let 
$$
N(g,x,k):=\#\{P\in\Pa_k:g(P_k(x)\cap U(g))\cap P\ne\es\}.
$$
Our first aim is to show that for every integer $k\ge 1$ large enough,
say $k\ge k(g)$, we have 
\beq\label{9.1.1}
N(g,x,k)\le  \pi(|g'(x)|+2)^2.
\eeq
Indeed, fix $x\in X$ and consider an integer $k\ge 1$ so large ($\ge
k(g)$) that both
$\Pa_k(x)\sbt U(g)$ and the Lipschitz constant of $g|_{\Pa_k(x)}$ does
not exceed $|g'(x)|+1$. Thus the 
set
$g(\Pa_k(x))$ is contained in the ball $B(g(x),(|g'(x)|+1)a_k)$. Therefore
if $g(\Pa_k(x))\cap P\ne\es$, then
$$
P\sbt B(g(x),(|g'(x)|+1)a_k+a_k)=B(g(x),(|g'(x)|+2)a_k).
$$
Hence
$N(g,x,k)\le \pi (|g'(x)|+2)^2a_k^2/a_k^2=\pi(|g'(x)|+2)^2$, and
\eqref{9.1.1} is proved. Let 
$$
N(g,x):=\sup_{k\ge k(g)} N(g,x,k).
$$
In view of \eqref{9.1.1} we get,
\beq\label{9.1.2}
N(g,x)\le \pi(|g'(x)|+2)^2.
\eeq
Now note that for every finite partition $\A$ consisting of Borel sets
we have,
\begin{eqnarray}\label{9.1.3}
\h_\mu(g,\A)&=&\lim_{n\to\infty}\frac 1{n+1}\H_\mu(\A^n)\nonumber\\
&=& \lim_{n\to\infty}\frac 1{n+1}  \Bigl(  \H_\mu(g^{-n}(\A)|\A^{n-1})
+\dots +\H_\mu(g^{-1}(\A)|\A)+\H_\mu(\A)  \Bigr)  \nonumber\\
&\le &
 \lim_{n\to\infty}\frac 1n  \Bigl(  \H_\mu(g^{-n}(\A)|g^{-(n-1)}(\A))
+\dots +\H_\mu(g^{-1}(\A)|\A)  \Bigr)
\nonumber\\&=&\H_\mu(g^{-1}(\A)|\A).  
\end{eqnarray}
Since
$$
\H{\mu_{\Pa_k(x)}}(g^{-1}(\Pa_k)|\Pa_k(x))\le
\log\#\{P\in\Pa_k: g^{-1}(P)\cap\Pa_k(x)\not=\es\}
=\log N(g,x,k),
$$
looking at Corollary~\ref{sem2cor}, we obtain 
$$
\aligned
\H_\mu(g)
&\le \limsup_{k\to\infty}\H_\mu(g^{-1}(\Pa_k)|\Pa_k)
=\limsup_{k\to\infty}\int_X\H_{\mu_{\Pa_k(x)}}(g^{-1}(\Pa_k)|\Pa_k(x))
\,d\mu(x)\\
&\le \limsup_{k\to\infty}\int_X\log N(g,x,k)\,d\mu(x) \\
&\le\int_X\log N(g,x)\,d\mu(x).  
\endaligned
$$
Applying this inequality to $g=f^n$, where $n\ge 1$ is an integer, and employing
\eqref{9.1.2} we get,
$$
\aligned
\H_\mu(f)
&={1\over n}\H_\mu(f^n)
 \le {1\over n}\int_X\log N(f^n,x)\,d\mu(x) 
=\int_X{1\over n}\log N(f^n,x)\,d\mu(x) \\
&\le \int_X{1\over n}\log \pi(|(f^n)'(x)|+2)^2\,d\mu(x).
\endaligned
$$
Since  
$$
0
\le {1\over n}\log (|(f^n)'(x)|+2)^2
\le 2(\log\(\sup_X |f'|)+1)
$$ 
and  
$$
\lim_{n\to\infty}
{1\over n}\log(|(f^n)'(x)|+2)=\max\{0,\chi_\mu(x)\},
$$ 
for $\mu$-a.e $x\in X$, applying Lebesgue's Dominated
Convergence Theorem we get that 
$$
\H_\mu(f)
\le \lim_{n\to\infty}\int_X{1\over n}\log (|(f^n)'(x)|+2)^2\,d\mu(x)
=\int_X\max\{0,2\chi_\mu(x)\}\,d\mu.
$$
The proof is complete. 
\epf

\subsection{Pesin's Theory}\label{PT}

\fr In this section we work in the same general setting of this chapter and we follow the same notation as in the previous sections. We present here a (very special) version of Pesin's Theory whose foundations were laid down in the fundamental works \cite{Pe1} and \cite{Pe2} of Yasha Pesin. Since then Pesin's Theory was extended, generalized, refined, and cited in numerous articles, research books, and textbooks. Our exposition here follows Section~11.2 of the \cite{PU2} book. 

\sp We begin with the following.

\blem\label{l9.2.1.} If $\mu$ is a Borel finite measure on $\R^d$, $d\ge
1$, $a$ is an arbitrary point in $\R^d$ and the function $\R^d\ni z\longmapsto\log\|z-a\|$ is $\mu$--integrable, then for every $C>0$ and every $0<t<1$,
$$
\sum_{n\ge 1}\mu(B(a,Ct^n))<+\infty.
$$
\elem

\bpf Since $\mu$ is finite and since given $t<s<1$ there exists
$q\ge 1$ such that $Ct^n\le s^n$ for all $n\ge q$, without loosing generality
we may assume that $C=1$. Recall that given $b\in\R^d$, and two
numbers $0\le r<R$, 
$$
R(b;r,R)=\{z\in\R^d:r\le|z-b|<R\}.
$$
Since
$-\log(t^n)\le -\log||z-a||$ for
every $z\in B(a,t^n)$ we get the following.
$$
\aligned
\sum_{n\ge 1}\mu(B(a, t^n))
&=\sum_{n\ge 1}n\mu(A(a;t^{n+1},t^n))
={-1\over \log t}\sum_{n\ge 1}-\log(t^n)\mu(A(a;t^{n+1},t^n)) \\
&\le {-1\over \log t}\int_{B(a,t)}-\log\|z-a\|\,d\mu(z) \\
&<+\infty.
\endaligned
$$
The proof is finished. 
\epf

\blem\label{l9.2.2.} 
Let $Y$ be a Riemann surface, $X\sbt Y$ be a compact set, and let $f\in\A(X)$. If $\mu$ is a Borel finite measure on $X$ and $\log|f'|$ is $\mu$ integrable, then the function 
$$
X\ni z\longmapsto\log|z-c|\in L^1(\mu)
$$ 
for every critical point $c$ of $f$. If additionally $\mu$ is
$f$--invariant, then also the function 
$$
X\ni z\longmapsto\log|z-f(c)|\in L^1(\mu).
$$
\elem

\bpf Since $X\sbt Y$ is compact, we may assume without loss of generality that $Y=\C$. That 
$$
X\ni z\longmapsto\log|z-c|\in L^1(\mu)
$$
follows from the fact that near $c$ we have 
$$
C^{-1}|z-c|^{q-1}\le |f'(z)|\le C|z-c|^{q-1},
$$
where $q\ge 2$ is the order of the critical point $c$ and $C\ge 1$ is a  constant independent of $z$, and since out of any neighbourhood of the set
of critical points of $f$, the derivative $|f'(z)|$ is uniformly bounded away from zero and infinity. 

\sp In order to prove the second part of the lemma,
consider a closed ray $R$ emanating from $f(c)$ such that $\mu(R)=0$ and a
disk $B(f(c),r)$ such that 
$$
f_c^{-1}:B(f(c),r)\sms R\lra \C
$$
the holomprphic inverse branch of $f$ sending $f(c)$ to $c$, is
well--defined. Let 
$$
D:= B(f(c),r)\sms R.
$$
We may additionally require $r>0$ to be so small that 
$$
|z-f(c)|\comp |f_c^{-1}(z)-c|^q
$$
for all $z\in f_c^{-1}\(B(f(c),r)\)$. It suffices to show that the integral 
$$
\int_D\log|z-f(c)|\,d\mu(z)
$$
is finite. And indeed, by $f$-invariance of $\mu$, we have
$$
\aligned
\int_D\log|z-f(c)|\,d\mu(z)
&=\int_X\!\!\1_D(z)\log|z-f(c)|\,d\mu(z) 
\comp_A \int_X\!\!\1_D(z)\log|f_c^{-1}(z)-c|^q\,d\mu(z) \\
&=\int_X(\1_D\circ f)(z)\log|z-c|^q\,d\mu(z)\\
&=\int_X\1_{f^{-1}(D)}\log|z-c|^q\,d\mu(z) 
\endaligned
$$
Notice here that the function 
$$
z\longmapsto\1_D(z)\log|f_c^{-1}(z)-c|^q
$$ 
]is indeed well--defined on $X$and that the comparability sign $\comp_A$, 
in the formula above, means an additive comparability. The finiteness
of the last integral follows from the first part of this lemma. 
\epf

\bthm\label{t9.2.3.} 
Let $(Z,\F,\nu)$ be a measure space with an ergodic measure preserving automorphism $T:Z\to Z$. Let $Y$ be a Riemann surface, $X\sbt Y$ be a compact set, and let $f\in\A(X)$.

Suppose that $\mu$ is an $f$--invariant ergodic measure on $X$ with positive Lyapunov exponent $\chi_\mu(f)$. Suppose also that $h:Z\to X$ is a measurable mapping such that 
$$
\nu\circ h^{-1}=\mu
\ \ {\rm and} \  \
h\circ T=f\circ h
$$ 
$\nu$-a.e..

Then for $\nu$-a.e. $z\in Z$ there
exists $r(z)>0$ such that the function $Z\ni z\mapsto r(z)$ is measurable
and the following is satisfied: 

\sp\fr For every $n\ge 1$ there exists $f_{x_n}^{-n}:B(x,r(z))\lra Y$,
a holomorphic inverse branch of $f^n$ sending $x:=h(z)$ to $x_n:=h(T^{-n}(z))$. Given in addition an arbitrary $\chi\in (-\chi_\mu(f),0)$, there exists a constant $K_\chi(z)\in(0,+\infty)$ independent of $n$ and $z$, such that
$$
|(f_{x_n}^{-n})'(y)|<K(z)\ep^{\chi n} \  \text{ and } \  {|(f_{x_n}^{-n})'(w)|
\over |(f_{x_n}^{-n})'(y)|} \le K
$$
for all $y,w\in B(x,r(z))$. $K$ as usually is in here the Koebe constant corresponding to the scale $1/2$.
\ethm

\bpf Suppose first that $\mu\(\bu_{n\ge 1}f^n(\Crit(f))\)>0$.
Since the measure $\mu$ is ergodic this implies that $\mu$ must be
concentrated on a 
periodic orbit of an element $w\in \bu_{n\ge 1}f^n(\Crit(f))$. This means
that $w=f^q(c)=f^{q+k}(c)$ for some $q,k\ge 1$ and $c\in \Crit(f)$, and
$$
\mu(\{f^q(c), f^{q+1}(c),\ld,f^{q+k-1}(c)\})=1.
$$
Since $\int\log|f'|\,d\mu >0$, we have that $|(f^k)'(f^q(c))|>1$. Thus the theorem is obviously true for the set
$$
h^{-1}(\{f^q(c), f^{q+1}(c),\ld,f^{q+k-1}(c)\})
$$ \

of $\nu$ measure 1.

So, suppose that 
$$
\mu\lt(\bu_{n\ge 1}f^n(\Crit(f))\rt)=0.
$$
Set 
$$
R:=\min\big\{1,\dist(X,Y\sms U(f))\big\}
$$ 
and fix $\l\in (\ep^{{1\over 4}\chi},1)$. Consider $z\in Z$ such that $x=h(z)\notin
\bu_{n\ge 1}f^n(\Crit(f))$,
$$
\lim_{n\to\infty}{1\over n}
\log|(f^n)'(h(T^{-n}(z))|=\chi_\mu(f),
$$
and $x_n=h(T^{-n}(z))$ belongs to
$B(f(\Crit(f)),R\l^n)$ for finitely many $n$'s only. We shall first
demonstrate that the set of points satisfying these properties is of
full measure $\nu$. Indeed, the first requirement is satisfied by our
hypothesis, the second is due to Birkhoff's Ergodic Theorem. In order
to prove that the set of points satisfying the third condition has $\nu$
measure 1, notice that
$$
\aligned
\sum_{n\ge 1} \nu\(T^n(h^{-1} (B(f(\Crit(f)), R\l^n)))\) 
&=\sum_{n\ge 1} \nu\(h^{-1} (B(f(\Crit(f)), R\l^n))\) \\
&=\sum_{n\ge 1} \mu (B(f(\Crit(f)), R\l^n)) \\
&< +\infty, 
\endaligned
$$
where the last inequality we wrote due to Lemma~\ref{l9.2.2.} and
Lemma~\ref{l9.2.1.}.  
The application of the Borel--Cantelli lemma finishes now the
demonstration. Fix now an integer $n_1=n_1(z)$ so large that 
$$
x_n=h(T^{-n}(z))\notin B(f(\Crit(f)),R\l^n)
$$ 
for all $n\ge n_1$. Notice that because of our choices there exists $n_2\ge n_1$ such that
$$
|(f^n)'(x_n)|^{-1/4}<\l^n
$$ 
for all $n\ge n_2$. Finally set
$$
S:=\sum_{n\ge 1}|(f^n)'(x_n)|^{-1/4},
$$
$$
b_n:={1\over 2}S^{-1}|(f^{n+1})'(x_{n+1})|^{{-1\over 4}},
$$
and
$$
\Pi:=\Pi_{n=1}^\infty(1-b_n)^{-1}.
$$
This infinite product converges since the series $\sum_{n\ge 1}b_n$ does. Choose
now $r=r(z)$ so small that $16r(z)\Pi KS^3\le R$, all the inverse branches
$f_{x_n}^{-n}:B\(x_0,\Pi r(z)\)\lra\oc$ are well-defined for all $n=1,2,\ld,
n_2$ and 
$$
\diam\(f_{x_{n_2}}^{-n_2}(B\(x_0,r\Pi_{k\ge n_2}(1-b_k)^{-1})\)\le
\l^{n_2}R. 
$$
We shall show by induction that for every
$n\ge n_2$ there exists an analytic inverse branch
$f_{x_n}^{-n}:B\(x_0,r\Pi_{k\ge n}(1-b_k)^{-1}\)\lra Y$, sending $x_0$ to
$x_n$ and such that
$$
\diam\(f_{x_n}^{-n}(B\(x_0,r\Pi_{k\ge n}(1-b_k)^{-1})\)\le \l^nR.
$$
Indeed, for $n=n_2$ this immediately follows from our requirements imposed
on $r(z)$. So, suppose that the claim is true for some $n\ge n_2$. Since 
$$
x_n=f_{x_n}^{-n}(x_0)\notin B(\Crit(f),R\l^n)
$$ 
and since $\l^nR\le R$, there
exists an inverse branch $f_{x_{n+1}}^{-1}:B(x_n,\l^nR)\lra Y$ sending
$x_n$ to $x_{n+1}$. Since 
$$
\diam\(f_{x_n}^{-n}(B\((x_0,r\Pi_{k\ge n}
(1-b_k)^{-1})\)\le \l^nR,
$$ 
the composition  
$$
f_{x_{n+1}}^{-1}\circ f_{x_n}^{-n}B\(x_0,r\Pi_{k\ge n}(1-b_k)^{-1})\lra Y
$$ 
is well--defined and forms the inverse branch of $f^{n+1}$ that sends $x_0$ to $x_{n+1}$. By Koebe's Distortion Theorem (Theorem~\ref{Euclid-I}), we now estimate
$$
\aligned
\diam\(f_{x_{n+1}}^{-(n+1)} &(B\(x_0,r\Pi_{k\ge n+1}(1-b_k)^{-1}))\)
\le \\
&\le 2r\Pi_{k\ge n+1}(1-b_k)^{-1}|(f^{n+1})'(x_{n+1})|^{-1}Kb_n^{-3} \\
&\le 16r\Pi KS^3|(f^{n+1})'(x_{n+1})|^{-1}|(f^{n+1})'(x_{n+1})|^{{3\over 4}} \\
&=16r\Pi KS^3|(f^{n+1})'(x_{n+1})|^{-{1\over 4}} \\
&\le R\l^{n+1}, 
\endaligned
$$
where the last inequality sign we wrote due to our choice of
$r$ and the number $n_2$. Putting $r(z)=r/2$
the second part of this theorem follows now as a combined application of
the equality $\lim_{n\to\infty}{1\over n}\log|(f^n)'(x_n)|=\chi_\mu(f)$ and
Koebe's Distortion Theorem (Theorem~\ref{Euclid-I}).
\epf

\sp\fr As an immediate consequence of Theorem~\ref{t9.2.3.} we get the following.

\bcor\label{c9.2.4} Assume the same notation and asumptions as in
Theorem~\ref{t9.2.3.}.Fix $\e>0$. Then there exist a measurable set $Z(\e)\sbt Z$, the
numbers $r(\e)\in (0,1)$ and $K(\e)\ge 1$ such that $\mu(Z(\e))>1-\e$,
$r(z)\ge r(\e)$ for all $z\in Z(\e)$ and with $x_n=h(T^{-n}(z))$, we have that
$$
K(\e)^{-1}
\le\exp(-(\chi_\mu+\e)n)
\le |(f_{x_n}^{-n})'(y)|
\le K(\e)\exp(-(\chi_\mu-\e)n) 
$$
and 
$$
K^{-1}\le  {|(f_{x_n}^{-n})'(w)|\over |(f_{x_n}^{-n})'(y)|}\le K
$$
for all $n\ge 1$, all $z\in Z(\e)$ and all $y,w\in B(x_0,r(\e))$. $K$
is here the Koebe constant corresponding to the scale $1/2$.
\ecor

\brem\label{9.2.5.} In our future applications the system $(Z,f,\nu)$ will
be usually given by the Rokhlin's natural extension (see
Theorem~\ref{t2ms33}) of the holomorphic system $(f,\mu)$. 
\erem

\sp\subsection{Two Auxiliary Partitions}

For the purpose of this section for every Lebesgue measurable
subset $A$ of $\R$ denote by $|A|$ the Lebesgue measure of $A$. We
start with the following technical fact which is however 
interesting on its own.

\blem\label{l7.1.7}
Every monotone increasing function $k:I\to \R$ defined on a bounded
closed interval $I\subset\R$ is Lipschitz continuous at Lebesgue
almost every point in $I$. In other words,
for every $\e>0$ there exist $L>0$ and a set $A\subset I$ such that
$|I\setminus A|<\e$ and the function $k:I\to \R$
is Lipschitz continuous at each point of $A$ with the Lipschitz
constant not exceeding $L$.
\elem

\bpf For every $y\in I$ let 
$$
L_y:=\sup\lt\{z\in I\sms\{y\}: {|k(z)-k(y)| \over |z-y|}\rt\}.
$$
Suppose on the contrary, that the set 
$$
B:=\{y\in I: L_y=+\infty\}
$$
has positive Lebesgue measure. Write $I=[a,b]$. Replacing $B$ by its
subset if necessary, we may assume without loosing generality that the set
$B$ is compact and contains
neither $a$ nor $b$. For every $y\in B$ choose $y'\in I\sms\{y\}$ such that
\beq\label{7.1.6'}
{|k(y)-k(y')| \over |y-y'|}>2{k(b)-k(a) \over |B|}.     
\eeq
Every such pair $y,y'$ replace now by a pair $x,x'$ such that  $I\spt
(x,x')\supset [y,y']$, and the points $x,x'$ are so close 
to respective points $y,y'$ that \eqref{7.1.6'} still holds with $x,x'$.
However, if (at least) one of the points $y$ or $y'$ happens to be an end-point
of $I$, we make no replacement. Now from the family of intervals
$(x,x')$ choose a finite family  $\cI$ covering our compact set $B$.
From this family it is in turn possible to choose a subfamily which is
a union of two subfamilies $\cI_1$ and $\cI_2$ each of which consist
of mutually disjoint open intervals. Using also monotonicity of the
function $k:I\to\R$ and formula \eqref{7.1.6'}, we get for $i=1,2$ that 
$$
k(b)-k(a)\ge \sum_{(w,w')\in {\cI}^i} (k(w')-k(w)) >
2{k(b)-k(a) \over
   |B|}
\sum_{(w,w')\in {\cI}^i} (w'-w) 
$$
Hence, taking into account that ${\cI}^1\cup{\cI}^2$ covers $B$,
we get
$$
2(k(b)-k(a)) >
2{k(b)-k(a) \over |B|}\!\!\! \sum_{(w,w') \in {\cI}^1\cup {\cI}^2}\!\!\!(w'-w)
\ge 2{k(b)-k(a) \over |B|}|B|=2(k(b)-k(a)),
$$
which is a contradiction. Thus, the proof is finished.
\epf

\bcor\label{c7.1.8}
 For every Borel probability measure $\nu$ on a
compact metric space
$(X,\rho)$ and for every $r>0$ there exists a finite partition
${\cP}=\{P_j\}_{j=1}^M$ of $X$ into Borel subsets of $X$, each of which with
positive measure $\nu$, that satisfies the following two properties.
\beq\label{7.1.7a} 
\diam(\cP) < r
\eeq
and there exists a constant $C>0$ such that for every $\eta>0$
\beq\label{7.1.7}
\nu(\partial_\eta\cP)\le C\eta,               
\eeq
where 
$$
\partial_\eta\cP:=\bigcap_{j=1}^M\Bigl(\bigcup_{k\not=j}B(P_s,\eta)\Bigr).
$$
\ecor 
\bpf Let $\{x_1,\dots ,x_N\}$ be a finite $r/4$-spanning set in $X$.
Fix  $\e\in (0,r/4N)$. For each monotone increasing function
$$
I:=[r/4,r/2]\ni t\longmapsto k_i(t):=\nu(B(x_i,t)),
$$
$i=1,2,\ld,N$, apply Lemma~\ref{l7.1.7}, and consider respective constants $L_i$ and sets $A_i$.
Let 
$$
L=\max\{L_i, i=1,\dots ,N\}
$$ 
and let
$$
A=\bigcap_{i=1,\dots ,N} A_i.
$$
The set $A$ has positive Lebesgue measure by the choice
of $\e$. So, we can choose its point $s$ different from $r/4$ and $r/2$.
Therefore, for all $0<\eta<\eta_0:=\min\{s-r/4,r/2-s\}$ and for all $i\in
\{1,2,\ld,N\}, $ we have that 
$$
\nu(B(x_i,r_0+\eta)\setminus B(x_i,s-\eta))\le 2L\eta.
$$
Hence, putting
$$
\Delta(a):=\bigcup_{i=1}^N \bigl( B(x_i,s+\eta)\setminus B(x_i,s-\eta)\bigr),
$$
we get $\nu(\Delta(\eta))\le 2LN\eta$. Put
$$
B^+(x_i,s):=B(x_i,s) \  \text{ and }  \ B^-(x_i,s):=X\setminus
B(x_i,s).
$$
Let $\Om$ be the set of all functions $\kappa:\{1,\dots ,N\}\to \{+,-\}$.
Define 
$$
\cP=\lt\{\bigcap_{i=1}^N B^{\kappa(i)}(x_i,s):\ka\in\Om\rt\}.
$$ 
After removing from $X$ of a set of measure 0, the partition $\cP$ covers $X$.
Since $s\ge r/4$, the balls $B(x_i,s)$, $i-1,2\ld,N$, cover $X$. Hence, for each
non-empty $P_j\in{\cP}$ at least one value of $\kappa$ is equal to $+$. Thus
$\diam(P_j)\le 2s<r$.
We shall show now that 
\beq\label{920130511}
\partial_\eta\cP\subset\Delta(a).
\eeq
Indeed, let $x\in\partial_\eta\cP$. Since $\cP$ covers $X$ there exists $j_0$
such that $x\in P_{j_0}$; so $x\notin P_j$ for all $j\not=j_0$. However,
since $x\in \bigcup_{j\not=j_0}B(P_t,\eta)$,
there exists $j_1\not=j_0$ such that
$\dist (x,P_{j_1})<a$. Let $B:=B(x_i,s)$ be such that 
$$
P_{j_0}\subset B^+
\  \  \  {\rm and} \  \  \  
P_{tj1}\subset B^-,
$$
or  vice versa. 

In the case when 
$$
x\in P_{t_0}\subset B^+,
$$ 
by the triangle inequality $\rho(x,x_i)>s-\eta$ and since
$\rho(x,x_i)<s$, we get that
$$
x\in\Delta(\eta).
$$

In the case when
$$
x\in P_{j_0}\subset B^-
$$ 
we have
$$
x\in B(x_i,s+\eta) 
\setminus B(x_i,s)\subset\Delta(\eta).
$$
Formula \eqref{920130511} is therefore proved.

We thus conclude that 
$$
\nu(\partial_\eta\cP)\le \nu(\Delta(a))\le 2LNa 
$$
for all $\eta<\eta_0$.
For $\eta\ge\eta_0$ it suffices to take $C\ge 1/\eta_0$. So the corollary
is proved, with $C:=\max\{2LN,1/\eta_0\}$. 
\epf

\sp\bcor\label{c7.1.9} 
Let $\nu$ be a Borel probability measure on a compact metric space
$(X,\rho)$ and let $T:X\lra X$ be a Borel endomorphism on $X$
preserving the measure $\nu$. 

Then for every
$r>0$ there exists a finite partition ${\cP}=\{P_j\}_{j=1}^M$
of $X$ into Borel sets of positive measure $\nu$ with $\diam(\cP)<r$ such
that for every $\d>0$ and $\nu$--a.e. $x\in X$ there exists an integer
$n_0=n_0(x)$ such that for every $n\ge n_0$,
\beq\label{7.1.8}
B\(T^n(x),\exp(-n\d)\)\subset {\cP}(T^n(x)).                
\eeq
\ecor
\bpf Let $\cP$ be the partition from Corollary \ref{c7.1.8}.
Fix an arbitrary $\d>0$. Then by Corollary \ref{c7.1.8}
$$
\sum_{n=0}^\infty
\nu\(\partial_{\exp(-n\d)}\cP\)\le\sum_{n=0}^\infty C\exp(-n\d)<\infty.
$$
Hence by the $T$-invariance of $\nu$, we obtain this,
$$
\sum_{n=0}^\infty \nu\(T^{-n}(\partial_{\exp(-n\d)}\cP\)\) <\infty.
$$
Applying now the Borel-Cantelli Lemma for the family
$\lt\{T^{-n}\(\partial_{\exp(-n\d)}\cP\)\rt\}_{n=1}^\infty$ we conclude that
for $\nu$-a.e $x\in X$ there exists $n_0=n_0(x)$ such that for every
$n\ge n_0$ we have that $x\notin
T^{-n}\(\partial_{\exp(-n\d)}\cP\)$. This means that
$T^n(x)\notin \partial_{\exp(-n\d)}\cP$. Hence, by the definition of
$\partial_{\exp(-n\d)}\cP$, if 
$T^n(x)\in P_j$ for some $j=1,2,\ld,M$, then
$$
T^n(x)\notin \bigcup_{k\not=j}B\(P_k,\exp(-n\d)\).
$$
Thus
$$
B\(T^n(x),\exp-n\d\)\subset \cP(T^n(x)).
$$
The proof is finished.
\epf

\

\subsection{Ma\~n\'e's Partition}

\fr In this subsection, basically following Ma\~n\'e's book
\cite{Ma1987}, we construct the so 
called Ma\~n\'e's partition which will play an
important role in the proof of a part of the Volume Lemma given in the next
section. We begin with the following elementary fact.

\sp\blem\label{l9.3.1.} If $x_n\in(0,1)$ for every $n\ge 1$ and
$\sum_{n=1}^\infty nx_n<+\infty$, then 
$$
\sum_{n=1}^\infty -x_n\log x_n<+\infty.
$$
\elem
\bpf Let $S:=\{n\ge 1:-\log x_n\ge n\}$. Then
$$
\sum_{n=1}^\infty -x_n\log x_n =\sum_{n\notin S}-x_n\log x_n +\sum_{n\in S}-
x_n\log x_n\le \sum_{n=1}^\infty nx_n +\sum_{n\in S}-x_n\log x_n.
$$
Since $n\in S$ means that $x_n\le e^{-n}$ and since $\log t\le 2\sqrt t$ for
all $t\ge 1$, we have
$$
\sum_{n\in S}x_n\log{1\over x_n} \le 2\sum_{n=1}^\infty x_n\sqrt{{1\over x_n}}
\le 2\sum_{n=1}^\infty e^{-{1\over 2}n} <\infty.
$$
The proof is finished. 
\epf

\sp The next lemma is the main and simultaneously the last result of this subsection.

\sp

\blem\label{l9.3.2.} If $\mu$ is a Borel probability measure supported on
a bounded subset $M$ of a Euclidean space and $\rho:M\lra (0,1]$ is a
measurable function such that $\log\rho$ is integrable with respect to $\mu$,
then there exists a countable measurable
partition, called \index{(N)}{partition Ma\~n\'e's} {\it Ma\~n\'e's 
  partition}, $\Pa$ of $M$ such that 
$\H_\mu(\Pa)<+\infty$ and
$$
\diam(\Pa(x))\le \rho(x)
$$
for $\mu$--almost every $x\in M$. 
\elem

\bpf Let $q$ be the dimension of the Euclidean space containing
$M$. Since $M$ is bounded, there exists a constant $C>0$ such that for every
$0<r<1$ there exists a partition $\Pa_r$ of $M$ of diameter $\le r$ and which
consists of at most $Cr^{-q}$ elements.
For every $n\ge 0$ put 
$$
U_n:=\{x\in M: e^{-(n+1)}<\rho(x)\le e^{-n}\}.
$$
Since $\log\rho$ is a non-positive integrable function, we have
\beq
\sum_{n=1}^\infty-n\mu(U_n)
\ge \sum_{n=1}^\infty\int_{U_n}\log\rho\,d\mu
=\int_M\log\rho\,d\mu >-\infty,
\eeq
so that
\beq\label{9.3.1}
\sum_{n=1}^\infty n\mu(U_n) <+\infty.
\eeq
Define now $\Pa$ as the partition whose atoms are of the form $Q\cap U_n$,
where $n\ge 0$ and $Q\in\Pa_{r_n}$, $r_n=e^{-(n+1)}$. Then
$$
\H_\mu(\Pa)=\sum_{n=0}^\infty\lt(-\sum_{U_n\spt P\in\Pa}\mu(P)\log\mu(P)\rt).
$$
But for every $n\ge 0$,
$$
\aligned
-\sum_{U_n\spt P\in\Pa}\mu(P)\log\mu(P) 
&= \mu(U_n)\sum_P-{\mu(P)\over\mu(U_n)} \log\({\mu(P)\over\mu(U_n)}\)-
\mu(U_n)\sum_P{\mu(P)\over\mu(U_n)} \log(\mu(U_n)) \\ 
&\le \mu(U_n)(\log C - q\log r_n) - \mu(U_n)\log\mu(U_n) \cr  &\le
\mu(U_n)\log C +q(n+1)\mu(U_n)-\mu(U_n)\log\mu(U_n). 
\endaligned
$$
Thus, summing over all $n\ge 0$, we obtain
$$
\H_\mu(\Pa)\le \log C +q+q\sum_{n=0}^\infty n\mu(U_n)+\sum_{n=0}^\infty-
\mu(U_n)\log\mu(U_n).
$$
Therefore looking at \eqref{9.3.1} and Lemma~\ref{9.3.1} we conclude
that $\H_\mu(\Pa)$ is 
finite.
Also, if $x\in U_n$, then the atom $\Pa(x)$ is contained in some
atom of $\P_{r_n}$ and therefore
$$
\diam(\Pa(x))\le r_n=e^{-(n+1)}<\rho(x).
$$
Now the remark that the union of all the sets $U_n$ is of measure 1 completes
the proof. 
\epf

\sp \subsection{Volume Lemmas and Hausdorff Dimension of Invariant
 Measures} 

This sections is entirely devoted to prove provide a closed formula for Hausdorff and packing dimensions of Borel probability measures invariant under a map from ${\mathcal A}(X)$. For historical reasons, partly justified, this formula is frequently referred to as a Volume Lemma. Its first forms can be traced back to the works \cite{Eggleston} and \cite {Billingsley} of Eggleston and Billingsley respectively. From the dynamical point of view a kind of breakthrough was the paper \cite{LSY0} by Lai Sang Young. Since then a multitude of papers appeared. We would would like to mention only some early ones: \cite{MM}, \cite{P1}, \cite{M2}. Also, early papers \cite{Bow2} and \cite{PUZI} shed some light on the nature of dimensions of measures. 

The main result of this subsection and of the entire section is this. 

\bthm[Volume Lemma]\label{t20120801}
Let $Y$ be a Riemann surface, $X\sbt Y$ be a compact set. If $f  
\in {\mathcal A}(X)$ and $\mu$ is a Borel probability $f$--invariant
ergodic measure on $X$ with $\chi_\mu(f)>0$, then
$$
\lim_{r\downto 0}\frac{\log\mu(B(x,r))}{\log r}
=\frac{\hmu(f)}{\chi_\mu(f)}
$$
for $\mu$--a.e. $x\in X$. In particular, the measure $\mu$ is dimensional exact and, by Proposition~\ref{p320191130}, 
$$
\HD(\mu)=\PD(mu)=\frac{\hmu(f)}{\chi_\mu(f)}.
$$
\ethm
\bpf 
In view of Theorem~\ref{t6.6.4.}, the second formula follows from
the first one and we therefore only need to prove the first one. Let us
prove first that
\beq\label{9.4.1}
\liminf_{r\to 0}{\log(\mu(B(x,r)))\over \log r}
\ge {\hmu(f)\over \chi_\mu(f)}  
\eeq
for $\mu$--a.e. $x\in X$. By Corollary~\ref{c7.1.9}
there exists a finite partition $\Pa$ 
such that for an arbitrary $\e>0$ and every $x$ in a set $X_o$
of full measure $\mu$ there exists $n(x)\ge 0$
such that for all $n\ge n(x)$.
\beq\label{9.4.2}
B(f^n(x),e^{-\e n})\sbt \Pa(f^n(x)). 
\eeq
Let us work from now on in the Rokhlin's natural extension (see
Theorem~\ref{t2ms33}) $(\^X,\^f,\^\mu)$. Let 
$\^X(\e)$ and $r(\e)$ be given by Corollary~\ref{c9.2.4},
i.e. $\^X(\e)=Z(\e)$. In view of  Birkhoff's
Ergodic Theorem (Theorem~\ref{Birkhoff}) there exists a measurable set $\^F(\e)\sbt\^X(\e)$ such that
$\^\mu(\^F(\e))=\^\mu(\^X(\e))$ and
$$
\lim_{n\to\infty}{1\over n}\sum_{j=1}^{n-1}\chi_{\^X(\e)}\circ\^f^n(\^x)=
\^\mu(\^X(\e))
$$
for every $\^x\in\^F(\e)$. Let $F(\e)=\pi(\^F(\e))$. Then 
$$
\mu(F(\e))
=\^\mu(\pi^{-1}(F(\e)) \ge
\^\mu(\^F(\e))=\^\mu(\^X(\e))
$$ 
converges to $1$ if $\e\downto
0$. Consider now an arbitrary point
$$
x\in F(\e)\cap X_o
$$ 
and take $\^x\in\^F(\e)$ such that
$x=\pi(\^x)$. Then by the above there exists an increasing sequence $\{n_k=
n_k(x):k\ge 1\}$ such that $\^f^{n_k}(\^x)\in\^X(\e)$ and
\beq\label{9.4.3}
{n_{k+1}-n_k\over n_k}\le \e
\eeq
for every $k\ge 1$. Moreover, we can assume that $n_1\ge n(x)$.
Consider now an integer $n\ge n_1$ and the ball 
$$
B\(x,Cr(\e)\exp(-(\chi_\mu+(2+\log\|f'\|)\e)n)\),
$$
where $0<C<K(\e)^{-1}$ is a
constant (possibly depending on $x$) so small that
\beq\label{9.4.4}
f^q\lt(B\lt(x,Cr(\e)\exp\Big(-(\chi_\mu+(2+\log\|f'\|)\e)n)\Big)\rt)\rt) \sbt P(f^q(x))
\eeq
for every $q\le n_1$ and $K(\e)\ge 1$ is the constant appearing
in Corollary~\ref{c9.2.4}.
Take now any $q$, $n_1\le q\le n$, and associate $k$ such that $n_k\le
q\le n_{k+1}$.
Since $\^f^{n_k}(\^x)\in\^X(\e)$ and since $\pi(\^f^{n_k}(\^x)) =f^{n_k}(x)$,
Corollary~\ref{c9.2.4} produces a holomorphic inverse branch 
$$
f_x^{-n_k}:B(f^{n_k}(x),r(\e))\lra\oc
$$ 
of $f^{n_k}$ such that $f_x^{-n_k}(
f^{n_k}(x))=x$. This corollary also yields
$$
f_x^{-n_k}\(B(f^{n_k}(x),r(\e))\)
\spt B\Big(x,K(\e)^{-1}r(\e)\exp\(-(\chi_\mu+\e)n_k\)\Big).
$$
Since 
$$
B\lt(x,Cr(\e)\exp\Big(-(\chi_\mu+(2+\log\|f'\|)\e)n\Big)\rt)
\sbt B\Big(x,K(\e)^{-1}r(\e)\exp\(-(\chi_\mu+\e)n_k\)\Big),
$$
it follows from Corollary~\ref{c9.2.4} that
$$
\aligned
f^{n_k}\Big(B\Big(x,Cr(\e)\exp\(&-(\chi_\mu+(2+\log\|f'\|)\e)n\)\Big)\Big)
 \sbt \\
&\sbt B\Big(f^{n_k}(x),CK(\e)r(\e)e^{-
\chi_\mu(n-n_k)}\exp(\e(n_k-(2+\log\|f'\|)n))\Big). 
\endaligned
$$
Since $n\ge n_k$ and since, by \eqref{9.4.3}, $q-n_k\le \e n_k$, we
therefore obtain 
$$
\aligned
f^q &\Big(B\Big(x,Cr(\e)\exp\(-(\chi_\mu+(2+\log\|f'\|)\e)n\)\Big)\Big)\sbt \\
&\sbt B\Big(f^q(x),CK(\e)r(\e)e^{-\chi_\mu(n-n_k)}\exp\(\e(n_k-(2+\log\|f'\|)n)\)\exp\((q-n_k)\log\|f'\|\)\Big) \\
&\sbt B\Big(f^q(x),CK(\e)r(\e)\exp\(\e(n_k\log\|f'\|+n_k-2n-n\log\|f'\|\)\Big) \\
&\sbt B\(f^q(x),CK(\e)r(\e)e^{-\e n}\) \\
&\sbt B(f^q(x),e^{-\e q}).
\endaligned
$$
Combining this, \eqref{9.4.2}, and \eqref{9.4.4}, we get
$$
B\lt(x,Cr(\e)\exp\Big(-(\chi_\mu+(2+\log\|f'\|)\e)n\Big)\rt) 
\sbt \bigvee_{j=0}^nf^{-j}(\Pa)(x).
$$
Therefore, applying Theorem~\ref{SMBTEC} (the ergodic case of the
Shannon--McMillan--Breiman Theorem), we have 
$$
\varliminf_{n\to\infty}-{1\over n}\log\mu\Big(B(x,Cr(\e)
\exp\(-(\chi_\mu+(2+\log\|f'\|) \e)n\)\Big)
\ge \hmu(f,\Pa)\ge \hmu(f) -\e.
$$
This means that denoting the number
$Cr(\e)\exp\(-(\chi_\mu+(2+\log\|f'\|)\e)n\)$ by $r_n$, we have
$$
\liminf_{n\to\infty}{\log\mu(B(x,r_n)\over \log r_n}\ge {\hmu(f) -\e \over
\chi_\mu(f)+(2+\log\|f'\|)\e}.
$$
Now, since $\{r_n\}$ is a geometric sequence and since $\e>0$ can be taken
arbitrarily small, we conclude that for $\mu$-a.e. $x\in X$
$$
\liminf_{n\to\infty}{\log\mu(B(x,r)\over \log r}\ge {\hmu(f)\over
\chi_\mu(f)}.
$$
This completes the proof of \eqref{9.4.1}.

\sp Now let us prove that
\beq\label{9.4.5}
\limsup_{r\to 0}{\log(\mu(B(x,r)))\over \log r}\le \hmu(f)/\chi_\mu(f)
\eeq
for $\mu$-a.e. $x\in X$.

\sp\fr In order to prove this formula we again work in the Rokhlin's
natural extension 
$(\^X,\^f,\^\mu)$ and we apply Pesin theory. In particular the sets
$\^X(\e)$, $\^F(\e)\sbt\^X(\e)$ and the radius $r(\e)$, produced in
Corollary~\ref{c9.2.4} have the same meaning as in 
the proof of \eqref{9.4.1}. 

To begin with notice that there exist two numbers $R>0$ 
and $0<Q<\min\{1,r(\e)/2\}$ such that the following two conditions are
satisfied.
\beq\label{9.4.6}
\mbox{If $z\notin B(\Crit(f),R)$, then $f|_{B(z,Q)}$ is injective.}
\eeq
and
\beq\label{9.4.7}
\mbox{If $z\in B(\Crit(f),R)$, then $f|_{B(z,Q\dist(z,\Crit(f)))}$ is
injective. }
\eeq
Observe also that if $z$ is sufficiently close to a critical point $c$, then
$f'(z)$ is of order $(z-c)^{q-1}$, where $q\ge 2$ is the order of critical
point $c$. In particular the quotient of $f'(z)$ and $(z-c)^{q-1}$ remains
bounded away from 0 and $\infty$ and therefore there exists a constant number
$B>1$ such that 
$$
|f'(z)|\le B\dist(z,\Crit(f)).
$$
So, in view of Lemma~\ref{l9.2.2.}, the logarithm of the function
$$
\rho(z):=Q\min\{1,\dist(z,\Crit(f))
$$ 
is integrable and consequently Lemma~\ref{l9.3.2.} 
applies. Let $\Pa$ be the Ma\~n\'e's partition produced by this lemma. Then
$B(x,\rho(x))\spt\Pa(x)$ for $\mu$--a.e. $x\in X$, say for a subset $X_\rho$ of
$X$ of measure 1. Consequently
\beq\label{9.4.8}
B_n(x,\rho)=\bi_{j=0}^{n-1}f^{-j}\(B(f^j(x),\rho(f^j(x)))\)\spt \Pa_0^n(x)
\eeq
for every $n\ge 1$ and every $x\in X_\rho$.
By our choice of $Q$ and the definition of $\rho$, the function $f$ is
injective on all balls $B(f^j(x),\rho(f^j(x)))$, $j\ge 0$, and therefore $f^k$ is injective on the set $B_n(x,\rho)$ for every $0\le k\le n-1$.
Now, let 
$$
x\in F(\e)\cap X_\rho
$$ 
and let $k$ be the greatest subscript such that
$q=n_k(x)\le n-1$. Denote by $f_x^{-q}$ the unique holomorphic inverse branch
of $f^q$ produced by Corollary~\ref{c9.2.4} which sends $f^q(x)$ to $x$. Of course
$$
B_n(x,\rho)\sbt f^{-q}(B(f^q(x),\rho(f^q(x))))
$$ 
and since $f^q$ is injective on $B_n(x,\rho)$, we even have that
$$
B_n(x,\rho)\sbt f_x^{-q}(B(f^q(x),\rho(f^q(x)))).
$$
By Corollary~\ref{c9.2.4}
$$
\diam\(f_x^{-q}(B(f^q(x),\rho(f^q(x))))\)\le K\exp(- 
q(\chi_\mu-\e)).
$$
Since by \eqref{9.4.3}, $n\le q(1+\e)$ we finally deduce that
$$
B_n(x,\rho)\sbt B\biggl(x,K\exp\biggl(-n{\chi_\mu-\e\over 1+\e}\biggr)\biggr).
$$
Thus, in view of \eqref{9.4.8},
$$
B\lt(x,K\exp\biggl(-n{\chi_\mu-\e\over 1+\e}\biggr)\rt) \spt \Pa_0^n(x).
$$
Therefore, denoting by $r_n$ the radius of the ball above, it follows
from Theorem~\ref{SMBTEC} (ergodic case of the
Shannon--McMillan--Breiman Theorem) that for $\mu$-a.e $x\in X$
$$
\limsup_{n\to\infty}-{1\over n}\log\mu(B(x,r_n)\le \hmu(f,\Pa)\le \hmu(f).
$$
So
$$
\limsup_{n\to\infty}{\log\mu(B(x,r_n)\over \log r_n}\le {\hmu(f) \over
\chi_\mu(f)-\e}(1+\e).
$$
Now, since $\{r_n\}_{n=1}^\infty$ is a geometric sequence and since $\e$ can be taken arbitrarily small, we conclude that for $\mu$--a.e. $x\in X$ we have that
$$
\limsup_{n\to\infty}{\log\mu(B(x,r)\over \log r}\le {\hmu(f)\over
\chi_\mu(f)}.
$$
This completes the proof of \eqref{9.4.5} and because of \eqref{9.4.1}
also the proof of the first part Theorem~\ref{t20120801}. The second part is now an immediate consequence of the first part and Proposition~\ref{p320191130}.
\epf

\section{General Notion of Conformal Measures}\label{generalconformalmeasures}

In this section we encounter for the first time the concept of conformal measures. As we will soon see below, it is closely related to the notion of quasi--invariance and it will be a central theme for the most of the rest of our book. Conformal measures were first defined and introduced by Samuel Patterson in in his seminal paper \cite{Pat1} (see also \cite{Pat2}) in the context of Fuchsian groups. Dennis Sullivan extended this concept to all Kleinian groups in \cite{Su1}--\cite{Su3}. He then, in the papers \cite{Su4} --\cite{Su6}, defined conformal measures for all rational functions of the Riemann sphere $\oc$; he also proved their existence therein. Both Patterson
and Sullivan came up with conformal measures in order to get an understanding of geometric measures, i.e. Hausdorff and packing ones. Although already Sullivan noticed that there are conformal measures for Kleinian groups that are not equal, nor even equivalent to any Hausdorff and packing (generalized) measure, the main purpose to deal with them is still to understand Hausdorff and packing measures but goes beyond. Next chapter, Chapter~\ref{Markov-systems}, Graph Directed Markov Systems, and Part~\ref{EFA}, Elliptic Functions A, and, especially, Part~\ref{EFB}, Elliptic Functions B, of our book, provide a good evidence. 

Conformal measures, in the sense of Sullivan have been studied in greater detail in \cite{DU2}, where, in particular, the structure of the set of their exponents was examined. This is the theme of the next section, i.e. Section~\ref{Sullivan_Conformal_Measures}, Sullivan's Conformal Measures. The general concept of conformal measures has been introduced in \cite{DU1} and this is the theme of the next section. 

Since then conformal measures in the context of rational functions have been studied in numerous research works. We list here only very few of them appearing in the early stages of the development of their theory: \cite{DU LMS}, \cite{DU3}, \cite{DU4}. Subsequently the concept of conformal measures, in the sense of Sullivan, has been extended to countable alphabet iterated functions systems in \cite{MU1} and to conformal graph directed Markov systems in \cite{MU2}. These are treated at length in Chapter~\ref{Markov-systems}, Graph Directed Markov Systems. It was furthermore  extended to transcendental meromorphic dynamics in \cite{KU1}, \cite{UZ1}, and \cite{MyU1}; see also \cite{UZ2}, \cite{MyU2}, \cite{BKZ1}, and \cite{BKZ1}. Last, the concept of conformal measures found its place also in random dynamics; we cite only \cite{MSU}.

\sp Consider first an arbitrary, quasi--invariant Borel measure $m$ for a Borel measurable map $T:X\lra X$, where $X$ is a metrizable space. Assume that $T$ is (at most) countable--to--one, i.e. 
$$
X=\bu_{j\in I}X_j,
$$
where $X_j$, $j\in I$, are measurable, pairwise disjoint sets, and for each $j\in I$, the map 
$$
T|_{X_j}:X_j\lra T(X_j)
$$ 
is a measurable isomorphism. Recall that 
$$
\L_m:L_m^1\lra L_m^1
$$ 
is the transfer operator corresponding to the quasi--invariant measure $m$, defined in Section~\ref{QEC} by formula \eqref{1_mu_2014_11_15}. We now can introduce the measure $\L_m^*m$ by the formula
\beq\label{2_mu_2014_11_15}
\L_m^*m(\psi)=m(\L_m\psi).
\eeq
Putting in the formula \eqref{1j157} $g:=\1_X$ and $f=\psi\in L_m^1$, we get from this formula that
\beq\label{3_mu_2014_11_15}
 \int_{X}\psi\,d\L_m^*m = \int_X\psi\,dm.
\eeq
This precisely means that
\beq\label{3a_mu_2014_11_15}
\L_m^*m=m.
\eeq
Denote 
\beq\label{5_mu_2014_11_15}
J_m^{-1}(T):=\frac{d(m\circ (T|_{X_j})^{-1})}{dm}.
\eeq
Obviously, for all $u\in L_m^1$ we have that
\beq\label{4_mu_2014_11_15}
\L_m(u)(x)=\L_{\log J_m^{-1}(T)} (x):=\sum_{y\in T^{-1}(x)}u(y)J_m^{-1}(T)(y),
\eeq
We say that $A\sbt X$ is a special set\index{(N)}{special set} if $A$ is a Borel set and the restriction $T|_A$ is injective.    
Assume now in addition that 

\centerline{$T:X\lra X$ is non--singular with respect to $m$,}

\fr meaning that 
$$
m\(T((J_m(T)^{-1}(+\infty))\)=0
$$ 
(note that $m\(T((J_m(T)^{-1}(0))\)=0$ merely by quasi--invariance). 
Hence, if $A\sbt X$ is a special set, then this formula along with \eqref{3_mu_2014_11_15} and \eqref{4_mu_2014_11_15}, yield
\beq\label{5_mu_2014_11_15B}
\begin{aligned}
\int_A J_m(T)\,dm
&=\int_X\1_AJ_m(T)\,dm
=\int_X\1_AJ_m(T)\,d\L_m^*m
=\int_X\L_m\(\1_AJ_m(T)\)\,dm \\
&=\int_X\lt(\sum_{y\in T^{-1}(x)}\1_A(y)J_m(T)(y)J_m^{-1}(T)(y)\rt)\,dm(x) \\
&=\int_{T(A)}\lt(\sum_{y\in T^{-1}(x)}\1_A(y)\rt)\,dm(x) \\
&=\int_{T(A)}\1\,dm 
=m(T(A)).
\end{aligned}
\eeq
We want to single out this property by saying that the non--singular measure $m$ is $J_m(T)$-conformal. This prompts us to introduce the following.

\sp\bdfn\label{general_conformal_measure}
Let $T:X\lra X$ be a Borel measurable map of a metric space $(X,\rho)$
and let $g:X\lra\R$ be a non--negative measurable function.
A Borel probability measure $m$ on $X$ is said to be $g$--{\it
conformal}\index{(N)}{conformal measure} for $T:X\lra X$ if and only if 
\beq\label{10.1.1}
m(T(A))=\int_A g\,dm 
\eeq
for every special set $A\sbt X$
\edfn

\sp\bobs
Observe that, on the other hand, if, in the above definition, $g>0$ then the $g$--conformal measure $m$ is  non-singular with respect to the map $T$ and
$$
J_m^{-1}(T)=1/g.
$$
\eobs

\

\fr Notice that even if $T$ is continuous, the operator $\L_m$ need not in general map $C(X)$ into $C(X)$. But this is the case if, in addition to being continuous, the map 
$T$ is also open. Nevertheless, if we just assume that $\L_m(C(X))\sbt
C(X)$, then the linear operator $\L_m:C(X)\to C(X)$ is bounded, and so is its dual
$\L_m^*:C^*(X)\to C ^*(X)$.  Then, under the above constrains: $T$ being 
countable -- to -- $1$ and $g$ being positive, we get the following.

\bprop\label{prop10.1.1}
A probability measure $m$ is $g$--conformal if and only if
$$
\L^*_{-\log g}(m)=m.
$$
\eprop

\fr Therefore, since we can have troubles with the operator $\L^*$ for the
maps $T$ that are not open,  rather than looking for theorems that would assure the existence of fixed points of $\L^*$, we shall provide another 
general method of constructing conformal measures, called
Patterson--Sullivan method. This construction will make use of the
following simple fact. For a sequence $\(a_n\)_{n=1}^\infty$ of reals, the
number
\beq\label{10.1.3}
c:=\limsup_{n\to\infty}\frac{a_n}{n}  
\eeq
will be called the transition parameter of the sequence $\(a_n\)_{n=1}^\infty$. It is uniquely determined by the property that
$$
\sum_{n\ge 1}\exp(a_n-ns)
$$
converges for $s>c$ and diverges for $s<c$. For $s=c$ the sum may converge or
diverge. By a simple argument one obtains the following.

\blem\label{10.1.2}
There exits a sequence  $\(b_n\)_{n=1}^\infty$ of positive reals such that
$$
\sum^\infty _{n=1}b_n \exp(a_n-ns) \begin{cases} <\infty & \text{if } s>c\\
=\infty & \text{if } s\le c \end{cases}
$$
and 
$$
\lim_{n\to\infty}\frac{b_n}{b_{n+1}}=1.
$$
\elem

\bpf If 
$$
\sum_{n\ge 1}\exp(a_n-nc)=\infty,
$$
put $b_n=1$ for every $n\ge 1$. If 
$$
\sum_{n\ge 1}\exp(a_n-nc)<\infty,
$$
choose a sequence $\(n_k\)_{k=1}^\infty$ of positive integers such that
$$
\lim_{k\to\infty}\frac{n_k}{n_{k+1}}=0 \  \  \text{ and } \  \
\e_k:=a_{n_k}n_k^{-1}-c\longrightarrow 0.
$$
Setting then
$$ 
b_n=\exp\lt(n\lt(\frac {n_k-n}{n_k-n_{k-1}}
\e_{k-1}+\frac{n-n_{k-1}}
{n_k-n_{k-1}}\e_k\rt)\rt)\ \ \ \ \text{for}\ n_{k-1}\le n<n_k,
$$
it is easy to check that the lemma follows. 
\epf

\sp Getting back to dynamics, let
$\{E_n\}_{n=1}^\infty$ be a sequence of finite subsets of $X$ such that
\beq\label{10.1.4}
T^{-1}(E_n)\subset E_{n+1}
\eeq
for every $n\ge 1$. Let $\phi:X\lra\R$ be an arbitrary bounded fucnction. Let
\beq\label{520181213}
a_n:= \log\lt(\sum_{x\in E_n}\exp (S_n \phi(x))\rt)
\eeq
where, we recall, 
$$
S_n \phi=\sum_{k=0}^{n-1} \phi\circ T^k.
$$
Denote by $c$ the transition parameter of this sequence. Choose
a sequence $(b_n)_1^\infty$ of positive reals as in Lemma~\ref{10.1.2} for the
sequence $(a_n)_1^\infty$. For every $s>c$ define
\beq\label{10.1.5}
M_s:=\sum_{n=1}^\infty b_n\exp(a_n-ns)  
\eeq
and the normalized measure
\beq\label{10.1.6}
m_s=\frac 1{M_s}\sum_{n=1}^\infty\sum_{x\in E_n}b_n
\exp(S_n\phi(x)-ns)\d_x, 
\eeq
where $\d_x$ denotes the Dirac $\d$ measure supported at the point $x\in X$.
Let $A$ be a special set. Using \eqref{10.1.4} and \eqref{10.1.6}, it follows that
\begin{eqnarray}\label{10.1.7}
m_s(T(A)) &= &\frac 1{M_s}\sum_{n=1}^\infty\sum_{x\in E_n\cap T(A)}
b_n\exp(S_n\phi(x)-ns)\nonumber \\
&=& \frac 1{M_s}\sum_{n=1}^\infty\sum_{x\in A\cap T^{-1}E_n}
b_n\exp(S_n\phi(T(x))-ns)\nonumber  \\
&=& \frac 1{M_s}\sum_{n=1}^\infty\sum_{x\in A\cap E_{n+1}}
b_n\exp[S_{n+1}\phi(x)-(n+1)s]\exp(s-\phi(x))-\nonumber \\
&& \  \  \  \  \  \  \  \  \  \  \  \ - \frac 1{M_s}\sum_{n=1}^\infty\sum_{x\in A\cap
(E_{n+1}\sms T^{-1}E_n)} b_n\exp(S_n\phi (T(x))-ns). 
\end{eqnarray}
Set
$$
\begin{aligned}
\Delta_A(s)
:&= \bigg|\frac 1{M_s}\sum_{n=1}^\infty\sum_{x\in A\cap E_{n+1}}
b_n\exp[S_{n+1}\phi(x)-(n+1)s]\exp(s-\phi(x))- \\ 
& \  \  \  \   \  \  \  \  \  \  \  \ \  \  \  \  \  \  \  \ -\int_A\exp(c-\phi)\,
dm_s\bigg|
\end{aligned}
$$
and observe that
$$
\aligned
\Delta_A(s)
={}&\frac 1{M_s}\Big|\sum_{n=1}^\infty\sum_{x\in A\cap E_{n+1}}
\exp[S_{n+1}\phi(x)-(n+1)s]\exp(-\phi(x))\bigl[b_ne^s-b_{n+1}e^c
\bigr] - \\
& \  \  \  \   \  \  \  \   \  \  \  \   \  \  \  \ -b_1\sum_{x\in A\cap E_1}e^{c-s}\Big|\le  \\
&\le \frac 1{M_s}\sum_{n=1}^\infty\sum_{x\in A\cap E_{n+1}}
\bigg|\frac {b_n}{b_{n+1}}-e^{c-s}\bigg| b_{n+1}\exp(s-\phi(x))\exp
[S_{n+1}\phi(x)-(n+1)s] + \\
&  \  \  \  \   \  \  \  \   \  \  \  \   \  \  \  \ +\frac 1{M_s}b_1\exp(c-s)\,\sharp (A\cap E_1) \le \\
&\le \frac 1{M_s}\sum_{n=1}^\infty\sum_{x\in E_{n+1}}
\bigg|\frac {b_n}{b_{n+1}}-e^{c-s}\bigg| b_{n+1}\exp(s-\phi(x))\exp
[S_{n+1}\phi(x)-(n+1)s] + \\
& \  \  \  \   \  \  \  \   \  \  \  \   \  \  \  \ + \frac 1{M_s}b_1\exp(c-s)\,\sharp  E_1.
\endaligned
$$
By Lemma~\ref{10.1.2} we have $\lim_{n\to\infty}b_{n+1}/b_n=1$ and $\lim_{
s\downto c}M_s=\infty$. Therefore
\beq\label{10.1.8}
\lim_{s\downto c}\Delta_A(s)=0  
\eeq
uniformly with respect to all special sets $A$.

\sp\bdfn\label{d120191130}
Any weak accumulation point, when $s\searrow c$, of the measures $(m_s)_{s>c}$, defined by \eqref{10.1.6}, will be called a limit measure (associated to the function $\phi$ and the
sequence $(E_n)_1^\infty$).
\edfn

\sp\fr In order to find conformal measures among the limit
measures, it is necessary to examine \eqref{10.1.7} in greater detail. To begin with, for a Borel set $D\sbt X$, consider the following condition:
\beq\label{10.1.9}
\lim_{s\downarrow c}\frac 1{M_s}\sum_{n=1}^\infty
\sum_{x\in D\cap (E_{n+1}\sms T^{-1}E_n)} b_n
\exp[S_n\phi(T(x))-ns]=0. 
\eeq
We will need the following definitions. A point $x\in X$ is said to be
singular for $T$ 
if at least one of the following two conditions is satisfied:
\beq\label{10.1.10}  \ \text{There is no open neighbourhood}\;\; U\;\;
\text{of}\;\; 
x \;\; \text{such that}\;\; 
T|_U \;\;\text{is injective}. 
\eeq
\beq\label{10.1.11}
\forall_{\e >0}\  \exists_{0<r<\e}\;\;\text{such that}\;\; T(B(x,r))\;\;
\text{is not an open subset of}\;\; X. 
\eeq
The set of all singular points is denoted by $\Sing(T)$, the
set of all points satisfying condition \eqref{10.1.10} is denoted by
$\Crit(T)$ and the
set of all points satisfying condition \eqref{10.1.11} is denoted by
$X_0(T)$.

\sp\fr It is easy to give examples where $X_0(T)\cap \Crit(T)\ne \es$.
If $T:X\lra X$ is an open map, no point satisfies condition
\eqref{10.1.11} that is $X_0(T)=\es$.

\blem\label{l10.1.3} 
Let $T:X\lra X$ be a Borel measurable map of a metric space $(X,\rho)$
and let $g:X\lra\R$ be a non--negative Borel measurable function.
Let $m$ be a Borel probability measure on
$X$ and let $\Ga$ be a compact set containing $\Sing(T)$.
If $g$ is integrable with respect to $m$ and \eqref{10.1.1} holds for every special set $A$ whose closure is disjoint from $\Ga$ and such that 
$m(\bd A)=m(\bd T(A))=0$, then \eqref{10.1.1} continues to hold for
every special set $A$ disjoint from $\Ga$. 
\elem

\bpf Let $A$ be a special set disjoint from $\Ga$.
Fix $\e>0$. Since on the complement of $\Ga$ the map $T$ is open, for each
point $x\in A$ there exists an open neighborhood $U(x)$ of $x$ such that
$T|_{U(x)}$ is a homeomorphism, 
$$
m(\bd U(x))=m(\bd T(U(x)))=0, \  \  \ov{U(x)}\cap\Ga=\es
$$ 
and such that
$$
\int_{\bu_{x\in A} U(x)\sms A}\,g\,dm <\e.
$$
Choose a countable family $\{U_k\}_{k=1}^\infty$ from $\{U(x)\}_{x\in
A}$ which covers $A$ and define recursively $A_1:=U_1$ and 
$$
A_n:=U_n\sms \bu_{k<n}U_k.
$$
By the hypotheses of the lemma, each set $A_k$ satisfies
\eqref{10.1.1}, and hence
$$
\aligned
m(T(A)) &=m\Big(\bu_{k=1}^\infty T(A\cap A_k)\Big)
\le\sum_{k=1}^\infty m(T(A_k)) \\
&=\sum_{k=1}^\infty\int_{A_k}\,g\,dm=\int_A\,g\,dm+
\sum_{k=1}^\infty \int_{A_k\sms A}\,g\,dm \\
&\le \int_A\,g\,dm+\e.
\endaligned
$$
If $\e\downto 0$, it follows that
\beq\label{1_mu_2014_11_17}
m(T(B))\le \int_B\,g\,dm
\eeq
for any special set $B$ disjoint from $\Ga$.
Using this fact, the lower bound for $m(T(A))$ is obtained
from the following estimate, if $\e \downto 0$:
$$
\aligned
m(T(A))  &=m\Big(\bu_{k=1}^\infty T(A\cap A_k)\Big)
= \sum_{k=1}^\infty m(T(A\cap A_k)) \\
&= \sum_{k=1}^\infty \lt(m(T(A_k))-m(T(A_k\sms A))\rt)
\ge \sum_{k=1}^\infty\int_{A_k}g\,dm- \int_{A_k\sms
A}g\,dm \\
& = \int_{\cup_{k\ge 1}A_k}g\,dm- \int_{\cup_{k\ge 1}A_k
\sms A}g\,dm \ge\int_Ag\,dm -\e. 
\endaligned
$$
Letting $\e\downto 0$, we thus get that $m(T(A)\ge \int_Ag\,dm$. Along with \eqref{1_mu_2014_11_17}, this gives that
$$
m(T(A)=\int_Ag\,dm.
$$
The lemma is proved. 
\epf

\sp\blem\label{l10.1.4} 
Let $T:X\lra X$ be a Borel measurable map of a metric space $(X,\rho)$.
Let $\phi:X\lra\R$ be a bounded Borel measurable function. Let $m$ be a limit measure of Definition~\ref{d120191130}, and let $\Ga$ be a compact set containing$\Sing(T)$. Assume that every special set $D\sbt X$ with
$m(\bd D)=m(\bd T(D))=0$ and $\bar D\cap \Ga=\es$ satisfies condition
\eqref{10.1.9}. Then 
$$
m(T(A))=\int_A\exp(c-\phi)\,dm
$$ 
for every special set $A$ disjoint from $\Ga$.
\elem

\bpf Let $D\sbt X$ be a
special set such that $\bar D\cap \Ga=\es$ and
$m(\bd D)=m(\bd T(D))=0$. It follows immediately from
\eqref{10.1.7}--\eqref{10.1.9} that
$$
m(T(D))=\int_D\exp(c-\phi)\,dm. 
$$
Applying now Lemma~\ref{l10.1.3} completes the proof. 
\epf

\blem\label{l10.1.5} Let $T:X\lra X$ be a Borel measurable map of a metric space $(X,\rho)$. Let $\phi:X\lra\R$ be a bounded Borel measurable function. Let $m$ be a limit measure of Definition~\ref{d120191130}. If condition \eqref{10.1.9} is satisfied for $D=X$, then
$$
m(T(A))\ge \int_A \exp(c-f)\,dm
$$
for every special set $A$ disjoint from $\Crit(T)$.
\elem

\bpf Suppose first that $A$ is compact and $m(\bd A)=0$. From
\eqref{10.1.7}, \eqref{10.1.8} and the assumption one obtains
$$
\lim_{s\in J}\Big| m_s(T(A))-\int_A\exp(c-\f)\,dm_s\Big| =0
$$
where $J$ denotes the subsequence along which $m_s$ converges to $m$.  Since
$T(A)$ is compact, this implies that
$$
m(T(A))\ge\liminf_{s\in J}m_s(T(A))=\lim_{s\in J}
\int_A\exp(c-\f)\,dm_s=\int_A\exp(c-\f)\,dm
$$
Now, drop the assumption $m(\bd A)=0$ but keep $A$ compact and assume
additionally that for some $\e>0$ the ball $B(A,\e)$ is also special. Choose a
descending sequence $A_n$ of compact subsets of $B(A,\e)$ whose intersection
equals $A$
and $m(\bd A_n)=0$ for every $n\ge 0$. By what has been already proved
$$
m(T(A))=\lim_{n\to\infty}m(T(A_n))\ge \int_{A_n}\exp(c-\f)\,dm=\int_A\exp(c-
\f)\,dm 
$$

The next step is to prove the lemma for $A$, an arbitrary open special set
disjoint from $\Crit(T)$, by partitioning it into countably many compact sets. Then one approximates from above special sets of sufficiently small diameters by special open sets and the last step is to partition an arbitrary special set disjoint from $\Crit(T)$ by sets of so small diameters that the lemma holds.
\epf

\blem\label{l10.1.6}  Let $T:X\lra X$ be a Borel measurable map of a metric space $(X,\rho)$. Let $\Ga$ be a compact subset of $X$ containing
$\Sing(T)$. Suppose that for every integer
$n\ge 1$ there are a continuous function $g_n:X\to X$ and a
measure $m_n$ on $X$ satisfying the following conditions.
\begin{itemize}
\item[(a)] \eqref{10.1.1} holds for $g=g_n$ and for every special
set $A\sbt X$ with
$$
\ov A\cap \Ga =\es.
$$
\item[(b)] 
$$
m_n(T(B))\ge \int_B g_n \,dm_n
$$
for any special set $B\sbt X$ such that $B\cap \Crit(T)=\es$.

\,

\item[(c)] The sequence $\{g_n\}_{n=1}^\infty$ converges uniformly
to a continuous function $g:X\to\R$. 
\end{itemize}
Then for any weak accumulation point $m$ of the sequence
$\{m_n\}_{n=1}^\infty$ we have 
\begin{itemize}
\item[(1)]
$$
m(T(A))=\int_A g\,dm     
$$
for all special sets $A\sbt X$ such that $A\cap \Ga =\es$ and
\item[(2)]
$$
m(T(B))\ge \int_B g\,dm   
$$
for all special sets $B\sbt X$ such that $B\cap \Crit(T)=\es$.
\end{itemize}
Moreover, if (a) is replaced by
\begin{itemize}
\item[(a')]
$$
\ov A\cap (\Ga\sms(\Crit(T)\sms X_0(T)))=\es, 
$$
then for any $x\in \Crit(T)\sms X_0(T)$ 
\end{itemize}
\begin{itemize}
\item[(3)]
$$
m(\{T(x)\})\le g(x)m(\{x\})\le q(x)m(\{T(x)\})  
$$
\end{itemize}
where $q(x)$ denotes the maximal number of preimages of single points under
the transformation $T$ restricted to a sufficiently small neighbourhood of $x$.
\elem

\sp\fr The proof of property (1) is a simplification of the proof of
Lemma~\ref{l10.1.4} and
the proof of property (2) is a simplification of the proof of
Lemma~\ref{l10.1.5}. The 
proof of (3) uses the same techniques and is left for the reader.

\sp\section{Sullivan's Conformal Measures}\label{Sullivan_Conformal_Measures}

This section is devoted to a short study of more special conformal measures, called Sullivan's conformal measures. An extended historical and motivational discussion of these measures and more general ones as well was given in the previous section, Section~\ref{generalconformalmeasures}, General Conformal Measures. We therefore start right now with actual mathematics. 

In this section as in Section~\ref{Volume_Lemma} let $Y$ be a Riemann surface and let $X$ be a compact subset of $Y$. Assume that $f   
\in {\mathcal A}(X)$. We now shall provide a construction aiming to
establish the existence of measures that will be called (Sullivan's) $t$--conformal
measures \index{(N)}{$t$-conformal measure}. In fact we will prove the existence of measures with some slightly
weaker properties. We will call them semi $t$--conformal Sullivan's
measures. To begin, fix  $z\in X$ and for all $n\ge 0$ set 
$$
E_n:=f^{-n}(z).
$$ 
Then
$$
E_{n+1}=f^{-1}(E_n)
$$ 
and therefore the sequence $\{E_n\}_{n=1}^\infty$ satisfies
\eqref{10.1.9} with all sets $D\sbt X$. Fix an arbtrary $t\ge 0$ and let $c(t)$ be the transition parameter associated to this sequence and the function 
$$
\phi:=-t\log|f'|
$$ 
according to \eqref{10.1.3} and \eqref{520181213}. As an immediate consequence of Lemma~\ref{l10.1.4} and Lemma~\ref{l10.1.5}, we get the following. 

\blem\label{l120120827} 
Let $Y$ be a Riemann surface and let $X$ be a compact subset of $Y$. If $f\in {\mathcal A}(X)$ and $t\ge 0$, then there exists a Borel
probability measure $m_t$ on $X$ such that
$$
m_t(f(A)\ge \int_A e^{c(t)}|f'|^t\,dm_t
$$
for every special set $A\sbt X$, and 
$$
m_t(f(A)=\int_A e^{c(t)}|f'|^t\,dm_t
$$
if in addition $A\cap X_0(f)=\es$.
\elem

\sp We now want to establish some useful properties of the function
$c(t)$ and to get a parameter $t\ge 0$ such that $c(t)=0$. Let
$$
\P(t):=\P(f,-t\log|f'|)
$$
be the topological pressure of the potential $-t\log|f'|$ with respect
to the dynamical system $f:X\lra X$. We shall prove the following.

\sp\blem\label{l10.3.4} 
For every $t\ge 0$ we have that $c(t)\le \P(t)$.
\elem

\bpf Since the map $f:X\to X$ has no critical points, it is
locally 1-to-1 at all points of $X$. Since $X$ is compact, this means that 
there exists $\d>0$ such that the map $f$ restricted to any set with
diameter $\le \d$ is 1-to-1. Consequently, all the sets $E_n$ are
$(n,\e)$-separated for all $0<\e<\d$. Hence, the required inequality
$c(t)\le \P(t)$ follows immediately from Theorem~\ref{ch9sem2thm9.12}. 
\epf

\sp The standard straightforward convexity arguments showing
continuity of topological pressure prove also the following.

\blem\label{l10.3.5} 
The function $[0,+\infty)\ni t\longmapsto c(t)\in\R$ is
continuous.
\elem

\fr Set
\beq\label{120181217}
s(f)=\inf\{t\ge 0: c(t)\le 0\}.
\eeq
We call the dynamical dimension \index{(N)}{dynamical dimension} of $X$ the number 
$$
\DD(X):=\sup\{\HD(\mu)\},
$$
where the supremum is taken over all Borel probability $f$--invariant ergodic
measures $\mu$ on $X$. We shall prove the following.

\blem\label{l10.3.6} 
Let $Y$ be a Riemann surface and let $X$ be a compact subset of $Y$. If $f\in {\mathcal A}(X)$, then 
$$
s(f)\le \DD(X)\le\HD(X)\le 2.
$$
\elem

\bpf Suppose on the contrary that $\DD(X)<s(f)$ and take
$0\le\DD(X)<t<s(f)$. From this choice and by Lemma~\ref{l10.3.4} we have that
$$
0<c(t)\le \P(t),
$$ 
and, by the Variational Principle,
Theorem~\ref{sem2thm9.13}, there exists an ergodic $f$-invariant Borel
probability measure $\mu$ on $X$ such that 
$$
\P(t)\le \hmu(f)-t\chi_\mu(f)+c(t)/2.
$$
Therefore, by Corollary~\ref{c10.3.2} we get
$\hmu(f)\ge c(t)/2>0$ 
and applying in adition Theorem~\ref{t9.1.1} (Ruelle's inequality),
$\chi_\mu(f)>0$. Hence, it follows from Theorem~\ref{t20120801} that
$$
t\le \HD(\mu)-{1\over 2}{c(t)\over \chi_\mu(f)}<\HD(\mu)\le \DD(X).
$$
This contradiction finishes the proof. 
\epf 

\

\fr Since $c(0)\ge 0$, as an immediate consequence of Lemma~\ref{l10.3.5}, we get that
\beq\label{120120827}
c(s(f))=0.
\eeq
Inserting this to Lemma~\ref{l120120827}, we get the following main
result of this section.

\blem\label{l220120827} 
Let $Y$ be a Riemann surface and let $X$ be a compact subset of $Y$. If $f\in {\mathcal A}(X)$ then there exists a Borel
probability measure $m$ on $X$ such that
$$
m(f(A)\ge \int_A|f'|^{s(f)}\,dm
$$
for every special set $A\sbt X$, and 
$$
m(f(A)=\int_A|f'|^{s(f)}\,dm
$$
if in addition $A\cap X_0(f)=\es$.
\elem

\bdfn\label{dsemiconformal}
Any measure with the former property above (with some parameter $t$ in place of $s(f)$) will be called (Sullivan's) semi $t$--conformal \index{(N)}{semi $t$--conformal measure}. If in addition the latter property (with $s(f)$ replaced by $t$) holds for all special sets $A\sbt X$, then the measure $m$ is called (Sullivan's) $t$--conformal \index{(N)}{$t$--conformal measure}.
\edfn

\section{Conformal Pairs of Measures}\label{conformal pairs of measures}

In this section we examine in detail how Sullivan's conformal
measures, or rather Sullivan's like conformal measures, behave under
transformations, particularly those having critical points. We will
eventually use them for a one given conformal measure but in order to
see clearly what is going on and what matters, it is most natural to
present these results in the setting of conformal and semi conformal
pairs of measures. The results of this section, quite technical, are
all taken from Section~2, Section~3, and Section~4 of
\cite{U1}. Motivated by the last definition in the previous section we
formulate the following.

\sp
\bdfn\index{(N)}{conformal pairs of measures}\label{1d20120909}
Fix $t\ge 0$. Let $G$ and $H$ be non-empty open subsets of $\oc$. Let
$f:G\lra H$ be a meromorphic map. A pair ($m_G$, $m_H$) of Borel finite
measures on $G$ and $H$ respectively is called 
spherical semi $t$--conformal\index{(N)}{spherical semi $t$--conformal
pair of measures} for the map $f:G\lra H$, if and only if
\beq\lab{semconf}
m_H(f(A))\ge \int_A|f^*|^t\,dm_G
\eeq
for every Borel set $A\sbt G$ such that $f|_A$ is injective, and this pair
is said to be a spherical $t$--conformal\index{(N)}{spherical
$t$--conformal measure} for $f:G\to H$ if and only if
\beq\lab{dc1.2}
m_H(f(A))=\int_A|f^*|^t\,dm_G
\eeq
for these sets $A$. 

\begin{itemize}
\sp\item[(a)] If $G\sbt\C$ and in the two formulas above the
spherical derivative $f^*$ is replaced by the Euclidean derivative
$f'$, the pair ($m_G$, $m_H$) is called Euclidean (semi) $t$-conformal
for the map $f:G\lra H$. 

\sp\item[(b)]  If it does not matter whether a spherical (semi)
conformal pair of measures or a Euclidean (semi) conformal pair of
measures is involved, we frequently simply speak about (semi)
conformal pair of measures. 

\sp\item[(c)]  If both measures $m_G$ and $m_H$ are restrictions of the same Borel finite measure $m$ defined on $G\cup H$, we refer to $m$ as a (semi) $t$--conformal measure for the map $f:G\lra H$.
\end{itemize}
\edfn

\sp As long as we are far from poles and infinity it does not really
matter weather our pair of measures is spherical (semi) conformal or
Euclidean (semi) conformal. This is explained by the following.

\sp\bobs\label{o120120908}
If ($m_G$, $m_H$) is a Euclidean (semi) $t$--conformal pair of measures
for a meromorphic map $f:G\lra H$, then ($m_G^*$, $m_H^*$) is a spherical
$t$--conformal pair of measures for $f$, where the measures $m_G^*$ and $m_H^*$ are respectively defined  by
\beq\label{1p6a}
\frac{dm_G^*}{dm_G}(z)=(1+|z|^2)^{-t}
\eeq
and 
\beq\label{1p6b}
\frac{dm_H^*}{dm_H}(z)=(1+|z|^2)^{-t}.
\eeq
Similar conversion holds if we start from a spherical $t$--conformal
pair of measures.
\eobs

\sp Notice that in the context of item (c) of Definition~\ref{1d20120909} the $\sigma$--finite  measure
$m_e$\index{(S)}{$m_e$} related to the spherical one  by the formula 
\beq\label{1p6}
\frac{dm_e}{dm_s}(z)=(1+|z|^2)^t
\eeq\index{(S)}{$\frac{dm_e}{dm_s}$}
has the property that
$$m_e(f(A))=\int_{A}|f'|^tdm_e$$
for  each special set $A$. It will be called  an Euclidean
$t$-conformal measure. 
In the  sequel we respect the convention that the spherical
conformal measure (or their weaker versions) are labeled with the
subscript '$s$'  whereas Euclidean  conformal  measures (and their
weaker versions) are labeled  with the subscript '$e$'. \,If no \,
subscript is used, the  conformal measure under consideration can be
spherical equally well as Euclidean. 

\sp\fr We provide in  this section some more technical facts taken from
Section~2, Section~3, and Section~4 of \cite{U1}.

\

\bdfn\lab{dncp12.3.} Given $t>0$, $r>0$, and $L>0$, a point $x\in {\mathbb
C}$ is said to be $(r,L)$--$t$--upper estimable \index{(N)}{point is
$(r,L)$--$t$--upper estimable} with respect to a finite Borel measure
$\nu$ if 
$$
\nu(B(x,r))\le Lr^t
$$ 
and the point $x$ is said to be 
$(r,L)$--$t$--lower estimable \index{(N)}{point is $(r,L)$--$t$--lower
estimable} with respect to $\nu$ if 
$$
\nu(B(x,r))\ge Lr^t.
$$
We will frequently abbreviate the notation writing
$(r,L)$--$t$--u.e. \index{(S)}{$(r,L)$--$t$--u.e.} for $(r,L)$--$t$--upper estimable
and $(r,L)$-$t$-l.e.\index{(S)}{$(r,L)$--$t$--l.e.} for $(r,L)$--$t$--lower
estimable. 

\sp We also say that the point $x$ is $t$--upper
estimable\index{(N)}{point is $t$--upper estimable} ($t$--lower
estimable\index{(N)}{point is $t$--lower estimable}) if it is
$(r,L)$--$t$--upper estimable {\rm (}$(r,L)$--$t$--lower estimable) for
some $L>0$ and all $r>0$ sufficiently small. 
\edfn

\sp

\bdfn\lab{dncp12.4.} 
Given $r>0$, $\sg>0$ and $L>0$ the point $x\in
X$ is said to be $(r,\sg,L)$--$t$--strongly lower
estimable\index{(N)}{point is $(r,\sg,L)$-$t$-strongly lower
estimable}, or shorter
$(r,\sg,L)$--$t$--s.l.e.\index{(S)}{$(r,\sg,L)$--$t$--s.l.e.} with respect
to a finite Borel measure $\nu$ if 
$$
\nu(B_e(y,\sg r))\ge Lr^t
$$ 
for every $y\in B_e(x,r)$. 
\edfn

\sp

\blem\lab{lncp12.5.} If $z$ is $(r,\sg,L)$--$t$--s.l.e. with respect
to a finite Borel measure $\nu$, then every
point $x\in B_e(z,r/2)$ is $(r/2,2\sg,2^tL)$--$t$--s.l.e. with respect
to this measure. 
\elem

\bpf Let $x\in B_e(z,r/2)$. Then $x\in B_e(z,r)$ and
therefore 
$$
\nu(B_e(x,2\sg({r/ 2})))=\nu(B_e(x,\sg r))\ge
Lr^t=2^tL(r/2)^t.
$$ 
\endpf

\sp\blem\lab{lncp12.6.} If $x$ is $(r,\sg,L)$--$t$--s.l.e. with respect
to a finite Borel measure $\nu$, then for
every $0<u\le 1$ it is $(ur,\sg/u,Lu^{-t})$--$t$--s.l.e. with respect
to this measure. 
\elem

\bpf If $y\in B_e(x,ur)$, then $y\in B_e(x,r)$ and
therefore $\nu(B_e(y,({\sg/ u})ur))=\nu(B_e(y,\sg r))\ge
Lr^t=Lu^{-t}(ur)^t.$  
\endpf

\sp We would like to finish this part, not involving transformations,
with the following obvious statement.

\sp\blem\lab{lncp12.7.} 
If $\nu$ is positive on nonempty open sets,
then for every $r>0$ there exists $E(r)\ge 1$ such that every point
$x\in X$ is $(r,E(r))$--$t$--u.e. and $(r,E(r)^{-1})$--$t$--l.e. with respect
to $\nu$. 
\elem

\sp The following lemma, the first one involving dynamics, is a
straightforward consequence of Theorem~\ref{Euclid-I}.

\sp\blem\lab{lncp13.1.} 
Let ($m_G$, $m_H$) be a Euclidean semi $t$-conformal pair of measures
for a univalent conformal map $f:G\lra H$. If $z\in B(z,2R)\sbt G$, then
for every $0\le r\le R$, we have that
$$
m_G(B_e(z,K^{-1}r|f'(z)|^{-1})) \le K^t|f'(z)|^{-t}m_H(B_e(f(z),r)).
$$
If, in addition, the pair ($m_G$, $m_H$) is $t$--conformal, then also
$$
 m_G(B_e(z,Kr|f'(z)|^{-1})) \ge K^{-t}|f'(z)|^{-t}m_H(B_e(f(z),r)).
$$
\elem

\sp\blem\lab{lncp13.2.} 
Let ($m_G$, $m_H$) be a Euclidean semi $t$-conformal pair of measures
for a univalent conformal map $f:G\lra H$. Let $z\in B(z,2R)\sbt G$.

If the point $H(z)$ is $(r,\sg,L)$--$t$--s.l.e. with respect
to $m_H$, where $r\le R/2$ and $\sg\le 1$, then the point $z$ is
$(K^{-1}|f'(z)|^{-1}r,K^2\sg,L)$--$t$--s.l.e. with respect
to $m_G$. 
\elem

\fr{\sl  Proof.} In this proof we apply Lemma~\ref{lncp13.1.}
several times without special indicating. Consider 
$$x
\in B_e(z,K^{-1}r|f'(z)|^{-1}).
$$ 
Then $f(x)\in B_e(f(z),r)$ and
therefore $m_H(B_e(f(x),\sg r))\ge Lr^t$. Since
$$
B_e(f(x),\sg r)\sbt B_e(f(z),2r)\sbt B_e(f(z),R),
$$ 
we have
$$
f_z^{-1}\(B_e(f(x),\sg r)\)
\sbt B_e(x,K\sg r|f'(z)|^{-1})
=B_e(x,K^2\sg( K^{- 1}|f'(z)|^{-1}r)).
$$
Thus 
$$\nu_e\(B_e(x,K^2\sg(K^{-1}|f'(z)|^{-1}r))
\ge K^{-t}|f'(z)|^{-t}Lr^t= L(K^{- 1}|f'(z)|^{-1}r)^t.
$$ 
The proof is finished. \endpf

\sp\blem\lab{lncp13.3.} 
Let ($m_G$, $m_H$) be a Euclidean semi $t$--conformal pair of measures
for an analytic map $f:G\to H$. If $0<r\le R(H,c)$ and $f(c)$ is
$(r,L)$--$t$--l.e. with respect to $m_H$, then $c$ is
$$
((A(c)r)^{1/p_c},(A(c))^{-2t}L)-t-{\rm l.e.}
$$
 with respect to $m_G$, where $A(c):=A(f,c)$ was defined just after
Definition~\ref{dcomp} and $p_c$  is the order of $f$
at the critical point $c$. 
\elem

\bpf By Definition~\ref{dcomp} we get
$B_e(f(c),r)=f(\Comp(c,f,r))$. If $x\in \Comp(c,f,r)$,
then  
$$
(A(c))^{-1}|x-c|^{p_c} \leq  |f(x)-f(c)| <r,
$$
which implies that $x\in B_e(c,(A(c)r)^{1/p_c})$. Thus $B_e(f(c),r)\sbt
f(B_e(c,(A(c)r)^{1/p_c})$, and therefore
$$\begin{aligned}
Lr^t &\le m_H(B_e(f(c),r))
\le m_H\(f(B_e(c,(A(c)r)^{1/p_c}))\)
\le \int_{B_e(c,(A(c)r)^{1/p_c})}|f'(z)|^t\,dm_G(z) \\
&\le \int_{B_e(c,(A(c)r)^{1/p_c})}(A(c))^t(|z-c|^{p_c-1})^t\,dm_G(z) \\
&\le (A(c))^t(Ar)^{{p_c-1\over p_c}t}m_G(B_e(c,A(c)r)^{1/p_c})).
\end{aligned}
$$
So, $m_G(B_e(c,(A(c)r)^{1/p_c})) \ge (A(c))^{-2t}L((A(c)r)^{1/p_c})^t$.
\endpf

\sp\blem\lab{lncp13.4.} 
Let ($m_G$, $m_H$) be a Euclidean semi $t$-conformal pair of measures
for an analytic map $f:G\lra H$. Let $c\in G$ be a critical point of
$f$ such that $m_G(c)=0$. If $0<r\le R(f,c)$ and $f(c)$ is 
$(s,L)$--$t$--u.e. with respect to $m_H$ for all $0<s\le r$, then $c$ is
$$
\(((A(c))^{-1}r)^{1/p_c},q(2(A(c))^2)^t(2^{t/p_c}-1)^{-1}L\)-{t-u.e.},
$$
with respect to $m_G$. 
\elem

\bpf Take any $s\le r$. Then
$f(B_e(c,(A(c))^{-1}s)^{1/p_c}))\sbt B_e(f(c),r)$. Therefore, recalling that
$A(c:a,b)=\{z:a\le |z-c|<b\}$, denoting
$$
A_n(c):A(c,2^{-1/p_c}((A(c))^{-1}s)^{1/p_c},((A(c))^{-1}s)^{1/p_c}),
$$ 
and using the decomposition of $B(c,((A(c))^{-1}s)^{1/p_c})$ described after
Definition~\ref{dcomp} we obtain,
$$
\begin{aligned}
Ls^t &\ge m_H(B_e(f(c),s))
     \ge m_H \(f(B_e(c,((A(c))^{-1}s)^{1/p_c}))\) \\
     &=p_c^{-1}\int_{B_e(c,((A(c))^{-1}s)^{1/p_c})}|f'(z)|^t\,dm_G(z) \\
&\ge p_c^{-1}\int_{A_n(c)}|f'(z)|^t\,dm_G(z) \\
&\ge p_c^{-1}A^{-t}(2^{-1}(A(c))^{-1}s)^{{p_c-1\over p_c}t}m_G(R(c)).
\end{aligned}
$$
So, 
$$
\begin{aligned}
m_G\(A(c;2^{-1/p_c}((A(c))^{-1}s)^{1/p_c},&((A(c))^{-1}s)^{1/p_c})\)\le \\
&\le p_c2^{t(1-{1\over p_c})}(A(c))^{2t}L(((A(c))^{-1}s)^{1/p_c})^t,
\end{aligned}
$$ 
and therefore
$$
\begin{aligned}
m_G \(B_e(c, &((A(c))^{-1}r)^{1/p_c})\)= \\
& =m_G \lt(\bu_{n=0}^\infty A\(c;2^{-{n+1\over p_c}}((A(c))^{-1}r)^{1/p_c},
2^{-{n\over p_c}}((A(c))^{-1}r)^{1/p_c})\rt) \\
&=\sum_{n=0}^\infty m_G \(A(c;2^{-{1\over
    p_c}}((A(c))^{-1}2^{-n}r)^{1/p_c},((A(c))^{-1}2^{- 
n}r)^{1/p_c})\)\\
&\le p_c(2^{1-{1\over p_c}}(A(c))^2)^tL\sum_{n=0}^\infty ((A(c))^{-1}2^{-
n}r)^{t/p_c}\\
&=p_c(2^{1-{1\over p_c}}(A(c))^2)^t{L\over 1-2^{-{t\over
      p_c}}}(((A(c))^{-1}r)^{1/p_c})^t \\ 
&=p_c(2(A(c))^2)^t(2^{t/p_c}-1)^{-1}L(((A(c))^{-1}r)^{1/p_c})^t.
\end{aligned}
$$
The proof is finished. \endpf

\sp

\blem\lab{lncp13.5.} 
Let ($m_G$, $m_H$) be a Euclidean $t$-conformal pair of measures
for an analytic map $f:G\lra H$. Let $c\in G$ be a critical point of an analytic
map $f:G\lra H$. If $0<r\le {1\over 3}R(H,c)$, $0<\sg\le
1$ and $f(c)$ is $(r,\sg,L)$--$t$--s.l.e with respect to $m_H$, then $c$ is
$$
(((A(c))^{-1}r)^{1/p_c},\^\sg,\^L)-t-{\rm s.l.e} 
$$
with respect to $m_G$, where
$\^\sg=(2^{p_c+1}K(A(c))^2\sg)^{1/p_c}$ and
$\^L=L\min\{K^{-t},((A(c))^2\sg)^{{1-p_c\over p_c}t}\}$,
\elem

\bpf Let $x\in B_e(c,((A(c))^{-1}r)^{1/p_c})$. If $\^\sg((A(c))^{-
1}r)^{1/p_c}\ge 2|x-c|$, then
$$\begin{aligned}
B_e(x,\^\sg((A(c))^{-1}r)^{1/p_c})
&\spt B_e(c,\^\sg((A(c))^{-1}r)^{1/p_c}/2)
=B_e(c,(2K)^{1/p_c}(A(c)\sg r)^{1/p_c}) \\
   & \spt B_e(c,(A(c)\sg r)^{1/p_c}).
\end{aligned}$$
It follows from our assumptions that $f(c)$ is $(\sg
r,\sg^{-t}L)$-l.e. with respect to $m_H$,
and therefore, in view of Lemma~\ref{lncp13.3.} the critical point
$c$ is $((A(c)\sg r)^{1/p_c},(A(c))^{-2t}\sg^{-t}L)$-l.e. with respect
to $m_G$. Thus 
\beq\lab{1042606}
\begin{aligned}
m_G\(B_e(x,\^\sg((A(c))^{-1}r)^{1/p_c})\)
&\ge (A(c))^{-2t}\sg^{-t}L(A(c)\sg r)^{t/p_c}\\
& =((A(c))^2\sg)^{{1- p_c\over p_c}t}L(((A(c))^{-1}r)^{1/p_c})^t.
\end{aligned}
\eeq
So, suppose that
\beq\lab{2042606}
\^\sg((A(c))^{-1}r)^{1/p_c}< 2|x-c|.
\eeq
Since $c$ is a critical point we have 
$$
|f'(x)|\ge
(A(c))^{-1}|x-c|^{p_c-1}
\ge (A(c))^{-1}\^\sg^{p_c- 1}((A(c))^{-1}r)^{p_c-1\over p_c}2^{1-p_c},
$$ which means that
\beq\lab{3042606}
\begin{aligned}
\^\sg((A(c))^{-1}r)^{1/p_c} 
& \ge (A(c))^{-1}\^\sg^p_c(A(c))^{-1}r2^{1-p_c}|f'(x)|^{-1}\\
& = 4K\sg r|f'(x)|^{-1} \\
&\ge K\sg r|f'(x)|^{-1}.
\end{aligned}
\eeq
In view of (\ref{2042606})
$$
|f(x)-f(c)|\ge (A(c))^{-1}|x-c|^{p_c}\ge
(A(c))^{-1}2^{-p_c}\^\sg^{p_c}(A(c))^{-1}r=2K\sg r\ge 2\sg r
$$
which implies that
\beq\lab{4042606}
f(c)\notin B_e(f(x),2\sg r).
\eeq
Since $|f(x)-f(c)|\le A(c)|x-c|^{p_c}\le R/3$, we have $B_e(f(x),2\sg
r)\sbt B_e(f(c),R)$. So, (\ref{4042606}) implies the existence of a
holomorphic inverse branch 
$$
f_x^{-1}:B_e(f(x),2\sg r)\lra G
$$ 
of $f$ which sends $f(x)$ to $x$. Since, by our assumptions $f(x)$ 
is $(\sg r,\sg^{-t}L)$--l.e. with respect to $m_H$, it follows from
Lemma~\ref{lncp13.2.} that $x$ is $(K\sg
r|f'(x)|^{-1},(K^2\sg)^{-t}L)$--l.e. with respect to $m_G$. Thus, using
(\ref{3042606}), we get
$$\begin{aligned}
m_G\(B_e(x,\^\sg((A(c))^{-1}r)^{1/p_c})\) 
\ge m_G \(B_e(x,K\sg r|f'(x)|^{-1}))
& \ge (K^2\sg)^{-t}L(K\sg r|f'(x)|^{-1})^t \\
&\ge K^{-t}Lr^t(A(c))^{-t}|x-c|^{(1-p_c)t}\\
& \ge K^{-t}L((A(c))^{-1}r)^t((A(c))^{-1}r)^{{1-p_c\over p_c}t} \\
&=K^{-t}L(((A(c))^{-1}r)^{1/p_c})^t.
\end{aligned}
$$
In view of this and (\ref{1042606}) the proof is completed.
\endpf

\chapter{Graph Directed Markov Systems}\label{Markov-systems}

\sp In this chapter we describe a powerful method to construct and study geometric and dynamical properties of fractal sets. This method is given by the theory of countable alphabet conformal iterated function systems, or more generally of countable alphabet conformal
graph directed Markov systems, as developed in the papers \cite{MU1}, \cite{MU5}, and  the book \cite{MU2}.
We present some elements of this theory now, primarily those related to conformal measures and a version of Bowen's Formula for the Hausdorff dimension of limit sets of such systems. In particular we will get an almost cost free, effective, lower estimate for the Hausdorff dimension of such limit sets. More about conformal graph directed Markov systems can be found 
in many papers and books; we bring up here some of them: \cite{MU3}--\cite{MU6}, \cite{MPU}, \cite{MSzU}, \cite{U3}, \cite{U4}, \cite{CTU}, \cite{HD_Spectrum}, \cite{CF Subsystems}, and \cite{CU-Porosity}. 
 
Afterwards, in Part 3 and especially in Part 4, we apply the techniques developed here to get a quite good, explicit estimate from below of the Hausdorff dimensions of the Julia sets of all elliptic functions and to explore stochastic properties of invariant versions of conformal measures for parabolic and subexpanding elliptic functions. 

\sp\section[Subshifts of Finite Type over Infinite Alphabets]{Subshifts of Finite Type over Infinite Alphabets; Topological Pressure}\label{SOFToIA;TP}

Let 
$$
\mathbb{N}:=\{1, 2, \ldots \}
$$  
be the set of all natural numbers, i.e. of all positive integers
and let $E$ be  a  countable  either finite or infinite, set, called
in the sequel  an alphabet. Let 
$$
\sg: E^\mathbb{N} \lra E^\mathbb{N}
$$
\index{(S)}{$\sg$} be the shift map\index{(N)}{shift
map}, i.e. the map cutting off the first coordinate. It is given by the
formula
$$  
\sg\( (\om_n)^\infty_{n=1}  \) :=  \( (\om_{n+1})^\infty_{n=1}  \).
$$
We also set
$$
E^*
:=\bigcup_{n=0}^\infty E^n.
$$\index{(S)}{$E^*$}
For every $\om \in E^*$, by $|\om|$  we mean the only integer
$n \geq 0$ such that $\om \in E^n$. We  call $|\om|$ the length of
$\om$. We make  a  convention that $E^0=\es$. If $ \om \in
E^\mathbb{N}$ and $n \geq 1$, we put
$$ \om_{|n}:=\om_1\ldots \om_n\in E^n.$$
Given $\om,\tau\in E^\infty$, we define $\omega\wedge\tau  \in
E^\infty\cup E^*$\index{(S)}{$\omega\wedge\tau$}  to be the longest
initial block common to both $\om$ and $\tau$. For each $\alpha >
0$, we define a $d_\alpha$
  metric $d_\alpha$\index{(S)}{$d_\alpha$}, on
$E^\infty$\index{(S)}{$E^\infty$}, by setting 
\beq\label{120191105}
d_\alpha(\om,\tau) := {\rm e}^{-\alpha|\om\wedge\tau|}.
\eeq
These metrics are all mutually equivalent, meaning that they
induce the same topology on $E^\infty$. Of course, they then define the same $\sg$--algebras of Borel sets. 

\sp A function defined on $E^\infty$ is uniformly 
continuous with respect to one of these metrics if and only if it is
uniformly continuous with respect to all of them. Also, a function is
H\"older continuous with respect to one of these metrics if and only if it is H\"older with respect to all of them; of course the H\"older exponent depends on the metric. If no metric is specifically mentioned, we take it to
be $d_1$.

\sp Now consider a 0--1 matrix 
$$
A:E \times E \lra \{0,1\}.
$$
Any such matrix is called an incidence matrix. \index{(N)}{incidence matrix}Set
$$
E^\infty_A:=\big\{\om \in E^\N: A_{\om_i\om_{i+1}}=1
\,\, \mbox{for all}\,\,   i \geq 1 \big\}.
$$\index{(S)}{$E^\mathbb{N}_A$}  
All  elements   of
$E^\infty_A$ are called $A$--admissible\index{(N)}{admissible}.
We also set
$$
E^n_A:=\big\{w \in E^\mathbb{N}:  \,\, A_{\om_i\om_{i+1}}=1  \,\, \mbox{for
  all}\,\,  1\leq i \leq  n-1\big\}, \,\, \ n \geq 1,
$$\index{(S)}{$E^n_A$} 
and
$$
E^*_A:=\bigcup_{n=0}^\infty E^n_A.
$$\index{(S)}{$E^*_A$}
The elements of these sets are also  called $A$-admissible. For
every  $\om \in E^*_A$, we put  
$$
[\om]:=\big\{\tau \in E^\infty_A:
\,\, \tau_{|_{|\om|}}=\om \big\},
$$
\index{(S)}{$[\om]$} and call this set the cylinder generated by $\om$. The  following fact is obvious.

\sp\bprop\label{p1j83} If $E$ is a countable set and $A:E \times E \lra \{0,1\}$ is an incidence matrix, then  $E^\mathbb{N}_A$ is a closed subset of $E^\infty$ invariant under the shift map $\sg: E^\infty\lra E^\infty$. The later means  that 
$$
\sg (E^\infty_A)\sbt E^\infty_A.
$$
\eprop

\sp The matrix $A:E \times E \lra \{0,1\}$ is said to be
\index{(N)}{finitely irreducible}{\it finitely irreducible} if there
exists a finite set $\Lambda\sbt E_A^*$ such that for all $i,j\in E$
there exists a path $\om\in \Lambda$ for which 
$$
i\om j\in E_A^*.
$$
Given a set $F\sbt E$ we  put
$$
F^\infty:=\{\om \in E^\N: \om_i\in F
\,\,\, \mbox{for all }\, \, \,  i \geq
1\}.
$$\index{(S)}{$F^\infty$}

\fr  Given $ F \sbt E$ and a function $f:F_A^\infty\lra
\mathbb{R}$, we define the standard {\it
n--th partition function}\index{(N)}{$n$-th partition function} by
$$
Z_n(F,f)\index{(S)}{$Z_n(F,f)$}
:=\sum_{\om \in F_A^n}\exp
\lt(\sup_{\tau\in [\om]}\sum_{j=0}^{n-1}f(\sg^j(\tau))\rt).
$$
Recall from Definition~\ref{1d20190927} that a sequence $\{a_n\}_{n=1}^\infty$ consisting of real numbers 
is said to be subadditive\index{(N)}{subadditive sequence} if 
$$
a_{n+m}\le a_n+a_n
$$
for all $m,n\ge 1$. 
We will need the following.

\sp\blem\lab{subadd2} 
The sequence $\cN\ni n\longmapsto \log Z_n(F,f)$ is
subadditive. 
\elem

\bpf We need to show that the sequence $\cN\ni n\longmapsto
Z_n(F,f)$ is submultiplicative\index{(N)}{submultiplicative
sequence}, i.e. that 
$$
Z_{m+n}(F,f)\le Z_m(F,f)Z_n(F,f)
$$ 
for all $m,n\ge 1$. And indeed,
$$
\aligned Z_{m+n}(F,f) 
&=\sum_{\om \in F_A^{m+n}}\exp\lt(\sup_{\tau\in
  [\om]}\sum_{j=0}^{mn-1}f(\sg^j(\tau))\rt) \\
&=\sum_{\om \in F_A^{m+n}}\exp\lt(\sup_{\tau\in
  [\om]}\lt\{\sum_{j=0}^{m-1}f(\sg^j(\tau))+
  \sum_{j=0}^{n-1}f(\sg^j(\sg^m(\tau)))\rt\}\rt) \\
&\le \sum_{\om\in F_A^{m+n}}\exp\lt(\sup_{\tau\in
  [\om]}\sum_{j=0}^{m-1}f(\sg^j(\tau))+
 \sup_{\tau\in [\om]_F} \sum_{j=0}^{n-1}f(\sg^j(\sg^m(\tau))\rt) \\
&\le \sum_{\om \in F_A^m}\sum_{\rho \in F_A^n}
\exp\lt(\sup_{\tau\in [\om]}\sum_{j=0}^{m-1}f(\sg^j(\tau))+
\sup_{\g\in [\rho]}\sum_{j=0}^{n-1}f(\sg^j(\g))\rt) \\
&=\sum_{\om \in F_A^m}\exp\lt(\sup_{\tau\in [\om]}
 \sum_{j=0}^{m-1}f(\sg^j(\tau))\rt)\cdot
 \sum_{\rho \in F_A^n}\exp\lt(\sup_{\g\in [\rho]}
 \sum_{j=0}^{m-1}f(\sg^j(\g))\rt) \\
&=Z_m(F,f)Z_n(F,f).
\endaligned
$$\endpf

\sp We can now define the \index{(N)}{topological pressure}
topological pressure of $f$ with respect to the shift map
$\sg:F_A^\mathbb{N}\lra F_A^\mathbb{N}$. Indeed, combining Lemma~\ref{subadd2}and Lemma~\ref{ch6subadditive}, we see that the folowing limit exists and the equality following it holds.
\beq\lab{2.1.1}
\P_F(f):=\lim_{n\to\infty}{1\over n}\log Z_n(F,f)
       = \inf_{n\in\cN}\lt\{{1\over n}\log Z_n(F,f)\rt\}.
\eeq
If $F=E,$ we suppress the subscript $F$ and write simply $\P(f)$ for
$\P_E(f)$ and $Z_n(f)$ for $Z_n(E,f)$.

\sp\bdfn\lab{d2.1.2}  A function $f:E_A^\infty\lra \mathbb{R}$ is said to be
\index{(N)}{acceptable function} {\it acceptable} provided it is
uniformly continuous and \index{(S)}{${\rm osc}(f)$}
$$
{\rm osc}(f):=\sup\big\{\sup(f|_{[e]})-\inf(f|_{[e]}):e \in E\big\}<+\infty.
$$
\edfn

\sp Note that an acceptable function need not be bounded. In fact, in what follows, it usually will not. We shall prove the following.

\sp\bthm\lab{t2.1.3} 
If $f:E^\infty_A\lra \mathbb{R}$ is acceptable and
$A$ is finitely irreducible, then
$$
\P(f)=\sup\{\P_F(f)\},
$$
where the supremum is taken over all finite subsets $F$ of $E$.
\ethm

\bpf The inequality 
$$
\P(f) \geq \sup\{\P_F(f)\}
$$ 
is obvious. In order to prove the converse let $\Lambda\sbt E_A^*$ witness finite irreducibly of the matrix $A$. We shall show first that 
$$
\P(f)<+\infty.
$$
Put 
$$
q:=\#\La, \ p:=\max\{|\om|:\om\in \Lambda\}, \  \text{ and } \
T:=\min\left\{\inf\lt\{\sum_{j=0}^{|\om|-1}f\circ \sg^j|_{[\om]}\rt\}:\om\in
\Lambda\right\}.
$$
Fix $\e>0$. By acceptability of the function $f^\infty_A\lra \mathbb{R}$, we have
$M:= \osc(f)<+\infty$ and there exists $l\ge 1$ such that 
$$
|f(\om)-f(\tau)|<\e
$$ 
whenever $\om|_l=\tau|_l$. Now, fix $k \geq l.$ By Lemma~\ref{subadd2},
${1\over k}\log Z_k(f) \ge \P(f)$. Therefore, there exists a finite
set $F\sbt I$ such that
\beq\lab{2.1.2}
{1\over k}\log Z_k(F,f) > \P(f)-\e.
\eeq
We may assume that $F$ contains $\Lambda$. Put
$$
\ov f:=\sum_{j=0}^{k-1}f\circ\sg^j.\index{(S)}{$\ov f$}
$$
Now, for every element $\tau=\tau_1, \tau_2,\ld,\tau_n\in F_A^k
\times \cdots \times F_A^k$ ($n$ factors) one can choose
elements $\a_1,\a_2,\ld,\a_{n-1}\in \Lambda$ such that
$$
\ov{\tau}:=\tau_1\a_1\tau_2\a_2\ld\tau_{n-1}\a_{n-1}\tau_n\in E^*.
$$
Notice that the defined in this way function $\tau\longmapsto \ov{\tau}$ is at most $u^{n-1}$--to--$1$, where $u$ is the number
of lengths of words composing $\La$). Then for every $n\ge 1$,
$$
\aligned q^{n-1}\sum_{i=kn}^{kn+p(n-1)}Z_i(F,f) 
&\ge \sum_{\tau\in (F_A^k)^n}
\exp\lt(\sup_{[{\ov\tau}]}\sum_{j=0}^{|\ov\tau|}f\circ\sg^j\rt) 
\ge \sum_{\tau\in (F_A^k)^n}
   \exp\lt(\inf_{[\ov\tau]}\sum_{j=0}^{|\ov\tau|}f\circ\sg^j\rt) 
     \\
&\ge \sum_{\tau\in (F_A^k)^n}\exp\lt(\sum_{i=1}^n\inf\ov
f|_{[\tau_i]} +T(n-1)\rt) \\
&=\exp(T(n-1))\sum_{\tau\in (F_A^k)^n}\exp\sum_{i=1}^n\inf\ov
f|_{[\tau_i]}\\
&\ge \exp(T(n-1))\sum_{\tau\in (F_A^k)^n}\exp\lt(
\sum_{i=1}^n(\sup\ov
f|_{[\tau_i]} -(k-l)\e-Ml)\rt) \\
&=\exp(T(n-1)-(k-l)\e n - Mln)\sum_{\tau\in (F_A^k)^n}\exp
\sum_{i=1}^n\sup\ov f|_{[\tau_i]}\\
&=\ep^{-T}\exp\(n(T-(k-l)\e-Ml)\)\lt(\sum_{\tau\in F_A^k
}\exp(\sup\ov f|_{[\tau]})\rt)^n.
\endaligned
$$
Hence, there exists $kn\le i_n\le (k+p)n$ such that
$$
Z_{i_n}(F,f)\ge {1\over pn}\ep^{-T}\exp\(n(T-(k-l)\e-Ml-\log
q)\)Z_k(F,f)^n
$$
and therefore, using (\ref{2.1.2}), we obtain
$$
\P_F(f) =\lim_{n\to\infty}{1\over i_n}\log Z_{i_n}(F,f) \ge
{-|T|\over k}-\e +{l\e\over k+p}-{Ml+\log p\over k}+\P(f)-2\e \ge
\P(f)-7\e
$$
provided that $k$ is large enough. Thus, letting $\e\downto 0$, the
theorem follows. The case $\P(f)=\infty$ can be treated similarly.
\endpf

\

\section{Graph Directed Markov Systems} 

\fr We start with the following definition of the main object of this chapter.

\bdfn\label{GDMS_def}
A Graph Directed Markov System (GDMS)\index{(N)}{graph directed
Markov systems}\index{(S)}{\rm{GDMS}} consists of 

\sp\begin{itemize}
\item[(a)] a directed multigraph $(E,V)$ with a countable set of edges $E$
  and a finite set of vertices $V$,

\sp\item[(b)] an incidence matrix $A:E\times E\lra\{0,1\}$,

\sp\sp\item[(c)] two functions $i,t:E\lra V$ such that $t(a)=i(b)$ whenever $A_{ab}=1$.

\sp\item[(d)] a family of non--empty compact metric spaces $\{X_v\}_{v\in V}$, 

\sp\item[(e)] a number $\b\in (0,1)$, and

\sp\item[(f)] a collection $\big\{\phi_e:X_{t(e)}\lra X_{i(e)}:e\in E\big\}$ of 1--to--1 contractions, all with a Lipschitz constant $\le \b$. 
\end{itemize}

\sp\fr Briefly, the set
$$
S:=\{\phi_e:X_{t(e)}\lra X_{i(e)}\}_{e\in E}
$$
is called a Graph Directed Markov System (GDMS). The corresponding
number $\b\in (0,1)$ is denoted by $\b_S$. 

A graph  directed  Markov system is called an iterated function system (IFS) if $V$\index{(N)}{iterated function system}, the set of vertices, is a singleton and $A(E\times E)=\{1\}$.
\edfn  

\sp The main object of interest in this chapter will be 
the limit set of the system $S$ and its geometric features.
We now define the limit set. For each $\om \in E^*_A$, say $\om\in
E^n_A$, we consider the map coded by $\om$:
$$
\phi_\om:=\phi_{\om_1}\circ\cdots\circ\phi_{\om_n}:X_{t(\om_n)}\lra
X_{i(\om_1)}.\index{(S)}{$\phi_\om$}
$$
For $\om \in E^\infty_A$, the sets
$\{\phi_{\om|_n}\(X_{t(\om_n)}\)\}_{n\ge 1}$ form a descending
sequence of non--empty compact sets and therefore $\bi_{n\ge
1}\phi_{\om|_n}\(X_{t(\om_n)}\)\ne\es$. Since for every $n\ge 1$,
$$
\diam\(\phi_{\om|_n}\(X_{t(\om_n)}\)\)\le
\b_S^n\diam\(X_{t(\om_n)}\)\le \b_S^n\max\{\diam(X_v):v\in V\},
$$ 
we conclude that the intersection
$$
\bi_{n\ge1}\phi_{\om|_n}\(X_{t(\om_n)}\)
$$
is a singleton and we denote its only element by $\pi(\om)$. In this
way we have defined the ``projection'' map $\pi$ 
$$
\pi:E^\infty_A\lra \du_{v\in V}X_v
$$
from the coding space $E_A^\infty$ to $\du_{v\in V}X_v$, the disjoint union of the
compact sets $X_v$. The set 
$$ 
J=J_S:=\pi(E^\infty_A)\index{(S)}{$J_S$}
$$
will be called the \index{(N)}{limit set}{\it limit set} of the GDMS
$S$. 

\brem\label{r120190212}
As a matter of fact we do not need all the maps $\phi_e:X_{t(e)}\lra X_{i(e)}$, $e\in E$, to be uniform contractions. it entirely suffices to know that there exist $\b\in(0,1)$ and some integer $q\ge 1$ such that all the maps $\phi_\om:X_{t(\om)}\lra X_{i(\om)}$, $\om\in E_A^q$ are contractions with a Lipschitz constant $\le \b$. 
\erem

\

\section{Conformal Graph Directed Markov Systems}\label{CGDMS}

\sp In order to be able anything meaningful about fractal properties of the limit set of a graph directed Markov system, the contracting maps $\phi_e$, $e\in E$, need to have some additional smooth properties, than merely continuity. There are basically two options: the contracting maps of the system $S$ are commonly assumed to be either affine maps or conformal maps. In this book we consider the latter case.

\bdfn\label{d120190419}
We call a GDMS\index{(S)}{\rm{GDMS}}  conformal graph directed Markov system\index{(N)}{conformal graph directed Markov system} (CGDMS)
\index{(S)}{\rm{CGDMS}} if the following conditions are satisfied.

\sp
\begin{itemize}
\item[(4a)] For every vertex $v\in V$, $X_v$ is a compact 
subset of a Euclidean space $\mathbb{R}^d$ (the dimension $d$ common
for all $v\in V$) and $X_v=\ov{\Int(X_v)}$.

\sp\item[(4b)] \index{(N)}{Open Set Condition}{\rm({\it Open Set
Condition})(OSC)}.\index{(S)}{\rm{OSC}} For all 
$a,b\in E$ with $a\ne b$, we have 
$$
\phi_a(\Int(X_{t(a)})\cap \phi_b(\Int(X_{t(b)})=\es.
$$
\item[(4c)]  For every vertex $v\in V$ there exists an open connected set
$W_v\spt X_v$ such that for every $\om\in E_A^*$ the map
$\f_\om$ extends to a $C^1$ conformal diffeomorphism from $W_{t(\om)}$ into
$\R^d$.

\sp\item[(4d)] There are two constants $L\ge 1$ and $\a>0$ such that
$$
\bigg|\frac{|\f_e'(y)|}{|\f_e'(x)|}-1\bigg|\le
L\|y-x\|^\a
$$
for every $e\in E$ and every pair of points $x,y\in W_{t(e)}$, where
$|\f_\om'(x)|$ means the norm of the derivative; equivalently this is
the scaling factor of the similarity map $\f_\om'(x):\cR^d\lra\cR^d$.
\end{itemize}
\edfn

\brem\label{r1_2018_02_14}
We would like to emphasize that unlike \cite{MU1} and \cite{MU2}, until Section~\ref{SOSC}, entitled The Strong Open Set Condition, we do NOT assume in this chapter any kind of cone condition or strong open set condition. 
\erem

\brem\label{r120190212B}
The requirement of uniform contractions of the maps $\phi_e$, $e\in E$, or rather, its weaker form explained in Remark~\ref{r120190212}, is now replaced by a slightly stronger condition requiring that
\beq\label{820190212}
\|\phi_\om'\|_\infty:=\sup\big\{|\phi_\om'(x)|:x\in W_{t(\om)}\big\}
\le\b=\b_S< 1
\eeq
for some integer $q\ge 1$ and all $\om\in E_A^q$. Passing to the (iterated) GDMS 
$$
\big\{\phi_\om:X_t(\om)\lra X_{i(\om)}:\om\in E_A^q\big\},
$$
we will frequently assume, without loss of generality that \eqref{820190212} holds with $q=1$.
\erem 

\brem\label{r220190212}
The condition (4c) follows from the following (stronger) condition
\begin{itemize}
\item[(4c')] For every vertex $v\in V$ there exists an open connected set
$W_v\spt X_v$ such that for every $e\in E$ with $t(e)=v$, the map
$\f_\om$ extends to a $C^1$ conformal diffeomorphism from $W_{t(e)}$ into
$W_{i(e)}$.
\end{itemize}
\erem

\sp\fr We start our investigations of conformal GDMSs by proving the following.

\blem\label{l1201200102}
If $S=\{\phi_e\}_{e\in E}$ is a CGDMS, then for every vertex $v\in V$, 
there exists an open, connected and bounded set $W_v'$ such that 
$X_v\sbt W_v'\sbt\overline{W_v'}\sbt W_v$. 
\elem
\bpf Fix $v\in V$. 
Let 
$$
\delta_v:=\mathrm{dist}\(X_v,\R^d\backslash W_v\)>0.
$$
The neighborhood
$$
B(X_v,\delta_v/2)=\bigcup_{x\in X_v}B(x,\delta_v/2)
$$ 
of $X_v$ is open and bounded but not necessarily connected. However, each open ball 
$$
B(x,\delta_v/2), \  \ x\in X_v,
$$
connected and thus every connected component of $B(X_v,\delta_v/2)$ contains at least
one such ball. Since each set $B(X_v,\delta_v/2)$, $v\in V$. is bounded and its connected components 
are mutually disjoint and each contain at least one ball of radius $\delta_v/2$, 
the neighborhood $B(X_v,\delta_v/2)$ has only finitely many connected components. 
Denote these components by 
$$
B_v^{(1)},B_v^{(2)},\ldots,B_v^{(k_v)}.
$$
For each $1\leq j\leq k_v$, choose
$$
x_v^{(j)}\in B_v^{(j)}\cap X_v.
$$
Since the set $W_v$ is path--connected, for every $1\leq j<k_v$ there is 
a piecewise $C^1$ curve $\gamma_v^{(j)}$ joining $x_v^{(j)}$ and $x_v^{(j+1)}$ in $W_v$.
Since $\gamma_v^{(j)}$ is compact and $W_v$ is open, there is $\delta_v^{(j)}>0$ 
such that 
$$
B(\gamma_v^{(j)},\delta_v^{(j)})\sbt W_v.
$$
So, setting
$$
W_v':=B(X_v,\delta_v/2)\cup\bigcup_{j=1}^{k-1}B(\gamma_v^{(j)},\delta_v^{(j)}/2)
$$ 
finishes the proof.
\epf

\sp\bprop\lab{p1.033101} If $d\ge 2$ and a family $S=\{\phi_e\}_{e\in
I}$ satisfies conditions (4a) and (4c), then it also satisfies
condition (4d) with $\a=1$ after replacing all the sets $W_v$, $v\in V$, by the corresponding sets $W_v'$, $v\in V$, produced in 
Lemma~\ref{l1201200102}. 
\eprop

\bpf 
Since the set of vertices $V$ is finite, this proposition is an immediate consequence of Theorem~\ref{tKoebe-1} in the case
$d=2$. In the case $d \geq 3$, which will not  be considered in Parts 3 and 4 of 
this book, the proposition is proved in \cite{{MU2}} and ultimately
relays on Liouville's Classification Theorem of conformal
homeomorphisms  of $\ov{\mathbb R^d}$, which asserts that every conformal map in $\R^d$, $d\ge 3$, is a composition of the inversion with respect to a sphere with radius $1$ (the center can be $\infty$) and a similarity map.
\epf

\sp

As a rather straightforward consequence of hypothesis (4d) we get the
following.

\blem\lab{l2.033101} If $S=\{\phi_e\}_{e\in E}$ is a CGDMS, then for
all $\om\in E^*$ and all $x,y\in W_{t(\om)}$, we have
$$
\bigl|\log|\phi_\om'(y)|-\log|\phi_\om'(x)|\bigr|\le {L\over
1-\b_S^\a}\|y-x\|^\a.
$$
\elem

\bpf For every $\om\in E^*$, say $\om\in E^n$, and every
$z \in W_{t(\om)}$ put
$z_k=\f_{\om_{n-k+1}}\circ\f_{\om_{n-k+2}}\circ\cdots
\circ\f_{\om_n}(z)$; put also $z_0=z$. In view of (4d) for any two
points $x,y\in  W_{t(\om)}$ we have
\beq\lab{Log1}
\aligned \bigl|\log(|\f_\om'(y)|)-\log(|\f_\om'(x)|)\bigr|
&=\left|\sum_{j=1}^n\log\left(1+{|\f_{\om_j}'(y_{n-j})|-|\f_{\om_j}'(x_{n-
j})|\over |\f_{\om_j}'(x_{n-j})|}\right)\right| 
\\
&\le \sum_{j=1}^n\frac{\left|\f_{\om_j}'(y_{n-j})|-
|\f_{\om_j}'(x_{n- j})|\right|}{|\f_{\om_j}'(x_{n-j})|}
\\
&\le \sum_{j=1}^nL\|y_{n-j}-x_{n-j}\|^\a 
\le L\sum_{j=1}^n\b_S^{\a (n-j)}\|y-x\|^\a 
\\
&\le {L\over 1-\b_S^\a}\|y-x\|^\a.
\endaligned
\eeq 
\endpf

\sp As a straightforward consequence of this lemma (and ultimately of (4d)), we get the following.

\sp

\begin{itemize}
\item[(4e)](Bounded Distortion Property)\index{(N)}{bounded distortion
    property}. There exists $K\geq 
1$ such that for all $\om\in E^*$ and all $x,y\in W_{t(\om)}$
$$
|\f_\om'(y)|\leq K|\f_\om'(x)|.
$$
\end{itemize}

\sp We shall now prove some basic geometric consequences of the
properties (4a)--(4e). Because of the Mean Value Inequality, for all finite words $\om\in E^*$, all convex subsets $C$ of $W_{t(\om)}$, all $x\in
X_{t(\om)}$ and all radii $0<r\le \dist(X_{t(\om)}, \bd W_{t(\om)})$,
we have
\beq\lab{4.1.7}
\diam(\f_\om(C))\le \|\f_\om'\|_\infty\diam(C), \, \,\,
\phi_\om(B(x,r))\sbt B(\phi_\om(x),\|\phi_\om'\|r).
\eeq
We now shall prove the following.

\sp\blem\label{l22013_03_12}
Replacing, if necessary, all the sets $W_v$, $v\in V$, by the corresponding sets $W_v'$, $v\in V$, produced in Lemma~\ref{l1201200102}, we will have that there exists a constant $D\ge 1$ such that
\beq\lab{4.1.9a}
\|\f_\om(y)-\f_\om(x)\|\le D\|\phi_\om'\|\cdot\|y-x\|
\eeq
for all finite words $\om\in E_A^*$ and all $x,y\in W_{t(\om)}$. 
\elem

\bpf Fix
$$
R:=\min\{\dist(W_v',\bd W_v):v\in V\}>0.
$$
Since $V$, the set of all vertices, is finite and since
each set $\ov{W_v'}$ is compact and connected, there exists an integer $D\ge
3$ such that each such set can be covered 
by finitely many balls $B(x_1,R/2),\ld, B(x_{D-1},R/2)$ with centers
$x_1,\ld,x_{D-1}$ in $W_v'$ and with the property that
$$
B(x_i,R/2)\cap B(x_{i+1},R/2)\ne\es
$$
for all $i=1,2,\ld,D-2$. Therefore, for every vertex $v\in V$ and all
points $x, y\in \ov{W_v'}$ there are $k\le D+1$ points $x=z_1,z_2,\ld,z_k=y$
in $B(W_v',R)$ such that $\|z_i-z_{i+1}\|\le \|x-y\|$ for all
$i=1,2,\ld, k-1$. Hence, using the Mean Value Inequality, we get
$$
\aligned
\|\phi_\om(x)-\phi_\om(y))\|
&\le \sum_{i=1}^{k-1}\|\phi_\om(z_i)-\phi_\om(z_{i+1})\|
\le \sum_{i=1}^{k-1}\|\phi_\om'\|_\infty\|z_i-z_{i+1}\| \\
&\le (k-1)\|\phi_\om'\|_\infty\|x-y\| \\
&\le D\|\phi_\om'\|_\infty\|x-y\|. 
\endaligned
$$
Thus, the proof is complete.
\epf

\sp Since the set $V$ is finite, replacing the number $D$ of Lemma~\ref{l22013_03_12} by 
$$
D\max\big\{1,\max\{\diam(W_v):v\in V\}\big\},
$$
as an immediate consequence of this lemma we get the following.

\sp\bcor\label{c23013_03_12}
For all finite words $\om\in E_A^*$ we have that
\beq\lab{4.1.9}
\diam\(\f_\om(W_{t(\om)})\)\le D \|\f_\om'\|_\infty.
\eeq
\ecor

\sp  Now we shall prove the following.

\sp

\blem\label{l42013_03_12}
For all finite words $\om\in E_A^*$, all $x\in X_{t(\om)}$ and all radii
$0<r<\dist(X_{t(\om)}, \bd W_{t(\om)})$, we have that
\beq\lab{4.1.8}
\f_\om(B(x,r))\spt B(\f_\om(x),K^{-1}\|\f_\om'\|r),
\eeq
\elem

\bpf 
First notice that $B(x,r)\sbt W_{t(\om)}$. Take also any $\om\in E^*$
and let $R\ge 0$ be the maximal radius such that
\beq\lab{1.051101}
B(\f_\om(x),R)\sbt \f_\om(B(x,r)).
\eeq
Then 
\beq\lab{120200103}
\bd\(B(\f_\om(x),R)\)\cap \bd\(\f_\om(B(x,r))\)\ne\es,
\eeq
and in view of the Mean Value Inequality along with the Bounded Distortion
Property (4e), we get that
$$
\f_\om^{-1}\(B(\f_\om(x),R)\) 
\sbt B(x,K\|\f_\om'\|_\infty^{-1}R).
$$
If $K\|\f_\om'\|_\infty^{-1}R<r$, then the set $\f_\om\(B(x,K\|\f_\om'\|_\infty^{-1}R)\)$ is well defined and 
$$
B(\f_\om(x),R)\sbt
\f_\om\(B(x,K\|\f_\om'\|_\infty^{-1}R)\).
$$
But then, both \eqref{1.051101} and \eqref{120200103} imply that $K\|\f_\om'\|_\infty^{-1}R\ge r$. This contradiction shows that  
$$
K\|\f_\om'\|_\infty^{-1}R\ge r.
$$ 
So, using (\ref{1.051101}) we obtain (\ref{4.1.8}), which completes
the proof. 
\epf

\sp Let $R_S$ be the radius of the largest open ball that can be
inscribed in all the sets $X_v$, $v \in V$. Let $x_v\in \Int(X_v)$, $v \in V$, be the
centers of respective balls. As an immediate
consequence of Lemma~\ref{l42013_03_12} we get, perhaps with a larger
constant $D$, the following.  
 
\sp

\blem\label{l52013_03_12}
For all finite words $\om\in E_A^*$ we have that
\beq\lab{4.1.10a}
\f_\om(\Int\(X_{t(\om)}\)\spt
B\(\phi_\om(x_{t(\om)}),K^{-1}\|\f_\om'\|R_S),
\eeq
and so, 
\beq\lab{4.1.10}
\diam(\f_\om(X_{t(\om)}))\ge D^{-1} \|\f_\om'\|_\infty.
\eeq
\elem

\

\section{Topological Pressure, $ \theta$--Number, and Bowen's
Parameter}\label{sec-top-press}

Let $S$  be a finitely irreducible conformal 
GDMS. For every $t\geq 0$ let
$$
Z_n(t):=\sum_{\om\in E^n_A} \|\phi_\om'\|_\infty^t.
$$
Since 
$$
\|\phi_{\om \tau}'\|_\infty \leq \|\phi_\om'\|_\infty\cdot\|\phi_{ \tau}'\|_\infty,
$$ 
we see that
$Z_{m+n}(t)\leq Z_m(t)Z_n(t)$ for all integers $m, n\ge 1$, and consequently, the sequence 
$$  
\N\ni n \longmapsto  \log Z_n(t)\in\R
$$
is subadditive. Thus,  the limit
$$ \lim_{n \to  \infty}  \frac{1}{n} \log Z_n(t)$$
exists and is equal to  
$$
\inf_{n \geq 1}\lt\{\frac{1}{n} \log Z_n(t)\rt\}.
$$
This limit is denoted by $\P(t)$, or, if we want to be more precise,
by $\P_E(t)$ or $\P_S(t)$. It
is called the topological pressure of the system $S$ \index{(N)}{topological pressure}
evaluated at the parameter $t$. Let $\zeta: E^\mathbb{N}_A \lra
\mathbb{R}$ be the function defined by the formula
$$ 
\zeta(\om):= t \log |\phi_{\om_1}'(\pi(\sg(\om))|.
$$
\index{(S)}{$\zeta(\om)$} As a straightforward consequence of  Lemma~\ref{l2.033101} we get the following.

\sp

\blem\label{l1j85} 
For every $t \geq 0$ the function $t \zeta:E^\mathbb{N}_A \lra
\mathbb{R}$ is H\"older continuous and $\P(\sg,t\zeta)=\P(t)$.
\elem

\sp\fr Let 
$$
F(S):=\{  t \geq 0:   \,\, \P(t)< \infty\}
$$ 
and let
$$ 
\theta(S):=\inf(F(S)).
$$
\index{(S)}{$\theta(S)$} Having bounded distortion (4e) and uniform
contraction of all generators of the systems $S$, the following facts
are easy to prove. 

\sp\bprop\label{p2j85} If $S$ is a finitely irreducible conformal
GDMS, then 
\begin{itemize}
\item [(a)] $\P(t)< +\infty \ \Leftrightarrow \  Z_1(t) < + \infty$.

\sp\item [(b)]  $\inf \{t\ge 0: Z_1(t)< + \infty\}=\theta(S)$.

\sp\item [(c)] The  topological pressure function $\P(t)$ is
\begin{itemize}

\,

\item[(c1)] non--increasing on $[0, + \infty)$, 

\,

\item[(c2)] strictly decreasing on $[0,+\infty)$ with $\lim_{t\to +\infty}\P(t)=-\infty$, 

\,

\item[(c3)] convex and continuous on $F(S)$
\end{itemize}

\sp\item [(d)] $\P(0)=+\infty$ if and only if $E$  is infinite.
\end{itemize}
\eprop

\sp\fr The number 
\beq\lab{120200108}
h_S:=h:= \inf\{t\geq 0:  \P(t)\leq 0\}
\eeq
is called the Bowen's
parameter of the  system $S$ \index{(N)}{Bowen's parameter}. It will turn out to be equal to the Hausdorff dimension of the limit set $J_S$. In view  of Proposition~\ref{p2j85}(c1) we have 
\beq\lab{220200108}
\P(h_S)\leq 0.
\eeq
The following useful observation is now obvious.

\sp\bobs\label{o120131007}
If $\P(t)=0$ for some $t\ge 0$, then $t=h$.
\eobs

\bdfn\label{d5_2017_11_18}
We say  that the system $S$ is 

\sp\begin{itemize}
\item regular if $\P(h_S)=0$, 

\sp \item strongly regular if 
there exists $t \geq 0$ such that $0< \P(t) <+\infty$, and 

\sp\item  co--finitely regular if $\P(\theta_S)=+\infty$. 
\end{itemize}
\edfn

\sp\fr It is easy to see that each co--finitely regular system
is strongly regular and each strongly regular one is regular. 
We shall prove the following.

\sp\bprop\label{cof-reg}
If a finitely irreducible conformal GDMS $S$ is  co--finitely regular, then
for every co--finite finitely irreducible subset $F\sbt E$, the system $S_F$ is also co--finitely regular, thus regular. 

In addition if $S$ is a conformal iterated function system (IFS), and its every co--finite subsystem is regular, then $S$ is  co--finitely regular.
\eprop

\bpf
First notice that $\th_F=\th_S$ for every co--finite subset $F$ of
$E$. Suppose now that $S$ is co--finitely regular. By virtue of
Proposition~\ref{p2j85}(a) this means that $Z_1(\theta)=+\infty$. But
then $Z_1(F,\th)=+\infty$ for every co--finite subset $F$ of
$E$. This however, by Proposition~\ref{p2j85}(a) again, means that
each such finitely irreducible system is co--finitely regular, thus regular. 

\sp For the converse suppose that $S$ is an  iterated function system and the system $S$ is not
co--finitely regular. By virtue of Proposition~\ref{p2j85}(a) used for the
third time, this means that $Z_1(\th)<+\infty$. But then there exists a
co-finite subset $F$ of $E$ such that $Z_1(F,\th)<1$. Hence, by the
definition of topological pressure, $\P_F(\th)<0$. But then $F$ is not
regular and we are done. 
\epf

\sp Let us record the following obvious fact.

\sp

\bdf\label{f1j87} If $S$ is a finitely irreducible conformal {\rm GDMS}, then $\theta(S)\leq h_S$. If $S$ is  strongly regular, in particular, if $S$ is co--finitely  regular, then $\theta(S)< h_S$. 
\edf

\section{Hausdorff Dimension and Bowen's Formula for
GDMS}\label{sec-HD-Bowen} 

In this section  we  prove a formula  for the Hausdorff dimension of
the  limit set of a finitely irreducible  conformal GDMS. It is
entirely expressed in dynamical terms and, because of its
correspondence to the formula obtained in the breakthrough Bowen's work \cite{Bow2} dealing with quasi--Fuchsian groups, we refer to it as Bowen's formula. In the case of finite alphabets $E$, $S$ being an IFS, and the maps $\phi_e$, $e\in E$, all being similarities, it was obtained in the seminal paper \cite{Hutchinson}, where Hutchinson also introduced the concept of the Open Set Condition. In the case when the alphabet $E$ is still finite, $S$ is an IFS, but the maps $\phi_e$, $e\in E$, are all assumed only to be conformal, Bowen's formula was obtained in \cite{Bedford}. In the case of arbitrary (countable alphabet) conformal IFSs this formula was obtained in \cite{MU1}. Finally it was proved in its final form of arbitrary (countable alphabet) conformal GMMSs in\cite{MU2}. This is its form we formulate and prove in the current section. It got also some extensions to conformal IFSs
acting in Hilbert spaces, see \cite{MSzU}.  

We start  with the  following simple geometrical observation following from the Open Set Condition.

\sp\blem\label{l1j81} If $S$ is a conformal GDMS, then for
all $0< \kappa_1 < \kappa_2< +  \infty$, for all $r >0$ and for all
$x \in X$, the cardinality  of any  collection of mutually 
incomparable words $\om \in E^*_A$ that satisfy the  condition
$$ 
B(x,r) \cap \phi_\om(X_{t( \om)})\neq \es
$$
and
$$ 
\kappa_1 r \leq \diam(\phi_\om( X_{t(\om)}))\le \kappa_2 r,
$$
is bounded  above by the number
$$
((1+\kappa_2)KD(R_S \kappa_1)^{-1})^d,
$$
where, we recall, $R_S$  is the radius of the largest open ball that
can be inscribed in all the sets $X_v$, $v \in V$.  
\elem

\bpf Recall that $\l_d$ is the Lebesgue measure on $\R^d$ and
that $V_d$ is the Lebesgue measure of the unit ball in $\R^d$.  
For every $v\in V$ let $x_v$ be the center of a
ball with radius $R_S$ which is contained in $\Int(X_v)$. 
Let $F$ be any collection of $A$--admissible words satisfying
the hypotheses of our lemma.  Then for  every $\om \in  F$, we have
$$ 
\phi_\om(X_{t(\om)})
\sbt  B\(x, r +\diam(\phi_\om(X_{t(\om)}))\) 
\sbt B(x, (1+ \kappa_2)r).
$$ 
Since, by the Open  Set  Condition, all the
sets $\{\phi_\om(\Int X_{t(\om)})\}_{\om \in F}$ are mutually
disjoint, applying Lemma~\ref{l42013_03_12}, we thus get 
$$ 
\begin{aligned}
V_d(1+\kappa_2)^dr^d
&= \l_d( B(x,(1+\kappa_2) r ))
 \geq \l_d\lt(\bigcup_{\om \in F} \phi_\om(X_{t(\om)})\rt)\\
&=  \sum_{\om \in F} \l_d(\phi_\om(\Int X_{t(\om)}))
\geq \sum_{\om \in F} \l_d\(\phi_\om(B(x_{t(\om)},R_S)\) \\
&\geq \sum_{\om \in F} \l_d \(B(\phi_\om(x_{t(\om)}),
K^{-1}R_S\| \phi_\om'\|_\infty)\)\\ 
& \geq \sum_{\om \in F} \l_d\(B(\phi_\om(x_{t(\om)}),
(KD)^{-1}R_S \diam(X_{i(\om)})\)\\ 
& \geq \sum_{\om \in F} \l_d\(B(\phi_\om(x_{t(\om)}),
(KD)^{-1}R_S \kappa_1 r )\)\\ 
&=\sharp F((KD)^{-1} R_S\kappa_1)^dV_dr^d \\
\end{aligned}
$$
Hence 
$$
\sharp F\leq ((1+\kappa_2)(KD) (R_S\kappa_1)^{-1})^d.
$$ 
We are done. 
\epf

\

Now assume that the alphabet $E$ is finite and keep the incidence
matrix $A$ (finitely) irreducible. Fix $t\ge 0$. Consider the
operator $\mathcal{L}_t$ given by the formula
\beq\label{1j89} 
\mathcal{L}_t g(\om):= \sum_{i:\, A_{i \om_1}=1}
g(i\om)| \phi_i'(\pi (\om)|^t, \quad  \om \in E^\infty_A,
\eeq
where $G$ is a real--valued function defined on $C(E^\infty_A)$. Note that  
$$
\mathcal{L}_t(C(E^\infty_A)) \sbt C(E^\infty_A)
$$ 
and the the linear operator $\mathcal{L}_t$ acts continuously on
$C(E^\infty_A)$, the Banach space of all real--valued continuous
functions on $E_A^\infty$ endowed with the supremum norm. In fact the
norm of $\mathcal{L}_t$ is bounded by $Z_1(t)$. A
straightforward inductive calculation gives that for all $n \geq 1$ we have
\beq\label{1j90} 
\mathcal{L}_t^n g(\om)= \sum_{\tau \in E^n_A, \tau \om \in 
E^\infty_A} g(\tau\om)|\phi_\tau'(\pi(\om))|^t.
\eeq

Let $\mathcal{L}_t^*: C^*(E^\infty_A)\lra C^*(
E^\infty_A)$
be the dual operator of $\L_t$. Denote by  $M_A$ the set of all Borel
probability measures on $E^{\N}_A$ considered as convex subset of the Banach space $C^*(E^\infty_A)$ via the canonical embedding
$$
M_A\ni m\lmt \lt(g\longmapsto \mu(g)=\int_{E^\infty_A}d\,dm\rt)\in C^*(E^\infty_A). 
$$
Consider the map
\beq\label{1j91}
M_A\ni m \longmapsto \frac{ \mathcal{L}_t^* m}{\mathcal{L}_t^* m( \1)}\in M_A.
\eeq
This map is  well defined  since
$$\mathcal{L}^*_t  m (\1)=m(\mathcal{L}_t
\1)=\int_{E^\mathbb{N}_A}  \sum_{i:\, A_{i \om_1}=1} |\phi_i'(\pi
(\om))|^tdm,
$$ 
and the  integrand  is everywhere  positive, whence the integral is positive. Since
this  map  is continuous in the weak$^*$--topology on $M_A$
and since $M_A$  is a compact (because $E$ is finite) convex subset
of the  locally convex 
topological vector space $C^*(E^\infty_A)$ endowed with the weak$^*$--topology, it follows  from the
Schauder--Tichonov Theorem  that the map  defined  in  (\ref{1j91})
has a fixed  point. Denote  it by  $\^m_t$ and put 
$$
\l_t=\mathcal{L}_t^*\^m_t(\1)>0.
$$ 
We then have
\beq\label{2j91}
\mathcal{L}_t^* \^m_t = \l_t \^m_t.
\eeq
Iterating (\ref{2j91}), we get for  every  $n \geq
1$ that
$$
\begin{aligned}
\l^n_t
&= \l_t^n\^m_t(\1)
=\mathcal{L}_t^{*n} \^m_t(\1)
 =\^m_t(\mathcal{L}^n_t \1)
=\int_{E^\infty_A}\sum_{\tau \in E^n_A: A_{\tau_n\om_1}=1} |\phi_\tau'(\pi
(\om))|^t\,d\^m_t\\
&\le \int_{E^\infty_A} \sum_{\tau \in E^n_A} \|\phi_\tau \|_\infty^td\,\^m_t
= \int_{E^\infty_A}  Z_n(t)\, d\^m_t\\
&=Z_n(t).
\end{aligned}
$$
Therefore, 
\beq\label{120130508}
\P(t)
=\lim_{n \to \infty}\frac{1}{n}\log Z_n(t)
\ge \lim_{n \to \infty}\frac{1}{n}\log\l^n_t
=\log \l_t.
\eeq
The formulas (\ref{1j89}) and  (\ref{1j90}) clearly  extend to all Borel  bounded  functions $g:E^\infty_A\lra \mathbb{R}$. The standard approximation arguments show then that
\beq\label{mtdual}
\mathcal{L}_t^{*n} \^m_t(g)=\l_t^n\^m_t(g)
\eeq
for all  Borel  bounded   functions  $g:E^\infty_A\lra
\mathbb{R}$. In particular, taking $g:=\1_{\phi_{\om(x)}} \circ \pi$ and
$\om \in E^*_A$, we get  with $n=|\om|$ that,
\beq\lab{1201200106}
\begin{aligned} \^m_t([\om])
& =\l_t^{-n}\mathcal{L}_t^{*n}\^m_t(\1_{[\om]})
  =\l_t^{-n}\^m_t\(\mathcal{L}_t^{n}(\1_{[\om]})\)\\
& =\l_t^{-n}\int_{E^\mathbb{N}_A}\sum_{\tau \in E^n_A, \,
A_{\tau_n\rho_1=1}}|\phi_\tau( \pi(\rho))|^t \1_{[\om]}( \tau \rho) \,d\^m_t(\rho)\\
& =\l_t^{-n}\int_{\rho \in E^\mathbb{N}_A}  \sum_{A_{\om_n
    \rho_1}=1}|\phi_\om'(\pi (\rho))|^t \, d\^m_t.  
\end{aligned}
\eeq
From this formula, we get
\beq\label{420130508}
\^m_t([\om])
\leq \l_t^{-n}\|\phi_\om'\|_\infty^t 
 \^m_t\(\{ \rho \in E^\infty_A:  A_{\om_n\rho_1}=1\}\)
\leq \l_t^{-n}\|\phi_\om '\|_\infty^t,
\eeq
and
\beq\label{2j92}
\begin{aligned}
 \^m_t([\om])
& \geq K^{-t}\l_t^{-n}\|\phi_\om'\|_\infty^t \^m_t\(\{ \rho \in E^\infty_A:  A_{\om_n\rho_1}=1\}\)
\\
& =K^{-t}\l_t^{-n}\|\phi_\om '\|_\infty^t \^m_t\lt(\bu_{ e\in E,\, A_{\om_n e=1}} [e]\rt)\\
& = K^{-t}\l_t^{-n}\|\phi_\om'\|_\infty^t \sum_{ e\in E, A_{\om_n e=1}}\^m_t ([e]).
\end{aligned}
\eeq
Now, since $\^m_t\(E^\infty_A\)=1$, there is $b \in E$ such that
$\^m_t([b])>0$. Let $\La$ be the finite subset of $E_A^*$ witnessing
finite irreducibility of the matrix $A$. Put
$$
M:=\min\{\|\phi_\tau'\|_\infty:\tau\in\La\}.
$$
Then for every $\om\in E^*_A$, there exists a word $\tau\in\La$ such
that  $\om \tau b\in E^*_A$.  We therefore get  from (\ref{2j92}) that
\beq\label{320130508}
\aligned
\^m_t([\om])
&\geq \^m_t ( [\om \tau]) 
\geq K^{-t}\l_t^{-n}\|\phi_{\om \tau}'\|_\infty^t \^m_t([b]) \\
&\geq K^{-2t}\^m_t([b])\|\phi_\tau'\|_\infty^t\l_t^{-n}\|\phi_\om'\|_\infty^t 
\\
&\geq K^{-2t}\^m_t([b])M^t\l_t^{-n}\|\phi_\om'\|_\infty^t \\
&>0.
\endaligned
\eeq
Rearranging terms we get,
\beq\label{520130508}
\|\phi_\om'\|_\infty^t\le K^{2t}(\^m_t([b])M^t)^{-1}\l_t^n\^m_t([\om]).
\eeq
Summing over all $A$-admissible words of length $n$, this gives
$$
Z_n(t)
\le K^{2t}(\^m_t([b])M^t)^{-1}\l_t^n\sum_{\om\in E_A^n}\^m_t([\om])
= K^{2t}(\^m_t([b])M^t)^{-1}\l_t^n.
$$
Therefore, 
$$
\P(t)
=\lim_{n \to \infty}\frac{1}{n}\log Z_n(t)
\le \lim_{n \to \infty}\frac{1}{n}\log\l_t^n
=\log\l_t.
$$
Along with \eqref{120130508} this yields,
\beq\label{220130508}
\log\l_t=\P(t).
\eeq
Since the system $S$ is finite, with $h$ defined in \eqref{120200108}, we have $\P(h)=0$, and therefore formulas
\eqref{420130508} and \eqref{520130508}, applied with $t=h$, yield
\beq\label{2j93}
M_h\|\phi_\om'\|_\infty^h\leq \^m_h([\om])\leq \| \phi_\om'\|_\infty^h
\eeq
with some positive constant $M_h$. Let
\beq\label{220131007}
m_h:=\^m_h\circ \pi^{-1}.
\eeq
and more generally
\beq\label{320131007}
m_t:=\^m_t\circ \pi^{-1}
\eeq
for all $t\ge 0$. We call $m_h$ the $h$--conformal measure for the system $S$. More generally, we call $m_t$ the $(t,e^{\P(t)})$--conformal measure for the system $S$. We shall prove the following.

\sp\bthm \label{t1j93} 
Suppose \,  $S=\{\phi_e\}_{e\in E}$ is
a finite  irreducible conformal  GDMS. Then
\begin{itemize}
\item[{\rm (a)}] There exists a constant $C \geq 1$  such that  for all
$x\in J_S$ and all $ 0<r \leq 1$, we have  that
$$ 
C^{-1}\leq  \frac{m_h(B(x,r))}{r^h }\leq C,
$$
\item[{\rm (b)}]  $0< \H_h(J_S), \Pi_h(J_S) <  +\infty$,

\sp\item[{\rm (c)}]  $\HD(J_S)=h$.
\end{itemize}
\ethm

\bpf Since the alphabet $E$ is finite, we have that
$$
\xi:=\inf\{\| \phi_e'\|_\infty: e \in  E\}>0.
$$
Fix  $x\in J_S$ and $0< r < \frac{1}{2}\min  \{  \diam(X_v):\, v \in
V\}$. Then   $x= \pi(\tau)$  for some $\tau \in E^\infty_A$. Let $
n=n(\tau) \geq 0$ be the least integer such that
$$
\phi_{\tau_{|n}}\(X_{t(\tau_n)}\)\sbt B(x,r).
$$
Then by (\ref{2j93}) we have the following.
$$
m_h(B(x,r))
\geq m_h(\phi_{\tau_{|n}}(X_{t(\tau_n)})
\geq\^m_h([\tau_{|n}]))
\geq  M_h K^{-h}\|\phi_{\tau_{|n}}'\|_\infty^h.
$$ 
By the definition of $n$ and Corollary~\ref{c23013_03_12}, we have 
$$
r \leq
\diam(\phi_{\tau_{|n-1}}(X_{t(\tau_{n-1})}))
\leq D \|\phi_{\tau_{|n-1}}'\|_\infty.
$$
Hence
\beq\label{1j95}
 m_h(B(x,r))\geq M_h(DK)^{-h} r^h.
\eeq
In order  to prove the opposite inequality, let $Z$  be
the family of all minimal (in  the sense  of length) words $\om \in
E^*_A$ such that
\beq\label{2j95}
 \phi_\om(X_{t(\om)})\cap B(x,r)
\neq \es \quad \mbox{and} \quad \phi_\om(X_{t(\om)})\sbt B(x,2r).
\eeq
Consider an arbitrary $\om \in Z$. Then
\beq\label{120130320}
\diam\(\phi_\om(X_{t(\om)})\)\le 2r
\eeq
and $\diam(\phi_{\om_{|_{|\om|-1}}}(X_{t(\om)}))\geq 2r$. Therefore,
making use of Corollary~\ref{c23013_03_12} twice, we get the following.
\beq\label{220130320}
\begin{aligned}
\diam\(\phi_\om(X_{t(\om)})\)
&\ge D^{-1}\|\phi_\om'\|_\infty
 \geq (DK)^{-1}\| \phi_{\om_{|_{|\om|-1}}}\|_\infty\cdot\|\phi_{\om_{|\om|}}\|_\infty
 \\
& \geq K^{-1}\xi(DK)^{-1}\diam\(\phi_{\om_{|_{|\om|-1}}}(X_{t(\om_{|\om|-1})})\)
\\
& \geq 2(KD^2)^{-1}\xi r.
\end{aligned}
\eeq
Hence, by virtue of formula \eqref{4.1.10a} of
Lemma~\ref{l52013_03_12}, we get that 
\beq\label{291210}
\phi_\om(X_{t(\om)}) \supset B\(\phi_\om(x_{t(\om)}),2(K^2DL)^{-1}R_S\xi r\)
\eeq
Since, by its very definition, the family $Z$ consists of mutually
incomparable words, Lemma~\ref{l1j81} along with \eqref{120130320} and
\eqref{220130320} imply that
\beq\label{3j95}
\sharp Z \leq \Ga:=((1+2(KD^2)^{-1})2KD(R_S)^{-1})^Q.
\eeq
Since 
$$
\pi^{-1}(B(x,r)) \sbt  \bigcup_{\om \in Z} [\om],
$$
we get from (\ref{2j95}), \eqref{2j93}, \eqref{l52013_03_12},
\eqref{120130320}, and (\ref{3j95}) that 
$$\begin{aligned}
m_h(B(x,r))& 
=\^m_h \circ \pi^{-1}(B(x,r)) 
\leq \^m_h\lt(\bigcup_{\om \in Z} [\om]\rt)
= \sum_{\om \in Z}\^m_h([\om]) 
\leq \sum_{\om \in Z}\|\phi_\om'\|_\infty^h\\
&\leq \sum_{\om \in Z}(D \diam(\phi_{\om}(X_{t(\om)})))^h \leq
(4D)^h \sum_{\om \in Z} r^h
=(4D)^h\sharp Z r^h\\
& \leq (4D)^h \Ga r^h.
\end{aligned}$$
Along with (\ref{1j95}) this proves item (a). Items (b) and (c)
are now an immediate consequence of Frostmann's Converse Lemmas, i. e. 
Theorem~\ref{t23hg53} and Theorem~\ref{t24hg59}. The proof is complete.
\epf 

\sp We are now ready to provide a short simple proof of the following
main theorem of this section proven in \cite{MU2}, comp. \cite{MU1}.

\sp\bthm[Bowen's Formula]\label{t1j97} 
If $S$ is a finitely irreducible conformal GDMS, then
$$
h_S=\inf\{ t\geq 0: \,\, \P(t)\leq 0 \}
=\HD(J_\mathcal{S}) 
=\sup \{\HD(J_F):  \, F \sbt E \, \, \mbox{and} \, \, F \,\,  \mbox{is
  finite}\,\}.  
$$
\ethm
\bpf The first equality is just the definition of $h_S$. Put 
$$
h_\infty:=\sup \{\HD(J_F):  \,\, F \sbt E \, \,
\mbox{and} \,\, \, F \,\,  \mbox{is finite}\,\}
$$ 
and 
$$
H:=\HD(J_S).
$$
Fix $t> h$. Then $\P(t)< 0$. Therefore, for all $n\geq 0$ large enough, we have 
$$  
\sum_{\om \in E^n_A}\|\phi_\om'|_\infty^t \leq \exp \lt( \frac{1}{2}\P(t)n\rt).
$$
Hence
$$
\sum_{\om \in E^n_A}\diam^t(\phi_\om(X_{t(\om}))
\leq D^t \sum_{\om \in E^n_A}\|\phi_\om'\|_\infty^t
\leq D^t \exp\lt( \frac{1}{2}\P(t)n\rt).
$$ 
Since the
family $\big\{\phi_\om (X_{t(\om)})\big\}_{\om \in E^n}$ forms a
cover  of $J_S$, letting $ n \to \infty$, we thus get  that $\H_t(J_S)=0$.
This implies  that $t \geq H$, and in consequence, $h_S\geq H$. Since
obviously, $h_\infty \leq H$, we thus have
$$  
h_\infty \leq H \leq h_S.
$$
We are left to show that $h_S\leq h_\infty$. If $F$  is a finite
subset of $E$, then  $h_F\leq h_\infty$, and by virtue  of
Theorem~\ref{t1j93}, $\P_F(h_\infty)\leq 0$. So, in view of
Theorem~\ref{t2.1.3} and Lemma~\ref{l1j85}, we have
$$
\P(h_\infty)=\sup\{\P_F(h_\infty): F\sbt E, \,
\mbox{and}\,\,F\,\, \mbox{is finite} \}\leq 0.
$$
Hence $h_\infty \geq h_S$ and the proof is complete. 
\epf

\sp Combining  this theorem and Fact~\ref{f1j87}, we get the following.

\sp\bthm\label{HD-theta}
If  $S$ is a finitely irreducible conformal GDMS,
then 
$$
\HD(J_S)=h_S\geq \theta(\mathcal{S}).
$$
If in addition $S$ is strongly regular, in particular  if $S$
is co--finitely regular, then 
$$
\HD(J_S)=h_S > \theta(S).
$$ 
\ethm

\section[Conformal and Invariant Measures for Conformal GDMSs]{Conformal and Invariant Measures for Conformal GDMSs}
Sticking to the previous section, 
$$
S=\{\phi_e:e\in E\}
$$ 
is assumed to be a finitely irreducible conformal graph directed Markov system with a countable alphabet $E$. Obviously, for all 
$$
t\in F(S)=\{t\ge 0:\P(t)<+\infty\}=\{t\ge 0:Z_1(t)<+\infty\},
$$
formula \eqref{1j89} defines a bounded linear operator from $C_b(E_A^\infty)$, into itself, where $C_b(E_A^\infty)$ is the Banach space of all bounded real (or complex)--valued continuous functions on $E_A^\infty$ endowed with the supremum norm. This operator is here also denoted by $\L_t$. Our primary goal in this section is to prove for every $t\in F(S)$ the existence of a Borel probability measure $\tilde m_t$ on $E_A^\infty$ satisfying equation 
\eqref{2j91} (and with $\l_t=e^{\P(t)}$). We first need the following auxiliary result. We will provide its short proof for
the sake of completeness and convenience of the reader. It is more natural and convenient to formulate it
in the language of directed graphs. Let us recall that a directed graph is
said to be strongly connected if and only if its incidence matrix is
irreducible. In other words, it means that every two vertices can be joined
by a path of admissible edges.

\blem\lab{l2.7.2}
If $\Ga=<E,V>$ is a strongly connected directed
graph, then there exists a sequence of strongly connected
subgraphs $<E_n, V_n>$ of $\Ga$ such that all the vertices $V_n\sbt V$
and all the edges $E_n$ are finite,
$\{V_n\}_{n=1}^\infty$ is an ascending sequence of vertices,
$\{E_n\}_{n=1}^\infty$ is an ascending sequence of edges,
$$
\bu_{n=1}^\infty V_n=V
\  \  \  {\rm and} \  \  \ 
\bu_{n=1}^\infty E_n=E.
$$
\elem

\bpf Let 
$$
V=\{v_n:n\ge 1\}
$$ 
be a sequence of all vertices of $\Ga$ and let 
$$
E=\{e_n:n\ge 1\}
$$ 
be a sequence of edges
of $\Ga$. We will proceed inductively to construct the sequences 
$$
\{V_n\}_{n=1}^\infty 
\  \  \  {\rm and} \  \  \
\{E_n\}_{n=1}^\infty.
$$
In order to
construct $<E_1,V_1>$ let $\a$ be a path joining $v_1$ and $v_2$
($i(\a)=v_1$, $t(\a)=v_2$) and let $\b$ be a path joining $v_2$ and $v_1$
($i(\b)=v_2$, $t(\b)=v_1$). These paths exist since $\Ga$ is strongly
connected. We define $V_1\sbt V$ to be the set of all vertices of
paths $\a$ and $\b$ and $E_1\sbt E$ to be the set of all edges from
$\a$ and $\b$ enlarged by $e_1$ if this edge is
among all the edges joining the vertices of $V_1$. Obviously
$<E_1,V_1>$ is strongly
connected and the first step of inductive procedure is
complete. 

Suppose now that a strongly connected graph $<E_n, V_n>$ has
been constructed. 

If $v_{n+1}\in V_n$, we set $V_{n+1}=V_n$ and
$E_{n+1}$ is then defined to be the union of $E_n$ and all the edges
from $\{e_1,e_2,\ld, e_n,e_{n+1}\}$ that are among all the edges
joining the vertices of $V_n$. 

If $v_{n+1}\notin V_n$, let $\a_n$ be a
path joining $v_n$ and $v_{n+1}$ and let $\b_n$ be a path joining $v_{n+1}$
and $v_n$. We define $V_{n+1}$ to be the union of $V_n$ and the set of all
vertices of $\a_n$ and $\b_n$. $E_{n+1}$ is then defined to be the
union of $E_n$, all the edges building the paths $\a_n$ and $\b_n$ and all
the edges from $\{e_1,e_2,\ld, e_n,e_{n+1}\}$ that are among all the edges
joining the vertices of $V_{n+1}$. Since $<E_n,
V_n>$ was strongly connected, so is $<E_{n+1}, V_{n+1}>$. The
inductive procedure is complete. It immediately follows from the
construction that 
$$
V_n\sbt V_{n+1}, \ E_n\sbt
E_{n+1}, \ \bu_{n=1}^\infty V_n=V,
\  \  \  {\rm and} \  \  \
\bu_{n=1}^\infty E_n=E.
$$
We are done. 
\epf

\sp The following theorem was proved in \cite{MU2}, comp. also Corollary~6.32 and Theorem~7.4 in \cite{CTU}, in a more general setting of arbitrary H\"older continuous summable potentials rather than merely $t\zeta$. We provide here its proof for the convenience of the reader and the sake of completeness.

\sp\bthm\label{t220131009} 
If $S$ is a finitely irreducible graph directed Markov system then, for all $t\in F(S)$ there exists a unique Borel probability measures $\tilde m_t$ on $E_A^\infty$ such that \eqref{2j91} holds, i.e.
\beq\label{2j91N}
\mathcal{L}_t^* \^m_t = e^{\P(t)} \^m_t.
\eeq
\ethm 

\bpf Without loss of generality we may assume that $E=\N$. Since
the incidence matrix $A$ is irreducible, it follows from
Lemma~\ref{l2.7.2} that we can reorder the set $\N$ such that there exists a sequence $\(l_n\)_{n\ge 1}$, increasing to infinity, and such that for every $n\ge 1$ the matrix 
$$
A|_{\{1,\ld,l_n\}\times
\{1,\ld,l_n\}}
$$ 
is irreducible. Given an integer $q\ge 1$ denote
$$
\N(q):=\{1,2,\ld, q\}\sbt E=\N.
$$
In view of the formula \eqref{2j91N}, holding for finite alphabets,
there exists an eigenmeasure $\^m_n$ of the operator $\L_n^*$, conjugate to the Perron--Frobenius operator\index{(N)}{Perron--Frobenius operator} 
$$
\L_n:C\(\N(l_n)_A^\infty\)\lra C\(\N(l_n)_A^\infty\)
$$
associated to the potential $t\zeta\big|_{\N(l_n)_A^\infty}$.
More precisely, for every function $g\in C\(\N(l_n)_A^\infty\)$,
$$
\L_ng(\om):= \sum_{i\in \N(l_n):\, A_{i \om_1}=1}
g(i\om)| \phi_i'(\pi (\om)|^t, \quad  \om \in \N(l_n)^\infty_A.
$$
A family of Borel probability measures $\cM$ in a topological space $X$ is called tight \index{(N)}{tight family of measures}, or uniformly tight if and only if for every $\e>0$ there exists a compact set $K_\e\sbt X$ such that 
$$
\mu(K_\e) > 1-\e
$$ 
for all $\mu \in \cM$. If $X$ is a completely metrizable space and $\cM$ is a tight family of Borel probability measures then Prohorov's Theorem, see e.g. \cite[Book II, Theorem~8.6.2]{bog:meas}, asserts that every sequence in $\cM$ contains a weakly convergent subsequence. 

We will prove the following. 

\sp {\sl Claim}~$1^0$. The sequence $\{\^m_n\}_{n\ge 1}$ is tight with all measures $\^m_n$, $n\ge 1$, treated as Borel probability measures on $E_A^\infty$. 

\bpf 
Let 
$$
\P_n(t):=\P_{\N(l_n)}(t).
$$
Obviously $\P_n(t)\ge \P_1(t)$ for all $n\ge 1$. For every $k\ge 1$ let $\pi_k:E_A^\infty\lra \N$ be the projection
onto the $k$-th coordinate, i.e. 
$$
\pi\(\{(\tau_u)_{u\ge 1}\}\):=\tau_k.
$$
By formula \eqref{220130508}, $e^{\P_n(t)}$ is the eigenvalue of $\L_n^*$ corresponding to the eigenmeasure $\^m_n$.

Therefore, applying (\ref{1201200106}), we obtain for every $n\ge
1$, every $k\ge 1$, and every $s\in\N$ that
$$
\aligned
\^m_n(\pi_k^{-1}(s))
&=\sum_{\om\in \N(l_n)_A^k:\om_k=s}\^m_n([\om])
\le \sum_{\om\in \N(l_n)_A^k:\om_k=s}\|\phi_\om'\|_\infty^te^{-\P_n(t)k}
  \\
&\le e^{-\P_n(t)k}\sum_{\om\in \N(l_n)_A^k:\om_k=s}\|\phi_{\om|_{k-1}}'\|_\infty^t\|\phi_s'\|_\infty^t
  \\
&\le e^{-\P_1(t)k}\lt(\sum_{i\in \N}\|\phi_i'\|_\infty^t\rt)^{k-1}\|\phi_s'\|_\infty^t.
\endaligned
$$
Therefore
$$
\^m_n\(\pi_k^{-1}([s+1,\infty))\)
\le e^{-\P_1(t)k}\lt(\sum_{i\in \N}\|\phi_i'\|_\infty^t\rt)^{k-1}
\sum_{j=s+1}^\infty\|\phi_s'\|_\infty^t.
$$
Fix now $\e>0$ and for every $k\ge 1$ choose an integer $n_k\ge 1$ so large that
$$
e^{-\P_1k}\lt(\sum_{i\in \N}\|\phi_i'\|_\infty^t\rt)^{k-1}
\sum_{j=n_k+1}^\infty\|\phi_s'\|_\infty^t 
\le \frac{\e}{2^k}.
$$
Then, for every $n\ge 1$ and every $k\ge 1$, 
$$
\^m_n\(\pi_k^{-1}([n_k+1,\infty))\)\le \e/2^k.
$$
Hence
$$
\^m_n\lt(E_A^\infty\cap\prod_{k\ge 1}[1,n_k]\rt)
\ge 1-\sum_{k\ge 1}\^m_n\(\pi_k^{-1}([n_k+1,\infty))\)
\ge 1-\sum_{k\ge 1} \frac{\e}{2^k}
=1-\e.
$$
Since $E_A^\infty\cap\prod_{k\ge 1}[1,n_k]$ is a compact subset of $E_A^\N$,
the tightness of the sequence $\{\^m_n\}_{n\ge 1}$ is therefore proved.
\epf

\sp Thus, in view of Prohorov's Theorem 
there exists $\^m$, a Borel probability measure on $E_A^\N$, which is a weak*--limit point of the sequence $\{\^m_n\}_{n\ge 1}$. Passing to a subsequence, we may assume that the sequence $\{\^m_n\}_{n\ge 1}$ itself converges weakly to the measure $\^m$. Let
$$
\L_{0,n}=e^{-\P_n(t)}\L_n \  \  \  {\rm and } \   \   \   \L_0=e^{-\P(t)}\L_t
$$
be the corresponding normalized operators. Fix $g\in
C_b(E_A^\infty)$ and $\e>0$. Let us now consider an integer $n\ge 1$ so large
that the following requirements are satisfied.
\beq\label{2.7.1}
\sum_{i=n+1}\|g\|_\infty \|\phi_i'\|_\infty^te^{-\P(t)} \le \frac{\e}{6},
\eeq
\beq\lab{2.7.2}
\sum_{i\le n}\|g\|_\infty\|\phi_i'\|_\infty^t\big|e^{-\P(t)}- e^{-\P_n(t)}\big|
\le \frac{\e}{6},
\eeq
\beq\label{2.7.3}
|\^m_n(g)-\^m(g)|\le \frac{\e}{3},
\eeq
and
\beq\label{2.7.4}
\left|\int_{E_A^\infty}\L_0(g)\,d\^m -\int_{E_A^\infty}\L_0(g)\,d\^m_n\right|\le \frac{\e}{3}.
\eeq
It is possible to make condition (\ref{2.7.2}) satisfied since, due to
Theorem~\ref{t2.1.3}, $\lim_{n\to\infty}\P_n(t)=\P(t)$. Let
$$
g_n:=g|_{E_{l_n}^\infty}.
$$
The first two observations are the following.
\beq\label{2.7.5}
\aligned
\L_{0,n}^*\^m_n(g)
&=\int_{E_A^\infty}\sum_{i\le n:A_{i\om_n}=1}g(i\om)\big|\phi_i'(\pi(\sg(\om)))\big|^t e^{-\P_n(t)}\,d\^m_n(\om) 
 \\
&=\int_{\N(l_n)_A^\infty}\sum_{i\le n:A_{i\om_n}=1}g(i\om)\big|\phi_i'(\pi(\sg(\om)))\big|^t e^{-\P_n(t)}\,d\^m_n(\om)
 \\
&=\int_{\N(l_n)_A^\infty}\sum_{i\le n:A_{i\om_n}=1}g_n(i\om)\big|\phi_i'(\pi(\sg(\om)))\big|^t e^{-\P_n(t)}\,d\^m_n(\om)
\\
&=\L_{0,n}^*\^m_n(g_n)
 =\^m_n(g_n),
\endaligned
\eeq
and
\begin{equation}\label{2.7.6}
\^m_n(g_n)-\^m_n(g)
=\int_{\N(l_n)^\N}(g_n-g)d\^m_n
=\int_{\N(l_n)^\N}0d\^m_n
=0.
\end{equation}
Using the triangle inequality we get the following.
\begin{equation}\label{2.7.7}
\aligned
\big|\L_0^*\^m(g)-\^m(g)\big|
&\le |\L_0^*\^m(g)-\L_0^*\^m_n(g)| + |\L_0^*\^m_n(g)-\L_{0,n}^*\^m_n(g)|+ \\
&\  \  \  \  \ +|\L_{0,n}^*\^m_n(g)-\^m_n(g_n)|+|\^m_n(g_n)-\^m_n(g)|
   +|\^m_n(g)-\^m(g)|.
\endaligned
\end{equation}
Let us look first at the second summand. Applying (\ref{2.7.2}) and
(\ref{2.7.1}) we get
\begin{equation}\label{2.7.8}
\aligned
|\L_0^*\^m_n(g) &-\L_{0,n}^*\^m_n(g)|= \\
&=\bigg|\int_{E_A^\infty}\sum_{i\le n:A_{i\om_n}=1}g(i\om)\(\big|\phi_i'(\pi(\sg(\om)))\big|^t e^{-\P(t)}-\big|\phi_i'(\pi(\sg(\om)))\big|^t e^{-\P_n(t)}\)\, d\^m_n(\om) 
\\
&\  \  \   \   \   \  \  \  \  \  \  \   \  \  \  + \int_{E_A^\infty}
\sum_{i> n:A_{i\om_n}=1}g(i\om)\big|\phi_i'(\pi(\sg(\om)))\big|^t e^{-\P(t)}\,d\^m_n(\om)\bigg| 
\\
&\le \sum_{i\le n}\|g\|_\infty\|\phi_i'\|_\infty
     \big|e^{-\P(t)}- e^{-\P_n}\big|+\sum_{i=n+1}\|g\|_\infty^t \|\phi_i'\|_\infty ^t e^{-\P(t)}
 \\
&\le \frac{\e}{6}+\frac{\e}{6}
=\frac{\e}{3}.
\endaligned
\end{equation}
Combining now in turn (\ref{2.7.4}), (\ref{2.7.8}), (\ref{2.7.5}), (\ref{2.7.6}) and (\ref{2.7.3}) we get
from (\ref{2.7.7}) that
$$
|\L_0^*\^m(g)-\^m(g)|\le \frac{\e}{3}+\frac{\e}{3}+\frac{\e}{3}=\e.
$$
Letting $\e\downto 0$ we therefore get $\L_0^*\^m(g)=\^m(g)$ or $\L_t^*\^m(g)
=e^{\P(t)}\^m(g)$. Hence 
$$
\L_t^*\^m=e^{\P(t)}\^m,
$$ 
and the proof is complete.
\epf

\sp\brem\label{r120200107}
Note that once we have an eigenmeasure of $\L_t^*$, then all the formulas from \eqref{2j91} to \eqref{120130508} hold for any countable alphabet $E$ regardless of whether it is finite or not. 
{1201200106}
\erem

\sp We will now introduce the concept of Gibbs states and shift--invariant Gibbs states for a parameter $t\in F(S)$. These play for us in this book only a very limited auxiliary role; they are used by us only to simplify notation, formulation of some results, and proofs. This concept can be however extended to the class of all H\"older continuous potential on $E_A^\infty$ leading to a rich meaningful and powerful theory whose systematic account can be found in \cite{MU2}, comp. also \cite{CTU} and \cite{MU5}.

If $S$ is a finitely irreducible graph directed Markov system with an incidence matrix $A$ and $t\in F(S)$, then a Borel probability
measure $\^m$ on $E_A^\infty$ is called a  Gibbs state \index{(N)}{Gibbs state} for $t$ if there exist constants $Q_t\ge 1$ and $\P_{\^m}\in\R$ such that for every $\om\in E_A^*$ and every $\tau\in [\om]$,
\begin{equation}\lab{2.2.1}
Q_t^{-1}
\le \frac{\^m([\om])}{\big|\phi_\om'\(\pi(\sg^{|\om|}(\tau)\)\big|^t\exp\(-\P_{\^m}|\om|\)}
\le Q_t.
\end{equation}
If additionally $\^m$ is shift--invariant, then $\^m$ is called an
invariant Gibbs state. \index{(N)}{invariant Gibbs state} Because of the Bounded Distortion Property (4e), formula \eqref{2.2.1} takes also on the following somewhat simpler form.

For every $\om\in E_A^*$ 
\begin{equation}\lab{2.2.1B}
Q_t^{-1}
\le \frac{\^m([\om])}{\big\|\phi_\om'\big\|_\infty^t\exp\(-\P_{\^m}|\om|\)}
\le Q_t.
\end{equation}

\bprop\lab{p2.2.2}
If $S$ is a finitely irreducible graph directed Markov system with an incidence matrix $A$ and $t\in F(S)$, then the following hold:
\begin{itemize}
\item[(a)] For every Gibbs state $\^m$ for $t$, $\P_{\^m}=\P(t)$.
\item[(b)] Any two Gibbs states for the function $t$ are boundedly equivalent
with Radon--Nikodym derivatives bounded away from zero and infinity.
\end{itemize}
\eprop

\bpf We shall first prove (a). Towards this end fix
$n\ge 1$ and sum up (\ref{2.2.1}) over all words
$\om\in E_A^n$. Since 
$$
\sum_{|\om|=n}\^m([\om])=1,
$$
we then get
$$
Q_t^{-1}e^{-\P_{\^m}n}\sum_{|\om|=n}\big\|\phi_\om'\big\|_\infty^t
\le 1
\le Q_te^{-\P_{\^m}n}\sum_{|\om|=n}\big\|\phi_\om'\big\|_\infty^t.
$$
Applying logarithms to all three terms of this formula, dividing all
the terms by $n$ and taking the limit as $n\to\infty$, we obtain
$$
-\P_{\^m}+\P(t)\le 0\le -\P_{\^m}+\P(t)
$$ 
which means that
$$
\P_{\^m}=\P(t).
$$
The proof of item (a) is thus complete.

\sp In order to prove part (b) suppose that $m$ and $\nu$ are two
Gibbs states of $t$. Notice now that part (a) implies the
existence of a constant $T\ge 1$ such that
$$
T^{-1}\le \frac{\nu([\om])}{m([\om])}\le T
$$
for all words $\om\in E_A^*$. Straightforward reasoning gives now that
$\nu$ and $m$ are equivalent and 
$$
T^{-1}\le \frac{d\nu}{dm}\le T.
$$
The proof of Proposition~\ref{p2.2.2} is complete.
\epf

\bthm\lab{t2.3.3}
If $S$ is a finitely irreducible graph directed Markov system with an incidence matrix $A$ and $t\in F(S)$, then the any eigenmeasure $\^m$ of the dual operator $\mathcal{L}_t^*: C^*(E^\infty_A)\lra C^*(
E^\infty_A)$ is a Gibbs state for $t$. In addition, its corresponding eigenvalue $\l$ is equal to $e^{\P(t)}$. 
\ethm

\bpf Let $\La$ be a minimal set which witnesses the finite
irreducibility of $A$. Let
$$
T:=\min\big\{\inf(|\phi_\a'|) e^{-\P(t)|\a|}:\a\in \La\big\}\in\
$$
Fix $\om\in E_A^n$ and put $n:=|\om|$. For every $\a\in\La$, let
$$
E_\a(\om):=\big\{\tau\in E_A^\infty:\om\a\tau\in E_A^\infty\big\}.
$$
By the definition of $\La$, $\bu_{\a\in\La}\a(\om)=E_A^\infty$. Hence, there exists $\g\in \La$ such that 
$$
\^m(E_\g)\ge (\#\La)^{-1}.
$$
Writing $p=|\g|$ and using the Bounded Distortion Property (4e), we therefore  get
\begin{equation}\lab{2.3.1227b}
\aligned
\^m([\om])
&\ge \^m([\om\g])
=\l^{-(n+p)} \int_{\rho\in E_A^\infty:A_{\g_p\rho_1}=1}
\big|\phi_{\om\g}'(\pi(\rho))\big|^t\,d\^m(\rho)
 \\
&=\l^{-(n+p)} \int_{\rho\in E^\infty:A_{\g_p\rho_1}=1}
\big|\phi_{\om}'(\pi(\g\rho))\big|^t\cdot\big|\phi_{\g}'(\pi(\rho))\big|^t
\, d\^m(\rho) 
 \\
&\ge \l^{-n}T\l^{-p} \int_{\rho\in E_A^\infty:A_{\g_p\rho_1}=1}\big|\phi_{\om}'(\pi(\g\rho))\big|^t\,d\^m(\rho) 
 \\
&=\l^{-n}T\l^{-p} \int_{E_\g(\om)}\big|\phi_{\om}'(\pi(\g\rho))\big|^t\,d\^m(\rho)
\ge T\l^{-p}K^{-t}\^m(E_\g(\om))\l^{-n}\big\|\phi_{\om}'\big\|_\infty^t
\\
&\ge T\l^{-p}K^{-t}(\#\La)^{-1}\l^{-n}\big\|\phi_{\om}'\big\|_\infty^t,
\endaligned
\end{equation}
Along with \eqref{420130508}, which holds because of Remark~\ref{r120200107}, this shows that $\^m$ is a Gibbs state for $t$. The equality $\l=e^{\P(t)}$ follows now immediately from Proposition~\ref{p2.2.2}. The proof of Theorem~\ref{t2.3.3} is complete.
\epf

\bthm\lab{t2.2.4} 
If $S$ is a finitely irreducible graph directed Markov system with an incidence matrix $A$, then for every $t\in F(S)$ there exists a unique shift--invariant Gibbs state $\^\mu_t$ of $t$. The shift--invariant Gibbs state $\^\mu_t$ is ergodic. In addition, if the incidence matrix $A$ is finitely primitive, then the Gibbs state $\^\mu_t$ is completely ergodic.
\ethm

\bpf Let $\^m$ be a Gibbs state for $t\in F(S)$. Fixing $\om\in
E_A^*$, using (\ref{2.2.1}), and Proposition~\ref{p2.2.2}(a) we get for every $n\ge 1$
\begin{equation}\lab{2.2.2}
\aligned
\^m(\sg^{-n}([\om]))
&=\sum_{\tau\in E_A^n:A_{\tau_n\om_1=1}}\^m([\tau\om]) 
\le \sum_{\tau\in E_A^n:A_{\tau_n\om_1}=1}
Q_t\|\phi_{\tau\om}\|_\infty \exp(-\P(t)(n+|\om|)
\\
&\le \sum_{\tau\in E_A^n:A_{\tau_n\om_1}=1}Q_t
\|\phi_\tau'\|_\infty e^{-\P(t)n}\|\phi_\om'\|_\infty e^{-\P(t)|\om|}
\\
&\le \sum_{\tau\in E_n^\om}Q_tQ_t\^m([\tau])Q_t\^m([\om])\\
&\le Q_t^3\^m([\om]).
\endaligned
\end{equation}
Let the finite set of words $\La$ witness the finite irreducibility of the incidence matrix $A$ and let $p$ be the
maximal length of a word in $\La$. Let
$$
T:=\min\big\{\inf(|\phi_\a'|) e^{-\P(t)|\a|}:\a\in \La\big\}\in(0,+\infty).
$$
For each $\tau,\om \in E_A^*$, let $\a = \alpha (\tau,\om) \in \La$ be such
that $\tau\alpha\om \in E_A^*.$ Then, we have for all $\om\in E_A^*$
and all integers $n\ge 1$ that
\begin{equation}\lab{2.2.3}
\aligned
 \sum_{i= n}^{n+p} \^m(\sg^{-i}([\om]))
&=\sum_{i= n}^{n+p}\sum_{\tau\in E_A^i:A_{\tau_i\om_1}=1}\^m([\tau\om]) 
\ge \sum_{\tau\in E_A^n} \^m([\tau\alpha (\tau,\om)\om]) 
  \\
&\ge \sum_{\tau\in E_A^n}
  Q_t^{-1}\inf\(\big|\phi_{\tau\alpha (\tau,\om)\om}'\big|^t\)
  \exp\(-\P(t)(|\tau|+|\alpha (\tau,\om)|+|\om|)\) 
  \\
&\ge Q_t^{-1}\sum_{\tau\in E_A^n}
 \inf\(\big|\phi_{\tau}'\big|^t\)\exp\(-\P(t)(|\tau|)\)\cdot \\
&\qquad\  \  \  \ \cdot\inf\(\big|\phi_{\alpha(\tau,\om)}'\big|^t\)
  \exp\(-\P(t)|\alpha (\tau,\om)|\) 
 \inf\(\big|\phi_{\om}'\big|^t\)\exp\(-\P(t)(|\om|)\)
 \\
&\ge Q_t^{-1}T\inf\(\big|\phi_{\om}'\big|^t\)\sum_{\tau\in E_A^n}\exp\(-\P(t)|\om|\)\inf\(\big|\phi_{\tau}'\big|^t\)\exp\(-\P(t)(|\tau|)\) 
 \\
&\ge Q_t^{-2}T\^m([\om])
\sum_{\tau\in E_A^n}\exp\(-\P(t)|\om|\)\inf\(\big|\phi_{\tau}'\big|^t\)\exp\(-\P(t)(|\tau|)\)
\\
&\ge Q_t^{-2}T\^m([\om])Q_t^{-1}\sum_{\tau\in E_A^n}\^m([\tau])
 \\
&= Q_t^{-3}T\^m([\om]).
\endaligned
\end{equation}
Let  
$$
l_B:l_\infty\lra\R
$$ 
be a Banach limit, see formula \eqref{220200106} and properties (a)--(g) following it. It is then not difficult to check that the formula 
$$
\^\mu(B)=l_B\((\^m(\sg^{-n}(B)))_{n\ge 0}\)
$$
defines a shift--invariant, finitely additive probability measure on
Borel sets of $E_A^\infty$ satisfying
\begin{equation}\lab{2.2.4}
{\frac {Q_t^{-3}T}{ p}}\^m(B)
\le \^\mu(B)
\le Q_t^3\^m(B),
\end{equation}
for every Borel set $B \subset E_A^\infty$. Since $\^m$ is a countably additive measure, we easily deduce that $\^\mu$ is also countably additive. It then immediately follows from \eqref{2.2.4} that $\tilde\mu$ is a Gibbs state for $t$.

\sp Let us prove the ergodicity of $\^\mu$. For each $\tau\in E_A^*$, as in \eqref{2.2.3} we getd:
\begin{equation}\lab{2.2.5}
 \sum_{i=n}^{n+p}\^\mu(\sg^{-i}([\tau])\cap [\om])
\ge \^\mu([\om\alpha(\omega,\tau)\tau]) 
\ge Q_t^{-3}T\^\mu([\tau])\^\mu([\om]).
\end{equation}

Take now an arbitrary Borel set $B\sbt E_A^\infty$ and  fix $\e>0$. Since the
nested family of sets $\{[\tau]:\tau\in E_A^*\}$ generates the Borel
$\sg$--algebra on $E_A^\infty$, for every $n\ge 0$ and every $\om\in E_A^n$
we can find a subfamily $Z$ of $E_A^*$ consisting of mutually incomparable
words such that 
$$
B\sbt \bu_{\tau\in Z}[\tau]
$$ 
and for all $n\leq i \leq n+p$,
$$
\sum_{\tau\in Z}\^\mu(\sg^{-i}([\tau])\cap [\om]) \le \^\mu\([\om]\cap
\sg^{-i}(B)\) +\e/p.
$$
Then, using \eqref{2.2.5} we get
\begin{equation}
\aligned
\e+\sum_{i=n}^{n+p}\^\mu\([\om]\cap\sg^{-I}(B)\)
&\ge \sum_{i=n}^{n+p}\sum_{\tau\in Z} \^\mu\([\om]\cap\sg^{-i}([\tau])\) 
\ge \sum_{\tau\in Z}Q_t^{-3}{T}\^\mu([\tau])\^\mu([\om]) \\
&\ge Q_t^{-3}{T}\^\mu(B)\^\mu([\om]).
\endaligned
\end{equation}
Hence, letting $\e\downto 0$, we get
$$
\sum_{i=n}^{n+p}\^\mu\([\om]\cap\sg^{-(i)}(B)\) \ge Q_t^{-3}e^T\^\mu(B)\^\mu([\om]).
$$
From this inequality we get
$$
\aligned
\sum_{i=n}^{n+p}\^\mu\(\sg^{-i}(E_A^\N\sms B)\cap [\om]\)
& = \sum_{i=n}^{n+p}\^\mu\([\om]\sms \sg^{-i}(B)\cap [\om]\) 
= \sum_{i=n}^{n+p}\^\mu([\om])-\^\mu\(\sg^{-i}(B)\cap [\om]\)\\
&\le (p-Q_t^{-3}T\^\mu(B))\^\mu([\om]).
\endaligned
$$
Thus, for every Borel set $B\sbt E_A^\infty$, for every $n\ge 0$, and
for every $\om\in E_A^n$ we have
\begin{equation}\lab{2.2.6}
\sum_{i=n}^{n+p}\^\mu(\sg^{-i}(B)\cap[\om]\)
\le (p-Q_t^{-3}T(1-\^\mu(B)))\^\mu([\om]).
\end{equation}
In order to conclude the proof of the ergodicity of $\sg$, suppose
that 
$$
\sg^{-1}(B)=B 
\  \  \  {\rm with} \  \  \
0<\^\mu(B)<1.
$$
Put 
$$
\g:=1-Q_g^{-3}e^T(1-\^\mu(B))/p.
$$
Note that \eqref{2.2.4} implies that $Q_g^{-3}T p^{-1} \leq 1$, hence $0<\g<1$.
In view of (\ref{2.2.6}), for every $\om\in E_A^*$ we get
$$\^\mu(B\cap[\om])= \^\mu(\sg^{-i}(B)\cap[\om]\) \le \g\^\mu([\om]).$$
Take now $\eta>1$ so small that
$\g\eta<1$ and choose a subfamily $R$ of $E_A^*$ consisting of mutually
incomparable words such that
$$
B\sbt \bu_{\om\in R}[\om]
\ \  \mbox{ and } \  \
\^\mu\left(\bu_{ \om \in R}[\om]\right) \le \eta\^\mu(B).
$$ 
Then
$$
\^\mu(B)
\le \sum_{\om\in R}\^\mu(B\cap[\om])\le \sum_{\om\in R}\g\^\mu([\om])
=\g\^\mu\lt(\bu_{\om\in R}[\om]\rt) 
\le \g\eta\^\mu(B
<\^\mu(B).
$$
This contradiction finishes the proof of the existence part.

\sp It follows from ergodicity of $\^\mu$ and Proposition~\ref{p2.2.2}(b) that any two Gibbs states are ergodic. So, the uniqueness of invariant Gibbs states follows immediately from Theorem~\ref{t1j74.1-2019-11-18} and Proposition~\ref{p2.2.2}(b) invoked again.

\sp Finally, let us prove the complete ergodicity of $\^\mu$ in the case when $A$ is
finitely primitive. Essentially, we repeat the argument just given.
Let $\La$ be a finite set of words all of length $q$
which witness the finite primitivity of $A$. Fix $r \in \N$. Let  $\om\in E_A^n$. For each $\tau\in E_A^*$, we get the following improvement
of (\ref{2.2.3}).
\begin{equation}\lab{2.2.7}
 \^\mu(\sg^{-(n+qr)}([\tau])\cap [\om])
\ge \sum_{\a\in\La^r\cap E^{qr}:A_{\om_n\a_1}=A_{\a_{qr}\tau_1}=1}
\^\mu([\om\a\tau]) 
\ge Q_t^{-3}{T^r}\^\mu([\tau])\^\mu([\om]).
\end{equation}
Take now an arbitrary Borel set $B\sbt E_A^\infty$. Fix $\e>0$. Since the
nested family of sets $\{[\tau]:\tau\in E_A^*\}$ generates the Borel
$\sg$--algebra on $E_A^\infty$, for every $n\ge 0$ and every $\om\in E_A^n$
we can find a subfamily $Z$ of $E_A^*$ consisting of mutually incomparable
words such that 
$$
B\sbt \bu_{\tau\in Z}[\tau]
$$ 
and
$$
\sum_{\tau\in Z}\^\mu\(\sg^{-(n+qr)}([\tau])\cap [\om]\) 
\le \^\mu\([\om]\cap
\sg^{-(n+qr)}(B)\) +\e.
$$
Then, using (\ref{2.2.7}) we get
$$
\e+\^\mu\([\om]\cap\sg^{-(n+qr)}(B)\)
\ge \sum_{\tau\in Z}Q_t^{-3}T^r\^\mu([\tau])\^\mu([\om]) 
\ge Q_t^{-3}T^r\^\mu(B)\^\mu([\om]).
$$
Hence, letting $\e\downto 0$, we get
$$
\^\mu\([\om]\cap\sg^{-(n+qr)}(B)\)
\ge \^Q(r)\^\mu(B)\^\mu([\om]),
$$
where $\^Q(r):=Q_t^{-3}\exp(rT)$. Note that it follows from this
last inequality that $\^Q(r) \leq 1.$  Also, from this inequality we find that
$$
\aligned
\^\mu\(\sg^{-(n+qr)}(E_A^\infty \sms B)\cap [\om]\) 
&=\^\mu\([\om]\sms \sg^{-n}(B)\cap [\om]\)  
\\
&=\^\mu([\om])-\^\mu\(\sg^{-(n+qr)}(B)\cap [\om]\)
 \\
&\le(1-\^Q(r)\,\^\mu(B))\^\mu([\om]).
\endaligned
$$ 
Thus, for every Borel set $B\sbt E_A^\infty$, for every $n\ge 0$, and for every $\om\in E_A^n$ we have
\begin{equation}\lab{2.2.8}
\^\mu(\sg^{-(n+qr)}(B)\cap[\om]\) \le \(1-\^Q(r)\,(1-\^\mu(B))\)\^\mu([\om]).
\end{equation}
In order to conclude the proof of the complete ergodicity of $\sg$ suppose
that 
$$
\sg^{-r}(B)=B
\  \  \ {\rm with} \  \  \
0<\^\mu(B)<1.
$$
Let
$$
{(E_A^r)}^\ast:=\bigcup_{k \in \N} {(E_A^r)}^k.$$
Put 
$$
\g:=1-\^Q(r)(1-\^\mu(B)).
$$
Note that $0<\g<1$. In view of (\ref{2.2.8}), for every $\om\in (E_A^r)^*$ we get
$$
\^\mu(B\cap[\om])= \^\mu(\sg^{-(|\om|+qr)}(B)\cap[\om]\) \le \g\^\mu([\om]).
$$
Take now $\eta>1$ so small that $\g\eta<1$ and choose a subfamily $R$ of $(E_A^r)^*$ consisting of mutually
incomparable words such that 
$$
B\sbt \bu_{\om\in R}[\om]
$$ 
and
$$
\^\mu\lt(\bu_{\om\in R}[\om]:\rt) \le \eta\^\mu(B).
$$ 
Then
$$
\^\mu(B)
\le \sum_{\om\in R}\^\mu(B\cap[\om])\le \sum_{\om\in R}\g\^\mu([\om]) 
=\g\^\mu\lt(\bu_{\om\in R}[\om]\rt) 
\le \g\eta\^\mu(B)
<\^\mu(B).
$$ 
This contradiction finishes the proof of the complete ergodicity of  $\^\mu$.
The proof of Theorem~\ref{t2.2.4} is complete.
\epf

\bthm\lab{t2.3.6}
If $S$ is a finitely irreducible graph directed Markov system with an incidence matrix $A$ and $t\in F(S)$, then the conjugate operator 
$$
e^{-\P(t)}\L_t^*:C_b(E_A^\infty)\lra C_b(E_A^\infty)
$$ 
fixes at most one Borel probability measure on $E_A^\infty$.
\ethm

\bpf Suppose that $\^m$ and $\^m_1$ are such two fixed points.
In view of Proposition~\ref{p2.2.2}(b) and Theorem~\ref{t2.3.3}, the
measures $\^m$ and $\^m_1$ are equivalent. Consider the Radon--Nikodym derivative 
$$
\rho:=\frac{d\^m_1}{d\^m}.
$$
Temporarily fix $\om\in E_A^*$, say $\om\in E_A^n$. Denote
$$
Z(\sg(\om)):=\big\{\tau\in E_A^\infty:A_{\om_n\tau_1}=1\big\}.
$$
Note that if $n\ge 2$, which we assume from now on, then
$$
Z(\sg(\om))=Z(\sg(\om))
$$
It then follows from (\ref{1201200106}) and Theorem~\ref{t2.3.3} that 
$$
\aligned
\^m([\om])
&=\int_{Z(\om)}
\big|\phi_{\om}'\big|^t\exp(-\P(t)n)\, d\^m 
\\
&=\int_{Z(\om)}\big|\phi_{\sg(\om)}'(\pi(\tau))\big|^t\exp\(-\P(t)(n-1)\)
\big|\phi_{\om_1}'(\pi(\sg(\om)\tau))\big|^t\exp(-\P(t)) \,d\^m(\tau)
 \\
&=\int_{Z(\sg(\om))}\big|\phi_{\sg(\om)}'(\pi(\tau))\big|^t\exp\(-\P(t)(n-1)\)
\big|\phi_{\om_1}'(\pi(\sg(\om)\tau))\big|^t\exp(-\P(t)) \,d\^m(\tau).
\endaligned
$$
Since, by the same token, also
$$
\^m([\sg(\om)])
=\int_{Z(\sg(\om))}\big|\phi_{\sg(\om)}'(\pi(\tau))\big|^t\exp\(-\P(t)(n-1)\) \,d\^m(\tau),
$$
we thus conclude that
$$
\bal
\inf\Big\{\big|\phi_{\om_1}'(\pi(\sg(\om)\tau))\big|^t&\exp(-\P(t)):\tau\in Z(\sg(\om))\Big\}\^m([\sg(\om)])\le \\
&\le\^m([\om])\le  \\
&\le\sup\Big\{\big|\phi_{\om_1}'(\pi(\sg(\om)\tau))\big|^t\exp(-\P(t)):\tau\in Z(\sg(\om))\Big\}\^m([\sg(\om)]).
\eal
$$
We therefore conclude that for every $\om\in E_A^\infty$
\begin{equation}\lab{2.3.5}
\lim_{n\to\infty} \frac{\^m\([\om|_n]\)}{\^m\([\sg(\om)|_{n-1}]\)}
=\big|\phi_{\om_1}'(\pi(\sg(\om)))\big|^t\exp(-\P(t)).
\end{equation}
Of course,  the same formula is true with $\^m$ replaced by $\^m_1$. Since 
the measures $\^m$ and $\^m_1$ are equivalent, by applying Theorem~\ref{t2.3.3} and Theorem~\ref{t2.2.4} along with Proposition~\ref{p2.2.2} (ergodicity),
we deduce that there exists a measurable $\sg$--invariant set $\Ga\sbt E_A^\infty$ such that the Radon--Nikodym derivative $\rho(\om)$ is defined for every $\om\in\Ga$. So, if $\om\in\Ga$ then the Radon--Nikodym derivatives $\rho(\om)$ and $\rho(\sg(\om))$ are both defined. Then from (\ref{2.3.5}) and its version for $\^m_1$ we obtain:
$$
\aligned
\rho(\om)
&=\lim_{n\to\infty}\lt(\frac{\^m_1\([\om|_n]\)}{\^m\([\om|_n]\)}\rt)
\\
&=\lim_{n\to\infty}\lt(\frac{\^m_1\([\om|_n]\)}{\^m_1\([\sg(\om)|_{n-1}]\)}
\cdot\frac{\^m_1\([\sg(\om)|_{n-1}]\)}{\^m\([\sg(\om)|_{n-1}]\)}\cdot
\frac{\^m\([\sg(\om)|_{n-1}]\)}{\^m\([\om|_n]\)}\rt) 
\\
&=\big|\phi_{\om_1}'(\pi(\sg(\om)))\big|^t\exp(-\P(t))\rho(\sg(\om))
\big|\phi_{\om_1}'(\pi(\sg(\om)))\big|^{-t}\exp(\P(t))
\\
&=\rho(\sg(\om)).
\endaligned
$$
But, since, according to Theorem~\ref{t2.2.4} again, the shift map $\sg:E_A^\infty\lra E_A^\infty$ is ergodic with respect to a shift--invariant measure equivalent with $\^m$, we conclude that the Radon--Nikodym derivative $\rho$ is $\^m$-almost everywhere constant. Since $\^m_1$ and $\^m$ are
both probability measures, we thus have that $\^m_1=\^m$. The proof of Theorem~\ref{t2.3.6} is complete.
\epf

\sp As an immediate consequence of Theorem~\ref{t220131009}, Theorem~\ref{t2.3.6}, Theorem~\ref{t2.3.3}, and Theorem~\ref{t2.2.4}, we get the following result summarizing what we did about the thermodynamic formalism.

\bcor\lab{c2.7.5}
If $S$ is a finitely irreducible graph directed Markov system with an incidence matrix $A$ and $t\in F(S)$, then

\begin{enumerate}
\item There exists a unique eigenmeasure $\^m_t$ of the conjugate
Perron--Frobenius operator \index{(N)}{Perron--Frobenius operator} $\L_t^*:C_b^*(E_A^\infty)\lra C_b^*(E_A^\infty)$ and the corresponding eigenvalue is equal to $e^{\P^(t)}$.

\,

\item The eigenmeasure $\^m_t$ is a Gibbs state \index{(N)}{Gibbs state} for $t$.

\,

\item There exists a unique shift--invariant Gibbs state $\^\mu_t$ for $t$.

\,

\item The Gibbs state $\tilde\mu_t$ is ergodic, equivalent to $\tilde m_t$ and $\log(d\tilde\mu_t/d\tilde m_t)$ is uniformly bounded.


\,

\item If the incidence matrix $A$ is finitely primitive, then the Gibbs state $\^\mu_t$ is completely ergodic. 
\end{enumerate}
\ecor

\sp\fr For $t=h$ this corollary results gives the following.

\sp\bprop\label{p320130508}
If $S$ is a finitely irreducible confomal GDMS, then $S$ is regular if and
only if there exists a Borel probability measure $\nu$ on $E_A^\infty$ such that 
\beq\label{120131005}
\L_h^*\nu=\nu. 
\eeq
Furthermore, if such a measure exists, then it is unique.

In addition, if $\L_t^*\nu=\nu$, for some $t\ge 0$ and some Borel probability measure $\nu$, then $t=h$ and $\nu=\tilde m_h$. In particular, formula \eqref{2j93} holds.
\eprop

\sp\section{Finer Geometrical Properties of CGDMSs}\label{FGPoCGDMSs}
In this section we continue our investigations of conformal graph directed Markov systems. We deal with pressure function, conformal measures, strongly regular systems, and we prove finiteness of the Hausdorff measure along with positivity of the packing measure. 

\sp Throughout the whole section $S$ is a finitely irreducible graph directed Markov system with an incidence matrix $A$ and $t\in F(S)$.

With the measures $\tilde m_t$ and $\tilde \mu_t$ coming from Corollary~\ref{c2.7.5}, we denote
$$
m_t:=\tilde m_t\circ\pi^{-1}.
$$
and
$$
\mu_t:=\tilde\mu_t\circ\pi^{-1}.
$$
In the Part III of the book, we will need the following three results, out of which the last two concern strongly regular confomal systems.

\sp\bprop\label{p120131009}
If $S_E=\{\phi_e:e\in E\}$ is a finitely irreducible confomal GDMS and $E'$ is a proper subset of $E$ such that the system $S_{E'}=\{\phi_e:e\in E'\}$ is also finitely irreducible, then
$$
\P_{E'}(t)<\P_E(t)
$$
for all $t\in F(S)$.
\eprop

\bpf
Fix $b\in E\sms E'$. Suppose for a contradiction that $\P_{E'}(t)=\P_E(t)$ for some $t\in F(S)$. Let $\tilde m_t'$ be the measure produced in Theorem~\ref{t220131009} for the system $S_{E'}$. It then follows from Remark~\ref{r120200107}, formulas \eqref{420130508}, \eqref{520130508},  \eqref{120130508}, and the third assertion of Corollary~\ref{c2.7.5} (4), that 
$$
\tilde \mu_t'([\om])\le C\tilde \mu_t([\om])
$$
for some constant $C>0$ and all $\om\in {E'_{A'}}^*$, where $A'=A|_{E'\times E'}$. Furthermore, this inequality holds for all $\om\in E_A^*$ since, if $\om\notin {E'_{A'}}^*$, then $\tilde \mu_t'([\om])=0$. Thus, the measure $\tilde \mu_t'$, considered as a Borel probability measure on $E_A^{\infty}$, is absolutely continuous with respect to the measure $\tilde \mu_t$. Since, by the first assertion of Corollary~\ref{c2.7.5} (4), both measures $\tilde \mu_t'$ and $\tilde \mu_t$ are ergodic, we conclude (see Theorem~\ref{t1j74.1-2019-11-18}) that these these two measures are equal. This however is a contradiction as $\tilde \mu_t'([b])=0$ while $\tilde \mu_t([b])>0$. We are done.
\epf

\sp\bcor\label{c420131009}
If $S_E=\{\phi_e:e\in E\}$ is a strongly regular finitely irreducible confomal GDMS and $E'$ is a proper subset of $E$ such that the system $S_{E'}=\{\phi_e:e\in E'\}$ is also finitely irreducible, then
$$
h_{S_{E'}}<h_{S_E}.
$$
\ecor

\bpf
Since the system $S_E$ is strongly regular, there exists $u<h_E$ such that $0<\P_E(u)<+\infty$. Hence,
$$
\P_{E'}(u)\le \P_E(u)<+\infty,
$$
and therefore the function $\P_{E'}$ restricted to $[u,h_E]$ is continuous. Since also, by Proposition~\ref{p120131009}, $\P_{E'}(h_{E})<\P_E(h_E)=0$, we therefore see that there exists $t\in [u,h_E)$ such that $\P_{E'}(t)<0$. Hence, it follows from Theorem~\ref{t1j97} that $h_{E'}\le t<h_E$. The proof is complete. 
\epf

\sp\bthm\label{GDMS_Finite_Entropy} 
If $S$ is a finitely irreducible strongly regular confomal GDMS, then
the metric entropy $\h_{\tilde\mu_{h_S}}(\sg)$ of the dynamical system
$(\sg:E_A^\infty\lra E_A^\infty)$ with respect to the $\sg$--invariant measure
$\tilde\mu_{h_S}$ is finite.
\ethm
\bpf
Let $\a$ be the partition of $E_A^\cN$ into initial cylinders of
length one, i.e. 
$$
\a=\{[e]\}_{e\in E}.
$$
Since the system $S$ is strongly regular, by virtue of
Proposition~\ref{p2j85} there exists $\eta>0$ such that
$Z_1(h-\eta)<\infty$. This means that 
$$
\sum_{e\in E}\|\f_e'\|_\infty^{h_S-\eta}<+\infty.
$$ 
Since $\|\f_e'\|_\infty^{-\eta}\ge -h_S\log\|\f_e'\|_\infty$ for all but, 
perhaps, fintely many $e\in E$, the series 
$$
\sum_{e\in E}-h_S\log\(\|\f_e'\|_\infty\)\|\f_e'\|_\infty^{h_S}
$$ 
converges too. Hence, by \eqref{2j93},
$$
\H_{\^\mu_{h_S}}(\a)=\sum_{i\in I}-\log(\^\mu_{h_S}([e])\^\mu_{h_S}([e])<+\infty
$$
Since the partition $\a$ is a metric (even topological) generator of
the dynamical system $\sg:E_A^\infty\lra E_A^\infty$, we have
$$
\h_{\^\mu_{h_S}}(\sg)=\H_{\^\mu_{h_S}}(\a)<+\infty,
$$
and the proof is complete.
\epf 

\sp\fr Since we always have $\P(h_S)\le 0$, formulas \eqref{520130508} and
\eqref{4.1.9}, immediately give the following.

\sp\bprop\label{p420130508}
If $S$ is a finitely irreducible confomal GDMS, then the Hausdorff measure $\H_{h_S}(J_S)<+\infty$. \index{(N)}{Hausdorff measure}
\eprop

\fr We need however a slightly stronger statement.

\sp\bthm\lab{t4.4.1} 
If $S$ is finitely irreducible regular confomal GDMS, then
the Hausdorff measure \index{(N)}{Hausdorff measure} $\H_{h_S}$ restricted to the limt set $J_S$ 
is absolutely continuous with respect to the conformal measure $m_h$
and 
$$
\|d\H_h/dm_h\|_\infty<+\infty.  
$$
In particular we get again that the Hausdorff measure $\H_{h_S}(J_S)$ is finite. 
\ethm
\bpf
Let $A$ be an arbitrary closed subset of $J_S$. For every integer $n\ge 1$ 
put
$$
A_n:=\{\om\in E_A^n:\f_\om(J_S)\cap A\ne\es\}.
$$
Then the sequence of sets,
$$
\lt(\bu_{\om\in A_n}\f_\om(X_{t(\om)})\rt)_{n=1}^\infty
$$ 
is descending and  
$$
\bi_{n\ge 1}\lt(\bu_{\om\in A_n}\f_\om(X_{t(\om)})\rt)=A.
$$
Therefore, using (\ref{4.1.9}) in Corollary~\ref{c23013_03_12}, and
\eqref{2j93}, 
$$
\aligned
\H_h(A) 
&\le \varliminf_{n\to\infty}\sum_{\om\in A_n}\(\diam(\f_\om(X_{t(\om)}))\)^h
\le D^h\varliminf_{n\to\infty}\sum_{\om\in A_n}\|\f_\om'\|^h \\
&\le D^hM_h^{-1}\varliminf_{n\to\infty}\sum_{\om\in A_n}\tilde m_h([\om])\\
&\le D^hM_h^{-1}\varliminf_{n\to\infty}
          \lt(m\(\bu_{\om\in A_n}\f_\om( X_{t(\om)})\)\rt)\\
&=D^hM_h^{-1}m_h(A).
\endaligned 
$$
Since $J_S$ is a metric separable space, the measure $m_h$ is regular, and 
therefore the  inequality 
$$
\H_t(A)\le D^hM_h^{-1}m_h(A)
$$
extends to all Borel subsets of $J_S$. The proof is finished. 
\epf

\sp A dual statement holds for packing measures. It however requires an additional mild hypothesis common in the theory of conformal graph directed Markov systems. Since this hypothesis has many more significant consequences, we formulta it in the next section, where the first result is about packing measures and other results then follow. 

\sp\section{The Strong Open Set Condition}\label{SOSC}
As it was indicated at the end of the previous section, in the present one we formulate the Strong Open Set Condition and derive several of its remarkable consequences. We want to emphasize already now that this is a mild natural condition and it is satisfied for ``most'' GDMSs. Here it is:

\sp\bdfn
A conformal GDMS $S=\{\phi_e\}_{e\in E}$ (satisfying the Open Set Condition) is said to satisfy the Strong Open Set Condition (SOSC) if
$$
J_S\cap \Int X\ne\es.
$$
\edfn 

\sp\bprop\label{p420130508B}
If $S$ is a finitely irreducible regular confomal GDMS satisfying the Strong
Open Set Condition (SOSC), then the packing measure \index{(N)}{packing measure}$\Pi_h(J_S)>0$.
\eprop
\bpf
By our (SOSC) hypothesis there exists $\tau\in E_A^\infty$ such that
$x:=\pi(\tau)\in\Int(X)$. There then exists $R>0$ such that
$B(x,R)\sbt\Int(X)$. Hence there exists $q\ge 1$ so large that
\beq\label{620130510}
\pi([\tau|_q])\sbt B(x,R/2).
\eeq
Since $\^\mu_h([\tau|_q])>0$ and the $\sg$--invariant mesure $\^\mu_h$ is ergodic, it follows from Theorem~\ref{t1_mu_2014_11-18} that $\^\mu_h(E_\tau^\infty)=1$, where
$$
E_\tau^\infty:=\big\{\om\in E_A^\infty: \sg^n(\om)\in [\tau|_q] \text{ for infinitely many } n'\text{s}\big\}.
$$
So, fix $\om\in E_\tau^\infty$ and $n\ge 1$ such that $\sg^n(\om)\in [\tau|_q]$. Then $$
B(\pi(\sg^n(\om)),R/2)\sbt B(x,R)\sbt\Int(X).
$$
It also follows from Lemma~\ref{l42013_03_12} that
$$
\phi_{\om|_n}(B(\pi(\sg^n(\om)),R/2)\)\spt B\(\pi(\om), (2K)^{-1}R\|\phi_{\om|_n}'\|_\infty\).
$$
Therefore, using also the last assertion of Proposition~\ref{p320130508}, we get
$$
\begin{aligned}
m_h\(B\(\pi(\om),(2K)^{-1}R\|\phi_{\om|_n}'\|_\infty\)\) 
&\le m_h\(\phi_{\om|_n}(B(\pi(\sg^n(\om)),R/2))\) 
\le m_h\(\phi_{\om|_n}(\Int(X))\) \\
&=\^m_h\([\om|_n]\) 
\le \|\phi_{\om|_n}'\|_\infty^h \\
&=(2K/R)^h\((2K)^{-1}R\|\phi_{\om|_n}'\|_\infty\)^h.
\end{aligned}
$$
Hence, denoting by $\N(\om)$ the infinite set of integers $n\ge 1$ such that $\sg^n(\om)\in [\tau|_q]$
$$
\varliminf_{r\to 0}\frac{m_h\(B\(\pi(\om),r)\)}{r^h}
\le \varliminf_{{n\to\infty}\atop {n\in\N(\om)}}\frac{m_h\(B\(\pi(\om),(2K)^{-1}R\|\phi_{\om|_n}'\|_\infty\)\)}
    {\((2K)^{-1}R\|\phi_{\om|_n}'\|_\infty\)^h}
\le (2K/R)^h
>0.
$$
So, applying Frostman's Converse Theorem for Packing Measures, i.e. Theorem~\ref{tncp12.2.} (1), we therefore get 
$$
\Pi_h\(J_S\)\ge \Pi_h\(\pi(E_\tau^\infty)\)>0.
$$
The proof is complete.
\epf

\sp A straightforward observation is that if 
$$
\phi_\om\(X_{t(\om)}\)\sbt \Int\(X_{i(\om)}\)
$$
for some $\om\in E_A^*$, then (SOSC) holds. Before we provide further significant consequences of (SOSC), we prove several technical ones. First, it directly follows from (OSC) that
\beq\label{1_2018_02_16}
\phi_\om\(X_{t(\om)}\)\cap \phi_\tau\(X_{t(\tau)}\)
=\phi_\om\(\bd X_{t(\om)}\)\cap \phi_\tau\(\bd X_{t(\tau)}\)
\eeq
for all incomparable words $\om,\tau\in E_A^*$. Second:

\blem\label{l2sc1}
If $S=\{\phi_e\}_{e\in E}$ is a conformal GDMS, then
$$
\sg(\pi^{-1}(\bd X))\sbt \pi^{-1}(\bd X).
$$
\elem

{\sl Proof.} 
Since $\phi_e\(X_{t(e)}\)\sbt\Int\(X_{i(e)}\)$ and $\bd X_v\cap\Int X_v=\es$ for every $e\in E$ and every $v\in V$, we conclude that
\beq\label{1sc1}
\phi_e^{-1}\(\bd X_{i(e)}\)\sbt \bd X_{t(e)}
\eeq
Now, if $\om\in \pi^{-1}(\bd X)$, then $\om)\in \bd X_{i(\om_1)}$. Since also $\pi(\om)=\phi_{\om_1}(\pi(\sg(\om)))$, using \eqref{1sc1}, we conclude that
$$
\pi(\sg(\om))\in \phi_{\om_1}^{-1}\(\bd X_{i(\om_1)}\)
\sbt\bd X_{t(\om_1)}.
$$
Hence, $\sg(\om)\in \pi^{-1}\(\bd X_{i(\om_1)}\)\sbt\pi^{-1}(\bd X)$, and the proof is complete.
\endpf

\blem\label{l1sc3}
If $S=\{\phi_e\}_{e\in E}$ is a conformal GDMS satisfying (SOSC), then for every $\om\in E_A^*$ we have that
$$
\sg^{|\om|}\(\pi^{-1}\(\phi_\om\(\bd X_{t(\om)}\)\)\sbt \pi^{-1}(\bd X).
$$
\elem
{\sl Proof.} Put $n=|\om|$. Let $\tau\in \sg^{|\om|}\(\pi^{-1}\(\phi_\om\(\bd X_{t(\om)}\)\)$. Then $\tau=\sg^n(\g)$ for some $\g\in\pi^{-1}\(\phi_\om\(\bd X_{t(\om)}\)\)$. So,
$$
\phi_{\g|_n}(\pi(\sg^n(\g)))=\pi(\g)\in \phi_\om\(\bd X_{t(\om)}\).
$$
It thus follows from (OSC) that $\pi(\tau)=\pi(\sg^n(\g))\in\bd X_{t(\g_n)}\sbt\bd X$. Hence, $\tau\in\pi^{-1}(\bd X)$, and the proof is complete.
\endpf

\sp Now, we are in position to prove the following first main result of this section. 

\bthm\label{t2sc3}
If $S=\{\phi_e\}_{e\in E}$ is a conformal GDMS satisfying {\rm(SOSC)} and $\mu$ is a Borel probability $\sg$--invariant ergodic measure on $E_A^\infty$ with full topological support, then
\begin{itemize}
\item[(a)] $\mu\circ\pi^{-1}(\bd X)=0$.

\,\item[(b)] $\mu\circ\pi^{-1}\(\phi_\om\(\bd X_{t(\om)}\)\)=0$ for each $\om\in E_A^*$.

\,\item[(c)] If $\om,\tau\in E_A^*$ are incomparable, then
$$
\mu\circ\pi^{-1}\(\phi_\om\(\bd X_{t(\om)}\)\cap \phi_\tau\(\bd X_{t(\tau)}\)\)=0.
$$
\end{itemize} 
\ethm
{\sl Proof.} Since $S$ satisfies (SOSC), there exists $\om\in E_A^\infty$ such that $\pi(\om)\in\Int X_{i(\om_1)}$. But then there exists an integer $n\ge 1$ such that 
$$
\phi_{\om|_n}\(X_{t(\om_n)}\)\sbt \Int X_{i(\om_1)}.
$$
Hence, using also the fact that $\supp(\mu)=E_A^\infty$, we get that
$$
\mu\circ\pi^{-1}(\Int X)
\ge \mu\circ\pi^{-1}\(\Int X_{i(\om_1)}\)
\ge \mu\circ\pi^{-1}\(\phi_{\om|_n}\(X_{t(\om_n)}\)\)
\ge \mu\circ\pi^{-1}(\om|_n)>0.
$$
Therefore, $\mu\circ\pi^{-1}(\bd X)<1$. So, it follows from ergodicity of the measure $\mu$ and Lemma~\ref{l2sc1} that $\mu\circ\pi^{-1}(\bd X)=0$. The proof of item (a) is complete. 

Prooving item (b), it follows from item (a) of Lemma~\ref{l1sc3} and $\sg$--invariance of the measure $\mu$, that
$$
0=\mu\circ\pi^{-1}(\bd X)
\ge \mu\(\sg^{|\om|}\(\pi^{-1}\(\phi_\om\(\bd X_{t(\om)}\)\)\)\)
\ge \mu\(\pi^{-1}\(\phi_\om\(\bd X_{t(\om)}\)\)\).
$$
The proof of item (b) is complete.

Now, we shall prove item (c). Since $\om$ and $\tau$ are incomparable, $\om=\g a\a$ and $\tau=\g b\b$ with some $\g\in E_A^*$, $\a,\b\in E_A^\infty$ and $a,b\in E$ with $a\ne b$. Then by (OSC),
$$
\begin{aligned}
\phi_\om\(X_{t(\om)}\)\cap \phi_\tau\(X_{t(\tau)}\)
&\sbt\phi_\g\(\phi_a\(X_{t(a)}\)\cap \phi_b\(X_{t(b)}\)\)
=\phi_\g\(\phi_a\(\bd X_{t(a)}\)\cap \phi_b\(\bd X_{t(b)}\)\)\\
&\sbt \phi_{\g a}\(\bd X_{t(a)}\).
\end{aligned}
$$
So, by item (b), 
$$
\mu\circ\pi^{-1}\(\phi_\om\(\bd X_{t(\om)}\)\cap \phi_\tau\(\bd X_{t(\tau)}\)\)
\le \mu\circ\pi^{-1}\( \phi_{\g a}\(\bd X_{t(a)}\)\)=0.
$$
Thus, the proof of item (c) is complete and we are done.
\endpf

\sp Since for every $t\in F(S)$ the invariant measure $\mu_t$ is of full topological support and since the measures $\mu_t$ and $m_t$ are equivalent, as an immediate consequence of Theorem~\ref{t2sc3}, we get the following.

\bthm\label{t1sc5}
If $S=\{\phi_e\}_{e\in E}$ is a finitely irreducible conformal GDMS satisfying {\rm (SOSC)} and $t\in F(S)$, then
\begin{itemize}
\item[(a)] $\mu_t(\bd X)=0$.

\,\item[(b)] $\mu_t\(\phi_\om\(\bd X_{t(\om)}\)\)=0$ for each $\om\in E_A^*$.

\,\item[(c)] If $\om,\tau\in E_A^*$ are incomparable, then
$$
\mu_t\(\phi_\om\(\bd X_{t(\om)}\)\cap \phi_\tau\(\bd X_{t(\tau)}\)\)=0.
$$
\end{itemize} 
Also, the same holds with $\mu_t$ replaced by $m_t$.
\ethm

Now, we are in position to prove the conformality property of the measure $m_t$. For every $e \in E$, let
$$
E_A^{+}(e):=\{ \om \in E^{\infty}_{A}: A_{e\om_1}=1\}.
$$

\bthm\label{t2sc5} 
If $S=\{\varphi_e\}_{e \in E} $ is a conformal GDMS satisfying {\rm (SOSC)} and $t\in F(S)$, then
$$ 
m_t(\varphi_\om(F))= e^{-\P(t)|\om|}\int_{F} |\varphi_\om'|^tdm_t
$$
for every $\om \in E^{*}_{A}$ and every Borel set $F \sbt \pi(E_A^+(\om_{|\om|}))$.
\ethm

\fr {\sl Proof.} Put $n:=|\om|$. Denote
$$  
Z_\om(F):=\{\gamma\in \pi^{-1}(F): \,\,  A_{\om_n\gamma_1}=1\}.
$$
Note that if $x \in F$ then $x=\pi(\gamma)$, where $\gamma\in E^{+}(\om_n)$. So, $ \gamma \in \pi^{-1}(F)$ and $A_{\om_n \gamma_1}=1$. Whence $ \gamma \in Z_\om(F)$. We have thus proved that 
\beq\label{(1sc6)}
\pi(Z_\om(F))=F.
\eeq
Now, if $\tau \in \pi^{-1}(F)\sms Z_\om(F)$, then $A_{\om_n\tau_1}=0$. Hence,
$$ 
\pi^{-1}(F)\sms Z_\om(F)\sbt \bigcup_{e\in E:A_{\om_ne}=0}[e].
$$
So, using (\ref{(1sc6)}), we  get
$$
\begin{aligned}
\pi(\pi^{-1}(F)\sms Z_\om(F)) 
& \sbt F \cap \pi\lt(\bigcup_{b\in E:A\om_nb=0} [e]\rt)\\ 
& \sbt \bigcup_{a\in E:A_{\om_n a}=1} \varphi_a\(X_{t(a)}\)\cap \bigcup_{b\in E:A_{\om_nb}=1} \varphi_b\(X_{t(b)}\)\\
&= \bigcup_{a\in E:A_{\om_n}a=1}\,\bigcup_{b\in E:A_{\om_n}b=1} \varphi_a\(X_{t(a)}\)\cap  \varphi_b\(X_{t(b)}\).
\end{aligned}
$$
Hence, it follows from Theorem~\ref{t1sc5} (c) that $m_t\(\pi(\pi^{-1}(F)\sms Z_\om(F))\)=0$. Therefore,
\beq\label{(2sc6)}
\tilde{m}_t(\pi^{-1}(F)\sms Z_\om(F))=0.
\eeq
Now, $\tau \in \pi^{-1}(\varphi_\om(F))$ if and only if $\pi(\tau) \in \varphi_\om(F)$, i.e.
 $\varphi_{\tau|n} (\pi(\sigma^n(\tau))) \in \varphi_\om(F)$. Hence,
$$
\om Z_\om(F) 
\in \pi^{-1}(\varphi_\om(F)) 
\sbt Z_\om(F) \cup \pi^{-1}\lt(\bigcup_{\tau \in E^n_A\sms \{\om\}} \varphi_\tau\(X_{t(\tau)}\)\cap \varphi_\om\(X_{t(\om)}\)\rt).
$$
It therefore  follows from Theorem~\ref{t1sc5} (c) that
$$ 
\^m_t(\om Z_\om(F))
=\tilde{m}_t(\pi^{-1}(\varphi(F))).
$$
So, remembering  that $\lambda_t= e^{\P(t)}$ and using formulas (\ref{mtdual}), (\ref{(2sc6)}), and  (\ref{(1sc6)}), we get that
$$
\begin{aligned}m_t(\varphi_\om(F))
&= \tilde{m}_t \circ \pi^{-1}(\varphi_\om(F))
 =\tilde{m}_t(\om Z_\om(F))
 =e^{-\P(t)n}{\mathcal L}^{*n}_t \tilde{m}_t\(\1_{[\om Z_\om(F)]}\)\\
 &=e^{-\P(t)}\tilde{m}_t\({ \mathcal L}^n_t \1_{[\om Z_\om(F)]}\)\\
 & = e^{-\P(t)n}\int_{E^{\mathbb N}_A} \sum_{\tau  \in E^n_A} \1_{[\om Z_\om(F)]}(\tau \gamma)| \varphi_\tau'( \pi(\gamma))|^t d\tilde{m}_t(\gamma)\\
 & =e^{-\P(t)n}\int_{Z_\om (F)}|\varphi_\om'(\pi(\gamma))|^t d\tilde{m}_t(\gamma)\\
 &  =e^{-\P(t)n}\int_{\pi^{-1}(F)}|\varphi_\om'(\pi(\gamma))|^t d\tilde{m}_t(\gamma)\\
 & =e^{-\P(t)n}\int_F|\varphi_\om'|^t dm_t.\\
 \end{aligned}
 $$
 The  proof is complete. \qed

\sp\section{Conformal Maximal Graph Directed Systems}

We already know the significance of the measure $m_h$ defined by formula \eqref{220131007}. In this section we want to give an intrinsic characterization of this measure, i. e. one that does not invoke any symbol dynamics. For this we however need to restrict our attention to a narrower class of graph directed Markov systems. In fact we will do it more generally, for all measures $m_t$, 
$t\in F(\cS)$.

\sp\bdfn\label{d2_2017_11_18}
A graph directed Markov system $S=\{\phi_e:e\in E\}$ is called maximal  \index{(N)}{maximal graph directed Markov system} if $A_{ab}=1$ whenever $i(b)=t(a)$. Equivalently, $A_{ab}=1$ if and only if $i(b)=t(a)$.
\edfn

\sp\fr In particular every iterated function system is a maximal graph directed Markov system. Now, let $S=\{\phi_e:e\in E\}$ be a conformal maximal graph directed Markov system. For every $t\in F(S)$ the (intrinsic) Perron-Frobenius operator $L_t:C(X)\to C(X)$ is defined as follows. If $v\in V$ and $x\in X_v$, then
\beq\label{320130107} 
L_t g(x)
:= \sum_{e\in E:t(e)=v}
g(\phi_e(x))|\phi_e'(x)|^t.
\eeq
The following theorem of this section was partly proved in \cite{MU1} and  \cite{MU2}.

\sp\bthm\label{t420131007} 
Let $S=\{\phi_e:e\in E\}$ be a finitely irreducible conformal maximal graph directed Markov system satisfying {\rm(SOSC)} and let $t\in F(S)$. If $m$ is a Borel probability measure on $X$, then the following conditions are equivalent.

\sp\begin{itemize}

\sp\item[(a)] $m=m_t$.

\sp\item[(b)] $L_t^*(m)=e^{\P(t)}m$.

\sp\item[(c)]  
$$
m(\phi_e(F))=e^{-\P(t)}\int_F|\phi_e'|^t\,dm
$$
for every $e\in E$ and evry Borel set $F\sbt X_{t(e)}$.

\sp\item[(d)] Item (c) holds and 
$$
m\(\phi_a(X_{t(a)})\cap \phi_b(X_{t(b)})\)=0
$$
whenever $a,b\in E$ and $a\ne b$. 

\sp\item[(e)] $L_t^*(m)=\g m$ for some $\g>0$.

\sp\item[(f)] 
$m(J_S)=1$ and
$$
m(\phi_e(F))=\g^{-1}\int_F|\phi_e'|^t\,dm
$$
for some $\g>0$, every $e\in E$ and evry Borel set $F\sbt X_{t(e)}$.
\item[(g)] Item (f) holds and 
$$
m\(\phi_a(X_{t(a)})\cap \phi_a(X_{t(a)})\)=0
$$
whenever $a,b\in E$ and $a\ne b$. 
\end{itemize}
\ethm

{\sl Proof.}
We will first establish a  close relation between the operators ${\mathcal L}_t$ and $L_t$, where $t \in F(S)$.

\sp\fr {\bf Claim 1.} {\it Let $t \in F(S)$. If $ n \geq 0$, then for every $g \in C(X)$, we have that
$$
{\mathcal L}^n_t( g \circ \pi) = L^n_t(g)\circ \pi
,$$
where, we recall, $ \pi : E^{\infty}_A\lra J_S$ is the canonical H\"older continuous projection from the symbol space onto the limit set of the GDMS $S$.}

\fr {\sl Proof.} The proof goes through direct calculation. For every $\om \in E^{\infty}_A$ we have
$$\begin{aligned}  
{\mathcal L}^n_t(g\circ \pi)(\om)
& =\sum_{\tau\in E^n_A:A_{\tau_n\om_1}=1} g \circ \pi( \tau \om)|\varphi_\tau'(\pi(\om))|^t\\
 &=\sum_{\tau\in E^n_A:t(\tau)=i(\om)} g (\varphi_\tau( \pi(\om)))|\varphi_\tau'(\pi(\om))|^t\\
&  = L^n_t(g)(\pi(\om)).
\end{aligned}
$$
The proof of Claim 1 is complete.
\qed

\sp The standard approximation argument, based on Theorem~2.4.6 in \cite{MU2}, which involves H\"older continuous functions only, gives the following.

\sp\fr {\bf Fact~9.8.}\label{f1mg2} {\it For every $k \in C_b(E^{\infty}_A)$ the sequence  $(e^{-\P(t)n}{\mathcal L}^n_t(k))^{\infty}_{n=0}$ converges uniformly to
$\tilde{m}_t(k)\tilde{\rho}_t$, where $ \tilde{\rho}_t=\frac{d \tilde{\mu}_t}{d\tilde{m}_t}$.}

\sp\fr {\bf Claim 2.} {\it  There exists a unique continuous function $\rho_t: J_S \lra (0, +\infty)$ such that
$$ 
\tilde{\rho}_t ={\rho}_t\circ \pi,
$$
and for every $g\in C_b(J_S)$, the sequence $ (e^{-\P(t)n} L^n_tg)_{n=0}^\infty$ converges  uniformly to the function $m_t(g)\rho_t$.}

\bpf  Applying Fact~\ref{f1mg2} to the function $k=\1$, along with Claim~1, applied to the function $g=\1$, we get for every $x \in J_S$ and every $\om \in  \pi^{-1}(x)$, that
$$
\begin{aligned} \tilde{\rho}_t(\om)
& =\lim_{n \to \infty}\(e^{-\P(t)n}{\mathcal L}^n_t(\1 \circ \Pi)\)(\om)
  =\lim_{n \to \infty}\( e^{-\P(t)n}{L}^n_t\1\)(\pi(\om))\\
 & =\lim_{n \to \infty} (e^{-\P(t)n}{ L}^n_t\1 )(x)
 \end{aligned}
 $$
and the convergence is uniform  with respect to $ \om \in E^\infty_A$ and $x \in J_S$. In particular $\tilde{\rho}_t$ is constant on $ \pi^{-1}(x)$. Since also $\tilde{\rho}_t:E^\infty_A \lra (0, +\infty)$ is  continuous, we thus conclude that there exists a continuous function $\rho_t: J_S \lra (0, 1)$ such that 
$$
\tilde{\rho}_t=\rho_t\circ \pi.
$$
Now, the uniqueness of $\rho_t$ follows from surjectivity of the projection $\pi: E^\infty_A \lra J_S$. Applying again Claim~1 with the function $g$ and Fact~\ref{f1mg2} with the function $k = g \circ \Pi$, we get that for every $x \in J_S$ and every $x \in \pi^{-1}(\om)$, that
$$\begin{aligned}\lim_{n \to \infty}(e^{-\P(t)n}L^n_t\1)(x)
&= \lim_{n \to \infty}(e^{-\P(t)n}L^n_tg)(\pi(\om))
= \lim_{n \to \infty}( e^{-\P(t)n}{\mathcal L}^n_t(g\circ \pi))(\om)\\
& = \tilde{m}_t ( g \circ \pi)\tilde{\rho}_t(\om)
={m}_t ( g ){\rho}_t( \pi(\om))\\
& =m_t(g)\rho_t(x)
\end{aligned}
$$
and the convergence is uniform with respect to $ \om\in E^\infty_A$ and $x \in J_S$. The proof of Claim~2 is complete.
\qed

\sp Now we shall prove that $\rm{(a)}\Rightarrow \rm{(b)}$. Indeed, by 
Claim~1, we have for every $g \in C(X)$ that
$$
\begin{aligned}L^{*}_t m_t(g)
& =m_t(L_tg)
=\tilde{m}_t\circ \pi^{-1}(L_t g)
=\tilde{m}_t((L_t g)\circ \pi) \\
&=\tilde{m}_t({\mathcal L}_t (g\circ \pi))
 = {\mathcal L}^{*}_t\tilde{m}_t(g \circ \pi) \\
 &= e^{\P(t)}\tilde{m}_t ( g\circ \pi) 
 =e^{\P(t)}\tilde{m}_t\circ \pi^{-1}(g) \\
 &=e^{\P(t)}{m}_t (g).
\end{aligned}
$$
Therefore, $L^{*}m_t =m_t$ and the implication $\rm{(a)}\Rightarrow \rm{(b)}$ is established. 

\sp Of course $\rm{(b)}\Rightarrow \rm{(e)}$. Now, we  shall prove that $\rm{(e)}\Rightarrow \rm{(a)}$. Because of Claim 2, we get for every $g \in C(X)$ that
$$
\begin{aligned} 
\lim_{n \to \infty} e^{-\P(t)n}L^{*n}_t m(g)
& = \lim_{n \to \infty}m(e^{-\P(t)n}L^{n}_t (g))\\
& =\lim_{n \to \infty}\int_{J_S}e^{-\P(t)n}L^{n}_t g(x)dm(x)\\
&=m_t(g)m(\rho_t).
\end{aligned}
$$
Invoking (e), we then obtain
$$ 
\lim_{n \to \infty}\(( \gamma e^{-\P(t)}\)^n m(g)\)= m_t(g)m(\rho_t).
$$
Taking $g=\1$ and noting that $m(\rho_t)>0$ (as $\rho_t$ is everywhere positive on $E^{\infty}_A$), we thus obtain that $\gamma=e^{\P(t)}$ and that
\beq\label{(1mg2)}
m(g)=m_t(g)m(\rho_t).
\eeq
Taking again $g=\1$, we  thus get $m(\rho_t)=\1$. Then (\ref{(1mg2)}) becomes $m(g)=m_t(g)$. This means that $m=m_t$ and the proof of the implication  $\rm{(e)}\Rightarrow \rm{(a)}$
is complete.

\sp Now we shall prove that $\rm{(a)}\Rightarrow \rm{(d)}$. The second part  of (d) holds because of (a) and Theorem~\ref{t1sc5}. Because of Theorem~\ref{t2sc5},  in order to prove  the first part of (d), {\it i.e.} item (c), it suffices to show  that
\beq\label{(1mg14)}
F\cap J_S\sbt \pi(E^+_A(e)).
\eeq
and
\beq\label{(1mg14B)}
m_t\((J_S\cap\varphi_e(F))\sms \varphi_e(F\cap J_S)\)=0.
\eeq
Indeed, proving \eqref{(1mg14)}, if $x \in F\cap J_S$, then there exists $\om\in E_A^\infty$ such that $x=\Pi(\om)$ and $i(\om_1)=t(e)$. But since our system $S$ is maximal, this means that $A_{e\om_1}=1$. Hence, $ \om \in E^+_A(e)$,. Thus, $x =\pi(\om) \in \pi(E^+_A(e))$. Thus (\ref{(1mg14)}) holds.

Proving \eqref{(1mg14B)}, let 
$$
x\in (J_S\cap\varphi_e(F))\sms \varphi_e(F\cap J_S).
$$
Thren $x=\varphi_e(z)$ for some $z\in F\sms J_S$ and $x=\pi(\om)$ for some $\om\in E_A^\infty$. So, $x=\varphi_{\om_1}\(\pi(\sg(\om))\)$. Since the map $\varphi_e$ is 1--to--1 and $\pi(\sg(\om))\in J_S$, we thus infer that $\om_1\ne e$. Since also $\in \varphi_e\(X_{t(e)}\)\cap\varphi_{\om_1}\(X_{t(\om_1}\)$, we thus conclude that
$$
(J_S\cap\varphi_e(F))\sms \varphi_e(F\cap J_S)
\sbt \bu_{a\ne b}\phi_a(X_{t(a)})\cap \phi_b(X_{t(b)}).
$$
In conjunction with Theorem~\ref{t1sc5}, this finishes the proof of \eqref{(1mg14B)}, and (d) is proved.

\sp Obviously $\rm{(g)}\Rightarrow \rm{(f)}$. 

\sp We will show now that $\rm{(f)}\Rightarrow \rm{(c)}$. This means that we are to show that $\gamma=e^{\P(t)}$. We have for  every $n \geq 1$  that
$$
\begin{aligned} 
1=m(J_S)
&\leq m\left(\bigcup_{\om \in E^n_A}\varphi_\om(X_{t(\om)})\right)
\leq \sum_{\om \in E^n_A}m(\varphi_\om(X_{t(\om)})) \\
&= \sum_{\om \in E^n_A} \gamma^{-n}\int_{X_t(\om)} |\varphi_\om'|^t dm\\
& \leq \gamma^{-n}\sum_{\om \in E^n_A} \|\varphi_\om'\|^t m (X_{t(\om)}) \\
      & \leq \gamma^{-n}\sum_{\om \in E^n_A}\|\varphi_\om'\|^t\\
      & = \gamma^{-n} Z_n(t).
\end{aligned}
$$
Hence $\gamma^n\leq  Z_n(t)$, and therefore
\beq\label{(1mg3)}
\log \gamma 
= \lim_{n \to \infty}\frac{1}{n} \log \gamma^n 
\leq \liminf_{n \to\infty} \frac{1}{n} \log Z_n(t) 
=\P(t).
\eeq
On the other hand, because we have already proved that $\rm{(a)}\Rightarrow \rm{(d)}$ (so (c) holds with $m$ replaced by $m_t$), we get for every $\om \in E^*_A$ that
\beq\label{(2mg3)}
m_t\(\varphi_\om(\Int X_{t(\om)})\)
=e^{-\P(t)|\om|} \int_{\Int X_{t(\om)}} |\varphi_\om'|^t dm_t 
\leq e^{-\P(t)|\om|} \|\varphi_\om'\|^t.
\eeq
Since the system $S$ is irreducible and satisfies (SOSC), it immediately  follows from (f) that $m(\Int X_v)>0$ for all $v \in V$.  Therefore
\beq\label{(3mg3)}
 M:= \min \{m(\Int X_v):v \in V \}>0.
\eeq
Therefore, we get for every $\om \in E^{*}_A$ that
$$ 
m(\varphi_\om(\Int X_{t(\om)})) 
= {\gamma}^{-|\om|}\int_{\Int X_{t(\om)}} |\varphi_\om'|^t dm
\geq  \gamma^{-|\om|}K^{-t} \|\varphi_\om'\|^t m(\Int X_{t(\om)})
 \geq M  K^{-t} \gamma^{-|\om|} \|\varphi_\om'\|^t.
$$
Combining this with (\ref{(2mg3)}), we get
\beq\label{(1mg4)}
m_t\(\varphi_\om(\Int X_{t(\om)})\) 
\leq M^{-1}K^{t} \gamma^{|\om|}e^{-\P(t)|\om|}m(\varphi_\om\(\Int X_{t(\om)})\),
\eeq
where the equality sign ``='' follows from Theorem~\ref{t1sc5} and since $m_t(J_S)=1$ while the last inequality sign ``$\le $'' follows from (OSC). Hence, we get for every integer $n \geq 1$ that
$$
1
=\sum_{\om \in E^n_A}m_t(\varphi_\om\(\Int X_{t(\om)})\)  
\leq M^{-1}K^{t} \left(\frac{\gamma}{e^{\P(t)}}\right)^n \sum_{\om \in E^n_A}m(\varphi_\om\(\Int X_{t(\om)})\)
\leq M^{-1}  K^{t} \left(\frac{\gamma}{e^{P(t)}}\right)^n.
$$
Letting $ n \to \infty$, we thus conclude that $\gamma \geq e^{\P(t)}$. Along with (\ref{(1mg3)}) this finishes the proof of equality  $ \gamma=e^{\P(t)}$, thus establishing (c).

\sp So, in order to conclude the proof of our theorem, it suffices to show that $\rm{(c)}\Rightarrow \rm{(a)}$. We will  do it now. We start with the following.

 \

 \fr {\bf Claim 4.} {\it If (c) holds, then $m_t$  is absolutely continuous with respect to $m$.} \\

\bpf Let $ G \sbt X$ be an open set. Thus for every $x \in J_S\cap G$ there exists at least one $\om \in E^{*}_A $  such that $ x \in \varphi_\om(X_{t(\om)})\sbt G$. Therefore, there exists $\mathfrak F \sbt E^{*}_A$, a family of  mutually incomparable finite words such that  $$J_S\cap G \sbt \bigcup_{\om \in \mathfrak F}  \varphi_\om( X_{t(\om)})\sbt G.$$
Hence, using also (\ref{(1mg4)}) with $\gamma =e^{\P(t)}$, Theorem~\ref{t1sc5} (b), and {\rm OSC}, we obtain
\beq\label{(1mg5)} 
\begin{aligned}
m_t(G)
&=m_t(J_S\cap G) 
\leq \sum_{\om \in \mathfrak F}m_t( \varphi_\om(X_{t(\om)}))
=\sum_{\om \in \mathfrak F}m_t(\varphi_\om(\Int X_{t(\om)})) \\
& \leq M^{-1}K^t\sum_{\om \in \mathfrak F}m( \varphi_\om(\Int X_{t(\om)}))\\
& = M^{-1}K^t m\left( \bigcup_{\om \in \mathfrak F} \varphi_\om(\Int X_{t(\om)})\right)\\
& \leq M^{-1}K^t m(G).
 \end{aligned}
 \eeq
 Now, if $\Gamma \sbt X$ is  a $G_\delta$  set, then $\Gamma=\bigcap_{n=1}^\infty G_n$, where $(G_n)^\infty_{n=1}$ is a descending sequence of open sets. So, using
 (\ref{(1mg5)}), we obtain
\beq\label{(2mg5)}
 m_t(\Gamma)
 =\lim_{n \to \infty}m_t(G_n)
 \leq M^{-1}K^t \lim_{n \to \infty}m_(G_n)=M^{-1}K^t  m(G).
 \eeq
 Finally, if $Y$ is an arbitrary Borel subset of $X$, then there exists $\hat{Y}$, a $G_\delta$ subset of $X$ such that $Y \sbt  \hat{Y}$, $ m_t(\hat{Y})= m_t({Y})$ and
 $ m(\hat{Y})= m({Y})$. So, using (\ref{(2mg5)}), we obtain
$$ 
m_t(Y)= m_t(\hat{Y}) \leq M^{-1}K^t m(\hat{Y})=  M^{-1}K^t m(Y).
$$
The proof of Claim~4 is complete. \qed

\sp We assume that (c) holds. Let
$$ 
\rho:= \frac{dm_t}{dm}.
$$
We then have for every $\om \in E^{*}_A$  and every Borel set  $F \in X_{t(\om)}$ that
$$  
m_t ( \varphi_\om(F))
= \int_{\varphi_\om(F)} \rho dm 
= e^{-\P(t)|\om|}\int_{F} \rho\circ\varphi_\om |\varphi_\om'|^t dm.
$$
On the other hand, because we have already proved that $\rm{(a)}\Rightarrow \rm{(d)}$ (so (c) holds with $m$ replaced by $m_t$), we get 
$$ 
m_t(\varphi_\om(F))
=e^{-\P(t)|\om|}\int_F|\varphi_\om'|^t dm_t
=e^{-\P(t)|\om|}\int_F \rho |\varphi_\om'|^t dm.
$$
Therefore,
$$ 
\rho\circ \varphi_\om(x)= \rho(x)
$$
for $m_t$--a.e. $x\in X_{t(\om)}$. So, for $\tilde{m}_t$--a.e. $\tau\in E^{\infty}_A$, we get 
$$ 
\rho\circ \pi(\sigma(\tau))
= \rho \circ\varphi_{\tau_1}(\pi( \sigma(\tau)))
= \rho\circ \Pi(\tau).
$$
This means that the function $\rho\circ \pi: E^{\infty}_A \lra \mathbb R$ is $ \sigma$--invariant, and ergodicity of the  measure $\tilde{\mu}_t$ with respect to $\sigma$, yield that the function $\rho\circ\pi$ is $\tilde{m}_t$ --a.e.constant. Thus, we  have also:

\sp\fr{\bf Claim 5}. 
{\it The function $\rho: J_S \lra \mathbb R$ is  $m_t$--a.e. constant.}

\sp The next step in the proof of (a) is to show the following. 

\sp\fr{\bf Claim 6}. {\it  $ \rho=1$  $m_t$--{\it a.e.}, and consequently $ m=m_t$.}

\sp\bpf Seeking contradiction assume that $\rho>1$. Then  $1=m_t(J_S)=\rho m(J_S)$. So, $m(J_S)=1/\rho< 1$. Then $m(X \sms J_S)=1-\frac{1}{\rho}>0$, whence the formula
$$  
\hat{m}(A):=\frac{\rho}{\rho-1}m(A\sms J_S)
$$
defines a Borel probability measure on $X$. For every $e \in E$ and every Borel set $F \sbt X_{t(e)}$, we thus have
$$ 
\begin{aligned} m(\varphi_e(F)\sms J_S)  
& = m(\varphi_e(F))-m(\varphi_e(F)\cap J_S)
  = e^{-\P(t)} \int_F|\varphi_e'|^t dm - \frac{1}{\rho}m_t(\varphi_e(F)\cap J_S)\\
 & = e^{-\P(t)} \int_F  | \varphi_e'|^t dm - \frac{1}{\rho}m_t(\varphi_e(F))\\
 & = e^{-\P(t)} \int_F  | \varphi_e'|^t dm - \frac{1}{\rho}e^{-\P(t)} \int_F  | \varphi_e'|^t dm_t \\
 & = e^{-\P(t)} \int_F  | \varphi_e'|^t dm - e^{-\P(t)} \int_{F\cap J_S}  | \varphi_e'|^t dm \\
  & = e^{-\P(t)} \int_{F\sms  J_S} | \varphi_e'|^t dm.
\end{aligned}
$$
Therefore,
$$ 
\hat{m}(\varphi_e(F))= e^{-\P(t)}\int_{F} |\varphi_e'|^t d\hat m.
$$
Hence, $ \hat{m}$ is a Borel probability measure on $X$ satisfying  condition (c). So, by  already proven Claim 4, $m_t$ is absolutely  continuous with respect to $\hat{m}$.
Consequently, $ m_t( J_S)=0$ as $\hat{m}(J_S)=0$. This contradiction shows that $\rho =1$ and $m=m_t$. The implication $\rm{(c)}\Rightarrow \rm{(a)}$ is thus established.

\sp By the just established implication $\rm{(c)}\Rightarrow \rm{(a)}$, we conclude that if (c) holds, then $m(J_S)=1$. Having this, the implication $\rm{(d)}\Rightarrow\rm{(g)}$ becomes obvious, and the proof of Theorem~\ref{t420131007} is complete. \qed

\sp\section{Conjugacies of Conformal Graph Directed Systems} 

\bdfn\label{d720180604}
We say that two conformal graph directed Markov systems $S_1=\{\phi_e:X_1\to X_1:e\in E\}$ and $S_2=\{\psi_e:X_2\to X_2:e\in E\}$ with the sames alphabet $E$, the same set of vertices, and the same incidence matrix $A:E\times E\to\{0,1\}$, are topologically conjugate if there exists a homeomorphism $H:J_{S_1}\lra J_{S_2}$ such that 
\beq\label{820180604}
H\circ \phi_e=\psi_e\circ  H
\eeq
for every $e\in E$.

We say that $S_1$ and $S_2$ are bi--Lipschitz conjugate if the map $H$ is bi--Lipschitz continuous. Then $H$ uniquely extends to a bi--Lipschitz map 
from $\ov{J}_{S_1}$ to $\ov{J}_{S_1}$ and \eqref{820180604} continues to hold for this extension. 

We say $S_1$ and $S_2$ are conformally conjugate if the map $H:J_{S_1}\lra J_{S_2}$ has a conformal extension to a map from $X_1$ to $X_2$. 
\edfn

We record the following two immediate observations.

\bobs\label{o1220180604}
Any two conformally conjugate conformal graph directed Markov systems are bi--Lipschitz conjugate and any two bi--Lipschitz conjugate systems are topologically conjugate. Also then \eqref{820180604} continues to hold for these extensions. 
\eobs

\bobs\label{o2220180604} 
If two conformal graph directed Markov systems $S_1=\{\phi_e:X_1\to X_1:e\in E\}$ and $S_2=\{\psi_e:X_2\to X_2:e\in E\}$ are bi--Lipschitz conjugate then the following hold:
\begin{enumerate}
\item $\P_{S_1}(t)=\P_{S_2}(t)$ for all $t\ge 0$.

\, \item $\th(S_1)=\th(S_2)$.

\, \item The system $S_2$ is regular, strongly regular or hereditarily regular respectively if and only if $S_1$ is regular, strongly regular or hereditarily regular. 
\end{enumerate}
\eobs

\part{Elliptic Functions A}\label{EFA} 

In this part of the book we first, in Chapter~\ref{elliptic-theory}, General Theory of Elliptic Functions; Selected Properties, present some part of a well known classical theory of elliptic functions, with especial emphasis in Weierstrass elliptic functions. Then we pass to the dynamics of elliptic functions. As a matter of fact, in Chapter~\ref{topological-picture}, Topological Picture of Iterations of (all!) Meromorphic
Functions, we provide a short account of topological dynamics of all meromorphic functions with emphasize on Fatou domains. Chapter~\ref{geometry-and-dynamics}, Geometry and Dynamics of Elliptic Functions, is about all elliptic functions. We prove in it good estimates of Hausdorff dimension of Julia sets of all elliptic functions and estimate Hausdorff dimension of escaping points. In Chapter~\ref{first-outlook}, Compactly Non--Recurrent Elliptic Functions: First Outlook, we introduce the concepts of non--recurrent and compactly non--recurrent elliptic functions which will be the main object of our interest in Part~\ref{EFB}, Elliptic Functions B, through the end of the book. The last chapter, Chapter~\ref{examples}, Examples of Compactly Non--Recurrent Elliptic Functions, of the current part of the book, Part~\ref{EFA}, Elliptic Functions A, is devoted to provide many (classes of) examples of compactly non--recurrent elliptic functions. Examples exhibiting a variety of topological and dynamical properties. These will be mainly Weierstrass elliptic functions, and this is where, the theory presented in Chapter~\ref{elliptic-theory}, General Theory of Elliptic Functions; Selected Properties, will be primarily used in our book.

\chapter{General Theory of Elliptic Functions; \\ Selected Properties}\label{elliptic-theory}

\sp In this chapter, using a variety of contemporary sources, we present some parts of a well known classical theory of elliptic functions, with especial emphasis in Weierstrass elliptic functions. For dynamical, topological and fractal properties, this chapter will be primarily used in Chapter~\ref{examples}, Examples of Compactly Non--Recurrent Elliptic Functions, devoted to provide many (classes of) examples of compactly non--recurrent elliptic functions that exhibit a variety of topological and dynamical properties.  We primarily follow here the classical books \cite{Du} and \cite{JS}. We would also like to bring reader's attention to the books \cite{AE} and \cite{Lang-Elliptic}.

\sp

\section{Periods, Lattices and Fundamental Regions}\label{periods}

Let $f$ be a function  defined on the complex plane $\mathbb{C}$. Then a  complex number $w$  is called a period of $f$\index{(N)}{period of function} if 
$$  
f(z+w)=f(z)
$$
for all $z\in \mathbb{C}$. The function $f$  is called periodic if it has a
period $w\neq 0$. For example, the maps 
$$
\C\ni z\longmapsto \sin(z)\in\C
\  \  \  {\rm and}  \  \  \ 
\C\ni z\longmapsto\cos(z)\in\C
$$ 
have period $2\pi $, the map 
$$
\C\ni z\longmapsto e^z\in\C
$$ 
has  period $2\pi i$, while for any complex number $w\neq 0$, the map 
$$
\C\ni z\longmapsto\sin(2\pi z/w)
$$
has period $w$. The set $\Lambda_f$ \index{(S)}{$\Lambda_f$} of all periods of a function $f$ has two important
properties: one algebraic, valid for  all $f$, and one topological,
valid for  non--constant meromorphic functions. These properties are
given in Theorem~\ref{t1jones} and Theorem~\ref{t2jones}.

\sp \bthm\label{t1jones} For every 
function $f:\mathbb{C}\lra\oc$, the set $\La_f$ is a subgroup of the
additive group $\mathbb{C}$. 
\ethm

\bpf Let $\a, \b\in \La_f$. Then
$$
f(z+(\a+\b))=f((z+\a)+\b)=f(z+\a)=f(z)
$$ 
and so $\a+\b\in
\La_f$. Moreover, $f(z-\a)=f((z-\a)+\a)=f(z)$, and so
$-\a\in\Lambda_f$. Finally $f(z+0)=f(z)$, and so $0\in \La_f$. Thus
$\La_f$ is  a subgroup of $\mathbb{C}$. \endpf

\sp

A subset $\Delta$ of a topological space is called
discrete\index{(N)}{discrete subset} if every pint $x \in \Delta$ has  a
neighbourhood $U$ such that  $U\cap \Delta=\{x\}$. For  example

\begin{itemize}

\sp\item [(a)] the integers $\mathbb{Z}$  form a  discrete subset of
$\mathbb{R}$;

\sp\item [(b)] any finite subset of ${\mathbb R}^n$  is discrete;

\sp\item [(c)] $\{ 1/n; \,  n \in {\mathbb Z},\, n \neq 0\}$  is a
discrete subset of $\mathbb R$.
\end{itemize}

\sp \fr However, the set $\{ 1/n;\,  n \in {\mathbb Z}\}\cup\{0\}$ is not discrete,
since every  neighbourhood  of $0$  contains a real number of the form
$1/n$, $n \in \mathbb Z$.

\sp Let $f:\mathbb{C}\to\oc$ be a meromorphic function. In particular $f$ is continuous, and therefore the set $\La_f$ is closed. Note that for every $\xi\in \C\cap f(\C)$ and for every $z\in f^{-1}(\xi)$, 
$$
z+\La_f\sbt f^{-1}(\xi).
$$
Since the set $z+\La_f$ is discrete if and only $\La_f$ is, and since zeros of any non--constant meromorphic function, form a discrete set, we get the following.

\bthm\label{t2jones} For every non--constant meromorphic function
$f:\mathbb{C}\lra\oc$, $\La_f$ is a closed discrete subset of $\mathbb{C}$. 
\ethm

\sp To summarize, the periods of  a non--constant  meromorphic function
form a discrete subgroup of $ \mathbb C$. We now will show that there
are exactly three  types of discrete  additive  subgroup of $\mathbb C$: isomorphic  to $\{0\}$, $\mathbb Z$, and $\mathbb Z \times \mathbb Z$  respectively. 

Given two complex number $w_1,w_2\in\C$, we denote
$$
[w_1,w_2]:=\{mw_1+nw_2: m, n \in \mathbb Z\}.
$$
Of course $[w_1,w_2]$ is an additive subgroup of $\C$ but we will have more to say about it. 
\bthm\label{t3jones} Let $\La_f$ be a discrete subgroup of $\mathbb
C$. Then one of the following  holds:
\begin{itemize}
\item [(1)]  $\La =\{0\}$.

\sp\item [(2)] There exists $w\in \C$ such that $\La=[w,w]= \{nw:n\in \mathbb Z\}$; in particular. $\La $ is isomorphic to $\mathbb Z$.

\sp\item [(3)] $\La=[w_1,w_2]$ for some 
$w_1, w_2 \in \mathbb C$ that are linearly
independent over $\mathbb R$, that is, $w_1\neq 0\neq w_2$ and
$w_1/w_2\notin \mathbb R$; in particular, the additive group $\La$ is isomorphic to
$\mathbb Z\times \mathbb Z$.
\end{itemize}
\ethm

\bpf Suppose that $\La\neq \{0\}$. Since $\La$ is discrete and closed subset of $\C$, there exists $w_1\in \La\sms\{0\}$ with a least value of $|w_1|$. Of course, $w_1$ is not unique: $-w_1$ for example will do equally  well. Now let  
$$
L=\{t w_1: t\in \mathbb R\}
$$ 
be the  line through $0$ and  $w_1$ in $\mathbb C$. Then 
$$
\{n w_1:n \in \mathbb Z \}\sbt\La\cap L.
$$
First suppose that $\La \sbt  L$. We then
claim that
$$
\La =\{ n w_1:n \in \La\},
$$
so that $\La$ is  of type (2). For if $\La  $ contains  some $w \neq nw_1$, then since $\La\sbt L$ we have $w =tw_1$ for some $t\in \mathbb R \sms \mathbb
Z $. So, $ n< t <n+1$ for some $n \in  \mathbb Z$. Since $\La $ is a
group containing  $w $ and $nw_1$,  it contains $w-nw_1
=(t-n)w_1$. However, $0\neq| (t- n)w_1|< |w_1|$, which contradicts
the minimality of $|w_1|$. Hence 
$$
\La =\{nw_1;\, n \in \mathbb Z\},
$$
and we are  in case (2).

Now suppose that $\La \nsubseteq L$. Since $\La \sms  L$ is a discrete and closed subset of $\C$ (as a closed subset of $L$ (note the relative topology on $L$ is discrete), the set $\La \sms  L$ is closed in $\C$). Hence, $\La \sms  L$ has an element $w_2$ with a least modulus. Then   
$$
\La\spt\Om:=\{m w_1+n w_2: m,n\in \mathbb Z\},
$$
and since  $w_2\notin
L$, the vectors $w_1$ and $w_2$ are  linearly independent over $\mathbb R$. So
$\Om$ consists of  the vertices of a tessellation of $\mathbb C$ by
congruent  parallelograms.

We now will show that 
$$
\La=\Om.
$$
If this is not the case, then  there exists $w\in \La\sms\Om$. Let $w=\l w_1+\mu w_2$  with $\l, \mu \in \mathbb R$. Then  by
adding suitable integral multiplies  of $w_1$  and $w_2$, we
may assume  that 
$$
 |\l|\leq \frac{1}{2}
 \  \  \  {\rm and} \  \  \
 |\mu|\leq \frac{1}{2}.
$$
If $\mu=0$ then $w=\l w_1\in L$, with $|w|=|\l w_1|< |w_1|$. By
minimality of $|w_1|$ we then have $w=0$, and hence $w\in \Om$, contrary to
our assumption. If $\l=0$ then $w=\mu w_2$, and again  $w=0$, this
time by minimality of $|w_2|$.  Hence $\l w_1$ and $\mu w_2$ are
non-zero and therefore  linearly independent  over  $\mathbb R$,
giving
$$ 
|w|< |\l w_1| +|\l w_2|
\leq \frac{1}{2}|w_1|+\frac{1}{2}|w_2|
\leq \frac{1}{2}|w_2|+ \frac{1}{2}|w_2|=|w_2|
$$
Now $w\in \La \sms L$ since $\mu\neq0 $. So, by
minimality of $|w_2|$, we have $w=0$, contradicting the fact that
$w\notin\Om$. Thus $\Om=\La$ and we are in the case (3).
\epf

\sp

\bdfn Let $f:\C\lra\oc$ be  a meromorphic function. 

\begin{itemize}

\sp\item If $f$ has its set $\La_f$\index{(S)}{$\La_f$}
of periods of type (2), then $f:\C\lra\oc$ is called simply
periodic\index{(N)}{simply periodic function}.

\sp\item If $\La_f$ is of type
(3), then $f$\index{(N)}{doubly periodic function} is called doubly
periodic.

\sp\item Any group $\La$ of type (3) is called a lattice\index{(N)}{lattice}, and any pair $\{w_1, w_2\}$ of complex numbers such that 
$$
\La=[w_1,w_2]
$$ 
is called a basis\index{(N)}{basis for a lattice}, or a generator, for the lattice $\La$. 
\end{itemize}
\edfn

\fr The above results show that a non--constant meromorphic periodic  function is either simply periodic  or doubly  periodic.

\sp We  now turn to the study of doubly periodic meromorphic functions, starting
with  a closer look at some of algebraic and geometric ideas
involved. Let $\La$ be a group of type (3). Write 
$$
\La=[w_1,w_2],
$$
where $\{w_1, w_2\}$ is some generator of $\La$. Of course, there
are other bases for $\La$ besides $\{w_1, w_2\}$. For instance
$\{w_1,w_1+ w_2\}$ is also a basis, for if $\l \in [w_1,w_2]$
then
$$ 
w=mw_1+nw_2=(m-n)\l_1+n(w_1+w_2),
$$
with $m-n\in \mathbb Z$. In general, if $w_1', w_2'\in [w_1,w_2]$ then
\beq\label{ej1}
    w_1'=aw_1+bw_2 \quad \mbox{and}\quad
    w_2'=cw_1+w_2,
\eeq
where $a,b,c,d$ are integers.

\sp\bthm\label{t4jones} 
Equations $(\ref{ej1})$ define a basis $\{w_1',
w_2'\}$ for $[w_1, w_2]$ if and only if 
$$
ad-bc=\pm 1.
$$
\ethm

\bpf It is convenient to write $(\ref{ej1})$  using
matrix notation

\beq\label{ej2}
\left(
  \begin{array}{c}
    w_1' \\
    w_2' \\
  \end{array}
\right)=A\left(
           \begin{array}{c}
              w_1 \\
              w_2 \\
           \end{array}
         \right),
\eeq where
$$
\left(\begin{array}{cc}
        a & b \\
        c & d\\
      \end{array}
      \right).
$$
If $ad-bc= \pm1$, then $A^{-1}$ has  integer  coefficients, and we
have
$$\left(
  \begin{array}{c}
    w_1 \\
    w_2 \\
  \end{array}
\right)=A^{-1}\left(
           \begin{array}{c}
               w_1' \\
               w_2' \\
           \end{array}
         \right)=\pm \left(\begin{array}{cc}
        d & - b \\
        - c & a \\
      \end{array}
      \right)\left(
  \begin{array}{c}
    w_1' \\
    w_2' \\
  \end{array}
\right).
$$ 
Thus $w_1, w_2 \in [w_1', w_2']$ and hence $[w_1, w_2]\sbt [w_1', w_2']$. The reverse inclusion is obvious, so
$$
[w_1, w_2]=[w_1', w_2']
$$ 
and hence  $\{w_1', w_2'\}$ is a
basis. Conversely,  suppose that equation $(\ref{ej2})$ define a
basis $\{w_1',w_2'\}$ for $[w_1, w_2]$. Expressing the
elements $w_1, w_2$ in therms  of the basis, we  have
$$ \left(
  \begin{array}{c}
    w_1 \\
    w_2 \\
  \end{array}
\right)=B\left(
           \begin{array}{c}
               w_1' \\
               w_2' \\
           \end{array}
         \right),
$$
for some fix matrix $B$ with  integer coefficients. So
\beq\label{ej3}
\left(
  \begin{array}{c}
    w_1 \\
    w_2 \\
  \end{array}
\right)=BA\left(
           \begin{array}{c}
               w_1 \\
               w_2 \\
           \end{array}
         \right).
\eeq
Since $w_1$ and $w_2$ form a basis of $\C$ considered as a 2-dimensional  vector space over $\mathbb R$, we thus have that $BA=\Id$. Therefore, 
$$  
\det(B)\cdot\det(A)=1.
$$
Since both $A$ and $B$  have  integer coefficients, their  determinants
are integers, and therefore $\det(A)=\pm 1$, that is, $ad-bc=\pm 1$. \endpf

\

If $\La$ is a lattice  and $a\in \mathbb{C}$, then we define 
$$
a\La:=\{aw: \, w \in \La\}.
$$
This is a lattice if and only if $a\neq 0$.

\sp\bdfn\lab{d12020200309}
Two lattices $\La $ and $\La'$ are called similar\index{(N)}{similar
lattices} if and only if 
$$
\La =a \La'
$$ 
for some $a\neq 0$. 
\edfn

\sp Keep $\{w_1',w_2'\}$ a generator of a lattice $\La$ and $\{w_1',w_2'\}$ a generator of a lattice $\La'$. Linear independence of $w_1$ and $w_2$ over
$\mathbb{R}$ means that $\im (w_1/w_2)\neq 0$. Interchanging $w_1$ and
$w_2$ if necessary, we may assume that 
$$
\im \lt(w_1/w_2\rt)>0.
$$ 

We define the {\it modulus}\index{(N)}{modulus of the  basis} of the basis
$\{w_1,w_2\}$ to be
$$ 
\tau:=w_1/w_2,
$$\index{(S)}{$\tau$}
\fr where the numbering is such that $ \im(\tau) >0$.
Each lattice $\La $ determines its set of  moduli, the moduli of its
all bases. Since $\mu w_1/\mu w_2=w_1/w_2$, similar lattices
have the same  set of moduli. Putting 
$$
\tau':=w_1'/w_2',
$$
the modulus of the basis $\{w_1',w_2'\}$ for $\La'$, we see that $\La $ and $\La'$ are similar if and only if
\beq\label{ej4}
\tau'=\frac{a\tau+b}{c\tau+d},
\eeq
where $a,b,c.d\in \mathbb{Z}$, and $ad-bc=\pm 1$. Now, both $\tau$
and $\tau'$, being moduli, lie in the upper half plane
$$
\mathbb H=\{z\in \mathbb{C}: \,  \im z>0\};
$$
and if $ad-bc=-1$ then the  M\"obius transformation
$$
\oc\ni z\lmt T(z)=\frac{az+b}{cz+d}\in\oc
$$
maps $\mathbb{H}$ onto the lower half-plane,
so we  must therefore have $ad-bc=1$. Conversely, since $a,b,c,d \in
\mathbb{Z}$ and $ad-bc=1$, we have that
\beq\label{ej5}
\begin{aligned}
  w_1' &= \mu(aw_2+bw_1) \\
  w_2' &= \mu(cw_2+dw_1)
\end{aligned}
\eeq 
gives basis $[w_1',w_2]$ for a lattice $\La'$ similar to $\La$.

\sp  The M\"{o}bius transformations\index{(N)}{M\"{o}bius transformation}
$$ 
\oc\ni z \longmapsto T(z):=\frac{az+b}{cz+d}\in\oc, \  \,\, a,b,c,d \in \mathbb{Z}, \, \, ad-bc=1,
$$ 
form a discrete group, called  the modular group \index{(N)}{modular group} $
\mathbb{Z}$\index{(S)}{$PSL(2, \mathbb{Z})$}.

\sp\bthm\label{t5jones} If $\La=[w_1, w_2]$ and $\La'=[w_1',
w_2']$ are lattices in $\mathbb C$, with respective moduli $\tau=w_2/w_1$ and
$\tau'=w_2'/w_1'$, both in $\mathbb{H}$, then the
lattices $\La$ and $\La'$ are similar if and only if 
$$
\tau'=T(\tau)
$$
for some $T \in PSL(2,\Z)$. 
\ethm

\sp Given a lattice $\La$, we define $z_1, z_2\in \mathbb C$ to be 
congruent \index{(N)}{congruent points}\!\!\! \!\!\!\!\! $\mod \La$, written 
$$
z_1\sim_\La z_2,
$$\index{(S)}{$\sim_\La $}
if and only if 
$$
z_1-z_2\in \La.
$$
The congruence mod
$\La$ is easy to be seen an equivalence relation  on $\mathbb
C$, and the equivalence classes are cosets $z+ \La$  of $\La$ in
the additive group $\mathbb C$. Alternatively, we  may regard $\La$
as acting on $\mathbb C$ as transformation  group. Each $w\in \La$
induces the transformation
$$  
\C\ni z\lmt T_w (z) =z+w\in\C
$$ 
of $\mathbb C$ onto itself. Since
$$ 
T_{w_1+w_2}= T_{w_1}\circ T_{w_2},
$$
we have a group  isomorphism 
$$
\La\ni w\longmapsto (T_w:\C\lra\C).
$$
Then  two points $z_1, z_2 \in \mathbb C$ are  congruent mod $\La$
if and only if  they lie  in the same orbit under this action of
$\La$.

\sp \bdfn A closed, connected  subset $\mathcal{R}$  of $\mathbb C$ is
called a fundamental  region\index{(N)}{fundamental region}
for a lattice $\La\sbt\C$ if and only if

\sp\begin{enumerate}
  \item   for each $z\in \mathbb C$, $\mathcal{R}$  contains at
  least one point in the same $\La$--orbit as $z$, i.e. every point
  $z\in \mathbb C$ is congruent to some  point  in $\mathcal{R}$,
  
  \sp\item  no two points in the interior of $\mathcal{R}$ are  in
  the same  orbit of $\La$, i.e. no two distinct points in the interior of
  $\mathcal{R}$ are congruent mod $\La$.
\end{enumerate}
\edfn

\sp If, as usually is the case in applications, $\mathcal{R}$  is also an Euclidean
 polygon, with a finite numbers of sides, then we  call
 $\mathcal{R}$ a {\it fundamental  parallelogram}\index{(N)}{fundamental  parallelogram} for $\La$.

\sp For example, the  parallelogram  $\mathcal{R}$  with vertices $0,
 w_1, w_2, w_1+w_2$ is a fundamental parallelogram for the
 lattice $[w_1, w_2]$. Conditions (1), (2) of the above definition ensure that if
 $\mathcal{R}$  is an  fundamental  regions for a lattice $\La$,
 then $\mathcal{R} $ and its image under action of  $\La$ (that is,
 its translates $ \mathcal{R}+w, \, w \in \La$) cover the plane
 $\mathbb C$ completely, overlapping only at their  boundaries. This
 type of covering is known as tessellation\index{(N)}{tessellation} of $\mathbb C$. By using (\ref{ej1})  we can obtain fundamental
 parallelograms  of  different shapes, and hence obtain different
 tessellations of $\mathbb C$.

 \sp  If $\mathcal{R}$  is any fundamental region  for $\La$, then  for
  a fixed  $t\in \mathbb C$, the set
$$ 
\mathcal{R}+t=\{ z+t: \,  z \in \mathcal{R}\}
$$
is also a fundamental region. This is  useful when we need to find a
fundamental  region containing or avoiding certain  specified
points, for  example, we can always find  a fundamental
parallelogram for $\La$  with 0 in its  interior.

\sp We have seen that a fundamental region for a lattice is not unique.
However the next theorem shows that, that its area is unique, and
may therefore be regarded as a function of the lattice alone. For notational convenience, we will write $T_w(X)$ for $X+w$; since  the  translation $T_w(z)=z+w$ is an isometry of $\mathbb  C$, we have 
$$
S(T_w(X))=S(X),
$$
where $S$ denotes planar Lebesgue measure on $\C$.

\sp\bthm\label{t6jones} If $\mathcal{R}_1$ and $\mathcal{R}_2$ are two fundamental regions for a lattice $\La$ whose boundaries are of planar Lebesgue measure zero, for example if $\mathcal{R}_1$ and $\mathcal{R}_2$ are
polygons, then
$$
S(\mathcal{R}_1)=S(\mathcal{R}_2).
$$
\ethm

\bpf Because of our hypotheses, we have that
$$
S(\mathcal{R}_1)=S(\Int\mathcal{R}_1)
\  \  \  {\rm and} \  \  \
S(\mathcal{R}_2)=S(\Int\mathcal{R}_2)
$$
We have 
$$
\mathcal{R}_1
\supseteq \mathcal{R}_1 \cap \bigcup_{w \in\La} T_w(\Int\mathcal{R}_2)
=\bigcup_{w \in\La}\mathcal{R}_1 \cap T_w(Q_2).
$$
As $\Int\mathcal{R}_2$ is the interior of a fundamental region, the sets
$\{\mathcal{R}_1\cap T_w(\Int\mathcal{R}_2):w\in\La\}$ are pairwise disjoint. Hence,
$$
 S(\mathcal{R}_1)
 \geq \sum_{w\in \La}S\( \mathcal{R}_1\cap
T_w(\Int\mathcal{R}_2)\)
    =\sum_{w \in \La} S\(T_{-w}(\mathcal{R}_1)\cap \Int\mathcal{R}_2\)
    =\sum_{w \in \La} S\(T_w(\mathcal{R}_1)\cap \Int\mathcal{R}_2\).
$$
Now as $\mathcal{R}_1$ is a  fundamental region, $\bu_{w \in \La} T_w(\mathcal{R}_1)=
\mathbb C$, and hence 
$$
\bu_{w \in \La} T_w(\mathcal{R}_1)\cap \Int\mathcal{R}_2
=\Int\mathcal{R}_2.
$$
Therefore,
$$
\sum_{w\in \La} S(T_w(\mathcal{R}_1)\cap \Int\mathcal{R}_2)\geq S\(\Int\mathcal{R}_2\)=S(\mathcal{R}_2),
$$
giving $S(\mathcal{R}_1)\geq S(\mathcal{R}_2).$ Interchanging
$\mathcal{R}_1$ and $\mathcal{R}_2$ we have $S(\mathcal{R}_2)\geq
S(\mathcal{R}_1)$. Thus 
$$
S(\mathcal{R}_1)=S(\mathcal{R}_2)
$$ 
and the theorem is proved. 
\endpf

\sp The proof of the next, well known, theorem can be found for instance in \cite{JS}, see Theorem~5.8.4 therein.

\sp

\bthm\label{t7jones}  
The set 
$$
{\mathfrak F}:=\lt\{z\in {\mathbb H}: 
|z|\ge 1\  \  \mbox{and}\  \ |\re z | \leq \frac{1}{2}\rt\}
$$ 
is a fundamental region for the modular group $PSL(2,\Z)$. 
\ethm

\sp\section{General Properties of Elliptic Functions}\label{elliptic}

\bdfn 
A meromorphic  function $f:\mathbb C \lra \oc$ is called
elliptic \index{(N)}{elliptic function} if and only if it is doubly periodic. \index{(N)}{doubly periodic} This means that $\La_f$, the set of periods of $f$, is a lattice in $\C$. 
\edfn

\sp

\fr Thus, if $f:\mathbb C \lra \ov{\mathbb C}$ is an elliptic function, then 
$$
f(z+w)=f(z)
$$ 
for all $z\in \mathbb C$ and all $w\in \La_f$.

\sp So far, the only elliptic  functions we  have met in this chapter
were the constant functions, and it is a substantial problem to
construct non-constant elliptic functions. Before doing this, we
will  examine some of the elementary properties which elliptic
functions must possess.

\sp Given a lattice $\La\sbt\C$, we denote
$$
\mT_\La:=\C/\La,
$$
the quotient space generated by the lattice $\La$. So, $\mT_\La$ is a compact connected Riemann surface and topologically it is the $2$-dimensional torus. Let
$$
\Pi=\Pi_\La:\C\lra\mT_\La
$$
be the corresponding quotient map.

\sp\fr Keep now $f:\mathbb C \lra \oc$, a non--constant elliptic function. We define
\beq\label{120140104}
\mT_f:=\mT_{\La_f}=\C/\La_f \  \  \text{ and }  \  \  \Pi_f:=\Pi_{\La_f}:\C\lra\mT_f.
\eeq
Denote
$$
\hat\mT_f:=\mT_f\sms \Pi_f(f^{-1}(\infty)).
$$
Of course the function $f:\mathbb C \lra \oc$ induces a unique continuous map 
$$
\hat f:\hat\mT_f\lra\mT_f
$$ 
such that the following diagram commutes:
\beq\label{120200408}\begin{tikzcd}
{\C\sms f^{-1}(\infty)} \arrow{r}{f} \arrow[swap]{d}{\Pi_f} & {\C}\arrow{d}{\Pi_f} \\
{\  \hat{\mathbb T}_f} \arrow{r}{\hat f} & {\ {\mathbb T}_f}
\end{tikzcd}
\eeq
Evidently $\hat f:\hat\mT_f\lra\mT_f$ is a holomorphic map. Actually even more, we have a holomorphic map 
\beq\label{120180816}
\^f:\mT_f\lra\oc
\eeq
such that the following diagram commutes:
\beq\label{220180816}
\begin{tikzcd}
{\hat{\mathbb T}_f}  \arrow{r}{\tilde f} \arrow[swap]{dr}{\hat f} 
    & \C\arrow{d}{\Pi_f} \\
 & {{\mathbb T}_f}.
\end{tikzcd}
\eeq
This means that
\beq\label{320180816}
\hat f=\Pi_f\circ\^f
\eeq
Fix $c \in \oc$. Then the solutions of the equation 
$$
f(z)=c
$$  
are isolated, and  each solution has finite  multiplicity, with congruent solutions having the same  multiplicity. Since the solutions are isolated, any
fundamental  polygon $\mathcal{R}$ for $\La_f$ contains only finitely
many  solutions (since $\mathcal{R} $ is compact), and, replacing
$\mathcal{R}$ by $\mathcal{R}+t$ with some appropriate $t\in \mathbb C$ if necessary,
we may assume that there are no solutions on the boundary
$\partial{\mathcal{R}}$ of $\mathcal{R}$. Let  
$$
\{z_1, z_2, \ldots, z_r\}:=\mathcal{R} \cap f^{-1}(c).
$$
Let 
$$
k_1,k_2,\ld,k_r
$$
denote respective multiplicities of solutions $z_1, z_2, \ldots, z_r$. Let 
$$
N:=k_1+k_2+\ldots+ k_r.
$$
Then we say that there are $N$ (counted with multiplicities) solutions to the equation $f(z)=c$. Of course, $N$ is also the number of solutions (counted with multiplicities) to the equation
$$
\hat f(z)=\Pi_f(c)
$$
on the torus $\hat\mT$, and 
$$
\hat f^{-1}(\Pi_f(c))=\Pi_f(f^{-1}(c))=\Pi_f(\{z_1, z_2, \ldots, z_r\}).
$$
\bdfn 
The order {\rm ord}$(f)$\index{(S)}{{\rm
ord}$(f)$}\index{(N)}{order of elliptic function} of an elliptic
function $f$ is the number of solutions to the equation $f(z)=\infty$, that is,
the sum of multiplicities of the congruence classes of poles of $f$. 
\edfn

\sp

This is analogous  to the degree  $\deg(f)$\index{(S)}{$\deg(f)$} of
a rational function $f$, which is equal to the number of solutions of
$f(z)=\infty$, counting multiplicities. 

\sp For the rest  of this
section, we  assume that $f$ is a non--constant elliptic function.
We denote ${\rm ord}(f)=N$, we denote by $\mathcal{R}$ a fundamental
parallelogram for $\La_f$ with vertices 
$$
t, \ t+w_1, \ t+w_2, \
t+w_1+w_2+w_3,
$$
where $[w_1, w_2]$ is a basis for $\La_f$, ad $t$ is
chosen so that $\partial{\mathcal{R}}$ contains no  zeros  or poles
of $f$.

\sp\bthm\label{t8jones} An elliptic function $f$ is constant if and only if $N=0$. In particular, a holomorphic elliptic function  must be constant. \ethm

\bpf  If $f$  is a constant and  meromorphic, then it has
no poles  in $\mathbb C$, so $N=0$. Conversely, suppose that $N=0$.
Then $f$ has no poles, so $f$ is holomorphic on $\mathbb C$. Now,
$\mathcal{R}$ is compact, and $f$ is continuous, so $f(\mathcal{R})$
is a compact subset of $\mathbb C$ and is therefore bounded. Since
$f(\mathbb C)=f(\mathcal{R})$, it follows  that $f$ is  bounded on
$\mathbb C$. Therefore, Liouville's Theore  implies that
$f$ must be constant. \endpf

\sp

\bthm\label{t9jones}  
The sum of the residues of any elliptic function $f:\C\lra\oc$ within its fundamental polygon $\mathcal R$ is equal to zero. 
\ethm

\bpf Since $f$ is meromorphic, and it is holomorphic on
$\partial{\mathcal R}$, we have that
$$ 
\frac{1}{2\pi i}\int_{\partial{\mathcal R}}f(z)\,dz
$$ is
equal to the sum of residues within $\mathcal R$. Now let
$\Gamma_1$, $\Gamma_2$, $\Gamma_3$ and $\Gamma_4$ be the sides of
$\mathcal R$ from $t$ to $t+w_1$, $t+w_1$ to $t+w_1+w_2$,
$t+w_1+w_2$ to $t+w_2$, and $t+w_2$ to $t$  respectively, so that
$$
\int_{\partial{\mathcal R}} f(z)\, dz
= \sum_{j=1}^4\int_{\Gamma_i}f(z)\, dz,
$$
where the direction  along $\Gamma_j$ is consistent with the
positive (i.e. anti-clock wise) orientation on ${\partial{\mathcal
R}}$. Now, since $w_2$ is a period of $f$ and since $\Gamma_3=-\Gamma_1+w_2$, where $-\Gamma_1$ is $\Ga_1$ with the reverse orientation, we can write
$$
\int_{\Gamma_3}f(z)\,dz
=\int_{-\Gamma_1+w_2}f(z)\,dz
=\int_{-\Gamma_1}f(z+w_2)\,dz
=\int_{-\Gamma_1}f(z)\,dz
=-\int_{\Gamma_1}f(z)\,dz.
$$ 
By the same token, $\int_{\Gamma_4}f(z)\,dz=-\int_{\Gamma_2}f(z)\,dz$. Hence, the
sum of the residues is equal to zero. \endpf

\

\bcor\label{c1jones} 
There are no elliptic  functions of order $N=1$. 
\ecor

\bpf  If $f$ were elliptic of order 1, it would have a
single pole of order 1 in $\mathbb R$, say, at $z=a\in \mathbb C$. Then
$$
f(z)=\sum_{j=-1}^\infty a_j(z-a)^j
$$ 
near $z=a$, with $a_{-1}\neq0$. Thus the sum  of the residues of $f$ within $\mathcal R$ would be equal to $a_{-1}$, which is non-zero, contradicting
Theorem~\ref{t9jones}.
\endpf

\

\bthm\label{t10jones}  
If $f$  has  order $N>0 $ then $f$  takes
each value $c\in \ov{\mathbb C}$ exactly $N$  times. 
\ethm

\bpf This is the definition of $N$ if $c=\infty$, so we
may assume that $c\in \mathbb C$. Replacing $f$ by $f-c$ (which is
has the same order as $f$), we may assume that $c=0$. Now $f'/f$ is
meromorphic, and since $\partial{R}$ contains  no poles or zeros of
$f$, the function $f'/f$ is  analytic on $\partial{R}$ . We may therefore
integrate $f'/f$,and applying the argument used in the proof of
Theorem~\ref{t9jones}, we get that
$$
\int_{\partial{R}}\frac{f'(z)}{f(z)}\,dz=0.
$$
This means that the sum of  the residue $f'/f$ within $\mathcal R$ is equal to zero.
Now $f'/f$ has poles at zeros and poles of $f$, and
nowhere else. Suppose that $f$  has zero of multiplicity  $k$ at
$ z=a\in \mathcal R$, so that 
$$
f(z)=(z-a)^kg(z)
$$ 
near $z=a$, where $g$ is analytic and $g(a)\neq 0$. Then
$$
f'(z)=k(z-a)^{k-1}g(z)+(z-a)^kg'(z)
$$  
near $z=a$, and so
$$ \frac{f'(z)}{f(z)}=\frac{k}{z-a}+\frac{g'(z}{g(z)}$$ near $z=a$,
so that $f'/f$ has residue $k$ at $z=a$. A similar argument, with
$f(z)=(z-a)^{-k}g(z)$, shows that $f'/f$ has residue $-k$ at each
pole of $f$ of multiplicity $k$. Since the sum of the residue of $f'/f$
is zero, the number of zeros of $f$ must be equal the number of
poles (counting multiplicities). Thus, the equation $f(z)=0$ has $N$
solutions, as required. The proof is complete.
\endpf

\sp

\bthm\label{t11jones}  
Let $f$ and $g$ be elliptic  functions  with
respect to some lattice $\La$. If $f^{-1}(\infty)=g^{-1}(\infty)$ and both $f$ and $g$ have the same principal parts at each pole, then $f(z)=g(z)+c$ for
some constant $c\in C$. 
\ethm

\bpf The function $f-g$ is elliptic and of order
0 since it has no poles. Thus, $f-g$ is constant by
Theorem~\ref{t8jones}. 
\endpf

\sp

\bthm\label{t12jones}  Let $f$ and $g$ be elliptic functions with the following properties.

\sp\begin{itemize}
\item $\La_f=\La_g$,

\sp\item $f^{-1}(0)=g^{-1}(0)$ and $f^{-1}(\infty)=g^{-1}(\infty)$,

\sp\item $f$ and $g$ have equal multiplicities at their zeros and their poles.
\end{itemize}
Then $f(z)=cg(z)$ for some constant $c\neq 0$. 
\ethm

\bpf Replace $f-g$ by $f/g$ in the proof of
Theorem~\ref{t11jones}. 
\endpf

\sp A rational function $f:\oc\lra\oc$ which is
not identically zero, must have finitely many zeros (say at
$a_1,\ldots, a_r$ with multiplicities  $k_1, \ldots, k_r$) and
finitely many poles (say at $b_1,\ldots, b_s$ with multiplicities
$l_1, \ldots, l_s$). Conversely, given any choice of points $
a_1,\ldots, a_r, b_1,\ldots, b_s\in \ov{\mathbb C}$ and
multiplicities $k_1, \ldots, k_r; l_1, \ldots, l_s\geq 1$, there
exists a rational function $f$ with these zeros and  poles, with
these  multiplicities, provided

\sp\begin{enumerate}
  \item  $k_1+ \ldots +  k_r=l_1+ \ldots +l_s$ (both must be equal to the
  degree of $f$), and

\sp\item  the sets $\{a_1,\ldots, a_r\}$ and $\{ b_1,\ldots, b_s\}$ are
 disjoint (zeros and poles  cannot overlap).
\end{enumerate}
We just take
$$ 
f(z):=\prod_j(z-a_j)^{k_j}\big/\prod_j(z-b_j)^{l_j},
$$
where these products range over all $j$ such that  $a_j, b_j\in
\mathbb C$, but  exclude factors where $a_j=\infty$ or $b_j=\infty$.

\sp If an elliptic function $f:\C\lra\oc$ is to
have its zeros and poles at the congruent classes $[a_1], \ldots,
[a_r]$ and $[b_1], \ldots, [b_s]$   with multiplicities $k_1, \ldots
,  k_r$ and $l_1,  \ldots, l_s$, then condition (1) is necessary by
Theorem~\ref{t10jones}, and corresponding to (2), we have  the
necessary condition

\sp\begin{itemize}
\item[{\rm(2')}] the sets $[a_1]\cup \ldots \cup
[a_r]$ and $[b_1]\cup \ldots\cup [b_s]$, are disjoint.
\end{itemize}

\sp

\fr The next result  shows that, in contrast with the situation for
rational  functions, these conditions are not {\it sufficient} for
the existence of $f$.

\sp

\bthm\label{t13jones}  Let the congruence classes  of zeros  and
poles of an elliptic function $f$ respectively be $[a_1],\ldots,[a_r]$ and
$[b_1], \ldots\cup [b_s]$ with respective multiplicities $k_1, \ldots, k_r$
and $l_1, \ldots, l_s$. Then
$$ 
\sum_{j=1}^r k_ja_j \sim_{\La_f}\sum_{j=1}^s l_j b_j.
$$
\ethm

\bpf As usually,  let $\mathcal R$ be a fundamental
parallelogram for $\La_f$, chosen so, that $f$ has no zeros or poles
on $\partial{\mathcal R}$. The effect of replacing any $a_j$ or $-
b_j$ by a congruent point is to add an element of $\La_f$ to 
$$
\sum_{j=1}^rk_ja_j-\sum_{j=1}^sl_j b_j,
$$  
and this does not affect the condition of the theorem, so we may assume that $a_j$, $b_j\in \mathcal R$ for all $j$. First we prove that 
$$ 
\sum_{j=1}^rk_j a_j-\sum_{j=1}^sl_jb_j=\frac{1}{2\pi i}
\int_{\partial{\mathcal R}} \frac{z f'(z)}{f(z)}\, dz.
$$
Indeed, the poles of $zf'(z)/f(z)$ are the zeros and poles of $f$, and if $f$ has a zero of
multiplicity $k$ at $z=a$, then $f(z)=(z-a)^kg(z)$ near $z=a$, with
$g$ analytic  and $g(a)\neq 0$. Then
$$\begin{aligned}
\frac{z f'(z)}{f(z)}&
=\frac{z}{(z-a)^kg(z)}(k(z-a)^{k-1}g(z)+(z-a)^kg'(z))\\
&=\frac{kz}{(z-a)}+ \frac{zg'(z)}{g(z)}
\end{aligned}
$$
near $z=a$, with $zg'/g$ analytic at $a$, so $zf'/f$ has  residue
$ka$ at $z=a$. Similarly, if $f$  has a pole  of multiplicity $l$ at
$z=b$, then $ zf'/f$ has residue $-lb$ at  $z=b$. Now the  zeros
and  poles of $f$ within $\mathcal R $ are $a_1, \ldots, a_r$ and
$b_1, \ldots, b_s$  with multiplicities $k_1, \ldots, k_r$ and $b_1,
\ldots, b_s$, so
$$\frac{1}{2\pi i} \int_{\partial{\mathcal R}} \frac{z
f'(z)}{f(z)}dz,$$ which is equal  to the sum of the residue of $z
f'/f$ , takes the value 
$$ 
\sum_{j=1}^r k_ja_j - \sum_{j=1}^s l_jb_j.
$$
We label  the sides of $\mathcal{R} $ as in the proof of
Theorem~\ref{t9jones}. Since $-\Gamma_4$, the path $\Ga_4$ with reverse orientation, is just the path $\Gamma_2-w_1$, we get
$$
\begin{aligned}
\int_{\Gamma_2} \frac{z f'(z)}{f(z)}\, dz
&=\int_{\Gamma_2} \frac{(z-w_1) f'(z)}{f(z)}\, dz
  + \int_{\Gamma_2} \frac{(w_1)f'(z)}{f(z)}\, dz\\
&=\int_{-\Gamma_4+w_1}\frac{(z-w_1) f'(z)}{f(z)}\, dz+w_1[\log
f(z)]_{\Gamma_2}
\end{aligned}
$$
Now we will calculate both summands separately. First, since $w_1$ is a period of both $f$ and $f'$, we get:
$$
\begin{aligned}
\int_{-\Gamma_4+w_1}\frac{(z-w_1) f'(z)}{f(z)}\, dz
=&\int_{-\Gamma_4} \frac{z f'(z+w_1)}{f(z+w_1)}\, dz \\
=&-\int_{\Gamma_4} \frac{z f'(z+w_1)}{f(z+w_1)}\, dz
=&-\int_{\Gamma_4} \frac{z f'(z)}{f(z)}\, dz.
\end{aligned}
$$
Secondly,
$$
w_1[\log f(z)]_{\Gamma_2}=2\pi n_1i w_1
$$
for some $n_1\in \mathbb{Z}$ since $\log f(z)$ changes its value by an integer
multiple of $ 2\pi i$ as  $z$ travels along $\Gamma_2$ from
$t+w_1$ to $t+w_1+w_2$ (as $f(t+w_1)=f(t+w_1+w_2))$. Thus,
$$
\int_{\Gamma_2} \frac{z f'(z)}{f(z)}\, dz
=-\int_{\Gamma_4} \frac{z f'(z)}{f(z)}\, dz + 2\pi n_1i w_1.
$$
Similarly,
$$
\int_{\Gamma_1} \frac{z f'(z)}{f(z)}dz
=-\int_{\Gamma_3} \frac{zf'(z)}{f(z)}dz +2\pi n_2i w_2 
$$ 
for some $n_2\in \mathbb{Z}$. Hence,
$$
\begin{aligned}
\sum_{j=1}^r k_ja_j - \sum_{j=1}^s l_jb_j
&=\frac{1}{2\pi i}
\int_{\partial{\mathcal R}} \frac{z f'(z)}{f(z)}\, dz
=\frac{1}{2\pi i}\sum_{j=1}^4 \int_{\Gamma_2} \frac{zf'(z)}{f(z)}\, dz\\
&= \frac{1}{2\pi i}(2\pi n_1i w_1  +2\pi n_2i w_2)\\
&=n_1i w_1  + n_2i w_2,
\end{aligned}
$$
which is an element  of $\La_f$ as required. The proof is complete.
\endpf

\

When we come to construct  elliptic functions, we shall  see that
the   conditions (1) and (2') together with the  condition 

\sp\begin{itemize}
\item[(3)] \
$\sum_{j=1}^rk_ja_j\sim_\La \sum_{j=1}^sl_jb_j$,
\end{itemize}

\sp\fr are sufficient for the
existence  of an elliptic function $f$ with prescribed  zeros and poles.

\sp

\section{Weierstrass $\wp$--Functions I}\label{weierstrass}

Let $\La=[w_1, w_2]$  be a lattice with a basis $\{w_1,w_2\}$ and
let $\mathcal R$  be a fundamental parallelogram for $\La$ with no
elements of $\La$ on $\partial{\mathcal R}$. 

Our goal in this section is to construct
non--constant functions $f$  which are elliptic with respect to
$\La$. By Theorem~\ref{t8jones}  we know that such a function cannot
be analytic and so it must have poles in $\mathcal R$.  By Corollary~\ref{c1jones} we know that $f$ cannot just have one simple pole in $\mathcal R$.
So the simplest, potential, non-constant elliptic functions must have order 2, with
either two simple poles or else  a single pole of order 2. 

\sp In this section we shall introduce the
Weierstrass $\wp(z)$--function which is elliptic of order  2 with
respect to $\La$ and has a single pole of order 2 in $\mathcal R$.
This will be our basic elliptic function in the sense that every
function which is elliptic with respect to $\La$ is a rational
function of $\wp$ and its  derivative $\wp'(z)$ (see
Theorem~\ref{t21jones}).

The sets 
$$
\Omega_r:=\big\{aw_1+bw_2: a, b \in \mathbb R\, \, \,
\mbox{ and}\, \, \,   \max(|a|, |b|=r)\big\},
$$ 
for integers $r\geq 1$, are similar parallelograms centered on 0. Defining 
$$
\La_r:=  \La\cap\Omega_r,
$$ 
we have
$$
\La_r=\{mw_1+nw_2: \,\,\, m, n \in \mathbb Z\, \, \,  \mbox{
and}\,\,\,  \max(|m|, |n|=r)\}.
$$ 
Now $\La$ is a disjoint union
$$
\La=\{0\}\cup \La_1\cup \La_2\cup\La_2\cup \ldots,
$$
and for each $r\geq 1$ we
have
\beq\label{ej6}
\#\La_r=8r.
\eeq
  We  can order the elements of $\La$ by starting at 0 and then
  listing the elements of $\La_1, \La_2, \ldots$ in turn, rotating
  around each $\La_r$ in the order $rw_1, rw_1+w_2, \ldots,
  rw_1-w_2$ so that the sequence  spirals  outwards from 0. If we
denote this ordering by $w^{(0)}, w^{(1)}, w^{(2)}, \ldots $ then
$$
w^{(0)}=0, \ w^{(1)}=w_1, \ w^{(2)}=w_1+w_2, \  w^{(3)}=w_2,\ldots, \
w^{(8)}=w_1-w_2,  w^{(9)}=2w_1, \ 
$$
and
$$
w^{(10)}=2w_1+w_2, \ldots.
$$ 
Clearly
$\lim_{k \to+\infty}|w^{(k)}|=+\infty$. By 
$$
\sum_{w\in \La} \  \  \text{ and }  \  \  \sum_{w\in \La}\!\!{}'
$$ 
we shall, naturally, mean the sum over all (respectively all non--zero) lattice-points $w$ taken in the above order. Thus 
$$
\sum_{w\in \La}h(w)=\sum_{k=0}^\infty h(w^{(k)})
$$
for any function $h$, and similarly, 
$$
\sum_{w\in\La}\!\!{}'h(w)=\sum_{k=1}^\infty h(w^{(k)}).
$$
By 
$$
\prod_{w\in \La} \  \  \text{ and }  \  \  \prod_{w\in \La}\!\!{}'
$$ 
we shall mean the product over all, respectively all non-zero, lattice-points, again, in the above order. For convenience we will  abbreviate the  notation  to $\sum, \sum^{'}$, etc., the  particular lattice $\La$ being understood. In
practise, the particular ordering of $\La$ will not often be
important, as  the sums, and products which concern us are  usually
absolutely convergent and hence  invariant under  rearrangements.

 \sp The convergence  properties  of the Weierstrass functions
depend  on the following result, which is  a 2-dimensional  analogue
of the fact that the series  $\sum_{n=1}^\infty n^{-s}$ defining the
Riemann zeta--function converges if and only if $s>1$. Its proof is straightforward and left for the reader.

\sp

\bthm\label{t14jones}  If $s\in \mathbb R $, then $\sum_{w\in \La}'
|w|^{-s}$ converges if and only if $s>2$. 
\ethm  

\sp

It is now easy to construct elliptic function of each order $N\geq
3$.

\sp

\bthm\label{t15jones} For each integer $N\geq 3$, the function
$$
\C\ni z\lmt F_N(z):= \sum_{w\in \La} ( z-w)^{-N}\in\oc
$$ 
is elliptic of order $N$ with respect to the lattice $\La$. \ethm

\bpf If $K$ is any compact  subset of $\mathbb C\sms
\La$ then the terms $(z-w)^{-N}$ are  analytic and therefore
bounded on $K$. Since $K$ is bounded, there exists a set $\Phi\sbt \La$ such that $\La\sms\Phi$ is finite and
$$
|w|\geq 2|z|
$$
for all $w\in\Phi$ and all $z\in K$. It follows that 
$$
|z-w|\geq |w|-|z|\geq \frac12|w|
$$ 
for all $z\in K$ and all $w\in \Phi$. Hence
$$
|(z-w)^{-N}|\leq 2^N|w|^{-N}
$$ 
for all $w\in \Phi$ and all $z\in K$. So, ff $N\geq 3$, then Theorem~\ref{t14jones}  and the comparison test imply that the series
$$
\sum_{w\in \La}(z-w)^{-N}
$$
converges absolutely uniformly on $K$. Since each term $(z-w)^{-N}$ is analytic   on $K$, this implies that $F_N(z)$ is analytic on $\mathbb C \sms \La$. Fix an element $\xi\in \La$. Since $\La$ is closed and discrete, there exists $r>0$ so small that $B(\xi,r)\cap B(\La\sms\{\xi\},r)=\es$. Then 
$$
F_N(z)=(z-\xi)^N+\sum_{w\in\La\sms\{\xi\}}(z-w)^{-N}
$$
The first term of this sum is meromorphic on $B(\xi,r)$ while the second one can be shown as above to be absolutely uniformly convergent on $B(\xi,r)$. Thus the function $F_N$ restricted to $B(\xi,r)$ is meromorhic as a sum of meromorphic and analytic function. Therefore the function $F_N$ is meromorphic on $\mathbb C$. Since, as we have shown, the series defining $F_N$ converges absolutely at each point of $\C\sms\La$, we can rearrange its terms arbitrarily and get for every $\xi\in \La$ that
$$
F_N(z+\xi)
 = \sum_{w\in \La} ((z+\xi)-w)^{-N}
 =\sum_{w'\in \La} (z+w')^{-N}
= F_N(z),
$$
as $w'=w\sms\xi$ ranges over $\La$ as $w$ does (though in
different order). The proof is complete. 
\endpf

\sp Clearly, this method fails to produce an elliptic  function $F_2(z)
$ of order 2, since Theorem~\ref{t14jones} cannot be used to prove
convergence of 
$$
\sum_{w\in \La} (z-w)^{-2}.
$$
In order to guarantee convergence, we make the terms of this series smaller, replacing $(z-w)^2$ by
 $(z-w)^2+w^{-2}$ for each $w\neq 0$. Thus the resulting series is: \index{(S)}{$\wp$}
\beq\label{ej7}
\wp=\wp_{\La} :=\frac{1}{z^2}+\sum_{w\in \La}^{}\!\!{'}\left(\frac{1}{(z-w)^2}-\frac{1}{w^2}\right).
\eeq

\sp

\bthm\label{weierstrass_convergence}
Given a lattice $\La\sbt\C$ the series $\wp_\La$ in formula \eqref{ej7} defines a meromorphic function from $\C$ to $\oc$. Moreover, this function 

\sp\begin{itemize}

\item is periodic with respect to the lattice $\La$; in fact, the set of periods of $\wp_\La$ coincides with $\La$, i.e. $\La_\wp=\La$.

\sp\item has poles of order $2$ at all elements of $\La$ and

\sp\item is holomorphic in $\C\sms\La$. 
\end{itemize}

\sp\fr This function is called the Weierstrass elliptic function \index{(N)}{ Weierstrass elliptic function} induced by the lattice $\La$. If we want to stress the importance of the lattice, we denote it by $\wp_\La$ \index{(S)}{$\wp_\La$} rather than merely by $\wp$ \index{(S)}{$\wp$}.
\ethm

\bpf
For every $w\in\La$, we have 
$$
\lt|\frac{1}{(z-w)^2}-\frac{1}{w^2}\rt|
=\lt|\frac{z(z-2w)}{(z-w)^2w^2}\rt|
$$
So, if $|w|\ge 2|z|$, then
$$
\lt|\frac{1}{(z-w)^2}-\frac{1}{w^2}\rt|
=\lt|\frac{z(z-2w)}{(z-w)^2w^2}\rt|
\le |z|\frac{3|w|}{|w|^4}
=3|z|\cdot|w|^{-3}.
$$
It therefore follows from  Theorem~\ref{t14jones} that for each point $z\in\C$, there exists $r>0$ so small that the series of \eqref{ej7} can be represented on $B(z,r)$ as a sum of its finitely many terms, which form a meromorphic function (even on $\C$) and an absolutely uniformly convergent series of analytic functions. Thus, the formula \eqref{ej7} defines a meromorphic function in $\C$. Checking periodicity of $\wp$, we rearange terms to get for all $z\in\C\sms\La$ and all $\xi\in\La$ that
$$
\begin{aligned}
\wp_\La(z+\xi)
&=\frac{1}{(z+\xi)^2}+\sum_{w\in \La}^{}\!\!{'}\left(\frac{1}{((z+\xi)-w)^2}
        -\frac{1}{w^2}\right) \\
&=\frac{1}{(z+\xi)^2}+\sum_{w\in \La\sms\{\xi\}}^{}\!\!\!\!\!\!\!{'}\,\left(\frac{1}
         {(z-(w-\xi))^2}-\frac{1}{(w-\xi)^2}\right)+\frac1{z^2}-\frac1{(-\xi)^2} \\
&=\frac1{z^2}+\frac{1}{(z+\xi)^2}+\!\!\!\sum_{w'\in \La}^{}\!\!{'}^{}\!\!{'}\left(\!\frac{1}{(z-w')^2}-\frac{1}{(w')^2}\!\right)-\frac{1}  {(z+\xi)^2}+\frac1{(-\xi)^2}-\frac1{(-\xi)^2} \\
&=\frac1{z^2}+\sum_{w'\in \La}^{}\!\!{'}^{}\!\!{'}\left(\frac{1}{(z-w')^2}-\frac{1}{(w')^2}\rt)\\
&=\wp_\La(z).     
\end{aligned}
$$
This means that 
$$
\La\sbt\La_\wp.
$$
The statement about poles and analyticity in $\C\sms\La$ is obvious; and this has also further consequences. Namely, since $0$ is a pole and $0\in\La_\wp$, we see that each element of $\La_\wp$ is a pole of $\wp$, i.e. $\La_\wp\sbt\wp^{-1}(\infty)$. But we already know that $\wp^{-1}(\infty)\sbt\La$. So, $\La_\wp\sbt\La$, and the equality $\La_\wp=\La$ is established. The proof is complete.
\epf

\sp

\bthm\label{t17jones} 
Any Weierstrass $\wp$--function has order 2, and its derivative $\wp'$ has
order 3.
\ethm

\bpf  $\wp_\La$ has a single congruence class of poles, namely the lattice $\La$. And since each element pf $\La$ is of multiplicity $2$, the function $\wp$ has order $2$. Similarly, $\wp_\La'=-2F_3$ has a single class of poles all of multiplicity $3$. Thus, $\wp_\La'$ has order 3. 
\endpf

\sp

We now shall derive an important equation connecting $\wp_\La(z)$ and $\wp_\La'(z)$ obtained from the Laurent series for $\wp_\La(z)$ near $z=0$. We start
with finding the Laurent series for the function
\beq\label{ej11}
\zeta(z)
:=\frac{1}{z}+ \sum_{w\in\La}^{}\!\!{'}\left(\frac{1}{z-w}+\frac{1}{w} + \frac{z}{w^2}\right).
\eeq
Let 
$$
m=\min\big\{|w|:w \in \La\sms \{0\}\big\}>0,
$$
Since
$$
\frac{1}{z-w}+\frac{1}{w} + \frac{z}{w^2}=\frac{z^2}{w^2(z-w)},
$$
we see, by comparison with  $\sum'|w|^{-3}$ that the series
$$ 
\sum_{w\in \La}^{}\!\!{'}\,\frac1{z-w}+\frac1w + \frac{z}{w^2}
$$ 
is absolutely convergent for each $z\in \mathbb
C\sms \La$. Moreover, for each $w\in \La \sms \{0\}$ the binomial
series
$$
\frac{1}{z-w}=-\frac{1}{w}-\frac{z}{w^2}-\frac{z^2}{w^3}-\ldots
$$
is absolutely convergent for $z \in B(0,m)$, so we may substitute this in
(\ref{ej11})  and  reverse  the order of summation  to  obtain
$$
\zeta(z)
=\frac{1}{z}+\sum_{w\in \La}^{}\!{'}\,\left(-\frac{z^2}{w^3}-\frac{z^3}{w^4}-\ldots\right)
=\frac{1}{z} - G_3z^2-G_4z^3-\ldots
$$
for all $z\in B(0,m)$, where
$$
G_k=G_k(\La):=\sum_{w\in \La}^{}\!\!{'}\, w^{-k}.
$$
These series $G_k$, $k\ge 3$, are called the Eisenstein series for $\La$; these are absolutely convergent by Theorem~\ref{t14jones}). For odd integers $k\ge 3$, the
terms $w^{-k}$ and $(-w)^{-k}$ cancel each other out, giving $G_k=0$. So, the
Laurent series of $\zeta$ near $0$ is
\beq\label{ej12}
    \zeta(z)=\frac{1}{z}- \sum_{n=2}^\infty G_{2n}z^{2n-1}.
\eeq
Hence,
\beq\label{ej13}
    \wp_\La(z)=-\zeta'(z)=\frac{1}{z^2}+\sum_{n=2}^\infty (2n-1) G_{2n}z^{2n-2}.
\eeq
This is thus the Laurent series for $\wp_\La(z)$, valid for all $z\in B(0,m)$. Now a
straightforward calculations give
 \begin{eqnarray*}
\wp_\La'(z)&=&-\frac{2}{z^3}+6G_4z +20G_6z^3+\ldots
\end{eqnarray*}
and
\begin{eqnarray*}
\wp_\La'(z)^2
&=&\frac{4}{z^6}-\frac{24G_4}{z^2}-80G_6 + z^2\phi_1(z) \\
4\wp_\La(z)^3 
&=& \frac{4}{z^6}+\frac{36G_4}{z^2}+60G_6z + z^2\phi_2(z)  \\
60G_4\wp_\La(z) 
&=&\frac{60G_4}{z^2} +z^2\phi_3(z ),
\end{eqnarray*}
where  $\phi_1(z), \phi_2(z)$ are some power series convergent in $B(0,m)$.
These  last three equations give
$$
\wp'(z)^2-4\wp_\La(z)^3+60G_4\wp_\La(z)+140G_6=z^2\phi(z),
$$
where $\phi(z)=\phi_1(z)- \phi_2(z)+\phi_3(z)$  is some power series
convergent in $B(0,m)$. As $\wp_\La$ and $\wp_\La'$ are periodic with respect
 to $\La$, the function
$$ 
f(z)=\wp_\La'(z)^2-4\wp_\La(z)^3+60G_4\wp_\La(z)+140G_6
$$ 
is also periodic with respect to $\La$.
Since $f(z)=z^2\phi(z)$ in $D$, with $\phi$ being analytic, the function $f$ vanishes
at $0$, and hence at all $w\in \La$. However, by its construction $f$
can have poles  only where $\wp$ or $\wp'$ have poles, that is,  at
the lattice--points $\La$. Therefore $f$ has no poles  and so, by
Theorem~\ref{t8jones} $f(z)$ is constant, which must be zero, since
$f(z)=0$. Thus we have proved the following.

\sp

\bthm~\label{t18jones}
$\wp_\La'(z)^2=4\wp_\La(z)^3-60 G_4 \wp_\La(z)-140 G_6$.
\ethm

\sp

\fr This is the  differential equation for $\wp(z)$. It is
customary to define
\beq\label{ej14a}
g_2=g_2(\La):= 60G_4= 60\sum_{w\in \La}^{}\!\!{'}\,w^{-4}
 \eeq\index{(S)}{$g_2$}
and 
\beq\label{ej14b}
g_3=g_3(\La):=140 G_6= 140\sum_{w\in \La}^{}\!\!{'}\, w^{-6},
\eeq\index{(S)}{$g_3$}
so that
\beq\label{ej15}
\wp_\La'(z)^2=4\wp_\La(z)^3-g_2 \wp_\La(z)-g_3.
\eeq
If we put $z=\wp(t)$ then this means that
$$
\left(\frac{dz}{dt}\right)^2=4z^3-g_2z-g_3
$$
So, locally, inverse branches of $\wp$-functions take on the form
$$ 
\wp^{-1}(z)=t=\int\frac{dz}{\sqrt{p(z)}},
$$ 
where $p(z)$ is the cubic polynomial $4z^3-g_2z-g_3$. This shows how the local inverses of elliptic  functions appear as indefinite integrals.

\sp\bthm~\label{t19jones} 
Let $\La=[w_1,w_2]$ and let $w_3=w_1+w_2$. If $\mathcal{R}$ is a fundamental
parallelogram for $\La$ with 
$$
\lt\{0, \frac{1}{2}w_1, \frac{1}{2}w_2, \frac{1}{2}w_3\rt\}\sbt \Int \mathcal{R}
$$ 
then 
$$
\frac{1}{2}w_1,\, \frac{1}{2}w_2, \, \frac{1}{2}w_3
$$ 
are the only critical points of $\wp_\La$ in $\mathcal R$, i.e. the only zeros of $\wp_\La'$ in $\mathcal R$. 
\ethm

\bpf By Theorem~\ref{t17jones}, $\wp_\La'(z)$ has order 3,
and hence has exactly three (counting with multiplicities) zeros in  $\mathcal{R}$.  
If $w\in \mathcal{R}$ then, because $\frac{1}{2}w\sim_\La -\frac{1}{2}w$, we have
$$
\wp_\La'\lt(\frac{1}{2}w\rt)= \wp_\La'\lt(-\frac{1}{2}w\rt).
$$
Since in addition $\wp_\La'(z)$ is an odd
function, we have $\wp_\La'(-\frac{1}{2}w)=- \wp_\La'(-\frac{1}{2}w)$, and hence
$\wp_\La'(-\frac{1}{2}w)=0$. Since also $\wp_\La^{-1}(\infty)=\La$ and $\{w_1/2,w_2/2,w_3/2\}\cap \La=\es$, we must therefore have that 
$$
\wp_\La'\lt(-\frac{1}{2}w_i\rt)=0
$$
for all $j=1,2,3$. The proof is complete.
\endpf

\sp We define
\beq\label{20200130}
e_j:=\wp_\La(w_j/2), \  \  j=1,2,3,
\eeq
and, as an immediate consequence of Theorem~\ref{t19jones}, we get the following.

\sp\bthm\label{t120200129}
If  $\La=[w_1,w_2]$ and $w_3=w_1+w_2$, then 
$$
\wp_\La\(\Crit(\wp_\La)\)=\{e_1,e_2,e_3\}=
\big\{\wp_\La(w_1/2), \wp_\La(w_2/2),\wp_\La(w_3/2)\big\}
$$
\ethm

\bcor\label{c2jones} For each  $c\in \ov{\mathbb C}\sms \{e_1, e_2,
e_3\}$ the equation $\wp(z)=c$ has two simple solutions while for $c= e_1,
e_2, e_3$ or $\infty$ this equation has one double solution. 
\ecor

\bpf  Since $\wp_\La$ is elliptic  of order 2, it takes each
value $c\in \ov{\mathbb C}$  twice by Theorem~\ref{t10jones} giving
either two simple solutions (some $z$ and $-z$  since $\wp$ is even) or one double solution. If
$c\in \mathbb C$ then $\wp_\La(z)=c$ has double solution if and only if
$\wp_\La'(z)=0$ giving $ z \sim_\La\frac{1}{2}w_j$, $ j=1,2,3$. So,
$c=e_j$. The pole  of order 2 at $z=0$
shows that  $\wp(z)=\infty$ has double solution. 
\endpf

\sp

\bthm\label{t20jones}  The complex numbers $e_1, e_2$ and  $e_3$ are mutually  distinct.
\ethm

\bpf  Let 
$$
f_j(z)=\wp_\La(z)-e_j
$$
for $j=1,2,3$. Like $\wp_\La$, the functions $f_j$ are of order $2$. As 
$$
f_j(w_j/2)=f_j'(w_j/2)=0,
$$
the function$f_j$ has double zeros on the $\La$ equivalence class $[w/2]$, and hence has  no others zeros. In particular, $f_j(w_k/2)\neq 0$ for $j\neq k$. Since $f_j(w_k/2)=\wp(w_k/2)-e_j=e_k-e_j$, it follows that $e_j\neq
e_k$ for $j\neq k$. 
\endpf

\sp We want to end this section by noting that by (\ref{ej15}), the polynomial
$$ 
p(z)=4z^3-g_2z-g_3
$$
has zero at $z=\wp_\La(t)$ where $\wp_\La'(t)=0$, so $p(z)$  has three
distinct zeros $z=e_1, e_2$ and $e_3$.

\sp\section{The Field of Elliptic Functions}\label{field}

In this section we consider a fixed lattice $\La$. An elliptic
function will mean a function which is elliptic with respect to
$\La$. If $f$ and $g$ are elliptic, then so are $f+g$, $f-g$ and
$fg$, and if $g$ is not  identically zero, then $1/g$ is elliptic.
Thus the set of all elliptic functions  is a field, which  we shall
denote by $\mathcal{E}(\La)$.\index{(S)}{$\mathcal{E}(\La)$} This
field contains a subfield
$\mathcal{E}_1(\La)$\index{(S)}{$\mathcal{E_1}(\La)$} consisting of
even elliptic functions. The constant functions form a subfield of
$\mathcal{E}_1(\La)$ isomorphic to $\mathbb C$, so we may regard
$\mathcal{E}(\La)$ and $\mathcal{E}_1(\La)$ as an extension field of
$\mathbb C$. Since $\mathcal{E}_1(\La)$  contains
$\wp(z):=\wp_\La(z)$, it contains all rational functions of $\wp$
with complex coefficients; these rational functions form a field
${\mathbb C}(\wp)$\index{(S)}{${\mathbb C}(\wp)$}, the smallest
field containing $\wp$ and the constant functions $\mathbb C$.
Similarly, $\mathcal{E}(\La)$ contains $\wp$ and $\wp'$, and hence
contains the field $\mathbb C(\wp, \wp')$\index{(S)}{${\mathbb
C}(\wp, \wp')$} of rational functions of $\wp$ and $\wp'$. This is
the smallest field containing $\wp, \wp'$ and $\mathbb C$. We start with the main result of this section.

\bthm\label{t21jones}
Let $\La\sbt\C$ be a lattice. Then

\sp\begin{enumerate}
\item If $f$ is an even elliptic function with respect to $\La$, then $f=R_1(\wp_\La)$ for some rational  function $R_1$; thus ${\mathcal E}_1(\La)=\mathbb C(\wp_\La)$.

\sp\item If $f$ is any elliptic function with respect to $\La$, then
$$
f=R_1(\wp_\La)+\wp'R_2(\wp_\La),
$$
where $R_1$ and $R_2$ are rational functions; thus ${\mathcal
E}(\La)=\mathbb C(\wp_\La, \wp_\La')$.
\end{enumerate}
 \ethm

 \bpf (1) Let $f$ be  an even elliptic function. The
 result is obvious  for  constant  functions, so suppose  that $f$
 has order $N>0$. If $k \in \mathbb C$, then $f(z)=k$ has  multiple
 roots only where $f'(z)=0$, and  this occurs  at only finitely many
  congruence  classes of points $z$. Thus $f(z)=k$ has its roots
  simple  for all but  finitely many  values of $k$. We  can
  therefore  choose two distinct  complex  numbers $c$ and $d$ so
  that the roots of $f(z)=c$ and of $f(z)=d$ are all simple, and so
  that  none of these roots are congruent  to 0 or $w_j/2$,
  $j=1,2,3$.  Since  $f$ is even, a complete set of roots  of
  $f(z)=c$ will have the form $\{a_1, -a_1, \ldots, a_n, -a_n\}$, these
  being  simple and  mutually  non-congruent, and similarly for the
  roots $\{b_1, -b_1, \ldots, b_n, -b_n\}$ of $f(z)=d$. Hence the
  elliptic  function
  $$ 
  g(z):=\frac{f(z)-c}{f(z)-d}
  $$ 
  has  simple zeros at $a_1, -a_1, \ldots, a_n,
  -a_n$ and simple poles $b_1, -b_1, \ldots, b_n, -b_n$.

 Now, Corollary~\ref{c2jones} implies that the equations $\wp(z)=\wp(a_i)$
    and $\wp(z)=\wp(b_i)$ have simple roots respectively at $z=\pm a_i$  and
    $z=\pm b_i$ for every $1\leq i \leq n$. So, the elliptic
    function
$$
 h(z):=\frac{(\wp(z)-\wp(a_1))(\wp(z)-\wp(a_2))\ldots
(\wp(z)-\wp(a_n))}{(\wp(z)-\wp(b_1))(\wp(z)-\wp(b_2))\ldots
(\wp(z)-\wp(b_n))}
$$ 
has the  same zeros and poles  as $g$, with the
same  multiplicities (all simple).  Hence Theorem~\ref{t12jones}
implies that $g=\mu h$ for some   constant $\mu\neq 0$. Solving
$$ \frac{f(z)-c}{f(z)-d}= \mu\frac{(\wp(z)-\wp(a_1))(\wp(z)-\wp(a_2))\ldots
(\wp(z)-\wp(a_n))}{(\wp(z)-\wp(b_1))(\wp(z)-\wp(b_2))\ldots
(\wp(z)-\wp(b_n))}$$ for $f(z)$, we see that  $f$ is rational
function with complex  coefficients $R_1(\wp)$ of $\wp$.

\sp (2)  If $f$ is odd, then $f/\wp'$ is even, so by (i) we have
$f=\wp'R_2(\wp)$ for some rational function $R_2$. In general, if
$f$ is any elliptic  functions then
$$ 
f(z)=\frac{1}{2}( f(z)+ f(-z)) +  \frac{1}{2}( f(z)-f(-z)),
$$
where $\frac{1}{2}( f(z)+ f(-z))$ is even  and elliptic  while
$\frac{1}{2}( f(z)- f(-z))$ is odd and elliptc, so by  the above
arguments we have
$f=R_1(\wp)+\wp'R_2(\wp),$ where  $R_1$ and $R_2$ are some rational functions. The proof is complete.
\endpf

\sp Using the differential equation $(\wp')^2= 4\wp^3-g_2\wp- g_3$, we
can reduce any rational function of $ \wp $ and $\wp'$ to  the form
$R_1( \wp)+\wp'R_2(\wp)$  by eliminating powers of $\wp '$; for
example
$$
\frac{\wp \wp'}{\wp +1}= \frac{\wp \wp'( \wp'-1)}{(\wp')^2
-1}=\frac{(4\wp^4 -g_2\wp^2- g_3\wp)-\wp' \wp }{4\wp^3 -g_2\wp-
g_3\wp-1} .$$ We can  view Theorem~\ref{t18jones} as an algebraic
equation between the functions $\wp$ and $\wp'$. We  now show that
any  two functions in $\mathcal E(\La)$ are connected by an
algebraic equation.

\sp

\bthm\label{t21ajones} 
If $f, g \in \mathcal E({\La})$  then there
exists a non-zero irreducible polynomial $\phi(x,y)$ with complex
coefficients, such that  $\phi(f,g)$ is identically  zero. 
\ethm

\bpf  If we choose  any polynomial in two variables $x,y$, say
$$
F(x,y)=\sum_{k=1}^m\sum_{l=1}^n \a_{kl}x^ky^l,\quad \a_{kl}\in\mathbb C,
$$ 
then the function 
$$
h(z)=F(f(z), g(z))
$$ 
is elliptic function with poles only at the poles  of $f$ or $g$. If $f$ and $g$
 have $M$ and $N$ poles  respectively, then  $h$ has   at most
 $mM+nN$  poles counting with multiplicities in each case.
 Therefore,  unless  $h$ is identically  zero, it has at most
 $mM+nN$ zeros by Theorem~\ref{t10jones}. We now show that  if $m$ and $n$ are
 large enough  then we  can choose  the coefficients $\a_{kl}$, so
 that   $h$ has more than  $mM+nN$ zeros and  hence
 $h(z)\equiv 0$.

    To do this, we let $z_1, \ldots, z_{mn-1}$ be $mn-1$ non-congruent
    points distinct from the poles of $f$ and $g$. Now we regard
\beq\label{ej15a}
h(z_j)=\sum_{k=1}^m \sum_{l=1}^n \a_{kl}f(z_j)^k g(z_j)^l=0,
     \,\, \  j=1, \ldots, mn-1
\eeq
as  a set of $mn-1$ homogeneous  linear  equations in the $mn$
unknowns $\a_{kl}$. As there are more unknowns than  equations, this
set  of equations has a non-trivial solution, that  is, there
exists $\a_{kl}$, not all zero, satisfying (\ref{ej15a}). Thus  $f(x, y)$ is
not identically zero, but 
$$
F(f(z), g(z))=h(z)=0
$$ 
at the points $z=z_1, \ldots, z_{mn-1}$. Now  for $m,n$ large  enough, $mn-1>mM+nN$, and so by choosing the coefficients $\a_{kl}$  as above we must  have
$h(z)$ identically equal to $0$, that is, $F(f,g)=0$. We  can factorize $F(x,y)$  within the polynomial ring $\mathbb C[x, y]$ as product
$$ F(x,y)=F_1(x,y)F_2(x,y)\ldots F_r(x,y)$$
of irreducible  polynomials. Thus $F_1(f,g)F_2(f,g)\ldots
F_r(f,g)=0$ within the  field $\mathcal E(\La)$, so some
$F_i(f,g)=0$, and we can take $\phi$ to be $F_i$. 
\endpf

\sp

The Weierstrass $\sigma$--function\index{(N)}{Weierstrass $\sigma$-function} is defined 
\beq\label{ej8}
\begin{aligned}
\sg(z)  :&= z\prod_{w\in \La}^{}\!\!{'}\,g(w,z), \  \text{ where } \\
 g(w,z) &:=\left(1-\frac{z}{w}\right) \exp\left(\frac{z}{w} + \frac{1}{2} \left(
 \frac{z}{w}\right)^2\right).
\end{aligned}
\eeq\index{(S)}{$ \sg(z)$}
The factor $(1-(z/w))$ is introduced in $g(w,z)$ to give $g(z)$ a
simple zero at each lattice-point $w\in \La$, while the exponential factor
is included  to guarantee convergence of the  infinite product. In order to see this convergence, consider $\in\La$ with large moduli, denote by $\log_0$ the holomorphic branch of logarithm defined on a neighborhood of $1$ and sending $1$ to $0$. Define then $\log_*$, a holomorphic branch of $\log g$ as follows
$$
\log_*g(w,z):=\log_0\left(1-\frac{z}{w}\right)+\frac{z}{w} + \frac{1}{2} \left(
 \frac{z}{w}\right)^2.
$$
Expanding now $\log_0\left(1-\frac{z}{w}\right)$ into its Taylor series near zero, we obtain that
$$
\log_*g(w,z)= \frac13\lt(\frac{z}{\om}\rt)^3 + \  \text{ higher terms}.
$$
This guaranties the convergence of the product defining $\sg(z)$. We have that
$$
\frac{\sg'(z)}{\sg(z)}
=\frac{d}{dz}\log(\sg(z))
=\frac1z+\sum_{w\in \La}^{}\!\!{'}\,\lt(\frac1{z-w}+\frac1w+\frac{z}{w^2}\rt)
=\zeta(z).
$$
After differentiating (term by term), we recover formula \eqref{ej13} from the previous section:
$$
\zeta'(z)=-\wp(z).
$$
The functions $\zeta$ and $\sigma$ which were just introduced are not elliptic, for,  as we shall show, that they are not invariant under the translations $\C\ni z\mapsto z+w\in\C$, $w \in
\La$. However, an examination of their  behavior  under translations will enable us to construct elliptic functions with prescribed properties.

\sp Since $\zeta'(z)=-\wp(z)$ we  have $\zeta'(z+w_j)=\zeta'(z)$ for
$j=1,2$, the  integration with respect to $z$ gives
$$
\zeta(z+w_j)= \zeta(z)+\eta_j, \, j=1,2,
$$
where $\eta_1,\eta_2$ are constants independent of $z$.  If
$w\in\La$ then $w=mw_1+nw_2$, where $m,n\in \mathbb Z$, and hence
\beq\label{ej16}
  \zeta(z+w)=\zeta(z)+\eta
\eeq
where
\beq\label{ej17}
\eta=m\eta_1+n\eta_2.
\eeq
Let $\mathcal R$ be a fundamental parallelogram for $\La$,
containing $0$ in its interior, with vertices $\xi$, $\xi+w_1$,
$\xi+w_1+w_2,$ $\xi+w_2$, $\xi\in\C$. Denote its edges by $\Gamma_1$, $\Gamma_2$, $\Gamma_3$, $\Gamma_4$, with anticlockwise orientation.
Since the function $\zeta:\C\lra\oc$ is
meromorphic and has a single pole in $\mathcal R$ at 0 with
residue equal to $1$, we have
$$ 
2\pi i 
= \int_{\partial{\mathcal R}} \zeta(z)\, dz
= \sum_{j=1}^4\int_{\Gamma_j}\zeta(z)\, dz.
$$ 
Now 
$$
\int_{\Gamma_3} \zeta(z)\, dz
= -\int_{\Gamma_1}(\zeta(z+w_2)\, dz
=-\int_{\Gamma_1}\(\zeta(z)+\eta_2\)\, dz
$$
So,
$$
\int_{\Gamma_1} \zeta(z)\, dz+ \int_{\Gamma_3} \zeta(z)\,dz
=\int_{\Gamma_1} \eta_2\, dz= -\eta_2 w_1.
$$
Similarly,
$$
\int_{\Gamma_2} \zeta(z)dz+\int_{\Gamma_4} \zeta(z)\, dz
=\eta_1 w_2,
$$
and hence
\beq\label{ej18}
\eta_1w_2-\eta_2w_1=2\pi i.
\eeq
This equation is usually  referred  to as  Legendre'a relation. It
implies that at least one of $\eta_1, \eta_2$ is non-zero, so
$\zeta(z)$ is not elliptic. To see how $\sigma(z)$ behaves under
translations, we use the formula
$$
\frac{\sigma'(z)}{\sigma(z)}=\zeta(z).
$$ From
this and  (\ref{ej16}) we get
$$
\frac{\sigma'(z+w)}{\sigma(z+w)}=\frac{\sigma'(z)}{\sigma(z)}+\eta
,$$
where $\eta=\eta_w:=m\eta_1+n \eta_2$ for  $w=mw_1+nw_2$. Integrating this,
we obtain 
$$
\log\sg(z+w)=\log\sg(z)+\eta z +c,
$$
$log$s are any fixed branches of logarithm and where $c$ is a constant depending only  $w$. Hence 
$$
\sg(z+w)=\sg(z)\exp(\eta z+c).
$$ 
We now will evaluate $c$. First suppose  that 
$$
w/2\notin \La,
$$
so  that $\sg(w/2)\neq 0$. Putting
$z=-(w/2)$ and using the fact that $\sg$ is an odd function,
we obtain
$$
\sg\lt(\frac{1}{2}w\rt)
=-\sg\lt(\frac{1}{2}w\rt)\exp\lt(-\frac{1}{2}\eta w +c\rt).
$$
So, canceling $\sg(w/2)$, we get 
$$\
\exp(c)=-\exp\lt(\frac{1}{2}\eta w\rt).
$$ 
Thus
\beq\label{ej19}
\sg(z+w)=- \sg(z)\exp\lt(\lt(\eta(z+\frac{1}{2}w\rt)\rt),
\eeq
provided $w/2\notin \La$. Repeating this, we get
$$ 
\sg(z+2w)
=\sg(z)\exp\lt(\eta\lt(z+\frac{1}{2}w\rt)\rt)\exp\lt(\eta\lt(z+\frac{3}{2}w\rt)\rt)=
\sg(z)\exp(2\eta(z+w)).
$$ 
Turning to the notation $\eta=\eta_w$, a straightforward induction then yields:
\beq\label{1_mu_2014_1122}
\sg(z+2kw)=\sg(z)\exp\(2k\eta_w(z+kw)\)
\eeq
for every integer $k\ge 1$ whenever $w/2\notin \La$. Assume now that 
$$
w/2\in\La\sms\{0\}.
$$
Then there exists a unique integer $l\ge 1$ such that 
$$
2^{-l}w\in\La
\  \  \  {\rm but } \  \  \
2^{-(l+1)}w\notin\La.
$$
Applying then \eqref{1_mu_2014_1122} with $w$ replaced by $2^{-l}w$ and $k=2^{l-1}$, we get
$$
\begin{aligned}
\sg(z+w)
&=\sg\(z+2\cdot2^{l-1}(w/2^l)\)
=\sg(z)\exp\(2^l\eta_{w/2^l}(z+2^{l-1}(w/2^l)\) \\
&=\sg(z)\exp\lt(\eta_w\lt(z+\frac{w}2\rt)\rt).
\end{aligned}
$$
Combining this along with (\ref{ej19}) and including the trivial case of $w=0$, we get that
\beq\label{ej21}
\sg(z+w)=\varepsilon\sg(z)\exp\lt(\eta_w\lt(z+\frac{1}{2}w\rt)\rt),
\eeq
where
$$
\varepsilon
:=\left\{\begin{array}{ll} +1\quad\mbox{if}\quad \frac{1}{2}w \in \La  \\
               -1 \quad\mbox{otherwise.}
\end{array}
\right.
$$
Of course $w/2\in \La$ if and only if both $m$ and $n$ are even in the representation $w=mw_1+nw_2$. So, $\varepsilon=(-1)^{mn+m+n}$.

\

We  now turn to the problem, posed in Section~\ref{elliptic}, of
finding an elliptic function $f\in \mathcal E(\La)$ with given zeros
and poles. We will show that the conditions (1), (2'), and (3), given in
Section~\ref{elliptic}, are not only necessary but also sufficient
for the existence of such $f$. This is in contrast with the situation for
rational functions on the sphere, where already the conditions
(1) and (2) of Section~\ref{elliptic} are necessary and sufficient.

\

\bthm\label{t22jones} 
Let $[a_1], \ldots, [a_r]$ and $[b_1], \ldots,
[b_s]$ be elements of $\mathbb C/ \La$ for a lattice $\La$, and let
$k_1, \ldots, k_r, l_1,\ldots,l_s$  be positive integers. If
conditions (1), (2'), and (3) of Section~\ref{elliptic} hold, then there exists an elliptic function $f \in {\mathcal E}(\La)$  with zeros of  multiplicity
$k_j$ at each $[a_j]$, poles of multiplicity $l_j$ at each $b_j$,
and no other zeros or poles. \ethm

\bpf Let $u_1, \ldots, u_n$  be the elements $a_1,
\ldots, a_r$, each $a_j$ being listed $k_j$ times, so that
$n=\sum_{j=1}^r k_j$. Similarly, let $v_1,\ldots, u_n$  be the elements $b_1,
\ldots, b_s$ counting with multiplicities $l_j$. Then (3) takes on the
form
$$
\sum_{j=1}^n u_j -\sum_{j=1}^n v_j=w\in \La,
$$
and there is no loss in replacing $u_1$ by $v_1+w$ so that (3) now reads
\beq\label{ej22}
\sum_{j=1}^n u_j =\sum_{j=1}^n v_j.
\eeq
Now consider  the function
$$ 
f(z):=\frac{\sg(z-u_1)\ldots \sg(z-u_n)}{\sg(z-v_1)\ldots \sg(z-v_n)}.
$$
Since  $\sg(z)$  is an analytic  function, $f(z)$ is meromorphic. By
(\ref{ej21})  we have 
$$
\sg(z-u_j+w_i)
=-\sg(z-u_j) \exp\lt(\eta_i\lt(z-u_j+\frac{1}{2}w_i\rt)\rt), \quad  j=1,\ldots, n, \quad i=1,2,
$$
and similarly  for $\sg( z-v_j+w_i)$. Hence, with the use of \eqref{ej22}, for $i=1,2$ we have
$$ \begin{aligned}
f(z+w_i)&=\frac{(-1)^n \exp \left(\sum_{j=1}^n
\eta_i(z-u_i+\frac{1}{2}w_i)\right)}{(-1)^n \exp \left(\sum_{j=1}^n
\eta_i(z-v_i+\frac{1}{2}w_i)\right)}\cdot f(z)\\
&=\exp \left(\eta_i \sum_{j=1}^n (v_j-u_j)\right) f(z)\\
&=\exp \left(\eta_i\cdot 0 \right)f(z)\\
&=f(z).
\end{aligned}$$
Thus, the function  $f$ is doubly periodic  with respect to $\La$, and hence $f$ is elliptic. Applying Theorem~\ref{3.8.7jones} to
 infinite product (\ref{ej8}) for  $\sg(z)$, we see that $\sg(z)$
has simple  zeros at the lattice-points  $z\in \La$ and that
$\sg(z)\neq0$ for  $z \notin \La$. Hence zeros and poles of $f(z)$
are at $ [a_1], \ldots, [a_r]$ and  $ [b_1], \ldots, [b_s]$  with
multiplicities $k_1, \ldots, k_r$ and $l_1, \ldots, l_s$
respectively. \endpf

\sp If $g$  is any other  elliptic  function with  the same  zeros and
poles as $f$,  then by Theorem~\ref{t12jones} $g(z)=cf(z)$ for some
constant $c\neq 0$.

\sp

\section{The Discriminant of a Cubic  Polynomial}

Let $\La\sbt\C$ be a lattice. It follows from Theorem~\ref{t18jones} that the Weierstrass elliptic
function $\wp_\La$ satisfies the equation $\wp_\La'=\sqrt{p}$, where $p$ is a
cubic polynomial of the form
\beq\label{ej23}
     p(z)=4z^3-c_2z-c_3, \quad  c_2, c_3\in \mathbb C, \quad z=\wp_\La.
\eeq
Any polynomial of the form (\ref{ej23})  is said to be in 
Weierstrass normal form.\index{(N)}{Weierstrass normal form of $p$}
By means of the substitution $\theta(z):=az+b$, $a,b \in \mathbb C$,
$a\neq 0$, any cubic polynomial may be  brought into this form. Now, the map
$\theta: \mathbb{C}\to \mathbb{C}$  is a bijection, preserving  the
multiplicities of the roots of polynomials. So, without loss of generality, we can
restrict our attention to cubic polynomial $p$ in Weierstrass normal
form. If $e_1, e_2, e_3 $\index{(S)}{$e_1, e_2, e_3$} are the roots  of the polynomial $p$ in (\ref{ej23}), then we  define the  discriminant \index{(S)}{$\Delta_p$} of $p$ to be
\beq\label{ej24}
    \Delta_p= 16(e_1-e_2)^2(e_2-e_3)^2(e_3-e_1)^2.
\eeq
Clearly, these roots are distinct if and only if $\Delta_p\neq
0$\index{(N)}{discriminant}. We shall prove the following.

\sp

\bthm~\label{t23jones} If $p(z)=4z^3-c_2z-c_3$, then
$$
\Delta_p=c_2^3-27c^2_3.
$$ 
\ethm
\bpf Since
\beq\label{ej25}
p(z)=4(z-e_1)(z-e_2)(z-e_3),
\eeq
by equating coefficients between  this and (\ref{ej23}), we  have
\begin{eqnarray}\label{ej26}
\begin{aligned}
  e_1+e_2+ e_3 &=&0, \\
  e_1e_2+e_2e_3+ e_3e_1 &=&-\frac{c_2}{4}, \\
  e_1e_2e_3 &=&\frac{c_3}{4}.
\end{aligned}
\end{eqnarray}
The remaining  symmetric functions of the roots may be obtained from
(\ref{ej26}), for example
$$ 
e_1^2+e_2^2+ e_3^2=(e_1+e_2+e_3)^2-
2(e_1e_2+e_2e_3+ e_3e_1)=\frac{c_2}{2},
$$ 
and
$$
e_1^2e_2^2+e_2^2e_3^2+ e_3^2e_1^2=(e_1e_2+e_2e_3+ e_3e_1)^2-2
e_1e_2e_3( e_1+e_2+ e_3)=\frac{c_2^2}{16}.
$$ 
Now, differentiating (\ref{ej23}) and (\ref{ej25}) at $z=e_1$, we have
$$ 
4(e_1-e_2)(e_1-e_3)=p'(e_1)=12e_1^2-c_2
$$
with similar expression for $p'(e_2)$ and $p'(e_3)$. Hence
$$\begin{aligned}
\Delta_p 
&= - \frac{1}{4}p'(e_1)p'(e_2)p'(e_3)\\
 &=- \frac{1}{4}\prod_{i=1}^3( 12 e^2_i-c_2)\\
 & = - \frac{1}{4}\(1728(e_1e_2e_3)^2
 -144c_2(e_1^2e_2^2+e_2^2e_3^2+e_3^2e_1^2)
+12c_2^2(e^2_1+e^2_2+e^2_3)-c_2^3\)\\
&= -\frac{1}{4} (108c_3^2-9c_2^3+6c_2^3-c_2^3)\\
&= c_2^3-27c_3^2.
\end{aligned}
$$
\endpf

\sp

\bcor~\label{c3jones} 
A polynomial $p(z)=4z^3-c_2z-c_3$ has distinct roots if and only if $c_2^3-27c_3^2\neq 0$. 
\ecor

\sp One can give a direct proof of Corollary~\ref{c3jones} without
introducing $\Delta_p$, by eliminating $z$  between the equation
$p(z)=0$ and $p'(z)=0$ thus giving a necessary and sufficiently
condition for $p$ and $p'$ to have a common  root. We have chosen
the above proof since the discriminant is an interesting function in
its own right.

\sp

By Theorem~\ref{t18jones}, (\ref{ej14a}) and  (\ref{ej14b}) the
Weierstrass function $\wp$ associated with a lattice $\La$ satisfies
$\wp_\La'=\sqrt{p(\wp_\La)}$, where $p$ is a cubic polynomial in Weierstrass
normal form
\beq\label{ej27}
    p(z)=4z^3-g_2z -g_3
\eeq
with
$$  
g_2=g_2(\La)=60\sum_{w\in \La}^{}\!\!{'}\,w^{-4}
$$\index{(S)}{$g_2(\La)$} 
and 
$$
g_3=g_3(\La)=140\sum_{w\in \La} ^{}\!\!{'}\,w^{-6}.
$$\index{(S)}{$g_3(\La)$}
If we write $\Delta(\La)$\index{(S)}{$\Delta(\La)$} for the
discriminant $\Delta_p$ of $p$, then by Theorem~\ref{t23jones}, we
have
\beq\label{ej28}
\Delta(\La)=g_2(\La)^3-27g_3(\La)^2.
\eeq
Theorem~\ref{t20jones} implies that  $p$ has  distinct roots, so
$\Delta(\La)\neq 0$ by Corollary~\ref{c3jones}, and hence we may
define a function $J: \La\lra J(\La)$ by
\beq\label{ej29}
J(\La):=\frac{g_2(\La)^3}{\Delta(\La)}=\frac{g_2(\La)^3}{\Delta(\La)^2-27g_3(\La)^2}.
\eeq
This function $J:\La\lra J(\La)$ is called the modular function associated to the lattice $\La$. For a similar lattice $\mu\La, \,\mu\neq 0$, we have
\beq\lab{120200306}
g_2(\mu\La)=60\sum_{w\in \La}^{}\!\!{'}\,(\mu w)^{-4}=\mu^{-4}g_2(\La)
\eeq
and
\beq\lab{220200306}
g_3(\mu\La)=140\sum_{w\in \La}^{}\!\!{'}\,(\mu w)^{-6}=\mu^{-6}g_3(\La),
\eeq
so that
$$
\Delta(\mu\La)=\mu^{-12}\Delta(\La).
$$
Thus,
\beq\label{ej30}
    J(\mu\La)=J(\La)
\eeq
for all $\mu \in \mathbb C \sms \{0\}$, so that similar lattices
determine the same values of $J$.
We therefore may and now will do regard $g_1, g_2$, $\Delta$ and $J$ as functions of the modulus $\tau \in \mathbb{H}$ by  evaluating them on the lattice 
$$
\La=[1, \tau]
$$
which has $\tau$ as one of its moduli. Thus
\beq\label{ej31}
\begin{aligned}
  g_2(\tau)&= 60\sum_{m,n}^{}\!{'}\, (m+n\tau)^{-4}, \\
  g_3(\tau)&= 140 \sum_{m,n}^{}\!{'}\, (m+n\tau)^{-6},
\end{aligned}
\eeq
where $\sum_{m,n}^{'} $ denotes summation over all $(m,n)\in \mathbb
{Z} \times \mathbb{Z}$ except for $(0,0)$. Then
\beq\label{ej32}
    \Delta(\tau)=g_2(\tau)^3-27g_3(\tau)^2
\eeq
and
\beq\label{ej33}
    J(\tau)=\frac{g_2^3(\tau)}{\Delta(\tau)}.
\eeq
In \cite{JS} it is proved (see formula (6.4.9)) that expansion of
$J(\tau)$ has  form
\beq\label{ej33a}
    J(\tau)=\frac{1}{1728}\left(\frac{1}{q}+744+196884q+\ldots\right)
\eeq
By Theorem~\ref{t5jones} for every $T\in PSL(3,\Z)$ the lattices $\La=[1,\tau]$ and $[1,T(\tau)]$ are similar, and hence $J(T(\tau))=J(\tau)$ by
(\ref{ej30}). So, we have the following.

\sp

\bthm\label{t24jones}  
$J(T(\tau))=J(\tau)$ for all $\tau \in PSL(2,\Z)$. \ethm

\sp

\fr This precisely means that $J(\tau)$ is invariant under the action of the modular group
$PSL(2,\Z)$. We now will show that the functions $g_2(\tau), g_3(\tau)$ and
$\Delta(\tau)$ come close to sharing this property. If  
$$
T: \oc\ni\tau \lra \frac{a\tau+b}{c\tau+d }\in\oc
$$ 
is an element of $PSL(3,\Z)$, then
$$
\begin{aligned}
  g_2(\tau)&= 60\sum_{m,n}^{}\!{'} \left(m+n\frac{a\tau+b}{c\tau+d }\right)^{-4} \\
& = 60(c\tau+d)^{-4}\sum_{m,n}^{}\!{'}(m(c\tau+d)+n(a\tau+b))^{-4}\\
&= 60(c\tau+d)^{-4}\sum_{m,n}^{}\!{'}((md+nb)+ (mc+na)\tau)^{-4}.
\end{aligned}
$$
Since $ad-bc=1$, the transformation 
$$
(m,n)\longmapsto (md+nb, mc+na)
$$ 
merely permutes the elements of the indexing set
$\mathbb{Z}\times\mathbb{Z}\sms\{(0,0)\}$. Hence, by Theorem~\ref{t4jones}
and by absolute convergence (implying unconditional convergence), we have
\beq\label{ej34}
 g_2(T(\tau))=
60(c\tau+d)^{-4}\sum_{m,n}^{}\!{'} (m
+n\tau+b)^{-4}=(c\tau+d)^{-4}g_2(\tau).
 \eeq
Similarly
\beq\label{ej35}
 g_3(T(\tau))=(c\tau+d)^{-6}g_3(\tau),
 \eeq
and hence
\beq\label{ej36}
 \Delta(T(\tau))=(c\tau+d)^{-12}\Delta(\tau),
\eeq
from which we immediately obtain an alternative  proof of
Theorem~\ref{t24jones}. In the special  case, where $a=b=d=1$ and
$c=0$, we have $T(\tau)=\tau+1$, giving

\sp

\bthm\label{t25jones} The functions $g_2(\tau)$, $g_3(\tau)$, and
$J(\tau)$ are periodic with respect to $\mathbb{Z}$. 
\ethm

\sp It is also useful  to determine the effect  on these functions of
the orientation reserving transformation of $\mathbb{H}$ of the
form
\beq\label{ej37}
    T(\tau)=\frac{a\ov{\tau}+b}{c\ov{\tau}+d}, \quad a,b,c,d \in
    \mathbb Z, \,\, ad-bc=-1.
\eeq
Calculations similar to  those  above give
\beq\label{ej38}
\begin{aligned}
   g_2(T(\tau))=&(c\ov{\tau}+d)^{-4}\ov{g_2(\tau)},\\
  g_3(T(\tau))=&(c\ov{\tau}+d)^{-6}\ov{g_3(\tau)},\\
\De(T(\tau))=&(c\ov{\tau}+d)^{-13}\ov{\Delta(\tau)},\\
 J(T(\tau))=&\ov{J(\tau)}.
\end{aligned}
\eeq

Now we shall prove the following.

\sp\bthm\label{t26jones} 
The functions $g_2$, $g_3$, $\Delta$ and
$J:\mathbb{H}\lra \mathbb{C}$ are analytic on $\mathbb{H}$. \ethm

\bpf Given any $\tau_0\in \mathbb{H}$, let
$$
\delta:=\frac{1}{2}\im (\tau_0),
$$
so that $\delta>0$, and let
$$
\ov{D(\tau_0, \delta)}:=\{\tau \in \mathbb{H}:  \,
|\tau-\tau_0|\leq \delta\}.
$$
Now, the functions 
$$
\mathbb{H}\ni\tau\longmapsto(m+n\tau)^{-4}\in\C
\  \  \  {\rm and} \  \  \
\mathbb{H}\ni\tau\longmapsto (m+n\tau)^{-6}\in\C
$$ 
are holomorphic for all $(m,n)\in
\mathbb{Z}\times\mathbb{Z}\sms\{(0,0)\}$, so if  we can show  that
the series (\ref{ej31}) defining the functions $g_2$ and  $g_3$ are
absolutely uniformly convergent on each set $\ov{D(\tau_0, \delta)}$, $ \tau_0 \in
\mathbb{H}$, then  these  two functions are analytic on
$\mathbb{H}$. For all $m,n \in \mathbb{Z}$ with  $n\neq 0$ we have $-m/n\in
\mathbb{R}$, and hence
$$
\left|\frac{m}{n}+\tau_0\right| \geq \im (\tau_0)=2\delta.
$$
Therefore  for all $m,n\in \mathbb{Z}$ (including n=0) and $\tau \in
\ov{D(\tau_0, \delta)}$ we have
$$
|(m+n\tau)-(m+n\tau_0)|=|n||\tau-\tau_0|
\leq |n|\delta
\leq \frac{1}{2}|m+n\tau_0|, 
$$
and so  the triangle inequality gives
$$
|m+n\tau|\geq |m+n\tau_0|-|(m+n\tau)-(m+n\tau_0)|\geq
\frac{1}{2}|m+n\tau_0|.
$$ 
Thus for any $r>0$ we have
$$
|m+n\tau|^{-2r}\leq 2^{2r}|m+n\tau_0|^{-2r}
$$
for all $\tau \in \ov{D(\tau_0, \delta)}$ and $(m,n)\neq (0,0)$. By
Theorem~\ref{t14jones} the series 
$$
\sum_{m,n}^{}\!{'}|m+n\tau_0|^{-2r}
$$
converges for each $r>1$, so, by the Weierstrass test, the series
$$
\sum_{m,n}^{}\!{'}|m+n\tau_0|^{-2r}
$$ 
converges absolutely uniformly on $\ov{D(\tau_0, \delta)}$.  Putting $r=2$ and
$r=3$, we thus see that the functions $g_2$ and $g_3$ are analytic  on
$\mathbb{H}$. It immediately follows from (\ref{ej32}) that $\Delta(\tau)$ is
analytic on $\mathbb{H}$, and since  $\Delta(\tau)\neq 0$ on
$\mathbb{H}$ by Theorem~\ref{t20jones} and Corollary~\ref{c3jones}, it follows from (\ref{ej33})  that the function $J$ is analytic on $\mathbb{H}$.
\endpf

\sp

\blem\label{l1jones}
We have the following two properties.
\begin{enumerate}
  \item If  $2 \re(\tau)\in \mathbb{Z}$ then $g_2$, $g_3$, $\Delta$ and
  $J$ are  all real.
\,

\item If $|\tau|=1$ then 
$$
g_2(\tau)=\tau^4\ov{g_2(\tau)}, \  \
  g_3(\tau)=\tau^6\ov{g_3(\tau)}, \  \
  \Delta(\tau)=\tau^{12}\ov{\Delta}(\tau) \  \
{\rm and}, \  \ J(\tau)=\ov{J(\tau)}.
$$
\end{enumerate}
\elem

\bpf (1) If $2 \re(\tau)=n \in \mathbb{Z}$ then $\tau $
is fixed by reflection $T(\tau)=n-\ov{\tau}$, which is of type
(\ref{ej37}) with $a=-1$, $b=n$ and $d=1$, so the result follows
from (\ref{ej38}).

(2)  If $|\tau|=1$ then $\tau$ is fixed  by  the inversion
$T(\tau)=1/\tau$ in the unit circle; putting $a=d=0$
and $b=c=1$ in (\ref{ej38}) we have
$$
g_2(\tau)=g_2(1/\ov{\tau})={\ov{\tau}}^{-4}\ov{g_2(\tau)}={\tau}^{4}\ov{g_2(\tau)},$$
and similarly for the other three functions. \endpf

\sp Theorem~\ref{t7jones} asserts that the modular group $PSL(2,\Z)$ has a fundamental region $$
{\mathfrak F}:=\lt\{z\in {\mathbb H}: 
|z|\ge 1\  \  \mbox{and}\  \ |\re z | \leq \frac{1}{2}\rt\}.
$$
We thus immediately have

\sp

\bcor\label{c4jones} $J(\tau)$ is real whenever $\tau$ is on the
imaginary axis or on the boundary $\partial{\mathcal{R}}$ of
$\mathcal{R}$.
 \ecor

\bcor\label{c5jones} $g_2(\epsilon)=g_3(i)=J(\epsilon)=0$ and $J(i)=1$,
where $\epsilon=e^{2\pi i/3}$. 
\ecor

\bpf Part (1) of Lemma~\ref{l1jones}  shows that $g_2$
and $g_3$ both take real  values at $i$ and $\epsilon$, while part (2)
shows that $g_2(\rho)=\rho \ov{g_2(\epsilon)}$ and $g_3(i)=-
\ov{g_3(i)}$. Thus $g_2(\rho)=g_3(i)=0$ so (\ref{ej33}) gives
$J(\epsilon)=0$ and $J(i)=1$.
\endpf

\sp

Let 
$$
L:=L_1\cup L_2\cup L_3,
$$ 
where
$$
\begin{aligned}
L_1&:=\lt\{\tau \in \mathbb{H}: |\tau| \geq 1\quad\mbox{and}\quad
\re(\tau)=-\frac{1}{2}\rt\}, \\
L_2&:=\lt\{\tau \in \mathbb{H}: |\tau|=1\quad\mbox{and}\quad
-\frac{1}{2}\leq \re(\tau) \leq 0\rt\}, \\
L_3&:=\{\tau \in \mathbb{H}: |\tau|=1\quad\mbox{and}\quad
 \re(\tau) =0\}.
\end{aligned}
$$
By Corollary~\ref{c4jones} we  have $J(L)\sbt \mathbb{R}$, but in
fact we can prove equality.

\

\bthm\label{t27jones} $J(L)=\mathbb{R}$. \ethm

\bpf If $\tau \in L_3$, then we have $\tau=iy$ with some
$y\geq 1$. Then 
$$
q=:e^{2\pi i\tau}=e^{-2\pi y}\lra 0
$$
through positive real numbers as $y \to +\infty$. So (\ref{ej33a})
gives 
\beq\label{120200109}
J(\tau)=(q^{-1}+744 +\ldots )/1728\xrightarrow[\ \, y \to +\infty \ ]{}  
+\infty.
\eeq
Similarly, on $L_1$ we have $\tau=-\frac{1}{2}+iy$ with some $y\ge 0$ and 
$$
q:=-e^{2\pi y}\xrightarrow[\ \, y \to +\infty \ ]{} 0
$$
and the convergence is through negative real numbers. Hence 
\beq\label{220200109}
\lim_{y \to +\infty}J(\tau)=-\infty.
\eeq
Since $J(L)\sbt\R$, $J|_L:T\lra\R$ is continuous as being analytic on $\mathbb H$ by Theorem~\ref{t26jones}, and since $L$ is connected, it follows from \eqref{120200109} and \eqref{220200109}
$J(L)=\mathbb R$.
\endpf

\sp

By Theorem~\ref{t24jones}, $J$ is constant on each orbit of $PSL(2,\Z)$
in $\mathbb{H}$. We shall prove more:

\sp

\bthm\label{t28jones} For each $c\in \mathbb C$ there is exactly one
orbit of $PSL(2,\Z)$ in $\mathbb H$ on which  $J$ takes the value $c$.
\ethm

\bpf Each orbit of $\Gamma$ meets the fundamental region
$\mathfrak{ F}$ either at a unique point in the interior of
$\mathfrak{ F}$, or else  at one or two equivalent  points on
$\partial{\mathfrak{ F}}$.

  First suppose that $ c \in\mathbb{C}\sms \mathbb{R}$. Since
  $J(\partial{\mathfrak{ F}})\subset \mathbb{R}$ by Corollary~\ref{c4jones}, it
  is sufficient to show  that there is a unique solution to the equation
  $J(\tau)=c$ in ${\rm Int}{\mathfrak{ F}}$. By Theorem~\ref{t26jones} and Corollary~\ref{c5jones}
  $J$ is analytic and not identically equal to $c$, so the function
  $$g(\tau)=\frac{J'(\tau)}{J(\tau)-c}$$
is meromorphic on $\mathbb{H}$. We can use (\ref{ej33a})  to
express $g(\tau)$ as a function of $q=e^{2\pi i\tau}$, meromorphic
at $q=0$ since $J(\tau)$ is. Hence $g(\tau)$ is analytic for
sufficiently small non-zero $|q|$, that  is, provided $\im \tau $ is
sufficiently large, say $\im \tau \geq K $ for some $K>1$. Thus the
poles of $g(\tau)$ in  $\mathfrak{ F}$ all lie in  the interior of
$$
G:=\{\tau\in \mathfrak{ F}:\im \tau \leq K\}.
$$
So, the sum of  the
residues of $g(\tau)$ in $\mathfrak{ F}$ (and hence the number  of
solutions, counting multiplicities, of $J(\tau)=c$ in $\mathfrak{ F}$)
is equal  to
\beq\label{ej39}
    \frac{1}{2\pi i}\int_{\partial{G}}g(\tau)d\tau,
\eeq
where the boundary $\partial{G}$ is given the positive orientation.

Now the sides $\re \tau=-1/2$ and  $\re \tau=1/2$ of
$G$ are equivalent under the  transformation 
$$
\tau\lmt \tau+1
$$ 
of $\Gamma$, so $J(\tau)$ and  hence $g(\tau)$ take the same values at
equivalent points  on these sides. Hence the integrals of $g(\tau)$
along these sides cancel in (\ref{ej39}), and similarly the integral
along the unit  circle from $\rho$ to $i$  cancels with the integral
from $i$ to $\rho+i$, using the transformation $\tau \to 1/ \tau$.
Hence,
 $$ 
 \int_{\partial{G}}g(\tau)d\tau=\int_{\gamma}g(\tau)d\tau,
 $$
where $\gamma$ is the side  $\im (\tau)=K$ of $G$ oriented from
$\frac{1}{2}+iK$ to $-\frac{1}{2}+iK$. Away from the poles  of
$g(\tau)$, each branch of the logarithm function of $J$ satisfies
$$
\frac{d}{d\tau }(\log J(\tau)-c))=\frac{J'(\tau)}{J(\tau)-c}=g(\tau).
$$
So
$$ 
\int_{\gamma}G(\tau)d\tau=[\log(J(\tau)-c)]_\gamma,
$$
the change in the  value of $\log(J(\tau)-c)$ arising from analytic
continuation along $\gamma$. As $\tau$ follows $\gamma$, $q$ winds
once (in the negative direction) around the circle $\delta$ given by
$|q|=e^{-2\pi K}$, starting and finishing at $-e^{-2\pi K}$. By
(\ref{ej33a}), $q(J(\tau)-c)$ is analytic and non-zero for $0\leq
|q|\leq e^{-2\pi K}$, and since this set is simply-connected, the
Monodromy Theorem~\ref{monodromy}  implies that
$$ 
[\log q(j(\tau)-c)]_\gamma=0.
$$
So 
$$[\log q(j(\tau)-c)]_\gamma
=[\log q(j(\tau)-c)-\log q]=[\log q]_\gamma=2\pi i.
$$ Hence  (\ref{ej39}) shows that the number of
solutions of $J(\tau)=c$ in $\mathfrak{ F}$ is  equal  to $(1/2\pi
i)=1$, as required.

 Finally, suppose that $ c \in \mathbb{R}$. By Theorem~\ref{t27jones} there is
 at least one orbit of $\Gamma$ on which $J$ takes the value $c$.
If there were more than one such orbit, there would be two
inequivalent solutions $\tau_1, \tau_2$ of $J(\tau)=c$, so by
choosing $c'\in  \mathbb C\sms \mathbb R$ sufficiently close to $c$
we would have two inequivalent solutions $\tau_1'$ and $\tau_2'$ of
$J(\tau)=c'$, close to $\tau_1$ and $\tau_2$ respectively. We have
already shown that this is impossible, so the orbit is unique.
\endpf

We record the following two consequences of this theorem.

\sp\bcor\label{c6jones}
If $c_2, c_3 \in \mathbb C$ satisfy $c_2^3-27c_3^2\neq 0$, then
there is a lattice $\La \sbt \mathbb C$ with $g_k(\La)=c_k$ for $k=2,
3$.
 \ecor
 \bpf First suppose that $c_2=0$, so that $c_3\neq 0$.
By Corollary~\ref{c5jones} $g_2(\rho)=0$ and hence $g_3(\rho)\neq 0$
since $g_2(\tau)^3-27c_3^2\Delta(\tau)$ does not vanish on
$\mathbb{H}$. We can therefore choose $\mu \in \mathbb C \sms
\{0\}$ such that $\mu^{-6}g_3(\rho)=c_3$, so putting
$$
\La :=\mu[1, \tau]= [\mu, \mu \tau]
$$
we have  $g_2(\La)=\mu^{-4}g_2(\rho)=0=c_2$ and hence
$g_3(\La)=\mu^{-6}g_3(\rho)=c_3$, as required.

Similarly, if $c_3=0$ then $c_2\neq 0$. We have $g_3(i)=0\neq
  g_2(i)$, so we  can choose $\mu \in \mathbb C\sms \{0\}$
  satisfying $\mu^{-4}g_2(i)=c_2$, and then $\La = [\mu, \mu i]$
  satisfies 
$$
  g_2(\La)=\mu^{-4}g_2(i)=c_2
\  \  \   {\rm and} \  \  \ 
g_3(\La)=0=c_3.
$$
 Finally  we consider  the general case, where $c_2\neq 0 \neq c_3$.
 By Theorem~\ref{t28jones} there exists $\tau \in \mathbb{H}$  such that
  \beq\label{ej40}
     J(\tau)=\frac{c_2^3}{c_2^3-27 c_3^2}.
  \eeq
For any $\mu\in \mathbb C \sms\{0\}$ the lattice $\La=[\mu, \mu
\tau]$ satisfies $g_2(\La)=\mu^{-4}g_2(\tau)$ and
$g_3(\La)=\mu^{-6}g_3(\tau)$, both  non-zero since

\beq\label{ej41}
\begin{aligned}
 J(\tau)&=\frac{c_2^3}{c_2^3-27 c_3^2},\\
   J(\La)&=\frac{g_2^3(\La)}{g_2^3(\La)-27 g_3^2(\La)}
\end{aligned}
\eeq
does not take the value $0$ or $1$  by (\ref{ej40}) and the fact
that $c_2\neq 0 \neq c_3$. We can therefore  choose $\mu \neq 0$ so
that
$$ \mu^2=\frac{c_2}{c_3}\frac{g_3(\tau)}{g_2(\tau)},$$
and hence
$$\frac{g_2(\tau)}{g_3(\tau)}=\frac{\mu^{-4}g_2(\tau)}{\mu^{-6}g_3(\tau)}=\frac{c_2}{c_3}.
$$
Thus $g_k(\La)=\l c_k$ (k=2,3) for some $\l \neq 0$,  so
substituting in (\ref{ej41}) and using (\ref{ej40}) we have
$$  J(\La)= \frac{c_2^3}{c_2^3-27 c_3^2}$$
and $$J(\La)=\frac{\l^3c_2^3}{\l^3c_2^3-27\l^2
c_3^2}=\frac{c_2^3}{c_2^3-27\l^{-1} c_3^2}. $$ Hence $\l=1$ and so
$g_k(\La)=c_k,$ $k=2,3$, as  required. \endpf

\sp

\bcor\label{hawkins}
The numbers $g_2$ and $g_3$ form a set of full invariants for lattices in $\C$. More precisely, if $\La, \La'$ are two lattices in $\C$, then
$$
\Lambda'=\Lambda
 \  \lekra \ 
\big[g_2(\Lambda')=g_2(\Lambda)
\  \  {\rm  and}  \  \
g_3(\Lambda')=g_3(\Lambda)\big]
$$
In addition, given $\xi\in\C\sms\{0\}$,
$$
\La'=\xi\La
 \ \lekra \ 
\big[g_2(\La')=\xi^4g_2(\La)
\  \  {\rm  and}  \  \
g_3(\La')=\xi^6g_2(\La)\big].
$$
\ecor

\sp\section{Weierstrass $\wp$--functions II}\label{weierstrassII}

In this section we collect and summarize, for the convenience of the reader, some properties of the Weierstrass $\wp$--functions which have been proved in the present chapter. We actually bring up only these properties which will be used in Chapter~\ref{examples} in our constructions of dynamically significant examples.

\sp Recall that given any lattice $\Lambda$ in $\C$, the 
Weierstrass elliptic function \index{(N)}{Weierstrass elliptic function} $\wp_\Lambda$ \index{(S)}{$\wp_\Lambda$}is defined by the following formula:
$$  
\C\ni z\longmapsto {\wp}_{\Lambda}(z)=\frac{1}{z^2}+ \sum_{w\in
\Lambda\sms\{0\}}\left(\frac{1}{(z-w)^2}-\frac{1}{w^2}\right)\in\oc.
$$
Replacing $z$ by $-z$ in the
definition we see that $\wp_{\Lambda}$ is an even function which is
analytic  in $\mathbb{C} \sms \La$ and has poles of order 2 at each
$w \in \La$ (see Theorem~\ref{weierstrass_convergence} and Theorem~\ref{t17jones}.

\sp

The derivative of the Weierstrass elliptic function $\wp_{\Lambda}$
is also an elliptic function, periodic with respect to
$\Lambda$, and is expressible by the series 
\beq\lab{220200317}
{\wp}_{\Lambda}'(z)=-2\sum_{w\in
\Lambda}\frac{1}{(z-w)^3}.
\eeq
The Weierstrass elliptic  function $\wp_{\Lambda}$ and its derivative
are related by the differential equation
\beq\lab{220200307}
{\wp}_{\Lambda}'(z)^2= 4 \wp_{\Lambda}(z)^3-g_2\wp_{\Lambda}(z) -g_3,    
\eeq
where (see (\ref{ej14a}), (\ref{ej14b}) and (\ref{ej15})), we recall,  
$$
g_2(\Lambda)=60\sum_{w\in \Lambda\sms \{0\}}w^{-4}
\  \  \  {\rm and}  \  \  \
g_3(\Lambda)=140\sum_{w\in \Lambda\sms \{0\}}w^{-6}.
$$

Recall that the zeros of $\wp_{\Lambda}'$ are called the critical points of ${\wp}_{\Lambda}$ and the set of all of them is denoted by $\Crit(\wp_{\Lambda})$. As an immediate consequence of 
Theorem~\ref{t19jones} we get the following.

\bthm\lab{t320200303}
If $\La=[w_1,w_2]\sbt\C$ is a lattice, then
$$
\Crit(\wp_{\Lambda})
=\lt(\frac{w_1}{2}+\Lambda\rt)\cup \lt(\frac{w_2}{2}+\Lambda\rt)\cup
\lt(\frac{w_3}{2}+\Lambda\rt)
$$ 
and
$$
\wp_{\Lambda}\(\Crit(\wp_{\Lambda})\)=\{e_1,e_2,e_3\},
$$
where 
$$
e_1:=\wp_{\Lambda}\left(\frac{w_1}{2}\right), \ \ 
e_2:=\wp_{\Lambda}\left(\frac{w_2}{2}\right), \  \text{ and } \
e_3:=\wp_{\Lambda}\left(\frac{w_3}{2}\right). 
$$
\ethm

\sp For any lattice $\La$, the Weierstrass elliptic function and its
derivative satisfy the following properties: for every $z\in\C$ and every
$s\in \mathbb{C}\sms\{0\}$, we have that
\beq\label{2}
\wp_{s\La}(sz)=\frac{1}{s^2}\wp_\La(z)\  \  \  -  \  \  \ \mbox{homogenity of
$\wp_\La$}
\eeq\index{(N)}{homogenity of
$\wp_\La$}
\beq\label{2a}
\wp_{s\La}'(sz)=\frac{1}{s^3}\wp_\La'(z)\  \  \  - \  \  \  \mbox{homogenity of
$\wp_\La'$}
\eeq\index{(N)}{homogenity of
$\wp_\La'$}
Verification of these properties can be done  by
substituting into the corresponding series definitions.

\sp We discern below some distinguished types of lattices.

\sp\bdfn\label{d120190216}
We define the following four classes of lattices:

\sp\begin{enumerate}

\item[(0)] A lattice $\La$ is is called real
if and only if 
$$
\Lambda=\ov{\Lambda}.
$$
\index{(N)}{real lattice}

\item A lattice $\Lambda$ is called {\em real rectangular}\index{(N)}{real rectangular lattice} if and only if there exist $\l_1, \l_2\in\C$ such that such that $\Lambda=[\l_1, \l_2]$ and
$$
\l_1\in\R
\  \  \  {\rm and}  \  \  \  \l_2\in i\R.
$$
Then (obviously) $\Lambda$ is real. Any lattice similar to a real rectangular lattice\index{(N)}{real rectangular lattice} is called rectangular. 

\,

\item A lattice $\Lambda$ is called real rhombic  \index{(N)}{real rhombic lattice} if and only if  
$$
\Lambda=[\l, \ov\l]
$$ 
with some $\l\in\C\sms\{0\}$. Then (obviously) $\Lambda$ is real. Any lattice similar to a real rhombic lattice is called rhombic. 

\,

\item A lattice $\Lambda$ is called a {\em square} lattice if and only if 
$$
i\Lambda=\Lambda.
$$
\index{(N)}{square lattice} Equivalently, $\Lambda$ is a square lattice if and only if it is similar to the lattice $[1,i]$. In addition, a square lattice is real if and only if it is of the form $[\l, \l i]$, $\lambda\in(0,+\infty)$.

\,

\item A lattice $\Lambda$ is called {\em triangular}\index{(N)}{triangular
lattice} if and only if  
$$
\Lambda=e^{2\pi i/3}\Lambda.
$$
\end{enumerate}
\edfn 

\sp In each of the cases (1)--(3) the numbers $\l_1$ and $\l_2$ can be  chosen so that the fundamental parallelogram with
vertices $0, \l_1, \l_2$, and $\l_3:=\l_1+\l_2$ is a rectangle, rhombus, or square respectively. In the  case
(4) the  fundamental parallelogram is comprised of two
equilateral  triangles.

\sp
\bdfn\label{d220190216} 
A meromorphic function $f:\C\lra\oc$ is called  {\it real}\index{(N)}{real function} if and only if   
$$
f(\bar{z})=\ov{f(z)}
$$
for all $z \in \mathbb C$.
\edfn

\sp The proof of the following proposition can be found in\cite{Du}.

\sp\bprop\label{th4.3in{HK1}} 
For a lattice $\La\sbt\C$ the following conditions are equivalent:
\ben
\item $\La$ is real.

\,

\item $\La$ is real rectangular or real rhombic.

\,

\item $g_2$ and $g_3$  are real.

\,

\item  The Weierstrass $\wp_\La$ function is real.   
\een
\eprop

As an immediate consequence of the definition of real matrices we get the following.

\bobs\label{o120200313}
If a lattice $\La\sbt\C$ is real, then so is the lattice $i\La$. 
\eobs

\bprop\lab{p120200306}
If $\La\sbt\C$ is a lattice, then the following three conditions are equivalent.
\ben
\item $\La$ is triangular.

\item 
$$
\La=\lt[\lambda, \lambda e^{\frac{2\pi i }{3}}\rt]
=\lt[\lambda e^{\frac{\pi i }{3}}, \lambda e^{-\frac{\pi i }{3}}\rt],
$$ 
with some $\lambda\in\C\sms\{0\}$.

\,

\item $g_2(\La)=0$. 
\een

\,

In addition, 
\ben 

\item [(4)] For every $g\in \C\sms\{0\}$ there exists a unique triangular lattice $\La\in\C$ such that
$$
g_3(\La)=g.
$$

\item [(5)] The (three) critical values of the Weierstrass $\wp_\La$ function coincide with the cubic roots of $g_3(\La)/4$. If $g_3(\La)\in\R$, then $e_3(\La)$ is the real cubic roots of $g_3(\La)/4$.

\,

\item [(6)] A triangular lattice $\La\in\C$ is real if and only if 
$$
\La=\lt[\lambda, \lambda e^{\frac{2\pi i }{3}}\rt]
=\lt[\lambda e^{\frac{\pi i }{3}}, \lambda e^{-\frac{\pi i }{3}}\rt],
$$ 
with some $\lambda\in (0,+\infty)\cup i(0,+\infty)$. 

\,

Then we have the following:
\begin{enumerate}
\item [(6a)] 
$$
g_3(\l[1,\epsilon])>0 \  \eqv  \  \l\in (0,+\infty)
\  \  \  {\rm and}  \   \   \
g_3(\l[1,\epsilon])<0 \  \eqv  \  \l\in i(0,+\infty).
$$

\item [(6b)] 
$$
\La\cap\R=\l\Z, \  \    \Crit\(\wp_\La\)\cap\R=\frac{\l}{2}+\l\Z,
\  \  \  {\rm and}  \   \
\wp''(z)>0
$$
for all $z\in\R\sms \l\Z$, and 

\,

\item [(6c)] For every $k\in\Z$, the real--valued function $\wp_\La|_{[\l k, \l(k+1)]}$
\begin{enumerate}

\,

\item[(6c1)] is continuous,

\,

\item [(6c2)] is strictly decreasing on $\lt[\l k, \l k+\frac{\l}{2}\rt]$, 

\,

\item [(6c3)] is strictly increasing on $\lt[\l k+\frac{\l}{2},\l(k+1)]\rt]$,

\,

\item [(6c4)] has a unique absolute minimum at the point $\l k+\frac{\l}{2}$ with the value $e_3(\La)$. 

\,

\fr In addition,

\,

\item [(6c5)] 
$$
\wp_\La\lt(\lt[\l k, \l k+\frac{\l}{2}\rt]\rt)
=[e_3(\La),+\infty)
=\wp_\La\lt(\lt[\l k+\frac{\l}{2},\l(k+1)\rt]\rt)
$$
\een
\een 
\een
\eprop

\bpf  
The implication (2)$\imp$(1) is obvious. 

The implication (1)$\imp$(3) directly follows from the formula \eqref{120200306} applied with $\mu=\epsilon$. 

The formula \eqref{220200306} yields
\beq\lab{120200307}
g_3([\mu,\epsilon\mu])=\mu^{-6}g_3([1,\epsilon]). 
\eeq
So, knowing that the lattices $[\mu,\epsilon\mu]$, $\mu\ne 0$, are all triangular, by varying $\mu$ all over $\C\sms\{0\}$, we conclude from Corollary~\ref{c6jones},  Corollary~\ref{hawkins}, and already proven implication (1)$\imp$(3), that item (4) holds, that (1) implies (2), and that (3) implies (2). 

Item (5) directly follows from \eqref{220200307} and item (3) of the current proposition. The proof is complete.

Proving the equivalence (first) part of item (6) and also item (6a), the fact that the lattices $[\lambda, \lambda \epsilon]
=[\lambda \epsilon, \lambda \epsilon]$, $\lambda\in(\R\cup i\R)\sms\{0\}$, are all real, is obvious. Since $g_3(\La)\in \R\sms\{0\}$ for every real lattice $\La$, by varying $\mu$ all over $(0,+\infty)\cup i(0,+\infty)$ in the formula \eqref{120200307}, we obtain the converse implication too and item (6a) as well. 

The first two formulas of item (6b) are immediate with the use Theorem~\ref{t320200303} and realness and triangularity of of $\La$. The last formula of item (6b), i.e $\wp_\La''(z)>0$ directly follows from the formula \eqref{220200317} after differentiating it. Items (6c1)--(6c4) are immediate consequences of item (6b). 
\epf

\sp 

\bprop\lab{p120200310}
If $\La\sbt\C$ is a lattice, then the following four conditions are equivalent.
\ben
\item $\La$ is a square lattice.

\item 
$$
\La=[\lambda, \lambda i]
$$ 
with some $\lambda\in\C\sms\{0\}$.

\,

\item $g_3(\La)=0$. 

\sp

\fr In addition: 

\,

\item [(4)] For every $g\in \C\sms\{0\}$ there exists a unique square lattice $\La\in\C$ such that
$$
g_2(\La)=g.
$$

\item [(5)] If $\La$ is a square lattice, then the (three) critical values of the Weierstrass $\wp_\La$ function are:
$$
e_1=\frac{1}{2}g_2(\La)^{\frac{1}{2}}, \  \  \  e_2=-e_1,
\  \  \  {\rm and } \  \  \
e_3=0.
$$ 
In particular, $e_3$ is a pole of $\wp_\La$. If in addition $\lambda\in (0,+\infty)$, then $g_2(\La)>0$ and $e_1$ is a real positive number. 

\,

\item [(6)] A square lattice $\La\in\C$ is real if and only if 
$$
\La=[\lambda,\lambda i]
$$ 
with some $\lambda\in (0,+\infty)\cup e^{\frac{\pi i}{4}}(0,+\infty)$.

\,

Then we have the following:
\begin{enumerate}
\item [(6a)] 
$$
g_2(\l[1,i])>0 \  \eqv  \  \l\in (0,+\infty)
\  \  \  {\rm and}  \   \   \
g_2(\l[1,i])<0 \  \eqv  \  \l\in i(0,+\infty).
$$

\item [(6b)] 
$$
\La\cap\R=\l\Z, \  \    \Crit\(\wp_\La\)\cap\R=\frac{\l}{2}+\l\Z,
\  \  \  {\rm and}  \   \
\wp''(z)>0
$$
for all $z\in\R\sms \l\Z$, and 

\,

\item [(6c)] For every $k\in\Z$, the real--valued function $\wp_\La|_{[\l k, \l(k+1)]}$
\begin{enumerate}

\,

\item[(6c1)] is continuous,

\,

\item [(6c2)] is strictly decreasing on $\lt[\l k, \l k+\frac{\l}{2}\rt]$, 

\,

\item [(6c3)] is strictly increasing on $\lt[\l k+\frac{\l}{2},\l(k+1)]\rt]$,

\,

\item [(6c4)] has a unique absolute minimum at the point $\l k+\frac{\l}{2}$ with the value $e_1(\La)$. 

\,

\fr In addition,

\,

\item [(6c5)] 
$$
\wp_\La\lt(\lt[\l k, \l k+\frac{\l}{2}\rt]\rt)
=[e_1(\La),+\infty)
=\wp_\La\lt(\lt[\l k+\frac{\l}{2},\l(k+1)\rt]\rt)
$$
\een
\een
\een
\eprop

\bpf 

The implication (2)$\imp$(1) is obvious. 

The implication (1)$\imp$(3) directly follows from the formula \eqref{120200306} applied with $\mu=i$. 

The formula \eqref{220200306} yields
\beq\lab{120200307B}
g_2([\mu,\mu i])=\mu^{-4}g_2([1,i]). 
\eeq
So, knowing that the lattices $[\mu,i\mu]$, $\mu\ne 0$, are all square, by varying $\mu$ all over $\C\sms\{0\}$, we conclude from Corollary~\ref{c6jones},  Corollary~\ref {hawkins}, and already proven implication (1)$\imp$(3), that item (4) holds, that (1) implies (2), and that (3) implies (2). 

Proving item (5), since, by item (3) of the current proposition $g_3(\La)=0$, the equation \eqref{220200307} takes on the form
$$
4e_i^2-g_2=0, \  \  \    i=1, 2, 3.
$$
So, we are only left to show that if $\lambda\in (0,+\infty)$ , then $g_2(\La)>0$ and $e_1>0$. But these follow by a direct calculation. The proof is complete.

Proving item (6), the fact that the lattices $[\lambda, \lambda i]$,
$\lambda\in(\R\cup e^{\frac{\pi i}{4}}\R)\sms\{0\}$, are all real, is obvious. Since $g_2(\La)\in \R\sms\{0]\}$ for every real lattice $\La$, by varying $\mu$ all over $(0,+\infty)\cup e^{\frac{\pi i}{4}}(0,\infty)$ in the formula \eqref{120200307B}, we obtain the converse implication too. 

The first two formulas of item (6b) are immediate with the use Theorem~\ref{t320200303} and realness and squareness of $\La$. The last formula of item (6b), i.e $\wp_\La''(z)>0$, directly follows from the formula \eqref{220200317} after differentiating it. Items (6c1)--(6c4) are immediate consequences of item (6b). \epf

\sp We call a lattice $\La\sbt\C$ real rhombic square \index{(N)}{real rhombic square} if and only if it is real rhombic and square. The following proposition is an immediate consequence of the previous one.

\sp\bprop\lab{p520200310} 
If $\La\sbt\C$ is a lattice, then the following three conditions are equivalent. 
\ben

\,

\item $\La\sbt\C$ real rhombic square lattice.

\,

\item $\La=ae^{\frac{\pi}{4}i}[i,1]=b[1+i,1-i]$ for some $a,b\in(0,+\infty)$.

\,

\item $g_3(\La)=0$ and $g_2(\La)<0$, and these two values determine $\La$ uniquely.

\sp

\fr In addition:

\,

\item  If $\La\sbt\C$ is a real rhombic square lattice, then the (three) critical values of the Weierstrass $\wp_\La$ function are:
$$
e_1(\La)=\frac{1}{2}g_2(\La)^{\frac{1}{2}}, \  \  \  e_2(\La)=-e_1(\La),
\  \  \  {\rm and } \  \  \
e_3(\La)=0.
$$ 
In particular, both $e_1(\La)$ and $e_2(\La)$ are pure imaginary numbers, $e_2(\La)=-e_1(\La)$, $e_1(\La)i\in(0,+\infty)$, and $e_3$ is a pole of $\wp_\La$. 
\een
\eprop

\sp

\chapter{Topological Picture of Iterations of (all!) Meromorphic
Functions}\label{topological-picture} 
 
In this chapter we provide a relatively short and condensed account of topological dynamics of all meromorphic functions with emphasize on Fatou domains including Baker domains that are exclusive for transcendental functions and do not occur for rational functions. We do this for all meromorphic functions and not merely for elliptic ones. In particular, we provide a complete proof of Fatou's classification of Fatou periodic components. We analyze the structure of these components and the structure of their boundaries in greater detail. Particularly, we provide a very detailed qualitative and quantitative description of the local behavior around rationally indifferent periodic points and the structure of corresponding Leau--Fatou flower petals including the Fatou Flower Petal Theorem. 

We also bring up the definitions of Speiser class $\cS$ and Eremenko--Lyubich class $\cB$, which play a seminal role in the recent development of the theory of iteration of transcendental meromorphic functions, and we prove some structural theorems about their Fatou components. 

In the last section of this chapter, Section~\ref{NiceSetsGeneral}, Nice Sets for Analytic Maps, we introduce and thoroughly study the objects related to the powerful concept of nice sets, which will be our indispensable tool in last part of the book, Part~\ref{EFB}, Elliptic Functions B. 

Up to our best knowledge there is no systematic book account of topological dynamics of transcendental meromorphic functions. Some results, with and without proofs, can be found in \cite{BKL1}--\cite{BKL1} and in \cite{Be}. 
Essentially all results in the current Chapter~\ref{topological-picture} of our book are supplied with proofs.

\section{Basic Iteration of Meromorphic Functions}\label{Basic iteration of meromorphic functions} 

In this section we define Fatou and Julia sets of meromorphic function. We also classify all periodic points of such functions. We prove some basic, rather elementary facts, of all of them.

\sp Let $f:\mathbb{C} \lra \oc$ be a transcendental
meromorphic function. As always in this book have $f^n$ is to be the $n$--th iterate of $f$\index{(N)}{iterate of function}, that is, 
$$
f^0(z)=z \  \text{ and } \ 
f^n(z)=f(f^{n-1}(z)) \  \text{ for all } \  n \geq 1.
$$
Then all values $f^n(z)$ are well--defined  for
all $z\in \mathbb{C}$ except for a countable set which consists of
the poles of $f$, $f^2$, $\ldots$,$f^{n-1}$. We call them prepoles
i.e. formally:

\sp\bdfn\lab{d120200311}
Given a natural number $n\ge 0$, we define the prepoles\index{(N)}{prepole} of order $n \geq 0$ as
$$
f^{-n}(\infty):=\{z\in {\mathbb C}:  f^n(z) \  \text{is well--defined and } \ f^n(z)=\infty\}.
$$
\index{(S)}{$f^{-n}(\infty)$} 
\edfn 
\fr Note that a pole is just an order $1$ prepole and $\infty$ is the sole prepole of order $0$. We shall prove the following.

\sp\blem\label{r12012} 
For every integer $n\ge 0$ the set of accumulation points of $f^{-n}(\infty)$ is contained in 
$$
\{\infty\}\cup \bigcup_{k=0}^{n-1} f^{-k}(\infty).
$$ 
\elem

\bpf Fix $n\in\mathbb N$. Of course $\infty$ is an accumulation points of $f^{-n}(\infty)$. Let $z\in \C$  be an  accumulation of $f^{-n}(\infty)$. Then there exists a sequence $(z_k)_{k=1}^\infty$ of mutually distinct in the set $f^{-n}(\infty)$ such that
$$
z=\lim_{k\to \infty}z_k.
$$
Seeking contradiction suppose that 
$$
z \notin \bigcup_{k=0}^{n-1}f^{-k}(\infty).
$$
Then there  exists a neighborhood $U$ of $z$  such that
$f^n$ restricted to $U$ is a meromorphic function. But this is
a contradiction since $f^n$ has poles at the points $z_k$. $k\ge$
and infinitely many of them belong to $U$. 
 
For the converse, of course $\infty$ is an accumulation points of $f^{-n}(\infty)$. Let 
$$
w\in \bigcup_{k=0}^{n-1} f^{-k}(\infty).
$$
Then there exists $0\le k\le n-1$ such that
$$
w\in f^{-k}(\infty).
$$
Since $f$ is an elliptic function, we then have for every $r>0$ that 
$$
f^{k+1}\lt(B(w,r)\sms \bu_{j=0}^k\rt)=\oc.
$$
Hence 
$$
B(w,r)\cap f^{-(k+1)}\(f^{-(n-(k+1))}(\infty)\)\ne \es.
$$
This means that
$$
B(w,r)\cap f^{-n}(\infty)\ne \es.
$$
So, $w$ is an accumulation points of $f^{-n}(\infty)$, and the proof is complete.
\endpf

\bdfn\lab{Fatou} The Fatou set\index{(N)}{Fatou set}
$F(f)$\index{(S)}{$F(f)$} of a meromorphic function $f: {\mathbb C}
\lra \oc$ is defined in exactly the same manner as for
rational functions: $F(f)$ is the set of all points $z \in {\mathbb C}$
for which all the iterates $f^n$ of $f$ are defined and form a normal family on some neighborhood of $z$. 
\edfn

\bdfn\label{Julia}
The {\em Julia set}\index{(N)}{Julia  set} $J(f)$\index{(S)}{$J(f)$} of a meromorphic function $f: {\mathbb C}\lra \oc$  is defined to be the complement of $F(f)$ in $\oc$, i.e $J(f)=\oc\sms F(f)$. 
\edfn

\fr Thus, the Fatou set $F(f)$ is open while the Julia set $J(f)$ is closed. 
We adopt the following definition.
 
\bdfn
Let $f: {\mathbb C}\lra \oc$ be a  meromorphic function.
\sp \begin{enumerate}
\item A number $\a\in \hat{\mathbb{C}}$ is called an asymptotic value of
$f$, if there exists a curve $\g:[0, 1)\lra \mathbb{C}$ such that
$\lim_{t \to 1} \g(t)=\infty$ and $\lim f(\g(t))=\a$. 

\,\item A complex number $\alpha\in \hat{\mathbb C}$ is called an
{\em omitted value}\index{(N)}{omitted value} of the meromorphic
function $f$, if $f(z)\neq \alpha$ for all $z \in \mathbb{C}$. Every
omitted value is {\em an asymptotic value}.\index{(N)}{asymptotic
value} 
\end{enumerate}
\edfn

We exclude from our considerations the class of meromorphic functions with exactly one  pole which is an  omitted value. 

\

\fr For all $z\in\ov{\mathbb{C}}$, we define 
$$
O^{+}(z):=\{f^{n}(z):n\ge 0\}
$$
\index{(S)}{$O^{+}(z)$} the {\em the forward orbit} of $z$ \index{(N)}{forward orbit}, and  
$$
O^{-}(z):=\bigcup_{n=0}^{\infty}f^{-n}(z),
$$ 
{\em the backward orbit} of $z$ \index{(N)}{backward orbit}. Then the above exclusion is
equivalent to  saying that 

\sp
\centerline{$O^{-}(\infty)$\index{(S)}{$O^{-}(z)$} is
an infinite set.}

\sp We note that it easily follows from Montel's
Criterion of normality that for considered  class of meromorphic
functions
\beq\label{density-prepoles}
 J(f)=\ov{\bu_{n\ge 0}f^{-n}(\infty)}.
\eeq
We now recall the basic properties of the Fatou set and Julia set. It is directly seen from the definitions that the Fatou set $F(f)$ is completely invariant while
$$
f^{-1}(J(f))\subset J(f) \  \text{ and }  \ f(J(f)\sms  \{\infty\})=J(f).
$$
We shall prove the following.

\bthm\label{duality}
If $f: {\mathbb C}\lra \ov{ \mathbb C}$ \index{(S)}{$J(f)$} is a  meromorphic function, then either $J(f)=\oc$ \ or $J(f)$ has empty interior (is nowhere dense).
\ethm

\bpf Suppose that $J(f)$ has non-empty interior. Denote this interior by
$ W$. Let $W_0= \bigcup_{n=0}^{\infty} f^n(W)$. If $\oc\sms W_0$
contains three distinct points, then Montel's Theorem yields that thefamily of iterates
$\{f^n|_W\}$ is normal, and therefore $W \sbt F(f)$. This is a 
contradiction, and thus $\oc\sms W_0$ contains at most two points. Hence, $J(f)\spt \ov{W_0}=\oc$. The proof is complete.
\endpf

\sp We say that $z\in\C$ is exceptional\index{(N)}{exceptional point} if
$O^{-}(z)$ is finite. The transcendental meromorphic  function has
at most two exceptional values. Again, Montel's Theorem implies
that, if $z$ is not exceptional and $z\in J(f)$, then
$$
J(f)=\ov{O^{-}(z)}.
$$ 
We recall that a set is perfect if it contains no isolated points. We shall prove the following.

\bthm\label{t_perfect_julia}
If $f: {\mathbb C}\lra \ov{ \mathbb C}$ $J(f)$\index{(S)}{$J(f)$} is a  meromorphic function, then the Julia set $J(f)$ is perfect.
\ethm

\bpf Fix $z \in J(f)$ arbitray. Let $U$ be an open
neighbourhood of $z$. As $O^{-}(\infty)$ is infinite, we can find three mutually distinct points $z_1, z_2, z_3 \in O^{-}(\infty)\sms O^{+}(z)$. Since the sequence $\{f^n_{|U}\}$ is not normal, $ z_j \in O^+(U)$ for at least one $j \in \{1, 2, 3\}$. Hence
$O^{-}(z_j) \cap (U\sms \{z\})\ne\es$. As $O^{-}(z_j)\sbt J$, this entails that $J\cap (U\sms \{z\})\ne\es$. Hence $z$ is not isolated in $J$. Thus
$J(f)$ is perfect, and the proof is complete.  
\endpf

\sp Now we provide the classical classification of periodic
points.\index{(N)}{periodic point} A point $\xi\in\C$ is called 
periodic if 
$$
f^p(\xi)=\xi
$$ 
for some $p\ge 1$. In this case the number $p$ is called a period of $\xi$, and the smallest $p$ with this property is called the minimal (or prime) period of $\xi$\index{(N)}{period of point}. If $p=1$, then $\xi$ is also called (naturally) a fixed point of $f$. We denote by 
$$
\Per(f), \  \  \Per_p(f), \  \  {\rm and }  \  \  \Per_p^*(f),
$$
respectively the set of all periodic points of $f$, all periodic points of $f$ of period $p$, and all periodic points of $f$ of period prime $p$. If $\xi$ is a periodic point of $f$ of prime period $p$, then the complex number 
$$
(f^p)'(\xi)
$$ 
is called the multiplier of $\xi$\index{(N)}{multiplier}. We classify periodic points of $f$ as follows.

\bdfn
Let $\xi$ be a periodic point of a meromorphic function $f:\C\to\oc$ with minimal period $p\ge 1$. The periodic point $\xi$ is called  

\begin{enumerate}
\item {\em attracting}\index{(N)}{attracting periodic point}, 

\,\item {\em supper attracting}\index{(N)}{supper attracting periodic point} (of course being supper attracting yields attracting), 

\,\item {\em indifferent} (or neutral), \index{(N)}{indifferent periodic point}  or

\,\item {\em repelling}\index{(N)}{repelling periodic point} 
\end{enumerate}

\fr respectively as the modulus of its multiplier is less than $1$, equal to $0$, equal to $1$, or greater than $1$. 
\edfn

\bdfn
\sp\fr\fr Writing the multiplier of an indifferent periodic point $\xi$ in the form $e^{2\pi i \alpha}$ where $0\leq \a <1$, we say that 

\begin{itemize}
\item[(a)] $\xi$ is {\em rationally indifferent (parabolic)} \index{(N)}{rationally indifferent (parabolic) periodic point} if $\a$ is rational  and 

\,\item[(b)] $\xi$ is {\em irrationally indifferent} if $\a$ is irrational. \index{(N)}{irrationally indifferent periodic point} 
\end{itemize}
\edfn

Since for every sufficiently small $r>0$, $f^p(B(\xi,r))\sbt B(\xi,r)$, we get the following.

\bthm
Each attracting periodic point of a meromorphic function $f:\C\to\oc$ belongs to the the Fatou set $F(f)$ of $f$.
\ethm

\fr For repelling periodic points just the opposite is true.

\bthm\label{repelling periodic}
Each repelling periodic point of a meromorphic function $f:\C\lra\oc$ belongs to the Julia set $J(f)$ of $f$.
\ethm

\bpf
Seeking contradiction suppose that a repelling periodic point $\xi$ of $f$ belongs to the the Fatou set $F(f)$. Denote the minimal period of $\xi$ by $p$. Let $U$ be such an open neighborhood of $\xi$ that the family $\(f^n|_U\)_{n=0}^\infty$ is normal. Then there exists a meromorphic function $g:U\to\oc$ to which some sequence $\(f^{pk_n}|_U\)_{n=0}^\infty$ converges uniformly on compact subsets of $U$. Then $g(\xi)=\xi$, and therefore, $g$ is holomorphic, and in particular, $g'(\xi)\in\C$. But on the other hand 
$$
|g'(\xi)|
=\lim_{k\to\infty}\big|\(f^{pk_n}\)'(\xi)\big|
=\lim_{k\to\infty}|(f^p)'(\xi)|^{k_n}
=+\infty
$$ 
as $|(f^p)'(\xi)|>1$. This contradiction finishes the proof.
\epf

\bthm
Each rationally indifferent periodic point of a meromorphic function $f:\C\lra\oc$ belongs to the Julia set $J(f)$ of $f$.
\ethm

\bpf
Changing coordinates by a translation we may assume without loss of generality that this periodic point is equal to $0$. Passing to a sufficiently high iterate, we my further assume that $0$ is a simple parabolic fixed point of $f$. Then the Taylor series expansion of $f$ about $0$ takes on the form
$$
f(z)=z+az^{p+1}+ \text{ higher terms of }z
$$
where $a\ne 0$ and $p\ge 1$ is an integer. We shall show by induction that
$$
f^n(z)=z+naz^{p+1}+ \text{ higher terms of }z.
$$
Indeed, this is of course true for $n=1$. Assuming its truth for some $n\ge 1$, and denoting higher than $k\ge 0$ terms of $w$ (a power series of $w$ starting with $w^{k+1}$)by $HT_k(w)$, we get
$$
\begin{aligned}
f^{n+1}(z)
&=f(f^n(z))
=f^n(z)+a(f^n(z))^{p+1}+ HT_{p+1}(f^n(z)) \\
&=f^n(z)+a(f^n(z))^{p+1}+ HT_{p+1}(z)\\
&=z+naz^{p+1}+HT_{p+1}(z) + a\(z+naz^{p+1}+HT_{p+1}(z)\)^{p+1}+HT_{p+1}(z)\\
&=z+(n+1)az^{p+1}+HT_{p+1}(z).
\end{aligned}
$$
The inductive proof is complete. Consequently,
$$
(f^n)^{(p+1)}(0)=a(p+1)!n.
$$
Therefore, $\lim_{n\to\infty}|(f^n)^{(p+1)}(0)|=+\infty$, and as $f(0)=0$, the proof can be now concluded in the same way as the proof of Theorem~\ref{repelling periodic}.
\epf

\bdfn\label{d-Siegel MU}
Let $\xi$ be a periodic point of a meromorphic function $f:\C\lra\oc$ with minimal period $p\ge 1$. The map $f^p$ is called {\em linearizable}\index{(N)}{linearizable map} near the periodic point
$\xi$ if and only if $f^p$ is topologically conjugate to its differential
$$
z\lmt g(z):=\xi+ (f^p)'(\xi)(z-\xi)
$$ 
in some (sufficiently small) neighborhood of $\xi$. 
\edfn

\bthm\label{Siegel MU}
An irrationally neutral periodic point $\xi$ of a meromorphic function $f:\C\lra\oc$ with minimal period $p\ge 1$, belongs to the Fatou set $F(f)$ if and only if
$f^p$ is linearizable near $\xi$. If this holds, the point $\xi$ is called a {\em Siegel periodic point}\index{(N)}{Siegel periodic point} of $f$. The corresponding topological conjugacy then also yields a holomorphic one.
\ethm

\bpf
Replacing $f$ by $f^p$ we may assume without loss of generality that $p=1$, i.e. that $\xi$ is a fixed point of $f$. Furthermore, changing coordinates by a translation we may assume without loss of generality that $\xi=0$. Write
$$
f'(0):=\g, \ \  |\g|=1.
$$
First, we assume that $f$ is linearizable near $0$. This means that 
$$
H\circ f\circ H^{-1}(z)= \g z, \  \  z\in D,
$$
where $D$ is a sufficiently small disk centered at $0$ and $H:D\to D$ is a homeomorphism. Iterating this equation we get for every $n\ge 0$ that
$$
H\circ f^n\circ H^{-1}(z)= \g^nz, \  \  z\in D,
$$
Equivalently
$$
f^n(z)=H^{-1}(\g^nH(z)), \  \  z\in D.
$$
In particular, $f^n(D)\sbt D$ for every $n\ge 0$. Therefore the family of iterates $\(f^n|_D\)_{n=0}^\infty$ is normal, so $0\in F(f)$.

\sp We now assume that $0\in F(f)$. Then, there is a
neighbourhood of $0$ on which the sequence  $(f^n)_{n=0}^\infty$ is equicontinuous and from  this, we  see that there exists some ball $B(0,r)$ ($0<r<1$), of $0$ such that for all $n\ge 0$ and all $z\in U$, we have that
\beq\label{e19.06.2012}
|f^n(z)| =  | f^n(z)- f^n(0)| < 1.
\eeq
Now, for every $n \geq 1$, define a function $T_n:B(0,r)\to\C$ by the formula
$$
T_n(z):=n^{-1}\(z + \g^{-1}f(z) + \g^{-2}f^2(z)\ldots +  \g^{-(n-1)}f^{n-1}(z)\).
$$
Note that as $|\g|=1$, we have that
\beq\label{e29.06.2012}
T_n(B(0,r))\sbt B(0,1)
\eeq
for every $n \geq 1$. A direct verification shows that the functions $T_n$ satisfy the following.
$$
(n/\g)T_n(f(z))+ z = ( n+1) T_{n+1}(z)= n T_n(z)+ \g^{-n}f^n(z),
$$
Hence,
$$
T_n(f(z))- \g T_n(z)=n^{-1}(\g^{1-n}f^n(z)-\g z).
$$
Since $|\g|=1$ and invoking (\ref{e19.06.2012}), we thus conclude that
\beq\label{e39.06.2012}
T_n(f(z))- \g T_n(z) \to 0
\eeq
uniformly on $B(0,r)$  as $n\to\infty$. Next (\ref{e29.06.2012})
implies that the sequence $\(T_n|_{B(0,r)}\)_{n=1}^\infty$ is normal on $U$ and it thus follows  that there exists $(k_n)_{n=1}^\infty$, an increasing sequence of positive integers, such that $T_{k_n}:B(0,r)\lra B(0,1)$  converges locally uniformly on $B(0,r)$  to some holomorphic function $H:B(0,r)\lra B(0,1)$. By (\ref{e39.06.2012}), satisfies  
$$ 
H\circ f(z)= \g H(z)
$$
for all $z\in U$. Since $T_n'(0)=1$, we also have  $H'(0)=1$ so $H$ is a homeomorphism on a sufficiently small neighborhood of $0$. This completes the proof. 
\epf

\bdfn
An irrationally neutral periodic point $\xi$ of a meromorphic function $f:\C\lra\oc$ with minimal period $p\ge 1$ belonging to $J(f)$ is called a {\em Cremer}\index{(N)}{Cremer  periodic point} periodic point of $f$. If $p\ge 1$ denotes the prime period of $\xi$, then near $\xi$ the map $f^{p}$ is not topologically conjugate to its differential (see Theorem~\ref{Siegel MU}).
\edfn

\fr From all the above we have the following.

\bthm
All attracting and Siegel periodic points of  a meromorphic function $f:\C\lra\oc$ are in the Fatou set of $f$, while  repelling, rationally indifferent, and Cremer periodic points are in the Julia set of $f$. 
\ethm

Also a point $\xi\in\oc$ is called preperiodic  if $f^n(\xi)$ is periodic \index{(N)}{preperiodic point} for some $n \geq 0$.

\sp If $f$ is transcendental meromorphic  function and $n \geq 2$, then $f$  has  infinitely many periodic points of mininimal period $n$. In  fact, $f$ has
 infinitely many  repelling periodic points of minimal period $n$ and the Julia
 set of $f$ is the closure of the set of all repelling periodic points of $f$. We shall now prove this result. Moreover, the Julia set of $f$ will turn out to be perfect.

\sp

\bthm\label{density}
If $f:\mathbb C  \to \ov{\mathbb C}$ is a transcendental  meromorphic function, then $J(f)$, the Julia sets of $f$, is the closure of the set of all repelling periodic points of $f$.
\ethm

\sp

In order to prove Theorem~\ref{density} we need the well-known {\it 'five island theorem'} of Ahlfors (see \cite{Ah}) from complex analysis.

\sp

\fr\bthm\label{Ahlfors}  {\rm (Ahlfors' Five Island
Theorem).}\index{(N)}{Ahlfors' Theorem} Let $f:\mathbb C \to
\ov{\mathbb C}$ be a transcendental meromorphic function, and let
$D_1,D_2, \ldots D_5$ be any five simply connected domains in $\mathbb C$
with mutually disjoint closures. Then there exists at least one $j\in \{1, 2, \ldots, 5\}$
and for every $R>0$ there exists a simply connected domain $G \sbt \{z \in \mathbb
C: \, |z| > R\}$ such $f|_G$ maps conformally $G$ onto $D_j$. If $f$ has only finitely many poles, then 'five'  may be replaced by 'three'. 
\ethm

\

\fr The next lemma follows from  Theorem~\ref{Ahlfors}.

\

\blem\label{l6in{Be}} Suppose that  $f:\mathbb C  \to \ov{\mathbb
C}$ is a transcendental  meromorphic function and that some five points $z_1,  z_2,
\ldots, z_5\in O^{-1}(\infty)\sms \{\infty\}$ are mutually distinct. Define
$n_j\ge 1$ uniquely by the property that $f^{n_j} (z_j)=\infty$. Then there exists  $j \in \{1,
\ldots, 5\}$ such that $z_j$ is a limit point of repelling periodic
points of $f$ with minimal period equal to $n_j+1$. If $f$  has only finitely
many poles then 'five'  may be replaced by 'three'.\elem

\bpf in order to deduce Lemma~\ref{l6in{Be}}  from
Theorem~\ref{Ahlfors} we choose the $D_j$ to be discs around $z_j$
where the radii are chosen so small that the $D_j$ do not contain
critical points of $f^{n_j}$ nor other poles of $f^{n_j}$ apart from $z_j$, and that their closures are pairwise disjoint. Since each $z_j$ is a pole of $f^{n_j}$,
there exists $R>0$ such that 
\beq\label{MU_2015_08_24}
\bu_{k=0}^{n_{j-1}}f^k(D_j)\sbt B(0,R/4) \  \  \text{ and } \  \ f^{n_j}\lt(\frac12D_j\rt)\supset \{z: |z|>R/2 \}\cup \{\infty\}
\eeq
for all $j=1, 2, \ld,5$. We pick $j$ and $G$  according to Ahlfors' Five Island Theorem, i.e. Theorem~\ref{Ahlfors}.
Then we can find an open connected set $H\sbt (1/2)D_j$ such that  $f^{n_j}(H)=G$. 
So, the map $f^{n_j+1}|_H:H \lra  D_j$ is a holomorphic surjection without critical points.
Since the open set $D_j$ is simply connected, the Modromy Theorem implies that the inverse map
$$
f_*^{-(n_j+1)}:D_j\lra H\sbt \frac12D_j
$$
is well-defined, and of course holomorphic. By Brouwer's Fixed Point Theorem (an alternative argument would be to use contraction of $f_*^{-(n_j+1)}$ in hyperbolic metric on $D_j$ and to apply Banach Contraction Principle) there exists a fixed point $\xi$ of $f_*^{-(n_j+1)}$ in $(1/2)D_j$, and by Schwarz Lemma this point is attracting. Hence, $\xi$ is a repelling fixed point of $f^{n_j+1}$. In other words, $\xi$ is a repelling  periodic point of $f$ of period
$n_j+1$, and by \eqref{MU_2015_08_24} its minimal period is equal to $n_j+1$. Because the disks $D_j$ can be chosen arbitrarily small, we therefore now conclude $z_j$ is the limit
of a sequence of repelling  periodic points with minimal period equal to $n_j+1$. This completes the proof of  the lemma. 
\endpf

\sp

\fr Since (see \eqref{density-prepoles}) for our transcendental meromorphic function $f$, the Julia set $J(f)$ is the closure of $O^{-1}(\infty)$, Theorem~\ref{density} directly follows from Lemma~\ref{l6in{Be}}.

\sp\section[Fatou Components of (General) Meromorphic Functions I] {Fatou Components of (General) Meromorphic Functions I; Classification of Periodic Components}

In this section having a meromorphic function $f:\C\lra\oc$, we want to describe and analyze in some detail the structure of connected components of the Fatou set $F(f)$. Let $U$ be such a component. Then for every $n\ge 0$, the open set $f^n(U)$ is contained in some connected component of $F(f)$, which we denote by $U_n$. 

\bdfn A connected component $U$ of $F(f)$ is called  preperiodic \index{(N)}{preperiodic component} if there exist $n> m\geq 0$ such that 
$$
U_n=U_m.
$$
In particular, if $m=0$, then  $U$ is called
periodic with period $n$.\index{(N)}{periodic component} The finite set 
$$
\{U_0, U_1, \ldots, U_{n-1}\}
$$ 
is then called  a (periodic)  cycle of
components. The smallest $n\ge 1$  with this property is called the
minimal period of $U$. In the case when $n=1$, that is, if 
$$
f(U)\subset U,
$$
the component $U$ is called $f$--invariant\index{(N)}{invariant component} or just invariant. 

\sp\fr A component of $F(f)$ which is not preperiodic is called a wandering
component.\index{(N)}{wandering components}
\edfn

\

\fr If $\xi$ is an attracting periodic point of $f$ with minimal period $p\geq1$, then for all sufficiently small $R>0$, we have that
\[
f^p(B(\xi, R))\subset B\left(\xi,\frac{1+|(f^p)'(\xi)|}{2}R\right) \subset B(\xi, R).
\]
Thus, by Montel's Theorem, it follows that $B(\xi,R)\subset F(f)$. So, if we define
\[
A(\xi):=\big\{z\in\C:\lim_{n\to\infty}f^{pn}(z)=\xi\big\},
\]
then $B(\xi, R)\subset A(\xi)\subset F(f)$, and, furthermore, if we denote by $A^*(\xi)$ the connected component of $A(\xi)$ that contains $\xi$, then $B(\xi, R)\subset A^*(\xi)$. Since all limit points of iterates of $f$ on $A(\xi)$ are constant functions (with values in $\{\xi, f(\xi), \ldots, f^{p-1}(\xi)\}$, we conclude that no point on the boundary of $A^*(\xi)$ may belong to the Fatou set $F(f)$. Thus the set $
A^*(\xi)$ is a connected component of the Fatou set $F(f)$. We collect these observations in the following theorem.

\bthm\label{t1ms113}
If $\xi\in\C$ is an attracting periodic point of a meromrphic function $f:\C\lra\oc$, then $A^*(\xi)$ and $A(\xi)$ are open sets,
\[
\xi\in A^*(\xi)\subset A(\xi) \subset F(f),
\]
and $A^*(\xi)$ is a connected component of the Fatou set $F(f)$.
\ethm

\fr The set $A(\xi)$ is called the \index{(N)}{basin of attraction}{\em basin of attraction to $\xi$}, while the set $A^*(\xi)$ is called the \index{(N)}{basin of attraction!immediate}{\em immediate basin of attraction to $\xi$}. If $p\geq 1$ is the minimal period of $\xi$, we call the sets
\begin{eqnarray}\label{1ms115}
A^*_p(\xi):=\bigcup_{j=0}^{p-1}A^*(f^j(\xi))\ \text{ and }\ A_p(\xi):=\bigcup_{j=0}^{p-1}A(f^j(\xi)),
\end{eqnarray}
which are both contained in the Fatou set $F(f)$, respectively the basin of immediate attraction and the basin of attraction to the periodic orbit $\{\xi, f(\xi), \ldots, f^{p-1}(\xi)\}$.

\sp Now we formulate and will prove the fundamental classification theorem of all periodic components of the Fatou set $F(f)$. This theorem is essentially due  to Cremer~\cite{Cre} and Fatou~\cite{F1}. Fatou~\cite{F1} [$\S$ 56, p. 249] proved that if $\{f^n|_V\}$  has only constant limit functions, then  $V$ is an immediate attractive basin or a  Fatou--Leau domain, provided $f$ is rational. His proof shows that shows that the only further possibility in the case of transcendental functions is that of  a Baker domain. Cremer~\cite{Cre} [p. 317] proved that if $\{f^n|_V\}$ has non--constant limit functions, then $V$  is a Siegel disk or  a Herman ring. Neither Fatou nor Cremer stated the full classification theorem, but T\"opfer remarks~\cite{T}  [p. 211] come fairly close to it.


\sp\bthm[Fatou Periodic Components]\label{Fatou Periodic Components}
Let $U$ be a periodic connected component of the Fatou set $F(f)$, of some period $p\ge 1$, of a meromorphic function $f:\C\lra\oc$. Then we have one of the following possibilities:

\sp \begin{enumerate}
\item $U$ contains an attracting  periodic point $\xi$ of period
$p$. Then 
$$
\lim_{n\to\infty}f^{np}(\xi)=\xi
$$ 
for all $z\in U$, and, we recall, $U=A^*(\xi)$ is called the immediate  attracting  basin\index{(N)}{immediate attracting  basin} of $\xi$.

\sp\item $\partial{U}$  contains a periodic point $\xi$ of period $p$
and 
$$
\lim_{n\to\infty}f^{np}(\xi)=\xi
$$ 
for all $z\in U$. Then $(f^{p})'(\xi)$ is a root of unity. In this case, $U$ is called a Leau--Fatou domain of $\xi$.\index{(N)}{Leau--Fatou domain}

\sp\item There exists an analytic homeomorphism $H:B(0,1)\lra U$ 
such that the following diagram commutes:
\[\begin{tikzcd}
B(0,1) \arrow{r}{R_\a} \arrow[swap]{d}{H} & B(0,1)\arrow{d}{H} \\
{U} \arrow{r}{f^p} & {U}
\end{tikzcd}
\]
i.e.
$$
H^{-1}(f^p(H(z)))=e^{2\pi i\a}z, \  \  z\in B(0,1),
$$  
for some $\a\in \mathbb{R}\sms \mathbb{Q}$, where $R_\a:B(0,1)\lra  B(0,1)$ is the rotation $R_\a(z)=e^{2\pi i\a}z$. In this case, $U$
is called a Siegel disk of $f$.\index{(N)}{Siegel disc}

\sp\item There exists an analytic homeomorphism $H:A(0;1,r)\lra U$, $r>1$, such that the following diagram commutes,:
\[\begin{tikzcd}
A(0;1,r) \arrow{r}{R_\a} \arrow[swap]{d}{H} & A(0;1,r)\arrow{d}{H} \\
{U} \arrow{r}{f^p} & {U}
\end{tikzcd}
\]
i.e.
$$
H^{-1}(f^p(H(z)))=e^{2\pi i \a}z, \  \  z\in A(0;1,r),
$$ 
for some $\a\in \mathbb{R}\sms \mathbb{Q}$. In this case, $U$ is called a Herman ring of $f$.\index{(N)}{Herman ring}

\sp\item There exists $\xi\in \partial{U}$ such that 
$$
\lim_{n\to\infty}f^{np}(z)=\xi
$$
for all $z\in U$ but there is no continuous extension of $f^p$ from $U$ to $\xi$. In this case $U$ is called a Baker domain of $f$.\index{(N)}{Baker domain}
\end{enumerate}
\ethm

\sp This classification theorem follows from Theorem~\ref{Limit_functions} below. The reader may also consult the books \cite{CG}, \cite{Bea}, \cite{steinmetz}, and \cite{Mi1} for the case of rational functions. 

We start with the following.

\sp\bprop\label{p1_2017_07-26}
If $f:\C\lra\oc$ is a meromorphic function, $V$ is a connected component of the Fatou set $F(f)$, and the sequence $\(f^n|_V\)_{n=0}^\infty$ has a non--constant limit function, then the component $V$ is preperiodic.
\eprop

\bpf 
It follows from the hypotheses of our proposition that there exists some  increasing sequence $(n_j)_{j=1}^\infty$ of positive integers such that the sequence $\(f^{n_j}\)_{j=1}^\infty$ converges to some non--constant analytic function 
$$
g:V\lra\C
$$ 
uniformly on each compact subset of $V$. Then, there exists a point $\xi\in V$ such that $g'(\xi)\neq 0$. Thus, there exists $r>0$ such that $B(\xi,2r)\sbt V$ and 
$$
g(z)\neq g(\zeta)
$$  
for all $z\in \bd $. Then, keeping such a point $z$, we will have for all integers $j\ge 1$ large enough that
$$ 
\begin{aligned}
\Big|\(f^{n_j}(z)-g(\xi)\)-\(g(z)-g(\xi)\)\Big| 
&=\Big|\(f^{n_j}(z)- g(z)\)\Big| 
< \inf\big\{|g(w)-g(\xi)|:w \in \bd B(\xi,r)\big\}\\
& \le |g(z)-g(\xi)|.
\end{aligned}
$$
Hence, by  Rouch\'e's Theorem, there exists a point $\zeta\in B(\xi,r)$ such that 
$$
f^{n_j}(\zeta)=g(\xi).
$$
But then also $f^{n_{j+1}}(\zeta)=g(\xi)$. So,
$$
f^{n_j}(V)\cap f^{n_{j+1}}(V)\ne\es
$$
Hence,
$$
f^{n_{j+1}}(V)\sbt f^{n_j}(V),
$$
and the proof is complete.
\epf

\blem\label{l120190625}
Let $f:\mathbb C\lra \oc$ be a non--affine (different from any map of the form $z\mapsto az+b$) meromorphic function. 
Suppose that $D$ is a forward $f$--invariant connected component of the Fatou set $F(f)$ of $f$. 

If there exists some non--constant limit function of the sequence $\(f^n|_D\)_{n=0}^\infty$, then the function $f|_D:D\lra D$ is a conformal homeomorphism, and the identity map $\Id_D$ on $D$ is a limit point of the sequence $\(f^n|_D\)_{n=0}^\infty$.
\elem

\bpf Let $g:D\lra\oc$ be a non--constant analytic limit function of the sequence $\(f^n|_D\)_{n=0}^\infty$; by our hypotheses at least one such a function exists. We first shall show that
\begin{equation}\label{7.2.1}
g(D)\sbt D.
\end{equation}
Indeed, by the definition of $g$ there exists $(n_j)_{j=1}^\infty$, an unbounded increasing sequence of positive integers such that 
$$
f^{n_j}|_D\lra g
$$ 
uniformly on compact subsets of $D$ as $j\to\infty$. Fix any $w \in D$. Since $g$ is non--constant, the zeros of the map 
$$
D\ni z \lmt g(z)-g(w)\in\C
$$ 
are isolated. Hence, there exists $r>0$ such that
$$
g(z)\ne g(w)
$$
for all $z\in \bd B(w,2r)$. Since $\bd B(w,r)$ is a compact subset of $D$, we therefore conclude that 
$$
\inf\big\{|g(x)- g(w)|:x\in\bd B(w,r)\big\}>0
$$
and 
$$ 
\big|(f^{n_j}(z)-g(w)) - (g(z)-g(w))\big|
=|f^{n_j}(z) - g(z)| 
< \inf\big\{|g(x)- g(w)|:x\in\bd B(w,r)\big\}
$$
for all sufficiently large $j\ge 1$, say $j\ge k\ge 1$, and all $z$ in $\bd B(w,r)$. So, by Rouch\'e's Theorem, the two functions 
$$
D\ni z\longmapsto g(z) -g(w)\in\C
\  \  \  {\rm and} \  \  \ 
D\ni z\longmapsto f^{n_j}(z) -g(w)\in\C
$$ 
have the same number of zeros for all $j\ge k$. As the first of these two functions vanishes at $w$, thus the second one must also have at least one zero in $B(w,r)$. Denote it by $y_j$. So, $y_j\in D$ and $f^{n_j}(y_j) -g(w)=0$. Thus
$$
g(w)=f^{n_j}(y_j)\in f^{n_j}(D)\sbt D,
$$
and so (\ref{7.2.1}) holds.

By passing to a subsequence and relabelling, we may assume that
$$ 
\lim_{j\to\infty}(n_j-n_{j-1})=+\infty.
$$
For every $j\ge 1$ put 
$$
m_j:=n_j-n_{j-1}.
$$
Since the family $\(f^{m_j}|_D\)_{j=1}^\infty$ is normal, there are some function $h:D\to\C$ and some infinite set $N\sbt \N$ such that the sequence
$\(f^{m_j}|_D\)_{j\in N}$ converges to $h$ uniformly on compact subsets of $D$.

Now take any $z \in D$. Since $\lim_{j\to+\infty}f^{n_j}(z)=g(z)$ with $g(z)\in D$ by (\ref{7.2.1}), there exists $s>0$ such that 
$$
B(g(z),2s)\sbt D
$$
and 
$$
f^{n_j}(z)\in \ov{B(g(z),s)}
$$
for all $j\ge 1$ large enough. Therefore,
$$
\lim_{N\ni j\to+\infty}f^{m_j}=h
$$
uniformly on $\ov{B(g(z),s)}$. In conclusion,
\begin{equation}\label{7.2.2}
h(g(z))
=\lim_{N\ni j\to+\infty} f^{m_j}(f^{n_{j-1}}(z))
=\lim_{N\ni j\to+\infty} f^{n_j}(z) 
=g(z).
\end{equation}
As $g$ is not constant, $h$ must be the identity map $\Id_D$ on $D$, and this proves the second claim of our lemma.

The fact that he function $f|_D:D\lra D$ is a conformal homeomorphism follows now easily. Indeed, $f(D)=D$ because $D$ is forward invariant under $f$. In order to prove that $f$ is one--to--one assume that $a, b\in D$ and  $f(a)=f(b)$. Then,  
$$
f^{m_j}(a)=f^{m_{j}-1}(f (a))=f^{m_{j}-1}(f (b)) =f^{m_j}(b)
$$
and letting $j \to \infty$ in $N$, we obtain, by already proved part of our lemma, that 
$$
a=\Id_D(a)=\Id_D(b)=b.
$$
This completes the proof of Lemma~\ref{l120190625}. \qed

\blem\lab{l220190801}
For every $r\in(0,1)$, the only conformal homeomorphisms of the annulus 
$$
A_r:=\{z\in\C:r<|z|<1\},
$$
that are of infinite order are Euclidean rotations by angles which are irrational multiples of $\pi$.
\elem

\bpf
Let 
$$
\ka:=\log(1/r).
$$
We use notation of Proposition~\ref{p1j220}. In particular, the map $g:\mathbb H\lra \mathbb H$ is given by the formula
$$
g(z):=kz,
$$
where 
$$
k:=\exp\lt(\frac{2\pi^2}{\ka}\rt)
$$
and
$$
\Pi_\ka:\mathbb H \lra A_r
$$ 
is the analytic covering map of $A_r$ given by the formula
\beq\lab{120200128}
\Pi_\ka(z):=\exp\lt(\frac{\ka}{\pi}i\log z\rt)
\eeq
Let 
$$
G:A_r\lra A_r
$$ 
be an infinite order conformal homeomorphism of $A_r$. Let 
$$
\tilde G:\mathbb H\lra \mathbb H
$$
be a lift of $G$, i.e. such a continuous map that
$$
G\circ\Pi_\ka=\Pi_\ka\circ \tilde G,
$$ 
i.e. the following diagram commutes.
\[\begin{tikzcd}
\mathbb H \arrow{r}{\tilde G} \arrow[swap]{d}{\Pi_\ka} & \mathbb H\arrow{d}{\Pi_\ka} \\
{\  \C^*} \arrow{r}{G} & {\  \C^*}
\end{tikzcd}
\]
Then for any two points $z, w\in \mathbb H$ with 
$$
\Pi_\ka(z)=\Pi_\ka(w),
$$
we have that
$$
\Pi_\ka(\tilde G(z))
=G\circ\pi(z)
=G\circ\pi(w)
=\Pi_\ka(\tilde G(w)).
$$
Therefore, because of Proposition~\ref{p1j220},
$$
\tilde G(z)=g^m\tilde G(w)
$$
with some integer $m$. It in turn follows from this that for  each $z \in \mathbb H$ there  is an integer $m(z)$ such that 
\beq\label{120190802}
\tilde G(kz)=k^{m(z)}\tilde G(z).
\eeq
So, the function
$$
\mathbb G\ni z\lmt k^{m(z)}
=\tilde G(kz)\(\tilde G(z)\)^{-1}\in\C
$$
is continuous. But in fact this function takes values in $(0,+\infty)$ and 
$$
m(z)=\frac1{\log k}\log\(k^{m(z)}\)
=\frac1{\log k}\log\lt(\frac{\tilde G(kz)}{\tilde G(z)}\rt).
$$
So, we conclude that the  function 
$$
\mathbb H\ni z\lmt m(z)\in\Z
$$
is continuous. Therefore, since the space $\mathbb H$ is connected, the function $\mathbb H\ni z\longmapsto m(z)\in\Z$ is constant. Denote its only value by $m$. Formula \eqref{120190802} takes then on the form
\beq\label{220190802}
\tilde G(kz)=k^m\tilde G(z),
\eeq
and by immediate induction:
\beq\label{320190802} 
\tilde G(k^nz)=k^{mn}\tilde G(z)
\eeq
for all $z\in\mathbb H$ and all $n\in\Z$. Recall that $\tilde G$ is a M\"obius map and suppose for a moment that $m=0$. Letting $n \to +\infty$, it would then follow from \eqref{320190802} that $\tilde G$ is a constant function whose only value is $\tilde G(\infty)$, contrary to the fact that $\tilde G$ is an invertible M\"obius map. Thus,
$$
m\ne 0.
$$
Recalling again that $\tilde H$ is a M\"obius map and letting $n \to +\infty$, we conclude from \eqref{320190802} that $\tilde G(\infty)$ is either $\infty$ (if $m>0$), or $0$ (if $m <0$). By letting $ n \to -\infty$, we conclude that $\tilde G(0)$ is either $0$ or $\infty$. Therefore, either $\tilde G$ fixes both $0$ and $\infty$, or interchange them. This means that $\tilde G$ is one of the form
$$
\oc\ni z \lmt az\in\oc
\  \  \  {\rm or} \  \  \
\oc\ni z \lmt b/z\in\oc,
$$
where $a>0$ and $b<0$. Using the formula \eqref{120200128}, we get in the former case that
$$
\begin{aligned}
G(\Pi_\ka(z))
&=\Pi_\ka\circ \tilde G(z)
=\Pi_\ka(az)
=\exp\lt(\frac{\ka}{\pi} i\log(az)\rt) 
\\
&=\exp\lt(\frac{\ka}{\pi} i\log(a)\rt)\exp\lt(\frac{\ka}{\pi} i\log(z)\rt) 
\\
&=\exp\lt(\frac{\ka}{\pi} i\log(a)\rt)\Pi_\ka(z)
\end{aligned}
$$
for all $z\in \mathbb H$. Likewise, in the latter case.
$$
G(\Pi_\ka(z))
=\frac{\exp\lt(\frac{\ka}{\pi} i\log(b)\rt)}{\Pi_\ka(z)}
$$
for all $z\in \mathbb H$. Hence,
$$
G(w)=e^{i\a}w 
\  \  \  {\rm and} \  \  \
G(w)=e^{i\b}/w 
$$
with $\a:=\frac{\ka}{\pi}{\log(a)}\in\R$ and $\b:=\frac{\ka}{\pi}{\log(b)}\in\R$, for all $w\in A_r$ respectively in the former and the latter case. But the map $\oc\ni w\mapsto e^{i\b}/w\in\oc$ is of order $2$, whence, as $H$ is of infinite order, it follows that
$$
G(w)=e^{i\a}w 
$$
for all $w\in A_r$. The proof is complete.
\epf

\sp Now, we shall prove the following, already announced, theorem.

\bthm\label{Limit_functions}
Let $ f:\mathbb C\lra \oc$ be a non--affine (different from $z\mapsto az+b$) meromorphic function and let $V\sbt\oc$ be a periodic connected component of the Fatou set $F(f)$ of $f$.
 
\begin{itemize}
\item [(a)] If all the limit functions of the sequence $\(f^n|_V\)^\infty_{n=0}$ are constant and, in the case when $\infty$ is an essential singularity of $f$, each of them is different from $\infty$, then $V$ is either the immediate basin of attraction of an attracting periodic point of $f$ or $V$ is a Leau domain of $f$.

\,

\item [(b)] If the limit functions of the sequence $\(f^n|_V\)^\infty_{n=0}$ contain a non-constant function, then $V$ is either a Siegel disk or a Herman ring.

\,

\item [(c)] If $\infty$ is an essential singularity of $f$, all the limit functions of the sequence $\(f^n|_V\)^\infty_{n=0}$ are constant and some of them is equal to $\infty$, then $V$ is a Baker domain.
\end{itemize}
\ethm

\bpf  We first assume that $V$ is $f$--invariant, i.e. $f(V)\sbt V$.

\sp Item (a). If $\infty$ is a limit function, then by our hypothesis $f$ is a rational function, and in what follows we actually entirely work on the Riemann sphere $\oc$. Since $V$ is a connected component of the Fatou set $F(f)$, the family 
\beq\label{520190802}
\big\{f^n|_V: V \lra\oc\big\}_{n=0}^\infty
\eeq
is normal. Hence, the set $L$ of all limit functions of this family, which are all by our hypotheses constant, is non--empty. We shall show that each element of $L$  is a fixed point of $f$. Indeed, let $\xi\in L$. Then there exists $(n_k)_{k=1}^\infty$, an unbounded increasing sequence, such that
$$
\lim_{k\to\infty}f^{n_k}(z)=\xi
$$
for all $z \in V$. So, take any $w \in V$. Then $f(w)\in f(V)\sbt V$ and
$$ 
f(\xi)
=f\(\lim_{k\to \infty}f^{n_k}(w)\)
=\lim_{k \to\infty}f^{n_k}(f(w))
=\xi.
$$
Now again take any point $w\in V$. Since $V$ is arcwise connected, there exists a continuous map (even a homeomorphic embedding) 
$$
p:[0,1] \lra V
$$  
such that 
$$
p(0)=w
\  \  \  {\rm and} \  \  \  p(1)=f(w).
$$
Next, extend this function to a function
$$
p_\infty:[0, +\infty)\lra\oc
$$
by the formula:
$$
p_\infty(t):=f^{[t]}(p(t-[t])).
$$
it is immediate from this definition that the function $p_\infty$ is continuous. Since the set $[0, +\infty)$, the  set of the accumulation points of $p_\infty$ is a closed (so compact) connected subset of $\oc$. Thus, it is either a singleton or contains uncountably many points. Note that the set of accumulation points of $p_\infty$ coincides with the set of accumulation points of of the sequence $\(f^n(p([0,1])\)_{n=0}^\infty$ which is equal to $L$ as $p([0,1])$ is a compact subset of $V$. Thus, $L$ is either a singleton or contains uncountably many points. 

Seeking contradiction, suppose that the set $L$ is uncountable. Consider a meromorphic function  $h: \mathbb C \to \oc$ defined as 
$$
h(z):=f(z)-z.
$$
Then  $h(\xi)=0$ for all $\xi\in L$. So, $h$ has an uncountable set of  zeros and therefore is identically equal to zero. Thus $f(z)=z$ for all $z\in\C$, contrary to our hypotheses. Therefore the set $L$ is a singleton, and in accordance with our considerations up to now, denote its only element by $\xi$. Of course
$$
\xi\in\ov V.
$$ 
Suppose first that 
$$
\xi \in V.
$$
Fix a number $R>0$ such that
$$
\ov{B(\xi,R)}\sbt V.
$$
By the definition of $\xi$ there exists an integer $m\ge 1$ such that 
$$
f^m(\ov{B(\xi,R)})\sbt B(\xi,R/2).
$$
So, by Schwarz Lemma,
$$ 
|f'(\xi)|^m = |(f^m)'(\xi)|<1.
$$
Hence, 
$$ 
|f'(\xi)|<1,
$$
and $\xi\in V$ is an attracting fixed point of $f$, and, as we already know,  the full sequence 
$$
\big\{f^n|_V: V \lra\oc\big\}_{n=0}^\infty
$$ 
converges to $\xi $ on $V$ uniformly on compact subsets of $V$. Thus $V$ is the basin of immediate attraction of an attracting fixed point $\xi\in V$.

Suppose in turn that  
$$
\xi \in \partial{V}.
$$
Since the full sequence 
$$
\big\{f^n|_V: V \lra\oc\big\}_{n=0}^\infty
$$ 
converges to $\xi $ on $V$ uniformly on compact subsets of $V$ and $f(\xi)=\xi$, it remains to prove that 
$$
f'(\xi)=1. 
$$
Obviously, 
$$
|f'(\xi)|\geq 1 
$$ 
since otherwise $\xi$ would lie in $F(f)$. Seeking contradiction, suppose that $$
|f'(\xi)|> 1.
$$
Then, since also fixing any $z\in V$, we know that $\lim_{n\to\infty}f^n(z)=\xi$, we conclude that
$$
f^n(z)=\xi
$$ 
for some integer $n\ge 0$. But as $\xi\in J(f)$, so also $z\in J(f)$. This contradicts the fact that $z\in V\sbt F(f)$, and get that
$$
|f'(\xi)|=1.
$$
Conjugating (changing coordinates) by an affine function ($z\mapsto az+b$, $a\ne0$), we may assume without loss of generality that 
$$
\xi=0
$$ 
Let 
$$
\lambda:= f'(0).
$$
Then
$$
f(z)
=\lambda z+a_2z^2+a_3z^3+\ldots
=\lambda z+z^2H(z),
$$ 
where, for every and $R>0$ small enough, so that 
$$
B(0,2R)\sbt V,
$$
and
$$
H:B(0,R)\lra\C
$$
is some bounded holomorphic function. Since $ f'(0)\neq 0$, we may require $R>0$ to be so small that the map  
$$
f|_{B(0,R)}: B(0,R)\lra \mathbb C
$$
is one--to--one. Since $ f|^n_{V}\lra 0$ uniformly on compact subsets of V as $n\to+\infty$, and since $B(0,2r) \sbt V$, there exist $ w \in V\sms\{0\}$ and $r>0$ such that
\beq\label{1fpc1}
B(w, 2r)\cap B(0, 2r)=\es,
\eeq
$$
B(w, 4r)\sbt B(0,R),
$$
and moreover 
$$
f^n(B(w,4r))\sbt B(0,R)
$$
for every integer $n\geq 0$. In particular, the function
$$
f|^n_{B(w,4r)}: B(w, 2r) \lra B(0,R)
$$
is holomorphic for every integer $n \geq 0$ and  
$$
0 \notin f^n(B(w,4r))
$$
for every integer $ n \geq 0$. Therefore, we can, for every integer $ n \geq 0$, define a function
$$g
_n :B(w, 4r) \lra  \mathbb C
$$
by the formula
$$
g_n(z):=\frac{f^n(z)}{f^n(w)}.
$$
Of course, all functions  $g_n$ are holomorphic, injective,
$$
g_n(w)=1,\ 0, \infty \notin g_n(B(w, 4r)),  
\  \  \  {\rm and} \  \  \ 
g^n (1)=\{w\}.
$$
It thus follows from $\frac{1}{4}$-Koebe's Distortion Theorem that
$$
\frac{1}{4} |g_n'(w)|\leq 1.
$$
So, it follows from Koebe's Distortion Theorem that
$$ 
g_n(B(w, 2r)) \sbt B\(1, K|g_n'(w)|\) \sbt B(1, 4K)  \sbt B(0, 1+4K).
$$
Thus the sequence of holomorphic functions
$$
\left(g_n|_{B(w, 2r)}\right)^\infty_{n=0}
$$
is normal. Let  $(n_k)^\infty_{k=1}$ be any increasing sequence of positive integers such that
the sequence  $\left( g_{n_k}\right)^\infty_{k=1}$ converges uniformly on compact subsets of $B(w,2r)$, and let
$$
g:=\lim_{k \to \infty} g_{n_k}.
$$
Of course  
$$ 
g: B(w, 2r)\to \mathbb C
$$
is a holomorphic function. Since, for every $n \geq 0$  and every $z \in B(w,r)$,
$$ 
\begin{aligned}
g_{n+1}(z)
& =\frac{f(f^n(z))}{f(f^n(w))}=\frac{\lambda f^n(z)+(f^n(z))^2 H(f^n(z))}{\lambda f^n(w)+(f^n(w))^2 H(f^n(w))}\\
&=\frac{f^n(z)\(\lambda +f^n(z) H(f^n(z))\)}{f^n(w)\(\lambda +f^n(w) H(f^n(w))\)}\\
& = g_n(z)\frac{\lambda +f^n(z) H(f^n(z))}{\lambda +f^n(w) H(f^n(w))}.
\end{aligned}
$$
Since the function H is bounded, $ \lambda \neq 0$, and $\lim_{n \to \infty}f^n(z)=0 $, we therefore conclude that the sequence $\(f_{n_k+1}\)^\infty_{k=1}$  also converges uniformly on compact subsets of  $B(w, 2r)$, and
\begin{equation}\label{1fpc2}
\lim_{k \to \infty}g_{n_k+1}=g.
\end{equation}
Also,
$$
\begin{aligned}
g_{n+1}(z)& =\frac{f^n(f(z))}{f^n(w)}\cdot  \frac{f^n(w)}{f^{n+1}(w)}
   = \frac{f^n(f(z))}{f^n(w)}\cdot \frac{f^n(w)}{f^n(w)(\lambda +f^n(w) H(f^n(w)))}\\
   &=g_n(f(z))\cdot \frac{1}{\lambda +f^n(w) H(f^n(w))}.
   \end{aligned}$$
Again, since $H$ is bounded, $\lambda\neq0$, and  $ \lim_{n \to \infty}f^n(w)=0$, using   (\ref{1fpc2}), we conclude that
$$
g(z)=\lim_{k \to \infty}g_{n_k+1}(z)=g(f(z))\frac{1}{\lambda},
$$
i.e.
\begin{equation}\label{3fpc2}
g\circ f=\lambda g
\end{equation}
on $B(w,2r)$. Iterating this equality, we get
\begin{equation}\label{2fpc2}
g\circ f^n=\lambda^n g
\end{equation}
on $B(w,2r)$. Seeking contradiction suppose that g is injective. Then (\ref{2fpc2}) is equivalent to $f^n(z)=g^{-1}(\lambda^n g(z))$
for every $ z \in B(w,2r)$. In particular,
$$
f^n(w)=g^{-1}(\lambda^ng(w))=g^{-1}(\lambda^n)\in g^{-1}(\partial{\mathbb D}) \sbt B(w,2r).
$$
Invoking also (\ref{1fpc1}), we thus conclude that
$$
\big\{f^n(w):\, n \geq 0\big\}\sbt \mathbb C \sms B(0,2r),
$$
contrary to the fact that $ \lim_{n \to \infty}f^n(w)=0$. So, the function $g$ is not injective, and, as $g(w)=1$, it follows from Hurwitz's Theorem that
$g=\1$  on $B(w,2r)$. Therefore, the sequence
$$
\({g_n}_{|B(w, 2r)}\)^\infty_{n=0}
$$
converges to the function $\1$ uniformly on compact subsets of $B(w,2r)$ . It thus follows from (\ref{3fpc2}) that $\lambda=1$.

Consequently, $V$ is a Leau domain corresponding to $\xi$.

\sp  Item (b). Since $J(f)$ is infinite, the unit disk $\D$ is the universal covering space of $V$. Denote the corresponding cover group (the group of deck transformations) by $\Gamma$, and let 
$$
\pi:\D\lra V
$$
be the corresponding covering projection map. Changing coordinates, we may assume without loss of generality that $0\in V$, and that 
$$
\pi(0)=0.
$$
For every integer $n\ge 0$ let 
$$
\tilde f^n:\D\lra \D
$$
be a lift of $f^n|_V:V\lra V$ via $\pi$, i.e. such a continuous, in fact holomorphic, map from $\D$ to $\D$ that
\beq\label{120190625}
\pi\circ\tilde f^n=f^n\circ\pi,
\eeq
i.e. the following diagram commutes:
\[\begin{tikzcd}
\D \arrow{r}{\tilde f^n} \arrow[swap]{d}{\pi} & \D\arrow{d}{\pi} \\
{V} \arrow{r}{f^n} & {V}
\end{tikzcd}
\]
Let $R>0$ be so small that the restriction 
$$
\pi_R:=\pi|_{B(0,2R)}:B(0,2R)\lra V
$$
is injective. By our assumptions and Lemma~\ref{l120190625} the set of limit functions of the sequence 
$$
\(f^n|_V\)^\infty_{n=0}
$$  
contains the identity map $\Id_V$ on $V$. This means that there exists an unbounded increasing sequence $(n_k)_{k=1}^\infty$ of positive integers such that
\beq\label{120190730}
\Id_V=\lim_{k\to\infty}f^{n_k}|_V
\eeq
uniformly on compact subsets of $V$. So, disregarding finitely many terms, and also since $\pi(B(0,R))$ is an open subset of $V$ containing $0$, we may assume without loss of generality that 
$$
f^{n_k}(\pi(B(0,R)))\sbt \pi(B(0,2R))
$$
for all integers $k\ge 1$. So,
$$
\pi\(\tilde f^{n_k}(B(0,R))\)
=f^{n_k}(\pi(B(0,R)))
\sbt \pi(B(0,2R)).
$$
Therefore, there exists a deck transformation 
$$
g_k:\D\lra \D 
$$
of $V$ such that
$$
g_k\circ (B(0,R))\sbt B(0,2R).
$$
Replacing now $\tilde f^{n_k}$ by $g_k\circ \tilde f^{n_k}$, we will have 
formula \eqref{120190625} still satisfied and 
\beq\label{220190730}
\tilde f^{n_k}(B(0,R))\sbt B(0,2R).
\eeq
It follows from \eqref{120190730} and \eqref{120190625} that 
$$
\lim_{k\to\infty} \pi\circ\tilde f^{n_k}=\pi
$$ 
with the convergence being uniform on all compact subsets of $\D$. Along with \eqref{220190730} and injectivity of $\pi_R$, this gives that 
\beq\label{120190731}
\lim_{k\to\infty}\tilde f^{n_k}|_{B(0,R)}=\Id_{B(0,R)}.
\eeq
In particular,
\beq\label{320190730}
\lim_{k\to\infty}\tilde f^{n_k}(0)=0.
\eeq
Since the sequence $\(\tilde f^{n_k}\)_{k=1}^\infty$ is normal, it follows from \eqref{120190731} that 
\beq\label{220190731}
\lim_{k\to\infty}\tilde f^{n_k}=\Id_{\D}, 
\eeq
with the convergence being uniform on all compact subsets of $\D$. There are two cases to consider. Namely:

\,

\begin{itemize}
\item [Case 1:] $\tilde f^{n_k}(0)=0$ for some $k\ge 1$, 
  
\, 

\fr and
  
\,
  
\item [Case 2:] $\tilde f^{n_k}(0)\ne 0$ for all $k\ge 1$.
\end{itemize}
First consider Case 1.  With $k\ge 1$ as therein, we have that
$$
f^{n_k}(0)=f^{n_k}(\pi(0))= \pi\(\tilde f^{n_k}(0)\)=\pi(0)=0,
$$
meaning that $0$ is a fixed point of $f^{n_k}$. So, since $0$ is in the Fatou set of $f^{n_k}$, it is either attracting or irrationally indifferent. But it cannot be attracting since then all the limit functions of the sequence 
$\(f^n|_V\)_{n=1}^\infty$ would be constant equal $0$, $f(0),\ld,f^{n_k-2}(0)$, or $f^{n_k}(0)$. Thus $0$ is an irrationally indifferent fixed point of $f^{n_k}$. It then follows from Theorem~\ref{Siegel MU} that $V$ is a Siegel disk for $f^{n_k}$. In particular $f^{n_k}|_V$ is conjugate to some irrational rotation; hence it can have only one fixed point, namely $0$. Consequently $0$ must be a fixed point of $f$, and moreover, irrationally indifferent one. So, applying Theorem~\ref{Siegel MU} again, finishes the proof in this case. 

\sp It remains to consider Case 2. It follows from the formula \eqref{220190731} and Hurwitz's Theorem that the holomorphic maps $\tilde f^{n_k}:\D\lra \D$ are conformal homeomorphisms for all $k\ge 1$ sufficiently large; in fact, it follows from Lemma~\ref{l120190625} that $f:V\to V$ is a conformal homeomorphism, so is each of its iterates $f^n:V\to V$, $n\ge 0$, and it is well known (and easy to prove) in algebraic topology that a lift of any hoemomorphism is a homeomorphism. So, disregarding finitely many terms (if using only Hurwitz's Theorem), we may assume without loss of generality that all holomorphic maps $\tilde f^{n_k}:\D\lra \D$, $k\ge 1$, are conformal homeomorphisms. Then the inverse maps
$$
\tilde f^{-n_k}:\D\lra \D,
$$
$k\ge 1$, are well defined and holomorphic, and form a normal family in the sense of Montel. Let $H:\D\lra\D$ be any limit function of this sequence, say 
\beq\label{320190731}
H=\lim_{j\to\infty}\tilde f^{-n_{k_j}}
\eeq
for some unbounded increasing sequence $(k_j)_{j=1}^\infty$. Fix a compact set $S\sbt\D$ and $r>0$ so small that $B(H(S),2r)\sbt \D$. It follows from \eqref{320190731} that 
$$
f^{-n_{k_j}}(S)\sbt \ov{B(H(S),r)}\sbt \D
$$
for all $j\ge 1$. Having this and noting that the set $\ov{B(H(S),r)}\sbt \D$ is compact, it follows from \eqref{320190731} and \eqref{220190731} that
$$
H|_S
=\lim_{j\to\infty}\Id_{\D}\circ \tilde f^{-n_{k_j}}|_S
=\lim_{j\to\infty}\tilde f^{n_{k_j}}\Big|_{\ov{B(H(S),r)}}\circ f^{-n_{k_j}}|_S
=\lim_{j\to\infty}\Id_S. 
$$
Thus $H=\Id_{\D}$ and it follows from \eqref{320190731} that
\beq\label{320190731B}
\lim_{k\to\infty}\tilde f^{-n_k}=\Id_{\D},
\eeq
with the convergence being uniform on all compact subsets of $\D$. 

Now, for every integer $k\ge 1$ both sets $\tilde f^{n_k}\Ga$ and $\Ga\tilde f^{n_k}$ constitute the collection of all lifts of $f^{n_k}$ to $\D$. In particular, 
$$
\tilde f^{n_k}\Ga=\Ga\tilde f^{n_k}.
$$
Therefore, for any element $\gamma\in\Gamma$,
$$
\tilde f^{-n_k}\circ\g\circ\tilde f^{n_k}\in\Ga.
$$
In addition, we conclude from \eqref{220190731} and \eqref{320190731B} that 
\beq\label{420190731}
\lim_{k\to\infty}\tilde f^{-n_k}\circ\g\circ\tilde f^{n_k}=\g,
\eeq
with the convergence being uniform on all compact subsets of $\D$. Since the $\Gamma$--orbit of any point in $D(0,1)$ cannot accumulate in $B\D$ (is discrete therein), we deduce that for all sufficiently large $k\ge 1$, say for  $k\ge k(\gamma)$, we have that
\beq\label{120190801}
\tilde f^{-n_k}\circ\g\circ\tilde f^{n_k}=\g.
\eeq
meaning that $\tilde f^{n_k}$ and $\g$ commute. 

It is convenient now to replace $\D$ (as the universal covering space of $V$) by the upper half--plane $\mathbb H$. The ``new'' cover groups (acting on $\mathbb H$ and depending on the choice of a covering map from $\mathbb H$ onto $V$) of $V$ are equal to
$$
h \Gamma h^{-1},
$$
where $h$ ranges over arbitrary M\"obius transformations from $\oc$ onto $\oc$ such that $h(\D)= \mathbb H$.  However, for simplicity for the simplicity of exposition, we continue  to use the notation $\Gamma$ and $\tilde f^{n_k}$, where now $\tilde f^{n_k}$, $k\ge 1$, are conformal homeomorphisms of $\mathbb H$ (and  $\tilde f^{n_k}$ is still a lift of $f^{n_k}$). By choosing the map $h$ suitably, we may assume that the group $\Gamma$ contains one the maps
$$ 
\mathbb H\ni z \lmt z+1\in \mathbb H 
\  \  \  {\rm or} \  \  \
\mathbb H\ni z \lmt kz\in \mathbb H,
$$
with some real number $k >1$. Assume first that $\Ga$ contains the former, i.e. the map $\g:\mathbb H\to\mathbb H$ given by the formula
$$
\g(z)=z+1
$$
belongs to $\Ga$. Let $\rho\in\Ga$ be arbitrary. Recall, see \eqref{120190801}, that for large $k\ge \max\{k(\gamma),k(\rho)\}$, the map $\tilde f^{n_k}$ commutes with both $\gamma$ and $\rho$. A straightforward computational reasoning  shows  that the only M\"obius  maps which  preserve
$\mathbb H$ and which commute with  a map $ z \mapsto z+t$, $t\in\R$, are all the maps of the form $z\mapsto z+s$, $s\in\R$. Thus,
$$
\tilde f^{n_k}(z)=z+t, \ z\in \mathbb H,
$$ 
with some $t\in\R$, and hence 
$$
\rho(z) = z +s, \ z\in \mathbb H,
$$ 
with some $s\in\R$. This means that $\Gamma$ is a discrete group  of real translations and so, up to conjugacy, is generated by the map $\mathbb H\ni z \mapsto z+1\in \mathbb H$. In this case the quotient space $\mathbb H/ \Gamma$ is, up to conformal conjugacy, the punctured unit disk
$$ 
\D_0:=\{ z \in \mathbb C:0< |z| < 1\}
$$
and the corresponding quotient map is given by the formula
$$
\mathbb H\ni z \lmt \exp(2 \pi i z)\in\D_0.
$$
It follows that there is a conformal homeomorphism 
$$
Q:\D_0\lra V
$$
from $\D_0$ onto $V$. Since the point $0\in\D$ is neither pole nor an essential singularity of $Q$, it must be a removable singularity. Therefore, $Q$ extends to a meromorphic map from $\D$ onto $V\cup\{Q(0)\}$ and $Q(0)$ is an isolated point of $\partial{V}$. Hence, $Q(0)$ is an isolated point of $J(f)$, contrary to Theorem~\ref{t_perfect_julia}. Thus, $\Gamma$ must contain an element $\g$ of the form
$$
\mathbb H\ni z \lmt \g(z)=kz\in \mathbb H,
$$
with real numbers $k >1$. The argument follows now essentially the same lines as above. The only M\"obius  maps which commute with $\gamma$ and preserve $\mathbb H$ are of the form
\beq\label{220190801}
\mathbb H\ni z \lmt tz,
\eeq
where $t>0$, and the maps of the form
\beq\label{320190801}
\mathbb H\ni z \lmt u/z,
\eeq
where $u<0$. Now, maps of the latter type of order two, so if $\Ga$ contains such element and $k\ge 1$ is large enough, then by \eqref{420190731}, 
$$
\tilde f^{2n_k}=\Id_{\mathbb H}.
$$
Thus,
$$
f^{2n_k}|_V=\Id_V.
$$
Since both functions $f^{2n_k}:\C\to\oc$ and $\Id_{\C}:\C\to\oc$ are meromorphic, this implies that
$$
f^{2n_k}=\Id_{\C},
$$
contrary to the fact that the degree of $f$ is larger than $1$. It follows that each map $\tilde f^{n_k}:\mathbb H\to\mathbb H$ with $k\ge 1$ large enough is of the form from \eqref{220190801}.
Now, for every $k\ge 1$ large enough, any element $\rho\in\Gamma$ commutes with $\tilde f^{n_k}$, so $\rho$ is one of the maps of \eqref{220190801} or \eqref{320190801}. If $\rho\in\Ga\sms\{\Id_{\mathbb H}\}$, then $\rho$ has no fixed points in $\mathbb H$. So, $\rho$ must be of the form from \eqref{220190801}, in addition with $t\ne 1$. Hence, $\Gamma$ is a discrete subgroup of the maps of \eqref{220190801}. Thus, it is a cyclic group generated by some map 
$$
\mathbb H\ni z \longmapsto kz\in \mathbb H,
$$
$k\in(0,+\infty)\sms\{1\}$. Write
$$
k=\frac{2\pi^2}{\ka}
$$
with some unique $\ka>0$. It therefore follows from Proposition~\ref{p1j220} that the quotient space ${\mathbb H}/\Gamma$ is conformally equivalent to the annulus $A(0;e^{-\ka},1)$. Since it is also conformally  equivalent to $V$, item (b) is now proved by a direct application of Lemma~\ref{l220190801}.

\sp Proving item (c), note first that if the family of \eqref{520190802} has a limit function which is constant, then the reasoning from the beginning of proving item (a), shows that the collection of all limit functions of this family is a singleton. But this is a contradiction since $\infty$ is such a function. Thus,
$$
\lim_{n\to\infty}f^{n}(z)=\infty
$$
for all $z\in V$ uniformly on all compact subsets of $V$. The proof of item (c) is complete, and simultaneously, the proof of the whole theorem in the case of forward invariant Fatou components $V$. 

\sp If $V$ is periodic with some period $p\ge 1$, then we pass to the $p$th iterate $f^p$ and apply the already proven case of invariant domains (i.e. of period $1$). One should only notice that the fact that $f^p$ need not be longer meromorphic (pre-poles of $f$ of orders $\le p-1$ become essential singularities for $f^p$) causes no problems. The proof of Theorem~\ref {Limit_functions} is complete.
\epf

\sp\section[Fatou Components of (General) Meromorphic Functions II] {Fatou Components of (General) Meromorphic Functions II; Properties of Periodic Components}

 The importance of periodic components of Fatou sets is multifolded. The present section is an evidence of this. But there are much more reasons.
Later, in Section~\ref{FC B+S}, we will see that any meromorphic function in Speiser Class $\cS$ (to be defined in Section~\ref{FC B+S}) has no wandering components, in other words, that every component of the Fatou set is preperiodic. 

\sp The periodic domains of the Fatou set $F(f)$ of $f$ are closely related to the  set of singular values of $f$. Indeed, we say that $y \in\C$ is a regular point \index{(N)}{regular point} of $f^{-1}$ if for every  $r>0$ small enough and every connected  component $C$ of $f^{-1}(B(y,r))$ the restriction 
$$
f|_C: C\lra B(y,r)
$$
is a homeomorphism from $C$ onto $B(y,r)$. Otherwise we say that $y$ is a singular point \index{(N)}{singular point} of $f^{-1}$ and we denote by $\Sing(f^{-1})$ \index{(S)}{$\Sing(f^{-1})$} the set of all such singular points. It is well known that $\Sing(f^{-1})$ consists of critical ($f(\Crit(f))$) and, finite asymptotic, values of $f$. We also set
$$
{\rm PS}(f):=\bigcup_{n=0}^\infty f^n(\Sing (f^{-1}))
$$
\index{(N)}{PS(f)}
and call it the postsingular set of $f$ \index{(N)}{postsingular set}. The following two theorem show the significance of the set $\Sing(f^{-1})$ when studying attracting components of the Fatou set. Its importance (see all subsequent chapters) goes far beyond.

\sp\bthm\label{t1ms123}
Let $f:\C\lra\oc$ be meromorphic function. If $\xi\in\C$ is an attracting periodic point of $f$, then $A_p^*(\xi)$, the basin of immediate attraction to the forward orbit of $\xi$, contains at least one singular point of $f^{-1}$ whose forward orbit under $f$ either does not contain $\xi$ or it coincides with the forward orbit of $\xi$. In particular, in the former case this orbit is infinite while in the later case the orbit of $\xi$ is supper attracting. 
\ethm

\bpf
Assume first that $\xi$ is an attracting fixed point of $f$. If $f'(\xi)=0$, we are done. So, suppose that
$$
f'(\xi)\ne0.
$$
Then there exists $U_0$, an open, connected, simply connected, neighborhood of $\xi$ such that
$$
f(U_0)\sbt U_0
$$
and the map 
$$
f|_{U_0}:U_0\to U_0
$$
is 1--to--1. Observe that then
$$
U_0\cap\bu_{n=0}^\infty f^{-n}(\xi)=\{\xi\}.
$$
Now, by the way of contradiction, suppose that all singular points of   $f^{-1}$ lying in $A^*(\xi)$ are eventually mapped onto $\xi$ under some iterates of $f$, i.e. that they belong to the set $\bu_{n=1}^\infty f^{-n}(\xi)$. We shall prove by induction that there exists an ascending sequence $(U_n)_{n=0}^\infty$ of connected, simply connected, open subsets of $A^*(\xi)$ such that:
\begin{itemize}
  \item [($a_n$)] $\xi\in U_n$ for all $n\geq0$.
  
\,\item [($b_n$)] For every $n\geq1$ there exists $f_n^{-1}:U_{n-1}\to U_n$, a unique surjective holomorphic inverse branch of $f$ such that $f_n^{-1}(\xi)=\xi$.

\,\item [($c_n$)]
$$
U_n\cap\bu_{n=0}^\infty f^{-n}(\xi)=\{\xi\}
$$
for all integers $n\geq0$.
\end{itemize}
Indeed, the set $U_0$ has been already defined satisfying ($a_0$) and ($c_0$), and, vaciously ($b_0$). For the inductive step, suppose that $n\geq0$ and that connected, simply connected subsets
$$
U_0\subset U_1\subset U_2\subset\cdots \subset U_n\subset A^*(\xi),
$$ 
satisfying conditions (a), (b), and (c), with appropriate sunscripts, have been defined. Since the open, connected, and simply connected set $U_n$ contains no singular values of $f|_{A^*(\xi)}$ other than $\xi$, the branch $f_n^{-1}$ can be analytically continued throughout $U_n$ and the Monodromy Theorem yields the existence of a holomorphic inverse branch
$$
f^{-1}_{n+1}:U_n\to A^*(\xi)
$$ 
uniquely determined by the requirement that $f_{n+1}^{-1}(\xi)=\xi$ or, equivalently, that $f_{n+1}^{-1}|_{U_{n-1}}=f_n^{-1}$. Set
\[
U_{n+1}:=f^{-1}_{n+1}(U_n)\subseteq A^*(\xi)
\]
Immediately from this definition we see that $U_{n+1}$ is an open connected, simply connected subset of $A^*(\xi)$, and 
$$
\xi=f^{-1}_{n+1}(\xi)\in f_{n+1}^{-1}(U_n)=U_{n+1}.
$$
In conclusion, item $(b_{n+1})$ holds. Now, we shall show that 
$$
U_n\subseteq U_{n+1}.
$$
Indeed, 
$$
U_n
=f_n^{-1}(U_{n-1})
=f_{n+1}^{-1}(U_{n-1})
\sbt f_{n+1}^{-1}(U_n)
=U_{n+1}.
$$
So we are only left to show that $(c_{n+1})$ holds. Indeed, suppose that 
$$
z\in U_{n+1}\cap\bu_{n=0}^\infty f^{-n}(\xi).
$$
Then, as $f(\xi)=\xi$ and $f(U_{n+1})=U_n$, we have that
$$
f(z)\in U_n\cap\bu_{n=0}^\infty f^{-n}(\xi).
$$
Hence, by ($c_n)$, $f(z)=\xi$. Thus,
$$
z=f_{n+1}^{-1}(\xi)=\xi,
$$
and item $(c_{n+1})$ is established The inductive construction of the sequence $(U_n)_{n=0}^\infty$ is thus complete.

It follows from all the above properties of the sequence $(U_n, f_{n+1})_{n=0}^\infty$ that
for all $n\geq1$, the composition
\[
f_\xi^{-n}:=f_n^{-1}\circ f^{-1}_{n-1}\circ \cdots \circ f_2^{-1}\circ f_1^{-1}:U_0\to U_n
\subset A^*(\xi)
\]
is well defined and holomorphic. Since $A^*(\xi)\subset F(f)$ and since the Julia
set $J(f)$ contains at least three points (because it is perfect, by Theorem
\ref{t_perfect_julia}), it follows from Montel's Theorem that the family
of maps $(f^{-n}_\xi:U_0\to A^*(\xi))_{n=0}^\infty$ is normal. This, however, produces a contradiction, since 
$$
\lim_{n\to\infty}(f^{-n}_\xi)'(\xi)=\lim_{n\to\infty}((f^n)'(\xi))^{-1}=\infty.
$$
We are thus done in the case when $\xi$ is a fixed point of $f$.

In general, if $\xi$ is a periodic point of $f$, say of minimal period $p\ge 1$, then, as we have just proved, the map $f^p:A^*(\xi)\lra A^*(\xi)$ contains a singular point of $f^{-p}$. Moreover, as 
$$
A^*_p(\xi)=\bigcup_{j=0}^{p-1}A^*(f^j(\xi))=\bigcup_{j=0}^{p-1}f^j(A^*(\xi))
$$ 
and $f^p=f\circ f\circ\ld\circ f$ ($n$-folded composition), we thus conclude that $A^*_p(\xi)$ contains a singular value of $f$ and the proof is complete.
\epf

\sp Later on in Section~\ref{FCGMFII}, after appropriate preparations, we shall prove as Theorem~\ref{t2ms135} an analogous result for basins of immediate attraction to rationally indifferent periodic points. Now, we shall prove a somehow analogous theorem for Siegel disks and Herman rings. We will need the following property.

\sp\bproper[Standard Property]\label{pr1j144} 
Let $Y$ be a complete Riemann surface with constant curvature $0$, i.e. $Y$ is either the complex plane $\C$, a complex torus $\mT_\La=\C/\La$, where $\La$ is a lattice on $\C$ or an infinite cylinder $\C/2\pi i\Z$.
Let $X$ be a non--empty open subset of $Y$. An analytic map 
$$
f:X\lra Y
$$ 
is said to have the Standard Property if and only if the set $(Y \sms X ) \cup  {\rm Per }(f)$ contains at least three points, where, as usually,
${\rm Per }(f)$ denotes the set of all periodic points of $f$.
\eproper

\sp Let $X$, $Y$ and $f: X \lra Y$ be as in Property~\ref{pr1j144}. In particular $f$ has the Standard Property. We say that a point $z \in X$  belongs \ to the Fatou set $F(f)$ of $f$ if there is  an open neighborhood $U$ of $z$ such that all the iterates 
$$
f|_U^n:U \lra  Y, \  \  \  n \geq 1, 
$$
are all well--defined and
contain a subsequence forming a normal family. The Julia set $J(f)$  is defined  as $Y \sms F(f)$. Clearly $J(f)$ is a closed subset of $Y$ and
$$
f^{-1}(J(f)) \sbt J(f), \  \  \  J(f)\cap X \sbt  f^{-1}(J(f)), \  \  \  J(f)\cap X \sbt J(f) 
$$
and
$$
J(f)\cap f(X) \sbt  f(J(f)).
$$

Given $w \in Y, \, R >0$, and  $W$, an open neighborhood of  $w$ compactly contained in $B(w,R)$, for every integer $n \geq 0$ we denote  by 
$$
\Comp_n^*(w,R)
$$
\index{(S)}{$\Comp_n^*$} the collection of all connected components $\Ga$ of $f^{-n}(B(w,R))$ such that the map 
$$
f^n|_{\^\Ga}: \^\Ga\lra B(w,2R)
$$ 
is a conformal homeomorphism from $\tilde{\Ga}$ onto $B(w,2R)$, where $\tilde{\Ga}$ is the (only) connected component of $f^{-n}(B(w,2R))$ containing $\Ga$. We denote by 
$$
f^{-n}_\Ga:B(w,2R)\lra\tilde{\Ga}
$$ 
its inverse $(f^n|_{\tilde{\Ga}})^{-1}$. We then put $n(\Ga)$ to be $n$.  Note that if $w \notin \ov{\rm PS (f)}$ and $ R\le \frac12\dist(w,\ov{\rm PS (f)})$, then  
$$
\Comp^*_n(w,R)=\Comp_n(w,R),
$$
\index{(S)}{$\Comp_n$} the collection of  all connected component of $f^{-n}(B(w,R))$.

\sp We start with the following general result which will be used now and frequently in the sequel.

\blem\label{l1j293}
Let $Y$ be a complete Riemann surface with constant curvature $0$, i.e. $Y$ is either the complex plane $\C$, a complex torus $\mT_\La=\C/\La$, where $\La$ is a lattice on $\C$ or an infinite cylinder $\C/2\pi i\Z$.
Let $X$ be a non--empty open subset of $Y$ and let 
$$
f:X \lra Y
$$ 
be an analytic map with the Standard Property. Let $Q\sbt Y$ be a set witnessing this property, i.e. $Q\sbt (Y\sms X)\cup\Per(f)$ and $Q$ has at least three elements. Let $w \in J(f)\sms O_+(Q)$. 

Then for every $u>0$ we have that \index{(S)}{$\Comp_n^*$}
\beq\label{1j293}
  \liminf_{n \to   \infty} \big\{ \inf\{ |(f^n)'(z)| :\,  z \in  V \}: V \in  \Comp_n^*(w,u)\big\} = +\infty.
\eeq
\elem

\bpf For  every $n \geq 0$ and every $V \in \Comp_n^*(w,u)$ recall that 
$$
f^{-n}_V:B(w,2u) \lra \hat{V}
$$
is the unique holomorphic inverse branch of $f^{-n}$ defined on $B(w,2u)$ and mapping it onto $V$. Let $\hat V$ be the only connected component of $f^{-n}_V(B(w,2u))$ containing $V$. By our choice of $Q$ and $w$, the family
$$
{ \mathfrak F}:= \big\{f^{-n}_V: \,  n\geq 1, \,  V \in \Comp_n^*(w,u)\big\}
$$ 
omits the set $Q$ which consists of at  least three points. Thus, by Montel's Theorem, the family $\mathfrak F$ is normal. 

Seeking contradiction suppose that (\ref{1j293}) fails. This means that there exist $\gamma >0$ and an infinite sequence $( n_k)_{1}^{\infty}$ such that  for every  $k \geq 1$
there exist $V_k \in  \Comp^*_{n_k}(w,u)$ and a point $\xi_k \in V_k$ such that   $ |(f^{n_k})'( \xi_k)|\leq  \gamma$. It then follows from Koebe's Distortion Theorem  that
\beq\label{2j293}
|(f^{n_k})'(z)|  \leq M\gamma
\eeq
for all  $z \in \hat V_k$ with some constant $M\ge 1$. Passing to a subsequence we may assume  without loss of generality that the sequence $(f^{-n_k}_{V_k})_{k=1}^{\infty}$ converges to an analytic function 
$$
g:B(w,2u) \lra Y
$$
uniformly on compact subsets of $B(w,2u)$. It follows from (\ref{2j293})  that  
$$
|g'(z)|\geq  (M\gamma)^{-1}
$$ 
for all $z\in B(w,u)$. In particular, $g: B(w,2u) \lra Y$ is not a constant function, and so, $g(B(w,u/2))$ is an open neighborhood of $g(w)$. Since $g(\ov{B(w,u/2)})$ is a compact subset of the open set $g(B(w,u))$, we thus  conclude that
$$
f^{n_k}(g(B(w,u/2)))\sbt B(w,u)
$$
for all $k\geq 1$ large enough. Hence the family $\(f^{n_k}_{| g(B(w,u/2))}\)^{\infty}_{k=1}$ is normal, and therefore the point $g(w)$ belongs to $F(f)$, the Fatou set of $f$. But, on the other hand, 
$$ 
g(w)=\lim_{n_\to \infty} f^{-n_k}_{V_k}(w)\in \ov{J(f)}=J(f).
$$
This contradiction finishes the proof. \qed

\bthm\label{t1sh11}
Let $f:\C\lra\oc$ be a meromorphic function. If $\{U_0, U_1, \ldots, U_{p-1}\}$ is a cycle of Siegel disks or Herman rings, then 
$$
\partial{U}_j \subset \ov{O^+(\Sing(f^{-1}))}
$$ 
for all $j\in \{1, \ldots, p-1\}$.
\ethm

\bpf
We aim to apply Lemma~\ref{l1j293}. We set,
$$ 
Y:= \mathbb C, \quad X:=\mathbb C \sms f^{-1}(\infty)
$$
and keep $f$  the same, i.e. to be really formally correct, when applying Lemma~\ref{l1j293}, we consider $f|_X$. We take as $Q$ an arbitrary finite subset of $\Per(f)(\sbt \mathbb C)$ having at least three elements. Seeking contradiction suppose that
\beq\label{1sh11}
\partial{U}_j \nsubseteq \ov{O^{+}(\Sing (f^{-1}))}
\eeq
for some $ j \in \{0,1, \ldots, p-1\}$. Denote 
$$
D :=\partial{U}_j.
$$
Let 
$$
\Delta:=B(0,1)
$$  
if $D$ is a Siegel disk and 
$$
\Delta:=A(0;1, r)
$$ 
with $ r >1$ coming from Theorem~\ref{Fatou Periodic Components} if $D$ is a Herman ring. Let 
$$
H:\Delta \lra D
$$ 
be the analytic homeomorphism  resulting from Theorem~\ref{Fatou Periodic Components} (3) and (4) respectively in the Siegel or Herman case. Because of (\ref{1sh11}) there exists a point 
\beq\label{1sh6}
\xi \in \partial{D}\sms  \ov{O^{+}(\Sing (f^{-1}))}.
\eeq
Using the fact that either $D$ is simply connected or doubly connected and
 $r \in (1, +\infty)$, we deduce that $ \xi$ is not an isolated point of $\partial{D}$, in fact $\partial{D}$ has no isolated points. Therefore, since $O^{+}(Q) \cup f^{-1}(\infty)$  is a countable set, we may assume in addition that 
$$
\xi \notin O^{+}(Q)\cup f^{-1}(\infty),
$$ 
so that Lemma~\ref{l1j293} applies. Fix $s>0$ so small that 
$$
B_e(\xi,3s)\cap {\rm PS}(f)=\es.
$$
Pick a point
\beq\label{2sh6}
w \in B(\xi, s)\cap D.
\eeq
Then
\beq\label{3sh6}
B(w,2s) \cap {\rm PS}(f)=\es.
\eeq
Let
\beq\label{4sh6}
F:=H\(\{ z \in \Delta:\,  |z|=|H^{-1}(w)|\}\).
\eeq
Note that $F$ is a compact set (homeomorphic to a circle) and $F \sbt D$. Since $w \in D$, for every $n \geq 1$ the intersection
$D \cap f^{-n}(w)$ is a singleton. Denoting its only element by $w_n$, we will have that
\beq\label{1sh7}
w_n\in F.
\eeq
By (\ref{3sh6}), for every $n \geq 1$ there exists a unique holomorphic branch
\beq\label{2sh7}
f^{-n}_n : B(w,2s)\lra \mathbb C
\eeq
of $f^{-n}$ sending $w$ to $w_n$. By (\ref{1sh6}), $ \xi \in J(f)$, whence $f^{-n}_n(\xi)\in J(f)$. Thus, $f^{-n}_n (\xi)\notin D$ and therefore, using also
(\ref{2sh7}) and (\ref{1sh7}), we  conclude that
$$
\diam_e\(f^{-n}_n (B(w,2s)))\geq \dist_e(w_n,f^{-n}_n(\xi)\)
\geq \dist_e(F, \mathbb C \sms D)>0.
$$
Therefore,
$$ 
\liminf_{n \to \infty}\diam_e\(f^{-n}_n (B(w,r))\)
\geq  \dist_e(F, \mathbb C \sms D)>0.
$$
It thus follows from Koebe's Distortion Theorem i.e. Theorem~\ref{Euclid-I}, that
$$
\liminf_{n \to \infty}|(f^n)'(w)| < +\infty,
$$
contrary to Lemma~\ref{l1j293}. The proof of Theorem~\ref{t1sh11} is thus finished.
\epf

\sp

\section[Local and Asymptotic Behavior of (General) Meromorphic Functions] {Local and Asymptotic Behavior of (General) Meromorphic Functions Locally Defined Around Rationally Indifferent (Parabolic) Fixed Points; Part I}\label{local-parabolic}

 In this section we want to  bring up some basic results about local behavior
of meromorphic functions about their parabolic (rationally indifferent) fixed points, or
somewhat more generally  about parabolic periodic points. As a
matter of fact our analysis is so local that all what we will be
assuming throughout this section is that our given analytic map
$\varphi$ is defined on some  neighborhood of its parabolic fixed point.

\sp In conformity to the previous sections we call a holomorphic map $\varphi$, defined around a point $\om\in\mathbb C$, a locally holomorphic map parabolic at $\om$ (or rationally indifferent at $\om$), if $\om$ is a periodic point of $\varphi$, meaning that $\varphi^p(\om)=\om$ for some $p\ge
1$, $(\varphi^p)'(\om)$ is a root of unity and no iterate of $\varphi$ is equal to the identity map. We then also call $\om$
parabolic or rationally indifferent. We call the locally holomorphic parabolic map $\varphi$, or the point $\om$, simple parabolic, \index{(N)}{simple parabolic periodic point} if 

\begin{itemize}
\item[(a)] $\varphi(\om)=\om$,

\,
 
\item[(b)] $\varphi'(\om)=1$, and

\,
 
\item[(c)] $\varphi$ is not the identity map. 
\end{itemize}

\fr Note that some sufficiently high iterate of any locally holomorphic parabolic map is simple. Therefore, in order to analyze the behavior of locally holomorphic parabolic maps, it essentially suffices to do this for simple parabolic maps
only. Thus, throughout this section the map $\varphi$ (and its fixed
point $\om$) is always assumed to be locally holomorphic simple parabolic. Then on a sufficiently small 
neighborhood of $\om$, the map $\phi$ has the following Taylor
series expansion:
\beq\label{2j101}
\varphi(z)=z+a(z-\om)^{p+1}+b(z-\om)^{p+2}+\cdot\cdot\cdot
\eeq
with some integer $p=p(\om)\ge 1$ and $a\in \mathbb C\sms \{0\}$. Being in
the circle of ideas related to Fatou's Flower Theorem (see [Al] for
extended historical information), we now want to analyze
qualitatively and especially quantitatively the behavior of $\varphi$
in a sufficiently small neighborhood of the parabolic point $\om$.
Let us recall that the rays coming out from $\om$ and forming the
set
$$
\{z\in\C:a(z-\om)^p<0\}
$$
are called attracting directions \index{(N)}{attracting directions} and the rays forming the set $$
\{z\in\C:a(z-\om)^p>0\}
$$
are called repelling directions. Fix an attracting direction, \index{(N)}{repelling  directions} say
$$
A:=\om+\sqrt[p]{-a^{-1}(0,\infty)},
$$
where $\sqrt[p]{\cdot}$ is a holomorphic branch of the $p$th radical defined on $\mathbb C\sms a^{-1}(0,\infty)$. In order to simplify our analysis let us change the system of coordinates with the help of the affine map
$$
\rho_A(z)=\sqrt[p]{-a^{-1}}z+\om.
$$
I.e. we define
$$
\varphi_{A,0}:=\rho_A^{-1}\circ\varphi\circ\rho_A
$$
We then get
\beq\label{120120728}
\varphi_{A,0}(z)
=\rho_A^{-1}\circ\varphi\circ\rho_A(z) 
=z-z^{p+1}+b\sqrt[p]{-a^{-1}}z^{p+2}+\cdot\cdot\cdot
=z-z^{p+1}+\sum_{n=2}^\infty a_nz^{p+n},
\eeq
and
$$
\varphi_{A,0}(0)=0, \  \  \(\varphi_{A,0}\)'(0)=1,
$$
i.e. $0$ is a simple parabolic fixed point of $\varphi_{A,0}$. In addition, $\rho^{-1}(A)=(0,\infty)$ is an attracting direction for
$\varphi_{A,0}$. 

We first want to analyze the behavior of $\varphi_{A,0}$ on sufficiently small neighborhoods of $0$.
In order to do this, similarly as in the previous section, we
conjugate $\varphi_{A,0}$ on $\mathbb C\sms (-\infty,0]$ to a map defined
"near" infinity. More precisely, we consider $\sqrt[p]{\cdot}$, the
holomorphic branch of the $p$th radical defined on $\mathbb C\sms
(-\infty,0]$ and leaving the point $1$ fixed. We will also frequently denote this branch by $z^{1/p}$. Then we define the map
$$
H(z)={1\over \sqrt[p]{z}}=z^{-\frac1p}.
$$
This map has a meromorphic inverse $H^{-1}:\oc\lra\oc$ defined by the formula
$$
H^{-1}(z)=z^{-p}.
$$
Consider the conjugate map
\beq\label{1j101}
\^\varphi_A
=H^{-1}\circ\varphi_{A,0}\circ H
=(\rho_A\circ H)^{-1}\circ\varphi\circ (\rho_A\circ H),
\eeq
defined on $V_\varphi\sms (-\infty,0]$, where $V_\varphi\sbt\C$ is a sufficiently small neighborhood of $0$ on which $\vp_0$ is defined. We shall prove the following technical but very useful result.

\sp\blem\label{l1j101} 
If $\varphi$ is a locally holomorphic simple parabolic map and $A$ is an attracting direction of $\varphi$, then shrinking the neighborhood $V_\varphi$ if necessary, there exists a holomorhic function $B:V_\varphi\lra\C$ such that for all $z \in \(\C\sms (-\infty,0]\)\cap H^{-1}(V_\varphi)$ we have 
$$
\^\varphi_A(z)=z+1+B(H(z)).
$$
and 
$$
\^\varphi_A'(z)=1+(B\circ H)'(z),
$$
with
$$
|B(H(z))|\le M|H(z)|=M|z|^{-\frac1p}
$$
and
$$
|(B\circ H)'(z)|
=|B'(H(z))|\cdot|H'(z))|
\le Mp^{-1}|z|^{-\frac{p+1}{p}},
$$ 
with some appropriate constant $M\in(0,+\infty)$. 
\elem

\bpf We have $H^{-1}(z)=z^{-p}$. For all $z \in \C\sms (-\infty,0]\cap H^{-1}(V_\varphi)$ we then get
\beq\label{1h41}
\begin{aligned}
\tilde{\varphi}_A(z) 
&= H^{-1}\( \varphi_{A,0} (H(z)) \) 
  = H^{-1} \( H(z)-H(z)^{p+1} + \sum_{n=2}^{\infty}
  a_n H(z)^{p+n} \)  \\ 
& = H^{-1} \lt(z^{-\frac1p}-az^{-\frac {p+1}{p}} +
\sum_{n=2}^{\infty} a_n z^{-\frac {p+n}{p}}\rt) \\ 
&= H^{-1} \lt(z^{-\frac1p}\lt(1-z^{-1} +
\sum_{n=2}^{\infty} a_n z^{-\frac {p+n-1}{p}} \rt)  \rt) \\ 
&= \frac {z}{\(  1-z^{-1} +  \sum_{n=2}^{\infty} a_n
  z^{-\frac {p+n-1}{p}} \)^{p}} 
\end{aligned}
\eeq
Set $w = H(z) =z^{-\frac {1}{p}}$ and put
\beq\label{1ath20111010}
G(w)=1-w^{p} +  \sum_{n=2}^{\infty} a_nw^{p+n-1}, \ \
w\in V_\varphi.
\eeq
Now 
$$
G(0)=1, \  \ \frac {\partial^{k}G(w)}{\partial w^{k}}\Big|_{0} 
=0 \  \text{ for}\ \text{all} \ k=1, 2, . . ., p-1,  \
\text{and}\  \frac {\partial^{p} G(w)}{\partial w^{p}}
\Big|_{0}=-p!.
$$
Therefore,
$$
G^{-p}(0)
=1, \  \ \frac {\partial^{k}\(G^{-p})(w)}{\partial w^{k}}\Big|_{0} 
=0 \  \text{ for}\ \text{all} \ k=1, 2, . . ., p-1,  \
\text{and}\  \frac {\partial^{p}(G^{-p})(w)}{\partial w^{p}}
\Big|_{0}
=p!
$$
and
\beq\label{2ath20111010}
G^{-p}(w) = 1+w^{p} +  \sum_{n=1}^{\infty} b_n w^{p+n}
\eeq
with some appropriate coefficients $b_{n}(\psi)$, $n\ge 1$, where the series
$$
\sum_{n=1}^{\infty} b_nw^{p+n}
$$
converges uniformly absolutely on some sufficiently small neighborhood of $0$; shrinking $V_\varphi$ if needed, we may indentify this neighborhood with $V_\varphi$. Then the series
$$
B(w):=\sum_{n=1}^{\infty} b_nw^n
$$
also converges uniformly absolutely on $V$, and represents a holomorphic function. Shrinking $V$ even more if needed, we will have a constant $M\in (0,+\infty)$ such that
\beq\label{f1_2015_09_08}
|B(w)|\le M|w| \  \  \text{ and } \  \  |B'(w)|\le M 
\eeq
for all $w\in V_\varphi$. Going back to the variable $\C\sms (-\infty,0]\cap H^{-1}(V_\varphi)\ni z=w^{-p}$ we get from
(\ref{1h41}), \eqref{1ath20111010}, and (\ref{2ath20111010}) that,
\beq\label{2ath141}
\begin{aligned}
\tilde\vp_A(z)
&=z G^{-p}(H(z))=z\(1+H(z)^p+\sum_{n=1}^{\infty} b_n(z) H(z)^{p+n}\)\\ 
&= z\lt(1+\frac{1}{z} + \frac{1}{z}\sum_{n=1}^{\infty} b_nH(z)^n\rt)\\
&=z+1+B(H(z)).
\end{aligned}
\eeq
It also follows from \eqref{f1_2015_09_08} that
\beq\label{2ath141B}
|B(H(z))|\le M|H(z)|=M|z|^{-\frac1p}
\eeq
and 
\beq\label{2ath141C}
|(B\circ H)'(z)|
=|B'(H(z))|\cdot|H'(z))|
\le Mp^{-1}|z|^{-\frac{p+1}{p}}.
\eeq
The proof is complete.
\epf

\sp Given now a point $x \in (0, \infty)$ and an angle $\a \in (0, \pi)$, let
$$ 
S(x,\a):=\{z\in \mathbb C: \, - \a <\arg(z-x)< \a
\}.
$$\index{(S)}{$S(x,\a)$} 
By Lemma~\ref{l1j101} for all $\kappa \in
(0,1)$ and all $\a \in (0,\pi)$ there exists $x(\a, \kappa) \in (0,
+\infty)$ such that for every $x \geq x(\a, \kappa)$, we have that
\beq\label{2j103}
    \^\varphi_A(S(x,\a))\sbt S(x+1-\kappa, \a) \sbt S(x, \a)
\eeq
and
\beq\label{1j105}
\re(\^\varphi_A(z)) \geq \re(z)+1-\kappa
\eeq
for all $z \in S(x, \a)$. Then by an immediate induction,
\beq\label{f2_2015_09_08}
\re(\^\varphi_A^n(z)) \geq \re(z)+(1-\kappa)n
\eeq
for all $z \in S(x,\a)$ and all integers $n\ge 0$. Summarizing, we have that the following two statements.

\sp
\bprop\label{p1j104} 
Let $\varphi$ be a locally holomorphic simple parabolic map and let $A$ be an attracting direction of $\varphi$.
Fix $\a\in (0, \pi)$, $\kappa \in (0,1)$, and $x \geq x(\a, \kappa)$.
Then all the iterates $\^\varphi_A^n: S(x,\a)\lra \mathbb C$, $n\ge 0$, are well defined and
$$
\^\varphi_A^n(S(x,\a)) \sbt S(x+(1-\kappa)n,\a)\sbt S(x,\a).
$$ 
In particular, $\^\vp_A^n(z)\to \infty$ uniformly on $S(x,\a)$.
\eprop

\sp
\blem\label{l1j105} 
If $\varphi$ is a locally holomorphic simple parabolic map and $A$ is an attracting direction of $\varphi$, then for all $\a \in (0, \pi)$ and $\kappa >0$, we have
$$
|\^\varphi_A^n(z)|
\geq \re (\^\varphi_A^n(z)) 
\geq \re z + (1-\kappa)n 
$$
for all $n\geq 0$ and all $z \in S(x(\a,\kappa),\a)$. 
\elem

\sp\fr We shall prove the following.

\sp
\blem\label{l2j105} 
If $\varphi$ is a locally holomorphic simple parabolic map and $A$ is an attracting direction of $\varphi$, then for all $\a\in (0, \pi)$, $\kappa \in (0,1)$, $n \geq 1$, and $  z\in  S(x(\a,\kappa),\a)$, we have
$$
\^\varphi_A^n(z)=z+n+O_{\re(z)}\(\max\{n^{1-\frac1p},\log(n+1)\}\),
$$
where $O_t$ is a big O symbol depending decreasingly on $t\in\R$ and converging to $0$ when $t\lra +\infty$.
\elem

\bpf  Fix $z \in S(x(\a,\kappa),\a)$. By
Lemmas~\ref{l1j101} and \ref{l1j105}, we  have for every $n \geq 1$ that
$$
\^\varphi_A^{n+1}(z)=\^\varphi_A^n(z)+1+O_{\re(z)}(n^{-\frac{1}{p}}).
$$
Therefore, by induction,
$$
\^\varphi_A^n(z)
=z+n+O_{\re(z)}\left(\sum_{k=1}^{n-1}k^{-\frac{1}{p}}\right)
=z+n+O_{\re(z)}(\max\{n^{1-\frac1p}, \log(n+1)\}).
$$ 
The proof is complete. \endpf

\sp\fr As an immediate consequence of this lemma, we get the following.

\sp
\blem\label{l1j107} If $\varphi$ is a locally holomorphic simple parabolic map and $A$ is an attracting direction of $\varphi$, then for all $\a\in (0, \pi)$,  $ \kappa \in (0,1)$
and $t\in\R$, there exists a constant $C=C(\a, \kappa, t)$ such that
for all $n \geq 1$ and all $z \in S(x(\a,\kappa),\a)$ with $\re(z)\ge t$, we
have
$$
C^{-1}n\le |\^\varphi_A^n(z)|\le Cn.
$$
\elem

\sp\fr For every $R>0$ let
$$
S(x,\a,R):=S(x,\a)\cap B(0,R).
$$
We shall prove the following. 

\sp
\blem\label{l2j107} 
Let $\varphi$ be a locally holomorphic simple parabolic map and let $A$ be an attracting direction of $\varphi$.
Fix $\a\in(0,\pi)$, $\kappa \in (0,1)$ and $R>0$. Then
$$
\begin{aligned}
0<\inf\big\{|(\^\varphi_A^{n})'(z)|:z\in \,&S(x(\a,\kappa),\a,R),\, n\geq 1\big\}\le \\
&\le \sup\{|(\^\varphi_A^{n})'(z)|:z\in S(x(\a,\kappa),\a,R), \, n
\geq 1\}< +\infty.
\end{aligned}
$$ 
Furthermore, for every $\g>1$, 
$$
\lim_{x\to+\infty}\sup\big\{|(\^\varphi_A^{n})'(z)-1|:\, z\in S(x,\a,\g x), \, n\geq 1\big\}=0
$$ 
and
$$
\lim_{x\to+\infty}\sup\big\{|(\^\varphi_A^{n})'(z)-1|:\, z\in S(x,\pi/2),\, n\geq 1\big\}=0.
$$ 
\elem

\bpf For every $z\in S(x(\a,\kappa),\a)$ let
$$
g(z)=\^\varphi_A'(z)-1.
$$
By the Chain Rule, we have for every $n \geq 1$ that
$$
(\^\varphi_A^n)'(z)
=\prod_{j=0}^{n-1}\^\varphi_A'(\^\varphi_A^j(z)))
=\prod_{j=0}^{n-1}(1+g(\^\varphi_A^j(z))).
$$
But by Lemma~\ref{l1j101} and Lemma~\ref{l1j105}, we have
$$
|g(\^\varphi_A^j(z))|\lek|\^\varphi_A^j(z)|^{-\frac{p+1}{p}}.
$$
For every $x>0$ let $r=r(x,\a)$ be the radius of the maximal ball centered at $0$ which is disjoint from $S(x,\a)$. Of course 
$$
r(x,\a)\le ux
$$
with some constant $u\in(0,+\infty)$ depending only on $\a$. Let $k_R\ge 0$ be the least integer such that 
$$
-R+(1-\ka)(k_R+1)\ge R.
$$
Then $-R+(1-\ka)k_R< R$, and therefore
$$
k_R<\frac{2R}{1-\ka}.
$$
Fix $x\ge x(\a,\ka)$ arbitrary. Because of Proposition~\ref{p1j104} we have for all $z\in S(x,\a)$ that 
$$
|\^\varphi_A^j(z)|\ge r(x,\a)\ge ux
$$
for all $0\le j\le k_R$. By virtue of Lemma~\ref{l1j105} we have for all $j\ge k_R+1$ and all $z\in S(x(\a,\kappa),\a)\cap B(0,R)$ that 
$$
|\^\varphi_A^j(z)|\ge R+(1-\ka)(j-(k_R+1)).
$$
Therefore,
$$
\begin{aligned}
\sum_{j=0}^\infty|g(\^\varphi_A^j(z))|
&=   \sum_{j=0}^{k_R}|g(\^\varphi_A^j(z))|+\sum_{j=k_R+1}^\infty|g(\^\varphi^j(z))| \\
&\lek \sum_{j=0}^{k_R}|\^\varphi_A^j(z)|^{-\frac{p+1}{p}}+\sum_{j=k_R+1}^\infty
      |\^\varphi_A^j(z)|^{-\frac{p+1}{p}} \\
&\lek u^{-\frac{p+1}{p}}x^{-\frac{p+1}{p}}\frac{2R}{1-\ka}+
      \sum_{j=k_R+1}^\infty\big|R+(1-\ka)(j-(k_R+1))\big|^{-\frac{p+1}{p}} \\
&\lek x^{-\frac{p+1}{p}}R+\sum_{j=0}^\infty\big|R+(1-\ka)j)\big|^{-\frac{p+1}{p}} \\
&\comp Rx^{-\frac{p+1}{p}}+\int_0^\infty \(R+(1-\ka)t\)^{-\frac{p+1}{p}}\, dt \\
&= Rx^{-\frac{p+1}{p}}+\frac1{1-\ka}\int_0^\infty (R+s)^{-\frac{p+1}{p}}\, ds \\
&\comp Rx^{-\frac{p+1}{p}}+\int_R^\infty s^{-\frac{p+1}{p}}\, ds \\
&\comp Rx^{-\frac{p+1}{p}}+R^{-\frac1p}.
\end{aligned}
$$
This proves the first assertion of our lemma. For the second assertion, putting $R:=\g x$ in the above formula, we obtain
$$
\sum_{j=0}^\infty|g(\^\varphi_A^j(z))|
\lek \g x^{-\frac{p+1}{p}}+\g^{-\frac1p}x^{-\frac1p}
\comp x^{-\frac{p+1}{p}}\lra 0.
$$
when $x\to+\infty$. This proves the second assertion. For the proof of the third assertion, we note that if $z\in S(x,\pi/2)$, then $\re(z)\ge x$ and, by virtue of Lemma~\ref{l1j105}, we have as above that
$$
\sum_{j=0}^\infty|g(\^\varphi_A^j(z))|
\le \sum_{j=0}^\infty\big|x+(1-\ka)j)\big|^{-\frac{p+1}{p}}
\comp x^{-\frac1p}.
$$
Since $\lim_{x\to+\infty}x^{-\frac1p}=0$, the third assertion is also proved. We are done.
\epf

\sp Now we can prove the following.

\sp
\blem\label{l111.3_3}
Let $\varphi$ be a locally holomorphic simple parabolic map and let $A$ be an attracting direction of $\varphi$. Fix $\a\in(0,\pi)$ and $\kappa \in (0, 1)$. If $z\in S(x(\a,\ka),\a)$ then the limit
$$
\^\varphi_{A,\infty}'(z):=\lim_{n\to\infty}(\^\varphi_A^n)'(z)
$$
\begin{enumerate}

\item exists, belongs to $\C\sms\{0\}$, and the convergence is uniform on every set $S(x(\a,\ka),\a)\cap\{w\in\C:\re(w)\ge t\}$, $t\in\R$, 

\item in particular on every set of the form 
$S(x(\a,\ka),\a,R)$, $R>0$. 

\item Therefore, the function $S(x(\a,\ka),\a)\ni z\longmapsto \^\varphi_{A,\infty}'(z)\in\C\sms\{0\}$ is holomorphic.

\item 
$$
\lim_{z\in S(x(\a,\ka),\a)\atop \re z\to+\infty} \^\varphi_{A,\infty}'(z)=1.
$$

\item So, if in addition $\a\in (0,\pi/2)$, then 
$$
\lim_{S(x(\a,\ka),\a)\ni z\to\infty} \^\varphi_{A,\infty}'(z)=1.
$$
\end{enumerate}
\elem

\bpf
Fix an arbitrary $\e>0$. By Lemma~\ref{l2j107} there exists $x\ge 0$ such that
\beq\label{1_2017_07_31}
\big|(\^\varphi_A^n)'(\xi)-1\big|<\e
\eeq
for all $\xi\in S(x,\pi/2)$ and all integers $n\ge 0$. By Lemma~\ref{l1j105} there exists $k\ge 0$ so large that $\Re\^\phi^n(z)\ge x$, i.e. $\^\phi^n(z)\in S(x,\pi/2)$, for all 
$$
z\in S(x(\a,\ka),\a)\cap\{w\in\C:\re(w)\ge t\}
$$ 
and all $n\ge k$. Then for every $i\ge k$ and every $j\ge 0$ we have
$$
\lt|\frac{(\^\varphi_A^{i+j})'(z)}{(\^\varphi_A^i)'(z)}-1\rt|
=\big|(\^\varphi_A^j)'(\^\phi^i(z))-1\big|
<\e.
$$
Therefore, the sequence $\((\^\varphi_A^n)'(z)\)_{n=0}^\infty$ is uniformly quotient Cauchy on  
$$
S(x(\a,\ka),\a)\cap\{w\in\C:\re(w)\ge t\}
$$ 
and the first assertion of our lemma is proved. Assertions (2) and (3) then follow immediately while (4) is a direct consequence of \eqref{1_2017_07_31} and (5) is a direct consequence of (4).
\epf

\sp It immediately follows from this lemma that 
\beq\label{1ch11.3_4}
\^\varphi_{A,\infty}'(z)=\^\varphi_{A,\infty}'(\^\varphi(z))\^\phi'(z)
\eeq
for all $z\in S(x(\a,\ka),\a)$. In particular, the holomorphic function $S(x(\a,\ka),\a)\ni z\longmapsto \^\varphi_\infty'(z)$, and its modulus as well, are not constant. 

\sp For every $x\in (0,\infty)$, $\a\in (0,\pi)$ and $R>0$ let 
\beq\label{120190627}
S_0(x,\a):=H(S(x,\a))\index{(S)}{$S_0(x,\a)$},
\eeq
\beq\label{{MU2}p208}
S_\varphi^A(x,\a):=\rho_A\circ H(S(x,\a))
=\rho_A(S_0(x,\a)).
\eeq\index{(S)}{$S_\phi^A(x,\a)$}
and 
\beq\label{220190627}
S_0(x,\a,R):=S_0(x,\a)\cap B^c(0,R)=H\(S(x,\a,R^{-p})\)
\eeq
while
\beq\label{320190627}
S_\varphi^A(x,\a,R):=S_\varphi^A(x,\a)\cap B^c(\om,R)=\rho_A\circ H\(S(x,\a,|a|^{-1}R^{-p})\).
\eeq
The regions $S_0(x,\a)$ and $S_\varphi^A(x,\a)$ look like flower petals
which are respectively symmetric about the rays $(0,\infty)$ and 
$A=\om+\sqrt[p]{-a^{-1}(0,\infty)}$, contain initial segments of these rays, and form with them two ``angles'' of measure $\a/p$ at the points $0$ and $\om$
respectively.

\sp

\fr The following proposition is an immediate consequence of
Proposition~\ref{p1j104}.

\sp

\bprop\label{p1j108} 
Let $\varphi$ be a locally holomorphic simple parabolic map at $\om\in\C$ and let $A$ be an attracting direction of $\varphi$. 
Fix $\a\in(0,\pi)$ and $\kappa \in (0, 1)$. Then for every $x\ge x(\a,\kappa),\a)$, all iterates $\varphi^n: S^A_\varphi(x,\a)\lra \mathbb C$, $n\ge 0$, are well defined,
$$
\varphi^n(S^A_\varphi(x,\a))\sbt
S^A_\varphi(x, \a).
$$ 
and the sequence $\(\vp^n(z)\)_{n=0}^{\infty}$ converges to $\om$ uniformly on $S_\varphi^A(x,\a)$.
\eprop

\sp

A direct elementary calculation based on \eqref{2j101} yields this.

\blem\label{dp4.3A}
If $\varphi$ is a locally holomorphic simple parabolic map at $\om\in\C$ and  $A$ is an attracting direction of $\varphi$, then
for every $0\le \a<\frac{\pi}{2p(\om)}$ there exists $R_\a(\om)>0$ such that
$$
|\varphi(z)-\om|<|z-\om| 
\  \  {\rm and} \  \
|\varphi'(z)|<1
$$
for all $z\in \(S_\phi^A(x,\a)\cap B(\om,R_\a(\om))\)\sms\{\om\}$. 
\elem

\sp We now shall prove the following.

\sp

\bprop\label{p1j109} 
Let $\varphi$ be a locally holomorphic simple parabolic map at $\om\in\C$ and let $A$ be an attracting direction of $\varphi$.
Fix $\a\in(0,\pi)$ and $\kappa \in (0, 1)$. Recall that $p=p(\om)$. For all
$z\in S^A_\varphi(x(\a,\kappa),\a)$, we  have
$$
\lim_{n \to \infty}\( n^{\frac{p+1}{p}}|\varphi^{n+1}(z)-
\varphi^{n}(z)|\)=|a|^{-\frac{1}{p}}
$$
and 
$$
\lim_{n \to\infty}\(n^{\frac{1}{p}}|\varphi^{n}(z)-\om|\)
=|a|^{-\frac{1}{p}}.
$$
In addition, in the two limits above, the convergence is uniform on
the set $S^A_\varphi(x(\a,\kappa),\a, R)$ for every $R>0$. 
\eprop

\bpf  We have by Lemma~\ref{l2j105} for all $\xi\in S_0(x(\a,\kappa),\a,R)$ that
$$
\begin{aligned}
|\varphi^n_{A,0}(\xi)|
&=| H \circ \^\varphi^n \circ H^{-1}(\xi)|\\
&=\left|\^\varphi^n(H^{-1}(\xi)))\right|^{-\frac{1}{p}}\\
&=\left|H^{-1}(\xi)+n+O_{-R^{-p}}\(\max\{n^{1-\frac1p},\log(n+1)\}\)\right|^{-\frac{1}{p}}.
\end{aligned}
$$
Therefore, as $|H^{-1}(\xi)|\le R^{-p}$, we get that
$$
\begin{aligned}
\lim_{n \to \infty}n^{\frac{1}{p}}|\varphi^n_{A,0}(\xi)|
=\lim_{n \to \infty}\Big|n^{-1}H^{-1}(\xi)+1 + 
  O_{-R^{-p}}\(\max\{n^{-\frac1p}, n^{-1}\log(n+1)\}\)\Big|^{-\frac{1}{p}} 
=1,
\end{aligned}
$$ 
and the convergence is uniform on $\xi\in S_0(x(\a,\kappa),\a,R)$. 
Since for all points $z\in S^A_\varphi(x(\a,\kappa),\a,R)$ we have
$$
\varphi^n(z)-\om 
=\rho(\varphi^n_{A,0}(\rho^{-1}(z)))-\om 
=\sqrt[p]{-a^{-1}}\varphi^n_{A,0}(\rho^{-1}(z)),
$$ 
the second formula of our proposition follows along with the appropriate uniform convergence. Turning attention to the first formula, we have by (\ref{2j101}), by the first formula of the proposition, and by that the last assertion of Proposition~\ref{p1j109}, that 
$$
\begin{aligned}
&\lim_{n\to\infty}n^{\frac{p+1}{p}}|\varphi^{n+1}(z)-\varphi^n(z) |\\
&=\lim_{n\to\infty}n^{\frac{p+1}{p}}\big|a(\varphi^n(z)-\om)^{p+1}+O(|\varphi^n(z)-\om|^{p+2})\big|\\
&=|a|\lim_{n \to \infty}(n^{\frac{1}{p}}|\varphi^n(z)-\om|)^{p+1}+
\lim_{n \to \infty}O\((n^{\frac{1}{p}}|\varphi^n(z)-\om|)^{p+1}|\varphi^n(z)-\om|\)\\
& =|a|\lim_{n \to \infty}(n^{\frac{1}{p}}|\varphi^n(z)-\om|)^{p+1}\\
&=|a||a|^{-\frac{p+1}{p}}
 =|a|^{-\frac{1}{p}}.
\end{aligned}
$$
where the big $O$ symbol represents an absolute constant depending only on the Taylor series expansion (\ref{2j101}) of $\phi$; moreover the convergence is uniform on $\xi\in S_0(x(\a,\kappa),\a,R)$. The proof is complete. \endpf

\sp For every $z\in S^A_\varphi(x(\a,\kappa),\a)$ put 
\beq\lab{4_2017_11_28}
\varphi_{a,\infty}'(z):=\^\varphi_{A,\infty}'(\rho\circ H^{-1}(z)).
\eeq
Our last result about the iterates of $\varphi$ is this.

\sp\bprop\label{p1j111} 
Let $\varphi$ be a locally holomorphic simple parabolic map at $\om\in\C$ and let $A$ be an attracting direction of $\varphi$.
Fix $\a \in (0, \pi)$ and $\kappa \in (0,1)$. Recall that $p=p(\om)$. Then for every $z\in S^A_\varphi(x(\a,\kappa),\a)$
\beq\label{1_2017_08_04}
\lim_{n\to+\infty}n^{\frac{p+1}{p}}|(\varphi^n)'(z)|
=|a|^{-\frac{p+1}{p}}|z-\om|^{-(p+1)}|\varphi_{A,\infty}'(z)|,
\eeq
and the convergence is uniform on $S^A_\varphi(x(\a,\kappa),\a,R)$ for every $R>0$. In addition, if $\a\in (0,\pi/2)$, then
\beq\label{2_2017_08_04}
\lim_{S^A_\varphi(x(\a,\kappa),\a)\ni z\to\om}\varphi_{A,\infty}'(z)=1.
\eeq
\eprop

\bpf By virtue of Lemma~\ref{l2j105} we get for all $z
\in S_0(x(\a,\kappa),\a)$ that
\beq\label{1_2015_10_10}
\begin{aligned}
|z&|^{p+1}n^{\frac{p+1}{p}}|(\varphi_{A,0}^n)'(z)|= \\
&=n^{\frac{p+1}{p}}|H'(\^\varphi_A^n(H^{-1}(z)))(\^\varphi_A^n)'(H^{-1}(z))(H^{-1})'(z)|  \\
&= |z|^{p+1}n^{\frac{p+1}{p}}\frac{1}{p}|\^\varphi_A^n(H^{-1}(z))|^{-\frac{p+1}{p}}|(\^\varphi_A^n)'(H^{-1}(z))|p|z|^{-(p+1)}\\
&=n^{\frac{p+1}{p}}|\^\varphi_A^n(H^{-1}(z))|^{-\frac{p+1}{p}}|(\^\vp_A^n)'(H^{-1}(z))|\\
&=n^{\frac{p+1}{p}}|(\^\vp_A^n)'(H^{-1}(z))|\Big|H^{-1}(z)+n +O_{\Re(H^{-1}( z))}\(\max\{n^{1-\frac1p}, \log(n+1)\}\)\Big|^{-\frac{p+1}{p}}  \\
&=|(\^\vp_A^n)'(H^{-1}(z))|\lt|\frac{H^{-1}(z)}{n}+1 +O_{\Re(z^{-p})}\(\max\{ n^{-\frac1p},n^{-1}\log(n+1)\}\)\rt|^{-\frac{p+1}{p}}  
\end{aligned}
\eeq
So, if $z\in S^A_\varphi(x(\a,\kappa),\a,R)$, then  formula \eqref{1_2015_10_10} takes then on the form
$$
\begin{aligned}
|z|^{p+1}&n^{\frac{p+1}{p}}|(\varphi_{A,0}^n)'(z)|= \\
&=|(\^\vp_A^n)'(H^{-1}(z))|\lt|\frac{H^{-1}(z)}{n}+1 +O_{-R^{-p}}\(\max\{ n^{-\frac1p},n^{-1}\log(n+1)\}\)\rt|^{-\frac{p+1}{p}},
\end{aligned}
$$
Therefore,
\beq\label{1_2017_02_14}
\begin{aligned}
n^{\frac{p+1}{p}}&|(\varphi^n)'(z)|
=n^{\frac{p+1}{p}}|(\varphi^n_{A,0})'(\rho^{-1}(z))|=\\
&=|\rho^{-1}(z)|^{-(p+1)}
\big|(\^\vp_A^n)'(\rho\circ H)^{-1}(z))\big|\cdot \\
& \ \  \  \  \  \  \  \  \  \  \  \  \  \  \  \   \  \  \  \  \cdot\lt|\frac{(\rho\circ H)^{-1}(z)}{n}+1 +O_{-R^{-p}}\(\max\{ n^{-\frac1p},n^{-1}\log(n+1)\}\)\rt|^{-\frac{p+1}{p}} \\
&=|a|^{-\frac{p+1}{p}}|z-\om|^{-p+1}
\big|(\^\vp_A^n)'(\rho\circ H)^{-1}(z))\big|\cdot \\
& \ \  \  \  \  \  \  \  \  \  \  \  \  \  \  \   \  \  \  \  \cdot\lt|\frac{(\rho\circ H)^{-1}(z)}{n}+1 +O_{-R^{-p}}\(\max\{ n^{-\frac1p},n^{-1}\log(n+1)\}\)\rt|^{-\frac{p+1}{p}}.
\end{aligned}
\eeq
Since also
$$
(\rho\circ H)^{-1}\(S^A_\varphi(x,\a,R)\)=S\(x,\a,(|a|R)^{-1/p}\),
$$
the proof is thus concluded by applying Lemma~\ref{l111.3_3}. We are done. 
\endpf  

\sp The next theorem establishes the the existence of functions called Fatou coordinates around simple parabolic fixed points. Its primary application in this book is used in the description of basins of attractions to parabolic points. This theorem is however of local character and we formulate and prove it here. For every $\a\in (0,\pi)$ put
$$
Q_\a:=S(x(\a,1/2),\a) \  \  \text{ and } \  \  Q_{\alpha,t}:=Q_\alpha \cap \{z\in \C: \re(z)\geq t\}.
$$

\bthm\label{t1ms127}
Let $\varphi$ be a locally holomorphic simple parabolic map at $\om\in\C$ and let $A$ be an attracting direction of $\varphi$. Then
for every $\a\in (0,\pi)$ there exists a unique, up to an additive constant, holomorphic function
$\widetilde{F}_\alpha: Q_\alpha \lra \C$ such that
\begin{eqnarray}\label{1ms127}
\widetilde{F}_\alpha\circ \widetilde{\vp}_A= \widetilde{F}_\alpha+1
\end{eqnarray}
and
\begin{eqnarray}\label{2ms127}
\lim_{z\to+\infty} \frac{\widetilde{F}_\alpha(z)}{z}=1
\end{eqnarray}
uniformly on $Q_{\alpha,t}$ for every $t\in \R$. In addition, for $t\in \R$ large enough, the map 
$$
\widetilde{F}_\alpha|_{Q_{\alpha, t}}:Q_{\alpha, t}\lra \C
$$ 
is univalent.
\ethm

\bpf First of all notice that, becuase of Proposition~\ref{p1j104} and \eqref{1j105}, we have that
\beq\label{2_2017_02_14}
\widetilde{\vp}_A(Q_\alpha)\sbt Q_\alpha \  \  \text{ and } \  \  \widetilde{\vp}_A(Q_{\alpha,t})\sbt Q_{\alpha,t}
\eeq
for all $t\in \R$. In particular, all iterates $\widetilde{\vp}_A^n$, $n\ge 0$, of $\widetilde{\vp}_A$ are well defined on $Q_\a$ and $ Q_{\alpha,t}$.
For uniqueness, suppose that a holomorphic function $\widetilde{F}_\alpha: Q_\alpha\to\C$ satisfies (\ref{1ms127}) and (\ref{2ms127}). Fix an arbitrary  $\varepsilon>0$. In view of the definition of $Q_\alpha$ and (\ref{2ms127}) there exist $\b\in (0,\a)$ small enough and $s>\max\{2,x(\a,1/2)\}$ large enough such that if $z\in Q_{\b,s}$, then
\[
\overline{B}(z, |z|/2)\subseteq Q_\alpha\ \text{ and }\ \left|\frac{{\widetilde{F}_\alpha}(z)}{z}-1\right|<\frac{\varepsilon}{3}.
\]
Cauchy's Integral Formula then yields,
\begin{eqnarray}\label{1ms129}
|\widetilde{F}_\alpha(z)-1|&=& \left|\frac{1}{2\pi i}\int_{\partial B(z, |z|/2)} \frac{\widetilde{F}_\alpha(\xi)}{(\xi-z)^2}\ d\xi-1\right|\nonumber=
\left|\frac{1}{2\pi i}\int_{\partial B(z, |z|/2)} \frac{\widetilde{F}_\alpha(\xi)-\xi}{(\xi-z)^2}\ d\xi\right|\nonumber\\&\leq&
\frac{1}{2\pi}\int_{\partial \overline{B}(z, |z|/2)} \frac{|\widetilde{F}_\alpha(\xi)-\xi|}{|\xi-z|^2}\ |d\xi|\nonumber\\&=&\frac{2}{\pi|z|^2}\int_{\partial \overline{B}(z, |z|/2)} |\widetilde{F}_\alpha(\xi)-\xi|\ |d\xi|\nonumber\\&\leq&
\frac{2}{\pi|z|^2}\int_{\partial \overline{B}(z, |z|/2)} \frac{\varepsilon}{3}|\xi|\ |d\xi|\nonumber\\&\leq&\frac{2}{\pi|z|^2}
\pi|z|\frac{\varepsilon}{3}\frac{3|z|}{2} = \varepsilon.
\end{eqnarray}
Now fix also $w\in  Q_{\b,t}$. Using (\ref{1ms127}), we obtain that
\begin{eqnarray*}
\widetilde{F}_\alpha(z)-\widetilde{F}_\alpha(w)
&=&\widetilde{F}_\alpha(\widetilde{\vp}_A^n(z))-\widetilde{F}_\alpha(\widetilde{\vp}_A^n(w))\\
&=&\int_{\Delta_n}\widetilde{F}_\alpha'(\xi)\ d\xi 
=\int_{\Delta_n}(1+(1-\widetilde{F}_\alpha'(\xi)))\ d\xi\\
&=&\widetilde{\vp}_A^n(z)-\widetilde{\vp}_A^n(w)+\int_{\Delta_n}(1-\widetilde{F}_\alpha'(\xi))\ d\xi,
\end{eqnarray*}
where $\Delta_n$ is the oriented segment from $\widetilde{\vp}_A^n(w)$ to $\widetilde{\vp}_A^n(z)$. So, applying (\ref{1ms129}), it follows that
\begin{eqnarray*}
\left|\widetilde{F}_\alpha(z)-\widetilde{F}_\alpha(w)-(\widetilde{\vp}_A^n(z)-\widetilde{\vp}_A^n(w))\right|
&\leq&
\left|\int_{\Delta_n}(1-\widetilde{F}_\alpha'(\xi))\ d\xi\right|\\&\leq& \int_{\Delta_n}\left|\widetilde{F}_\alpha'(\xi)-1\right|\ \left|d\xi\right|\\&\leq&\varepsilon \left|\Delta_n\right|
=\varepsilon \left|\widetilde{\vp}_A^n(z)-\widetilde{\vp}_A^n(w)\right|.
\end{eqnarray*}
Thus,
\[
\left|\frac{\widetilde{F}_\alpha(z)-\widetilde{F}_\alpha(w)}{\widetilde{\vp}_A^n(z)-\widetilde{\vp}_A^n(w)}-1\right|\leq \varepsilon.
\]
Hence, the limit of the sequence $(\widetilde{\vp}_A^n(z)-\widetilde{\vp}_A^n(w))_{n=0}^\infty$ exists and
\[
\widetilde{F}_\alpha(z)-\widetilde{F}_\alpha(w)=\lim_{n\to\infty}(\widetilde{\vp}_A^n(z)-\widetilde{\vp}_A^n(w)).
\]
The uniqueness part is thus established. This argument also gives us a hint as to how to prove the existence part. Keep $t\in\R$ fixed and also fix $w\in Q_{\b,s}$ as above. Consider the functions
\begin{eqnarray}\label{2ms129}
Q_\a\ni z\longmapsto\psi_n(z)
:=\widetilde{\vp}_A^n(z)-\widetilde{\vp}_A^n(w), \  \  n\geq0.
\end{eqnarray}
We then have that
\[
\psi_{n+1}(z)-\psi_n(z)
=(\widetilde{\vp}_A(\widetilde{\vp}_A^n(z))-(\widetilde{\vp}_A^n(z)+1))-(\widetilde{\vp}_A(\widetilde{\vp}^n(w))-(\widetilde{\vp}_A^n(w)+1)).
\]
So, making use of the Mean Value Inequality, in virtue of Lemma~\ref{l1j101}, formula \eqref{2_2017_02_14} and \eqref{f2_2015_09_08}, we obtain
\begin{eqnarray*}
\left|\psi_{n+1}(z)-\psi_n(z)\right|
&=& \left|\widetilde{\vp}_A(\widetilde{\vp}_A^n(z))-(\widetilde{\vp}_A^n(z)+1))-(\widetilde{\vp}_A(\widetilde{\vp}_A^n(w))-(\widetilde{\vp}_A^n(w)+1)\right|\\
&=&\left|B(H(\widetilde{\vp}_A^n(z)))-B(H(\widetilde{\vp}_A^n(w)))\right|\\
&=&\lt|B\circ H(\widetilde{\vp}_A^n(z))-B\circ H(\widetilde{\vp}_A^n(w))\rt|\\
&\leq& Mp^{-1}\sup\{|\xi|^{-\frac{p+1}p}:\xi\in\Delta_n\}\left|\widetilde{\vp}_A^n(z)-\widetilde{\vp}_A^n(w)\right|\\
&\leq& C_t n^{-\frac{p+1}p}\left|\widetilde{\vp}_A^n(z)-\widetilde{\vp}_A^n(w)\right|\\
&=& C_t\left|\psi_n(z)\right|n^{-\frac{p+1}p}.
\end{eqnarray*}
with some constant $C_t\in (0,+\infty)$ independent of $z$ and $n$ but depending (in general) on $t$. Equivalently,
\begin{eqnarray}\label{1ms131}
\left|\frac{\psi_{n+1}(z)}{\psi_n(z)}-1\right|\leq C_t n^{-\frac{p+1}p}.
\end{eqnarray}
Since the series $\sum_{n=1}^\infty n^{-\frac{p+1}p}$ converges, the estimate (\ref{1ms131}) implies that the sequence $(\psi_n)_{n=1}^\infty$ converges uniformly on $Q_{\a,t}$. So, the limit
\begin{eqnarray}\label{3ms131}
\widetilde{F}_\alpha=\lim_{n\to\infty} \psi_n
\end{eqnarray}
is a holomorphic function on $Q_{\alpha, t}$ and letting $t\to-\infty$, we conclude that \eqref{3ms131} defines a holomorphic function from $Q_\a$ to $\C$. Moreover, due to (\ref{3ms131}) and Lemma~\ref{l1j101} combined with \eqref{f2_2015_09_08}, for all $z\in Q_\a$ we obtain the following equality: 
\begin{eqnarray*}
\widetilde{F}_\alpha\left(\widetilde{\vp}_A(z)\right)
&=&\lim_{n\to\infty}\psi_{n}\left(\widetilde{\vp}_A(z)\right)
=\lim_{n\to\infty}\(\widetilde{\vp}_A^{n+1}(z)-\widetilde{\vp}_A^n(w)\)\\
&=&\lim_{n\to\infty}\left(\left(\widetilde{\vp}_A^{n+1}(z)-\widetilde{\vp}_A^{n+1}
  (w)\right)+\left(\widetilde{\vp}_A^{n+1}(w)-\widetilde{\vp}_A^{n}(w) \right)\right) \\
&=& \lim_{n\to\infty}\psi_{n+1}(z)+1=\widetilde{F}_\alpha(z)+1.
\end{eqnarray*}
This proves (\ref{1ms127}).

In order to obtain (\ref{2ms127}), we must do a little bit more work. Keep again $t\in\R$ fixed. By immediate induction, based on the first and on the third assertion of Lemma~\ref{l1j101}, we deduce that
\begin{eqnarray}\label{1ms130}
\left|\widetilde{\vp}_A^n(z)-(z+n)\right|\leq L_{t,n}|z|^{-1/p},
\end{eqnarray}
where $L_{n,t}\geq1$ is some constant independent of $z\in Q_{\a,t}$ but may depend on $t$ and $n$. Fixing now $\varepsilon>0$, in virtue of (\ref{3ms131}) and \eqref{1ms131}, there exists $k\geq1$ so large that for all $z\in Q_{\a,t}$ we have that,
\[
\left|\frac{\widetilde{F}_\alpha(z)}{\psi_k(z)}-1\right|<\frac{\varepsilon}{2}.
\]
Equivalently,
\[
\left|\frac{\widetilde{F}_\alpha(z)/z}{\psi_k(z)/z}-1\right|<\frac{\varepsilon}{2}.
\]
But, by (\ref{1ms130}) along with the definition \eqref{2ms129}, we have that $\lim_{Q_t\ni z\to+\infty}\frac{\psi_k(z)}{z}=1$. Hence,
\[
\left|\frac{\widetilde{F}_\alpha(z)}{z}-1\right|<\varepsilon
\]
for all $z\in Q_{\a,t}$ with sufficiently real parts. Therefore,
\[
\lim_{Q_t\ni z\to+\infty}\frac{\widetilde{F}_\alpha(z)}{z}=1,
\]
and formula (\ref{2ms127}) is established. 

It only remains to show that the function 
$$
\widetilde{F}_\alpha|_{Q_{\alpha, t}}:Q_{\a,t}\to\C
$$ 
is univalent for all $t\in \R$ large enough. As the main ingredient of our proof of this statement, we shall prove the following.

\sp{\bf Claim:} The map $\widetilde{\vp}_A|_{Q_{\alpha, t}}:Q_{\alpha,t}\lra\C$ is univalent for all $t\in \R$ large enough. 

\bpf Indeed, fix  $t>4$ so large that
$$
2M(t-2)^{-\frac1p}<1,
$$
where $M>0$ is the constant coming from Lemma~\ref{l1j101}. So, if $\widetilde{\vp}(y)=
\widetilde{\vp}(x)$ for some $x, y\in Q_{\a,t}$, then it follows from this lemma that
\beq\label{f1_2015_09_09}
|y-x|\le 2Mt^{-\frac1p}<1.
\eeq
On the other hand if $|z-x|=2$, then it also follows from Lemma~\ref{l1j101} that
$$
|\widetilde{\vp}_A(z)-(z+1)|
\le M|z|^{-\frac1p}
\le M(|x|-2)^{-\frac1p}
\le M(t-2)^{-\frac1p}
<1
<|z+1|,
$$
where the last inequality follows from the fact that $|z+1|\ge |z|-1\ge |x|-2-1=|x|-3\ge \Re(x)-3\ge t-3>1$. Thus Rouch\'{e}'s Theorem asserts that the function $z\mapsto \widetilde{\vp}_A(z)$ has the same number of zeros in $B(x,2)$ as the function $z\mapsto z+1$. Therefore, the map $z\mapsto \widetilde{\vp}(z)$ is 1-to-1 on $B(x,2)$. Along with 
\eqref{f1_2015_09_09}, this yields $y=x$, and the proof of the claim is complete.
\epf

\sp In view of this claim and of the uniform convergence of the sequence $\(\psi_n\)_{n=1}^\infty$ on $Q_{\a,t}$, by applying Hurwitz's Theorem we see that the function $\widetilde{F}_\alpha|_{Q_{\alpha,t}}$ is also univalent for all $t\in \R$ large enough. This completes the proof of our theorem.
\epf

\sp The functions $\widetilde{F}_\a$ are called the Fatou coordinates of the parabolic map $\widetilde{\vp}_A$. We shall now establish one important property of the functions $\widetilde{F}_\alpha: Q_\alpha \lra \C$.

\blem\label{l1ms130}
Let $\varphi$ be a locally holomorphic simple parabolic map at $\om\in\C$ and let $A$ be an attracting direction of $\varphi$.
For every $\alpha\in [\pi/2, \pi)$ and every $t\in\R$ there exists $s:=s(t)>0$ such that
\[
\widetilde{F}_\alpha(Q_{\alpha,t})\supseteq \{z\in \C:\re(z)>s\}.
\]
\elem

\bpf
By virtue of formula (\ref{2ms127}), there exists $s\geq \max\{0, 4(x(\a,1/2)+t)\}$ so large that
\[
\left|\frac{\widetilde{F}_\alpha(z)}{z}-1\right|\leq \frac18
\]
whenever $\re(z)>s/4$. Take now an arbitrary number $w\in \C$ with $\Re(w)>s$. Observe that
\beq\label{1_2017_02_16}
B(w,2|w|/3)\sbt Q_{\alpha,t}.
\eeq
If $z\in \partial B(w,|w|/2)$, then
\begin{eqnarray*}
\left|(\widetilde{F}_\alpha(z)-w)-(z-w)\right|&=&\left|\widetilde{F}_\alpha(z)-z\right|=|z|\left|\frac{\widetilde{F}_\alpha(z)}{z}-1\right|\\&\leq& \frac{|z|}{8}\leq \frac{|w|+|w|/2}{8}\leq\frac3{16}|w|<\frac14|w|.
\end{eqnarray*}
Hence,
\[
|w-z|=\frac{|w|}2>\big|\widetilde{F}_\alpha(z)-w-(z-w)\big|.
\]
Thus Rouch\'{e}'s Theorem asserts that the function $z\mapsto \widetilde{F}_\alpha(z)-w$ has the same number of zeros in $B(w,|w|/2)\subseteq Q_{\alpha,t}$ as the function $z\mapsto z-w$. But the latter function has exactly one zero, namely $z=w$. Thus, there exists $z\in B(w, |w|/2)\subseteq Q_{\alpha,t}$ such that $\widetilde{F}_\alpha(z)=w$. This means that 
$$
\widetilde{F}_\alpha(Q_{\alpha,t})\supseteq \{z\in \C:\Re(z)>s\}.
$$ 
The proof is complete.
\epf

\sp Fix $\a \in (0, \pi)$, $\kappa=1/2$, and denote by 
\beq\label{1_2017_29A}
S^1_a(\om, \a),\ldots,S^{p(\om)}_a(\om,\a)
\eeq
the corresponding attracting sectors\index{(N)}{attracting sectors} for $\vp$ defined in (\ref{{MU2}p208}) with $x:=x(\a,1/2)$. Denote also by 
\beq\label{1_2017_29B}
S^1_r(\om, \a),\ldots,S^{p(\om)}_r(\om,\a)
\eeq
the analogous attracting sectors for $\vp^{-1}$. We will also frequently call them the repelling sectors\index{(N)}{repelling sectors} for $\vp$ (at the parabolic point $\om$). As an immediate consequence of their definition, we get the following.

\blem\label{l120190617}
If $\varphi$ is a locally holomorphic simple parabolic map at $\om\in\C$, then for every $\alpha\in (0,\pi)$, both collections 
$$
\Big\{S^1_a(\om, \a),\ldots,S^{p(\om)}_a(\om,\a)\Big\} 
\  \  \  {\rm and} \  \  \
\Big\{S^1_r(\om, \a),\ldots,S^{p(\om)}_r(\om,\a)\Big\} 
$$
consist of mutually disjoint sets.
\elem

Putting
\begin{eqnarray}\label{1ms133}
F_\alpha:=\widetilde{F}_\alpha\circ\rho\circ H^{-1},
\end{eqnarray}
and passing to the original parabolic function $\vp$, as an immediate consequence of Proposition~\ref{p1j104}, Theorem \ref{t1ms127}, and Lemma \ref{l1ms130}, we obtain the following result.

\bthm\label{t1ms135}
If $\varphi$ is a locally holomorphic simple parabolic map at $\om\in\C$, then for every $\alpha\in (0,\pi)$, we have that
\begin{eqnarray}\label{1ms135}
\vp(S^j_a(\om,\a))\subset S^j_a(\om,\a)
\end{eqnarray}
for all $j=1, \ldots, p(\om)$ and
\beq\label{3_2017_02_17}
F_\alpha\circ \vp=F_\alpha +1
\eeq
on $S^j_a(\om,\a)$. Also, 
\begin{enumerate}

\item With approriate normalizations (fixing the values of $F_\a$ at one point), if $0<\beta,\gamma<\pi$, then 
$$
F_{\gamma}|_{S^j_a(\om,\b)\cap S^j_a(\om,\g)}=F_\beta|_{S^j_a(\om,\b)\cap S^j_a(\om,\a)}.
$$
\item  Increasing the point $x(\a,1/2)\in\R$ appriopriately, we will have that the map   $F_\alpha|_{S^j_a(\om,\a)}$ is univalent.

\sp\item If $\a\in[\pi/2,\pi)$, then there exists $s>0$ such that 
\[
F_\alpha(S^j_a(\om,\a))\supset\{z\in \C:\re(z)>s\}.
\] 
\end{enumerate}
\ethm 

\sp

\section[Fatou Components of (General) Meromorphic Functions  II;
Leau--Fatou] {Fatou Components of (General) Meromorphic Functions  II; \\ Leau--Fatou Flower Petals}\label{FCGMFII}

In this section we continue our study of parabolic fixed and periodic points undertaken in the previous section. However, up to this point our considerations have had entirely local character in the sense that the map $f$ (denoted so far by $\vp$) was defined only on a small neighborhood of its simple rationally indifferent fixed point $\om$. Throughout this section $f$ is instead global, that is 
$$
f:\C\lra\oc
$$ 
is a meromorphic function  and it has a simple rationally indifferent fixed point $\om$. We define in this section the corresponding Leau--Fatou Flower Petals and study their dynamical and topological structure in greater detail.

It follows from Proposition~\ref{p1j108} that for every $j\in\{1, \ldots, p(\om)\}$, $S^j_a(\om,\a)$ is a subset of the Fatou set $F(f)$ of $f$. We then define the set $A_j^{*}(\om)$ to be the connected component of the Fatou set $F(f)$ that contains $S^j_a(\om,\a)$ for some (equivalently, for every), $\alpha\in (0, \pi)$. Let $\d_f>0$ be so small that there exists 
$$
f^{-1}_\om:B(\om,\d_f)\lra \C,
$$
a unique holomorphic inverse branch of $f$ that sends $\om$ to $\om$. We require in addition that the ma p
$$
f|_{B(\om,2\d_f)}
$$ 
is injective. Recall that the sets $S^j_r(\om,\a)$, $j\in\{0,\ldots, p(\om)-1\}$, are the attracting sectors for the map $f^{-1}_\om:B(\om,\d_f)\lra\C$. Invoking Proposition~\ref{p1j108} again, we obtain the following.

\blem\label{l1_2015_09_10}
If $f:\C\lra\oc$ is a meromorphic function and $\om$ is its simple rationally indifferent fixed point, then 
$$
f(A_j^{*}(\om))=A_j^{*}(\om)
$$ 
for every $j\in\{1, \ldots, p(\om)\}$ and
$$
\lim_{n\to\infty}f^n(z)=\om
$$
for every $z\in \bu_{j=0}^{p(\om)-1}A_j^{*}(\om)$.
\elem

\fr We further define
\[
A_j(\om):=\bigcup_{n=0}^\infty f^{-n}(A_j^{*}(\om)).
\]
The open set $A_j(\om)$ is called the basin of attraction to $\om$ in the direction determined by $j$, and $A^*_j(\om)$ is called the immediate basin of attraction to $\om$ in this direction. We shall now prove the following.

\blem\label{t3ms134A}
If $f:\C\lra\oc$ is a meromorphic function and $\om$ is a simple parabolic fixed point of $f$, then for every $\alpha\in (0, \pi)$, we have that
\[
A^*(\om):=\bu_{j=1}^{p(\om)} A^*_j(\om)\sbt \bu_{n=0}^\infty f^{-n}\lt(\bu_{j=1}^{p(\om)} S^j_a(\om,\a)\rt).
\]
\elem

\bpf
Suppose, by way of contradiction, that this inclusion does not hold. This means that there exist some $i\in\{1, \ldots, p(\om)\}$ and some
\beq\label{2_2015_10_06}
z\in A^*_i(\om)\setminus \bigcup_{n=0}^\infty f^{-n}\lt(\bu_{j=1}^{p(\om)}S^j_a(\om,\a)\rt).
\eeq
Since $f(A_j^{*}(\om))=A_j^{*}(\om)$ and since, by Lemma~\ref{l1_2015_09_10}, $\lim_{k\to\infty} f^k(z)=\om$, passing to a sufficiently large iterate of $z$, we may assume without loss of generality that
\beq\label{1_2015_10_06}
f^k(z)\in B(\om,\d_f)
\eeq
for all integers $k\ge 0$. Again since $\lim_{k\to\infty} f^k(z)=\om$, we conclude from \eqref{2_2015_10_06} that 
\beq\label{3_2015_10_06}
f^n(z)\in\bigcup_{j=1}^{p(\om)}S^j_r(\om,\pi-\a).
\eeq
for all $n\ge 0$ large enough, say $n\ge q\ge 1$.
But, in light of Proposition~\ref{p1j108}, for all $k\ge 0$ large enough, we have that
\beq\label{4_2015_10_06}
f_\om^{-k}\lt(\bigcup_{j=1}^{p(\om)}S^j_r(\om,\pi-\a)\rt)\subset B(\om,|z-\om|/2).
\eeq
Since by virtue of \eqref{1_2015_10_06}, we have that $z=f_\om^{-k}(f^k(z))$ for all $k\ge 0$, we thus obtain from \eqref{3_2015_10_06} and \eqref{4_2015_10_06} that $z\in B(\om, |z-\om|/2)$. This contradiction finishes the proof.
\epf

\sp\fr We call $A^*(\om)$ the immediate basin of attraction to the parabolic fixed point $\om$.

\bthm\label{t3ms134}
If $f:\C\to\oc$ is a meromorphic function and $\om$ is a simple parabolic fixed point of $f$, then the sets $\(A^*_j(\om)\)_{j=1}^{p(\om)}$, called Leau--Fatou petals or Leau domains, are mutually disjoint and, moreover, for every $j\in\{1, \ldots, p(\om)\}$ and every $\alpha\in (0, \pi)$, we have that
\[
A^*_j(\om)\sbt \bigcup_{n=0}^\infty f^{-n}\(S^j_a(\om,\a)\).
\]
\ethm

\bpf
Since $\lt(\bu_{n=0}^\infty f^{-n}\(S^j_a(\om,\a)\)\rt)_{j=1}^{p(\om)}$ is a finite collection of mutually disjoint open sets, and since all the sets $\(A^*_j(\om)\)_{j=1}^{p(\om)}$ are connected, making use of Lemma~\ref{t3ms134A}, it follows that each set $A^*_i(\om)$ is contained in exactly one set of the form $\bu_{n=0}^\infty f^{-n}\(S^j_a(\om,\a)\)$. Since, $A^*_i(\om)\cap \bu_{n=0}^\infty f^{-n}\(S^i_a(\om,\a)\)\ne\es$, we thus conclude that 
$$
A^*_i(\om)\sbt \bu_{n=0}^\infty f^{-n}\(S^i_a(\om,\a)\).
$$
This also implies that all the sets $\(A^*_j(\om)\)_{j=0}^{p-1}$ are mutually disjoint and the proof is complete.
\epf

Coming back to Theorem \ref{t1ms135}, if the map $\vp$ is global, that is, if $\vp:\C\to\oc$ is a meromorphic function,  we can now say more. It turns out, as we shall prove momentarily, that all the maps $F_\alpha$ extend holomorphically to the entire petals $A_j^*(\om)$.

\bthm\label{t1ms134}
If $f:\C\lra\oc$ is a meromorphic function of degree at least equal to 2 and if $\om\in \C$ is a simple rationally indifferent fixed point of $f$, then for every $j\in\{1, \ldots, p(\om)\}$, there exists a holomorphic function $F:A_j^*(\om)\lra \C$ such that
\begin{eqnarray}\label{1ms134}
F\circ f=F+1
\end{eqnarray}
and for all $\alpha\in (0,\pi)$ the function $F|_{S^j_a(\om,\a)}$ is univalent, and in addition, if $\a\in[\pi/2,\pi)$, then there exists $s>0$ such that
\beq\label{1_2015_10_07}
F(S^j_a(\om,\a))\supseteq \{z\in\C:\Re(z)>s\}.
\eeq
\ethm
\bpf
Fix $z\in A^*_j(w)$. By virtue of Theorem~\ref{t3ms134}, there exists $k\geq 0$ such that $f^k(z)\in S^j_a(\om,\a)$. We then define
\begin{eqnarray}\label{1ms134.a}
F(z):=F_\alpha(f^k(z))-k.
\end{eqnarray}
To see that this definition is independent of $k$, notice that Theorem \ref{t1ms135} implies that for every $j\geq0$,
\[
F(z)=F_\alpha(f^{k+j}(z))-(k+j).
\]
So,
\begin{eqnarray}\label{2ms134.a}
F(z):=\lim_{n\to\infty}F_\alpha(f^n(z))-n
\end{eqnarray}
is independent of $k$, and (\ref{1ms134}) holds. Formula \eqref{1_2015_10_07} is an immediate consequence of Theorem~\ref{t1ms135} (3).  Since formula (\ref{1ms134.a}) directly implies that the map $F$ is holomorphic, the proof is thus complete.
\epf

\sp Being more general, let $\om\in \C$ be a rationally indifferent periodic point of $f$. Then there exists a least integer $\ell\geq1$ such that $\om$ is a simple parabolic fixed point of $f^\ell$. We define $A^*(\om)$ to be the immediate basin of attraction to $\om$ under the iterates of $f^\ell$. We also define $A_j^*(\om)$, $j\in\{1,1,\ld,p(\om)\}$, to be the Leau--Fatou petals of $\om$ considered as a simple parabolic fixed point of $f^\ell$. We further define
\[
A(\om):=\bigcup_{j=0}^\infty f^{-j}(A^*(\om)), \  \  A^*_p(\om):=\bigcup_{j=0}^{\ell-1}A_*(f^j(\om))
\]  
and 
\[ 
A_p(\om):=\bigcup_{j=0}^{\ell-1}A(f^j(\om)).
\]
Now we are in the position to prove the following theorem.

\bthm\label{t2ms135}
If $f:\C\lra\oc$ is a meromorphic function of degree at least equal to 2 and $\om\in \C$ is a rationally indifferent periodic point of $f$, then each Leau--Fatou petal of $\om$, contains at least one singular point of $f^{-1}$. 
\ethm

\bpf
The first part of the proof parallels the corresponding part of the proof of Theorem \ref{t1ms123}, which concerns attracting periodic orbits. The argument generating a contradiction is, however, more difficult. So, assume first that $\om$ is a simple rationally indifferent fixed point of $f$. Fix 
$$
j\in\{1,1,\ld,p(\om)\}.
$$
Assume for a contradiction that $A_j^*(\om)$ contains no singular point of $f^{-1}$. We shall prove by induction that there exists an ascending sequence $(U_n)_{n=0}^\infty$ of open, connected, simply connected subsets of $A_j^*(\om)$ such that
\begin{itemize}
  \item [(a)] $\om\in \overline{U_0}$.
 
\, \item [(b)] For every $n\geq1$ there exists $f^{-1}_n:U_{n-1}\lra U_n$, a surjective holomorphic inverse branch of $f$.

\,  \item [(c)] For every $n\geq1$, $f^{-1}_{n+1}|_{U_{n-1}}=f^{-1}_n$.
\end{itemize}
Indeed, set
\[
U_0:=S^j_a(\om,\pi/2)\sbt A_j^*(\om).
\]
By virtue of formula (\ref{1ms135}) in Theorem~\ref{t1ms135}, we have that
\[
U_0\sbt f^{-1}_\om(U_0)\sbt A_j^*(\om),
\]
where, we recall, $f^{-1}_\om$ is the unique local holomorphic inverse branch of $f$ that sends $\om$ to $\om$. Therefore, setting $f_1^{-1}:=f_\om^{-1}|_{U_0}$, as $U_0$ is simply, we have that $U_1:=f_1^{-1}(U_0)\spt U_0$ is also simply connected, and the base case of the induction is complete.

For the inductive step, suppose that $n\geq1$ and 
$$
U_0\subseteq U_1\subseteq \cdots \subseteq U_n,
$$
an ascending sequence of open, connected, simply connected sets has been constructed such that properties (a) and (b) are satisfied for all $1\leq j\leq n$, whereas (c) is satisfied for all $1\leq j\leq n-1$. Since $A_j^*(\om)\cap f^{-1}(U_n)$ contains no singular points of $f^{-1}$, since $f:A_j^*(\om)\to A_j^*(\om)$ is a holomorphic surjective map, and since $U_n\supset U_{n-1}$, there exists $f^{-1}_{n+1}: U_n\to A_j^*(\om)$, a holomorphic inverse branch of $f$ such that
\begin{eqnarray}\label{1ms137}
f^{-1}_{n+1}|_{U_{n-1}}=f^{-1}_n.
\end{eqnarray}
So, setting
\begin{eqnarray}\label{2ms137}
U_{n+1}:= f^{-1}_{n+1}(U_{n-1}),
\end{eqnarray}
we only need to check that $U_{n+1}\supseteq U_n$. Invoking (\ref{2ms137}), (\ref{1ms137}) and property (b) above, we infer that
\[
U_{n+1}= f^{-1}_{n+1}(U_{n-1})\spt f^{-1}_n(U_{n-1})=U_n.
\]
The inductive construction is thus finished. Since, by Lemma~\ref{t3ms134A}, $\bu_{n=0}^\infty U_n=A_j^*(\om)$, conditions (a)--(c) along with ascendence of the sequence $(U_n)_{n=0}^\infty$, yield the existence of a holomorphic bijective function 
$$
f_*^{-1}:A_j^*(\om)\lra A_j^*(\om)
$$ 
such that 
$$
f_*^{-1}|_{U_{n-1}}=f_n^{-1}
$$ 
for every $n\ge 1$. Since also $f\circ f_*^{-1}=\Id_{A_j^*(\om)}$, it follows that the function $f:A_j^*(\om)\to A_j^*(\om)$ is also injective, and since we already know that it is surjective, we finally conclude that the holomorphic function 
$$
f:A_j^*(\om)\lra A_j^*(\om)
$$
is bijective. Iterating formula \eqref{1ms134} of Theorem~\ref{t1ms134} yields
\beq\label{2_2015_10_07}
F\circ f^n=F+n
\eeq
for all $n\ge 0$. Taking now any two distinct points $w, z\in A_j^*(\om)$ we find $n\ge 0$ so large that both points $f^n(w)$ and $f^n(z)$ belong to $S^j_a(\om,\pi/2)$. It then follows from bijectivity of the map $f^n:A_j^*(\om)\to A_j^*(\om)$, formula \eqref{2_2015_10_07}, and the middle placed assertion of Theorem~\ref{t1ms134} that $F(w)\ne F(z)$. In consequence, the holomorphic function 
$$
F:A_j^*(\om)\lra \C
$$
(produced in Theorem~\ref{t1ms134}) is univalent. Formula \eqref{2_2015_10_07} also yields
\[
F\circ f^{-n}_*=F-n.
\]
It therefore follows from formula \eqref{1_2015_10_07} of Theorem~\ref{t1ms134} that
$$
\begin{aligned}
F(A_j^*(\om))
&\supseteq F(U_n)=F(f^{-n}_*(U_0))=F(U_0) -n\\
&\supseteq\{z\in \C:\re(z)>s\}-n \\
&= \{z\in \C:\re(z)>s-n\}.
\end{aligned}
$$
Taking the union over all $n\geq1$, we therefore obtain that 
$$
F(A_j^*(\om))=\C.
$$
Since we already know that the map $F:A_j^*(\om)\to \C$ is univalent, we thus have that the inverse
$F^{-1}:\C \to A_j^*(w)$ is holomorphic. Since $A_j^*(\om)\cap \mathcal{J}(f)\neq \emptyset$ and $J(f)$ being perfect (see Theorem \ref{t_perfect_julia}), contains at least three distinct points, we therefore conclude from Picard's Theorem that the function $F^{-1}:\C \to A_j^*(\om)$ is constant. This contradiction finishes the proof in the case when $\om$ is a simple parabolic fixed point of $f$.

The general case follows from the preceeding one (simple parabolic fixed point) by applying it to $f^l$, where $l\ge 1$ is the least integer making $\om$ a simple parabolic fixed point of $f^l$. The proof is complete.\epf

\sp

\section{Fatou Components of (General) Meromorphic Functions III; \\ Simple Connectedness}\label{Fatou Components IV}

In this section we provide several sufficient conditions for a Fatou component (especially periodic) of of meromorphic function to be simply connected. We start with the following.

\bprop\label{p1ex1}
Let $f:\mathbb C\lra \bar{\mathbb C}$ be a meromorphic function. If $W$  is a periodic  connected  component of the Fatou set $F(f)$ of $f$ which is neither a Herman ring nor a Baker domain and if $W$ contains at most one singular value of $f$, then $W$ is simply connected.
\eprop

\fr {\sl Proof.}  Passing to an iterate of $f$ (even thought such iterate may fail to be meromorphic) we may assume  without loss of generality that $W$ is a forward invariant connected component of $F(f)$, i.e. $f(F(f)) \sbt F(f)$. The assertion of our theorem is immediate if $W$ is a Siegel disk because of Definition~\ref{d-Siegel MU}, Theorem~\ref{Siegel MU} and Theorem~\ref{Fatou Periodic Components} (3). It is true for all Siegel disks whatsoever.

Assume now  that $W$ is the basin of immediate attraction to an attracting fixed point. Denote it by $\xi$. We proceed somewhat similarly to the proof of Theorem~\ref{t1ms123}. Similarly as therein fix $U_0$, an open connected simply connected neighborhood of $\xi$, such that 
\beq\label{1ex1}
f(\ov{U_0})\sbt U_0
\eeq
and $\partial{U}_0$ intersects no forward orbit of any singular value of $f$. We shall construct inductively a sequence $(U_n)_{n=0}^\infty$ of  open subsets of $\mathbb C$ with the following properties:

\,

\begin{itemize}
\item[\rm($a_n$)]  $ U_{n+1}\supset U_n$,

\,

\item[\rm($b_n$)]  $U_n\sbt A^{*}( \xi) =W$,

\,

\item[\rm($c_n$)]  $U_{n+1}$ is a  connected  component of $f^{-1}(U_n)$,

\,

\item[\rm($d_n$)]  $U_n$  is  simply connected.
\end{itemize}

\,

\fr Indeed, suppose that the open sets $U_0, U_1, \ldots ,U_n$ have been defined such that all conditions ($a_k$)--($d_k$) hold for all $0\leq k \leq n-1$ along with $(b_n)$ and $(d_n)$. Look at $f^{-1}(U_n)$. By $(a_n)$ and $(c_{n-1})$,
$$
f^{-1}(U_n)\supset f^{-1}(U_{n-1})\supset  U_{n}.
$$
Let $U_{n+1}$  be the connected component of $f^{-1}(U_n)$ containing $U_n$. Then immediately, the conditions $(a_n)$ and $(c_n)$  hold. Conditions $(b_n)$ and $(d_n)$ hold  because of our inductive assumption. By ($a_n$) and ($b_n$), we have that 
$$
U_{n+1}\cap A^{*}(\xi)\neq \es.
$$
Since $U_{n+1}\sbt F(f)$ and $A^{*}(\xi)$ is a connected  component of $F(f)$, this implies that $U_{n+1}\sbt A^{*}(\xi)$, meaning that $b_{n+1}$  holds. Item $d_{n+1}$ follows directly from Corollary~\ref{General Riemann Hurwitz Simply Connected_2} and Corollary~\ref{General Riemann Hurwitz Simply Connected_3}. Let 
$$
U_\infty := \bigcup_{n=0}^\infty U_n.
$$
Then  $U_\infty$ is open and also connected  and simply connected  as being the union of an ascending sequence of open  connected simply connected sets. Thus, in order to finish the proof in the  current case, it is suffices, to show  that
$$
U_\infty=W.
$$
By conditions $(b_n)$, we  have that 
$$
U_\infty\sbt W.
$$
Seeking a contradiction, suppose that 
$$ 
U_\infty\varsubsetneq W.
$$
Since both sets, $W$ and $U_\infty$ are connected, this gives that $\partial{ }_W U_\infty\neq \es$. So, there exists a point
\beq\label{1ex2}
 z \in \partial{ }_W U_\infty.
 \eeq
By ($c_n)$  we have that $f(U_{n+1})=U_n$, whence $f(U_\infty)=U_\infty$.
Therefore,  
$$
f(\ov{U}_\infty)\sbt \ov{f(U_\infty)}=\ov{U}_\infty.
$$
Hence, for every integer $k \geq 0$,
\beq\label{2ex2}
f^k(\partial{ }_W U_\infty)\sbt  \ov{U}_\infty\cap W.
\eeq
Since $ z \in W=A^{*}(\xi)$, there exists an integer $ l \geq 0$ such that $ f^l(z) \in U_0$. In particular  $f^l(z)\in U_\infty$.
Let $ n \geq 0$ be the least integer with this property. It follows from (\ref{1ex2}) that $n \geq 1$ and, so, with the help of (\ref{2ex2}),
$f^{n-1}(z) \in \partial{ }_W U_\infty$. So, replacing $z$ by  $f^{n-1}(z)$, we may assume that
\beq\label{3ex2}
z \in \partial{ }_W U_\infty \quad \text{and} \quad  f(z) \in U_\infty.
\eeq
Hence,  there exists $n \geq 0$ such that
\beq\label{4ex2}
 f(z) \in U_n.
\eeq
Since the map $f$ is continuous, there thus exists an  $\varepsilon >0$ such that
\beq\label{5ex2}
 B(z, \varepsilon)  \sbt f^{-1}(U_n)\sbt F(f).
\eeq
Since the sequence $(U_n)^\infty_{n=0}$ is ascending, it follows from the first part of (\ref{3ex2}),
that 
\beq\label{1ex3}
B(z, \varepsilon)\cap U_{s+1}\neq \es
\eeq
for some $ s \geq n$. Using again the fact that the sequence $(U_n)^\infty_{n=0}$ is ascending,
it follows from (\ref{5ex2}), that $B(z, \varepsilon)  \sbt f^{-1}(U_s)$ . Along with ($c_s$), this yields
$$
B(z, \varepsilon)\cup U_{s+1}  \sbt f^{-1}(U_s).
$$
Since both sets, $B(z, \varepsilon)$ and  $U_{s+1}$, are connected, it follows from (\ref{1ex3})  that the set $B(z, \varepsilon)\cup U_{s+1}$ is connected.
Combining this,(\ref{5ex2}), and ($c_s$) again, we conclude that 
$$
B(z, \varepsilon)\sbt U_{s+1}.
$$
Hence, $z\in U_\infty$ contrary to the first part of (\ref{3ex2}) and openness of $U_\infty $ in $W$. We are done in this case.

Because of Theorem~\ref{Fatou Periodic Components} there is only left the case of $\xi$  being a rationally indifferent fixed point of $f$. Passing to yet larger iterate, we may assume without loss of generality that $\xi$ is simple, i.e. $f'(\xi)=1$. Then 
$$
W=A^{*}_j(\xi)
$$ 
with some $ j \in \{1,2, \ldots, p(\xi)\}$. Recall that $S^j_a(\xi, \alpha)$, $\alpha \in (0, \pi)$ is the sector defined by the formula \eqref{1_2017_29A}. Taking 
$$
U_0:=S^j_a(\xi, \alpha)
$$
and making use of Lemma~\ref{l1_2015_09_10} and Theorem~\ref{t3ms134}, we may proceed in exactly the same way as in the case of an attracting
fixed point, to conclude that the set $W=A^{*}_j(\xi)$  is simply connected. The proof of Proposition~\ref{p1ex1} is complete. 
\qed

\sp As a fairly immediate consequence of Proposition~\ref{p1ex1} we get the following.

 \bcor\label{c1ex3} Let $f: \mathbb C  \lra \oc$ be a meromorphic function. If $W$ is a periodic connected component of the Fatou set $F(f)$ of f which is not a Herman ring nor a Baker domain and each connected component of the set 
$$
W_\infty:= \bigcup_{n=0}^\infty f^{-n}(W)
$$
contains at most one singular value of f, then each of these components is simply connected.
\ecor

\bpf For each connected component $\Gamma$  of $W_\infty$ let $N_\Gamma \geq 0$ be the least integer such that
$$
f^{N_\Gamma}(\Gamma) \sbt W.
$$
We proceed by induction with respect to $N_\Gamma$. If $N_\Gamma=0$, then $ \Gamma=W$, and $\Gamma$ is simply connected because of Proposition~\ref{p1ex1}. So,
fix $n \geq 0$ and suppose that the assertion of Corollary~\ref{c1ex3} holds for all components G of $W_\infty$  for which $N_G \leq n$. Let $\Gamma$
 be a connected component of $W_\infty$ for which $N_\Gamma=n+1$ . Then $ f(\Gamma) \sbt W_\infty$ and $N_{f(\Gamma)}=n$. So, $f(\Gamma)$ is simply connected. Since it contains at most one singular value of $f$ and since $\Gamma$  is a connected component of $f^{-1}(f(\Gamma))$, it follows directly from Corollary~\ref{General Riemann Hurwitz Simply Connected_2} and Corollary\ref{General Riemann Hurwitz Simply Connected_3}that $\Gamma$ is simply connected, and the proof of Corollary~\ref{c1ex3} is complete. 
\epf

\section[Fatou Components of (General) Meromorphic Functions V; Baker] {Fatou Components of (General) Meromorphic Functions V; Baker Domains}

In this section we deal with Baker domains. Indeed, except for Baker domains all the periodic components of Fatou sets described in Theorem~\ref{Fatou Periodic Components} appear already in the realm of the dynamics of rational functions. The first example of an entire transcendental function having a Baker domain was already given by Fatou in \cite{F}, who considered the entire function function $f:\C\lra\C$ given by the formula
$$
f(z)=z+1+e^{-z},
$$
and proved that 
$$
\lim_{n\to\infty}f^n(z)=+\infty
$$ 
for every $z\in\C$ with $\re(z)>0$. This implies that the right half-plane is contained in an invariant Baker domain. For ergodic and fractal properties of the function $f$ see \cite{KU4}. An
example of a Baker domain of higher period was given in \cite{BKL3},
where it was shown that the function 
$$
f(z):=z^{-1}-e^z
$$ 
has a cycle $U_0, U_1$ of Baker domains such that 
$$
f^{2n}_{|U_0} \to\infty 
\  \  \  {\rm and}  \  \  \
f^{2n}_{|U_1} \to 0
$$ 
as $ n \to \infty$.

\sp Now we explore some  general properties of Baker domains. Let $\{U_0, U_1, \ldots, U_{p-1}\}$ be a periodic cycle of Baker domains, and denote by $\xi_j$ the limit corresponding to $U_j$, that is  
$$
\lim_{n\to\infty}f^{np}(z)=\xi_j
$$ 
for all $z \in U_j$. Clearly 
$$
f(\xi_j)=\xi_{j+1}
$$ 
if $z_{j}\neq \infty$, where
$\xi_p:=\xi_0$ and there exists at least one $j\in \{0, 1,
\ldots, p-1\}$ such that  $\xi_j=\infty$, and for all $j\in \{0, 1,
\ldots, p-1\}$ there exists $l=l(j) \in \{0, 1, \ldots, p-1\}$ such
that 
$$
f^l(\xi_j)=\infty.
$$
We claim that each set $U_j$ contains a (topological) curve $\g_j$ converging to $\xi_j$ such that
$$
f^p(\g_j)\sbt \g_j \  \  \text{ and } \  \ f^p(z)\to \xi_j
$$ 
as $ z \to \xi_j$ in $\g_j$.  To see this we fix $w\in U_0$ and then a closed (compact) topological arc $\sigma
 \sbt  U_0$ that joins $w_0$ and $f^p(w)$. We define
$$
\g_0:=\bigcup_{n=0}^\infty f^{np}(\sigma)\  \  \  {\rm and } \  \  \ \g_j:=f^j(\g_0)
$$
for $j \in \{1, 2, \ldots, p-1\}$. Then the curves $\g_j$ have the desired properties. Moreover, 
$$
\lim_{\g_j\ni z\to \xi_j} f(z)= \xi_{j+1}.
$$ 
We deduce that if $\xi_j=\infty$, then $\xi_{j+1}$ is an asymptotic value
of $f$, the asymptotic path being contained in $U_j$. We have thus proved the following theorem.
 
\sp
\bthm\label{th13in{Be}} Assume that $f:\C\lra\oc$ is a meromorphic function having 
$$
\{U_0, U_1, \ldots, U_{p-1}\},
$$
a periodic cycle of Baker domains. Denote by $\xi_j$, $j=0, 1, \ld,p-1$, the limit corresponding to $U_j$, and define $\xi_p:=\xi_0$. Then 
$$
\xi_j \in
\bigcup_{n=0}^{p-1}f^{-n}(\infty)
$$ 
for all $ j \in \{0, 1, \ldots,
p-1\}$, and  $\xi_j=\infty$ for at least  one $j\in \{0, 1, \ldots,
p-1\}$. If $\xi_j=\infty$, then $\xi_{j+1}$ is an asymptotic value of
$f$. 
\ethm

The two following consequences of this theorem are immediate.

\bcor\label{c1in{Be}}
If a meromorphic function  $f:\C\lra\oc$ has a periodic cycle $\{U_0, U_1, \ldots, U_{p-1}\}$ of Baker domains such 
$$
f^n_{| U_0}\lra \infty \  \  {\rm as } \  \ n\to+\infty,
$$
then $\infty$ is an asymptotic value of $f$. In particular, this is the case if $f$ has a Baker domain which is $f$--invariant .
\ecor

\bcor\label{c2in{Be}}
If a meromorphic function $f:\C\lra\oc$ has a cycle $\{U_0, U_1,\ldots, U_{p-1}\}$ of Baker domains such 
$$
f^n_{| U_0}\not\!\!\lra \infty \  \  {\rm as } \  \ n\to+\infty,
$$
then $f$ has a finite asymptotic value. 
\ecor

\section{Fatou Components of Meromorphic Functions IV; Class ${\mathcal B}$ and ${\mathcal S}$}\label{FC B+S}

\sp It is fairly common in transcendental dynamics to consider the following two classes of transcendental meromorphic functions. Firstly: 
$$
{\mathcal S}:=\big\{f:{\mathbb C}\lra \oc: f \ {\rm is \ meromorphic \ and } \
\Sing(f^{-1})  \mbox{ is finite} \big\}
$$ 
This collection of functions is usually refered to as Speiser class ${\mathcal S}$. Secondly:
$$
{\mathcal B}:=\big\{f:{\mathbb C}\lra \ov{\mathbb C}: f \ {\rm is \ meromorphic \ and } \ \Sing(f^{-1})  \mbox{ is bounded} \big\}.
$$
\index{(S)}{$\mathcal S$} \index{(S)}{$\mathcalB$} \! This collection of functions is usually called Eremenko--Lyubich class ${\mathcal B}$. It was introduced in \cite{ELclasB} and widely used since then. Of course
$$
{\mathcal S}\sbt {\mathcal B}.
$$
As an immediate consequence of Theorem~\ref{t1ms123} and Theorem~\ref{t2ms135}, we get the following.

\bthm\label{t120200404}
Any meromorphic function in Speiser class $\mathcal S$ has only finitely many attracting and rationally indifferent periodic points.
\ethm

The fundamental theorem about Speiser Class $\cS$ is the following.

\bthm\label{Sullivan} 
No functions in Speiser class $\mathcal S$ have wandering domains. 
\ethm

This theorem has been conjectured by P. Fatou in \cite{F0} and \cite{F1}for rational functions and was proved in this form by D. Sullivan in \cite{Su0}. In fact it holds for all meromorphic functions in Speiser class $\cS$ and was proved by I.N. Baker, J. Kotus, and Y. L\"u in \cite{BKL4}. 
See this paper and references therein for more historical and bibliographical information. Examples of wandering components for analytic self--maps of $\mathbb C^*$ were provided by J. Kotus in \cite{K2} and by I.N.~Baker, J.~Kotus and Y. L\"u  in \cite{BKL2}. We provide in Appendix~B its original proof stemming from \cite{BKL4}. 

\sp For every integer $n\ge 1$, we introduce an additional auxiliary class of functions. For a meromorphic function $f:\C\lra\oc$, let
$$
S_n(f):=\bigcup_{j=0}^{n-1}f^j(\Sing(f^{-1})\sms f^{-j}(\infty)).
$$
We define 
\beq\label{B_n} 
\cB_n:=\{f:\C\to\oc: f \mbox{ is meromorphic transcendental with}\, S_n(f) \,
\mbox{bounded}\}.
\eeq
Of course $S_1(f)=\Sing(f^{-1})$ and $\cB_1=\cB$ is class $\cB$ introduced above. Also
$$
\cS\sbt\bi_{n=1}^\infty\cB_n
$$
For every $R>0$ we define
$$
B_\infty(R):=\{z\in \ov{\mathbb C}:|z| >R\}
\  \ {\rm and} \ \
B_\infty^*(R) :=\{z\in {\mathbb C} |z|>R\}.
$$
\index{(S)}{$\ov B_R$}\index{(S)}{$\ov B_R^*$} Aiming to prove Theorem~\ref{baker-domain} about non--existence of Baker domains for functions from Speiser class $\cS$, wee need first the following three  auxiliary results.

\blem\label{l2.1in {RS}} 
If $f \in \cB_n$ and $S_n(f) \sbt B(0, R)$, then each connected component of $f^{-n}(B_\infty(R))$ in $\mathbb C$ is simply connected. 
\elem

\bpf Let $V$ be a connected component of $f^{-n}(B_R(\infty))$. For every $z\in V$ let 
$f_z^{-n}$ be the germ of holomorphic inverse branches of $f^n$ defined on some neighborhood of $f^n(z)$ and sending $f^n(z)$ to $z$. This germ is well-defined since $S_n(f)\cap B_\infty)(R)=\es$. Let $h_{z,z_k}$, $z\in V$, $k\in\Z$, with $e^{z_k}=f^n(z)$, be the germs equal to $f_z^{-n}\circ\exp$, defined on sifficiently small neighborhoods of 
$$
z_k\in H_R:=\{w\in\C: \Re(w) > \log R\}.
$$
These germs can be analytically continued along any arc in $H_R$. Fix one such germ. Since $H_R$ is an open connected simply connected subset of $\C$, we thus conclude that by the Monodromy Theorem (Theorem~\ref{monodromy}), there is $h:H_R\to\C$, an analytic continuation of this germ, and 
\beq\label{1MU20150825}
h(H_R)=V\sms f^{-n}(\infty).
\eeq 
Note that
\beq\label{1_2017_07_24}
f^n\circ h=\exp 
\eeq
on $H_R$. 

There are two cases to consider. First, if the function $h:H_R\to V$ is 1-to-1, then $V\sms f^{-n}(\infty)=h(H_R)$ is open simply connected as a homeomorphic image of an open simply connected set. Since, in addition, both sets $V$ and $h(H_R)$ are open, we thus conclude that $V\cap f^{-n}(\infty)\sbt h(H_R)$. But as $f^{-n}(\infty)\cap h(H_R)=\es$, we further conclude $V\cap f^{-n}(\infty)=\es$. Thus, finally $V=h(H_R)$ and we are done in this case. 

\sp So, suppose that $h:H_R\lra V$ is not 1--to--1. This means that there are two points $z, w\in H_R$ with $z\ne w$ such that 
$$
h(z)=h(w).
$$
But then $f^n\circ h(z)=f^n\circ h(w)$, meaning, because of \eqref{1_2017_07_24}, that $e^z=e^w$. Hence, there exists $m\in\Z$ such that 
$$
w=z+2\pi im.
$$
Exchanging $z$ with $w$ if needed, we may assume without loss of generality that $m>0$, and furthermore that $m>0$ is minimal with this property. It now follows from the Open Mapping Theorem that for every $\e>0$ there exists $\d>0$ such that for every $\xi\in B(z,\d)$ there exists $\xi'\in B(w,\e)$ such that 
$$
h(\xi')=h(\xi).
$$ 
But then by the same token as above there exists $k\in \Z\sms\{0\}$ such that $\xi'=\xi+2\pi ik$. Since $B(z,\d)$ is uncountable while $\Z$ is countable, there thus exist $q\in\Z\sms\{0\}$ (assuming that $\e>0$ and $\d>0$ are small enough) and a sequence $\xi_j\in B(z,\d)\sms\{z\}$, $j\ge 1$, converging to $z$, such that the sequence $\xi_j+2\pi iq$, $j\ge 1$, converges to $w$ and  $h(\xi_j+2\pi iq)=h(\xi_j)$ for every $j\ge 1$. It follows that 
$$
h(\xi+2\pi iq)=h(\xi)
$$
for all $\xi\in H_R$. Also, letting $j\to\infty$ and looking at the sequences $\xi_j$ and  $\xi_j+2\pi iq$, $j\ge 1$, we conclude that $q=m>0$. So, $h:H_R\to V$ is periodic with period $2\pi im$ and (by minimality of $m$) is univalent on each horizontal strip of with $2\pi m$. Therefore, we can represent $h:H_R\to V$ in the form
$$
h(s)=\phi(e^{s/m}),
$$
where $\phi:B_\infty^*(R^{1/m})\to V\sbt\C$ is a holomorphic univalent map. Therefore, the Laurent series expansion of $\phi$ in $B_\infty^*(R^{1/m})$ takes on the form
$$
\phi(s)=a_1 s +a_0+a_{1}s^{-1}+a_{2}s^{-2}+\ldots \,.
$$
By \eqref{1_2017_07_24} we have that $f^n(\phi(e^{s/m}))=f^n(h(s))=e^s$. So, writing $z=\phi(e^{s/m})\in V$, we have $e^{s/m}=\phi^{-1}(z)$, and further
$$
f^n(z)=e^s=(\phi^{-1}(z))^m
$$
for all $z\in \phi\(B_\infty^*(R^{1/m})\)\sbt V$. Now, if $a_1\neq 0$, then $\phi\(B_\infty^*(R^{1/m})\)$ contains a deleted neighbourhood of $\infty$ and
\beq\label{2.1in{RS}}
\lim_{z\to\infty}\frac{f^n(z)}{a_1^{-m} z^m}=1.
\eeq
But this implies that $\infty$ is a pole of $f^n$ whereas it in fact is its essential singularity. Thus $a_1=0$ and so, invoking \eqref{1MU20150825}, we get
$$
(V\sms f^{-n}(\infty))\cup \{a_0\}
=\{a_0\}\cup \phi\(B_\infty^*(R^{1/m})\)
=\phi\(B_\infty(R^{1/m})\)
$$
is an open simply connected subset of $\mathbb C$. Thus $a_0\in V$ and, consequently, $V\sms (f^{-n}(\infty)\sms\{a_0\})$ is open and simply connected. Hence ($f^{-n}(\infty)\sms\{a_0\}$ is a countable set), $V$ itself is 
simply connected, and the proof is complete.
\endpf

\sp 

\blem\label{l2.2in{RS}} If $f:\C\lra\oc$ is a transcendental meromorphic
function, then there exists $R_f$ such that, if $R >R_f$, $S_n(f) \sbt
B(0,R)$ and $|z|, |f^n(z)| > R^2,$  then
$$
|(f^n)'(z)|> \frac{|f^n(z)|\log |f^n(z)|}{16 \pi |z|}.
$$
\elem

\bpf First fix $\xi\in J(f)$, a periodic point of $f$ and take
$R_f\ge 0$  so large that $|f^n(\xi)|< R_f$ for each $ n\in \mathbb N$.
Now take any $R > R_f$, any $z\in\C$, and suppose that 
$$
S_n(f) \sbt B(0,R)  \  \  \text{ and }  \   \ |z|, |f^n(z)| >R^2.
$$
Let $V$ be a connected component of $f^{-n}(B_R(\infty))$. For every $z\in V$ let 
$f_z^{-n}$ be the germ of holomorphic inverse branches of $f^n$ defining on some neighborhood of $f^n(z)$ and sending $f^n(z)$ to $z$. This germ is well-defined since $S_n(f)\cap B_\infty(R)=\es$.
Since $\xi\notin V$ and since, by  Lemma~\ref{l2.1in {RS}}, the open set $V$ is simply connected, there exists $\log_\xi:V\to\C$, a holomorphic branch of the logarithm of the function $V\ni z\mapsto z-\xi$  (of course there are infinitely coutably many of such functions). With,
$$
H_R=\{w\in\C: \re(w) > \log R\},
$$
being the same as in the previous proof, the colection
$$
\Ga:=\(\log_\xi\circ f_{w_k}^{-n}\circ\exp\)_{w\in H_R},
$$
where $w_k$ ranges over the set $V\cap f^{-n}(e^w)$, forms a germ of holomorphic functions which can be analitically continued along any arc in $H_R$.  Fix one such germ. Since $H_R$ is an open connected simply connected subset of $\C$, we thus conclude that by the Monodromy Theorem (Theorem~\ref{monodromy}), there exists
$$
G:H_R\lra\C,
$$
an analytic continuation of this germ. It immediately follows from the definition of $\Ga$ and $G$ that 
$G(H_R)$ does not contain any disc of radius greater than $\pi$. Thus, by Bloch's Theorem (see \cite{C}), we get that
$$
|G'(w)|\leq \frac{\pi}{B(\re(w) - \log R)},
$$
for $t\in H$, where $B>0$ is the Bloch's constant. Developing the derivative of $G$, we thus get that
$$
\left|\frac{(f_z^{-n})'(e^w) e^w}{f_z^{-n}(e^w)-\xi}\right| \leq \frac{\pi}{B( \re(w)
-\log R)},
$$ 
for any $w\in H_R$, and any $z\in V\cap f^{-n}(e^w)$ such that $f^n(z)=e^w$. Hence 
\beq\label{2.2in{RS}}
\left|\frac{f^n(z)}{(z-\xi)(f^n)'(z)}\right|=
\left|\frac{(f_z^{-n})'(f^n(z))f^n(z)}{z-\xi}\right|\leq \frac{\pi}{B( \log
|f^n(z)|  -\log R)}.
\eeq
The lemma follows by using $|z-c|\leq |z|+ |c| < 2 |z|$, $
\log|f^n(z)|> 2 \log R$ and $B>1/4$. 
\endpf

\sp

\bthm\label{thAin{RS}} 
If $f \in \cB_n$, then there is no connected component $U$ of the Fatou set $F(f)$ such that $f|_U^{mn}\to \infty$ uniformly on compact subsets of $U$ as $m\to \infty$. 
\ethm

\bpf If $f \in \cB_n$, then there exists $R> \max(e^{16},
R_f)$ with $S_n(f) \sbt  B(0, R)$.  If $F(f)$ has a  component $U$
such that $f^{mn}(z)\to \infty$ uniformly on compact subsets of $U$ as $m \to \infty$, then  there exists $p \in  \mathbb N$, $ w \in F(f)$ and  $r >0$ such that
\beq\label{f1_2015_08_29}
f^{pn}(U)\spt \ov{B(w,r)} \  \  \text{ and }  \  \  |f^{mn}(z)| >R^2
\eeq
for each  $z\in \ov B(w,r)$ and all integers $m\ge 0$.
Now  let $V_n$ be the connected component of $f^{-n}(B_R)$ containing 
$U_m:=f^{mn}(B(w,r))$. Taking $\xi$ to be the same  periodic
point as in the proof of Lemma~\ref{l2.2in{RS}}. We know from Lemma~\ref{l2.1in {RS}} that $V_m$ is simply connected. So, as in the proof of Lemma~\ref{l2.2in{RS}}, there exists $\log_\xi:V_n\to\C$, a holomorphic branch of the logarithm of the function $V_n\ni z\mapsto z-\xi$. Next, put 
$$
T_m:=\log_\xi(U_m)
$$ 
and
$$
F_m(t):=\log_\xi\circ f^n\circ(\xi+\exp):T_m\to\C,
$$
so that $T_{m+1}=F_m(T_m)$. It now  follows from (\ref{2.2in{RS}}) and from \eqref{f1_2015_08_29} that if $w \in T_m$, then for $z=\xi+e^w\in U_m$
$$
\begin{aligned} 
|F_m'(w)| 
&=\left| \frac{(f^n)'(\xi+e^w)e^w}{f^n(\xi+e^w)-\xi}\right| 
=\left| \frac{(f^n)'(z)(z-c)}{f^n(z)-\xi}\right|
 \geq \left| \frac{(f^n)'(z)(z-c)}{2f^n(z)}\right|\\
 & \geq \frac{B}{2\pi}( \log|f^n(z)|-\log R)
 \geq  \frac {B\log R}{2 \pi} \geq 2,
\end{aligned}
$$
Hence, using the Chain Rule we obtain
$$
|(F_m\circ \ldots \circ F_1)'(w)| \geq 2^m
$$
for all $w\in T_1$ and all $m\ge 1$. 
Thus, by Bloch's Theorem, each set $T_m$ contains a disc of radius $\const\cdot 2^m$ diverging to infinity as $m \to \infty$. This, is however a contradiction
since $T_m \sbt\log_\xi(V_n)$ which is contained in a horizontal strip of height $2\pi$, thus contains no disk of radius greater than $\pi$. We are done.
\endpf

\sp

Since $\cS\sbt\bi_{n=1}^\infty\cB_n$, as an immediate consequence of Theorem~\ref{thAin{RS}}, we get the following remarkable result. 

\bthm\label{baker-domain} 
No functions in Speiser class $\mathcal S$ have Baker domains. 
\ethm

\sp As an immediate consequence of Corollary~\ref{c1ex3} and Theorem~\ref{baker-domain}
we get the following. 

\bcor\label{c1ex4} 
Let $f:\mathbb C \lra \hat{\mathbb C}$ be a meromorphic function in Speiser class $\mathcal S$. If $W$ is a periodic connected component of the Fatou set $F(f)$ of $f$ which is not a Herman ring and each connected component of
$$
W_\infty:=\bigcup_{n=0}^{\infty}f^{-n}(W)
$$
contains at most one singular value of f, then each of these components is simply connected.
\ecor

Since the Julia set is, by definition, the complement of the Fatou set, as an immediate consequence of this corollary and Theorem~\ref{herman}, we get the following.

\bthm\label{t2ex4}
Let $f:\mathbb C \lra \hat{\mathbb C}$ be an even elliptic function in Speiser class $\mathcal S$. If each connected component of the Fatou set $F(f)$ of $f$ contains at most one critical value of $f$, then each connected component of $F(f)$ is simply connected, and (therefore) the Julia set $J(f)$ of $f$ is connected.
\ethm

In subsequent sections we will be dealing with the $\omega$-limit set
and  $\alpha$-limit set of points. For a given $z \in \mathbb C$ the
$\omega$--limit set\index{(N)}{$\omega$-limit set}
$\om(z)$\index{(S)}{$\om(z)$} is  the set of accumulation points of
$O^+(z)$ in $\ov{\mathbb C}$. Analogously the $\a$-limit
set\index{(N)}{$\a$-limit set} $\a(z)$\index{(S)}{$\a(z)$} is the
set of accumulation points of $O^-(z)$ in $\ov{\mathbb C}$. We end this section with a concept which will play an important role in part of the book devoted to elliptic functions, particularly to Hausdorff dimension of their Julia sets.
\beq\label{3_2017_09_25}\index{(S)}{$I_{\infty}(f)$}
I_{\infty}(f):=\lt\{z\in {\mathbb C} :z\in \bu_{n\ge 0}f^{-n}(\infty)
\text{ or } \lim_{n\to \infty} f^n(z)= \infty \rt\},
\eeq
The set $I_\infty(f)$ is called the set of points escaping to $\infty$ under iterates of $f$, or just the escaping set of $f$\index{(N)}{the set of points escaping to $\infty$}\index{(N)}{the escaping set}. As an immediate consequence of Theorem~\ref{Sullivan}, Theorem~\ref{baker-domain}, and Theorem~\ref{Fatou Periodic Components}, we get the following. 

\bthm\label{t2_2017_09_25}
If $f\in\mathcal S$, then set $I_\infty(f)\sbt J(f)$. 
\ethm

For a further background concerning the topological dynamics of transcendental
meromorphic functions  the reader is referred  to \cite{BKL1},
\cite{BKL2},  \cite{BKL3}, \cite{BKL4} and \cite{Be}.

\sp \section[Local and Asymptotic Behavior of Meromorphic Functions; Part II] {Local and Asymptotic Behavior of (General) Meromorphic Functions Around Rationally Indifferent Periodic Points; \\ Part II: 
Fatou Flower Theorem and Fundamental Domains}\label{local-parabolic II} 

In this section $f:{\mathbb C}\lra \oc$ is an entirely
arbitrary (global) meromorphic function; in particular any elliptic one is
allowed. The results we formulate here (with proofs) are of fairly
classical nature (except those concerning semi--conformal measures) and
are scattered widely in the literature; you may find them for example in
\cite{ADU}, \cite{DU3}, and \cite{DU4}. 

Let $\Om(f)$\index{(S)}{$\Om(f)$}  denote the set of rationally
indifferent periodic points of $f$. Throughout this section, unless otherwise stated, $\om$\index{(S)}{$\om\in\C$} is assumed to be a simple parabolic fixed
point  \index{(N)}{simple parabolic fixed point} of $f$, that is, we recall, 
$$
f(\om)=\om  \  \  {\rm and }  \  \  f'(\om)=1.
$$
We denote the set of all simple parabolic fixed
points of $f$ by $\Om_0(f)$. \index{(S)}{$\Om_0(f)$}
On a sufficiently small neighborhood $V$  of $\om$ a unique
holomorphic inverse branch $f^{-1}_\om: V \to \mathbb  C$, which
sends $\om $ to $\om $, is well--defined. Thus, all the results of
Section~\ref{local-parabolic} apply with both $\varphi=f^{-1}_\om$ or
$\varphi=f$. Fix $\a \in (0, \pi)$ ($\kappa=1/2$), and recall that 
$$ 
S^1_a(\om, \a),\ldots,S^{p(\om)}_a(\om, \a)
$$
are the corresponding attracting sectors\index{(N)}{attracting sectors} for $f$ defined in (\ref{1_2017_29A}) and 
$$
S^1_r(\om, \a),\ldots, S^{p(\om)}_r(\om, \a)
$$
are the corresponding repelling sectors\index{(N)}{repelling sectors} defined in (\ref{1_2017_29B}). Equivalently, these are the
attractive sectors for $f^{-1}_\om$.  Since the family of iterates
of $f$ is not normal on any  neighbourhood of any point in the Julia
set, as fairly immediate consequence of Proposition~\ref{p1j108}, we get the
following celebrated classical result. 

\

\bthm\lab{lFf}{\rm (Fatou's Flower Theorem).}\index{(N)}{Fatou's
Flowers Theorem}
Let $f:\C\lra\oc$ be  a meromorphic function. If $\om\in\Om(f)$, i. e. if $\om$ is a rationally indifferent periodic points of $f$, then
for every $\a\in(0,\pi)$ there exists $\th_\a(\om)\in(0,\th_f)$ such that
$$
J(f)\cap B(\om,\th_\a(\om)) 
\sbt S^1_r(\om,\a)\cup\ld\cup S^{p(\om)}_r(\om,\a).
$$
\ethm

\bpf
Passing to a suffficiently high iterate of $f$ we may assume without loss of generality that $\om$ is a simple parabolic fixed point of $f$. 
Seeking contradiction suppose that there exists a sequence
$(z_n)_{n=1}^\infty$ of points in $J(f)\sms\cup_{j=1}^{p(\om)}S_r^j(\om,\a)$ such that  
$$
\lim_{n\to+\infty}|z_n-\omega|=0. 
$$
Then, fixing $\b\in(0,\a)$, for every $n\ge 1$ large
enough there would exist $i\in\{1,2,\ld,p(\om)\}$ such that $z_n\in
S^i_a(\om,\pi-\b)$. Fix one such $n$ and $i$. Then, by 
Proposition~\ref{p1j108}, $S^i_a(\om,\pi-\b)$ is a subset of the Fatou
set of $f$. Thus $z_n\notin J(f)$. This contradiction
finishes the proof.
\epf

\sp Since the Julia set $J(f)$ is fully invariant ($f^{-1}(J(f))\sbt
J(f)$ and $f(J(f)\sms \{\infty\})=J(f)$), we infer from Fatou's
Flower Theorem (Theorem~\ref{lFf}), Proposition~\ref{p1j108}, and
\eqref{120120728} (actuallly the derivative of $\phi_0$) that for every $\om\in\Om_0(f)$ there exists
$\th(\om)\in\(0,\th_{\pi/4}(\om)\)$ such that for every
$0<R\le\th(\om)$, we have 
\beq\label{1070506}
f_\om^{-1}(J(f)\cap B(\om,R)))\sbt J(f)\cap B(\om,R).
\eeq
Thus all iterates 
\beq\label{1070506B}
f_\om^{-n}:J(f)\cap B(\om,R)\lra J(f)\cap B(\om,R),
\eeq
$n=0,1,2,\ld$ are well defined. From
Theorem~\ref{lFf} and (\ref{1070506}) we obtain for all $\a>0$ and all
$j\in\{1,2,\ld,p(\om)\}$, that
\beq\lab{dp4.4}
 f_\om^{-1}\(J(f)\cap S^j_r(\om,\a)\)\sbt J(f)\cap S^j_r(\om,\a).
\eeq
By Proposition~\ref{p1j108}
\beq\lab{dp4.4D}
 f_\om^{-1}\(S^j_r(\om,\a)\)\sbt S^j_r(\om,\a),
\eeq
and therefore all iterates 
\beq\lab{dp4.4E}
 f_\om^{-n}:S^j_r(\om,\a)\lra S^j_r(\om,\a),
\eeq
$n\ge 0$, are well defined. Furthermore, it follows from Proposition~\ref{p1j108} that 
\beq\lab{dp4.6}
\lim_{n \to \infty}f_\om^{-n}(z)=\om
\eeq
uniformly on $z\in S^j_r(\om,\a)$, while from Lemma~\ref{dp4.3A} that
\beq\lab{dp4.3}
|f_\om^{-1}(z)-\om| < |z-\om|  \  \  \and  \  \ |(f_\om^{-1})'(z)|< 1
\eeq
for every $z\in \(S^j_r(\om,\a)\cap B(\om,R_\a(\om)\)\sms\{\om\}$ if $\a<\frac{\pi}{2p(\om)}$.
Put
\beq\lab{dp4.5}
\th=\th(f) := \min\big\{\min\{\th(\om),R_{\pi/4}(\om)\}:\om\in \Om(f)\big\}.
\eeq

\fr As an immediate consequence of Fatou's
Flower Theorem (Theorem~\ref{lFf}) and \eqref{1070506B}, we also get the following.

\blem\lab{ldp4.4} 
Let $f:\C\lra\oc$ be a meromorphic function. If $\tau >0$ is sufficiently small (with $\th>0$ decreased if necessary), then for every $\om\in\Om_0(f)$ and every $z\in J(f)\cap B(\om,\th)$ there exists $j\in\{1,2,\ld,p(\om)\}$ such that
$$
B_e(z,2\tau|z-\om|)\sbt B(\om,\th)\cap S^j_r(\om,\pi/2).
$$
In addition, all the local holomorphic inverse branches
$$
f_{\om}^{-n}:B(z,2\tau|z-\om|)\lra\C
$$
are well defined for all integers $n\ge 0$. 
\elem

\fr As an immediate consequence of Proposition~\ref{p1j109} and Fatou's
Flowers Theorem (Theorem~\ref{lFf}), we get this. 

\blem\lab{l1042606} 
If $f:\C\lra\oc$ is a meromorphic function, then for every $\om\in\Om_0(f)$, every $j\in\{1,2,\ld,p(\om)\}$, and every $z\in S^j_r(\om,\a)$, in particular for every $z\in J(f)\cap B(\om,\th)$, with $p:=p(\om)$, we have that,
$$
\lim_{n \to \infty}( n^{\frac{p+1}{p}}|f_\om^{-(n+1)}(z)-
f_\om^{-n}(z)|)=|a|^{-\frac{1}{p}},
$$
and 
$$
\lim_{n \to\infty}(n^{\frac{1}{p}}|f_\om^{-n}(z)-\om|)
=|a|^{-\frac{1}{p}}.
$$
In addition, in these two limits, the convergence is uniform on
the set $S^j_r(\om,\a)\cap B^c(\om,R)$, in particular on the set $J(f)\cap B(\om,\th)\cap B^c(\om,R)$, for every $R>0$.
\elem

\sp\fr As an immediate consequence of Proposition~\ref{p1j111} and Fatou's
Flowers Theorem (Theorem~\ref{lFf}), we get the following.

\bprop\label{p1j111DL} 
If $f:\C\lra\oc$ is a meromorphic function, then for every $\om\in\Om_0(f)$
there exists a holomorphic function 
$$
\(f_\om^{-\infty}\)':\bu_{i=1}^{p(\om)}A_j^{*}(\om)\longrightarrow\C\sms\{0\}
$$
with the following properties. 
Fix $\a \in (0,\pi)$ and $j\in\{1,2\ld,p(\om)\}$. Then for every $z\in S^j_r(\om,\a)$, in particular for every $z\in J(f)\cap B(\om,\th)$, with $p:=p(\om)$, we have that
\beq\label{1_2015_10_09}
\lim_{n\to+\infty}
n^{\frac{p+1}{p}}|(f_\om^{-n})'(z)|
= |a|^{-\frac{p+1}{p}}|z-\om|^{-(p+1)}\(f_\om^{-\infty}\)'(z)
\eeq
and the convergence is uniform on $S^j_r(\om,\a)\cap B^c(\om,R)$ for every $R>0$. In addition, if $\a\in (0,\pi/2)$, then
\beq\label{2_2017_08_04L}
\lim_{S^j_r(\om,\a)\ni z\to\om}\(f_\om^{-\infty}\)'(z)=1.
\eeq
\eprop 

As the last result of the local behavior of $f$ around $\om$, we will prove the following.

\bprop[local expansiveness at rationally indifferent parabolic points]\label{p1_2017_10_30}
Let $f:\C\lra\oc$ be a meromorphic function. If $\om\in\Om_0(f)$, $z\in J(f)$, and $f^n(z)\in B(\om,\th)$ for all integers $n\ge 0$, then $z=\om$.
\eprop

{\sl Proof.} Since $z\in J(f)$, we have that $f^n(z)\in J(f)$ for every $n\ge 0$. It therefore follows from Theorem~\ref{lFf} (Fatou's Flower Theorem) and \eqref{dp4.6} (the uniform convergence does matter), that
$$
\om=\lim_{n\to\infty}f_\om^{-n}(f^n(z))=\lim_{n\to\infty}z=z.
$$
\endpf  

We end this section with the following two lemmas which are interesting on their own and will be substantially used in ``parabolic'' sections of Chapter~\ref{invariant-p.s.n.r.}.

\sp Let  
$$
S_r(\om, \alpha):=\bigcup_{j=1}^{p(\om)} S_r^j(\om, \alpha).
$$

\blem~\label{l1ch11.4} If $f: \mathbb C \lra \oc$ is a meromorphic function, $\om \in \Om_0(f)$, and  $ \alpha \in(0, \pi)$, then the set $J(f)\cap \(S_r(\om, \alpha) \sms f^{-1}_\om(S_r(\om, \alpha))\)$ is a fundamental domain for the map $f^{-1}_\om$ acting on $J(f)$ near $\om$. More precisely,
\begin{itemize}
\item [(a)] If $(i,k),  (j,l) \in \{1, 2, \ldots, p(\om)\}\times \mathbb N_0$  and $(i,k)\neq (j,l)$,  then
$$ 
f^{-k}_\om\(S_r^i(\om, \alpha) \sms f^{-1}_\om(S_r(\om, \alpha))\)\cap f^{-l}_\om \(S_r^j(\om, \alpha) \sms f^{-1}_\om(S_r(\om, \alpha))\)=\es$$
and
\item [(b)]  
$$
\om \in \Int_{J(f)}\lt(J(f) \cap \Big(\{\om\}\cup \bigcup_{n=0}^\infty f^{-n}_\om \(S_r(\om, \alpha) \sms f^{-1}_\om(S_r(\om, \alpha))\)\Big)\rt).
$$   
\end{itemize}
\elem

\fr{ \sl Proof.} We start with proving (a). Assume first that $i = j$ but $ k \neq l$ and furthermore, without loss of generality that $ k < l$. Using (\ref{{MU2}p208}) and (\ref{1j101}), we get
$$\begin{aligned}
f^{-k}_\om&\left(S_r^i(\om, \alpha) \sms f^{-1}_\om(S_r(\om, \alpha))\right)\cap f^{-l}_\om \left(S_r^i(\om, \alpha) \sms f^{-1}_\om(S_r(\om, \alpha))\right)=\hspace{5cm}\\
&=f^{-k}_\om \left(\rho \circ H \(S(x(\alpha, \kappa), \alpha)\) \sms  f^{-1}_\om \(\rho\circ H (S(x(\alpha, \kappa), \alpha))\)\right)\cap  \hspace{5cm} \\
   &\hspace{5cm} \cap f^{-l}_\om \left(\rho\circ H (S(x(\alpha, \kappa), \alpha))\sms  f^{-1}_\om \(\rho\circ H (S(x(\alpha, \kappa), \alpha))\)\right)\\
&=f^{-k}_\om \left(\rho\circ H \(S(x(\alpha, \kappa), \alpha)\)\sms \rho\circ H \(\tilde{f}^{-1}_\om (S(x(\alpha, \kappa), \alpha))\)\right)\cap \hspace{5cm}\\
&\hspace{5cm} \cap f^{-l}_\om \left(\rho\circ H \(S(x(\alpha, \kappa), \alpha)\)\sms \rho\circ H \( \tilde{f}^{-l}_\om ( S(x(\alpha, \kappa), \alpha))\)\right)\\
&=f^{-k}_\om \left(\rho\circ H \(S(x(\alpha, \kappa), \alpha)\)\sms\tilde{f}^{-1}_\om \(S(x(\alpha, \kappa), \alpha)\)\right)\cap \hspace{5cm}\\
&\hspace{5cm} \cap f^{-l}_\om \left(\rho\circ H \(S(x(\alpha, \kappa), \alpha)\) \sms  \tilde{f}^{-l}_\om \(S(x(\alpha, \kappa), \alpha)\)\right)\\
&=\rho\circ H \left(\tilde{f}^{-k}_\om \(S(x(\alpha, \kappa), \alpha)\) \sms \tilde{f}^{-1}_\om \(S(x(\alpha, \kappa), \alpha)\)\right)\cap \hspace{5cm}\\
&\hspace{5cm} \cap \rho\circ H \left(\tilde{f}^{-l}_\om \(S(x(\alpha, \kappa), \alpha)\) \sms \tilde{f}^{-1}_\om \(S(x(\alpha, \kappa), \alpha)\)\right)\\
&=\rho\circ H \left(\tilde{f}^{-k}_\om \(S(x(\alpha, \kappa), \alpha)\) \sms \tilde{f}^{-1}_\om \(S(x(\alpha, \kappa), \alpha)\)\right)\cap \hspace{5cm}\\
&\hspace{3cm} \cap \tilde{f}^{-l}_\om \left(S(x(\alpha, \kappa), \alpha) \sms \tilde{f}^{-1}_\om \(S(x(\alpha, \kappa), \alpha)\)\right)\\
&=\rho\circ H \circ \tilde{f}^{-k}_\om \left(S(x(\alpha, \kappa), \alpha) \sms \tilde{f}^{-1}_\om \(S(x(\alpha, \kappa), \alpha)\)\right)\cap \hspace{5cm}\\
&\hspace{5cm} \cap\tilde{f}^{-(l-k)}_\om \left(S(x(\alpha, \kappa), \alpha) \sms \tilde{f}^{-1}_\om \(S(x(\alpha, \kappa), \alpha)\)\right)\\
&\sbt \rho\circ H \circ \tilde{f}^{-k}_\om \left(S(x(\alpha, \kappa), \alpha) \sms \tilde{f}^{-1}_\om (S(x(\alpha, \kappa), \alpha))\right)\cap \tilde{f}^{-(l-k)}_\om (S(x(\alpha, \kappa), \alpha)).
\end{aligned}
$$
Now, since $l-k\geq 1$, employing (\ref{2j103}), we continue as follows:
$$
\begin{aligned}
 f^{-k}_\om\Big(S_r^i(\om, \alpha) &\sms f^{-1}_\om(S_r(\om, \alpha))\Big) \cap  f^{-l}_\om \Big(S_r^i(\om, \alpha) \sms f^{-1}_\om(S_r(\om, \alpha))\Big)\\
 &\sbt  \rho\circ H \circ \tilde{f}^{-k}_\om \left(S(x(\alpha, \kappa), \alpha) \sms \tilde{f}^{-1}_\om (S(x(\alpha, \kappa), \alpha))\right)\cap
\tilde{f}^{-1}_\om \(S(x(\alpha, \kappa), \alpha)\)\\
& =\rho\circ H \circ \tilde{f}^{-k}_\om (\es)=\es.
\end{aligned}
$$
We are thus done with (a) in this case.

Now assume that $ i \neq j$. Then we are immediately done by virtue of Lemma~\ref{l120190617} and formula (\ref{1ms135}) of Theorem~\ref{t1ms135}.

\sp Part (b) is  an immediate consequence of the formula~\ref{dp4.4D} and Fatou's Flower Theorem~3, i.e. Theorem~\ref{lFf}. The proof is complete.\qed 

\sp

\blem\label{l1ch11.5} If  $f: \mathbb C\lra \oc$ is a meromorphic function, $\om \in \Omega_0(f)$,  and $ \alpha \in (0, \pi)$,  then there exists
$u>0$ such that denoting for every $j =1, 2, \ldots, p(\om)$,
$$ \Delta_j(\om, \alpha):= \Delta_j(\om, \alpha, u)=\left( S^j_r(\om, \alpha)\sms f^{-1}_\om(S^j_r(\om, \alpha))\right)\cap \(\mathbb C \sms B(w, u)\)
$$
and
$$ 
\Delta(\om)
:=\Delta(\om , \alpha)
:= \bigcup_{j=1}^{p(\om)}\Delta_j(\om, \alpha),
$$
we have that the set $ \Delta(\om , \alpha)$ is, a fundamental domain for the map $f^{-1}_\om$ acting on $J(f)$ near  $\om$. More precisely,
\begin{itemize}
\item [(a)]	If $(i, k), (j,l) \in \{1,2, \ldots, p(\om)\}\times \mathbb N_0$ and $(i,k)\neq (j,l)$ then 
$$
f^{-k}_\om (\Delta_j\(\om, \alpha))\cap f^{-l}_\om ( \Delta_j(\om, \alpha)\)=\es
$$
and

\item [(b)]	$\om \in\Int_{J(f)}\lt(J(f) \cap \Big(\{\om\}\cup\bigcup_{n=0}^\infty f^{-n}_\om (\Delta(\om, \alpha)\Big)\rt)$.
\end{itemize}
\elem

\bpf Item (a) is an immediate consequence of Lemma~\ref{l1ch11.4}.  We now shall prove (b). Since $f^{-1}_\om$ is a $C^1$-diffeomorphism, there is $ u \in (0, \theta_\alpha(\om))$ such that  
$$
B(\om, u)\sbt f^{-1}_\om(B(\om, \theta_\alpha(\om))).
$$
It then follows from Fatou's Flowers Theorem, i.e  Theorem~\ref{lFf}, that
 \beq\label{1Ch11_5}
 J(f) \cap B(\om, u)\sbt \bigcup_{j=1}^{p(\om)} f^{-1}_\om (S_r^j(\om, \alpha)) =f^{-1}(S_r(\om, \alpha)).
 \eeq
Because of Lemma~\ref{l1ch11.4} (b), it is enough to show that
$$
J(f)\cap \left(\bigcup_{n=0}^{\infty} f^{-n}_\om\Big(S_r(\om, \alpha)\sms f_\om^{-1}(S_r(\om, \alpha))\Big)\right)
\sbt \bigcup_{n=0}^\infty  f^{-n}_\om ( \Delta(\om, \alpha)).
$$
So, let $z$ belong to the left--hand side of this inclusion. Then there exists $ n \geq 0$ such that
$$ 
z \in J(f)\cap f^{-n}_\om \(S_r(\om, \alpha)\sms f^{-1}_\om (S_r(\om, \alpha))\).
$$
But then, because of (\ref{1Ch11_5}), we get
$$
\begin{aligned}
f^n(z) 
&\in J(f)\cap \(S_r(\om, \alpha)\sms f^{-1}_\om (S_r(\om, \alpha))\) \\
&\sbt J(f)\cap \(S_r(\om, \alpha)\sms f^{-1}_\om (S_r(\om,\alpha))\) \cap \(\mathbb C \sms B(\om, u)\)\\
&\sbt \Delta (\om, \alpha).
\end{aligned}
$$
Hence $ z \in f^{-n}(\Delta (\om, \alpha))$ and the proof is complete.\qed

\section[Local and Asymptotic Behavior of Meromorphic Functions: Part III] {Local and Asymptotic Behavior of (General) Meromorphic Functions Around Rationally Indifferent Periodic Points; \\ Part III: 
Conformal Measures}\label{local-parabolic III} 

\sp We now, in this section, pass to deal with the local behavior of conformal measures and its generalizations around parabolic points. Recall that in Definition~\ref{1d20120909} we have introduced the concept of 
semi--conformal and conformal measures. Passing to deal with them, we shall prove the following.  

\sp
\blem\lab{ldp4.5} 
Let $m$ be a semi $t$--conformal measure for a meromorphic map
$f:\C\lra\oc$ defined on some neighborhood of a simple parabolic fixed point $\om$ of $f$. Then for every $R>0$ there exists a constant $C=C(t,\om,R)\ge 1$ such that for every $0<r\le R$
$$
{m(B(\om,r)\sms\{\om\})\over r^{\a_t(\om)}} \le C,
$$
where $\a_t(\om):=t+p(\om)(t-1)$\index{(S)}{$\a_t(\om)$}. If $m$ is
$t$--conformal, then in addition
$$
{m(B(\om,r)\sms\{\om\})\over r^{\a_t(\om)}} \ge C^{-1}.
$$
\elem

\bpf Take $R>0$ so small, for example $R=\th/2$, that
$\overline{B_e(\om,R)}\sbt B(\om, \th)$ and let
\beq\label{320190603}
P := J(f)\cap \Big\{z\in\C: R (2\|f_\om'\|)^{-1}\le |z-\om| \le R \Big\},
\eeq
where 
$$
\|f_\om'\|_\infty:=\sup\Big\{|f'(z)|: z\in \ov{B(\om, \th)}\cap J(f)\Big\}.
$$
Fix $\tau>0$ so small as needed in Lemma~\ref{ldp4.4}. Let
$$
\d = \tau \inf\{|z-\om|: z\in P\}> 0.
$$
Since $P$ is compact, there are finitely many points $z^1,\ld,z^q$
in $P$, such that
$$
P\sbt J(f)\cap \(B(z^{1},\d)\cup\ld\cup B(z^{q},\d)\),
$$
and we may assume that $\d$ is so small that
$$
f_{\om}^{-n}(B(z^i,\d)) \cap B(z^i,\d) = \es
$$
for  $i=1,\ld,q$ and  $n=1,2,\ld$. Denote
$$
p:=p(\om).
$$
For every $n\ge 1$ define:
$$
P_n := \lt\{z\in B(\om, \th)\cap J(f): \frac12|a|^{\frac1p}n^{-1/p}\le
|f_\om^{-n}(z)-\om|\le 2|a|^{-\frac1p}n^{-1/p}\rt\}.
$$
By the local behavior of $f$ around a  parabolic point we conclude
that for every 
$$
z\in (B(\om,R )\sms \{\om\})\cap J(f)
$$ 
there exists an integer $l\ge 0$ such that
$$ 
R(2\|f_\om'\|_\infty)^{-1}< |f^l(z)-\om| < R,
$$
i.e. $f^l(z) \in P$. Therefore, the set
$$
J(f)\cap \big\{z:  R (2\|f_\om'\|_\infty)^{-1}< |z-\om|<R\big\}
$$ is non--empty ($B_e(\om,R)\cap (J(f) 
\sms\{\om\})$ is non--empty since $J(f)$ is perfect). Moreover, since
it is open in $J(f)$, we deduce that for some $1\le j\le q$ the set
$B_e(z^j,\d)\cap P$ has non-empty interior in $J(f)$. Hence
$$
M:=m_e(B_e(z^j,\d)\cap P)>0.
$$
By Lemma~\ref{l1042606} there is an integer $n_0\ge 1$ such that $
f_\om^{-n}(z) \in P_n$ for every $n\ge n_0$ and and every $z\in P$. In other words this means that $P_n\supset f_\om^{-n}(P)$ for $n\ge n_0$.
Thus,
\beq\label{120190603}
B\(\om,2|a|^{-\frac1p}n^{-1/p}\)
\supset \bu_{k=n}^\infty P_k\supset
\bu_{k=n}^\infty f_\om^{-k}(P)\supset \bu_{k=n}^\infty \bu_{i=1}^q
f_{\om}^{-n}(B(z^i, \d)\cap P).
\eeq
On the other hand, because of Proposition~\ref{p1_2017_10_30}, for every $z\in J(f)\cap(B(\om, R)\sms\{\om\})$ there exists an integer $l\ge 0$ such that $f^l(z)\in P$. Let
$l(z)\ge 0$ be such smallest integer. Take
$n_1\ge n_0$ so large that if $z\in B(\om,2|a|^{-\frac1p} n_1^{-1/p})$, then
$l(z)\ge n_0$. Consider now any 
$$
z\in J(f)\cap B(\om,2|a|^{-\frac1p}n^{-1/p})\sms\{\om\}
$$ 
with $n\ge n_1$. Since $l(z)\ge n_0$ and
$f^{l(z)} (z)\in P$ we see that
$z=f^{-l(z)}_\om(f^{l(z)}(z))\in P_{l(z)}$. Therefore
$\frac12|a|^p l(z)^{-1/p} \le 2|a|^{-\frac1p}n^{-1/p}$ and consequently $l(z)\ge 4^{-p}|a|^2$. Hence,
\beq\label{220190603}
\begin{aligned}
J(f)\cap B(\om,2|a|^{-\frac1p}n^{-1/p})
&\sbt \{\om\}\cup \bu_{l\ge 4^{-p}|a|^2n} f_\om^{-l}(P) \\
&= \{\om\}\cup \bu_{i=1}^q \bu_{l\ge  4^{-p}|a|^2n}
f_{\om}^{- l}(B(z^i,\d)).
\end{aligned}
\eeq
Since for every $1\le j\le q$ the sets $\{f_{\om}^{-n}(J(f)\cap B(z^j,\d))\}_{n=0}^\infty$, are mutually disjoint, it follows from
Theorem~\ref{Euclid-I}, Proposition~\ref{p1j111DL} and semi-conformality of the measure $m$that
$$ 
\begin{aligned}
m\(B(\om, 2|a|^{-\frac1p}n^{-1/p})\sms \{\om\} \) 
& \le m\Bigl(\bu_{i=1}^q\bu_{l\ge 4^{-p}|a|^2n} f_{\om,i}^{-l}(B(z^i,\d))\Bigr)\\
& \le qK^tC^t \sum_{l\ge 4^{-p}|a|^2n}l^{-\frac{p+1}pt} \\
&\le C'(n^{-1/p})^{\a_t(\omega)}
\end{aligned}
$$
where $C>0$ and $C'$ are some constants. If, in addition, $m$  is
$t$-conformal we have
$$\begin{aligned}
m\(B(\om,2|a|^{-\frac1p}n^{-1/p})\sms \{\om\} \) &\ge \sum_{k=n}^\infty
m_e(f_\om^{-k}(B_e(z^j,\d)\cap P) \\ 
&\ge \sum_{k=n}^\infty K^{-t}C^{-t}(k^{-\frac{p+1}p})^tM\\
&\ge MK^{-t}C^{-t}\sum_{k=n}^{2n} (k^{-\frac{p+1}p})^t \\
&\ge MK^{-t}C^{-t}n((2n)^{-\frac{p+1}pt}) \\
& =2^{-\frac{p+1}pt} MK^{-h}C^{-t} (n^{-1/p})^{\a_t(\om)}.
\end{aligned}
$$
The proof is now concluded by observing that
$\lim_{n\to\infty}\frac{(n+1)^{-1/p}}{n^{-1/p}} = 1$. 
\endpf

\

\fr As a subproof of this proof we obtained the following.

\sp\blem\label{l2.7a} 
Let $m$ be a semi $t$--conformal measure for a meromorphic map
$f:\C\lra\oc$ defined on some neighborhood of $\Om(f)$.
Then for all $R>0$ sufficiently small there exists a constant
$C(R)>0$ such for all $\om  \in \Omega(f)$ and for all $k \geq 0$, we
have
$$
m(f^{-q(\om)k}_\om (B(\om, R) \cap J(f)\sms\{\om\})
\comp C(R)( k+1)^{1-\frac{p(\om)+1}{p(\om)}t},
$$
where $q(\om)$ is the smallest integer $\ge 1$ turning $\om$ into
simple parabolic point, and $p(\om)$ is the number of its petals.
\elem

\sp\fr As an immediate consequence of formulas \eqref{120190603} and \eqref{220190603}, we get the following.

\sp\blem\label{l120190603}
If $f:\C\lra\oc$ is a meromorphic map defined on some neighborhood of a simple parabolic fixed point $\om$ of $f$, then for every integer $n\ge 0$, we have that
$$
\Int_{J(f)}\Big(J(f)\cap\Big(\{\om\}\cup\bu_{k=n}^\infty f_\om^{-k}(P)\Big)\Big)\ne\es,
$$
where the set $P$ is defined by formula \eqref{320190603}.
\elem

\sp\blem\label{confmeasparabsector} 
If $m$ is a semi $t$--conformal measure for a meromorphic map
$f:\C\lra\oc$, defined (i.e. $m$) on some neighborhood of $\Om(f)$, then \nl 
$\forall \, \beta >0 \  \  \exists \, C_\b\ge 
1 \  \  \forall \, \om\in\Om(f) \  \  \forall \, z\in J(f)$,
$$
m(B(z,\beta|z-\om|)) \ge C_\b|z-\om|^{\a_t(\om)}.
$$
\elem

\bpf Since the set $\Om(f)$ is finite, it suffices to show the lemma
for some fixed $\om\in\Om$ and every $0<\beta<1$. Furthermore, passing to a sufficiently high iterate of $f$, we may assume that $\om$ is simple. Denote
$$
p:=p(\om).
$$
For any $\xi\in\C$ let
$$
A(\xi):= \big\{x\in\C:(1-\beta)|\xi-\om| \le |x-\om| \le (1+\beta)|\xi-\om| \big\}.
$$
Observe that there exists $0<\a\in(0,\pi)$ such that, if
$j\in\{1,2,\ld,p(\om)\}$ and $\xi\in S_r^j(\om,\a)$, then
$$
A(\xi)\cap S_r^j(\om,\a)\sbt B(\xi,2\beta|\xi-\om|).
$$
Now apply the construction of the proof of Lemma~\ref{ldp4.5}. We use the same
notation and, due to Theorem~\ref{lFf}, may assume in addition that 
$$
P\sbt B(\om,\th_\a(\om))
$$
and that the radius $\d$ of balls $B(z^1,\d),\ld, 
B(z^{q},\d)$ is so small that each of these balls $B(z^i,\d)$ is 
contained in exactly one sector $S_r^j(\om,\a)$, $j\in\{1,\ld,p\}$. Choose $\e>0$ so small that 
$$
s:=\left(\frac{e+\e}{(1+\beta)(e-\e)}\right)^p <1 \qquad 
\text{and} \qquad 
u:=\left(\frac{e-\e}{(1-\beta)(e+\e)}\right)^p>1, 
$$
where $e=|a|^{-1/p(\om)}$. With these definitions, by virtue of
Lemma~\ref{l1042606} there exists $n_0\ge 1$ so large that for every
$n\ge sn_0$ and every $y\in P$,
$$
[un]-([sn]+1)\ge \frac{u-s}2 n,
$$
where $[t]$ denotes the integer part of $t$, and
$$
(e-\e)n^{-1/p} \le |f_\om^{-n}(y)-\om| \le (e+\e)n^{-1/p},
$$
where $p=p(\om)$. Let for every $n\ge 1$,
$$
A_n:= \big\{x\in\C:(1-\beta)(e+\e)n^{-1/p} \le |x-\om| \le (1+\beta) (e-
\e)n^{-1/p}\big\}.
$$
By Theorem~\ref{lFf}, for every $\xi\in J(f)\cap B(\om,\th_\a(\om))\sms
\{\om\}$ there  exists exactly one $1\le j(\xi)\le p$ such that $\xi\in S_r^{j(z)}(\om,\a)$. 
Moreover, $j\(f_\om^{-n}(\xi)\)=j(\xi)$ for every $\xi\in J(f)\cap B(\om,\th_\a(\om)) \sms\{\om\}$ and every integer $n\ge 0$. Hence, for every $y\in P$ and every $n\ge n_0$ we have 
$$
A_n\cap S_r^{j(y)}(\om,\th_\a(\om)) \sbt B\(f_\om^{-n}(y),2\beta|z_n-\om|\).
$$
Now consider an arbitrary 
$$
z\in \bu_{n\ge n_0}f_\om^{-n}(P).
$$
Then $z\in f_\om^{-n}(P)$ for some $n\ge n_0$. Let $x:=f^n(z)$ and choose $1\le 
i\le q$ such that $x\in B(z^i,\d)$. Having the radius $\d>0$ is small enough, it follows from Theorem~\ref{lFf} that for every $l\ge 0$, we have that
$$
f_\om^{-l}(J(f)\cap B(z^i,\d)) \sbt S_r^{j(x)}(\om,\a).
$$
If $y\in P$ and $sn_0\le sn\le l\le un$ then, by the choice of $\e$ and 
$n_0$, we have that
$$
|f_\om^{-l}(y)-\om|
\le (e+\e)l^{-1/p}
\le (e+\e)s^{-1/p}n^{-1/p}
= (1+\beta)(e-\e)n^{-1/p}
$$
and
$$
|f_\om^{-l}(y)-\om|
\ge (e-\e)l^{-1/p}\ge (e-\e)u^{-1/p}n^{-1/p}
= (1-\beta)(e+\e)n^{-1/p}. 
$$
Thus $f_\om^{-l}(P)\sbt A_n$, whence
$$ 
\aligned
\bu_{sn\le l\le un}f_\om^{-l}(J(f)\cap B(z^i,\d)) 
&\sbt A_n\cap S_r^{j(x)}(\om,\a)  \\
&\sbt B\(f_\om^{-n}(x), 2\beta\big|f_\om^{-n}(x)-\om\big|\) \\
&= B(z,2\beta|z-\om|).
\endaligned
$$
Since the sets $f_\om^{-l}(J(f)\cap B(z^i,\d))$, $l=[sn]+1,\ld, 
[un]+1$, are mutually disjoint, it follows from
Proposition~\ref{p1j111DL} and Koebe's Distortion Theorem
(Theorem~\ref{Euclid-I}) that with some constant $C_4>0$,
$$
m\lt(\bu_{sn\le l\le un}f_\om^{-l}(J(f)\cap B(z^i,\d))\rt)\ge 
\sum_{l=[sn]}^{[un]}C_4l^{-\frac{p+1}ph}
$$ 
and therefore, 
$$
\aligned
m(B(z,2\beta|z-\om|)) 
&\ge C_4\frac{u-s}2nu^{-\frac{p+1}ph}n^{-\frac{p+1}ph}
=C_4\frac{u-s}2u^{-\frac{p+1}ph}(n^{-1/p})^{h+p(h-1)} \\ 
&\ge C_4u^{-\frac{p+1}ph}\frac{u-s}{2(e+\e)^{\a(\om)}}
     \big|f_\om^{-n}(x)-\om\big|^{\a(\om)}\\
&= C_4u^{-\frac{p+1}ph}\frac{u-s}{2(e+\e)^{\a(\om)}}|z-\om|^{\a(\om)}. 
\endaligned
$$
This proves the lemma for all points in the set 
$$
G:=J(f)\cap\bu_{n\ge n_0}f_\om^{-n}(P).
$$
Since for $\om$ itself there is nothing to prove, the 
statement of the lemma holds on the set $J(f)\cap (\{\omega\}\cup G)$  which is, by Lemma~\ref{l120190603}, a neighborhood of $\omega$ in $J(f)$. Since $m$ is positive on non--empty open sets, for $ z\notin
\{\om\}\cup G$ we have that
$$
m(B(z,2\beta|z-\om|)) \ge \const>0,
$$
finishing the proof of our lemma. 
\endpf

\sp\blem\lab{lncp14.7.} 
Suppose $m$ is a semi $t$--conformal measure for the meromorphic map
$f:\C\lra\oc$ defined on some neighborhood of $\Om(f)$. Then for every
$\om\in\Om$, every $R>0$, and every $0<\sg\le 1$ 
there exists $L=L(\om,R,\sg)>0$ such that for every $0<r\le R$ the
point $\om$ is $(r,\sg,L)$-$\a_t(\om)$-s.l.e. with respect
to the measure $m$. 
\elem

\bpf Let $z\in B_e(\om,r)$. If $\sg r\ge 2|z-\om|$, then
$B_e(z,\sg r)\spt B_e\(\om,{\sg\over 2}r\)$ and therefore by
Lemma~\ref{ldp4.5}
\beq\label{1070806}
\nu_e(B_e(z,\sg r))\ge C(R/2)\lt({\sg\over
2}r\rt)^{\a_t(\om)}=\lt({\sg\over
2}\rt)^{\a_t(\om)}C(R/2)r^{\a_t(\om)}.
\eeq
In order to deal with the opposite case first notice that always
$$
B_e(z,\sg r)\spt B_e(z,\sg|z-\om|)
$$ and therefore, by Lemma~\ref{confmeasparabsector},
we have 
$$
\nu_e(B_e(z,\sg r))\ge C^{-1}(\sg)|z-\om|^{\a_t(\om)}.
$$ 
As $\sg r< 2|z-\om|,$ this implies that 
$$
\nu_e(B_e(z,\sg r))
\ge C^{-1}(\sg)\({\sg/2}\)^{\a_t(\om)}r^{\a_t(\om)}.
$$ 
So, putting
$$
L(\om,R,\sg)=({\sg/2})^{\a_t(\om)}\min\{C(R/2),C^{-1}(\sg)\}
$$
finishes the proof. \endpf

\sp 

\section{Nice Sets for Analytic Maps}\label{NiceSetsGeneral}

This section of this chapter is somewhat different than the previous ones. It deals with a powerful tool of nice sets. Indeed, in this section  we introduce the concept of pre--nice sets and nice sets \index{(N)}{pre--nice sets} \index{(N)}{pre--nice sets} \index{(N)}{nice sets} and prove their  existence. Nice sets naturally  merged
in dynamical  systems in  the   context   of self--maps of an interval, where  their existence was  sort of  obvious.
They were adapted,  although in a somewhat  obscure way, to holomorphic endomorphisms of the  Riemann sphere, by  Juan Rivera-Letelier
in \cite{Ri}. A much clearer proof of the existence of nice sets in a more general  context of maps that are meromorphic, transcendental or
rational alike, from the complex  plane  to the Riemann sphere was provided by N. Dobbs in \cite{Do}. Nice sets are  a powerful  tool
indeed. They formed a central theme in the fairly  complete  treatment of  Collect--Eckmann rational functions given
in \cite{PR}, as well in later papers of Feliks Przytycki and Rivera--Letelier. 

In this section we  define pre--nice sets and nice sets  for holomorphic  maps of the Riemann surfaces (one of which is an open subset of the other).
We  then prove their existence. We need  such generality in order to  deal with  projected  maps:
$$
\hat f: \mathbb{T}_f \sms \Pi_f(f^{-1}(\infty))\lra \mathbb{T}_f.
$$

Nice sets  give  naturally  rise to  conformal
infinite iterated function systems systematically explored in Chapter~\ref{Markov-systems}. 
In this way  many problems  of transcendental and rational  holomorphic dynamics  can be  successfully treated (or even reduced to)
by means of  the theory  of conformal iterated  function systems, or  their generalization  formed by  graph directed Markov systems, presented, as just said above, in Chapter~\ref{Markov-systems}. 

In this book we apply the technique of nice sets, via means of iterated  function systems,
to demonstrate finitenes of Krengel's entropy for  the dynamics generated by  an elliptic  function (and invariant measure equivalent to the conformal one), and, most notably, to prove refined stochastic laws in various cases of finite invariant measures.

\sp We start with somewhat long and involved topological preparations.

\bprop\label{p1fp29.3}
Let $W$ be a  subset of a Hausdorff topological space $Y$. Let $V$ be a subset of $W$ such that $W\sms V$ is connected. Let $G$ be an open subset of a Hausdorff topological space $X$ such that $\overline{G} \sbt X $ is compact. If $f:X \to Y$ is an open continuous map such that
$$
f(\partial{G})\sbt V
\  \  \  {\rm and} \  \  \
f(G)\nsupseteq W\sms V,
$$
then 
$$
f(G)\sbt V.
$$
\eprop

\bpf Since $G$ is open, $ \overline{G}$ is compact, and $f$ is continuous and open, we  have that  $\partial{f(G)} \sbt f(\partial{G})\sbt V$. Therefore,
$$ 
\partial_{W \sms V}{\(f(G)\cap (W \sms V)\)}
\sbt (W \sms V)\cap \partial{f}(G)\sbt (W\sms V)\cap f(\partial{G})
\sbt (W \sms V)\cap V =\es.$$
Thus,
$$
\partial_{W\sms V}{\(f(G)\cap (W \sms V)\)}=\es.
$$
Since the set $W\sms V$ is connected, this implies that either
$$
f(G)\cap (W\sms V)=W\sms V
\  \  \  {\rm or} \  \  \ 
f(G) \cap (W\sms V)=\es.
$$
If the first part of this alternative holds, then $f(G) \supset W\sms V$,  contrary to our  hypothesis. So, the other part of this alternative holds, meaning that $f(G)\sbt V$. The proof is complete. \qed

\sp Since every non-constant holomorphic mapping between Riemann surfaces is  open, as an immediate consequence of this proposition, we get the following.

\bcor\label{c1fp29.4}
Let $W$ be a subset of a boundaryless Riemann surface $Y$. Let $V$ be a subset of $ W$   such $W\sms V$ is connected. Let $G$ be an open subset of a Riemann surface $X$ such that  $\overline{G} \sbt X$ is compact. If $f: X \to Y$ is a holomorphic map, $f(\partial{G}) \sbt V$ and $f(G) 
\not\varsupsetneq W\sms V$, then 
$$
f(G)\sbt V.
$$
\ecor

\bprop\label{p2j144}
Let $X$ and $Z$  be two  boundaryless Riemann surfaces. 

\begin{enumerate}
\item Let $D\sbt X$ be an open set  conformally  equivalent to the unit disc $\mathbb D$. 

\,

\item Let $W$ be a subset  of $Z$  and let $V$  be  a subset of $W$ such that $W\sms V$ is connected. 

\,

\item Let $G \sbt D$ be an open  connected set such that $\bar{G} \sbt \mathbb D$ is compact. 

\,

\item Let $\Gamma $ be the only connected  component of $D\sms G$ which is not compact. 

\,

\item Let $\hat{G}:=D\sms \Gamma \supset G$. 
\end{enumerate}

Then $\hat{G}$ is open connected simply connected and $\partial{\hat{G}}\sbt \partial{G}$. Furthermore, if $f:D \to  Z$ is a holomorphic map,  $f(\partial{G})\sbt V$ and $f(D)\not\varsupsetneq W\sms V$, then 
$$
f(\hat{G}) \sbt V.
$$
\eprop

\sp \fr {\sl Proof.} The set $\hat G$ is open because $D$ is open and $\Ga$ is closed in $D$. 
As an immediate consequence of the (very general) Theorem~5, page 140, in \cite{Kur2}, we conclude that $\hat G$ is connected. Simple connectivity of $\hat G$ follows from Theorem~4.4, page 144, in \cite{Newman}. 

 By the definition of  $\Gamma$  we have that $ G  \sbt \bar{G} \sbt \overline{\hat{G}} \sbt  D$. Since also both sets $G$ and $\hat{G}$ are open, we get that 
   \beq\label{1fp29.5}
    \partial{G}=\partial_D{G}\quad \text{and} \quad \partial{\hat{G}}= \partial_D\hat{G}.
    \eeq
Because of Theorem~3, page 238 in \cite{Kur2} (which applies since $X$ is locally connected), we have that
$$  \partial_D{\hat{G}}=\partial_D{(D \sms \Gamma)} =\partial_D{\Gamma}\sbt \partial_D{(D\sms G)}=\partial_D{G}.$$
Along  with (\ref{1fp29.5}) this gives that $\partial{\hat{G}}\sbt \partial{G}$ and the proof of the first part of our  proposition is complete. Therefore,
$$f(\partial{\hat{G}})\sbt f(\partial{G})\sbt V.$$
Since $\hat{G}\sbt D $ and $f(D) \not\varsupsetneq W \sms V$, we have also that $f(\hat{G})\not\varsupsetneq W\sms V$. Thus  a direct application of Corollary~\ref{c1fp29.4}  gives that $f(\hat{G})\sbt V$. The proof of Proposition~\ref{p2j144} is complete.\qed

\sp  Let $Y$ be a complete Riemann surface with constant curvature $0$, i.e. $Y$ is either the complex plane $\C$, a complex torus $\mT_\La=\C/\La$, where $\La$ is a lattice on $\C$ or an infinite cylinder $\C/2\pi i\Z$. Let $X$ be a non-empty open subset of $Y$. Let 
$$
f:X \lra Y
$$ 
be an  analytic map. As in the case of meromorphic functions, we say that $y \in Y$ is a regular point \index{(N)}{regular point} of $f^{-1}$ if for every  $r>0$ small enough and every connected component $C$ of  $f^{-1}(B(y,r))$ the restriction of 
$$
f|_C: C\lra B(y,r)
$$ 
is a (conformal) homeomorphism from $C$ onto $B(y,r)$. Otherwise we say that $y$ is a singular point \index{(N)}{singular point} of $f^{-1}$ and we denote by $\Sing(f^{-1})$ \index{(S)}{$\Sing(f^{-1})$}the set of all such singular points. We also set
$$
{\rm PS} (f):=\bigcup_{n=0}^\infty f^n(\Sing (f^{-1}))
$$
with the convention that  
$$
f(\{z\}) =\es
$$  
if $z\in Y\sms X$.

\sp\bdfn\label{d1j291-}
A non--empty set $V\sbt Y$ is said to be a pre--nice set \index{(N)}{pre--nice sets} for the analytic  map $f: X \to Y$ if the following conditions are satisfied:

\sp\begin{itemize}
\item [(a)] $\ov{V}$ is compact.
 
\sp\item [(b)] $V$ has finitely many connected  components.

\sp\item [(c)] Every connected component of $V$ is simply connected.

\sp \item [(d)] $V \cap \bigcup_{n=0}^\infty f^n(\partial{V})=\es$.
\end{itemize}

\edfn

\sp\bdfn\label{d1j291}
A non--empty open set  $V \sbt Y$  is said to be a nice set \index{(N)}{nice sets} for the analytic  map $f: X \to Y$  if the following  conditions  are satisfied:
 
\sp\begin{itemize}
\item [(a)] $\ov{V}$ is compact.
 
\sp\item [(b)] $V$ has  finitely many connected  components.

\sp\item [(c)] If $\Ga$ is a connected component of $V$, then $\Ga$ is simply connected and  there exists $\hat{\Ga}$, an open  connected simply connected subset of $Y$ containing $\ov\Ga$ such that
$$
\ov{\hat{\Ga}}\cap \ov{\rm PS(f)} =\es. 
$$
\item [(d)] $V \cap \bigcup_{n=0}^\infty f^n(\partial{V})=\es$.

\sp\item [(e)] There exists $\l>1$ such that 
$$
|(f^n)'(z)|\ge \l \  \  
$$
for every $n\ge 1$ and all $z\in V\cap f^{-n}(V)$.
\end{itemize}
\edfn

\

\fr Of course, every pre--nice set is nice. Before establishing the existence of nice sets and pre--nice sets we shall draw the  most characteristic and particularly useful, for us, properties of such sets. We start with the following.

\sp

\bprop\label{p1j290} 
Suppose that $V$ is a  pre--nice set for a holomorphic map $f :X \to Y$. Let $U$ and $W$  be two  distinct components  of $V$. If $j, k \geq 0$
are  two  integers and  $A$ and $B$  are connected  components  respectively  of $f^{-j}(U)$ and $f^{-k}(W)$, then  either

\sp\begin{itemize}
\sp \item [(a)]   $A \cap B = \es $
\sp\item [(b)]  $A \sbt B$

or  

\sp\item [(c)]  $B \sbt A$.
\end{itemize}
\eprop

\bpf Assume without loss of generality  that $j \leq k$. Seeking  contradiction assume that $A \cap B \neq \es$ but  neither $ A \sbt B$ nor $ B \sbt A$. Then the  first
and the third  of theses properties yield
$ B \cap  \partial {A}\neq \es$.  As  $f^j(\partial{ A})= \partial{U}$, we thus get
$$ 
W \cap f^{k-j}( \partial{U})= f^k(B) \cap f^k (\partial  A) \supset f^k( B \cap \partial{A}) \neq \es.
$$
Since $\partial{U} \sbt \partial{V}$ and $ W \sbt V$, this  gives $V \cap f^{k-j}(\partial{V})\neq\es $, contrary   to condition (d) of Definition~\ref{d1j291-}. The  proof is finished.\qed

\sp Now, given a nice set $V$ let $\mathcal C_1^\infty(V)$ \index{(S)}{$\mathcal C_1^\infty(V)$}be the family of all connected components of the set $ V \cap \bigcup_{n=1}^\infty  f^{-n}(V)$. Proposition~\ref{p1j290}   entails the following.

\sp\bprop\label{p2j290}
Suppose that $V$ is a pre--nice set for a holomorphic map $f: X \to Y$. Then the following hold.
\begin{itemize}
\item [(a)] If $\Ga\in \mathcal C_1^\infty(V)$ \index{(S)}{$\mathcal C_1^\infty(V)$} then there exists a unique  connected component $W$ of $V$ and a unique integer $n(\Ga)\geq 1$ such that $\Ga$ is a connected component of $f^{-n(\Ga)}(W)$.

\, 

\item [(b)] If $W$ is a connected component of $V$, $n\geq 1$, and $\Ga$ is a connected component of $f^{-n}(W)$ such that 
$$
\Ga\cap V \neq \es 
\  \  \  {\rm and} \  \  \
f^k(\Ga) \cap V =\es
$$ 
for all  $1 \leq  k <n$, then $\Ga \in{\mathcal C_1^\infty}(V) $ and  $n(\Ga)=n$.
\end{itemize}
\eprop

\sp\fr This proposition implies that if $V$ is a nice set, then for every $U \in {\mathcal C_1^\infty}(V)$ there exists a unique connected component $\Ga^*$ of $V$ and a unique holomorphic inverse branch $f^{-n(\Ga)}_{{\hat \Ga}^*}: \hat{\Ga}^*\to X$ of $f^{n(\Ga)}$ such that  
$$ 
f^{-n(\Ga)}_{\Ga}(\Ga^*)=\Ga.
$$
We therefore  obtain the following fundamental consequence of the existence of nice sets.

\sp

\bthm~\label{t1j290a} If $f: X \lra Y$ is an analytic  map and $V$  is a nice set for $f$, then 
$$
{\mathcal S}_V=\big\{ f^{-n(\Ga)}_\Ga:\hat{\Ga}^*\lra X\big\}_{\Ga \in {\mathcal C_1^\infty}(V)}
$$
is a maximal conformal graph directed Markov system \index{(N)}{maximal conformal graph directed Markov system} in the sense of Chapter~\ref{Markov-systems}, maximality of a system being defined in Definition~\ref{d2_2017_11_18}.
\ethm
The main technical result of this subsection concerns the existence of pre--nice sets:

\sp\bthm\label{t2j291}   
Let $Y$ be a complete Riemann surface with constant curvature $0$, i.e. $Y$ is either the complex plane $\C$, a complex torus $\mT_\La=\C/\La$, where $\La$ is a lattice on $\C$ or an infinite cylinder $\C/2\pi i\Z$.

Fix $R>0$ so small that every open ball in $Y$, with radius $R$, is simply connected, i.e. isometrically equivalent $B(0,R)\sbt\C$. Fix $\kappa\in(1,2]$.

Let $X$ be a non-empty open subset of $Y$ and let $f:X \lra Y$ be an analytic map. 

\fr Fix 

\begin{itemize}

\,

\item $F$, a finite subset of $J(f)\sms\ov{{\rm PS}(f)}$, 

\,

\item a collection $\{U_0(b)\}_{b\in F}$ of open subsets of $Y$, 

\fr and 

\,

\item a vector $\bold r=\(r_b:b\in F\)$ with the following properties:
\end{itemize}

\, 

\begin{enumerate}  
\item
$$
U_0(b)\sbt B(b,r_b)
$$
for all $b\in F$.

\, \, \item
$F$ can be represented as a disjoint union:
$$
F=F_0\cup F_1
$$
such that
\begin{enumerate}

\, \, \item $f(b)=b$ for every  $b\in F_0$.

\, \, \item  For each $b\in F_0$, the map $f|_{B(b,6r_b)}$ is 1--to--1 and the holomorphic inverse branch $f_b^{-1}:B(b,6r_b)\lra X$, sending $b$ to $b$,  is well defined.

\, \, \item  $f_b^{-1}\(B(b,3r_b)\)\sbt B(b,6r_b)$.

\, \, \item $b\in \bd U_0(b)$ for all $b\in F_0$,

\, \, \item For each $b\in F_0$ there exists an integer $p_b\ge 1$ such that $U_0(b)$ has exactly $p_b$ (open) connected components $U_0(b,j)$, $j=1,2,\ld, p_b$, each of which is simply connected and its closure is a closed topological disk.

\, \, \item For each $b\in F_0$ and each $j=1,2,\ld, p_b$:
$$
f_b^{-1}(U_0(b,j))\sbt U_0\(b,j)
$$
and

\, \, \item The sequence 

$$
\Big(f_b|_{U_0(b,j)}^{-n}:U_0(b,j)\lra X\Big)_{n=0}^\infty
$$ 
converges uniformly to the constant function which assignes to each point in $U_0(b,j)$ the value $b$.

\, \, \item $B(b,r_b/2)\sbt U_0(b)$ for all $b\in F_1$,

\, \, \item For each $b\in F_1$ the set $U_0(b)$ is connected, simply connected, and its closure is a closed topological disk. We then put $p_b:=1$ and denote also $U_0(b)$ by $U_0(b,1)$.

\, \, \item  For each $b \in F_0$ and $s>0$ there exists $s^{-}_b \in (0, r_b]$ such that if
$$ 
z \in B(b, s^{-}_b)\sms \bigcup_{j=1}^{p_b}U_0(b, j),
$$
then 
$$
f^n(z) \in B(b,s)
$$ 
for all $n\geq 0$.
\end{enumerate}

\, \, \item
$$
r_b\in \lt(0,\frac{1}{6}\min\Big\{R,\min\big\{\rho(b,c):c \in F\sms \{b\}\big\}\Big\}\rt)
$$
for all $b\in F_0$.

\, \, \item
$$
r_b\in \lt(0,\frac{1}{6}\min\Big\{R,\min\big\{\rho(b,c):c \in F\sms \{b\}\big\}, \dist\(b, \ov{\rm PS}(f)\)\Big\}\rt)
$$
for all $b\in F_1$.

\, \, \item 
Suppose that $a, b\in F$ and $n\ge 1$ is an integer. Assume that either $b\in F_1$ or there exists a point 
\beq\label{120180824}
w\in f^{-n}(U_0(b))
\eeq
such that 
\beq\label{220180824B}
f^{n-1}(w)\notin U_0(b),
\eeq
then
 
\begin{enumerate}
\, \, \item 

the holomorphic inverse branch $f_w^{-n}:B(b,6r_b)\to X$ of $f^n$, sending $f^n(w)$ to $w$, is well--defined (this is only an extra hypothesis if $b\in F_0$; if $b\in F_1$, this follows from (4)).

\, \, \item 
If, in addition, the connected component of $f^{-n}(B(b,2r_b))$ containing $w$ intersects $B(a,2r_a)$, and 
\beq\label{220180824}
f^j(w)\notin \bu_{c\in F}U_0(c)
\eeq
for all $j=1,2,\ld, n-1$, then 
\beq\label{320180824}
|(f^{n})'(z)|\ge \frac{8\kappa}{\kappa-1}\cdot \frac{r_b}{r_a}
\eeq
for all $z\in f_w^{-n}(B(b,2r_b))$. 
\end{enumerate}

\, \, 

\item For every $b \in F_0$ and every $j \in \{1,2, \ldots, p_b\}$ there exists, $V_j(b)$, an open neighborhood of $b$ in $Y$ such that:

\,

\begin{itemize}
\item [(a)] The set $\overline{V_j(b)}\cap \overline{U_0(b, j)}$ is a closed topological disk,

\,

\item [(b)] The set $V_j(b)\sms \overline{U_0(b, j)}$ is  connected (and, in consequence, the set $\overline{V_j(b)}\sms \overline{U_0(b, j)}$ is connected too), 

\,

\item [(c)] $\overline{V_j(b)}\sbt B(b, r_b'/4)$ with some $r_b'\in(0, (r_b)_b^-)$.
\end{itemize}

\, \, 

\item $|f'(b)|=1$ for every $b\in F_0$.

\, \, 

\item  
$$ 
\Crit(f)\cap \bigcup_{w \in F_0} B(w, 6r_w) =\es.
$$

\, \, 

\item $\forall b\in F_0 \, \, \, \, \om(\Crit(f))\cap B(b,6r_b) \sbt \{b\}$.
\end{enumerate}

\sp Then there exists a pre--nice set \index{(N)}{pre--nice sets} $U=U_{\bold r}$ (although we do not indicate it here this set does depend also at least on $F$ and the sets $U_0(\xi)$, $\xi\in F_1$, and $U_0(\xi,i)$, $\xi\in F_0$, $j=1,2,\ld, p_\xi$), with the following additional properties:
\begin{itemize}
\item [(A)] For all $b\in F$:
$$ 
U_0(b)\sbt U  \sbt \bu_{b\in F}B(b,\kappa r_b)
\sbt\bu_{b\in F}B(b,2\kappa r_b),
$$
\item [(B)] For all $b\in F_1$:
$$ 
U_0(b)\sbt U\sbt \bu_{b\in F}B(b,\kappa r_b)
\sbt\bu_{b\in F}B(b,2\kappa r_b)
\sbt\bu_{b\in F}B(b,3\kappa r_b)
\sbt Y\sms  \ov{{\rm PS}(f)},
$$
\item [(C)] If $W$ is a connected component of $U$, then there exist either a unique $b\in F$ and a unique $j\in\{1,2,\ld, p_b\}$ such that
$$
U_0(b,j)\sbt W.
$$ 
We then denote $W$ by $U(b,j)$.

\,

\item [(D)] 
$$ 
U(b,j)\cap B(b, r_b')=U_0(b,j)\cap B(b, r_b')
$$
for all $b\in F_0$ and all $j\in \{1,2, \ldots, p_b\}$.

\sp

\item [(E)] 
$$ 
f^{-1}_b( U(b,j))\sbt U(b,j)\sbt B(b,r_b)
$$
for all $b\in F_0$ and all $j\in \{1,2, \ldots, p_b\}$.

\sp

\item [(F)] For all $b\in F_0$ and all $j\in \{1,2, \ldots, p_b\}$ there exists an open connected set $W(b, j)$ such that 
\beq\label{220190204}
\ov{U(b,j)}\sms\{b\}\sbt W(b,j)
\eeq
and for every integer $n\geq 0$ and every $z\in f^{-n}(b)$, the holomorphic branch 
$$
f_z^{-n}: W(b,j)\lra X 
$$
of $f^{-n}$, sending $f^n(z)$ to $z$, is well defined.

\sp

\item [(G)] For every $b\in F_0$, every $j\in\{1,2,\ld, p_b\}$, and every $u\in(0,r_b'/4)$ small enough there exist $V_j(b,u)\sbt V_j(b)$, an open neighborhood of $b$, and an open connected set $W(b,j;u)\sbt W(b,j)$ such that
\beq\label{2fp29.3.c}
 \overline{U(b, j)}\sms V_j(b,u)\sbt W(b,j;u).
\eeq
and the maps 
$$
f^{-n}_b:W(b,j;u)\lra X, \  \  n\ge 0,
$$ 
converge uniformly to the constant function whose range is equal to $\{b\}$
\end{itemize}
\ethm

{\sl Proof.}
For all $w \in F$ and all $1\le i\le p_w$, define $(U_n(w,i))_{0}^{\infty}$, an ascending sequence of open connected sets, as follows:
For all integers $n\geq 1$, the set $U_n(w,i)$ is defined to be the connected component of the union
$$
U_0(w,i)\cup\bigcup_{\xi\in F}\bigcup_{k=0}^n\bigcup_{V\in\Comp_k(U_0(\xi,i))}V 
$$  
that contains $w$. Set
$$
U_0(w):=\bu_{i=1}^{p_w}U_0(w,i).
$$
We will  prove by induction  that
\beq\label{1j299}
U_n(w) \sbt B(w,\kappa r_w)
\eeq
for all $w \in F$ and all $n\geq 0$. For $n=0$ this is immediate as, by hypothesis (1), 
$$
U_0(w)\sbt B(w,r_w)\sbt B(w,\kappa r_w).
$$
For the inductive step suppose that (\ref{1j299}) is true for all $ 0\leq j <n$ with  same $n \geq 1$. Fix $w\in F$, $1\le i\le p_w$, and let $C$ be a connected component of 
$$
U_n(w,i)\sms \ov{B(w,r_w)}.
$$
Then 
\beq\label{120180905}
\ov{C} \cap \ov{B(w,r_w)}\neq \es
\eeq
and there exist $\xi^*\in F$, $s\in\{1,2,\ld,p_{\xi^*}\}$, a minimal integer
 $1\leq k \leq n$, and a connected component $W\in\Comp_k\(U_0(\xi^*,s)\)$ such that  
\beq\label{f1_2015_08_31}
W \cap C \neq \es. 
\eeq
By the definitions of $C$, $W$, and $k$, we have that
$$ 
C\sbt\bigcup_{d\in F}\bu_{l=1}^{p_d} \bigcup_{j=k}^n \bigcup_{V\in\Comp_j(U_0(d,l))}V,
$$
and therefore,
$$ 
f^k(C)\sbt f^k\left(\bigcup_{d\in F}\bu_{l=1}^{p_d} \bigcup_{j=k}^n \bigcup_{V\in\Comp_j(U_0(d,l))}V \right) \\
\sbt \bigcup_{d\in F}\bu_{l=1}^{p_d} \bigcup_{j=0}^{n-k} \bigcup_{V\in\Comp_j(U_0(d,l))}V.
$$
Since the set $f^k(C)$ is connected and since, by \eqref{f1_2015_08_31}, $f^k(C)\cap U_0(\xi^*,s)\ne\es$, we therefore conclude that  
\beq\label{f2_2015_08_31}
f^k(C)\sbt U_{n-k}(\xi^*,s).
\eeq
Therefore, since $0\le n-k <n$, the inductive assumption yields
\beq\label{2j299}
f^k(C)
\sbt U_{n-k}(\xi^*,s)
\sbt B(\xi^*,\kappa r_{\xi^*}) 
\sbt B(\xi^*,2r_{\xi^*}).
\eeq
But then, $W\cup C$ is a connected set (by \eqref{f1_2015_08_31}) with   
\beq\label{3j299}
f^k(W\cup C) 
\sbt U_0(\xi^*,s)\cup U_{n-k}(\xi^*,s)
=U_{n-k}(\xi^*,s)
\sbt B\(\xi^*,\kappa r_{\xi^*}\)
\sbt B(\xi^*,2r_{\xi^*}).
\eeq
Consider two cases. First assume that either 

\begin{itemize}
\, \item $\xi^*\in F_1$ or 

\, \item $k\ge 2$ or 

\, \item $k=1$, $\xi^*\in F_0$ and $w\ne\xi^*$. 
\end{itemize}

\, 

\fr Then, as, by minimality of $k$, formula \eqref{220180824} holds for all $z\in C$ (so \eqref{220180824B} also holds), it follows from hypothesis (5) that the holomorphic inverse branch 
$$
f^{-k}_W:B(\xi^*,6r_{\xi^*})\lra X,
$$
satisfying $f^{-k}_W(B(\xi^*,6r_{\xi^*})\)\spt W$, is well defined. Consequently, using also the fact that $W\cup C$ is connected, we get that
\beq\label{3j299C}
W\cup C\sbt f^{-k}_W\(B(\xi^*,2r_{\xi^*})\),
\eeq
and. Also, with the help of \eqref{120180905}, we have that
$$
B(w,2r_w)\cap f^{-k}_W\(B(\xi^*,2r_{\xi^*})\)\ne\es.
$$
Therefore, we get from hypothesis (5b) that
$$
\big|(f^{-k}_W)'(x)\big|
\le \frac{\kappa-1}{8\kappa}\cdot \frac{r_w}{r_{\xi^*}}
\le \frac{\kappa-1}{2\kappa}\cdot \frac{r_w}{r_{\xi^*}}
$$ 
for all $x\in B\(\xi^*,\kappa r_{\xi^*}\)$. Consequently, recalling also (\ref{3j299C}), we get that
$$ 
\begin{aligned}
\diam (C) 
& \leq \diam\(f^{-k}_W \(B (\xi^*, \kappa r_{\xi^*})\)\)
< \frac{\kappa-1}{8\kappa}\cdot \frac{r_w}{r_{\xi^*}}\diam\(B(\xi^*,\kappa r_{\xi^*})\)\\
& = \frac{\kappa-1}{2 \kappa}\cdot \frac{r_w}{r_{\xi^*}}\cdot 2 \kappa r_{\xi^*} 
= \frac{\kappa-1}4 r_w.
\end{aligned} 
$$
Because of \eqref{120180905}, we thus conclude that
\beq\label{120180907} 
C\sbt \ov B\lt(w,r_w+\frac{\kappa-1}4r_w\rt)
=\ov B\lt(w,\frac{\kappa+3}4r_w\rt)
\sbt B(w,\kappa r_w).
\eeq
Suppose now that the remining case holds, meaning that $\xi^*\in F_0$, $w=\xi^*$ and $k=1$. It then follows from \eqref{3j299} that 
\beq\label{520180907}
f(W\cup C)\sbt B(w,2r_{w}).
\eeq
So, 
\beq\label{620180907}
f(\ov W\cup \ov C)\sbt B(w,3r_{w}).
\eeq
Since, by, \eqref{f1_2015_08_31}, the set $\ov W\cup \ov C$ is connected, since, by \eqref{120180905}, $\ov W\cup \ov C$ intersects $\ov B(w,r_w)$, since, by (2b), the map $f|_{B(w,6r_w)}$ is 1--to--1, and since, by (2c), $f_w^{-1}(B(w,3r_w))\sbt B(w,6r_w)$, we conclude, with the help of \eqref{620180907}, that 
$$
W\sbt \ov W\cup \ov C \sbt f_w^{-1}(B(w,3r_{w}).
$$
But, since $W\in\Comp_1(U_0(w,s))$ and $U_0(w,s)\sbt B(w,r_{w})\sbt B(w,3r_{w})$, we thus conclude that 
$$
W=f_w^{-1}(U_0(w,s)).
$$
So, by (2f) and by (1): 
$$
W\sbt U_0(w,s)\sbt B(w,r_w).
$$
This however contradics \eqref{f1_2015_08_31} and the very definiition of $C$. Hence, the considered case is ruled out and \eqref{120180907} holds.

Thus, taking the union over all $i\in\{1,2,\ld, p_w\}$ and over all connected components of $U_n(w,i)\sms \ov{B(w,r_w)}$ along with $U_0(w)\sbt B(w,r_w)$, we thus see that the  inductive proof of (\ref{1j299}) is complete. 

Now, for every $w \in F$ and all $1\le i\le p_w$, define
$$
U'(w,i):=\bigcup_{n=0}^\infty U_n(w,i).
$$
From our construction and from (\ref{1j299}), $U'(w)$ is thus an open  connected set such that
\beq\label{220180825}
B(w,r_w) \sbt U'(w,i) \sbt B(w,\kappa r_w)
\eeq
for all $w\in F$ and all $1\le i\le p_w$. 

Let $\Ga(w,j)$ be the connected component of $\ov B(w,4r_w)\sms U'(w,i)$ containing $\bd B(w,4r_w)$. It immediately follows from \eqref{220180825} that
\beq\label{120181008}
\Ga(w,j) \spt \ov B(w,4r_w)\sms B(w,\kappa r_w).
\eeq
Let 
\beq\label{220181008}
U(w,i):=\ov B(w,4r_w)\sms \Ga(w,j).
\eeq
Then
\beq\label{1j301}
U_0(w,i)\sbt U'(w,i) \sbt U(w,i) \sbt B(w,\kappa r_w),
\eeq
and we have:

\bobs\label{o120181008}
For every $w \in F$ and all $1\le i\le p_w$, the following hold:
\begin{enumerate}
\item [(1*)] For every $\g\in[\ka r_w,6r_w]$,
$$
U(w,i)
=\ov B(w,\g)\sms \Ga(w,j)
=B(w,\g)\sms \Ga(w,j).
$$
\item[(2*)]  $U(w,i)$ is an open subset of $B(w,\kappa \g)$.

\, \item[(3*)] The set $U(w,i)$ is connected.

\, \item[(4*)] The set $U(w,i)$ is simply connected.

\, \item[(5*)] $\bd U(w,i)\sbt \bd U'(w,i)$.

\end{enumerate}
\eobs

\bpf Item (1*) is obvious. Item (2*) follows immediately from (1*) and the fact that the set $\Ga(w,j)$ is closed. Item (3*) is an immediate consequence of the (very general) Theorem~5, page 140, in \cite{Kur2}. Item (4*) follows from Theorem~4.4, page 144, in \cite{Newman}. Because of (1*) and Theorem~3, page 238 in \cite{Kur2} (which applies since $Y$ is locally connected), we have that
$$
\begin{aligned}
\bd U(w,i)
&\sbt \ov B(w,2r_w)\cap\(\bd \Ga(w,j)\cup \bd B(w,6r_w)\)
=\ov B(w,2r_w)\cap\bd \Ga(w,j) \\
&\sbt \ov B(w,2r_w)\cap\bd\(B(w,4r_w)\sms U'(w,j)\)
 \sbt \ov B(w,2r_w)\cap(\bd B(w,4r_w)\cup \bd U'(w,j)\) \\
&=\ov B(w,2r_w)\cap\bd U'(w,j) \\
&=\bd U'(w,j).
\end{aligned}
$$
\epf

For every $b\in F$ set
$$
U(b):=\bu_{i=1}^{p_b}U (p,i)
$$
and
$$
U:=\bigcup_{w \in F}U(w)=\bigcup_{w \in F}\bu_{i=1}^{p_w}U (w,i).
$$
Our goal now is to show that for each $b\in F_0$ the sets

\sp (*) \!\!\!\!\!\!\!\!\!\!\!\! \centerline{$U(b,i)$, $1\le i\le p_b$, are pairwise disjoint.}

\sp\fr Having this done, Property (c) of Definition~\ref{d1j291-} follows from (4*) while its property (a) follows from the inclusion $\ov{U(w,i)} \sbt \ov{B(w, \kappa r)}$ and since the latter set is compact. Properties (A) and (B) of Theorem~\ref{t2j291} directly follow from (\ref{1j301}) and hypothesis (4) of this theorem. Property (C) of Theorem~\ref{t2j291} also holds since all the sets $U(w,i)$, $w \in F$, $1\le i\le p_w$, are connected and simply connected, since 
$$
U(w)\sbt B(w,\kappa r_w )\sbt B(w,2r_w),
$$
and since the sets $B(w,2r_w)$, $w\in F$, are pairwise disjoint. 

\sp\fr So, in order to have all the claims of the preceding paragraph established, we now focus on proving (*). Toward this end, we shall first prove the following:

\sp\fr {\sc Claim~$1^0$.}  For every $b\in F_0$ and every $ i \in \{1,2, \ldots, p_b\}$  all holomorphic  iterates
$$ f^{-n}_b: U'(b, i)\lra X, $$
$n \geq 0$, are well--defined, and
$$f^{-n}_b(U'(b, i))\sbt U'(b, i) \sbt B(b, \kappa r_b).$$

{\sl Proof.} We proceed by induction. For $n=0$, the claim follows immediately  from formula~\ref{220180825}. Suppose it holds for some $n \geq 0$. Then,  by (2c) (and as $\kappa \leq 2$),
$$ f^{-(n+1)}_b:=f^{-1}_b\circ f^{-n}_b: U'(b,i)\lra X$$ is well--defined. Since $U_0(b,i) \sbt U'(b,i)$ and since $f^{-1}_b(U_0(b,i))\sbt U_0(b,i)$, we  conclude from the construction
of the  set $U'(b,i)$ that  $ f^{-(n+1)}_b(U'(b,i))\sbt U'(b,i)$. Since also, by (2c), $U'(b,i)\sbt B(b, \kappa r_b)$, the proof of Claim~$1^0$ is complete.

\sp\fr As an immediate consequence of this claim, connectedness of the sets $U'(b,i)$, and the hypothesis (2g), we get the following:

\sp\fr {\sc Claim~$2^0$.} For every $b \in F_0$ and every  $ i \in \{1,2, \ldots, p_b\}$, the sequence $$\left({f^{-n}_b|}_{U'(b,i)} : U'(b, i) \lra  U'(b,i)\right)_{n=0}^\infty$$
converges uniformly on compact subsets of $U'(b, i)$ to  the constant function whose range is equal to $\{b\}$.

Now we are in position to prove that the property (*) holds. Indeed, because of Claim~$1^0$, we  have that
  $$ f^{-1}_b( U'(b, i))\sbt U'(b,i)\sbt U(b,i).$$
  So, since $U(b,i)$ is conformally equivalent to the unit disc $\mathbb{D}$, it follows from Proposition~\ref{p2j144} and (\ref{1j301}) that
  \beq\label{1fp28}
    f^{-1}_b( U(b, i))\sbt U(b,i)\sbt B(b,r_b).
  \eeq
So, item (E) is proved. Furthermore, all iterates
   $$f^{-n}_b: U(b, i)\lra U(b,i), \quad n \geq 0,$$
are well-defined and form a normal family. It then follows  from the hypothesis (2g) that:

\sp\fr {\sc Claim~$3^0$.}  For every $b \in F_0$ and every  $ i \in \{1,2, \ldots, p_b\}$, the sequence 
$$
\left({f^{-n}_b|}_{U(b,i)} : U(b, i) \lra  U(b,i)\right)_{n=0}^\infty
$$
converges uniformly on compact subsets of $U(b, i)$ to  the constant function whose range is equal to $\{b\}$.

\sp We conclude from this claim and property (2h) that for every $z \in U(b,i)$ there exists $n \geq 0$  such that
  $$f^{-n}_b(z) \in \bigcup_{j=1}^{p_b} U_0(b,j).$$
So, if for every $j \in \{1,2, \ldots, p_b\}$
$$ U_\infty(b, i; j):=\big\{ z \in U(b,i):\,  \exists n \geq 0 \, \,  f^{-n}_b(z)\in U_0(b,j)\big\},$$
then
\beq\label{1fp29}
 \bigcup_{j=1}^{p_b} U_\infty (b,i;j):= U(b,i).
  \eeq
Obviously, all the sets $U_\infty (b,i;j), \, j=1, \ldots, p_b$, are open, and, because of (2e) and (2f), there  are pairwise disjoint. Since the set $U(b,i)$ is connected, we conclude that only one term in the union of (\ref{1fp29}) is non--empty. Since $U_\infty(b, i; i) \supset U_0(b, i) \neq \es $, we further  conclude that
$$ 
U(b,i)= U_\infty(b,i;i).
$$

Since by (2e), applied again, the sets $\{U_\infty(b,i,;i)\}, \, i=1,\ldots, p_b$, are  pairwise disjoint, the property (*) follows.

\sp We thus conclude that the sets $U(w,i)$, $w \in F$, $i=1,2,\ld,p_b$, coincide with the collection of all connected components of $U$. Hence, the number of connected components of $U$ is finite, which means that condition (b) of Definition~\ref{d1j291-} holds.

\sp We shall now show  that condition (d) of Definition~\ref{d1j291-} holds.  Indeed, seeking contradiction suppose that 
$$
U \cap f^n(\partial{U})\neq \es
$$    
with some $n \geq 1$. Consider a minimal $n\ge 1$ with this property. Then  there exists 
$$
x\in\partial{U} 
\  \  \ {\rm such \ that} \  \  \ 
f^n(x)\in  U.
$$
Assume first that  
$$
f^n(x)\in U',
$$
where 
$$
U':=\bigcup_{w\in F}\bu_{i=1}^{p_w}U'(w,i),
$$ 
and assume only that $x \in \partial{U'}$ (this is a weaker requirement than $x\in\bd U$ as $\partial{U} \sbt \partial{U'}$). Therefore there exist $\xi\in F$, $i\in\{1,2,\ld p_{\xi}\}$, an integer $k\geq 0$, and $W \in  \Comp_k\(U_0(\xi,i)\)$ such that $f^n(x) \in W$. But then $ x \in \tilde{W}$, the connected component of $f^{-n}(W) $ containing $x$. We also immediately see that 
$$ 
\tilde{W} \in \Comp_{n+k}\(U_0(\xi,i)\).
$$
Furthemore  $ \tilde{W} \cap U'\neq \es$ as $ x \in \partial{U}'$. By our construction  of the set $U'$, this implies that $ \tilde{W} \sbt U'$. Hence $x \in U'$  and this contradiction  rules out the considered case.

\sp Now, continuing the general case, there exist $w\in F$ and $1\le i\le p_w$ such that 
$$
f^n(x)\in U(w,i).
$$
It therefore follows from \eqref{1j301}, the formula $\ka\le 2$, and minimality of $n$ that formula \eqref{220180824B} holds. Hence, by hypothesis (5a), there exists $f^{-n}_x:B(w,6r_w) \lra X$, a unique holomorphic inverse branch of $f^n$ defined on $B(w,6r_w)$ and sending $f^n(x)$ back to $x$. If 
$$
f^{-n}_x(U'(w,i))\cap U'(y,j)\neq \es
$$
for some $y \in F$ and $j\in\{1,2,\ld, p_y\}$, then by the previous  case, 
$$
f^{-n}_x (U'(w,i))\sbt  U'(y,j)\sbt U(y,j).
$$ 
By virtue of Proposition~\ref{p2j144}, we thus then get $f^{-n}_x(U(w,i))\sbt U(y,j)$. In particular $x =f^{-n}_x(f^n(x))\in U(y,j)\sbt U$. This contradiction yields
\beq\label{1j305}
U'\cap f^{-n}_x (U'(w,i))=\es.
\eeq
But
\beq\label{120180929}
x\in\partial{U}(\xi,k)\cap f^{-n}_x(U(w,i)) 
\sbt \partial{U}'(\xi,k)\cap f^{-n}_x(U(w,i))
\eeq
for some $\xi\in F$ and $k\in\{1,2,\ld, p_\xi\}$, where the inclusion part of this formula holds because of Observation~\ref{o120181008}(5*). Hence,
$$
U'(\xi,k) \cap f^{-n}_x (U(w,i))\neq \es.
$$
Because of Observation~\ref{o120181008}(5*), applied this time to $U(w,i)$, it follows from this and \eqref{1j305} that 
\beq\label{2j305}
U'(\xi,k) \sbt f^{-n}_x (U(w,i)).
\eeq
But then, using \eqref{1j301}, we get that
$$
f^{-n}_x(B(w, \kappa r_w))\cap B(\xi,\kappa r_\xi)\neq \es.
$$
Therefore, remembering also about minimality of $n$, we see that hypothesis (5) of our theorem is satisfied, and so, it follows from \eqref{320180824} with the help of hypothesis (1) and \eqref{1j301}, that 
\beq\label{3j305}
\diam\(f^{-n}_x(U(w,i))\)
\le \frac{\kappa-1}{8\kappa}\cdot\frac{r_\xi}{r_w}\diam (U(w,i))
\le \frac{\kappa-1}{8\kappa}\cdot \frac{r_\xi}{r_w} 2\kappa r_w
=\frac{\kappa-1}{4}r_\xi
\le r_\xi/4.
\eeq
On the other hand, if $\xi\in F_1$, then it follows from (\ref{2j305}) and the hypothesis (2h) that 
$$
\diam\(f^{-n}_x(U(w,i))\)
\ge \diam(U'(\xi,k))
\ge r_\xi/2.
$$
This however contradicts \eqref{3j305} and finishes the proof of condition (d) of Definition~\ref{d1j291-} in the case when $\xi\in F_1$. 

\sp So, suppose that $\xi\in F_0$. It then follows from \eqref{2j305} and (2d) that 
$$
\xi\in \ov{U'(\xi,k)}\sbt f^{-n}_x\(\ov{U(w,i)}\).
$$
Then, by  virtue of (2a) we get that
       $$\xi=f^n(\xi)\in f^n(f^{-n}_x(\overline{U(w, i)}))=\overline{U(w, i)}.$$
It thus follows from the hypothesis (3) and (\ref{1j301}) (remember also that $\kappa \leq 2$) that $\xi=w$. So, $ w \in f^{-n}_x(\overline{U(w,i)})$, and  as $f^n(w)=w$ and $w \in\overline{U(w,i)} \sbt B(w,2r_w)$, we  conclude that $f^{-n}_x(w)=w$, and $f^{-n}_x=f^{-n}_w$.  We  deduce from this, (\ref{2j305}), the first  inclusion of (\ref{1fp28}), and the property (*), that $k=i$. It eventually follows from (\ref{120180929}) and (\ref{1fp28}) that $x\in \partial{U}(w, i)\cap U(w,i)$. Since the set $U(w,i)$ is open this is a  contradiction and the proof of condition (d) of Definition~\ref{d1j291-} is complete.

\sp We conclude the proof of this theorem with proving item (D). We start with the following.

\sp\fr {\sc Claim~$4^0$.} For all $b \in F_0$ and all $j\in\{1,2, \ldots, p_b\}$ we have that
$$ 
U'(b, j) \cap B(b, r_b') = U_0(b, j)\cap B(b,r_b').
$$

\sp\bpf Of course 
$$
U_0(b, j) \cap B(b, r_b') \sbt U'(b, j)\cap B(b,r_b').
$$
In order  to prove the opposite inclusion, take any point
$$
z \in U'(b, j) \cap B(b,r_b').
$$
If $ z \notin U_0(b,j)$, then $f^n(z)\in B(b, r_b)$ for all $ n \geq 0$ by virtue of condition (*), hypothesis (2j) and since $r_b'\le (r_b)_b^-$.  Since also $z \in U'(b, j)$, by using property (*) and Claim~$1^0$, we thus conclude that $f^k(z)\in U_0(B,j)$ for some $k\geq 0$. It then follows from Clam~1, used again, that
$$z=f^{-k}_b(f^k(z))\in f^{-k}_b(U_0(b, j))\sbt U_0(b, j)).$$
 This contradiction shows that $z \in U_0(b,j)$, finishing the proof of Claim~$4^0$.
 
\sp Fix $b\in F_0$ and $ j \in \{1,2, \ldots, p_b\}$. By our hypothesis (6) the set  $V_j(b)\sms \overline{U_0(b, j)}$ is connected.

\sp Seeking contradiction suppose that $$G:=U(b,j)\cap ( V_j(b)\sms \overline{U_0(b, j)})\neq \emptyset.$$
Then
\beq\label{1fp29.2}
 \partial{G}\sbt \overline{V_j(b)}\cap (\partial{U}(b, j)\cup \partial{U}_0(b,j)).
  \eeq
So, using item (5*) of Observation~\ref{o120181008} and Claim~$4^0$, we get
 $$\partial{G}\sbt \overline{V_j(b)}\cap (\partial{U}'(b, j)\cup \partial{U}_0(b,j))\sbt ( \overline{V_j(b)}\cap \overline{U'(b,j)})\cup \overline{U_0(b,j)}=\overline{U_0(b,j)}.$$
 Hence,
 $$\partial{}_{V_j(b)\sms \overline{U_0(b,j)}}G\sbt  (V_j(b)\sms \overline{U_0(b,j)})\cap \partial{G}=\es.
$$
Thus, because of (6b),  
\beq\label{120290204} 
{\rm either} \quad  G=\es\quad \text{{\rm or}} \quad G=V_j(b)\sms \overline{U_0(b,j)}.
\eeq
 If $G=V_j(b)\sms \overline{U_0(b,j)}$, then $V_j(b)\sms \overline{U_0(b,j)}\sbt U(b,j)$ and therefore, $V_j(b)\sbt \overline{U(b,j)}$. So, if $s_1>0$ is so small that $B(b,s_1)\sbt V_j(b)$, then $$\overline{B(b, s_1)}\sbt \Int(\overline{U(b, j)}) =U(b,j).$$
It therefore follows from Claim~$3^0$ that
$$f^{-q}_b(\overline{B(b, s_1)})\sbt B(b, s_1/2)$$
for all (just one  suffices) $q\geq 1$ large enough. But then the Schwarz Lemma yields $|(f^{-q}_b)'(b)|< 1$. Therefore,  $|f'(b)|> 1$. This  however contradicts our hypothesis (7) and proves that $G =\es$. Therefore, keeping also in mind  that the set $U_0(b,j)$ is open and $\overline{U_0(b,j)}$ is a closed topological disk, we get
$$\begin{aligned}
U(b, j)\cap (V_j(b)\sms U_0(b, j)) &= U(b,j)\cap (V_j(b)\cap \partial{U}_0(b,j))\\
                              &=U(b,j)\cap (V_j(b)\cap\partial{U}_0(b,j))\\
                              & \sbt U(b,j)\cap(V_j(b)\cap  \overline{(Y \sms \overline{U_0(b,j)})}).
                              \end{aligned}$$
 Seeking contradiction, suppose that this set is not empty. Since $U(b,j)\cap V_j(b)$ is open, we would  then have that
 $$ U(b, j)\cap V_j(b)\cap (Y\sms \overline{U_0(b,j)})\neq \es.$$
 Equivalently
 $$ U(b, j)\cap (V_j(b)\sms  \overline{U_0(b,j)})\neq \es,$$
  meaning that $G\neq \es$. This contradiction shows that
 $$U(b,j)\cap(V_j(b)\sms U_0(b, j))=\es,$$
 yielding
 \beq\label{1fp29.3.0a}
 U(b, j)\cap V_j(b)=U_0(b,j)\cap V_j(b).
  \eeq
Since $V_j(b)$ is an open set containing $b$, there exists $\b>0$ such that $B(b,\b)\sbt V_j(b).$
Consequently, for all $\a \in (0, \b]$, we  have that
\beq\label{1fp29.3}
U(b, j)\cap B(b,\a)= U_0(b, j)\cap B(b,\a)
\eeq
and the proof of item D of Theorem~\ref{t2j291} is complete.

\sp Now, assuming condition (8) and (9) we shall prove (F) and (G). 
Let $\xi \in \partial{U}(b, j)\sms \{b\}$ for some $b \in F_0$ and $j \in \{1, 2, \ldots, p_b\}.$
 Fix some $s>0$ so small that
 $B(\xi, s) \sbt B(b, 3r_b)$. Then, because of condition (9) there exists an integer $N_\xi \geq 0$ such that
 $$\lt(\bigcup_{n=N_\xi}^\infty f^n(\Crit(f))\rt)\cap B(\xi, s)=\es.$$
 Equivalently,
 \beq\label{1fp29.3.a}
 \Crit(f) \cap  \bigcup_{n=N_\xi}^\infty f^{-n}(B(\xi, s))=\es.
\eeq
It immediately follows from (8), (\ref{1fp28}), and continuity of $f^{-1}_b $ that for every $n \geq 0$ there exists $r_n(\xi)\in (0, +\infty)$ such that
\beq\label{2fp29.3.a}
\Crit(f) \cap  \Comp(f^{-n}_b(\xi), k, 8r_n(\xi))=\es
\eeq
for every $k \in \{0,1, \ldots, n\}$. Put
$$ r_\xi:=\min \{r_{N_\xi}, s/8\}.$$
It then follows from (\ref{1fp29.3.a}) and (\ref{2fp29.3.a}) that
$$\Crit (f)\cap \bigcup_{n=0}^\infty \Comp ( f^{-n}_b(\xi), n, 8r_\xi)=\es.$$
Hence, for all $n \geq 0$ we have that
$$\Crit(f^n)\cap \Comp ( f^{-n}_b(\xi), n, 8r_\xi)=\es, $$
and therefore there exists a unique holomorphic branch $ f^{-n}_{b,\xi}: B(\xi, 8r_\xi)\lra X$ of $f^{-n}$ of $f^n$ sending  $\xi$ to $f_b^{-n}(\xi)$. So, taking
$$
W(b,j):=U(b,j)\cup\bigcup_{\xi\in \partial{U}(b,j)\sms\{b\}} B(\xi, 8r_\xi)
$$
and invoking also condition (5a), we conclude that item (F) holds.

\

We now shall prove the following.

\sp\fr {\sc Claim~$5^0$}. For every $\xi  \in \partial{U}(b,j) \sms \{b\}$, the sequence $ {f^{-n}_{b,\xi}}|_{B(\xi, 4r_\xi)}:B(\xi,4r_\xi) \lra X$
converges uniformly to the constant function whose  range is equal to $\{b\}$.

\sp\fr {\sl Proof.} Since $ \xi \in \partial{U}(b,j)$, we have that $B(\xi, r_\xi)\cap U(b,j)\neq \es$. Fix a point $y$ in this  intersection. Take $\rho>0$ so small that
$$
B(y,8\rho)\sbt B(\xi,r_\xi)\cap U(b,j).
$$
By virtue of Claim~$3^0$, the sequence $ f^{-n}_{b,\xi}|_{B(y, 4\rho)}:B(y,4\rho) \lra X$ converges uniformly to the constant function whose range is equal to $\{b\}$. Consequently, 
$$ 
\lim_{n \to \infty}\big\|({f^{-n}_{b,\xi}}|_{B(y,4\rho)})'\big\|_{\infty}=0.
$$
It follows from this and Koebe's Distortion Theorem that
$$ 
\lim_{n \to \infty}\big\|({f^{-n}_{b,\xi}}|_{B(\xi,4r_\xi)})'\big\|_{\infty}=0. 
$$
Hence  
$$ 
\lim_{n \to \infty}\diam (f^{-n}_{b,\xi}\( B(\xi,4r_\xi))\)=0.
$$
and the claim follows. \qed

\sp By condition (6), for every $u\in (0, r_b')$, small enough, there exists  an open set  $V_j(b,u)\sbt B(b,u)\cap V_j(b)$, containing $b$, with the following properties:
\begin{itemize}
\item [(1*)]  The set $\overline{V_j(b,u)} \cap \overline{U_0(b,j)} $ is a closed topological disk,

\,

\item [(2*)]  The set $V_j(b,u) \sms \overline{U_0(b,j)} $ is  connected,

\,

\item [(3*)]  $V_j(b,u) \sbt B(b, r_b'/4)$,

\,

\item [(4*)]  $\partial{ }\left(\overline{V_j(b,u)} \cap \overline{U_0(b,j)}\right)=\left(\overline{V_j(b,u)} \cap \partial{U}_0(b,j)\right)\cup \left(\overline{U_0(b,u)} \cap \partial{V}_j(b,j)\right)$,

\,

\item [(5*)] Both  sets $\overline{V_j(b,u)} \cap \partial{U}_0(b,j)$ and $\overline{U_0(b, j)} \cap \partial{V}_j(b,u)$ are compact topological arcs,

\,

\item [(6*)] The intersection  
$$
\left(\overline{V_j(b,u)} \cap \partial{U}_0(b,j)\right) \cap \left(\overline{U_0(b,u)} \cap \partial{V}_j(b,j)\right)
$$
consisits of two distinct points, which we denote by $x_j(b,u)$  and $y_j(b,u)$.
\end{itemize}

\sp We  shall prove the following.

\sp\fr {\sc Claim~$6^0$}. The set $\overline{U(b, j)} \sms V_j(b,u)$  is connected.

\sp\fr {\sl Proof.} Fix two points $w,z\in U(b,j)\sms V_j(b,u)$. Since $U(b,j)$ is an open connected subset of the Riemann surface $Y$, it is arcwise connected. Therefore, there  exists a homeomorphic embedding 
$$
\gamma: [0,1]\to  U(b,j)
$$ 
such that $\gamma(0)=w$ and $\gamma(1)=z$. If $\gamma([0,1])\sbt U(b,j)\sms V_j(b,u)$, we are done. So, suppose that 
$$
\gamma([0,1]) \cap V_j(b,u)\neq \es.
$$
But, because of (\ref{1fp29.3.0a}),
\beq\label{1fp29.3.c1}
 \gamma([0,1])\cap V_j(b,u)\sbt  U_0(b,j).
 \eeq
Denote:
$$
t_w:= \inf\{t\in[0,1]: \, \, \gamma(t) \in V_j(b,u)\}
$$ 
and  
$$
t_z:= \sup\{t\in[0,1]: \, \, \gamma(t) \in V_j(b,u)\}.
$$
Then, keeping also in mind that the set $V_j(b,u)$ is open, we have that
$$ 
\gamma([0, t_w])\cap V_j(b,u)=\es
\  \  \  {\rm and}\  \  \  
\gamma([t_z, 1])\cap V_j(b,u)=\es,
$$
and, because of (\ref{1fp29.3.c1}), 
$$
\gamma(t_w), \gamma(t_z)\in \overline{U_0(b,j)}.
$$
We have also that
$$
\gamma(t_w), \gamma(t_z)\in \partial{V}_j(b,u).
$$
In conclusion 
$$
\gamma(t_w), \gamma(t_z)\in \overline{U_0(b,j)}\cap \partial{ V}_j(b,u).
$$
Therefore, by (5*), there exists a compact topological  arc $\Delta\sbt \overline{U_0(b,j)}\cap \partial{V}_j(b,u)$, with $\gamma(t_w)$ and  $\gamma(t_z)$ being its endpoints, which is entirely contained  in  $U_0(b,j)\cap \partial{V}_j(b,u)$ except, perhaps, the endpoints $\gamma(t_w)$ and  $\gamma(t_z)$. But both $\gamma(t_w)$ and  $\gamma(t_z)$ belong to $U(b,j)$. Thus,
$$ \Delta \sbt U(b,j)\sms V_j(b,u).$$  Hence
  $[\gamma(0),\gamma(t_w)]\cup \Delta \cup [\gamma(t_z), \gamma(1)] $
is a compact topological arc joining $\gamma(0)$ and $\gamma(1)$, entirely contained in  $U(b,j)\sms V_j(b,u)$. So, any two points in $U(b,j)\sms V_j(b,u)$ belong to some  connected subset of $U(b,j)\sms V_j(b,u)$. Therefore the set $U(b,j)\sms V_j(b,u)$ is also  connected. Since also
$$ \overline{U(b,j)}\sms V_j(b,u)\sbt \overline{U(b,j)\sms V_j(b,u)},$$
we thus conclude that the set $\overline{U(b,j)}\sms V_j(b,u)$ is also connected. The  proof of Claim~$6^0$ is complete.\qed

\sp Obviously,
\beq\label{1fp29.3.c}
\left((\overline{U(b, j)} \sms V_j(b,u)\right)\sms \bigcup_{\xi \in \partial{U}(b,j)\sms V_j(b,u)}B(\xi, 4r_\xi)\sbt U(b,j)\sms V_j(b,u)
=U_0(b,j)\sms V_j(b,u)
\eeq
where the equality sign was written due to (\ref{1fp29.3})
 and the  choice of $u$. Since $U_0(b,j)$ is an open set, for every $\xi\in U_0(b,j)\sms V_j(b,u)$ there exists $r_\xi>0$ so small that 
$$
B(\xi, 8r_\xi)\sbt U_0(b,j).
$$
Since by (\ref{1fp29.3.c}), the collection
$$
\big\{B(\xi,4r_\xi):\,  \xi\in\left( \partial{U}(b, j)\sms V_j(b,u)\right)\cup \left( U_0(b, j)\sms V_j(b,u)\right)\big\}
$$
is an open cover of the compact set $\overline{U(b, j)}\sms V_j(b,u)$, there exists a finite set
$$ 
Z \sbt \left(\partial{U}(b, j)\sms V_j(b,u)\right)\cup \left( U_0(b,j)\sms V_j(b,u) \right)
$$ 
such that
 \beq\label{2fp29.3.cB}
 \overline{U(b,j)}\sms V_j(b,u)\sbt W(b,j;u):=\bigcup_{\xi \in Z} B(\xi,4r_\xi).
 \eeq
Obviously, the set $W(b, j;u)$ is open. It is connected because by Claim~$6^0$, the set $\overline{U(b,j)}\sms V_j(b,u)$
is connected, each ball $B(\xi, 4r_\xi)$, $ \xi \in Z$, is connected, and each such ball intersects $\overline{U(b,j)}\sms V_j(b,u)$. Since  $W(b,j;u)\sbt W(b,j)$, for every $n \geq 0$ and every $z\in f^{-n}(b)$, the holomorphic branch $f^{-n}_z:W(b,j;u)\lra X$ of $f^{-n}$ is well defined. Also, the maps $f^{-n}_b:W(b,j;u)\lra X$ converge uniformly to the constant function whose range is equal to $\{b\}$  because of Claim~$5^0$, the hypothesis (2g), and since the set $Z$ is finite. The proof of (G) is complete. 

\sp\fr Thus, the proof of Theorem~\ref{t2j291} is complete. 
\qed

\sp We remark that the hypotheses about the behavior of $f$ around points $b\in F_0$ are modeled on those when $b$ is a parabolic fixed point. These will be addressed in depth and length in the following chapters. 

\sp We shall now derive several consequence of Theorem~\ref{t2j291}. We start with the following.

\sp\bcor\label{c120180827A}   
Let $Y$ be a complete Riemann surface with constant curvature $0$, i.e. $Y$ is either the complex plane $\C$, a complex torus $\mT_\La=\C/\La$, where $\La$ is a lattice on $\C$ or an infinite cylinder $\C/2\pi i\Z$.
 
 Fix $R>0$ so small that every open ball in $Y$, with radius $R$ is simply connected, i.e. isometrically equivalent $B(0,R)\sbt\C$. Fix $\l>1$ and $\kappa\in(1,2]$.

Let $X$ be a non-empty open subset of $Y$ and let $f:X \lra Y$ be an analytic map. 

\fr Fix $F$, a finite subset of $J(f)\sms\ov{{\rm PS}(f)}$, and a vector $\bold r=\(r_b:b\in F\)$ with the following properties:
\, \begin{enumerate}  
\item
$$
r_b\in \lt(0,\frac{1}{6}\min\Big\{R,\min\big\{\rho(b,c):c \in F\sms \{b\}\big\}, \dist\(b, \ov{\rm PS}(f)\)\Big\}\rt)
$$
for all $b\in F$.

\, \, \item 
If $a, b\in F$, $n\ge 1$ is an integer, and $w\in f^{-n}(B(b,2r_b))$ is such that the connected component of $f^{-n}(B(b,2r_b))$ containing $w$ intersects $B(a,2r_a)$, and 
\beq\label{220180824J}
f^j(w)\notin \bu_{c\in F}B(c,r_c)
\eeq
for all $j=1,2,\ld, n-1$, then 
\beq\label{320180824J}
|(f^{n})'(z)|\ge \max\lt\{\frac{8\kappa}{\kappa-1}\cdot \frac{r_b}{r_a},\l\rt\}
\eeq
for all $z\in f_w^{-n}(B(b,2r_b))$.  
\end{enumerate}

\sp\fr Then there exists a nice set \index{(N)}{nice sets} $U=U_{\bold r}$ with the following properties:
\begin{itemize}
\item [(a)] For all $b\in F$:
$$ 
B(b,r_b)\sbt U \sbt \bu_{b\in F}B(b,\kappa r_b)
\sbt B(b,3\kappa r_b)\sbt Y\sms \ov{{\rm PS}(f)},
$$
\item [(b)] If $W$ is a connected component of $U$, then $W \cap F$ is a singleton. 
\end{itemize}
\ecor

\bpf
Setting $F:=F_1$ and $U_0(b):=B(b,r_b)$ for all $b\in F$, this corollary is an immediate consequence of Theorem~\ref{t2j291} except perhaps condition (e) of Definition~\ref{d1j291}. But this one follows  immediately from formula \eqref{320180824J} of condition (2) of our corollary. 
\epf

\sp\fr As an immediate consequence of this corollary and Theorem~\ref{t1j290a}, we get the following. 

\sp\bthm\label{t2j290e}  
The system ${\mathcal S}_{U_{\bold r}}$ resulting from this corollary and Theorem~\ref{t1j290a} is a maximal conformal graph directed Markov systems in the  sense of Chapter~\ref{Markov-systems}.
\ethm

\sp\fr As an immediate consequence of Corollary~\ref{c120180827A} we get the following.

\sp\bcor\label{c120180827}   
Let $Y$ be a complete Riemann surface with constant curvature $0$, i.e. $Y$ is either the complex plane $\C$, a complex torus $\mT_\La=\C/\La$, where $\La$ is a lattice on $\C$ or an infinite cylinder $\C/2\pi i\Z$.

Fix $R>0$ so small that every open ball in $Y$, with radius $R$ is simply connected, i.e. isometrically equivalent $B(0,R)\sbt\C$. Fix $\l>1$ and $\kappa>1$.

Let $X$ be a non--empty open subset of $Y$ and let $f:X \lra Y$ be an analytic map. 

\fr Fix $F$, a finite subset of $J(f)\sms\ov{{\rm PS}(f)}$, and a radius $r>0$ with the following properties:
\, \begin{enumerate}  
\item
$$
r\in \lt(0,\frac{1}{6}\min\Big\{R,\min\big\{\rho(a,b):a,b \in F, a\ne b\}\big\}, \dist\(F, \ov{\rm PS}(f)\)\Big\}\rt)
$$

\, \, \item 
If $a, b\in F$, $n\ge 1$ is an integer, and $w\in f^{-n}(B(b,2r))$ is such that the connected component of $f^{-n}(B(b,2r))$ containing $w$ intersects $B(a,2r)$, and 
\beq\label{220180824A}
f^j(w)\notin \bu_{c\in F}B(c,r)
\eeq
for all $j=1,2,\ld, n-1$, then 
\beq\label{320180824A}
|(f^{n})'(z)|\ge \max\lt\{\frac{8\kappa}{\kappa-1},\l\rt\}
\eeq
for all $z\in f_w^{-n}(B(b,2r))$. 
\end{enumerate}

\sp Then there exists a nice set $U=U_r$ with the following properties:
\begin{itemize}
\item [(a)] For all $b\in F$:
$$ 
B(b,r)\sbt U  \sbt \bu_{b\in F}B(b,\kappa r)
\sbt B(b,3\kappa r)\sbt Y\sms \ov{{\rm PS}(f)},
$$
\item [(b)] If $W$ is a connected component of $U$, then $W \cap F$ is a singleton. 
\end{itemize}
\ecor

\sp\fr As an immediate consequence of this corollary and Theorem~\ref{t1j290a}, we get the following. 

\sp\bthm\label{t2j290f}  
The system ${\mathcal S}_{U_r}$ resulting from this corollary and Theorem~\ref{t1j290a} is a maximal conformal graph directed Markov systems \index{(N)}{maximal conformal graph directed Markov system }in the  sense of Chapter~\ref{Markov-systems}.
\ethm

We shall now provide some sufficient conditions for the hypotheses of Corollary~\ref{c120180827} to be satisfied. We shall prove the following. 
 
\sp\bthm~\label{t2j291C}   
Let $Y$ be a complete Riemann surface with constant curvature $0$, i.e. $Y$ is either the complex plane $\C$, a complex torus $\mT_\La=\C/\La$, where $\La$ is a lattice on $\C$ or an infinite cylinder $\C/2\pi i\Z$.

\begin{itemize}
\item Fix $R>0$ so small that every open ball in $Y$, with radius $R$ is simply connected, i.e. isometrically equivalent $B(0,R)\sbt\C$. Fix $\l>1$ and $\kappa>1$.

\,

\item Let $X$ be a non--empty open subset of $Y$ and let 

\,

\item $f:X \lra Y$ be an analytic map with the Standard Property. 

\,

\item Let $Q\sbt Y$ be a set witnessing this property, i.e. $Q\sbt (Y\sms X)\cup\Per(f)$ and $Q$ has at least three elements. 

\,

\item Fix $F$, a finite subset of 
$$
\bal
J(f)\sms \lt(Q\cup \ov{{\rm PS}(f)}\cup\bu_{n=0}^\infty f^{-n}(Y\sms X)\rt)
&=J(f)\cap\bi_{n=0}^\infty f^{-n}(X)\sms \big(Q\cup \ov{{\rm PS}(f)}\big) 
\\
&=J(f)\cap\bi_{n=0}^\infty f^{-n}(Y)\sms \big(Q\cup \ov{{\rm PS}(f)}\big)
\eal
$$
such that  
\beq\label{520180827}
F\cap \bigcup_{n=1}^\infty  f^n(F)=\es.
\eeq
\end{itemize}

Then for every 
$$
r\in \lt(0,\frac{1}{6}\min\Big\{R,\min\big\{\rho(a,b):a,b \in F, a\ne b\}\big\}, \dist\(F, Q\cup\ov{\rm PS}(f)\)\Big\}\rt)
$$
small enough there exists a nice set $U=U_r$ with the following properties:
\begin{itemize}
\item [(a)] $ B(F,r)\sbt U  \sbt  B(F, \kappa r)\sbt  B(F,2\kappa r)\sbt J(f)\sms  \ov{{\rm PS}(f)}$,

\sp\item [(b)] If $W$  is a connected  component  of $U$, then  $W \cap F$  is  a singleton. Denoting this singleton by $b$, we then set $\hat W:=B(b,\ka r)$.
\end{itemize}
\ethm

\bpf
We shall prove that for every 
$$
r\in \De:=\lt(0,\frac{1}{6}\min\Big\{R,\min\big\{\rho(a,b):a,b \in F, a\ne b\}\big\}, \dist\(F, Q\cup\ov{\rm PS}(f)\)\Big\}\rt)
$$
small enough the hypotheses of Corollary~\ref{c120180827} are satisfied. This in fact means that we are supposed to check property (2) of this corollary for such radii $r$. Indeed, by Lemma~\ref{l1j293} there exists $N\ge 1$ so large that 
\beq\label{320180824M}
|(f^{n})'(z)|\ge \max\lt\{\frac{8\kappa}{\kappa-1}, \lambda\rt\}
\eeq
for all $n\ge N$ and all $z\in f^{-n}\(B(F,\De)$. On the ther hand, since the set $F$ is finite, it follows from \eqref{520180827} that there exists $\De_1\in(0,\De/4)$ so small that 
\beq\label{620180827}
B(F,2\De_1)\cap \bigcup_{n=1}^N  f^n\(B(F,2\De_1)\)=\es.
\eeq
So, taking any $r\in (0,\De_1)$ finishes the proof.
\epf

\sp\fr As an immediate consequence of this theorem and Theorem~\ref{t1j290a}, we get the following. 

\sp\bthm\label{t2j290g}  
For all $r>0$ small enough the systems ${\mathcal S}_{U_r}$ resulting from Theorem~\ref{t2j291C} and Theorem~\ref{t1j290a} are maximal conformal graph directed Markov systems in the  sense of Chapter~\ref{Markov-systems}.
\ethm

\sp We will also prove the following result, close to the current subject, and interesting in itself. 

\sp\blem\label{l1j295} 
Let $Y$ be a complete Riemann surface with constant curvature $0$, i.e. $Y$ is either the complex plane $\C$, a complex torus $\mT_\La=\C/\La$, where $\La$ is a lattice on $\C$ or an infinite cylinder $\C/2\pi i\Z$.

\begin{itemize}

\item Fix $R>0$ so small that every open ball in $Y$, with radius $R$ is simply connected, i.e. isometrically equivalent $B(0,R)\sbt\C$. 

\,

\item Let $X$ be a non--empty open subset of $Y$ with finite area.

\,

\item Let $f:X \lra Y$ be an analytic map. 

\,

\item Let $w \in X$, and let $r\in(0,R]$ be so small that $B(w,2r)\sbt X$. 
\end{itemize}

\,

\fr Then for every $ n \geq 1$,
$$
\sup\lt\{\big|(f^{-n}_V)'(z)\big|: \,   V \in  \bigcup_{k=1}^{n}\Comp^*_k(w,r), \,  z \in V\rt\}  < + \infty
$$
and for  every  $ \vep >0$,
$$ 
\# \lt\{V \in  \bigcup_{k=1}^{n}\Comp^*_k(w,r):\big|(f^{-n}_V)'(z)\big| \geq \vep \mbox{ for some z } \in  V\rt\} < + \infty.
$$
\elem

\bpf  Obviously, it is enough to prove that
\beq\label{1j295}
  \sup\big\{\big|(f^{-n}_V)'(z)\big|:V \in  \Comp^*_n(w,r), \,  z \in V\big\}< + \infty
\eeq
for all $ n \geq 1$ and
\beq\label{2j295}
\sharp \lt\{V \in  \Comp^*_n(w,r): \big|(f^{-n}_V)'(z)\big| \geq \vep \ \mbox{ for some z } \in V\rt\} < + \infty
\eeq
for all $ n \geq  1$ and all $ \vep>0$. Since the components $V\in \Comp^*_n(w,r)$ are mutually disjoint, by applying  Koebe's Distortion Theorem, we get that
\beq\label{1j297}
+\infty> {\rm Area}(X)  
\geq \sum_{V \in  \Comp^*_n(w,r)} {\rm Area}(V) 
\gtrsim  r^2 \sum_{V \in  \Comp^*_n(w,r)}\big|(f^{-n}_V)'(w)\big|^2.
\eeq
Therefore, using Koebe's Distortion Theorem again
$$
\begin{aligned} \sup\{\big|(f^{-n}_V)'(z)\big|: & \, \,   V \in \Comp^*_n(w,r), z\in W \}\\
  &\leq K  \sup\lt\{|(f^{-n}_V)'(z)|:V \in \Comp^*_n(w,r) \rt\}< +\infty.
\end{aligned}
$$
Denote the set  of components $V$ involved in (\ref{2j295}) by $C_n(\vep)$ and the corresponding  points $z \in V$ by $z_V$.
Using (\ref{1j297}) and  employing  Koebe's Distortion Theorem, we get that
$$
\begin{aligned}
+\infty> {\rm Area}(X) 
&\gek r^2\sum_{V \in C_n(\vep)} |f^{-n}_V)'(w)|^2 
\geq K^{-2}r^2 \sum_{V \in C_n(\vep)}\big||f^{-n}_V)'(z_V)\big||^2\\
& \geq K{-2}r^2\vep^2 \sharp C_n(\vep).
\end{aligned} 
$$
Hence, $\sharp C_n(\vep) \lesssim K^2 r^{-2} \vep^{-2} {\rm Area}(X) < +\infty$. The proof is complete. 
\qed

\sp We end this section with the following theorem whose proof requires the full power of Theorem~\ref{t2j291}.

With the notation of Theorem~\ref{t2j291} and its proof, for every $b \in F_1$, set
\beq\label{1gns30A}
X_b:= X(b,1;u):=\overline{U(b,1)} 
\quad \text{and} \quad 
W(b,1;u):=B(b,4r_b).
\eeq
For every $b \in F_0$ and $j\in \{1,2, \ldots, p_b\}$, set
 \beq\label{4gns30}
 X^*(b,j):= \overline{U(b,j)}\sms f^{-1}_b(U(b,j))
 \eeq
For every $r\in(0, r_b')$ as small as required in Theorem~\ref{t2j291} (G), set
\beq\label{1gns30}
 X(b,j;u):=X^*(b,j)\sms V_j(b,u).
\eeq
For every $A=X(a,j;u)$, $ a \in F$, $j\in \{1, \ldots, p_a\}$, put 
$$
U_A:=U(a,j;u)
\  \  \  {\rm and}  \  \  \
W_A:=W(a,j;u).
$$
We form a system $\cS_{U}=\cS_{U_{{\bold r}}}$ as follows. The pairs 
\beq\label{720190211}
\big\{\(X(b,j;u), W(b,j;u)\), \, b \in F,\, \, j \in \{1, 2,\ldots, p_b\}\big\}
\eeq
form the domains of $\cS_{U}$. In order to define the maps of $\cS_{U}$,
look at all sets
\beq\label{520190211}
A, B \in\big\{X(b,j;u): \, b \in F\, \, \text{ and}\,\,  j \in \{1, 2,\ldots, p_b\}\big\},
\eeq
and all integers $n \geq 1$ such that
 \beq\label{3gns30}
 \Int(A) \cap f^{-n}(\Int (B))\neq \es.
 \eeq
Because of the hypothesis (5) and item (F) of Theorem~\ref{t2j291}, for every $\xi \in \Int (A) \cap f^{-n}(\Int (B))$ there exists 
\beq\label{120190211}
f^{-n}_\xi: W_B \lra X,
\eeq
a unique holomorphic branch of $f^{-n}$, defined on $W_B$ and sending $f^n(\xi) \in \Int(B) \sbt W_B$ to $\xi$. We declare that the map $f^{-n}_\xi: W_B \lra X$ belongs to $\cS_U$ if 
 \beq\label{5gns30}
  f^s\(f^{-n}_\xi(\Int(B))\) \cap \bigcup_{w \in F}\bigcup_{j=1}^{p_w} \Int(X(w,j;r))=\es
 \eeq
 for all integers $0< s< n$. It follows from \eqref{3gns30} that

\beq\label{120190209}
U_A\cap f^{-n}_\xi(U_B)\ne\es.
\eeq

\sp\bthm\label{t2j291B} Let $Y$ be a complete Riemann surface with constant curvature $0$, i.e. $Y$ is either the complex plane $\C$, a complex torus $\mT_\La=\C/\La$, where $\La$ is a lattice on $\C$ or an infinite cylinder $\C/2\pi i\Z$.

Fix $R>0$ so small that every open ball in $Y$, with radius $R$, is simply connected, i.e. isometrically equivalent $B(0,R)\sbt\C$. Fix $\kappa\in(1,2]$.

Let $X$ be a non-empty open subset of $Y$ and let $f:X \lra Y$ be an analytic map with the Standard Property. Let $Q\sbt Y$ be a set witnessing this property, i.e. $Q\sbt (Y\sms X)\cup\Per(f)$ and $Q$ has at least three elements. 

\fr Fix 

\begin{itemize}

\,

\item $F$, a finite subset of $J(f)\sms O_+(Q)$, 

\,

\item a collection $\{U_0(b)\}_{b\in F}$ of open subsets of $Y$, 

\fr and 

\,

\item a vector $\bold r=\(r_b:b\in F\)$ with the following properties:
\end{itemize}

\, 

\begin{enumerate}  
\item
$$
U_0(b)\sbt B(b,r_b)
$$
for all $b\in F$.

\, \, \item
$F$ can be represented as a disjoint union:
$$
F=F_0\cup F_1
$$
such that
\begin{enumerate}

\, \, \item $f(b)=b$ for every  $b\in F_0$.

\, \, \item  For each $b\in F_0$, the map $f|_{B(b,6r_b)}$ is 1--to--1 and the holomorphic inverse branch $f_b^{-1}:B(b,6r_b)\lra X$, sending $b$ to $b$,  is well defined.

\, \, \item  $f_b^{-1}\(B(b,3r_b)\)\sbt B(b,6r_b)$.

\, \, \item $b\in \bd U_0(b)$ for all $b\in F_0$,

\, \, \item For each $b\in F_0$ there exists an integer $p_b\ge 1$ such that $U_0(b)$ has exactly $p_b$ (open) connected components $U_0(b,j)$, $j=1,2,\ld, p_b$, each of which is simply connected and its closure is a closed topological disk.

\, \, \item For each $b\in F_0$ and each $j=1,2,\ld, p_b$:
$$
f_b^{-1}(U_0(b,j))\sbt U_0\(b,j)
$$
and

\, \, \item The sequence 

$$
\Big(f_b|_{U_0(b,j)}^{-n}:U_0(b,j)\lra X\Big)_{n=0}^\infty
$$ 
converges uniformly to the constant function which assignes to each point in $U_0(b,j)$ the value $b$.

\, \, \item $B(b,r_b/2)\sbt U_0(b)$ for all $b\in F_1$,

\, \, \item For each $b\in F_1$ the set $U_0(b)$ is connected, simply connected, and its closure is a closed topological disk. We then put $p_b:=1$ and denote also $U_0(b)$ by $U_0(b,1)$.

\, \, \item  For each $b \in F_0$ and $s>0$ there exists $s^{-}_b \in (0, r_b]$ such that if
$$ 
z \in B(b, s^{-}_b)\sms \bigcup_{j=1}^{p_b}U_0(b, j),
$$
then 
$$
f^n(z) \in B(b,s)
$$ 
for all $n\geq 0$.
\end{enumerate}

\, \, \item
$$
r_b\in \lt(0,\frac{1}{6}\min\Big\{R,\min\big\{\rho(b,c):c \in F\sms \{b\}\big\}\Big\}\rt)
$$
for all $b\in F_0$.

\, \, \item
$$
r_b\in \lt(0,\frac{1}{6}\min\Big\{R,\min\big\{\rho(b,c):c \in F\sms \{b\}\big\}, \dist\(b, \ov{\rm PS}(f)\)\Big\}\rt)
$$
for all $b\in F_1$.

\, \, \item 
Suppose that $a, b\in F$ and $n\ge 1$ is an integer. Assume that either $b\in F_1$ or there exists a point 
\beq\label{120180824R}
w\in f^{-n}(U_0(b))
\eeq
such that 
\beq\label{220180824RB}
f^{n-1}(w)\notin U_0(b),
\eeq
then
 
\begin{enumerate}
\, \, \item 

the holomorphic inverse branch $f_w^{-n}:B(b,6r_b)\lra X$ of $f^n$, sending $f^n(w)$ to $w$, is well--defined (this is only an extra hypothesis if $b\in F_0$; if $b\in F_1$, this follows from (4)).

\, \, \item 
If, in addition, the connected component of $f^{-n}(B(b,2r_b))$ containing $w$ intersects $B(a,2r_a)$, and 
\beq\label{220180824R}
f^j(w)\notin \bu_{c\in F}U_0(c)
\eeq
for all $j=1,2,\ld, n-1$, then 
\beq\label{32018082R}
|(f^{n})'(z)|\ge \frac{8\kappa}{\kappa-1}\cdot \frac{r_b}{r_a}
\eeq
for all $z\in f_w^{-n}(B(b,2r_b))$. 
\end{enumerate}

\, \, 

\item For every $b \in F_0$ and every $j \in \{1,2, \ldots, p_b\}$ there exists, $V_j(b)$, an open neighborhood of $b$ in $Y$ such that:

\,

\begin{itemize}
\item [(a)] The set $\overline{V_j(b)}\cap \overline{U_0(b, j)}$ is a closed topological disk,

\,

\item [(b)] The set $V_j(b)\sms \overline{U_0(b, j)}$ is  connected (and, in consequence, the set $\overline{V_j(b)}\sms \overline{U_0(b, j)}$ is connected too), 

\,

\item [(c)] $\overline{V_j(b)}\sbt B(b, r_b'/4)$ with some $r_b'\in(0, (r_b)_b^-)$.
\end{itemize}

\, \, 

\item $|f'(b)|=1$ for every $b\in F_0$.

\, \, 

\item  
$$ 
\Crit(f)\cap \bigcup_{w \in F_0} B(w, 6r_w) =\es.
$$

\, \, 

\item $\forall b\in F_0 \, \, \, \, \om(\Crit(f))\cap B(b,6r_b) \sbt \{b\}$.
\end{enumerate}

\sp\fr Then the system $\cS_U$ generated by the domains of \eqref{720190211} and by the maps of \eqref{120190211} is a maximal conformal graph directed Markov system in the sense of Definition~\ref{d120190419} and Definition~\ref{d2_2017_11_18}.
\ethm

\bpf 
The number $u\in(0, r_b')$ will be further required to be sufficiently small in the course of the proof. For every integer $n \geq 1$, every $b \in F_1$ and every $j\in \{1, \ldots, b\}$, we have
$$
\begin{aligned} 
f^{-n}_b(X(b,j;u))& \sbt f^{-n}_b(X^*(b, j))
 \sbt f^{-n}_b\(\overline{U(b,j)}\sms f^{-1}_b(U(b,j))\)\\
&\sbt f^{-n}_b(\overline{U(b,j)})\sbt f^{-1}_b(\overline{U(b,j)}).
\end{aligned}
$$
Hence
  \beq\label{2gns30}
\begin{aligned}
 \Int\(X(b,j;u))\cap f^{-n}_b(X(b,j;u))\)
 & \sbt \Int\lt(\overline{U(b,j)}\sms f^{-1}_b(U(b,j))\rt)\cap  f^{-1}_b(\overline{U(b,j)})\\
 & = \left(\Int(\overline{U(b,j)})\sms \overline{f^{-1}_b(U(b,j))}\right)\cap  f^{-1}_b(\overline{U(b,j)})\\
& \sbt \left(\Int(\overline{U(b,j)})\sms \overline{f^{-1}_b(U(b,j))}\right)\cap  \overline{f^{-1}_b(U(b,j))}\\
&=\es.
\end{aligned}
\eeq
Fix two sets $A, B$ as in formula \eqref{520190211}. Suppose that the formula \eqref{3gns30} holds and fix 
$$
\xi \in \Int (A) \cap f^{-n}(\Int (B)).
$$ 
We consider several cases.

\sp\fr Case~$1^0$. $A=X_a$ and $B=X_b$ for some $a,b \in F_1$. It then follows from Theorem~\ref{t2j291}, Proposition~\ref{p1j290}, and formula \eqref{120190209} that 
\beq\label{1gns33}
 f^{-n}_\xi(X(b,j;u))
 = f^{-n}_\xi(\overline{U(b,1)})
\sbt\overline{f^{-n}_\xi(U(b,1))}
\sbt \overline{U(a,1),
}=X_a=X(a,1;u)
\eeq
and we are done in this case.

\sp\fr Case~$2^0$.  $A=X_a$ for some $a \in F_1$ and $B=X(b,j;u)$ for some $b \in F_0$ and $j \in \{1,2, \ldots, p_b\}$. Then, similarly as in the previous case, it follows from Theorem~\ref{t2j291}, Proposition~\ref{p1j290}, and formula \eqref{120190209} that 
\beq\label{2gns33.1}
 f^{-n}_\xi(X(b,j;u))\sbt f^{-n}_\xi(\overline{U(b,j)})\sbt
\overline{f^{-n}_\xi(U(b,j))}\sbt \overline{U(a,1)}=X_a=X(a,1;u),
\eeq
and we are done in this case too.

\sp\fr Case~$3^0$. $A=X(a,j;r)$ for some $a\in F_0$ and $j \in \{1,2, \ldots, p_a\}$ and $B=X_b$ for some $ b \in F_1$. 
Seeking contradiction suppose that 
$$
f^{-n}_\xi ( U(b,1))\nsubseteq U(a,j)) \sms f^{-1}_a(\overline{U(a, j)}).
$$ 
But by  Theorem~\ref{t2j291}, Proposition~\ref{p1j290}, and formula \eqref{120190209},
$$
f^{-n}_\xi ( U(b,1))\subset U(a,j)).
$$ 
So,
$$
f^{-n}_\xi ( U(b,1))\cap f^{-1}_a(\overline{U(a, j)})\neq \es.
$$
Hence, $f^{-n}_\xi ( U(b,1))\cap f^{-1}_a(U(a, j))\neq \es.$
Thus $f \circ f^{-n}_\xi ( U(b,1))\cap U(a, j)\neq \es.$
Let then $k \geq 1$ be the largest integer such that
\beq\label{1gns34}
f^s \circ f^{-n}_\xi ( U(b,1))\cap U(a, j)\neq \es.
\eeq
for all  $1\leq s \leq k$. Such integer $k$ exists and $k <n$ because \beq\label{320190209}
U(b,1)\cap U(a,j)=\es  \  \  \  ({\rm as} \  a \neq b). 
\eeq
It also follows from (\ref{1gns34}) and (\ref{5gns30}) that
$$ \begin{aligned}
\es \neq f^k\circ f^{-n}_\xi(U(b,1))\cap U(a, j)&\sbt U(a,j)\sms\left(\Int(\overline{U(a, j)})\sms \overline{f^{-1}_a(U(a,j))}\right)\\
   & \sbt U(a,j)\cap \overline{f^{-1}_a(U(a,j))}.
   \end{aligned}$$
   Since the set $f^k\circ f^{-n}_\xi(U(b,1))$ is open, this implies that
$$   f^k \circ f^{-n}_\xi ( U(b,1))\cap f^{-1}_a(U(a, j))\neq \es.$$
Hence
$$   f^{k+1} \circ f^{-n}_\xi ( U(b,1))\cap U(a, j)\neq \es, $$
contrary to the definition of $k$. We have thus proved that
\beq\label{1gns35}
f^{-n}_\xi(U(b,1))\sbt U(a,j)\sms f^{-1}_a(\overline{U(a,j)}).
\eeq
Since  the sets $\big\{U(c,k);\, c \in F, \, k\in \{1, \ldots, p_c\}\big\}$ are pairwise disjoint, it follows from (\ref{1gns35}) and item (E) of Theorem~\ref{t2j291} that
$$ 
f^{-n}_\xi(U(b,1))\cap\bigcup_{c \in F} \bigcup_{k=1}^{p_c} f^{-1}_c(U(c,k))=\es.
$$
Hence 
\beq\label{1gns35.1}
f \circ f^{-n}_\xi(U(b,1))\cap \bigcup_{c \in F} \bigcup_{k=1}^{p_c} U(c,k)=\es.
\eeq
Now, seeking contradiction suppose that 
$$ 
f^{-n}_\xi(U(b,1))\nsubseteq X(a,j;u).
$$
It then follows  from (\ref{1gns35}), (\ref{4gns30}), and (\ref{1gns30}), that
\beq\label{2gns35}
f^{-n}_\xi(U(b,1))\cap V_j(a,u) \neq \es.
\eeq
Now, it is the moment to fix $u\in(0, r_b')$ sufficiently small. We assume, on top of all others assumptions the number $r$ to be so small that
$$ 
V_i(w,u)\cap f(V_i(w,u))\cup f^2(V_i(w,u))\sbt B(w, r_w')
$$
for all $w \in F_0$ and  $i \in \{1, \ldots, p_w\}$. It then follows from (\ref{2gns35}) that

\beq\label{5gns35}
f\circ f^{-n}_\xi(U(b,1))\cap B(a,r_a') \neq \es.
\eeq
This along with (\ref{1gns35.1}) and hypothesis (2j) (and since $r_b'\le (r_b)_b^-$ and $U_0(c, k)\sbt U(c,k)$ for all $c \in F$ and $k \in \{1, \ldots, p_c\}$), yields
\beq\label{6gns35}
\es \neq f^l\left(f\circ f^{-n}_\xi(U(b,1))\cap B(a,r_a')\right) \sbt B(a,r_a)
\eeq
for all $l \geq 0$. In particular, taking $l=n-1\geq 0$, we get that $U(b,1)\cap B(a, r_a)\neq \es$.
Since $a\neq b$, this contradicts the hypotheses (3) and (4). The inclusion
$$ 
f^k\circ f^{-n}_\xi(U(b,1))\sbt X(a,j;u)
$$ is thus  proved. This in turn, by minimality of $n$, implies that $k=0$, whence,
\beq\label{4gns35}
 f^{-n}_\xi(U(b,1))\sbt X(a,j;u).
 \eeq
Therefore, as the set $X(b,j;u)$ is closed,
\beq\label{3gns35}
f^{-n}_\xi(B)=f^{-n}_\xi(X_b)=f^{-n}_\xi(\overline{U(b,1)}) \sbt \overline{f^{-n}_\xi(U(b,1))}\sbt X(a,j;u)=A.
\eeq

\sp\fr Case~$4^0$. $ A= X(a,i;u)$ and $B=X(b,j;u)$ for some $a \in F_0$, $ i \in \{1, \ldots, p_a\}$, $ b \in F_0 \sms \{a\}$ and $j \in \{1, \ldots, p_b\}$.

Since $U(b,j)\sbt W(b,j;u)$ (so $f^{-n}_\xi$ is well defined on $U(b,j)$, we can proceed in this case in entirely the same way as in Case~$3^0$, with only $U(b,1)$ replaced by $U(b,j)$.
We then end up with the formula
$f^{-n}_\xi(U(b,j))\sbt X(a,i;u)$, corresponding to formula  (\ref{4gns35}) of Case~$3^0$. Then
\beq\label{1gns36}
f^{-n}_{\xi}(X(b,j;u))\sbt f^{-n}_{\xi}(\overline{U(b,j)}) \sbt \overline{f^{-n}_{\xi}(U(b,j))}\sbt X(a,i;u).
\eeq

\sp Consider in turn

\sp\fr Case~$5^0$.
$ A=X(c,i;u)$ and $B=X(c,j;r)$ for some $c \in F_0$, and  $ i\neq j \in  \{1, \ldots, p_c\}$. 

We proceed exactly as in Case~$3^0$ with obvious replacements $X(a,j;u)$
by $X(c,i;u)$ and $U(b,1)$ by $U(b,j)$. Also now the formula 
$$
U(a,i)\cap U(a,j)=\es
$$
replaces \eqref{320190209} and holds not because $a\ne b$ but because $i\ne j$. We then obtain the formula
$$ 
f_\xi^{-n}(U(c,j)) \sbt X(c,i;u),
$$
corresponding to the formula (\ref{4gns35}) of Case~$3^0$. Hence,
\beq\label{420190209}
f^{-n}_{\xi}(X(c,j;u))
\sbt f^{-n}_{\xi}(\overline{U(c,j)}) \sbt \overline{f^{-n}_{\xi}(U(c,j))}\sbt X(c,i;u).
\eeq
 
\sp We now consider: 

\sp\fr Case~$6^0$. Eventually assume that $ A=X(a,j;u)$ and $B=X(a,j;u) $  with some $c \in F_0$ and  $j \in  \{1, \ldots, p_c\}$. If
$$
f^u(f^{-n}_\xi(U(c,j)) \cap B(c,6r_c))\neq \es
$$  
for all integers $u=0,1, \ldots, n$, then $f^{-n}_\xi=f^{-n}_c$, and so, this is ruled out by the formula (\ref{2gns30}). Hence
\beq\label{2gns37}
f^u(f^{-n}_\xi(U(c,j))) \cap B(c,6r_c))=\es
\eeq
for some integer $u\in\{0,1, \ldots, n\}$. We now proceed again exactly as in Case~$3^0$ with obvious replacements $X(a,j;u)$
by $X(c,j;u)$ and $U(b,1)$ by $U(c,j)$. But now the fact that the integer $k \geq 1$, involved in the formula (\ref{1gns34}), is strictly smaller than $n$ follows not from $a\ne b$ but from (\ref{2gns37}).  As in the previous case (Case~$5^0$), we then obtain the formula
$$ 
f_\xi^{-n}(U(c,j)) \sbt X(c,j;u),
$$
corresponding to the formula (\ref{4gns35}) of Case~$3^0$. Hence, 
\beq\label{420190209B}
f^{-n}_{\xi}(X(c,j;u))
\sbt f^{-n}_{\xi}(\overline{U(c,j)}) \sbt \overline{f^{-n}_{\xi}(U(c,j))}\sbt X(c,j;u).
\eeq
Because of this, Lemma~\ref{l1j293}, and Remark~\ref{r120190212} Condition (f) of the Definition~\ref{GDMS_def} of graph directed Markov systems is thus satisfied. 

\sp Let us now prove condition (4b) of the Definition~\ref{d120190419} of conformal graph directed Markov systems, i.e. the Open Set Condition. So, suppose that $A, B, C$ and $D$ are in the set of \eqref{520190211},  
$$
\xi \in \Int (A) \cap f^{-m}(\Int (B)), \  \  \  \zeta \in \Int (C) \cap f^{-n}(\Int (D))
$$
formula \eqref{5gns30} holds for $\xi$ and $B$ and for $\zeta$ and $D$, and
\beq\label{620190211}
f_\xi^{-m}(\Int (B))\cap f_\zeta^{-n}(\Int (D))\ne\es
\eeq
for some integers $m, n\ge 1$. Because of, already proven, Condition (f) of the Definition~\ref{GDMS_def}, we conclude from (\ref{620190211}) that $A\cap C\ne\es$. Therefore, $A=C$. Suppose without loss of generality that $m\le n$. It also follows from (\ref{620190211}) that 
$$
\Int (B)\cap f^m\(f_\zeta^{-n}(\Int (D)\)\ne\es. 
$$
Along with \eqref{5gns30} this implies that $m=n$. But then it follows from (\ref{620190211}) that $B=D$. Applying then (\ref{620190211}) once more and remembering that $U_A=U_B$ is an open connected simply connected 
set, we conclude that $f_\xi^{-m}=f_\zeta^{-n}$. The Open Set Condition, i.e. condition (4b) of the Definition~\ref{d120190419} is satisfied.

\sp Since condition (4c) of the Definition~\ref{d120190419} of conformal graph directed Markov systems is satisfied because of \eqref{120190211}, the proof of Theorem~\ref{t2j291B} is complete.
\epf

\brem\label{r120190211}
Instead of assuming in Theorem~\ref{t2j291B} that $f$ satisfies the Standard Property and $F\cap O_+(Q)$, we could have have assumed that 
$r_a=r_b$ for all $a, b\in F$ and use this assumption rather than Lemma~\ref{l1j293} to prove Condition (f) of the Definition~\ref{GDMS_def}. 
\erem



\chapter{Geometry and Dynamics of Elliptic Functions}\label{geometry-and-dynamics}

Through this whole chapter we deal with general non--constant elliptic functions, i. e. impose no assumptions on a given non--constant elliptic function. We first give some basic preliminary facts about such functions. Then, following the paper \cite{KU2}, by associating to a given elliptic function an infinite alphabet conformal iterated function systems, and heavily utilizing its $\th$ number, we provide a strong, somewhat surprising, lower bound for the Hausdorff dimension of the Julia sets of all non--constant elliptic functions. In particular, this estimate shows that the Hausdorff dimension of the Julia sets of any non--constant elliptic function is strictly larger than $1$. We also provide a simple closed formula for the Hausdorff dimension of the set of points escaping to infinity of an elliptic function.
In the last section of this chapter we prove that no conformal measure of an elliptic function charges the set of escaping points.

\section{Selected Preliminaries}\label{Selected Preliminaries}

As  indicated in the title, throughout this chapter $f:{\mathbb
C}\lra \oc$ is  a non--constant elliptic function. Every
such function is doubly periodic and meromorphic. Let 
$$
\La_f\sbt\C
$$ 
be the set of all priods of $f$. We know from the first two sections of
Chapter~\ref{elliptic-theory} that then
there exist unique two vectors  $\l_1,\, \l_2$, $\im(\l_1/\l_2)\neq
0$, such that  
$$
\La_f=[\l_1, \l_2]= \{m\l_1 + n\l_2: \; m,n\in {\mathbb Z}\; \}.
$$
In particular,
$$
f(z)=f(z+m\l_1 + n\l_2)
$$
for all $z\in {\mathbb C}$ and all $n,m \in {\mathbb
Z}$. Recall that the set $\La_f$ is called the lattice of the elliptic function $f$. We say that
$$
z\sim_{\La_f} w
$$
if $w-z\in \La_f$. Let
$$\index{(S)}{${\mathcal R}_f$}
{\mathcal R}_f:=\{t_1 \l_1+t_2\l_2:  0\leq t_1, t_2 \leq 1\}
$$
be the basic fundamental parallelogram\index{(N)}{fundamental
parallelogram} of $f$. It follows from the periodicity of $f$ that
$$
f({\mathbb C})=f({\mathcal R}_f).
$$
Therefore the set $f({\mathbb C})$, is simultaneously compact, whence
closed, and open in $\oc$. Since $\oc$ is connected, the set $f({\mathbb C})$ is equal to $\oc$. This means that each elliptic function is surjective, and a bit more:
\beq\label{220190913}
f({\mathbb C})=f({\mathcal R}_f)=\oc
\eeq
But even more is true. Having this,  the following observation is immediate.

\bobs\label{o420200303}
If $f:\C\lra\oc$ is a non--constant elliptic function, then $\Sing(f^{-1})=f(\Crit(f))$.
\eobs
Of course we have the following.

\bobs\label{o220190913}
If $f:\C\lra\oc$ is a non--constant elliptic function, then
$$
J(f)=\ov{\bu_{n=1}^\infty f^{-n}(\infty)}.
$$
\eobs

\fr A profound extension of \eqref{220190913} is this.

\bprop\label{p120190913}
Each non--constant elliptic function $f:\C\lra\oc$ is topologically exact \index{(N)}{topologically exact} in the sense that if $U\sbt \C$ is an open set intersecting the Julia set $J(f)$, then there exists an integer $l\ge 1$ such that
$$
f^l(U)=\oc. 
$$
\eprop
\bpf
It follows from Observation~\ref{o220190913} that for some integer 
$l\ge 2$, the image $f^{l-1}(U)$ contains an open neighborhood of $\infty$ in $\oc$. Thus, it contains at least one (in fact
infinitely many) congruent copy of the fundamental parallelogram $\mathcal R_f$ of $f$. Consequently, by \eqref{220190913}, we get that
$$
f^l(U)\spt f(\mathcal R_f)=\oc.
$$
The proof is complete. 
\epf

\sp It also follows from the periodicity of $f$ that
$$
f^{-1}(\infty)=\bu_{m,n\in {\mathbb Z}}\({\mathcal R}_f\cap
f^{-1}(\infty) +m\l_1+ n\l_2\).
$$
For every pole $b$ of $f$ let $q_b$\index{(S)}{$q_b$} denote its
multiplicity\index{(N)}{multiplicity of the pole}. We define
$$
q=q_{\max}(f):=\sup\{q_b: b\in  f^{-1}(\infty)\}
  =\max \{q_b: b\in  f^{-1}(\infty) \cap {\mathcal R}_f\}.
$$\index{(S)}{$q$}\index{(S)}{$q_{\max}(f)$}
For every $R>0$ we have defined in Section~\ref{FC B+S} the following two sets:
$$
B_\infty(R)=\{z\in \ov{\mathbb C}: |z| >R\} \  \
{\rm and} \ \
B_\infty^*(\infty)=\{z\in {\mathbb C}: |z| >R\}.
$$
Given $b\in f^{-1}(\infty)$ let 

\sp \centerline{$B_b(R)$\index{(S)}{$B_b(R)$} be the
connected component of $f^{-1}(B_\infty(R))$ containing $b$} 

\fr and let

\, \centerline{$B_b^*(R):=B_b(R)\sms\{b\}$}\index{(S)}{$B_b^*(R)$}.

\,\fr More generally, given $k\geq 1$ and $b \in f^{-k}(\infty)$, let
$B^k_b(R)$\index{(S)}{$B^k_b(R)$} be the connected component of
$f^{-k}(B_R) $ containing $b$, and $B^{*k}_b(R):=B^k_b(R)\sms\{b\}$. Fix $T=T(f)\geq 1$  so large that all
components $B_b(T)$, $ b \in {f}^{-1}(\infty),$ are mutually
disjoint.

\sp Since each elliptic function is meromorphic, all considerations of Chapter~\ref{topological-picture}, Topological Picture of Iterations of (all!) Meromorphic Functions, apply to them. Recall  that $\Crit(f)$ is the set of critical points of $f$, i.e.
$$
\Crit(f):=\{z\in\C:\, f'(z)=0\}.
$$  
Its image, $f(\Crit(f))$, has been called the set of critical 
values of $f$. Since ${\mathcal R}_f\cap\Crit(f)$ is finite and since
$f(\Crit(f))= f({\mathcal R}_f\cap\Crit(f))$, the set of critical
values $f(\Crit(f))$ is also finite. Since each elliptic function has obviously no asymptotic values, this gives us the following obvious observation.

\bobs\label{o1_2017_09_25}
All non--constant elliptic functions belong to Speiser class ${\mathcal S}$. 
\eobs

As an immediate consequence of this observation and Theorem~\ref{t120200404} we get the following.

\bthm\label{t220200404}
Any elliptic function has only finitely many attracting and rationally indifferent periodic points.
\ethm

Because of this same observation, all considerations of Section~\ref{FC B+S}, Fatou Components of (General) Meromorphic Functions II; Class ${\mathcal B}$ and ${\mathcal S}$, apply to them. Therefore, as immediate consequence of Theorem~\ref{baker-domain} and Theorem~\ref{Sullivan}, we get the following.

\bthm\label{baker-domain+Sullivan for elliptic} 
No elliptic functions have Baker or wandering domains. 
\ethm 

\sp Now  we shall prove the following. 

\bthm\label{herman} 
No even elliptic function has a cycle of Herman rings. 
\ethm
 
 \bpf Suppose that $f$ has a cycle of Herman rings
 $\{U_0, U_1, \ldots ,U_{p-1}\}$ of some period $p\ge 1$. Then for any $i=0, 1,
 \ldots, p-1$, the iterate $f^p:U_i\to U_i$ is conjugate to an irrational
 rotation of an annulus $A(0;1,R)$ with some $R>1$, and thus it is bijective and $\infty\notin U_0\cup U_1\cup\ld\cup U_{p-1}$. The preimages under these conjugacies from $A(0;1,R)$ to $U_i$, $i=0, 1, \ldots, p-1$, of the circles $\{z\in\C:|z|=r\}$, $ 1<r <R$, foliate the rings $U_i$ with $f^p$ forward invariant leaves on which $f$ is bijective. Let $\g$ be an $f^p$ forward invariant leaf of $U_0$, and let
 $B_\g$ denote the bounded component of the complement of $\g$; remember that $\infty\notin U_0$, so $\infty\notin \g$. If $B_\g$ contained no pre-pole of $f^p$, then the standard Maximal Modulus Theorem argument would yield 
$$
f^{pk}(B_\g)\sbt B_\g
$$ 
for all $k\ge 0$. Hence, $B_\g$ would be contained in the Fatou set of $f$, contrary to the definition of a Herman ring. Thus $B_\g$ contains a pre--pole of $f^p$. Therefore, there is a smallest integer $n\ge 0$ such $f^n(\g)$ has a lattice point $\om$ in $B_{f^n(\g)}$. Let $U_j$, $j=0,1,\ldots,p-1$, denote the Herman ring $f^n(U_0)$.

Since the function $f$ is even, $-U_j$ is also a connected component of the Fatou set $F(f)$. Since $f$ is $\La$-invariant, so is the set $-U_j+2\om$. In addition $\om \in -B_{f^n(\g)}+2\om$. If $\om$ belongs to one of the sets $U_j$ or $-U_j+2\om$, then it also belongs to the other, and so
\beq\label{f1_2015_09_03}
U_j=-U_j+2\om,
\eeq
as both these sets are comonents of the Fatou set $F(f)$. Otherwise, $\om\notin U_j\cup (-U_j+2\om)$. Hence, $\om$ belongs to the bounded connected components of both the complement of $U_j$ and of $-U_j+2\om$. Therefore, supposing that $U_j$ and $-U_j+2\om$ are disjoint, we conclude that one of the sets $U_j$ or $-U_j+2\om$ is contained in the bounded compontent of the complements of the other. This however is a contradiction since $\sup\{|z-\om|:z\in U_j\}=\sup\{|z-\om|:z\in -U_j+2\om\}$. Since any two components of the Fatou set of $f$ are either disjoint or equal this yields \eqref{f1_2015_09_03} in this case too. But then for every $z\in U_j$, $-z+2\om$ is also in $U_j$, the points $z$ and $-z+2\om$ (if only $z\ne\om$) are different and $f(z)=f(-z+2\om)$. This contradicts bijectivity of $f$ on $U_j$ and finishes the proof.
\endpf

\brem 
For elliptic functions $f:\C\to\oc$ we slightly modify the definition of the Julia set  $J(f)$, i.e.  $J(f)$ \index{(N)}{Julia set of elliptic function} is the complement of the Fatou set $F(f)$ in ${\mathbb
 C}$. Thus $\infty \notin J(f)$. This is a prerequisite  for  our considerations  concerning conformal measures and geometric measures supported on Julia sets.
\erem

\, Now, we analyze in greater detail the behavior of elliptic functions near poles. Since the set $f(\Crit(f))$ is finite, if $R>0$ is large enough, say $R\ge R_0\ge 4$, then 
\beq\label{1_2017_10_02}
B_\infty(R/4)\cap f(\Crit(f))=\es,
\eeq
i.e. $B_{R/4}$ contains no critical values of $f$. In addition, the sets $B_b(R)$\index{(S)}{$B_b(R)$}, $b\in f^{-1}(\infty)$, are connected, simply connected, mutually disjoint, and there exists $A_1=A_1(f)\ge 1$\index{(S)}{$A_1(f,b)$} such that for all $b\in f^{-1}(\infty)$ and all $z\in
B_b(R)$, we have that
\beq\lab{u1}
A_1^{-1} |z-b|^{-q_b} \leq |f(z)| \leq  A_1 |z-b|^{-q_b}.
\eeq
The fact that the number $A_1=A_1(f)$ can be taken here independently of $b\in f^{-1}(\infty)$ follows from the fact that there are only finitely many equivalence classes of poles mod $\La_f$. 

\sp\fr Keep $b\in f^{-1}(\infty)$. If $U\sbt B_\infty^*(R)$ is an open connected simply connected set, then all the holomorphic inverse branches
$f_{b,U,1}^{-1},\ld,f_{b,U,q_b}^{-1}$ of $f$ are well-defined on
$U$. In addition, if $R_0>0$ is large enough, then there exists $A_2=A_2(f)\ge 1 $\index{(S)}{$A_2(f, b)$} such that for every $1\le j\le q_b$ and all $z\in U$ we have
\beq\lab{1110301}
A_2^{-1} |z|^{-{q_b+1\over q_b}} \leq |(f_{b,U,j}^{-1})'(z)|\leq A_2
|z|^{-{q_b+1\over q_b}}.
\eeq
The number $A_2$ is independent of $b\in f^{-1}(\infty)$ for the same reasons as $A_1$ was. As an immediate consequence of \eqref{1110301} we get for all $z\in U$ the following.
\beq\lab{u2}
\begin{aligned}
 (2A_2)^{-1} {|z|^{{q_b-1\over q_b}}\over |b|^2}\leq (2 A_2)^{-1} {|z|^{{q_b-1\over q_b}}\over 1+|b|^2}
\leq |(f_{b,U,j}^{-1})^*(z)|
    \leq 2 A_2 {|z|^{{q_b-1\over
q_b}}\over 1+|b|^2} \leq 2A_2 {|z|^{{q_b-1\over q_b}}\over |b|^2},
\end{aligned}
\eeq
where the  first left  and  the  second  right  inequality signs were
written assuming in addition that $|b|$ is large enough, i.e. that $R_0>0$ is large enough.
\beq\lab{1042506}\index{(S)}{$A(f,b)$}
A(f):=\max\{ A_1(f),\, A_2(f)\}.
\eeq
We want to emphasize once more that this constant is independent of poles $b$ and simply connected domains $U\sbt B_\infty^*(R)$.
Also based on congruency of poles, a slightly more general straightforward observation from the local behavior around  poles is  that  for every $k \geq 1$
there exist  constants $L_k\ge 1$ and $R_k>0$, both monotone increasing with $k$, such that for all
$b\in f^{-k}(\infty)$ and all $R\ge R_k$, we have
\beq\lab{u4}
\aligned L^{-1}_kR^{-{1\over q_b}} \le &\diam_e(B_b^k(R))
\le L_kR^{-{1\over q_b}}, \\
L^{-1}_kR^{-{1\over q_b}}(1+|b|^2)^{-1} \le &\diam_s(B_b^k(R)) \le
L_kR^{-{1\over q_b}}(1+|b|^2)^{-1}.
\endaligned
\eeq
Frequently, we will write  $L$ for $L_1$.

\

\section[Hausdorff Dimension of Julia sets]{Hausdorff Dimension of
Julia Sets of (General) Elliptic Functions} 

In this  section  we apply  the results of Section~\ref{sec-HD-Bowen} 
(Bowen's Formula (Theorem~\ref{t1j97}) and Theorem~\ref{HD-theta}) to provide  a strong, somewhat surprising,
lower bound for the  Hausdorff  dimension of the Julia sets of all
non--constant elliptic functions. The idea is to associate to each
elliptic function an iterated function system and to apply the above mentioned theorems. We also provide in the next section a closed formula for the Hausdorff dimension of $I_\infty(f)$, the set of points escaping to infinity under iteration of $f$. These two  estimates in particular show that $I_\infty(f)$  is a very thin subset of the Julia set $J(f)$.

\sp\bthm\label{thm:julia} 
If $f: \mathbb{C} \lra \ov{\mathbb{C}}$ is an elliptic function, then
$$
{\rm HD}(J(f))> \frac{2q_{\max}(f)}{q_{\max}(f)+1}\ge 1.
$$
\ethm

\fr {\sl  Proof.} We may assume that $R_1\ge R_0$ is so large that
\beq\lab{u5a}
LR_1^{-{1\over q_b}}<R_0
\eeq
for all $b\in f^{-1}(\infty)$.

Given two poles $b_1,b_2\in B_\infty(2R_1)$, denote by
$f_{b_2,b_1,j}^{-1}:B(b_1,R_0)\to\mathbb{C}$\index{(S)}{$f_{b_2,b_1,j}^{-1}$}, $1\le j\le q_b$, all the holomorphic inverse branches $f_{b_2,B(b_1,R_0),j}^{-1}$, $1\le j\le q_b$. It
follows from (\ref{u4}) and (\ref{u5a}) that
\beq\lab{u6a}
f_{b_2,b_1,j}^{-1}\(B(b_1,R_0)\) \sbt B_{b_2}(2R_1-R_0) \sbt
B_{b_2}(R_1) \sbt B(b_2,R_0).
\eeq
Fix a pole $a\in B_\infty(2R_1)$ with $q_a=q=q_{\max}(f)$. For every pole $b\in
B_\infty(2R_1)\cap f^{-1}(\infty)$ with $q_b=q$, fix the inverse branches
$$
f_{b,a,1}^{-1}:\ov B(a,R_0)\to\mathbb{C} \  \text{ and } \ 
f_{a,b,1}^{-1}:\ov B(b,R_0)\to\mathbb{C}
$$
of $f$. In view of (\ref{u6a}),
$$
f_{b,a,1}^{-1}\(\ov B(a,R_0)\)\sbt \ov B(b,R_0) \  \text{ and } \
f_{a,b,1}^{-1}\(\ov B(b,R_0)\)\sbt \ov B(a,R_0).
$$
The family
$$
S=\big\{f_{a,b,1}^{-1}\circ f_{b,a,1}^{-1}:   \,\,\,\ov B(a,R_0)\to\ov
B(a,R_0)\big\}_{b\in B_{2R_1}\cap f^{-1}(\infty)}
$$
thus forms a conformal infinite iterated function
system.\index{(N)}{conformal infinite iterated function system}
Given $t\ge 0$ we consider the function
$$
Z_1(t)=\sum_{b\in B_{2R_1}\cap f^{-1}(\infty)}||\phi_b'||^t
$$
and the number
$$
\th_S=\inf \{t\ge 0:\,  Z_1(t)<\infty\}.\index{(S)}{$\th_S$} 
$$
Our proof is based on demonstrating that $\th_S={2q\over q+1}$ and
$Z_1(\th_S)=+\infty$. In view of (\ref{1110301}), we can write:
$$
Z_1(t) \comp\sum_{b\in B_{2R_1}\cap f^{-1}(\infty)}|a|^{-{q+1\over
q}t}|b|^{-{q+1\over q}t} \\
\comp\sum_{b\in B_{2R_1}\cap f^{-1}(\infty)}|b|^{-{q+1\over q}t}
$$
But the series $\sum_{b\in B_{R_2}\cap
f^{-1}(\infty)}|b|^{-{q+1\over q}t}$ converges if and only if
$t>{2q\over q+1}$, and therefore the formulas
\beq\label{120120831}
\th_S={2q\over q+1}>1 \  \  \text{ and }  \  \ Z_1(\th_S)=\infty
\eeq
are proved. Applying now the second assertion of Theorem~\ref{HD-theta}
the desired  result follows. \endpf 

\sp

\fr As  an immediate consequence of this theorem we get the following.

\bcor \label{cor:sup} If $\La$ is a lattice in $\mathbb{C}$ and
${\mathcal K_\La}$ is the field of all elliptic functions with respect
to $L$, then 
$$
\sup\{{\rm HD}(J(f)):f\in {\mathcal K}_\La\}=2. 
$$
\ecor

\section[Hausdorff Dimension of Escaping Set]{Hausdorff
  Dimension of Escaping Sets of Elliptic Functions} 

We prove in this section that the Hausdorff dimension of the set of points escaping to infinity under the action of an elliptic function $f$ is precisely equal to 
$$
\frac{2q_{\max}(f)}{q_{\max}(f)+1}.
$$
We first provide an upper bound for this dimension and then the lower bound.
Along with the estimate of the previous section this
shows that $I_\infty(f)$, the set of points escaping to infinity, (unlike the case of exponential functions for example (see \cite{McM1}), is a very thin subset of the Julia set. We now recall, from formula \eqref{3_2017_09_25} of Section~\ref{FC B+S} the formal definition of $I_\infty(f)$. 
$$\index{(S)}{$I_{\infty}(f)$}
I_{\infty}(f)=\lt\{z\in {\mathbb C} :z\in \bu_{n\ge 0}f^{-n}(\infty)
\  \text{ or } \  \lim_{n\to \infty} f^n(z)= \infty \rt\},
$$  
As an immediate consequence of Theorem~\ref{t2_2017_09_25} and Theorem~\ref{o1_2017_09_25}, we get the following. 

\bthm\label{t3_2017_09_25}
If $f: \mathbb{C} \to \ov{\mathbb{C}}$ is a non-constant  elliptic
function, then $I_\infty(f)\sbt J(f)$. 
\ethm

\fr As announced, we start with the upper bound.

\blem\label{thm:infinity} 
If $f: \mathbb{C} \to \ov{\mathbb{C}}$ is a non-constant  elliptic
function, then 
$$
{\rm HD}(I_{\infty}(f))\leq \frac{2q_{\max}(f)}{q_{\max}(f)+1}.
$$ 
\elem

\bpf As in the previous section, put 
$$
q:=q_{\max}(f).
$$
Keep $R_1\ge R_0$ the same as in the proof of Theorem~\ref{thm:julia}, i.e. satisfying (\ref{u5a}). For every real number $R>0$, set
$$
I_R(f):=\{z\in\mathbb{C}: \forall_{n\ge 0}\, |f^n(z)|>
R\}.\index{(S)}{$I_R(f)$}
$$
Since the series $\sum_{b\in f^{-1}(\infty)\sms \{0\}}|b|^{-s}$
converges for all $s>{2q\over q+1}$, given $t>{2q\over q+1}$ there exists
$R_*\ge R_1$ such that
\beq\lab{u3a}
q(4A_2)^t\sum_{b\in B_{2R_*}\cap f^{-1}(\infty)}|b|^{-{q+1\over q}t}\le 1.
\eeq
Fix $R\ge 4R_*$. Put
$$
E:=f^{-1}(\infty) \cap B_\infty(R/2).
$$
It follows from (\ref{u4}), (\ref{u5a}), and \eqref{u6a} that for every $l\ge 1$ the family 
$$
\left\{f_{b_l,b_{l-1},j_l}^{-1}\circ
f_{b_{l-1},b_{l-2},j_{l-1}}^{-1} \ldots \circ
f_{b_1,b_0,j_1}^{-1}\(B_{b_0}(R/2)\):\,\, b_i\in E_R:1\le j_i\le
q_{b_i}, \,\,i=0,1,\ld,l\right\},
$$
denoted by $W_l$, is well-defined and covers $I_R(f)$. Applying \eqref{u2}, the second part of \eqref{u4}, and noting that $B_e(b,R_0)\sbt B_e(0,2|b|)$ for every $b\in E$, we may now estimate as follows:
$$
\aligned
&\Sg_l(R):=\index{(S)}{$\Sg_l$} \\
&=\sum_{b_l\in E}\sum_{j_l=1}^{q_{b_l}} \ld
  \sum_{b_1\in E}\sum_{j_1=1}^{q_{b_1}}
  \sum_{b_0\in E}
  \diam_s^t\lt(f_{b_l,b_{l-1},j_l}^{-1}\circ
  f_{b_{l-1},b_{l-2},j_{l-1}}^{-1} \ldots \circ f_{b_1,b_0,j_1}^{-1}\(B_{b_0}(R/2)\)\rt)\\
  &\le \sum_{b_l\in E}\sum_{j_l=1}^{q_{b_l}} \ld
  \sum_{b_1\in E}\sum_{j_1=1}^{q_{b_1}}
  \sum_{b_0\in E}
  \big\|\(f_{b_l,b_{l-1},j_l}^{-1}\circ
  f_{b_{l-1},b_{l-2},j_{l-1}}^{-1} \ldots \circ f_{b_1,b_0,j_1}^{-1}\)^\sh\big|_{B_{R_0}}\big\|_\infty^t\times
  \diam_s^t\(B_{b_0}(R/2)\)\\
  &\le \sum_{b_l\in E}\sum_{j_l=1}^{q_{b_l}} \ld
  \sum_{b_1\in E}\sum_{j_1=1}^{q_{b_1}}
  \sum_{b_0\in E}\lt(2A_2{(2|b_{l-1}|)^{q_{b_l}-1
  \over q_{b_l}}\over |b_l|^2}\rt)^t\cdot \lt(2A_2{(2|b_{l-2})|^{q_{b_{l-1}}-1
  \over q_{b_{l-1}}}\over |b_{l-1}|^2}\rt)^t\ld  \lt(2A_2{(2|b_0|)^{q_{b_1}-1
  \over q_{b_1}}\over |b_1|^2}\rt)^t \times L^t\lt({R\over 2}\rt)^{-{t\over q_{b_0}}} {1\over  |b_0|^{2t}} \\
 &\le L^t\lt({2\over R}\rt)^{{t\over q}}(4A_2)^{lt}
  \sum_{b_l\in E}\sum_{j_l=1}^{q_{b_l}} \ld \sum_{b_1\in I}\sum_{j_1=1}^{q_{b_1}}\sum_{b_0\in E}
  |b_l|^{-2t}\(|b_{l-1}|^{-{q+1\over q}t}\ld |b_0|^{-{q+1\over
  q}t}\)\\
&\le L^t\lt({2\over R}\rt)^{{t\over q}}(4A_2)^{lt}
  \sum_{b_l\in E}\sum_{j_l=1}^{q_{b_l}} \ld
  \sum_{b_1\in E}\sum_{j_1=1}^{q_{b_1}}\sum_{b_0\in E}
  \(|b_l|^{-{q+1\over q}t}|b_{l-1}|^{-{q+1\over q}t}\ld
  |b_0|^{-{q+1\over q}t}\)\\
 &\le L^t\lt({2\over R}\rt)^{{t\over q}}(4A_2)^{lt}\lt(\sum_{b\in
  E}|b|^{-{q+1\over q}t}\rt)^lq^l \\
&\le L^t\lt({2\over R}\rt)^{{t\over q}}\lt(q(4A_2)^t\sum_{b\in
B_{2R_*}\cap f^{-1}(\infty)}|b|^{-{q+1\over q+}t}\rt)^l.
\endaligned
$$
Applying (\ref{u3a}) we therefore get $\Sg_l\le L^t(2/R)^{t/q}$.
Since the diameters (in the spherical metric) of the sets of the
covers $W_l$ converge uniformly to $0$ as $l\downto\infty$, we
therefore infer that $\H_s^t(I_R(f))\le L^t(2/R)^{t/q}$, where the
subscript $s$ indicates that the Hausdorff measure is considered
with respect to the spherical metric. Consequently $\HD(I_R(f))\le
t$ and if we put
$$
I_{R,e}(f):=\lt\{z\in\mathbb{C}:\liminf_{n\to\infty}|f^n(z)|>R\rt\}
=\bu_{k\ge 1}f^{-k}(I_R(f)),
$$\index{(S)}{$I_{R,e}(f)$}
then also $\HD(I_\infty(f))\le \HD(I_{R,e}(f))= \HD(I_R(f))\le
t$. Letting now $t\downto {2q\over q+1}$ finishes the proof.
\endpf

\

\fr As an immediate consequence of this theorem  we obtain the
following.

\bcor \label{cor:lebesgue} 
If $f: \mathbb{C} \to \ov{\mathbb{C}}$ is
an elliptic function, $h:=\HD(J(f))$, then $\H^h(I_{\infty}(f))=0$,
and consequently $S(I_{\infty}(f))=0$, $S$ being planar Lebesgue measure on $\C$. 
\ecor

This corollary and the previous lemma show that the escaping set $I_{\infty}(f)$ is a fairly small subset of the Julia set $J(f)$. Now we shall prove the opposite inequality.

\blem\label{thm:infinityB} 
If $f: \mathbb{C} \to \ov{\mathbb{C}}$ is a non-constant  elliptic
function, then 
$$
\HD(I_{\infty}(f))\ge \frac{2q_{\max}(f)}{q_{\max}(f)+1}.
$$ 
\elem

\fr{\sl Proof.}
As in the previous proof, keep
$$
q:=q_{\max}(f). 
$$
Let $f^{-1}_q(\infty)$ be the set of all poles of $f$  of order $q$. Fix $R_0$ so  large as required in formulas \eqref{1_2017_10_02}, (\ref{u1}), (\ref{u2}) and  (\ref{1110301}). Then fix 
\beq\label{1_2017_10_04}
R \geq 4\max\{R_0,\diam(\mathcal R_F)\}.
\eeq
Our goal is to apply Theorem~\ref{p1_mcmullen_HD_estimate}. We perform an inductive  construction, required  by this  theorem, as  follows. Fix $\xi \in f^{-1}_q(\infty)$.  Let
$$ 
E_0:= \hat{\mathbb C}.
$$
As an inductive  assumption, suppose  that  for some $ n \geq 1$  a collection $ E_0, \ldots, E_{n-1}$ of mutually disjoint compact connected subsets of $E_0$ has been defined  with the  following properties:

\begin{itemize}
\item [(a)] For every $F\in E_n$, there exists a unique set $F_{-}\in E_{n-1}$ such that $F\sbt F_{-}$. The set $F_{-}$ will  be referred  to as  the parent of $F$  and $F$ as a  child of $F_{-}$.

\,\item [(b)]  If $F\in  E_n$, then  there exists a unique pole in $(\xi+\Lambda_f)\cap (\bar{B}_\infty(2^{n+1}A_1R^q)\sms B_\infty(2^{n+2}A_1R^q))$, denoted by $b_F$, such that $f^n(F)=\bar{B}(b_F,R^{-1})$ and there exists a compact  connected set $\hat{F}$  containing $F$  such that

\,\item[(c)] $f^n(\hat{F})=\bar{B}(b_F, 2R^{-1})$  and
   
\, \item[ (d)] $f^n|_{\hat{F}}:\hat{F} \to \bar{B}(b_F,2R^{-1})$ is 1-to-1 (and consequently a homeomorphism).
\end{itemize}

\, \fr As the inductive step, for every $F\in E_n$ define $E_{n+1}(F)$ to be the collection of all the sets of the form
\beq\label{1j325}
F_b:=(f^n|_F)^{-1}(f^{-1}_{b_F, B(b,1), 1}(\bar{B}(b, R^{-1})))\sbt F
\eeq
and
\beq\label{3j325}
\hat{F_b}:=(f^n|_{\hat{F}})^{-1}(f^{-1}_{b_F, B(b,1), 1}(\bar{B}(b, 2 R^{-1})))\sbt \hat{F},
\eeq
where $b$ ranges over all elements of the set 
\beq\label{2_2017_10_04}
(\xi+\Lambda_f)\cap (\bar{B}_\infty(2^{n+2}A_1R^q)\sms B_\infty(2^{n+3}A_1R^q)).
\eeq
Obviously all  elements of $E_{n+1}(F)$ are  compact and mutually disjoint. Let  
$$
\Ga:=f(\bar{B}(b_F, R^{-1})).
$$
Note that $\Gamma= f(\bar{B}(a,R^{-1}))$ for all $a \in (\xi+\Lambda_f)$ and, by (\ref{u1}),
$$
\Gamma 
\supset B_\infty \(A_1(R^{-1})^{-q}\) 
= B_\infty (A_1 R^q)
\supset \bar{B}(b,2R^{-1}) 
$$  
for every 
$$
b \in (\xi+\Lambda_f)\cap \bar{B}_\infty(2A_1R^q)
\spt(\xi+\Lambda_f)\cap (\bar{B}_\infty(2^{n+2}A_1R^q)\sms B_\infty(2^{n+3}A_1R^q)).
$$
This along with (a) implies that for such $b$s, we have that $f^{n+1}(F_b)=\bar{B}(b, R^{-1})$,  $f^{n+1}(\hat{F}_b)=\bar{B}(b, 2R^{-1})$, and, obviously, the map $f^{n+1}|_{\hat{F}_b} $  is 1--to--1. Thus defining
$$
E_{n+1}:= \bigcup_{F\in E_n}\big\{F_b: \, b \in ( \xi+\Lambda_f) \cap
(\bar{B}_\infty(2^{n+2}A_1R^q)\sms B_\infty(2^{n+3}A_1R^q)),
$$
finishes our inductive construction so that (a), (b), (c) and (d) hold.

\,

Fix $n \geq 1$ and $F\in E_n$. The conditions (a)--(d) imply that Koebe's  Distortion Theorem (Theorem~\ref{Euclid-I})
applies to the map $f^{n-1}|_{\hat{F}_{-}}$ and
$f^{n}|_{\hat{F}}$. We therefore  get from (\ref{1j325})  and (\ref{1110301}) that
\beq\label{2j325}
\begin{aligned}
\frac{\diam_e(F)}{\diam_e(F_{-})}& =
\frac{\diam_e\((f^{n-1}|_{F_{-}})^{-1}(f^{-1}_{b_{F_{-}}, B(b_{F},1),1} (\bar{B}(b_F, R^{-1})))\)}{\diam_e\((f^{n-1}|_{F_{-}})^{-1} \(\bar{B}\(b_{F_{-}}, R^{-1}\)\)\)}\\
& \leq   K \frac{\diam_e\(f^{-1}_{b_{F_{-}}, B(b_{F},1),1} (\bar{B}(b_F, R^{-1}))\)}{\diam_e\(\bar{B}\(b_{F_{-}}, R^{-1})\)} \\
& \leq   K  A_2 (|b_F|-R^{-1})^{-\frac{q+1}{q}}R^{-1}\\
& \leq K A_2(2^nA_1R^q)^{-\frac{q+1}{q}}R^{-1}\\
&=\g_12^{-\frac{q+1}{q}n},
\end{aligned}
\eeq
where
$$
\g_1:=KA_2A_1^{-\frac{q+1}{q}}R^{-(q+2)}.
$$
So, an immediate induction yields
\beq\label{3j325B}
\diam_e(F)\leq \gamma^n_1 \prod_{k=1}^n 2^{-\frac{q+1}{q}k}.
\eeq 
Hence, with the notattion of Theorem~\ref{p1_mcmullen_HD_estimate}, we get that
\beq\label{2j327}
d_n=\max\{\diam_e(F):  F \in E_n\} \leq \gamma^n_1 \prod_{k=1}^n 2^{-\frac{q+1}{q}k}.
\eeq
So, in conclusion, $\{E_n\}_{n=1}^\infty$ is a McMullen's sequence, and Theorem~\ref{p1_mcmullen_HD_estimate} applies indeed. Still keeping $F\in E_n$, by the same arguments as those generating (\ref{2j325}), we  get
\beq\label{1j327}
\begin{aligned}
\frac{S(F)}{S(F_{-})} 
& = \frac{S\((f^{n-1}|_{F_{-}})^{-1}\(f^{-1}_{b_{F_{-}}, B(b_{F},1),1} (\bar{B}(b_F, R^{-1}_0))\)\)}{S\((f^{n-1}|_{F_{-}})^{-1} (B(b_{F_{-}},  R^{-1}))\)}\\
& \geq  K^{-2}\frac{S\(f^{-1}_{b_{F_{-}}, B(b_{F},1),1} (\bar{B}(b_F, R^{-1}))\)}{S\(B(b_{F_{-}},R^{-1})\)}\\
& \geq K^{-2} A_2^{-2} (|b_F|+R^{-1})^{-2\frac{q+1}{q}}\\
& \geq K^{-2}A_2^{-2}( 2^{n+3}A_1R^q)^{-2 \frac{q+1}{q}}\\
&= \gamma^{-2}_2 4^{-\frac{q+1}{q}n},
\end{aligned}
\eeq
where $\gamma_2:=KA_2A_1^{2\frac{q+1}{q}}R^{2(q+1)}64^{\frac{q+1}{q}}$ and $S$, as always, denotes planar Lebesgue measure on $\C$. Since  with some  $\gamma_3>0$ there are at least $\gamma_34^n$ (see the formulas \eqref{1_2017_10_04} and \eqref{2_2017_10_04}) children of each element in $E_n$, we thus get from (\ref{1j327}) that for each  $G \in E_n$, we have
$$ \Delta_n(G)
:= \frac{S\(\bigcup_{H\in E_{n+1 (G)}} H\)}{S(G)} 
\geq \gamma_3 \gamma_2^{-2} 4^n4^{-\frac{q+1}{q}n}=
\gamma_3 \gamma_2^{-2}4^{-\frac{n}{q}}.
$$
Therefore,
$$ 
\Delta_n:= \min\{ \Delta_n(G): G \in E_n\}
\geq \gamma_3\gamma_2^{-2}4^{-\frac{n}{q}}.
$$
Combining this with (\ref{2j327}), we thus  get
$$ 
\varlimsup_{n \to  \infty}\frac{\sum_{k=1}^{n-1} \log \Delta_k}{\log d_n} 
\leq\varlimsup_{n\to\infty}\frac{(n-1)\log(\gamma_3\gamma_2^{-2})-\frac{1}{q}\log 4\sum_{k=1}^{n-1} k}{n\log \gamma_1-\frac{q+1}{q}\log 2 \sum_{k=1}^{n-1} k}
=\frac{\frac{1}{q}\log4}{\frac{q+1}{q}\log2}=\frac{2}{q+1}.
$$
Thus
\beq\label{3j327}
2 -\varlimsup_{n \to  \infty} \frac{\sum_{k=1}^{n-1} \log \Delta_n}{\log d_n}
\geq  2 - \frac{2}{q+1}= \frac{2q}{q+1}.
\eeq
As in Theorem~\ref{p1_mcmullen_HD_estimate}, denote by $E_\infty$ the set 
$$
\bigcap_{n=1}^\infty\bigcup_{F\in E_n} F.
$$
It follows from Theorem~\ref{p1_mcmullen_HD_estimate}, applied with planar Lebesgue measure $S$ on $\C$, and (\ref{3j327}) that
$\HD( E_\infty)\geq \frac{2q}{q+1}$. Since also, by (b) ($|b_F| \geq 2^{n+1}A_1R^q$ for all $F\in E_n$) and by (c), we have that $E_\infty\sbt I_\infty(f)$, we  thus conclude that 
$$
\HD(I_\infty(f)) \geq \frac{2q}{q+1}.
$$
The proof is complete.
\endpf

\sp Along with Lemma~\ref{thm:infinity} this gives the following main result of the current section. 

\bthm\label{thm:infinityC} 
If $f: \mathbb{C} \lra\oc$ is a non--constant  elliptic
function, then 
$$
\HD(I_{\infty}(f))=\frac{2q_{\max}(f)}{q_{\max}(f)+1}\in [1,2).
$$ 
\ethm

\sp

\section[Conformal Measures of Escaping Sets]{Conformal Measures
  of Escaping Sets of Elliptic Functions} 

In Chapter~\ref{conformal-measure}, we will start presenting a systematic account of the theory of conformal measures for non--recurrent elliptic functions. In this section we still deal with arbitrary elliptic functions and we will "only" prove that any conformal measure for any elliptic function vanishes on its set of escaping points. This fact is interesting on its own and will be needed in the proof of ergodicity and conservativity of conformal measures (see Theorem~\ref{tmaincm}). Its proof is similar to the proof of Theorem~\ref{thm:infinity}. Let 
\beq\label{120131007}
I_-(f):=\bu_{n\ge 1}f^{-n}(\infty)\index{(S)}{$I_-(f)$}.
\eeq
Rcall that conformal measures have been introduced in Definition~\ref{1d20120909}. 
The result announced to be proved in this section is the following: 

\blem\lab{liinfty} 
If $m$ is a $t$--conformal measure for an elliptic function $f:{\mathbb
C}\lra \oc$, then 
$$
m(I_\infty(f)\sms I_-(f))=0.
$$ 
Furthermore, there exists $R>0$ such that
$$
m(\{z\in\C:\liminf_{n\to\infty}|f^n(z)|>R\})=0.
$$
\elem

\bpf  As in the previous sections we denote
$$
q=q_{\max}(f).
$$
It suffices to prove the lemma for the spherical
measure $m_s$. Let $b$ be a pole of $f:{\mathbb C}\to \ov{\mathbb
C}$. We shall obtain first an upper estimate on $m_s(B_b(R))$
similar to the second inequality in (\ref{u4}). And indeed, considering two
open connected simply connected sets
$$
B_R^+=\{z\in B_\infty^*(R):\im z>0\} \ \text{ and } \
B_R^-=\{z\in B_\infty^*(R)\sms:\im z<1\}
\index{(S)}{$B_R^+$}\index{(S)}{$B_R^1$}
$$
for every $j=1,2,\ld,q_b$, we obtain
\beq\label{1_2017_09_27}
1\ge m_s\lt(f_{b,B_R^{\pm},j}^{-1}\(B_R^{\pm}\)\rt)
=\int_{B_R^{\pm}}\big|(f_{b,B_R^{\pm},j}^{-1})^*\big|^tdm_s.
\eeq
Using now also (\ref{u2}), we obtain
\beq\label{2_2017_09_27}
1\ge \int_{B_R^{\pm}}\big|(f_{b,B_R^{\pm},j}^{-1})^*\big|^tdm_s 
\ge\int_{B_R^{\pm}}\lt({(2A_2)^{-1}\over 1+|b|^2}|z|^{{q_b-1\over q_b}}\rt)^tdm_s(z)
={(2A_2)^{-t}\over (1+|b|^2)^t}\int_{B_R^{\pm}}|z|^{{q_b-1\over q_b}t}dm_s(z) 
\eeq
and
\beq\label{3_2017_09_27}
\aligned 
\int_{B_R^{\pm}}\big|(f_{b,B_R^{\pm},j}^{-1})^*\big|^tdm_s 
&\le\int_{B_R^{\pm}}\lt({2A_2\over 1+|b|^2}|z|^{{q_b-1\over q_b}}\rt)^tdm_s(z)
 ={(2A_2)^t\over (1+|b|^2)^t}\int_{B_R^{\pm}}|z|^{{q_b-1\over q_b}t}dm_s(z) \\
&\le (2A_2)^t(1+|b|^2)^{-t}\int_{B_R^+}|z|^{{q-1\over q}t}dm_s(z)\\
&\le (2A_2)^t|b|^{-2t}\int_{B_R^+}|z|^{{q-1\over q}t}dm_s(z).
\endaligned
\eeq
Taking any pole $b$ of maximal multiplicity, i.e. with $q_b=q$, we see from \eqref{2_2017_09_27} that 
$$
\int_{B_R^\pm}|z|^{{q-1\over q}t}dm_s(z)
\le (2A_2)^t(1+|b|^2)^t<+\infty.
$$ 
Put
$$
\Sg_R:=\max\lt\{\int_{B_R^+}|z|^{{q-1\over q}t}dm_s(z),
\int_{B_R^-}|z|^{{q-1\over q}t}dm_s(z)\rt\}\in [0,+\infty).
$$
Since $B_R^*=B_R^+\cup B_R^-$, we obtain
$$
B_b^*(R)
=\bu_{j=1}^{q_b}f_{b,B_R^{\pm},j}^{-1}\(B_R^+\)\cup \bu_{j=1}^{q_b}f_{b,B_R^{\pm},j}^{-1}\(B_R^-\).
$$
Hence, using also \eqref{1_2017_09_27} and \eqref{3_2017_09_27}, we obtain for any pole $b$ of $f$ that
\beq\lab{diinfty.2}
m_s(B_b^*(R)) 
\le \sum_{j=1}^{q_b}\int_{B_R^+}\big|(f_{b,B_R^+,j}^{-1})^*\big|^tdm_s
 + \sum_{j=1}^{q_b}\int_{B_R^-}\big|(f_{b,B_R^-,j}^{-1})^*\big|^tdm_s 
\le 2(2A_2)^tq\Sg_R|b|^{-2t}.
\eeq
Now the argument goes essentially in the same way as in the proof of Lemma~\ref{thm:infinity}. We present it here for the sake of completeness. 
Keep $R_1\ge R_0$ the same as in the proof of Theorem~\ref{thm:julia}, i.e. satisfying (\ref{u5a}). Recall that for every real number $R>0$, we put in the proof of Lemma~\ref{thm:infinity}:
$$
I_R(f)=\{z\in\mathbb{C}: \forall_{n\ge 0}\, |f^n(z)|>
R\}.\index{(S)}{$I_R(f)$}
$$
Since the series $\sum_{b\in f^{-1}(\infty)\sms \{0\}}|b|^{-s}$
converges for all $s>{2q\over q+1}$, given $t>{2q\over q+1}$ there exists
$R_*\ge R_1$ such that
\beq\lab{u3aB}
q(4A_2)^t\sum_{b\in B_{2R_*}\cap f^{-1}(\infty)}|b|^{-{q+1\over q}t}\le 1/2.
\eeq
Fix $R\ge 4R_*$. Put
$$
E:=f^{-1}(\infty) \cap B_{(R/2)}.
$$
It follows from (\ref{u4}), (\ref{u5a}), and \eqref{u6a} that for every $l\ge 1$ the family 
$$
\left\{f_{b_l,b_{l-1},j_l}^{-1}\circ
f_{b_{l-1},b_{l-2},j_{l-1}}^{-1} \ldots \circ
f_{b_1,b_0,j_1}^{-1}\(B_{b_0}(R/2)\):\,\, b_i\in E_R:1\le j_i\le
q_{b_i}, \,\,i=0,1,\ld,l\right\},
$$
denoted by $W_l$, is well-defined and covers $I_R(f)$. Applying \eqref{u2}, (\ref{diinfty.2}), and noting that $B_e(b,R_0)\sbt B_e(0,2|b|)$ for every $b\in E$, we may now estimate as follows:
$$
\aligned
&m_s(I_R(f))\le \\
&\le\sum_{b_l\in E}\sum_{j_l=1}^{q_{b_l}} \ld
  \sum_{b_1\in E}\sum_{j_1=1}^{q_{b_1}}
  \sum_{b_0\in E}
  m_s\lt(f_{b_l,b_{l-1},j_l}^{-1}\circ
  f_{b_{l-1},b_{l-2},j_{l-1}}^{-1}\circ\ld\circ
  f_{b_2,b_1,j_2}^{-1}\circ f_{b_1,b_0,j_1}^{-1}\(B_{b_0}(R/2)\)\rt)\\
&\le \sum_{b_l\in E}\sum_{j_l=1}^{q_{b_l}} \ld
  \sum_{b_1\in E}\sum_{j_1=1}^{q_{b_1}}
  \sum_{b_0\in E}
  \big\|\(f_{b_l,b_{l-1},j_l}^{-1}\circ
  f_{b_{l-1},b_{l-2},j_{l-1}}^{-1}\circ\ld\circ
  f_{b_2,b_1,j_2}^{-1}\circ
  f_{b_1,b_0,j_1}^{-1}\)^*\big|_{B_{b_0}(R/2)}\big\|_\infty^t\times \\
&  \  \  \   \   \   \   \  \  \   \   \   \   \  \  \   \   \   \   \  \  \   \   \   \   \times m_s\(B_{b_0}(R/2)\) \\
&\le \sum_{b_l\in E}\sum_{j_l=1}^{q_{b_l}} \ld
  \sum_{b_1\in E}\sum_{j_1=1}^{q_{b_1}}
  \sum_{b_0\in E}(4A_2)^{lt}\lt({|b_{l-1}|^{q_{b_l}-1
  \over q_{b_l}}\over |b_l|^2}\rt)^t\cdot \lt({|b_{l-2}|^{q_{b_{l-1}}-1
  \over q_{b_{l-1}}}\over |b_{l-1}|^2}\rt)^t\ld \lt({|b_0|^{q_{b_1}-1
  \over q_{b_1}}\over |b_1|^2}\rt)^t\times \\
& \  \  \   \   \   \   \  \  \   \   \   \   \  \  \   \   \   \   \  \  \   \   \   \   \times 2q(2A_2)^t\Sg_R|b_0|^{-2t} \\
&= 2q(2A_2)^t\Sg_R(4A_2)^{lt}
  \sum_{b_l\in I}\sum_{j_l=1}^{q_{b_l}} \ld
  \sum_{b_1\in I}\sum_{j_1=1}^{q_{b_1}}\sum_{b_0\in E}
  |b_l|^{-2t}\(|b_{l-1}|^{-{q+1\over q}t}\ld |b_0|^{-{q+1\over q}t}\) \\
&\le 2q(2A_2)^t\Sg_R(4A_2)^{lt}
  \sum_{b_l\in I}\sum_{j_l=1}^{q_{b_l}} \ld
  \sum_{b_1\in I}\sum_{j_1=1}^{q_{b_1}}\sum_{b_0\in E}
  \(|b_l|^{-{q+1\over q}t}|b_{l-1}|^{-{q+1\over q}t}\ld
  |b_0|^{-{q+1\over q}t}\) \\
&\le 2q(2A_2)^t\Sg_R(4A_2)^{lt}\lt(\sum_{b\in E}|b|^{-{q+1\over q}t}\rt)^lq^l\\
&\le 2q(2A_2)^t\Sg_R\lt(q(4A_2)^t\sum_{b\in B_{2R_*}\cap f^{-1}(\infty)}
|b|^{-{q+1\over q}t}\rt)^l.
\endaligned
$$
Applying (\ref{u3aB}) we therefore get $m_s(I_R(f)) \le 2q(2A_2)^t\Sg_R2^{-l}$. Letting $l\to\infty$ we thus get $m_s(I_R(f))=0$. Since $m_s\circ f^{-1}$ is absolutely continuous with respect to $m_s$ and since
$$
\big\{z\in\C:\liminf_{n\to\infty}|f^n(z)|>R\big\}=\bu_{j=0}^\infty
f^{-j}(I_R(f)),
$$ 
we conclude that
$$
m_s\(\{z:\liminf_{n\to\infty}|f^n(z)|>R\}\)=0.
$$ 
The proof is complete. \endpf

\chapter{Compactly Non--Recurrent Elliptic Functions: First Outlook}\label{first-outlook}

In this chapter we define the class of non--recurrent and, more notably, the  class of compactly non--recurrent elliptic functions. This is the class of elliptic functions which will be dealt with by us since now through the end of the book in greatest detail. Its history goes back to the papers \cite{U1} and \cite{U2} by the second named author, and \cite {KU3}. One should also mention the paper \cite{CJY}. Similarly, as in all these papers, our treatment of non--recurrent elliptic functions is based on, in fact possible at all, due to an appropriate version of the breakthrough Ma\~ne's Theorem proven in \cite{M1} in the context of rational functions. In our setting of elliptic functions this is Theorem~\ref{mnt6.3a}. The first section of the current chapter is entirely devoted to proving this theorem, its first most fundamental consequences, and some other results surrounding it. The next two sections of this chapter, also relying on Ma\~ne's Theorem, provide us with further refined technical tools to study the structure of Julia sets and holomorphic inverse branches. 

The last section of this chapter has somewhat different character. It systematically defines and describes various subclasses of, mainly compactly non--recurrent, elliptic functions we will be dealing with in Part~4 of the book. Mostly, these classes of elliptic functions will be defined in terms of how strongly expanding these
functions are. We would like to add that while in the theory of rational functions such classes pop up in a natural and fairly obvious way, and for example metric and topological definitions of expanding rational functions describe the same class of functions, in the theory of iteration of transcendental meromorphic functions such classification is by no means obvious, topological and metric analogs of rational function realm concepts do not usually coincide, and the definitions of expanding, hyperbolic, topologically hyperbolic, subhyperbolic, etc, functions vary from author to author. Our definitions seem to us pretty natural and fit well for our purpose of detailed investigation of dynamical and geometric properties of elliptic functions they define.  

\section[Fundamentals of Non--Recurrent Elliptic Functions; Mane's Theorem] {Fundamental Properties of Non--Recurrent Elliptic Functions; \\ Mane's Theorem}\label{FPoN-REF MT}

In this section we prove results forming fundamental dynamical and topological properties of non--constant elliptic functions. They stem from the breakthrough Ma\~ne's Theorem proven in \cite{M1} in the context of rational functions. Our results in this section are based on an appropriate version of Ma\~ne's Theorem for elliptic functions, i. e. Theorem~\ref{mnt6.3a}. This theorem, roughly speaking, asserts that the connected components of sufficiently small balls centered at points of Julia sets, different from rationally indifferent periodic points, and not belonging to the omega limit set of recurrent critical points, have small diameters. It has an enormously powerful consequences for the structure of such components, gives some sort of hyperbolicity, and it is, along with its technical consequences, an indispensable tool for us throughout the rest of the book. In particular, see Theorem~\ref{r071708}, and this is not the most significant consequence of Ma\~ne's Theorem, this theorem rules out the existence of Siegel disks, Herman rings, and Cremer points for all elliptic functions.

Denote by $\Crit_r(f)$ \index{(S)}{$\Crit_r(f)$} the set of all {\em recurrent critical points} \index{(N)}{recurrent critical points} of an elliptic function $f:\C\lra\oc$, i.e. the set of such critical points $c$ of $f$ that
$$
c\in\om(c).
$$
Note that except for periodic attracting points, all non-recurrent critical of $f$ are contained in the Julia set $J(f)$ of $f$. Essentially, from now on throughout the entire book, we will deal with the class of 
non--recurrent elliptic functions. Here is their definition.

\bdfn\label{pseudo-non-recurrent_B} 
We say that a non--constant elliptic
function $f:\C\lra\oc$ is {\em non--recurrent} \index{(N)}{non--recurrent} (NR) \index{(S)}{(NR)} if and only if
$$
\Crit_r(f)\cap J(f)=\es.
$$
\edfn

\fr Directly from Theorem~\ref{baker-domain+Sullivan for elliptic} and Theorem~\ref{Fatou Periodic Components}, we obtain the following.

\sp\bobs\label{o1_2017_09_27} 
If $f:\C\lra\oc$ is a non--recurrent elliptic function and if $c \in \Crit(f) \cap F(f)$, then  there exists  either
an attracting periodic point  $\om$ of $f$ or a  rationally
indifferent periodic point $\om$ of $f$, in either case such that 
$$
\om(c)\sbt \{f^{n}(\om): n \geq 0\}.
$$
\eobs 

\sp\fr The class of non--recurrent elliptic function has several important subclasses. These are particularly transparent in the realm of rational functions on the Riemann sphere $\oc$. In the context of meromorphic transcendental functions, in particular, elliptic functions, their analogues are not uniquely obvious and definitions vary from author to author and from paper to paper. We will primarily deal with the subclass of compactly non--recurrent elliptic functions and some of its subclasses such as various subclasses of subexpanding
and parabolic elliptic functions, which will be defined in the next section.

\sp As was already said, similarly as in the  paper \cite{KU3}, the basic indispensable technical tool in the present book for our dealing with geometry and dynamics of non-recurrent elliptic functions, is an appropriate version of  Ma\~ne's
Theorem\index{(N)}{Ma\~ne's Theorem} originally proved by him, in \cite{M1} for the class of rational functions. We start with several preparatory results, then we prove the elliptic version of  Ma\~ne's Theorem, and then we derive many of its consequences for recurrent elliptic functions. 

First observe that if $f:\C\to\oc$ is elliptic, then for every $\e>0$ there exists $\d>0$ such that for every $z\in\oc$ every connected component of $f^{-1}(B_s(z,\d))$ has the Euclidean diameter smaller than $\e$. Thus, by an obvious induction, we get the following.

\bobs\label{o1_2017_10_04}
If $f:\C\lra\oc$ is an elliptic function, then 

\,

\fr $\forall\,\e>0 \  \forall\, n\ge 0\  \exists\, \d>0\ \forall\, z\in\oc\  \forall (0\le k\le n)$ 

\centerline{every connected component of $f^{-k}(B_e(z,\d))$ has the Euclidean diameter smaller than $\e$.}

\, 

and

\,

\fr $\forall\, z\in\oc\  \forall (1\le k\le n)$ 

\centerline{every connected component of $f^{-k}(B_s(z,\d))$ has the Euclidean diameter smaller than $\e$.}
\eobs

From now on if $f:\C\lra\oc$ is an elliptic function, then we fix $\eta_0(f)>0$ to be so \index{(S)}{$\eta_0(f)}small that the following two conditions are satisfied:
\begin{enumerate}
\item
$$
\min\big\{|w-z|:w,z\in \Crit(f)\cup f^{-1}(\infty): w\ne z\big\}>2\eta_0(f),
$$
\item For every point $z\in\oc$, the Euclidean diameter of each connected component of $f^{-1}\(B_s(z,8\eta_0(f))\)$ is finite.

\, \item If $U$ is an open connected simply conntected subset of $\oc$ with $\diam_s(U)\le \eta_0(f)$, then each connected component of $f^{-1}(U)$ is simply connected.
\end{enumerate}

\sp\fr Having (2) above, the simple observation (3) follows from Corollary~\ref{General Riemann Hurwitz Simply Connected_2} and Corollary~\ref{General Riemann Hurwitz Simply Connected_3} since the set $f(\Crit(f))$ is finite. It will turn out to be useful many times in the sequel. We will frequantly make use of it without explicite invoking. 

Recall that for any set $A\sbt  {\mathbb C}$, we denote \index{(S)}{$O_+(A)$}
$$
O_+(A)=\bu_{n \geq 0}f^n(A).
$$

Now, we shall prove an appropriate version of Przytycki's Lemma, stated and proved in the realm of rational functions, in \cite{P2}. This is an important ingredient of the proof of Theorem~\ref{mnt6.3a}. Note that there are places in the several proofs we provide below, where one has to proceed subtler than in the case of rational functions. Let 
$$
{\rm Attr}(f)
$$
denote \index{(S)}{Attr$(f$)} the set of all attracting periodic points of an elliptic function $f$.

\blem\lab{felman} 
If $f:\C\lra\oc$ is an elliptic function, then 

\,

\fr $\forall N\in\N\ \forall \lambda\in(0,1)\ \forall\varepsilon >0\ \forall\ka>0 \ \exists\, \delta_0>0\ \forall\delta\in(0,\delta_0]$ if
$$
x\in {\mathbb  C}\sms B_e\({\rm Attr}(f)\cup\Om(f),\ka\),
$$ 
then for every integer $n \geq 0$ and every connected
component $W$ of $f^{-n}(B_e(x,\delta))$ such that 

\,\,
\centerline{$f^n|_W$ has at most $N$ critical points,}
\,\,
\fr each connected component $W'$ of $f^{-n}(B_e(x,\lambda\delta))$ contained in $W$, satisfies
$$ 
\diam_e(W')\leq \varepsilon 
$$
and 
$$
\lim_{n\to\infty}\diam_e(W') = 0
$$ 
uniformly with respect to $x$ and $W'$. 
\elem

\bpf Because of Obervation~\ref{o1_2017_10_04}, it suffices to prove the second, i.e. the last assertion of this lemma. Proving it, suppose, for a contrary, that there exist 
\begin{enumerate}
\item a sequence $\{x_n\}_{n=1}^\infty$ of points belonging to ${\mathbb  C}\sms B_e\({\rm Attr}(f)\cup\Om(f),\ka\)$,

\, \item a sequence $\d_n\downto 0$, 

\, \item a sequence $\{k_n\}_{n=1}^\infty$ of positive integers diverging to infinity,

\,\item a sequence of connected components $W_n$ contained in $f^{-k_n}(B_e(x_n, \d_n))$ such that the number of critical  points (counted with  multiplicities) of each map $f^{k_n}$ on $W_n$ is bounded above by $N$ and 

\, \item a sequence $\{W_{n}'\}_{n=1}^\infty$ of connected components of $f^{-n}(B_e(x_n,\lambda\delta_n))$ contained in $W_n$ such that the following limit exists and 
$$
\lim_{n\to\infty}\diam_e(W_n')>0.
$$
\end{enumerate}
Then, for each integer $n\ge 1$, there exists an integer $L=L(n)\in [0, N]$,  such that the annulus 
$$
A(n):=B_e\left(x_n,\delta_n\left(\lambda+(1-\lambda)\frac{L+1}{K+1}\right)\right)
\sms
B_e\left(x_n,\delta_n\left(\lambda+(1-\lambda)\frac{L}{K+1}\right)\right)
$$
contains no  critical value of $f^{k_n}_{|W_n}$. We may assume, without loss of generality, that all the components $W_n'$ intersect the fundamental region $\mathcal {R}_f$ of $f$. Let $W^{(1)}_n$ and $W^{(2)}_n$ be the connected components, respectively of 
$$
f^{-k_n}\lt( B_e\lt(x_n,\delta_n\lt(\lambda+(1-\lambda)\frac{L}{K+1}\rt)\rt)\rt)
$$
and
$$
f^{-k_n}\lt( B_e\left(x_n,\delta_n\left(\lambda+(1-\lambda)\frac{L+1}{K+1}\right)\rt)\right)
$$
containing $W'_n$. Put 
$$
A_n:=W^{(2)}_n\sms  W^{(1)}_n,
$$
and  for all $0\leq m\leq k_n$ and  $i=1,2$, let
$$ 
W^{(i)}_{n,m}:=f^{k_n-m}(W^{(i)}_n),  \quad  
A_{n,m}:=f^{k_n-m}(A_n)=W^{(2)}_{n,m}\sms  W^{(1)}_{n,m}.
$$
For each integer $n\ge 0$, let  $0\le m=m(n)\leq k_n$  be the  least  integer  such
that
$$
{\rm diam}_e(W^{(1)}_{n,m})\geq \inf \big\{{\rm dist}_e(c_1, c_2); \,\,
c_1, c_2\in \Crit(f), \,c_1\neq c_2\big\}>0.
$$
Therefore, for every integer $0\leq t < m(n)$, at each step back, by $f^{-1}$, from $A_{n,t-1}$ to $A_{n,t}$, there is at most one branch point for 
$f^{-1}$ acting from $W^{(i)}_{n, t-1}$ to $W^{(i)}_{n, t}$, $i=1,2$. Hence, for every integer $0\leq t < m(n)$, the set $A_{n,t}$ is a topological
annulus.

Now, all the annuli $A_{n,m(n)-1}$ have moduli  bounded  below by
$2^{-N}(1-\lambda)(L+1)^{-1}$. Since, in addition, all the connected
components $W_n'$ intersect the fundamental region $\mathcal R_f$, it
follows that there exist an unbounded subsequence of integers $\{n_s\}_{s=1}^\infty$ and a  topological (maybe  not geometric) annulus $A$ contained in all annuli $A_{n_s,m(n_s)-1}$'s and not contractible to a point in any of them. Denote the bounded connected component of $\C\sms A$ by $D$. So 
$$
D\subset  W_{n_s,m(n_s)-1}^{(2)}.
$$
Hence,
\beq\label{5_2017_09_29}
f^{m(n_s)-1}(D) \subset B_e(x_n,\delta_n).
\eeq
Thus the family of functions $\{f^{m(n_s)-1}:D\to {\mathbb C}\}_{s=1}^{\infty}$ is normal, and, consequently $D$, cannot intersect the Julia set
$J(f)$. If $D$ were contained in a preimage of a Siegel disk or a
Herman ring, the limit of diameters of iterates $f^{m(n_s)-1}(D)$
would be positive contrary to item (2) and \eqref{5_2017_09_29}. Thus, because of  Theorem~\ref{baker-domain+Sullivan for elliptic} and Theorem~\ref{Fatou Periodic Components}, $D$ is contained in the basin of attraction
to an attracting periodic orbit or a parabolic periodic orbit. In
either case, as $s\to\infty$, the sets $f^{m(n_s)-1}(D)$ would have some limit point being either an attracting periodic point or a parabolic
periodic point. Because of \eqref{5_2017_09_29} and item (2), the sequence $\{x_n\}_{n=1}^\infty$ would have also some limit point being either an attracting periodic point or a parabolic periodic point. This would however contradict item (1). The proof is complete.  
\epf

\sp Our second preparatory result is the folowing obvious observation. 

\bobs\label{o1_2017_10_06}
Let $f:\C\lra\oc$ be an elliptic function. If  $x\in
J(f)\sms \om(\Crit_r(f))$, then there exists $\eta_1(f,x)\in (0,\eta_0(f)]$ \index{(S)}{$\eta_1(f,x)} such that
\begin{itemize}
\item [(1)] there is no critical point $c$ of $f$ for which there
exist two integers $0< n_1 \leq n_2$ satisfying
$$
|f^{n_1}(c)- c| \le \eta_1(f,x) \quad {\rm and }\quad |f^{n_2}(c)- x| \le
\eta_1(f,x),
$$

\,

\fr If $\Crit_r(f)=\es$, i.e. if the elliptic function $f$ is non-recurrent, then we set
$$
\eta_1(f):=\min\Big\{\dist_e\(J(f),{\rm Attr}(f)\),\min\big\{\dist_e(c,O_+(f(c)): \,c\in J(f)\cap\Crit(f)\big\}\Big\},
$$
and item (1) takes on the following stronger form
\item[(1')]
$$
|f^{n}(c)- c|> \eta_1(f)
$$
for all $c\in J(f)\cap\Crit(f)$ and all integers $n\ge 1$.
\end{itemize}
\eobs

\fr Now, let $N_f$ be \index{(S)}{$N_f$} the number of equivalence classes of the relation
$\sim_f$ between critical points of an elliptic function $f:{\mathbb C}\to \ov{\mathbb C}$. In other words
$$
N_f=\#(\Crit(f)\cap \mathcal R)
$$
for any fundamental domain $\mathcal R$ of $f$ whose boundary contains no critical points of $f$.

\,

Our last preparatory result, interesting on its own, is this.

\sp\blem\lab{lmane2} 
Let $f:\C\lra\oc$ be an elliptic function. Let $x\in
J(f)\sms \om(\Crit_r(f))$. If $U\subset {\mathbb C} $ is an open connected simply connected
neighborhood of $x$, $n\ge 0$ is an integer and $V$ is a connected component of $f^{-n}(U)$ for which
$$
{\rm diam}_e(f^i(V))\leq \eta_1(f,x), 
$$
for all $0\leq i  \leq n$, where $\eta_1(f,x)$ comes from Observation~\ref{o1_2017_10_06}, then 

\begin{enumerate}

\,

\item Each set $f^i(V)$, $0\leq i \leq n$, is simply connected and contains at most one critical point of $f$.

\, 

\item The equivalence class of the equivalence relation $\sim_f$ of each critical point of $f$ intersects at most one of the sets $f^i(V)$, $0\leq i \leq n$.

\,

\item If $f^i(V)\cap\Crit(f)=\es$, for some $0\leq i \leq n$, then the map $f|_{f^i(V)}:f^i(V)\to\C$ is 1--to--1.

\,

\item Consequently,
$$
\deg\(f^n|_V)\le N_f^*:=\prod_cp_c,
$$
where the product is taken over any (fixed) selector of the relation $\sim_f$ between critical points of $f$. 

\,

\item Hence, 
$$
\#\(V\cap\Crit(f^n)\) \leq N_f^*-1.
$$
\end{enumerate}
\elem

\bpf 
Item (1) is immediate from the definition of $\eta_0(f)$ and since $\eta_1(f,x)\le \eta_0(f)$.

In order to prove item (2) suppose for a contradiction that
$$\begin{aligned}
c_1\in\Crit(f)\cap f^{k_1}(V), \ \ 
   c_2\in \Crit(f)\cap f^{k_2}(V)
\end{aligned}
$$ 
and $c_1\sim_f c_2$, where $0\le k_1<k_2\le n$. But then also
$$
f^{k_2-k_1}(c_2)=f^{k_2-k_1}(c_1) \in
f^{k_2}(V).
$$
So,
\beq\label{1_2017_10_09}
|f^{k_2-k_1}(c_2)-c_2|\le\eta_1(f).
\eeq
On the other hand,
\beq\label{2_2017_10_09}
f^{n-k_1}(c_2)
=f^{n-k_1}(c_1)
\in f^{n-k_1}\(f^{k_1}(V)\)
=V.
\eeq
So, $|f^{n-k_1}(c_2)-x|\le\diam(V)\le\eta_1(f)$. Since also $k_2-k_1>0$ and $n-k_1\ge k_2-k_1$, both, \eqref{1_2017_10_09} and \eqref{2_2017_10_09}, contradict Observation~\ref{o1_2017_10_06}, and the proof of item (2) is finished.
Item (3) now follows immediately from the first part of item (1) and Corollary~\ref{General Riemann Hurwitz Simply Connected_2}. 

Item (4) is an immediate consequence of items (1)--(3) along with Corollary~\ref{General Riemann Hurwitz Simply Connected_2}, Corollary~\ref{General Riemann Hurwitz Simply Connected_3}, and multiplicity of degree.

Item (5) in turn is an immediate consequence of formula \eqref{2_2017_10_12} from Corollary~\ref{General Riemann Hurwitz Simply Connected_1}.
\endpf

\sp\bthm[Ma\~ne's Theorem for Elliptic Functions]\lab{mnt6.3a} 
Let $f:\C\lra\oc$ be an elliptic function and, as always, let $\Om(f)$\index{(S)}{$\Om(f)$} denote the set of all rationally indifferent periodic points of $f$. If 
$$
x\in J(f)\sms (\Om(f)\cup \om(\Crit_r(f))
$$
then 

\,

$\forall\, \varepsilon >0\  \exists\, \d>0$ 
 such that
\begin{itemize}
\item [(a)] For all integers $n\geq 0$, every  connected component of
  $f^{-n}(B_e(x,\d))$ has Euclidean diameter $\leq \varepsilon$.

\, \item [(b)] $\deg\(f^{n}_{|V}\)\le N_f^*$ for all integers $n\geq 0$ and for every connected component $V$ of $f^{-n}(B_e(x,\d))$.
\end{itemize}
\ethm

{\sl  Proof.} The core of this theorem is item (a). Then, item (b) follows immediately from item (a) and Lemma~\ref{lmane2}.

Given an open set $U\subset \ov{\mathbb C}$, denote, in this proof, by 
\begin{itemize}
\item[(1)]
$$
c(U,n)
$$
\end{itemize}
the collection of all connected components of $f^{-n}(U)$. Of course if $V\in c(U,n)$, then $f^j(V) \in c(U, n-j)$  for all $0\leq j \leq n$.

Recall that given any real number $\a>0$ and an open ball $B=B(z,r)$, we denote by $\a B$ the ball $B(z,\a r)$. 

\,

If $B$ is an open ball with radius $r$, then denote  by
${\mathcal L}(B)$ the  family of all open balls contained in
$\frac32B\sms B$ with radii equal to $r/4$. Denote further by ${\mathcal
L}^*(B)$ the family of all squares $\frac32D$  with $D\in
{\mathcal L}(B)$. 

\,

Take $\eta_2(f,x)\in (0,\eta_1(f,x)]$ so small that
\begin{itemize}
\item[(2)]
$$
\dist_e(x,\Om(f)\cup {\rm Attr}(f))> 10 \eta_2(f,x). 
$$
\end{itemize}
Because of Observation~\ref{o1_2017_10_04}, applied with $n=1$, for every $\varepsilon >0$, there exists $\varepsilon_1 >0 $ satisfying the following two conditions.
 
\begin{itemize}
\item [(3)] $0< \varepsilon_1 < \min\lt\{\frac{\varepsilon}{10},\frac{\eta_2(f,x)}{10}\rt\}$.

\,\item [(4)] If $U$ is an open connected  set with ${\rm  diam}_e(U) \leq
2\varepsilon_1$, then ${ \rm diam}_e(W) \leq \eta_2(f,x)$ for all $W\in
c(U,1)$.
\end{itemize}
Let $\delta>0$ be  given by

\, \begin{itemize}
\item [(5)] $\delta:=\min\lt\{\frac{\eta_2(f,x)}{10}, \delta_0\rt\}$, where
$\d_0>0$ is the number produced in Lemma~\ref{felman} for 
$N=N_f^*$, $\l=2/3$, $\varepsilon={\varepsilon_1}/{(20N_f^*)}$, and $\ka=\eta_2(f,x)$.
\end{itemize}

\,

Let $B_0$ be the square with center $x$ and  radius $\delta$. Suppose
that Theorem~\ref{mnt6.3a}(a) fails. Then there  exists an integer
$n\ge 0$ and a set $V\in c(B_0, n)$  with $\rm diam_e(V) \geq \varepsilon
\geq 10 \varepsilon_1$. Hence there exists an
integer $n_0\geq 0$ such that  there exists $V_0\in c\lt(\frac32B_0,n_0\rt)$  satisfying

\, \begin{itemize}
\item [(6)] $\diam_e(f^{-(n_0-i)}(B_0)\cap f^i(V_0))\leq \varepsilon_1 $
for all $1\leq i \leq n_0$, and

\, \,

\item [(7)] $\diam_e(f^{-n_0}(B_0)\cap V_0)>  \varepsilon_1 $.
\end{itemize}
Since, by (6), ${\rm diam}_e(B_0) = 2\delta < \varepsilon_1$, it follows that $n_0\ge 1$. Now, starting  with $B_0$ we shall construct inductively a sequence of squares $\(B_j\)_{j=0}^\infty$ and a monotone decreasing sequence of  positive integers $\(n_j\)_{j=0}^\infty$ satisfying the following

\, \begin{enumerate}
\item [(8)] $B_{j+1}\in {\mathcal L}^*(B_j)$ and

\, 

\item [(9)]   there exists $V_j \in c\lt(\frac32B_j,n_j\rt)$ such that
$$
{\rm diam}_e(f^{-(n_j-i)}(B_j)\cap f^i(V_j))\leq \varepsilon_1
$$
for all  $1\leq i \leq n_j$ and
$$ {\rm diam}_e(f^{-n_j}(B_j)\cap V_j)> \varepsilon_1.
$$
\end{enumerate}
If we construct  such a sequence of squares  and integers, then  Theorem~\ref{mnt6.3a} (a) will be proved as follows. The sequence $\(n_j\)_{j=0}^\infty$ must stabilize, i.e. $n_j=n_i$ for all $j\geq i$ and some
$i\ge 0$. But the definition of the operation ${\mathcal L}^*$ implies that the radius of $B_j$ is equal to $(3/8)^j\d$. In particular  
\beq\label{1_2017_10_06}
\lim_{j\to\infty}\diam_e (B_j)=0.
\eeq
On the other hand, by (9),
$$
V_j\in c\lt(\frac32B_j,n_j\rt)= c\lt(\frac32B_j, n_i\rt)
$$
and
$$
{\rm diam}_e(f^{-n_i}(B_j)\cap V_j)
={\rm diam}_e(f^{-n_j}(B_j)\cap V_j)
>\varepsilon_1.
$$
It therefore follows that 
$$
\varliminf_{j\to\infty}\diam_e\lt(c\lt(\frac32B_j, n_i\rt)\rt)
\ge \varepsilon_1>0.
$$
This however contradicts \eqref{1_2017_10_06} and Observation~\ref{o1_2017_10_04}.

\,

The sequences $\(B_j\)_{j=0}^\infty$ and $\(n_j\)_{j=0}^\infty$ will be constructed  by induction starting with $B_0$. It follows from (6) and  (7) that $B_0$  satisfies (9). For the inductive step fix an integer $j\ge 0$ and suppose  that $B_i$ and $n_i$ are constructed  for all $0\leq i \leq j$ satisfying (8) and (10). In order to find $B_{j+1}$ and $n_{j+1}$, we begin by observing that by (8) and from (a) and the definition of the operation ${\mathcal L}^*$, we have for all $z \in B \in {\mathcal L}^*(B_j)$, that
$$
|z-x|\leq \sum_{i=0}^{j+1} {\rm diam}_e (B_i) 
= \sum_{i=0}^{j+1}
\lt(\frac{3}{8}\rt)^i {\rm diam}_e(B_0) =
2\sum_{i=0}^{j+1}\lt(\frac{3}{8}\rt)^i \delta \leq
4\delta.
<\eta_2(f,x).
$$
This means that $B\sbt B(x,\eta_2(f,x))$. Hence, by invoking (2) of Observation~\ref{o1_2017_10_06}, we get that
\begin{enumerate}
\item [(10)] $ \dist_e\(B,\Om(f)\cup {\rm Attr}(f)\)>9\eta_2(f,x)$. 
\end{enumerate}

\sp\fr For the induction step (i.e., the construction of $B_{j+1}$ and
$n_{j+1}$), we  we shall prove the following.

\sp{\bf Claim~$1^0$:}
There exists a ball $B\in {\mathcal L}(B_j)$ for which there exist an integer $1\le n\leq n_j$ and a set $V\in c(B,n)$ such that
$$
{\rm diam}_e(V) \geq\frac{\varepsilon_1}{10N_f^*}. 
$$
{\sl Proof.} 
Seeking contradition, suppose that the claim is false. Then, invoking also (9), we see that for all $1\leq i \leq n_j$, we get
$$\aligned
{\rm diam}_e(f^i(V_j)) 
& \leq {\rm diam}_e\(f^{-(n_j-i)}(B_j)\cap f^i(V_j)\)+ \\
& \  \  \  \  \  \  \  \  \  \  \  \  +\sup\big\{{\rm diam}_e(W): B \in {\mathcal L}^*(B_j) \ {\rm and} \  W\in c(B, n_j-i)\big\} \\
 &\leq   \varepsilon_1+ \frac{\varepsilon_1}{10 N_f^*}
 \leq  2 \varepsilon_1.
\endaligned
$$
From this inequality, applied to $i=1$, and  from property (4), we get that
$$
{\rm diam}_e(V_j) \leq \eta_2(f,x).
$$
In addition, since $2\varepsilon_1\leq \eta_2(f,x)$ (by (3)), we have 
$$
{ \rm  diam }_e(f^i(V_j)) \leq \eta_2(f,x)
$$
for all $ 1\leq i \leq n_j$, hence, in consequence, for all $0\leq i \leq n_j$. By
Lemma~\ref{lmane2}, this proves that 
\beq\label{1_2017_10_17}
\deg\(f^{n_j}|_{V_j}\) \leq N_f^*.
\eeq
Hence,
\beq\label{1_2017_10_11}
\# \{W \in c(B_j, n_j): W \subset V_j \} \leq \deg\(f^{n_j}|_{V_j}\) \leq N_f^*.
\eeq
Also, since $V_j\in  c\lt(\frac23B_j,n_j\rt)$, it follows from
(10), \eqref{1_2017_10_17}, Lemma~\ref{lmane2} (5) and Lemma~\ref{felman}, that 
\beq\label{2_2017_10_11}
[W\in c(B_j,n_j) \ \ {\rm and} \  \ W \subset  V_j]\  \imp \ {\rm diam}_e(W) \leq
\frac{\varepsilon_1}{20 N_f^*},
\eeq

\beq\label{3_2017_10_11}
[S\in {\mathcal L}(B_j) \ \ {\rm and} \  \ G \in c(B, n_j)]\ \imp  \ {\rm  diam}_e(G) 
\leq\frac{\varepsilon}{20 N_f^*}.
\eeq
and

Now observe that $V_j$ is the union of all sets $W \in c(B_j,
n_j)$, where $W \subset V_j$, along with the sets $G\in c(B,n_j)$, where $G
\subset V_j$ and $B\in {\mathcal L}(B_j)$. Denote the former family of these sets by $\cF_1$ and the latter by $\cF_2$. Then, for any two sets $W'$, $W''$ in
$\cF_1\cup\cF_2$ there exist an integer $k\ge 0$ and mutually distinct sets 
$$
W'=W_0, W_1, \ldots, W_k=W''
$$
alternately belonging to either $\cF_1$ or $\cF_2$ and such that
$$
\ov W_l\cap \ov W_{l+1}\ne\es
$$
or all $0\leq l<k$. Then $k\le 2N_f^*$ by \eqref{1_2017_10_11}, and, by \eqref{2_2017_10_11} and \eqref{3_2017_10_11}, we get
$$
{\rm diam}_e(V_j)  
\leq (2N_f^*+1) \frac{\varepsilon_1}{10 N_f^*}
\le \frac{3}{10}{\varepsilon_1},
$$
This contradicts  the last inequality in condition (9), and the proof of Claim~$1^0$ is finished.  
\endpf

\sp Now take $V$ produced in Claim~$1^0$. Denote by $\tilde{V}$ the only set in   
$c\lt(\frac32 B,n\rt)$ containing $V$. If $\deg\(f^n|_{\tilde{V}}\)\)\leq N_f^*$, then, by  formula \eqref{2_2017_10_12} of Corollary~\ref{General Riemann Hurwitz Simply Connected_1}, by Lemma~\ref{felman}, and by condition (5), we get that
$$
{\rm diam}_e(V) \leq \frac{\varepsilon_1}{20N_f^*},
$$
since $V\in c\lt(\frac23\lt(\frac32 S\rt),n\rt)$ and it is contained in
$\tilde{V}$. This however contradicts Claim~$1^0$ and proves that
$$
\deg\(f^n|_{\tilde{V}}\)\)\geq N_f^*+1.
$$
It then follows from Lemma~\ref{lmane2} that
$$
\diam_e(f^l(\tilde{V})) > \eta_2(f,x)
$$
for some $0 \leq l \leq n$. Now we define 
$$
B_{j+1}:=\frac32B\in {\mathcal L}^*(B_j).
$$
Then $f^i(\tilde{V})\in c(B_{j+1}, n-l)$ and $ {\rm diam}_e(f^l(
\tilde{V})) > \delta_0 \geq 10 \varepsilon_1$. Moreover, ${\rm
diam}_e(B_{j+1})\leq  2  \delta < \varepsilon_1$. Hence, there exists
$0 \leq n_{j+1} \leq n-l\leq n_j-l$ and $V_{j+1} \in \frac32c(B_{j+1},
n_{j+1}) $ such that
$$ 
{\rm diam}_e(f^{-n_{j+1}}(B_{j+1})\cap V_{j+1}) >\varepsilon_1
$$
and
$$
{\rm diam }_e(f^{-n_{j+1}+i}(B_{j+1})\cap f^i( V_{j+1}))\leq
\varepsilon_1
$$
for all $1\le i\le n_{j+1}$.
Observe that $n_{j+1}>0$ since $\rm diam_e\(B_{j+1}\) <
\varepsilon_1$. This completes the construction of the sequences
$\{B_j\}_{j=0}^\infty$ and $\{n_j\}_{j=0}^\infty$, and simultaneously the proof of entire Theorem~\ref{mnt6.3a}. 
\endpf

\sp

As a kind of maximal consequence of this theorem, we shall prove the following.

\sp\bthm\lab{mnt6.3A} 
Let $f:\C\lra\oc$ be a non-constant elliptic function. If 
$$
X\sbt J(f)\sms \(\Om(f)\cup\om(\Crit_r(f))\)
$$
is compact (remember that now $\infty\in J(f)$), then for every $\e>0$ there exists $\d>0$ such that  
$$
\sup\big\{\diam_e\(\Comp_s(z,f^n,\d)\):n\ge 1,\, z\in f^{-n}(X)\big\}\le \e,
$$
where the subscript $s$ above indicates that this is a connected component of $f^{-n}(B_s(z,f^n,\d))$. 
\ethm

\bpf Assume first that 

\sp \fr Case $1^0$:
$$
X\sbt J(f)\sms \(\Om(f)\cup\om(\Crit_r(f))\cup\{\infty\}\).
$$
Then, by virtue of Theorem~\ref{mnt6.3a}, for every $x\in X$ there exists $\d_x>0$ such that for every integer $n\ge 0$ all the connected
components of $f^{-n}(B_s(x,\d_x))$ have Euclidean diameters $\le \e$. Let $2\d>0$ be a Lebesgue number of the cover $\{B_e(x,\d_x)\}_{x\in X}$. Then for every $x\in X$ there exists $y\in X$ such that 
$$
B_s(x,\d)\sbt B_e(y,\d_y).
$$
So, for every integer $n\ge 0$ each connected component of $f^{-n}(B_s(x,\d))$ is contained in a (unique) connected component of $f^{-n}(B_e(y,\d_y))$, whence its 
Euclidean diameter is $\le\e$. 

For the general case let 
$$
\De:=\dist_e\(\Om(f)\cup\om(\Crit_r(f)), \, f^{-1}(\infty)\)>0.
$$
In view of (\ref{1110301}) and (\ref{u4}) there exists $R>0$ so
large that if $|f(z)|\ge R/2$, then for some  $b\in f^{-1}(\infty)$,
$z\in B_b(R/2)$
\beq\lab{2110301}
|f'(z)|\ge 2 \  \  {\rm and } \  \  \diam_e(B_b(R/2))\le \De/2.
\eeq
Consider now the compact set 
$$
Y:=X\cup \Big(J(f)\sms \Big(B_e\(\Om(f)\cup\om(\Crit_r(f)), \, \De/2\)\cup B_\infty(R)\Big)\Big)
$$ 
Let $0<\d_1\le\min\{\e,R/2\}$ be the number $\d$ ascribed to the set $Y$ and the number $\min\{\e,R/2,\d_1\}$ according to Case $1^0$. Let $\d_2\in(0,\d_1]$ be the number $\d>0$ ascribed to $\min\{\e,R/2,\d_1\}$ and the integer $n=1$ according to Observation~\ref{o1_2017_10_04}.
Finally, let $0<\d\le\d_2$ be ascribed to the set $Y$ and the number $\d_2$ according to Case $1^0$.

\sp Because of Case $1^0$, in order to complete the proof of our theorem, it suffices to prove the following . 

\sp {\bf Claim~$1^0$:} If $x\in B_\infty(R)$, then for each integer $n\ge 1$ the
Euclidean diameter of each connected component $V$ of
$f^{-n}(B_s(x,\d))$ does not exceed $\e$. 

\,

\bpf 
Fix $w\in f^{-n}(x)\cap V$ and let $0\le k\le n$ be the least
integer such that $f^{n-k}(w)\notin B_\infty(R)$ provided it exists.
Otherwise, set $k=n$. Since $f^n(w)=x\in B_\infty(R)$, we have $k\ge 1$. We shall show by induction that
\beq\lab{3110301}
\diam_e\(f^{n-j}(V)\)\le \d_1
\eeq
for every $1\le j\le k$. 

For $j=1$, this formula is true since $\d\le\d_2$ and because of the second part of Observation~\ref{o1_2017_10_04}. 

For the inductive step suppose that \eqref{3110301} holds for some $1\le j\le
k-1$. Since $f^{n-j}(w)\in B_\infty(R)$ and since $\diam_e
(f^{n-j}(V)\le \d_1\le R/2$, we conclude that
\beq\lab{4110301}
f^{n-j}(V)\sbt B_\infty(R/2).
\eeq
It therefore follows from the first part of formula (\ref{2110301})
that
$$
\diam_e\(f^{n-(j+1)}(V)\) 
\le{1\over
2}\diam_e\(f^{n-j}(V)\) 
\le \frac12\d_1
\le \d_1.
$$
This proves formula (\ref{3110301}). 

In the case when $k=n$, Claim~$1^0$ follows from (\ref{3110301}) since $\d_1\le\e$. Otherwise, i. e. if $1\le k\le n-1$, note first (\ref{4110301}) also holds for $j=0$, and then that it follows from (\ref{4110301}), applied with $j=k-1$, and from the second part of formula (\ref{2110301}) that
$$
f^{n-k}(w)\in{\mathbb C}\sms B_e\(\Om(f) \, \cup \om(\Crit_r(f)), \,\De/2\).
$$ 
Since we also
know that $f^{n-k}(w)\notin B_\infty(R)$, we conclude that $f^{n-k}(w)\in
Y$. It thus follows from (\ref{3110301}) and, already proven, Case $1^0$ that 
$$
\diam_e(V)\le \min\{\e,R/2,\d_1\}\le \e.
$$
We are done.
\endpf 

\sp\fr As a fairly immediate consequence of this this theorem, we get the following.

\sp\bthm\lab{mnt6.3} 
Let $f:\C\lra\oc$ be a non--recurrent elliptic function. If 
$$
X\sbt J(f)\sms\Om(f)
$$
is compact (remember that now $\infty\in J(f)$), then for every $\e>0$ there exists $\d>0$ such that
\begin{enumerate}
\item 
$$
\sup\big\{\diam_e\(\Comp_s(z,f^n,\d)\):n\ge 1, z\in f^{-n}(X)\big\}\le \e,
$$
and 

\item
$$
\lim_{n\to\infty}\sup\big\{\diam_e\(\Comp_s(z,f^n,\d)\):z\in f^{-n}(X)\big\}=0,
$$
where the subscript $s$ above indicates that this is a connected component of the set $f^{-n}\(B_s(z,f^n,\d)\)$. 
\end{enumerate} 
\ethm

{\sl Proof.} Item (1) is indeed an immediate consequence of Theorem~\ref{mnt6.3A}. For item (2), note that the number $\eta_1(f)>0$ of Observation~\ref{o1_2017_10_06} is now well define and we may assume without loss of generality that 
$$
\d\le\e<\eta_1(f).
$$
Now, since $\Om(f)$ is $f$--invariant, we have that
$$
f^{-1}(X)\sbt J(f)\sms\Om(f),
$$
and, as $f^{-1}(X)$ is closed and $\Om(f)$ is finite, we have that
$$
\dist_e(f^{-1}(X),\Om(f))>0.
$$
Therefore, both Lemmas~\ref{lmane2} and \ref{felman}, the latter used for the set of points $x\in f^{-1}(X)$, apply to yield item (2).
\endpf

\sp Let us now take first, basic, fruits of Theorem~\ref{mnt6.3}. From now on throughout this chapter and in fact throughout the entire book, $f:{\mathbb C}\to \ov{\mathbb C}$ is assumed to be a non--recurrent elliptic function. As an immediate consequence of Theorem~\ref{mnt6.3} and Theorem~\ref{Euclid-I}, we get the following. 

\sp\blem\label{l120180306}
Let $f:\C\lra\oc$ be a non--recurrent elliptic function. If 
$$
X\sbt J(f)\sms\PC(f)
$$
is compact (remember that now $\infty\in J(f)$, so, in particular, both $X$ and $\PC(f)$ may contain infinity), then for every $\d\in(0,\dist_e(X,\PC(f))$
we have that
\begin{enumerate}
\item 
$$
\sup\big\{\diam_e\(f_z^{-n}\(B_e(f^n(z),\d)\)\):n\ge 1, z\in f^{-n}(X)\big\}<+\infty,
$$
\item
$$
\lim_{n\to\infty}\sup\big\{\diam_e\(f_z^{-n}\(B_e(f^n(z),\d)\)\):z\in f^{-n}(X)\big\}=0,
$$
\item  
$$
\inf\big\{|(f^n)'(z)|:n\ge 1, z\in f^{-n}(X)\big\}>0,
$$
\fr and

\item  
$$
\lim_{n\to\infty}\inf\big\{|(f^n)'(z)|:z\in f^{-n}(X)\big\}=+\infty.
$$
\end{enumerate} 
\elem

\brem\label{r120190427}
We could have alternatively easily deduced this lemma from Lemma~\ref{l1j293} and Lemma~\ref{l1j295}, using the latter after projecting on the torus.
\erem

In \eqref{dp4.5} the number $\th(f)$ was defined. We now want to utilize it. However we need a (possibly) smaller number, so we redefine $\th(f)$ to be
\beq\lab{d5.2}\index{(S)}{$\th$}
\th=\th(f):=\min\lt\{\min\{\th(f,\om):\om\in\Om(f)\},\, {1\over
2}\dist_e(\Om(f),\Crit(f))\rt\} >0.
\eeq
We denote
\beq\lab{d5.3}\index{(S)}{$A(f)$}
A=A(f):=\max \{ A(f,c),\, A(f,b):\, c\in \Crit(f), \, b\in
f^{-1}(\infty)\, \},
\eeq
where $A(f,c)$ was defined just after Definition~\ref{dcomp} while $A(f,
b)$ was defined in (\ref{1042506}). Recall from
Chapter~\ref{elliptic-theory} that two points $z$ and $w$  are
equivalent and  write 
$$
z\sim_f w \  \  {\rm if} \ \ w-z\in \La_f,
$$
the lattice associated to the elliptic function $f$. Obviously, $z\sim_f w$
implies that 
$$
O_+(z)=O_+(w) \  \ {\rm and} \  \ \om(z)=\om(w).
$$
Since the number $N_f$ of equivalnce classes of critical points with respect to the relation $\sim_f$ is finite, each of the following four numbers below is positive.
$$
\aligned
&\eta_1(f), \th/2, \\
&\min\{(A(f,c)R(f,c))^{1/p_c}:\, c\in \Crit(f)\} \\
&\min\{|c-c'|:\, c, c'\in \Crit(f) \and c\ne c'\},
\endaligned
$$
where $p_c=p(f,c)$ is the order of the critical point $c$ of $f$.
Both, $p_c$  and $R(f,c)$, were  defined just after
Definition~\ref{dcomp}. Fix a positive  constant
\beq\lab{1042206}\index{(S)}{$\beta$}
\beta=\b_f
\eeq
smaller than all these four numbers.
It immediately follows from Theorem~\ref{mnt6.3} that there exists
$0<\g=\g_f<1/4$ such that if $n\ge 0$ is an integer, $z\in J(f)$ and
$f^n(z)\notin B_e(\Om(f),\th)$, then
\beq\lab{dp5.4}
\diam_e\(\Comp(z,f^n,2\g)\)<\b_f.
\eeq \index{(S)}{$\Comp(z,f^n,2\g)$}
From now on fix also 
\beq\label{1_2017_11_02}
0<\tau<\th^{-1}\min\{\b,\g\}
\eeq
so small as required in Lemma~\ref{ldp4.4} for every $\om\in\Om(f)$ and so small
that for every $z\in J(f)$
\beq\lab{2012102}
\diam_e\(\Comp(z,f,\th\tau)\)<\min\{\b,\g\}.
\eeq

As an immediate consequence of Lemma~\ref{lmane2}, we get the following.

\sp\blem\lab{ld5.1}  
Let $f:\C\lra\oc$ be a non--recurrent elliptic function. If $n\ge 0$ is an integer, $\eta>0$, $z\in J(f)$
and for every $k\in\{0,1,\ld,n\}$
$$
\diam_e\(\Comp(f^k(z),f^{n-k},\eta)\)\le\b_f,
$$
then 

\begin{enumerate}

\,

\item Each connected component $\Comp(f^k(z),f^{n-k},\eta)$, $k\in\{0,1,\ld,n\}$, is simply connected and contains at most one critical point of $f$.

\, 

\item The equivalence class of the equivalence relation $\sim_f$ of each critical point of $f$ intersects at most one of the sets $\Comp(f^k(z),f^{n-k},\eta)$, $k\in\{0,1,\ld,n\}$.

\,

\item If $\Comp(f^k(z),f^{n-k},\eta)\cap\Crit(f)=\es$, for some $k\in\{0,1,\ld,n\}$, then the map $f$ restricted to $\Comp(f^k(z),f^{n-k},\eta)$, is 1--to--1. 

\,

\item Consequently,
$$
\deg\lt(f^n|_{\Comp(z,f^n,\eta)}\rt)\le N_f^*-1.
$$

\,

\item Hence, 
$$
\#\(\Comp(z,f^n,\eta)\cap\Crit(f^n)\) \leq N_f^*-1.
$$
\end{enumerate}
\elem

The argument provided in the proof below has actually already appeared in the proof of Lemma~\ref{felman}. Now, we single it out and make it more transparent.

\sp

\blem\lab{lmod} 
Let $f:\C\lra\oc$ be a non-recurrent elliptic function. If $z\in J(f)$, $n\ge 0$ is an integer, and $f^n(z)\notin B_e(\Om(f),\th)$, then the set $\Comp(z,f^n,2\g)\sms \Comp(z,f^n,\g)$ is a topological annulus and
$$
\Mod\(\Comp(z,f^n,2\g)\sms \Comp(z,f^n,\g)\) 
\ge \frac{\log 2}{(N_f^*)^2}. 
$$
\elem

\bpf The first assertion follows immediately from Lemma~\ref{ld5.1} (1) and the definition of $\g$. By Lemma~\ref{ld5.1} (5) there exists a geometric
annulus
$$
R\sbt B_e(f^n(z),2\g)\sms B_e(f^n(z),\g)
$$
centered at $f^n(z)$ of modulus $\log 2/N_f^*$ such that
$$
R_n\cap \Crit(f^n)=\es, 
$$
where $R_n:=f^{-n}(R)\cap\Comp(z,f^n,2\g)$ is, also due to Lemma~\ref{ld5.1} (1), a topological annulus. Hence, applying Corollary~\ref{c1j225} (monotonicity of modulus) and Theorem~\ref{t1j235} along with Lemma~\ref{ld5.1} (4), we conclude that
\beq\lab{1050206}
\Mod \(\Comp(z,f^n,2\g)\sms
\Comp(z,f^n,\g)\)
\ge \text{Mod}(R_n)
\ge \frac{\log 2}{N_f^*}(N_f^*)^{-1}  
=\frac{\log 2}{(N_f^*)^2}.
\eeq
\endpf

\sp\fr As an immediate consequence of this lemma and Theorem~\ref{Euclid-II} we get the following.

\blem\lab{lncp23.2} 
Let $f:\C\lra\oc$ be a non--recurrent elliptic function. Suppose that $z\in J(f)$ and  $f^n(z)\notin
B_e(\Om(f),\th)$. If $0\le k\le n$ and
$$
f^k:\Comp(z,f^n,2\g)\lra \Comp(f^k(z),f^{n-k},2\g)
$$ 
is univalent, then
$$
{|(f^k)'(y)|\over |(f^k)'(x)|}\le w\lt(\frac{\log 2}{(N_f^*)^2}\rt)
$$
for all $x,y\in \Comp(z,f^n,\g)$, where $w:(0,+\infty)\to [1,+\infty)$ is the function produced in Theorem~\ref{Euclid-II}.
\elem

\blem\lab{lncp23.3} Let $f:\C\lra\oc$ be a non-recurrent elliptic function. Suppose that $z\in J(f)$ and  $f^n(z)\notin
B_e(\Om(f),\th)$. Suppose also that $Q^{(1)}\sbt Q^{(2)}\sbt
B(f^n(z),\g)$ are connected sets. If $Q_n^{(2)}$ is a connected
component of $f^{-n}(Q^{(2)})$ contained in $\Comp(z,f^n,\g)$
and $Q_n^{(1)}$ is a connected component of $f^{-n}(Q^{(1)})$
contained in $Q_n^{(2)}$, then
$$
{\diam_e\(Q_n^{(1)}\)\over \diam_e\(Q_n^{(2)}\)} \gek
{\diam_e\(Q^{(1)}\)\over \diam_e\(Q^{(2)}\)}.
$$
\elem

\bpf Let $1\leq n_1\leq \ldots \leq n_u\leq n$  be all
the integers $k$ between  1 and $n$  such that
$$ \Crit(f) \cap \Comp(f^{n-k}(z),f^k, 2\g)\neq \es.$$
Fix  $1\leq i \leq u $. If $ j \in [ n_i, n_{i+1}-1]$ (we set
$n_{u+1}=n-1$), then  by Lemma~\ref{lncp23.2} there exists a
universal constant $T>0$ such that
\beq\label{1020907}
\frac{\diam_e(Q_{j}^{(1)})}{\diam_e(Q_{j}^{(2)})} \geq T
\frac{\diam_e(Q_{n_i}^{(1)})}{\diam_e(Q_{n_i}^{(2)})}.
\eeq
Since, in view of  Lemma~\ref{ld5.1}, $ u \leq \sharp( \Crit(f)\cap
{\mathcal R})$, in order   to conclude  the proof  it is enough  to
show the existence  of a universal   constant $E >0$ such that  for
every $ 1  \leq i \leq u$
$$\frac{\diam_e(Q^{(1)}_{n_i})}{\diam_e(Q^{(2)}_{n_i})}\geq E
\frac{\diam_e(Q^{(1)}_{n_i-1})}{\diam_e(Q^{(2)}_{n_i-1})}.
$$
Indeed, let $c$  be the  critical point in $\Comp(f^{n-n_i}(z),
f^{n_i}, 2\g)$ and let $p_c\ge 2$ be, as always, its  order. Since  both
sets $Q^{(1)}_{n_i}$ and  $Q^{(1)}_{n_i}$  are connected, we get for
$j=1,2$ that
$$\begin{aligned}
\diam_e(Q_{n_i-1}^{(j)})
&\comp \diam_e(Q_{n_i}^{(j)})\sup\{|f'(x)|:x\in Q_{n_i}^{(j)}\} \\
&\comp \diam_e(Q_{n_i}^{(i)})\Dist_e(c, Q_{n_i}^{(i)}).
\end{aligned}
$$
Hence
$$\begin{aligned}
\frac{\diam_e(Q_{n_i}^{(1)})}{\diam_e(Q_{n_i}^{(2)})} 
\comp\frac{\diam_e(Q_{n_i-1}^{(1)})}{\Dist_e(c, Q_{n_i}^{(1)})} \cdot
\frac{\Dist_e(c, Q_{n_i}^{(2)})}{\diam_e(Q_{n_i-1}^{(2)})} 
\ge \frac{\diam_e(Q_{n_i-1}^{(1)})}{ \diam_e(Q_{n_i-1}^{(2)})}.
\end{aligned}
$$
The proof is finished. \endpf

\sp An important consequence of Theorem~\ref{mnt6.3} is the following result which along with Theorem~\ref{baker-domain+Sullivan for elliptic}  sheds a lot of light on the structure of Fatou and Julia sets of non--recurrent  elliptic functions.

\sp\bthm\label{r071708} 
If $f:{\mathbb C}\lra\oc$ is a non--recurrent  elliptic function, then $f$ has no Siegel disks or Herman rings nor Cremer periodic points.
\ethm

{\sl Proof.} 
Seeking contradiction suppose that the function $f$ has a Siegel disk or a Herman ring. Denote it by $D$.
Let 
$$
\Delta:=B(0,1)
$$  
if $D$ is a Siegel disk and 
$$
\Delta:=A(0;1, r)
$$ 
with $ r >1$ coming from Theorem~\ref{Fatou Periodic Components} if $D$ is a Herman ring. Let 
$$
H:\Delta \lra D
$$ 
be the analytic homeomorphism  resulting from Theorem~\ref{Fatou Periodic Components} (3) and (4) respectively in the Siegel or Herman case. Because of (\ref{1sh11}) there exists a point 
\beq\label{1sh6B}
\xi \in \partial{D}\sms  \ov{O^{+}(\Sing (f^{-1}))}.
\eeq
Using the fact that either $D$ is simply connected or doubly connected and
 $r \in (1, +\infty)$, we deduce that $ \xi$ is not an isolated point of $\partial{D}$, in fact $\partial{D}$ has no isolated points. Therefore  we  may assume in addition that
$$ 
\xi \notin \Omega(f).
$$
Fix $ \delta >0$ so small as required in Theorem~\ref{mnt6.3} (2) and fix a point $ w \in D \cap B(\xi, \delta)$. Let
\beq\label{4sh12}
F:=H\(\{ z \in \Delta: \, |z|=|H^{-1}(w)|\}\).
\eeq
Note that $F$ is a compact set (homeomorphic to a circle) and 
$$
F \sbt D.
$$
Since $w \in D$, for every $n \geq 1$ the intersection $D\cap f^{-n}(w)$ is a  singleton. Denoting its only element by $w_n$, we will have that
\beq\label{1sh12}
w_n \in F.
\eeq
By (\ref{3sh6}), for every $n\geq 1$ there exists a unique holomorphic branch
\beq\label{2sh13}
f^{-n}_n: B(w, 2s)\lra \mathbb C
\eeq
of $f^{-n}$ sending  $w$ to $w_n$. By (\ref{1sh6B}), $ \xi \in J(f)$, whence $f^{-n}_n(\xi) \in J(f)$.  Thus, $ f^{-n}_n(\xi) \notin D$, and therefore,
 using also (\ref{2sh13}) and
(\ref{1sh12}), we conclude that
$$\diam_e\(f^{-n}_n(B(w,s))\) 
\geq \dist_e(w_n,f^{-n}_n(\xi))
\geq \dist(F, \mathbb C \sms D)>0.
$$
Therefore
$$
\liminf_{n \to \infty}\diam_e\(f^{-n}_n(B(w, r))\)
\geq \dist_e(F, \C \sms D)>0,
$$
contrary to Theorem~\ref{mnt6.3} (2). Siegel disks and Herman are thus ruled out.

\sp Now, proceeding again by the way of contradiction, suppose that $f$  has a periodic Cremer point. Call it $\xi$ and denote by $q \geq 1$ one of its periods. So, we have
\beq\label{1sh13}
\xi \in J(f), \  \ f^q(\xi)=\xi, \  \ |(f^q)'(\xi)|=1,
\eeq
and $(f^q)'(\xi)$ is not a root of unity. In particular 
$$ 
\xi \notin \Omega(f).
$$
Take $ \varepsilon >0$ so small that the map
$$
f^q|_{B(\xi,4\varepsilon)}: B(\xi, 4 \varepsilon) \lra \mathbb C
$$  
is 1-to-1. Let then $ \delta>0$ be ascribed to this $\varepsilon >0$ and the set $ X =\{\xi\}$ according to Theorem~\ref{mnt6.3}. Because of this theorem we have that 
$$
\diam_e\(\Comp_e(\xi, f^{qn}, \delta)\) \leq \varepsilon
$$ 
for every integer $n \geq 1$. In particular 
$$
\Comp_e\(\xi, f^{qn}, \delta\)\sbt B(\xi, \varepsilon).
$$
It then follows by immediate induction that the restriction
$$ 
f^{qn}|_{\Comp_e\(\xi, f^{qn}, \delta\)}: \Comp_e\(\xi,f^{qn},\delta\)\lra B(\xi, \delta)
$$
is a 1-to-1 map. So, its inverse
$$
f^{-qn}_\xi: B(\xi, \delta) \lra \Comp_e\(\xi, f^{qn}, \delta\)\sbt \mathbb  C
$$ 
is a well--defined holomorphic injective map sending $\xi$ to $\xi$. So, it follows from $\frac{1}{4}$--Koebe's Distortion Theorem, i.e Theorem~\ref{one-quater}, that
$$
f^{-qn}_{\xi} (B(\xi, \delta))
 \supset B\lt(\xi, \frac{1}{4}\delta|(f^{-qn}_\xi)'(\xi)|\rt)
 =B(\xi,\delta/4).
$$ 
Hence, 
$$
f^{qn}(B\(\xi,\delta/4)\)\sbt B(\xi, \delta)
$$ 
for every integer $n\geq 1$. Thus, the family of maps 
$$
\big\{f^{qn}|_{B(\xi,\delta/4)}\big\}_{n=1}^\infty
$$ 
is normal. Hence $\xi \in F(f^q)=F(f)$ contrary to the first item of the formula (\ref{1sh13}). The proof of Theorem~\ref{r071708} is finished. 
\endpf

\sp As an immediate consequence of this theorem, Theorem~\ref{t220200404}, Theorem~\ref{baker-domain+Sullivan for elliptic}, and Theorem~\ref {Fatou Periodic Components}, we get the following.

\bthm\label{t220200404B}
Any non--recurrent elliptic function has only finitely many Fatou periodic components and all of them are either basins of attraction to
attracting or rationally indifferent periodic points.
\ethm

\sp

\section[Compactly Non--Recurrent Elliptic Functions]{Compactly Non--Recurrent Elliptic Functions: \\ Definition, Partial order in
$\Crit_c(J(f))$, and \\ Stratification of Closed Forward-Invariant
Subsets of $J(f)$}

 From now on throughout the book we will need a stronger property than mere non--recurrence. Indeed, from now on all chapters deal almost exclusively with compactly non--recurrent elliptic functions and some of their subclasses such as regular, subexpanding, and parabolic; see Section~\ref{DDCoEF}, Dynamically Distinguished Classes of Elliptic Functions, for their definitions and some basic properties. In this section we lay down the foundations of such functions which will be used throughout the reminder of the book. 

Compact non--recurrence is given by the following definition.

\sp\bdfn\label{pseudo-non-recurrent} We say that a non--constant elliptic
function $f:{\mathbb C}\lra\oc$ is {\em compactly non--recurrent} (CNR) \index{(N)}{compactly non--recurrent} \index{(S)}{(CNR)}if and only if whenever $c \in \Crit(f)\cap J(f)$, then either
\begin{enumerate}
\item[(1)]\,
$\om$-limit set $\om(c)$ is a compact subset of $\mathbb C$ (i.e.
$\infty\notin \om(c)$) and $c\notin \om(c)$ or

\,\item[(2)]\,
 $c\in \bigcup_{ n\geq 1}f^{-n}(\infty)$, or
 
\, \item[(3)]\,
$ \lim_{n \to\ \infty}f^n(c)=\infty$. 
\end{enumerate}
\edfn

\fr Of course every compactly non--recurrent elliptic function is  non--recurrent.



\bdfn\label{d1071805} 
Given an arbitrary non--constant elliptic
function $f:{\mathbb C}\lra\oc$, the set of critical points captured by items (1), (2) or (3) in Definition~\ref{pseudo-non-recurrent} will  be respectively refereed  to as
$\Crit_c(f)$\index{(S)}{$\Crit_c(f)$},
$\Crit_p(f)$\index{(S)}{$\Crit_p(f)$}  and
$\Crit_\infty(f)$\index{(S)}{$\Crit_\infty(f)$}. 
\edfn

\fr Let us record the following immediate observation.

\sp\bobs\lab{o1_2017_11_18}
If $f:{\mathbb C}\lra\oc$ is an arbitrary elliptic function, then
$$
\Crit_p(f)\cup \Crit_\infty(f)\sbt J(f).
$$
\eobs

{\sl Proof.}
$\Crit_p(f)\sbt J(f)$ because all poles of $f$ are contained in $J(f)$ and $f^{-1}(J(f))\sbt J(f)$, while the inclusion $\Crit_\infty(f)\sbt J(f)$ immediately follows from Theorem~\ref{baker-domain+Sullivan for elliptic} and Theorem~\ref{Fatou Periodic Components}. 
\endpf

\sp\fr Keeping $f:{\mathbb C}\lra\ov{\mathbb C}$, an arbitrary elliptic function, for every $c\in \Crit_p(f)$ we in fact have that
$$
c\in \bigcup_{n\geq 2}f^{-n}(\infty),
$$
and let $n(c)\geq 2$  be the only integer such $f^j(c)$ is well  defined  for all
 $0 \leq j \leq n(c)$  and $f^{n}(c)=\infty$. Set
\beq\lab{1042706}
\begin{aligned}
  {\rm PC}_c(f)&:=\bigcup_{c\in  \Crit_c(f)}\{ f^j(c): j \geq 1\}\index{(S)}{${\rm PC}_c(f)$}\\
  {\rm PC}^0_c(f)&:=\Crit_c(f)\cup{\rm PC}_c(f)\index{(S)}{${\rm PC}^0_c(f)$},\\
  {\rm PC}_p(f)&:=\bigcup_{c\in \Crit_p(f)}\{ f^j(c):1\leq j \leq n(c)-1  \}\index{(S)}{${\rm PC}_p(f)$},\\
{\rm PC}^0_p(f)&:=\Crit_p(f)\cup{\rm PC}_p(f)\index{(S)}{${\rm PC}^0_p(f)$},\\
{\rm PC}_\infty(f)&:=\bigcup_{c\in \Crit_\infty(f)}\{ f^j(c):   j
\geq 1 \}\index{(S)}{${\rm PC}_\infty(f)$},\\
{\rm PC}^0_\infty(f)&:=\Crit_\infty(f)\cup{\rm PC}_\infty(f)\index{(S)}{${\rm PC}^0_\infty(f)$},\\
 {\rm PC}(f)&:=\bigcup_{c\in  \Crit(f)}\{ f^j(c): j \geq 1\}\index{(S)}{${\rm PC}(f)$}\\
  {\rm PC}^0(f)&:=\Crit(f)\cup{\rm PC}(f)\index{(S)}{${\rm PC}^0(f)$},\\
\end{aligned}
\eeq

Throughout this section we assume that an elliptic function $f:\C\lra\oc$ is compactly non--recurrent. Unless otherwise explicitly stated all the  distances and all closures are understood with respect to Euclidean metric and topology on the Euclidean complex plane $\mathbb C$. In particular if
$\lim_{ n \to \infty} f^n(z)=\infty$ or $z\in
\bigcup_{n=1}^{\infty}f^{-n}(\infty)$, then $\om(z)=\es$. Also
$\dist(A, \es)=0$. 

\sp We record the following immediate.

\bobs\label{120200120}
If $f:\C\lra\oc$ is a compactly non--recurrent elliptic function, then
$$
\bal
{\rm PC}(f)&:={\rm PC}_c(f)\cup {\rm PC}_p(f)\cup  {\rm
PC}_\infty(f)\index{(S)}{${\rm PC}(f)$},\\
{\rm PC}^0(f)&:={\rm PC}^0_c(f)\cup {\rm PC}^0_p(f)\cup  {\rm
PC}^0_\infty(f)\index{(S)}{${\rm PC}^0_\infty(f)$}.
\eal
$$
\eobs

In this section we introduce an order in the set of critical points and a stratification of the Julia set. Both of these are crucial for inductive proofs in forthcoming sections and chapters. The results and proofs provided in the present section pretty closely follow those from Section~2.4 of \cite{KU3}. 

Set
$$
\Crit_c(J(f)):=\Crit_c(f)\cap J(f).\index{(S)}{$\Crit_c(J(f))$}
$$

\blem\lab{ld5.2} 
If $f:\C\lra\oc$ is a compactly non--recurrent elliptic function, then the set $\om(\Crit_c(J(f)))$ is nowhere dense in $J(f)$. 
\elem

\bpf Suppose that the interior (relative to $J(f)$) of
$\om(\Crit_c(J(f)))$ is nonempty. Then there exists $c\in
\Crit_c(J(f))$ such that $\om(c)$ has nonempty interior. But then, because of Proposition~\ref{p120190913}, there would exist $n\ge 0$ such that $f^n(\om(c))=J(f)$ and, consequently, $\om(c)=J(f)$. This, however, is a contradiction, as $c\notin \om(c)$. 
\epf

\sp Since $\ov{O_+\(J(f)\cap\Crit(f)\)}$ is the union of $\om(\Crit_c(J(f)))$ and the forward orbit of $J(f)\cap\Crit(f)$, which is a countable set, as an immediate consequence of this lemma and Baire Category Theorem, we get the following. 

\bprop\label{p120190914}
If $f:\C\lra\oc$ is a compactly non--recurrent elliptic function, then the set 
$$
\ov{O_+\(J(f)\cap\Crit(f)\)}
$$ 
is nowhere dense in $J(f)$.
\eprop

\sp Now we introduce in $\Crit_c(J(f))$ a relation $<$ which, in view
of Lemma~\ref{ld5.3} below, is an ordering relation, by putting
\beq\lab{d5.5}
c_1<c_2 \  \iff  \  c_1\in\om(c_2).\index{(S)}{$c_1<c_2$}
\eeq
Since $c_2\sim_f c_3$ implies $\om(c_2)=\om(c_3)$, we have that if $c_1<c_2$
and $c_2\sim_f c_3$, then $c_1<c_3$

\blem\lab{ld5.3} Let $f:\C\lra\oc$ be a compactly non--recurrent elliptic function. If $c_1<c_2$ and $c_2<c_3$, then $c_1<c_3$. \elem

\bpf Indeed, we have $c_1\in\om(c_2)\sbt\om(c_3)$.
\endpf

\sp

\blem\lab{ld5.3a} If $f:\C\to\oc$ is a compactly non--recurrent elliptic function, then there is no infinite, linear subset of the partially ordered set $(\Crit_c(J(f)),<)$.
\elem

\bpf Indeed, suppose on the contrary that $c_1<c_2<\ld$
is an infinite, linearly ordered subset of $\Crit_c(J(f))$. Since
the number of equivalency classes of relation $\sim_f$ is equal to
$\#(\Crit_c(J(f))\cap\mathcal R)$ which is finite, there exist two
numbers $1\le i<j$ such that $\om(c_i)=\om(c_j)$. But this implies
that $c_i\in \om(c_j)=\om(c_i)$ and we get a contradiction. The
proof is finished. \endpf

\sp The following observation is a reformulation of the condition that
$J(f)$ contains no recurrent critical points.

\blem\lab{ld5.4} Let $f:\C\lra\oc$ be a compactly non--recurrent elliptic function. If $c\in \Crit_c(J(f))$, then it is not the case
that $c<c$. 
\elem

Now define inductively a sequence
$$
\{Cr_i(f)\}_{i=0}^\infty
$$\index{(S)}{$Cr_i(f)$} of subsets of $\Crit_c(J(f))$ by
setting $Cr_0(f)=\es$ and
\beq\lab{d5.6}
Cr_{i+1}(f):=\lt\{c\in \Crit_c(J(f))\sms\bu_{j=0}^iCr_j(f): \tif c'<c, \text{ then }
c'\in Cr_0(f)\cup\ld\cup Cr_i(f)\rt\}.
\eeq

\blem\lab{ld5.5}
Let $f:\C\to\oc$ be a compactly non--recurrent elliptic function. Then:
\begin{enumerate}
\item[(a)] If $c\in Cr_i(f)$ and $c'\sim_f c$, then $c'\in Cr_i(f)$.

\, \item[(b)] The sets $\{Cr_i(f)\}$ are mutually disjoint.

\, \item[(c)] $\exists_{p\ge 1} \  \forall_{i\ge p+1} \  \  Cr_i(f)=\es$.

\, \item[(d)] $Cr_0(f)\cup\ld\cup Cr_p(f)=\Crit_c(J(f))$.

\, \item[(e)] If $\Crit_c(J(f))\ne\es$, then $Cr_1(f)\ne\es$.
\end{enumerate}
\elem

\bpf Part (a) follows immediately from the definition of
the sets $Cr_i$ and the fact that two equivalent points have the
same $\om$-limit sets. By definition $Cr_{i+1}(f)\cap
\bu_{j=1}^iCr_j(f)=\es$, so disjointness in (b) is clear. As the
number $N_f$ of equivalence classes of the relation $\sim_f$ is is finite, (a) and (b) imply (c). 

Take $p$ to be the minimal number satisfying (c) and suppose
that 
$$
c\in \Crit_c(J(f))\sms\bu_{j=1}^pCr_j(f).
$$
Since $Cr_{p+1}(f)=\es$, there exists $c'\notin\bu_{j=1}^pCr_j(f)$ such
that $c'<c$. Iterating this procedure, we would obtain an infinite
sequence of critical points $c_1=c>c_2=c'>c_3>\ld$. But this
contradicts Lemma~\ref{ld5.3a} proving (d). 

Now, if $Cr_1=\es$, then it would follow from (\ref{d5.6}) by a straightforward induction that $Cr_i=\es$ for every $i\ge 0$. Then $\Crit_c(J(f))=\es$ by (d). But this would contradict our hypothesis that $\Crit_c(J(f))\ne\es$, whence proving (e). \endpf

\sp As an immediate consequence of the definition of the sequence
$\{Cr_i(f)\}_{i=0}^p$, we get the following simple lemma.

\blem\lab{ld5.6} 
Let $f:\C\lra\oc$ be a compactly non--recurrent elliptic function. If $c,c'\in Cr_i(f)$, then it is not the case that $c<c'$. 
\elem

For each point $z\in J(f)$, define the set
$$
\Crit_c(z):=\{c\in \Crit_c(J(f)):c\in \om(z)\}.
$$\index{(S)}{$\Crit_c(z)$}

\blem\lab{ld5.7} 
Let $f:\C\lra\oc$ be a compactly non--recurrent elliptic function. If $z\in J(f)\sms I_\infty(f)$, then either

\,

\begin{enumerate}
\item $z\in\bu_{n\ge 0}f^{-n}(\Om(f))$ or

\,

\item $\om(z)\sms \{\infty\}$ is not contained in $\ov{O_+(f(\Crit_c(z))}\cup \Om(f)$. 
\end{enumerate}
\elem

\bpf Suppose that $z\notin\bu_{n\ge 0}f^{-n}(\Om(f))\cup
I_\infty(f)$. Then $\om(z)\sms \{\infty\}\ne\es$ and, by Proposition~\ref{p1_2017_10_30}, the set $\om(z)\sms \{\infty\}$
is not contained in $\Om(f)$. So, if we assume that
\beq\lab{d5.7}
\om(z)\sms \{\infty \} \sbt\ov{O_+(f(\Crit_c(z))} \cup \Om(f),
\eeq
then, as $\om(z)\sms \{\infty\}\ne\es$, we conclude that
$\Crit(z)\ne\es$. Let $c_1\in \Crit_c(z)$.  This means that
$c_1\in\om(z)$; and as $c_1\notin\Om(f)$, it follows from
(\ref{d5.7}) that there exists $c_2\in \Crit_c(z)$ such that either
$c_1\in \om(c_2)$ or $c_1=f^{n_1}(c_2)$ for some $n_1\ge 1$.
Iterating this procedure, we obtain an infinite sequence
$\{c_j\}_{j=1}^\infty$ such that for every $j\ge 1$, either $c_j\in
\om(c_{j+1})$ or $c_j=f^{n_j}(c_{j+1})$ for some $n_j\ge 1$.
Consider an arbitrary block $c_k,c_{k+1},\ld,c_l$ such that
$c_j=f^{n_j}(c_{j+1})$ for every $k\le j\le l-1$ and suppose that
$$
l-(k-1)\ge N_f.
$$ 
Then there are two indexes $k\le a<b\le l$ such that $c_a\sim_f c_b$. Then
$$
f^{n_a+n_{a+1}+\ld+n_{b-1}}(c_a)=f^{n_a+n_{a+1}+\ld+n_{b-1}}(c_b)=c_a;
$$
and consequently, as $$n_a+n_{a+1}+\ld+n_{b-1}\ge b-a\ge 1,$$
 $c_a$
is a super-attracting periodic point of $f$. Since $c_a\in J(f)$,
this is a contradiction; and in consequence, the length of the block
$c_k,c_{k+1},\ld,c_l$ is bounded above by $N_f$. Hence there exists an infinite subsequence $\{n_k\}_{k\ge
1}$ such that $c_{n_k}\in \om(c_{n_k+1})$ for every $k\ge 1$ or,
equivalently, $c_{n_k}<c_{n_{k+1}}$ for every $k\ge 1$. This,
however, contradicts Lemma~\ref{ld5.3a} and we are done.
\endpf

\sp Recall that the integer $p$ was defined in Lemma~\ref{ld5.5}. Define
for every $i=0,1,\ld,p$,
\beq\lab{1071303}
S_i(f)=Cr_0(f)\cup\ld\cup Cr_i(f)
\eeq\index{(S)}{$S_i(f)$}
and for every $i=0,1,\ld,p-1$, consider $c'\in\bu_{c\in
Cr_{i+1}(f)}\om(c)\cap \Crit_c(J(f))$. Then there exists $c\in
Cr_{i+1}(f)$ such that $c'\in\om(c)$, which  means that $c'<c$.
Thus, by (\ref{d5.6}), we get $c'\in S_i(f)$. So
\beq\lab{d5.8}
\bu_{c\in Cr_{i+1}(f)}\om(c)\cap (\Crit_c(J(f))\sms S_i(f))=\es.
\eeq
Therefore, since the set $\bu_{c\in Cr_{i+1}(f)}\om(c)\sbt {\mathbb
C}$ is compact and $\Crit_c(J(f))\sms S_i(f)$ has no accumulation
points in $\mathbb C$,
\beq\lab{d5.9}
\d_i\index{(S)}{$\d_i$}:=\dist_e\lt(\bu_{c\in Cr_{i+1}(f)}\om(c),
\Crit_c(J(f))\sms S_i(f)\rt)>0.
\eeq
Set
\beq\label{1_2017_11_07}
\rho:=\frac{1}{2}\min\Big\{\min\{\d_i:i=0,1,\ld,p-1\}, \, \dist_e
\left(O_+(\Crit_c(f)), \Crit_p(f)\cup\Crit_\infty(f)\right) \Big\}.
\eeq
Fix a closed forward--invariant subset $E$ of $J(f)$, and for every
$i=0,1,\ld,p$, define
$$
E_i(f):=\big\{z\in E:\dist_e\(O_+(z), \Crit_c(J(f))\sms
S_i(f)\)\ge\rho\big\}\index{(S)}{$E_i(f)$}.
$$
Let us  now  prove the following two lemmas concerning the sets $E_i(f)$.\index{(S)}{$E_i(f)$}, $ i=0, \ldots, p$.

\blem\lab{ld5.8} 
If $f:\C\lra\oc$ is a compactly non--recurrent elliptic function, then
$$
E_0(f)\sbt E_1(f)\sbt\ld\sbt E_p(f)=E(f). 
$$
\elem

\bpf Since $S_{i+1}(f)\spt S_i(f)$, the inclusions
$E_i(f)\sbt E_{i+1}(f)$ is obvious. Since $S_p(f)=\Crit_c(J(f))$, it holds
$E_p=E$. We are done. \endpf

\

\blem\lab{ld5.9} 
If $f:\C\to\oc$ is a compactly non--recurrent elliptic function, then there exists $l=l(f)\geq 1 $ such that for every $i=0,1,\ld,p-1$, we have that
$$
\bu_{c\in Cr_{i+1}(f)}\om(c)\sbt\ov{O_+(f^l(Cr_{i+1}(f)))}\sbt
\ov{{\rm PC}_c(f)_i}\index{(S)}{$\ov{{\rm PC}_c(f)_i}$} \sbt
\ov{{\rm PC}^0_c(f)_i}.\index{(S)}{$\ov{{\rm PC}^0_c(f)_i}$}$$ \elem

\bpf The left-hand inclusion is obvious regardless
whatever $l(f)\geq 1$ is. In order to prove the right-hand
inclusion, fix $i\in\{0,1,\ld,p-1\}$. By the definition of
$\om$-limit sets, there exists $l_i\ge 1$ such that for every $c\in
Cr_{i+1}(f)$ we have
$$\dist_e\left(O_+(f^{l_i}(c)), \bu_{c\in
Cr_{i+1}(f)}\om(c)\right)<\d_i/2.$$ Thus, by (\ref{d5.9}),
$$\dist_e\(\ov{O_+(f^{l_i}(c))}, \Crit_c(J(f))\sms S_i(f)\)>\d_i/2.$$
Since $\rho\le\d_i/2$ and since for every $z\in
\ov{O_+(f^{l_i}(c))}$ also $O_+(z)\sbt \ov{O_+(f^{l_i}(c))}$, we
therefore get $$\ov{O_+(f^l(Cr_{i+1}(f)))}\sbt \ov{{\rm
PC}^0_c(f)_i}.$$ So, by putting $l(f)=\max\{l_i:i=0,1,\ld,p- 1\}$ the
proof is completed.  
\endpf

\sp

\section{Holomorphic Inverse Branches}

\sp \fr This section has a technical character The main (technical) result, Proposition~\ref{p1071305}, of this section concerns compactly non--recurrent elliptic functions and provides us with abundance of holomorphic inverse branches of iterates of elliptic functions. It will be used many times in the sequel. 

However, at the beginning of this section our considerations are fairly general. So, let $f:\C\lra\oc$ be an elliptic function. Set
$$
\Sing^-(f)
:=\bu_{n\ge 0}f^{-n}\(\Om(f)\cup \Crit(J(f))\cup f^{-1}(\infty)\)\index{(S)}{$\Sing^-(f)$},
$$ 
and recall that 
$$ 
I_-(f)=\bu_{n\ge 1}f^{-n}(\infty).
$$
Consider any $T\ge R_0$, where $R_0>0$ comes from \eqref{1_2017_10_02}. 
For every $b \in  f^{-1}(\infty)$ and every $w\in B_{b}(2T)$ let 
$$
f^{-1}_{b,w}:B_e(f(w),T)\longrightarrow B_b(T) 
$$ 
be the inverse branch of $f$ sending $f(w)$ to $w$. It follows from \eqref{u1},
 (\ref{1110301}), and (\ref{u4}) that there exists a constant $L \geq 1$ so large that for every $b\in f^{-1}(\infty)$,
$$
B(b,L^{-1})\sbt B_b(R_0), 
$$
and the following properties are satisfied.
\beq\label{1052405}
\diam_e (B_b(T))\leq LT^{-1/q_b}  \quad \mbox{and}\quad B_{b}(2T)
\supset B_e(b,L^{-1}T^{-1/q_b})
\eeq
and for every $R\in (0,T]$ sufficiently small
\beq\label{2052405}
\begin{aligned}
B_e\left(w ,L^{-1}R|f(w)|^{-\frac{q_b+1}{q_b}}\right) & \subset
f^{-1}_{b,w}(B_e
(f(w),R))\subset B_e\left(w,LR|f(w)|^{-\frac{q_b+1}{q_b}}\right) \\
& \subset B_e(w,R),
\end{aligned}
\eeq
where the last   inclusion  was written   assuming  that $|f(w)|\geq
L^{\frac{q_b+1}{q_b}}$. Since there are  only finitely many
equivalence
 classes   of  the relation $\sim_f$  generated by the  poles of $f$, there exists $R_1>0 $ so  small  that $w\in
B_e({f}^{-1}(\infty), R_1)$  then $|f(w)| \geq L$.

\sp Using now  (\ref{2052405})  and the right-hand side of (\ref{1052405}) with  $T$ replaced  by $2T$, a straightforward induction  gives the following.

\blem\label{l1052405} 
If $f:\C\lra\oc$ is a compactly non--recurrent elliptic function, then there exists   
$$
R_2\in \lt(0, \min\lt\{T,
(2LT)^{-1}, R_1, \frac{1}{2}\dist_e(f^{-1}(\infty), \Crit(f))\rt\}\rt)
$$
so small  that if $z\in {\mathbb C}$, $n \geq 1$ and if 
$$
\{f^j(z):0\leq j \leq n-1\} \sbt  B_e({f}^{-1}(\infty), R_2),
$$
then there exists a unique holomorphic inverse branch 
$$
f^{-n}_z: B_e(f^n(z),R_2)\lra B_e(z, R_2)
$$ 
sending  $f^n(z)$ to $z$. 
\elem

\sp Now we shall prove   the following.

\blem\label{l2052405} 
If $f:\C\lra \oc$ is an elliptic function, then for every $\varepsilon >0$ there exists $a=a(\varepsilon) \geq 1$ such that  if $z\in {\mathbb C} \sms B_e({f}^{-1}(\infty), \varepsilon)$  then 
$$
z\notin
\bigcup_{j=a+1}^{\infty} B_e(f^j( \Crit_{\infty}(f)), 5).
$$
\elem

\bpf  Suppose on the  contrary  that there exists $
\varepsilon
>0$  for every $a \geq 1$  there exists $z_a\in
\bigcup_{j=a+1}^{\infty} B_e(f^j(\Crit_{\infty}(f)), 5)\sms
B_e({f}^{-1}(\infty), \varepsilon)$. Since the sets
$f^j(\Crit_\infty(f)) $   converge to $ \infty$ when $j\to \infty$,
it follows that $\lim_{a \to \infty}{z_a}=\infty$. But then  $z_a\in
B_e({f}^{-1}(\infty), \varepsilon)$ for all $a\geq 1$ large  enough.
This contradiction finishes the proof. 
\endpf

\sp We now pass to deal with compactly non--recurrent elliptic functions. Since for all such functions the sets ${\rm PC}_c(f)$ and ${\rm PC}_p(f)$ are bounded, the number
$$
D:=\frac{1}{2}\Dist_e({\rm PC}_c(f)\cup {\rm PC}_p(f),0)
$$ 
is finite. The main result of this section is the following.

\bprop\label{p1071305} 
Let $f:\C\lra\oc$ be a compactly non--recurrent elliptic function. If $z\in J(f)\sms \Sing^{-}(f)$ then there exist:

\,

\begin{enumerate}
\item[(a)] $ \eta(z)>0$. 

\, 

\item[(b)] $\{n_j\}_{j=1}^\infty$, an
increasing   sequence  of positive  integers. 

\,

\item[(c)] a sequence
$\{x_j(z)\}_{j=1}^\infty \sbt  J(f)\sms \(\Om(f) \cup\om
(\Crit_c(z))\)$ with  the following properties:

\,

\begin{enumerate}
\item[(1)]$ \Comp(z,f^{n_j}, \eta(z))\cap \Crit(f^{n_j})=\es$.

\, 

\item[(2)] $\lim_{j \to \infty}|f^{n_j}(z)-x_j(z)|=0$.

\,

\item [(3)] If $|x_j(z)|\geq 2D$  for all $j\geq 1$, then $\eta(z) \geq \min\{2,R_2\}$.

\,

\item [(4)] If  the sequence  $\{x_j(z)\}_{j=0}^{\infty}$ is bounded, then it is constant.

\,

\item [(5)] If $\lim_{j\to \infty}  x_j(z)=\infty$  and $z \notin I_{\infty}(f)$ then $x_j(z)\sim_f  x_k(z)$ for all $j,k\geq 1$.
\end{enumerate}
\end{enumerate}
\eprop

\bpf If $z\in I_{\infty}(f)$, equivalently, if $\lim_{n
\to \infty}\dist_e(f^n(z),{f}^{-1}(\infty))=0$, then  in view of
Lemma~\ref{l1052405}
 we are done by
setting $n_j=j+u$  and $x_j(z)=f^{j+u}(z)$ with some $u \geq 0$
large  enough, so that  $\dist_e(f^{j+u}(z),{f}^{-1}(\infty))< R_2$.
So suppose that  there exists  an $\varepsilon >0$  such that  the
set
$$
S=\big\{k \geq 0:\dist_e(f^k(z), {f}^{-1}(\infty))>\varepsilon\big\}
$$
is infinite.

 Now,  suppose that $A^d_S$, the set  of  limits points of $A_S=\{f^k(z):\,\,  k \in S\}$
is unbounded. There exists  $w\in A^d_S$  and $a(\varepsilon)\geq 1
$ ($a(\varepsilon)$ comes from  Lemma~\ref{l2052405})  such that
\beq\label{3052405}
B_e(w,5)\cap
\left(\bigcup_{j=0}^{a(\varepsilon)}f^j(\Crit_\infty(f))\cup {\rm
PC}_c(f) \cup {\rm PC}_p(f) \right)=\es
\eeq
and $|w|\geq 4D$. There  also exists an  increasing  sequence
$\{n_j\}_{j=1}^\infty\sbt S$ such that 
$$
\lim_{j \to \infty}f^{n_j}(z)=w.
$$
Disregarding  finitely many  elements of this
sequence, we may assume  without  loss of generality that
$f^{n_j}(z)\in B_e(w, \min\{1, D\}) $ for all $j \geq 1$. In view of
the definition  of $S$, Lemma~\ref{l2052405}, and (\ref{3052405}),
we see that for every $j \geq 1$ there exists a holomorphic inverse
branch $f^{-n_j}_z:B_e(f^{n_j}(z), 4) \mapsto \mathbb C$ sending
$f^{n_j}(z)$ to $z$. Because of  the same  premises $ w \notin
\om(\Crit(f))$ and we are done in this case  by setting $x_j(z)=w$
for all $j \geq 1$ and $\eta(z)=4$.

\sp So, suppose that the set  $A^d_S$ is  bounded. Assume first that
$$\liminf_{n \to \infty} \dist_e(f^n(z),f^{-1}(\infty))=0.$$
Then there exists $\{k_j\}_{j=1}^{\infty}$, an increasing sequence
of positive integers, such
\beq\label{1071405}
\dist_e (f^{k_j+1}(z), {f}^{-1}(\infty))> \varepsilon
\eeq
(i.e.  $k_j+1\in S$, $j \geq 1$) and we  require $f^{k_j}(z)$ to be
so close ${f}^{-1}(\infty)$ (assuming  $\varepsilon >0$  to be
sufficiently  small) that
\beq\label{1071305}
|f^{k_j+1}(z)|\geq 4D +2
\eeq
for all $j \geq 1$. Passing to a subsequence, we may assume without
loss of generality  that  the  sequence $\{\Pi
(f^{k_j+1}(z))\}_{j=1}^{\infty}$ on the   torus $\mathcal T$
converges to the same point $y\in \mathcal T$, where $\Pi:{\mathbb
C} \mapsto {\mathcal T}= {\mathbb C}/\sim_f $ is the  canonical
projection from $\mathbb C$ onto the torus $\mathcal T$. Clearly,
there exists a sequence $\{x_j(z)\}_{j=1}^{\infty}$ such that
$\lim_{j \to \infty} |f^{k_j+1}(z)-x_j(z)|=0$ and $\Pi (x_j(z))=y$
for all $ j \geq 1$.

\sp Assume  first  that   the sequence $\{f^{k_j+1}(z)\}_{j=1}^\infty $ is unbounded. Passing  to a subsequence we
may assume without  loss  of generality that $\lim_{j\to
\infty}f^{k_j+1}(z)=\infty.$  Then, applying (\ref{1071405}),
Lemma~\ref{l2052405}, (\ref{1071305}) and  the definition of $D$, we
are  done with $n_j=k_j+1$ and $\eta(z)=2$. So, assume   that  the
sequence  $\{f^{k_j+1}(z)\}_{j=1}^\infty$  is   bounded. We already
know that
$$B_e(f^{k_j+1}(z), 2)\cap \left(\ov{{\rm PC}_c(f)}\cup \ov{{\rm PC}_p(f)}\cup \ov{\bigcup_{j=a(\varepsilon)+1}^\infty
f^j(\Crit_\infty(f))}\right)=0.$$ Since the   second  component  of
this intersection  is forward-invariant, we  conclude  that   no
accumulation  point  of the  sequence
$\{f^{k_j-a(\varepsilon)}(z)\}_{j=1}^\infty$ belongs to
$$\ov{{\rm PC}_c(f)}\cup \ov{{\rm PC}_p(f)}\cup \ov{{\rm PC}_\infty(f)}=\ov{{\rm PC}(f)}.$$  If the  sequence
$\{f^{k_j-a(\varepsilon)}(z)\}_{j=1}^\infty$ is unbounded, we  may
complete  the argument  in exactly the same  way as above with
$k_j+1$ replaced  by $k_j-a(\varepsilon).$ If the  sequence
$\{f^{k_j-a(\varepsilon)}(z)\}_{j=1}^\infty$ is bounded, we  are
immediately  done  by passing to a  converging subsequence.

\sp So assume that
$$
\liminf_{n \to \infty} \dist_e(f^n(z),f^{-1}(\infty))>0.
$$
Then $\liminf_{n \to \infty}|f^n(z)|<\infty$ and the $\om$-limit set
is compact in the plane $\mathbb C$.  In view of Lemma~\ref{ld5.7}
there exists $x\in \om(z)\sms (\Om(f)\cup\ov{O_+(f(\Crit_c(z))}\cup
\{\infty\})$. The number
$$
\eta:=\frac{1}{2}\dist_e\(x,\Om(f)\cup\ov{O_+(f(\Crit_c(z))}\)
$$ 
is positive since $\om(\Crit_c(z))$ is a compact subset of $\mathbb C$
and $\Om(f)$ is finite. Then there exists an infinite increasing
sequence $\{m_j\}_{j\ge 1}$ such that
\beq\lab{d6.1}
\lim_{j\to\infty}f^{m_j}(z)=x
\eeq
and
\beq\lab{d6.2}
B_e(f^{m_j}(z),\eta)\cap\bu_{n\ge 1}f^n(\Crit_c(z))=\es.
\eeq
Now we claim that there exists $\eta(z)$ such that for every $j\ge
1$ large enough
\beq\lab{d6.3}
\Comp(z,f^{m_j},\eta(z))\cap \Crit(f^{m_j})=\es.
\eeq
Otherwise we would find an increasing to infinity subsequence
$\{m_{j_i}\}$ of $\{m_j\}$ and a decreasing to zero sequence of
positive numbers $\eta_i$ such that $\eta_i<\eta$ and
$$
\Comp(z,f^{m_{j_i}},\eta_i)\cap \Crit(f^{m_{j_i}})\ne
\es.
$$
Let $\^c_i\in \Comp(z,f^{m_{j_i}},\eta_i)\cap
\Crit(f^{m_{j_i}})$. Then there exists $c_i\in \Crit(f)$ such that
$f^{p_i}(\^c_i)=c_i$ for some $0\le p_i\le m_{j_i}-1$. Since the set
$f^{-1}(x)$ is at a positive distance from $\Om(f)$ and since
$\eta_i\to 0$, it follows from Theorem~\ref{mnt6.3} that
$\lim_{i\to\infty}\^c_i=z$. Since $z\notin \bu_{n\ge
0}f^{-n}(\Crit(f))$, it implies that $\lim_{i\to\infty}p_i=\infty$.
But then using Theorem~\ref{mnt6.3} again and the formula
$f^{p_i}(\^c_i)=c_i$ we conclude that the set of all accumulation
points of the sequence $\{c_i\}$ is contained in $\om(z)$. Hence,
passing to a subsequence, we may assume that the limit
$c=\lim_{i\to\infty}c_i$ exists. But since $c\in \om(z)$, since
$\om(z)$ is a compact subset of $\mathbb C$ and since $\infty$ is
the only accumulation point of $\Crit(f)$, we conclude that the
sequence $c_i$ is eventually constant. Thus, dropping some finite
number of initial terms, we may assume that this sequence is
constant. This means that $c_i=c$ for all $i=1,2,\ld$. Since
$c=f^{p_i}(\^c_i)$, we get
$$
|f^{m_{j_i}}(z)-f^{m_{j_i}-p_i}(c)|
=|f^{m_{j_i}}(z)-f^{m_{j_i}}(\^c_i)| <\eta_i.
$$
Since $\lim_{i\to\infty}\eta_i=0$ and since $\om(z)$ is a compact
subset of $\mathbb C$, we conclude that
$$\lim_{i\to\infty}|f^{m_{j_i}}(z)-f^{m_{j_i}-p_i}(c)|=0.$$ Since $c\in\Crit(z)$,
in view of (\ref{d6.2}) this implies that $m_{j_i}-p_i\le 0$ for all
$i$ large enough. So, we get a contradiction as $0\le p_i\le
m_{j_i}-1$ and (\ref{d6.3}) is proved. 
\endpf

\sp In particular, with  the  notation of this
proposition, for every integer $j \geq 1$ there exists the holomorphic inverse branch 
$$
f^{-n_j}_z:B_e(x_j(z), \eta(z))\lra \mathbb C.
$$
For every integer $j \geq 1$, let $T_j:{\mathbb C} \to {\mathbb
C}$ be a translation given by the formula
$$
T_j(w)=w+x_j(z).
$$

We shall  prove  the following.

\blem\label{l1052504} With hypotheses and notation of
Proposition~\ref{p1071305} the family 
$$
{\mathcal F}_z:=\big\{f^{-n_j}_z\circ T_j: B_e(0, \eta(z)) \lra {\mathbb C}\, \big| \, j
\geq 1 \big\}
$$ 
is normal and all its limit functions are constant. 
\elem

\bpf Decreasing $\eta(z)>0$ if necessary, we can always
find a periodic orbit of $f$ of period $\geq$ 3 disjoint  from all
the balls $B_e(x_j(z), \eta(z))$. Then this  orbit
 is also  disjoint  from  all   the sets
 $$f^{-n_j}_z\circ T_j(B_e(0, \eta(z))= f^{-n_j}_z\circ T_j(B_e(x_j(z),
 \eta(z))).$$
Hence, by Montel's Theorem, ${\mathcal F}_z$ is a normal family. If
there were   non-constant limit functions of the family ${\mathcal
F}_z$, then  there would exists a  radius $r>0$  and  an increasing
subsequence $\{n_{j_k}\}_{k=1}^\infty$ such that $$T_{j_k}^{-1}\circ
f^{n_{j_k}}(B_e(z, r)) \subset  B_e(0, \eta(z))$$ or equivalently
$$f^{n_{j_k}}(B_e(z, r))\sbt T_{j_k}(B_e(0, \eta(z))).$$
Passing  yet to another  subsequence, we may assume  without loss of
generality that
$$\oc\sms \bigcup_{k=1}^{\infty} T_{j_k}(B_e(0, \eta(z)))$$
has a non-empty  interior, and consequently contains at least three
points.  Thus  the family 
$$
\{f^{n_{j_k}}:B_e(z,r)\lra {\mathbb C}\}_{k=1}^\infty
$$ 
would  be  normal, contrary   to the   fact
$z\in J(f)$. We are  done. 
\endpf

\sp As  an immediate  consequence  of this  lemma, we get  the
following.

\bcor\label{c1041806} Let $f:\C\lra\oc$ be a compactly non--recurrent elliptic function. If $z\in J(f)\sms \Sing^{-}(f)$ and  the
sequence  $\{n_j\}_{j=1}^\infty$ is taken from Proposition~\ref{p1071305}, then
$$ 
\lim_{k \to \infty}|(f^{n_{j_k}})'(z)|=+\infty.
$$
\ecor

\sp\section{Dynamically Distinguished Classes of Elliptic Functions}\label{DDCoEF} 

In this section we recall already defined and define new classes of dynamically significant elliptic functions. These will form essentially all major classes of elliptic functions dealt with in this book. Mostly, they will be defined in terms of how strongly expanding these functions are. 

We would like to add that while in the theory of rational functions such classes pop up in a natural and fairly obvious way, and for example metric and topological definitions of expanding rational functions describe the same class of functions, in the theory of iteration of transcendental meromorphic functions such classification is by no means obvious, topological and metric analogs of rational function realm concepts do not usually coincide, and the definitions of expanding, hyperbolic, topologically hyperbolic, subhyperbolic, etc, functions vary from author to author. Our definitions seem to us pretty natural and fit well for our purpose of detailed investigation of dynamical and geometric properties of elliptic functions they define.  

We start with the following.

\bobs\label{o2_2017_09_27}
If $f:\C\lra\oc$ is an elliptic function, then for every point $z\in \C\sms \bigcup_{ n\geq 1}f^{-n}(\infty)$, we have that
$$
\lim_{n\to\infty}f^n(z)=\infty 
\  \ {\rm\text{if and only if}} \ \
\limsup_{n\to\infty}\dist_e(f^n(z),f^{-1}(\infty))=0.
$$
\eobs

\fr Let
\beq\label{4_2017_11_08}
\d(f^{-1}(\infty)):=\min\big\{|b-a|:a,b\in f^{-1}(\infty) \ {\rm and} \ a\ne b\big\}>0.
\eeq
Thus, for every $c\in \Crit_{\infty}(f)$ and all integers $n\ge 1$ large enough there exists a unique pole $b_n$ of $f$ such that 
\beq\label{3_2017_11_08}
|f^n(c)-b_n|<\d(f^{-1}(\infty))/2.
\eeq
Let
\beq\label{q_c}
q_c=\limsup_{n \to \infty} q_{b_n}\index{(S)}{$q_c$},
\eeq
where, we recall, $q_{b_n}$ is the multiplicity the pole $b_n$. The integer $p_c\ge 2$, i.e. the multiplicity of the critical point $c$, was defined at the beginning of Section~\ref{Local Properties of Critical Points of Holomorphic Functions}. Let 
\beq\lab{120200120B}
l_{\infty}:=\max\{p_cq_c:  c\in  \Crit_{\infty}(f)\}
\eeq
\index{(S)}{$l_{\infty}$}
(if $\Crit_{\infty}(f)=\es$, then $l_{\infty}=0$). We say that the elliptic function $f:{\mathbb C} \lra \oc$ is {\em regular}\index{(N)}{regular elliptic function} if and only if 
\beq\lab{3121605}
h:=\HD(J(f))>\frac{2 l_{\infty}}{l_{\infty}+1}.
\eeq
Regularity of elliptic functions will be our standing hypothesis from Section~\ref{Elliptic Regularity} onward essentially through the end of the book. 

\sp We also recall that an elliptic function $f:\mathbb C \lra \oc$ is  called non--recurrent or compactly non--recurrent if and only if respective definitions~\ref{pseudo-non-recurrent_B} or \ref{pseudo-non-recurrent} are satisfied. Of Course:

\sp\bobs\label{o120200120}
Every compactly non--recurrent elliptic function is non--recurrent.
\eobs

We will now provide the definitions of all those classes of elliptic functions $f$ dealt with in this section, for which $\Om(f)=\es$, where, we recall, $\Om(f)$ is the set of all rationally indifferent periodic points of $f$. Later in this sections we will deal with classes of elliptic functions for which $\Om(f)\ne \es$.

\bdfn\label{d1cef2} We  say that  an elliptic function $f: \mathbb C \lra \oc$ is normal \index{(N)}{normal} if and only if 
\begin{equation}\label{120200118}
\(f^{-1}(\infty)\cup \Crit (f)\) \cap \bigcup_{n=1}^{\infty}f^n(\Crit(f))=\es.
\end{equation}
\edfn

\bdfn\label{d2cef2}  We say that  an elliptic function $f: \mathbb C \lra \oc$ is of finite character  \index{(N)}{finite character} if and only if 
$$
\Crit_\infty(f)=\es.
$$
\edfn

\bdfn\label{d3cef2}  
We say that an elliptic function $f: \mathbb C \lra \oc$ is weakly semi--expanding  \index{(N)}{weakly semi--expanding} if and only if it is non--recurrent and $\Om(f)=\es$. 
\edfn

\bdfn\label{d4cef2}  We  say that  an elliptic function $f: \mathbb C \lra\oc$ is  semi--expanding  \index{(N)}{semi--expanding} if and only if it is compactly  non-recurrent and $\Om(f)=\es$. 
\edfn

\bdfn\label{d5cef2}  
We  say that  an elliptic function $f: \mathbb C \lra\oc$ is subexpanding \index{(N)}{subexpanding} if and only if it is semi--expanding and
\begin{equation}\label{220200118}
\om\(f(\Crit (f)\) \cap J(f)\cap \Crit(f)=\es.
\end{equation}
\edfn

\bdfn\label{d6cef2}  
We say that an elliptic function $f: \mathbb C \lra\oc$ is expanding \index{(N)}{expanding} if and only
\begin{equation}\label{320200118}
J(f)\cap \ov{{\rm PC}^0(f)}=\es,
\end{equation}
where, we  recall ${\rm PC}^0(f)=\bigcup_{n=0}^\infty f^n(\Crit(f)).$\edfn

We shall now bring up some inclusion relations between these classes of elliptic functions. All of them except one are obvious. We start with the following.

\bobs\label{o1cef2.1} 
Every elliptic function  of finite character is regular. \index{(N)}{regular}
\eobs

\bprop\label{p1cef3} Each expanding elliptic function is normal subexpanding of finite character.\eprop

\bpf Let $f: \mathbb C \lra\oc$ be an expanding elliptic function. It is compactly non--recurrent because $\Crit(f)\cap J(f)=\es$. It follows from Theorem~\ref{t2ms135}and formula (\ref{320200118}) that $\Om(f)=\es$. Hence, $f$ is semi--expanding. It is subexpanding  since $J(f)\cap \Crit(f)=\es$. By the same token $f$ is normal and of finite character.\qed

\sp So, we get the following.

\bprop\label{p6cef3} Referring to elliptic functions the following inclusions hold:
$$
\bal
Expanding \, &\sbt\, Subexpanding\, \sbt\, Semi-expanding \, \sbt \,Compactly \,  non-recurrent \\
&\sbt\, Non-recurrent
\eal
$$ 
and 
$$
Weakly\, semi-expanding \, \sbt\, Semi-expanding.
$$
\eprop

We shall now provide  same useful  characterizations  of some of these classes of elliptic function. We will  first  prove the following fact, which is interesting
by itself and which will be heavily used in this and  the next section. In the theory of iteration of rational functions of the Riemann sphere this fact is commonly
referred to as  the exponential shrinking property. We follows this tradition. 

\sp\bthm\label{t1j237} {\rm (Exponential Shrinking \index{(N)}{exponential shrinking} on $ \mathbb C$)} If $f:\C\lra\oc$ is a weakly semi--expanding
elliptic function, then  there exist $ 0<
\hat{\gamma} \leq \gamma_f $ and $M_f>0$ such that
$$
\diam_e(\Comp (z, f^n, \hat{\gamma})) \leq 8 e^{-M_f n}
$$ 
for all integers $n \geq 0$ and all $ z \in f^{-n}( J(f))$.
\ethm

\bpf In view of the last assertion of Lemma~\ref{felman},
there exist $ 0< \hat{\gamma} \leq \gamma $ and $l\geq 1$ such that
\beq\label{1j239}
\diam_e(\Comp(z,f^k,\hat{\gamma}))< \hat{\gamma}/2
\eeq
for  all $ k \geq l$ and all $ z \in f^{-k}(J(f))$. For every $k\geq
1$ let
\beq\label{3j239}
B_k(z):=\Comp (z, f^{kl}, \hat{\gamma}) \quad \mbox{and}\quad
A_k(z):= B_k(z)\sms \ov{B}_{k+1}(z).
\eeq
Observe that $f^{kl}(B_{k+1}(z))$ is a  connected subset of $\mathbb
C$ containing the point $f^{kl}(z)$, and $f^l(f^{kl}(B_{k+1}(z)))=B(f^{(k+1)l}
(z), \hat{\gamma})$. So,
 $$f^{kl}(B_{k+1}(z)) \sbt \Comp(f^{kl}(z), f^{l}, \hat{\gamma}) \sbt
 B(f^{kl}(z), \hat{\gamma}/2).$$
Hence
\beq\label{2j239}
f^{kl}(\ov{B}_{k+1}(z))\sbt \ov{B}(f^{kl}(z), \hat{\gamma}/2) \sbt
B(f^{kl}(z), \hat{\gamma}).
\eeq
Since $\ov{B}_{k+1}(z)$ is  a connected subset containing $z$, we
thus conclude that
\beq\label{4j239}
\ov{B}_{k+1}(z) \sbt B_k(z).
\eeq
Also, by (\ref{2j239}), 
$$ 
A_k(z)  \supset \Comp(z, f^{kl},
\hat{\gamma})\sms \ov{\Comp}(z, f^{kl}, \hat{\gamma}/2).
$$ 
Let $A_k^*(z)$  be the connected component  of $A_k(z)$ separating $z$
from $\ov{\mathbb C} \sms B_k(z)$. Then $A_k^*(z)$   is a
topological  annulus  containing  $(\Comp(z, f^{kl},\g)\sms
\ov{\Comp}(z, f^{kl},\g/2))_{*}$, and, by Corollary~\ref{c1j225} and
Lemma~\ref{lmod}, we get that
\beq\label{5j239}
\Mod(A^*_k(z)) \geq \Mod ((\Comp(z, f^{kl}, \hat{\gamma})\sms
\ov{\Comp}(z, f^{kl}, \hat{\gamma}/2))_* )\geq M_f,
\eeq
where  $M_f:=\log 2/(N_f^*)^2$. By virtue of (\ref{3j239}), (\ref{4j239}) and the definition of $A_k^*(z)$, $ k \geq 1$, we get that $(A_k^*(z))_{k
\geq 0}$ is a sequence of mutually disjoint annuli separating  $z$
from $\partial{B}(z, \hat{\gamma})$. It also  follows  from
(\ref{3j239}) and (\ref{4j239}) that
$$ A^{*}_k(z) \sbt A_k(z) \sbt B(z, \hat{\gamma}) \sms
\ov{B}_{k+1}(z) \sbt B(z, \hat{\gamma}) \sms \ov{B}_n (z)$$ for all
$k=0, 1, \ldots, n-1$. So $ A^{*}_k(z) \sbt (B(z, \hat{\gamma})\sms
\ov{B}_n(z))_*$. So,  $(A^{*}_k(z))_{k=0}^{n-1}$ is  a sequence of
mutually  disjoint annuli essentially contained in the annulus
$((B(z,r)\sms \ov{B}_n(z))_{*}$. It therefore follows from
Theorem~\ref{t1j223}  and (\ref{5j239}) that
 $$\Mod ((B(z, \gamma)\sms \ov{B}_n)_*)\geq
 \sum_{k=0}^{n-1}\Mod(A^{*}_k(z))\geq M_f n.$$
 So, apply  Theorem~\ref{t1j233} we get
 $$ \diam_e(\ov{B}_n) \leq 2 r_z (\ov{B}_n)  \leq 8e^{-M_fn}.$$ We
 are done. \endpf

\sp\fr As an immediate  consequence of this theorem we get the
following.

\sp\bcor\label{c1j249} {\rm (Exponential Shrinking on $\mathbb T_f$)} \index{(N)}{exponential shrinking} If
$f:\C\lra\oc$ is a weakly semi--hyperbolic
elliptic function, then there exist $ 0\leq
\hat{\g} \leq  \g$ and $M_f>0$ such that 
$$ 
\diam(\Comp(z,\hat f^n,\hat{\g})) \leq 8 e^{-M_fn}
$$ 
for all $n \geq 0$ and all $ z \in \hat{f}^{-n}(J(\hat{f}))$. 
\ecor

Now we shall  prove the following characterization of expandingness  which fully  justifies the name.

\bthm\label{t2cef3}  
An  elliptic function $f: \mathbb C \to \ov{\mathbb C}$ is expanding if and only if \nl 
$\exists q\geq 1$\,  $\forall \xi \in J(f)$\, $\forall z \in f^{-q}(\xi)$
\begin{equation}\label{1cef3}
|(f^q)'(z)|\geq 2.
\end{equation}
\ethm

\bpf Assume first that $f: \mathbb C \lra\oc$ is expanding. It then  follows  from  Definition~\ref{d1cef2} that $\Om(f)=\es$  and $J(f)\cap \Crit(f)=\es$. It therefore  follows from  Theorem~\ref{Fatou Periodic Components}, Theorem~\ref{baker-domain+Sullivan for elliptic} and finiteness of the set $f(\Crit(f))$, that
\begin{equation}\label{2cef3}
\om(\Crit (f))\sbt A \sbt F(f),
\end{equation}
where $A$ is the finite set of all attracting points of $f$. Along with the definition expandingness, this implies that
\begin{equation}\label{3cef3}
R=\frac{1}{2}\dist_e\(O^{+}(\Crit(f)/),  J(f)>0.
\end{equation}
Hence, for all $n \geq 1$, and all $\xi \in J(f)$, and all $z \in f^{-n}(\xi)$, there exists 
$$
f^{-n}_z: B_e(\xi, 2R)\to \mathbb C,
$$
a unique holomorphic branch of $f^{-n}$ sending $\xi $ to $z$. The proof (\ref{1cef3}) is now concluded by applying Theorem~\ref{t1j237} along with  Koebe's Distortion Theorem. 

For  the opposite direction, suppose that (\ref{1cef3}) holds. Then $J(f)\cap \Crit(f)=\es$ and $J(f)\cap \Om(f)=\es$.  So, exactly as in the first part of the proof, we conclude that formula (\ref{2cef3}) holds, and furthermore, (\ref{3cef3}) holds. In, particular, $\ov{O^{+}(\Crit(f))}\cap J(f)=\es$, meaning that  $f: \mathbb C \lra\oc$ is expanding. The proof of Theorem~\ref{t2cef3} is complete.
\qed

\

Now we  shall  prove a similar characterization  of subexpanding elliptic  functions.

\bthm\label{t3cef3}
 An  elliptic function $f:\lra\oc$ is subexpanding if and only if there exists $q\geq 1$ such that  for all  $ \xi \in J(f)\cap \om(\Crit(f))$ and all $z \in f^{-q}(\xi)\cap J(f)\cap \om(\Crit(f))$, we have that
 \begin{equation}\label{2cef4}
|(f^q)'(z)|\geq 2
\end{equation}
and for  every point $c \in J(f)\cap \Crit(f)$ either $\om(c)$  is a  compact subset of $\mathbb C$, $ c \in \bigcup_{n=1}^\infty f^{-n}(\infty)$  or the sequence $(f^n(c))_{n=1}^\infty$ converges  to $\infty$ (meaning that $c\in\Crit_\infty(f))$.
\ethm

\fr{\sl Proof} Assume first that $f: \mathbb C \lra\oc$ is expanding. So, $f$ is compactly non--recurrent, and we only are to show that formula (\ref{2cef4})  holds. Since
the set $\Crit(f)$  has only finitely many  congruent classes  modulo $\Lambda_f$, since  $\dist_e(\Crit(f), f^{-1}(\infty))>0$, and  since $f$ is compactly non--recurrent, it follows from  formula (\ref{220200118}) that
\begin{equation}\label{1cef4}
R=\frac{1}{2}\dist_e\(\om(\Crit (f))\cap J(f),\Crit(f)\)>0.
\end{equation}
Since, also $\Om(f)=\es$,  it follows  from Theorem~\ref{mnt6.3}, applied with 
$$
X:=\(\om(\Crit(f))\cap J(f)\)\cup\{\infty\} 
$$ 
and $\varepsilon:=\frac{1}{2}R$, that there exists $\delta >0$ such that for all $n \geq 1$, also $0\leq k \leq n$, all $\xi \in J(f)\cap \om(\Crit(f))$, and all $z \in f^{-n}(\xi)$, we  have that
$$
\diam_e\(f^k(\Comp(z,f^n, 2 \delta)\)<R.
$$
So, if addition $z \in J(f)\cap \om(\Crit(f))$, then
$$
f^k\(\Comp(z,f^n,2 \delta)\)\sbt B_e\(J(f)\cap\om(\Crit(f)),R\)
$$ 
for all $ 0\leq k \leq n$. Along with (\ref{1cef4}), this gives that 
$$
f^k\(\Comp(z,f^n,2 \delta)\)\cap \Crit(f)=\es
$$ 
for all $0\leq k\leq n$. Hence, there exists
$$
f^{-n}_z: B_e(\xi, 2 \delta) \lra \mathbb C,
$$
a unique holomorphic branch of $f^{-n}$ sending $\xi$ to $z$. The proof of (\ref{2cef4})  is now concluded by applying Theorem~\ref{t1j237} along with Koebe's Distortion Theorem. 

For the opposite implication, suppose that its hypotheses hold. We are then only to show that formula (\ref{220200118}) holds and $\Omega(f)=\es$, but both of them follow immediately from (\ref{2cef4}) and the fact that $\Om(f)\sbt \om(\Crit(f))\cap J(f) $  (which in turn follows  from Theorem~\ref{t2ms135}). The proof of Theorem~\ref{t3cef3} is complete. 
\qed

\sp We now  shall prove the following

\bthm\label{t1cef5} If $f:\mathbb C \lra \oc$ is a normal elliptic function of finite character, then $f$ is subexpanding if and only if  $\ov{{\rm PC}(f)}\cap J(f)$ is a compact subset of $\mathbb C $, $ \Om(f)=\es$, and
\begin{equation}\label{2cef5}
\dist_e\(f^{-1}(\infty)\cup \Crit(f), {\rm  PC}(f)\cap J(f)\)>0.
\end{equation}
\ethm

\fr{\sl Proof}. Suppose  first that $f: \mathbb C \lra\oc$ is subexpanding. Then $\Om(f)=\es$. Since $f$ is normal and of finite character, we have that
$$
\Crit_p(f)=\Crit_\infty(f)=\es.
$$
Hence, since $f$ is  compactly non--recurrent, it follows that $\om(\Crit(f))\cap J(f)=\om_c(\Crit(f)\cap J(f))$ is a compact subset of $
\mathbb C$ and
\begin{equation}\label{1cef5}
f^{-1}(\infty)\cap \om(\Crit(f))\cap J(f)=\es.
\end{equation}
Thus ${\rm PC}(f)\cap J(f)$ is a bounded subset of $\mathbb C$, whence $\ov{{\rm PC}(f)}\cap J(f)$ is a compact subset of $\mathbb C$. Since 
$\om(f(\Crit(f))\cap J(f))=\om(f(\Crit(f))\cap J(f)$ is  a compact subset of $\mathbb C$ and the set $\Crit(f)\cup f^{-1}(\infty) \sbt \mathbb C$ is closed, it follows from (\ref{220200118}) and (\ref{1cef5})
$$
\dist_e\(f^{-1}(\infty)\cup \Crit(f), \om(f(\Crit(f))\cap J(f)\)>0.
$$ 
Having  this and invoking  normality of $f$, we deduce that  formula (\ref{2cef5}) holds.

Suppose in turn  that $\ov{\rm PC}(f)\cap J(f)$ is a compact subset of $\mathbb C$, $ \Om(f)=\es$, and formula (\ref{2cef5}) holds.  Then immediately formula (\ref{220200118})  holds, $f$ is compactly non--recurrent, and, moreover, semi--expanding. This means that $f$ is sub--expanding. The proof of Theorem~\ref{t1cef5} is complete. 
\qed
 
\sp As an immediate consequence of this theorem we get the following.

\bthm\label{t2cef5} If  $f:\mathbb C \lra\oc$ is a normal elliptic function of finite character, then $\ov{{\rm  PC}(\hat{f})}\cap J(\hat f)$ is a compact subset of $\hat\mT_f$ and 
$$
\dist_e\(B(\hat{f}), \ov{{\rm  PC}(\hat{f})}\cap J(\hat f)\)>0.
$$
\ethm

\sp The last class of elliptic functions we will dealing with in this book, more precisely in Sections~\ref{PEM-I}, \ref{parabolicfiniteinvariantmeasures}, and \ref{Darling_Kac-elliptic}, are parabolic elliptic functions\index{(N)}{parabolic elliptic function} i.e. all elliptic functions $f:\mathbb {C} \lra\oc$ for which
$$
\Crit(f)\cap J(f)=\es 
\ \ {\rm and } 
\  \  \Om(f)\ne\es.
$$
As an immediate of the first half of this definition we get the following observation. 

\sp
\bobs\label{o120200125}
Every parabolic elliptic function is regular normal compactly non--recurrent of finite character.
\eobs 

\sp We furthermore discern in the Sections~\ref{PEM-I}, \ref{parabolicfiniteinvariantmeasures}, and \ref{Darling_Kac-elliptic}, mentioned above two disjoint subclasses: of finite class and of infinite class, respectively depending on whether the invariant measure $\mu_h$ of Theorem~\ref{t120190619} is finite or infinite.

\chapter{Various Examples of Compactly Non--Recurrent Elliptic
Functions}\label{examples}

The purpose of this chapter is to provide examples of elliptic
functions with prescribed properties of the orbits of critical
points (an values). We are primarily focused on constructing examples of various classes of compactly non--recurrent elliptic functions discerned in Section~\ref{DDCoEF}. All these examples are either Weierstrass 
$\wp_\La$ elliptic functions or their modifications. The dynamics of such functions depends heavily on the lattice $\La$ and varies drastically from $\La$ to $\La$. 

The first three sections of this chapter have a preparatory character and respectively describe basic dynamical and geometric properties of all Weierstrass $\wp_\La$ elliptic functions, the ones generated by square lattices, and triangular lattices. 

In Section~\ref{Simple Examples of Elliptic Functions} we provide simple constructions of many classes of elliptic functions discerned in Section~\ref{DDCoEF}. We essentially cover all of them. All these examples stem from 
Weierstrass $\wp$ functions.

We then, starting with Section~\ref{connectedjulia}, provide also some different, interesting by their own, and historically first examples of various kinds of Weierstrass $\wp$ elliptic functions and their modifications. These come from the series of papers \cite{HK1}, \cite{HK2}, \cite{HK3}, \cite{HKK}, \cite{HL} by Jane Hawkins and her collaborators. Throughout the whole current chapter we use the notation and terminology introduced in Chapter~\ref{elliptic-theory}, particularly in Sections~\ref{weierstrass}--\ref{weierstrassII}.

\section{The Dynamics of Weierstrass Elliptic Functions: Some Selected General Facts} 
 
As an immediate consequence of Theorem~\ref{baker-domain+Sullivan for elliptic}, we get the following.

\sp\bthm\label{t120200303} 
If $\wp:{\mathbb C}\lra\oc$ is a non--constant Weierstrass elliptic function, then $f$ has no Baker or wandering domains. 
\ethm

As an immediate consequence of this theorem, Theorem~\ref{Fatou Periodic Components} (Fatou Periodic Components), Theorem~\ref{t1ms123}, Theorem~\ref{t2ms135}, Theorem~\ref{t320200303}, and Observation~\ref{o420200303}, we get the following. 

\sp\bthm\label{t220200316} 
If $\wp:{\mathbb C}\lra\oc$ is a non--constant Weierstrass elliptic function that has no Siegel disks or Herman rings, then 

\begin{enumerate} 

\item Either the Julia set $J(\wp)$ is the whole sphere $\oc$ or else all Fatou connected components are basins of attractions to (super) attracting periodic points or rationally indifferent periodic points. 

\,

\item There are at most three cycles of periodic components of the Fatou set $F(\wp)$ and each of them contains a critical value of $\wp$ which is not preperiodic.
\end{enumerate}
\ethm

Invoking in addition Theorem~\ref{t1sh11}, we get the following. 

\sp\bthm\label{t120200316} 
If $\wp:{\mathbb C}\lra\oc$ is a non--constant Weierstrass elliptic function that has three periodic orbits each of which is either (super) attracting or rationally indifferent, then the collection of Fatou periodic connected components of $\wp$ consists of immediate basins of attraction to these three periodic orbits and the Fatou set of $\wp$ is the union of basins of attraction to these three periodic orbits.
\ethm

As an immediate consequence of Theorem~\ref{t220200316} and Theorem~\ref{r071708} we get the following.

\sp\bthm\label{t120200303B} 
If $\wp:{\mathbb C}\lra\oc$ is a non--recurrent Weierstrass elliptic function, then 

\begin{enumerate} 
\item $\wp$ has no Baker or wandering domain, no Siegel disks or Herman rings nor Cremer periodic points. 

\,

\item Either the Julia set $J(\wp)$ is the whole sphere $\oc$ or else all Fatou connected components are basins of attractions to (super) attracting periodic points or rationally indifferent periodic points. 

\,

\item There are at most three periodic components of the Fatou set $F(\wp)$ and each of them contains a critical value of $\wp$ which is not preperiodic.
\end{enumerate}
\ethm

\sp Since each Weierstarss  $\wp$ function is even, as an immediate consequence of Theorem~\ref{herman} we get the following.

\bthm\label{hermanB} 
For any lattice $\La$, the Weierstarss $\wp_\La$ function has no cycle of Herman rings. 
\ethm

\sp Since each Weierstrass elliptic function is even, as an immediate consequence of Theorem~\ref{t2ex4}, we get the following.

\sp\bthm\label{th3.1in{HK3}} If $\wp:\mathbb C\lra \oc$ is a Weierstrass elliptic function such that each connected component of the Fatou set $F(\wp_\Lambda)$ contains at most one critical value of $\wp$, then the Julia set $J(\wp)$ is connected.
\ethm

\sp \section[The Dynamics of Square Weierstrass Elliptic Functions] {The Dynamics of Square Weierstrass Elliptic Functions: Some Selected Facts}

Square Weierstrass elliptic functions will play an important role in creating many of our examples. We will need several of its properties. The first one is the following.

\bprop\label{symmetry}
If $\La$ be a square lattice, then  the Fatou set  $F(\wp_\La)$ and the Julia set $J(\wp_\La)$ of the Weierstrass $\wp_\La$
function are invariant under the multiplication by $i$, or, in geometric terms, under rotation about the angle $\pi/2$, i.e.
$$
i J(\wp_\La)=J(\wp_\La)
\  \  \  {\rm and} \  \  \
i F(\wp_\La)=F(\wp_\La).
$$
In addition,
\beq\label{Aiterates}
\wp_\La^n (i z)= i\wp_\La^n (z)
\eeq
for all integers $n \ge 0$ and all $z\in\C$ whenever $\wp_\La^n(z)$ is well defined.  
\eprop
 \bpf Formula \eqref{Aiterates} directly  follows, by a straightforward induction, from the homogeneity property (\ref{2}. 
Take $z \in F(\wp_\La)$ and an open neighborhood $U$  of $z$ such that the family  $\( \wp_\La^n|_U\)_{n=1}^\infty$ is normal. By (\ref{Aiterates}), we have
$$
\wp_\La^n(iU) =i\wp_\La^n(U)
$$
for all $n \geq 1$. Thus, the family$\( \wp_\La^n|_{iU}\)_{n=1}^\infty$ is normal. So, $iz\in F(\wp_\La)$. Hence $iF(\wp_\La\sbt F(\wp_\La)$. The proof of the converse inclusion is analogous. Thus the second assertion of our proposition is established. The first one then follows immediately since $J(\wp_\La)=\C\sms F(\wp_\La)$. The proof is complete.
\qed

\sp We now shall prove the following.

\sp\bthm\label{St620200303C} 
Let $\wp:\mathbb C\lra \oc$ be a square Weierstrass elliptic function. 

If $\xi\in\C$ is either a superattracting, attracting, or rationally indifferent periodic point, then 

\ben \item $i\xi$ is also respective superattracting, attracting, or rationally indifferent periodic point and both $\xi$ and $i\xi$ lie on the same periodic orbit.
 
\item
$$
\wp'(\xi)=\wp'(i\xi).
$$ 



\item If $\xi$ is a superattracting periodic point of $\wp$, then

\,

\ben
\item there is exactly one superattracting periodic orbits of $\wp$, namely the one $\xi$.

\item there is exactly one periodic connected components of the Fatou set $F(\wp)$, namely the basin of immediate attraction of $\wp$ to $\xi$.

\,

\item all connected components of the Fatou set $F(\wp)$ are  basins of attraction of $\wp$ to $\xi$.
\een

\,

\item 
If $\xi$ is an attracting (but not superattracting) periodic point of $\wp$, then

\ben

\item there is exactly three attracting (not superattracting) periodic orbits of $\wp$, namely the one $\xi$.

\item there is exactly one periodic connected components of the Fatou set $F(\wp)$, namely the basin of immediate attraction of $\wp$ to $\xi$.

\,

\item all connected components of the Fatou set $F(\wp)$ are  basins of attraction of $\wp$ to $\xi$. 
\een

\,

\item If $\xi$ is a rationally indifferent periodic point of $\wp$, then

\,

\ben 

\item there is exactly one rationally indifferent periodic cycles of $\wp$, namely the one $\xi$.

\item there is exactly one periodic connected components of the Fatou set $F(\wp)$, namely the basin of immediate attraction of $\wp$ to $\xi$.

\,

\item all connected components of the Fatou set $F(\wp)$ are  basins of attraction of $\wp$ to $\xi$. 
\een
\end{enumerate}
\ethm

\bpf   
Since $\wp$ is a square function, there is a square lattice $\La\sbt\C$ such that 
$$
\wp=\wp_\La.
$$
We first will show that if $f$ has one (super) attracting or rationally indifferent periodic point $\xi\in\C$ with some period $p\ge 1$, then $i\xi$ is also respective (super) attracting or rationally indifferent periodic points and 
$$
\wp'(\xi)=\wp'(i\xi).
$$ 
Furthermore, the periodic orbits of $\xi$ and $i\xi$ coincide. 

So, let $\xi\in\C$ be such a periodic point with period $p\ge 1$. Because of \eqref{Aiterates} we get,
$$
\wp_\La^p(i \xi)
=i\wp_\La^p(\xi)
$$
meaning that $i\xi$ is also periodic points of $\wp_\La$ with period $p$. Differentiating now both equations in \eqref{Aiterates} and using the Chain Rule, we get that
$$
i\(\wp_\La^p\)'(\xi)
=\(\wp_\La^p\circ i\)'(\xi)
=i\(\wp_\La^p\)'(i\xi).
$$
Hence,
$$
\(\wp_\La^p\)'(i\xi)=\(\wp_\La^p\)'(\xi).
$$
Now, since, by Proposition~\ref{p120200310} (6), $e_1(\La)=-e_2(\La)$ and since the elliptic function $\wp_\La$ is even, we have that $\wp_\La\(e_1(\La)\)=\wp_\La\(e_2(\La)\)$. Since, also by Proposition~\ref{p120200310} (6), the third  (and the last) critical value $e_3(\La)=0$, 
of $\wp_\La$ is a pole of $\wp_\La$, we directly deduce all the remining assertions of our theorem from Theorem~\ref{Fatou Periodic Components} (Fatou Periodic Components), Theorem~\ref{t1ms123}, Theorem~\ref{t2ms135}, Theorem~\ref{t320200303}, and Observation~\ref{o420200303}. The proof is complete.
\epf

\sp As an immediate consequence of this theorem, Theorem~\ref{t120200303B}, and Theorem~\ref{r071708}, we get the following.

\sp\bthm\label{t120200330} 
If $\wp:\mathbb C\lra \oc$ is a square Weierstrass non--recurrent elliptic function, then exactly one of the following holds 

\ben 
\item $J(\wp)=\C$,

\item $\wp$ has exactly one (super) attracting periodic orbit. Denote one of its points by $\xi$. Then

\ben 

\,

\item there is exactly one periodic connected components of the Fatou set $F(\wp)$, namely the basin of immediate attraction of $\wp$ to $\xi$.

\,

\item all connected components of the Fatou set $F(\wp)$ are basins of attraction of $\wp$ to $\xi$.
\een

\,

\item $\wp$ has exactly one rationally indifferent periodic orbit. Denote one of its points by $\xi$. Then

\ben 

\item there is exactly one rationally indifferent periodic cycles of $\wp$, namely the one $\xi$.

\item there is exactly one periodic connected components of the Fatou set $F(\wp)$, namely the basin of immediate attraction of $\wp$ to $\xi$.

\,

\item all connected components of the Fatou set $F(\wp)$ are  basins of attraction of $\wp$ to $\xi$. 
\een
\end{enumerate}
\ethm

\sp We recall from Definition~\ref{d120190216} (0) that a lattice $\Lambda\sbt\C$ is called real if and only if
$$
\ov\La=\La.
$$
We need the following properties of square lattices and properties of  Weierstrass $\wp_\La$ functions generated by such lattices. 

\sp\fr\bprop\label{Scor4.6in{HK1}}
Let $\La\sbt\C$ be a square lattice. 
\begin{itemize}
\item [(1)] Then the critical values of the Weierstrass $\wp_\La$ function
are 
$$
e_1=\frac{1}{2}g_2(\La)^{\frac{1}{2}}, \  \  \  e_2=-e_1,
\  \  \  {\rm and } \  \  \
e_3=0.
$$ 
In particular, $e_3$ is a pole of $\wp_\La$.

\,

\item [(2)] If the lattice $\La\sbt\C$ is real, then 
$$
{\rm PC}(\wp_\La)\sbt \{-e_1\}\cup\{0\}\cup [e_1, +\infty]\sbt\R\cup\{+\infty\}.
$$
\end{itemize}
\eprop
\bpf Item (1) is a reformulation of item (6) of Proposition~\ref{p120200310}. 

\sp Proving item (2), it follows from Proposition~\ref{p120200310} (5) and (6c5) that

\beq\label{120190216B}
\wp_\La([e_1,+\infty])\sbt [e_1,+\infty]\sbt\R\cup\{+\infty\}.
\eeq
Therefore,
$$
O_+(e_1)=\{\wp_\La^n(e_1):n\ge 0\}\sbt [e_1, +\infty]
$$  
Since also, Proposition~\ref{p120200310} (5) and evenness of $\wp_\La$, $e_3=0$, $e_2=-e_1$, and $\wp_\L(e_2)=\wp_\L(e_2)$, item (2) follows.
The proof is complete.
\qed

\sp Using this proposition we get the following.

\sp\bthm\label{t220200330} 
If $\wp:\mathbb C\lra \oc$ is a a real square Weierstrass  elliptic function, then exactly one of the following holds 

\ben 
\item $J(\wp)=\C$,

\,

\item $\wp$ has exactly one (super) attracting periodic orbit. Denote one of its points by $\xi$. Then

\ben 

\,

\item there is exactly one periodic connected components of the Fatou set $F(\wp)$, namely the basin of immediate attraction of $\wp$ to $\xi$.

\,

\item all connected components of the Fatou set $F(\wp)$ are basins of attraction of $\wp$ to $\xi$.
\een

\,

\item $\wp$ has exactly one rationally indifferent periodic orbit. Denote one of its points by $\xi$. Then

\ben 

\item there is exactly one rationally indifferent periodic cycles of $\wp$, namely the one $\xi$.

\item there is exactly one periodic connected components of the Fatou set $F(\wp)$, namely the basin of immediate attraction of $\wp$ to $\xi$.

\,

\item all connected components of the Fatou set $F(\wp)$ are  basins of attraction of $\wp$ to $\xi$. 
\een
\end{enumerate}
\ethm

\bpf 
The same arguments as those for Theorem~\ref{t120200330} work, so all what we need to  show is that $\wp_{\La}$ has no Siegel disks or Herman rings. We will do it now. So, seeking a contradiction suppose that $U$ is either a Siegel disk or Herman ring of the function $\wp_{\La}$. By Proposition~\ref{Scor4.6in{HK1}} and Theorem~\ref{t1sh11},
\beq\lab{S120200314}
\bd U\sbt\{-e_1\}\cup\{0\}\cup [e_1, +\infty]\sbt\R\cup\{+\infty\},
\eeq
where $e_1\in(0,+\infty)$. Since $\bd U\sbt J(\wp)$ and $J(\wp)$ has no isolated points, we further conclude that
\beq\lab{120200331}
\bd U\sbt [e_1, +\infty]\sbt\R\cup\{+\infty\},
\eeq
Since the point $0$ is a pole of $\wp_{\La}$, it belongs to the Julia set $J(\wp_{\La})$, whence it does not belong to $U$.
But by \eqref{120200331} and since $e_1\in(0,+\infty)$, the origin $0$ does not belong to the boundary $\bd U$ either. Therefore,
\beq\lab{220200331}
0\in \C\sms (U\cup \bd U).
\eeq
On the other hand, Since $U$ is not empty and open in $\C$, it is not contained in $[e_1, +\infty]$. So, there exists a point 
\beq\lab{320200331}
\xi \in U\sms [e_1, +\infty].
\eeq
Let $[0,\xi]$ be the closed line segment joining $0$ and $\xi$. Then, on the one hand,
\beq\lab{420200331}
[0,\xi]\sbt \C\sms [e_1, +\infty]\sbt \C\sms \bd U.
\eeq
On the the other hand, by \eqref{320200331} and \eqref{220200331}, we have that
$$
[0,\xi]\cap U\ne\es
\  \  \  {\rm and}  \  \  \
[0,\xi]\cap (\C\sms U)\ne\es.
$$
Along with \eqref{420200331}, this produces contradiction since the segment $[0,\xi]$ is connected. We are done.
\epf

\sp

\section[The Dynamics of Triangular Weierstrass Elliptic Functions] {The Dynamics of Triangular Weierstrass Elliptic Functions: Some Selected Facts}

\sp We recall from Definition~\ref{d120190216} (4), that a lattice $\La\sbt\C$ is called triangular if and only if
\beq\lab{120200304}
\epsilon\Lambda=\Lambda.
\eeq
where, we recall, 
$$
\epsilon=e^{\frac{2\pi i}{3}}.
$$

Triangular Weierstrass elliptic functions will play an important role in creating many of our examples. We will need several of its properties. The first one is the following.

\bprop\label{symmetryT}
If $\La$ be a triangular lattice, then  the Fatou set  $F(\wp_\La)$ and the Julia set $J(\wp_\La)$ of the Weierstrass $\wp_\La$
function are invariant under the multiplication by $\epsilon=e^{\frac{2\pi i}{3}}$, or, in geometric terms, under rotation about the angle $2\pi/3$, i.e.
$$
\epsilon^2 J(\wp_\La)=\epsilon J(\wp_\La)=J(\wp_\La)
\  \  \  {\rm and} \  \  \
\epsilon^2 F(\wp_\La)=\epsilon F(\wp_\La)=F(\wp_\La).
$$
In addition,
\beq\label{iterates}
\wp_\La^n (\epsilon z)= \epsilon\wp_\La^n (z)\quad \text{and} \quad \wp_\La^n (\epsilon^2 z)= \epsilon^2\wp_\La^n (z)
\eeq
for all integers $n \ge 0$ and all $z\in\C$ whenever $\wp_\La^n(z)$ is well defined.  
\eprop
 \bpf Formula \eqref{iterates} directly  follows, by a straightforward induction, from the homogeneity property (\ref{2}. 
Take $z \in F(\wp_\La)$ and an open neighborhood $U$  of $z$ such that the family  $\( \wp_\La^n|_U\)_{n=1}^\infty$ is normal. By (\ref{iterates}), we have
$$
\wp_\La^n(V) =\wp_\La^n(\epsilon U)=\epsilon  \wp_\La^n(U)
$$
for all $n \geq 1$. Thus, the family
$\( \wp_\La^n|_{\epsilon U}\)_{n=1}^\infty$ is normal. So, $\epsilon z\in F(\wp_\La)$. Hence $\epsilon F(\wp_\La\sbt F(\wp_\La)$. The proof of the converse inclusion is analogous. Thus the second assertion of our proposition is established. The first one then follows immediately since $J(\wp_\La)=\C\sms F(\wp_\La)$. This immediately implies that 
$$
\epsilon^2 F(\wp_\La)=F(\wp_\La)
\  \  \  {\rm  and}  \  \  \   
\epsilon^2 J(\wp_\La)=\epsilon J(\wp_\La)=J(\wp_\La).
$$
The proof is complete.
\qed

\sp We now shall prove the following.

\sp\bthm\label{t620200303C} 
Let $\wp:\mathbb C\lra \oc$ be a triangular Weierstrass elliptic function. 

If $\xi\in\C$ is either a superattracting, attracting, or rationally indifferent periodic point, then 

\ben \item $\epsilon\xi$ and $\epsilon^2\xi$ are also respective superattracting, attracting, or rationally indifferent periodic points,
 
\item
$$
\wp'(\xi)=\wp'(\epsilon\xi)=\wp'(\epsilon^2\xi).
$$

\item Either some two of the three periodic orbits of $\xi$, $\epsilon\xi$ or $\epsilon^2\xi$ intersect, or all three of them are mutually disjoint. 
\een

Assume the latter case holds. Then:

\begin{enumerate}


\item [(A)] If $\xi$ is a superattracting periodic point of $\wp$, then

\,

\ben
\item there are exactly three superattracting periodic orbits of $\wp$ and these are equal to the periodic orbits of $\xi$, $\epsilon\xi$, and $\epsilon^2\xi$,

\item the periodic connected components of the Fatou set $F(\wp)$ are respectively the three basins of immediate attraction of $\wp$ to the periodic orbits of $\xi$, $\epsilon\xi$, and $\epsilon^2\xi$,

\,

\item all connected components of the Fatou set $F(\wp)$ are respectively the three basins of attraction of $\wp$ to the periodic orbits of $\xi$, $\epsilon\xi$, and $\epsilon^2\xi$.
\een

\,

\item [(B)] If $\xi$ is an attracting (but not superattracting) periodic point of $\wp$, then

\ben

\item there are exactly three attracting (none of which is superattracting) periodic orbits of $\wp$ and these are equal to the periodic orbits of $\xi$, $\epsilon\xi$, and $\epsilon^2\xi$

\,

\item the periodic connected components of the Fatou set $F(\wp)$ are respectively the three basins of immediate attraction of $\wp$ to the periodic orbits of $\xi$, $\epsilon\xi$, and $\epsilon^2\xi$,

\,

\item all connected components of the Fatou set $F(\wp)$ are respectively the three basins of attraction of $\wp$ to the periodic orbits of $\xi$, $\epsilon\xi$, and $\epsilon^2\xi$.

\,

\item the sets of derivatives of $\wp$ for each of these three attracting periodic orbits of $\xi$, $\epsilon\xi$, and $\epsilon^2\xi$ are equal. 
\een

\,

\item [(C)] If $\xi$ is a rationally indifferent periodic point of $\wp$, then

\,

\ben 

\item there are exactly three rationally indifferent periodic cycles of $\wp$,

\,

\item the periodic connected components of the Fatou set $F(\wp)$ are respectively the three basins of immediate attraction of $\wp$ to the periodic orbits of $\xi$, $\epsilon\xi$, and $\epsilon^2\xi$,

\,

\item all connected components of the Fatou set $F(\wp)$ are respectively the three basins of attraction of $\wp$ to the periodic orbits of $\xi$, $\epsilon\xi$, and $\epsilon^2\xi$.

\,

\item the sets of derivatives of $\wp$ for each of these three attracting periodic orbits of $\xi$, $\epsilon\xi$, and $\epsilon^2\xi$ are equal.  
\een
\end{enumerate}
\ethm

\bpf 
Since $\wp$ is triangular, there is a triangular lattice $\La\sbt\C$ such that 
$$
\wp=\wp_\La.
$$
We first will show that if $f$ has one (super) attracting or rationally indifferent periodic point $\xi\in\C$ with some period $p\ge 1$, then $\epsilon\xi$ and $\epsilon^2\xi$ are also respective (super) attracting or rationally indifferent periodic points,
the derivative of $\wp_\La$ at each of these three periodic points is the same, and that if any two of the three periodic orbits of $\xi$, $\epsilon\xi$ or $\epsilon^2\xi$ intersect, then all these three periodic orbits coincide. 

So, let $\xi\in\C$ be such a periodic point with period $p\ge 1$. Because of \eqref{iterates} we get,
$$
\wp_\La^p(\epsilon \xi)
=\epsilon\wp_\La^p(\xi)
=\epsilon\xi
\  \  \  {\rm and} \  \  \
\wp_\La^p(\epsilon^2 \xi)
=\epsilon^2\wp_\La^p(\xi)
=\epsilon\xi
$$
meaning that $\epsilon\xi$ and $\epsilon^2\xi$ are also periodic points of $\wp_\La$ with period $p$. Differentiating now both equations in \eqref{iterates} and using the Chain Rule, we respectively get that
$$
\epsilon\(\wp_\La^p\)'(\xi)
=\(\wp_\La^p\circ\epsilon\)'(\xi)
=\epsilon\(\wp_\La^p\)'(\epsilon\xi).
$$
Hence,
$$
\(\wp_\La^p\)'(\epsilon\xi)=\(\wp_\La^p\)'(\xi).
$$
Likewise,
$$
\(\wp_\La^p\)'(\epsilon^2\xi)=\(\wp_\La^p\)'(\xi).
$$
Now, suppose that some two of the three the three periodic orbits of $\xi$, $\epsilon\xi$ or $\epsilon^2\xi$ intersect. We may assume without loss of generality that these are orbits of $\xi$ and $\epsilon\xi$. But any two periodic orbits that intersect, coincide. Thus, in particular, the orbits of 
$\xi$ and $\epsilon\xi$ are equal. But then there exists an integer $1\le p$  such that
$$
\wp^k(\xi)=\epsilon\xi.
$$
Having this and using \eqref{iterates}, we get 
$$
\wp^k(\epsilon\xi)
=\epsilon\wp^k(\xi)
=\epsilon^2\xi.
$$
Therefore, all the orbits of $\xi$, $\epsilon\xi$ or $\epsilon^2\xi$, coincide, and we have proved what we claimed. 

Having this and assuming that the orbits of $\xi$, $\epsilon\xi$ and $\epsilon^2\xi$ do not intersect, we are immediately done by virtue of Theorem~\ref{t120200316}. The proof is complete.
\epf

\sp The case in the above theorem when the three periodic orbits are mutually distinct is typical while the case when they coincide is exceptional. It is not hard to calculate that such exceptional phenomenon occurs for example for the triangular lattice with $g_3=5.5\epsilon$. This lattice is not real. It is illustrated on Figure~1 below. The only attracting, in fact superattracting, periodic orbit with period $3$ is marked in black. This lattice is a rotation of the real triangular lattice with $g_3=5.5$ having instead $3$ super attracting fixed points.

\begin{figure}[ht]\label{Figure 1}
\begin{center}
\includegraphics[scale=0.2]{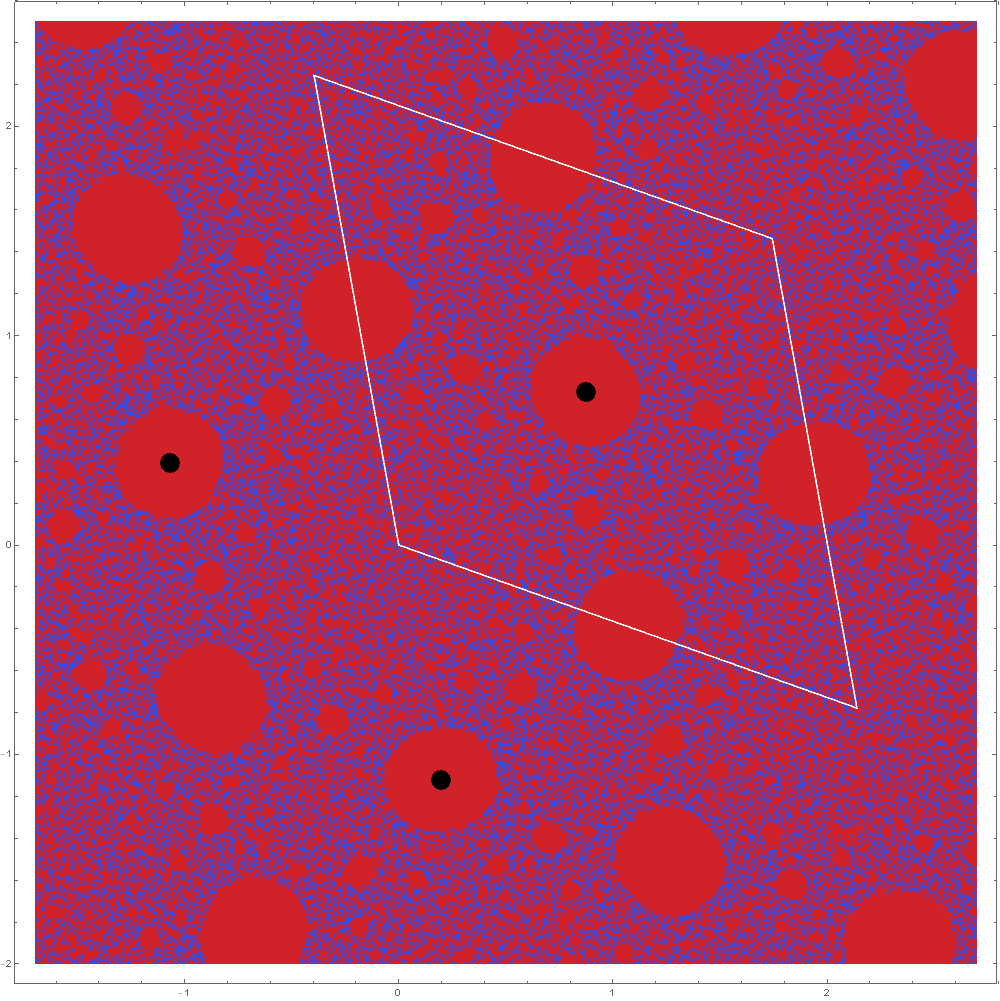}

\caption{$J(\wp_\La)$ for the triangular lattice $\La$ with $g_3(\La)=5.5\epsilon$. The only attracting, in fact superattracting, periodic orbit with period $3$ is marked in black.}
\end{center}
\end{figure}

\sp As an immediate consequence of the previous theorem we get the following.

\sp\bthm\label{t120200305} 
Let $\wp:\mathbb C\lra \oc$ be a triangular non--recurrent Weierstrass elliptic function. 

\begin{enumerate} 

\item If $\wp$ has a (super) attracting fixed point $\xi\in\C$, then 

\,

\ben

\item $\wp$ has exactly three (super)attracting fixed points $\xi$, $\epsilon\xi$, and $\epsilon^2\xi$, 

\,

\item the periodic connected components of the Fatou set $F(\wp)$ are the three basins of immediate attraction of $\wp$ to these fixed points $\xi$, $\epsilon\xi$, and $\epsilon^2\xi$, 

\,

\item all connected components of the Fatou set $F(\wp)$ are the three basins of attraction of $\wp$ to these fixed points $\xi$, $\epsilon\xi$, and $\epsilon^2\xi$. 

\,

\item In consequence, the elliptic function $\wp$ is expanding, so normal compactly non--recurrent of finite character. 
\een

\,

\item If $\wp$ has a rationally indifferent fixed point, then 

\,

\ben

\item $\wp$ has exactly three rationally indifferent fixed points $\xi$, $\epsilon\xi$, and $\epsilon^2\xi$, 

\,

\item the periodic connected components of the Fatou set $F(\wp)$ are the three basins of immediate attractions of $\wp$ to these fixed points $\xi$, $\epsilon\xi$, and $\epsilon^2\xi$,

\,

\item all connected components of the Fatou set $F(\wp)$ are the three basins of attraction of $\wp$ to these fixed points $\xi$, $\epsilon\xi$, and $\epsilon^2\xi$.

\,

\item In consequence, the elliptic function $\wp$ is parabolic.
\een
\end{enumerate}
\ethm

\sp We recall from Definition~\ref{d120190216} (0) that a lattice $\Lambda\sbt\C$ is called real if and only if
$$
\ov\La=\La.
$$
We need the following properties of triangular lattices and properties of  Weierstrass $\wp_\La$ functions generated by such lattices.

\sp\fr\bprop\label{cor4.6in{HK1}}
Let $\La\sbt\C$ be a triangular lattice.
\begin{itemize}
\item [(1)] Then the critical values of the Weierstrass $\wp_\La$ function
are the cubic roots of $g_3(\La)/4$. So,
$$ 
e_3\in\C\sms\{0\}, \  \
e_1=e^{4\pi i /3} e_3, \  \  {\rm and} \  \
e_2=e^{2\pi i /3} e_3.
$$
In particular, if $g_3(\La)=4$ (then $\La$ is real by Proposition~\ref{p120200306}) then the  critical values of  the Weierstrass $\wp_\Lambda$  function are the cubic roots of unity i.e.
$$
e_1=\epsilon^2, \  \  e_2=\epsilon \ \ {\rm and} \  \
e_3=1,
$$
where, we recall, $\epsilon:=e^{\frac{2\pi i}{3}}$.

\,

\item [(2)] If the lattice $\La\sbt\C$ is real, then the  postcritical set ${\rm PC}(\wp_\La)$ of the  Weierstrass  $\wp_\La$ function is contained in the union of the following three $\wp_\La$--forward invariant rays: 
$$
[e_3, +\infty]\sbt\R\cup\{\infty\}, \  \  \epsilon[e_3, +\infty], \  \  
{\rm and} \  \  \epsilon^2[e_3, +\infty].
$$
\end{itemize}
\eprop
\bpf Item (1) is a reformulation of item (6) of Proposition~\ref{p120200306}. 

\sp Proving item (2), it follows from Proposition~\ref{p120200306} (5) and (6c5) that

\beq\label{120190216}
\wp_\La([e_3,+\infty])\sbt [e_3,+\infty]\sbt\R.
\eeq
Therefore,
$$
O_+(e_3)=\{\wp_\La^n(e_3):n\ge 0\}\sbt [e_3, +\infty]
$$  
Since $ \epsilon \Lambda = \Lambda$, the homogeneity property (\ref{2}) implies that
$$
\wp_\La (\epsilon z)
=\frac{1}{\epsilon^2}\wp_\La(z)
=\epsilon \wp_\La(z)
$$  
for all $z\in\C$. Similarly (or consequently),
$$
\wp_\La (\epsilon^2 u)=\epsilon^2\wp_\La(u)
$$ 
for all $z\in\C$. These two above displayed formulas along with \eqref{120190216}, yield
\beq\label{220190216}
\wp_\La(\epsilon[e_3,+\infty])\sbt \epsilon[e_3,+\infty]
\  \  \  {\rm and} \  \  \
\wp_\La(\epsilon^2[e_3,+\infty])\sbt \epsilon^2[e_3,+\infty].
\eeq
Thus 
$$
\wp_\La(e_1)
=\wp_\La(\epsilon^2)
=\epsilon^2\wp_\La(1)\in \epsilon^2[e_3, +\infty]
$$
and 
$$
\wp_\La(e_2)
=\wp_\La(\epsilon)
=\epsilon\wp_\La(1)\in \epsilon[e_3, +\infty]
$$
Along with \eqref{120190216} and \eqref{220190216}, these imply that
$$
{\rm PC}(\wp_\La)
=O_+(e_1)\cup O_+(e_2)\cup O_+(e_3)
\sbt [e_3, +\infty]\cup \epsilon[e_3, +\infty]\cup \epsilon^2[e_3, +\infty].
$$
The proof is complete.
\qed


\sp Now we shall prove one more theorem in the form of Theorem~\ref{t620200303C} and Theorem~\ref{t120200305}. 

\sp\bthm\label{t120200323} 
If $\La$ is a (real) triangular lattice with $g_3(\La)>0$, then the Weierstrass elliptic $\wp_\La$ function has no Siegel disk or Herman ring. Moreover, if $\xi\in\C$ is either a superattracting, attracting, or rationally indifferent periodic point, then 

\ben \item $\epsilon\xi$ and $\epsilon^2\xi$ are also respective superattracting, attracting, or rationally indifferent periodic points,
 
\item
$$
\wp_\La'(\xi)=\wp_\La'(\epsilon\xi)=\wp_\La'(\epsilon^2\xi).
$$

\item Either some two of the three periodic orbits of $\xi$, $\epsilon\xi$ or $\epsilon^2\xi$ intersect, or all three of them are mutually disjoint. 
\een

Assume the latter case holds. Then:

\begin{enumerate}


\item [(A)] If $\xi$ is a superattracting periodic point of $\wp_\La$, then

\,

\ben
\item there are exactly three superattracting periodic orbits of $\wp_\La$ and these are equal to the periodic orbits of $\xi$, $\epsilon\xi$, and $\epsilon^2\xi$,

\item the periodic connected components of the Fatou set $F\(\wp_\La\)$ are respectively the three basins of immediate attraction of $\wp_\La$ to the periodic orbits of $\xi$, $\epsilon\xi$, and $\epsilon^2\xi$,

\,

\item all connected components of the Fatou set $F\(\wp_\La\)$ are respectively the three basins of attraction of $\wp$ to the periodic orbits of $\xi$, $\epsilon\xi$, and $\epsilon^2\xi$.
\een

\,

\item [(B)] If $\xi$ is an attracting (but not superattracting) periodic point of $\wp_\La$, then

\ben

\item there are exactly three attracting (none of which is superattracting) periodic orbits of $\wp_\La$ and these are equal to the periodic orbits of $\xi$, $\epsilon\xi$, and $\epsilon^2\xi$

\,

\item the periodic connected components of the Fatou set $F_\La(\wp_\La\)$ are respectively the three basins of immediate attraction of $\wp_\La$ to the periodic orbits of $\xi$, $\epsilon\xi$, and $\epsilon^2\xi$,

\,

\item all connected components of the Fatou set $F\(\wp_\La\)$ are respectively the three basins of attraction of $\wp$ to the periodic orbits of $\xi$, $\epsilon\xi$, and $\epsilon^2\xi$.

\,

\item the sets of derivatives of $\wp_\La$ for each of these three attracting periodic orbits of $\xi$, $\epsilon\xi$, and $\epsilon^2\xi$ are equal. 
\een

\,

\item [(C)] If $\xi$ is a rationally indifferent periodic point of $\wp_\La$, then

\,

\ben 

\item there are exactly three rationally indifferent periodic cycles of $\wp_\La$,

\,

\item the periodic connected components of the Fatou set $F\(\wp_\La\)$ are respectively the three basins of immediate attraction of $\wp_\La$ to the periodic orbits of $\xi$, $\epsilon\xi$, and $\epsilon^2\xi$,

\,

\item all connected components of the Fatou set $F\(\wp_\La\)$ are respectively the three basins of attraction of $\wp_\La$ to the periodic orbits of $\xi$, $\epsilon\xi$, and $\epsilon^2\xi$.

\,

\item the sets of derivatives of $\wp_\La$ for each of these three attracting periodic orbits of $\xi$, $\epsilon\xi$, and $\epsilon^2\xi$ are equal.  
\een
\end{enumerate}
\ethm

\bpf
The same arguments as those for Theorem~\ref{t120200303B} and Theorem~\ref{t120200305} work, so all what we need to   show is that $\wp_{\La}$ has no Siegel disks or Herman rings. We will do it now. So, seeking a contradiction suppose that $U$ is either a Siegel disk or Herman ring of the function $\wp_{\La}$. By Proposition~\ref{cor4.6in{HK1}} and Theorem~\ref{t1sh11},
\beq\lab{120200314}
\bd U\sbt[e_3, +\infty]\cup \epsilon[e_3, +\infty]\cup \epsilon^2[e_3, +\infty],
\eeq
where $e_3\in(0,+\infty)$. Since the point $0$ is a pole of $\wp_{\La}$, it belongs to the Julia set $J(\wp_{\La})$, whence it does not belong to $U$. But by \eqref{120200314} and since $e_3\in(0,+\infty)$, the origin $0$ does not belong to the boundary $\bd U$ either. Therefore,
\beq\lab{220200314}
0\in \C\sms (U\cup \bd U).
\eeq
On the other hand, Since $U$ is not empty and open in $\C$, it is not contained in $[e_3, +\infty]\cup \epsilon[e_3, +\infty]\cup \epsilon^2[e_3, +\infty]$. So, there exists a point 
\beq\lab{320200314}
\xi \in U\sms \([e_3, +\infty]\cup \epsilon[e_3, +\infty]\cup \epsilon^2[e_3, +\infty]\).
\eeq
Let $[0,\xi]$ be the closed line segment joining $0$ and $\xi$. Then, on the one hand,
\beq\lab{420200314}
[0,\xi]\sbt \C\sms \([e_3, +\infty]\cup \epsilon[e_3, +\infty]\cup \epsilon^2[e_3, +\infty]\)
\sbt \C\sms \bd U.
\eeq
On the the other hand, by \eqref{220200314} and \eqref{320200314}, we have
$$
[0,\xi]\cap U\ne\es
\  \  \  {\rm and}  \  \  \
[0,\xi]\cap (\C\sms U)\ne\es.
$$
Along with \eqref{420200314}, this produces contradiction since the segment $[0,\xi]$ is connected. We are done.
\epf

\sp As an immediate consequence of this theorem, we get the following.

\sp\bthm\label{t120200305B} 
Let $\La$ be a (real) triangular lattice with $g_3(\La)>0$. Then the following hold.

\begin{enumerate} 

\item If $\wp_\La$ has a (super) attracting fixed point $\xi\in\C$, then 

\,

\ben

\item $\wp_\La$ has exactly three (super)attracting fixed points $\xi$, $\epsilon\xi$, and $\epsilon^2\xi$, 

\,

\item the periodic connected components of the Fatou set $F\(\wp_\La\)$ are the three basins of immediate attraction of $\wp$ to these fixed points $\xi$, $\epsilon\xi$, and $\epsilon^2\xi$, 

\,

\item all connected components of the Fatou set $F(\wp)$ are the three basins of attraction of $\wp_\La$ to these fixed points $\xi$, $\epsilon\xi$, and $\epsilon^2\xi$. 

\,

\item In consequence, the elliptic function $\wp_\La$ is expanding, so normal compactly non--recurrent of finite character. 
\een

\,

\item If $\wp_\La$ has a rationally indifferent fixed point, then 

\,

\ben

\item $\wp_\La$ has exactly three rationally indifferent fixed points $\xi$, $\epsilon\xi$, and $\epsilon^2\xi$, 

\,

\item the periodic connected components of the Fatou set $F\(\wp_\La\)$ are the three basins of immediate attractions of $\wp_\La$ to these fixed points $\xi$, $\epsilon\xi$, and $\epsilon^2\xi$,

\,

\item all connected components of the Fatou set $F\(\wp_\La\)$ are the three basins of attraction of $\wp$ to these fixed points $\xi$, $\epsilon\xi$, and $\epsilon^2\xi$.

\,

\item In consequence, the elliptic function $\wp_\La$ is parabolic.
\een
\end{enumerate}
\ethm

\section{Simple Examples of Dynamically Different Elliptic Functions}\label{Simple Examples of Elliptic Functions}


Let  $\Lambda\sbt \mathbb C$ be a lattice. Let $\gamma  \in \mathbb C \sms \{0\}$. Finally, let  $ s  \in \mathbb C \sms \{0\}$. Then, by the homogeneity property (\ref{2}), for any $ z \in  \mathbb C \sms \{0\}$ we  have that
$$ 
\wp_\Lambda(sz)=sz\quad  \Longleftrightarrow  \quad \frac{\wp_\Lambda(z)}{z}=s^3,
$$
while  by the  homogeneity property (\ref{2}), we  have that
$$ 
\wp_\Lambda'(sz)=\gamma \quad  \Longleftrightarrow  \quad \gamma^{-1} \wp_\Lambda'(z)=s^3.
$$
Therefore, we got the following.

\blem\label{l1ex7} If $\Lambda \sbt \mathbb C$ is a lattice and $\gamma \in \mathbb C\sms\{0\}$, then there exists $s  \in \mathbb C \sms \{0\}$ and $ \xi
\in \mathbb C \sms (\Lambda \cup \wp_\Lambda ^{-1}(0))$, such that
$$ \wp_{s\Lambda}(s \xi )= s \xi \quad \text{and}   \quad \wp_{s\Lambda}'(z\xi )=\gamma$$
if and only if the equation
$$ 
\frac{z \wp_\Lambda'(z)}{\wp_\lambda(z)}=\gamma
$$
has solution in $ \mathbb C \sms (\Lambda \cup \wp_\Lambda^{-1}(0))$.
\elem

Let  $F_\Lambda: \mathbb  C \lra \hat{\mathbb C}$ be  the  meromorphic function defined  by the formula
$$
F_\Lambda(z):=\frac{z \wp_\lambda'(z)}{\wp_\Lambda (z)}.
$$ 
Our first result is this.

\blem\label{l2ex7} If $ \Lambda \sbt \mathbb C$ is a lattice and $t\in (0, +\infty)$, then there exists $G_t$, a non--empty open (in the relative topology) subset of $\{
z \in \mathbb C:|z|=t\}$ such that for every $w \in G_t$ the equation
\begin{equation}\label{1ex7}
\frac{z \wp_\Lambda'(z)}{\wp_\Lambda (z)}=w
\end{equation}
has a solution in $ \mathbb C \sms (\Lambda \cup \wp_\Lambda^{-1}(0))$.
\elem

\bpf Since the set $\mathbb C$ is open and connected and the function $F_\Lambda: \mathbb C \to \hat{\mathbb C}$ is meromorphic, we get the following

\

Claim~$1^0$:  The set $F_\Lambda (\mathbb C)$ is open and connected.

\sp

\fr  Obviously,

\sp

Claim~$2^0$:   $F_\Lambda(\Crit(\wp_\Lambda))=\{0\}$ and $F_\Lambda(\Lambda \sms \{0\})=\{\infty\}$.

\sp Since the set $\{z \in \mathbb C:|z|=t \}$ separates $0$  and $\infty \in \hat{\mathbb C}$, it follows  from Claim~$1^0$ and Claim~${2^0}$ that $\{z \in \mathbb C:|z|=t\}\cap F_\Lambda(\mathbb C)$ is a non--empty open set in the relative topology on $\{z \in \mathbb C:|z|=t\}$. Now, there are at least two arguments for the set
$$
\big\{z \in \mathbb C:|z|=t\big\}\cap F_\Lambda(\mathbb C \sms (\Lambda \cup \wp_\Lambda^{-1}(0))
$$
to be non--empty. The first one is to observe that the set  
$$
\big\{z \in \mathbb C:|z|=t\big\}\cap F_\Lambda(\mathbb C \sms (\Lambda \cup \wp_\Lambda^{-1}(0))
$$ 
is a countable $F_\sigma$ subset of 
$$
\big\{z \in \mathbb C:|z|=t\big\}\cap F_\Lambda(\mathbb C)
$$
and to apply the Baire  Category Theorem. The second one is to note that the set $\{z \in \mathbb C:|z|=t\}\cap F_\Lambda(\mathbb C)$ is uncountable. We are done. \epf

\sp As an immediate consequence of this lemma and Lemma~\ref{l1ex7}, we get the following:

\sp

\blem \label{l1ex8} If  $\Lambda \sbt \mathbb C$ is a lattice, $ t \in (0,+\infty)$, and $ w \in G_t$, the set produced in Lemma~\ref{l2ex7}, then there exists $s\in \mathbb C\sms \{0\}$ and
$\xi \in \mathbb C \sms (\Lambda \cup \wp_\Lambda^{-1}(0))$, such that
$$ \wp_\Lambda(s\xi )=s\xi   \quad \text{and}\quad {\wp_{\xi\Lambda}'(s \xi)}=w.$$
\elem

\sp Let's take the first this remarkable  fruits of this lemma and our theorems about triangular lattices. Indeed, as an immediate consequence of this lemma and Theorem~\ref{t620200303C}, we get the following four theorems.

\sp

\bthm\label{t2ex8} For every $t \in (0,1)$ there exists a triangular lattice $\Lambda \sbt \mathbb C$
such that
\begin{itemize}
\item [(1)] There exist three attracting fixed points $\xi_1, \xi_2,\xi_3\in \mathbb C$ of $\wp_\Lambda$  such that 
$$
 \wp_{\lambda}'(\xi_1)=\wp_{\lambda}'(\xi_2)=\wp_{\lambda}'(\xi_3)
 \  \  \  {\rm and} \  \  \ |\wp_\Lambda'(\xi_1)|=t.
$$

\item [(2)] The periodic Fatou components of $\wp_{\Lambda}$ consists  of the three  basins of immediate attraction to $\xi_1, \xi_2, \xi_3$.

\,

\item [(3)] The Fatou set $F(\wp_\Lambda)$ of $\wp_\Lambda$ is the union of basins of attraction to $\xi_1$, $\xi_2$  and $\xi_3$.

\,

\item [(4)] The Weierstrass elliptic function $\wp_\Lambda$ is expanding, thus normal  compactly non-recurrent of finite character.
\end{itemize}.
\ethm

\bthm\label{t3ex8} There exists a triangular lattice  $\Lambda \sbt \mathbb C$ such that
\begin{itemize}

\item[(1)] There exist three rationally  indifferent fixed points $\xi_1, \xi_2,\xi_3\in \mathbb C$ of $\wp_\Lambda$ such that
$$
\wp_{\Lambda}'(\xi_1)=\wp_{\Lambda}'(\xi_2)=\wp_{\Lambda}'(\xi_3).
$$

\item[(2)]  The periodic Fatou components of $\wp_\Lambda$ consist of the three basins of immediate attraction to
$\xi_1, \xi_2$  and $\xi_3$.

\,

\item[(3)]  The Fatou set $F(\wp_\Lambda)$ of $\wp_\Lambda$ is the union of basins of attraction to  $\xi_1, \xi_2$ and $\xi_3$.

\,

\item [(4)]  The Weierstrass elliptic  function $ \wp_\Lambda$ is parabolic.
\end{itemize}
\ethm

\bthm\label{t1ex9} 
For each lattice $\Lambda \sbt \mathbb C$ there exists a lattice $\Gamma\sbt \mathbb C$, similar to $\Lambda$, such that the Weierstrass elliptic function $\wp_\Ga$ has at  least  one Siegel disk whose center is a fixed point of $\wp_\Ga$.
\ethm

\sp

\bthm~\label{t3ex9} For each lattice $\Lambda \sbt \mathbb C$ there exists a lattice $\Gamma\sbt \mathbb C$, similar to $\Lambda$, such that the Weierstrass elliptic function $\wp_\Ga$ has at least one Cremer fixed point.\ethm

We would like  to have  a little  bit more than Theorem~\ref{t3ex8}, namely we would like to have a parabolic triangular Weierstrass elliptic function with specific, explicitly  known  value of derivative of its rationally indifferent fixed points. We shall prove the following.

\bthm\label{t2ex9}  For every (real) triangular lattice $\Lambda$ with $g_3(\Lambda)>0$ there  exists a real triangular lattice $\Gamma$, similar to $\Lambda$, such that the Weierstrass elliptic function $\wp_\Gamma$ is parabolic  and has three rationally  indifferent  fixed points $\xi_1, \xi_2$ and $\xi_3$ such that
$$
\wp_\Gamma'(\xi_1)=\wp_\Gamma'(\xi_2)=\wp_\Gamma'(\xi_3)=1.
$$
\ethm

\bpf Because of Proposition~\ref{p120200306} (6), 
$$
\Lambda =\lambda[\epsilon, \epsilon^{-1}]
$$
with some $\lambda \in(0, +\infty)$. Fix arbitrary $k \in \mathbb N$. By Proposition~\ref{p120200306} (8), we  have
 \begin{equation}\label{1ex9}
F_\Lambda\left(\frac{\lambda}{2}\right)=0.
\end{equation}
Since $\wp_\Lambda$ has a pole of order $2$ at $\lambda$, we have that
\begin{equation}\label{2ex9}
\wp_\Lambda(z)=A(z-\lambda)^{-2}+ O(|z-\lambda|^{-1})
\end{equation}
and
$$ \wp_\Lambda'(z)=-2A(z-\lambda)^{-3}+ O(|z-\lambda|^{-2}),$$
with same $A\in \mathbb R$. Hence,
$$
\begin{aligned}
\lim_{\mathbb R \ni x \upto \lambda} F_\Lambda(x)
& = \lambda \lim_{\mathbb R \ni x \upto \lambda}\frac{\wp_\Lambda'(x)}{\wp_\Lambda (x)}
=\lambda \lim_{\mathbb R \ni x \upto \lambda}\frac{-2A(x-\lambda)^{-3} + O(|x-\lambda|^{-2})}{A(x-\lambda)^{-2}+ O(|z-\lambda|^{-1})}\\
& =\lambda\lim_{\mathbb R \ni x \upto \lambda}\frac{-2A(x-\lambda)^{-1}+O(1)}{A}\\
&=O(1)-2 \lambda \lim_{\mathbb R \ni x \upto \lambda}(x -\lambda)^{-1}=+\infty.
\end{aligned}
$$
It follows from  this, formula (\ref{1ex9}), continuity of the function
$$
F_{\Lambda}\, |_{[\lambda/2, \lambda)}:[\lambda/2, \lambda) \lra \mathbb R,
$$
and the Intermediate Value Theorem, that
$$
F_\Lambda((\lambda/2, \lambda))\supset (0, +\infty).
$$
Hence, there exists $x  \in (\lambda/2, \lambda)$ such that $ F_\Lambda(x)=1$. Therefore, a direct application of Lemma~\ref{l1ex7} completes the proof.
\epf


\

Let $\Lambda=\lambda[1,i]$ be a real square  lattice with $\lambda >0$. For every $\alpha \in \mathbb C$, let 
$$
h_\alpha:=\wp_\Lambda+\alpha.
$$
We shall prove the following.

\bthm\label{t1ex10} 
If $\Lambda=2 \lambda[1,i]$ is a real  square  lattice   with some $ \lambda \in (0,+\infty)$, then for every $ \e\in (0,1)$ there exists $M_\e>0$ such that if $(m_k)_{k=1}^\infty$ is any sequence of positive integers  such that $m_k\geq \frac{M_\e}{2\lambda}$, then there exists $ \beta \in (-\epsilon,\epsilon)$ such that
$$ 
h_\beta^k(\lambda)\in [2 \lambda m_k, 2 \lambda(m_k+1)]
$$
for every integer $k \geq 2$.
\ethm

\bpf For  every integer $k \geq 1$ define the function $G_k:\mathbb C \lra \hat{\mathbb C}$ by the formula 
$$
G_k:=h^k_\alpha(\lambda+\lambda i)
$$ 
with the convention that $h_\alpha(\infty)=\infty$. We then have 
$$
G_{k+1}(\alpha)=\wp_\Lambda \circ G_k(\alpha)+\alpha
$$
and
$$
G_1(\alpha)=\wp_\Lambda(\lambda+\lambda i)+\alpha=\alpha, \quad G_2(\alpha)=\wp_\Lambda(\alpha)+\alpha, \quad \text{and}\quad  G_2(0)=\wp_1(0)=\infty.
$$
So, $G_2$ is meromorphic (and takes on the value $\infty$) on some open ball $B_{\mathbb C}(0,2 \delta)$ with $ \delta \in (0,\xi)$. In addition, since $\wp_\Lambda$ is real,
$$
G_2((-2\delta, 2 \delta))\sbt \mathbb R,
$$
and 
\begin{equation}\label{1ex11}
G_2((-\delta, \delta))\supset [M_\e, +\infty) \sbt \mathbb R
\end{equation}
with some $M_\e>0$. We shall now define inductively a descending sequence $\{[a_k, b_k]\}_{k=1}^\infty$ of compact intervals in $\mathbb R$ such that for every $k \geq 2$:
\begin{itemize}
\item [(1)]  $[a_k, b_k] \sbt [-\delta, \delta] \sbt (-\e, \e)$
\end{itemize}
and
\begin{itemize}
\item[(2)] 
$$ 
G_k(a_k)=2 \lambda m_k, \  \  G_k(b_k)=2\lambda(m_k+1)
$$ 
and
$$
G_k([a_k, b_k])=[2 \lambda m_k, 2 \lambda(m_k+1)]\sbt \mathbb R.
$$
\end{itemize}
By our hypothesis, 
$$ 
[2 \lambda m_2, 2 \lambda(m_2+1)]\sbt [M_\epsilon, +\infty).
$$
So, it follows from (\ref{1ex11}) and continuity of $G_2$, that there exists closed interval $[a_2,b_2]\sbt [-\delta, \delta]$ such that 
$$
G_2(a_2)=2\lambda m_2, \  \  G_2(b_2)=2\lambda (m_2+1)
$$ 
and 
$$ 
G_2([a_2,b_2])=[2 \lambda m_2, 2 \lambda(m_2+1)].
$$
We are thus done with the base of induction. For the inductive step,
suppose that (1) and (2) hold for some integer $k \geq 2$. Then by (2) there exists $c_k\in[a_k, b_k]$ such that $G_k(c_k)=2 \lambda m_k +\lambda$. Hence,
$$
\begin{aligned} G_{k+1}(c_k)
&=\wp_\Lambda\circ G_k(c_k)+c_k=\wp_\Lambda(2 \lambda m_k+\lambda) +c_k\\ &=\wp_\Lambda (\lambda)  +c_k=e_1(\Lambda)+c_k
\leq e_1(\Lambda)+\e \\
& \leq e_1(\Lambda)+1.
\end{aligned}
$$
Also, by (2),
$$  
G_{k+1}(a_k)=\wp_\Lambda\circ G_k(a_k)+a_k=\wp_\Lambda(2 \lambda m_k) +a_k=\infty.
$$
Thus, using also Proposition~\ref{p120200310} (6c), we get
$$
G_{k+1}([a_k, b_k])\supset [e_1(\Lambda)+1, +\infty).
$$
Therefore, there exists an interval $[a_{k+1}, b_{k+1}]\sbt [a_k, b_k]$ such that
$$
G_{k+1}(a_{k+1})=2 \lambda m_{k+1}, \  \  
G_{k+1}(b_{k+1})=2\lambda(m_{k+1}+1),
$$ 
and 
$$
G_{k+1}([a_{k+1}, b_{k+1}])=[2 \lambda m_{k+1}, 2 \lambda(m_{k+1}+1)].
$$
So, our inductive construction is complete. Since $\{[a_k,b_k]\}_{k=1}^\infty$ is a descending sequence of compact intervals contained in $[-\delta, \delta]\sbt (-\e, \e)$, we get that
$$ 
\es\neq \Gamma :=\bigcap_{k=1}^\infty [a_k, b_k]\sbt [-\delta, \delta]\sbt (-\e, \e).
$$
Therefore, for any $ \beta \in \Gamma$, we have for all $k \geq 2$ that
$$
h^k_\beta(\lambda) =G_k (\beta) \in G_k([a_k, b_k])=[2 \lambda m_k, 2 \lambda ( m_k+1)].
$$
The proof is complete.
\epf

\sp We can  now prove one of the main results of this chapter.

\sp\bthm\label{t1ex12}
If $\Lambda=2 \lambda [1, i]$ is a real  square lattice with some $\lambda \in (0,+\infty)$ which has an attracting periodic point, then for every $\e\in (0, 1)$ sufficiently small there exists $\beta \in (-\e, \e)$ such that the elliptic function 
$$
h_\beta=\wp_\Lambda+\beta: \mathbb C \to \hat{\mathbb C}
$$  
has the following properties:
\begin{itemize}
\item [(1)]  $h_\beta$ has  exactly one attracting periodic orbit.

\,

\item [(2)]  $\beta$  is a critical  value of  $h_\beta$ and its orbit $\{h_\beta(\beta):\,  n\geq 0\}$ is bounded and infinite.

\,

\item [(3)] The Julia set $J(h_\beta)$ is a proper nowhere dense subset of $\mathbb C$.

\,

\item[(4)] $h_\beta$ is a normal subexpanding (in particular compactly non-recurrent) elliptic function of finite character, and it is non--expanding.
\end{itemize}
\ethm

\bpf By Theorem~\ref{t220200330} and item (6) of Proposition~\ref{p120200310}, $\wp_\Lambda$ has exactly one periodic orbit and both critical values  $e_1(\Lambda)$ and  $e_2(\Lambda)$  belong to its basin of attraction. Since attracting periodic orbits are stable under locally  smart perturbation, if $\e\in (0,1)$ is small enough, then  each function
$h_\alpha$, $ \alpha \in ( -\e, \)$, has an attracting periodic orbit and both its critical values   $e_1(\Lambda)+\alpha$ and
$e_2(\Lambda)+\alpha$  belong to its basin of attraction. So, (3) and a part of (1) are proved. Now we take $(m_k)_{k=1}^\infty$ a bounded sequence of positive integers  which is not eventually periodic with  $m_2\geq \frac{M_\e}{2\lambda}$. Let then  $\beta \in (-\e, \e)$  be the corresponding  number  produced in Theorem~\ref{t1ex10}.  Then (2) is automatically satisfied. Item (1) then follows since $h_\beta$ has exactly three critical values. Item (4) follows as well. The proof is complete.
\epf

\sp In order to complete the picture, we will now describe a large, uncountable, class of real square lattices $\Lambda$ satisfying  the hypothesis of the above  Theorem~\ref{t1ex12}.

\sp

\bthm\label{t2ex10} Let $\Lambda= 2 \lambda [1,i]$ be a real square lattice with some $ \lambda \in (0, +\infty)$. Let $m$ be any odd integer and let
$k:=\left(\frac{e_1}{m \lambda}\right)^{\frac{1}{3}}$, the real cubic root of $\frac{e_1}{m \lambda}$. 

If $\Gamma:=k\Lambda$ (a real square lattice), then the corresponding Weierstrass elliptic function $\wp_\Gamma: \mathbb C \lra \hat{\mathbb C}$ has a superattracting  periodic  point $mk\lambda$.

In consequence, by Theorem~\ref{t1ex10}, there exists $\beta \in\mathbb R \sms \{0\}$ such that the elliptic function $\wp_\Gamma+\beta: \mathbb C \lra \hat{\mathbb C}$ has all the properties (1)--(4) of Theorem~\ref{t1ex10}.
\ethm

\bpf By the homogeneity equation \eqref{2a}, we have that 
$$
\wp_\Gamma'(mk \lambda)=
\wp_\Gamma'( k \lambda)=k^{-3}\wp_\Gamma'(\lambda)=k^{-3}\cdot 0=0.
$$
By the homogeneity equation \eqref{2}, we have that
$$\wp_\Gamma(mk \lambda)=
\wp_\Gamma( k \lambda)=k^{-2}\wp_\Gamma(\lambda)=k^{-2}e_1 = k^{-2}k^3m\lambda=mk\lambda.$$
So, $mk\lambda$ is indeed a superattracting periodic point of $ \wp_\Lambda$ and the proof is complete.
\epf

\sp This theorem/example stemmed form Lemma~7.2 in \cite{HK2}. Theorem~9.3 (with $n=0$) and Theorem~9.4 in this paper are other sources of examples in Theorem~\ref{t1ex10}.

\section[Expanding Weierstrass Functions--Nowhere Dense Connected Julia Sets]
{Expanding (thus Compactly Non--Recurrent) \\ Triangular Weierstrass Elliptic Functions \\ with \\ Nowhere Dense Connected Julia Sets} \label{connectedjulia} 

In this section we shall provide several examples of compactly non--recurrent elliptic with various behavior of critical points (and critical values) functions, all of them having connected Julia sets. All these examples are motivated by the work of J. Hawkins and her collaborators, primarily by \cite{HK1}. The first following example was published therein.

\sp\fr\bthm\label{th8.4in{H K1}}  
Let $\Om=[w_1, w_2]$, where $w_1=\epsilon w_2$ with $\epsilon=e^{2\pi i/3}$, be the triangular lattice associated with the invariants $g_2=0$ and $g_3=4$; see Proposition~\ref{p120200306} for its existence and uniqueness. In particular, by Proposition~\ref{p120200306} (6), $\Om$ is real. Let $m$ be a negative odd integer.

If
$$
\g_1:=\lt({\frac{2 w_1^2\epsilon^2}{m}}\rt)^{\frac{1}{3}},  \  \  \
\g_2:=\g_1\frac{w_2}{w_1},
\  \  \  {\rm and} \  \  \
\Lambda:=[\g_1, \g_2],
$$
then 

\begin{enumerate}

\item The corresponding Weierstrass function $\wp_{\Lambda}:\C\lra\oc$ has exactly three superattracting fixed points.

\, 

\item The periodic connected components of the Fatou set $F(\wp)$ are the three basins of immediate attraction of $\wp$ to these fixed points, and all connected components of the Fatou set $F(\wp)$ are the three basins of attraction of $\wp$ to these fixed points. 

\, 

\item The Weierstrass function $\wp_{\Lambda}:\C\lra\oc$ is expanding, thus normal compactly non--recurrent of finite character.

\, 

\item The Julia set $J(\wp_\Lambda)$ is a proper nowhere dense subset of $\C$.

\, 

\item The Julia set $J(\wp_\Gamma)$ is connected. 
\end{enumerate}
\ethm

\fr {\sl Proof.} First, we will verify that the lattice $\Lambda$ is triangular. Indeed, since $w_1=\epsilon w_2$  and $\g_1=\g_2\frac{w_1}{w_2}$, we have that $\g_1=\epsilon \g_2$. Thus 
\beq\lab{120200131}
\epsilon \La = \La,
\eeq
meaning that $\La$ is triangular. It then follows from Theorem~\ref{t120200129}, formula \eqref{20200130} and Proposition~\ref{cor4.6in{HK1}} (1) that
$$
\wp_\Om(w_1/2)=e_1=\epsilon^2.
$$
Then the homogeneity equation~(\ref{2}) gives that
\beq\label{fixed}
 \begin{aligned}
\wp_\La\left(\frac{\g_1}{2}\right)
=\wp_{\frac{\g_1}{w_1}\Om}\lt(\frac{\g_1}{w_1}\cdot\frac{w_1}{2}\rt)
= \frac{1}{\left(\frac{\g_1}{w_1}\right)^2}\wp_\Om\left(\frac{w_1}{2}\right)
=\frac{1}{\left(\frac{\g_1}{w_1}\right)^2}\epsilon^2
=\frac{w_1^2\epsilon^2}{\left(\frac{2w_1^2\epsilon^2}{m}\right)^{\frac{2}{3}}}
= \frac{m \g_1}{2}.
\end{aligned}
\eeq
Since $m$ is an odd integer, we have that $\wp_\La(\frac{m\g_1}{2})=\wp_\La(\frac{\g_1}{2})$. Thus,
\beq\label{fixedB}
\wp_\La\lt(\frac{m\g_1}{2}\rt)=\frac{m\g_1}{2}
\eeq
meaning that $\frac{m\g_1}{2}$ is a fixed point of $\wp_\La$. Also,
\beq\lab{320200131}
\wp_\La'\lt(\frac{m\g_1}{2}\rt)=\wp_\La'\lt(\frac{\g_1}{2}\rt)= 0,
\eeq
whence $\frac{m\g_1}{2}$ is a superattracting fixed point of $\wp_\La$. Therefore, the items (1), (2), and (3) directly follow from Theorem~\ref{t120200305B}. Item (5) follows now directly from Theorem~\ref{th3.1in{HK3}}. Item (4) follows from Theorem~\ref{duality} since $J(\wp)\ne\C$.
\qed

\sp Figure~2  shows  the Julia  set $J(\wp_\La)$  of  the triangular Weierstrass $\wp_\La$  function defined in
Theorem~\ref{th8.4in{H K1}} for the case when $m=-1$. This lattice $\Lambda$  is generated by
 The corresponding  Weierstrass $\wp_\La$ function is expanding and has
three superattracting fixed points at 1.1382, -1.1382$\epsilon$, and
1.1382$\epsilon^2$.  A fundamental region is shown in white (so we see it is made of 2 equilateral triangles), the 3 distinct fixed points are shown in black (so we see their symmetry around the origin), and each attracting basin is a different color.  The origin is marked, but it is tiny in this format. It is the left vertex of the region.

\begin{figure}[ht]\label{Figure 2}
\begin{center}
\includegraphics[scale=0.2]{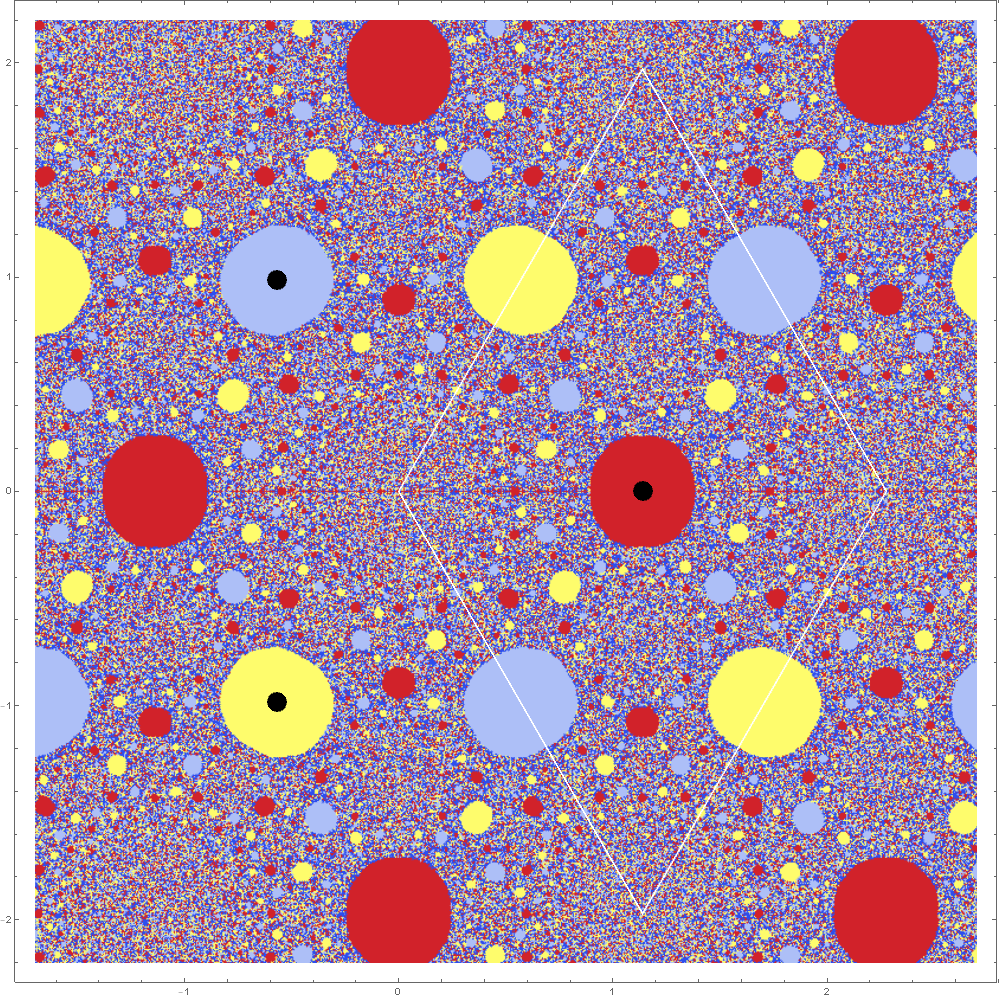}

\caption{$J(\wp_\La)$ where $\La=[\l_1,\l_2]$ with  $\l_1\approx1.1382+1.9741 i$,  $\l_2\approx 1.1382-1.9714i$. $\wp_\La$ is triangular, expanding, and has three superattracting fixed points.}
\end{center}
\end{figure}

\sp The following theorem has a proof analogous to the proof of Theorem~\ref{th8.4in{H K1}} and its proof will be omitted.

\sp\bthm[\cite{HK2}, Theorem~8.3]\label{th8.3in{HK2}} 
Let $\Omega=[\lambda, \lambda e^{\frac{2\pi i}{3}}]$, $\lambda>0$, be the  triangular lattice associated with the invariants $g_2=0$, $g_3=4$.
For any $m, n \in \mathbb Z$, if 
$$
k=\left((\l/2)+m \l + n
\l e^{\frac{2\pi i}{3}}\right)^{-1/3},
$$
then for 
$$
\Lambda:=k \Omega,
$$
the invariant  
$$
g_3(\Lambda)=4(\l/2 +m \l +n \l e^{\frac{2\pi
i}{3}})^2.
$$
Furthermore, 

\begin{enumerate}

\item The corresponding Weierstrass function $\wp_{\Lambda}:\C\lra\oc$ has exactly three superattracting fixed points.

\, 

\item The periodic connected components of the Fatou set $F(\wp)$ are the three basins of immediate attraction of $\wp$ to these fixed points, and all connected components of the Fatou set $F(\wp)$ are the three basins of attraction of $\wp$ to these fixed points. 

\, 

\item The Weierstrass function $\wp_{\Lambda}:\C\lra\oc$ is expanding, thus normal compactly non--recurrent of finite character.

\, 

\item The Julia set $J(\wp_\Lambda)$ is a proper nowhere dense subset of $\C$.

\, 

\item The Julia set $J(\wp_\Gamma)$ is connected.
\end{enumerate}
\ethm

\sp Figure~3 illustrates   the  Julia set of the triangular expanding  Weierstrass
$\wp_{\Lambda}$ function associated to the triangular lattice $\Lambda$  with the invariant $g_3(\La)\approx 5.67 + 2.08i$. In this case, there are three attracting cycles of period 3. One of these cycles occurs at approximately
$$
\xi=\big\{1.139+0.134i, 0.989+0.11i, 1.131-0.068i\big\}
$$
while the second cycle is located at $e^{\frac{2\pi i}{3}}p$ and the third one at $e^{\frac{4\pi i}{3}}p$.  Both
pictures are the same, with the fundamental period outlined.  Each
attracting basin is colored differently and the period 3 orbits are
marked with points.

\begin{figure}[ht]\label{Figure3}
\begin{center}
\includegraphics[scale=0.35]{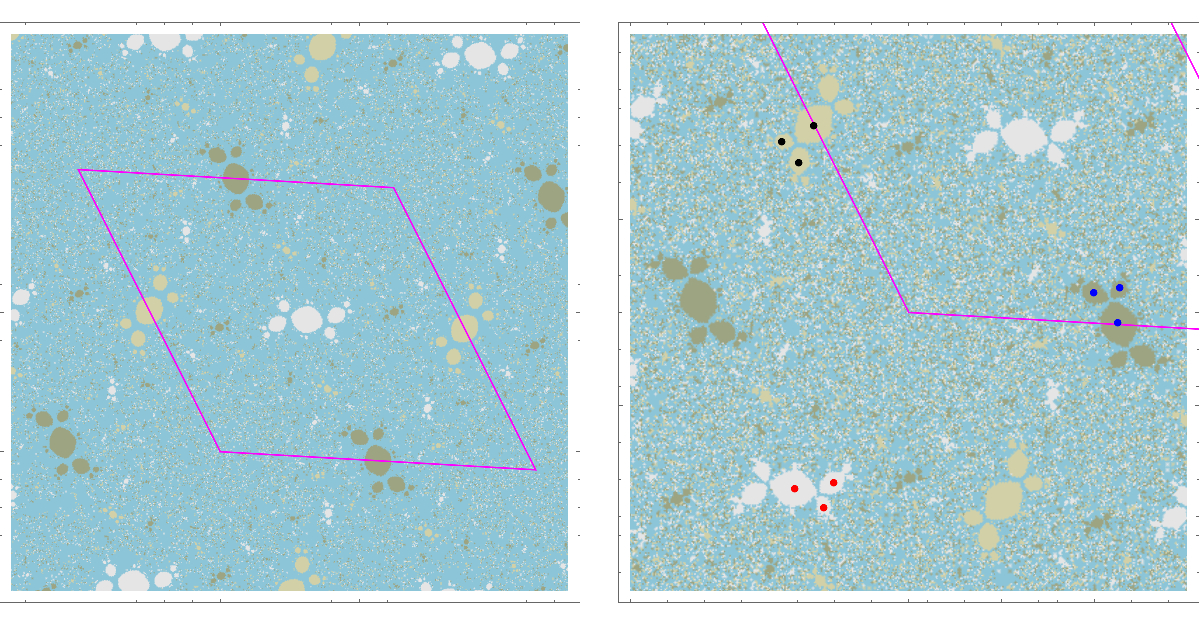}

\caption{$J(\wp_\La)$ where $g_3(\La)\approx 5.67 + 2.08i$. $\wp_\La$ is triangular, expanding, and has three attracting cycles of period $3$.
}
\end{center}
\end{figure}

\sp

\section[Weierstrass Elliptic Functions whose Critical Values are Preperiodic] {Triangular Weierstrass Elliptic Functions \\ whose Critical Values are Preperiodic; \\ Thus being Subexpanding}\label{preperiodic}

The main result of this section is the following theorem proved in \cite{HK1} as Theorem~8.6.

\sp\fr\bthm\label{th8.6in{HK1}} There exists a triangular  lattice
$\Lambda$ such that the critical values of the Weierstrass $\wp_\La$ function are preperiodic. More precisely,
$\wp_{\Lambda}(e_1)$,  $\wp_{\Lambda}(e_2)$  and $\wp_{\Lambda}(e_3)$ are
repelling fixed points of $\wp_\La$. 

In consequence, the elliptic function $\wp_{\Gamma}:\C\lra\oc$ is normal subexpanding of finite character (in particular, compactly non--recurrent),
and  non--expanding. Also $J(\wp_{\Gamma})=\oc$. 
\ethm

\bpf Denote again the real triangular lattice with $g_3=4$ (see Proposition~\ref{p120200306} for its existence and uniqueness) by $\Om$.
By Proposition~\ref{p120200306},
$$
\Om =[w_1,w_2],
$$
where
$$
w_1=t e^{\frac{\pi i}{3}}
\  \  \  {\rm and }  \  \  \  
w_2=t e^{-\frac{\pi i}{3}}
$$ 
with some $t>0$. Then $w_3=w_1+w_2$ is a real period of the Weierstrass $\wp_\Om$ function. Using the tables in \cite{MT} we have $w_3\approx 2.42 \ldots$, and so 
$$
2<w_3<3.
$$
Given $k\in\C\sms\{0\}$, let 
$$
\La:=k\Om=[kw_1, kw_2]=[\l_1, \l_2],
$$
As usually, denote 
$$
\l_3:=\l_1+\l_2.
$$
Then, by Theorem~\ref{t320200303}, $\l_3/2$ is a critical point of $\wp_\La$, i.e.
$$
\wp_\La'(\l_3/2)=0.
$$
The homogeneity equation~(\ref{2}) and the property  $\wp_\Om\left(\frac{w_1+w_2}{2}\right)=1$ (see Proposition~\ref{cor4.6in{HK1}} (1)) imply  that
$$
e_{3, \La}=\wp_\La\left(\frac{\l_3}{2}\right)=
\frac{1}{k^2}\wp_\Om\left(\frac{w_1+w_2}{2}\right)=\frac{1}{k^2}. 
$$
We are looking for a value of $k\in(0,+\infty)$ such that 
$$
\wp_\La (e_{3, \La})=\wp_\La(k^{-2})
$$ 
is a fixed point of the function $\wp_\La$ .
First, we will show that if there exists $k\in\R\sms\{0\}$ such that
\beq\label{war}
\wp_\La\left(-\frac{1}{k^2}+2\l_3\right)=-\frac{1}{k^2}+2 \l_3,
\eeq
then
\beq\label{fix}
\wp_\La\left(\wp_\La\left(\frac{1}{k^2}\right)\right) = \wp_\La\left(\frac{1}{k^2}\right).
\eeq
Indeed, using $\La$--periodicity and symmetry (i.e. eveness) of the function $\wp_\La$, we get from \eqref{war} that
\beq\label{left}
\wp_\La\left( \wp_\La\left(-\frac{1}{k^2}+2\l_3\right)\right)=\wp_\La\left( \wp_\La\left(-\frac{1}{k^2}\right)\right)=\wp_\La\left( \wp_\La\left(\frac{1}{k^2}\right)\right)
\eeq
and
\beq\label{right}
\wp_\La\left(-\frac{1}{k^2}+2 \l_3\right)
=\wp_\La\left(-\frac{1}{k^2}\right) 
=\wp_\La\left(\frac{1}{k^2}\right).
\eeq
So, (\ref{fix}) follows. Applying the homogeneity equation~(\ref{2}), we  get
$$
\wp_\La\left(-\frac{1}{k^2}+2\l_3\right)=\wp_{k \Om}\left(k
\left(-\frac{1}{k^3}+2w_3\right) \right)=
\frac{1}{k^2}\wp_\Om\left(-\frac{1}{k^3}\right),
$$ 
so we can  rewrite (\ref{war}) as
$$
\frac{1}{k^2}\wp_\Om\left(-\frac{1}{k^3}\right)
=-\frac{1}{k^2}+ 2\l_3
=-\frac{1}{k^2}+2k w_3,
$$
or equivalently as
\beq\label{war2}
\wp_\Om\left(-\frac{1}{k^3}\right)=-1+2k^3w_3.
\eeq
In order to find  $k \in(0,+\infty)$ satisfying (\ref{war2}), we consider two auxiliary real--valued functions 
$$
f,g:[r,s]:=\left[\left(\frac{1}{w_3}\right)^{\frac{1}{3}}, \left(\frac{2}{w_3}\right)^{\frac{1}{3}}\right] \lra \R,
$$  
defined as
$$
g(k):=\wp_\Om(-1/k^3)
\  \  \  {\rm and}  \  \  \
f(k):=-1+2k^3w_3.
$$  
Since $2<w_3<3$, we have that $r>1/2$ and $s<1$. Note that $g(r)=+\infty$, $g'(s)=0$, and $g$ is monotone decreasing on the interval $[r,s]$. Hence, there exists $\d>0$ such that 
$$
r+\d<s<1
\  \  \  {\rm and}  \  \  \ 
g(r+\d)>10.
$$
Thus $f(r+\d)<-1+ 2 w_3<-1+6=5$ and  
$$
g(r+\d)-f(r+\d)>0.
$$
By Proposition~\ref{cor4.6in{HK1}} (1), we get that
$$
g(s)=\wp_\Om(-w_3/2)=\wp_\Om(w_3/2).
$$
Hence,
$$
g(s)-f(s)<0.
$$
Since both $g$ and $f$ are continuous functions on the interval $[r+\delta,s]$, the Intermediate Value Theorem implies that there exists a number $k \in(r+\d,s)$ such that 
$$
g(k)=f(k).
$$ 
Hence, the condition (\ref{war2}) is satisfied, and, in consequence, (\ref{war}) holds. So, \eqref{fix} too. Thus, we have a triangular
lattice  $\La$  such that $e_{3, \La}=k^{-2}$ is preperiodic. Therefore, the  homogeneity properties  (\ref{iterates})  and  (\ref{2a}) imply that 
$$
\wp_\La^n(\epsilon z)
=\epsilon \wp_\La^n(z)\quad \text{and}\quad
\wp_\La^n(\epsilon^2 z)
=\epsilon^2 \wp_\La^n(z)
$$
for all integers $n\ge 0$, and 
$$
\wp_\La'(\epsilon z)=\epsilon \wp_\La'(z)\quad \text{and}\quad \wp_\La'(\epsilon ^2 z)=\epsilon^2 \wp_\La'(z). 
$$ 
Therefore, using Proposition~\ref{cor4.6in{HK1}} (1), we conclude that the critical values $e_1(\La)$ and $e_2(\La)$ are also preperiodic. Since  $\l_3/2\in (-\infty,0)$ and $k^{-2}\in(0,+\infty)$, we get that
\beq\lab{420200307}
\l_3/2\ne e_{3,\La}.
\eeq
Seeking contradiction, suppose that 
$$
\wp_\La\(e_{3,\La}\)=\l_3/2.
$$
Then \eqref{fix} would imply that
$$
e_{3,\La}=\wp_\La(\l_3/2)=\l_3/2,
$$
contrary to \eqref{420200307}. Thus, the critical point $\l_3/2$ is strictly preperiodic. Then, by virtue of Theorem~\ref{t120200305}, Theorem~\ref{t1ms123}, and Theorem~\ref{t2ms135}, we deduce that the periodic point 
$\wp_\La\(e_{3,\La}\)$, as well the points $\wp_\La\(e_{1,\La}\)$ and $\wp_\La\(e_{2,\La}\)$, are all repelling. This means that the first assertion of our theorem is proved. The second one is then an immediate consequence of the definitions elliptic functions being expanding, subexpandig, normal, compactly non--recurrent, and of finite character. The third one follows from it and Theorem~\ref{t120200303}.
\epf

\sp \section[Weierstrass Functions whose Critical Values are Poles or Prepoles] {Weierstrass Elliptic Functions whose Critical Values are \\ Poles or Prepoles \\ Thus being Subexpanding, thus Compactly Non--Recurrent}

\sp We shall provide in this section examples of elliptic functions whose critical values are prepoles of arbitrarily high orders. These come from \cite{HKK}.

\sp\fr\bthm\label{th8.9in{HK2}} 
There  exists a triangular lattice $\Gamma\sbt\C$  such that the critical values of the Weierstrass function $\wp_\Gamma:\C\lra\oc$ are the lattice points (i.e. belong to $\Ga$), thus poles. 

In consequence, the Weierstrass elliptic function $\wp_{\Gamma}:\C\lra\oc$ is subexpanding (thus compactly non--recurrent) of finite character, non--expanding, and not normal. In addition, $J(\wp_{\Gamma})=\oc$.  
\ethm

\bpf Let $\Lambda=[\lambda, \lambda \epsilon$, $\lambda >0$, be the real triangular lattice associated with the invariants $g_2=0$
and $g_3=4$; see Proposition~\ref{p120200306} for its existence and uniqueness. Consider the real triangular lattice  
$$
\Gamma:=k\Lambda:=[k\l, k\l\epsilon] 
$$
where  
$$
k:=(m\l+n \l \epsilon)^{-1/3}
$$ 
and $m, m$ are some fixed non--zero integers. Denote:
$$
\g:=k\l.
$$
The homogeneity equation~(\ref{2}) and Proposition~\ref{cor4.6in{HK1}} (1) give
\beq\label{prepoles}
\wp_\Ga \left(\frac{\g}{2}\right)
=\wp_{k\Lambda}\left(\frac{k\lambda}{2}\right)
=\frac{1}{k^2}\wp_\Lambda\left(\frac{\lambda}{2}\right)
=\frac{1}{k^2}
=(m \gamma+ n \gamma  \epsilon),
\eeq
where the last equality in (\ref{prepoles})  follows  from
$$\begin{aligned}
\frac{1}{k^2}&=\(m \lambda + n\lambda  \epsilon\)^{\frac{2}{3}}=
\(m \lambda + n\lambda  \epsilon\)\(m \lambda + n\lambda \epsilon\)^{-\frac{1}{3}}\\
&=k \(m \lambda + n\lambda  \epsilon\)=m \gamma+ n \gamma  \epsilon.
\end{aligned}$$
Thus $e_{1, \Gamma}=\wp_\Gamma (\gamma/2)= m  \gamma+ n \gamma  \epsilon$  is  a  lattice point and so, a pole of $\Gamma$.
Since $\Gamma$ is triangular, by applying (\ref{iterates}), we get that
$$
\wp_\La\lt(\epsilon \frac{\g}{2}\rt)
= \epsilon \wp_\La\lt(\frac{\g}{2}\rt)
= \epsilon \( m \gamma+ n \gamma  \epsilon\)
$$ 
and
$$
\wp_\La \lt(\epsilon^2 \frac{\g}{2}\rt)
= \epsilon^2 \wp_\La\lt(\frac{\g}{2}\rt)
= \epsilon^2(m \gamma+ n \gamma \epsilon).
$$
But $ \epsilon \( m \gamma+ n \gamma  \epsilon\)$ and $\epsilon^2 \(m \gamma+ n \gamma  \epsilon\)$ are also lattice points of $\Ga$ and so, poles of $\wp_{\Ga}$.
Hence
$$
\frac{\g}{2},\, \epsilon \frac{\g }{2},\, \epsilon^2 \frac{\g }{2} \in \wp_{\Ga}^{-2}(\infty).
$$
Therefore, all three critical values of the Weierstrass function $\wp_\Gamma$ are prepoles. The second assertion of our theorem immediately follows from it. The fact that $J(\wp_\Ga)=\oc$ now  directly follows from the second assertion and Theorem~\ref{t120200303}. 
\qed

\sp We now pass to such examples based on square, and later rhombic lattices. We recall that a lattice $\Lambda\sbt\C$ is called a {\em square} lattice if and only if
$$
i\Lambda=\Lambda.
$$

\sp Now we will give some examples of elliptic functions whose critical values are poles.

\bthm[\cite{HK1}, Theorem~8.2]\label{th8.2in{HK1}} 
There exists a real square lattice $\Gamma\sbt\C$ such that all the critical values of the Weierstrass function $\wp_\Gamma:\C\lra\oc$
are the lattice points (i.e. belong to $\Ga$), thus poles. 
\fr In consequence, the elliptic function $\wp_{\Gamma}:\C\lra\oc$ is subexpanding (thus compactly non--recurrent) of finite character, non--expanding, and not normal. In addition, $J(\wp_{\Gamma})=\oc$.
\ethm

\bpf Let $\Delta=[\delta, \delta i]$, $\d\in(0,+\infty)$, be the real square lattice associated with the invariants $g_2=4$ (and $g_3=0$); see Proposition~\ref{p120200310} for its existence and uniqueness. Consider the lattice  
$$
\Gamma:=[\g, \g i], 
$$
where $\g=\(\d^2/m\)^{\frac {1}{3}}$ and $m \in \mathbb N$. The homogeneity equation~(\ref{2}) and Proposition~\ref{p120200310} (6) give
$$
\frac{1}{(m\d)^{\frac {2}{3}}}\wp_\Ga \left(\frac{\g}{2}\right)
=\wp_\Delta\left(\frac{\delta}{2}\right)
=e_1(\De)
=1.
$$
Hence,
$$
\wp_\Ga\left(\frac{\g}{2}\right)=(m\d)^{\frac {2}{3}}=m\g.
$$
But $m\g$ is a lattice point of $\wp_{\Ga}$ and thus, a pole of $\wp_{\Ga}$.
The homogeneity property (\ref{2}) implies that
\beq\label{120190219}
\wp_\La (i z)=\frac{1}{i^2}\wp_\La(z)= - \wp_\La(z).
\eeq
Therefore,
$$
\wp_\La (i\g/2)= -\wp_\La(\g/2)=-m\g,
$$
which is also a lattice point  nd so, a pole of $\wp_{\Ga}$. In addition, it follows from Proposition~\ref{p120200310} (6) that $e_3=0$, also a pole of $\wp_{\Ga}$. Hence,
$$
\{e_1, e_2, e_3\}\sbt \wp_{\Ga}^{-2}(\infty).
$$
The second assertion of our theorem immediately follows from it. 
The fact that $J(\wp_\Ga)=\oc$ now  directly follows from the second assertion and Theorem~\ref{t120200303}. 
\epf

\sp We recall, see Section~\ref{weierstrassII} that any lattice $\Lambda$ of the form $[w, \ov w]$, $w\in\C\sms\{0\}$ is called real rhombic and any lattice similar to it is called rhombic. We also recall from this section that any lattice that is real rhombic and square is called real rhombic square.

\sp The  proofs of the next two results are analogous  to the proof of
Theorem~\ref{th8.2in{HK1}} (see Proposition~4.1  in \cite{HKK},
and Theorem~8.9 in \cite{HK2}).

\sp

\bthm\label{prop4.1in{HKK}}  
There exists a real rhombic square lattice $\Gamma\sbt\C$  such that all the critical values of the Weierstrass function $\wp_\Gamma:\C\lra\oc$
are the lattice points (i.e. belong to $\Ga$), thus poles. 

In consequence, the Weierstrass elliptic function $\wp_{\Gamma}:\C\lra\oc$ is subexpanding (thus compactly non--recurrent) of finite character, non--expanding, and not normal. In addition, $J(\wp_{\Gamma})=\oc$.  
\ethm

\bpf Let 
$$
\Lambda:=[2b+2bi, 2b-2bi]=2b[1+i,1-i],  \  \  b>0,
$$
be the real rhombic square lattice associated with the invariants $g_2=-4$
and $g_3=0$, see Proposition~\ref{p520200310} for its existence and uniqueness. It follows from Proposition~\ref{p520200310} that for critical points $b+bi$, $b-bi$, and $2b$ of $\wp_\La$, we have that
\beq\label{220190226}
\wp_\La(b+bi)=-i, \  \
\wp_\La(b-bi)=i, \  \
{\rm and}  \  \
\wp_\La(2b)=0.
\eeq
Note that for any square lattice $\La$, $p\in\C$ is a pole of $\wp_\Lambda$ if and only if  $\pm p i $ is a pole. Thus, all the points $2bj$, $j\in\Z$, are the poles of $\wp_\La$. Fix $j\in\N$. Let 
\beq\label{220190221}
k:=\left(2bj\right)^{-\frac{1}{3}}.
\eeq
Let 
$$
\Gamma:=k\Lambda.
$$
We will show that for critical points $k(b+bi)$, $k(b-bi)$, and $2kb$ of $\wp_\Ga$, we have that
\beq\label{220190219}
\wp_\Ga^2\(\{k(b+bi), k(b-bi), 2kb\}\)=\{\infty\}. 
\eeq
The homogeneity property~(\ref{2}) along with \eqref{220190226} yield
\beq\label{120190221}
\wp_\Ga
 \left(k\(b+bi)\right)=\wp_{k \Lambda}\left(k(b+bi)\right)=\frac{1}{k^2}\wp_\Lambda(b+bi)=-\frac{i}{k^2}=-2jkbi.
\eeq
and
\beq\label{120190221B}
\wp_\Ga
 \left(k\(b-bi)\right)
 =\wp_{k \Lambda}\left(k(b-bi)\right)
 =\frac{i}{k^2}=2jkbi.
\eeq
Since $\pm 2jbi \in \Lambda$, we have that $\pm 2jkbi \in \Gamma$, whence $\wp_\Gamma(\pm 2jkb i)=\infty$. Hence, 
\beq\label{Ga-poles}
\wp_\Ga^2\(\{k(b+bi),k(b-bi)\}\)=\infty.
\eeq
Analogously, we get that
\beq\label{920190223}
\wp_\Ga
  \left(2bk\right)
  =\wp_{k \Lambda}\left(2kb\right)
  =\frac{1}{k^2}\wp_\Lambda(2b)=0.
\eeq
Thus,
\beq\label{zero-poles}
\wp_\Ga^2(2kb)=\infty.
\eeq
The formulas \eqref{Ga-poles} and \eqref{zero-poles} entail \eqref{220190219}. Therefore, all three critical values of the Weierstrass function $\wp_\Gamma$ are prepoles. The second assertion of our theorem immediately follows from it. The fact that $J(\wp_\Ga)=\oc$ now directly follows from the second assertion and Theorem~\ref{t120200303}. 
\epf

\section{Compactly Non--Recurrent Elliptic Functions with Critical Orbits
Clustering at Infinity}\label{clustering}

Theorem~\ref{cor4.8in{HKK}}, the first one we prove below, apart from being interesting in itself is also the main ingredient in the construct of the two classes of examples which will follow it.

\sp We recall, see Definition~\ref{d120200311}, that the prepoles of order $n \geq 0$ of a meromorphic function $f:\C\lra \oc$ were defined as
$$
f^{-n}(\infty)=\big\{z\in {\mathbb C}: \,\,   f^n(z) \,\, \text{is well--defined and}\, f^n(z)=\infty\big\}.
$$
In particular, the poles coincide with order $1$ prepoles. For a Weierstrass function $\wp_\Gamma:\C\lra\oc$, we have the following immediate fact which was used several time above:
$$
f^{-1}(\infty)=\Ga.
$$
Fix a lattice
$$
\Ga=[\g_1,\g_2]\sbt\C.
$$
For any $\a\in {\mathbb C}\sms \{0\}$ set
\beq\label{g}
g_\a:= \a \wp_{\Gamma}.
\eeq
Of course the poles of  $g_\a$ are the same as for $\wp_\Gamma$, i.e. equal to $\Ga$. 
The critical points of $g_\a$ are also obviously the same
as for $\wp_\Ga$, i.e.: 
$$
\Crit\(\wp_{\Ga}\)=(c_1+\Ga)\cup (c_2+\Ga)\cup (c_3+\Ga),
$$ 
but the critical values of $g_\a$ are, in general, different from critical values of $\wp_\Gamma$ since
$$
g_\a\(\Crit(g_\a)\)=\a\wp_{\Ga}\(\Crit\(\wp_{\Ga}\)\), 
$$
where 
$$
c_1=\frac{\g_1}{2},  \  \  c_2=\frac{\g_2}{2},  \  \  c_3=\frac{\g_1+\g_2}{2}.
$$
Thus, we can denote all the critical values of $g_\a$ respectively as
$$
e_1(\a)=g_\a(c_1)=\a e_1, \, \  e_2 (\a)=g_\a(c_2)=\a e_2, \, \
e_3(\a)=g_\a(c_3)=\a e_3.
$$
Recall that
$$
B_R(\infty)=\{z\in \oc: |z|>R \}.
$$ 
The first main technical result of this section, needed in the next ones, is the following. 

\sp\bthm\label{cor4.8in{HKK}} 
Let $\Ga\sbt\C$ be a lattice and let $\wp_\Ga:\C\lra\oc$ be the corresponding Weierstrass function. 

If $q\ge 2$ is an integer and $\xi\in\wp_\Ga^{-q}(\infty)$, then for every $\e\in(0,1)$ there exist a sequence $(\b_k)_{k=q}^\infty\sbt B(1,\e)\sbt\C$ and a sequence $(\e_k)_{k=}^\infty$ of real positive numbers such that

\begin{enumerate}

\, 

\item [(0)]
$$
\b_q=1 
\  \  \  {\rm and}  \  \  \
\e_q=\e,
$$

\item 
$$
\b_k\in B(\b_{k-1},\e_{k-1})\sbt B(1,\e)
$$ 
for all integers $k\ge q+1$.

\,

\item 
$$
\xi\in g_{\b_k}^{-k}(\infty)
$$
for every $k\ge q+1$,

\, \, \,

\item 
$$
\ov B\(\b_{k+1},\e_{k+1}\)\sbt \ov B\(\b_k,\e_k\)\sbt \ov B(1,\e)
$$
for every $k\ge q+1$,

\sp and 

\, \, \,

\item
$$
\liminf_{k\to\infty}\inf\big\{\big|g_\a^k(\xi)\big|:\a\in \ov B(\b_k,\e_k)\big\}=+\infty.
$$
\end{enumerate}
\ethm

\bpf We shall prove Theorem~\ref{cor4.8in{HKK}} by induction on $k\ge q+1$, i.e. we define inductively the sequences 
$$
(\b_k)_{k=q}^\infty\sbt \ov B(1,\e)
\  \  \  {\rm and}  \  \  \ 
(\e_k)_{k=q}^\infty
$$ 
so that the conditions (1)--(4) are satisfied. First, 
$$
\b_q:=1 
\  \  \  {\rm and}  \  \  \
\e_q:=\e,
$$
Now, the base of induction, i.e. $k=q+1$ (yes, $k=q+1$ rather than $k=q$). We define the meromorphic function
$$
G_q:B(1,\varepsilon) \longrightarrow \oc
$$
by
$$
G_q(\alpha):= g_{\alpha}^q(\xi).
$$
Therefore
$$
G_q(1)=g_1^q(\xi)=\wp_\Ga^q(\xi)=\infty.
$$
It is also immediate from the definition of $G_q$ that $G_q\(B(1,\e)\sms\{1\}\)\sbt\C$ if $\varepsilon\in(0,1)$ is small enough. Therefore, for such $\varepsilon$'s, the function $G_q$ has a a unique pole at $\alpha=1$. Consequently, there exists an $R_q > 0$ such that 
$$
B_{R_q}(\infty) \subset G_q(B(1,\varepsilon)).
$$ 
Thus, fixing a point $\gamma_q\in \Ga\cap B_{R_q}(\infty)$, there exists a parameter $\b_{q+1}\in B(1,\varepsilon)$ such that
\[  
g_{\b_{q+1}}^q(\xi)
=G_q(\b_{q+1}) 
= \gamma_q.
\]
Therefore, 
$$
g_{\b_{q+1}}^{q+1}(\xi)=\infty,
$$
meaning that
$$
\xi\in g_{\b_{q+1}}^{-(q+1)}(\infty).
$$
So, 
we have items (1) and (2) established for $k=q+1$.

\sp For the inductive step, assume that for some $k\ge q+1$ the numbers $\b_{q+1}, \b_{q+2},\ld,\b_k\in\C$ and $\e_{q+1},\ld,\e_{k-1}\in(0,1)$ have been defined such that the following hold:

\begin{itemize}

\,

\item item (1) has been established for all integers $q+1\le j\le k$, meaning that 
$$
\b_j\in B(\b_{j-1},\e_{j-1})\sbt B(1,\e)
$$ 
for all integers $q+1\le j\le k$,

\,

\item item (2) has been established for all integers $q+1\le j\le k$, 

\,

\item item (3) has been established for all integers $q+1\le j\le k-2$, 

\,

\item item (4'): 
$$
\inf\big\{\big|g_\a^j(c_1)\big|:\a\in \ov B(\b_j,\e_j)\big\}\ge j
$$
for all integers $q+1\le j\le k-1$,

\,

\item 
and also a non--constant meromorphic functions 
\beq\label{220190222}
G_j:B(1,\varepsilon) \setminus \cup_{3\le i< j}G_i^{-1}(\infty)
\longrightarrow \oc
\eeq
have been defined for all integers $q+1\le j\le k$ by the formula
$$
G_k(\alpha):=g_{\alpha}^k(\xi).
$$
\end{itemize}

\sp Define the meromorphic function
\beq\label{220190222B}
G_{k+1}:B_k(\e):=B(1,\varepsilon) \setminus \cup_{3\le i< k+1}G_i^{-1}(\infty)
\longrightarrow \oc
\eeq
analogously, i.e. by the formula
$$
G_{k+1}(\alpha):=g_{\alpha}^{k+1}(\xi ).
$$
Condition (2) for $k$ means that $G_k(\b_k)=\infty$. Since the function $G_k$ is continuous, this, along with condition (1) (for $k$), implies that there exists $\e_k>0$ so small that 
$$
B(\b_k,\e_k)\sbt B(\b_{k-1},\e_{k-1})
$$ 
and 
$$
\inf\big\{\big|g_\a^k(c_1)\big|:\a\in \ov B(\b_k,\e_k)\big\}\ge k.
$$
This means that conditions (3) and (4') are established respectively for $k-1$ and for $k$. Furthermore, there exists a number $R_k>0$ such that
$$
G_k\(B(\b_k,\e_k)\)\supset B_{R_k}(\infty).
$$
Hence, fixing a point $\gamma_k\in \Ga\cap B_{R_k}(\infty)$, there exists a parameter 
\beq\label{520190223}
\b_{k+1}\in B(\b_k,\e_k)
\eeq 
such that
\[  
g_{\b_{k+1}}^k(\xi)
=G_k(\b_{k+1}) 
= \gamma_k.
\]
Therefore, 
$$
g_{\b_{k+1}}^{k+1}(\xi)=\infty,
$$
meaning that
$$
\xi\in g_{\b_{k+1}}^{-(k+1)}(\infty).
$$
This along with \eqref{520190223} means that conditions (1) and (2) are established for $k+1$. The inductive step of our construction is complete.

\sp Since item (4) follows directly form item (4'), the  proof of Theorem~\ref{cor4.8in{HKK}} is finished.
\epf

\sp We now directly pass to the announced examples. We first briefly recall the examples of  elliptic functions considered in the proof of Theorem~\ref{prop4.1in{HKK}}. Let 
$$
\Lambda =[2b+2bi, 2b-2bi], \  b >0,
$$
be the real rhombic square lattice introduced in its proof. 
We choose $j=4$ as the positive even integer appearing in the proof of Theorem~\ref{prop4.1in{HKK}}. As in this proof,
$$
\Ga=k\La,
$$
where, as in \eqref{220190221} with $j=4$, 
$$
k:=\left(8b\right)^{-\frac{1}{3}}.
$$
We denote
\begin{gather*}\label{720190223}
c_1:=k(b+bi), \  \  c_2:=k(b-bi), \  \  c_3:=2kb 
\  \  {\rm and} \  \  \\
e_1=\wp_\Ga(c_1), \  \  e_2=\wp_\Ga(c_2), \  \  e_3=\wp_\Ga(c_3)
\end{gather*}
Then, recalling \eqref{120190221} and \eqref{Ga-poles}, we have that 
\beq\label{320190221}
\wp_\Ga(c_1) =-(8b)^{\frac{2}{3}}i =e_1\in \Gamma
\  \  {\rm and} \  \ 
\wp_\Ga^2(c_1)
=\wp_\Ga(e_1)
=\infty.
\eeq

\sp\bthm[\cite{HKK},Theorem~4.9]\label{th4.9in{HKK}} 
Let $\Ga$ and $c_1$, $c_2$, and $c_3$ be as above. If 
$$
g_\a=a\wp_\Gamma, \  \  \alpha\in \mathbb C \sms\{0\},
$$
is the family of elliptic functions defined in formula \eqref{g}, then for every $\e\in(0,1)$ there exists a parameter $\b\in B(1,\e)\sbt\C$ such that 
\beq\label{620190223}
\lim_{n \to \infty}g^n_\b(c_1)
=\lim_{n \to \infty}g^n_\b\(\b e_1\)
=\infty
=\lim_{n \to \infty}g^n_\b(c_2)
=\lim_{n \to \infty}g^n_\b\(\b e_2\)
\eeq
and
\beq\label{820190223}
e_3=0,
\eeq
i.e. the forward iterates of the two critical values $\b e_1$ and $\b e_2$ of $g_\b$ diverge to $\infty$ under the action of $g_\b$ while the third one, i.e. $e_3=0$, is its pole. 

\fr In consequence, the elliptic function $g_\b:\C\lra\oc$ is non--expanding, subexpanding, thus compactly non--recurrent. In addition, $J(g_\b)=\oc$. 
\ethm

\bpf Fix an arbitrary $\varepsilon\in(0,1)$ and, looking at \eqref{320190221}, apply Theorem~\ref{cor4.8in{HKK}} with $\xi:=c_1$ and $q=2$ to get the sequences $(\b_k)_{k=2}^\infty\sbt B(1,\e)$ and a sequence $(\e_k)_{k=2}^\infty\sbt(0,1)$. Because of condition (3) of this theorem we have that
$$
\bi_{k=3}^\infty\ov B(\b_k,\e_k)\ne\es.
$$
Fix any $\b$ in this intersection. Then $\b\in B(1,\e)$ and it follows from item (4) of Theorem~\ref{cor4.8in{HKK}} that the first two equality signs of formula \eqref{620190223} hold. 

Since $c_2=c_1i$, we can show, analogously as in (\ref{iterates}), that $g_\b^n(c_2)=g_\b^n(c_1)$ for all integers $n\geq 2$. Thus, the remaining two equality signs of \eqref{620190223} also follow. Formula \eqref{820190223} follows directly from \eqref{720190223} and \eqref{920190223}. Hence, the elliptic function $g_\b:\C\lra\oc$ is non--expanding, subexpanding, thus compactly non--recurrent. 
Furthermore, because of Theorem~\ref{t120200303}, Theorem~\ref{r071708}, Theorem~\ref{Fatou Periodic Components} (Fatou Periodic Components), and because $\Crit(g_\b)=\{\b e_1, \b e_2, e_3\}$, we have that $J(g_\b)=\oc$. 
The proof is complete.
\epf

\sp

We finish this section with a class of examples of non--recurrent elliptic functions whose all critical values diverge to infinity; compare Theorem~4.10 in \cite{HKK}).

\bthm\label{th4.10in{HKK}} 
If $\Gamma\sbt\C$ is a triangular lattice such that the critical value $e_3$ (see Proposition~\ref{cor4.6in{HK1}} (1)) 
of the Weierstrass function $\wp_\Gamma:\C\lra\oc$ belongs to $\wp_\La^{-q}(\infty)$ with some integer $q\ge 2$, (see Theorem~\ref{th8.9in{HK2}} for such examples,
then for every $\e\in(0,1)$ there exists $\b\in B(1,\e)$ such that for the elliptic function $g_\beta:=\b \wp_\Ga:\C\lra\oc$, we have that 
$$
\lim_{n \to \infty}g^n_\b\(\b e_i\)
=\infty
$$
for all $i=1, 2, 3$,  i.e, the forward iterates of all critical
values of $g_\b$ diverge to $\infty$. 

\fr In consequence, the elliptic function $g_\b:\C\lra\oc$ is non--expanding, subexpanding, thus compactly non--recurrent, and $J(g_\b)=\oc$.
\ethm

\bpf Let $c_3\in\C$ be such critical point of $\wp_\Ga$ that $\wp_\Ga(c_3)=e_3$. By out hypotheses $c_1$ is a prepole of $\wp_\Ga$, say of order $q\ge 1$. So, may we apply Theorem~\ref{cor4.8in{HKK}}, in the same way as in the proof of the previous theorem (Theorem~\ref{th4.9in{HKK}}), with $\xi:=c_3$ and $q$ being $q$, to conclude that there exists $\beta \in B(1,\e)$ such that  
$$
\lim_{n \rightarrow \infty}g_{\beta}^n(c_3) 
=\lim_{n \to \infty}g^n_\b\(\b e_3\)
= \infty.
$$
By Proposition~\ref{cor4.6in{HK1}} and (\ref{iterates}), we get
$$
g_{\beta}^n (e_1) = \epsilon^2 g_{\beta}^n(e_3)
\  \  \ {\rm and} \  \  \
g_{\beta}^n (e_2) = \epsilon g_{\beta}^n(e_3).
$$
Therefore,  
$$
\lim_{n \rightarrow \infty} g_{\beta}^n(c_1) = \infty
\  \  \ {\rm and} \  \  \
\lim_{n \rightarrow \infty}
g_{\beta}^n(c_2) = \infty.
$$
In conclusion,
$$
\lim_{n \rightarrow \infty} g_{\beta}^n\(\b e_1\) = \infty
\  \  \ {\rm and} \  \  \
\lim_{n \rightarrow \infty}
g_{\beta}^n\(\b e_2\) 
= \infty.
$$
The proof now can be completed in exactly the same way as the one of the previous theorem.
\epf

\sp\section{Further Examples of Compactly Non--Recurrent Elliptic Functions}

Finally we will provide some examples of non-recurrent elliptic functions with non--empty Fatou set. All of them are non--hyperbolic maps.

\sp In Theorem~\ref{prop4.1in{HKK}} we  considered the Weierstrass  function $\wp_\Gamma:\C\lra\oc$ induced by a real rhombic square lattice
$$
\Gamma=2kb[1+i,1-i],
$$
where 
$$
b\in(0,+\infty), \  \   
k=(2bj)^{-\frac{1}{3}},
\  \  {\rm and} \  \  
j\in\N.
$$ 
We know that
$$
c_{1}=k(b+bi), \  \  c_{2}= k(b-bi),
\  \  {\rm and} \  \ 
c_{3}=2kb
$$
are critical points of $\wp_\Ga$ and every other critical point of $\wp_\Ga$ is congruent mod $\Ga$ with one of them. Also, 
$$
e_i=\wp_\Ga(c_i), \  \  i=1, 2, 3,
$$
are all critical values of $\wp_\Ga$ and that they all are poles of $\wp_\Ga$. In addition $e_3=0$ which is also a pole of $\wp_\Ga$.

We fix $j=1$, so that $k=(2b)^{-\frac{1}{3}}$, and we define 
\beq\label{120190226}
h:=\wp_\Gamma + 2kbi:\C\lra\oc.
\eeq 
The function $h$ is elliptic with poles of order two at each point of $\Gamma$. We shall prove the following.

\bthm\label{prop4.14in{HKK}} 
If $h:\C\lra\oc$ is the elliptic function defined by the formula \eqref{120190226}, then  

\begin{enumerate}

\,

\item $0$, $4kbi$, and $2kbi$ are its all critical values,

\,

\item $0$ and $4kbi$ are also poles of $h$, 

\,

\item $h(2kb)=2kbi$ and $2kbi$ is a superattracting fixed point of $h$. 
\end{enumerate}

\fr In consequence, the elliptic function $h:\C\lra\oc$ is non--expanding, subexpanding, thus compactly non--recurrent, of finite character, and the Julia set $J(h)$ is a proper, nowhere dense, subset of $\oc$.
\ethm

\bpf By formulas \eqref{120190221} and \eqref{120190221B}, we have that
\beq
h\left(k(b + bi)\right)
=\wp_\Gamma \left(k(b + bi)\right)+2kbi= -2kbi + 2kbi=0.
\eeq
and
\beq
h\left(k(b- bi)\right) =
\wp\left(k(b-bi)\right)+2kbi=k2bi +2kbi=4kbi.
\eeq
By \eqref{920190223}, we have that
\beq\label{b-poles}
h(2kb)=\wp_\Ga(2kb)+2kbi=0+2kbi = 2kbi.  
\eeq
Thus, item (1) follows. Item (2) is obvious since $4kbi=2kb(1+i)-2kb(1-i)\in\Ga$.

By \eqref{920190223} and  using also the fact that $\wp_\Lambda(iz)=-i\wp_\Lambda(z)$ (holding because $\Lambda$ is a real rhombic lattice), we get that
\beq\label{b-polesB}
h(2kbi)
=\wp_\Ga(2kbi)+2kbi
=-i\wp_\Ga(2kb)+2kbi
=0+2kbi 
= 2kbi.
\eeq
Thus $2kbi$ is a fixed point of $h$. Since, $\Ga$ is a square lattice, by the the homogeneity equation~(\ref{2a}), we get
$$
h'(2kbi)
=\wp_\Ga'(2kbi)
=\wp_{i\Ga}'(2kbi)
=i^{-3}\wp_\Ga'(2kb)
=i\wp_\Ga'(c_3)
=0.
$$
Therefore, $2kbi$ is a superattracting fixed point of $h$. All other assertions of the theorem follow now immediately and the proof is complete.
\epf

\sp Figure~4 illustrates  the  Julia set  of
  $\wp_{\Lambda}$ with $g_2(\Lambda)\approx 26.5626$ and
  $g_3(\Lambda)\approx -26.2672$ ($\Lambda$ is real rectangular).
  For these values $\wp_{\Lambda}$ has  an attracting  fixed point
  $\xi\approx 1.5566$. A fundamental region is (precisely) a 2 by 1 horizontal real rectangle
  The outline of one fundamental region is in magenta.  The critical values are shown in red, and one is a pole (at around
-2.97). The attracting fixed point is shown in black, and the other two
critical values (red) are in its immediate attracting basin.  Only the darkest blue points are in the Julia set.

\begin{figure}[ht]\label{Figure 4}
\begin{center}
\includegraphics[scale=0.2]{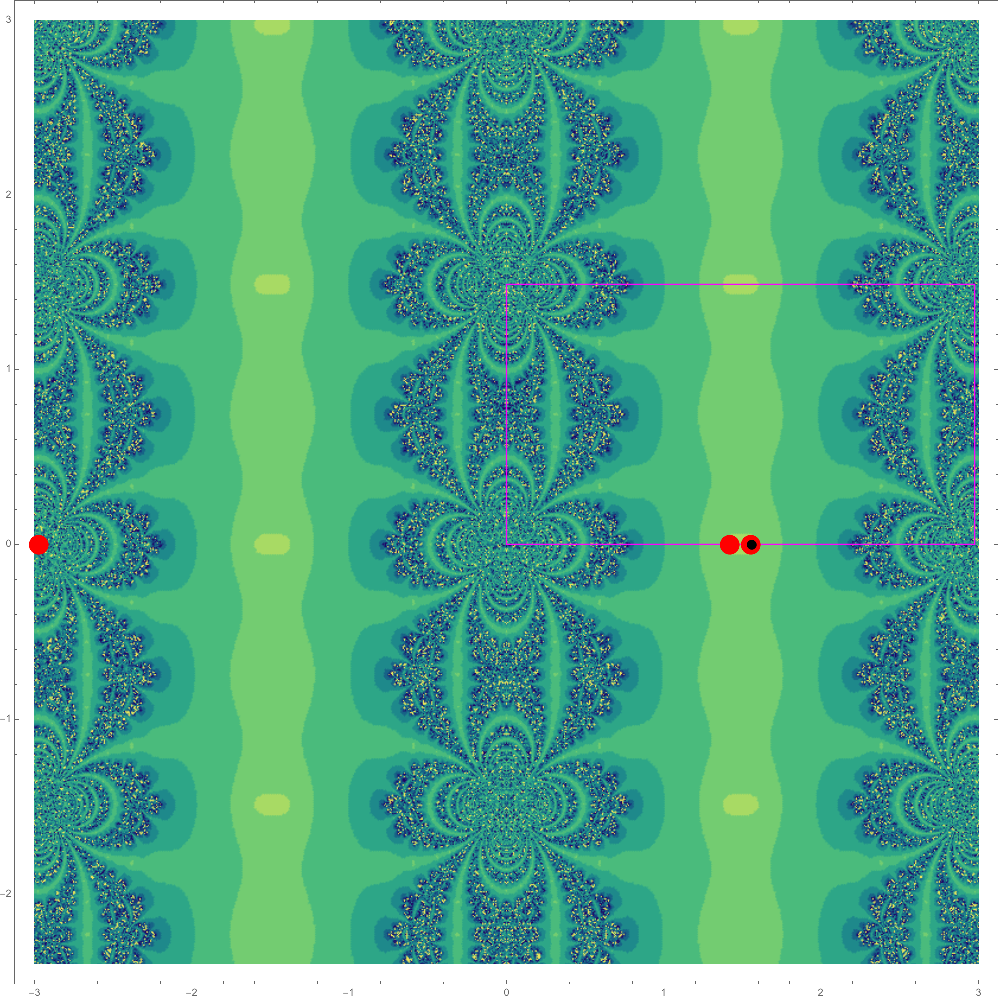}

\caption{$J(\wp_\La)$ with $g_2(\Lambda)\approx 26.5626$.}
\end{center}
\end{figure}

\sp \bthm\label{th4.15in{HKK}} 
Let $h:\C\lra\oc$ be the elliptic function defined by the formula \eqref{120190226} and for every $\alpha\in \mathbb C \sms \{0\}$ let $g_\alpha:=\alpha h$. Then for every $\e\in(0,1)$ there exists $\b\in B(1,\e)$ such that 

\begin{enumerate}
\item The critical value $\b 2kbi$ of $g_\b$ is attracted to
some attracting fixed point of $g_\b$, 

\,

\item The critical value $0$ of $g_\b$ is its pole, and  

\,

\item The critical value $\b 4kbi$ of $g_\b$  diverges to infinity. 
\end{enumerate}

\fr In consequence, the elliptic function $g_\b:\C\lra\oc$ is compactly non--recurrent, non--expanding, and the Julia set $J(g_\b)$ is a proper, nowhere dense, subset of $\oc$. 
\ethm

\bpf Since attracting fixed points are stable under perturbations that are small near them and since, by Theorem~\ref{prop4.14in{HKK}} (3), $2kbi$ is a superattracting fixed point of $h$, there exists $\e\in(0,1)$ so small that for every $\a\in B(1,\e)$ the elliptic function $g_\a$ has an attracting fixed point $\xi_\a$ such that 
$$
\lim_{\a\to 1}\xi_\a=2kbi
$$
and $\a 2kbi$, a critical value of $g_\a$, belongs to the basin of immediate attraction to $\xi_\a$. So, (1) holds for every $\a\in B(1,\e)$.

By Theorem~\ref{prop4.14in{HKK}} (1) and the definition of $g_\a$, $0$ is a critical value of $g_\a$. By Theorem~\ref{prop4.14in{HKK}} (2) it is a pole of $g_\a$, meaning that item (2) of our theorem holds. 
 Applying now Theorem~\ref{cor4.8in{HKK}} in the same way as it was applied in the proof of Theorem~\ref{th4.9in{HKK}}, we produce $\b\in B(1,\e)$ such that
$$
\lim_{n\to \infty}g^n_\b(\b 4kbi)=\infty.
$$
The proof is complete.
\epf

\part{Elliptic Functions B}\label{EFB}

\chapter[Conformal Measures for c.n.r.r. Elliptic Functions]{Sullivan's $h$--Conformal Measures for Compactly Non--Recurrent Elliptic Functions}\label{conformal-measure}

 In this chapter we deal systematically with one of the primary concepts of the book, namely that of (Sullivan's) $h$--conformal (as always $h=\HD(J(f))$) measures for compactly non--recurrent elliptic functions. We will prove their existence for this class of elliptic functions. In Section~\ref{Elliptic Regularity} we will introduce an important class of regular compactly non--recurrent elliptic functions. For this class of elliptic functions we will prove uniqueness and atomlessness of $h$--conformal measures along with their first basic stochastic properties such as ergodicity and conservativity.  We will then assume an elliptic function to be compactly non--recurrent regular throughout the entire book unless explicitly stated otherwise. 
 
We met the concept of conformal measures already in Section~\ref{generalconformalmeasures}, General Notion of Conformal Measures, and in  Section~\ref{Sullivan_Conformal_Measures}, Sullivan's Conformal Measures. We there treated them in a very general setting of the former of these two sections, while in the latter in a setting and spirit quite  close to the one we will be dealing with in the current section. In particular we will heavily use the results of these two sections in the current one. 

We gave in Section~\ref{generalconformalmeasures} quite extended historical account of the concept of conformal measures, particularly Sullivan's ones. We repeat a part of it here for the sake of completeness and convenience of the reader.

Conformal measures were first defined and introduced by Samuel Patterson in in his seminal paper \cite{Pat1} (see also \cite{Pat2}) in the context of Fuchsian groups. Dennis Sullivan extended this concept to all Kleinian groups in \cite{Su1}--\cite{Su3}. He then, in the papers \cite{Su4} --\cite{Su6}, defined conformal measures for all rational functions of the Riemann sphere $\oc$; he also proved their existence therein. Both Patterson
and Sullivan came up with conformal measures in order to get an understanding of geometric measures, i.e. Hausdorff and packing ones. Although already Sullivan noticed that there are conformal measures for Klenian groups that are not equal, nor even equivalent to any Hausdorff and packing (generalized) measure, the main purpose to deal with them is to understand Hausdorff and packing. Chapter~\ref{Markov-systems}, Graph Directed Markov Systems, and Part~\ref{EFA}, Elliptic Functions A, and, especially, the current Part~\ref{EFB}, Elliptic Functions B, of our book, provide a good evidence. 

Conformal measures, in the sense of Sullivan have been studied in in the context of rational functions in greater detail in \cite{DU2}, where, in particular, the structure of the set of their exponents was examined.

Since then conformal measures in the context of rational functions have been studied in numerous research works. We list here only very few of them appearing in the early stages of the development of their theory: \cite{DU LMS}, \cite{DU3}, \cite{DU4}. Subsequently the concept of conformal measures, in the sense of Sullivan, has been extended to countable alphabet iterated functions systems in \cite{MU1} and to conformal graph directed Markov systems in \cite{MU2}. These were treated at length in Chapter~\ref{Markov-systems}, Graph Directed Markov Systems. It was furthermore  extended to some transcendental meromorphic dynamics in \cite{KU1}, \cite{UZ1}, and \cite{MyU1}; see also \cite{UZ2}, \cite{MyU2}, \cite{BKZ1}, and \cite{BKZ1}. Our current construction fits well in this line of development.

Last, the concept of conformal measures found its place also in random dynamics; we cite only \cite{MSU}.
 
\sp

\section[Conformal Measures for Compactly Non--Recurrent Elliptic Functions] {Existence of Conformal Measures for Compactly Non--Recurrent Elliptic Functions}

In this section we prove the existence of $h$--conformal measures for compactly non--recurrent elliptic functions. We also locate their potential atoms.

Let $H:U_1\to U_2$ be an analytic map of open subsets $U_1$, $U_2$ of the complex plane $\mathbb C$. We recall from Definition~\ref{1d20120909} that given $t\ge 0$, a Borel measure $\nu_e$, finite on bounded sets of $\mathbb C$, is an Euclidean semi $t$--conformal measure \index{(N)}{semi $t$-conformal measure} if and only if
$$
\nu_e(H(A))\ge \int_A|H'|^t \,d\nu_e
$$
for every Borel subset $A$ of $U_1$ such that $H|_A$ is one--to--one
and is called $t$--conformal \index{(N)}{$t$--conformal measure} if the ``$\ge$'' sign can be replaced by an ``$=$'' sign. We observed in Section~\ref{conformal pairs of measures} that if $m_s$ is a spherical semi $t$--conformal measure for $f:J(f)\lra J(f)\cup\{\infty\}$, then the measure $m_e$, defined by the requirement that
$\frac{dm_{e}}{dm_s}(z)=(1+|z|^2)^t$ is Euclidean semi $t$--conformal, i.e.
$$
m_{e}(f(A))\ge \int_A|f'|^tdm_{e}  
$$
for every Borel set $A\sbt J(f)$ such that $f|_A$ is $1$--to--$1$. If
$m_s$ is $t$--conformal, then so is $m_e$ in the obvious sense. We spoke in Section~\ref{conformal pairs of measures} of measure $m$ as being either Euclidean or spherical without specifying which one, and respective 
measures $m_{e}$ and $m_s$ being its Euclidean and spherical versions. Obviously $m_{e}$ is equivalent to $m_s$ and is finite on bounded subsets of
$\mathbb C$. 

\sp Assume now that $f:\C\lra\oc$ is an elliptic function and that $m$ is some $t$-conformal measure for $f$. Let $\La$ be the corresponding lattice. Then for every  $w \in \La$ we get
that
$$
\int_A|f'|^tdm_e=m_e(f(A))=m_e(f(A+w))=\int_{A+w}|f'|^tdm_e.
$$
Since the  derivative $f'$  is periodic with respect to the lattice
$\Lambda$, we  thus get the following.

\sp

\bprop\label{p1071405} 
If $f:\C\lra\oc$ is an elliptic function and $m$ is some $t$--conformal measure for $f$, then the  Euclidean $t$--conformal measure $m_e$ is $T_w$-invariant for every $w\in \Lambda$, where $T_w: {\mathbb C} \longmapsto {\mathbb C}$ is the translation  about  the  vector $w$ given by the formula $T_w(z)=z+w$\index{(S)}{$T_w$}. 
\eprop

\sp\fr As an immediate consequence of this proposition, we get  the
following.

\sp
\bcor\label{c2071405} With the hypotheses of the previous proposition, for every $r>0$ we have,
$$
M(t,r):=\inf\{m_e(B_e(z,r)): z \in J(f)\}>0.
$$\index{(S)}{$M(t,r)$} 
\ecor

\sp Throughout the entire chapter $f:\C\lra\oc$ is assumed to be a non--constant elliptic function. Usually and, as ever before
$$
h:=\HD(J(f)).
$$
We prove the existence conformal measures now starting with the following lemma which is also interesting on its own and will be used several times in the book. 

\sp\blem\lab{lhlet} 
Let $f:\C\lra\oc$ be a compactly non--recurrent elliptic function. If $t\ge 0$ and $\nu$ is a $t$--conformal measure, either Euclidean
or spherical, then $t\ge \HD(J(f))$ and $\H^t|_{J(f)}$ is
absolutely continuous with respect to $\nu$ i.e. $ \H^t|_{J(f)}\abs
\nu$. 
\elem

\bpf Since the  measures $\nu_e$ and $\nu_s$ are equivalent
as well   as $\H^t_e$ and $\H^t_s$ are, it suffices to prove the
lemma for Euclidean measures $\nu_e$ and $\H^t_e$.  Now, fix $z\in
J(f)\sms \Sing^-(f)$. Let the sequence $\{n_j\}_{j=1}^{\infty}$,
associated to $z$, come from Proposition~\ref{p1071305}. It follows
from Koebe's Distortion Theorem, I (Euclidean  version)  that
$$ 
f^{-n_j}_z(B_e( f^{n_j}(z), \eta(z)))
\sbt B_e(z, 4^{-1}K \eta(z)|(f^{n_j})'(z)|^{-1}).
$$
Applying Koebe's Distortion Theorem, I (Euclidean  version),  again,
Corollary~\ref{c2071405} and conformality of the measure $m$, we
thus  get
$$
\begin{aligned}
\nu_e\(B_e(z, 4^{-1}K \eta(z)|(f^{n_j})'(z)|^{-1})\)
&\geq K^{-t} |(f^{n_j})'(z)|^{-t}\nu_e (B_e(f^{n_j}(z), \eta(z)/4))\\
& \geq K^{-t} |(f^{n_j})'(z)|^{-t} M(t, \eta(z)/4)\\
&=M(t, \eta(z)/4) K^{-2t}4^t \eta(z)^{-t}
(4^{-1}K \eta(z)|(f^{n_j})'(z)|^{-1})^t.
\end{aligned}
$$
Thus, 
$$
\limsup_{r\to 0}\frac{\nu_e(B_e(z,r))}{r^t}
\geq M(t, \eta(z)/4)(4^{-1}K^2 \eta(z))^{-t},
$$
and we are done because of Theorem~\ref{tncp12.1.}\,(1) and since the set $\Sing^-(f)$ is countable. 
\endpf

\sp Now, we shall prove the main result of this section.

\sp\blem\label{l1071805} 
If $f:\C\lra\oc$ is a compactly non--recurrent elliptic function, then there exists an $h$--conformal measure $m_h$ for $f$ (remember that its spherical version, as all spherical conformal measures considered in this book, is finite) whose all atoms are contained in
$$
\Crit(J(f))\cup \bigcup_{n=1}^{\infty}f^{-n}(\infty).
$$ 
\elem

\bpf First of all we show that it suffice to prove that
all atoms of $m_h$ (if already known to be $h$--conformal) are contained in
$$
Z_f:=\bigcup_{n=0}^{\infty}f^{-n}\(\Crit(J(f))\) \cup \bigcup_{n=1}^{\infty}f^{-n}(\infty).
$$ 
Indeed, assume that $m_h(\Crit(J(f))=0$, $z\in J(f)$ and $f^n(z)\in \Crit(J(f))$ for some integer $n\ge 0$. Then, let $0\le k\le n$ be the least integer such that $f^k(z)\in \Crit(J(f))$. Since $m_{h,e}(f^k(z))=0$, we then conclude, by conformality of $m_{h,e}$, that $m_{h,e}(z)=0$.

We thus pass to the the proof of the existence of an $h$--conformal measure $m_h$
for $f$ whose all atoms are contained in $Z_f$. With the help of
Section~\ref{Sullivan_Conformal_Measures} , particularly 
Lemma~\ref{l220120827}, we shall construct an $h$-conformal measure with
required properties by utilizing the method of $K(V)$ sets \index{(S)}{$K_J(V)$} developed in \cite{DU2}, comp. \cite{KU5}. In order to begin, we
call 
$$
Y \sbt \{\infty\}\cup \Om(f)\cup \bu_{n=1}^\infty f^n(\Crit(J(f)))
$$ 
a crossing set\index{(N)}{crossing set} if
$Y$ is finite and  the following four conditions are  satisfied.

\, \begin{enumerate}

\, \item [(y1)] $ \infty \in Y$.

\, \item [(y2)] $Y \cap \{f^n(x): n \geq 1\} $ is a singleton  for all
$x \in \Crit(J(f))$.

\, \item [(y3)] $Y\cap \Crit(f)=\es$.

\, \item [(y4)] $\Omega(f) \sbt Y$.
\end{enumerate} 

\, \fr  Since $f(\Crit(f))$ is finite, crossing sets do exists. Let $V
\sbt \ov{\mathbb{C}}$ be an open neighborhood of $Y$  such that
\beq\label{e1013007}
\Crit(J(f))\cap  \partial{V}=\es.
\eeq
Define 
\beq\label{KV1}
K_J(V):=J(f)\cap \bigcap_{n\geq 0} f^{-n}(\oc\sms
V)=\{z\in J(f):\forall n \geq0\,\,  f^n(z)\notin V\, \}.
\eeq
These sets will play an important role throughout the book and will be more systematically treated in Section~\ref{CKS}, Closed Invariant Subsets, $K(V)$ Sets, and Summability Properties, where also some historical outlook will be given. 

Obviously \index{(S)}{$K_J(V)$}
$$
f(K_J(V))\sbt K_J(V).  
$$
Since $f: \mathbb{C} \to \oc$ is continuous  and $V$ is  open,  we see  that $K_J(V)$ is a closed subset of $\mathbb{C}$. Since $V \cap K_J(V)=\es$ and
$\infty \in V$, the set $K_J(V)$ has no critical points and is bounded, whence compact. Since 
$$
K_J(V)_0(f)\sbt \bd V,
$$
applying Lemma~\ref{l220120827} with $X=K_J(V)$ and $U(f)=\mathbb{C}$, we directly obtain a number 
\beq\label{220120831}
s(V)\index{(S)}{$s(V)$}:=s\(f|_{K_J(V)}\) \in [0, h]
\eeq 
and a Borel probability measure $m_V $ supported on $K_J(V)$ such that
$$
m_V(f(A))\geq \int_A|f^*|^{s(V)}dm_V
$$  
for every special set $A \sbt K_J(V)$ and
$$
m_V(f(A))= \int_A|f^*|^{s(V)}dm_V
$$
for every special  set $A \sbt K_J(V)\sms \bd V$. Treating the measure
$m_V $ as supported on $J(f)$, with a direct calculation, we get the following.
\beq\label{10.3.7}
m_V(f(A))\geq \int_A|f^*|^{s(V)}dm_V
\eeq
for every special set $A \sbt J(f)$ and
\beq\label{10.3.8}
m_V(f(A))= \int_A|f^*|^{s(V)}dm_V
\eeq
for every special  set $A \sbt J(f)\sms \ov{V}$. 

\sp From now on throughout the entire chapter, we fix a crossing set $Y$ and we consider an open neighborhood $\hat V\sbt \oc$ of $Y$ such that the
closure of $\hat V$ is disjoint from at least one fundamental
parallelogram of $f$. All other neighborhoods of $Y$ considered in
this chapter  will be always assumed to be contained in this
set $\hat V$. 

\sp In order to continue with the proof of
Lemma~\ref{l1071805} we will now single out in a separate definition the
property \eqref{10.3.8} of the measures $m_V$ which will be important for the current proof and also in further investigations of geometrical properties of the
$h$-conformal measure constructed in Lemma~\ref{l1071805}. We will
then proof one lemma about measures with such properties, and afterwads, we will come back to the direct proof of Lemma~\ref{l1071805}.

\sp

\bdfn\label{defsemiconformal}
If $f:\C\lra\oc$ is an elliptic function, then a Euclidean  semi $t$--conformal measure $m_e$ is said to be almost
$t$--conformal\index{(N)}{almost $t$--conformal measure} if and only if there exists a neighborhood $V_{m_e}\sbt\hat V$ of $Y$ such that
$$
m_e(f(A))= \int_A|f'|^tdm_e
$$
\fr for every Borel set $A\sbt J(f)$ such that $f|_A$ is $1$--to--$1$
and $A\cap \ov V_{m_e}=\es$. 
\edfn

\

\fr Hence for every Borel set $A$ such that
$f|_A$ is $1$-to-$1$ and $A\cap \ov V_{m_e}=\es$ and for every $w\in \La$,
we have
\beq\label{1022507}
\int_A|f'|^tdm_{e} =m_{e}(f(A)) =m_{e}(f(A+w)) \ge
\int_{A+w}|f'|^tdm_{e},
\eeq
and the last inequality sign becomes an equality either if in
addition $(A+w)\cap \ov V_{m_e}=\es$ or if $m_e$ is a $t$-conformal
measure, and we assume only that $f|_A$ is $1$-to-$1$. Since $f'$ is
periodic with respect to a lattice $\La$ of $f$, all the above statements
and assumptions lead to the following.

\sp

\blem\lab{l1012502} 
If $f:\C\lra\oc$ is an elliptic function and $\La$ is a lattice of $f$, then for every $w\in\La$, every Borel set $A\sbt
{\mathbb C}$ such that $A\cap\ov V_{m_e}=\es$ and every almost
$t$-conformal measure $m$, we have
\beq\label{2022507}
 m_{e}(A+w)\le m_{e}(A).
\eeq
If either in addition $(A+w)\cap \ov V_{m_e}=\es$ or if $m$ is
$h$-conformal and we assume only that $f|_A$ is $1$-to-$1$, then
this inequality becomes an equality. For every $r>0$ there exists
$M(r)\in (0,\infty)$ independent of any almost $t$-conformal measure
$m$ such that
\beq\lab{1110402}
m_{e}(F)\le M(r)
\eeq
for every Borel set $F\sbt {\mathbb C}$ with the diameter $\le r$.
If in addition $m$ is $h$-conformal, then for every $R>0$ there
exist constants $Q(R)$ and $Q_h(R)$ such that
\beq\lab{2110402}
m_e(B_e(x,r))\ge Q(R)r^2\ge Q_h(R)r^h
\eeq
for all $x\in J(f)$ and all $r\ge R$. \elem

\bpf Inequality (\ref{2022507}) as well as its equality
counterpart are an immediate application of (\ref{1022507}). Formula
(\ref{1110402}) follows directly from (\ref{2022507}) and the fact
that $V$ is disjoint  from at least one fundamental parallelogram.
The second part of formula
 (\ref{2110402}) is  clear. In order  to prove  the first one, fix  a
 fundamental parallelogram $\mathcal R$ of $f$ and notice that
 $$ T(R):=\inf\{m_e(B(z,R)): \, z \in J(f) \cap {\mathcal R} \}>0.$$
Hence, if $$R\leq r\leq 4 \diam_e({\mathcal R}),$$ then  for any
$x\in J(f)$,
$$ m_e(B_e(x,r))\geq m_e(B_e(x,R))\geq T(R)= \frac{T(R)}{r^2}r^2\geq
\frac{1}{16}\frac{T(R)}{\diam^2_e({\mathcal R})}r^2,$$ and we are
done in this case. So suppose  that $r \geq 4\diam ({\mathcal R})$.
Then each ball $B_e(x,r)$  contains at least
$$ \left(\frac{  \sqrt{2}r}{2 \diam_e ({\mathcal R})}\right)^2=
\frac{r^2}{2\diam^2_e({\mathcal R})}$$ non-overlapping
$\La$-congruent copies of $ \mathcal R$. Therefore,
$$ m_e(B_e(x,r))\geq \frac{r^2}{2\diam^2_e({\mathcal R})}
m_e(R)=\frac{m_e({\mathcal R})}{2\diam_e({\mathcal R})}r^2.$$ We are
done. 
\endpf 

\sp Coming back to the proof of  Lemma~\ref{l1071805}, we see that
taking the neighborhood  $V\sbt\hat V$ with sufficiently small
diameter, the limit set $J_S$ of the iterated function system $S$ 
defined in the proof of Theorem~\ref{thm:julia} is contained in
$K_J(V)$. We are free to require the point $z$ used in the construction of 
conformal measures in Section~\ref{Sullivan_Conformal_Measures}, to belong to the set $J_S$. Then, with  notation of Chapter~\ref{Markov-systems} and
the proof of Theorem~\ref{thm:julia}, we get that 
$$
c(t)
\ge \limsup_{n\to\infty}\frac1{2n}\log\sum_{|\om|=n}\|\phi_\om'\|^t
=\frac12\P_S(t).
$$
Thus, by virtue of Theorem~\ref{t1j97} and of the first part of
\eqref{120120831}, we obtain
$$
s(V) \geq h_*:=\HD(J_S)\ge \th_S>1.
$$
Along with \eqref{220120831} this gives that
\beq\label{1J35}
1< h_*\le s(V) \leq h.
\eeq 
Fix from now on $(V_n)_{n=1}^\infty$, a descending  sequence of
neighborhood of $Y$  satisfying (\ref{e1013007})  and such that
$\diam_s(V_n)\le 1/n$.  In view of (\ref{1J35}), passing to
a subsequence, we may  assume without   loss of generality that the
sequence $(s(V_n))_{n=1}^\infty$  converges.  Denote its limit by
$s(Y)$. We then   have
\beq\label{1J37}
1< h_*\leq  s(Y) \leq h.
\eeq
Passing yet to  another subsequence, we may assume that  the
sequence $(m_{V_n})_{n=1}^\infty$,  treated as consisting of probability measures
on the compact  space $\ov{\mathbb{C}}$,  converges  weakly  to a
Borel probability  measure $m_Y$ on  $\ov{\mathbb{C}}$. For the sake
of the proof of Lemma~\ref{l1071805} we shall
prove now the following lemma.

\sp

\blem\label{sl1J37} 
The limit measure $m_Y$ enjoys the following properties.

\, \begin{enumerate}
\item [(a)] $m_Y(f(A))\geq \int_A|f^*|^{s(Y)}dm_Y$  for every  special set
$A\sbt J(f)$.

\, \item [(b)] $m_Y(f(A))= \int_A|f^*|^{s(Y)}dm_Y$  for every  special set
$A\sbt J(f)\sms Y $.

\, \item [(c)] $m_Y(\infty)=0$.

\, \item [(d)] $m_Y(\Om(f))=0$.
\end{enumerate}
\elem

\bpf Let $\Sing(f)$ be the singular set of $f$ as defined
by \eqref{10.1.10} and \eqref{10.1.11}. Note that the sequence of
continuous functions $g_n:=|f^*|^{s(V_n)}$, $n\ge 1$, defined on $J(f)$, converges uniformly to the continuous function 
$$
g:=|f^*|^{s(Y)}.
$$
Let $A$ be a special subset of
$J(f)$ such that 
\beq\label{10.3.9}
\ov A \cap (\Sing(f)\cup Y)=\es. 
\eeq
Then, one can find a compact set $\Ga\sbt J(f)$ disjoint from $\ov A$ and
such that $\text{Int}(\Ga)\spt \Sing(f)\cup Y$. So, invoking (y3) and
disregarding finitely many $n$s if necessary, we have for all $n\ge 1$ that
\beq\label{10.3.10}
\ov{V_n}\sbt \Ga\qquad\text{and}\qquad\ov{V_n}\cap
\Crit(f)=\es.
\eeq
Therefore, by \eqref{10.3.7} and \eqref{10.3.8}, we conclude that
Lemma~\ref{l10.1.6} applies to the sequence of measures
$(m_n)_{n=1}^\infty$ and the sequence of functions $(g_n)_{n=1}^\infty$.
Hence, property (a) of our lemma is satisfied for any special
subset of $J(f)$ and property (b) is satisfied for all our sets $A$ for which
\eqref{10.3.9} holds. So, since any special subset of $J(f)$ disjoint
from $\Sing(f)\cup Y$ can be expressed as a disjoint union
of special sets satisfying \eqref{10.3.9}, an immediate computation shows
that (b) is satisfied for all special sets disjoint from $\Sing(f)\cup Y$.
Therefore, in order to finish the proof of (b), it is enough to show
that the second 
requirement of the lemma is satisfied for every point of the set
$\Sing(f)\sms Y$. First note that by \eqref{10.3.10} and
\eqref{10.3.8}, formula (a') in Lemma~\ref{l10.1.6} holds. 
As $f:J(f)\to J(f)$ is an open map, the set $J(f)_0(f)$, generally defined right after formula \eqref{10.1.11}, is empty, and
$\Sing(f)=\Crit(f)$. Consequently formula (3) of Lemma~\ref{l10.1.6} 
is satisfied for every critical point $c\in J(f)$ of $f$. Since
$g(c)=|f^*(c)|^{s(Y)}=0$, this formula implies that $m(\{f(c)\})\le 0$. Thus
$$
m(\{f(c)\})=0=|f^*(c)|^{s(Y)}m(\{c\}).
$$
The proof of (b) is finished. 

\sp We now shall prove that (c) holds. For every $n\ge 1$ set
$$
m_n^s=m_{V_n}, \quad n \geq 1
$$ 
and 
$$m_n^e=(m_{V_n})_e, \quad s_n=s(V_n).
$$ 
For every $k\geq 0$, consider $S_k$,
the square centered  at the origin whose edges  are parallel to the
coordinates axises and are of length $2^k$. Since by \eqref{10.3.8}
and by (\ref{1J35}), each measure $m_n^s$ is almost
conformal with an  exponent ($=s_n$) in $[1, h]$, and since each 'annulus'
$A_k=S_{k+1}\sms S_k $ is a union of a $3\times 4^k$ unit squares,
it follows   from  (a) and (b),  which say that $m_n^s$ is almost
$s_n$-conformal, and from (\ref{1J35}) that for all $k\geq 1$ and
all $n \geq 1$
$$ 
m_n^e(A_k) \leq 3 M(1) 4^k.
$$
Consequently
$$
m_n^s(A_k)\leq \frac{3 M(1) 4^k}{(1+ 4^k)^{s_n}}
\leq \frac{3 M(1) 4^k}{4^{k s_n}} 
\leq \frac{3 M(1) 4^k}{ 4^{k
h_*}}=3M(1)4^{(1-h_*)k}.
$$ 
Since $h_*> 1$ (see again (\ref{1J35})), we thus get  for all $j\geq 1$ and all
$n\geq 1$ that
$$ 
\aligned
m_n^s(S_{j+1})
&=m_n^s\lt(\bu_{k=j}^\infty A_k\rt)
=\sum_{k=j}^\infty m_n^s(A_k) \\
&\leq \sum_{k=j}^\infty3M(1) 4^{(1-h_*)k}\\
&= 3 M(1)(1-4^{1-h_*})^{-1} 4^{(1-h_*)j}.
\endaligned
$$
Hence $m_Y(\infty)=0$ and we are done with item (c). 

\sp Proving item (d) keep the same measures $m_n$ as in the proof
of item (c). 
Since the measures $m_n$ are semi $s_n$-conformal (see
\eqref{10.3.7}), it follows from Lemma~\ref{ldp4.5} that there exists
a constant $C>0$ such that for every $\om\in \Om(f)$ and every
$r\in(0,1)$, we have 
$$
m_n^s(B(\om,r))
\le Cr^{s_n+p(\om)(s_n-1)}
\le Cr^{h_*+p(\om)(h_*-1)}.
$$
Hence, passing to the limit.
$$
m_Y^s(B(\om,r))
\le Cr^{h_*+p(\om)(h_*-1)}.
$$
Since $h_*>1$, letting $r\downto 0$, this implies that $m_Y(\om))=0$. The proof of Lemma~\ref{sl1J37} is complete. 
\endpf

\sp Now, we are in position to complete easily the proof  of
Lemma~\ref{l1071805}. Let
$$
m_h:=m_Y.
$$
Since the measure $m$ is finite (bounded above by $1$), it directly follows from the definition of $Y$ (in particular the fact that $Y\cap \Sing^{-}(f)=\es$),
Lemma~\ref{sl1J37}, and Corollary~\ref{c1041806} that $m_Y(Y)=0$.
Therefore, since in addition, $f(\Om(f))=\Om(f)$, in order to prove
$s(Y)$--conformality of the measure $m$, it suffices to show that
$$
m_h(f(Y \sms \Om(f)))=0. 
$$ 
But if $y \in Y\sms (\Om(f)\cup\{\infty\})$, then by our  definition  of $Y$, $y \notin
\Sing^{-}(f)$;  and the formula $m_h(f(y))=0$  immediately follows
from Corollary~\ref{c1041806} and the formula
$$
1\ge m_h(f^n(f(y)))\geq |(f^n)^*(f(y))|^{s(Y)}m_h(f(y)),
$$ 
Thus  the $s(Y)$-conformality of $m_h$  is proved; and, in addition, all the atoms  of $m_h$   must  be  contained  in $J(f)\sms \Om(f)$.
In view of \eqref{1J37} and Lemma~\ref{lhlet}, $s(Y)=h$.
Applying now Lemma~\ref{liinfty} and Corollary~\ref{c1041806} we see
that all atoms of $m_h$ must  be contained in $$I_{-}(f) \cup \bu_{n
\geq 0}f^{-n}( \Crit(f)).$$ The proof of Lemma~\ref{l1071805} is complete. 
\epf 

\sp

\section[holomorphic inverse branches]{Conformal Measures for Compactly Non--Recurrent Elliptic Functions and Holomorphic Inverse Branches} 

\sp In this section we keep $f:\C\lra\oc$, a compactly non--recurrent elliptic function. Let $m$ be an almost $t$-conformal measure and let $m_e$ be its
Euclidean version. The upper estimability and strongly lower
estimability will be considered in this section with respect to the
measure $m_e$. When we speak about lower estimability we will assume
more, that the measure $m$ is $t$--conformal. Since the number of
parabolic points is finite, passing to an appropriate iteration, we
assume without loss of generality, in this and the next section, that
all parabolic periodic points of $f$ are simple. 

Consider a closed forward $f$-invariant subset $E$ of $\mathbb C$ such that
$$
\|f'\|_E:=\sup\{|f'(z)|:z\in E\}< +\infty.
$$\index{(S)}{$\mid \mid f' \mid \mid_E$}
\fr\!\!\!\! Such sets will be called
$f$-pseudo-compact\index{(N)}{pseudo-compact}. Obviously, each
$f$-invariant compact subset $E$ of $\mathbb C$ is
$f$-pseudo-compact. Recall that $\th=\th(f)>0$ was defined in (\ref{d5.2}), $\b=\b_f>0$ was defined in \eqref{1042206}, $\a_t(\om)$ in Lemma~\ref{ldp4.5}, and that $\tau>0$ is so small as required in Lemma~\ref{ldp4.4}.

\sp The  proofs of Proposition~4.15 and Proposition~4.16  from
\cite{KU3} carry out verbatim to our current case. We present them now.

\sp

\bprop\lab{pncp16.3} 
Let $f:\C\to\oc$ be a compactly non-recurrent elliptic function. Fix an $f$-pseudo-compact subset $E$ of $J(f)$.
Let $z\in E$, $\l>0$ and let $0<r\le \tau\th\min\{1,\|f'\|_E^{-1}\}\l^{-1}$ be
a real number. Suppose that at least one of the following two
conditions is satisfied:
$$
z\in E\sms\bu_{n\ge 0}f^{-n}(\Crit(J(f))
$$
or
$$
z\in E \  \  \and  \  \ r>
\tau\th\min\big\{1,\|f'\|_E^{-1}\}\l^{-1}\inf\{|(f^n)'(z)|^{-1}:n=1,2,\ld\big\}.
$$

Then, there exists an integer $u=u(\l,r,z)\ge 0$ such that
$$r|(f^j)'(z)|\le \l^{-1}\th\tau$$  for all $0\leq j \leq u$ and the
following four conditions are satisfied.
\beq\lab{ncp1a}
\diam_e\(\Comp(f^j(z),f^{u-j},\l r|(f^u)'(z)|)\)\le\b=\b_f
\eeq
for every $j=0,1,\ld,u$.

For every $\eta>0$ there exists a continuous function $[0,\infty)\ni t\longmapsto
B_t=B_t(\l,\eta)>0$, (independent of $z$, $n$,
and $r$) and  such that if $f^u(z)\in B_e(\om,\th)$ for some
$\om\in\Om(f)$, then
\beq\lab{ncp1b}
f^u(z) \   \text{ is }  \  (\eta
r|(f^u)'(z)|,B_t)-\a_t(\om)\text{-u.e.}
\eeq
and there exists a function $W_t=W_t(\l,\eta):(0,1]\to(0,1]$
(independent of $z$, $n$, and $r$) such that if $f^u(z)\in
B_e(\om,\th)$ for some $\om\in\Om(f)$, then for every $\sg\in(0,1]$
\beq\lab{ncp1c}
f^u(z) \  \text{ is }   \  (\eta r|(f^u)'(z)|,\sg,W_t(\sg))-
\a_t(\om)\text{{\rm-s.l.e.}}
\eeq
with respect the almost $t$-conformal measure $m$.
If $f^u(z)\notin B_e(\Om(f),\th)$,\, then formulas\, (\ref{ncp1b})\,
and (\ref{ncp1c})\,  are\, also\, true\, with\,
\beq\lab{ncp1d}
\a_t(\om)\, \, \text{replaced by $t$.}
\eeq
\eprop

\bpf Suppose first that 
\beq\label{1_2017_11_01}
\sup\{\l r|(f^j)'(z)|:j\ge 0\}> \th\tau\min\{1,\|f'\|_E^{-1}\}
\eeq
and let $n=n(\l,z,r)\ge 0$ be the least integer for which
\beq\lab{ncp16.4}
\l r|(f^n)'(z)|> \th\tau \min\{ 1,\|f'\|_E^{-1}\}.
\eeq
Then $n\ge 1$ (due to the assumption imposed on $r$),
\beq\lab{ncp16.5A}
\l r|(f^j)'(z)|\le \th\tau \min\{ 1,\|f'\|_E^{-1}\}
\eeq
for all $0\le j\le n-1$, and also
\beq\lab{ncp16.5}
\l r|(f^n)'(z)|\le \th\tau.
\eeq

If $f^n(z)\notin B_e(\Om(f),\th)$, set $u=u(\l,r,z):=n$. Then the items
(\ref{ncp1b}), (\ref{ncp1c}) and (\ref{ncp1d}) are obvious in view
of \eqref{ncp16.4} and \eqref{ncp16.5}, while \eqref{ncp1a} follows from \eqref{dp5.4} and \eqref{1_2017_11_02} along with \eqref{ncp16.5}. Thus, we are done in this subcase.

\sp So, suppose that 
\beq\label{2_2017_11_02}
f^n(z)\in B_e(\Om(f),\th),
\eeq
say $f^n(z)\in B_e(\om,\th)$, $\om\in\Om(f)$. Let $0\le k=k(\l,z,r)\le n$ be the
least integer such that $f^j(z)\in B_e(\Om(f),\th)$ for every
$j=k,k+1,\ld,n$. Consider all the numbers
$$
r_i:=|f^i(z)-\om||(f^i)'(z)|^{-1},
$$
where $i=k,k+1,\ld,n$. Put 
$$
||f'||^+_E:= \max\{1,\|f'\|_E\}.
$$
By (\ref{ncp16.4}) we have
$$
r_n
=|f^n(z)-\om||(f^n)'(z)|^{-1}
\le \th||f'||_E^+\th^{-1}\tau^{-1}\l r
=\|f'\|_E^+\tau^{-1} \l r
$$
and, therefore, there exists a minimal $k\le u=u(\l,r,z)\le
n$ such that $r_u \le \|f'\|_E^+\tau^{-1}\l r$. In other words
\beq\lab{ncp16.6}
|f^u(z)-\om| \le \|f'\|_E^+\tau^{-1}\l r|(f^u)'(z)|.
\eeq

Now suppose that
\beq\label{2_2017_11_01}
\sup\{\l r|(f^j)'(z)|:j\ge 0\}\le\th\tau\min\{ 1,||f'||_E^{-1}\}.
\eeq
Then, it follows from Corollary~\ref{c1041806} and our hypotheses that 
$$
z\in \bu_{j\ge 0}f^{-j}(\Om(f)).
$$
Define then the three numbers $u(\l,z,r)$, $k(\l,z,r)$, and $n(\l,z,r)$ to be all equal to the least integer $j\ge 0$ such that $f^j(z)\in \Om(f)$. Denote:
$$
\om=f^u(z).
$$
Notice that in this case formulas (\ref{ncp16.5}) and
(\ref{ncp16.6}) are also satisfied. Our further considerations are
valid in both cases \eqref{1_2017_11_01} with \eqref{2_2017_11_02}, and  \eqref{2_2017_11_01}. First note that by (\ref{ncp16.6}) we have
\beq\lab{1012102}
B_e(f^u(z),\eta r|(f^u)'(z)|)\sbt
B_e(\om,(1+||f'||_E^+\tau^{-1}\eta^{-1}\l)\eta r|(f^u)'(z)|)
\eeq
and, in view of Lemma~\ref{ldp4.5} along with (\ref{ncp16.5A}) and (\ref{ncp16.5}), we get
$$\aligned
m_e\(B_e(f^u(z), \eta r|(f^u)'(z)|)\)\le
 C(1+||f'||_E^+\tau^{-1}\eta^{-1}\l)^{\a_t(\om)}(\eta
r|(f^u)'(z)|)^{\a_t(\om)}.
\endaligned
$$
So, item (\ref{ncp1b}) is proved. Also applying (\ref{ncp16.6}),
Lemma~\ref{lncp14.7.}, Lemma~\ref{lncp12.5.} and (\ref{ncp16.5}) we
see that the point $f^u(z)$ is
$$
\(\|f'\|_E^+\tau^{-1}\l r|(f^u)'(z)|,\sg\tau\|f'\|_E^{-1}\eta\l^{-1}
,2^{\a_t(\om)}L(\om,
2\|f'\|_E^+\th,\sg\tau(2\|f'\|_E^+)^{-1}\eta\l^{-1})\)-\a_t(\om)\text{{\rm-s.l.e.}}
$$
So, if $\|f'\|_E^+\tau^{-1}\l \ge\eta$, then by Lemma~\ref{lncp12.6.}, $f^u(z)$ is
$$
\(\eta r|(f^u)'(z)|,\sg,(2||f'||_E^+\tau^{-1}\l\eta^{-
1})^{\a_t(\om)}L(\om,2\|f'\|_E^+\th,\sg\tau(2||f'||_E)^{-1}
\eta)\l^{-1}\)-\a_t(\om)\text{{\rm-s.l.e.}}
$$
If instead $||f'||_E\tau^{-1}\l \le\eta$, then again it follows from
(\ref{ncp16.6}), Lemma~\ref{lncp14.7.}, Lemma~\ref{lncp12.5.} and
(\ref{ncp16.5}) that the point $f^u(z)$ is 
$$
\(\eta r|(f^u)'(z)|,\sg,2^{\a_t(\om)}L(\om,2\th\tau\l^{-1}\eta,\sg/2) \) - \a_t(\om)\text{{\rm-s.l.e.}}
$$
So, part (\ref{ncp1c}) is also proved.

\sp In order to prove  (\ref{ncp1a}) suppose first that $u=k$. In
particular this is the case if $z\in \bu_{j\ge 0}f^{-j}(\Om(f))$. If $k=0$, we are done since $\l r\le \tau\th$ by our hypotheses, while $\tau\th\le \b_f$ by \eqref{1_2017_11_02}. So, suppose $k\ge 1$. Since $0\le u\le n$, it then follows from \eqref{ncp16.5A} and \eqref{ncp16.5} that
$$
\Comp(f^{k-1}(z),f,r|(f^u)'(z)|)\sbt
\Comp(f^{k-1}(z),f,\th\tau), 
$$
and by the choice of $k$ and (\ref{dp4.3}), we have that $f^{k-1}(z)\notin
B_e(\Om(f),\th)$. Therefore, (\ref{ncp1a}) follows from  (\ref{dp5.4})) and \eqref{1_2017_11_02}.

\sp\fr  If $u>k$ (so we are in the case of \eqref{1_2017_11_01} and \eqref{2_2017_11_02}), then $r_{u-1}>\|f'\|_E^+\tau^{-1}\l r\ge \|f'\|_E\tau^{-1}\l r$, and, using also (\ref{dp4.3}), we get
$$
r_u= {|f^u(z)-\om|\over |f^{u-1}(z)-\om|} |f'(f^{u-1}(z))|^{-1}
r_{u- 1}\ge\|f'\|^{-1}_E r_{u-1} \ge \tau^{-1}\l r.
$$
Hence, $\l r|(f^u)'(z)|\le \tau |f^u(z)-\om|$ and, applying
Lemma~\ref{ldp4.4} and (\ref{dp4.3}) $u-k$ times, we conclude that
for every $k\le j\le u$
$$
\diam_e\(\Comp(f^j(z),f^{u-j},\l
r|(f^u)'(z)|)\)\le\th\tau<\b_f.
$$
And now for all $j=k-1,k-2,\ld,1,0$, the same argument as in the
case of $u=k$ finishes the proof. 
\endpf

\sp

\bprop\lab{pncp16.4} 
Let $f:\C\lra\oc$ be a compactly non-recurrent elliptic function. Fix an $f$--pseudo--compact subset $E$ of $J(f)$. Let $\varepsilon$ and $\l$ be both positive numbers\,  such that $\varepsilon <\l\min\{1,\tau^{-1},\th^{-1}\tau^{-1}\g\}$. If
$0<r<\tau\th\min\{1,\|f'\|_E^{- 1}\}\l^{-1}$ and $z\in E\sms \Crit(J(f))$,
then there exists an integer $s=s(\l,\varepsilon,r,z)\ge 1$ with the
following three properties.
\beq\lab{ncp1aa}
|(f^s)'(z)|\ne 0.
\eeq
If $u=u(\l,r,z)$ of Proposition~\ref{pncp16.3} is well--defined, then $s\le u(\l,r,z)$.

\fr If either $u$ is not defined or $s<u$, then there exists a
critical point $c\in Crit(f)$ such that
\beq\lab{ncp1bb}
|f^s(z)-c|\le\varepsilon r|(f^s)'(z)|.
\eeq
In any case
\beq\lab{ncp1cc}
\Comp\(z,f^s,(KA^2)^{-1}2^{-N_f}\varepsilon r|(f^s)'(z)|\) \cap \Crit(f^s)=\es,
\eeq
where  $A$ was defined in (\ref{d5.3}). 
\eprop

\bpf Since $z\notin \Crit(J(f))$ and in view of
Proposition~\ref{pncp16.3}, there exists a minimal number
$s=s(\l,\varepsilon,r,z)$ for which at least one of the following
two conditions is satisfied
\beq\lab{ncp16.7}
|f^s(z)-c|\le \varepsilon r|(f^s)'(z)|
\eeq
for some $c\in \Crit(J(f))$ or
\beq\lab{ncp16.8}
u(\l,r,z) \  \text{ is well--defined }  \  \  \and  \
s(\l,\varepsilon,r,z)=u(\l,r,z).
\eeq
Since $|(f^s)'(z)|\ne 0$, the parts (\ref{ncp1aa}) and
(\ref{ncp1bb}) are proved.

\sp\fr In order to prove (\ref{ncp1cc}) notice first that no matter which
of the two numbers $s$ is, in view of Proposition~\ref{pncp16.3} we
always have
\beq\lab{ncp16.9}
\varepsilon r|(f^s)'(z)|\le \varepsilon \l^{-1}\th\tau.
\eeq
Let us now argue that for every $0\le j\le s$
\beq\lab{ncp16.10}
\diam_e\(\Comp(f^{s-j}(z),f^j, \varepsilon r|(f^s)'(z)|)\)\le
\b_f.
\eeq
Indeed, if $s=u$, it follows immediately from
Proposition~\ref{pncp16.3} and (\ref{ncp1a}) since $\varepsilon \le
\l$. Otherwise $|f^s(z)-c|\le \varepsilon r|(f^s)'(z)| \le
\varepsilon \l^{-1}\th\tau <\th$ and therefore, by (\ref{d5.2}),
$f^s(z)\notin B_e(\Om(f),\th)$. Thus (\ref{ncp16.10}) follows from
(\ref{dp5.4}).

\sp Now by (\ref{ncp16.10}) and Lemma~\ref{ld5.1}, there exist
$0\le p\le N_f$, an increasing sequence
of integers $1\le k_1<k_2<\ld<k_p\le s$, and mutually distinct
critical points $c_1,c_2,\ld,c_p$ of $f$ such that
\beq\lab{ncp16.11}
\{c_l\}=\Comp\(f^{s-k_l}(z),f^{k_l},\varepsilon
r|(f^s)'(z)|\)\cap \Crit(f)
\eeq
for every $l=1,2,\ld,p$ and if $j\notin \{k_1,k_2,\ld,k_p\}$, then
\beq\lab{ncp16.12}
\Comp\(f^{s-j}(z),f^j,\varepsilon r|(f^s)'(z)|\)\cap
\Crit(f)=\es.
\eeq
Setting $k_0=0$ we shall show by induction that for every $0\le l\le
p$
\beq\lab{ncp16.13}
\Comp\(f^{s-k_l}(z),f^{k_l},(KA^2)^{-1}2^{-l}\varepsilon
r|(f^s)'(z)|\)\cap \Crit(f^{k_l})=\es.
\eeq
Indeed, for $l=0$ there is nothing to prove. So, suppose that
(\ref{ncp16.13}) is true for some $0\le l\le p-1$. Then by
(\ref{ncp16.12})
$$
\Comp\(f^{s-(k_{l+1}-1)}(z),f^{k_{l+1}-1},(KA^2)^{-1}2^{-l}\varepsilon
r |(f^s)'(z)|\)\cap \Crit(f^{k_{l+1}-1})=\es.
$$
So, if
$$
c_{l+1}\in \Comp(f^{s-k_{l+1}}(z),f^{k_{l+1}},(KA^2)^{-1}2^{-
(l+1)}\varepsilon r|(f^s)'(z)|)
$$
then by Lemma~\ref{lncp12.11.} applied for holomorphic maps $H=f$,
$Q=f^{k_{l+1}-1}$ and the radius $R=(KA^2)^{-1}2^{-(l+1)}\varepsilon
r|(f^s)'(z)|<\g$ we get
$$
\aligned |f^{s-k_{l+1}}(z)-c_{l+1}| &\le
KA^2|(f^{k_{l+1}})'(f^{s-k_{l+1}}(z))|^{-
1}(KA^2)^{-1}2^{-(l+1)}\varepsilon r|(f^s)'(z)| \\
&= 2^{-(l+1)}\varepsilon r|(f^{s-k_{l+1}}(z))'| \\
&\le \varepsilon r|(f^{s-k_{l+1}})'(z)|
\endaligned
$$
which contradicts the definition of $s$ and proves (\ref{ncp16.13})
for $l+1$. In particular it follows from (\ref{ncp16.13}) that
$$
\Comp\(z,f^s,(KA^2)^{-1}2^{-N_f}\varepsilon r|(f^s)'(z)|\)\cap \Crit(f^s)=\es.
$$
The proof is finished. 
\epf

\sp We will also need  the following  similar result.

\sp

\blem\lab{l1022506} 
Let $f:\C\lra\oc$ be a compactly non-recurrent elliptic function. Assume that $\Om(f)=\es$. Then there exist  two constants $a, \xi>0$ such  that the following  holds. Suppose that 
$$
z\in J(f)\sms \bu_{n=0}^{\infty}f^{-n}(\{\infty\}\cup \Crit(f)).
$$
Suppose also that  $r\in(0,\g(a\xi)^{-1})$, where $\g>0$ was defined in (\ref{dp5.4}). 

Then there exists an integer $s\geq 0$ with the following  properties:

\, \begin{enumerate}
 \item[(a)]
$ r a \xi |(f^s)'(z)|\geq \g$ or

\, \item[(b)]
$ r a \xi |(f^s)'(z)|<  \g$\\ and

\, \item[(c)]
there exists  a critical point $c\in \Crit(J(f))$  such that
$|(f^s)(z)-c |< r\xi |(f^s)'(z)|$   or

\, \item[(d)]
there exists  a pole  $b\in f^{-1}(\infty)$  such that $|(f^s)(z)-b
|< r\xi |(f^s)'(z)|$.
\end{enumerate}
In either case
$$  
\Comp\(z,f^s, 2\xi r |(f^s)'(z)|\)\cap \Crit(f^s)=\es.$$
\elem

\bpf Put $a=2KA^2 2^{N_f}$,
where  $A$ was defined  in  (\ref{d5.3}). Fix $\rho\in (0, 1/2)$ so
small  that for every $w\in {\mathbb C}\sms(\Crit(f)\cup
f^{-1}(\infty)) $, the map  $f$ restricted  to the set
$$
B_e\(w, 2\rho\, \dist_e \(w, \Crit(f)\cup f^{-1}(\infty)\)\)
$$ 
is one-to-one. Set $\xi= 2^{-4}\rho$. Take $\l>0$ in Proposition~\ref{pncp16.4} such
that $\varepsilon>0$ appearing  there  can be taken  to be equal to
$a \xi$.  In  view of Corollary~\ref{c1041806} there exists a least
integer  $n \geq 0$  such that  $r a \xi |(f^n)'(z)|\geq \g$. Since $ r < \g(a\xi)^{-1}$, we see that $ n\ge 1$. If there exists  an integer $0\leq j \leq n-1$ satisfying  (c) or (d), take $s$ to be the least  one. Otherwise
take $s=n$. By the definition of $n$, it follows  from (\ref{dp5.4}) that
$$ 
\diam_e(\Comp (z,f^k, 2\xi r |(f^k)'(z)|))< \beta_f
$$
for all $k=0, \ldots, n-1$.  Thus, we see that (\ref{ncp16.10}) is
satisfied  if $  s \leq n-1$ and the proof of the last  formula in
our  lemma  is complete  by verbatim  repetition of the  fragment of
the proof of the Lemma~\ref{pncp16.4}  from "Now  by
(\ref{ncp16.10})"  till  its end.  If $s=n$, the same argument shows
that
\beq\lab{11022506}
\Comp\(z, f^{n-1}, 2\xi r |(f^{n-1})'(z)|\)\cap
\Crit(f^{n-1})=\es.
\eeq
By the choice  $\xi$  and the definition  of $n$ we   also  know
that  the map $f^{n-1}$ restricted  to the ball  $B_e(f^{n-1}(z), 16
\xi r| (f^{n-1})'(z)|)$ is injective. Thus  by Koebe's
$\frac{1}{4}$-Theorem
$$
f(B_e(f^{n-1}(z),  16 \xi r| (f^{n-1})'(z)|))\supset B_e(f^{n}(z), 4 \xi r| (f^{n})'(z)|),
$$
and therefore
$$
\Comp\(f^{n-1}(z),f,  2\xi r |(f^{n-1})'(z)|\)\sbt
B_e(f^{n-1}(z), 16\xi r|(f^{n-1})'(z)|).
$$ 
Combining  this with (\ref{11022506}) and injectivity of $f$ restricted  to
$$
B_e(f^{n-1}(z), 16\xi r|(f^{n-1})'(z)|),$$ we  conclude  that
$$
\Comp\(z,f^{n}, 2\xi r |(f^{n})'(z)|\)\cap \Crit(f^{n})=\es.
$$
We are done. \endpf

\sp

\section[Conformal Measure are Atomless, Unique, Ergodic, and
Conservative] {Conformal Measures for Compactly Non--Recurrent Regular Elliptic Functions: Atomlessness, Uniqueness, Ergodicity and Conservativity}\label{Elliptic Regularity} 

In this section we continue dealing with conformal measures. We alreday have their existence, and our goal now is to prove their uniqueness, atomlessness, ergodicity and conservativity. This will require stronger hypotheses than mere compact non--recurrence. In fact it will require one more hypothesis. This hypothesis is the regularity of a compactly non--recurrent elliptic function introduced at the beginning of Section~\ref{DDCoEF}. Firstly, it is needed for us in order to be able to show that the $h$--conformal 
measure constructed in Lemma~\ref{l1071805} is atomless. This in turn is a
prerequisite for, essentially  all, our considerations concerning
geometric measures (Hausdorff and packing) and measurable dynamics
with respect to the measure class generated by the  conformal
measure $m_h$. The place in the book from which on we do need regularity
is the proof of Lemma~\ref{ncp1l7.12.}. Let us record the following immediate observation. 

\bobs\label{o120131012}
Every compactly non--recurrent elliptic function $f:{\mathbb C} \lra\oc$ with $\Crit_{\infty}(f)=\es$ is regular.
\eobs

\fr This simple observation starkingly indicates that the class of all regular
non--recurrent elliptic functions is large indeed, comp. also the entire Section~\ref{examples} devoted to examples of non--recurrent elliptic functions. As an immediate consequence of Observation~\ref{o120131012}, we have the following corollary.

\bcor\label{c1_2017_09_29}
Every expanding and parabolic elliptic function is regular.
\ecor 

\fr Another, sufficient condition, immediately following from Theorem~\ref{thm:julia} for a non--recurrent elliptic function to be regular is this.

\bprop\label{p2_2017_09_29}
If $f:\C\lra\oc$ is a compactly non--recurrent elliptic function and
$$
\frac{2q_{\max}(f)}{q_{\max}(f)+1}>\frac{2 l_{\infty}}{l_{\infty}+1},
$$
then $f$ is regular. 
\eprop

\fr Now, we derive from \eqref{3121605} a technical condition, \eqref{2121705}, which will be directly needed in our considerations involving continuity of conformal measures. It immediately follows from \eqref{3121605} that for every $c\in \Crit_\infty(f)$, $ h>\frac{2 p_cq_c}{p_c q_c+1}.$ Hence
$$ \frac{p_c-1}{p_c}h < (q_c+1)h-2q_c.$$
\fr So there   exists $h_{-}\in(1,h)$\index{(S)}{$h_{-}$} such that
\beq\label{h_{-}}
\frac{p_c-1}{p_c}h_{-} < (q_c+1)h_{-}-2q_c,
\eeq
 and therefore  there
exists $\ka_c>0$\index{(S)}{$\ka_c$} such that
\beq\label{3121705}
  \frac{p_c-1}{p_c}h_{-} < \ka_c< (q_c+1)h_{-}-2q_c.
\eeq
The right-hand  side of this formula  is equivalent to the following
\beq\label{2121705}
  \left(\frac{h_{-} -\ka_c}{2-\ka_c}\right)\left(  \frac{q_c+1}{q_c}\right)>1.
\eeq

We now pass to more general considerations. Let $m_s$ be a Borel probability measure on $\mathbb C$ and let $m_e$ be its Euclidean version, i.e.
$$
\frac{dm_e}{dm_s}(z)=(1+|z|^2)^t.
$$ 
We shall prove the following.

\sp

\blem\lab{l1012602} 
If $z\in \C$, $r_n\downto 0$ and there are two constants $\un M\le ov M$ such
$$
\underline{M}
\leq \liminf_{n\to\infty}\frac{m_e(B_e(z,r_n))}{r^t_n}
\leq \limsup_{n\to\infty}\frac{m_e(B_e(z,r_n))}{r^t_n}
\leq \ov{M},
$$ 
then
$$
\limsup_{n\to\infty}{m_s\(B_s(z,(2(1+|z|^2))^{-1}r_n)\) \over
((2(1+|z|^2))^{-1}r_n)^t} \le 2^t\ov{M}
$$
and
$$
\liminf_{n\to\infty}{m_s\(B_s(z,2(1+|z|^2)^{-1}r_n)\) \over
(2(1+|z|^2)^{-1}r_n)^t} \ge 2^{-t}\underline{M}.
$$
\elem

\bpf Since for every $r>0$ sufficiently small
$$
B_e(z,2^{-1}(1+|z|^2)^{-1}r)\sbt B_s(z,r)\sbt B_e(z,2(1+|z|^2)r)
$$
and since
$$
\lim_{r\downto 0}{m_e(B_e(z,r))\over m_s(B_e(z,r))}=(1+|z|^2)^t,
$$
we get
$$
\limsup_{n\to\infty}{m_s\(B_s(z,(2(1+|z|^2))^{-1}r_n)\) \over
((2(1+|z|^2))^{-1}r_n)^t} \le \lim_{n\to\infty}{m_s(B_e(z,r_n))\over
2^{-t}(1+|z|^2)^{-t}r_n^t}=2^t\ov{M}
$$
and
$$
\liminf_{n\to\infty}{m\(B_s(z,2(1+|z|^2)^{-1}r_n)\) \over
(2(1+|z|^2)^{-1}r_n)^t} \ge \lim_{n\to\infty}{m_s(B_e(z,r_n))\over
2^t(1+|z|^2)^{-t}r_n^t}=2^{-t}\underline{M}.
$$
We are done. \endpf

\sp

Assuming that the compactly non-recurrent elliptic function $f:\C\to\oc$ is regular, our first goal is to show that the $h$--conformal measure $m$
proven to exist in Lemma~\ref{l1071805} is atomless and that
$$
\H^h_s(J(f))=0
$$  
whenever  $h<2$. The regularity assumption will be needed only from Lemma~\ref{ncp1l7.12.} on. We now will consider for $f$ almost $t$-conformal
measures $\nu$ with $t\ge 1$. The notion of upper estimability
introduced in Definition~\ref{dncp12.3.} is considered with respect
to the Euclidean almost $t$--conformal measure $\nu_e$. Recall that
$l=l(f)\ge 1$ is the integer produced in Lemma~\ref{ld5.9} and put
\beq\lab{2071303}
\aligned R_l(f)\index{(S)}{$R_l(f)$}
&:=\inf\{R(f^j,c):c\in \Crit(f) \and 1\le j\le l(f)\} \\
&=\min\{R(f^j,c):c\in \Crit(f)\cap {\mathcal {R}} \and 1\le
j\le l(f)\} <\infty
\endaligned
\eeq
and
\beq\lab{3071303}
\aligned A_l(f)\index{(S)}{$A_l(f)$}
&:=\sup\{A(f^j,c):c\in \Crit(f) \and 1\le j\le l(f)\} \\
&=\max\{A(f^j,c):c\in \Crit(f)\cap {\mathcal {R}}\and 1\le
j\le l(f)\},
\endaligned
\eeq
where the numbers $R(f^j,c)$ and $A(f^j,c)$ are defined just after
Definition~\ref{dcomp}. Since
$$\ov{O_+(f(\Crit_c(J(f))))}$$ is a
compact $f$-invariant subset of $\mathbb C$ (so disjoint from
$f^{-1}(\infty)$) and since $$\ov{{\rm
PC}_c^0(f)}=\ov{O_+(\Crit_c(J(f)))}=\Crit_c(J(f))\cup
\ov{O_+(f(\Crit_c(J(f))))},$$ we have the following straightforward
but useful fact.

\sp\blem\lab{l1071903} 
If $f:\C\lra\oc$ is a compactly non--recurrent elliptic function, then the set $\ov{{\rm PC}^0_c(f)}$ is $f$--pseudo--compact. 
\elem

\sp Recall for the needs of the next two lemmas that the sequence
$\{Cr_i(f)\}_{i=1}^p$ was defined inductively by the formula (\ref{d5.6})
and the sequence $\{S_i(f)\}_{i=1}^p$ was defined by the formula
(\ref{1071303}), while the number $p$, here and in the sequel in this section, comes from Lemma~\ref{ld5.5} (c).

\sp Since the number $N_f$ of equivalence classes of the relation $\sim_f$ between critical points of an elliptic function $f:\C\lra\oc$, is
finite, looking at Lemma~\ref{ld5.9} and Lemma~\ref{l1012502}, the
following lemma follows immediately from Lemma~\ref{lncp13.4.}. 

\sp\blem\lab{ncp1l7.8.} 
Let $f:\C\lra\oc$ be a compactly non--recurrent elliptic function. Fix an integer $1\le i\le p-1$. If $R_{i}^{(u)}>0$ is a positive constant and $t\longmapsto C_{t,i}^{(u)}\in (0,\infty)$, $t\in [1,\infty)$, is a
continuous function such that all points $z\in \ov{{\rm PC}^0_c(f)_i}$ are $(r,C_{t,i}^{(u)})$--$t$--u.e. with respect to some
Euclidean almost $t$--conformal measure $\nu_e$ (with $t\ge 1$) for
all $0<r\le R_{i}^{(u)}$, then there exists a continuous function
$t\longmapsto\^C_{t,i}^{(u)}>0$, $t\in [1,\infty)$, such that all
critical points $c\in Cr_{i+1}(f)$ are
$(r,\^C_{t,i}^{(u)})$-$t$-u.e. with respect to the measure $\nu_e$ for all $0<r\le A_l^{-1}R_{i}^{(u)}$.
\elem

\sp

In the above lemma the superscript $u$ stands for "upper". In the
lemma below it has the same connotation. The number $u$ is also used
to denote the value of the function $u(\l,r,z)$ defined in
Proposition~\ref{pncp16.3}. This should not cause any confusion.

\sp

\blem\lab{ncp1l7.10.} 
Let $f:\C\lra\oc$ be a compactly non--recurrent elliptic function. Fix an integer $1\le i\le p$. If $R_{i,1}^{(u)}>0$ is a positive constant
and $t\longmapsto C_{t,i,1}^{(u)}\in (0,\infty)$, $t\in [1,\infty)$, is
a continuous function such that all critical points $c\in S_i(f)$
are $(r,C_{t,i,1}^{(u)})$--$t$--u.e. with respect to some Euclidean
almost $t$--conformal measure $\nu_e$ (with $t\ge 1$) for all $0<r\le
R_{i,1}^{(u)}$, then there exist  a continuous function $t\longmapsto
\^C_{t,i,1}^{(u)}>0$, $t\in [1,\infty)$, and $\^R_{i,1}^{(u)}>0$
such that all points $z\in \ov{{\rm PC}^0_c(f)_i}$ are
$(r,\^C_{t,i,1}^{(u)})$--$t$--u.e. with respect to the measure $\nu_e$ (with $t\ge 1$) for all $0<r\le\^R_{i,1}^{(u)}$. 
\elem

\bpf 
Put
$$
\varepsilon:=2K(KA^2)2^{N_f},
$$
where $A\ge 1$ was defined  in (\ref{d5.3}). Then fix $\l>0$ so large that
\beq\lab{ncp17.4}
\varepsilon
<\l\min\bigl\{1,\tau^{-1},\th^{-1}\tau^{-1}\min\{\rho,R_{i,1}^{(u)}/2\}\bigr\},
\eeq
where $\rho$ was defined in \eqref{1_2017_11_07}. We shall show that one can take 
$$
\^R_{i,1}^{(u)}
:= \min\lt\{\tau\th \l^{-1}\min\lt\{1,\|f'\|_{\ov{{\rm PC}^0_c(f)_i}}^{-1}\rt\}, R_{i,1}^{(u)}, 1\rt\}
$$
and  
$$ 
\^C_{t,i,1}^{(u)}:= \max\{K^22^tC_{t,i,1}^{(u)}, K^{2t}B_t\},
$$
where $B_t=B_t(\l,\eta)>0$ comes from Proposition~\ref{pncp16.3} with
$\eta=2K$.

Consider $0<r\le \^R_{i,1}^{(u)}$ and $z\in \ov{{\rm PC}^0_c(f)_i}$.
If $z\in \Crit(J(f))$, then $z \in\Crit_c(J(f))$ and  $z\in S_i(f)$,
and we are therefore done. Thus, we may assume that
$z\notin\Crit(J(f))$. Let $s=s(\l,\varepsilon,r,z)$. By the
definition of $\varepsilon$,
\beq\lab{ncp17.5}
2Kr|(f^s)'(z)|=(KA^2)^{-1}2^{-N_f}\varepsilon r|(f^s)'(z)|.
\eeq

Suppose first that $u(\l,r,z)$ is well defined and $s=u(\l,r,z)$.
Then by  item (\ref{ncp1b}) in Proposition~\ref{pncp16.3} or
 by item  (\ref{ncp1d}) in Proposition~\ref{pncp16.3}, we see that the point $f^s(z)$ is $(2K r|(f^s)'(z)|,B_t)$-$t$-u.e. Using (\ref{ncp17.5}),  we obtain  from item (\ref{ncp1cc}) in Proposition~\ref{pncp16.4} and
Lemma~\ref{lncp13.1.} that the point $z$ is $(r,K^{2h}B_t)$-$t$-u.e.. \nl

  If either $u$ is not defined or $s< u(\l,r,z)$, then in view of item (\ref{ncp1cc})
 in Proposition~\ref{pncp16.4}, there exists a
critical point $c\in \Crit_c(J(f))$ such that $|f^s(z)-c|\le
\varepsilon r|(f^s)'(z)|$. Since $s\le u$, by
Proposition~\ref{pncp16.3} and (\ref{ncp17.4}) we get
\beq\lab{ncp17.6}
2Kr|(f^s)'(z)|\le \varepsilon r|(f^s)'(z)| <
\min\{\rho,R_{i,1}^{(u)}/2\}.
\eeq
Since $z\in \ov{{\rm PC}^0_c(f)_i}$, this implies that $c\in
S_i(f)$. Therefore using (\ref{ncp17.6}), the assumptions of
Lemma~\ref{ncp1l7.10.}, and (\ref{ncp17.5}) and then applying
 item  (\ref{ncp1cc})  in Proposition~\ref{pncp16.4} (remember that by Lemma~\ref{l1071903}
the set $\ov{{\rm PC}^0_c(f)}$ is $f$-pseudo-compact) and
Lemma~\ref{lncp13.1.}, we conclude that $z$ is
$(r,K^22^tC_{t,i,1}^{(u)})$-$t$-u.e. The proof is complete. 
\endpf

\sp 

\fr Given an arbitrary integer $k \geq 1$ recall that for any pole $b$ of $f^k$, the number $q_b$ denotes its multiplicity and $B_b^k(R)$\index{(S)}{$B_b^k(R)$}
is the connected component of $f^{-k}(B_\infty^*(R))$ containing $b$. We  have proved Lemma~4.21  in \cite{KU3}  with  no constrained imposed on the elliptic function $f$. In fact, the following more general lemma is true (with  the same  proof), where $f^{-1}$ is replaced by $f^{-k}$.

\sp

\blem\lab{l1071702} 
Let $f:\C\lra\oc$ be a non--constant elliptic function. Fix an integer $k \geq 1$ and a point $b\in f^{-k}(\infty)$. 

If $\nu_e$ is an Euclidean almost $t$-conformal measure with $t>{2q_b\over
q_b+1}$ such that $\nu_e(b)=0$, and if $m$ is the $h$-conformal
measure proven to exist in Lemma~\ref{l1071805}, then
$$
\nu_e(B_b^k(R))\lek R^{2-{q_b+1\over q_b}t}
$$
and
$$
m_e(B(b,r))\gek r^{(q_b+1)h-2q_b}
$$
for all sufficiently small radii $0<r\le 1$.
\elem

\bpf It follows from  Lemma~\ref{l1012502} that
$m_e(\{z\in {\mathbb C}:R\leq |z|< 2R \})\comp R^2$ and
$\nu_e(\{z\in {\mathbb C}: \; R \leq |z| < 2R\})\lek R^2$  for all
$R>0$  large enough. It therefore follows from  (\ref{1110301}) that
\beq\label{1021007}
m_e\(B_b^k(R)\sms \ov{B_b^k(2R)}\)\comp R^2 R^{-\frac{q_b+1}{q_b}h}
\eeq
and
\beq\label{2021007}
\nu_e\(B_b^k(R)\sms \ov{B_b^k(2R)}\)\lek R^2 R^{-\frac{q_b+1}{q_b}t}.
\eeq
Now fix $r >0$ so small that $R=(r/L_k)^{-q_b}$ is large  enough for
formulas (\ref{1021007}) and (\ref{2021007}) to hold. Using
(\ref{u4}) and (\ref{2021007}), we get
\beq\label{2_2017_11_07}
\begin{aligned}
\nu_e(B_b^k(R))&= \nu_e \left( \bu_{j\geq 0}(B_b^k(2^jR)\sms
\ov{B^k_b(2^{j+1}R)})\right)
=\sum_{j= 0}^\infty \nu_e(B_b^k(2^jR)\sms \ov{B^k_b(2^{j+1}R)})\\
& \lek \sum_{j=0}^\infty (2^jR)^2(2^jR)^{-\frac{q_b+1}{q_b}t}\\
&= R^{2-\frac{q_b+1}{q_b}t}\sum_{j=0}^\infty
2^{j(2-\frac{q_b+1}{q_b}t)}\\
&=L_k^{q_b(2-\frac{q_b+1}{q_b}t)}
r^{(q_b+1)t-2q_b}\sum_{j=0}^\infty2^{j(2-\frac{q_b+1}{q_b}t)}\\
&\comp r^{(q_b+1)t-2q_b},
\end{aligned}
\eeq
where the last comparability  sign  holds since $\frac{q_b+1}{q_b}t
>2$. We  are done with  the first  part of our lemma.  

Now replace $\nu_e$ by $m_e$ and $t$ by $h$ (which is  greater  than
$\frac{2q_b}{q_b+1}$ because of Theorem~\ref{thm:julia}) in the above
formula. In this case, the '$\lek $' sign in \eqref{2_2017_11_07}, can, by virtue of (\ref{1021007}), be replaced by the comparability sign '$\comp$'. Since also the first equality sign in \eqref{2_2017_11_07} becomes '$\geq$' (we have not not ruled out the possibility that $m_e(b)>0$ yet) and $m_e(B(b,r))\geq m_e(B_b(R))$,
we are  also done in this case. 
\endpf

\sp
 From now onward, in all our  considerations in this chapter, we assume $f:\C\lra\oc$ to be a compactly non--recurrent regular elliptic function.

\sp

\fr We shall prove now the following.

\sp

\blem\lab{ncp1l7.12.} 
If $f:\C\lra\oc$ is a compactly non--recurrent regular elliptic function, then the $h$--conformal measure $m_h$, for $f:J(f)\lra J(f)\cup\{\infty\}$, proven to exist in Lemma~\ref{l1071805}, is atomless. 
\elem

\bpf By induction on $i=0,1, \ldots, p$ (remember that $p$ comes from Lemma~\ref{ld5.5} (c)), it follows from
Lemma~\ref{ncp1l7.10.} (this lemma provides the base of induction
as $S_0(f)=\es $ and simultaneously contributes to the inductive
step), Lemma~\ref{ncp1l7.8.} and Lemma~\ref{ld5.8} that there exists
a continuous function $t\longmapsto C_t\in (0, \infty)$,\,  $t \in [1,
\infty)$, such that if $\nu$  is an arbitrary almost $t$--conformal
measure on $J(f)$, then
\beq\label{e1020207}
 \nu_e(B(x,r))\leq C_tr^t
\eeq
for all $x\in \ov{{\rm PC}_c^0(f)}$  and all  $r \leq r_0$  for some
$r_0>0$ sufficiently small. Consider now  the almost $s_n$-conformal
measure 
$$
m_n^s:=m_{V_n}
$$
and its Euclidean version
$$
m_n^e:=\(m_{V_n})_e,
$$
both introduced at the beginning of the proof of Lemma~\ref{sl1J37}, where $s_n=s(V_n)$ also comes from the proof of Lemma~\ref{sl1J37}. Letting $n \to \infty$ and recalling that $m_{h,s}$ is a weak limit of measures $m_n^s$, we see
from formula~(\ref{e1020207}) that
\beq\label{e2020207}
 m_{h,e}(B(x,r))\leq C_hr^h
\eeq
for all $x\in \ov{{\rm PC}_c^0(f)}$  and all  $r \leq r_0$. It now follows from Lemma~\ref{l1012602} that
\beq\label{e3020207}
\limsup_{r\searrow 0} \frac{m_s(B(x,r)}{r^h}\leq 2^h C_h
\eeq
for all $x\in \ov{{\rm PC}_c^0(f)}$. In particular,
\beq\label{e4020207}
m_{h,s}(\Crit_c(f))=0.
\eeq
Now  fix $k\geq 1$, $b \in f^{-k}(\infty)$ and $u\in
\(\frac{2q_b}{q_b+1}, h\)$. Consider all integers $n \geq 1$ so
large that $s_n\geq u$. Since $m_n^e(f^{-k}(\infty))\leq m_n^e
(f^{-k}(V_n))=0$, it follows  from  Lemma~\ref{l1071702}  that
$$m_n^e(B^k_b(R))\lek R^{2-\frac{q_b+1}{q_b}s_n}\leq
R^{2-\frac{q_b+1}{q_b}u}.
$$
Hence $m_{h,e}(b)=0$. Since $m_{h,s}$ and $m_{h,e}$ are equivalent  on
$\mathbb{C}$, this gives $m_{h,s}(b)=0$. Consequently
\beq\label{1_2017_11_08}
m_{h,s}\lt(\bu_{n\geq 1} f^{-n}(\infty)\rt)=0.
\eeq
In particular
\beq\label{2_2017_11_08}
m_{h,s}(\Crit_p(f))=0. 
\eeq

We now pass to deal with the set $\Crit_{\infty}(f)$. Since $s_n\nearrow
h$ and since $h_{-}< h$ ($h_{-}$ was defined in (\ref{h_{-}})),
disregarding finitely many $j's$, we may assume without loss  of
generality that
\beq\lab{4121605}
s_j > h_{-}
\eeq
for all  $j\geq 1$.  

Fix $c\in \Crit_{\infty}(f)$. Fix  also $j\geq
1$ and put 
$$
t:=s_j.
$$
Since $\lim_{n\to \infty}f^n(c)=\infty$, there exists an integer $k \geq 1$ such that $q_{b_n}\leq q_c$ (where $b_n\in f^{-1}(\infty)$, defined in \eqref{3_2017_11_08}, is near $f^n(c)$, and $q_c$ was defined in (\ref{q_c})) and
\beq\lab{1121605}
 |f^n(c)|> \max\big\{1, 2\Dist_e(0, f(\Crit(f))\big\}
\eeq
for all $n \geq k$. We ma  need in the course
of the proof the integer $k\geq 1$ to be appropriately bigger. Put  
$$
a:=f^k(c).
$$ 
We recall that $\kappa_c$ was defined in (\ref{3121705}). We shall  prove the  following.

\sp {\bf  Claim 1.}\; {\it  There exists a constant $c_1\geq 1$,
independent of $j$, such that
$$
m_j^e(B_e(a,r))\leq  c_1r^{\kappa_c}
$$
for all $r>0$ small enough independently of $j$.} 

\sp\bpf Put $q=q_c$. In view  of (\ref{1121605}) for every
$n \geq 1$ there exists a unique holomorphic inverse branch
$$
f^{-1}_n: B_e\(f^n(a), \frac{1}{2}|f^n(a)|\)\longrightarrow \mathbb C
$$ 
of $f$ sending $f^n(a)$ to  $f^{n-1}(a)$. Then, by Lemma~\ref{lncp12.9.} and \eqref{3_2017_11_08}, we have for every $n\ge k$ that
$$\begin{aligned}
f^{-1}_n\left(B_e\left(f^n(a),\frac{1}{4}|f^n(a)|\right)\right)& \subset B_e\left(f^{n-1}(a),\frac{K}{4}|f^n(a)||f'(f^{n-1}(a))|^{-1}\right)  \\
& \sbt B_e(f^{n-1}(a), C|f^n(a)||f^n(a)|^{- \frac{q+1}{q}})\\
& = B_e(f^{n-1}(a), C|f^n(a)|^{- \frac{1}{q}})\\
& \sbt B_e\left(f^{n-1}(a), \frac{1}{2}|f^{n-1}(a)|\right)
\end{aligned}
$$
with some constant $C>0$, where the last inclusion was written  assuming
that $|f^{n-1}(a)|\geq 2C |f^n(a)|^{- \frac{1}{q}}$, which holds if the integer $k$ is taken large enough. So, the composition  
$$
f^{-n}_a=f^{-1}_1\circ f^{-1}_2\circ\ldots \circ
f^{-1}_n: B_e\left(f^n(a), \frac{1}{4}|f^n(a)|\right)\longrightarrow {\mathbb C}
$$ 
sending $f^n(a)$ to $a$ is well--defined  and forms a holomorphic branch of $f^{-n}$. Take  $ 0<r< \frac{1}{16}|a| $ and let $n+1\geq 1$ be the least integer such that
$$ 
r|(f^{n+1})'(a)|\geq \frac{1}{16}|f^{n+1}(a)|. 
$$
Such integer  exists since $ |f'(z)|\comp
|f(z)|^{\frac{q_{b}+1}{q_{b}}}$  if $z$ is near  a  pole $b$. By
definition $n \geq 0$, and since $r < \frac{1}{16}|a|$, we have
$$  r |(f^n)'(a)|< \frac{1}{16}|(f^n)(a)|.$$
Then, by $\frac{1}{4}$-Koebe's Distortion Theorem, we have that
\beq\lab{6121605}
B_e(a,r)\sbt f^{-n}_a\(B_e(f^n(a), 4r |(f^n)'(a)|)\).
\eeq
Now we consider three  cases  determined  by the value  of
$r|(f^n)'(a)|$.

\sp {\bf Case 1.}  $\d(f^{-1}(\infty))\leq r |(f^n)'(a)|<
\frac{1}{16}|f^n(a)|$, where $\d(f^{-1}(\infty))$ comes from \eqref{4_2017_11_08}.

\sp In view of (\ref{6121605}) and Koebe's Distortion Theorem  along
with almost conformality of the measure $m_j^e$, we get that
\beq\lab{5121605} 
\begin{aligned}
m_j^e(B_e(a, r))& \leq  K^t|(f^n)'(a)|^{-t}  m_j^e(B_e(f^n(a), 4r|(f^n)'(a)|)) \\
  & \lek K^t|(f^n)'(a)|^{-t}(4r|(f^n)'(a)|)^2\\
& \comp r^2 |(f^n)'(a)|^{2-t}.
\end{aligned}
\eeq
Put $$q_n:=q_{b_n}.$$ Since $t>h_{-}$ (see  formula (\ref{4121605}))
and $q_n\leq  q_c$,  it follows   from (\ref{2121705}) that
$$\lt(\frac{t-\ka_c}{2-\ka_c}\rt)\lt(\frac{q_n+1}{q_n}\rt)>1.$$
Hence,
$$
|f^n(a)| <   |f^n(a)|^{\frac{t-\ka_c}{2-\ka_c}\frac{q_n+1}{q_n}}
         \comp   |f'(f^{n-1}(a))|^{\frac{t-\ka_c}{2-\ka_c}}
         \lek   |(f^n)'(a)|^{\frac{t-\ka_c}{2-\ka_c}}
          = |(f^n)'(a)||(f^n)'(a)|^{\frac{t-2}{2-\ka_c}}.
$$
Combining this and the Case 1 assumption, we get
$$ 
r < \frac{1}{16}|(f^n)'(a)|^{-1}|f^n(a)| \lek |(f^n)'(a)|^{\frac{t-2}{2-\ka_c}}.
$$
Therefore  $r^{2-\ka_c}\lek  |(f^n)'(a)|^{t-2}$, or equivalently,
$r^2|(f^n)'(a)|^{2-t}\lek  r^{\ka_c}$. Together with
(\ref{5121605}), we obtain
\beq\label{5_2017_11_08}
m_j^e(B(a,r))\lek r^{\ka_c}.
\eeq

{\bf  Case 2.}  
$|f^n(a)-b_n|\leq 32
A^{\frac{q_{min}+1}{q_{min}}}r |(f^n)'(a)|
<{32}A^{\frac{q_{min}+1}{q_{min}}}\d(f^{-1}(\infty))$, where $A>0$ was
defined in (\ref{d5.3}) and $q_{\min}$ is the minimal order of all critical points and poles.

\sp Put $\a:=32 A^{\frac{q_{\min}+1}{q_{\min}}}$. Then 
$$
B_e(f^n(a), 4r|(f^n)'(a)|)
\sbt B_e( b_n, (4+c)r |(f^n)'(a)| )   
\sbt B_e(b_n,(4+\a)\d(f^{-1}(\infty)))
$$  
and it follows from Lemma~\ref{l1071702} that
$$  
m_j^e(B_e(f^n(a), 4r|(f^n)'(a)|)) \lek (4r |(f^n)'(a)|)^{(q_n+1)t-2q_n}.
$$
Thus
$$\begin{aligned}
m_j^e(B_e(a,r))& \leq  K^t|(f^n)'(a)|^{-t}(4r |(f^n)'(a)|)^{(q_n+1)t-2q_n}\\
  & \comp r^{(q_n+1)t-2q_n}|(f^n)'(a)|^{(t-2)q_n}\\
& \leq r^{(q_n+1)t-2q_n}.
\end{aligned}
$$
But, as  $q_n \leq q_c$  and $t>h_{-}$, it  follows  from
(\ref{3121705}) that $$ (q_n+1)t -2q_n \geq (q_n+1)t -2q_c > \ka_c$$
and therefore
$$ m_j^e(B(a,r))\leq r^{\ka_c}.$$

\

\fr It remains to consider

\

\fr {\bf  Case 3.} \,\, $ r|(f^n)'(a)|<  \frac{1}{32}
A^{-\frac{q_{\min}+1}{q_{\min}}}|f^n(a)-b_n|$.

\

\fr But then
$$\begin{aligned}
 r |(f^{n+1})'(a)| & =  r |(f^{n})'(a)|  |f'(f^{n}(a))|
                   < \frac{1}{32}   A^{-\frac{q_{\min}+1}{q_{\min}}}|f^n(a)-b_n|(A |f^{n+1}(a)|)^{\frac{q_n+1}{q_n}}\\
               & \leq \frac{1}{32}A^{-\frac{q_{\min}+1}{q_{\min}}}A^{\frac{1}{q_n}+1} |f^{n+1}(a)|\\
                & \leq \frac{1}{32}|f^{n+1}(a)|\\
              & \leq \frac{1}{16}|f^{n+1}(a)|
\end{aligned}
$$
contrary to the definition  of $n$. So, Claim 1  is proved. 
\endpf

\sp The last step  of our  proof is  to demonstrate  the following.

\sp{\bf  Claim 2.} {\it There  exist $c_2\geq  $ and $ R>0$, both
independent of $j$, such that
    $$ m_j^e(B_e(c,r))\leq c_2 r^{p_c\ka_c+ h(1-p_c)}$$
for all $j\geq 1$ and for all  $r\leq R$, where $p_c$ is  the order  of critical  point $c$ of  the map $f^k$.}

\sp{\sl Proof.}  Let $p:=p_c  \geq 2$. There exists $R>0$ so small
that $$f^k(B_e(c), R)\sbt B_e(f^k(c), 2^{-4}|f^k(c)|)$$ and that
there exists $M\geq 1$ such that
$$M^{-1}|z-c|^p\leq |f^k(z)-f^k(c)|\leq M |z-c|^p$$
and
$$ M^{-1}|z-c|^{p-1} \leq |(f^k)'(z)| \leq  M|z-c |^{p-1}$$
for all $z\in B_e(c, R)$. Thus, for all $k\geq 0$ and  all $r\leq R$
$$ 
f^k(A(c;2^{-(l+1)}r, 2^{-l}r)) 
\sbt A\(f^k(c); M^{-1} r^p 2^{-p(l+1)}, Mr^p 2^{-pl}\).
$$
Since the map $f^k_{|B_e(c, R)}$  is  $p$--to--one, using  almost
conformality of the measure $m_j^e$ and the right--hand side
 of (\ref{3121705}), we get   that
 $$\begin{aligned}
m_j^e(A(f^k(c); &\,  M^{-1} r^p 2^{-p(l+1)}, M r^p2^{-pl}))\\
& \geq
\frac{1}{p} M^{-h}( 2^{-(l+1)}r)^{t(p-1)}m_j^e(A(c;2^{-(l+1)}r, 2^{-l}r))\\
& \geq p^{-1}M^{-h}(2^{-(l+1)}r)^{h(p-1)}m_j^e(A(c;2^{-(l+1)}r,
2^{-l}r)).
\end{aligned}$$
Applying Claim 1, we   therefore  get

$$\begin{aligned}
m_j^e(B_e(\, c, \, r)) 
&= \sum_{l=0}^{\infty} m_j^e(A(c,2^{-(l+1)}r, 2^{-l}r))\\
&\leq pM^hr^{h(1-p)}\sum_{l=0}^{\infty}2^{h(p-1)(l+1)} m_j^e\(A\(f^k(c); M^{-1} r^p 2^{-p(l+1)}, M r^p2^{-pl}\)\)\\
& \leq p M^h c_1 2^{h(p-1)} r^{h(1-p)}\sum_{p=0}^{\infty}2^{h(p-1)l} (M r^p2^{-pl})^{\ka_c}\\
 & = p2^{h(p-1)} c_1 M^{h+\ka_c}r^{h(1-p)+p\ka_c}\sum_{l=0}^{\infty}2^{(h(p-1)-p\ka_c)l}\\
 & = p2^{h(p-1)} c_1 M^{h+\ka_c}(1-2^{h(p-1)-p\ka_c})^{-1}r^{p\ka_c+h(1-p)},
\end{aligned}
$$
where writing  the last equality  sign  we  used the fact that
$p\ka_c+h(1-p)>0$ equivalent to  the left--hand side of
(\ref{3121705}). Claim~2 is thus proved. 
\endpf

\sp Repeating again that $p\ka_c+h(1-p)>0$, Claim~2 implies that $m_h(c)=0$. So,
\beq\label{e4020207B}
m_{h,s}(\Crit_\infty(f)))=0.
\eeq
Along with \eqref{e4020207}, \eqref{1_2017_11_08}, \eqref{2_2017_11_08}, and Lemma~\ref{l1071805}, this shows that the measure $m_h$ is atomless and the proof of Lemma~\ref{ncp1l7.12.} is complete.
\endpf 

\sp

The argument from the beginning of the proof  of this lemma, based
on Lemma~\ref{ncp1l7.10.} and Lemma~\ref{ncp1l7.8.} gives the
following,

\sp

\blem\lab{l5042706a} 
If $f:\C\lra\oc$ is a compactly non--recurrent regular elliptic function, then the set $\ov{{\rm PC}^0_c(f)}$ is uniformly $h$--upper estimable with respect to the measure $m_h$ constructed in Lemma~\ref{l1071805}. 
\elem

\sp Denote by $\Tr(f)\index{(S)}{$Tr(f)$}\sbt J(f)$ the set of all
transitive points of\index{(N)}{set of transitive points} $f$, that
is the set of points in $J(f)$ such that $\om(z)=J(f)$.

\sp

\bthm\lab{tmaincm} 
If $f:\C\lra\oc$ is a compactly non--recurrent regular elliptic function, then

\begin{enumerate}
\item There exist a unique $t\ge 0$ and a unique atomless probability spherical $t$--conformal measure $m_s$ for $f:J(f)\to J(f)\cup\{\infty\}$. Then $t=h$.

\, \item The spherical $h$--conformal measure $m_s$ is weakly metrically exact, in particular ergodic\index{(N)}{ergodic measure} and conservative\index{(N)}{conservative measure}, 

\,

\item All other conformal measures are purely atomic, supported on $\Sing^-(f)$ with exponents larger than $h$. 

\, 
\item $m_s(\Tr(f))=1$. 
\end{enumerate} 

\fr In the sequel the $h$--conformal measure $m$, either spherical $m_s$ or its Euclidean version $m_e$, will be denoted by $m_h$. Following the convention of this book the spherical and Euclidean versions of $m_h$ will be respectively denoted by $m_{h,s}$ and $m_{h,e}$.
\ethm

\bpf In view of Lemma~\ref{ncp1l7.12.} there exists an atomless $h$--conformal measure $m_h$ for $f:J(f)\to J(f)\cup\{\infty\}$. So, the existence part of (1) is done. 

Continuing the proof, let $R>0$ be so large that the ball $B_e(0,R)$ contains a fundamental domain of $F$. For evey $w\in\C$, fix $w'\in B(0,R)$ such that
$$
w\sim_f w'.
$$
Suppose that $\nu_e$ is an arbitrary Euclidean $t$-conformal measure for $f$ and some $t\ge 0$. By Lemma~\ref{lhlet}, $t\ge h$. For each 
$$
z\in J(f)\sms \Sing^-(f)
$$
let $(x_k(z))_{k=1}^\infty$ be the sequence produced in Proposition~\ref{p1071305}.   Define for every $l\ge 1$
$$
Z_l:=\big\{z\in J(f)\sms \Sing^-(f):\eta(z)\ge 1/l\big\}.
$$
Fix $l\ge 1$ and assume that $z\in Z_l$. Disregarding finitely terms if needed, assume without loss of generality that
\beq\label{1_2017_11_10}
\big|f^{n_k}(z)-x_k(z)\big|<\frac{1}{32Kl}
\eeq
for all $\ge 1$. Then for all $k\ge 1$, 
$$
B\lt(f^{n_k}(z),\frac{1}{2l}\rt)\sbt B\lt(x_k(z),\frac{1}{l}\rt)
$$
and the holomorphic inverse branch 
$$
f_z^{-n_k}:B_e\left(f^{n_k}(z),{1\over 2l}\right)\longrightarrow\C
$$
sending $f^{n_k}(z)$ to $z$ is well--defined. Using conformality of the measure $\nu$ along with Koebe's $\frac{1}{4}$--Theorem (Theorem~\ref{one-quater}), Koebe's Distortion Theorem I, Euclidean Version (Theorem~\ref{Euclid-I}), and Proposition~\ref{p1071405}, we get the following:
\beq\label{grA}
\begin{aligned}
\nu_e\lt(B_e\lt(z,\frac{1}{16l}|(f^{n_k})'(z)|^{-1}\rt)\rt)
&\le \nu_e\lt(f_z^{-n_k}\lt(B_e\lt(f^{n_k}(z),{1\over 4l}\rt)\rt)\rt) \\
&\le K^t|(f^{n_k})'(z)|^{-t}\nu_e\lt(B_e\lt(f^{n_k}(z),{1\over 4l}\rt)\rt)\\
&\le K^t|(f^{n_k})'(z)|^{-t}\nu_e\lt(B_e\lt(x_k(z),{1\over 2l}\rt)\rt)\\
&=   K^t|(f^{n_k})'(z)|^{-t}\nu_e\lt(B_e\lt(x_k'(z),{1\over 2l}\rt)\rt)\\
&\le K^t\nu_e(B_e(0,R+1))|(f^{n_k})'(z)|^{-t}.
\end{aligned}
\eeq
Likewise, using Lemma~\ref{lncp12.9.}, Koebe's Distortion Theorem I, Euclidean Version (Theorem~\ref{Euclid-I}), and Corollary~\ref{c2071405}, we get the following:
\beq\label{grB}
\begin{aligned}
\nu_e\lt(B_e\lt(z,\frac{1}{16l}|(f^{n_k})'(z)|^{-1}\rt)\rt)
&\ge \nu_e\lt(f_z^{-n_k}\lt(B_e\lt(f^{n_k}(z),{1\over 16Kl}\rt)\rt)\rt) \\
&\ge K^{-t}|(f^{n_k})'(z)|^{-t}\nu_e\lt(B_e\lt(f^{n_k}(z),{1\over 16Kl}\rt)\rt)\\
&\ge K^{-t}|(f^{n_k})'(z)|^{-t}\nu_e\lt(B_e\lt(x_k(z),{1\over 32Kl}\rt)\rt)\\
&=   K^{-t}|(f^{n_k})'(z)|^{-t}\nu_e\lt(B_e\lt(x_k'(z),{1\over 32Kl}\rt)\rt)\\
&\ge K^tM\lt(t,{1\over 32Kl}\rt)|(f^{n_k})'(z)|^{-t},
\end{aligned}
\eeq 
where the constant $M\lt(t,{1\over 32Kl}\rt)$ comes from Corollary~\ref{c2071405}.
Summarizing \eqref{grA} and \eqref{grB}, we obtain
\beq\lab{gr}
B(\nu_e,l)^{-1}|(f^{n_k})'(z)|^{-t} 
\le \nu_e\lt(B_e\lt(z,\frac{1}{16l}|(f^{n_k})'(z)|^{-1}\rt)\rt)
\le B(\nu_e,l)|(f^{n_k})'(z)|^{-t},
\eeq
where $B(\nu_e,l)\ge 1$ is some constant depending only on $R$, $\nu_e$, and $l$.

Fix now $E$, an arbitrary bounded Borel set contained in $Z_l$. Since $m_{h,e}$
is outer regular, for every $x \in E$ there exists a radius $r(x)>0$ of
the form from (\ref{gr}) such that
\beq\lab{gc}
m_{h,e}\lt(\bigcup_{x \in E}B_e(x,r(x))\sms E\rt)< \varepsilon.
\eeq  
Now, by the  Besicovitch's Covering Theorem, i. e. Theorem~\ref{th:6.5.1}, we can choose a countable subcover 
$$
\{ B_e(x_i,r(x_i))\}_{i=1}^{\infty},
$$
$r(x_i)\le \varepsilon$, from the cover $\{B_e(x,r(x))\}_{x \in E}$
of $E$, of multiplicity bounded by some constant $C \ge 1$,
independent of the cover. Therefore by (\ref{gr}) and (\ref{gc}), we
obtain
\beq\lab{bes}
\aligned \nu_e (E) &\le \sum_{i=1}^{\infty}  \nu_e(B_e(x_i,r(x_i)))
 \le B(\nu_{h,e},l)\sum_{i=1}^{\infty}r(x_i)^t \\
&\le B(\nu_e,l)B(m_{h,e},l) \sum_{i=1}^{\infty}r(x_i)^{t-h}m_{e,h}(B_e(x_i,r(x_i)))\\
&\le B(\nu_e,l)B(m_{h,e},l)C\varepsilon^{t-h}m_{h,e}(\bigcup_{i=1}^{\infty}B_e(x_i,r(x_i)))\\
&\le CB(\nu_e,l)B(m_{h,e},l)\varepsilon^{t-h}(\varepsilon +m_{h,e}(E)).
\endaligned
\eeq
In the case when $t>h$, letting $\varepsilon \downto 0$ we obtain
$\nu_e(Z_l)=0$. Since 
$$
J(f)\sms \Sing^-(f)=\bu_{l=1}^\infty Z_l,
$$ 
we therefore get
$$
\nu_e(J(f)\sms\cup \Sing^-(f))=0,
$$ 
which means that $\nu_e(\Sing^-(f))=1$. Thus item (3) of our theorem is proved. 

\sp Suppose now that $t=h$. Then, letting $\varepsilon \downto 0$, formula \eqref{bes} takes on the form:
\beq\lab{bes_B}
\nu_e (E) \le  CB(\nu_e,l)B(m_{h,e},l)m_{h,e}(E).
\eeq
Since this holds for every integer $l\ge 1$, we thus wconclude that 
$$
\nu_e|_{J(f)\sms \Sing^-(f)}\abs m_{h,e}|_{J(f)\sms
\Sing^-(f)}\comp m_{h,s}|_{J(f)\sms \Sing^-(f)}.
$$ 
Reversing the roles of $m_{h,e}$ and $\nu_e$, we infer that
\beq\label{2_2017_11_09}
\nu_{e}|_{J(f)\sms\Sing^-(f)} \   \comp \   m_{h,s}|_{J(f)\sms \Sing^-(f)}.
\eeq
Suppose that $\nu_e(\Sing^-(f))>0$. Then there exists $$y\in
\Crit(J(f))\cup\Om(f) \cup f^{-1}(\infty)$$ such that $m_s(y)>0$.
But then
$$
\sum_{\xi\in y^-}|(f^{n(\xi)})^*(\xi)|^{-h}<+\infty, 
$$
where $y^-=\bu_{n\ge 0}f^{-n}(y)$ and for every $\xi\in y^-$,
$n(\xi)$ is the least integer $n\ge 0$ such that $f^n(\xi)=y$.
Hence,
$$
\nu_y:={\sum_{\xi\in y^-}|(f^{n(\xi)})^*(\xi)|^{-h}\d_{\xi}\over
\sum_{\xi\in y^-}|(f^{n(\xi)})^*(\xi)|^{-h}}
$$
is a spherical $h$--conformal measure supported on $y^-\sbt
\Sing^-(f)$. This contradicts the, already  proven (see \eqref{2_2017_11_09}), fact that the measures
$\nu_y|_{J(f)\sms \Sing^-(f)}$ and $m_{h,s}|_{J(f)\sms \Sing^-(f)}$ are
equivalent and $m_{h,s}(J(f)\sms \Sing^-(f))=1$. Thus $\nu_e$ and $m_{h,s}$
are equivalent.

\sp Let us now prove that any probability spherical $h$--conformal measure $\nu_s$ is ergodic.
Indeed, suppose to the contrary that $f^{-1}(G)=G$ for some Borel
set $G\sbt J(f)$ with $0<\nu_s(G)<1$. But then the two conditional
measures $\nu_G$ and $\nu_{J(f)\sms G}$
$$
\nu_G(B):={\nu_s(B\cap G)\over \nu_s(G)}  
\  \  {\rm and } \  \  
\nu_{J(f)\sms G}(B):={\nu_s(B\cap (J(f)\sms G))\over \nu_s(J(f)\sms G)}
$$
would be $h$--conformal and mutually singular; a contradiction. 

\sp If now $\nu_s$ is again an arbitrary probability spherical $h$--conformal
measure, then by a simple computation based on the definition of
conformal measures we see that the Radon--Nikodym derivative
$\phi:=d\nu_s/dm_{h,s}$ is constant on grand orbits of $f$. Therefore by
ergodicity of $m_{h,s}$ we conclude that $\phi$ is constant $m_{h,s}$--almost
everywhere. As both $m_{h,s}$ and $\nu_s$ are probability measures, this
implies that $\phi=1$ a.e., hence $\nu_s=m_{h,s}$. Thus, item (1) of our theorem is established. 

\sp Let us show now that the probability spherical $h$--conformal measure $m_{h,s}$ is conservative. We shall prove first that $E$, any forward invariant $(f(E)\sbt E)$ Borel subset of $J(f)$, is either of measure $0$ or $1$. Indeed, suppose to the contrary that 
$$
0<m_{h,s}(E)<1.
$$
Let 
$$
\hat E:=\La_f+E=\{w+y:w\in \La_f,\, y\in E\}
$$
Then, the set $\hat E$ is $\La_f$--translation invariant, i.e.
\beq\label{120190914}
w+\hat E=\hat E
\eeq
for all $w\in \La_f$. Furthermore,
$$
E\sbt \hat E, \  m_{h,s}(\hat E)>0,
$$
and
$$
f(\hat E)=f(E)\sbt E\sbt \hat E.  
$$
Since $m_{h,s}(E)<1$ and since $f$ maps the sets of measure $m_{h,s}$ equal to zero into sets of measure $m_{h,s}$ equal to zero, it follows from this that 
$$
m_{h,s}(\hat E)<1.
$$
Since 
$$
m_{h,s}(\Sing^-(f))=0,
$$ 
in order to get a contradiction, it suffices to show that
$$
m_{h,s}(\hat E\sms \Sing^-(f))=0.
$$
Fix an arbitrary point $x\in J(f)$ and an arbitrary radius $R>0$. Seeking contradiction suppose that 
$$
m_{h,e}(B_e(x,R)\sms \hat E)=0.
$$
Then also 
$$
m_{h,s}(B_e(x,R)\sms \hat E)=0.
$$
By conformality of $m_{h,s}$, we have that $m_{h,s}(f(Y))=0$ for all Borel sets $Y\sbt\C$ such that $m_{h,s}(Y)=0$. Hence, using also the fact that 
\beq\label{520190913}
f^n(B_e(x,R)\sms \hat E)\spt f^n(B_e(x,R))\sms f^n(\hat E),
\eeq 
we get that
\beq\lab{2.inv}
\aligned 0&=m_{h,s}\(f^n(B_e(x,R)\sms \hat E)\)
  \ge m_{h,s}\(f^n(B_e(x,R))\sms f^n(\hat E)\) \\
 &\ge m_{h,s}\(f^n(B_e(x,R))\sms \hat E\)
  \ge  m_{h,s}\(f^n(B_e(x,R)\) - m_{h,s}(\hat E)
\endaligned
\eeq
for all $n\ge 0$. By virtue of Proposition~\ref{p120190913} there exists an integer $l\ge 1$ such that $f^l(B_e(x,R))=\oc$. In particular 
$$
m_{h,s}\(f^l(B_e(x,R))\)=1.
$$
Then (\ref{2.inv}) implies that $0\ge
1-m_{h,s}(\hat E)$ which is a contradiction. Consequently 
\beq\label{320190913}
m_{h,e}(B_e(x,R)\sms \hat E)>0.
\eeq

Denote by $Z$ the Borel set of all points $z\in E\sms
(I_\infty(f)\cup\Sing^-(f))$ such that
\beq\lab{1.inv}
\lim_{r\to 0}{m_{h,e}\(B(z,r)\cap (\hat E\sms
(I_\infty(f)\cup\Sing^-(f)))\)\over m_{h,e}(B(z,r))}=1.
\eeq
In view of Lebesgue's Density Theorem, i.e., of Theorem~\ref{th:6.5.4}, we have that $m_{h,s}(Z)=m_{h,s}(\hat E)$. Since $m_{h,s}(E)>0$, there exists
at least one point $z\in Z$. Since
$$
z\in J(f)\sms (I_\infty(f)\cup\Sing^-(f)),
$$ 
Proposition~\ref{p1071305} applies, and let $(x_j(z))_{j=1}^\infty$, $\eta(z)>0$, and an increasing sequence $(n_j)_{j=1}^\infty$ be given by this proposition. Put
$$
\d=\eta(z)/8.
$$ 
It then follows from \eqref{320190913} and Proposition~\ref{p1071305} that for every $j\ge 1$ large enough we have that,
\beq\label{220190914}
m_{h,e}\(B_e(x_j(z),\d)\sms \hat E\)>0.
\eeq
Therefore, as $f^{-1}(J(f)\sms E)\sbt
J(f)\sms E$, the standard application of Theorem~\ref{Euclid-I} and
Lemma~\ref{lncp13.1.} shows that
\beq\label{320190914}
\limsup_{r\to 0}{m_{h,e}(B(z,r)\sms \hat E)\over m_{h,e}(B(z,r))}>0
\eeq
which contradicts (\ref{1.inv}). Thus either 
\beq\label{420190913}
m_{h,s}(E)=0
\  \  \  {\rm or} \  \  \
m_{h,s}(E)=1.
\eeq

Now conservativity is straightforward. One needs to prove that for
every Borel set $B\sbt J(f)$ with $m_{h,s}(B)>0$, one has $m_{h,s}(G)=0$, where
$$
G:=\lt\{x\in J(f): \sum_{n\ge 0}\1_B(f^n(x))<+\infty\rt\}.
$$
Indeed, suppose that $m(G)>0$, and, for all $n\ge 0$, let
$$
G_n:=\lt\{x\in J(f): \sum_{k\ge n}\1_B(f^n(x))=0\rt\}
=\big\{x\in J(f):f^k(x)\notin B \  \text{ for all } \  k\ge n\big\}.
$$
Since 
$$
G=\bu_{n\ge 0}G_n,
$$
there exists $k\ge 0$ such that
$m_{h,s}(G_k)>0$. Since all the sets $G_n$ are forward invariant we get from \eqref{420190913} that 
$$
m_h(G_k)=1.
$$
But, on the other hand, all the sets
$f^{-n}(B)$, $n\ge k$, are of positive measure and are disjoint from
$G_k$. This contradiction finishes the proof of conservativity of
$m_{h,s}$. The item (2) is established. Because of (2) and since $\supp(m_{h,s})=J(f)$, we have $m_{h,s}(\Tr(f))=1$, i.e. (4). 
\endpf 

\chapter[Hausdorff and Packing Measures of c.n.r.r. Elliptic Functions]{Hausdorff and Packing Measures of Compactly Non-Recurrent Regular Elliptic Functions}
\label{Hausdorff-and-packing}

 From now on, throughout this chapter, and, in fact, throughout the entire book, $\H^t_e\index{(S)}{$\H^t_e$}$ stands for  the $t$--dimensional
Hausdorff measure on $\mathbb{C}$ with respect to the Euclidean
metric whereas $\H^t_s$\index{(S)}{$\H^t_s$} refers  to its
spherical counterpart. The same convention is applied to the packing
measures $\Pi^t_e\index{(S)}{$\Pi^t_e$}$ and
$\Pi^t_s$\index{(S)}{$\Pi^t_s$}. Note that the  measures $\H^t_e$
and $\H^t_s$ as well as $\Pi^t_e$ and $\Pi^t_s$ are equivalent with
Radon-Nikodym derivative bounded away from zero and $\infty$ on
compact  subsets of $\mathbb{C}$. In particular  the Hausdorff
dimension of any  subset $A$  of $\mathbb{C}$  has the same value no
matter whether calculated with respect to the Euclidean or spherical
metric; it will  be denoted in the sequel simply by $\HD(A)$. If
$\H^t$ or $\Pi^t$ will be  endowed neither with  the subscript '$e$'
nor '$s$', it will refer simultaneously to both Euclidean as well as spherical measures. As in the previous chapters we keep
$$
h=\HD(J(f)).
$$

The goal of this chapter can be viewed as two folded. The first one is to provide a geometrical characterization of the $h$--conformal measure $m_h$, which, with the absence of parabolic points, turns out to be a normalized packing measure, and the second one is to give a complete description of geometric, Hausdorff and packing, measures of the Julia sets $J(f)$. All of this is contained in the following theorem.

\sp
\bthm\label{t1020907} 
Let $f:\C\lra\oc$  be a compactly non--recurrent regular elliptic function. If $h=\HD(J(f))=2$, then $J(f)=\mathbb C$. If $h <2$, then
\begin{enumerate}
\item [(a)] $\H_s^h(J(f))=0$.

\, \item [(b)] $\Pi^h_s(J(f))>0$.

\, \item [(c)] $\Pi^h_s(J(f))<+\infty$ if and only if $\Om(f)=\es$.
\end{enumerate}  
In either the case, (c) or if $\HD(J(f))=2$ the unique spherical $h$--conformal measure $m_h$ coincides with the normalized packing measure $\Pi^h_s/\Pi^h_s(J(f))$ restricted to the Julia set $J(f)$.
\ethm

\sp This theorem has an interesting story, for expanding rational functions $f$ we always have, essentially due to \cite{Bow2}, that 
$$
0<\H^h(J(f)), \Pi^h(J(f))<+\infty
$$
and these two measures coincided up to a multiplicative constant. Their probability version is then the unique $h$--conformal measure. If $f$ is still a rational function but parabolic, or more generally, non--recurrent, then (see \cite{DU3} and \cite{U1} respectively):

\ben 
\item [(a)] $\H_h(J(f))<+\infty$ and $\Pi_h(J(f))>0$. 

\,

\item [(b)] $\H_h(J(f))=0$ if and only if $h<1$ and $\Om(f)\ne\es$. 

\,

\item [(c)] $\Pi_h(J(f))=+\infty$ if and only if $h>1$ and $\Om(T)\ne\es$. 
\een

So, the description of Hausdorff and packing measures in both cases of non--recurrent rational functions and compactly non--recurrent regular elliptic function coincide except that in the latter case $\H_h(J(f))\le 1$ never holds. For other transcendental meromorphic and entire functions, even hyperbolic (expanding), the situation is generally less clear and varies from case to case. See for example \cite{UZ1} and \cite{MyU3}.

\sp As an immediate consequence of Theorem~\ref{t1020907}, we  get the following.

\sp\bcor\label{c1020907} 
Let $f:{\mathbb C} \lra\oc$  be a compactly non-recurrent regular elliptic function. If $\Om(f)=\es$, then the Euclidean
$h$-dimensional packing  measures $\Pi^h_e$\index{(S)}{$\Pi^h_e$} is
finite on each bounded  subset of $J(f)$. 
\ecor

\sp

\section{Hausdorff Measure}\label{Hausdorff}

\fr We start with the proof of the first part of Theorem~\ref{t1020907}. Our first preparatory result is the following.

\blem\label{l1_2017_11_10}
If $f:{\mathbb C}\lra\oc$  is a compactly non-recurrent elliptic function, then
$$
\bu_{j=1}^{\infty} f^{-j}(\infty)\sms \ov{O_+(\Crit(f))}\ne\es.
$$
\elem

{\sl Proof.} Seeking contradiction, suppose that
$$
f^{-1}(\infty)\sbt \ov{O_+(\Crit(f))}.
$$
So, for each $b\in f^{-1}(\infty)$ there exists $c_b\in \Crit(f)\cap J(f)$ such that $b\in \ov{O_+(c_b))}$. It then follows from Definition~\ref{pseudo-non-recurrent} that its items (1) and (3) are ruled out for $c_b$, whence item (2) must hold. We then conclude that \beq\label{5_2017_11_1}
b\in O_+(f(c_b)). 
\eeq
Since then $O_+(f(c_b))$ is a finite set and since $f(\Crit(f))$ is also a finite set, we conclude that 
$$
\bu_{b\in f^{-1}(\infty)}O_+(f(c_b)) 
$$
is a finite set. But, \eqref{5_2017_11_1} implies that 
$$
f^{-1}(\infty)\sbt \bu_{b\in f^{-1}(\infty)}O_+(f(c_b)).
$$
Since $f^{-1}(\infty)$ is infinite, we arrived at a contradiction, and we are thus done.
\endpf  

\sp{\bf The proof of part (a) of Theorem~\ref{t1020907}.} 
By Lemma~\ref{l1_2017_11_10} there exists $b \in f^{-1}(\infty)\sms \ov{O_+(\Crit(f))}$,
say
$$
b \in f^{-1}(\infty)\sms \ov{O_+(\Crit(f))}.
$$
Hence, there exists
$\ka>0$ such that
\beq\lab{1071802}
B_e(b,3\ka)\cap O_+(\Crit(f))=\es.
\eeq 
Consider an arbitrary point $z\in \Tr(f)$. Then, there exists an infinite increasing sequence
$\{n_j\}_{j=0}^\infty$ such that
\beq\lab{2071702}
\lim_{j\to\infty}f^{n_j}(z)=b \  \  \text{ and } \  \ |f^{n_j}(z)-b|<\ka/2
\eeq
for every $j\ge 1$. It follows from this and (\ref{1071802}) that
for every $j\ge 1$ there exists a holomorphic inverse branch
$$
f_z^{-n_j}:B_e(f^{n_j}(z),2\ka)\lra\mathbb C
$$ 
of $f^{n_j}$ sending $f^{n_j}(z)$ to $z$. Let $m_h$ be the unique $h$-conformal atomless measure proven to exist in Theorem~\ref{tmaincm}. Using now Theorem~\ref{Euclid-I}, Lemma~\ref{lncp12.9.}, Lemma~\ref{lncp13.1.}, and Lemma~\ref{l1071702}, we conclude that 
$$
\aligned 
m_{h,e}\big(B_e\big(z,2K|(f^{n_j})'(z)|^{-1}&|f^{n_j}(z)-b|\big)\big)\ge \\
&\ge m_{h,e}\big(f_z^{-n_j}\(B_e\(f^{n_j}(z), 2|f^{n_j}(z)-b|\)\)\) \\
&\ge K^{-h}m_{h,e}\(B_e\(f^{n_j}(z), 2|f^{n_j}(z)-b|\)\)|(f^{n_j})'(z)|^{-h} \\
&\ge K^{-h}m_{h,e}\(B_e(b,|f^{n_j}(z)-b|)\)|(f^{n_j})'(z)|^{-h} \\
&\gek |f^{n_j}(z)-b|^{(q_b+1)h-2q_b}|(f^{n_j})'(z)|^{-h} \\
&= \(2K|(f^{n_j})'(z)|^{-1}|f^{n_j}(z)-b|\)^h(2K)^{-h}
  |f^{n_j}(z)-b|^{q_b(h-2)}.
\endaligned
$$
Since $h<2$, using (\ref{2071702}), this implies that
$$
\ov{\lim_{r\to 0}}\frac{m_{h,e}(B_e(z,r))}{r^{h}}
\ge \ov{\lim_{j\to\infty}}\frac{m_{h,e}\big(B_e\big(z,2K|(f^{n_j})'(z)|^{-1}|f^{n_j}(z)-b|\big)\big)}{\(2K|(f^{n_j})'(z)|^{-1}|f^{n_j}(z)-b|\)^h}
=+\infty.
$$ 
Hence $$\H^h_e(\Tr(f))=0$$ in view of Theorem~\ref{tncp12.1.}\,(1). Since by
Theorem~\ref{tmaincm} $m_{h,e}(J(f)\sms \Tr(f))=0$, it follows from
Lemma~\ref{lhlet} that $\H^h_e(J(f)\sms \Tr(f))=0$. In conclusion,
$$
\H^h_e(J(f))=0,
$$ 
which completes the proof. 
\endpf

\sp

\section{Packing Measure I} 

In this section we shall prove Proposition~\ref{ncpp17.5}, stated just below. Its item (3) is just item (b) of Theorem~\ref{t1020907}, while item (1) contributes towards the last assertion of this theorem. We shall also prove Lemma~\ref{ncp1l7.7.} which establishes one side of item (c).

\sp\bprop\lab{ncpp17.5} 
If $f:\C\lra \oc$ is a compactly non-recurrent regular elliptic function, then

\begin{enumerate}
\item The $h$-conformal measure $m_h$ is absolutely
continuous with respect to the packing measure $\Pi^h$ and moreover,

\, \item The Radon-Nikodym derivative $dm_s/d\Pi_s^h$ is uniformly bounded away
from infinity. In particular:

\item
$$
\Pi^h(J(f))>0.
$$ 
\end{enumerate}
\eprop 

\bpf Since $$J(f)\cap\om\(\Crit(f)\sms
\Crit(J(f))\)=\Om(f),$$ we conclude from Lemma~\ref{ld5.2} that
there exists $y\in J(f)$ at a positive distance, denote it by $8\eta$, from
$O_+(\Crit(f))$. Fix $z\in\Tr(f)$.  Then there exists an infinite
sequence $(n_j)_{j=1}^\infty$ of increasing positive integers such that $f^{n_j}(z)\in B_e(y,\eta)$ for every $j\ge 1$. Hence,
$$
B_e(f^{n_j}(z),4\eta)\cap O_+(\Crit(f))=\es.
$$ 
Consequently,
$$
\Comp\(z,f^{n_j},4\eta\)\cap \Crit(f^{n_j})=\es.
$$
Hence, it follows from Lemma~\ref{lncp12.9.} and
Lemma~\ref{lncp13.1.} that
$$
\liminf_{r\to 0}{m_{h,e}(B_e(z,r))\over r^h}\le B
$$
for some constant $B\in (0,\infty)$ and all $z\in \Tr(f)$. Applying
Lemma~\ref{l1012602} we therefore get that
$$
\liminf_{r\to 0}{m_{h,s}(B_s(z,r))\over r^h}\le 2^hB.
$$
Therefore, by Theorem~\ref{tncp12.2.} (1), the measure $m_{h,s}|_{\Tr(f)}$ is
absolutely continuous with respect to $\Pi^h_s|_{\Tr(f)}$. Since, by
Theorem~\ref{tmaincm}, $m_{h,s}(J(f)\sms \Tr(f))=0$, we are done.
\endpf

\sp

\blem\lab{ncp1l7.7.} 
If $f:\C\lra \oc$ is a compactly non-recurrent regular elliptic function and $\Om(f)\ne\es$, then
$$
\Pi^h_s(J(f))=+\infty.
$$
\elem

\bpf Fix $\xi\in\Om(f)$. Since the set 
$$
\bu_{n\ge 0}f^{-n}(\xi)
$$
is dense in $J(f)$ and, since, by Lemma~\ref{ld5.2}, $\om(\Crit(f))$ is
nowhere dense in $J(f)$, there exist an integer $s\ge 0$, a real number
$\eta>0$, and a point 
$$
y\in f^{- s}(\xi)\sms B_e\lt(\bu_{n\ge 0}f^n(\Crit(f)),\eta\rt).
$$ 
Since by Theorem~\ref{thm:julia}, $h>1$,
it follows from Lemma~\ref{ldp4.5} and Lemma~\ref{lncp13.4.} ($y$
may happen to be a critical point of $f^s$!) that
\beq\lab{ncp17.2}
\liminf_{r\to 0}{m_e(B_e(y,r))\over r^h}=0.
\eeq
Consider now a transitive point $z\in J(f)$, i.e. $z\in \Tr(f)$.
Then there exists an infinite increasing sequence $n_j=n_j(z)\ge 1$, $j\ge 1$,
of positive integers such that
$$
\lim_{j\to\infty}|f^{n_j}(z)-y|=0 \  \  \  \and  \  \  \
r_j=|f^{n_j}(z)- y|<\eta/7
$$
for every $j\ge 1$. By the choice of $y$, for all $j\ge 1$ there
exist holomorphic inverse branches 
$$
f_z^{-n_j}:B_e(f^{n_j}(z),6r_j)\lra\C
$$ 
of $f^{n_j}$ sending $f^{n_j}(z)$ to $z$. So, applying Lemma~\ref{lncp12.9.} and
Lemma~\ref{lncp13.1.} with $R=3r_j$, we conclude from
(\ref{ncp17.2}) that
$$
\liminf_{r\to 0}{m_{h,e}(B_e(z,r))\over r^h}=0.
$$
Applying Lemma~\ref{l1012602}, we conclude that the same formulas
remain true with $m_{h,e}$ replaced by $m_{h,s}$ and $B_e(z,r)$ by
$B_s(z,r)$. Therefore, it follows from Theorem~\ref{tmaincm}
($m_{h,s}(\Tr(f))=1$) and Theorem~\ref{tncp12.2.} (1) that
$\Pi^h_s(J(f))=+\infty$. We are done. 
\endpf

\sp

\section{Packing measure II}

\sp As before, from now on throughout this section $m_h$ denotes the unique atomless $h$-conformal measure  proven to exist in Theorem~\ref{tmaincm}. Our aim in this
section is to show that in the absence of parabolic periodic points the
$h$-dimensional Euclidean  packing measure is finite on bounded subsets of $J(f)$ and that $\Pi^h_s(J(f))< +\infty$. This will complete the item (c) of Theorem~\ref{t1020907}.

\sp Recall that the numbers $R_l(f)$ and $A_l(f)$ have been defined by formulas
(\ref{2071303}) and (\ref{3071303}) respectively. 

\sp Recall for the needs of this section that the sequence
$\{Cr_i(f)\}_{i=1}^p$ was defined inductively by the formula (\ref{d5.6})
and the sequence $\{S_i(f)\}_{i=1}^p$ was defined by the formula
(\ref{1071303}), while the number $p$, here and in the sequel in this section, comes from Lemma~\ref{ld5.5} (c).

\sp 

Since the number $N_f$ of equivalence classes of the relation $\sim_f$ between critical points of an elliptic function $f:\C\lra \oc$, is
finite, looking at Lemma~\ref{ld5.9} and Lemma~\ref{l1012502}, the
following lemma follows immediately from Lemma~\ref{lncp13.5.}. 

\sp\blem\lab{ncp1l7.9.} 
Let $f:\C\lra \oc$ be a compactly non--recurrent regular elliptic function. Fix $0\le i\le p-1$. If $C_{i}^{(l)}>0$, $0<R_{i}^{(l)}\le
R_l(f)/3$, and $0<\sg\le 1$ are three real numbers such that all
points $z\in \ov{{\rm PC}^0_c(f)}_i$ are
$(r,\sg,C_{i}^{(l)})$-$h$-s.l.e. with respect to the measure $m_{h,e}$,
then there exists $\^C_{i}^{(l)}>0$ such that all critical points
$c\in Cr_{i+1}(f)$ are $(r,\^\sg,\^C_{i}^{(l)})$-$h$-s.l.e. with
respect to the measure $m_{h,e}$ for all $0<r\le
A_l(f)^{-1}R_{i}^{(l)}$, where $\^\sg$ was defined in
Lemma~\ref{lncp13.5.}. 
\elem

\sp Let us prove the following. 

\sp \blem\lab{ncp1l7.11.} 
Let $f:\C\lra \oc$ be a compactly non--recurrent regular elliptic function. Suppose that $\Om(f)=\es$. Fix $0\le i\le p$. Assume that
$C_{i,1}^{(l)}>0$, $R_{i,1}^{(l)}>0$ and $0<\sg\le 1$ are three real
numbers  such that all critical points $c\in S_i(f)$ are
$(r,\sg,C_{i,1}^{(l)})$-$h$-s.l.e. with respect to the measure $m_{h,e}$
for all $0<r\le R_{i,1}^{(l)}$. Then there exist
$\^C_{i,1}^{(l)}>0$, $\^R_{i,1}^{(l)}>0$ such that all points
$z\in \ov{{\rm PC}^0_c(f)_i}$ are $(r,8K^3A^22^{N_f}
\sg,\^C_{i,1}^{(l)})$-$h$-s.l.e. with respect to the
measure $m_{h,e}$ for all $0<r\le \^R_{i,1}^{(l)}$, where $A>0$ was
defined  in (\ref{d5.3}). 
\elem

\bpf Recall that by Lemma~\ref{l1071903} the set
$\ov{{\rm PC}^0_c(f)}$ is $f$-pseudo-compact. We shall show that
this time one can take
$$
\^R_{i,1}^{(l)}:=\min\big\{\tau\th\min\{1,\|f'\|_i^{-1}\}\l^{-1},R_{i,1}^{(l)},1\big\}
\  \and \  
\^C_{i,1}^{(l)}:= \(8(KA^2)2^{N_f}\)^hC_{i,1}^{(l)},
$$
where $\|f'\|_i:=\|f'\|_{\ov{{\rm PC}^0_c(f)}_i}$. Indeed, take
$\varepsilon:=4K(KA^2)2^{N_f}$ and then choose $\l>0$ so large that
\beq\lab{ncp17.7}
\varepsilon <\l\min\left\{1,\tau^{-1},\th^{-1}\tau^{-
1}\min\{\rho,R_{i,1}^{(l)}/2\}\right\}.
\eeq
Consider $0<r\le \^R_{i,1}^{(l)}$ and $z\in {\ov{{\rm
PC}^0_c(f)}_i}$. If $z\in \Crit_c(J(f))$, then $z\in S_i(f)$ and we
are done. Thus, we may assume that $z\notin \Crit_c(J(f))$, then
$z\notin \Crit(J(f))$.

\sp Let $s=s(\l,\varepsilon,r,z)$. By the definition of $\varepsilon$,
\beq\lab{ncp17.8}
4Kr|(f^s)'(z)|=(KA^2)^{-1}2^{-N_f}\varepsilon r|(f^s)'(z)|.
\eeq
Suppose first that $u(\l,r,z)$ is well defined and $s=u(\l,r,z)$.
Then by item  (\ref{ncp1c}) in Proposition~\ref{pncp16.3}, applied
with $\eta=K$, we see that the point
$$
f^s(z) \text{ is }
(Kr|(f^s)'(z)|,\sg/K^2,W_h(\sg/K^2))-h\text{-s.l.e.}.
$$
Using (\ref{ncp17.8}) it follows from item  (\ref{ncp1cc}) in
Proposition~\ref{pncp16.4}  and Lemma~\ref{lncp13.2.} that the point
$z$ is $(r,\sg,W_h(\sg/K^2))$-$h$-s.l.e.. If either $u$ is not
defined or $s\le u(\l,r,z)$, then in view of item  (\ref{ncp1bb}) in
Proposition~\ref{pncp16.4}, there exists a critical point $c\in
\Crit(f)$ such that $$|f^s(z)-c|\le \varepsilon r|(f^s)'(z)|.$$
Since $s\le u$, by Proposition~\ref{pncp16.4} and (\ref{ncp17.7}) we
get
\beq\lab{ncp17.9}
|f^s(z)-c |\le \varepsilon r|(f^s)'(z)| <
\min\{\rho,R_{i,1}^{(l)}/2\}.
\eeq
Since $z\in {\ov{{\rm PC}^0_c(f)}_i}$ it implies that $c\in S_i(f)$.
Therefore, by the assumptions of Lemma~\ref{ncp1l7.11.} and by
(\ref{ncp17.9}) we conclude that $c$ is $(2\varepsilon
r|(f^s)'(z)|,\sg,C_{i,1}^{(l)})$-$h$-s.l.e.. Consequently, in view
of Lemma~\ref{lncp12.5.}, the point $f^s(z)$ is $(\varepsilon
r|(f^s)'(z)|,2\sg,2^h C_{i,1}^{(l)})$-$h$-s.l.e.. So, by
Lemma~\ref{lncp12.6.} this point is
$$
(K r|(f^s)'(z)|,2\sg \varepsilon /K,(2\varepsilon
K^{-1})^hC_{i,1}^{(l)})-h\text{-s.l.e.}
$$
Using now formula (\ref{ncp17.8}) and  item (\ref{ncp1cc}) in
Proposition~\ref{pncp16.4}  along with  the fact that
$K\varepsilon^{-1}< 1$ we have from  Lemma~\ref{lncp13.2.} that the
point $z$ is $(r,2K \varepsilon \sg,(2 \varepsilon
K^{-1})^hC_{i,1}^{(l)})$-$h$-s.l.e..
The proof is complete. \endpf

\sp\fr As a fairly  straightforward  consequence  of these  two lemmas
we get the following.

\sp\blem\lab{l1042006} 
If $f:\C\lra \oc$ is a compactly non--recurrent regular elliptic function, then with some $R >0$ and some $G>0$, each point of $\ov{{\rm PC}^0_c(f)}$ (in particular  each point of $\Crit_c(f)$) is $(r,1/2,G)$-$h$-s.l.e. with respect to the measure $m_{h,e}$ for every $r\in[0,R]$.
\elem

\bpf Since $S_0(f)=\es $, starting  with $\sigma> 0$ as
small as we wish, it immediately follow from Lemma~\ref{ncp1l7.11.},
Lemma~\ref{ncp1l7.9.} and Lemma~\ref{ld5.8}   by induction on
$i=0,1,\ldots, p$ that all  the points  of $S_i(f)$ and  $\ov{{\rm
PC}_c^0(f)}_i$ are $(r, 1/2,G)$-$h$-s.l.e. with  same $G,R>0$ and
all $r\in [0, R]$. We are done. \endpf

\sp\fr This lemma  and  Lemma~\ref{l1071702}, taken  together, yield the
following.

\sp\blem\lab{l2022506} 
If $f:\C\lra \oc$ is a compactly non--recurrent regular elliptic function, then every point of the set $\Crit(J(f)) \cup
f^{-1}(\infty)$ is $h$-s.l.e. with respect to the measure $m_{h,e}$ with  $\sigma\in (0,1)$  arbitrary. 
\elem

\sp Fix $c\in\Crit_{\infty}(f)$. Since $\lim_{n\to \infty}f^n(c)=\infty$, there exists an integer $k \geq 1$ such that $q_{b_n}\leq q_c$ (where $b_n\in f^{-1}(\infty)$, defined in \eqref{3_2017_11_08}, is near $f^n(c)$, and $q_c$ was defined in (\ref{q_c})) and
\beq\lab{0042106}
 |f^n(c)|> \max\{1, 2\Dist_e(0, f(\Crit(f)))\}
\eeq
for all  $n \geq k$. Put  
\beq\lab{1_2017_11_11}
a:=f^k(c)
\eeq
(we may  need  in the  course of the proof $k\geq 1$ to be bigger). We shall prove the  following.

\sp\blem\lab{2042006} 
If $f:{\mathbb C}\to\ov{\mathbb C}$ is a compactly non-recurrent regular elliptic function and $c\in\Crit_{\infty}(f)$, then there exists a constant $C_1\geq 1$, such that
$$
m_{h,e}(B_e(a,r))\geq  C_1^{-1}r^{h}
$$
for all radii $r>0$ small enough, where $a$ is defined by \eqref{1_2017_11_11}. 
\elem

\bpf Put 
$$
q:=q_c.
$$
In view   of (\ref{0042106}) for every $n \geq 1$ there is a well--defined  holomorphic inverse branch
$$
f^{-1}_n: B_e\lt(f^n(a), \frac{1}{2}|f^n(a)|\rt)\lra {\mathbb C}
$$ 
of $f$ sending $f^n(a)$ to $f^{n-1}(a)$. Let $ b_n\in f^{-1}(\infty)$ be the unique
pole (assuming $k\neq 1$ is large  enough) such that 
$$ 
|f^n(a)-b_n|\leq \d(f^{-1}(\infty))\lek 1,
$$
where $\d(f^{-1}(\infty))$ comes from \eqref{4_2017_11_08}. Then, by Theorem~\ref{Euclid-I}
$$\begin{aligned}
f^{-1}_n\left(B_e\left(f^n(a),\frac{1}{4}|f^n(a)|\right)\right)&
\subset
B_e\left(f^{n-1}(a), \frac{K}{4}|f^n(a)||f'(f^{n-1}(a))|^{-1}\right)  \\
& \sbt B_e(f^{n-1}(a), C|f^n(a)||f^n(a)|^{- \frac{q+1}{q}})\\
& = B_e(f^{n-1}(a), C|f^n(a)|^{- \frac{1}{q}})\\
& \sbt B_e\left(f^{n-1}(a), \frac{1}{4}|f^{n-1}(a)|\right),
\end{aligned}
$$
where $C\in(0,+\infty)0$ is a constant and the last inclusion was written
assuming that 
$$
|f^{n-1}(a)|\geq 4c |f^n(a)|^{- \frac{1}{q}},
$$ 
which we can assume to hold for all  $n\geq k$ if $k$ is large enough.
So, the composition 
$$
f^{-n}_a:=f^{-1}_1\circ f^{-1}_2\circ\ldots
\circ f^{-1}_n: B_e\left(f^n(a), \frac{1}{4}|f^n(a)|\right)\lra
{\mathbb C},
$$ 
sending $f^n(a)$ to $a$, is well-defined and forms a
holomorphic branch of $f^{-n}$. Take  $0<r<8K/|a|$ and
let $n+1\geq 1$  be the least integer such that
$$ 
r|(f^{n+1})'(a)|\geq \frac{K}{8}|f^{n+1}(a)|. 
$$
Such integer  exists since $ |f'(z)|\comp
|f(z)|^{\frac{q_b+1}{q_b}}$ if $z$ is near a pole $b$. By its
definition $n \geq 0$ and, since $r<8K/|a|$, we have
$$  
r |(f^n)'(a)|< \frac{K}{8}|(f^n)(a)|.
$$
Then by Theorem~\ref{Euclid-I}, we have that
\beq\lab{11042106}
B_e(a,r)\supset f^{-n}_a\(B_e(f^n(a), K^{-1} |(f^n)'(a)|\).
\eeq

\sp Now we consider three  cases  determined   by the value  of
$r|(f^n)'(a)|$.

\sp\fr {\bf Case 1.}  $\d(f^{-1}(\infty))\leq r |(f^n)'(a)|<
\frac{K}{8}|f^n(a)|.$

\sp\fr In view of (\ref{0042106}) and Theorem~\ref{Euclid-I} along
with almost conformality of  the  measure $m_{h,e}$, we get that
\beq\lab{2042106}
\begin{aligned}
m_{h,e}(B(a, r))& \geq  K^{-h}|(f^n)'(a)|^{-h}  m_{h,e}(B_e(f^n(a), 4r|(f^n)'(a)|)) \\
  & \gek K^{-h}|(f^n)'(a)|^{-h}(4r|(f^n)'(a)|)^2\\
& \gek  |(f^n)'(a)|^{-h}(4r|(f^n)'(a)|)^{h}\\
 & = 4^h r^h,
\end{aligned}
\eeq
and we are done in this case.

\sp\fr {\bf  Case 2.}  $|f^n(a)-b_n|\leq 32
A^{\frac{q_{min}+1}{q_{min}}}r |(f^n)'(a)|<
{32}A^{\frac{q_{min}+1}{q_{min}}}\d(f^{-1}(\infty))$, {\it where $A$ was
defined  in (\ref{d5.3})}.

\sp\fr It then follows from Lemma~\ref{l1071702} that
$$  
m_{h,e}(B_e(f^n(a), K^{-1}|(f^n)'(a)|)) \gek (K^{-1} |(f^n)'(a)|)^{h}.
$$
Thus
$$\begin{aligned}
m_{h,e}(B_e(a,r)) \geq  K^{-h}|(f^n)'(a)|^{-h}(K^{-1} |(f^n)'(a)|)^{h}
\comp r^h.
\end{aligned}
$$
And we are done in this case.

\sp It remains to consider:

\sp\fr {\bf  Case 3.} \,\, $ r|(f^n)'(a)|<  \frac{1}{8}K
A^{-\frac{q_{\min}+1}{q_{\min}}}|f^n(a)-b_n|$.

\sp\fr But then
$$\begin{aligned}
r |(f^{n+1})'(a)| & =  r |(f^{n})'(a)|  |f'(f^{n}(a))|
 < \frac{K}{8}A^{-\frac{q_{\min}+1}{q_{\min}}}|f^n(a)-b_n|(A |f^{n+1}(z)|)^{\frac{q_n+1}{q_n}}\\
& \leq \frac{K}{8}A^{-\frac{q_{\min}+1}{q_{\min}}}A^{\frac{1}{q_n}+1} |f^{n+1}(a)|\\
& \leq \frac{K}{8}|f^{n+1}(a)|
\end{aligned}
$$
contrary to the definition of $n$. So this case is ruled out and
Lemma~\ref{2042006} is proved. \endpf

\sp We are ready to prove the following.

\sp\bthm\lab{t1042206} Let $f:\C\lra \oc$ be a compactly non-recurrent regular elliptic function. If $\Om(f)=\es$ then the $h$-dimensional
packing measure $\Pi^h_e$ of  every bounded  Borel subset of $J(f)$
is finite and $\Pi^h_s(J(f))<+\infty$. \ethm

\bpf Consider arbitrary  point
$$
z\in J(f)\sms\bu_{n=0}^{\infty} f^{-n}(\{\infty\}\cup \Crit(f)) 
$$ 
and a radius $r\in (0, \g(a \xi)^{-1} )$, where $\g>0$ was defined in (\ref{dp5.4}) while $a$ and $\xi$ come from Lemma~\ref{l1022506}. Let $s\ge 0$ be associated to the point $z$ and the radius $r/\xi$ according to Lemma~\ref{l1022506}. If
the case (a) of this lemma holds, then we get from Lemma~\ref{lncp12.9.} and Lemma~\ref{lncp13.1.} that
\beq\lab{5042006}
\begin{aligned}
m_{h,e}(B_e(z,r))
&\ge K^{-h}|(f^s)'(z)|^{-h} m_{h,e}\(B_e\(f^s(z), K^{-1}r |(f^s)'(z)|\)\)\\
&  \gek K^{-h}|(f^s)'(z)|^{-h} ( K^{-1}r|(f^s)'(z)|)^2 \\
&\comp  r^h(r|(f^s)'(z)|)^{2-h}\\
&\gek r^h.
 \end{aligned}
\eeq
If the case (b) of  Lemma~\ref{l1022506} holds, then applying  this lemma along
with Lemma~\ref{l2022506} (with $\sigma\leq K^{-1}\xi)$, we get that
$$
\begin{aligned}
m_{h,e}(B_e(z,r))
&\ge K^{-h}|(f^s)'(z)|^{-h} m_{h,e}(B_e(f^s(z), K^{-1}r |(f^s)'(z)|))\\
&  \gek K^{-4}|(f^s)'(z)|^{-h} (K^{-1}r|(f^s)'(z)|)^h \\
&\comp r^h.
 \end{aligned}
$$
Combining this and (\ref{5042006}), completes the proof of the first
part because of Theorem~\ref{tncp12.2.} (a). Since
$\Pi^h_e(A)=\Pi^h_e(A+w)$ for every $\om \in \La$ and since
$\frac{d\Pi^h_s}{d \Pi^h_e}(z)=(1+|z|^2)^{-h}$, we  get with $R=4
\diam ({\mathcal R})$, where $\mathcal R$ is an arbitrary fixed fundamental parallelogram,  that
$$\begin{aligned}
\Pi^h_s(J(f))=&\sum_{k=0}^\infty \Pi^h_s(A(0; 2^kR,2^{k+1}R))+
\Pi^h_s(B_e(0,R))\\
 &  \lek \Pi^h_e(B_e(0,R))+ \sum_{k=0}^\infty \Pi^h_e(A(0;2^kR,2^{k+1}R))R^{-2h}4^{-hk}\\
 & \lek \Pi^h_e(B_e(0,R))+ \sum_{k=0}^\infty (2^kR)^2 R^{-2h}4^{-hk}\\
&=\Pi^h_e(B_e(0,R)) + R^{(2(1-h)}\sum_{k=0}^\infty 4^{(1-h)k} <
 +\infty,
\end{aligned}
$$
since $h> 1$. We are done. \endpf

\sp \bprop\label{p2J39} Let $f:\C\lra \oc$ be a compactly non--recurrent regular elliptic function. If $\HD(J(f))=2$, then $J(f)=\mathbb C$.\eprop

\bpf Since $\Pi^2_e$ and  $S$, two-dimensional
Lebesgue measure on $\mathbb C$, coincide up to a  multiplicative
constant, it follows from  (already proved) Theorem~\ref{t1020907}
(b) that if $h=2$, then $S(J(f))>0$. So, in order  to prove   our
proposition it suffices to show that if $J(f)\varsubsetneq \mathbb
C$, then $S(J(f))=0$. So, suppose that $J(f)\neq  \mathbb C$. We are to
show that
$$
S(J(f)\sms \Sing^{-}(f)))=0.
$$ 
Let for any integer $l\geq 1$ the set $Z_l$ have exactly the same meaning as in the proof of Theorem~\ref{tmaincm}. Let $\La\sbt\C$ be a lattice associated to the elliptic function $f$. Since $J(f)$ is a $\La$--invariant nowhere dense subset of
$\mathbb C$, there exists $\varepsilon>0$ such that  for every  $y
\in \C$ there exists $y_\varepsilon\in B_e(y,\frac{1}{2l})$ such that
\beq\label{1J39}
 B_e(y_\varepsilon, \varepsilon)  \sbt B_e\left(y,\frac{1}{2l}\right)
\sms J(f).
\eeq
Keep notation from the proof of Theorem~\ref{tmaincm}. Fix arbitrary point $z \in Z_l$. By Theorem~\ref{Euclid-I}, Koebe's $\frac{1}{4}$-Theorem and (\ref{1J39}), we have
$$\begin{aligned}
f^{-n_k}_z(B_e(f^{n_k(z)}(z), \varepsilon))& \sbt
f^{-n_k}_z\(B_e\(f^{n_k(z)}(z),(2l)^{-1}\)\sms J(f)\)\\
 & \sbt B_e\(z,K|(f^{n_k})'(z)|^{-1} (2l)^{-1}\)\sms J(f) \end{aligned}
$$ and
$$
f^{-n_k}_z(B_e(f^{n_k}(z)_{\varepsilon}, \varepsilon)) 
\supset B_e\lt(f^{-n_k}_z(f^{n_k}(z)_{\varepsilon}),\frac{1}{4}\varepsilon|(f^{n_k(z)})'(z)|^{-1}\rt).
$$
Therefore, we see that
$$
\frac{S\(B_e\(z,K|(f^{n_k})'(z)|^{-1}(2l)^{-1}\)\sms
J(f)\)}{S\(B_e\(z,K|(f^{n_k})'(z)|^{-1}(2l)^{-1}\)\)}
\geq\left(\frac{\varepsilon l}{2K}\right)^2>0.
$$ 
So, $z$ is not a
Lebesgue's density point for the set $Z_l$, and therefore $S(Z_l)=0$. Hence 
$$
S(J(f))
=S\(J(f)\sms \Sing^{-}(f)\)
=S\lt(\bu_{l=1}^{\infty}Z_l\rt)
=\sum_{l=1}^{\infty}S(Z_l)
=0.
$$
The proof of Theorem~\ref{t1042206} is complete. \endpf

\sp Theorem~\ref{t1020907}  is now a logical  consequence of
Section~\ref{Hausdorff}, Proposition~\ref{ncpp17.5},
Lemma~\ref{ncp1l7.7.} and Theorem~\ref{t1042206}.

\sp\chapter[Conformal Invariant Measures for c.n.r.r. functions]{Conformal Invariant Measures for Compactly Non--Recurrent Regular Elliptic Functions}\label{invariant-p.s.n.r.}

Throughout this whole chapter $f:\C\lra\oc$ is assumed to be a compactly non--recurrent regular elliptic function. This chapter is in a sense a core of our book. Taking fruits of what has been done in all previous chapters,  we prove in it the existence and uniqueness, up to a multiplicative constant, of a $\sg$--finite $f$--invariant measure $\mu_h$ equivalent to the $h$-- conformal measure $m_h$ which in turn was proven to exist in Theorem~\ref{tmaincm}. Furthermore, still heavily based on what has been done in all previous chapters, particularly nice sets, first return map techniques, Young's towers, we provide here a systematic account of ergodic and refined stochastic properties of the dynamical system $(f,\mu_h)$ generated by such subclasses of compactly non--recurrent regular elliptic functions as normal subexanding elliptic functions of finite character and parabolic elliptic functions. 

\section[Conformal Invariant Measures for CNR Regular Elliptic Functions] {Conformal Invariant Measures for Compactly Non--Recurrent Regular Elliptic Functions: the Existence, Uniqueness, Ergodicity/Conservativity, and Points of Finite Condensation}\label{conformal-invariant}

In this section we deal with $\sg$--finite invariant measures
equivalent to the conformal measure $m_h$ proven to exist in Theorem~\ref{tmaincm}. We prove their existence, uniqueness up to a multiplicative constant, metric exactness, implying ergodicity, and  conservativity. We also study at length the points of their finite and infinite condensation giving the first outlook of the location of such points. In the context of rational functions the results of such kind were obtained for example in \cite{ADU}, \cite{DU LMS}, \cite{DU Forum}, \cite{U2}. In the context of quite general hyperbolic/expanding transcendental meromorphic and entire functions see for example  \cite{UZ2},  \cite{MyU1}, and  \cite{MyU1}. 

\sp The first result of this chapter is the following.

\sp \bthm\lab{tinv} 
If $f:\C\lra\oc$ is a compactly non--recurrent regular elliptic function, then there exists a unique, up to a multiplicative constant, $\sg$--finite $f$--invariant measure $\mu_h$ that is absolutely continuous with respect to the $h$--conformal measure $m_h$. In addition, 

\, \begin{enumerate}

\item $\mu_h$ is equivalent to $m_h$, 

\, \item $\mu_h$ is metrically exact, hence ergodic and conservative, and 

\, \item $\mu_h$ is given, according to Theorem~\ref{t1h75}, by formulas (\ref{5.9a})--(\ref{eq:muequ})

\, \item $\mu_h\(J(f)\sms \Tr(f)\)=0$ . 
\end{enumerate}
\ethm

\bpf Let $\xi\in {\mathbb C}$ be a periodic point of $f$
with some period $p\ge 3$. We put
$$
P_3(f):=\ov{O_+\(f(J(f)\cap\Crit(f))\)} \cup\{\xi,f(\xi),\ld,f^{p-1}(\xi)\}.
$$
Since, by Proposition~\ref{p120190914}, $\ov{O_+\(f(J(f)\cap\Crit(f)\)}$ is a forward--invariant nowhere-dense subset of $J(f)$ and since the $h$--conformal measure $m_{h,s}$ is
positive on nonempty open subsets of $J(f)$, it follows from
ergodicity and conservativity of $m_{h,s}$ (see Theorem~\ref{tmaincm})
that
$$
m_{h,s}(\ov{O_+(f(\Crit(f)))})=0.
$$
Since $m_{h,s}$ has no atoms (see Theorem~\ref{tmaincm}) we therefore
obtain that
\beq\label{1p115}
m_{h,s}(P_3(f))=0.
\eeq
Aiming to apply Theorem~\ref{t1h75},
we shall now construct a Marco Martens cover of Definition~\ref{d:mmmap} for the map $f:\mathbb J(f)\to J(f)\cup\{\infty\}$. For further purposes, our construction will be more involved and more specific than the one which would be needed just
for the sake of having a Marco Martens cover. Fix an integer $q\ge 1$ and for every integer $k\geq 1$, set
\beq\lab{1_2017_11_15}
Z_k=Z_k(q):=\lt(\ov{B}_e(0,qk)\sms B\lt(P_3(f),\frac{1}{qk}\rt)\rt)\cap J(f).
\eeq
Obviously each set $Z_k$ is compact. We shall inductively construct the sets
$$
(\hat{X}_n)_{n\geq0}
$$ 
as follows. 

First, by compactness of $Z_1$, we
can cover it by finitely many sets $\hat{X}_0, \hat{X}_1, \ldots,
\hat{X}_{n_1},$ each of which is an open ball centered at a point
$Z_1$ with radius equal to $1/4$.

For the inductive step fix $j\geq 1$  and suppose that the open balls 
$$
\hat{X_0}, \ldots,
\hat{X}_{n_1}, \hat{X}_{n_1+1},\ldots, \hat{X}_{n_1+n_2},\ldots,
\hat{X}_{n_1+\ldots +n_j}
$$ 
have been constructed with the following properties.
\begin{itemize}
  \item [(a)] $\forall_{ 0\leq i  \leq j}$  $Z_i \sbt
  \hat{X}_{n_{i-1}+1}\cup
  \hat{X}_{n_{i-1+2}} \cup  \ldots \cup  \hat{X}_{n_{i}}$.
 \, \item  [(b)] $\forall_{ 0\leq i  \leq j}$ the sets  $\hat{X}_{n_{i-1}+1}, \ldots,
  \hat{X}_{n_i}$ are open balls centered at points $Z_i$ with radii all
  equal to $4^{-i}$.
\,   \item [(c)] $\forall_{0 \leq i < k \leq j}$   $Z_i\cap
   (\hat{X}_{n_{k-1}+1}\cup\ldots \cup {X}_{n_k})=\es$.
\end{itemize}
We now look   at the family \index{(S)}{$\mathcal{B}_{j+1}$}
$\mathcal{B}_{j+1}$  of all open balls centered  with  points at
$Z_j$  and radii equal to $4^{-(j+1)}$. Since  $Z_{j+1}$  is a
compact  set we  can find \index{(S)}{$\mathcal{B}_{j+1}^1$}
$\mathcal{B}_{j+1}^1$, a finite subfamily  of $\mathcal{B}_{j+1}$
which covers $Z_{j+1}$. Observe that if $b \in \mathcal{B}_{j+1}^1$
and if $B\cap Z_i\neq \es$ for some $0\leq i \leq j$, then there
exists $n_{i-1}+1 \leq k \leq n_i$ such that $B \sbt \hat{X}_i$. So,
if  we  remove  from ${\mathcal B}_j^1$ all such  sets, and we label
them as $\hat{X}_{n_{j-1}+1}, \hat{X}_{n_{j-1}+2}, \ldots,
\hat{X}_{n_j}$, then the  conditions (a), (b), and (c) will be satisfied
with $j$ replaced by $j+1$. Our inductive construction of the
sequence $(\hat{X}_n)_{n=0}^\infty$ is thus complete. 

We set,
\beq\lab{120181221}
X_n:= \hat{X}_n \cap J(f), \quad n \geq 0.
\eeq
Condition (1) of Definition~\ref{d:mmmap} is  obviously satisfied
and condition (2) follows from (\ref{1p115}) and since
$\bigcup_{n=0}^\infty \hat{X}_n\supset \bigcup_{n=0}^\infty
Z_k=J(f)\sms P_3(f)$.  In order to demonstrate  condition (3)  it
suffices to note (remember that $f$  is an elliptic function) that
if $z\in J(f)$  and  $r>0$, then there exists $n\ge 1$ such that
$$
f^n (B_e(z,r))=\oc.
$$

Let us check condition (4). Indeed, by our construction we have that
$$
2\hat{X}_n\sbt \mathbb C \sms P_3(f)\sbt \mathbb C \sms \PC(f).
$$
Hence, all holomorphic inverse branches of $f^{-n}$ for all $n\geq
1$ are well-defined on $2X_n$, and in view of Theorem~\ref{Spher-I},
there exists a constant $\hat{K}\ge 1$ such that if $i \geq 0$ and
$$
f_*^{-n}: X_i\lra {\mathbb C}
$$ 
is a holomorphic
branch of $f^{-n}$, then for  all $x,y\in X_i$, we have
\beq\lab{1.dist}
{|(f_*^{-n})^*(y)|\over |(f_*^{-n})^*(x)|}\le \hat{K}.
\eeq
We therefore obtain for all Borel sets $A,B\sbt X_i$ with
$m_{h,s}(B)>0$, and all $n\ge 0$, that
$$
{m_{h,s}(f_*^{-n}(A))\over m_s(f_*^{-n}(B))} 
={\int_A|(f_*^{-n})^*|^hdm_{h,s}\over \int_A|(f_*^{-n})^*|^hdm_{h,s}}
\le {\sup_{A_k}\{|(f_*^{-n})^*|^h\}m_{h,s}(A)\over
inf_{A_k}\{|(f_*^{-n})^*|^h\}m_{h,s}(B)}
\le \hat{K}^h{m_{h,s}(A)\over m_{h,s}(B)}.
$$
Hence 
$$
\begin{aligned} 
m_{h,s}(f^{-n}(A))=&\sum_{*}m_{h,s}(f_*^{-n}(A))\leq
\hat{K}^h\frac{m_{h,s}(A)}{m_{h,s}(A)}\sum_{*}m_{h,s}(f_*^{-n}(B))\\
=&\hat{K}^h\frac{m_{h,s}(A)}{m_{h,s}(B)}m_{h,s}(f^{-n}(B)). 
\end{aligned}
$$ 
So, condition (4) of Definition~\ref{d:mmmap} is satisfied, and we have
shown that $(X_i)_{i=0}^\infty$ is a Marco Martens  cover. Since, item (2) of Theorem~\ref{tmaincm}, the map $f:(J(f),m_{h,s})\lra (J(f),m_{h,s})$ is ergodic and conservative, we see that Theorem~\ref{t1h75} applies, and its application yields items (1) and (3). Item (2) of our current theorem also follows from item (2) of Theorem~\ref{tmaincm}, this time from weak metrical exactness of $m_{h,s}$, from Proposition~\ref{p220190918} and Proposition~\ref{p120190918}. Item (4) follows now from item (4) of Theorem~\ref{tmaincm}.  The proof is complete. 
\endpf

\sp We want to record one particular property of the cover
$(X_n)_{n=0}^\infty$, which follows immediately from our construction:

\blem\label{l1j137} 
If $f:\C\lra\oc$ is a compactly non--recurrent regular elliptic function, then for every compact set $F \sbt \mathbb C \sms
P_3(f)$, the  set $\{n\geq 0; \, \,  X_n \cap F\neq \es \}$  is
finite. 
\elem

\fr The following lemma was established in Theorem~\ref{t1h75}.

\blem\lab{lAnfin} 
If $f:\C\lra\oc$ is a compactly non--recurrent regular elliptic function, then for every $n\ge 0$ we have 
$$
0<\mu_h(X_n)<\infty.
$$
In particular, $\mu_h \in {\mathcal M}^\infty_f$, \index{(S)}{${\mathcal M}^\infty_f$} meaning that it is a $\sg$--finite $f$--invariant measure.. 
\elem

 The following lemma immediately follows from results of Section~\ref{abstract}.

\blem\lab{ncp2l5.1.}
If $f:\C\lra\oc$ is a compactly non--recurrent regular elliptic function, then $J_{\mu_h}(\infty):=J_\mu(f)(\infty)$\index{(S)}{$J_{\mu_h}(\infty)$}, the
set of points of infinite condensation of the measure $\mu_h$, is
\begin{enumerate}

\, \item closed, 

\, \item $f(J_{\mu_h}(\infty))\sbt J_{\mu_h}(\infty)$,  and

\, \item $\mu(J_{\mu_h}(\infty))=0$. 
\end{enumerate}
\elem

\blem\lab{ncp2l5.2.} 
If $f:\C\lra\oc$ is a compactly non--recurrent regular elliptic function, then 
$$
J_{\mu_h}(\infty)\sbt\ov{O_+(\Crit(f))}\cup\{\infty\}. 
$$
\elem

\bpf Since  all the sets $X_n$, $n\geq 0$, are open, it
follows from Lemma~\ref{lAnfin} that $J_{\mu_h}(\infty) \sbt P_3(f)\cup\{\infty\}$.
Replacing $\xi, \ldots, f^{p-1}(\xi)$ by another periodic orbit with
period $\geq 3$ and using uniqueness  of the measure  $\mu_h$
(Theorem~\ref{tinv}), we thus conclude that
$J_{\mu_h}(\infty)\sbt\ov{O_+(\Crit(f))}\cup\{\infty\}$. The proof is
finished. 
\endpf

\sp Now, we further investigate in greater detail the structure of the set $J_{\mu_h}(\infty)$ of points of infinite condensation of the measure $\mu_h$. Fix  $$
w \in J(f)\sms B_e(\Om(f),\th(f)),
$$
where $\th(f)>0$ comes from \eqref{d5.2}, and an open Jordan domain 
$$
Q\sbt B_e(w, 2\g_f),
$$
where $\g_f$ is defined by \eqref{dp5.4}. A sequence $\{Q_n\}_{n=0}^\infty$ of connected components of inverse images of $f^{-n}(Q)$, $n\ge 1$, is called $w$--nested\index{(N)}{$w$-nested sequence} if and only if
$$
f(Q_{n+1})=Q_n
$$ 
for all integers $n\geq 0$. The following lemma follows immediately from Lemma~\ref{lmane2}, the definition of $\b_f$, and \eqref{dp5.4}.  

\sp
\blem\lab{l1050106} 
Let $f:\C\lra\oc$ be a non--recurrent elliptic function. Let $w \in  J(f)\sms B_e(\Om(f),\th(f))$ and let $Q\sbt B_e(w,\g_f)$ be a Jordan domain. If $\{Q_n\}_{n=0}^\infty$ is a $w$--nested  sequence  of connected
components of the sets $f^{-n}(Q)$, then all the sets $Q_n, n\geq 0$, are open Jordan domains.
\elem

\sp Let $\Crit_h(f)$\index{(S)}{$\Crit_h(f)$} be the set of all
critical points of $f$ that are $h$-upper estimable with respect to the $h$--conformal measure $m_h$. Because of Lemma~\ref{l5042706a}, we have this.
\beq\lab{1_2017_11_13}
\Crit_c(f)\sbt \Crit_h(f).
\eeq
The key fact for what will follow in this section is this.

\sp\blem\lab{l2050106} 
Let $f:\C\lra\oc$ be a compactly non--recurrent regular elliptic function. Let $w \in  J(f)\sms B_e(\Om(f), \th)$, let $Q\sbt B_e(w, \g)$  be a Jordan domain, and let $\{Q_n\}_{n=0}^\infty$  be a $w$--nested sequence of connected
components of the sets $f^{-n}(Q)$. For every $n\geq 0$  let $W_n$
be the connected component of $f^{-n}(B_e(w, 2 \g))$ containing
$Q_n$. If
$$ 
(\Crit(f)\sms \Crit_h(f))\cap \bu_{n=1}^\infty W_n =\es,
$$
then 
$$
m_{h,e}(Q_n) \lek \frac{m_{h,e}(Q)}{\diam_e^h(Q)} \diam^h_e(Q_n)
$$ for
all $n\geq 0$. 
\elem

\bpf For $n=0$ the required inequality is trivial. Fix
$k\geq 0$ and  $n\geq 0$. Suppose that $W_{k+n}$ contains no
critical points of $f^n$. It then follows from
Theorem~\ref{Euclid-II} and (\ref{1050206}) (note that also
$Q\sbt(B_e(w, \g)$) that
\beq\lab{3050106}
\begin{aligned}
m_{h,e}(Q_{k+n})
& \leq  \sup\big\{|(f^n)'(z)|^{-h}:\,\,  z \in Q_{k+n}\big\} m_{h,e}(Q_{k})\\
&\leq  K^h_* \inf\big\{|(f^n)'(z)|^{-h}:\,\,  z \in Q_{k+n}\big\} m_{h,e} (Q_{k})\\
     &   \leq K^h_* \frac{\diam^h_e(Q_{k+n})}{\diam^h_e(Q_{k})} m_{h,e} (Q_{k})\\
     &  = K^h_* \frac{m_{h,e} (Q_{k})}{\diam^h_e(Q_k)} \diam^h_e(Q_{k+n}),
\end{aligned}
\eeq
with an appropriate  universal constant $K^h_*\geq 1$. Now suppose
that $W_{k+1}$ contains a critical point $c$ of $f$.
 By (\ref{dp5.4}) and Lemma~\ref{ld5.1}, $c$ is the only critical point of  $f$ in $W_{k+1}$.
 Suppose first that
$$\dist_e(f(c), Q_k)\geq  4 \diam_e( Q_k).$$ Fix  $z\in Q_k$. Then $Q_k
\sbt  B_e(z, \diam_e(Q_k))$, $$ f(\Crit(f))\cap B_e(z, 2
\diam(Q_k))=\es$$ (assuming that $\g, \eta> 0$ sufficiently small),
which  makes other (finitely many) critical values  lying
sufficiently  far apart from $f(c)$. Hence denoting by 
$$
f^{-1}_*:B_e(z,2 \g \diam_e(Q))\lra {\mathbb C}
$$ 
the holomorphic inverse branches of $f$ whose range covers
$Q_k$, using Theorem~\ref{Euclid-II}, we estimate similarly as above
\beq\lab{2050106}
\begin{aligned}
m_{h,e}(Q_{k+1}) 
& \leq  \sup \big\{|(f^{-1}_*)'(x)|^h:\, x \in Q_{k+1}\big\} m_{h,e} (Q_{k})\\
&\leq  K^h\inf\big \{|(f^{-1}_*)'(x)|^h:\, x \in Q_{k+1}\big\} m_{h,e} (Q_{k})\\
     &   \leq K^h \frac{\diam^h_e(Q_{k+1})}{\diam^h_e(Q_k)}m_{h,e} (Q_{k})\\
     &   \leq K^h \frac{m_{h,e}(Q_k)}{\diam^h_e(Q_k)} \diam^h_e(Q_{k+1}).
\end{aligned}
\eeq
 Now, assume  that
\beq\lab{10050106}
\dist_e(f(c), Q_k)\leq 4 \diam_e (Q_k).
\eeq
We thus get  that $Q_k \sbt  B_e(f(c), 5\diam_e(Q_k))$, and
therefore
$$
Q_{k+1} \sbt B_e\(c, A(5 \diam_e(Q_k))^{1/p_c}\).
$$
Hence, making use of $h$-upper estimability of the point $c$, we get
that
\beq\lab{1050106}
 \begin{aligned}
 m_e(Q_{k+1}) \leq L\(A(5\diam_e(Q_k))^{1/p_c}\)^h.
\end{aligned}
\eeq
It follows from  (\ref{10050106})  that
$$\begin{aligned}
\Dist_e(c, Q_{k+1})
&\leq  A(\Dist_e(f(c), Q_k))^{1/p_c}\\
& \leq A(\dist_e(f(c), Q_k)+ \diam_e(Q_k))^{1/p_c}\\
&q\le A(5\diam_e(Q_k))^{1/p_c}.
\end{aligned}
$$
Therefore
$$\begin{aligned}
\diam_e (Q_k)& \leq \diam_e(Q_{k+1}) A (\Dist_e(c,Q_{k+1}))^{p_c-1}\\
             & \leq A^25^{1/p_c}\diam_e (Q_{k+1}) \diam^{\frac{p_c-1}{p_c}}_e(Q_k).
\end{aligned}
$$
Thus    $$\diam^{1/p_c}_e(Q_k) \leq A^2 5^{1/p_c}\diam_e
(Q_{k+1}).$$ Inserting this to (\ref{1050106}), we get that
$$m_e(Q_{k+1})\leq L(25)^{h/p_c}A^{3h/p_c}\diam^h_e(Q_{k+1}).$$
Applying this,  (\ref{3050106}), (\ref{2050106}) and making use of
Lemma~\ref{ld5.1}  along with (\ref{dp5.4}), a straightforward
inductive argument  yields  that for  every $j\geq 1$
$$ 
m_e(Q_j)\leq\max\big\{(L(25 A^3)^{1/p_c},K,K_*)^h\big\}^{N_f}\diam^h_e(Q_j).
$$
The proof is complete. 
\endpf

\sp\fr As an immediate consequence  of this lemma, Lemma~\ref{lncp23.3}
and Theorem~\ref{Euclid-II} we obtain the following.

\sp\blem\lab{l1050206} 
Let $f:\C\lra\oc$ be a compactly non--recurrent regular elliptic function. Having $\th\in(0,\th(f)$ and $\g\in(0,\g_f)$, let $w \in J(f)\sms B_e(\Om(f),\th)$, let $V\sbt B_e(w, \g)$  be a Jordan domain, and let  $U$ be a Jordan domain   contained in $V$.  Let $\{V_n \}_{n=0}^{\infty}$ be  a $w$--nested  sequence of connected  components  of $f^{-n}(V)$ and let $\{U_n\}_{n=0}^{\infty}$, with $U_n\sbt V_n$, be a $w$--nested  sequence of connected components of $f^{-n}(U)$.
 For every  $n \geq  0 $ let $W_n$ be the connected  component of
$f^{-n}(B_e(w, 2\g))$  containing $V_n$. Suppose  that
$$ 
(\Crit(f)\sms \Crit_h(f)) \cap \bu_{n=1}^{\infty}W_n=\es
$$
and that  there exists a Jordan domain $\^U$ such that $\ov{U} \sbt
\^U \sbt V$ and $\^U \cap {\rm PC}(f) =\es $. Then
$$
m_{h,e}(V_n) \lek \frac{m_{h,e}(V)}{m_{h,h}(U)}m_{h,e}(U_n),
$$
and the same inequality remains true (perhaps with  a smaller
constant on the right-hand side) with $m_{h,e}$ replaced  by $m_{h,s}$
 since  the diameters  of all the  sets $V_n$ are bounded above by $\b_f$.
\elem 

\sp\fr Now we can take the first fruit of this lemma.

\sp\bprop\lab{p2050206} 
All the  points of the set $\ov{{\rm
PC}_c^0(f)}\sms \Om(f)$  are of finite condensation with respect to
the invariant measure $\mu_h$. 
\eprop

\bpf Keep the Marco Martens cover $(X_n)_{n=0}^\infty$ as defined by \eqref{120181221}. We keep working within the framework of 
Theorem~\ref{t1h75} and Definition~\ref{d:mmmap} set up in the proof of Theorem~\ref{tinv}. Take an arbitrary point  $w \in \ov{{\rm
PC}_c^0(f)}\sms \Om(f)$. Assuming $\th>0 $  to be small enough, we
will  have $w \notin B_e(\Om(f), \th)$. Fix 
$$
V \sbt B_e(w, \g),
$$
an open centered at $w$ and disjoint from ${\rm PC}_p(f)\cup
{\rm PC}_\infty(f)$. Since, by Proposition~\ref{p120190914}, $\ov{{\rm PC}_c^0(f)}$  is a nowhere
dense subset of $J(f)$, by taking sufficiently large $q\ge 1$ in \eqref{1_2017_11_15}, we may assume without loss of generality
that $2X_0\sbt V$. Invoking \eqref{1_2017_11_13}, it immediately follows from  Lemma~\ref{l1050206} (note that $2X_0\cap \PC(f)=\es$)  that
$$
m_{h,s}(f^{-n}(X_0)\cap V_n)\gek m_{h,s}(V_n)
$$ 
for every integer $n\geq 0$, where $V_n$ is a connected component  of
$f^{-n}(V)$. Therefore, summing up over all  connected  components
$V_n$ of $f^{-n}(V)$, we obtain
$$
m_{h,s}(f^{-n}(X_0)) =m_{h,s}(f^{-n}(A_0)\cap f^{-n}(V))\gek m_{h,s}(f^{-n}(V)).
$$
Consequently
$$
m_{h,s,n}(V)=\frac{\sum_{n=0}^k m_{h,s}( f^{-n}(V))}{\sum_{n=0}^k m_{h,s}
(f^{-n}(X_0))}\lek 1,
$$
where $m_{h,s,n}$ are given by formula \eqref{5.9a}, and it follows from the  formula (\ref{eq:muequ}) of Theorem~\ref{t1h75} that $\mu_h(V)<+\infty$, which finishes the proof. \endpf 

\sp\blem\lab{l3050206} 
If $-nested$ is a compactly non--recurrent regular elliptic function, then all  the points of the set ${\rm PC}^0_p(f)\cup {\rm PC}^0_\infty(f)$ are of finite condensation with respect to the $f$--invariant measure $\mu_h$. \elem

\bpf We still keep the Marco Martens cover $(X_n)_{n=0}^\infty$ as defined by \eqref{120181221} and are working within the framework of 
Theorem~\ref{t1h75} and Definition~\ref{d:mmmap} set up in the proof of Theorem~\ref{tinv}. Fix a point $w\in {\rm PC}^0_p(f) \cup {\rm
PC}^0_\infty(f)$. There exists an integer $j\geq 0$  so large  that
$$
f^{-j}(w)\cap \(\ov{O_+( \Crit(f))}\cup \Om(f)\)=\es.
$$
Therefore, taking into account 
Corollary~\ref{General Riemann Hurwitz Simply Connected_3}, we see that with $\th>0$ and $\g>0$ small enough, there exists an open disk $V$ centered at $w$ with  the following properties:
\begin{itemize}
\item [(a)] For every $z\in f^{-j}(w)$, $\dist_e(z, \Om(f))>\th$.

\,\item [(b)] For every $z\in f^{-j}(w)$,  if $V_z$ is the connected component of $f^{-j}(V)$
 containing $z$, then $V_z$ is a Jordan domain and $V_z\sbt B_e(z,\g)$.

\, \item [(c)] $\bu_{z\in f^{-j}(w)}B_e(z,2\g)\cap \ov{{\rm PC}(f)}=\es $.
\end{itemize}
\fr We may assume  without loss of generality that $2X_0\sbt V$.  It
follows from  condition (c) that
$$
\Crit(f)\cap \bu_{z\in f^{-j}(w)}\bu_{n=0}^{\infty}f^{-n}(B_e(z,2\g))=\es.
$$
\fr So, we may apply  Lemma~\ref{l1050206} to the pairs $(U_z,V_z)$,
$z\in f^{-j}(w)$,  where $U_z$ are  the connected  components of
$f^{-j}(X_0)$ contained in $V_z$, to get,  similarly as in the proof
of Proposition~\ref{p2050206}, that  for every $z\in f^{-j}(w)$ and
every $n\geq 0$,
$$
\sum_{i=0}^n m_{h,s}( f^{-i}(V_z))\lek  \sum_{i=0}^n m_{h,s}( f^{-i}(U_z)). 
$$
Summing over all $z\in f^{-j}(w)$, we thus get
$$
\sum_{i=j}^{j+n} m_{h,s}( f^{-i}(V))
= \sum_{i=0}^{n} m_{h,s}( f^{-i}(f^{-j}(V)))
  \lek  \sum_{i=j}^{j+n} m_{h,s}( f^{-i}(X_0)).
$$
Since in addition both $\sum_{i=0}^{j-1} m_{h,s}( f^{-i}(X_0))$ and
$\sum_{i=0}^{j-1} m_{h,s}( f^{-i}(V))$ are finite, we thus get that
$m_{h,s,n}(V) \lek 1$, where $m_{h,s,n}$ are given by formula \eqref{5.9a}, for all $n \geq 0$. It therefore follows from
formula (\ref{eq:muequ}) quoted in Theorem~\ref{t1h75} that
$\mu_h(V)<+\infty$. The proof is complete. \endpf

\sp\fr As an immediate consequence of this  lemma,
Proposition~\ref{p2050206} and Lemma~\ref{ncp2l5.2.} we get the
following.

\bthm\lab{t4050206} 
If $f:\C\lra\oc$ is a compactly non--recurrent regular elliptic function, then 
$$
J_{\mu_h}(\infty) \sbt \Om(f)\cup \{\infty\}.
$$
\ethm

\sp Now we shall deal with the point $\infty$. Recall that the set $\Crit_h(f)$ was defined just above formula \eqref{1_2017_11_13} We shall prove the following. 

\sp\bprop\lab{p5050206} If $f:\C\lra\oc$ is a compactly non--recurrent regular elliptic function and $\Crit_{\infty}(f)\sbt \Crit_h(f)$, \index{(S)}{$\Crit_h(f)$} then $\infty$ is a point of  finite condensation of the $f$--invariant measure
$\mu_h$. 
\eprop

\bpf As above, we keep the Marco Martens cover $(X_n)_{n=0}^\infty$ as defined by \eqref{120181221} and are working within the framework of Theorem~\ref{t1h75} and Definition~\ref{d:mmmap} set up in the proof of Theorem~\ref{tinv}. Recall that $\g_f>0$ $\g_f$ is determined by \eqref{dp5.4}, while $\th(f)>0$ was defined in \eqref{d5.2}. Fix $\g\in(0,\g_f))$ and $\th\in (0,\th(f)$ so small that
\beq\lab{2_2017_11_15}
\Big(B_e\(\Om(f),\th(f)\)\cup \ov{{\rm PC}^0_c(f)}\Big)\cap\, \bu_{b\in f^{-1}(\infty)}\!\!\!\! B(b,2\g)=\es
\eeq
and 
\beq\lab{3_2017_11_15B}
\Big(\bu_{c\in\Crit_p(f)}O_+(c)\Big)\cap\Big(\bu_{b\in f^{-1}(\infty)\sms \bu_{c\in\Crit_p(f)}O_+(c)}B(b,2\g)\Big)=\es,
\eeq
Let $R>0$ be so large that 
\beq\lab{3_2017_11_15}
B_b(R)\sbt B(b,\g)
\eeq
for every $b\in f^{-1}(b)$, where, we recall, $B_b(R)$ 
\index{(S)}{$B_\infty(R)$} is the connected component of $f^{-1}(B_\infty(R))$ containing $b$. By taking sufficiently large $q\ge 1$ in the formula \eqref{1_2017_11_15}, we may assume without loss of generality that 
$$
2X_0 \sbt B_\infty(R).
$$
Of course,
$$
\xi:= \inf\{\diam_e (B_b(R)): b \in f^{-1}(\infty)\}>0.
$$
For every $b\in  f^{-1}(\infty)$ let $X_0^b$  be a connected component of $f^{-1}(X_0)$ contained in $B_b(R)$ and let $2X_0^b$ be the connected component of $f^{-1}(2X_0)$  containing $X_0^b$. Since the set 
$$
P_c:= f^{-1}(\infty)\cap\bu_{c\in\Crit_p(f)}O_+(c)
$$ 
is finite, proceeding as in the  second part of the proof of
Lemma~\ref{l3050206}, with the pair $(X_0,V)$ replaced  by
$(X_0^b,B_b(R))$, we see that there exists $j\geq 0$ (with analogous
meaning as in the proof of Lemma~\ref{l3050206}) such that for all
$n\geq 0$
\beq\lab{7050106}
\sum_{b\in P_c} \sum_{i=j}^{j+n} m_{h,s}(f^{-i}(B_b(R)))
\leq \sum_{b\in P_c} \sum_{i=j}^{n} m_{h,s}(f^{-i}(X^b_0)).
\eeq
Since $\Crit_\infty(f)\sbt \Crit_h(f)$, invoking \eqref{1_2017_11_15}--\eqref{3_2017_11_15}, it directly follows
from Lemma~\ref{l1050206} that for all $j\geq 0$,
$$ 
m_{h,s}(V)
\lek \frac{m_{h,s}(B_b(R))}{m_{h,s}(X^b_0)} m_{h,s}(f^{-i}(X^b_0)\cap V) 
$$
for all $b\in P_2:= f^{-1}(\infty)\sms P_1$ and all connected components $V$ of $f^{-i}(B_b(R))$. But 
$$
m_{h,e}(X^b_0)
\comp \dist(0,X_0)^{-\frac{q_b+1}{q_b}h}m_{h,e}(X_0)
\comp 1,
$$
so, applying Lemma~\ref{l1071702} and Theorem~\ref{thm:julia}, we get
$$ 
\frac{m_{h,s}(B_b(R))}{m_{h,s}(X^b_0)}
\comp  \frac{m_{h,e}(B_b(R))}{m_{h,e}(X^b_0)}
\comp R^{2-\frac{q_b+1}{q_b}h}
\lek 1.
$$
Therefore,
$$ 
m_{h,s}(V)\lek m_{h,s}(f^{-i}(X^b_0)\cap V).
$$
Summing this inequality over all connected components $V$ of $f^{-i}(B_b(R))$, we  thus get that $ m_{h,s}(f^{-i}(V_{b}))\lek m_{h,s}(f^{-i}(X^b_0))$. Hence, for all $n\geq 0$,
$$
\sum_{b\in P_2} \sum_{i=0}^{n} m_{h,s}(f^{-i}(B_b(R)))
\lek \sum_{b\in P_2} \sum_{i=0}^{n}m_{h,s}(f^{-i}(X^b_0)).
$$
Adding  this inequality and (\ref{7050106}) side by side, we get
that 
$$
\sum_{i=0}^{n}m_{h,s}(f^{-i}(B_R))-F \lek  \sum_{i=0}^{n} m_{h,s}(f^{-i}(X_0))-G
$$
with some positive numbers $F$ and $G$ independent of $n$, resulting from
the fact that we sum in (\ref{7050106}) from $i=j$ and not from $i=0$. Thus
$$
m_{h,s,n}(B_\infty(R))\lek 1,
$$
where, for all $n \geq 1$, the measures $m_{h,s,n}$ are given by formula \eqref{5.9a}. It therefore follows from formula (\ref{eq:muequ}) of Theorem~\ref{t1h75} that $\mu_{h,s}(B_\infty(R))<+\infty$. We are done. \endpf 

\sp\fr As an immediate consequence of Theorem~\ref{t4050206} and
Proposition~\ref{p5050206} we get the following.

\bcor\lab{l6050206} 
If $f:\C\lra\oc$ is a compactly non--recurrent regular elliptic function and $\Crit_{\infty}(f)\cup \Om(f)=\es$, then  the invariant measure
$\mu_h$, equivalent to the conformal measure $m_{h,s}$ (which in this case
coincides, up to a multiplicative constant, with the packing measure $\Pi^h$), is finite. We will then always assume that $\mu_h$ is a probability measure.
\ecor

\sp We also get the following.

\bcor\lab{l7050206} 
If $f:\C\lra\oc$ is a compactly non-recurrent elliptic  function whose Julia set is equal to the entire complex plane $\mathbb C$, then there exists a unique Borel probability $f$--invariant measure $\mu_2$ equivalent
to planar Lebesgue measure on $\mathbb C$. 
\ecor

\bpf The function $f$ is regular since $h=2$. Hence, since $\Om(f)=\es $ and $\Crit_2(f)=\Crit(f)$, the existence of $\mu_2$ follows immediately from Theorem~\ref{t4050206} and Proposition~\ref{p5050206}. Uniqueness is guaranteed by
Theorem~\ref{tinv}. 
\endpf

\sp\brem\label{r1drukp124} In fact we have shown in the proof of
Proposition~\ref{p5050206} that for all $ R>0$ large enough
$m_{h,s,n}(B_\infty(R))\lek 1$, where the measures $m_{h,s,n}$ are given by formula \eqref{5.9a}
\erem

\sp\section[Real Analyticity of the Radon--Nikodym
  Derivative] {Real Analyticity of the Radon--Nikodym
  Derivative $\frac{d\mu_h}{dm_h}$}

Throughout this section we keep $f:\C\lra\oc$, a compactly non--recurrent regular elliptic function and notation from the previous section; in particular $m_h$ is the unique $h$--conformal measure for $f$ and $\mu_h$ is the $\sg$--finite $f$--invariant measure produced in Theorem~\ref{tinv}. The goal of this section is to show that the Radon--Nikodym
derivative $\frac{d\mu_h}{dm_{h,e}}$ has a real--analytic extension to a
neighborhood of $J(f)\sms \ov{{\rm PC}(f)}$ in $\mathbb C$. In the context of conformal iterated systems and conformal graph directed Markov systems such result has appeared in \cite{MPU}, see also \cite{MU2}. The proof we provide in this chapter stems from the one given in \cite{MPU}.

The first step towards this end is to work with the analytic map $\hat f:\hat\mT_f\lra\mT_f$ defined at the beginning of Section~\ref{elliptic} so that the diagram \eqref{120200408} commutes. We repeat it here.  
\beq\label{620140104} 
\begin{tikzcd}
{\C\sms f^{-1}(\infty)} \arrow{r}{f} \arrow[swap]{d}{\Pi_f} & {\C}\arrow{d}{\Pi_f} \\
{\  \hat{\mathbb T}_f} \arrow{r}{\hat f} & {\ {\mathbb T}_f}
\end{tikzcd}
\eeq
Let 
$$
\hat{J}\index{(S)}{$\hat{J}$}=\hat J(f):=\index{(S)}{$\hat{J(f)}$}\Pi_f(J(f))\sbt {\mathbb T}_f.
$$
Define the Borel probability measure $\hat{m}_h$\index{(S)}{$\hat{m}$} on
$\hat{J}$ by the formula 
\beq\label{720140104}
\hat{m}_h\index{(S)}{$\hat{\mu}$}(A):=m_{h,e}(\Pi_f^{-1}(A)\cap {\mathcal R}),
\eeq
whre $\mathcal R$ is a fundamental domain of $f$. This definition is in fact
independent of the choice of the fundamental region $\mathcal R$,
and we clearly have the following.

\bprop\label{p1_2018_08_07}
If $f:\C\lra\oc$ is a compactly non-recurrent regular elliptic function, then the Borel probability measure $\hat{m}_h$ is $h$--conformal with respect to the map $\hat{f}:\hat\mT_f\lra\mT_f$. 
\eprop

\fr Because the diagram \eqref{620140104} commutes, as an immediate consequence of Theorem~\ref{tinv}, we get the following.

\bprop\label{p2_2018_08_07}
If $f:\C\lra\oc$ is a compactly non--recurrent elliptic regular function, then for the Borel $\sg$--finite measure
\beq\lab{5_2017_11_18}
\hat{\mu}_h:= \mu_h\circ \Pi_f^{-1}
\eeq
we have that 
\begin{enumerate}
\item $\hat{\mu}_h$ is $\hat{f}$--invariant,

\item $\hat{\mu}_h(J(\hat f))=1$,

\item $\hat{\mu}_h$ is equivalent to $\hat{m}_h$, 

\item the dynamical system $\(\hat f:J(\hat f)\to J(\hat f),\hat\mu_h\)$ is metrically exact, hence ergodic,

\item $\hat\mu_h\(J(\hat f)\sms \Tr(\hat f)))=\hat m_h\(J(\hat f)\sms \Tr(\hat f)\)=0$.
\end{enumerate}
\eprop

\fr We shall  prove the following.

\blem\label{l1pj3} 
If $f:\C\lra\oc$ is a compactly  non--recurrent regular elliptic function, then the Radon--Nikodym derivative
$\hat{\rho}\index{(S)}{$\hat{\rho}$}_h:=\frac{d\hat{\mu}_h}{d\hat{m}_h}$ \,
has a real--analytic extension to a neighborhood of
$J(\hat f)\sms \Pi_f\(\ov{{\rm PC}(f)}\)$ in $\mathbb T_f$. 
\elem

\bpf 
Since the measure $m_{h,e}$  is  ergodic and
conservative (Theorem~\ref{tmaincm}), so is the measure $\hat{m}_h$. Since, by Koebe's Distortion Theorem, i.e. Theorem~\ref{Euclid-I}, we have bounded distortion on the complement of $\Pi_f(\ov{{\rm PC}(f)})$, we see that the assumptions of Theorem~\ref{t1h75} are satisfied for the dynamical system $\hat{f}:J(\hat f) \lra J(\hat f)$ and the conformal measure $\hat{m}_h$. Therefore,
\beq\label{1pj5}
\hat{\rho}_h(z)=\lim_{n \to \infty}a^{-1}_n \sum_{k=0}^n \sum_{\xi\in
\hat{f}^{-k}(z)}|(\hat{f}^{k})'(\xi)|^{-h}
\eeq
for every $z\in J(\hat f)\sms \Pi-f\(\ov{{\rm PC}(f)}\),$ where
$$
a_n=\sum_{k=0}^n \hat{m}_h(\hat{f}^{-k}(A_0))
$$ 
with some set $A_0\sbt J(\hat f)\sms \Pi_f\(\ov{{\rm PC}(f)}\)$ as required in Theorem~\ref{t1h75}. Fix  such an
arbitrary  point $z\in J(\hat f)\sms \Pi_f\(\ov{{\rm PC}(f)}\)$  and take
$r=r(z)>0$ so small that 
$$
B(z,2r)\cap \Pi\( \ov{{\rm PC}(f)}\)=\es.
$$
We can assume without loss of generality that 
$$
A_0 \sbt B(z,r).
$$
For every $k \geq 0$  and every $\xi
 \in\hat{f}^{-k}(z)$, let 
$$
\hat{f}^{-k}_{\xi}: B(z,2r) \lra {\mathbb  T}_f
$$
be the unique  holomorphic inverse branch  of $\hat{f}^{k}$ defined on $B(z,2r)$ and determined by the requirement that $\hat{f}^{-k}_{\xi}(z)=\xi$. Now  embed $\mathbb C $ into ${\mathbb C}^2$ by the formula:
$$
\C\ni x+iy \lmt (x,y) \in {\mathbb C}^2.
$$
For every  $\xi \in \hat{f}^{-k}(z)$ consider  the map $g_\xi: B(z, 2r) \to \mathbb C$
defined as  follows
$$
g_\xi(w):= {(\hat{f}_{\xi}^{-k})'(w)\over (\hat{f}_{\xi}^{-k})'(z)}.
$$\index{(S)}{$g_\xi(w)$}
Since the  ball  $B(z,2r)$ is simply connected, since the Jacobian
$g_\xi$  nowhere vanishes  on $B(z,2r)$ and since $g_{\xi}(z)=1$,
there exists  
$$
\log g_{\xi}:B(z,2r)\lra {\mathbb C},
$$
a unique holomorphic branch of logarithm $g_\xi$, such that $\log
g_{\xi}(z)=0$. By Theorem~\ref{Euclid-I} and Koebe's Distortion
Theorem~\ref{th17.46in{H}} (see formula~(\ref{17.4.18b})), there
exists a constant $\hat{K}$ such that $|\log g_\xi|\leq \hat{K}$
throughout $B(z,r)$. Expand $\log g_\xi$ into its Taylor series:
$$
\log g_{\xi} = \sum_{n=0}^\infty u_n(w-z)^n.
$$
By Cauchy's estimates
\beq\lab{1jp7}
|u_n|\le \hat{K}/r^n
\eeq
for all integers $n \geq 0$. For every   point $x+iy\in B(z,2r),$ we can write
$$
\begin{aligned}
\re\log g_{\xi}(x+iy)
& = \re\left(\sum_{n=0}^\infty u_n
\left((x-\re z)+i(y-\im z)\right)^n\right)\\
&=  \sum_{p,q=0}^\infty \re\left(u_{p+q} {p+q \choose q} i^q\right)
(x-\re z)^p(y-\im z)^q \\
 & =\sum c_{p,q}(x-\re z)^p(y-\im z)^q.
\end{aligned}
$$
In view of (\ref{1jp7}) we have 
$$
|c_{p,q}|\le {\hat{K}
r^{-(p+q)}2^{p+q}}.
$$ 
Hence, $\re\log g_{\xi}$ extends, by the same
power series expansion 
$$
{\mathbb D}_{{\mathbb C}^2}(z,r/3)\ni (x,y)\longmapsto\sum c_{p,q}(x-\re z)^p(y-\im z)^q\in\C,
$$ 
to the polydisc ${\mathbb D}_{{\mathbb C}^2}(z,r/3)$ and its modulus is
bounded above by $4\hat{K}$. Denote this extension by
$\widehat{\re} \log g_\xi$. Now, for every $n\ge 0$ consider the
function $b_n: B(z,2r)\to {\mathbb C}$ given by the formula
$$
b_n(w):=a_n^{-1}\sum_{k=0}^n\sum_{\xi \in \hat{f}^{-k}(z)}
|(\hat{f}_{\xi}^{-k})'(w)|^h.
$$
Each function $b_n$ extends to a holomorphic  function $B_n:
{\mathbb D}_{{\mathbb C}^2}(z, 2r)\lra {\mathbb C}$ as  follows
$$
B_n:= a_n^{-1}\sum_{k=0}^n\sum_{\xi\in \hat{f}^{-k}(z)}
|(\hat{f}^{k})'(\xi)|^{-h}\exp(h\widehat{\re}\log g_\xi).
$$
Since $A_0 \sbt B(z,2r)$, it follows from  (\ref{1pj5}) and Koebe's
Distortion Theorem, I (Euclidean  version), i.e. Theorem~\ref{Euclid-I},  that $L:=\sup_{n \geq 0}\{b_n(z)\}< + \infty$. Therefore, for every $w \in {\mathbb D}_{{\mathbb C}^2}(z,r/3)$, we get
$$
\begin{aligned}
|B_n(z)| &\leq a_n^{-1}\sum_{k=0}^n\sum_{\xi\in \hat{f}^{-k}(z)}
  |(\hat{f}^{k})'(\xi)|^{-h}|\exp(h\widehat{\re}\log g_{\xi}(w))| \\
&\le a_n^{-1} \sum_{k=0}^n\sum_{\xi\in \hat{f}^{-k}(z)}
   |(\hat{f}^{k})'(\xi)|^{-h} \exp(h |\widehat{\re}\log g_{\xi}(w)|) \\
&\le \exp(4h \hat{K})a_n^{-1}\sum_{k=0}^n\sum_{\xi \in
\hat{f}^{-k}(z)}
  |(\hat{f}^{k})'(\xi)|^{-h}\\
  &  = \exp(4h\hat{K})b_n(z)\\
  & \le L  \exp(4h\hat{K}).
\end{aligned}
$$
Hence, applying  Cauchy's Integral  Formula (in ${\mathbb
D}_{{\mathbb C}^2}(z,r/2)$), we see that the family $\{
B_n\}_{n=0}^{\infty}$  is equicontinuous on ${\mathbb D}_{{\mathbb
C}^2}(z,r/4)$. Thus, we can choose from $\{B_n\}_{n=0}^{\infty}$ a
subsequence uniformly convergent on ${\mathbb D}_{{\mathbb
C}^2}(z,r/5)$. Its limit  function 
$$
G_z: {\mathbb D}_{{\mathbb C}^2}(z,r_z/5)\lra {\mathbb C}
$$ 
is analytic and
$${G_z}_{|\hat{J}\cap {\mathbb D}_{{\mathbb
C}^2}(z,r_z/5)}=\hat{\rho}_{|\hat{J}\cap {\mathbb D}_{{\mathbb
C}^2}(z,r_z/5)}.$$ So ${G_z}_{|B(z,r_z/5)}$ is a real--analytic
extension  of $\hat{\rho}_{|\hat{J}\cap {\mathbb D}_{{\mathbb
C}^2}(z,r_z/5)}$. Now, if $$B(z, r_z/10)\cap B(z', r_z'/10)\neq \es,
$$
$z, z'\in \hat{J}\sms \Pi\( \ov{{\rm PC}(f)}\),$ then  choose a point
$$v  \in B(z, r_z/10)\cap B(z', r_z'/10).$$ We may assume without
loss of generality  that $r_z\leq r_z'$. Then
$$z\in B\left(z',
\frac{r_z}{10} +\frac{r_z'}{10}\right) \sbt B(z', r_z'/5).$$ So, $
z\in B(z, r_z/5)\cap B(z', r_z'/5)$; in particular:
$$
\hat{J}\cap B(z, r_z/5)\cap B(z', r_z'/5)\neq \es.
$$ 
Since  this intersection is not contained  in any real--analytic curve (its Hausdorff dimension is larger than 1), we thus  conclude that
$${G_z}_{|B(z,r_z/5)\cap B(z',r_z'/5)}= {G_{z'}}_{|B(z,r_z/5)\cap
B(z',r_z'/5)}.$$ In particular,
$${G_z}_{|B(z,r_z/10)\cap B(z',r_z'/10)}= {G_{z'}}_{|B(z,r_z/10)\cap
B(z',r_z'/10)}.$$ So,  the formula 
$$
G(w):=G_z(w)
$$ 
if $z\in \hat{J}\sms \Pi( \ov{\rm {PC(f)}})$ and $w \in B(z, r_z/10)$
provides a well-defined real--analytic function on 
$$
\bu_{z\in\hat{J}\sms \Pi( \ov{\rm{PC(f)}})} B(z, r_z/10)
$$ 
which  coincides
with $\hat{\rho} $ on $\hat{J}\sms  \Pi(\ov{{\rm PC}(f)})$. The proof is complete. 
\endpf

\bdfn\label{przeciwobrazy} 
Let $f:\C\lra\oc$ be a non--constant elliptic function. For  every $z \in \hat{\mathbb C}$ let
$f^{-1}_0(z)$ be a maximal subset of points from $f^{-1}(z)$ that are mutually incongruent with respect to the equivalence relation $\sim_f$, i.e. modulo the lattice $\Lambda_f$ of $f$. We also fix $R>0$ so large that 
$$
f(B(0,R))=\oc
$$ 
and we require in addition that 
$$
f_0^{-1}(z)\sbt B(0,R)
$$ 
for all $z\in\oc$.
\edfn

Now, we can prove the following main result of this section.

\bthm\label{t1pj11}
If $f:\C\lra\oc$ is a compactly  non--recurrent regular elliptic function, then the Radon--Nikodym  derivative $\rho_h:=\frac{d\mu_h}{dm_{h,e}}$  has a real--analytic extension to a neighborhood of $J(f)\sms \ov{{\rm PC}(f)}$ in $\oc$. In particular, for every $r>0$, the function
$$
\rho_h\big|_{J(f)\cup\{\infty\}\sms B(\ov{{\rm PC}(f)},r)}
$$
is uniformly continuous with respect to the spherical mertic on $\oc$. 
\ethm

\bpf Fix a point $z\in  J(f)\sms \ov{{\rm PC}(f)}$ and
put $R_z:= \frac{1}{2} \dist_e(z, \ov{{\rm PC}(f)})>0$. Then for every
$\xi\in f^{-1}(z)$ and every Borel set $A\sbt B_e(z,R_z)$, we have
\beq\label{620190918}
\hat{m}_h\(\Pi_f(f^{-1}_\xi(A))\)
=m_{h,e}\Big(\Pi_f^{-1}\(\Pi_f(f^{-1}_\xi(A))\)\cap \mathcal R_f\Big)
=m_{h,e}(f^{-1}_\xi(A)), 
\eeq
where, we recall $\mathcal R_f$ is a fundamental domain of $f$. Also,
$$
\begin{aligned}
\mu_h(A)
=&  \mu_h(f^{-1}(A))
= \mu_h\lt(\bigcup_{\xi \in f^{-1}_0(z)}
    \bigcup_{w \in \Lambda}w + f^{-1}_\xi (A)\rt) \\
&=\mu_h\lt(\Pi_f^{-1}\lt(\Pi_f\lt(\bigcup_{\xi \in f^{-1}_0(z)}
f^{-1}_\xi(A)\rt)\rt)\rt)\\
& = \hat{\mu}_h\lt(\Pi_f\lt( \bigcup_{\xi \in f^{-1}_0(z)}
    f^{-1}_\xi(A)\rt)\rt)\\
&=\sum_{\xi \in f^{-1}_0(z)}\hat{\mu}_h
(\Pi_f(f^{-1}_\xi(A))). 
\end{aligned}
$$ 
Therefore, using this and \eqref{620190918}, we get that
$$ \begin{aligned}
\frac{\mu_h(A)}{m_{h,e}(A)}
&=\sum_{\xi \in f^{-1}_0(z)}\frac{\hat{\mu}_h(\Pi_f(f^{-1}_\xi(A)))}{m_{h,e}(A)}=
      \sum_{\xi \in f^{-1}_0(z)}\frac{\hat{\mu}_h(\Pi_f(f^{-1}_\xi(A)))}
{\hat{m}_h(\Pi_f(f^{-1}_\xi(A)))}\cdot 
\frac{\hat{m}_h(\Pi_f(f^{-1}_\xi(A)))}{m_{h,e}(A)}\\
&= \sum_{\xi \in
f^{-1}_0(z)}\frac{\hat{\mu}_h(\Pi_f(f^{-1}_\xi(A)))}
{\hat{m}_h(\Pi_f(f^{-1}_\xi(A)))}\cdot\frac{m_{h,e}(f^{-1}_\xi(A))}{m_{h,e}(A)}.   
 \end{aligned}
$$ 
Hence, passing to Radon--Nikodym derivatives, we get
\beq\label{1j181}
\frac{d\mu_h}{dm_{h,e}}(w)
=\sum_{\xi\in f^{-1}_0(z)}\frac{d\hat{\mu}_h}{d\hat{m_{h,e}}}(\Pi_f(f^{-1}_\xi(w)))
|(f^{-1}_\xi)'(w)|^h
\eeq
for all $w \in B(z,R_z)\cap J(f)$. By Lemma~\ref{l1pj3}, the function 
$$
B(z,
\hat{R_z})\ni w \longmapsto \frac{d\hat{\mu}}{d\hat{m}}(\Pi_f(f^{-1}_\xi(w)))
|(f^{-1}_\xi)'(w)|^h\in \mathbb R
$$
is real--analytic on some sufficiently smal ball $B(z,\hat{R_z})$, $0 < \hat{R}_z \leq R_z$. Hence, since the number of terms in the series \eqref{1j181} is finite, bounded above by the (finite) number of elements of the lattice $\La_f$, this series represents a real--analytic function from $B(z,\hat R_z)$ to $\R$. We thus conclude, exactly as in the last part of the proof of Lemma~\ref{l1pj3}, that the formula (\ref{1j181}) gives a real--analytic extension of the Radon--Nikodym derivative $\frac{d\mu_h}{dm_{h,e}}$ to the open  set 
$$
\bu_{z\in J(f)\sms \ov{{\rm PC}(f)}}B(z,\hat{R}_z).
$$ 
We are done.  
\endpf

\sp 

\section[Finite and Infinite Condensation]{Finite and Infinite Condensation of Parabolic Periodic Points with Respects to the Invariant Conformal Measure $\mu_h$} 

Throughout this section we keep $f:\C\lra\oc$, a compactly non--recurrent regular elliptic function and notation from the previous section; in particular $m_h$ is the unique $h$--conformal measure for $f$ and $\mu_h$ is the $\sg$--finite $f$--invariant measure produced in Theorem~\ref{tinv}.

We continue the detailed study of the points of their finite and infinite condensation originated in the end of Section~\ref{conformal-invariant}. More precisely, we provide, under some additional mild assumptions, precise characterizations of location of such points. In the context of rational functions the results of such kind were obtained for example in \cite{ADU}, \cite{DU LMS}, \cite{DU Forum}, and \cite{U2}. 

\sp We again utilize the Marco Martens cover $(X_n)_{n=0}^\infty$ defined by \eqref{120181221} and keep working throughout the whole section within with notation and concepts of Theorem~\ref{t1h75} and Definition~\ref{d:mmmap} as these were used in the proof of Theorem~\ref{tinv}. Let $\om\in\C$ be a simple parabolic fixed point of $f$, i.e.
$$
f(\om)=\om \  \  {\rm and} \ \ f'(\om)=1.
$$
Let $p=p(\om)\ge 1$ be the corresponding number of petals. Fix any $\a\in(0,\pi)$. Let
$$ 
\De(\om,\a):= \bigcup_{j=1}^p \De_j(\om,\a)
$$\index{(S)}{$\De_j(\om,\a)$} 
be the corresponding fundamental domain, defined in Lemma~\ref{l1ch11.5}. Let
$$ 
S_r(\om):=\{\om\}\cup \bigcup_{j=1}^p S^j_r(\om,\a)
$$\index{(S)}{$S_R(\om)$} 
be the union of corresponding sectors, defined by \eqref{1_2017_29B}. In particular, Proposition~\ref{p1j108} holds, i.e.:
$$
f^{-1}_\om (S_r(\om)) \sbt S_r(\om).
$$
For every $y \in f^{-1}(\om)$ let, $f^{-1}_y(B(\om,\th_f))$  be the
connected component of $f^{-1}(B(\om,\th_f))$ containing  $y$  and for
every set $A\sbt B(\om,\th_f)$ let
$$
f^{-1}_y(A):=f^{-1}(A)\cap f^{-1}_y(B(\om,\th_f)).
$$
Keep $m_h$ to be the $h$-conformal measure for $f$. For every Borel set $A\sbt S_r(\om)$ and every integer $n \geq 0$, we  can write,
\beq\label{1j125}
\begin{aligned}
S_n m_{h,s}(A)
:&=\sum_{j=0}^n m_{h,s}\(f^{-j}(A)\) \\
&= m_{h,s}(f^{-n}_\om(A))
+\sum_{y\in f^{-1}(\om)\sms\{\om\}} \sum_{k=0}^{n-1}S_{n-(k+1)}
   m_{h,s}(f^{-1}_y(f^{-k}_\om(A))).
\end{aligned}
\eeq

\fr We shall first prove the following.

\bprop\lab{p120181222}
Let $f:\C\lra\oc$ be a compactly non--recurrent regular elliptic function. If $\om\in\Omega(f)\sms \(\ov{\PC(f)}\cap J(f)\)$, is a simple parabolic fixed point of $f$, then 
\beq\lab{620181222}
\mu_h\(f_\om^{-l}(F)\)\comp (l+1)^{1-\frac{p(\om)+1}{p(\om)}h}m_h(F)
\  \  {\rm and} \  \
\frac{d\mu_h}{dm_h}\bigg|_{f_\om^{-l}(\De(\om))}\comp l+1
\eeq
for every integer $l\ge 0$ and every Borel sets $F\sbt \De(\om,\a)$.
\eprop

\bpf
Let 
$$
A:=f^{-l}_\om(F),
$$ 
where $F$ is a Borel set contained in $\Delta_j(\om, \alpha)$ for some $ 1\leq j \leq p$. Fix an integer $n\ge 1$. It then follows from Proposition~\ref{p1j111DL}, conformality of the measure $m_{h, s}$ and Koebe's Distortion Theorem, that for all $y\in f^{-1}(\om)\sms\{\om\}$ and all integers $0\le k \leq n-1$, we have that
$$ 
S_{n-(k+l)}m_{h,s}\(f_y^{-1}(f_\om^{-k}(f_\om^{-l}(F))\)
\comp (k+l)^{-\frac{p+1}{p}h}S_{n-(k+1)}m_{h,s}(f^{-1}_y(F)).
$$
So, substituting this to (\ref{1j125}), we obtain
\beq\label{1fpn32}
\begin{aligned}
m_{h,s,n}&(f_\om^{-l}(F)) =\\
&=\frac{S_nm_{h, s}(f_\om^{-l}(F))}{S_nm_{h, s}(X_0)}\comp \\
&\comp \frac{m_{h, s}(f_\om^{-n}(f_\om^{-l}(F)))}{S_nm_{h, s}(X_0)}+\\
& \  \ + \sum_{y \in f^{-1}(\om) \sms \{\om \}} \sum_{k=0}^{n-1} (k+l)^{-\frac{p+1}{p}h} \frac{S_{n-(k+1)}m_{h,s}(f^{-1}_y(F))}{S_{n-(k+1)}m_{h, s}(X_0)}
 \frac{S_{n-(k+1)}m_{h, s}(X_0)}{S_nm_{h, s}(X_0)}\\
  &\comp \frac{m_{h, s}(f_\om^{-n}(f_\om^{-l}(F)))}{S_nm_{h, s}(X_0)}+\\
  & \  \  +\sum_{y \in f^{-n}(\om) \sms \{\om \}} \sum_{k=0}^{n-1} (k+l)^{-\frac{p+1}{p}h} \frac{ S_{n-(k+1)}m_{h,s}(f^{-1}_y(X_0))}{S_{n-(k+1)}m_{h, s}(X_0)}\frac{m_{h,s}(F)}{m_{h,s}(X_0)}\frac{S_{n-(k+1)}m_{h, s}(X_0)}{S_nm_{h, s}(X_0)}\\
&\comp \frac{m_{h, s}(f_\om^{-n}(f_\om^{-l}(F)))}{S_nm_{h, s}(X_0)}+\\
& \  \ +\sum_{k=0}^{n-1} (k+l)^{-\frac{p+1}{p}h} \frac{\sum_{y \in f^{-n}(\om) \sms \{\om \}} S_{n-(k+1)}m_{h,s}(f^{-1}_y(F))}{S_{n-(k+1)}m_{h, s}(X_0)}\frac{m_{h,s}(F)}{m_{h,s}(X_0)}\frac{S_{n-(k+1)}m_{h, s}(X_0)}{S_nm_{h, s}(X_0)}\\
&\comp \frac{m_{h, s}(f_\om^{-n}(f_\om^{-l}(F)))}{S_nm_{h, s}(X_0)}+ \\
&\  \  +\sum_{k=0}^{n-1} (k+l)^{-\frac{p+1}{p}h}\frac{S_{n-(k+1)}m_{h,s}\lt(\bigcup_{y \in f^{-1}(\om) \sms \{\om \}}f^{-1}_y(X_0)\rt)}{S_{n-(k+1)}m_{h, s}(X_0)}\frac{S_{n-(k+1)}m_{h, s}(X_0)}{S_nm_{h, s}(X_0)}\frac{m_{h, s}(F)}{m_{h, s}(X_0)}
\end{aligned}
\eeq
Now,
\beq\label{120190424}
\begin{aligned}
\frac{S_{n-(k+1)}m_{h,s}\lt(\bigcup_{y \in f^{-1}(\om) \sms \{\om \}}f^{-1}_y(X_0)\rt)}{S_{n-(k+1)}m_{h, s}(X_0)}
& \leq  \frac{S_{n-k}m_{h, s}(X_0)}{S_{n-(k+1)}m_{h, s}(X_0)} 
=\frac{m_{h, s}(f^{-(n-k)}(X_0))}{S_{n-(k+1)}m_{h, s}(X_0)}+1\\
& \leq \frac{1}{m_{h, s}(X_0)}+1.
\end{aligned}
\eeq
Also, since there exists an integer $j\geq 0$ such that
$$
m_{h,s}\lt(X_0\cap f^{-j}\lt(\bigcup_{y \in f^{-1}(\om) \sms \{\om \}}f^{-1}_y(X_0)\rt)\rt)>0,
$$
by making use of item (4) of Definition~\ref{d:mmmap}, we conclude that
$$
\frac{S_{n-(k+1)}m_{h,s}\lt(\bigcup_{y \in f^{-1}(\om) \sms \{\om \}}f^{-1}_y(X_0)\rt)}{S_{n-(k+1)}m_{h, s}(X_0)}
\succeq 1.
$$
Taking  this and  (\ref{120190424}) together, we obtain
$$
\frac{S_{n-(k+1)}m_{h,s}\lt(\bigcup_{y \in f^{-1}(\om) \sms \{\om \}}f^{-1}_y(X_0)\rt)}{S_{n-(k+1)}m_{h, s}(X_0)}\comp 1.
$$
Inserting this to (\ref{1fpn32}), we get
\beq\label{1fpn33}
\begin{aligned}
m_{h,s,n}(f^{-l}_\om(F)) & = \frac{m_{h,s}(f^{-n}_\om(f^{-l}_\om(F)))}{S_n m_{h, s}(X_0)}
 +\sum_{k=0}^{n-1} (k+l)^{-\frac{p+1}{p}h} \frac{S_{n-(k+1)}m_{h,s}(X_0)}{S_n m_{h, s}(X_0)} + m_{h,s}(F)\\
 &\comp o(1) +\sum_{k=0}^{n-1} (k+l)^{-\frac{p+1}{p}h} \frac{S_{n-(k+1)}m_{h,s}(X_0)}{S_n m_{h, s}(X_0)} m_{h,s}(F).
\end{aligned}
\eeq
Therefore, on the one hand, 
$$ 
m_{h,s,n}(f^{-l}_\om(F)) 
\preceq o(1) +\sum_{k=0}^{n-1} (k+l)^{-\frac{p+1}{p}h} m_{h,s}(F) 
\leq o(1) +(l+1)^{1-\frac{p+1}{p}h}m_{h,s}(F),
$$
where the $o(1)$ symbol is with respect to $n\to \infty$. So, it follows from  formula (\ref{eq:muequ}) that
\beq\label{2fpn33}
\mu_h (f^{-l}_\om(F))\preceq (l+1)^{1-\frac{p+1}{p}h}m_{h,s}(F).
\eeq
Hence
$$
\frac{\mu_h(f^{-l}_\om(F))}{m_{h,s}(f^{-l}_\om(F))}\preceq \frac{(l+1)^{1-\frac{p+1}{p}h}m_{h,s}(F)}{(l+1)^{-\frac{p+1}{p}h}m_{h,s}(F)}=(l+1).
$$
Therefore,
\beq\label{3fpn33}
\frac{d \mu_h}{dm_{h,s}}\Big|_{ f^{-l}_\om(\Delta(\om))}
\preceq l+1.
\eeq
On the other hand, there exists $q \geq 1$ so large that
$$
\sum_{k=0}^q (k+l)^{-\frac{p+1}{p}h}
\geq \frac{1}{2}\lt(\frac{p+1}{p}h -1\rt)(l+1)^{1-\frac{p+1}{p}h}.
$$
Inserting  this into (\ref{1fpn33}), we get for every $n\geq q+1$ that
\beq\label{1fpn34}
\begin{aligned}
m_{h,s,n}(f^{-l}_\om(F)) 
&\succeq o(1) +\sum_{k=0}^{q} (k+l)^{-\frac{p+1}{p}h} \frac{S_{n-(k+1)}m_{h,s}(X_0)}{S_nm_{h, s}(X_0)}m_{h,s}(F)\\
&\succeq o(1)+ \frac{S_{n-(q+1)}m_{h,s}(X_0)}{S_n m_{h, s}(X_0)}\sum_{k=0}^{q} (k+l)^{-\frac{p+1}{p}h} m_{h,s}(F)\\
&\succeq  o(1)+l^{1-\frac{p+1}{p}h}\frac{S_{n-(q+1)}m_{h,s}(X_0)}{S_n m_{h, s}(X_0)}m_{h,s}(F)\\
&\succeq  o(1)+l^{1-\frac{p+1}{p}h}\frac{S_{n}m_{h,s}(X_0)-q}{S_n m_{h, s}(X_0)}m_{h,s}(F)\\
&\succeq o(1)+(l+1)^{1-\frac{p+1}{p}h}m_{h,s}(F).
\end{aligned}
\eeq
Since $f^{-l}_\om(F)$ is contained in a finite union of the sets $X_j$, $ j\geq 0$, it follows from (\ref{120190618}) and (\ref{1fpn34}) that
$$
\mu_{h}(f^{-l}_\om(F))\succeq l^{1-\frac{p+1}{p}h}m_{h,s}(F). 
$$
Hence
$$ \frac{\mu_n(f^{-l}_\om(F))}{m_{h,s}(f^{-l}_\om(F))} \succeq \frac{ (l+1)^{1-\frac{p+1}{p}h} m_{h,s}(F)} {(l+1)^{-\frac{p+1}{p}h} m_{h,s}(F)}=l+1.$$
Therefore,
$$
\frac{d\mu_h}{dm_{h,s}}\Big|_{f^{-l}_\om(\Delta(\om))}\succeq l+1.
$$
Combining theses two formulas with (\ref{2fpn33}) and (\ref{3fpn33}) we get that
$$ \mu_{h}(f^{-l}_\om(F))\comp (l+1)^{1-\frac{p+1}{p}h}m_{h,s}(F)$$
and
$$ 
{\frac{d\mu_n}{dm_{h}}}\Big|_{f^{-l}_\om(\delta(\om))} 
\comp l+1.
$$
The proof of Proposition~\ref{p120181222} is complete. 
\epf

\sp\fr As an immediate consequence of this proposition we get the following.

\bcor\label{c120190425}
Let $f:\C\lra\oc$ be a compactly non--recurrent regular elliptic function. If $\om\in\Omega(f)\sms \(\ov{\PC(f)}\cap J(f)\)$, is a simple parabolic fixed point of $f$, then for every $\a\in(0,\pi)$ and every integer $k\ge 0$, we have that
$$
\mu_h\lt(\bu_{j=k}^\infty f_\om^{-j}(\De(\om))\rt)
\comp (k+1)^{2-\frac{p(\om)+1}{p(\om)}h}
$$
if $h \ge \frac{2p(\om)}{p(\om)+1}$ and 
$$
\mu_h\lt(\bu_{j=k}^\infty f_\om^{-j}(\De(\om,\a))\rt)=+\infty
$$
if $h \leq \frac{2p(\om)}{p(\om)+1}$.
\ecor

Passing to an appropriate iterate if $\om$ is a parabolic periodic point, as an immediate consequence of this proposition, we get the following main result of this section.

\bthm\label{t1j125} 
Let $f:\C\lra\oc$ be a compactly non--recurrent regular elliptic function. If $\om\in\Omega(f)\sms \(\ov{\PC(f)}\cap J(f)\)$, then 
\nl \centerline{$\om \in J_{\mu_h}(\infty)$ \  \ if and only if \  \  $h \leq \frac{2p(\om)}{p(\om)+1}$.}
\ethm 

\sp Observe that if $f$ is just compactly non--recurrent, then we can always find a point $y\in f^{-1}(\om ) \sms \{\om \}$ such that  
$$
\PC(f)\cap
f^{-1}_y(B(\om, 2R)))=\es
$$ 
for all $R>0$ sufficiently small. We then get $\gek$ part of Proposition~\ref{p120181222}. In conclusion, we obtain the following.

\bthm\label{t1j131} 
Let $f:\C\lra\oc$ be a compactly non--recurrent regular elliptic function. If $\om \in \Omega(f)$ and $h \leq
\frac{2p(\om)}{p(\om)+1}$, then $\om \in J_{\mu_h}(\infty))$.
\ethm

We end this section with a proposition showing  that  critical
points eventually landing  at parabolic points make these latter to
be more likely of infinite condensation. We need the following.

\bthm\label{l2j131}
Let $f:\C\lra\oc$ be a compactly non-recurrent regular elliptic function. If $\om \in \Omega(f)$ and $c\in J(f)$ is  a critical point of $f$s such that $f^l(c)=\om$ for some integer $l \geq 1$,
then 
$$ 
m_{h,s}(B(c,r))\comp r^{h+q_cp(\om)(h-1)}.
$$ 
\ethm

\bpf The proof of this theorem follows a simple
"integration" method originated in [\cite{DU3}, Lemma 4.8] and since then frequently used in similar contexts. Recall that given $x\in \mathbb
C$ and $0<r_1<r_2$, we denoted
$$
A(x;r_1,r_2)=\{z\in\mathbb C: r_1<|z-x|\le r_2\}.
$$  
Set $q:=q_c$. Using then Lemma~\ref{ldp4.5}, conformality of the measure$m_{h,s}$ and its atomlesness, we can write for every $r>0$ sufficiently small, as follows.
$$
(r^q)^{h+p(h-1)} 
\comp m_{h,s}\(A(f^l(c),r^q,(2r)^q)\) 
\comp m_{h,s}\(f^l(A(c,r,2r))\) 
\comp m_{h,s}\(A(c,r,2r)\)r^{(q-1)h},
$$
and therefore $m_{h,s}\(A(c,r,2r)\)\comp r^{h+qp(h-1)}$. Thus
$$
\begin{aligned}
m_{h,s}(B(c,r)) 
&=\sum_{n=1}^\infty m_{h,s}\(A(c,2^{-n}r,2^{-(n-1)}r)\)
\comp\sum_{n=1}^\infty r^{h+qp(h-1)}2^{-n(h+qp(h-1))} \\
&= r^{h+qp(h-1)}\sum_{n=1}^\infty 2^{-n(h+qp(h-1))}.
\end{aligned}
$$
Since $m_{h,s}(B(c,r))$ is finite, it first follows from this formula that
$h+qp(h-1)>0$, or equivalently that $h>{qp\over 1+qp}$, and then
that 
$$
m_{h,s}(B(c,r))\comp r^{h+qp(h-1)}.
$$
We are done. 
\endpf

\bprop\label{p3j131} 
Let $f:\C\lra\oc$ be a compactly non--recurrent elliptic function with $\Crit_\infty(f)=\es$. If $\om \in \Omega(f)$, $l \geq 1$ is an integer, $c\in f^{-l}(\om) \sms \ov{\PC(f)}$ is a critical point of $f^l$, and $h \leq \frac{2 q_c p(\om)}{q_cp(\om)+1}$, then $\om \in J_{\mu_h}(\infty)$. 
\eprop

\bpf There exists $R_0>0$ such that  $\PC(f)\cap B(\om,
2R_0)=\es$. Fix $ 0< R \leq R_0$. We may assume without loss of
generality  that $X_0 \sbt B(c, R)$ is the set defined by formula \eqref{120181221} with $n=0$. As in the previous proof, set
$$
q:=q_c.
$$
Let $f_c^{-l}(B(\om,2R_0))$ be the connected component of $f_c^{-l}(B(\om,2R_0))$ containing $c$ and for every set $G\sbt B(\om,2R_0)$, let
$$
f_c^{-l}(G)=f_c^{-l}(B(\om,2R_0))\cap f^{-l}(G).
$$
By virtue of Lemma~\ref{l2j131}, assuming that $R_0>0$ is small enough, we get for every $k \geq 0$, that
$$ 
m_{h,s}\(f^{-l}_c\(f^{-kl}_\om(J(f)\cap B(\om, R))\)\)
\comp\left((kl)^{-\frac{1}{pq}}\right)^{h+pq(h-1)}
=(kl)^{1-h \frac{pq+1}{pq}}\comp k^{1-h\frac{pq+1}{pq}}.
$$
Applying Theorem~\ref{Spher-I}, we therefore conclude that for all integers
$k,n \geq 1$
\beq\label{1j133}
\begin{aligned}
m_{h,n}(f^{-l}_c(f^{-kl}_\om(J(f)\cap B(\om, R))))
&\comp\frac{m_{h,s}\(f^{-l}_c\(f^{-kl}_\om(J(f)\cap B(\om,R))\)\)}{m_s(X_0)}\\
&\comp m_{h,s}\(f^{-l}_c\(f^{-kl}_\om(J(f)\cap B(\om, R))\)\)\\
& \comp k^{1-h\frac{pq+1}{pq}},
\end{aligned}
\eeq
where $m_{h,n}:=(m_{h,s})_n$ are the measures defined by formula  \eqref{5.9a}. Since the sequence 
$$
\Big(m_{s,n}\(f^{-l}_c(f^{-kl}_\om (J(f))\cap B(\om,R))\)\Big)_{n=0}^{\infty}
$$
is bounded, and, since by virtue of Lemma~\ref{l1j137},
$$
f^{-l}_c(J(f)\cap B(\om, R)) \cap
\bigcup_{j=s}^\infty Y_j=\es 
$$
for all $s \geq 0$ large enough, the formulas (\ref{120190618}) and (\ref{1j133}) yield
$$
\mu_h\(f_c^{-l}(f^{-kl}_\om(J(f)\cap B(\om,R)))\)
=l_B\((m_{h,n}(f^{-l}_c(f^{-kl}_\om (J(f)\cap B(\om,R))))\)_{n=0}^\infty)
\gek k^{1-h \frac{pq+1}{pq}}.
$$ 
But, since the measure $\mu$  is $f^l$-invariant, we therefore get for every integer $k\geq 0$ that
$$
\begin{aligned}
\mu_h\(f^{-kl}_\om&(J(f)\cap B(\om,R))\)=\\
&=\mu_h\(f^{-l}(f^{-kl}_\om(J(f)\cap B(\om,R)))\)\\
&=\mu_h\(f^{-(k+1)}_\om(J(f)\cap B(\om, R))\)+\sum_{y \in f^{-l}(\om)\sms
\{\om\}} \mu_h\(f^{-l}_y(f^{-k}_\om(J(f)\cap B(\om,R)))\)\\
&\geq \mu_h\(f^{-(k+1)l}(J(f)\cap B(\om,R))\) 
+\mu_h\(f^{-l}_c(f^{-k}_\om(J(f)\cap B(\om, R)))\)\\
&\geq \mu_h\(f^{-(k+1)l}(J(f)\cap B(\om, R))\)+
B(k+1)^{1-h\frac{pq+1}{pq}},
\end{aligned}
$$ 
with some  constant $B\in(0,+\infty)$ independent  of $k \geq 0$. Therefore, by induction:
$$ 
\mu_h(B(\om, R))
\geq B\sum_{k=0}^\infty (k+1)^{1-h\frac{pq+1}{pq}}=B\sum_{k=1}^\infty
k^{1-h\frac{pq+1}{pq}}=+\infty,
$$
as, by our assumptions, $1-h\frac{pq+1}{pq}\leq -1$. Since this
holds for all $R>0$ small enough, we are done. 
\endpf

\sp\section{Closed Invariant Subsets, $K(V)$ Sets, and Summability Properties}\label{CKS}

\sp In this section we deal with compact forward invariant subsets of non--recurrent elliptic functions. For the sake of future applications we will be mostly occupied  with projected maps on the tori. We will define and characterize expanding and pseudo--expanding invariant sets. We will get good upper estimates of their Hausdorff and box--counting dimensions. Next, similarly as in the proof of Lemma~\ref{l1071805}, but even more thoroughly and actually working on the projected torus $\mT_f$, we will push forward in the current section the method of $K(V)$ sets \index{(S)}{$K(V)$} developed in \cite{DU2}, comp. \cite{KU5}. Its ultimate goal will be Lemma~\ref{l1j251A}, a summability result which will constitute the main ingredient in Sections~\ref{Pre YT}, \ref{Young's Towers for Subexpanding}, and  \ref{parabolicfiniteinvariantmeasures} for proving strong regularity of conformal graph directed Markov systems resulting from nice sets. This in turn will be instrumental, in the same sections, for  proving statistical properties of dynamical systems generated by elliptic functions and finiteness of corresponding Kolmogorov--Sinai metric entropy. 

Recall also that $\hat f:\hat\mT_f\lra\mT_f$ is the projection of an elliptic function $f:\C\lra\oc$ to the torus $\mT_f$ defined by the diagram \eqref{620140104}. Denote
\beq\label{1_2017_09_20}
B_\infty(\hat{f}):=\Pi_f(f^{-1}(\infty)).
\eeq
Note that $B_\infty(\hat{f})$ is a finite set. Recall that
\beq\lab{120200116}
\Crit(\hat f^l)
=\big\{x\in T^{(l)}_f:(\hat f^l)'(x)=0\big\}
=\bu_{j=0}^{l-1}\hat f^{-j}\(\Crit(\hat f)\).
\eeq
Set 
$$
\Om(\hat f):=\Pi_f(\Om(f)).
$$
Obviously $\Om(\hat f)$ is the set of all rationally indifferent periodic points of $\hat f:\hat\mT_f\lra\mT_f$.

\sp From now on throughout this section, let $f:\C\lra \oc$ be a non--recurrent elliptic function.

\sp\bdfn\label{d320200116}
Fix an integer $ l\geq 1$. We say that a closed (equivalently compact) set $ X \sbt \mathbb T^{(l)}_f\cap J(f)$ in $\mT_f$ is expanding \index{(N)}{expanding} for $\hat{f}^l: \mathbb T_f^{(l)}\lra  \mathbb T_f$ if and only if

\begin{itemize}
\item[(a)] $\hat{f}^l(X)\sbt X$ i.e.  $X$ is forward  invariant under $\hat{f}^l$

\,

and

\, 

\item[(b)] $X \cap\(\Crit(\hat f^l)\cup \Om(\hat{f})\)=\es$.

\,

\,

\fr If instead of (b), we have merely

\,

\, 

\item [(c)] $\overline{(X \sms \Om(f))}\cap \(\Crit(\hat f^l)\cup \Om(\hat{f})\)=\es$,

\,

\fr then  $X$ is called quasi--expanding \index{(N)}{quasi--expanding} for $\hat{f}^l$. 
\end{itemize}
\edfn

\sp\fr Of course, each expanding set is quasi--expanding.  We now shall prove some basic properties of expanding and quasi--expanding sets. In particular  we will prove some characterizations of  such sets
that will fully justify  the name expanding and quasi--expanding. We start with the following. Because of (a),
\beq\lab{120200115}
X\cap B_\infty(\hat{f})=\es.
\eeq
We therefore immediately get the following.

\bobs\label{o120200115}
Let $f: \mathbb C \lra \hat{\mathbb C}$ be a non--recurrent elliptic function. Fix an integer $ l\geq 1$. If $X \sbt \mathbb T_f^{(l)}\cap J(\hat{f})$ is a quasi--expanding set for $\hat{f}^{l}$, then
$$
X\cap \ov{\Crit(\hat f^l)}=\es.
$$
\eobs

\sp\fr In addition, looking up at the right--hand side of formula \eqref{120200116}, we get the following.

\bobs\label{o120200116}
Let $f: \mathbb C \lra \hat{\mathbb C}$ be a non--recurrent elliptic function. Fix an integer $l \geq 1$. Then,
\ben
\item If $X \sbt \mathbb T_f^{(l)}\cap J(\hat{f})$ is an expanding set for $\hat{f}^{l}$, then
$$
\lt(\bu_{j=0}^{l-1}f^j(X)\rt)\cap \ov{\Crit(\hat f)}=\es.
$$
\item If $X \sbt \mathbb T_f^{(l)}\cap J(\hat{f})$ is a quasi--expanding set for $\hat{f}^{l}$, then
$$
\lt(\bu_{j=0}^{l-1}f^j(X)\rt)\cap \ov{\Crit(\hat f)}=\es.
$$
\een
\eobs

\sp We also record the following immediate.

\bobs\label{o1kvs1} Let $f: \mathbb C \lra \hat{\mathbb C}$ be a non--recurrent elliptic function.  Fix an integer $l \geq 1$. If $ X \sbt \mathbb T_f^{(l)}\cap J(f)$ is a closed set in $\mT_f$ invariant under $\hat{f}^l$,  then  the following conditions are equivalent:
 \begin{itemize}
 \item [(1)]  $X$ is quasi--expanding for $\hat{f}^l$.
 
 \,
 
 \item [(2)] $ X \cap \Crit(\hat f^l) = \es $  and $ X \sms \Om(\hat{f}^l) \sbt \mathbb T_f$ is a closed set.
 \end{itemize}
 \eobs

\fr We shall prove the following.

\blem\label{l2kvs1} Let $f: \mathbb C \lra \hat{\mathbb C}$ be a non--recurrent elliptic function. Fix an integer $ l\geq 1$. If $X \sbt \mathbb T_f^{(l)}\cap J(\hat{f})$ is a quasi--expanding set for $\hat{f}^{l}$, then there exists $\eta >0$ such that if $ x \in X$, $n \geq 0$ is an integer, and $\hat{f}^{ln}(x)\notin \Om (\hat{f})$, then
\beq\label{1kvs1}
\Crit( \hat{f}^{ln}) \cap \Comp(x, \hat{f}^{ln}, 2 \eta)=\es.
\eeq
In particular, the holomorphic branch $\hat{f}^{-ln}_x: B(\hat{f}^{-ln}(x), 2 \eta)\lra \mathbb T_f$ of $ \hat{f}^{ln}$ sending $\hat{f}^{ln}(x)$ back to $x$ is well--defined.
\elem

\bpf Since the set $ X \sms \Om (\hat{f})$ is compact,  so is the set
$$
\tilde{X}
:=\Pi_f^{-1}(X \sms \Om(\hat{f}))\cup \{\infty\}\sbt J(f)\sms \Om(f).
$$
Also, because of Observation~\ref{o120200115},
$$ 
\varepsilon 
:= \dist\(\Crit(\hat f^l), X\)>0.
$$
It thus follows Theorem~\ref{mnt6.3} that there exists $\eta >0$ such that
$$
\diam_e\(\hat{f}^{lj} (\Comp(x,\hat{f}^{ln}, 2 \eta))\)< \varepsilon
$$
for every $n \geq 0$, every $0\leq j \leq n$ and every $x \in X$ such that $ \hat{f}^{ln}(x) \notin \Om(\hat{f})$. Therefore,
$$ 
\Crit(\hat{f}^{l})\cap \hat{f}^{lj} \(\Comp (x,\hat{f}^{ln}, 2 \eta)\)=\es
$$
for every $0 \leq j \leq n$. Formula (\ref{1kvs1})  thus holds and Lemma~\ref{l2kvs1} is proved. \qed

\blem\label{l1kvs2} Let $f: \mathbb C \lra \hat{\mathbb C}$ be a non--recurrent elliptic function. Fix an integer $l\geq 1$. If $X \sbt \mathbb T^{(l)}_f \cap J(\hat{f})$ is a quasi--expanding set for  $\hat{f}^{l}$, then there exist constants $C\in [1, +\infty)$ and $\beta >0$ such that  
$$
\big|(\hat{f}^{ln})'(z)\big| \geq Ce^{\beta n}
$$  
for every integer $n \geq 0$ and every $z \in X$ such that $\hat{f}^{ln}(z)\notin \Om (\hat{f})$.
\elem

\fr {\sl Proof.}  Applying Theorem~\ref{mnt6.3} in the same way as in the proof of Lemma~\ref{l2kvs1}  but this time also with the help of this lemma itself and Koebe's Distrortion Theorem, we conclude that there exists an integer $N \geq 1$ such that
\beq\label{1kvs2}
\big| (\hat{f}^{lk})'(x)\big| \geq e
\eeq
for every integer $k \geq N$ and every $x \in X $ such that $\hat{f}^{lk}(x)\notin \Om(\hat{f})$. If $ n \geq 0$ is an integer, $ z \in X$  and $\hat{f}^{ln}(z) \notin
\Om(\hat{f})$, then 
$$ 
\hat{f}^{lj}(z) \in X \sms \Om(f)
$$ 
for every integer $0 \leq j \leq n$.  Write uniquely $ln:= qlN +r$, where $ q, r\geq 0$ are integers with $ r \leq lN-1$. Applying then (\ref{1kvs2}) $qN$  times, we get
$$
\big|(\hat{f}^{ln})'(z)\big|
=\big|(\hat{f}^{qlN})'(z)\big|\cdot\big|(\hat{f}^{r})'(\hat{f}^{qlN}(z))\big|
\geq e^q  |(\hat{f}^{r})'(\hat{f}^{qlN}(z))| 
\geq  e^{-r}M^{lN} e^{\frac{1}{N}n},
$$
where $M:=\min\big\{1, \inf\{|(\hat{f}')(x)|:x \in X\}\big\}>0$ since, because of Observation~\ref{o120200116}, 
$$ 
\dist\lt(\bu_{j=0}^{l-1}f^j(X), \Crit(\hat f)\rt)>0. 
$$
The proof is complete. 
\qed

\sp As a fairly immediate  consequence of this lemma, we get the following.

\bcor\label{c2kvs2}  Let $f: \mathbb C \lra \hat{\mathbb C}$ be a non--recurrent elliptic function. Fix an integer $l\geq 1$. A closed $\hat{f}^l$--invariant set $ X \sbt J(f)\cap \mathbb T^{(l)}_f$ is expanding for $\hat{f}^{l}$ if and only if there exists constants $C\in [1, +\infty)$ and $ \beta >0$ such that
\beq\label{2kvs2}
\big|(\hat{f}^{ln})'(z)\big| \geq  Ce^{\beta n}
\eeq
for every integer $n \geq 0$ and every $z \in X$.
\ecor

\bpf The implication $(\Rightarrow)$ is a direct consequence of Lemma~\ref{l1kvs2}. Proving $(\Leftarrow)$, formula (\ref{2kvs2}) implies both that 
$$
X \cap \bigcup_{j=0}^{l-1}\hat{f}^{-j}(\Crit (\hat{f})) =\es
\  \  \  {\rm and}  \  \  \
X\cap \Om(\hat{f})=\es. 
$$
Looking up at the right--hand side of formula \eqref{120200116}, the proof is thus complete. 
\qed

\sp\fr We have also the following.

\bcor\label{c3kvs2} Let $f: \mathbb C \lra \hat{\mathbb C}$ be a non--recurrent elliptic function. Fix an integer $l\geq 1$. A closed $\hat{f}^l$-invariant set $ X \sbt J(f)\cap \mathbb T^{(l)}_f$ is quasi--expandng for $\hat{f}^l$ if and only if there exist constants  $C \in [1, +\infty)$ and $ \beta >0$ such that
\beq\label{1kvs3}
\big|(\hat{f}^{ln})'(z)\big| \geq  Ce^{\beta n}
\eeq
for every iteger $n \geq 0$ and every $z \in X$  such that $\hat{f}^{ln}(z) \notin \Om(\hat{f})$.
 \ecor

\bpf The implication $(\Rightarrow)$ is a direct consequence of Lemma~\ref{l1kvs2}. Proving $(\Leftarrow)$, formula (\ref{1kvs3}), as in the proof of the previous corollary, implies that 
$$
X \cap \bigcup_{j=0}^{l-1}\hat{f}^{-j}(\Crit (\hat{f})) =\es. 
$$
By virtue of the right--hand side of formula \eqref{120200116} this means that
$$
X \cap\Crit(\hat f^l)=\es.
$$
Seeking contradiction suppose that
\beq\label{2kvs3}
\overline{(X \sms \Om(\hat{f}))}\cap \Om(\hat{f}) \neq  \es.
\eeq
Passing to an appropriate  integral multiple of $l$ we may  assume  that each  element of $\Om(\hat{f})$  is a simple parabolic fixed point of $\hat{f}^{l}$. Fix  an integer $n \geq 1$. Then there exists $ \varepsilon >0$  such that
\beq\label{3kvs3}
\big|(\hat{f}^{ln})'(x)\big| \leq  2
\eeq
for all $x \in B(\Om(\hat{f}),\varepsilon)$. By virtue of (\ref{2kvs3})
and continuity properties of $\hat{f}^{l}$ there exists $x_n\in (X \sms \Om(\hat{f}^{l}))\cap B(\Om(\hat{f}^{l}), \varepsilon)$ such that $\hat{f}^{ln}(x) \notin \Om(\hat{f}^{l})$. It then follows from (\ref{3kvs3}) that $| (\hat{f}^{ln})'(x_n)| \leq  2$. Hence, 
$$
\limsup_{ n \to \infty}\big| (\hat{f}^{ln})'(x_n)\big| \leq  2.
$$
On the other hand, it follows  from (\ref{1kvs3}), that 
$$
\liminf_{n\to\infty}\big|(\hat{f}^{ln})'(x_n)\big|=+\infty.
$$
This contradiction  finishes the proof of Corollary~\ref{c3kvs2}. 
\qed

\sp We will also need the following.

\bobs\label{o1kvs3}
Let $f: \mathbb C \lra \hat{\mathbb C}$ be a non--recurrent elliptic function.  Fix an integer $l\geq 1$.  If  $ X \sbt \mathbb T^{(l)}_f \cap J(\hat{f})$ is quasi--expanding for $\hat{f}^l$,  then  $\hat{f}^l \in \mathcal A(X)$ in the sense of Chapter\ref{CMHDIM}.
\eobs

We shall now provide a natural way to construct expanding and quasi--expanding sets for elliptic functions. Let $f:\C\lra\oc$ be a non--constant elliptic function.  
Recall from \eqref{KV1} that if $V$ is an open subset of $\oc$, then
$$
K_J(V)=J(f)\cap \bigcap_{n\geq 0} f^{-n}(\oc\sms
V)
$$
Given an open set $\hat V\sbt \mT_f$, the set $K(\hat V)$ has an analogous meaning as the set $K(V)$ introduced above:
\beq\label{KV2}
K(\hat V):=\bigcap_{n\geq 0} \hat f^{-n}(\mT_f\sms \hat V)
\eeq
More generally, given an integer $l\ge 1$ and an open set $\hat V\sbt \mathbb T_f$, let
$$
K_l(\hat V):=\bigcap_{n\geq 0} f^{-ln}(\mT_f\sms \hat V)
$$
Of course, $K_1(\hat V)=K(\hat V)$. We collect the most basic properties of the sets $K_l(\hat V)$ in the following.

\bobs\label{o1kvs3.1-2020-01-16}
Let $f: \mathbb C \lra \oc$ be a non--recurrent elliptic function. Fix an integer $l\geq 1$ and an open set $\hat{V} \sbt \mathbb T_f$. Then:
\begin{itemize}

\item [(a)] The map $\hat f^l$ is well defined on $K_l(\hat{V})$ and
$$
\hat f^l\(K_l(\hat{V})\)\sbt K_l(\hat{V}).
$$
\item [(b)]
$$
K_l(\hat{V})\cap \bu_{n=0}^\infty\hat f^{-n}(B_\infty(\hat f))=\es.
$$
\item [(c)] $K_l(\hat{V})$ is a closed subset of $\mathbb T_f^{(1)}=\mathbb T_f\sms B_\infty(\hat f)$.

\,

\item [(d)] If $B_\infty(\hat{f}) \sbt \hat V$, then $K_l(\hat{V})$ is a closed, thus compact, subset of $\mathbb T_f^{(1)}$.

\, 

\item [(e)] If $B_\infty(\hat{f})\cup \Crit(\hat f^l) \sbt \hat V$ and $(\hat{V}\cup \Om(\hat f))\cap J(\hat f)$ is an open subset relative to $J(\hat{f})$, then the closed ${f}^l$--invariant set $K_l(\hat{V})$ is quasi--expanding for $\hat{f}^l$. 

\,

\item [(f)] If instead of the second hypothesis of (e) we assume more, namely that $\hat V\spt \Om(\hat f)$, then $K_l(\hat{V})$ is expanding for $\hat{f}^l$.
\end{itemize}
\eobs

\sp Sometimes, in order to be more precise, we will write $K_f(V)$ or $K_{\hat f}(\hat V)$ to indicate whether we mean the subsets of $\oc$ or $\mT_f$. 

\sp We shall prove the following.

\sp\blem\label{l1j251AA}
Let $f:\mathbb C \lra \oc$ be a non--recurrent elliptic function. If $X \sbt J(\hat f)$ is a closed (equivalently compact) $\hat f$--forward invariant expanding set, then 
$$ 
\ov{\BD}(X)<\HD(J(f)). 
$$
\elem

\bpf 
For every $j\ge 1$ put
$$
\hat W_j:=B\(B_\infty(\hat{f})\cup \Crit(\hat f)\cup \Om(\hat f),1/j\)\sbt \mT_f.
$$
Because of Observation~\ref{o1kvs3.1-2020-01-16} and Observation~\ref{o1kvs3}
the map $\hat f|_{K(\hat W_j)}:K(\hat W_j)\lra \hat K(W_j)$ belongs to $\A(K(\hat W_j))$ in the sense of Section~\ref{Volume_Lemma} (with the ambient Riemann surface $Y$ equal to $\mathbb T_f$). By Lemma~\ref{l10.3.6}, we have that
\beq\label{520131009}
s_j:=s\lt(\hat f|_{K(\hat W_j)}\rt)\le \HD(K(\hat W_j))\le \HD(J(f))=h.
\eeq
Seeking contradiction assume that $s_j=h$. Let $m_j$ be the Borel probability measure produced in Lemma~\ref{l220120827} with $X=K(\hat W_j)$ for the map $\hat f$. Since, by Observation~\ref{o1kvs3.1-2020-01-16}, the set $\hat K(W_j)$ is expanding for $\hat f$, it follows from Lemma~\ref{l2kvs1} that for every $j\ge 1$ there exists $\eta_j>0$ such that if $z\in K(W_j)$ and $n\ge0$, then there exists a unique holomorphic inverse branch 
$$
\hat f_z^{-n}:B(\hat f^n(z),2\eta_j)\lra \mT_f
$$ 
of $f^n$ defined on $B(\hat f^n(z),2\eta_j)$ that sends $\hat f^n(z)$ to $z$. So, applying in a standard way Koebe's Distortion Theorem~\ref{Euclid-I} along with Lemma~\ref{l220120827} and Lemma~\ref{l3123006}, we deduce that there exists a constant $C\in (0,+\infty)$ such that 
$$
m_j(A)\le C\hat m_h(A)
$$
for every Borel set $A\sbt K(\hat W_j)$. Thus, 
$$
\hat m_h(K(\hat W_j))\ge C^{-1}m_j(K(\hat W_j))=1/C>0,
$$
contrary to the last item of Proposition~\ref{p2_2018_08_07} as $K(\hat W_j)$ is a nowhere dense subset of $J(\hat f)$. So,
\beq\label{120131010}
s_j<h.
\eeq
Fix $u\ge 1$ so large that $X\sbt K(\hat W_u)$. Let $\hat m$ be a weak limit of the sequence $\(m_j\)_{j\ge u}$ treated as Borel probability measures on $\mT_f$, i. e.
\beq\label{1j248a}
\hat m:=\lim_{j\to\infty}m_{k_j}
\eeq
for some strictly increasing sequence $(k_j)_{j=u}^\infty$ of positive integers. Substituting the set $B(\hat f)\cup \Om(f)$ for $Y$, the considerations leading to Lemma~\ref{sl1J37} yield, with $s:=s(Y)$, the formula
\beq\label{1j248b}
\hat m(\hat f(A))=\int_A|\hat f'|^{s}d\hat m 
\eeq
for every special set $A\sbt J(\hat f)\sms \(B(\hat f)\cup \Om(\hat f)\)$, i. e. formula (b) of Lemma~\ref{sl1J37} will hold; (a), (c), and (d) will hold too. Now, since the map $f:\C\lra\oc$ is non--recurrent, it cannot have arbitrarily long chains of inverse images consisting of critical points of $\hat f$. We therefore deduce from \eqref{1j248b} that 
$$
\supp(\hat m)\spt \bu_{n=0}^\infty \hat f^{-n}(\supp(\hat m)).
$$
So, since $\supp(\hat m)\ne\es$ and since $\ov{\bu_{n=0}^\infty \hat f^{-n}(z)}= J(\hat f)$ for every $z\in J(\hat f)$, we conclude that
\beq\label{2j248b}
\supp(\hat m)=J(\hat f).
\eeq
Since $X$ is compact, we can cover it by some finitely many open balls $\{B(z,\eta_u/2)\}_{z\in E}$, where $E$ is some finite subset of $X$. 
Now fix $x\in X$, arbitrary, and a raduius $r\in(0,\eta_u]$. Since $X\sbt K(\hat W_u)$ is compact, $\hat f$--forward invariant, and disjoint from $\Om(\hat f))$, with the help of formulas \eqref{120180816}, \eqref{220180816}, and \eqref{320180816}, we deduce from  Theorem~\ref{mnt6.3} (and our choice of $\eta_u$) that there exists a least integer $n\ge 0$ such that 
\beq\label{420131014}
\Comp(x,\hat f^n,\eta_u)\sbt B(x,r).
\eeq
Then $n\ge 1$ and $\Comp(x,\hat f^{n-1},\eta_u)\not\sbt B(x,r)$. Therefore, using Koebe's Distortion Theorem (Theorem~\ref{Euclid-I}), we get
\beq\label{120131014}
r
\le \diam\(\Comp(x,\hat f^{n-1},\eta_u)\) 
\le K\eta_u|(\hat f^{n-1})'(x)|^{-1} 
\le K\eta_uD|(\hat f^n)'(x)|^{-1},
\eeq
where the number $D:=\|\hat f'|_{X}\|_\infty$ is finite since $X$ is compact and disjoint from $B(\hat f)$. Because of \eqref{2j248b} and \eqref{1j248a}, and the definition of $\e_u$, there exists $i\ge u$ (in fact all but finitely many) such that
$$
M:=\min\big\{m_{k_i}(B(y,\eta_u/2)):y\in E\big\}>0
$$
and 
$$
\hat f^k\(\Comp(x,\hat f^n,2\eta_u)\)\cap K(W_{k_i})_0=\es
$$ 
for all $k=0, 1,\ld, n$, where the set $K(W_{k_i})_0$ is understood in the sense of Section~\ref{generalconformalmeasures}. It then follows from Lemma~\ref{l220120827} that
$$
m_{k_i}\(\Comp(x,\hat f^n,\eta_u)\)
=\int_{B(\hat f^n(x),\eta_u)}|(\hat f_x^{-n})'|^{s_i}dm_{k_i}.
$$
Hence, using \eqref{420131014} and Theorem~\ref{Euclid-I} (Koebe's Distortion Theorem), we further obtain 
\beq\label{520131014}
m_{k_i}(B(x,r))
\ge m_{k_i}\(\Comp(x,\hat f^n,\eta_u)\) 
\ge K^{-s_{k_i}}\big|\(\hat f^n\)'(x)\big|^{-s_{k_i}}m(B(\hat f^n(x),\eta_u)).
\eeq
Since $f^n(x)\in X$, there exists $y\in E$ such that $\hat f^n(x)\in B(y,\eta_u/2)$. Therefore, 
$$
m_{k_i}(B(\hat f^n(x),\eta_u))\ge m_{k_i}(B(y,\eta_u/2))\ge M.
$$
Hence, \eqref{520131014} and \eqref{120131014} yield,
\beq\label{620131014}
m_{k_i}(B(x,r))
\ge MK^{-s_{k_i}}\big|\(\hat f^n\)'(x)\big|^{-s_{k_i}}
\ge M(K^{2}\eta_uD)^{-s_i}r^{s_i}.
\eeq
At this moment we can invoke Proposition~\ref{p720131014} to conclude that $\ov{\BD}(X)\le s_{k_i}$. Along with \eqref{120131010}, applied for $j=k_i$, this finally yields $\ov{\BD}(X)<h=\HD(J(f))$. The proof is complete.
\epf
 
\sp As an immediate consequence of Lemma~\ref{l1j251AA} and Observation~\ref{o1kvs3.1-2020-01-16}, we get the following.

\sp\blem\label{l1j251}
Let $f:\mathbb C \lra \ov{\mathbb C}$ be a non--recurrent elliptic function. If $\hat V\sbt \mathbb T_f$ is an open neighborhood of
$B_\infty(\hat{f})\cup \Crit(\hat f)\cup \Om(\hat f)$, then 
$$ 
\ov{\BD}\big(J(\hat f)\cap K(\hat V)\big)<\HD(J(f)). 
$$
\elem

\sp For the sake of future dealing with parabolic elliptic functions, we will need a  slightly stronger result and its slightly more general consequence. We shall prove the following. 

\sp\blem\label{l1j251P}
Let $f:\mathbb C \lra \ov{\mathbb C}$ be a non--recurrent elliptic function. If $\hat V\sbt \mathbb T_f$ is an open neighborhood of
$B_\infty(\hat{f})\cup\Crit(\hat{f}^l)\cup \Om(\hat f)$, then 
$$ 
\ov{\BD}\(J(\hat f)\cap K_l(\hat V)\)<\HD(J(f)). 
$$
for every integer $l\ge 1$.
\elem
\bpf
Because of Observation~\ref{o1kvs3.1-2020-01-16} the set 
$$
X:=J(\hat f)\cap \bu_{j=0}^{l-1}f^{lj}(K_l(\hat V))\sbt J(\hat f).
$$
is compact and $f$--forward invariant. Since $K_l(\hat V)\cap B_\infty(\hat{f})=\es$ and $K_l(\hat V)$ is $f^l$--forward invariant, we have that $X\cap B_\infty(\hat{f})=\es$. Since 
$$
K_l(\hat V)\cap \bu_{j=0}^{l-1}f^{-j}(\Crit(\hat f))=\es,
$$
we have that 
$$
X\cap \Crit(\hat f))=\es.
$$
Finally, since $\hat f(\Om(\hat f))=\Om(\hat f)$, we conclude that 
$$
X\cap\Om(\hat f)=\es.
$$
Therefore, by Observation~\ref{o1kvs3.1-2020-01-16}, the set $X$ is expanding for $f$ and it also satisfies all the requirements of Lemma~\ref{l1j251AA}. Noting that $\HD(X)=\HD\(J(\hat f)\cap K_l(\hat V)\)$ and applying this lemma, completes the proof. 
\epf

\sp As a fairly immediate consequence of this lemma, we get the following.

\sp\blem\label{l1j251R}
Let $f:\mathbb C \lra \oc$ be a non--recurrent elliptic function.
Fix an integer $l\ge 1$. If $\hat V\sbt \mathbb T_f$ is an open set containing $B_\infty(\hat{f})\cup\Crit(\hat{f}^l)$ such that $\(\hat V\cup\Om(\hat f)\)\cap J(\hat f)$ is an open neighborhood of $\Om(\hat f)$ relative to $J(\hat f)$, then 
$$ 
\HD(J(\hat f)\cap K_l(\hat V)\)<\HD(J(f)). 
$$
\elem
\bpf
By our hypotheses there exists an open set $\hat V^*\sbt \mT$ such that
$$
\hat V^*\cap J(\hat f)=\(\hat V\cup\Om(\hat f)\)\cap J(\hat f)
$$
Then $\hat V^*$ satisfies all the hypotheses of Lemma~\ref{l1j251P} and
$$
J(\hat f)\cap K_l(\hat V)
\sbt \(J(\hat f)\cap K_l(\hat V^*)\)\cup\bu_{n=0}^\infty \hat  f^{-n}(\Om(\hat f)).
$$
Therefore, remembering that the set $\bu_{n=0}^\infty \hat  f^{-n}(\Om(\hat f))$ is countable and applying Lemma~\ref{l1j251P}, we get that
$$
\begin{aligned}
\HD\(J(\hat f)\cap K_l(\hat V)\)
&\le \max\lt\{\HD\(J(\hat f)\cap K_l(\hat V^*)\),\HD\(\bu_{n=0}^\infty \hat  f^{-n}(\Om(\hat f))\)\rt\} \\
&=\HD\(J(\hat f)\cap K_l(\hat V^*)\) 
\le \ov{\BD}\(J(\hat f)\cap K_l(\hat V^*)\) \\
&<\HD(J(f)).
\end{aligned}
$$
The proof is complete.
\epf

\sp Sticking to $\hat V$, an open subset of $\mathbb T_f$, and an integer $l\ge 1$, let 
$$
{\mD}_l(\hat V)
$$ 
\fr \index{(S)}{${\mD}_l(\hat V)$}\fr denote the family of all connected components of $\mathbb T_f \sms K_l(\hat V)$. If $\hat \Ga$ is a connected component of $\hat V$, denote by $\hat \Ga_+$ the the only element of $\mD_l(\hat V)$ containing $\hat \Ga$. Let $\hat V_+$ be the union of all such components $\hat \Ga_+$. Note that
\beq\label{120181213}
K_l(\hat V_+)=K_l(\hat V)
\eeq
and 
\beq\label{220181213}
\mD_l(\hat V_+)=\mD_l(\hat V).
\eeq
We call the open set $\hat V$ dynamically maximal if 
$$
\hat V_+=\hat V.
$$
For every $W\in {\mD}_l(\hat V)$, let $k_{\hat V}(W)\geq 0$ be the
least integer $n\geq 0$ such that 
$$
\hat{f}^{ln}(W)\cap \hat V\ne\es.
$$
Equivalently:
$$
\hat{f}^{ln}(W)\cap \hat V_+\ne\es,
$$
and then 

\sp
\centerline{$\hat{f}^{lk_{\hat V}(W)}(W)$ is a connected component of $\hat V_+$.}

\sp\fr Given a component $\hat H$ of $\hat V_+$, denote
$$
\mD_l\(\hat V,\hat H\):=\big\{W\in \mD_l\(\hat V):\hat{f}^{lk_{\hat V}(W)}(W)=\hat H\big\}.
$$
We call an open connected set $\hat G\sbt\mT_f$ horizontal \index{(N)}{horizontal} of type $1$ if there exist a positive radius $R_{\hat G}>0$ and a point $z_{\hat G}\in \hat V$ such that 
$$
\hat V\sbt B\(z_{\hat G},R_{\hat G}\)
\  \  {\rm and} \  \
B\(z_{\hat G},2R_{\hat G}\)\cap {\rm PC}(\hat f)=\es.
$$
We call an open connected set $\hat V\sbt\mT_f$ horizontal of type $2$ if each connected component of $\hat f^{-1}(\hat V)$ whose closure is disjoint from $\Om(\hat f)$ is horizontal of type $1$.
Such sets exist for example if
$$
\om\notin \ov{\cO_+\(J(\hat f)\cap \Crit(\hat f)\)}.
$$
We call an open connected set $\hat G\sbt\mT_f$ horizontal if it is either horizontal of type $1$ or of type $2$. Finally, we call an open set contained in $\mT_f$ horizontal if all its connected components are horizontal.

Keeping $\hat V\sbt \mathbb T_f$ an open set, if $b\in \hat V \sms K_l(\hat V)$, then we denote by $\hat V(b)$ the connected component of $\hat V$ containing $b$. Then 
$$
\hat V_+(b):=\hat V(b)_+
$$ 
is the connected component of $\mathbb T_f \sms K_l(\hat V)$ containing $\hat V_b$ or, equivalently, containing just $b$. The following observation is immediate from the definition of nice sets defined and extensively dealt with in Section~\ref{NiceSetsGeneral}. 

\sp\bobs\label{o1_2017_11_20}
Let $f:\mathbb C \lra \oc$ be an elliptic function. If $l\ge 1$ is an integer and $\hat V \sbt \mathbb T_f$ is a pre--nice set (in particular of it is a nice set) for the map $\hat{f^l}: \hat{\mathbb T}_f^{(l)}\lra \mathbb T_f$, then $\hat V$ is a dynamically maximal open set, meaning that 
$$
\hat V_+=\hat V.
$$
Conversely, if a dynamically maximal open set $\hat V \sbt \mathbb T_f$ satisfies all the conditions (b) and (c) of Definition~\ref{d1j291-}, then it is a pre--nice set for the map $\hat{f^l}: \hat{\mathbb T}_f^{(l)}\lra \mathbb T_f$.
\eobs

\sp We record the following immediate consequence of Observation~\ref{o1kvs3.1-2020-01-16}.

\bobs\label{o1kvs3.1}
Let $f: \mathbb C \lra \oc$ be a non--recurrent elliptic function. Fix an integer $l\geq 1$ and an open set $\hat{V} \sbt \mathbb T_f$. Then:
\begin{itemize}

\item [(a)] The map $\hat f^l$ is well defined on $K_l(\hat{V})$ and
$$
\hat f^l\(K_l(\hat{V})\)\sbt K_l(\hat{V}).
$$
\item [(b)]
$$
K_l(\hat{V})\cap \bu_{n=0}^\infty\hat f^{-n}(B_\infty(\hat f))=\es.
$$
\item [(c)] $K_l(\hat{V})$ is a closed subset of $\mathbb T_f^{(1)}=\mathbb T_f\sms B_\infty(\hat f)$.

\,

\item [(d)] If $B_\infty(\hat{f}) \sbt \hat V^+$, then $K_l(\hat{V})$ is a closed, thus compact, subset of $\mathbb T_f^{(1)}$.

\, 

\item [(e)] If $B_\infty(\hat{f})\cup \Crit(\hat f^l) \sbt \hat V^+$ and $(\hat V^+\cup \Om(\hat f))\cap J(\hat f)$ is an open subset relative to $J(\hat{f})$, then the closed ${f}^l$--invariant set $K_l(\hat{V})$ is quasi--expanding for $\hat{f}^l$. 

\,

\item [(f)] If instead of the second hypothesis of (e) we assume more, namely that $\hat V^+\spt \Om(\hat f)$, then $K_l(\hat{V})$ is expanding for $\hat{f}^l$.
\end{itemize}
\eobs

\sp Although technical, the ultimate result of this section is the one formulated below. We want to add that its proof we present below is entirely different from the Przytycki and Rivera--Latelier's proof of an analogous result in \cite{PR}.

\sp\blem\label{l1j251A} 
Let $f:\mathbb C \lra \ov{\mathbb C}$ be a non--recurrent elliptic function. Fix an integer $l\ge 1$. Let $\hat V\sbt \mathbb T_f$ be an open set containing $B_\infty(\hat{f})\cup\Crit(\hat{f}^l)$ such that $\(\hat V\cup\Om(\hat f)\)\cap J(\hat f)$ is an open neighborhood of $\Om(\hat f)$ relative to $J(\hat f)$.  

If $\hat \Ga$ is a connected component of $\hat V$ such that that $\hat \Ga^+$ is horizontal and $J(\hat f)\cap\hat \Ga^+\ne\es$, then there exists a number $t \in \(\HD\(J(\hat f)\cap K_l(\hat V)\),\HD(J(\hat f))\)$, such that
$$
\sum_{W\in {\mD}_l(\hat V,\hat \Ga^+)}\diam^t(W) < +\infty.
$$
\elem

\bpf 
Since $J(\hat{f})\cap \hat{\Gamma}^{+}\neq \es$, there exists $\xi\in J(\hat{f})\cap \hat{\Gamma}^{+}$ which is a repelling periodic point  of $\hat{f}^l$. Furthermore, there exists a non--empty open set $G \sbt \hat{V}$  such that 
$$ 
G \cap \{ \hat{f}^{lj}(\xi): \,  j \geq 0\}= \es 
$$ 
and $ (G \cup \Om(\hat{f})) \cap J(\hat{f})$ is an open neighbourhood of $\Om(\hat{f})$ in $J(\hat{f})$. Then, 
$$
K_l(\hat{V})\sbt K_l(G), \  \xi \in K_l(G),
$$
and, by Observation~\ref{o1kvs3.1}, $K_l(\hat{V})\cap J(\hat{f})$ is a closed expanding $\hat{f}^l$--invariant subset of $J(\hat{f})$.  For every integer $n \geq 0$ define
$$
E_n:=\(\hat{f}^l|_{K_l(G)}\)^{-n}(\xi) \sbt K_l(G)\sms \Om(\hat{f}).
$$ 
By Observation~\ref{o1kvs3}, $ \hat{f}^l\in {\mathcal A}(K_l(G)) $ in the sense of Chapter~\ref{CMHDIM}, see particularly Section~\ref{Volume_Lemma}. For every $t \geq 0$ let $c(t)$ be the corresponding transition parameter  as defined in the first  paragraph of Section~\ref{Sullivan_Conformal_Measures}. We shall prove the following.

\sp\fr{\sc Claim~$1^0$.} The function $[0, +\infty) \ni  t \lra c(t)$ is strictly decreasing.

\,

\bpf  For every $\alpha \geq 0$ and $ \gamma \in \mathbb R$, let
$$ 
\sum(\alpha,\gamma)
:= \sum_{n=0}^\infty \sum_{x \in E_n}\big|(\hat{f}^{ln})'(x)\big|^{-\alpha}e^{-\gamma n}.
$$
Now fix $0\leq s <t$ and $u >c(s) +s-t$. It then follows from Lemma~\ref{l1kvs2} that for  every integer $n \geq 0$ and every point $ x \in E_n$, we have that
$$ 
\big|(\hat{f}^{ln})'(x)\big|^{-t}
=\big|(\hat{f}^{ln})'(x)|^{-s}|(\hat{f}^{ln})'(x)\big|^{s-t}\leq  C^{s-t}e^{\beta(s-t)n}\big|(\hat{f}^{ln})'(x)\big|^{-s}.
$$
Hence
$$ 
\begin{aligned}
\sum(t,u)
&= \sum_{n=0}^\infty \sum_{x \in E_n}|(\hat{f}^{ln})'(x)|^{-t}e^{-u n}
\leq C^{s-t}\sum_{n=0}^\infty \sum_{x \in E_n} e^{-\beta(t-s+u)n}\big|(\hat{f}^{ln})'\big|^{-s}\\
& = C^{s-t}\sum (s,u+t-s)
< +\infty,
\end{aligned} 
$$
the last inequality holding because $\mu+t-s> c(s)$. Therefore $c(t) \leq \mu$, whence $c(t) \leq c(s)+s-t< c(s)$. The proof of Claim~$1^0$ is complete.\qed

\sp Now, it follows from Lemma~\ref{l1j251R},  Lemma~\ref{l10.3.6}, and formula \eqref{120120827} that
 \beq\label{1kvs4}
 s_l \leq \HD(J(\hat f)\cap K_l(G)< \HD(J(\hat f))
\eeq
and
\beq\label{2kvs4}
c(s_l)=0,
\eeq
where
$s:=s(\hat{f}^l|_{K_l(G)})$ is defined by formula \eqref{120181217}. Fix an arbitrary 
$$
t \in \(\HD (J(f)\cap K_l(G)), \HD(J(f)\).
$$ 
It follows from Claim~$1^0$, and formulas (\ref{1kvs4}) and
(\ref{2kvs4}) that $c(t)< 0$. So, we conclude directly from the definition of $c(t)$ that
\beq\label{3kvs4}
 \sum_{n=0}^\infty \sum_{x \in E_n}\big|(\hat{f}^{ln})'(x)\big|^{-t} < +\infty.
\eeq
Now, if $ n\geq 1$ and $W\in {\mathcal D}_l(\hat{V}, \hat{\Gamma}^{+})$  is such that $k_{\hat{V}}(W)=n$, then, since $\hat{\Gamma}^{+}$ is horizontal, there exists a unique holomorphic branch 
$$
\hat f^{-lk_{\hat{V}}(W)}_W: B(z, 2 R) \lra \mathbb T_f 
$$ 
of $\hat{f}^{lk_{\hat{V}}(W)}$ such that  
$$
f^{-lk_{\hat{V}}(W)}(\hat{V}^{+})= \hat{\Gamma}^{+},
$$
where $z = z_{\hat{\Gamma}^{+}}$ and $R_{\hat{G}^{+}}$  are the parameters witnessing horizontality  of  $\hat{\Gamma}^{+}$. Let 
$$
\xi_w:=\hat{f}^{-lk_{\hat{V}}(W)}_W(z)= \hat{f}^{-ln}_W(z).
$$
Note that $\xi_W \in K_l(G)$. So, denoting
$$
{\mathcal D}^{(n)}_l(\hat{V}, \hat{\Gamma}^{+})
:=\big\{W \in {\mathcal D}_l(\hat{V}, \hat{\Gamma}^{+}):  \, k_{\hat{V}}(W)=n\big\}, 
$$
we have produced a 1--to--1 map
$$
{\mathcal D}^{(n)}_l(\hat{V}, \hat{\Gamma}^{+}) \ni W \longmapsto \xi_W\in E_n.
$$
It therefore follows from (\ref{3kvs4}) and Koebe's Distortion Theorem that
$$\begin{aligned}
\sum_{W\in {\mathcal D}_l(\hat{V},\hat{\Gamma}^{+})}\diam^t(W)
& =\sum_{n=0}^\infty \sum_{W\in {\mathcal D}^{(n)}_l(\hat{V}, \hat{\Gamma}^{+})}\diam^t(W)
 \leq K^t\sum_{n=0}^\infty \sum_{W \in {\mathcal D}^{(n)}_l(\hat{V}, \hat{\Gamma}^{+})}|(\hat{f}^{ln})'(\xi_W)|^{-t}
 \\
& \leq K^t \sum_{n=0}^{\infty}\sum_{x \in E_n}|(\hat{f}^{ln})'(x)|^{-t} \\
&< +\infty
\end{aligned} 
$$
The proof is complete.
\endpf

\sp\section[Subexpanding Functions; Stochastic Properties and Metric Entropy] {Normal Subexpanding Elliptic Functions of Finite Character; Stochastic Properties and Metric Entropy; Young Towers and Nice Sets Techniques}\label{finiteinvariantmeasures}

Throughout this section, unless stated otherwise, we assume that $f:\mathbb C \lra \oc$ is a normal subexpanding elliptic function of finite character; see Definitions~\ref{d1cef2}, \ref{d2cef2}, \ref{d5cef2}, and Theorem~\ref{t1cef5}.
\index{(N)}{normal subexpanding elliptic function of finite character} 
By virtue of Observation~\ref{o1cef2.1} we have the following.

\sp
\bobs\lab{o620200121}
Every normal subexpanding elliptic function of finite character is regular. 
\eobs  

Since the function $f$ is compactly non--recurrent, Theorems~\ref{tmaincm} and \ref{tinv} apply. Employing these theorems, Observation~\ref{o620200121}, Theorem~\ref{t1020907}, and Corollary~\ref{l6050206}, we get the following.

\sp 

\bthm\lab{t520200121} 
If $f:\C\lra\oc$ is a normal subexpanding elliptic function of finite character, then

\begin{enumerate}
\item There exist a unique $t\ge 0$ and a unique spherical $t$--conformal atomless probability measure $m_{t,s}$ for $f:J(f)\lra J(f)\cup\{\infty\}$. Then $t=h$.

\,

\item The spherical $h$--conformal measure $m_{h,s}$ coincides with a (constant) multiple of the (finite and positive) packing measure $\Pi^h_s$ on $J(f)$.

\, 

\item There exists a unique probability $f$--invariant measure $\mu_h$ absolutely continuous with respect to $m_h$. In addition, $\mu_h$ is atomless and equivalent to $m_h$.

\,

\item The $f$--invariant Borel probability measure $\mu_h$ is metrically exact, in particular ergodic\index{(N)}{ergodic measure} and conservative\index{(N)}{conservative measure}, 

\,

\item All other conformal measures are purely atomic, supported on $\Sing^-(f)$ with exponents larger than $h$. 

\,
 
\item $\mu_h(\Tr(f))=1$. 
\end{enumerate} 
\ethm

\sp 
This section, based on previous sections, and heavily using the concept and results about Nice Sets of Section~\ref{NiceSetsGeneral}, is devoted  to explore refined  stochastic  properties
of the dynamical system $\(f|_{J(f)},\mu_h\)$. Employing the existence of Nice Sets and their giving rise to conformal maximal graph directed Markov systems of Chapter~\ref{Markov-systems} (all done in Section~\ref{NiceSetsGeneral}), we will show that the Young's 
tower abstract approach, described in Chapter~\ref{finite-measure},  Probability (Finite) Invariant Measures; Finer Properties, applies, and we will derive lots of dynamical and stochastic consequences from it such as the Central Limit Theorem, the Law of Iterated Logarithm and Exponential Decay of Correlations.  

\sp \subsection{Geometric and Dynamical Preparations for Young's Towers of Finite Character Elliptic Functions}\label{Pre YT}

This subsection is motivated by the seminal paper \cite{PR}. The
absolutely crucial and necessary concept for us in this section is the one of Nice Sets from Section~\ref{NiceSetsGeneral}. We will work with the torus $\mathbb T_f$ phase space, and eventually we will transfer the fruits of this work back to $\C$ again. In particular, we will prove strong regularity of iterated function systems induced by nice sets. As the last preparation for applications of Young's towers results and techniques (see Chapter~\ref{finite-measure}) we will prove that the measures of points returning in time $n$ to a nice set decay exponentially fast with $n$. At the very end of this section we will
also prove that the Kolmogorov--Sinai metric entropy $\h_{\mu_h}$ is finite.

\sp Given an elliptic function $f:\C \lra \oc$, let 
\beq\lab{120200122}
B(f):= \Crit(f)\cup f^{-1}(\infty)
\eeq
and 
\beq\lab{220200116}
B(\hat f)
:=\Pi_f(B(f))
=\Pi_f\(f^{-1}(\infty)\cap\Crit(f)\)
=B_\infty(\hat f)\cup \Crit(\hat f),
\eeq 
where $B_\infty(\hat f)=\Pi_f(f^{-1}(\infty))$ has been defined in \eqref{1_2017_09_20}. From now on throughout this section we assume that $f:\C \lra\oc$ is a normal subexpanding elliptic function of finite character subexpanding\index{(N)}{subexpanding}. According to Theorem~\ref{t2cef5}
\beq\lab{3_2017_11_18}
R_f:=\frac{1}{2} \dist_{\mathbb T}\(B(\hat{f}), J(\hat f)\cap {\rm PC}(\hat{f})\)>0,
\eeq

\sp As an immediate consequence of Lemma~\ref{l1j251A}, we get the following. 

\sp\blem\label{l1j251M} 
Let $f:\mathbb C \lra \oc$ be a normal subexpanding elliptic function of finite character. If $\hat V \sbt \mathbb T_f$, a neighborhood of
$B(\hat{f})$, is a nice set for the map $\hat{f}: \hat{\mathbb T}_f\lra \mathbb T_f$, then there exists a number $t \in \(\ov{\BD}\(J(\hat f)\cap K(\hat V)\),\HD(J(f))\)$, such that
$$
\sum_{W \in {\mathcal D}(\hat V)}\diam^t(W) < +\infty,
$$
where, we recall, ${\mathcal D}(\hat V)$ denotes the family of all connected components of $\mathbb T_f \sms K(\hat V)$.
\elem

\sp Recall that the concept of nice sets introduced in Section~\ref{NiceSetsGeneral} concerned all analytic maps of Riemann surfaces whose range is either a Euclidean torus or a complex plane, and we have proved there the existence (especially see Corollary~\ref{c120180827}) of nice sets for all such maps satisfying some further hypotheses. Now, with 
$$
Y:=\mathbb T_f \  \  \text{ and }  \  \ X:=\hat\mT_f,
$$
we want to apply this theorem to the map 
$$
\hat{f}:\hat{\mathbb T}_f\lra\mathbb T_f,
$$
and the set 
$$
F:= B(\hat{f}).
$$ 

We shall prove the following.

\bthm\label{t120180829}
If $f:\C\lra\oc$ is a normal subexpanding elliptic function of finite character, then for $F=B(\hat{f})$ and every radius $r>0$ sufficiently small fulfilling condition (1) of Corollary~\ref{c120180827}, there exists a corresponding nice set 
$$
\hat V:=U_r
$$
for the quotient map $\hat{f}:\hat{\mathbb T}_f\lra\mathbb T_f$. 
\ethm

\bpf
We just are two check that condition (2) (with a fixed $\l$ and $\ka$ larger than $1$) of Corollary~\ref{c120180827} holds for all $r\in(0,\hat\g]$ sufficiently small, particularly satisfying condition (1) of Corollary~\ref{c120180827}, where $\hat\g>0$ comes from Corollary~\ref{c1j249} (Exponential Shrinking on $\mathbb T_f$). Indeed, if $r\in(0,\hat\g]$ if sufficiently small, then for every $z\in B\(B(\hat f),2r\)$, we have that 
$$
\big|\hat f'(z)\big|\ge r^{-1}.
$$
Hence, for every integer $n\ge 1$, every $\xi\in\hat f^{-n}(B(\hat f))$ such that
$$
\Comp\(\xi,\hat f^n,2r\)\cap B(\hat{f})\ne\es,
$$
and evey $z\in\Comp\(\xi,\hat f^n,2r\)$, we get by virtue of Corollary~\ref{c1j249} and Koebe's Distortion Theorem, that
$$
\big|(\hat f^n)'(z)|
=\big|\hat f'(z)\big|\cdot\big|(\hat f^{n-1})'(\hat f(z))|
\ge (8Kr)^{-1}e^{M_fn}
\ge (8Kr)^{-1}.
$$
Since 
$$
(8Kr)^{-1}\ge K\max\lt\{\frac{2 \kappa}{\kappa-1}, \lambda,\rt\}
$$
for all $r>0$ small enough, we are done.
\epf

\sp With having $\ka, \l>1$ fixed, we say that the, resulting from 
Corollary~\ref{c120180827}, nice set 
$$
\hat V:=U_r \sbt \mathbb T_f
$$
for the analytic map $\hat{f}:\hat{\mathbb T}_f\lra\mathbb T_f$ is small \index{(N)}{small nice set} (so in particular $r>0$ satisfies formula (1) from the hypotheses of Corollary~\ref{c120180827}) if and only if  
\beq\lab{2_2017_11_20}
\max\{ \diam(\hat V_b): \, b \in B(\hat{f})\} <
\frac{1}{6}\min\{\g_f, \b_f, R_f, \eta\}, 
\eeq
where $\b=\b_f>0$ was defined in \eqref{1042206}, $\g_f$ is determined by \eqref{dp5.4}, while $\eta>0$ is taken so small that the projection map $\Pi_f:\mathbb C \lra \mathbb T_f$ is one--to--one if restricted  to any ball $B(z,2\eta)$, $ z \in \mathbb C$. The sets $\hat V_b$ denote here the respective connected components of $V$ containing points $b$. Let 
$$
J_{\hat V}:=J_{\mathcal S_{\hat V}}
$$ 
be the limit set of the maximal conformal graph directed Markov system $\cS_{\hat V}$, produced by Theorem~\ref{t1j290a}. Clearly, $\cS_{\hat V}$ is finitely irreducible and $B(\hat{f})$ is the set of its vertices. We denote the set of all edges of the system $\mathcal S_{\hat V}$ be $\hat E$. For every $e\in\hat E$ we denote by 
\beq\label{1_2018_0515}
N(e)
\eeq
the unique integer $\ge 1$ such that 
$$
\hat\phi_e:\hat V_{t(e)}\lra \hat V_{i(e)},
$$
the element of $\cS_{\hat V}$ corresponding to the edge $e$, is a local holomorphic branch of $f^{-N(e)}$. With the terminology of Proposition~\ref{p2j290} (a), $N(e)=n(\hat V_{i(e)})$. We further define 
\beq\label{120140102}
N(\om):=\sum_{j=1}^{|\om|}N(\om_j).
\eeq
In other words, $N(\om)$ is the unique integer $\ge 1$ such that $\phi_\om:\hat V_{t(\om)}\lra \hat V_{i(\om)}$ is a local holomorphic branch of $f^{-N(\om)}$. First, we shall prove the following.

\sp\blem\label{l120131004} 
Let $f:\C\lra\oc$ be a normal subexpanding elliptic function of finite character.
If $\hat V\sbt\mT$ is a small nice set, produced in Theorem~\ref{t120180829},  for the map $\hat{f}: \hat{\mathbb T}_f\lra \mathbb T_f$, then  
$$
\hat m_h(J_{\hat V})>0 \  \  \  \text{ and }  \  \  \  \HD(J_{\hat V})=\HD(J(f)), 
$$
where, we recall, the measure $\hat m_h$ has been defined in \eqref{720140104}.
Furthermore, the maximal graph directed Markov system $\cS_{\hat V}$ is regular (in the sense of Definition~\ref{d5_2017_11_18}) and $\hat m_h/\hat m_h(J_{\hat V})$ is the $h$--conformal measure (in the sense of Theorem~\ref{t420131007} with $t=h$) for the system $\cS_{\hat V}$. 
\elem

\bpf
By the construction of the system $\cS_{\hat V}$, provided in Theorem~\ref{t2j291}, and by the definition of the limit set $J_{\hat V}$, the latter contains the set of all transitive points of the map $\hat f:\hat J(f)\cap\hat\mT_f\to \hat J(f)$, denoted by $\Tr(\hat f)$, that belong to $\hat V$. Since $\hat m_h$ is of full topological support and since $\Pi(\Tr(f))\spt\Tr(\hat f)$, it therefore follows from Theorem~\ref{tmaincm} that 
$$
\hat m_h(J_{\hat V})\ge \hat m_h(\hat V\cap\Tr(\hat f))=\hat m_h(\hat V)>0.
$$
Since the measure $\hat m_h$ is $h$-conformal for the map $\hat{f}: \hat{\mathbb T}\to \mathbb T$, it also satisfies condition (f) of Theorem~\ref{t420131007} with $t=h$ and $\g=1$. It therefore follows from this theorem that $\P(h)=0$ and $\hat m_h/\hat m_h(J_{\hat V})$ is the $h$--conformal measure for the conformal graph directed system $J_{\hat V}$. Consequently, the system $\cS_{\hat V}$ is regular and (the same theorem) $h=h_{\cS_{\hat V}}$. But, by Theorem~\ref{t1j97}, $h_{\cS_{\hat V}}=\HD(J_{\hat V})$. The proof is thus complete.
\epf 

\sp Recall that the concept of strongly regular conformal graph directed systems was defined in Section~\ref{sec-top-press}, Topological Pressure, $ \theta$-Number, and Bowen's Parameter. The crucial, technically  involved, intermediate fact that opens up the doors for the stochastic properties of the dynamical system $\(f|_{J(f)},\mu_h\)$, is this.

\sp\blem\label{l1j257} 
Let $f:\C\lra\oc$ be a normal subexpanding elliptic function of finite character. If $\hat V \sbt  \mathbb T_f$, a neighborhood of $B(\hat{f})$, is a small nice set, produced in Theorem~\ref{t120180829}, for the analytic map $\hat{f}:\hat{\mathbb T}_f
\lra \mathbb T_f$, for example the one proved in Theorem~\ref{t120180829}, then  
$$
{\mathcal S}_{\hat V}=\{\phi_e\}_{e \in \hat E},
$$ 
the corresponding maximal conformal graph directed Markov system, is strongly regular in the sense of Definition~\ref{d5_2017_11_18}.   
\elem

\bpf  As always let $ q\geq 1 $ be the maximum of the orders of all poles of $f$. Fix an arbitrary number $t\in
\(\ov{\BD}(J(\hat f)\cap K(\hat V)),\HD(J(f))\)$ produced in 
Lemma~\ref{l1j251M}. Since any number in $(t,\HD(J(f)))$ is also good for this lemma, we may assume, because of Theorem\ref{thm:julia}, that
\beq\label{120181220}
t\in\lt(\max\Big\{\ov{\BD}\(J(\hat f)\cap K(\hat V)\),\frac{2q}{q+1}\Big\}, \HD(J(f))\rt).
\eeq
For every $ a \in B(\hat{f})$ let
$$ 
\hat E_a:=\{e \in \hat E: i(e)=a\}=\{e \in \hat E: \hat \phi_e(\hat V_{t(e)})\sbt \hat V_a\}.
$$
Of course  all the sets  $\hat E_a$, $a \in \hat E$, are mutually  disjoint, and their union is equal to $\hat E$. If  $b \in
\Pi_f(f^{-1}(\infty))$ and $e \in \hat E_b$, then
$$ 
\hat{f}(\hat \phi_e(\hat V_{t(e)})) \in {\mathcal D}(\hat V)
$$
and 
$$ 
\diam(\hat \phi_e(\hat V_{t(e)})) 
\leq \diam(\hat{f}(\hat \phi_e(\hat V_{t(e)})))
$$
provided that the diameters of the elements of the nice set $\hat V$ are sufficiently small, so that the sets $\hat \phi_e(\hat V_{t(e)})$, $e \in \hat E_b$, lie sufficiently close to the pole $b$. For every  $ b \in \Pi_f(f^{-1}(\infty))$ the map
$$ 
\hat E_b\ni e \longmapsto  \hat{f}\(\hat \phi_e(\hat V_{t(e)})\)\in {\mathcal D}(\hat V) \
$$
is at most $q$--to--$1$. So,
\beq\label{1j257}
\sum_{b \in \Pi_f(f^{-1}(\infty))}\sum_{e \in \hat E_b}\diam^t(\hat \phi_e(\hat V_{t(e)}))
\leq \!\!\!\sum_{b \in\Pi(f^{-1}(\infty))}\sum_{e \in \hat E_b}
\diam^t(\hat{f}(\hat \phi_e(\hat V_{t(e)}))) 
\leq q\!\!\! \sum_{W \in {\mathcal D}(\hat V)}\!\!\!\! \diam^t(W),
\eeq
the last  estimate  being very crude but sufficient for  us. By
Corollary~\ref{c1j249} and \eqref{3_2017_11_18}, there exists $0 < \d_R \leq  \hat{\g}$, the latter coming from this corollary, such that 
\beq\label{1j259}
 \diam(\Comp(z, \hat{f}^n,\d_R))< R_f/2  
 \eeq
for all $n \geq 0$ and all $z \in \hat{f}^{-n}(\hat J(f))$,
and the function  $\hat{f}$ is injective  on every ball $B(z,\d_R)$ if $ z \in \rm{PC}(\hat{f})$. Now, by Lemma~\ref{l1j251M} for example, there exists an integer $l \geq 1$ so large that $\diam(W)< \d_R/4$ for all $W \in {\mathcal D}(\hat V)$ with $n(W)\geq l$. Since $\ov{\hat V}\cap \ov{\rm PC}(\hat{f})=\es$ and since the set $\ov{\rm PC}(\hat{f})$ is forward invariant, meaning that $f(\ov{\rm PC}(\hat{f}))\sbt \ov{\rm PC}(\hat{f})$, we have that
$$
\ov{\rm PC}(\hat{f})\cap \bigcup\{\ov{W}: \, W \in {\mathcal D}(\hat V)\}= \es.
$$ 
Since also the family $\{W \in {\mathcal D}(\hat V):n(W)=l\}$ is finite, we conclude that
$$  
\a := \dist\left(\ov{\rm PC}( \hat{f}), \bigcup\{\ov{W}: \, W \in {\mathcal D}(\hat V)\, \, \mbox{and}\,\, n(W)=l\}\right)
$$ 
is positive. Let 
$$
\kappa:= \frac{1}{4}\min\{\b_f, \hat\g,\a, R_f,\delta_R\}.
$$ 
Now fix $ c \in \Crit (\hat{f})$ and  $ e \in  \hat E_c$
with $ N(e) \geq l+1$. Since $ \a \ge\ka$ and since the set $\ov{\rm PC}(\hat{f})$ is forward invariant, there exists a least integer
$ 1 \leq k_e \leq N(e)-l$ such that
\beq\label{2j259}
\dist(\hat{f}^{k_e}(c), \hat{f}^{k_e}(\hat \phi_e(\hat V_{t(e)})))\geq \kappa.
\eeq
Fix a  arbitrary point  $\xi_e \in \hat \phi_e(\hat V_{t(e)})$ and put
$$
w_e:=\hat{f}^{k_e}(\xi_e).
$$
Consider the ball $B(w_e,\d_R)$. Since
$\hat{f}^{k_e}\(\hat \phi_e(\hat V_{t(e)})\)\in {\mathcal  D}(\hat V)$  and
$n(\hat{f}^{k_e}(\hat \phi_e(\hat V_{t(e)})))=N(e)- k_e\geq l,$ by the  choice
of  $l$, we have that $\hat{f}^{k_e}\(\hat\phi_e(\hat V_{t(e)})\) \sbt B(w_e,\d_R/4)$. Consequently,
\beq\label{1j261}
\hat\phi_e(V_{t(e)}) \sbt \Comp   ( \xi_e, w_e, \hat{f}^{k_e}, \d_R/4).
\eeq
It follows from  the definition of $R_f$, (\ref{1j259}) and the choice
of $k_e$, that for all $1\leq k \leq k_e-1$, we  have $$
\begin{aligned}
\dist( \Crit(\hat{f}), &\hat{f}^k( \Comp( \xi_e, \hat{f}^{k_e},
\d_R)))\\
     &\geq \dist( \Crit(\hat{f}),
\hat{f}^k(c))-\dist(\hat{f}^k(c), \hat{f}^k(\Comp(\xi_e, w_e,
\hat{f}^{k_e}, \d_R)))\ge \\
& \geq 2R_f - \dist(\hat{f}^n(c), \hat{f}^k(\hat \phi_e(\hat V_{t(e)})))-\diam(
\hat{f}^k(\Comp( \xi_e, w_e, \hat{f}^{k_e}, \d_R)))\\
&> 2R_f-\kappa -R_f/2 
\geq \frac{5}{4}R_f>0.
\end{aligned}
$$ 
In particular 
$$\Crit(\hat{f}) \cap
\hat{f}^k( \Comp(\xi_e, \hat{f}^{k_e}, \d_R))=\es
$$ for all
$1\leq k \leq k_e-1$. Thus, there exists  a well--defined unique
holomorphic inverse branch 
$$
\hat{f}^{-(k_e-1)}_e: B(w_e, \d_R) \lra
\mathbb C
$$ 
of $\hat{f}^{k_e-1}$ sending $w_e$ to $ \hat{f}(\xi_e)$.
It then  directly follows  from Koebe's Distortion Theorem, i. e. Theorem~\ref{Euclid-I}, applied to the map  $\hat{f}^{-(k_e-1)}_e: B(w_e, \d_R) \to
\mathbb C$, (\ref{1j261}), and
Corollary~\ref{c1j249} (remember that $\d_R< \hat{\g})$ that
\beq\label{2j261} 
\diam\(\hat{f}^{k_e}(\hat \phi_e(\hat V_{t(e)}))\)\cdot|(\hat{f}^{k_e-1})'(
\hat{f}(\xi_e))|^{-1}
\comp_C\diam\(\hat{f}(\hat \phi_e(\hat V_{t(e)}))\)
\leq 8 e^{-M_f(k_e-1)},
\eeq
with some comparability constant $C \geq 1$ independent of $e$. Since $\kappa < \d_R/4$, we  have that $ B(w_e, \kappa) \sbt B(w_e,
\d_R/4)$, and it therefore follows from Koebe's
$\frac{1}{4}$-Theorem, i .e. Theorem~\ref{one-quater}, that
$$
\hat{f}^{-(k_e-1)}_e(B(w_e, \kappa))
\supset B\lt(\hat{f}(\xi_e),
\frac{1}{4}|(\hat{f}^{k_e-1})'( \hat{f}( \xi_e))|^{-1}\rt).
$$ 
Since also, by (\ref{2j259}), $\hat{f}^{k_e}(c)\notin B(w_e, \kappa)$, we
conclude that 
$$
\hat{f}(c) \notin B\lt(\hat{f}(\xi_e), \frac{1}{4}|
(\hat{f}^{k_e-1})'( \hat{f}( \xi_e))|^{-1}\rt).
$$ 
This means that $|\hat{f}(c)-\hat{f}( \xi_e))| \geq
\frac{1}{4}|(\hat{f}^{k-1}_e)'( \hat{f}( \xi_e))|^{-1}$. Since
$\xi_e$ was an arbitrary point in $\phi_e(V_{t(e)})$, this implies
that
 $$
 \dist(\hat{f}(c), \hat{f}( \hat \phi_e(\hat V_{t(e)})))\geq \frac{1}{4}
|(\hat{f}^{k_e-1})'( \xi_e)|^{-1}.
$$ So,
$$ | \hat{f}'(z)| \geq A^{-1}\left(
\frac{1}{4}|(f^{k_e-1})'(\hat{f}(\xi_e))|\right)^{\frac{q_c-1}{q_c}}=A^{-1}{4}^{\frac{1-q_c}{q_c}}|(f^{k_e-1})'(\hat{f}(\xi_e))|^{\frac{1-q_c}{q_c}}$$
for all $z \in \hat  \phi_e(\hat V_{t(e)})$. Combining  this and
(\ref{2j261}), we get  that, 
\beq\label{2j263}
 \begin{aligned}
 \diam(\hat \phi_e( \hat V_{t(e)}))&\leq
 \diam( \hat{f}(\hat \phi_e( \hat V_{t(e)}))) A
4^{\frac{q_c-1}{q_c}}|(f^{k_e-1})'(\hat{f}(\xi_e))|^{\frac{q_c-1}{q_c}}\\
 & \leq A 4^{\frac{q_c-1}{q_c}} \diam
 (\hat{f}^{k_e-1}(\phi_e(V_{t(e)})))
|(f^{k_e-1})'(\hat{f}(\xi_e))|^{-1 +\frac{q_c-1}{q_c}}\\
& \leq 4AC |(f^{k_e-1})'(\hat{f}(\xi_e))|^{
-\frac{1}{q_c}}\diam(\hat{f}^{k_e-1}(\hat \phi_e(\hat V_{t(e)})))\\
&\leq  C_1 e^{-\frac{M_f}{q_c}k_e}\diam(\hat{f}^{k_e-1}(\hat \phi_e(
\hat V_{t(e)})))
\end{aligned}
\eeq 
with some constant $C_1\geq
1$ independent of $ e \in \hat E$. Let 
$$
\hat E^l_c:=\{e \in \hat E_c: \,  N(e)\geq l+1\}.
$$ 
For all $ e \in \hat E_c^l$,  put 
$$
W_e:= \hat{ f}^{k_e-1}(\hat\phi_e\(\hat V_{t(e)})\).
$$
We shall  prove the following.

\sp\bcla  
For every $ c \in  \Crit(\hat{f})$, the  function  $\hat E^l_c \ni
e \longmapsto ( k_e, W_e)$  is at most  $q_c$--to--$1$. 
\ecla 

\bpf Fix $c \in \Crit(\hat{f})$ and suppose for a
contrary that there are two elements $a, b \in \hat E_e$ such that $ a
\neq b$ and $ k:=k_a=k_b$, $W_a=W_b$. Suppose also that
$\hat{f}\(\hat\phi_a(\hat V_{t(a)})\)\neq \hat{f}\(\hat\phi_b(\hat V_{t(b)})\)$. Since $W_a=W_b$ there thus exists  a least $ 2 \leq j \leq k$ such
\beq\label{1j263}
\hat{f}^j\(\hat\phi_a(\hat V_{t(a)})\)= \hat{f}^j\(\hat\phi_b(\hat V_{t(b)})\).
\eeq
Then
\beq\label{1j265}
\hat{f}^{j-1}\(\hat\phi_a(\hat V_{t(a)})) \cap \hat{f}^{j-1}(\hat\phi_b(\hat V_{t(b)})\)
=\es.
\eeq
By the definition of $k_a$ and $k_b$, and since
$n(\hat{f}^{j-1}(\phi_a(\hat V_{t(a)}))),
n(\hat{f}^{j-1}(\phi_b(\hat V_{t(b)})))\geq l$  we have
\beq\label{2j265}
\begin{aligned}
\hat{f}^{j-1}(\hat\phi_a&(\hat V_{t(a)})) \cup
\hat{f}^{j-1}(\hat\phi_b(\hat V_{t(b)})) \sbt \\
&\sbt B( \hat{f}^{j-1}(c), \kappa+\max\{
\diam(\hat{f}^{j-1}(\hat\phi_a(\hat V_{t(a)}))),
\diam(\hat{f}^{j-1}(\hat\phi_b(\hat V_{t(b)})))\})\\
& \sbt B(\hat{f}^{j-1}(c), \kappa+ \d_R/4) \\
&\sbt B(\hat{f}^{j-1}(c),\d_R/2)\sbt B(\hat{f}^{j-1}(c),\d_R).
\end{aligned}
\eeq
But, by the choice  of $\d_R$, the function $\hat{f}$ restricted to
the ball $B(\hat{f}^{j-1}(c), \d_R)$ is one to one.  This  however
contradicts (\ref{2j265}), (\ref{1j265}), and
 (\ref{1j263})  taken together. Along with  the observation  that the
function $\hat E_c \ni e  \mapsto \hat{f} (\hat\phi_e(V_{t(e)}))$ is at most
$q_c$--to--$1$, this completes the proof  of the claim.
\endpf

\sp Let 
$$ 
q := \max \{ q_c:  c \in \Crit(\hat{f})\}.
$$ 
Applying this claim and (\ref{2j263}), we act for any $c\in\Crit(\hat{f})$ that

\beq\label{3j265} 
\begin{aligned}
\sum_{e \in \hat E_c^l} \diam^t(\hat\phi_e(\hat V_{t(e)}))
& \leq C_1^t \sum_{e \in \hat E^l_c} e^{-\frac{M_f}{q}tk_e} \diam^t(W_e)
\leq q_c C_1^t \sum_{W \in {\mathcal D}(\hat V)} \diam^t(W)
\sum_{k=1}^{\infty}e^{-\frac{M_f}{q}tk}\\
& \leq  C_1^tq Q_t \sum_{W \in {\mathcal D}(\hat V)} \diam^t(W).
\end{aligned}
\eeq
Since the set $\bigcup_{c \in \Crit(f)}\{e \in \hat E_c: N(e) \leq l\}$
is finite, the sum
$$
\sum{}_l^{(t)}
:= \sum_{c \in \Crit(\hat{f})}\sum_{e\in \hat E_c: N(e)
\leq l} \diam^t(\hat\phi_e(\hat V_{t(e)}))
$$ 
is also finite. Combining this along with (\ref{3j265}), (\ref{1j257}), and finally with Lemma~\ref{l1j251M}, we get the following.
$$
\begin{aligned} 
\sum_{e \in \hat E} & \diam^t(\hat\phi_e(\hat V_{t(e)})) = \\
&= \sum_{b \in \Pi_f(
 f^{-1}(\infty))}\sum_{e \in \hat E_b} \diam^t\(\hat\phi_e(\hat V_{t(e)})\)+\sum{}_l^{(t)}      
  +  \sum_{c \in \Crit(\hat{f})}\sum_{e\in \hat E_c}
 \diam^t\(\hat\phi_e(\hat V_{t(e)})\)   \\
& \leq q \sum_{W \in {\mathcal D}(\hat V)} \diam^t(W)
+\sum{}_l^{(t)}+\sum_{c \in \Crit(\hat{f})} C^t_1q Q_t\sum_{W \in
{\mathcal D}(\hat V)} \diam^t(W)\\
&=q(1+C^t_1Q_t)\sum_{W \in {\mathcal D}(\hat V)} \diam^t(W)+
\sum{}_l^{(t)} < + \infty.
\end{aligned}
$$
Since $t < \HD(J(f))=\HD(J_r)$, this implies that our system is
strongly regular. The proof of Lemma~\ref{l1j257} is complete. 
\endpf

\sp Now we want to form an appropriate strongly regular conformal maximal graph directed Markov system on $\mathbb C$ and to obtain from it the
stochastic laws such as the law of iterated logarithm, and, by
Young's tower  construction, the exponential decay of correlations,
and the Central Limit Theorem. It will also allow us to show that the metric entropy $\h_{\mu_h}(f)$ is finite.

\sp Keep $\hat V\sbt \mT_f$, a small nice set, produced in Theorem~\ref{t120180829}, for the map $\hat f:\hat\mT_f\lra\mT_f$, satisfying \eqref{2_2017_11_20}. As the set of vertices take $\cV$, a (finite) selector of the partition of 
$$
\Pi_f^{-1}(B(\hat f))=f^{-1}(\infty)\cup (\Crit(f)\cap J(f))
$$ 
into equivalence classes of the equivalence relation $\sim_f$, i.e. choose exactly one point from each such equivalence class intersected with $f^{-1}(\infty)\cup (\Crit(f)\cap J(f))$. For every $v\in\cV$ let
$V_v$ be the connected component of $\Pi_f^{-1}(\hat V_v)$ containing $v$.
Because of \eqref{2_2017_11_20} for each $v\in\cV$ the holomorphic map 
$$
\Pi_v:=\Pi_f\big|_{V_v}:V_v\lra \hat V_v
$$ 
is 1--to--1. As the set of edges, set
$$
E:=\big\{e\in\hat E:\exists v\in\cV \ \ {\rm s. t. } \ f^{N(e)-1}\(\^f(\hat V_{t(e)})\)=V_v\big\}, 
$$
where, we recall, the map $\^f:\mT_f\to\oc$ is given by \eqref{120180816}.
For each $e\in E$ let $t(e)$ be the only element of $\cV\cap\Pi_f^{-1}(t(e))$, where the latter $t(e)$ is understood as associated with the graph directed Markov systems $\cS_{\hat V}$. Likewise, $i(e)$ is the only element of $\cV\cap\Pi_f^{-1}(i(e))$, where the latter $i(e)$ is, as before, understood as associated with the graph directed Markov systems $\cS_{\hat V}$.
For each $e\in E$ define
$$
\phi_e:=\Pi_{i(e)}^{-1}\circ\hat\phi_e\circ\Pi_{t(e)}:V_{t(e)}\lra V_{i(e)}.
$$
Obviously 
$$
\cS_V:=\{\phi_e:V_{t(e)}\lra V_{i(e)}\}_{e\in E}
$$
is a conformal maximal graph directed Markov system with eadges $E$ and vertices $\cV$. We shall easily prove the following. 

\bthm\label{t120180509} 
Let $f:\C\lra\oc$ be a normal subexpanding elliptic function of finite character.
If $\hat V\sbt\mT$, a neighborhood of $B(\hat{f})$, is a small nice set, produced in Theorem~\ref{t120180829}, for the map $\hat{f}: \hat{\mathbb T}_f\to \mathbb T_f$, then  
$$
m_h(J_V)>0 \  \  \  \text{ and }  \  \  \  \HD(J_V)=\HD(J(f)), 
$$
where, we recall, $m_h$ is the $h$--conformal measure for the dynamical system $f:\C\to\oc$.

Furthermore, the maximal graph directed Markov system $\cS_V$ is strongly  regular (in the sense of Definition~\ref{d5_2017_11_18}) and $$
m_h/m_h(J_V)
$$ 
is the $h$--conformal measure (in the sense of Theorem~\ref{t420131007} with $t=h$) for the system $\cS_V$. 
\ethm 

\bpf
All but strong regularity follows in the same way as in the proof of Theorem~\ref{l120131004} once we realize that, by the construction of $\cS_V$, $J_V\spt \Tr(f)\cap V$. 

Strong regularity of the system $\cS_V$ is an immediate consequence of the facts that the system  $\cS_{\hat V}$ is strongly regular (see  Lemma~\ref{l1j257}), $\P_{S_V}(h)=\P_{S_{\hat V}}(h)=0 $, $\P_{S_V}(t)\le \P_{S_{\hat V}}(t)$ for all $t\ge 0$, the latter holding because $E\sbt \hat E$ and all the maps $\Pi_v$, $v\in V$, are isometries, and both maps $\(\th_{S_{\hat V}},+\infty\)\ni t\longmapsto \P_{S_V}(t), \P_{S_{\hat V}}(t)$ are strictly decreasing.
\epf

\sp The first consequence of Theorem~\ref{t120180509} is the following. 

\sp\blem\label{l1j267} 
Let $f:\C\lra\oc$ be a normal subexpanding elliptic function of finite character. 
If $\hat V\sbt \mathbb T_f$, a neighbourhood of
$B(\hat{f})$, is a small nice set, produced in Theorem~\ref{t120180829}, for $ \hat{f}: \hat{\mathbb T}_f\lra
\mathbb T_f$, then there exist $C>0$ and $\a>0$ such that 
$$
m_h\left( \bigcup_{N(e)\geq n}\phi_e\(V_{t(e)}\)\right)
\leq C e^{-\a n}
$$ 
for all $n \geq 1$, where, we recall, $m_h$ is the unique $h$-conformal
measure for the map $f:J(f)\to J(f)$. 
\elem

\bpf  As $ \HD(J_r)=h$, Theorem~\ref{t120180509} produces a number
$ t \in (0, h)$ such that
$$
Z_1(t)=\sum_{e\in E}\|\phi_e'\|^t< +\infty.
$$
So, employing Theorem~\ref{t1j237} and Koebe's Distortion Theorem, i. e. Theorem~\ref{Euclid-I}, we  get for every $n \geq 1$
that
$$ 
\begin{aligned}
m_h\left( \bigcup_{N(e)\geq n}\phi_e\(V_{t(e)}\)\right)
&=\sum_{N(e)=n}  m_h(\phi_e(V_e)) 
\comp \sum_{N(e)=n}\|\phi_e'\|^h\\
&=\sum_{N(e)=n}\|\phi_e'\|^t\|\phi_e'\|^{h-t} \\
&\lek \sum_{N(e)=n}\|\phi_e'\|^t e^{-M_f(h-t)n}\\
&=e^{-M_f(h-t)n}\| \sum_{N(e)=n}\|\phi_e'\|^t  \\
&\leq Z_1(t)e^{-M_f(h-t)n}.
\end{aligned}
$$
We are done. 
\endpf

\sp\subsection{Young's Towers for Subexpanding Elliptic Functions of Finite Character}\label{Young's Towers for Subexpanding}

\sp In this subsection we take fruits of the results obtained in the
previous subsection and in Section~\ref{Young-Abstract} in order to establish
the exponential decay of correlations, the Central Limit Theorem, and
the Law of Iterated Logarithm for a normal subexpanding elliptic function of finite character, with respect to the invariant measure $\mu_h$
equivalent to the $h$-conformal measure $m_h$.

\sp With the setting of the previous section, we define the map $F:J_V\to J_V$ by the formula:
$$
F_V(\phi_e(z))=z \  \text{ if }  \  e\in E \ 
 \text{ and }  \  z\in J_V.
$$
Our goal is to show that the system $\(F,\nu_h:=m_h|_{J_V}\))$ can be embeded as the base into an  appropriately defined Young's Tower fitting into the abstract 
framework of Section~\ref{Young-Abstract}. Then to check that the
hypotheses of Theorems \ref{lsyoung0}, \ref{lsyoung}, and \ref{LILII} are  
satisfied for this tower. Finally to show that the original map $f:J(f)\to J(f)$ is a measure--theoretic factor of the tower dynamical system. The stochastic laws will then follow almost automatically. 
For every $b\in B(\hat f)$ set
$$
J_b:=V_b\cap J_V.
$$
Our Young tower $\cY_f$ is constructed as follows.

\begin{enumerate}
\, \item The space $\Delta_0$ is now $J_V$, the limit set of the iterated
function system ${\mathcal S}_V$. 

\, \item The partition $\mathcal{P}_0$ consists of the
sets $\Ga_e:=\phi_e(J_{t(e)})$, $e\in E$. 

\, \item The measure $m$ is $m_h$, the $h$-conformal measure $m_h$, as a matter of fact, restricted $J_V$. 

\, \item The map $T_0:\Delta_0\to\Delta_0$ is, in our setting, just the map $F_V$. 

\, \item The function $R$, the return
time, is, naturally, defined as  
$$
R|_{\Ga_e}:=N(e),
$$
\end{enumerate}

\sp Fix $e\in E$ arbitrary and then two arbitrary points $x, y\in \Ga_e
=\phi_e(J_V)$. This means that $x=\phi_e(x')$ and
$y=\phi_\e(y')$ with some $x', y'\in J_V$. Since
$t^h-1=h(t-1)+O(|t-1|^2)$, and because of  Theorem~\ref{tKoebe-1}, 
there exist respective constants $C_1, C_2>0$ such that
we have,
\beq\label{120121207}
\aligned
\lt|\frac{Jac_{m_h}F_V(y)}{Jac_{m_h}F_V(x)}-1\rt|
&=\lt|\frac{|F_V'(y)|^h}{|F_V'(x)|^h}-1\rt|
=\lt|\frac{|\phi_e'(x')|^h}{|(\phi_e'(y')|^h}-1\rt|\\
&\le C_1\lt|\frac{|\phi_e'(x')|}{|\phi_e'(y')|}-1\rt|\\
&\le C_1C_2|y'-x'|.
\endaligned
\eeq
Now write $x'=\pi_V(\a)$ and $y'=\pi_V(\g)$
with appropriate $\a,\g\in E_A^\infty$. Put $\om:=\a\wedge\g$ and
$k:=|\om|$. Write also $x'':=\pi_V(\sg^k(\a))$ and
$y'':=\pi_V(\sg^k(\g))$. Then we get that 
\beq\label{120121208}
\aligned
|y'-x'|
&=   |\phi_\om(y'')-\phi_\om(x'')|
 \le |\phi_\om'(x'')||y''-x''| \\
&\le \diam_e(V_{t(\om)})\|\phi_\om'\| 
 \le \diam_e(V_{t(\om)})\b_V^{|\om|} \\
&=   \diam_e(V_{t(\om)})\b_V^{s(x',y')}\\
&=   \diam_e(V_{t(\om)})\b_V^{s(F(x),F(y))}\\
&\le \|V\|\b_V^{s(F(x),F(y))},
\endaligned 
\eeq
where, we recall, $\b_V\in(0,1)$ is the contracting factor of the system $\mathcal S_V$ and $\|V\|:=\sup\{\diam_e(V_v):v\in B(\hat{f})\}$. Since $s(F_V(x),F_V(y))=s(x,y)-1$ (as we know that $s(x,y)\ge 1$), along with \eqref{120121207}, the formula \eqref{120121208} gives that
$$
\lt|\frac{Jac_{m_h}F_V(y)}{Jac_{m_h}F_V(x)}-1\rt|
\le\|V\|\b_V^{s\(F_V(x),F_V(y)\)}
= \b_V^{-1}\|V\|\b_V^{s(x,y)}.
$$
So, \eqref{youngcond0} is established in our context. The fact that
the partition $\mathcal{P}_0$ is generating follows either from the
contracting property of graph directed Markov systems, or, more directly, from
Theorem~\ref{t1j237}. The Big Images Property holds because the alphabet of the GDMS $\cS_V$ is finite. The last assumption in Theorem~\ref{lsyoung0} is that the map $T:\Delta\to\Delta$, denoted now by $T_V$, is topologically mixing. We will prove it now, introducing on the way some concepts which will be also needed further in the main proof.

\blem\label{l120180515}
The dynamical system $T_V:\Delta\to\Delta$ is topologically mixing in the sense of Section~\ref{Young-Abstract}. 
\elem

\bpf
For every $\tau\in E_A^*$, set
$$
N(\tau):=\sum_{j=1}^{|\tau|}N(\tau_j).
$$
By our construction of the tower $\cY_f$, we have for every $e\in E$ that
$$
T_V^{R(\Ga_e)}\(\Ga_e\times\{0\}\)
=T_V(\Ga_e\times\{R(\Ga_e)\}\)
=F_V(\Ga_e)\times\{0\}
=J_{t(e)}\times\{0\}.
$$
So, by an immediate induction, we get for every $\tau\in E_A^*$ that
\beq\label{620180531}
T_V^{N(\tau)}\(\phi_\tau(J_{t(\tau)})\times\{0\}\)
=J_{t(\tau)}\times\{0\}.
\eeq
Now fix two arbitrary elements $a, b\in E$. Then there exists $s\ge 0$ such that
$$
f^u\(\phi_a(V_{t(a)})\)\spt \phi_b(V_{t(b)})
$$
for all $u\ge s$. Hence, for all such $u$ there exists a holomomorphic branch of $f^u$ mapping $\phi_b(V_{t(b)})$ into $\phi_a(V_{t(a)})$. By our construction of the conformal graph directed Markov system $\cS_V$, this holomorphic branch is an admissible composition of elements of $\cS_V$. This means that it is equal to 
$$
\phi_{\tau(u)}:V_{t(\tau(u))}\lra V_{i(\tau(u))}
$$ 
for some $\tau(u)\in E_A^*$ with $t(\tau(u))=i(b)$ and $i(\tau(u))=i(a)$. Then, applying \eqref{620180531}, we get for every $0\le k \le R(\Ga_a)-1$, every $0\le l \le R(\Ga_b)-1$, and every $n\ge s+R(\Ga_b)-1$, that 
$$
\begin{aligned}
T_V^n\(\Ga_a\times\{k\}\)
&=T_V^{n+k}\(\Ga_a\times\{0\}\)
=T_V^{l+n+k-l}\(\Ga_a\times\{0\}\) \\
&\spt T_V^l\(T_V^{n+k-l}\(\phi_{\tau(n+k-l)}(J_{t(\tau(n+k-l))})\times\{0\}\)\) \\
&=T_V^l\(J_{t(\tau(n+k-l))}\times\{0\}\)
\spt T_V^l\(\phi_b(J_b)\times\{0\}\) \\
&=\phi_b(J_b)\times\{l\} \\
&=\Ga_b\times\{l\}.
\end{aligned}
$$
The proof of Lemma~\ref{l120180515} is complete.
\epf

\sp\fr Having all the above and invoking also Lemma~\ref{l1j267}, which yields
$m_h(R^{-1}([n,\infty))\le Ce^{-\a n}$, we conclude from Observation~\ref{120200124} and Theorems~\ref{lsyoung0}, \ref{lsyoung}, and \ref{LILII}, that

\begin{itemize}
\item \beq\label{tildemphifinite}
\tilde m_h(\De)<+\infty,  
\eeq
where $\tilde m_h$ is derived out of $m_h$ restricted to $J_V$ according to formula \eqref{520180515},

\, \item The map $T_V:\Delta\lra\Delta$ admits a  probability $T$--invariant measure $\nu_h$ which is absolutely continuous with respect to $\tilde m_h$, 

\, \item Given $\g\in (0,1)$, for the dynamical system $T_V:(\Delta,\nu_h)\lra(\Delta,\nu_h)$ the following hold: 

\begin{itemize}

\, \item The Exponential Decay of Correlations in the form of (\ref{cov}) holds, 

\, \item The Central Limit Theorem and is true for all functions $g\in C_\g(\Delta)$ that are not cohomologous to a constant in $L^2(\nu_h)$. 

\, \item The Law of Iterated Logarithm is true for all functions $g\in C_\g(\Delta)$ that are not cohomologous to a constant in $L^2(\nu_h)$. 
\end{itemize}
\end{itemize}

\sp Now consider $H:\Delta\lra\mathbb{C}$, the natural
projection from the abstract tower $\Delta$ to the complex plane $\mathbb{C}$ given by the formula
$$
H(z,n)=f^n(z).
$$
Then
\beq\label{eq120110623B}
H\circ T =f\circ H,
\eeq
$$
\tilde m_h|_{\Delta_0}\circ H^{-1}=m_h|_{J_V},
$$ 
and  
$$
\tilde m_h|_{\Ga_e\times\{n\}} \circ H^{-1} 
=m_h|_{\Ga_e\times\{0\}}\circ f^{-n}
=m_h|_{\Ga_e}\circ f^{-n} 
$$ 
for all $e\in E$ and all $0\le n\le N(e)-1$.
So, $\tilde m_h|_{\Ga_e\times\{n\}}\circ H^{-1}$ is absolutely continuous
with respect to $m_h$, with the Radon--Nikodym derivative equal to
$$
J_{e,n}(z)
:=Jac_{m_h}\(f^n|_{\Ga_e}\)^{-n}(z)
=\Big|\(f^n\)'\(\(f^n|_{\Ga_e}\)^{-n}(z)\)\Big|^{-h}
$$ 
for all $z\in f^n(\Ga_e)$ and zero elsewhere in $J(f)$. Therefore, using \eqref{tildemphifinite}, we get that
$$
\int_{J(f)}\sum_{e\in E}\sum_{n=0}^{N(e)-1}J_{e,n}dm_h
= \sum_{e\in E}\sum_{n=0}^{N(e)-1}\!\!\int\!\! J_{e,n}dm_h 
=\tilde m_h\circ H^{-1}(J(f)) 
=\tilde m_h(\De) 
<+\infty.
$$ 
Thus, the function 
$$
\sum_{e\in E}\sum_{n=0}^{N(e)-1}J_{e,n}
$$ 
is integrable with respect to the measure $m_h$. This implies immediately that
the measure $\tilde m_h\circ H^{-1}$ is absolutely continuous with
respect to the measure $m_h$ with the Radon--Nikodym derivative equal to 
$$
\sum_{e\in E}\sum_{n=0}^{N(e)-1}J_{e,n}.  
$$
Hence, the measure $\nu_h\circ H^{-1}$ is also absolutely continuous
with respect to $m_h$. Since $\nu$ is $F_V$--invariant and $H\circ
T_V =f\circ H$, the measure $\nu_h\circ H^{-1}$ is $f$--invariant. 
But the measure $\mu_h$ is $f$--invariant, ergodic, and
equivalent to the conformal measure $m_h$. Hence, $\nu_h\circ H^{-1}$
is absolutely continuous with respect to the ergodic
measure $\mu_h$. In conclusion.  
\begin{lem}\label{l120110623}
If $f:\C\lra\oc$ is a normal subexpanding elliptic function of finite character, then
$$
\nu_h\circ H^{-1}=\mu_h.
$$
\end{lem} 

\sp We are now in position to prove the following.

\begin{thm}\label{t220110627}
If $f:\C\lra\oc$ is a normal subexpanding elliptic function of finite character and if
$\mu_h$ is the corresponding probability $f$-invariant measure
equivalent to the $h$-conformal measure $m_h$, then the dynamical
system $\(f|_{J(f)},\mu_h\)$ satisfies the following. If $g:J(f)\to \mathbb{R}$ is a bounded function, H\"older continuous with respect to the Euclidean metric on $J(f)$, then
\begin{enumerate}
\item{} For every bounded measurable  function $\psi:J(f)\to\mathbb{R}$, we have that
$$
\bigg|\int\psi\circ f^n \cdot g d\mu_h
-\int g d\mu_h\int\psi d\mu_h\bigg|
=O(\theta^n)
$$
for some $0<\theta<1$ depending on $\alpha$.

\sp\item{} The Central Limit Theorem holds for every H\"older continuous
bounded function $g:J(f)\to  \mathbb{R}$ that is not
  cohomologous to a constant in $L^2(\mu_h)$, i.e. for which there is no
  square integrable function $\eta$ for which $g={\rm
    const}+\eta\circ f-\eta$. More recisely, there exists $\sigma>0$ such that
$$
\frac{\sum_{j=0}^{n-1}g\circ f^j-n\int g d\mu_h}{\sqrt{n}} \lra \mathcal N(0,\sigma)
$$ 
in distribution, where, as usually, $\mathcal N(0,\sigma)$ denotes the Gauss (normal) distribution centered at $0$ with covariance $\sg$. 

\sp\item{} The Law of Iterated Logarithm holds for every H\"older continuous
bounded function $g:J(f)\to \mathbb{R}$ that is not cohomologous to a constant
in $L^2(\mu_h)$. This means that there exists a real positive
constant $A_g$ such that such that $\mu_h$ almost everywhere 
$$
\limsup_{n\to\infty}\frac{S_{n}g-n\int gd\mu_h}{\sqrt{n\log\log n}}=A_g. 
$$
\end{enumerate}
\end{thm}

{\sl Proof.}
Let $g:J(f)\lra\mathbb{R}$ and $\psi:J(f)\lra\mathbb{R}$ be as in the
hypotheses of our theorem. Define the functions 
$$
\tilde g:=g\circ H:\Delta\lra\mathbb{R} \  \  \text{ and } \ \
\tilde\psi:=\psi\circ H:\Delta\lra\mathbb{R}.
$$
We shall prove the following.

\sp\fr {\bf Claim~1:} The function $\tilde g:\Delta\lra\mathbb{R}$
belongs to the space $C_\beta$ for an appropriate exponent $\beta\in (0,1)$. 

\sp\fr Indeed, consider two arbitrary points $(x,k),(y,l)\in\De$. We
treat separately two cases depending on whether $s((x,k),(y,l))=0$ or
not. If $s((x,k),(y,l))=0$, then we get
\beq\label{220121208}
\aligned
|\^g(y,l)-\^g(x,k)|
&  =|g(H(y,l))-g(H(x,k))|
   =|g(f^l(y))-g(f^k(x))|\\
&\le|g(f^l(y))|+|g(f^k(x))| \\
&\le2||g||_\infty \\
&=  2||g||_\infty\b^{s((x,k),(y,l))}
\endaligned
\eeq
regardless of what the value of $\b\in(0,1)$ is, which will be specified in the next case. Indeed, if $s((x,k),(y,l))>0$, then $k=l$, $k<R(x)=R(y)$, and
\beq\label{320121208}
|\^g(y,l)-\^g(x,k)|
=|g(f^k(y))-g(f^k(x))|
\le H_g|f^k(y)-f^k(x)|^\g,
\eeq
where $H_g\ge 0$ and $\g>0$ are respectively the H\"older constant and
H\"older exponent of $g$. Moreover, $x, y\in\phi_\tau(V_{t(\tau)}$, with some $\tau\in E_A^{s((x,k),(y,l))}$, and then $f^k(x),f^k(y)\in f^k\(\phi_\tau(V_{t(\tau)})\)$ and the set $f^k\(\phi_\tau(V_{t(\tau)})\)$ is contained in a connected component of 
$$
f^{-(R(x)-k)}(F^{-(s((x,k),(y,l)))}\(V_{t(\tau)})\)
$$
whose diameter is, by Theorem~\ref{t1j237}, generously bounded above by 
$$
\exp\(-M_f\(s((x,k),(y,l))+R(x)-k\)\)
\le \exp\(-M_f(s((x,k),(y,l))\).
$$ 
In conjunction with \eqref{320121208} (and \eqref{220121208}), this
finishes the proof of Claim~1 by taking $\b=\exp(-\g M_f)$.

\sp\fr{\bf Claim~2:} The function $\tilde g$ is not
cohomologous to a constant in $L^2(\nu)$.

\sp\fr Indeed, assume without loss of generality that $\mu_h(g)=0$. By
virtue of Lemma~\ref{l320121210} the fact that  $g:J(f)\to \mathbb{R}$ is not 
a coboundary in $L^2(\mu_h)$ equivalently means that the sequence
$\(S_n(g)\)_{n=1}^\infty$ is not uniformly bounded in
$L^2(\mu_h)$. But because of Lemma~\ref{l120110623},  $\|S_n(\tilde
g)\|_{L^2(\nu)} = \|S_n(g)\|_{L^2(\mu_h)}$. So, the sequence $\(S_n(\tilde
g)\)_{n=0}^\infty$ is not uniformly bounded in $L^2(\nu)$. Thus, by
Lemma~\ref{l320121210} again, it is not a coboundary in $L^2(\nu)$. 

\sp\fr Having these two claims, all items, (1), (2), and (3), now
follow immediately from Theorem~\ref{lsyoung} with the use of
Lemma~\ref{l120110623} and formula \eqref{eq120110623B}. The proof is
finished. 
\endpf

\sp\subsection{Metric Entropy}

In this miniature subsection, by taking fruits of what has been done so
far, we will easily prove that the metric entropy of the dynamical system
$\(f|_{J(f)},\mu_h\)$ is finite.

\begin{thm}\label{t120121211}
If $f:\C\lra\oc$ is a normal subexpanding elliptic function of finite character and if $\mu_h$ is the corresponding Borel probability $f$--invariant measure
equivalent to the $h$--conformal measure $m_h$, then $\h_{\mu_h}(f)<+\infty$.
\ethm 

{\sl Proof.} It is a direct consequence of Theorem~\ref{t120180509} and
Theorem~\ref{GDMS_Finite_Entropy} to have finite entropy of the induced
system $F_V:J_V\lra J_V$ with respect to the probability $F_V$--invariant
measure $\mu_V=(\mu_h(J_V))^{-1}\mu_h|_{J_V}$. Then, as $F_V:J_V\lra J_V$ is the first return map of $f$ from $J_V$ to $J_V$, Abramov's Formula
(Theorem~\ref{abramovabstract}) yields $\h_{\mu_h}(f)<+\infty$. 
\endpf

\sp For the case when the Julia set is the whole complex plane $\C$,
we get from this the following.

\sp\begin{cor}\label{c220121211}
If $f:\C\lra\oc$ is a normal subexpanding elliptic function of finite character with $J(f)=\C$, then $\h_{\mu_2}(f)<+\infty$,
where $\mu_2$ is the (unique) Borel probability $f$--invariant
measure on $\C$ equivalent to planar Lebesgue measure on $\C$ .
\ecor

\sp\section[Parabolic Maps: Nice Sets, GDMSs, Conformal and Invariant Measures] {Parabolic Elliptic Maps: Nice Sets, Graph Directed Markov Systems, Conformal and Invariant Measures, \\ Metric Entropy}\label{PEM-I}

In this section we deal with parabolic elliptic functions\index{(N)}{parabolic elliptic function} $f:\mathbb {C} \lra \oc$, i.e. with all such elliptic functions $f:\mathbb {C} \lra \oc$ for which
$$
\Crit(f)\cap J(f)=\es \ \
{\rm and } 
\  \  \Om(f)\ne\es. 
$$ 
As an immediate of the first half of this definition we get the following observation, already included in Chapter~\ref{DDCoEF} as Observation~\ref{o120200125}. 

\bobs\label{o120200123}
Every parabolic elliptic function is regular normal compactly non--recurrent of finite character.
\eobs 

In consequence. all what we have proved so far for such functions applies to parabolic elliptic functions. 
 
Our aim in this section is to construct appropriate Pre--Nice and Nice Sets dealt with at length in Section~\ref{NiceSetsGeneral}, to construct, as it also was done Section~\ref{NiceSetsGeneral}, and to study the corresponding conformal maximal graph directed Markov systems in the sense of Chapter~\ref{Markov-systems}, showing particularly their strong regularity. Having all these tools we also prove that the Kolmogorov--Sinai entropy of the measure $\mu_h$ is always finite even if the measure $\mu_h$ itself is infinite. Lastly, at the end of the section, we shortly deal with Radon--Nikodym derivatives of invariant measures $\mu_h$ with respect to the conformal ones $m_h$. 

\sp 

We start with defining appropriate nice families for a parabolic map $f:\C\lra\oc$, actually for some of its iterates. Indeed, let $l=l_f\geq 1$ be such integer that each element of $\Omega(f)$ is a simple parabolic point of $f^l$. Our first goal is  to construct a nice set for the holomorphic map
$$ 
\hat{f}^l: \mathbb{T}_f^{(l)}=\mathbb{T}_f \sms \bigcup_{j=0}^{l-1}\hat{f}^{-j}(B_\infty(\hat{f}))\lra \mathbb T_f,$$
where, we recall,
$$
B_\infty (\hat{f})= \Pi_f( f^{-1}(\infty)). 
$$ 
We  want to  construct such a nice set with properties suitable to estimate from above the first entrance time to it sufficiently well so that the machinery of Young's Towers set up in Section~\ref{Young-Abstract}, Lai-Sang Young Towers, and applied in the preceding Section~\ref{Young's Towers for Subexpanding}, Young's Towers for Subexpanding Elliptic Functions of Finite Character, will be now applicable for those parabolic elliptic maps for which the measure $\mu_h$ is finite. For our construction of nice sets, we apply Theorem~\ref{t2j291B}. Sticking to its notation we set
\beq\label{120190417}
F_0
:=\Omega(\hat{f})
=\Pi_f(\Omega(f))
= \Omega( \hat{f}^l)
\quad \text{and} \quad  
F_1:= B_\infty(\hat{f}).
\eeq
Of course the finite set $\Om(\hat{f}^l)$ consists of all rationally indifferent  fixed points of $\hat{f}^l$.  For  all $b \in F_1$, we set 
$$
U_0(b):=B(b,r)
$$ 
with $r>0$ so small that
\beq\label{1fpn14}
B(b, 6r)\cap \({\rm PC}(f)\cup \Om(\hat{f})\)=\es,
\eeq
and as required further in the course of our construction. For every $\om \in F_0$, a fixed $\a\in(\pi/2,\pi)$, a fixed $\ka\in(0,1)$, and every $j\in \{1, \ldots, p(\om)\}$, set 
\beq\label{1_2019_04_16}
U_0(\om, j):=\Pi_f\(S^j_r(\om,\a)\)
\eeq
where the petals $S^j_r(\om,\a)$ are defined by formula \eqref{1_2017_29B}.
With
$$ 
X:=\mathbb T_f^{(l)}
:=\mathbb T_f\sms \bigcup_{j=0}^{l-1} \hat{f}^{-j}(B_\infty(\hat{f})), \   \ \  Y=\mathbb T_f, 
\  \   {\rm and} \  f \ {\rm being} \ \hat{f}^l,
$$
all the hypotheses  of Theorem~\ref{t2j291B} up to (2d) (included) are evidently satisfied for all $b \in F$ and the corresponding radii $r_b >0$, small enough. 

The hypothesis \rm{(2e)} follows directly from the definition \eqref{1_2019_04_16}. 

The hypotheses (2f) and (2g) follow immediately  from Proposition~\ref{p1j108}. 

The hypotheses (2h), (2i), (6), and (7), are also  immediate from  our definitions. 

\,

In order to verify hypothesis (2j) fix any $s>0$ and then any $x\ge x(\a,\ka)$ so large that 
$$
\bu_{j=1}^{p(\om)}\(\Pi_f(S^j_r(\om,\a))\cup \Pi_f(S^j_a(\om,\a))\)\sbt B(\om,s).
$$
Since $\a\in(\pi/2,\pi)$, we have that
$$
\Int\lt(\{\om\}\cup\bu_{j=1}^{p(\om)}\(\Pi_f(S^j_r(\om,\a))\cup \Pi_f(S^j_a(\om,\a))\}\rt)
\ne\es.
$$
Therefore, there exists $s_\om^->0$ such that 
$$
B(\om,s_\om^-)\sbt \{\om\}\cup\bu_{j=1}^{p(\om)}\(\Pi_f(S^j_r(\om,\a))\cup \Pi_f(S^j_a(\om,\a))\).
$$
So, if 
$$
z\in B(\om,s_\om^-)\sms \bu_{j=1}^{p(\om)}U_0(\om, j),
$$
then 
$$
z\in \bu_{j=1}^{p(\om)}\Pi_f(S^j_a(\om,\a)).
$$
Therefore, 
$$
f^n(z)\in \bu_{j=1}^{p(\om)}\Pi_f(S^j_a(\om,\a))\sbt B(\om,s)
$$
for all integers $n\ge 0$. Thus, the hypothesis (2j) is satisfied too. 

The hypothesis (3) clearly holds for all $r_\xi >0$, $ \xi \in F_0$, small  enough, and (4) holds for all $r_\xi>0$, $ \xi \in F_1$  small  enough, because  of (\ref{1fpn14}). 

The hypothesis (5a) holds for all $r_\xi>0$, $ \xi \in F$, small enough since
  \beq\label{1fpn15}
  \ov{{\rm PC} (\hat{f})} \cap J(\hat{f})= \Om (\hat{f}).
\eeq

We are left to show  that (5b) can be  satisfied  with an appropriate  choice of $r_\xi, \xi\in F$, small enough. By (\ref{220180824R}), being a hypothesis of (5b), we have that
$$
\inf\Big\{\big|\hat{f}^{l(n-1)}\circ\hat{f}^{-ln}_{a,b}(b)-b\big|:  b \in \Om(\hat{f}), a \in \Om(\hat{f})\sms \{b\}, n\ge 1\Big\} >0.
$$
We look for radii $r_\xi$, $\xi \in F$, being equal, i.e.  
$$ 
r:=r_a=r_b
$$ 
for all $a, b\in F$. Since $\hat{f}^{-1} (\Om(\hat{f})) \cap \Crit(\hat{f})= \es $, since the points of $\hat{f}^{-1}\(\Om(\hat{f})\)$ cluster towards $B_\infty(\hat{f})$, and since $\ov{\Pi_f^{-1}\(\hat{f}^{-1}\(\Om(\hat{f})\)\sms\Om(\hat{f})\)}$ is a compact subset of $J(f)$ disjoint from $\ov{{\rm PC}(f)}$, it follows from Lemma~\ref{l120180306} (3) that
\beq\label{1fpn15.1}
M:= \min\Big\{1, \inf\big\{|\hat{(f^l)}'(z)|:  z \in \hat{f}^{-l}( \Om(\hat{f})\big\}\Big\}>0.
\eeq
Let $ Q \sbt \mathbb T_f^{(l)} $ be a periodic orbit of $\hat{f}^l: \mathbb T_f^{(l)} \lra \mathbb T_f$ disjoint  from $\Om( \hat{f})$; it is  of course disjoint from $B_\infty( \hat{f})$. So, our analytic map $\hat{f}^l: \mathbb T_f^{(l)}\lra \mathbb T_f$ enjoys the  Standard Property. Since the set $B_\infty(\hat{f})$ is finite  and disjoint from $Q \cup\ov{\PC(f)}$ (see ( \ref{1fpn15})), it follows from Lemma~\ref{l120180306} (4) that there exists an integer $N_1\geq 1$ such that 
\beq\label{1fpn17}
| (f^n)'(\xi) |\geq M^{-1}K^2\frac{8 \kappa}{\kappa -1}
\eeq
 for every $n \geq N_1$ and every $\xi \in \hat{f}^{-n}(B_\infty(\hat{f}))$, where $K \geq 1$ is, as usually, the Koebe's constant corresponding to the scale $1/2$.

 Since $\Om(\hat{f})$ is a finite set consisting of fixed points of $\hat{f}$ only and it is disjoint from  the finite set $B_\infty(\hat{f})$, there exist $R_1  >0$ such that
\beq\label{2fpn17}
    B\(B_\infty(\hat{f}), 4R_1\) \cap \lt(\ov{\PC(f)}\cup \bigcup_{k=0}^{N_1 l} \hat{f}^k\( B_\infty (\Om(\hat{f}), 4R_1)\)\rt)=\es.
\eeq
Then, looking  up at (\ref{1fpn17}) and using Koebe's Distortion Theorem, we conclude that
\beq\label{2fpn17.1}
| (f^n)'(z) |\geq M^{-1}K\frac{8 \kappa}{\kappa -1}
\eeq
for every $n\geq N_1$ and every $ z\in f^{-n}\(B(B_\infty(\hat{f}), 2R_1)\)$.
Since the set
 \beq\label{3fpn17}
 Z:= \big(\hat{f}^{-l}(\Om(\hat{f}))\sms \Om(\hat{f})\big)\sms B(B_\infty(\hat{f}), R_1)
 \eeq
is finite and disjoint from $ Q \cup \ov{\rm PC} (f)$, it  again follows from Lemma~\ref{l120180306} (4) that  there exists an integer $N_2\geq N_1$ such that
\beq\label{4fpn17}
| (f^n)'(\xi) |\geq M^{-1}K\frac{8 \kappa}{\kappa -1}
\eeq
for every $n \geq N_2$ and every  $\xi \in f^{-n}(Z)$. As before, since $\Om(\hat{f})$ is a finite set consiting of fixed points of $\hat{f}$ only and it is  disjoint from $Z$, there exists $R_2\in(0, R_1]$ such that
\beq\label{5fpn17}
B(Z, 4R_2)\cap \bigcup_{k=0}^{N_2l}\hat{f}^k\(B(\Om(\hat{f}), 4R_2)\)=\es.
\eeq
Take any
\beq\label{6fpn17}
r \in (0, R_2].
\eeq
Dealing still with hypothesis (5b), consider several cases. 

\sp\fr {\sc  Case 1}: $ a\in F_0, b\in F_1$. Then $n \geq N_1$ by (\ref{2fpn17}). Hence, \eqref{32018082R} holds because of (\ref{6fpn17}), (\ref{1fpn15.1}) ($M\le 1$), (\ref{2fpn17.1}),
and Koebe's Distortion Theorem. 

\sp\fr {\sc  Case 2}: $a, b\in F_0.$ Then
$$
\hat{f}^{l(n-1)}(\hat{f}^{-ln}_{a,b}(b)) \in \hat{f}^{-l}(\Om(\hat{f}))\sms \Om (\hat{f}).
$$
Suppose first that 
$$
\hat{f}^{l(n-1)}( \hat{f}^{-ln}_{a,b}(b))\in Z.
$$
Then ${f}^{-ln}_{a,b}(b)\in \hat{f}^{-l(n-1)}(Z).$ So, it follows  from (\ref {5fpn17})  that $ n-1\geq N_2$. So, using (\ref{4fpn17}) and (\ref{1fpn15.1}), we get that
$$
|(\hat{f}^{ln})'( \hat{f}^{-ln}_{a,b}(b))| 
=|(\hat{f}^{l(n-1)})'(\hat{f}^{-ln}_{a,b}(b))|\cdot|(\hat{f}^l)'(\hat{f}^{l(n-1)}(\hat{f}^{-ln}_{a,b}(b)))|
\geq M^{-1}K\frac{8 \kappa}{\kappa -1}M
=K\frac{8 \kappa}{\kappa -1}.
$$
Invoking now Koebe's Distortion Theorem (and equality $r=r_a=r_b$),  formula  \eqref{32018082R} thus follows in this case.

\sp So, suppose that 
$$
\hat{f}^{l(n-1)}( \hat{f}^{-ln}_{a,b}(b))\in  B(B_\infty(\hat{f}), R_1).
$$
It then follows from (\ref{2fpn17}) that $n-1\geq N_1$. So using (\ref{2fpn17.1}) and (\ref{1fpn15.1}) we get that
$$ 
|(\hat{f}^{ln})'( \hat{f}^{-ln}_{a,b}(b))| 
=|(\hat{f}^{l(n-1)})'(\hat{f}^{-ln}_{a,b}(b))|\cdot|(\hat{f}^l)'(\hat{f}^{l(n-1)}(\hat{f}^{-ln}_{a,b}(b)))|
\geq M^{-1}K\frac{8 \kappa}{\kappa -1}M
=K\frac{8 \kappa}{\kappa -1}.
$$
Invoking, as  before,  Koebe's Distortoon Theorem and equality $r=r_a=r_b$,  formula \eqref{32018082R} thus follows in this case too.

\sp\fr {\sc Case 3}: $a\in F_1$, $ b\in F$. Assume first that 
\beq\label{1fpn18}
n\ge N_2+1.
\eeq
Consider three subcases:

\sp\fr {\sc Case 3A}: $b\in F_1$. Then the formula \eqref{32018082R} holds because of (\ref{2fpn17.1}).

\sp\fr {\sc Case 3B}: $ b\in F_0$. Then 
$$
\hat{f}^{l(n-1)}(\hat{f}^{-ln}_{a,b}(b))
\in \hat{f}^{-l}(\Om(\hat{f}))\sms \Om(\hat{f}),
$$
and, assume that
$$
\hat{f}^{l(n-1)}( \hat{f}^{-ln}_{a,b}(b))\in Z.
$$
Then $\hat{f}^{-ln}_{a,b}(b)\in \hat{f}^{-l(n-1)}(Z).$
Since by (\ref{1fpn18}), $n-1\geq N_2$, using (\ref{4fpn17}) and (\ref{1fpn15.1}), we get that
$$ 
|(\hat{f}^{ln})'( \hat{f}^{-ln}_{a,b}(b))| 
=|(\hat{f}^{l(n-1)})'(\hat{f}^{-ln}_{a,b}(b))|\cdot|(\hat{f}^l)'(\hat{f}^{l(n-1)}(\hat{f}^{-ln}_{a,b}(b)))|
\geq M^{-1}K\frac{8 \kappa}{\kappa -1}M
=K\frac{8 \kappa}{\kappa -1}.
$$
Invoking  now  Koebe's Distortion Theorem, formula \eqref{32018082R} thus follows in this case.

\sp\fr{\sc Case 3C}: $b\in F_0$  and $\hat{f}^{l(n-1)}(\hat{ f}^{-ln}_{a,b}(b))\in  B(B_\infty(\hat{f}),R_1).$ Then, using  (\ref{2fpn17.1}) and (\ref{1fpn15.1}), we get that
$$ 
|(\hat{f}^{ln})'(\hat{f}^{-ln}_{a,b}(b))| 
=|(\hat{f}^{l(n-1)})'(\hat{f}^{-ln}_{a,b}(b))|\cdot|(\hat{f}^l)'(\hat{f}^{l(n-1)}(\hat{f}^{-ln}_{a,b}(b)))|
\geq M^{-1}K\frac{8 \kappa}{\kappa -1}M
=K\frac{8 \kappa}{\kappa -1}.
$$
Invoking, as before, Koebe's Distortion Theorem, formula \eqref{32018082R} thus follows in this case too.

\sp So, finally assume that
$$ 
n \leq N_2.
$$
Since $\ov{\Pi_f^{-1}(B_\infty(\tilde f))}$ is a compact subset of $J(f)$ disjoint from $\ov{{\rm PC}(f)}$, the same argument as the one giving (\ref{1fpn15.1}), gives us that
$$
D:= \inf\Big\{\big|(\hat{f}^{lk})'(z)\big|:0\le k\le N_2,\, z \in \hat{f}^{-lk}(F)\Big\}>0.
$$
Also, by taking $r>0$ small enough so that $\hat{f}^{-ln}_{a,b}(b)$ (with all $a$, $b$, and and $n$ as in the considered case) is sufficiently close to  $B(\hat{f})$, we will have
$$
\big|(\hat{f}^l)'(\hat{f}^{-ln}_{a,b}(b))\big| 
\geq D^{-1}K \frac{8 \kappa}{\kappa -1}.
$$
Hence, 
$$ 
\big|(\hat{f}^{ln})'( \hat{f}^{-ln}_{a,b}(b))\big| 
=\big|(\hat{f}^l)'(\hat{f}^{-ln}_{a,b}(b))\big|\cdot\big|\(\hat{f}^{l(n-1)}\)'(\hat{f}^l(\hat{f}^{-ln}_{a,b}(b)))\big|
\geq D^{-1}K\frac{8 \kappa}{\kappa -1}D
=K\frac{8 \kappa}{\kappa -1}.$$
 Thus, a direct application of  Koebe's Distortoon Theorem,  completes the proof of the formula \eqref{32018082R}. The  hypothesis (5b) of Theorem~\ref{t2j291B} is therefore verified. Its hypotheses (8) and (9) also obviously hold for all $r\in(0,R_2$ small enough. So, Theorem~\ref{t2j291B} and Theorem~\ref{t2j291} as well, applies and we have the following.

\bthm\label{t1fpn20}
If $f:\C\lra\oc$ is a parabolic elliptic map, then the construction started with \eqref{120190417} leads, by applying Theorem~\ref{t2j291}, to 
a pre--nice $U$ for the map $\hat{f}^l: \mathbb T_f^{(l)} \lra \mathbb T_f$, which in turn leads, via Theorem~\ref{t2j291B}, to a  maximal conformal graph directed Markov system $\hat\cS_U$ (acting on $\mathbb T_f^{(l)}$) in the sense of Definition~\ref{d120190419} and Definition~\ref{d2_2017_11_18}. Denote the edges of this system by $\hat E$ and the corresponding contractions $\hat\phi_e:X_{t(e)}\lra X_{i(e)}$.
\ethm

We need a few more properties concerning nice sets produced in the above theorem.

\blem\label{l2fpn20} Let $f:\C\lra\oc$ be a parabolic elliptic map. If $U$ is the pre--nice set coming from Theorem~\ref{t1fpn20}, then for every $\om \in \Om(\hat{f})$ and every $j\in \{1, 2, \ldots, p_\om\}$, we have that
$$ 
J(\hat{f})\cap \left( \{\om\}\cup\bigcup_{j=1}^{p_\om}\bigcup_{n=0}^\infty  \hat{f}^{-ln}(X(\om,j))\right)
$$
is a neighborhood of $\om$ in $J(\hat{f})$, where the sets $X(\om,j)$ are the  members of the maximal conformal graph directed Markov system $\hat\cS_U$ from
Theorem~\ref{t1fpn20}; see Theorem~\ref{t2j291B}, Theorem~\ref{t2j291}, and  the formula \eqref{1gns30} for their definition (with the parameter $u$ being skipped).
\elem

\bpf Since, by Lemma~\ref{l1ch11.4},
 $J(\hat{f})\cap \left( \{\om\}\cup\bigcup_{j=1}^{p_\om} U_0(\om, j)\right)$
 is an open  neighborhood of $\om \in J(\hat{f})$, so it the set
$$
H_\om:= J(\hat{f})\cap \left( \{\om\}\cup\bigcup_{j=1}^{p_\om} U(\om, j)\right).
$$
Hence, there exists $u>0$ so small that
\beq\label{1fpn21}
J(\hat{f})\cap \left( \bigcup_{j=1}^{p_\om}V_j(\om,u)\cup  \bigcup_{j=1}^{p_\om} f^l( V_j(\om,u))\right) 
\sbt H_\om.
\eeq
By virtue of Proposition~\ref{p1_2017_10_30}, for every $z \in H_\om\sms \{\om\}$, there exists $k\geq 1$ such that

\beq\label{2fpn21}
f^{kl}(z) \notin H_\om.
\eeq
Let $k$  be the minimal integer with this property. Then  $f^{(k-1)l}(z) \in U(\om, j) \sms f^{-l}_\om(U(\om, j))$ with
some $j \in \{1, \ldots, p_\om\}$.  It  then follows from  (\ref{1fpn21}) and  (\ref{2fpn21}) that
$f^{(k-1)l}(z) \notin \bigcup_{i=1}^{p_\om} f^{-l}(V_j(b,u))$. Hence
$$ 
f^{(k-1)l}(z) \in \(U(\om,j)\sms  f^{-l}_\om(U(\om, j))\)\sms V_j(\om,u)\sbt X(\om,j;u).
$$
So,
$$ 
H_\om\sbt J(\hat{f})\cap \left( \{\om\} \cup  \bigcup_{j=1}^{p_\om} \bigcup_{n=0}^\infty f^{-ln}_\om(X(\om,j;u))\right)
$$
and, with $X(\om, j):= X(\om,j;u)$, the proof is complete. 
\epf

\sp\fr As an immediate consequence of Lemma~\ref{l1j251A}, we  get the following.

\blem\label{l1fpn27}
Let $f: \mathbb C \lra \hat{\mathbb  C}$ be a parabolic elliptic function and let $\hat V$ be an open subset of $\mathbb T_f$ such that
\begin{enumerate}
\item $B_\infty(\hat{f})\sbt \hat V$  

\,

\item $\(\hat V\cup\Om(\hat{f})\)\cap J(\hat f)$ is an open neighborhood of 
$\Om(\hat{f})$ in $J(\hat f)$.

\,

\item the open set $\hat V$ has finitely many connected components,

\item the open set $\hat{V}^+ $ is horizontal.
\end{enumerate}

\,

Then for every $k\geq 1$ there exists $t \in \(\ov{\BD}(J(\hat f)\cap K_k(\hat V)), \HD (J(f))\)$ such that, with the  notation of Section~\ref{CKS}, we have that 
$$ 
\sum_{Q \in {\mD}_k(\hat V)} \diam^t(Q)< +\infty.
$$
\elem

\sp\fr As an immediate consequence of this lemma along with we Lemma~\ref {l2fpn20} and Observation~\ref{o1_2017_11_20}, we get the following.

\blem\label{l3fpn27} 
If $f:\mathbb C \lra \hat{\mathbb C}$ is a parabolic elliptic function and $U$ is a pre--nice set of Theorem~\ref{t1fpn20} for the map $\hat{f}^l: \mathbb T_f^{(l)} \lra \mathbb T_f$, then there exists $t\in \(\ov{\BD}(K_k(U)),\HD(J(f))\)$ such that, with the  notation of Section~\ref{CKS}, we have that  
$$
\sum_{Q \in {\mD}_l(U)} \diam^t(Q)< +\infty.
$$
\elem

\sp\fr We are now able to  prove the following. 

\blem\label{l2fpn27}
If $f:\mathbb C \lra \hat{\mathbb C}$ is a parabolic  elliptic function and $U$ is a  pre--nice set produced in Theorem~\ref{t1fpn20} for the map $ \hat{f}^l: \mathbb T_f^{(l)} \lra \mathbb T_f$, then the corresponding maximal graph directed Markov system $\hat{\mathcal S}_{U}$, given also by Theorem~\ref{t1fpn20}, for the map $ \hat{f}^l:\mathbb T_f^{(l)} \lra \mathbb T_f$ is strongly regular. 
\elem

\bpf Let $ q \geq 1$ be the maximal multiplicity of a pole of $f$. Fix the number $ t\in (\ov{\BD}(K_k(V)), \HD (J(f)))$ produced in Lemma~\ref{l3fpn27}. Since any number in $(t, \HD(J(f)))$ is  also good for this lemma, we may assume, because of Theorem~\ref{thm:julia}, that
\beq\label{1fpn27}
 t \in \left( \max\lt\{ \ov{\BD}(J(\hat f)\cap K_l(U)),  \frac{2q}{q+1} \rt\}, \HD(J(f))\right).
\eeq
For all $a, b \in F$,   let  
$$
E_{a,b}:=\big\{e \in E:  p_1(i(e))=a \  \  \text{and} \  \  p_1(t(e))=b\big\}.
$$
Since  for every  $e \in E(0,1)$, we have that 
$$
\hat{f}^l\( \hat{\phi}_e(U_{t(e)})\) \in {\mD}_l(U)
$$ 
and the map 
$$
E(0,1)\ni e \longmapsto  \hat{f}^l\( \hat{\phi}_e(U_{t(e)}))\) \in {\mD}_l(U)
$$
is one--to--one (assuming that the radius of the pre--nice set $U$ is sufficiently small so that $\hat{f}^l_{|B(b, 6r)}$ is one--to--one for all $ b \in F_0=\Om(\hat{f})$, we get from Lemma~\ref{l3fpn27} that
 \beq\label{1fpn28}
 \sum_{e \in E(0,1)} \diam^t\( \hat{f}^l( \hat{\phi}_e(U_{t(e)}))\) <+ \infty.
 \eeq
  Since  for parabolic maps, we  have
\beq\label{4fpn28}
M:= \inf\big\{|\hat{f}'(z)|: \,  z \in J(\hat{f})\big\}>0, 
\eeq
with the help  of Koebe's Distortion Theorem, we  conclude that for all $e \in E(0,1)$, we have that
\beq\label{5fpn28}
\diam^t(\hat{f}^l\( \hat{\phi}_e(U_{t(e)}))\) 
\geq \Gamma M^{lt}\|\hat{\varphi}_e'\|^t
\eeq
with some constant $\Gamma \in(0,+\infty)$ independent of $e\in E(0,1)$. Inserting this to (\ref{1fpn28}), we get  that
\beq\label{2fpn28}
\sum_{e\in E(0,1)}\|\hat{\varphi}_e'\|^t< +\infty.
\eeq
Passing to dealing with the set $E(1,1)$, note that, similarly as in  the previous case,
\beq\label{3fpn28}
{\mD}_{1,1}
:= \big\{\hat{f}^l( \hat{\phi}_e(U_{t(e)})):e \in E(1,1)\big\} \sbt {\mD}_l(U).
\eeq
For every $ D \in {\mD}_{1,1}$ let
$$
E(D):=\big\{e \in E(1,1):\hat{f}^l( \hat{\phi}_e(U_{t(e)}))=D\big\}.
$$
Since $t>\frac{2q}{q+1}$, we  have that 
$$ 
\sum_{\xi\in \hat{f}^{-l}(z)}|(\hat{f}^l)'(\xi)|^{-t}
\comp \left(\sum_{\xi\in \Lambda_f}|x|^{-\frac{q+1}{q}t}\right)^l
<+\infty,
$$
for every $z \in J(f)$ with the comparability constant independent of $ z \in J(f)$. Hence, using also Koebe's Distortion Theorem, we conclude that
\beq\label{3fpn29}
\sum_{e\in  E(D)}\|\hat{\varphi}_e'\|^{t}\comp  \left(\sum_{\xi\in \Lambda_f}|x|^{-\frac{q+1}{q}t}\right)^l \diam^t(D)\comp \diam^t(D)
\eeq
holds for  all $D \in \mD_{1,1}$. Because of this, (\ref{3fpn28}), and Lemma~\ref{l3fpn27}, we  conclude that
\beq\label{1fpn29}
\sum_{e\in  E (1,1)}\|\hat{\varphi}_e'\|^{t}
= \sum_{D\in  {\mathcal D}_{1,1}}\sum_{e \in  E(D)}\|\hat{\varphi}_e'\|^{t} \comp \sum_{D\in {\mathcal D}_{1,1}}\diam^t(D)<+\infty.
\eeq

Let us now deal with the set $E(0,0)$. For every $e \in E(0,0)$, let $k(e)\geq 0$ be  the least integer such that
$$
\hat{f}^{lk(e)} (\hat{\varphi}_e(U_{t(e)})) \sbt U_{t(e)}
$$
for every integer $k(e)\le s\le N_U(e)$. For every $n \geq 0$ let
$$
E_n(0,0):= \{ e \in E(0,0): \, N_U(e)-k(e)=n\}.
$$ 
Then
\beq\label{1fpn29.1}
\hat{f}^{lk(e)}(\hat{\varphi}_e(U_{t(e)}))= \hat{f}^{-ln}_{p_1(t(e))}(U_{t(e)}).
\eeq
Let $\psi_e: B(p_1(t(e)), 6r) \lra \mathbb T_f$  be a unique holomorphic branch of $\hat{f}^{-lk(e)}$ such that
\beq\label{4fpn29}
\hat{\varphi}_e= \psi_e \circ \hat{f}^{-ln}_{p_1(t(e))}.
\eeq
Then
\beq\label{2fpn29}
\hat{f}^l\circ \psi_e(U_{t(e)})\in \mD_l(U).
\eeq
Invoking the assumption of the previous case that the radius of the pre--nice set $U$ is sufficiently small so that  $\hat{f}^l_{|B(b,6r)}$ is one--to--one for all $b \in F_0=\Om(f)$, we  conclude that  the map
$$ 
E_n(0,0)\ni e \lra  \hat{f}^l\circ \psi_e (U_{t(e)})\in \mD_l(U)
$$
is one--to--one. Therefore, using Koebe's Distortion Theorem, formula (\ref{1fpn29.1}), (\ref{4fpn29}), Proposition~\ref{p1j111DL}, and Lemma~\ref{l3fpn27}, we get that 
\beq\label{4fpn30}
\begin{aligned}
\Sg_n(0,0)&:= \sum_{e \in E_n(0,0)}\diam^t(\hat{f}^l( \hat{\phi}_e(X_{t(e)})))\\
& \comp \sum_{e \in E_n(0,0)}\diam^t(\hat{f}^l \circ \psi_e(U_{t(e)})) \frac{\diam^t(\hat{f}^{-ln}_{p_1(t(e))}(X_{t(e)}))}{\diam^t(U_{t(e)})}\\
& \comp \sum_{e \in E_n(0,0)}\diam^t(\hat{f}^l \circ \psi_e(U_{t(e)})) (n+1)^{- \frac{ p(p_1(t(e)))+1}{p(p_1(t(e)))}t}\\
&=  (n+1)^{- \frac{ p(p_1(t(e)))+1}{p(p_1(t(e)))}t}\sum_{e \in E_n(0,0)}\diam^t(\hat{f}^l \circ \psi_e(U_{t(e)}))\\
&\comp  (n+1)^{- \frac{ p(p_1(t(e)))+1}{p(p_1(t(e)))}t}.
\end{aligned}
\eeq
Therefore,
\beq\label{1fpn30}
\sum_{e \in E(0,0)}\diam^t(\hat{f}^l \circ \hat{\varphi}_e(X_{t(e)}))
=\sum_{n=0}^\infty\Sg_n(0,0) 
\comp \sum_{n=0}^\infty (n+1)^{- \frac{ p(p_1(t(e)))+1}{p(p_1(t(e)))}t}
< +\infty,
\eeq
where the last inequality sign was written since, by (\ref{1fpn27}), $t>1$. Since in our present case formula (\ref{5fpn28}) is also true, we  conclude from it  and (\ref{1fpn30}) that
\beq\label{2fpn30}
\sum_{e\in E(0,0)}\|\hat{\varphi}_e'\|^t< +\infty.
\eeq

We now deal with the last case: the set $E(1,0)$. This  is the case  where  the difficulties of both cases $E(1,1)$ and $E(0,0)$ come together. It goes in a sense by combining together these two cases. For every $e \in E(1,0)$, we define $0\le k(e)\le N_U(e)$ in exactly the same way as in the case of $E(0,0)$. Similarly as in this case let
$$ 
E_n(1,0):= \big\{ e \in E(1,0): \,  N_U(e)-k(e)=n\big\}.
$$
Again, as in  the  $E(0,0)$ case, formula (\ref{1fpn29.1}) holds, the  holomorphic map 
$$
\psi_e: B(p_1(t(e)), 6r) \lra \mathbb T_f
$$ 
is defined in the same way, formula (\ref{3fpn29}) and (\ref{2fpn29}) hold, and the map
$$
E_n(1,0)\ni e \longmapsto \hat{f}\circ \psi_e (U_{t(e)})\in \mD_l(U)
$$
is one--to--one. As in  the case of $E(1,1)$, we define $\mD_{1,0}(n)$ by the following formula corresponding to (\ref{3fpn28}). For every $n \geq 0$
$$
{\mD}_{1,0}(n)
:= \big\{\hat{f}^l( \hat{\psi}_e(U_{t(e)})):e \in E_n(1,0)\big\} 
\sbt {\mD}_l(U).
$$
Furthermore, analogously, for every $ D \in {\mD}_{1,0}(n)$, we  let
$$
E_n(D)
:= \big\{e \in E_n(1,0):\hat{f}^l({\psi}_e(U_{t(e)}))=D\big\}.
$$
Then, by the same reasoning as the one leading to (\ref{3fpn29}), we get
$$ \sum_{e \in E_n(D)}\diam^t(\psi_e(U_{t(e)})) \comp \diam^t(D)$$
for all $ n \geq 0$ and all $D \in \mD_{1,0}(n)$. Now, similarly as (\ref{4fpn30}), we get that 
$$ 
\begin{aligned}
\Sg_{n}(1,0)
:=& \sum_{e \in E_n(1,0)}\diam^t( \hat{\phi}_e(X_{t(e)}))\\
= & \sum_{D \in \mD_{1,0}(n)}\sum_{e \in E_n(D)}\diam^t(\hat{\phi}_e(X_{t(e)}))\\
 & \comp \sum_{D \in \mD_{1,0}(n)}\sum_{e \in E_n(D)}\diam^t(\psi_e(U_{t(e)})))\frac{\diam^t(\hat{f}^{-ln}_{p_1(t(e))}(X_{t(e)}))}{\diam^t(U_{t(e)})}\\
 & \comp \sum_{D \in \mD_{1,0}(n)}\sum_{e \in E_n(D)}\diam^t(\psi_e(U_{t(e)}))(n+1)^{- \frac{ p(p_1(t(e)))+1}{p(p_1(t(e)))}t}\\
&\comp(n+1)^{- \frac{ p(p_1(t(e)))+1}{p(p_1(t(e)))}t}\sum_{D \in \mD_{1,0}(n)}\sum_{e \in E_n(D)}\diam^t(\psi_e(U_{t(e)}))\\
&\comp(n+1)^{- \frac{ p(p_1(t(e)))+1}{p(p_1(t(e)))}t}\sum_{D \in \mD_{1,0}(n)}\diam^t(\psi_e(U_{t(e)}))\\
&\comp(n+1)^{- \frac{ p(p_1(t(e)))+1}{p(p_1(t(e)))}t}.
\end{aligned}
$$
Now, similarly as (\ref{1fpn30}), we get the following  
$$
\sum_{e \in E(1,0)}\diam^t(\hat{\varphi}_e(X_{t(e)}))
=\sum_{n=0}^\infty \Sg_n(1,0)\comp \sum_{n=0}^\infty (n+1)^{- \frac{ p(p_1(t(e)))+1}{p(p_1(t(e)))}t}< +\infty.
$$
The proof of Lemma~\ref{l2fpn27} is  complete. 
\qed

\sp Now we want to form an appropriate strongly regular conformal maximal graph directed Markov system on $\mathbb C$ and, by Young's tower  construction, to obtain from it the stochastic laws such as the exponential decay of correlations, and the Central Limit Theorem. It will also allow us to show that the metric entropy $\h_{\mu_h}(f)$ is finite.

\sp Keep $U$, the pre--nice set, produced in Theorem~\ref{t1fpn20}, for the map $\hat f^l:\hat\mT_f^{(l)}\lra\mT_f$ and
$$
\hat\cS_U=\big\{\hat\phi_e\big\}_{e\in\hat E},
$$ 
the corresponding conformal maximal graph directed Markov system also produced in Theorem~\ref{t1fpn20}. Set 
$$
F_0^*:=\Om(f)=\Om(f^l) 
$$
and define $F_1^*$ to be a (finite) selector of the partition of $\Pi_f^{-1}(B_\infty(\hat f))=f^{-1}(\infty)$ into equivalence classes of the equivalence relation $\sim_f$, i.e. choose exactly one point from each such equivalence class intersected with $f^{-1}(\infty)$. 
As the set of vertices take 
$$
\cV:=\big\{(\xi,j):\xi\in F_0^*\cup F_1^*, \ j=\{1,2,\ld,p_{\Pi_f(\xi)}\}\big\}.
$$
If the radius $r>0$ appearing in Theorem~\ref{t1fpn20} is small enough then all the maps
$$
\Pi_{\xi}:=\Pi_f|_{B(\xi,6r)}:B(\xi,6r)\lra \mT_f, \  \  \xi\in F_0^*\cup F_1^*,
$$
are $1$--to--$1$. Keeping $v=(\xi,j)\in\cV$, put
$$
\Pi_v:=\Pi_\xi,
$$
and set
$$
W_v^*:=\Pi_v^{-1}\(W(\xi,j)\) 
\  \  {\rm and} \  \
X_v^*:=\Pi_v^{-1}\(X(\xi,j)\).
$$
If $\xi\in F_1^*$, we will frequently simply write $X_\xi^*$ and $W_\xi^*$ for $X_v^*$ and $W_v^*$ respectively. As the set of edges, set
$$
E:=\big\{e\in\hat E:\exists v\in\cV \ \ {\rm s. t. } \ f^{N_{\cU}(e)-1}\(\^f(W(t(e))\)=W_v^*\big\},
$$
where, we recall, the map 
$$
\^f:\mT_f\lra\oc
$$ 
is given by \eqref{120180816}. For each edge $e\in E$ write $t(e)=(b,j)$, $b\in F_0\cup F_1$, $j\in\{1,2,\ld, p_b\}$, and let $t(e):=(\xi,j)$, where $\xi$ is the only element of $(F_0^*\cup F_1^*)\cap\Pi_f^{-1}(b)$. Likewise for $i(e)$.
For each $e\in E$ define
$$
\phi_e:=\Pi_{i(e)}^{-1}\circ\hat\phi_e\circ\Pi_{t(e)}:W_{t(e)}^*\lra W_{i(e)}^*.
$$
Obviously 
\beq\label{120181123}
\cS_U:=\big\{\phi_e:W_{t(e)}^*\lra W_{i(e)}^*\big\}_{e\in E}
\eeq
is a conformal maximal graph directed Markov system with eadges $E$, vertices $\cV$, and the sets $X_v^*$, $v\in \cV$. Denote by $J_{U}$ its limit set of $\cS_U$. Let 
\beq\label{220181123}
X^*:=\bu_{v\in \cV}X_v^*.
\eeq
and for every $ x \in U$ let $U(x)$  be the connected component of $U$  containing $x$. 

\sp Because of Observation~\ref{o120200123}
and since $\Crit_\infty(f)=\es$ for all parabolic elliptic functions, as an immediate consequence of Theorem~\ref{tmaincm}, Theorem~\ref{tinv}, Theorem~\ref{t1020907}, and Proposition~\ref{p5050206}, we get the following.

\sp\bthm\lab{t120190619} 
If $f:\C\lra\oc$ is a parabolic elliptic function, then

\begin{enumerate}

\,

\item $\H_s^h(J(f))=0$ and $\Pi^h_s(J(f))=+\infty$.

\,

\item There exist a unique $t\ge 0$ and a unique spherical $t$--conformal atomless probability measure $m_{t,s}$ for $f:J(f)\lra J(f)\cup\{\infty\}$. Then $t=h$.

\,

\item All other conformal measures are purely atomic, supported on $\Sing^-(f)$ with exponents larger than $h$. 

\,

\item There exists a unique, up to a multiplicative factor, $\sg$--finite $f$--invariant measure $\mu_h$ absolutely continuous with respect to $m_h$. 

\,

\item The measure $\mu_h$ is atomless on $\oc$,

\,

\item The measure $\mu_h$ is equivalent to $m_h$, 

\,

\item The measure $\mu_h$ is metrically exact, in particular ergodic\index{(N)}{ergodic measure} and conservative\index{(N)}{conservative measure}, 

\,

\item $\mu_h(\Tr(f))=1$.

\,

\item The measure $\mu_h$ is given, according to Theorem~\ref{t1h75}, by formulas (\ref{5.9a})--(\ref{eq:muequ}). 
\end{enumerate}
\ethm

\sp We first of all need the following.

\sp\bthm\label{l120131004B}
Let $f:\mathbb C \lra \oc$ be a parabolic elliptic function. 
If $U$ is the pre--nice set produced in Theorem~\ref{t1fpn20} for the map $\hat f^l:\hat\mT_f^{(l)}\lra\mT_f$ and $\cS_U$ is the conformal maximal graph directed Markov system, defined by formula \eqref{120181123}, then  
$$
m_h(J_{U})>0 \  \  \  \text{ and }  \  \  \  \HD(J_{U})=\HD(J(f)), 
$$
where, we recall, $m_h$ is the $h$--conformal measure for $f$ on $J(f)$.

Furthermore, the maximal conformal graph directed Markov system $\cS_{U}$ is strongly regular (in the sense of Definition~\ref{d5_2017_11_18}) and $m_h/m_h(J_{U})$ is the $h$--conformal measure (in the sense of Theorem~\ref{t420131007} with $t=h$) for the system $\cS_{U}$. 
\ethm

\bpf
By the construction of the system $\cS_{U}$, and by the definition of the limit set $J_{U}$, the set $J_{U}$ contains the set of all transitive points of the map $f^l:J(f)\to J(f)$, denoted, as usually, by $\Tr(f^l)$, that belong to $X^*$. Since $m_h$ is of full topological support, it therefore follows from Theorem~\ref{tmaincm} that 
$$
m_h(J_{U})\ge m_h(X^*\cap\Tr(f^l))=\hat m_h(X^*)>0.
$$
Since the measure $m_h$ is $h$-conformal for the map $f:\mathbb C \to \ov{\mathbb C}$, it also satisfies condition (f) of Theorem~\ref{t420131007} with $t=h$ and $\g=1$. It therefore follows from this theorem that $\P(h)=0$ and $m_h/m_h(J_{U})$ is the $h$--conformal measure for the conformal maximal graph directed Markov system $J_{U}$. Consequently, the system $\cS_{U}$ is regular and (the same theorem) $h=h_{\cS_{U}}$. But, by Theorem~\ref{t1j97}, $h_{\cS_{U}}=\HD(J_{U})$. 

So, we only need to show that the system $\cS_{U}$ is strongly regular. Exactly as the proof of strong regularity in Theorem~\ref{t120180509}, the proof of t strong regularity of the system $\cS_V$ is an immediate consequence of the facts that the system $\hat\cS_{U}$ is strongly regular (see  Lemma~\ref{l2fpn27}), $\P_{\cS_{U}}(h)=\P_{\hat\cS_{U}}(h)=0 $, $\P_{\cS_{U}}(t)\le \P_{\hat\cS_{U}}(t)$ for all $t\ge 0$, the latter holding because $E\sbt \hat E$ and all the maps $\Pi_v$, $v\in V$, are isometries, and both maps $\(\th_{S_{\hat V}},+\infty\)\ni t\longmapsto \P_{\cS_{U}}(t), \P_{\hat\cS_{U}}(t)$ are strictly decreasing. The proof is thus complete.
\epf 

\sp We define the map $F_{U}:J_{U}\lra J_{U}$ by the formula:
$$
F_{U}(\phi_e(z))=z \  \text{ if }  \  \ e \in E \ 
 \text{ and }  \  z\in J_{U}.
$$
As an immediate consequence of $f$--invariance of measure $\mu_h$, item (4) of Theorem~\ref{tmaincm}, item (1) of Theorem~\ref{tinv}, the fact that every transitive point of $f$ enters (infinitely often) $J_{U}$, and of Theorem~\ref{t1j64}, we get the following. 

\blem\label{l120181220}
If $f:\C\lra\oc$ is a parabolic elliptic function and $U$ is the pre--nice set produced in Theorem~\ref{t1fpn20} for the map $\hat f^l:\hat\mT_f^{(l)}\lra\mT_f$, then the induced
system $F_{\cU}:J_{U}\to J_{U}$ preserves the Borel probability conditional measure 
$$
\mu_{U}:=(\mu_h(J_{U}))^{-1}\mu_h|_{J_{U}}
$$
on $J_{U}$.
\elem

Now, as a fairly immediate consequence of the Theorem~\ref{l120131004B}, especially its strong regularity part, we shall prove the following.

\begin{thm}\label{t120121211B}
If $f:\C\lra\oc$ is a parabolic elliptic function and if $\mu_h$ is the corresponding Borel $\sg$--finite $f$--invariant measure
equivalent to the $h$--conformal measure $m_h$, then $\h_{\mu_h}(f)<+\infty$.
\ethm 

\bpf 
It is a direct consequence of Theorem~\ref{t120180509} and
Theorem~\ref{GDMS_Finite_Entropy} to have finite entropy of the induced
system $F_{U}:J_{U}\lra J_{U}$ with respect to the probability $F_{U}$--invariant
measure $\mu_{U}=(\mu_h(J_{U}))^{-1}\mu_h|_{J_{U}}$. Then, as $F_{U}:J_{U}\lra J_{U}$ is the first return map of $f$ from $J_{U}$ to $J_V$, Abramov's Formula
(Theorem~\ref{abramovabstract}) yields $\h_{\mu_h}(f)<+\infty$. 
\epf

\sp We will need the following two technical facts. The first one is this.

\blem\label{l12018_08_06}
Let $f :\mathbb C \lra \oc$ be a parabolic elliptic function. If $U$ is the pre--nice set of Theorem~\ref{t1fpn20}, then there exists a constant $C\in(0,+\infty)$ such that 
$$
\lt\|\(f^k\circ \phi_\tau\)'\rt\|_\infty\le Ce^{-\b_U|\tau|}
$$
for all $\tau\in E_A^*$ and all integers $0\le k\le |\tau_1|$, where $\b_U\in(0,1)$ is the contracting factor of the graph directed Markov system $\mathcal S_{U}$. 
\elem 

\fr {\sl Proof.} Denoting the infimum of Lemma~\ref{l120180306} (3) by $C_1$, we get for every $x\in X_{t(\tau)}$, that
$$
\big|\(f^k\circ \phi_{\tau}\)'(x)\big|
=\big|\(f^k\circ \phi_{\tau_1}\)'\(\phi_{\sg(\tau)}(x)\)\big|
\big|\(\phi_{\sg(\tau)}\)'(x)\big|
\le C_1e^{-\b_U|\sg(\tau)|}
=C_1e^{\b_U}e^{-\b_U|\tau|}.
$$
So, taking $C:=C_1e^{\b_U}$ finishes the proof.
\qed

\sp Passing to the second technical fact, let 
$$
\rho_s:=\frac{d\mu_h}{dm_{h,s}}
$$
\index{(S)}{$\rho_s$} be the Radon--Nikodym derivative of the measure  $\mu_h$ with respect to $m_{h,s}$ and let 
$$
\rho_e:=\frac{d\mu_h}{dm_{h,e}}
$$
\index{(S)}{$\rho_e$} be the Radon--Nikodym derivative of the measure $\mu_h$ with respect to $m_{h,e}$. We need to know the behavior of these derivatives near infinity. We shall prove the following..

\blem\label{l1j183} 
If $f:\C\lra\oc$ is a parabolic elliptic function, then there exist constants $R_0>0$ and $M\geq 1$ such that  
$$
M^{-1} |z|^{-\frac{q+1}{q}h}\leq \rho_e(z) \leq M |z|
^{-\frac{q+1}{q}h}
$$ 
for all $z\in J(f)\cap B_\infty^*(R_0)$.
\elem

\bpf Because of (\ref{u1}) there exist $R_1 >0$  and
$A \geq 1$ such that
\beq\label{1j183}
A^{-1}|f(z)|^{\frac{q+1}{q}} \leq | f'(z)| \leq A | f(z)
|^{\frac{q+1}{q}}
\eeq
for all $b \in f^{-1}(\infty)$ with  $q_b=q$ and all $z \in B_b(R_1)$. Because of Lemma~\ref{l1pj3}, the function 
$$
\hat{\rho}_h : J(f\sms \Pi (\ov{\PC(f)})\lra (0,
+\infty)
$$ 
is continuous. Hence, it follows from (\ref{1j181}) that there exists $R_0\ge R_1$ such that for all $z\in J(f)\cap B_\infty^*(R_0)$, we have
$$
\frac{1}{2}\sum_{b \in f^{-1}_0(\infty)}\sum_{w \in f^{-1}(z) \cap B_b(R_1)} 
\!\!\!\!\!\hat{\rho}_h (\Pi(b))|f'(w)|^{-h} 
\le \rho_e(z)
\le 2 \sum_{b \in f^{-1}_0(\infty)}\sum_{w \in f^{-1}(z) \cap
B_b(R_1)} \!\!\!\!\!\hat{\rho}_h (\Pi(b))|f'(w)|^{-h},
$$
where $f^{-1}_0(\infty)$ is given by Definition~\ref{przeciwobrazy}. Putting
\index{(S)}{$M_{-}$} 
$$
M_{-}:= \min\{ \hat{\rho}_h(\Pi(b)):  \, b \in
f^{-1}(\infty)\}\ \text{ and } \ \index{(S)}{M_{+}} 
M_{+}:= \max\{\hat{\rho}_h(\Pi(b)):  \, b \in f^{-1}(\infty)\},
$$
we further get
$$
\frac{M_{-}}{2}
\sum_{b \in f^{-1}_0(\infty)}\sum_{w \in f^{-1}(z)\cap B_b(R_1)}|f'(w)|^{-h}
\leq \rho_e(z)
\leq 2M_{+}\sum_{b \in f^{-1}_0(\infty)}\sum_{w \in f^{-1}(z)\cap B_b(R_1)}
|f'(w)|^{-h}.
$$ 
Invoking  finally (\ref{1j183}), we get
$$ 
\frac{M_{-}A^{-h}}{2}\#(f^{-1}_0(\infty))|z|^{-\frac{q+1}{q}h}
\leq \rho_e(z)
\leq 2 M_{+}\#(f^{-1}_0(\infty))|z|^{-\frac{q+1}{q}h},
$$ 
whenever $|z|\geq R_0$. We are done. \endpf

\sp\section[Parabolic Elliptic Maps with Finite Invariant Conformal Measures]  {Parabolic Elliptic Maps with Finite Invariant Conformal Measures; Statistical Laws; \\ Young Towers and Nice Sets Techniques}\label{parabolicfiniteinvariantmeasures}

In this section we continue dealing with parabolic elliptic functions $f:\C\lra\oc$. As an immediate consequence of Theorem~\ref{t4050206}, Proposition~\ref{p5050206}, and Theorem~\ref{t1j125} we get the following simple characterization of parabolic elliptic functions for which the invariant measure $\mu_h$ of Theorem~\ref{t120190619}is finite (and so, also for which it is infinite):

\bthm\lab{t2_2017_09_21} 
If $f:\C\lra\oc$ is a parabolic elliptic function, then the $f$--invariant measure $\mu_h$ of Theorem~\ref{t4050206} is finite if and only if 
$$
h=\HD(J(f))>\frac{2p_{\max}(f)}{p_{\max}(f)+1},
$$
where $p_{\max}(f):=\max\{p(\om):\om\in\Om(f)\}$.
\ethm

We assume throughout this section that the measure $\mu_h$ of Theorem~\ref{t120190619} is finite. We call such elliptic functions of finite class. \index{(N)}{parabolic elliptic functions of finite class} \index{(N)}{parabolic elliptic functions with finite measure $\mu_h$} We then always normalize the measure $\mu_h$ so that
$$
\mu_h(J(f))=1.
$$

Similarly as in Section~\ref{Young's Towers for Subexpanding}, we now want to apply the Young's Tower Machinery of Section~\ref{Young-Abstract}
too. Again, the primary  issue is to estimate well enough the measure $m_h$ of the set of  points in $J_{U}$ whose return time to $J_{U}$ is greater than or equal to $n$. However, now our approach is different than in Section~\ref{Young's Towers for Subexpanding}; it is closer to that  of  \cite{StratU1} and  \cite{StratU2}. In addition, it is much more convenient with this method, we  could  actually say that it is indispensable,  to work  with the invariant measure $\mu_h$ rather than the conformal one, $m_h$. For the sake of this method we bring up in the upcoming section the concept of the first entrance time to the sets $J_U$; it was in fact already introduced in Section~\ref{Young-Abstract} in the context of Young's Towers. We keep the setting and notation of the previous section, i.e. Section~\ref{PEM-I}. 

\sp\subsection{The First Entrance and the First Return Times} The first entrance time to the set $J_\cU$ is the function $E_{f^l,J_{U}}: J(f) \lra [0, \infty]$ defined as
$$ 
E_{f^l,J_{U}}(z):= \min\big\{ n \in \{0, 1, 2, \ldots,+\infty\}: \, f^{ln}(z) \in J_{U}\big\}.
$$
On a a set of full measure $m_h$, equivalently $\mu_h$, this function coincides with the function $E_{f^l,X}: J(f) \to [0, \infty]$, given by the formula: 
$$ 
E_{f^l,X^*}(z):= \min\big\{ n \in \{0, 1, 2, \ldots,+\infty\}: \, f^{ln}(z) \in X^*\big\}.
$$
We recall that 
$$
p_{\max}=\max \{ p_{f^l}(\om):\om \in \Omega(f)\}.
$$

\sp The main technical result of this section is the following.

\blem\label{l1fp5}
Let $f:\mathbb C \lra \ov{\mathbb C}$ be a parabolic elliptic function of finite class. If $X^*$ is the set defined by formula \eqref{220181123}, then there exists a constant $C \in(0, \infty)$ such that
$$ 
\mu_h\(\{z\in J(f):E_{f^l,J_{U}}(z) \geq n \}\) 
=\mu_h\(\{z\in J(f):E_{f^l,X^*}(z) \geq n \}\) 
\leq Cn^{2 - \frac{ p_{\max}+1}{p_{\max}}h}
$$  
for every integer $n \geq 1$.
\elem
 
\fr {\sl Proof.} Fix $b\in F_1^*$. Since $f^l\(\Int(X_b^*)\)\sbt \ov{\mathbb C}$ is an open set containing $\infty$ ($f^l\(\Int(X_b^*)\)=\ov{\mathbb C}$ if $l\ge 2$), there exists a fundamental domain ${\mathcal R} \sbt \mathbb C$ for  $f$ and an open ball $ B \supset {\mathcal R}$ such that
$$ 
f^l\(\Int(X_b^*)\) \supset B \quad \mbox{and} \quad 2B \cap \ov{\PC(f)}=\es.
$$
Then there exists a holomorphic inverse branch $f^{-l}_b:B \lra \mathbb C$ of $f$ such that
\beq\label{3fp5}
f^{-l}_b(B)\sbt\Int(X_b^*)
\eeq
and 
\beq\label{1fp5}
f(B) =\ov{\mathbb C}.
\eeq
Since the map $f:\mathbb C \lra \oc$ is parabolic, we  can partition the set
$$
J(f)\sms \left(\Omega(f)\cup \bigcup_{\om \in \Omega(f)} \bigcup_{j=1}^{p(\om)} \bigcup_{n=0}^\infty f^{-ln}_\om \(X_{(\om,j)}^*\)\right)     
$$
into countably many Borel sets $\{B_k\}_{k=1}^\infty$, such that
\beq\label{220181128}
\diam(B_k)\le 1
\eeq
for all $k\ge 1$,
\beq\label{2fp5}
\bigcup_{k=1}^\infty B\(B_k, \diam(B_k)\)\cap \ov{\PC(f)}=\es,
\eeq
and 
\beq\label{120181128}
\lim_{k \to \infty}\dist (0, B_k)=+\infty.
\eeq
In particular, all the  holomorphic branches of $f^{-ln}$, $n \geq 0$, are well defined on all sets $B\(B_k, \diam(B_k)\), \, k \geq 1$, and so,
Koebe's Distortion Theorems applies to all these holomorphic branches restricted to the sets $B_k,\, k \geq 1$. Now, fix an integer $n \geq 0$ and  $k\geq 1$. Put
$$
\Gamma_k:= B_k \cap E^{-1}_{f^l,X^*}([n,\infty]).
$$
Because of (\ref{1fp5}) there exists a holomorphic branch $f^{-l}_k: B\(B_k,\diam(B_k)\) \lra \mathbb C$ such that
\beq\label{1fp6}
f^j\(f^{-l}_k\(B(B_k,\diam(B_k))\)\)\sbt B
\eeq
for every integer $j=0, 1, 2,\ld, l-1$. Hence, looking up at this formula with $j=0$, we see that the composition
$$ 
f^{-l}_b \circ f^{-l}_k: B\(B_k, \diam(B_k)\)\lra \mathbb C
$$
is well defined, is holomorphic, and it follows from (\ref{3fp5}) and (\ref{1fp6}) that
\beq\label{2fp6}
f^{-l}_b \circ f^{-l}_k\(B\(B_k, \diam(B_k)\)\)\sbt\Int(X_b^*).
\eeq
Fix a point $\xi_k\in \Gamma_k$. Because of Koebe's Distortion  Theorem I, Euclidean Version, we have
\beq\label{3fp6}
\begin{aligned}
\mu_h\(f^{-l}_b \circ f^{-l}_k ( \Gamma_k)\)
& = \int_{\Gamma_k} \frac{ \rho_e\(f^{-l}_b \circ f^{-l}_k(z)\)}{\rho_e(z)}\big|(f^{-l}_b \circ f^{-l}_k)'(z)\big|^hd\mu_h(z)\\
& \geq C_1C_2\big|(f^{-l}_k)'(\xi_k)\big|^h\int_{\Gamma_k}\frac{1}{\rho_e(z)}d\mu_h(z),
\end{aligned}
\eeq
where
$$ 
C_1:=\inf\{\rho_e(x): \, x \in \Int(X_b^*)\} >0 
\quad \mbox{and} \quad 
C_2:= \inf\{|(f^{-l}_b)'(x)|:\,  x \in B\} >0.
$$
It follows from \eqref{120181128} and \eqref{220181128} that the sets $f^{l-1}\circ f^{-l}_k(B_k)$ accumulate towards poles in $B$. We can, in fact, choose the holomorphic branches $f_k^{-l}$ such that
$$
\lim_{k \to \infty}\Dist \(a, f^{l-1}\circ f_k^{-l}(B_k)\) = 0 
$$
and
$$
\ov{\bigcup_{k=1}^\infty f^{l-1}\circ f^{-l}_k(B_k)}\cap \left(f^{-1}(\infty)\sms\{a\}\right) =\es,
$$
where $ a \in f^{-1}(\infty) \cap B $  is such a pole that $q_a=q_{\max}$. Then, it follows from the behavior of meromorphic functions around poles and using also \eqref{1fp6}, we conclude that there exists a constant $C_3\in (0, \infty)$ such that
\beq\label{4fp6}
\big|\(f^{-l}_k\)'(\xi_k)| \geq C_3|\xi_k|^{-\frac{q+1}{q}}
\eeq
for all integers $k\ge 1$, where, as usually, we  denote $q_{\max}$ by $q$. On the other hand, it follows from Lemma~\ref{l2fpn20} that the set
\beq\label{120181206}
V_{\Omega(f)}:= \Omega(f) \cup  \bigcup_{\om \in \Omega} \bigcup_{j=1}^{p(\om)}  \bigcup_{n=0}^\infty f_\om^{-ln}\(X_{(\om,j)}^*\)
\eeq
is a neighbourhood  of $\Omega(f)$. Since all sets $B_u$, $u\ge 1$, are disjoint from this neighbourhood, we deduce from Theorem~\ref{t1pj11} and Lemma~\ref{l1j183} that there exists a constant $C_4\in (0, \infty)$ such that
\beq\label{120181126}
\rho_e(z) \leq C_4|\xi_u|^{-\frac{q+1}{q}h}
\eeq
for all integers $u\ge 1$ and all $z \in B_u$. Coming back to our fixed integer $k\ge 1$, inserting \eqref{120181126} to (\ref{3fp6}), and using (\ref{4fp6}), we get that
$$
\begin{aligned} 
\mu_h\(f^{-l}_b \circ f^{-l}_k( \Gamma_k)\)
& \geq C_1C_2C_3^h |\xi_k|^{-\frac{q+1}{q}h} C_4^{-1}|\xi_k|^{\frac{q+1}{q}h}\mu_h(\Gamma_k)\\
& =C_1C_2C_3^hC_4^{-1}\mu_h(\Gamma_k).
\end{aligned}
$$
Since, because of (\ref{2fp6}) $E_{f^l,X^*}(f^{-l}_b \circ f^{-l}_k( \Gamma_k)) \sbt \{0\},$ we thus have that
$$
\begin{aligned} 
\mu_h\lt(f^{-2l}( \Gamma_k) \cap E_{f^l,X^*}^{-1}([n+2,+\infty])\rt) 
& \leq \mu_h\(f^{-2l}( \Gamma_k)\)-\mu_h\(f^{-l}_b\circ f^{-l}_k(\Gamma_k)\)\\
& \leq \mu_h ( \Gamma_k) - C_1C_2C_3^hC_4^{-1}\mu_h(\Gamma_k)\\
&= (1- C_1C_2C_3^hC_h^{-1})\mu_h(\Gamma_k).
\end{aligned}
$$
So, setting 
$$ 
\alpha:=1- C_1C_2C_3^hC_h^{-1} \in [0,1],
$$
we have that
$$ 
\mu_h\lt(f^{-2l}( \Gamma_k) \cap E_{f^l,X^*}^{-1}([n+2,+\infty])\rt) 
\leq \alpha\mu_h(\Gamma_k).
$$
Recalling that $\{B_k\}_{k=1}^\infty$ is a partition of $J(f)\sms V_{\Omega(f)}$, by summing up the above inequality over all $k \geq 1$, we obtain
\beq\label{1fp7}
\mu_h\lt( f^{-2l}(J(f)\sms V_{\Omega(f)}) \cap E_{f^l,X^*}^{-1}([n+2,+\infty])\rt) 
\leq \alpha \mu_h\lt((J(f)\sms V_{\Omega(f)}) \cap E_{f^l,X^*}^{-1}([n,+\infty])\rt).
\eeq
Now suppose that $z\in V_{\Omega(f)}\sms \Omega(f)$. This means that there exists $\om \in \Omega(f)$, $ j \in \{1, 2, \ldots, p(\om)\}$  and $ i \geq 1$ such that
$$ 
z \in f^{-li}_\om\(X_{(\om,j)}^*\).
$$
Then $f^{li}(z) \in X_{(\om,j)}^* \sbt X^*$. Therefore,
$$ 
E_{f^l,X^*}(z) \leq i. 
$$
Thus 
$$
\(V_{\Om(f)}\sms\Omega(f)\)\cap E^{-1}_{f^l,X^*}([n, \infty])  
\sbt \bigcup_{i \geq n}f^{-li}_\om\(X_{(\om,j)}^*\).
$$
So using Proposition~\ref{p1j111DL}, applied to the function $f^l$, and Proposition\ref{p120181222}, we conclude that there exist constants
$C_5, C_6\in(0, \infty)$ such that
$$
\begin{aligned}
\mu_h \lt(V_{\Omega(f)} \cap E_{f^l,X^*}^{-1}([n,+\infty])\rt) 
&= \mu_h\lt(V_{\Omega(f)}\sms \Omega(f)) \cap E_{f^l,X^*}^{-1}([n, \infty])\rt)\\
&\leq \sum_{\om\in\Omega(f)}\sum_{j=1}^{p(\om)} \sum_{i\geq n}\mu_h\(f^{-li}_\om\(X_{(\om,j)}^*\)\)  \\
& \le C_5\sum_{ \om \in \Omega} \sum_{j=1}^{p(\om)} \sum_{i \geq n} i^{1-\frac{p(\om)+1}{p(\om)}h }\\
 &  \leq C_5\# \Omega p \sum_{i \geq n} i^{1-\frac{p+1}{p}h },
\end{aligned}
$$
where $p:=p_{\max}$. Since, by our hypotheses, $1-\frac{p+1}{p}h < -1$, we further obtain, with some constant $C_6\in (0,\infty)$, that
$$
\mu_h\(V_{\Omega(f)} \cap E^{-1}_{f^l,X^*}([n, \infty])\) 
\leq C_6n^{2-\frac{p+1}{p}h }.
$$
Since the measure $\mu_h$ is $f$--invariant, we therefore get
$$
\begin{aligned}
\mu_h\(f^{-2}\(V_{\Omega(f)}\) \cap E_{f,U}^{-1}([n+2, \infty])\) 
&\leq \mu_h\(f^{-2}\(V_{\Omega(f)} \cap E^{-1}_{f,U}([n, \infty])\) \\
&= \mu_h\(V_{\Omega(f)} \cap E^{-1}_{f,U}([n, \infty])\) \\
&\leq C_7n^{2-\frac{p+1}{p}h}.
\end{aligned}
$$
Combining this with (\ref{1fp7}), we  get
$$  
\mu_h\(E^{-1}_{f^l,X^*}([n+2, \infty])\) 
\leq \alpha\mu_h\(E^{-1}_{f,U}([n,\infty])\)+C_6n^{2-\frac{p+1}{p}h}.
$$
So, there exists a constant $C_7\in [C_6, \infty)$ such that
\beq\label{1fp8}
\mu_h\(E^{-1}_{f,U}([n+2, \infty])\) 
\leq \alpha\mu_h\(E^{-1}_{f^l,X^*}([n, \infty])\) +C_7(n+2)^{2-\frac{p+1}{p}h }.
\eeq
Fix an integer $s\geq 2$ so large that
\beq\label{2fp8}
\beta:=\alpha\left(\frac{s+2}{s}\right)^{\frac{p+1}{p}h-2 }<1.
\eeq
Put
$$ 
C:=\max\Big\{C_7(1-\beta)^{-1}, \max\big\{ k^{\frac{p+1}{p}h}\mu_h\(E^{-1}_{f^l,X^*}([k, \infty])\):   \,  1\leq k \leq s\big\}\Big\}.
$$
We shall prove by induction that
\beq\label{2fp9}
\mu_h\(E^{-1}_{f^l,X^*}([n, \infty])\) \leq  C n^{2-\frac{p+1}{p}h}
\eeq
for every integer $n\ge 1$. Indeed, for $ 1\leq n\leq s$ this is immediate from the definition of $C$. So suppose for the inductive step, that $n \geq s$ and (\ref{2fp9}) holds for all $ 1\leq k \leq n$. Then, using also (\ref{1fp8}) and (\ref{2fp8}), we  get
$$
\begin{aligned}
\mu_h\(E^{-1}_{f^l,X^*}([n+1, \infty])\) 
& \leq \alpha\mu_h\(E^{-1}_{f^l,X^*}([n-1,\infty])\)+C_7 (n+1)^{2-\frac{p+1}{p}h }\\
& \leq \alpha C(n-1)^{2-\frac{p+1}{p}h }+C_7(n+1)^{2-\frac{p+1}{p}h }\\
& \leq \beta  C(n+1)^{2-\frac{p+1}{p}h }+(1-\beta) C (n+1)^{2-\frac{p+1}{p}h }\\
& = C (n+1)^{2-\frac{p+1}{p}h }.
\end{aligned} 
$$
The inductive proof of (\ref{2fp8})  is thus  complete, and we are done.
\qed

\sp Let $N_U: J_{U} \to [1, \infty]$ be the first return time to $J_{U}$, {\it i.e.} $\forall z \in J_{U}$
$$ 
N_{U}(z):=\min\{n \geq 1:f^n(z) \in J_{U}\}.
$$
As an almost immediate consequence of Lemma~\ref{l1fp5}, we get the following result, which we will not truly need but which is interesting on its own and we will use it in a ``negative'' way.

\blem\label{l1fp6}
Let $f:\mathbb C \lra \oc$ be a parabolic elliptic function of finite class. If $U$ is the pre--nice set of Theorem~\ref{t1fpn20}, then there exists a constant $C\in (0,+\infty))$ such that 
$$
\mu_h\(\big\{z\in J_{U}:N_{U}(z)\geq n\big\}\) \leq  C n^{2-\frac{p_{\max}+1}{p_{\max}}h}
$$
for every integer $n \geq 1$.
\elem

\bpf Since
$$ 
\big\{ z\in U \cap J(f):N_U(z)\geq n\big\} 
\sbt \big\{z \in J(f): \, \, E_{f^l,J_{U}}(f^l(z))\geq n-1\}
= f^{-l}(E^{-1}_{f^l,J_{U}}([n-1, +\infty])),
$$
and since the measure $\mu_h$ is $f^l$-invariant, we get from Lemma~\ref{l1fp5} for all $n \geq 2$ that
$$
\mu_h\(\{ z \in J_{U}:N_U(z)\geq n\}\) 
\leq  C (n-1)^{2-\frac{p_{\max}+1}{p_{\max}}h} 
\leq  C'n^{2-\frac{p_{\max}+1}{p_{\max}}h},
$$
with some constant $C'\in(0,\infty)$, and the lemma  follows.
\epf

\sp Since the measures $\mu_h$ and $m_h$ restricted to our nice set $U$ are  equivalent with log bounded Radon--Nikodyn derivatives, as an immediate consequence of Lemma~\ref{l1fp6}, we get the following.

\blem\label{l1fp10} 
Let $f :\mathbb C \lra \oc$ be a parabolic elliptic function of finite class.
If $U$ is the pre--nice set of Theorem~\ref{t1fpn20}, then there exists a constant $ C \in (0,+\infty])$  such that
$$
m_h\(\big\{z\in J_{U}:N_{U}(z)\geq n\big\}\)\leq C n^{2-\frac{p_{\max}+1}{p_{\max}}h}
$$
for  every integer $n \geq 1$.
\elem

\sp\subsection{Young's Towers for Parabolic Elliptic Functions with Finite Invariant Measure $\mu_h$; Statistical Laws}

In this subsection we take fruits of the results obtained in the
previous subsection and in Section~\ref{Young-Abstract} in order to establish
a polynomial decay of correlations and the Central Limit Theorem for all parabolic elliptic functions $f:\C\lra\oc$ of finite class (i.e. with finite invariant measure $\mu_h$), and all H\"older continuous bounded observables. This subsection is very similar to Subsection~\ref{Young's Towers for Subexpanding}. We however present it in full for the sake of completeness, convenience of the reader, and since there is one subtlety in the proofs which needs to be taken care of differently than in Subsection~\ref{Young's Towers for Subexpanding}.

We keep in this section all hypotheses and notation from the previous one. In particular, $f:\C\lra\oc$ is a parabolic elliptic function of finite class. Furthermore, 
$$
U
$$
is the pre--nice set produced in the previous section, i.e. the one of Theorem~\ref{t1fpn20}. Let $\cS_{U}$ be the corresponding conformal graph directed Markov system also produced in Theorem~\ref{t1fpn20}. Keep $E$ to denote the set of its edges and $\cV$ to denote the set of vertices of $\cS_{U}$.

We recall that the map $F_{U}:J_{U}\lra J_{U}$ has been defined by the formula: 
$$
F_{U}(\phi_e(z))=z \  \text{ if }  \  \ e \in E \ 
 \text{ and }  \  z\in J_{U}.
$$
Our goal is to show that the system $(F_{\cU},m_{h})$ (in fact the measure $m_h$ is treated here as restricted to  $J_{U}$) fits into the framework of Section~\ref{Young-Abstract} and to check that the hypotheses of Theorems \ref{lsyoung0}, \ref{lsyoung}, and \ref{LILII} are satisfied for this (induced) system $(F_{U},m_{h})$. The stochastic laws will then automatically follow.

\begin{enumerate}
\, \item The space $\Delta_0$ is now $J_{U}$, the limit set of the iterated
function system ${\mathcal S}_{\cU}$. 

\, \item The partition $\mathcal{P}_0$ consists of the
sets $\Ga_e:=\phi_e(J_{U})$, $e\in E$. 

\, \item The measure $m$ is the measure $m_{h}$ restricted to $J_{U}$; it is positive because of Theorem~\ref{l120131004B}.

\, \item The map $T_0:\Delta_0\lra\Delta_0$ is, in our setting, just the map $F_{U}$. 

\, \item The function $R$, the return time, is, naturally, defined as  
$$
R|_{\Ga_e}:=N_{U}|_{\Ga_e}.
$$
\end{enumerate}
We also write
$$
N_U(e):=N_{U}|_{\Ga_e}
$$
Fix $e\in E$ arbitrary and then two arbitrary points $x, y\in \Ga_e
=\phi_e(J_{U})$. This means that $x=\phi_e(x')$ and
$y=\phi_e(y')$ with some $x', y'\in J_{U}$. Since
$t^h-1=h(t-1)+O(|t-1|^2)$, and because of  Theorem~\ref{tKoebe-1}, 
there exist respective constants $C_1, C_2>0$ such that
we have,
\beq\label{120121207_2}
\aligned
\lt|\frac{Jac_{m_h}F_{U}(y)}{Jac_{m_h}F_{U}(x)}-1\rt|
&=\lt|\frac{|F_{U}'(y)|^h}{|F_{U}'(x)|^h}-1\rt|
=\lt|\frac{|(\phi_e)'(x')|^h}{|(\phi_e)'(y')|^h}-1\rt|\\
&\le C_1\lt|\frac{|\phi_e)'(x')|}{|(\phi_e)'(y')|}-1\rt|\\
&\le C_1C_2|y'-x'|.
\endaligned
\eeq
Now write $x'=\pi_U(\a)$ and $y'=\pi_U(\g)$
with appropriate $\a,\g\in E$. Put $\om:=\a\wedge\g$ and
$k:=|\om|$. Write also $x'':=\pi_{U}(\sg^k(\a))$ and
$y'':=\pi_{U}(\sg^k(\g))$. Then we get that 
\beq\label{120121208_2}
\aligned
|y'-x'|
&=   |\phi_\om(y'')-\phi_\om(x'')|
 \le |\phi_\om'(x'')||y''-x''| \\
&\le \diam_e\(X_{t(\om)}^*\)\|\phi_\om'\|_\infty
 \le \diam_e\(X_{t(\om)}^*\)\b_{U}^{|\om|} \\
&=   \diam_e\(X_{t(\om)}^*\)\b_{U}^{s(x',y')}\\
&=\diam_e\(X_{t(\om)}\)\b_{U}^{s(F_U(x),F_U(y))}\\
&\le D_{U}\b_U^{s(F_{U}(x),F_{U}(y))},
\endaligned 
\eeq
where, we recall, $\b_U\in(0,1)$ is the contracting factor of the graph directed Markov system $\mathcal S_{\cU}$ and 
$$
D_{U}:=\max\big\{\diam_e(U_v):v\in F\big\}<+\infty.
$$
Since $s(F_{U}(x),F_{U}(y))=s(x,y)-1$ (as we know that $s(x,y)\ge 1$), along with \eqref{120121207_2}, the formula \eqref{120121208_2} gives that
$$
\lt|\frac{Jac_{m_h}F_{U}(y)}{Jac_{m_h}F_{U}(x)}-1\rt|
\le D_{U}\b_{U}^{s(F_U(x),F_{U}(y))}
= \b_{U}^{-1}D_U\b_{U}^{s(x,y)}.
$$
So, \eqref{youngcond0} is established in our context. The fact that
the partition $\mathcal{P}_0$ is generating follows from the
contracting property of the graph directed Markov system $\cS_U$. The Big Images Property holds because the alphabet of the GDMS $\cS_U$ is finite. The last assumption in Theorem~\ref{lsyoung0} is that the map $T:\Delta\lra\Delta$, which we now denote by $T_U$, is topologically mixing. The proof of this fact is virtually the same as the proof of Lemma~\ref{l120180515}  for subexpanding elliptic maps. We provide it here for the sake of completeness and convenience of the reader.

\blem\label{l120180515B}
The dynamical system $T_U:\Delta\lra\Delta$ is topologically mixing in the sense of Section~\ref{Young-Abstract}. 
\elem

\bpf
For every $\tau\in E_A^*$, set
$$
N_U(\tau)
:=\sum_{j=1}^{|\tau|}N_U(\tau_j)
=\sum_{j=1}^{|\tau|}R\(\Ga_{\tau_j}\).
$$
By our construction of the tower $\cY_f$, we have for every $e\in E$ that
$$
T_U^{R(\Ga_e)}\(\Ga_e\times\{0\}\)
=T_U(\Ga_e\times\{R(\Ga_e)\}\)
=F_U(\Ga_e)\times\{0\}
=J_{t(e)}\times\{0\}.
$$
So, by an immediate induction, we get for every $\tau\in E_A^*$ that
\beq\label{620180531B}
T_U^{N_\cU(\tau)}\(\phi_\tau(J_{t(\tau)})\times\{0\}\)
=J_{t(\tau)}\times\{0\}.
\eeq
Now fix two arbitrary elements $a, b\in E$. Then there exists $s\ge 0$ such that
$$
f^{lu}\(\phi_a(X_{t(a)}^*)\)\spt \phi_b(X_{t(b)}^*)
$$
for all $u\ge s$. Hence, for all such $u$ there exists a holomomorphic branch of $f^{lu}$ mapping $\phi_b\(X_{t(b)}^*\)$ into $\phi_a\(X_{t(a)}^*\)$. By our construction of the conformal graph directed Markov system $\cS_U$, this holomorphic branch is an admissible composition of elements of $\cS_U$. This means that it is equal to 
$$
\phi_{\tau(u)}:X_{t(\tau(u))}^*\lra X_{i(\tau(u))}^*
$$ 
for some $\tau(u)\in E_A^*$ with $t(\tau(u))=i(b)$ and $i(\tau(u))=i(a)$. Then, applying \eqref{620180531B}, we get for every $0\le p \le R(\Ga_a)-1$, every $0\le q \le R(\Ga_b)-1$, and every $n\ge s+R(\Ga_b)-1$, that 
$$
\begin{aligned}
T_U^n\(\Ga_a\times\{p\}\)
&=T_U^{n+p}\(\Ga_a\times\{0\}\)
=T_U^{q+n+p-q}\(\Ga_a\times\{0\}\) \\
&\spt T_U^q\(T_U^{n+p-q}\(\phi_{\tau(n+p-q)}(J_{t(\tau(n+p-q))})\times\{0\}\)\) \\
&=T_U^q\(J_{t(\tau(n+p-q))}\times\{0\}\)
\spt T_U^q\(\phi_b(J_b)\times\{0\}\) \\
&=\phi_b(J_b)\times\{q\} \\
&=\Ga_b\times\{q\}.
\end{aligned}
$$
The proof of Lemma~\ref{l120180515} is complete.
\epf 

\sp Since the measures $\mu_h$ and $m_h$, restricted to our nice set $U$, are  equivalent with log bounded Radon--Nikodyn derivatives, as an immediate consequence of Theorem~\ref{tmaincm} (particularly its item (2), the finiteness of the measure $\mu_h$, and Theorem~\ref{t1j69} (a) (Kac's Lemma), we get that
\beq\label{1_2018_08_03}
\int_{\De_0}R\,dm_h
\comp \int_{\De_0}R\,d\mu_h
<+\infty. 
\eeq
Therefore, as an immediate consequence of Observation~\ref{120200124} and Theorem~\ref{lsyoung0}, we get the following. 

\sp
\begin{thm}\label{lsyoung0_2018_08_04}
If $f:\C\lra\oc$ is a parabolic elliptic function of finite class, then 
\begin{itemize}
\item 
[(a)]
$\tilde m_h(\De)<+\infty$,  
where $\tilde m_h$ is derived out of $m_h$ restricted to $J_U$ according to formula \eqref{520180515},

\,

\item[(b)] There exists a unique probability $T_{U}$--invariant measure $\nu_h$, absolutely continuous with respect to $\tilde m_h$. 

\, 

\item[(c)]  The Radon--Niokodym derivative $d\nu_h/d\tilde m_h$ is bounded
from below by a positive constant. 

\, \item[(d)]  The dynamical system $\(T_U,\nu_h\)$ is metrically exact, thus ergodic.
\end{itemize}
\end{thm}

\sp Now consider $H:\Delta\lra \mathbb{C}$, the natural
projection from the abstract tower $\Delta$, generated by the above, to the complex plane $\mathbb{C}$ given by the formula
\beq\label{3_2018_08_04}
H(z,n)=f^{ln}(z).
\eeq
Then
\beq\label{eq120110623}
H\circ T_U =f^l\circ H,
\eeq
$$
\tilde m_h|_{\Delta_0}\circ H^{-1}=m_h|_{J_U},
$$ 
and  
$$
\tilde m_h|_{\Ga_e\times\{n\}} \circ H^{-1} 
=m_h|_{\Ga_e\times\{0\}}\circ f^{-ln}
=m_h|_{\Ga_e}\circ f^{-ln}
$$ 
for all $e\in\^E$ and all $0\le n\le N_U(e)$.
Now, $\tilde m_h|_{\Ga_e\times\{n\}}\circ H^{-1}$ is absolutely continuous
with respect to $m_h$, with the Radon--Nikodym derivative equal to
$$
J_{e,n}:=Jac_{m_h}\lt(\lt(f^{ln}|_{\Ga_e}\rt)^{-1}\rt)
$$ 
in $f^n(\Ga_e)$ and zero
elsewhere in $J(f)$. Therefore, using \eqref{1_2018_08_03}, we get that
$$
\int_{J(f)}\sum_{e\in E}\sum_{n=0}^{N_U(e)-1}J_{e,n}dm_h
= \sum_{\tau\in E}\sum_{n=0}^{N_U(e)-1}\!\!\int_{J(f)}\!\! J_{e,n}dm_h 
=\tilde m_h\circ H^{-1}(J(f)) 
=\tilde m_h(\De) 
<+\infty.
$$ 
Thus, the function 
$$
\sum_{e\in E}\sum_{n=0}^{N_U(e)-1}J_{e,n}
$$ 
is integrable with respect to the measure $m_h$. This implies immediately that
the measure $\tilde m_h\circ H^{-1}$ is absolutely continuous with
respect to the measure $m_h$ with the Radon--Nikodym derivative equal to 
$$
\sum_{e\in E}\sum_{n=0}^{N_U(e)-1}J_{e,n}.  
$$
Hence, the measure $\nu_h\circ H^{-1}$ is also absolutely continuous
with respect to $m_h$. Since $\nu_h$ is $F_U$--invariant and $H\circ
T_U =f^l\circ H$, the measure $\nu_h\circ H^{-1}$ is $f^l$--invariant. 
But the measure $\mu_h$ is $f$--invariant, ergodic, and
equivalent to the conformal measure $m_h$. Hence, $\nu_h\circ H^{-1}$
is absolutely continuous with respect to the ergodic
measure $\mu_h$. In conclusion, we get the following. 
 
\begin{lem}\label{l120110623B}
If $f:\C\lra\oc$ is a parabolic elliptic function of finite class, then 
$$
\nu_h\circ H^{-1}=\mu_h.
$$
\end{lem} 
Recall that the function $E_{T_{U}}$ is given by the formula \eqref{2_2018_08_04}. Using also \eqref{eq120110623} and \eqref{3_2018_08_04},
we get for every point $(z,k)\in\De$ that
$$
E_{T_{U}}(z,l)=E_{f^l,J_U}(f^{lk}(z))=E_{f^l,J_U}(H(z,k)).
$$
Therefore, using also Lemma~\ref{l120110623B}, Theorem~\ref{lsyoung0_2018_08_04}, and, at the end, Lemma~\ref{l1fp5}, we get 
\beq\label{4_2018_08_04}
\begin{aligned}
\^m_h\(E_{T_{U}}^{-1}(n,+\infty])\)
&\lek \nu_h\(E_{T_{U}}^{-1}(n,+\infty])\)
= \nu_h\(H^{-1}\(E_{f^l,J_U}^{-1}(n,+\infty])\)\)\\
&=\mu_h\(E_{f^l,J_U}^{-1}(n,+\infty])\) \\
&\lek n^{2-\frac{p_{\max}+1}{p_{\max}}h}.
\end{aligned}
\eeq 
Therefore, we conclude from Theorem~\ref{lsyoung} the following. 

\bthm\label{t32018_8_04}
If $f:\C\lra\oc$ is a parabolic elliptic function of finite class, then, for the dynamical system $T_U:(\Delta,\nu_h)\lra(\Delta,\nu_h)$, the following hold:

\begin{itemize}
\, \item[(a)] The Polynomial Decay of Correlations in the form of (\ref{cov_poly}) with the parameter $\a:=\frac{p_{\max}+1}{p_{\max}}h-2$, 

\, \item[(b)] If $h>\frac{3p_{\max}}{p_{\max}+1}$, then for all $\g\in (0,1)$ the Central Limit Theorem is true for all functions $g\in C_\g(\Delta)$ that are not cohomologous to a constant in $L^2(\nu_h)$. 
\end{itemize}
\ethm

\brem\label{r1_2018_08_06B}
Note that the hypothesis of item (b) has a chance to be satisfied only if $p_{\max}=1$. Then this hypothesis says that $h>3/2$. In order to verify that such inequality holds for particular examples, one may for example use Theorem~\ref{thm:julia}.

We did not formulate The Law of Iterated Logarithm as in order to establish it via Theorem\ref{LILII} and Lemma~\ref{l1fp10}, we would have to know that $h>2$, which, of course is never satisfied.
\erem

\sp We are now in position to prove the following.

\begin{thm}\label{t220110627B}
If $f:\C\lra\oc$ is a parabolic elliptic function of finite class, then the dynamical system $\(f^l|_{J(f)},\mu_h\)$ satisfies the following. If $g:J(f)\lra \mathbb{R}$ is a bounded function, H\"older continuous with respect to the Euclidean metric on $J(f)$, then

\, \begin{enumerate}
\item{} (The Polynomial Decay of Correlations) For
every bounded measurable function $\psi:J(f)\lra\mathbb{R}$, we have that
$$
\bigg|\int\psi\circ f^{ln} \cdot g d\mu_h
-\int g d\mu_h\int\psi d\mu_h\bigg|
\lek n^{2-\frac{p_{\max}+1}{p_{\max}}h}.
$$

\sp\item{} (The Central Limit Theorem) If in addition $g:J(f)\to\mathbb{R}$ is not cohomologous to a constant in $L^2(\mu_h)$ with respect to $f^l$, i.e. if there is no square integrable function $\eta$ for which $g={\rm const}+\eta\circ f^l-\eta$, then the The Central Limit Theorem holds for $g$. More precisely,
 there exists $\sigma>0$ such that
$$
\frac{\sum_{j=0}^{n-1}g\circ f^{lj}-n\int g d\mu_h}{\sqrt{n}} \xrightarrow[\ \, n\to\infty \ ]{} \mathcal N(0,\sigma)
$$
in distribution, where, as usually, $\mathcal N(0,\sigma)$ denotes the Gauss (normal) distribution centered at $0$ with covariance $\sg$.
\end{enumerate}
\end{thm}

{\sl Proof.} 
Keep
$$
U
$$
the pre--nice set produced in Theorem~\ref{t1fpn20} and $\cS_{U}$ the corresponding conformal graph directed Markov system also produced in Theorem~\ref{t1fpn20}.

Let $g:J(f)\lra\mathbb{R}$ and $\psi:J(f)\lra\mathbb{R}$ be as in the
hypotheses of our theorem. Define the functions 
$$
\tilde g:=g\circ H:\Delta\lra\mathbb{R} \  \  \text{ and } \ \
\tilde\psi:=\psi\circ H:\Delta\lra\mathbb{R}.
$$

We shall prove the following.

\sp\fr {\bf Claim~1:} The function $\tilde g:\Delta\lra\mathbb{R}$
belongs to the space $C_\beta$ for an appropriate exponent $\beta\in (0,1)$. 

\sp\fr Indeed, consider two arbitrary points $(x,k),(y,u)\in\De$. We
treat separately two cases depending on whether $s((x,k),(y,u))=0$ or
not. If $s((x,k),(y,u))=0$, then we get
\beq\label{220121208B}
\aligned
|\^g(y,u)-\^g(x,k)|
&  =|g(H(y,u))-g(H(x,k))|
   =|g(f^{lu}(y))-g(f^{lk}(x))|\\
&\le|g(f^{lu}(y))|+|g(f^{lk}(x))| \\
&\le2\|g\|_\infty \\
&=  2\|g\|_\infty\b^{s((x,k),(y,u))}
\endaligned
\eeq
regardless of what the value of $\b$ is, which will be specified in the
next case. Indeed, if $s((x,k),(y,u))>0$, then $k=u$, $0\le k<R(x)=R(y)$, and
\beq\label{320121208B}
|\^g(y,u)-\^g(x,k)|
=|g(f^{lk}(y))-g(f^{lk}(x))|
\le H_g|f^{lk}(y)-f^{lk}(x)|^\g,
\eeq
where $H_g$ and $\g$ are respectively the H\"older constant and
H\"older exponent of the function $g$. Moreover, $x,y\in\phi_\tau\(J_{t(\tau)}\)$, with some $\tau\in E_A^{s((x,k),(y,u))}$, and then 
$$
f^{lk}(x),f^{lk}(y)\in f^{lk}\(\phi_\tau\(J_{t(\tau)}\)\),
$$ 
$0\le k\le|\tau_1|$, and $f^{lk}\(\phi_\tau\(J_{t(\tau)}\)\)$ is contained in a connected component of 
$$
f^{lk}\circ \phi_\tau\(W_{t(\tau)}^*\)
$$ 
whose diameter is, by Lemma~\ref{l12018_08_06}, bounded above by 
$$
C\exp\(-\b_U\(s((x,k),(y,u))\).
$$
In conjunction with \eqref{320121208B} and \eqref{220121208B}, this
finishes the proof of Claim 1 by taking $\b:=e^{-\g\b_U}$.

\sp\fr{\bf Claim~2:} The function $\tilde g$ is not
cohomologous to a constant in $L^2(\nu)$.

\sp\fr Indeed, assume without loss of generality that $\mu_h(g)=0$. By
virtue of Lemma~\ref{l320121210} the fact that $g:J(f)\to \mathbb{R}$ is not 
a coboundary with respect to the map $f^l$ in $L^2(\mu_h)$ equivalently means that the sequence $\(S_n^{(l)}(g)\)_{n=1}^\infty$ is not uniformly bounded in
$L^2(\mu_h)$, where $S_n^{(l)}$ refers to the $n$th Birkhoff's sum with respect to the map $f^l$. But because of Lemma~\ref{l120110623B}, 
$$
\|S_n^{(l)}(\tilde g)\|_{L^2(\nu_h)} = \|S_n^{(l)}(g)\|_{L^2(\mu_h)}.
$$
So, the sequence $\(S_n^{(l)}(\tilde
g)\)_{n=0}^\infty$ is not uniformly bounded in $L^2(\nu_h)$. Thus, by
Lemma~\ref{l320121210} again, it is not a coboundary in $L^2(\nu)$. The proof of Claim~2 is finished.

\sp\fr Having these two claims, all items, (1), (2), and (3), now
follow immediately from Theorem~\ref{lsyoung} with the use of
Lemma~\ref{l120110623B} and formula \eqref{eq120110623}. The proof of Theorem~\ref{t220110627B} is finished. 
\epf

\sp As a fairly immediate consequence of Theorem~\ref{t220110627B}, we get the following. 

\begin{thm}\label{t220110627C}
If $f:\C\lra\oc$ is a parabolic elliptic function of finite class, then the dynamical system $\(f|_{J(f)},\mu_h\)$ satisfies the following. If $g:J(f)\lra \mathbb{R}$ is a bounded function, H\"older continuous with respect to the Euclidean metric on $J(f)$, then

\, \begin{enumerate}
\item{} (The Polynomial Decay of Correlations) For
every bounded measurable function $\psi:J(f)\lra\mathbb{R}$, we have that
$$
\bigg|\int\psi\circ f^n \cdot g d\mu_h
-\int g d\mu_h\int\psi d\mu_h\bigg|
\lek n^{2-\frac{p_{\max}+1}{p_{\max}}h}.
$$

\sp\item{} (The Central Limit Theorem) If in addition $g:J(f)\lra\mathbb{R}$ is not cohomologous to a constant in $L^2(\mu_h)$, i.e. if there is no square integrable function $\eta$ for which $g={\rm const}+\eta\circ f-\eta$, then the The Central Limit Theorem holds for $g$. More precisely,
there exists $\sigma>0$ such that
$$
\frac{\sum_{j=0}^{n-1}g\circ f^j-n\int g d\mu_h}{\sqrt{n}} \xrightarrow[\ \, n\to\infty \ ]{}  
 \mathcal N(0,\sigma)
$$
in distribution, where, as usually, $\mathcal N(0,\sigma)$ denotes the Gauss (normal) distribution centered at $0$ with covariance $\sg$.
\end{enumerate}
\end{thm}

\bpf
Of course each function $ \psi\circ  f^r: J(f) \lra \mathbb R$ , $ r >0$, is bounded, measurable, and $\|\psi\circ f^r\|_\infty \leq \|\psi\|_\infty$.
Let $C_r \in (0, +\infty)$ , $ r \geq 0$, be the constant witnessing  formula (1) in Theorem~\ref{t220110627B} for the  functions
$\psi\circ f^r$ and $g$. Let 
$$
C:=\max \{ C_r:\,  0 \leq r \leq l-1\}.
$$
Given  $ n \geq l$  write uniquely  $n=lk+r$, where  $r$ and $k$ are non--negative integers and $ r \in\{0,1,\ldots, l-1\}$. Then, by virtue of Theorem~\ref{t220110627B}, we get that
 $$\begin{aligned}
 \left| \int\psi\circ f^n g d\mu_h- \int gd\mu_h \int\psi d\mu_h\right|&=\left|  \int (\psi\circ f^r)\circ f^{lk} g d\mu_n- \int gd\mu_h \int\psi \circ f^r d\mu_h\right|\\
 & \leq C_r k^{2- \frac{p_{\max}+1}{p_{\max}}h} 
 \leq C_r\left(\frac{n}{2l}\right)^{2- \frac{p_{\max}+1}{p_{\max}}h} \\
 &\leq C'n^{2- \frac{p_{\max}+1}{p_{\max}}h}
 \end{aligned}$$
 for all $n \geq 0$ large enough, where  $ C'=C(2l)^{\frac{p_{\max}+1}{p_{\max}}-2}.$
Item (1)  is thus proved.

\sp Passing to item (2), assume without  loss of generality that $\mu_h(g)=0$. Now, given  the integer $n \geq 0$, write uniquely $n=lk+r, \, k=k_n \geq 0$,
$0\leq r \leq l-1$.
Then
\beq\label{1fpn22}
 S_ng=S^{(l)}_k(S_lg) + \sum_{j=n-r}^{n-1} g\circ f^j= S^{(l)}_k(S_lg)+S_r(g\circ f^{n-r}).
 \eeq
  Hence
 $$ \frac{S_ng}{\sqrt{n}}=\frac{S^{(l)}_k(S_lg)}{\sqrt{k}} \sqrt{\frac{k}{n}} +\frac{S_r(g\circ f^{n-r})}{\sqrt{n}},$$
so
 \beq\label{2fpn22}
 \left| \frac{S_ng}{\sqrt{n}} - \frac{1}{l}\frac{S^{(l)}_k(S_lg)}{\sqrt{k}}\right| \leq \left| \frac{S^{(l)}_k(S_lg)}{\sqrt{k}}\right| \left|  \sqrt{\frac{k}{n}} - \frac{1}{l}\right| + \left|\frac{S_r(g\circ f^{n-r})}{\sqrt{n}}\right|
\eeq
Likewise,
 $$ \frac{S_ng}{\sqrt{n}}\sqrt{\frac{n}{k}}=\frac{S^{(l)}_k(S_lg)}{\sqrt{k}}  +\frac{S_r(g\circ f^{n-r})}{\sqrt{k}}.$$
Therefore,
\beq\label{1fpn24}
\left|  \frac{S_ng}{\sqrt{n}}l- \frac{S^{(l)}_k(S_lg)}{\sqrt{k}}\right| \leq  \left|  \frac{S_ng}{\sqrt{n}}\right| \left|  \sqrt{\frac{n}{k}} - {l}\right| +
\left|\frac{S_r(g\circ f^{n-r})}{\sqrt{n}} \right|\sqrt{\frac{n}{k}}.
\eeq
Fix now a non--degenerate  open interval $(a,b)$  and $\varepsilon >0$.  Since the function $g$ is uniformly bounded, there exists $N \geq 1$  such that
\beq\label{3fpn24}
\left|\frac{S_r(g\circ f^{n-r})}{\sqrt{n}} \right|\sqrt{\frac{n}{k}} \leq \varepsilon, \, \, \, \left|  \sqrt{\frac{n}{k}} - {l}\right| \leq \varepsilon \, \,\, \text{and}\, \,\,
\left|  \sqrt{\frac{k}{n}} -\frac{1}{l}\right|\leq \varepsilon
\eeq
for every  $n \geq N$. So, if $ z \in J(f)$ and
$$
\frac{S_ng(z)}{\sqrt{n}}\in (a, b),
$$
then it follows  from (\ref{1fpn24}) and  (\ref{3fpn24}) that
$$  
\frac{S^{(l)}_{k_n}(S_lg(z))}{\sqrt{k_n}} 
\in \( la- (\max\{|a|, |b|\}+1)\varepsilon, lb+ (\max\{|a|, |b|\}+1)\varepsilon\).
$$
Therefore, applying Theorem~\ref{t220110627B} (2), which is applicable because of Corollary~\ref{c1fpn25}, we get that
$$
\begin{aligned} &\varlimsup_{ n \to \infty} \mu_h\lt( \{ z \in J(f): \frac{S_ng(z)}{n} \in (a, b) \}\rt)\le \\
& \leq \varliminf_{ n \to \infty} \mu_h\left( \left\{ z \in J(f):  \frac{S_k^{(l)}(S_lg(z))}{\sqrt{k}} \in \( la- (\max\{|a|, |b|\}+1)\varepsilon, lb+ (\max\{|a|, |b|\}+1)\varepsilon\)\right\}\right)
\\
& = \mathcal N (0, \sigma_l)\Big(\( la- (\max\{|a|, |b|\}+1)\varepsilon, lb+ (\max\{|a|, |b|\}+1)\varepsilon\)\Big),
\end{aligned}
$$
where $\sigma_l$ is the covariance  of Theorem~\ref{t220110627B} (2) ascribed to the dynamical system  $(f^l, \mu_h)$ and observable  $S_lg$. Since  the normal distribution has no atoms, by letting $\varepsilon\to 0$, it follows from  this formula that
\beq\label{2fpn24}
\varlimsup_{n \to \infty} \mu_h\lt(\Big\{ z \in J(f): \frac{S_ng(z)}{n} \in (a, b) \Big\}\rt) 
\leq \mathcal N (0, \sigma_l)\(( la, lb)\)
=\mathcal N (0,\sigma_l/l)\((a, b)\).
\eeq
Likewise if $  z \in J(f)$ and  
$$
\frac{1}{l}\frac{S_{k_n}^{(l)}(S_lg(z))}{\sqrt{k_n}} \in \( a + (\max\{|a|, |b|\}+1)\varepsilon, b- (\max\{|a|, |b|\}+1)\varepsilon\)
$$
 then it follows from  (\ref{2fpn22}) and (\ref{3fpn24}) that
$$
\frac{S_ng(z)}{\sqrt{n}} \in (a, b). 
$$ 
Therefore, by applying  Theorem~\ref{t220110627B} (2), we get that
$$
\begin{aligned}
&\varliminf_{ n \to \infty} \mu_h \lt(\Big\{ z \in J(f):  \frac{S_ng(z)}{n} \in (a, b) \Big\}\rt)\ge \\
&\geq \varliminf_{ n \to \infty} \mu_h\left( \left\{ z \in J(f): \,   \frac{S_k^{(l)}(S_lg)(z)}{\sqrt{k}} \in \( la + l(\max\{|a|, |b|\}+1)\varepsilon, lb-l(\max\{|a|, |b|\}+1)\varepsilon\) \right\} \right)\\
& = \mathcal N (0, \sigma_l)\Big(\( la+  l (\max\{|a|, |b|\}+1)\varepsilon, lb -l (\max\{|a|, |b|\}+1)\varepsilon\)\Big).
\end{aligned}
$$
Hence, letting $ \varepsilon \to 0$, we get that
$$
\varliminf_{ n \to \infty} \mu_h\lt( \Big\{ z \in J(f): \frac{S_ng(z)}{n} \in (a, b) \Big\}\rt) 
\geq \mathcal N (0, \sigma_l)\((la, lb)\)
= \mathcal N(0,\sigma_l/l)\((a, b)\).
$$
Along with (\ref{2fpn24}), this gives that,
$$
\lim_{ n \to \infty} \mu_h\lt( \Big\{ z \in J(f): \frac{S_ng(z)}{n} \in (a, b) \Big\}\rt)
=\mathcal N(0,\sigma_l/l)\((a, b)\),
$$
and the proof of Theorem~\ref{t220110627C} (2) is complete. Thus, the proof of entire Theorem~\ref{t220110627C} is complete.
\epf

\sp\section[Infinite Conformal Invariant Measures for Parabolic Functions] {Infinite Conformal Invariant Measures; Darling--Kac Theorem for Parabolic Elliptic Functions}~\label{Darling_Kac-elliptic}

In this section, similarly as in the two previous sections, we deal with parabolic elliptic functions\index{(N)}{parabolic elliptic function}  i.e. with all such elliptic functions $f:\mathbb {C} \lra\oc$ for which
$$
\Crit(f)\cap J(f)=\es \ \
{\rm and } 
\  \  \Om(f)\ne\es.
$$
We keep the notation of these two previous sections. Similarly as in the previous section, i.e. Section~\ref{parabolicfiniteinvariantmeasures} our main focus is on the invariant measure $\mu_h$. However now, we do not assume any longer that the measure $\mu_h$ is finite, indeed all results in this section concerning measure $\mu_h$ are non--trivial only if this measure is infinite.

Our ultimate goal in this section is to prove the Darling--Kac Theorem of Section~\ref{D-K AV} (Theorems~\ref{t1j163} and \ref{t1j165}) for the class of parabolic elliptic functions $f:\mathbb {C} \lra\oc$ with infinite invariant measure $\mu_h$; we call such functions of infinite class. \index{(N)}{parabolic elliptic functions of infinite class} We want to apply Theorem~\ref{t1j165}, and therefore, taking the
iterate $l \geq 1$ so large that all rationally indifferent periodic
points of $f$ become simple fixed points of $f^l$, all we need to do is to
find a Borel set $Y$  satisfying the hypotheses of
Theorem~\ref{t1j165}. Up to the very last theorem of this section, i.e. Darling--Kac Theorem for Parabolic Elliptic Functions of Infinite Class, Theorem~\ref{t2j203}, our considerations are valid regardless of whether the $f$--invariant measure $\mu_h$ is infinite or finite. 

To begin with, let
\beq\lab{2_2017_11_24}
Y:=J(f)\sms \left(\{\om\}\cup\bigcup_{\om \in
\Om}\bigcup_{n=1}^\infty f^{-ln}_\om (\De(\om))\right),
\eeq
where, we recall, given $\a\in(0,\pi)$, the sets 
$$
\De_j(\om,\a), \  \om\in\Om, 1\le j\le p(\om),
$$
are defined in Lemma~\ref{l1ch11.5}, 
$$
\De(\om)=\De(\om,\a)=\bu_{j=1}^{p(\om)}\De_j(\om,\a),
$$
and 
$$
\De=\bu_{\om\in\Om(f)}\De(\om).
$$
It follows from Theorem~\ref{t4050206} that
$$
0< \mu_h(Y) < +\infty.
$$
As in previous sections  let $m_{h,e}$ be the $h$--conformal  Euclidean
measure and let $m_{h,s}$ be the $h$--conformal spherical measure. We
know that $m_{h,s}$ is  finite and it follows from (\ref{1p6}) that
$$ 
\frac{dm_{h,s}}{dm_{h,e}}(z)= \frac{1}{(1+|z|)^{2h}}.
$$
Recall that $\rho_s:=\frac{d\mu_h}{dm_{h,s}}$\index{(S)}{$\rho_s$} is the
Radon--Nikodym derivative of the measure  $\mu_h$ with respect to
$m_{h,s}$ and $\rho_e:=\frac{d\mu_h}{dm_{h,e}}$\index{(S)}{$\rho_e$} is
the Radon--Nikodym derivative of the measure $\mu_h$ with respect to
$m_{h,e}$. Let
$$
{\mathcal L}_s\index{(S)}{${\mathcal L}_s$}:
C_b(J(f))\index{(S)}{$C_b(J(f))$}\lra C_b(J(f))
$$ 
and 
$$
{\mathcal L}_e\index{(S)}{${\mathcal L}_e$}: L_1(m_{h,e})\lra L_1(m_{h,e})
$$ 
be the transfer (Perron--Frobenius) operators \index{(N)}{transfer operator} \index{(N)}{Perron--Frobenius operator} operator ascribed respectively to the
quasi-invariant measures $m_{h,s}$ and $m_{h,e}$ as of \eqref{1_2017_11_21}. They are  given respectively by the formulas
$$ {\mathcal L}_s(g)(z)=\sum_{w\in
f^{-1}(z)}|f'(w)|^{-h}_sg(w) $$ and  $${\mathcal
L}_e(g)(z)=\sum_{w\in f^{-1}(z)}|f'(w)|^{-h}g(w).
$$
The Perron--Frobenius operator ${\mathcal
L}_0\index{(S)}{${\mathcal L}_0$}: C_b(J(f))\lra C_b(J(f))$ ascribed
to the measure $\mu_h$ is equal  to
$$ 
\begin{aligned}
{\mathcal L}_0(g)(z)
=& \rho^{-1}_s(z) \mathcal {L}_s( g\rho_s)(z)
=\rho^{-1}_s(z)\sum_{w \in f^{-1}(z)} g(w) \rho_s(w)| f'(w)|^{-h}_s
=\rho^{-1}_e (z) \mathcal {L}_e( g\rho_e)(z)\\
=&\rho^{-1}_e(z)\sum_{w \in f^{-1}(z)} g(w) \rho_e(w)|f'(w)|^{-h}.
\end{aligned}
$$
For every point $\om \in \Om$ and $\xi \in f^{-l}(\om) \sms \Om$,
let 
$$
f_\xi^{-l}:B(\om, \theta) \lra \mathbb C
$$ 
be the unique holomorphic  inverse branch  of $f^{-l}$, defined  on $B(\om,
\theta)$ with some  $\theta\in (0,\th(f)]$ sufficiently  small, and sending
$\om$ to  $\xi$. 

Now, fix $z \in \De(\om)\cap J(f)$ and $n_0 \geq 1$ so large that 
$$
f^{-ln}_\om(z)\in B(\om, \theta)
$$ 
for all $n \geq
n_0$. According to Proposition~\ref{p1j111DL} and Theorem~\ref{t1pj11}, we have that
\beq\label{1j173}
\begin{aligned} 
&\lim_{n\to + \infty}(n+1)^{\frac{p(\om)+1}{p(\om)}h}
\big|(f^{-l}_\xi \circ f^{-ln}_\om )'(z)\big|^h_s\frac{\rho_s(f^{-l}_\xi\circ
f^{-ln}_\om(z))}{\rho_s(z)}= \\
&=\lim_{n\to + \infty}n^{\frac{p(\om)+1}{p(\om)}h}\big|(f^{-ln}_\om)'(z)\big|^h
\big|(f^l)'(\xi)\big|_s^{-h}\frac{(1+|z|^2)^h}{(1+|\om|^2)^{h}}\frac{\rho_s(\xi)}{\rho_s(z)}\\
&=|a_\om|^{-\frac{p(\om)+1}{p(\om)}h}
\big|(f^l)'(\xi)\big|_s^{-h}(1+|\om|^2)^{-h}\rho_s(\xi)\big|\(f_\om^{-1}\)_\infty'(z)\big|\rho_s^{-1}(z)(1+|z|^2)^h|z-\om|^{-h(p(\om)+1)}
\end{aligned}
\eeq 
uniformly on $\De(\om)\cap J(f)$, where $\big|\(f_\om^{-1}\)_\infty'(z)\big|$ comes from Proposition~\ref{p1j111DL} and formula \eqref{4_2017_11_28}. Denote the product of the first four factors in the last line of \eqref{1j173} by $Q_\om(\xi)$ and the product of the last four factors by $G(z)$. With the sets $Y_n(f^l)$, defined in \eqref{5j164} with $T=f^l$, formula (\ref{1j173}) then implies that
\beq\label{2j173}
\lim_{n\to +\infty}n^{\frac{p(\om)+1}{p(\om)}h}{\mathcal
L}_0^{ln}\(\1_{Y_n(f^l)}\)(z)
=G(z)\sum_{\xi\in
f^{-l}(\om)\sms \Om}Q_\om(\xi)
\eeq
uniformly on $\De(\om)\cap J(f)$. Now  let
$$ 
p:= \min \{p(\om):\om \in \Om \}
$$
and let
$$
\begin{aligned}
D(z)
:=\left\{\begin{array}{ll}
0 & \mbox{ if $z \in J(f)\sms \De$}\\
G(z)\sum_{\xi \in f^{-l}(\om)\sms \Om }Q_\om(\xi) &
\mbox{ if $z \in \De(\om)\cap J(f)$} \mbox{ and $p(\om)=p$}\\
0 & \mbox{ if $z\in \De(\om) \cap J(f)$} \mbox{ and $p(\om) < p$}.
\end{array}\right. 
\end{aligned}
$$
Since, ${\mathcal L}^{ln}_0\(\1_{Y_n(f^l)}\)(z)=0$ for all $n \geq
2$ and all $z \in J(f)\sms\{\De\}$, it follows from (\ref{2j173}) that
\beq\label{1j175}
\lim_{n \to \infty} n^{\frac{p+1}{p}} {\mathcal
L}^{ln}_0\(\1_{Y_n(f^l)}\)(z)=D(z)
\eeq
uniformly on $J(f)$. Integrating this over $J(f)$, we thus
get
\begin{eqnarray}\label{1j172}
\begin{aligned} 
\lim_{n \to\infty}n^{\frac{p(\om)+1}{p(\om)}h}{\mu_Y}(\1_{Y_n(f^l)})
&=\frac{1}{\mu_h(Y)}\lim_{n \to\infty}n^{\frac{p+1}{p}h} \mu_h(Y_{f^l,n})
 \\
 &=\frac{1}{\mu_h(Y)}\lim_{n \to
 \infty}n^{\frac{p+1}{p}h}\int_{J(f)}
{\mathcal L}^{ln}_0 (\1_{Y_{ f^l,n}})(z) d \mu_h(z)\\
& =\frac{1}{\mu_h(Y)}\lim_{n \to
\infty}\int_{J(f)}n^{\frac{p+1}{p}h}
{\mathcal L}^{ln}_0 (\1_{Y_{ f^l,n}})(z) d \mu_h\\
&=\frac{1}{\mu_h(Y)}\int_{J(f)} D(z) d\mu_h(z).
\end{aligned}
\end{eqnarray}
Along with (\ref{1j175}) this yields that
\beq\lab{3_2017_11_28}
\lim_{n \to \infty}\frac{1}{\mu_h(Y_n(f^l))}{\mathcal L}^{ln}_0
(\1_{Y_n(f^l)})(z)=\frac{\mu_h(Y)D(z)}{\int_{J(f)} D(z) d
\mu_h(z)}
\eeq
uniformly on $J(f)$. So, Lemma~\ref{l1j168} applies with the function
\beq\lab{4_2017_11_28B}
\hat{H}\index{(S)}{$\hat{H}$}
:=\frac{\mu_h(Y)D}{\int_{J(f)} D(x)d \mu_h(x)}.
\eeq
We then have
\beq\label{1j185}
H:= \mu_h(Y)^{-1}\hat{H}= \frac{D}{\int_{J(f)}D(z)d\mu_h(z)}.
\eeq
Obviously both functions $\hat H$ and $H$ are supported on $Y$. In order to show that the function 
$$
H:J(f)\lra [0,+\infty)
$$ 
is uniformly sweeping on $Y$, we will need a more refined knowledge about the
behavior of the values of the Perron--Frobenius operator acting
on characteristic functions of sufficiently small neighborhoods of
poles with maximal order. Recall that
$$
q=q_f:=\max\{q_b:b\in f^{-1}(\infty)\}.
$$
We start with the following. 

\sp\blem\lab{l1j185} 
If $f:\C\lra\oc$ is a parabolic elliptic function, then, increasing perhaps $R_0>0$ of Lemma~\ref{l1j183}, there exists $\eta>0$ such that
$$
{\mathcal L}_0(\1_{B_b(R)})(z) \geq \eta
$$ 
whenever $b\in f^{-1}(\infty)$ with $q_b=q$, $R \geq R_0$, and $z\in J(f)\cap B_\infty^*(R_0)$. 
\elem

\bpf Keep $R_0>0$ so large as required in Lemma~\ref{l1j183} and also, invoking Theorem~\ref{t1pj11}, so that
$$
\frac12\le \frac{\rho_e(w)}{\rho_e(b)}\le 2
$$
whenever $b\in J(f)\cap f^{-1}(\infty)$ and $w\in J(f)\cap B_\infty^*(R_0)$. Take any $R\geq R_0$. By virtue of this inequality and Lemma~\ref{l1j183} we thus have for all $z\in J(f)\cap B_\infty^*(R_0)$, that
$$
\begin{aligned} 
{\mathcal L}_0 ( \1_{B_b(R)})(z) 
&=  \rho_e^{-1}(z) \sum_{w \in f^{-1}(z)\cap B_b(R)} \rho_e(w)|f'(w)|^{-h} \\
& \geq M^{-1}|z|^{\frac{q+1}{q}h} A^{-h} \sum_{w \in f^{-1}(z) \cap
    B_b(R)}|z|^{-\frac{q+1}{q}h}\rho_e(w)\\
& =(MA^h)^{-1}\sum_{w \in f^{-1}(z) \cap B_b(R)}\rho_e(w) \\
& \geq (2MA^h)^{-1} \rho_e(b).
\end{aligned}
$$
So, since, invoking Theorem~\ref{t1pj11} again, we see that 
$$
\eta:=(2MA^h)^{-1}\min\{\rho_e(b):b\in f^{-1}(\infty)\}
$$ 
is positive, the proof is complete. 
\epf

\sp We  shall prove the following.

\sp\blem\label{l1j191} 
If $f:\C\lra\oc$ is a parabolic elliptic function, then there exists $R_1\geq R_0$ so large that if $R\geq R_1$ and $b \in f^{-1}(\infty)\cap B_\infty(2R)$ with $q_b=q$, then
$$
{ \mathcal L}^n_0( \1_{B_b(R)})(z) \geq \eta ^{n-1} {\mathcal
L}_0(\1_{B_b(R)})(z)
$$
for all integers $n \geq 1$ and all $z \in J(f)$, where $\eta>0$  comes
from Lemma~\ref{l1j185}.  
\elem

\bpf Starting induction on $n\ge 1$, the lemma is obviously true for $n=1$. The next step is to prove it for $n=2$. Towards this end, select $R_1\geq R_0$  so large that if $ R\geq R_1$ and $ a \in f^{-1}(\infty)\cap B_\infty(2R)$, then $|w|>R$ whenever $w\in B_a(R)$. Applying Lemma~\ref{l1j185}, we can then estimate as follows with $b \in f^{-1}(\infty)\cap B_\infty(2R)$:
  $$\aligned
           {\mathcal L}^2_0(\1_{B_b(R)})(z)
           &= \sum_{w \in f^{-1}(z)} \rho^{-1}_e(z) |f'(w)|^{-h}\rho_e(w)
  {\mathcal L}_0(\1_{B_b(R)})(w)\\
  & \geq     \sum_{w \in f^{-1}(z)\cap B_b(R)} \rho^{-1}_e(z) |f'(w)|^{-h}\rho_e(w)
  {\mathcal L}_0(\1_{B_b(R)})(w)\\
   & \geq \eta   \sum_{w \in f^{-1}(z)\cap B_b(R)} \rho^{-1}_e(z) |f'(w)|^{-h}\rho_e(w)  \\
   &= \eta {\mathcal L}_0 (\1_{B_b(R)})(z).
   \endaligned  $$
So, our lemma  is also proved for $n=2$. Proceeding further  by
induction, suppose that the lemma is true for some $n\geq 2$ (or $1$, it does not matter). Using monotonicity of the operator  ${\mathcal L}_0$ we then get
$$
\begin{aligned} {\mathcal L}^{n+1}_0 ( \1_{B_b(R)})
& = {\mathcal L}_0( {\mathcal L}_0^n (\1_{B_b(R)}))
\geq{\mathcal L}_0( \eta^{n-1}{\mathcal L}_0(\1_{B_b(R)}))\\
&=\eta^{n-1}{\mathcal L}_0^2 (\1_{B_b(R)})\\
& \geq \eta^n {\mathcal L}_0 (\1_{B_b(R)}).
\end{aligned}
$$
The inductive proof is complete. \endpf

\sp\fr Now, we record the following.

\sp\blem\label{l2j191} 
If $f:\C\lra\oc$ is a parabolic elliptic function, then for every $n\geq 1$, we have that
$$
{\mathcal L}^n_0(\1_{B_b(R)})(z)\geq \eta^n
$$ 
for all $R\geq R_1$, all $b\in f^{-1}(\infty)\cap B_\infty(2R)$ with $q_b=q$, and all $z\in J(f)\cap B_\infty^*(R_0)$, where $\Delta >0$ comes from Lemma~\ref{l1j185}.
\elem

\bpf It follows from Lemma~\ref{l1j191} and Lemma~\ref{l1j185}  that
for all $n \geq 1$
\beq\label{3j193}
{\mathcal L}_0^n (\1 _{B_b(R)})(z) 
\geq \eta^{n-1}{\mathcal L}_0(\1_{B_b(R)})(z)
\geq \eta^n.
\eeq
This completes the proof.
\endpf

\sp The next lemma is this.

\sp\blem\label{l1_2017_11_25}
If $f:\C\lra\oc$ is a parabolic elliptic function, then for all $r>0$, all $R\geq R_1$ and all $b\in f^{-1}(\infty)\cap B_\infty(2R)$ with $q_b=q$, we have that
$$
\Ga(r):=\inf\big\{{\mathcal L}^2_0(\1_{B_b(R)})(z):z\in J(f)\sms B(\Om(f),r)\big\}>0.
$$
\elem

\bpf 
Because of Lemma~\ref{l2j191} it suffices to show that
$$
\inf\big\{{\mathcal L}^2_0(\1_{B_b(R)})(z):z\in (J(f)\cap \ov B(0,R_0))\sms B(\Om(f),r)\big\}>0.
$$
In order to do it, fix a Euclidean ball 
$$
B\sbt f(B_b(R)),
$$ 
centered at a point in $J(f)$, such that
$$
f(B)=\oc.
$$
So, if
$$
G:=B_b(R)\cap f^{-1}\(B\cap f^{-1}\((J(f)\cap \ov B(0,R_0))\sms B(\Om(f),r)\)\),
$$
then 
$$
f^2(G)=J(f)\cap \(\ov B(0,R_0))\sms B(\Om(f),r)\)
$$
and 
$$
(\ov G\cup \ov{f(G)})\cap (\Om(f)\cup f^{-1}(\infty))=\es.
$$
Therefore,
$$
\ka:=\sup\big\{\big|(f^2)'(w)\big|:w\in G\big\}<+\infty,
$$
and also, invoking Theorem~\ref{t1pj11}, there exists $M\ge 1$ such that
$$
\rho_e(z)\le M
$$
for all $z\in J(f)\cap \(\ov B(0,R_0))\sms B(\Om(f),r)\)$, and
$$
\rho_e(w)\ge 1/M
$$
for all $w\in G$. Hence, for every $z\in J(f)\cap \(\ov B(0,R_0))\sms B(\Om(f),r)\)$ the exists at least one $\xi_z\in G$ such that $f^2(\xi_z)=z$ and
$$
{\mathcal L}^2_0(\1_{B_b(R)})(z)
=\sum_{w \in f^{-2}(z)}\rho^{-1}_e(z) |(f^2)'(w)|^{-h}\rho_e(w)
\ge \frac{\rho_e(\xi_z)}{\rho_e(z)|(f^2)'(\xi_z)|^h}
\ge M^{-2}\ka^{-h}>0
$$
The proof is complete.
\endpf 

\sp As a fairly straightforward consequence of this lemma, we can prove
the following.

\sp\blem\label{l2_2017_11_25}
If $f:\C\lra\oc$ is a parabolic elliptic function, then for all $r>0$, all $R\geq R_1$, all $b\in f^{-1}(\infty)\cap B_\infty(2R)$ with $q_b=q$, and all integers $n\ge 1$, we have that
$$
\Ga_n(r):=\inf\big\{{\mathcal L}^{2n}_0(\1_{B_b(R)})(z):z\in J(f)\sms B(\Om(f),r)\big\}>0.
$$
\elem

\bpf We proceed by induction on $n\ge 1$. For $n=1$ this is just Lemma~\ref{l1_2017_11_25}. So Suppose that Lemma~\ref{l2_2017_11_25} is true for some $n\ge 1$. Fix $r>0$. There then exists $r'>0$ such that
$$
f^{-2}\(J(f)\sms B(\Om(f),r))\sbt J(f)\sms B(\Om(f),r')
$$
Using Lemma~\ref{l1_2017_11_25} and our inductive hypothesis, we then have for all $z\in J(f)\sms B(\Om(f),r))$ that
$$
\begin{aligned}
{\mathcal L}^{2{(n+1)}}_0(\1_{B_b(R)})(z)
&=\sum_{w\in f^{-2}(z)}\rho^{-1}_e(z)|(f^2)'(w)|^{-h}\rho_e(w){\mathcal L}^{2n}_0(\1_{B_b(R)})(w) \\
&\ge \Ga(r')\sum_{w\in f^{-2}(z)}\rho^{-1}_e(z)|(f^2)'(w)|^{-h}\rho_e(w) \\
&=\Ga(r'){\mathcal L}^{2n}_0(\1_{B_b(R)})(z) \\
&\ge\Ga(r')\Ga_n(r).
\end{aligned}
$$
The proof is complete.
\endpf 

\sp As a fairly straightforward consequence of this lemma, we can prove
the following.

\sp\blem\label{l1j195} Let $f:\C\lra\oc$ be a parabolic elliptic function. If $U \sbt J(f)$ is a non-empty bounded subset of $J(f)$ such that $\ov U\cap\Om(f)=\es$, then there exists an integer $s=s_U\geq 2$ such that  
$$
\inf\big\{{\mathcal L}^n_0 (\1_{U})(z): \, z \in Y\big\}>0,
$$ 
for all integers $n\ge s$, where, we recall, the set $Y\sbt J(f)$ is defined by \eqref{2_2017_11_24}. 
\elem

\bpf Fix $b \in f^{-1}(\infty)$ and $R \geq R_1$ as
required in Lemma~\ref{l1j191}. Since $J(f)=\ov{\bigcup_{n=1}^\infty f^{-n}(\infty)}$, there exists an even integer $k\geq 2$ such that 
\beq\lab{5_2017_11_25}
f^k(U)=J(f).
\eeq
Since $U$ is bounded and its closure is disjoint from $\Om(f)$, it follows from 
Theorem~\ref{t1pj11} that there exists a constant $Q\ge 0$ such that 
\beq\lab{3_2017_11_24}
\rho_e(\xi)\ge Q
\eeq
for all $\xi\in U$. Fix now an arbitrary point $a\in J(f)\cap U$ such that $f^k(a)=b$; in particular $\{a, f(a),\ld, f^k(a)\}\cup \(\{\infty\}=\es$. Then there exists an open ball $B(a,r)\sbt C$ such that $J(f)\cap B(a,r)\sbt U$
$$
\{\infty\}\cup\ov{\bu_{l=0}^kf^l(B(a,r))}=\es.
$$
Therefore, there exists a constant $T\in(0,+\infty)$ such that
\beq\label{2j195}
|(f^k)'(\xi)|\leq T
\eeq
for all $\xi\in B(a,r)$, and, Since $f^k((B(a,r))$ is an open neighborhood of $b$, there exists $R\ge R_1$ as above and so large that
$$
f^k(B(a,r))\spt B_b(R).
$$
Hence, for every $w\in B_b(R)$ the exists at least one $\xi_w\in B(a,r)$ such that $f^k(\xi_w)=w$. Using Lemma~\ref{l1j183} for $\rho^{-1}_e(w)$ below, formula \eqref{3_2017_11_24}, and (\ref{2j195}), and Lemma~\ref{l2_2017_11_25}, we can write for all $z\in Y$ and all even integers $j\ge 2$, that
$$\begin{aligned}
{\mathcal L}_0^{k+j} (\1_U)(z) 
&= \sum_{w \in f^{-j}(z)}\rho^{-1}_e(z) |(f^j)'(w)|^{-h}\rho_e(w)
 ({\mathcal L}_0^{k}\1_U)(w)\\
& \geq  \sum_{w \in f^{-j}(z)\cap B_b(R)}
   \rho^{-1}_e(z) |(f^j)'(w)|^{-h}\rho_e(w){\mathcal L}^k_0( \1_U)(w)\\
& \geq  \sum_{w \in f^{-j}(z)\cap B_b(R)} \rho^{-1}_e(z)
   |(f^j)'(w)|^{-h}\rho_e(w)\(\rho^{-1}_e(w)|(f^k)'(\xi_w)|^{-h}\rho_e(\xi_w)\)\\
& \geq M^{-1}|w|^{\frac{q+1}{q}h}QT^{-h}\sum_{w \in f^{-j}(z)\cap B_b(R)} \rho^{-1}_e(z)|(f^j)'(w)|^{-h}\rho_e(w)\\
& \ge QR^{\frac{q+1}{q}h}(MT^h)^{-1} {\mathcal L}^j_0(\1_{B_b(R)})(z)\\
& \geq  QR^{\frac{q+1}{q}h}(MT^h)^{-1}\Ga_{j/2}(r),
\end{aligned}
$$
where $r>0$ is such that $Y\sbt J(f)\sms B(\Om(f),r)$. So, we have proved our lemma for all even integers $n\ge 2$ large enough. If we take an odd integer $k\ge 1$ satisfying \eqref{5_2017_11_25}, then exactly the same reasoning will prove the lemma for all odd integers $n\ge 3$ large enough. We are done. \endpf 

\sp As a main consequence of this lemma, we  now are able to prove
easily the following.

\sp\bprop\label{p1j197} 
If $f:\C\lra\oc$ is a parabolic elliptic function and
the set $Y\sbt J(f)$ is defined by \eqref{2_2017_11_24}, then the function 
$$
H:J(f)\lra[0,+\infty),
$$
defined by the formula (\ref{1j185}), is uniformly sweeping on $Y$ for the map $f^l:J(f)\lra J(f)$ with respect to the invariant measure $\mu_h$. 
\eprop

\bpf Aiming to apply Lemma~\ref{l1j195}, take as $U\sbt J(f)$ any open bounded subset of $J(f)$ such that $\ov U\cap\Om(f)=\es$, and $U\sbt V$. Fix an integer $n \geq 1$ so large that $ln\ge s_U$, the number produced in Lemma~\ref{l1j195}. Since $I:=\inf\(H|_V\)>0$ and since ${\mathcal L}_0$ is a monotone operator, we get from Lemma~\ref{l1j195} for all $z\in Y$, that
$$
{\mathcal L}_0^{ln}(H)(z) 
\geq {\mathcal L}_0^{ln}(I \1_U)(z)
=I{\mathcal L}^{ln}_0 (\1_ U)(z)  
\geq I\inf \{{\mathcal L}^{ln}_0(\1_U)(w):w\in Y\}
>0.
$$ 
We are done.
\endpf

\sp Our goal now is to determine the wandering rates 
$$
w_n(Y):=w_n(f^l,Y),
$$
which, we recall, are defined by \eqref{1201407}, and  the set $Y\sbt J(f)$ is defined by \eqref{2_2017_11_24}. In order to do this fix $\a >0$. The standard elementary integral test gives that
\beq\label{1j1999}
\lim_{n \to  \infty} \left(n^{\a-1}\sum_{k=N}^n k^{-\a}\right)=
\frac{1}{1-\a}\quad  \mbox{if}\,\, \a <1 \, \, \mbox{for every} \,\,
N\geq 1,
\eeq

\beq\label{2j199}
\lim_{n \to  \infty} \left(n^{\a-1}\sum_{k=n}^\infty
k^{-\a}\right)=\frac{1}{\a-1}\quad  \mbox{if}\,\, \a >1,
\eeq
and
\beq\label{3j1999}
\lim_{n \to  \infty} \frac{1}{\log n}\sum_{k=N}^n k^{-1}
= \lim_{n\to \infty} \frac{1}{\log n}\sum_{k=n}^\infty k^{-1}=1.
\eeq
Let
$$ 
\ka:= \frac{p+1}{p}h>1, 
$$
the inequality holding because obviously $\frac{p+1}{p}>1$, while $h>1 $ because of Theorem~\ref{thm:julia}. We will calculate the asymptotic behavior of  the wandering rates
$(w_n(Y))_{n-1}^\infty$ in the following two lemmas.

\sp\blem\label{l1j199} 
If $f:\C\lra\oc$ is a parabolic elliptic function, and the set $Y\sbt J(f)$ is defined by \eqref{2_2017_11_24}, then there exists a constant $A_1\in0,+\infty)$ such that
$$
\lim_{n \to \infty} n^{\ka-2} n \mu_Y ( \tau_Y\geq n)=A_1,
$$
where, we recall, $\tau_Y:Y\lra\{1,2,\ld\}$ is the first return time to the set $Y$ under the action of the map $f^l$.
\elem

\bpf Put 
\beq\lab{1_2017_11_28}
C:=\mu_h(Y)^{-1}\int_Y D\, d\mu_h.
\eeq
Fix an arbitrary  $\varepsilon >0$. By virtue of  (\ref{1j175}) there exists $N_\varepsilon \geq 1$ such that
\beq\label{4j199}
(1+\varepsilon)^{-1}C 
\leq k^\ka\mu_Y (\tau_Y=k)\leq
(1+\varepsilon) C
\eeq
for all $k \geq N_\varepsilon$. Hence, for every $n\geq
N_\varepsilon$,
$$ (1+\varepsilon)^{-1}C \sum_{k=n}^\infty k^{-\ka} 
\leq  \sum_{k=n}^\infty \mu_Y (\tau_Y=k)
+ \mu_Y(\tau_Y\geq n)
\leq  (1+\varepsilon) C \sum_{k=n}^\infty k^{-\ka}.
$$ 
It therefore follows  from formulas (\ref{2j199}) that
$$ 
(1+\varepsilon)^{-1}\frac{C}{\ka-1} \liminf_{n\to \infty} n^{\ka-1}  \mu_Y (\tau_Y\geq n)
\leq \limsup_{n\to \infty}n^{\ka-1}  \mu_Y (\tau_Y\geq n)
\leq (1+\varepsilon) \frac{C}{\ka-1}.
$$ 
Letting $\varepsilon \searrow 0$, this  gives  that
$$
\limsup_{n\to \infty}n^{\ka-2}n  \mu_Y (\tau_Y\geq n)
=\frac{C}{\ka-1}$$ and the proof is complete. 
\endpf

\sp\blem\label{l1j201} If $f:\C\lra\oc$ is a parabolic elliptic function, and the set $Y\sbt J(f)$ is defined by \eqref{2_2017_11_24}, then there exists a constant $A_2\in (0,+\infty)$  such that
$$ \lim_{n\to \infty} n^{\ka-2} \sum_{k=1}^{n-1}
 k \mu_Y
 (\tau_Y=k)=A_2$$
if $\ka <2$  and
$$ \lim_{n\to \infty}\frac{1}{ \log n}  \sum_{k=1}^{n-1}
 k \mu_Y(\tau_Y=k)=A_2$$
if $ \ka=2$. \elem

\bpf Fix $\varepsilon >0$. With the constant $C$ defined in \eqref{1_2017_11_28} and $N_\varepsilon$ as introduced in the proof of Lemma~\ref{l1j199}, by virtue of (\ref{4j199}), for every $k\geq N_\varepsilon$ we  have
$$ 
C(1+\varepsilon)^{-1}k^{1-\ka}
\leq k \mu_Y( \tau_Y=k)\leq C(1+\varepsilon)k^{1-\ka}.
$$ 
Put 
$$
S_\varepsilon:=\sum_{k=1}^{n-1}k \mu_Y(\tau_Y=k).
$$
Then we have for all $n\geq N_\varepsilon +1$ that
$$
C(1+\varepsilon)^{-1}\sum_{k=N_\varepsilon}^{n-1}k^{1-\ka}\leq \sum_{k=1}^{n-1} k \mu_Y( \tau_Y=k)
\leq S_\varepsilon+C(1+\varepsilon) \sum_{k=N_\varepsilon}^{n-1}k^{1-\ka}.
$$
It therefore follows from
(\ref{1j1999}) and (\ref{3j1999}) that
$$ \frac{C}{2-\ka}(1+\varepsilon)^{-1}\liminf_{n \to \infty}  n^{\ka-2} \sum_{k=1}^{n-1} k \mu_Y(\tau_Y=k) \leq
\limsup_{n \to \infty}
 k \mu_Y(\tau_Y=k)\leq \frac{C}{2-\ka}(1+\varepsilon)
$$ 
if $\ka <2$ and
$$ 
(1+\varepsilon)^{-1}\liminf_{n \to \infty} \frac{1}{\log n} \sum_{k=1}^{n-1} k \mu_Y(\tau_Y=k) 
\leq\limsup_{n \to \infty} \frac{1}{\log n} \sum_{k=1}^{n-1}k\mu_Y(\tau_Y=k)
\leq C(1+\varepsilon)
$$      
if $\ka=2$. Letting $\varepsilon \searrow 0$ the derived result follows. \endpf

\

\fr Since, by Lemma~\ref{l1j160}
 $$w_n( f^q, Y)= \mu(Y) \left(\sum_{k=1}^{n-1} k \mu_Y( \tau_Y=k) + n
 \mu_Y( \tau_Y \geq n)\right),$$
as a direct consequence  of Lemma~\ref{l1j199} and
Lemma~\ref{l1j201} , we obtain  the following.

\

\bthm\label{t1j203} If $f:\mathbb C \lra  \ov{\mathbb C}$ is a parabolic elliptic function, the set $Y\sbt J(f)$ is defined by \eqref{2_2017_11_24}, and $\ka = \frac{p+1}{p}h\leq 2$  (which precisely means that the $f$--invariant measure $\mu_h$ is infinite), then
$$ 
\lim_{n \to \infty}\frac{w_n(f^l, Y)}{n^{2-\b}}\in (0, \infty)
\quad \mbox{if } \quad \ka <2
$$ and
$$ \lim_{n \to \infty}\frac{w_n(f^l, Y)}{\log n}\in (0, \infty)
\quad \mbox{if } \quad \ka =2.$$ In particular, the sequence
$(w_n(f^l, Y))^\infty_1 \in {\mathcal R}_{2-\ka}.$ 
\ethm

\sp\fr Now we can take the fruits. As  a direct consequence of Theorem~\ref{t1j165}, Lemma~\ref{l1j168}, formula \eqref{3_2017_11_28} (see also \eqref{4_2017_11_28B}), Proposition~\ref{p1j197}, and Theorem~\ref{t1j203}, we get the following.

\sp\bthm[Darling--Kac Theorem for Parabolic Elliptic Functions of infinite class]\label{t2j203} 
If $f:\mathbb C \lra \oc$ is a parabolic elliptic function of infinite class, meaning that the $f$--invariant measure $\mu_h$ is infinite, \index{(N)}{parabolic elliptic functions of infinite class} then for every function $ g \in L^1( \mu_h)$ with $\int
g d\mu_h \neq 0$, the sequence
$$ 
\frac{1}{b_n}\sum_{j=0}^{n-1} g \circ f^j
$$
converges strongly, with respect to  the measure $\mu_h$, to the Mittag--Leffler  distribution 
$$
\lt(\int g d\mu\rt) {\mathcal M}_{\b-1},
$$
where the numbers $b_n$, $n\ge 1$, defined by \eqref{2_2017_11_28}, satisfy
$$ 
\lim_{n \to \infty} \frac{b_n}{n^{\ka-1}}
= \frac{ l^{2-\ka}}{\Gamma(\b) \Gamma( 3-\b)},
$$
and 
$$
\ka:=\frac{p+1}{p}h\leq 2
$$ 
which, by virtue of Theorem~\ref{t2_2017_09_21}  means that the $f$--invariant measure $\mu_h$ is infinite.  
\ethm

As we mentioned, and discussed at some length, in Section~\ref{D-K AV}, there are many recent deep results on the stochastic laws coming out of measure--preserving transformation with $\sg$--finite infinite measure. These are very sophisticated and exciting direction of research but because our book is already quite long, we will not discuss them here. 

\chapter[Dynamical rigidity of c.n.r.e. functions]{Dynamical Rigidity of Compactly Non--Recurrent Regular Elliptic Functions} 

\sp In this chapter we deal with the problem of dynamical rigidity of compactly non--recurrent regular elliptic functions.
The issue at stake is whether a weak conjugacy such as Lipschitz one on Julia sets can be promoted to much better one such as affine on the whole complex plane $\C$. This topic has a long history and goes back at least to the seminal paper \cite{Su3} of Dennis Sullivan treating among others the dynamical rigidity of conformal expanding repellers. Its systematical account can be found in \cite{PU2}. Many various, in many contexts, smooth and conformal dynamical rigidity theorems then followed. The literature abounds. 

Our approach in this chapter stems from the original Sullivan's article\cite{Su3}. It is also influenced by \cite{PU1}, where the case of tame rational functions has been actually done, and \cite{SU}, where the equivalence of the statements (4) and (1) of Theorem~\ref{t1042706} was established for all tame transcendental meromorphic functions. Being tame meaning that the closure of the postsingular set does not contain the whole Julia set; in particular each non--recurrent elliptic function is tame. We would however like to emphasize that, unlike \cite{SU}, we chose in our book the approach which does not make use of dynamical rigidity results for conformal iterated function
systems proven in \cite {MPU}. Such approach is therefore on the one hand somewhat more economical in tools but on the other hand uses the existence of conformal measures $m_h$ and their invariant versions $\mu_h$ which are in general not available in the general context of tame transcendental meromorphic functions dealt with in \cite{SU}. Our approach in particular  provides us with one more way of proving dynamical rigidity, and it can be possibly used in some other contexts. We would also like to add that doing dynamical rigidity for transcendental functions, like in \cite{SU} or in the current chapter, is substantially more involved than in the case of rational functions of \cite{PU1}; the first much bigger difficulty being caused by infinite degree of transcendental functions.

\sp Since in this chapter we simultaneously deal with two different compactly 
non--recurrent regular elliptic functions and our considerations heavily depend on Theorem~\ref{tmaincm} and Theorem~\ref{tinv}, and the conformal measures $m_h$ and their invariant versions $\mu_h$ introduced therein, in order to clearly distinguish between the measures for different maps, we now denote the measures $m_h$ and $\mu_h$, ascribed by these theorems, to a compactly non--recurrent regular elliptic function $f:\C\lra\oc$ respectively by $m_f$ and $\mu_f$.

Our ultimate goal  in this chapter
is to prove the following rigidity theorem. 

\sp\bthm\lab{t1042706} 
Suppose that $f:\C\lra\oc$ and $g:\C\lra\oc$ are two  compactly non--recurrent regular elliptic functions. Let $h:J(f)\lra J(g)$ be a homeomorphism conjugating $f$ to $g$ on their respective Julia sets, namely 
$$
h\circ f = g\circ h
$$
on $J(f)$. Then the following conditions (1)--(6) are equivalent.

\, \begin{enumerate}
\item[(1)] $h$ extends to an affine  conjugacy  from $\oc$  to
$\oc$ between $f:{\mathbb C} \lra\oc$ and
$g:{\mathbb C} \lra\oc$.

\,\item[(2)] $h$ extends to a conformal homeomorphism conjugating $f$
and $g$ on some respective neighborhoods of $J(f)$ and $J(g)$ in ${\mathbb C}$.

\,\item[(3)] $h$ extends to a real--analytic diffeomorphism
conjugating $f$ and $g$ on some respective neighborhoods of $J(f)$ and $J(g)$ in ${\mathbb C}$.

\,\item[(4)] The homeomorphisms  $h:J(f)\lra  J(g)$  and
$h^{-1}: J(g)\lra J(f)$ are Lipschitz continuous.

\,\item[(5)] For every periodic point $z\in J(f)$ of $f$, say of period $p$,
$|(f^p)'(z)|=|(g^p)'(h(z))|$.

\,\item[(6)] The measure class of $m_f$
is transported under $h$ to the measure class of $m_g$.
\end{enumerate}
\ethm 

\section[No Compactly Non-Recurrent Regular Function is Esentially Linear] {No Compactly Non-Recurrent Regular Elliptic Function is Esentially Linear}

\fr Let $f:\C\lra\oc$ be a compactly non--recurrent regular elliptic function. Let
$$
D_{\mu_f} \index{(S)}{$D_\mu$}:=\frac{d\mu_f\circ f}{d\mu_f}
$$ 
be the Jacobian of the elliptinc function $f:\C\to\oc$ with respect to the $f$-invariant measure $\mu_h$. Since 
$$
\frac{d\mu_f\circ f}{d\mu_f}
=\lt(\frac{d\mu_f}{dm_{f,e}}\circ f\rt)|f'|^h\left(\frac{d\mu_f}{dm_{f,e}}\right)^{-1},
$$ 
as an immediate consequence of Theorem~\ref{t1pj11}, we get  the following.

\sp\bcor\label{c2pj11} 
If $f:\C\lra\oc$ is a compactly non--recurrent regular elliptic function, then the Jacobian  $D_{\mu_f} f=\frac{d\mu_f\circ
f}{d\mu_f}$ has a real--analytic extension on a neighborhood
of $J(f)\sms (\ov{{\rm PC}(f)}\cup f^{-1}(\infty))$ in $\mathbb C$.
\ecor

\sp\fr We recall that  the  sets  ${\rm PC}(f)$ and ${\rm PC}^0(f)$ were
defined in (\ref{1042706}). Let 
$${\rm PS}^0(f)\index{(S)}{${\rm
PS}^0(f)$}:={\rm PC}^0(f)\cup f^{-1}(\infty),\quad  {\rm
PS}_{-}(f)\index{(S)}{${\rm PS}_{-}(f)$}:=\bu_{n=0}^\infty
f^{-n}\({\rm PS}^0(f)\).
$$

If $Y$  is a  subset  of $\mathbb C$, then we say that
\index{(S)}{$u$} $u:Y\lra S^1=\bd\D=\{z\in\C:|z|=1\}$ is an $f$--invariant line
field\index{(N)}{invariant line fields} on $Y$ if
$$
u(f(x))=\left(\frac{f'(x)}{|f'(x)|}\right)^2u(x)
$$ 
for all $x\in Y\cap f^{-1}(Y)$.

\sp\bthm\label{t3042706} If $f:\C\lra\oc$ is a compactly
non--recurrent regular elliptic function, then no of the following statements is true.

\begin{enumerate}
\item[(a)] The Jacobian $D_{\mu_f}: J(f)\sms \ov{{\rm PS}^0(f)}\lra
(0,+\infty)$ is locally constant.

\, \item[(b)] The  function
$\log|f'|:J(f)\sms \ov{{\rm PS}^0(f)}\lra {\mathbb R}$ is
cohomologous\index{(N)}{cohomologous function} to a locally constant
function on $J(f)\sms \ov{{\rm PS}^0(f)}$ in the class of continuous
functions on $J(f)\sms \ov{{\rm PS}^0(f)}$.

\, \item[(c)] There exists a  continuous $f$--invariant line field on
$J(f)\sms \ov{{\rm PS}^0(f)}$.

\,

\,

\item[(d)] For  every $ n \geq 1$
and  every point $z\in J(f)\sms {\rm PS}_{-}(f)$,
$$
{\rm det}\(\nabla(D_{\mu_f}\circ f^n)(z),\nabla(D_{\mu_f}\)(z)=0.
$$
\end{enumerate}
\ethm

\bpf The structure of the proof is to establish
the following  implications (a) $\Longrightarrow$ (b), (d)
$\Longrightarrow$(a)$\vee$(c) and to show  that each item (b) and (c) leads to a contradiction.

\sp {\bf (a)$\bf \Longrightarrow$(b)}. As before, let
$$
\rho_f\index{(S)}{$\rho$}:=\frac{d\mu_f}{dm_{f,e}}.
$$
Since
$D_{\mu_f}=\rho_f\circ f|f'|^h\rho_f^{-1}$, we  get that  
$$
\log D_{\mu_f} =\log \rho_f \circ f+ h \log |f'|- \log \rho_f.
$$ 
Hence
$$
\log |f'|= h^{-1}\log D_{\mu_f}+ h^{-1}\log \rho_f-h^{-1}\log \rho_f \circ f.
$$
Since the function $\rho_f: J(f)\sms \ov{{\rm PS}^0(f)}\lra (0, +\infty)$ is, by
Theorem~\ref{t1pj11}, continuous, and since $\log D_{\mu_f}$  is
locally constant by (a), the item (b) follows. 

\sp {\bf (d) $\bf \Longrightarrow$ (a) $\vee$ (c)}. Suppose first
that $\nabla D_{\mu_f}\equiv 0$ on $J(f)\sms \ov{{\rm  PS}^0(f)}$.
This of course precisely means that $ D_{\mu_f}$ is locally
constant on $J(f)\sms \ov{{\rm  PS}^0(f)}$, yielding (a). So, suppose
that there exists $ v\in J(f)\sms \ov{{\rm PS}^0(f)}$ such that
$$
\nabla D_{\mu_f}(v)\neq 0.
$$
Since the gradient $\nabla D_{\mu_f}$ is locally real--analytic, there  thus exists a simply-connected neighborhood $W \sbt {\mathbb C} \sms \ov{{\rm PS}^0(f)}$ of $v\in{\mathbb C}$ on which the gradient nowhere vanishes. Then there
exists a continuous function 
$$
l:W \lra S^1
$$  
such that $l(z)$ is orthogonal to $\nabla D_{\mu_f}(z)$ at every point $z\in W$. Now, for every $z\in J(f)\sms \ov{{\rm PS}^0(f)}$ there exists $n \geq 0$
and $\xi \in W\cap f^{-n}(z)$. Then define:
$$
l(z):=(f^n)'(\xi)l(\xi).
$$
We  want to show first  that  the function  $l$  is well defined
on $J(f)\sms \ov{{\rm PS}^0(f)}$ i.e. that  if $\zeta \in
f^{-m}(z)\cap W$, $m \geq 0$, then
\beq\label{2071406}
(f^n)'(\xi)l(\xi)=(f^m)'(\zeta)l(\zeta).
\eeq
Suppose on the contrary that (\ref{2071406}) fails  with  some $z,
\xi, \zeta$ as above. Then there exists  a point $ x\in (J(f)  \sms
\ov{{\rm  PS}^0(f)})\cap W$, $ k \geq 0$, and $ w\in f^{-k}(x)$ so
close to $z$  that  there  are two
 points  $y_1\in f^{-n}(w)$ and $y_2\in f^{-m}(w)$  so close  respectively to $\xi$ and $\zeta$ that
$$(f^n)'(y_1)l(y_1)\neq (f^m)'(y_2)l(y_2).$$
Hence,
$$(f^{k+n})'(y_1)l(y_1)\neq (f^{k+m})'(y_2)l(y_2).$$
So,  either
\beq\label{3071406}
 (f^{k+n})'(y_1)l(y_1)\neq l(x)\quad \text{or}\quad  (f^{k+m})'(y_2)l(y_2)\neq l(x).
\eeq
Suppose without loss of generality  that the first  of these two
inequalities  holds. Consider now gradients  as horizontal vectors
and vectors  parallel  to $l$ as vertical ones. The  standard  inner
product   becomes then the  product of matrices (horizontal or
vertical). Let  $t$ be a unit vector  parallel  to  $l(x)$. Since, by the Chain Rule,
$$
\nabla (D_{\mu_f}\circ f^{k+n})(y_1)=\nabla D_{\mu_f}(f^{k+n}(y_1))(f^{k+n})'(y_1)
=\nabla D_{\mu_f}(x)(f^{k+n})'(y_1) 
$$
and since the matrix $((f^{k+n})'(y_1))^{-1}$ is symmetric, we get
$$
\begin{aligned}
<\nabla (D_{\mu_f}\circ  f^{k+n})(y_1),((f^{k+n})'(y_1))^{-1}t>
&= <\nabla (D_{\mu_f} \circ f^{k+n})(y_1)((f^{k+n})'(y_1))^{-1},t> \\
&=<\nabla D_{\mu_f}(x),t> =0.
\end{aligned}
$$
Combining  this  and (\ref{3071406})  we see that   $l(y_1)$ is  not
perpendicular  to $$\nabla (D_{\mu_f}\circ f^{k+n})(y_1).$$ This
means that $\nabla D_{\mu_f}(y_1)$ and $\nabla (D_{\mu_f}\circ
f^{k+n})(y_1)$ are not parallel, or equivalently,
$$
\det (\nabla ( D_{\mu_f}\circ f^{k+n})(y_1),\nabla D_{\mu_f}(y_1))\neq 0,
$$ 
contrary to (d) since $y_1\notin  {\rm PS}_{-}(f)$. Thus the function $l:J(f)\sms  \ov{{\rm  PS}^0(f)}\to S^1$ is well defined, and the invariance of the  line field
$$
u(z):=l^2(z), \  \  z\in J(f)\sms \ov{{\rm  PS}^0(f)}
$$
is immediate from the definition of $l$. The implication (d) $\Longrightarrow$(a)$ \vee$ (c)  is thus
established.

\sp {\bf Item (c) leads to a contradiction.} By item (c) there
exists a continuous function $l: J(f) \sms \ov{{\rm PS}^0(f)} \to
S^1$ such that
\beq\label{1110806}
l(f(z)) =  l(z)\left( \frac{f'(z)}{|f'(z)|}  \right)^2
\eeq
for all $z\in J(f) \sms (\ov{{\rm PS}^0(f)}\cup f^{-1} (\ov{{\rm
PS}^0(f)}))$. Fix a pole $b \in f^{-1}(\infty)$. Let $q \geq 1$ be
the order  of $b$.  Taking $R >0$  sufficiently small, there
exists a holomorphic function 
$$
A:B(b, r) \lra \C\sms\{0\}
$$ 
such that
$$
f(z)=A(z)(z-b)^{-q}
$$ for all $z\in B(b, R)$. Consequently
$$ 
f(z) =A(z-w)(z-w-b)^{-q}
$$
for all $w\in \La_f$ and for all $z\in B(b+w,R)$. So,  $$f'(z)=(z-w -
b)^{-q-1}(A'(z-w)( z-w-b)-qA(z)).$$ Therefore,
\beq\label{2110806}
\frac{f'(z)}{|f'(z)|}=\frac{|z-w-b|^{q+1}} {(z-w-b)^{q+1}}
\frac{A'(z-w)(z-w-b)-qA(z-w)}{|A'(z-w)(z-w-b)-qA(z-w)|}.
\eeq
Now note that for every $w\in \La_f$ with sufficiently large
modulus, there exists $z_w\in B(w+b,R)$ being a fixed point of $f$,
i.e.
\beq\label{4110806}
f(z_w)=z_w.
\eeq
Note  also that $\lim_{w\to \infty} |z_w-(b+w)|=0$. Therefore,
\beq\label{3110806}
\lim_{w\to\infty}\frac{A'(z_w-w)(z_w-w-b)-qA(z_w-w)}{|A'(z_w-w)(z_w-w-b)-qA(z_w-w)|}
=-\frac{A(b)}{|A(b)|}.
\eeq
Since the elliptic function $f$ is compactly non--recurrent, if $|w|$ is large enough,  then 
$$
z_w \in J(f)  \sms (\ov{{\rm PS}^0(f)}\cup f^{-1} (\ov{{\rm PS}^0(f)})),
$$ 
and it therefore follows   from (\ref{1110806}) that
$$
\left(f'(z_w)/|f'(z_w)|\right)^2=1.
$$ 
Hence, combining (\ref{2110806}) and  (\ref{3110806}), we get that
$$
\lim_{w\to\infty}\left(\frac{(z_w-w-b)^{q+1}} {|z_w-w-b|^{q+1}}\right)^2
=\left(\frac{A(b)}{|A(b)|}\right)^2.
$$
Thus,
\beq\label{5110806}
\lim_{w\to \infty}\frac{(z_w-w-b)^{2(q+1)}} {|z_w-w-b|^{2(q+1)}}=
\left(\frac{A(b)}{|A(b)|}\right)^{2}.
\eeq
By (\ref{4110806}), $A(z_w-w)(z_w-w-b)^{-q}=z_w$, or equivalently,
$(z_w-b-w)^q=A(z_w-w)z_w^{-1}$. Hence,
$$
\frac{(z_w-w-b)^{2q(q+1)}} {|z_w-w-b|^{2q(q+1)}}=
\left(\frac{A(z_w-w)}{|A(z_w-w)|}\right)^{2(q+1)}\left(\frac{|z_w|}{z_w}\right)^{2(q+1)}.
$$
Since $\lim_{w\to \infty}(z_w-w)=b$,  inserting this to
(\ref{5110806}), we  get  that 
$$
\lim_{w\to\infty}\left(\frac{z_w}{|z_w|}\right)^{2(q+1)}
=\left(\frac{A(b)}{|A(b)|}\right)^2.
$$ 
This is a contradiction  since
the  set of accumulation   points  of the sequence
$\left(\frac{z_w}{|z_w|}\right)_{w\in \La}$ is the entire  unit
circle $S^1$. We are done.

\sp{\bf Item (b) leads to a contradiction}. Let 
$$
N:=\#f(\Crit(f)).
$$ 
Since, by Theorem~\ref{thm:julia}, $\HD(J(f))>1$, in fact the Hausdorff dimension
of every non-empty open subset  of $J(f)$ is equal  to
$\HD(J(f))>1$, we conclude  that there  are closed polygonal arcs (homeomorphic to the unit interval $[0,1]$ in $\R$)
$$
\g_1, \hat{\g}_1, \ldots, \g_{N},  \hat{\g}_{N}, { \g}_{N+1},
\hat{\g}_{N+1}
$$ 
consisting of finitely many straight line segments in $\oc$ with the following properties: 
\begin{enumerate}
\item [(i)] There exists  $x\in {\mathbb C}\sms \bu_{n\geq 0}^\infty
f^{-n}(\ov{{\rm PS}^0(f)})$ such that for all different indeces $i, j\in
\{1, \ldots, N+1\}$, we have that 
$$
\{x\}
\sbt\g_i\cap\g_j
=\hat\g_i\cap\hat\g_j
=\hat\g_i\cap\g_j
=\hat\g_i\cap\hat\g_i
\sbt\{x,\infty\}.
$$

\, \item [(ii)] The polygonal arcs $\g_1, \hat{\g}_1, \ldots,\g_{N}, \hat{\g}_{N}$ are all contained in $\C$ and compact.

\, \item [(iii)] \, The arcs ${\g}_{N+1}$ and $\hat{\g}_{N+1}$ have as one of their endpoints $\infty$, as the other $x$, and 
$$
{ \g}_{N+1} \cap\hat{\g}_{N+1}=\{x,\infty\}.
$$

\item [(iv)] For every $1 \leq j \leq N$ the arcs $\g_j$ and $\hat{\g}_j$ have the same  endpoints, one of which is $x$, the other belongs to $f(\Crit(f))$, and the intersection   $\g_j\cap \hat{\g}_j$ is a doubleton. In particular, $\g_j \cup \hat{\g}_j$ is a  closed  Jordan curve.
\, \item [(v)] If $Q=\bigcup_{j=1}^{N+1} \g_j$ and  $\hat{Q}=\bigcup_{j=1}^{N+1}\hat{\g}_j$,
 then each connected  component  of ${\mathbb C} \sms (Q\cup \hat{Q})$ intersects $J(f)$.
 
\, \item [(vi)] $f(\Crit(f))$ is   contained in  the set of  endpoints  of $\g_j,  j=1, \ldots, N$ (and  also
$\hat{\g}_j,  j=1, \ldots, N$).
\end{enumerate}
As a matter of fact we now fix the arcs $\g_1, \g_2,\ld\g_N, \g_{N+1}$ as above and we specify what the arcs $\hat\g_1, \hat\g_2,\ld\hat\g_N, \hat\g_{N+1}$ are in Claim~4 below. These will be small perturbations of the arcs $\g_1, \g_2,\ld\g_N, \g_{N+1}$.

Take a repelling fixed (with respect to $f$) point $w\in J(f)\sms \ov{{\rm
PS}^0(f)}$   (note that there are infinitely many  of such points).
Fix  also an arbitrary   point $$\xi_0\in J(f)\sms (Q\cup \ov{{\rm
PS}^0(f)})$$ and an arbitrary  radius $R>0$  so large, say  $R \geq
R^*$, that $B(\xi_0,R)\sms Q$ and $B(\xi_0,R)\sms \hat{Q}$   are
open topological disks. Then, for every $\xi_1 \in f^{-1}(\xi_0)$  there
exists  a unique holomorphic inverse branch 
$$
f^{-1}_{\xi_1}:B(\xi_0,R)\sms Q \lra {\mathbb C}
$$ 
of  $f$ sending  $\xi_0$ to $\xi_{-1}$. Note that  all the sets $f^{-1}_{\xi_1}(B(\xi_0,R)\sms Q)$ are uniformly bounded, have piecewise  smooth boundaries, and
$$
\dist\( f^{-1}_{\xi_1}(B(\xi_0,R)\sms Q), f^{-1}(\infty)\)>0.
$$
Recalling  also that ${\rm PC}_c(f)$ is bounded  and  ${\rm
PC}_p(f)$ is bounded (even finite), we  thus deduce  that for all
$\xi_{1}\in f^{-1}(\xi_0)$ with sufficiently  large   modulus
(depending on $R$),
\beq\lab{1082406}
\ov{f^{-1}_{\xi_1}(B(\xi_0,R)\sms Q)}\cap \ov{{\rm PS}^0(f)}=\es.
\eeq
Denote by $V$ the connected component of the set  
$$
\oc\sms \ov{f^{-1}_{\xi_1}(B(\xi_0,R)\sms Q)}
$$ 
that contains $\infty$. Obviously $V$ is an open simply connected  set whose  boundary  is a piecewise smooth Jordan curve contained in $
\partial{(f^{-1}_{\xi_1}(B(\xi_0,R)\sms Q))}$. Fixing  $\xi_1$  with
sufficiently  large  modulus, we will also have, $$\ov{{\mathbb C}
\sms V}\cap \ov{{ \rm PC}(f)}=\es.$$ Then there exists $r>0$ such
that, for every $s\in (0,r]$, 
$$
W_s:=B({\mathbb C} \sms V, s) 
\supset B(f^{-1}_{\xi_1}(B(\xi_0,R)\sms Q), s)
$$ 
is a topological disk disjoint from ${\rm PC}(f)$. Extend now $\xi_0$ and $\xi_1$ to a sequence $\xi:=\{\xi_n\}_0^\infty$ such that $f(\xi_{n+1})=\xi_n$
for all $n \geq 0$ and
\beq\lab{2082406}
\lim_{n\to \infty}\xi_n=w.
\eeq
For every $n\geq 2$, let $f^{-(n-1)}_n:W_r \lra {\mathbb C}$, be the
unique holomorphic inverse branch of $f^{n-1}$ sending $\xi_1$ to
$\xi_n$. For every $z\in W_r$ and every $ n \geq 2$, put 
$$
z_n:= f^{-(n-1)}_n(z).
$$
We shall show that the series
\beq\lab{2071206}
z\longmapsto\sum_{n=2}^{\infty}\(\log|f'(\xi_n)|-\log|f'(z_n)|\)
\eeq
of harmonic functions converges uniformly on $W_{r/2}$. Indeed, take arbitrary  $
1\leq k \leq l$. Then by the Chain Rule
\beq\lab{2_2017_11_30}
\left|\sum_{j=k+1}^{l}\(\log|f'(\xi_j)|-\log|f'(z_j)|\)\right|
=\left|\log\left|\frac{\(f^{-(l-k)}_{k,l}\)'(z_k)}{\(f^{-(l-k)}_{k,l}\)'(\xi_k)
}\right|\right|,
\eeq
where $f^{-(l-k)}_{k,l}: f^{-(k-1)}_{k}(W_s)\lra
{\mathbb C}$, the unique  holomorphic  inverse branch  sending
$\xi_k$ to $\xi_l$, is equal to  $f^{-(l-1)}_{l}\circ f^{k-1}$.
Since for all integers $k \geq 1$  large enough, the map $f^{-(l-k)}_{k,l}$
extends univalently (and holomorphically) to the ball
$$
B\lt(\xi_k,\frac{1}{2}\dist(w, {\rm PC}(f))\rt)
$$
and since $\lim_{n \to\infty}|z_n - \xi_n|=0$ uniformly on $W_{r/2}$, we  conclude from \eqref{2_2017_11_30} and from Theorem~\ref{Euclid-II} that
\beq\lab{1071206}
\lim_{k \to \infty}\sup_{z\in  W_{r/2}}\sup_{l\geq
k}\left\{\left|\sum_{j=k+1}^{l}
\(\log|f'(\xi_j)|-\log|f'(z_j)|\)\right|\right\}=0.
\eeq
This means  that  the sequence  of   partial sums of the series
(\ref{2071206})  is uniformly Cauchy (fundamental), and it therefore
converges   uniformly  to a  harmonic function.
 Thus the function
\beq\lab{4071206}
u_R(z)
:= u(\xi_0)+\log |f'(\xi_1)| -\log|f'(f^{-1}_{\xi_1}(z))|
  + \sum_{n=2}^{\infty}\left(\log |f'(\xi_n)|
    -\log|f'((f^{-1}_{\xi_1}(z))_n)|\right)
\index{(S)}{$u_R$}
\eeq
is well defined and harmonic on $B(\xi_0,R)\sms Q$. By the
assumption (b)  there exist a continuous  function $u: J(f)\sms
\ov{{\rm PS}^0(f)}\lra {\mathbb R}$ and locally constant   function
$c: J(f)\sms \ov{{\rm PS}^0(f)}\lra {\mathbb R}$ such that
\beq\lab{3071206}
\log|f'(z)|= c(z)+u(z)-u(f(z))
\eeq
for all   $z\in J(f)\sms ({\rm PS}^0(f)\cup f^{-1}({\rm PS}^0(f)))$.

Consider now the  set $E \sbt J(f)$ consisting of all those points
$y$ for which $f^{-1}(y)$ is not a subset of
$\ov{{\rm PS}^0(f)}$. Note that $J(f)\sms E\sbt f(\Crit(f))$. Fix 
$$
x\in f^{-1}(y)\sms\ov{{\rm PS}^0(f)}.
$$
Making use of  (\ref{3071206}), for points $z$
near $x$, we deduce that $u$ extends  continuously  to $E$. Noting
also that if $ z\in J(f)\sms \ov{{\rm PS}^0(f)}$, then $f(z)\in E$,
we  further deduce   that (\ref{3071206})  holds for  all $z \in
J(f)\sms \ov{{\rm PS}^0(f)}$, i.e.
\beq\lab{3072606}
\log|f'(z)|= c(z)+u(z)-u(f(z)).
\eeq
Now using (\ref{1082406}) and  (\ref{2082406}) we  conclude  also
that there exists $\hat{R}_Q>0$ so small that  the function
$c:J(f)\sms \ov{{\rm PS}^0(f)} \to \mathbb R$ is constant on  each
set
$$f^{-(n-1)}_n(B(\xi_1,\hat{R}_Q))\cap(J(f)\sms \ov{{\rm PS}^0(f)})
$$ 
and on the set 
$$
B(\xi_1, \hat{R}_Q)\cap (J(f)\sms \ov{{\rm PS}^0(f)}).
$$ 
Taking $R_Q= R_{Q,\xi} \in (0, R)$ so small that 
$$
\ov{B(\xi_0,R_Q)}\sbt B(\xi_0,R)\sms Q 
\quad {\rm and } \quad
f^{-1}_{\xi_1}(B(\xi_0,R_Q)) \sbt B(\xi_1,\hat{R}_Q),
$$ 
and recalling that
$$
\lim_{n\to\infty}u((f^{-1}_{\xi_1}(z))_n)=u(w)
= \lim_{n\to\infty}u(\xi_n)
$$ 
for all $z \in B(\xi_0, R)\sms Q$, uniformly on
$B(\xi_0,R_Q)$ (apply (\ref{1082406}) along with
Theorem~\ref{Euclid-I}  and  the standard normality argument), we
conclude from (\ref{3072606}) and (\ref{4071206}) that
\beq\label{5071206}
u(z)=u_R(z)
\eeq
for all  $z\in B(\xi_0, R_Q)\cap (J(f)\sms \ov{{\rm PS}^0(f)})$ and
$u(z)- u_R(z)$ is locally   constant throughout $$(B(\xi_0, R)\sms
Q)\cap (J(f)\sms \ov{{\rm PS}^0(f)}).$$ Suppose now that $0 <
R^*\leq R_1 \leq R_2$. Since, by Theorem~\ref{thm:julia},
$\HD(J(f))>1$ and, by Proposition~\ref{p120190914}, $\ov{{\rm PS}^0(f)}$ is a nowhere dense subset of
$J(f)$, $J(f)\sms \ov{{\rm PS}^0(f)})$ is not contained in any
countable union of real--analytic curves. It therefore follows
from (\ref{5071206}) that the function $u_{R_2}$, restricted  to $B(\xi_0,
R_1)\sms Q$ coincides with $u_{R_1}$. Thus, we can define a
harmonic function 
$$
\^u_{Q} :{\mathbb C}\sms Q \lra \mathbb R
$$ 
by the formula \index{(S)}{$\^u_Q$} 
$$ 
{\^u_Q}(z)=u_{|z|+1}(z).
$$
Then,
\beq\label{1062606}
\^u_{Q}(z)=u(z)
\eeq
for all $z\in B(\xi_0,R_Q)\cap (J(f)\sms \ov{{\rm PS}^0(f)})$, and we have the following. \sp

{\bf Claim~1.} The function $({\mathbb C} \sms Q) \cap (J(f)\sms \ov{{\rm PS}^0(f)}) \ni z\longmapsto u(z)- \^u_Q(z)$ is locally constant.\sp

\fr Using, in addition, condition (v), we thus conclude the following. \sp

{\bf  Claim~2.}  The function $\^u_{\hat{Q}}-\^u_Q:{\mathbb C} \sms (Q\cup \hat{Q})\lra \R$ is constant on each connected component of ${\mathbb C} \sms (Q\cup \hat{Q})$.\sp

\fr Since $f^{-1}(Q)\spt \Crit(f)\cup f^{-1}(\infty)$, the function $\log|f'|$ is well defined on $\C\sms f^{-1}(Q)$. Hence, the function
$$
\^u_Q-\^u_Q\circ f - \log|f'|:{\mathbb C} \sms \(Q\cup  f^{-1}(Q)\)\lra \R.
$$
is well defined and harmonic. Combining Claim~1 again and (\ref{3072606}), we conclude that the restriction of the function $\^u_Q-\^u_Q\circ f - \log|f'|$ to the set 
$\({\mathbb C} \sms (Q\cup  f^{-1}(Q))\)\cap \(J(f)\sms \ov{{\rm PS}^0(f)}\)$
is locally constant. Using also the fact (see Proposition~\ref{p120190914}), that $\ov{{\rm PS}^0(f)}$ is nowhere dense in $J(f)$, we get the following. 

\sp {\bf  Claim~3.} The function 
$$
\^u_Q-\^u_Q\circ f - \log|f'|:{\mathbb C} \sms \(Q\cup  f^{-1}(Q)\)\lra \R
$$
is constant on each connected component  of ${\mathbb C} \sms (Q \cup
f^{-1}(Q))$ that intersects  $J(f)$.\\

\fr Our  nearest goal   now is  to extend  this  claim  to all connected
components of ${\mathbb C} \sms (Q \cup f^{-1}(Q))$. And indeed,
consider $S$, a connected  component  of $${\mathbb C} \sms
f^{-1}(Q)=f^{-1}({\mathbb C} \sms Q).$$ Two connected  components
$S_1$ and $S_2$ of $${\mathbb C} \sms ( Q\cup f^{-1}(Q))$$ contained
in $S$ are called equivalent if $\ov{S_1}\cap \ov{S_2}$ is a
non-degenerate segment of $Q$ (since $S$ is simply-connected, the
other possibilities  are that either $\ov{S_1}\cap \ov{S_2}=\es$ or
$\ov{S_1}\cap \ov{S_2}$  is a singleton contained in $Q$). A
connected component $S'$ of ${\mathbb C}  \sms (Q\cup f^{-1}(Q))$ is
called tame if the function 
$$
\^u_Q-\^u_Q\circ f - \log|f'|:{\mathbb C} \sms \(Q\cup  f^{-1}(Q)\)\lra \R
$$ 
is constant on $S'$. We shall prove the following.

\sp{\bf Claim~4.} Suppose $S$ is a bounded connected component of
${\mathbb C} \sms f^{-1}(Q)$.  If $S_1$ and $S_2$   are two
arbitrary equivalent connected components   of ${\mathbb C}\sms (Q
\cup f^{-1}(Q))$
contained in  $S$ and $S_1$ is tame, then so  is $S_2$.

\sp\bpf  Let $\Delta=\bar{S_1}\cap {\bar S_2}\sbt Q$.
Perturbing $\Delta$ a little  bit we can replace  it by a closed
segment $\hat{\Delta}$ with the following properties:
\begin{itemize}
\item  [-] $\hat{\Delta}$ has  the same endpoints as  $\Delta$.

\, \item  [-] $\Int \hat{\Delta}\sbt S_2$.

\, \item  [-] If $\hat{Q}$   is obtained from $Q$
by replacing $\Delta$ by $\hat{\Delta}$, then there exists $S_3$, a
connected component  of ${\mathbb C} \sms (\hat{Q}\cup
f^{-1}(\hat{Q}))$ that has non-empty intersections with $S_1$ and
$S_2$.
\end{itemize}
Consider $S_{3,1}$, a connected component of $S_1\cap S_3$. Then
$$
S_{3,1}\cup f(S_{3,1})\sbt {\mathbb C} \sms(Q\cup \hat{Q}),
$$ and
it follows from  Claim~2 that $\^u_{\hat{Q}}- \^u_{Q}$ and
$\^u_{\hat{Q}}\circ f - \^u_{Q}\circ f$ are both constant on
$S_{3,1}$. Since $S_1$ is tame, we thus conclude that the harmonic function
$$
\^u_{\hat{Q}}-\^u_{\hat{Q}}\circ f - \log|f'|:{\mathbb C} \sms \(Q\cup  f^{-1}(Q)\)\lra \R
$$ 
is constant on $S_{3,1}$, and hence, on $S_3$. Let $S_{3,2}$ be a connected
component of $S_3\cap S_2$. As above, $\^u_{Q}- \^u_{\hat{Q}}$ and
$\^u_{Q}\circ f - \^u_{\hat{Q}}\circ f$ are both constant on
$S_{3,2}$. Therefore 
$$
\^u_{Q} - \^u_{Q}\circ f - \log |f'|:{\mathbb C} \sms \(Q\cup  f^{-1}(Q)\)\lra \R
$$ 
is constant on $S_{3,2}$ and hence, $S_2$ is tame. The proof of
Claim~4 is complete. \endpf\sp

Since each connected   component $S$   of $${\mathbb C} \sms
f^{-1}(Q)= f^{-1}({\mathbb C} \sms Q)$$ intersects the   Julia set
$J(f)$, the component $S$  contains  at least   one  tame, connected
component of ${\mathbb C} \sms ( Q\cup f^{-1}(Q))$. Since  in
addition, any two connected components $S'$ and $S''$ of ${\mathbb
C} \sms( Q\cup f^{-1}(Q))$ contained in $S$ can be  connected   by a
sequence $S'=S_1, S_2, \ldots, S_k=S''$ of components  of ${\mathbb
C} \sms( Q\cup f^{-1}(Q))$ contained in $S$ such that any two
consecutive are equivalent, we thus conclude from Claim~4 the following.

\sp{\bf Claim~5.} The function 
$$
\^u_Q-\^u_Q\circ f - \log|f'|:{\mathbb C} \sms \(Q\cup  f^{-1}(Q)\)\lra \R
$$  
is constant on each connected component of ${\mathbb C} \sms(
Q\cup f^{-1}(Q))$.\sp

Define $\Sing_+(\^u_Q)$ and $\Sing_-(\^u_Q)$,  the sets of
 singularities  of $\^u_Q$ as follows.
$$
\Sing_+(\^u_Q):=\{ w \in {\mathbb C}: \limsup_{z \to  w}
  \^u_Q(z)=+\infty\}.\index{(S)}{$\Sing_+$}
$$
and
$$
\Sing_-(\^u_Q):=\{ w \in {\mathbb C}:  \liminf_{z \to w}
  \^u_Q(z)=-\infty\}.\index{(S)}{$\Sing_-$}
$$
Since the function $\^u_Q$ is defined on $\C\sms Q$ while $\^u_{\hat Q}$ is defined on $\C\sms \hat Q$, it immediately follows from Claim~2 that
\beq\label{1083006}
\Sing_+(\^u_Q) \cup \Sing_-(\^u_Q)\sbt \{x\}\cup f(\Crit(f)).
\eeq
Since  the family  of closures   of connected  components  of
${\mathbb C} \sms (Q\cup f^{-1}(Q))$ is   locally finite, as  an
immediate consequence  of Claim~5, we get the  following.
\beq\label{2083006}
f((\Sing_+(\^u_Q)\sms f^{-1}(\infty)) \cup ( \Crit(f)\sms
\Sing_-(\^u_Q))) \sbt \Sing_+(\^u_Q)
\eeq
and
\beq\label{2091106}
f^{-1}(\Sing_+(\^u_Q))  \sbt \Sing_+(\^u_Q)  \cup  \Crit(f).
\eeq
Now, since  the set  $f(\Crit(f))$ is finite, for every point
$w \in f( \Crit(f))$ there exists
$$
z \in \Big(\Crit(f) \sms\(\{x\}\cup
f(\Crit(f))\)\Big)\cap   \bu_{n=1}^\infty f^{-n}(w).
$$ 
It then follows from (\ref{1083006}) that $z \notin \Sing_-(\^u_Q)$. Hence, using
(\ref{2083006}), we see  that $f(z)\in \Sing_+(\^u_Q)$. Since there
is an integer $n\geq 1$ such that $w=f^n(z)$, applying (\ref{2083006}) $n-1$
times, we conclude that $ w \in \Sing_+(\^u_Q).$  So, we have proved
the following.
\beq\label{1091106}
f( \Crit(f)) \sbt \Sing_+(\^u_Q).
\eeq

\sp{\bf  Claim~6.}  The  set $\bu_{n=0}^\infty f^n(f( \Crit(f)))$ is
finite.\\

\bpf  Suppose on the contrary   that this  union  is
infinite. Since the   set $f(\Crit(f))$ is finite, there  thus
exists  $w \in f(\Crit(f))$   such that    all the points $f^n(w),
\, n\geq 0$, are mutually   distinct. If then  follows   from
(\ref{1091106})  and  (\ref{2083006})  that
$$\{f^n(w)\}_{n=0}^\infty \sbt \Sing_+(\^u_Q).$$ Since, by
(\ref{1083006}), the  set $\Sing_+(\^u_Q)$ is finite, we  get a
contradiction which finishes the  proof.

\sp Due to Theorem~\ref{baker-domain+Sullivan for elliptic}, Theorem~\ref{r071708}, Theorem~\ref{t1ms123},
Theorem~\ref{t2ms135}, and the classification of periodic connected components of the Fatou set $F(f)$, provided in Theorem~\ref{Fatou Periodic Components}, as an immediate consequence of this claim, we get the following.

\sp{\bf Claim~7.} 
Either the elliptic function $f:\C\lra\oc$ has a superattracting periodic orbit or else $J(f)={\mathbb C}$.

\sp\fr Consider now two cases. 

\sp\fr{\bf Case~1.} The elliptic function $f:\C\lra\oc$ has a
critical periodic point $c$.

\sp\fr Then, by (\ref{1091106}) and (\ref{2083006}), both points $c$ and
$f(c)$ belong  to $\Sing_+(\^u_Q)$ (the whole forward orbit of $c$
does). Since the set $\Sing_+(\^u_Q) \cup \Sing_-(\^u_Q)$ is finite by
(\ref{1083006}), and since  $[c]$, the equivalence class of $c$ with
respect to  the relation $\sim_f$, is infinite,  there exists $$ w \in
[c]\sms (\Sing_+(\^u_Q)\cup \Sing_{-}(\^u_Q)). 
$$ 
Then, there exist two constants $C_1, C_2\in {\mathbb R}$, $\l \in \Lambda$, and a
sequence $\{w_n\}_{n=1}^\infty$ with the following properties:

\, \begin{enumerate}
\item [($\a$)] $\forall n \geq 1$ $ w_n\notin \Sing_+(\^u_Q)\cup
\Sing_{-}(\^u_Q)$,

\sp \item [($\b$)] $ \lim_{ n\to \infty}w_n =w$ \ and \ $\lim_{ n\to \infty}
(w_n +\l)=c$,

\sp \item[($\gamma$)]   $\^u_Q(w_n)-\^u_Q(f(w_n))-\log|f'(w_n)|=C_1$, (Claim~5)

\sp \item[($\d$)] $\^u_Q(w_n+\l)-\^u_Q(f(w_n+\l))-\log|f'(w_n+\l)|=C_2$,
(Claim~5)

\sp \item[($\varepsilon$)]  $\lim_{n \to \infty}\^u_Q(w_n+\l)=+\infty.$
\end{enumerate}

\sp\fr Since $f(w_n+\l)=f(w_n)$  and $f'(w_n+\l)=f'(w_n)$,  ($\gamma$) and ($\delta$)
 imply that $\^u_Q(w_n+\l)=\^u_Q(w_n)+ C_2-C_1$. Since $w
 \notin\Sing_+(\^u_Q)$, we conclude that
$$
\ov{\lim}_{n \to \infty}\^u_Q(w_n+\l)
=C_2-C_1+\ov{\lim}_{n \to\infty}\^u_Q(w_n)<+\infty, 
$$ 
contrary to ($\varepsilon$). The Case~1 is thus ruled out.

\sp\fr{\bf  Case~2.} $J(f)={\mathbb C}$

\,\fr By  Claim~6, the set ${\mathbb C} \sms \ov{{\rm PS}^0(f)}$ is
connected, and therefore,  the locally  constant  function  $ c:
{\mathbb C} \sms \ov{{\rm PS}^0(f)}\lra {\mathbb R}$ is constant,
say equal to $c$. From Claim~6 and Banach Contraction Principle
(applied for univalent holomorphic inverse branches of $f$ mapping
fundamental regions into themselves) there exists a sequence
$\{z_n\}_{1}^\infty$ of points in ${\mathbb C} \sms \ov{{\rm
PS}^0(f)}$ with the following properties
$$
\lim_{n \to \infty} z_n = \infty
$$ 
and $f(z_n)=z_n$ for all $n \geq 1$. Since then $\lim_{n \to \infty} f( z_n)=\infty$, it follows that
$$
\lim_{n \to \infty} \dist( z_n, f^{-1}(\infty))=0,
$$ 
and therefore
$\lim_{n\to \infty}\log|f'(z_n)|=+\infty$. This is a  contradiction
with (\ref{3072606}) since this formula implies that
$\log|f'(z_n)|=c(z_n)=c$ for all $n\ge 1$. The Case~2 is thus also
ruled out. In conclusion, (b) leads to contradiction and the proof of
Theorem~\ref{t3042706} is complete.
\endpf

\sp Any function satisfying any, equivalently all, conditions of
Theorem~\ref{t3042706} is commonly referred to as essentially non--linear.

\sp\section{Proof of the Rigidity Theorem}

\fr We will prove in this section the chain of implications $\bf (1)\Rightarrow (2)
\Rightarrow (3) \Rightarrow (4)\Rightarrow (5) \Rightarrow (6)
\Rightarrow (1)$ from Theorem~\ref{t1042706}, 
which will complete its proof.

\sp We  will frequently apply the following easy fact.

\sp\blem\label{l1j13} Every homeomorphism $h: X\lra Y$, where $X$ and
$Y$ are  closed  subsets  of $\mathbb C$, is uniformly   continuous
with respect to the spherical metric  on $\mathbb C$. \elem

\bpf If  one of the sets $X$ or $Y$ is compact, then so
is the other, and we are  done. So, we may assume that neither $X$
nor $Y$ are  compact. It   suffices  to show  that the map 
$$
\hat{h}:X\cup \{\infty\} \lra  Y\cup \{\infty\},
$$
determined  by the requirements  
$$
\hat{h}_{|X}=h
\  \  \  {\rm and} \  \  \
\hat{h}(\infty)=\infty
$$
\index{(S)}{$\hat{h}(\infty)$}, is
continuous at $\infty$. To prove this  suppose for the contrary that
$\hat{h}$ is not continuous at $\infty$. This means that there
exists  a sequence $(x_n)_{n=0}^\infty\sbt X$ converging to $\infty$
such that $( \hat{h}(x_n))_{n=0}^\infty $ does not converge to
$\infty$. Passing to a subsequence, we may assume without  lost of
generality that the sequence $( \hat{h}(x_n))_{n=0}^\infty $  is
bounded. But then its closure $\(\ov{{h}(x_n)}\)_{n=0}^\infty $ is
compact. So, the set $h^{-1}\((\ov{{h}(x_n)})_{n=0}^\infty \)$ is
compact, thus bounded. Since $(x_n)_{n=0}^\infty  \sbt
h^{-1}({h}(x_n)_{n=0}^\infty) $, the sequence $(x_n)_{n=0}^\infty$
is bounded, contrary to  the assumption that it converges to
$\infty$. \endpf

\sp\fr Using this lemma we easily get  the following.

\sp\bprop\label{p2j13}
If  a homeomorphism  $h:J(f)\lra J(g)$ satisfying condition~(5) of
Theorem~\ref{t1042706} conjugates two elliptic functions  $f$ and $g$, then
$$
h(\Crit(J(f)))=\Crit(J(g)),\  \ h(\Om(f))=\Om(g),  \ \
h(f^{-1}(\infty))= g^{-1}(\infty),
$$
and 
$$
{\rm PC}(f)={\rm PC}(g).
$$
\eprop

\sp\fr The implications  $\bf (1)\Rightarrow (2)$, $\bf (2)\Rightarrow
(3)$ and  $\bf (3) \Rightarrow (4)$ 
of Theorem~\ref{t1042706} are  obvious. For the  sake  of completeness
we provide  now an 
 easy proof of the  implication  $\bf (4) \Rightarrow (5)$ from this
 theorem.  Suppose that 
 $x \in J(f)$ is a periodic  point of period $p$ and
 $$|(f^p)'(x)|\neq |(g^p)'(h(x))|.$$
 Without losing  generality we
 may assume that $$|(f^p)'(x)|< |(g^p)'(h(x))|.$$
 Fix $$|(f^p)'(x)|< \mu < \l < |(g^p)'(h(x))|.
 $$
 Let $U$ be a neighborhood  of $x$ such that  both continuous inverse
 branches 
 $$
 f^{-p}_x: U \lra U
 $$  
 of $f^p$ and 
 $$
 g^{-p}_{h(x)}: h(U\cap J(f)) \lra h(U\cap J(f))
 $$ 
 of $g^p$ sending   respectively  $x$ to $x$ and $h(x)$ to
 $h(x)$ are well defined. Taking $U$ sufficiently  small, we may
 assume that 
$$
|f^{-pn}_x(z)-x|\geq  \mu^{-n} 
\quad {\rm and} \quad
|g^{-pn}_{h(x)}(w)-h(x)|\leq  \l^{-n}
$$ 
for all $n \geq 1$, $z\in U$  and $w\in h(U)$. Hence
$$ 
\frac{|f^{-pn}_x(z)-x|}{|h(f^{-pn}_x(z))-h(x)| }  
\geq \frac{\mu^{-n}}{\l^{-n}}
= \left(\frac{\l}{\mu}\right)^n \lra\infty
$$
as $n \to \infty$. The implication $(4) \Rightarrow (5)$  is
proved.

\sp In order to show that $\bf (5)\Rightarrow(6)$ in
Theorem~\ref{t1042706}, we need first the
following version  of the closing lemma (or shadowing lemma). 

\sp\blem\label{l1112506} 
Let $f:\C\lra\oc$ be a non--recurrent elliptic function. Fix $ s>0$. Then for all $0< \rho_2< s$ there exist $ \rho_1>0$ and an integer $n_1\geq 1$ such that for every integer $n\geq n_1$ if 
$$
f^n(x)\in J(f)\sms B_s( \ov{{\rm PC}(f)},s)
\quad {\rm and} \quad
f^n(x) \in B_s(x, \rho_1),
$$
then there exists $y\in  J(f)$ such
that 
$$
f^n(y)=y, \quad |y-f^n(x)|_*< \rho_2,
$$
and 
$$
|f^j(x)- f^j(y)|_* \leq  \rho_2
$$ 
for all integers $0\leq j \leq n - n_1$. 
\elem

\bpf It follows from lemma~\ref{l120180306}
that 
$$
\lim_{n\to \infty} \sup \{ \diam_s(P_n)\}=0,
$$
where $P_n$ ranges over all connected components of 
$f^{-n}(\ov{B_s(z, \rho_2)})$, $ z\in J(f)\sms B_s(\ov{{\rm
PC}(f)},s)$. Take $n_1\ge 1$ so large that 
$$
\diam_s(P_n)< \rho_2/2
$$ 
for all integers $n \geq n_1$. Then take $\rho_1 < \rho_2/2$. Let $B_n$, $n\geq n_1$, be the connected component of $f^{-n}(\ov{B_s( f^n(x), \rho_2)})$
containing $x$. Let $f^{-n}_x:B_s( f^n(x),s)\lra B_n$ be
the holomorphic inverse branch of $f^n$ sending $f^n(x)$ to $x$.
We then have
$$f^{-n}_x(\ov{B_s( f^n(x), \rho_2)}) 
\sbt \ov{B_s( x, \rho_2/2)} 
\sbt \ov{B_s( f^n(x),
\rho_1+\rho_2/2)} 
\sbt \ov{B_s( f^n(x), \rho_2)}.
$$ 
Hence, by Brower's Fixed Point Theorem there exists $y\in  \ov{B_s( f^n(x),
\rho_2)}$ such that $f_x^{-n}(y)=y$ ,which implies that $f^n(y) =y$.
Finally note  that
$$|f^{-j}_x(f^n(x))- f^{-j}_x(f^n(y))|_*= |f^{-j}_x(f^n(x))-
f^{-j}_x(f^n(y))|_*\leq \rho_2/2 \leq \rho_2$$  for all $j\geq n_1$. The proof is complete.
\endpf

\sp By topological exactness of an elliptic function $f:\C\lra\oc$, the set  of  transitive  points of
$f$ is dense  in $J(f)$. Choose  one such  a point, say $x$. For
every $z\in J(f)\sms (\Crit(f)\cup f^{-1}(\infty))$  define
$$ 
\eta(z):= \log |g'(h(z))|- \log |f'(z)|. 
$$  
Note that $x \notin\Crit(f)\cup f^{-1}(\infty)$ and for every integer $n \geq 1$, set
\beq\label{1112506}
u(f^n(x)):=\sum_{j=0}^{n-1}\eta( f^j(x)).
\eeq

\blem\label{l2112506} Let $f:\C\lra\oc$ be a non--recurrent elliptic function. Suppose  that condition~(5) of Theorem~\ref{t1042706} holds. Then for  every $s >0$ the function $u$ restricted to the set 
$$
\(J(f)\sms B_s( \ov{{\rm PC}(f)},
s)\)\cap \{f^n(x):n\geq 0\}
$$ 
is uniformly  continuous with respect to the spherical metric. 
\elem

\bpf By Lemma~\ref{l1j13}  and Proposition~\ref{p2j13}
there exists $s'>0$  such that  $$h(J(f)\sms B_s(\ov{{\rm
PC}(f)},s)) \sbt J(g)\sms  B_s(\ov{{\rm PC}(g)},s').$$ Fix $\d\in
(0, \min\{s, s'\})$. By Lemma~\ref{l1j13}  there  exists $0< \rho_2<
\d$ so small that 
$$
|h(z)-h(w)|_*\leq \d
$$ 
whenever $|z-w|_*\leq\rho_2$. Choose $\rho_1$ and $n_1$ according  to
Lemma~\ref{l1112506} applied to the function $f$. Consider two
points $f^n(x)$ and $f^m(x)$  in $J(f)\sms B_s( \ov{{\rm PC}(f)},
s)$ such that  $n \geq m$ and 
$$
|f^n(x)-f^m(x)|_* < \rho_1.
$$
Then, in view of Lemma~\ref{l1112506}, there exists  a point $y \in J(f)$ such
that $f^{n-m}(f^m(y))=f^m(y)$,
$$ |f^{m+j}(x)- f^{m+j}(y)|_*< \rho_2$$
for all $j=0,1, \ldots, n-m-n_1$,  and $ |f^{n}(x)- f^{n}(y)|_*<
\rho_2$. Then $\sum_{j=m}^{n-1}\eta(f^j(y))=0$, and we get
$$\begin{aligned}
u(f^n(x))- u(f^m(x))=& \sum_{j=m}^{n-1}\eta(f^j(x))=
\sum_{j=m}^{n-1}\(\eta(f^j(x))- \eta(f^j(y))\)\\
=&\sum_{j=m}^{n-1}(\log|g'(h(f^j(x)))|-\log|g'(h(f^j(y)))|)\\
&  -  (\log| f'(f^j(x))|  - \log| f'(f^j(y))|)\\
=& \log\left|\frac{(g^{n-m})'(h(f^m(x)))}{(g^{n-m})'(h(f^m(y)))}
\right|- \log\left|\frac{(f^{n-m})'(f^m(x))}{(f^{n-m})'(f^m(y))}
\right|.
\end{aligned}
$$
We want $u(f^n(x))$ and $u(f^m(x))$  to  be  close to each
other if $f^n(x)$  and $f^m(x)$ are. For this  it suffices  to know
that both terms
$$
\log\left|\frac{(g^{n-m})'(h(g^m(x)))}{(g^{n-m})'(h(g^m(y)))}\right|
\quad {\rm and} \quad
\log\left|\frac{(f^{n-m})'(f^m(x))}{(f^{n-m})'(f^m(y))} \right|
$$ 
are small. But for  the later the term this follows from
Theorem~\ref{Spher-I} since $|f^m(x)-f^m(y)|_* < \rho_2< \delta$ and
the inverse branch $f^{-(n-m)}_{f^m(x)}$ sending $f^n(x)$  to
$f^m(x)$  is defined  on $B_s(f^n(x), s)$. For the former term,
since $|f^m(x)-f^m(y)|_* < \rho_2$ we have
$$|g^m(h(x))-g^m(h(x))|_*=|h(f^m(x))-h(f^m(y))|_*\leq \d, $$ and note
that the inverse branch  $g^{-(n-m)}_{g^m(h(x))}$ sending $g^n(x)$
to $g^m(x)$  is defined  on $B_s(g^n(x), s')$. The proof is
finished. \endpf

\sp\fr Consequently the function $u$ extends  continuously   to each  set
$$
J(f)\sms B_s(\ov{{\rm PC}(f)}, s),\quad s >0,
$$ 
and, therefore, to the set
$J(f)\sms \ov{{\rm PC}(f)}$.

\sp\blem\label{l1122606} 
If $f:\C\lra\oc$ is a non--recurrent elliptic function, then the functions 
$$
\log|f'(z)|
\quad {\rm and} \quad
\log|g'(h(z))|
$$ 
are cohomologous  in the class  of continuous
functions  on $ (J(f)\sms \ov{{\rm PC}(f)})\cap  f^{-1}(J(f)\sms
\ov{{\rm PC}(f)})$. More precisely there exists a  continuous
function $u: J(f)\sms \ov{{\rm PC}(f)} \lra \mathbb R$ such that
$$ \log|g'(h(z))|-\log|f'(z)|=u(f(z))-u(z)$$
 for all $z \in (J(f)\sms \ov{{ \rm PC}(f)})\cap f^{-1}(J(f)\sms {\ov{\rm
 PC}(f)})$.
\elem

\bpf From  (\ref{1112506})  we  get
$$\eta(f^n(x))=u(f^n(x))-u(f^{n-1}(x)).$$
Since  the set $\{f^n(x): n\geq 1\}$ is dense  in $$(J(f)\sms
\ov{{\rm PC}(f)})\cap( J(f)\sms f^{-1}(\ov{{\rm PC}(f)}))$$ and all
the functions $\eta$, $u$, $u \circ f$ are continuous in this set,
the lemma is proved.
\endpf

\sp\fr{\bf Proof of the implication $\bf (5)\Rightarrow(6)$ of
Theorem~\ref{t1042706}.} In the proof of this implication, we do assume that $f:\C\lra\oc$ is a compactly non--recurrent regular elliptic function.

\sp\fr Let 
$$
Q(f)\index{(S)}{$Q(f)$}:=\big\{z\in J(f): \limsup_{j \to \infty} \dist_s(\f^j(z), \ov{{\rm PC}(f)}\)> 0\big\}
$$
and
$$ 
Q_n(f)\index{(S)}{$Q_n(f)$}
:=\big\{z\in J(f): \limsup_{j \to \infty}\dist_s(f^j(z), \ov{{\rm 
 PC}(f)})> 2/n\big\},\,\, n \in {\mathbb N}.
$$
Then $Q_n(f) \sbt Q(f).$ By Lemma~\ref{l1j13} and
 Proposition~\ref{p2j13} there  exists $k_n 
\geq 1$  such that  $\dist_s(x, \ov{{\rm PC}(f)})\geq 2/n$ entails
$\dist_s(h(x), \ov{{\rm PC}(f)})\geq 2/{k_n}$. In particular:
$$
h(Q_n(f))\sbt Q_{k_n}(g).
$$ 
Now fix $z\in Q_n(f)$. Then there exists an infinite sequence $n_j=n_j(z)$, $j\ge 1$, such that
$$
\dist_s(f^{n_j}(z),\ov{{\rm PC}(f)})\geq 2/n
$$ 
for all integers $j \geq 1$. Applying now Lemma~\ref{l2112506} and Lemma~\ref{l1122606}, we see that there exists a constant $K_n \geq 1$  such that
\beq\label{3122706}
K_n^{-1} \leq \frac{|(g^{n_j})'(h(z))|}{|(f^{n_j})'(z)|} \leq K_n.
\eeq
In view of uniform  continuity of the homeomorphism $h$, there exists $\gamma_n \in (0,1/n)$ such that  
$$
h(J(f)\cap B_s(x,\g_n))\sbt  B_s(h(x), 1/k_n)
$$
for all $x\in J(f)$. By the  choice of  the sequence $n_j=n_j(z)$
for every $j \geq 1$  there exists a  holomorphic inverse  branch
$$
f^{-n_j}_z:B_s(f^{n_j}(z), 2/n)\lra\oc
$$  
of $f^{n_j}$, sending $f^{n_j}(z)$ to $z$. Set
$$
r_j(z):=K^{-1}|(f^{-n_j}_z)^*(f^{n_j}(z))|\g_n.
$$
Then, by Lemma~\ref{lncp12.9s.}, we have that
$$ 
B_s(z, r_j(z)) \supset f^{-n_j}_z ( B_s(f^{n_j}(z),K^{-2} \g_n)),
$$ 
and therefore
\beq\label{2122706}
\begin{aligned}
m_f(B_s(z, r_j(z)))& \geq K^{-h_f}|(f^{-n_j}_z)^*(f^{n_j}(z))|^{h_f}
m_f(B_s(f^{n_j}(z), K^{-2}\g_n))\\
& \geq M_n K^{-h_f}|(f^{-n_j}_z)^*(f^{n_j}(z))|^{h_f},
\end{aligned}
\eeq
where $M_n := \inf\{ m_f( B_s(x, K^{-2}\g_n)): x\in J(f)
\}>0$. Similarly, by Lemma~\ref{lncp12.9s.}, $B_s(z, r_j(z)) \sbt
f^{-n_j}_z(B_s(f^{n_j}(z), \g_n))$. Hence,
$$ \begin{aligned}
h(J(f)\cap B_s(z, r_j(z))) & \sbt h (f^{-n_j}_z(B_s(f^{n_j}(z),
\g_n)\cap J(f)))\\
& = g^{-n_j}_{h(z)}(h(B_s(f^{n_j}(z),
\g_n)\cap J(f)))\\
& \sbt g^{-n_j}_{h(z)}(B_s(h(f^{n_j}(z)),
1/k_n))\\
& = g^{-n_j}_{h(z)}(B_s(g^{n_j}(h(z)),
1/k_n)).\\
\end{aligned}
$$
Therefore using (\ref{3122706}), (\ref{2122706}) and
Theorem~\ref{Spher-II}, we get
$$
\begin{aligned}
m_g\circ h(B_s(z, r_j(z)))&\leq
m_g(g^{-n_j}_{h(z)}(B_s(g^{n_j}(h(z)), 1/k_n)))\\
&  \leq K^{h_g}|(g^{-n_j}_{h(z)})^*(g^{n_j}(h(z)))|^{h_g}
m_g(B_s(g^{n_j}(h(z)), 1/k_n))\\
& \leq K^{h_g}|(g^{n_j})^*(h(z))|^{-h_g}\\
&\leq  K^{h_g}K_n^{h_g}|(f^{n_j})^*(z)|^{-h_f +(h_f-h_g)}\\
&\leq M_n^{-1} K^{h_f+h_g}K_n^{h_g}m_f(B_s(z, r_j(z)))|(f^{n_j})^*(z)|^{h_f -h_g}.\\
\end{aligned}
$$
If $h_f-h_g<0$, then  it would follow  from Lemma~\ref{l3123006} and
the fact that $\lim_{j \to \infty}|(f^{n_j})'(z)|=+\infty$ (it
implies that $\lim_{j \to \infty}|(f^{n_j})^*(z)|=+\infty$), that
$m_g\circ h(Q_n(f))=0$  for every  $ n \geq 1$. Since
$$
\bu_{n \geq 1}Q_n(f)=Q(f),
$$  
it would imply that $m_g(h(Q(f)))=0$. But since $h( \ov{{\rm PC}(f)})=\ov{{\rm PC}(g)}$, using uniform continuity of $h$ and $h^{-1}$, we  conclude that
$h(Q(f))=Q(g)$. Hence $m_g(Q(g))=0$ which would contradict
Theorem~\ref{tmaincm}. Thus for every $ n \geq 1$ and every $z\in Q_n(f)$
$$ 
m_g\circ h(B(z, r_j(z))) \leq M_n ^{-1} K^{2h_f}K_n^{h_g} m_f( B_s(z,
r_j(z))).
$$ 
Therefore, applying Lemma~\ref{l3123006} we conclude
that $m_g\circ h|_{Q_n(f)}$ is  absolutely  continuous  with respect
to $m_f|_{Q_n(f)}$  for every  $n \geq 1$. Since $$\bu_{n \geq
1}Q_n(f)=Q(f),$$ this  implies that $m_g \circ h|_{T(f)}$
is absolutely continuous with respect  $m_f|_{Q(f)}$. Since 
$$
m_g\circ h(\ov{{\rm PC}(f)})=m_g(\ov{{\rm PC}(f)}) =0
$$ 
and 
$$m_f(\ov{{\rm PC}(f)})=0,
$$
we obtain that  $m_g\circ h$ is absolutely continuous with respect to
$m_f$. By symmetry $m_f \circ h^{-1}$ is   absolutely continuous
with respect  to $m_g$ and consequently  the measure $m_g\circ h $
and $m_f$  are equivalent. The proof  of the implication
(5)$\Rightarrow$(6) is finished. \endpf

\sp We  are thus left to establish  the  implication $\bf (6)\Rightarrow
(1)$ of Theorem~\ref{t1042706}. The proof is quite long. As the first step, we shall  prove the following.

\sp\blem\label{l1j15} 
Let $f:\C\lra\oc$ be a compactly non--recurrent regular elliptic function. Let $f:\C\lra\oc$ be a compactly non--recurrent elliptic function. If condition~(6) of Theorem~\ref{t1042706} is satisfied, then the homeomorphism $ h:J(f)\lra J(g)$ extends to a real--analytic endomorphism from a neighborhood of  $J(f)\sms \ov{{\rm
PC}(f)}$ onto a neighborhood  of $J(g)\sms \ov{{\rm PC}(g)}$.
\elem

\bpf  In view of  Theorem~\ref{t3042706}~(d) there exist $n \geq 1$ and $z\in J(g)\sms {\rm PS}_{-}(g)$ such that   
$$
\det\(\nabla(D_{\mu_g}\circ g^n)(z),\nabla(D_{\mu_g})\)(z))\neq 0.
$$
Therefore, using Corollary~\ref{c2pj11}, we  conclude  that there
exists  an  open set  $W \sbt {\mathbb C} \sms  \bu_{j=0}^n
g^{-j}(\infty)$ containing $z$ and such that
\beq\label{1j15}
\det\(\nabla(D_{\mu_g}\circ g^n)(w),\nabla (D_{\mu_g})(w)\)\neq 0
\eeq
for all $w \in W$. Since the measures $m_f\circ h^{-1}$ and
$m_g$ are equivalent, the ergodic measures $\mu_f \circ h^{-1}$
and $\mu_g$ coincide up to a multiplicative constant. Thus
\beq\label{2j15}
D_{\mu_f}\circ f^j= D_{\mu_g}\circ g^j\circ h
\eeq
throughout $J(f)\sms \bu_{i=0}^j f^{-i}(\ov{{\rm PS}^0(f)})$ for
every integer $j\geq 0$. Making use of Corollary~\ref{c2pj11}, define
$$ 
F(x):=\(D_{\mu_f}(x), D_{\mu_f}\circ f^n(x)\)
$$ 
and
$$
G(y):=\(D_{\mu_g}(y), D_{\mu_g}\circ g^n(y)\)
$$
for $x\in U$, an open neighborhood of $J(f)\sms \bu_{i=0}^n
f^{-i}(\ov{{\rm PS}^0(f)})$ and $y\in V$, an open neighborhood of $J(g)\sms \bu_{i=0}^ng^{-i}(\ov{{\rm PS}^0(g)})$. Shriking $W$ if necessary, we may assume without loss of generality that $W\sbt U$ and, because of (\ref{1j15}), that $G$  is invertible on $W$. By (\ref{2j15}), $F(h^{-1}(z))=G(z)$, and therefore, there exists $U_z\sbt U$, an open neighborhood of $h^{-1}(z)$, such that $F(U_z) \sbt
G(W)$. Hence, the map $G^{-1}\circ F$ is well defined on $U_z$,
and, because of $(\ref{2j15})$ again,
$$
G^{-1} \circ F(x)=h(x)
$$
for all $x\in J(f)\cap U_z$. 

Now consider $\xi$, an arbitrary
point in $J(f)\sms \ov{{\rm PC}(f)} $. By topological exactness of the the elliptic map $f:\C\lra\oc$, there exists an integer $k\geq 0$ and a point
$\hat{\xi}\in U_z\cap f^{-k}(\xi)$. Then there exists $r_\xi>0$  so
small that  
$$
r_\xi < \frac{1}{2} \dist(\xi, \ov{{\rm PC}(f)})
\  \  \  {\rm and} \  \  \
f^{-k}_{\hat\xi}(B(\xi, r_\xi)) \sbt U_z,
$$
where $f^{-k}_{\hat\xi}(B(\xi,2r_\xi)) \lra\C$ is the unique holomorphic branch of $f^{-1}$ sending $\xi$ to $\hat\xi$. Hence, the map 
$$
g^k\circ (G^{-1}\circ F)\circ f^{-k}_{\hat{\xi}}
$$ 
is well defined on $B(\xi, r_\xi)$, real--analytic, and, since $h$ conjugates $f$ 
and $g$, the map 
$$ 
g^k\circ (G^{-1}\circ F) \circ f^{-k}_{\hat{\xi}},
$$
restricted to $J(f)\cap B(\xi,r_\xi)$, coincides with $h$. Now,
since  no open subset of $J(f)$ is contained  in a countable  union
of real-analytic curves, the same argument as in the end of  the
proof of Lemma~\ref{l1pj3}, shows  that all  the maps $$g^k\circ
(G^{-1}\circ F) \circ f^{-k}_\xi$$  glue together  on the balls
$B(\xi, r_\xi/2)$ to form a real--analytic map from  
$$
\bu_{z\in J(f)\sms \ov{{\rm PC}(f)}}B(z, r_z/2)
$$ 
onto an open neighborhood of $J(g)\sms \ov{{\rm PC}(g)}$. This map restricted to $J(f)$ coincides with $h$. The proof is complete.
\endpf

\sp Now we pass to the next step.

\sp\blem\label{l1j19} 
Let $f:\C\lra\oc$ be a compactly non--recurrent elliptic function. If a topological conjugacy $h:J(f)\lra J(g)$ has a real--analytic  extension  on an  open  neighborhood of $J(f)\sms \ov{{\rm PC}(f)}$ in $\mathbb C$, then it  has a holomorphic conformal extension on an open neighborhood of
$J(f)\sms \ov{{\rm PC}(f)}$.  
\elem

\bpf  Let $ H: U \lra \mathbb C$ be a real--analytic
extension of $h$ on an open neighborhood  $U$ of $J(f)\sms
\ov{{\rm PS}^0(f)}$. Hence, the complex dilatation 
$$
\mu_H:=\frac{{\ov{\partial}{H}}}{\partial{H}}
$$ 
is well defined throughout $U$ (shrinking it if necessary). Since $H\circ f = g\circ H$, we get that $H=g\circ H\circ f^{-1}_z$ on some ball $B(f(z), r_z)$ and all
$$
z\in J(f)\sms f^{-1}(\ov{{\rm PC}(f)}).
$$ 
Hence, since the maps $f$ and $g$ are conformal,
\beq\label{1j19}
 \mu_H(f(z))=\mu_{g\circ
H}(z) \left( \frac{f'(z)}{|f'(z)|}\right)^2= \mu_{ H}(z)
 \left(\frac{f'(z)}{|f'(z)|}\right)^2.
\eeq
It follows  from  this  equation that if $\mu_H(w)=0$ at some
point $w\in J(f)\sms \ov{{\rm PC}(f)}$, then $\mu_H$ vanishes
everywhere  on $f^{-1}(w)$ and $f^{-1}(w)\cap \ov{{\rm
PS}^0(f)}=\es.$ So, by induction, $\mu_H$ vanishes on
$\bu_{n=0}^\infty f^{-n}(w)$. Since, by transitivity of  $f$, this
set  is dense  in $J(f)$  and since $\mu_H$ is continuous  on
$J(f)\sms  \ov{{\rm PC}(f)}$, we  thus conclude that $\mu_H$
vanishes throughout $J(f)\sms  \ov{{\rm PC}(f)}$. So, if
$\mu_H(z)\neq 0$ at some point $z\in J(f)\sms \ov{{\rm PS}^0(f)}$,
then $\mu_H$ does not vanish anywhere on $J(f) \sms \ov{{\rm
PC}(f)}$. It then firstly follows from (\ref{1j19}) that the
modulus of $\mu_H$  is constant on backward orbits of $f$, and secondly, that
$$
\frac{\mu_H}{| \mu_H|}(f(z))=
\frac{\mu_{H}}{| \mu_H|}(z)\left( \frac{f'(z)}{|f'(z)|}\right)^2,
\quad  z\in J(f)\sms (\ov{{\rm PC}(f)}\cup  f^{-1}(\ov{{\rm
PC}(f)})).
$$ 
Thus $\frac{\mu_H}{|\mu_H|}$ is a continuous
invariant line field on $J(f)\sms \ov{{\rm PC}(f)}$, contrary to
Theorem~\ref{t3042706}~(c). So, $\mu_H(z)=0$ for all $z\in J(f) \sms
\ov{{\rm PC}(f)}$. Since $\mu_H$  is real--analytic on $U$ and
since no  non-empty open subset of $J(f) \sms \ov{{\rm PC}(f)}$
is contained in a  countable  union  of real-analytic  curves, we
conclude  that $\mu_H$   vanishes on some neighborhood of $J(f)\sms
 \ov{{\rm PC}(f)}$, meaning  that $H$  is conformal on this
neighborhood. We are done. \endpf

\sp \blem\label{l1j21} 
Let $f:\C\lra\oc$ be a compactly non--recurrent elliptic function. 
Suppose that a topological conjugacy 
$$
h:J(f)\lra J(g)
$$ 
has a holomorphic conformal extension on some neighborhood of $J(f) \sms \ov{{\rm PC}(f)}$ and 
$$
h^{-1}: J(g)\lra J(f)
$$ 
has a holomorphic conformal extension on some neighborhood of $J(g)\sms \ov{{\rm PC}(g)}$. 

Then $h$ extends to an affine map ($z\mapsto az
+b$) on $\mathbb C$ and this map is a conjugacy between $f:\C\lra\oc$ and $g: \C\lra\oc$.
\elem

\bpf Let $U_f$ and $U_g$ be the respective neighborhoods of $J(f)
\sms \ov{{\rm PC}(f)}$ and $J(g) \sms \ov{{\rm PC}(g)}$ coming form hypotheses of Lemma~\ref{l1j21}. Let 
$$
H: U_f \lra {\mathbb C}
$$  
be the holomorphic extension of the conjugacy $h: J(f) \lra J(g)$ coming from the hypothesis of the lemma. Because of Lemma~\ref{l1j13}, shrinking $U_f$ if necessary, we may assume  without loss of generality that $H$ maps bounded subsets of $U_f$  onto bounded subsets of $\mathbb C$.

Now, for every $z\in J(f)\sms \ov{{\rm PC}(f)}$ and every $w\in J(g)\sms \ov{{\rm PC}(g)}$ there exist radii $r_z>0$ and $r_w>0$ such that 
$$
B_e(z,2r_z)\sbt U_f\sms \ov{{\rm PC}(f)}
\  \  \  {\rm and} \  \  \
B_e(w,2r_w)\sbt U_g\sms \ov{{\rm PC}(g)}.
$$
Let 
$$
U_f^*:=\bu_{z\in J(f)\sms \ov{{\rm PC}(f)}}B_e(z,r_z)
\  \   \   {\rm and} \  \  \
U_g^*:=\bu_{w\in J(g)\sms \ov{{\rm PC}(g)}}B_e(w,r_w).
$$
Then for every $z\in U_f^*$ and every $w\in U_g^*$ there are points $z'\in J(f)\sms \ov{{\rm PC}(f)}$ and $w\in J(g)\sms \ov{{\rm PC}(g)}$ such that $z\in B(z',r_{z'})$ and $w\in B(w',r_{w'})$. Consequently, setting $r_z^*:=r_{z'}$ and $r_w^*:=r_{w'}$, we have that
\beq\lab{1_2017_12_05}
J(f)\cap B_e(z,r_z^*)\neq \es, \quad B_e(z,r_z^*) \sbt U_f\sms \ov{{\rm PC}(f)},
\eeq
and 
\beq\lab{2_2017_12_05}
J(g)\cap B_e(w,r_w^*)\neq \es, \quad B_e(w,r_w^*) \sbt U_g\sms \ov{{\rm PC}(g)},
\eeq

We shall prove the following.

\sp\fr{\bf Claim~$1^0$.} {\it If $B$  is a round open ball contained  in
$U_f^*$ and $f^{-n}_*:B \to {\mathbb C}$ is a holomorphic inverse
branch of $f^n$ such that $f^{-n}_*(B) \sbt U_f^*$, then 
$$ 
g^n\circ H \circ f^{-n}_*: B \lra \hat{\mathbb C}
$$   
coincides with $H|_B$}. 

\bpf  Fix $z\in B$. We have $w:=f^{-n}_*(z) \in U_f^*$ and, by \eqref{1_2017_12_05}, $B_e(z,r_z^*) \cap \ov{{\rm PC}(f)}=\es$. Thus, there exists a
unique holomorphic inverse branch $f^{-n}_w :B_e(z, r_z^*) \lra
{\mathbb C}$ of $f^n$ sending $z$ to $w$. Since $f^{-n}_* $ and
$f^{-n}_w$ agree on $B\cap B_e(z,r_z^*)$, they glue together  to a
holomorphic function 
$$
F: B\cup B_e(z, r_z^*) \lra {\mathbb C}. 
$$
By \eqref{1_2017_12_05} there exists a point $\xi\in J(f)\cap B_e(z,r_z^*)$. Then, by \eqref{1_2017_12_05} again,
$$
f^{-n}_w(\xi)\in J(f)\sms f^{-1}\(\ov{{\rm PC}(f)}\)
\sbt J(f)\sms \ov{{\rm PC}(f)}
$$
Hence, keeping also in mind that $J(f)\sms \ov{{\rm PC}(f)}\sbt U_f^*$, there exists $\d>0$ such that $B_e(\xi,\d)\sbt B_e(z,r_z^*)$ and
$$
f^{-n}_w(B_e(\xi,\d))\sbt U_f^*\sms \ov{{\rm PC}(f)}.
$$
Therefore, the map $g^n\circ H \circ f^{-n}_w|_{B_e(\xi,\d)}$ is well defined and 
$$
g^n\circ H \circ f^{-n}_w|_{B_e(\xi,\d)\cap J(f)}
= h|_{B_e(\xi,\d)\cap J(f)}
= H|_{B_e(\xi,\d)\cap J(f)}.
$$ 
Since also $B_e(\xi,\d)\cap J(f)\sms \ov{{\rm PC}(f)}$ is an uncountable set (as $\xi\in J(f)\sms \ov{{\rm PC}(f)}$ and $\ov{{\rm PC}(f)}$ is, by Proposition~\ref{p120190914}, a nowhere dense subset of $J(f)$), we therefore conclude that
$$
g^n\circ H \circ f^{-n}_w|_{B_e(\xi,\d)}= H|_{B_e(\xi,\d)}.
$$ 
This means that
$$
g^n\circ H \circ F|_{B_e(\xi,\d)}= H|_{B_e(\xi,\d)}.
$$ 
Thus
$$
g^n\circ H \circ F|_{B\cup B_e(\xi,\d)}= H|_{B\cup B_e(\xi,\d)}.
$$ 
Hence,
$$ 
g^n\circ H \circ f^{-n}_*
=g^n\circ H \circ F|_B
=H|_B,
$$
and the proof of Claim~$1^0$ is complete. 
\endpf

\sp Now, continuing the proof of Lemma~\ref{l1j21}, fix $u\in J(f)\sms \ov{{\rm PC}(f)}$ and $r_u>0$ so small that
\beq\label{1_2017_12_09}
B_e(u,2r_u)\sbt U_f^*.
\eeq
Our strategy is to apply Kuratowski--Zorn Lemma. In order to do this consider the 
family $\cF$ of all open subsets $W$ of $\C$
containing $U_f^*$ for which there exists a holomorphic function 
$H_W:W\to\oc$ with the following two properties.

\begin{itemize}
\item[(a)] $H_W|_{U_f^*}=H$.

\, \item[(b)] If $z\in U_f^*$ and $f^n(z)\in W$, then $g^n\circ H(z)=
  H_W\circ f^n(z)$.
\end{itemize}

\sp\fr The family $\cF$ is partially ordered by inclusion and, by 
Claim~$1^0$, it contains $U_f^*$, so $\cF$ is not empty. If $\cC$ is a 
linearly ordered subset of $\cF$, then $H_{W_2}|_{W_1}=\tilde
H_{W_1}$ whenever $W_1,W_2\in\cC$ and $W_1\sbt W_2$. This is so since
$W_1\spt U_f^*$ and (a) holds. Thus putting 
$$
W:=\bu\{G:G\in\cC\}
$$
and defining 
$$
\tilde H_W(z):=H_G(z)
$$ 
if $z\in G\in\cC$, we see that
$H_W:W\to\oc$ is a well defined holomorphic function satisfying
the requirements (a) and (b). So, $W$ is an upper bound of $\cC$. We
therefore conclude from Kuratowski--Zorn Lemma that $\cF$ contains a
maximal element, and we denote it by $G$. We claim that
\beq\label{1pdns10}
G=\C.
\eeq
Indeed, seeking contradiction suppose that $G\ne\C$. Then there exists a point
$w\in\bd G$. Assume that 
\beq\lab{2_2017_12_09}
w\notin \Om(f).
\eeq
Then there exist an integer $s\ge 0$ and a point $\xi\in U_f^*$ such that $f^s(\xi)=w$. Take $R>0$ smaller than $\d>0$, ascribed to the set $X:=\{\xi\}$ and $\varepsilon:=r_u$, according to Theorem~\ref{mnt6.3A}, and also so small that 
\beq\label{2dns148} 
\Comp(\xi, f^s,R)\cap f^{-s}(w)=\{\xi\}, 
\eeq 
\beq\label{1dns148} 
\Comp(\xi, f^s,R)\sbt U_f^*, 
\eeq 
and the map $f^s|_{\Comp(\xi, f^s,R)}:\Comp(\xi, f^s,R)\lra B(w,R)$ has no other critical points except possibly $\xi$. Let $l$ be an arbitrary closed line
segment joining $w$ and $\bd B(w,R)$. Then, there exists
$f_l^{-s}:B(w,R)\sms l\to\C$, a holomorphic branch of $f^{-k}$ such that
\beq\label{5_2017_12_06}
f_l^{-s}\(B(w,R)\sms l\)\sbt \Comp(\xi,f^s,R).
\eeq
Define the holomorphic map $\tilde H_l:B(w,R)\sms l\to\C$ by
$$
\tilde H_l:=g^s\circ H\circ f_l^{-s}.
$$
Because of \eqref{5_2017_12_06}, \eqref{1dns148}, and item (b) applied to $G$, we have that
$$
\tilde H_l|_{G\cap(B(w,R)\sms l)}=H_G|_{G\cap(B(w,R)\sms l)},
$$
Therefore, if $q$ is another closed line segment joining $w$ and
$\bd B(w,R)$, then $\tilde H_l$ and $\tilde H_q$ coincide on the
uncountable set $G\cap(B(w,R)\sms(l\cup q))$. Hence, they glue
together to a single holomorphic map 
$$
\tilde H_w:B(w,R)\sms\{w\}\lra\C.
$$
By virtue of \eqref{2dns148}, $\lim_{z\to 
w}f_l^{-s}(z)=\xi$ and $\lim_{z\to w}f_q^{-s}(z)=\xi$. Therefore,
$$
\lim_{z\to w}\tilde H_w(z)=g^s(H(\xi)).
$$
Consequently $\tilde H_w$ extends holomorphically to a function from
$B(w,R)$ to $\C$. Since $\tilde H_w$ and $H_G$ coincide on
$G\cap B(w,R)$, they glue together to a single holomorphic function
$H_w:G_w\to\C$, where $G_w:=G\cup B(w,R)$. We shall prove that
$$
G_w\in \cF.
$$
Indeed, $G_w$ is an open connected subset of $\C$ containing
$U_f^*$. Moreover, property (a) holds since $G_w\spt G\spt U$ and 
$$
H_w|_{U_f^*}=H_G|_{U_f^*}=H.
$$
We shall show that (b) holds too. In
order to prove it, consider an integer $n\ge 0$ and $C$, a 
connected component of $U_f^*\cap f^{-n}(G_w)$. If $f^n(C)\cap G\ne\es$,
then (b) holds for $G_w$ because it holds for $G$. So, we may assume
without loss of generality that 
$$
f^n(C)\cap G=\es,
$$
in particular 
\beq\label{1dns153}
f^n(C)\sbt B(w,R).
\eeq
Let $\g$ be a compact topological arc in $B(w,R)$ joining $f^n(C)$ and
$G$, and disjoint from $\bu_{j=1}^\infty f^j(\PC(f))$. Fix a point 
$$
y\in f_*^{-n}(\g\cap G).
$$
Since the elliptic map $f:\C\lra\oc$ is, by Proposition~\ref{p120190913}, topologically exact, there exist a point $x\in B(u,r_u)$ and an integer $k\ge 0$ such that $f^k(x)=y$. Let 
$$
V\sbt B(w,R)
$$ 
be an open connected simply connected neighborhood of
$\g$ disjoint from $\bu_{j=0}^{n+k} f^j(\PC(f))$. Let
$f_*^{-n}:V\lra\C$ be a holomorphic inverse branch of $f^n$
defined on $V$ and satisfying: 
\beq\lab{5_2017_12_07}
C\cap f_*^{-n}(V\cap f^n(C))\ne\es.
\eeq
Let
$f_x^{-(n+k)}:V\lra\C$ be a unique holomorphic inverse branch of $f^{n+k}$
defined on $V$ and sending $f^n(y)$ to $x$. Because of Theorem~\ref{mnt6.3A} and the choices of $r_u$ and of $R$, we have that
\beq\lab{1_2017_12_07}
f_x^{-(n+k)}(V)\sbt B(u,2r_u)\sbt U_f^*.
\eeq
The immediate observations are that 

\begin{itemize}
\item[(1)] $f^k\circ f_x^{-(n+k)}=f_*^{-n}|_V$,

\, \item[(2)] $f_x^{-k}:=f_x^{-(n+k)}\circ f^n:f_*^{-n}(V)\lra U_f^*$ \newline 
is the unique holomorphic inverse branch of $f^k$ defined on $f_*^{-n}(V)$ and
sending $y$ to $x$.
\end{itemize}

\,\fr Furthermore,

\begin{itemize}
\item[(3)] $U_f^*\cap f_x^{-(n+k)}(G\cap V)\ne\es$ \  \  \
as  \  \  \  $x\in U_f^*\cap f_x^{-(n+k)}(G\cap V)$
\end{itemize}

\,\fr and

\begin{itemize}
\item[(4)] $U_f^*\cap f_x^{-k}(U_f^*\cap f_*^{-n}(V))\ne\es$  \newline
as, by \eqref{5_2017_12_07},

\, \begin{itemize}
\item[(4a)]
$$
U_f^*\cap f_*^{-n}(V)
\spt C\cap f_*^{-n}(V\cap f^n(C)))\ne\es,
$$
and, by \eqref{1_2017_12_07}, 

\, \item[(4b)] 
$$
f_x^{-k}(U_f^*\cap f_*^{-n}(V))
\sbt f_x^{-k}\circ f_*^{-n}(V)
=f_x^{-(n+k)}(V)\sbt U_f^*.
$$
\end{itemize} 
\end{itemize}

\sp\fr By (3) there exists $D$, an (open) connected component of $G\cap V$, such that $U_f^*\cap f_x^{-(n+k)}(D)\ne\es$. It then foolows from \eqref{1_2017_12_07} and (b) that
$$
H_w|_{D}=g^{n+k}\circ H\circ f_x^{-(n+k)}|_{D},
$$
whence
$$
H_w|_V=g^{n+k}\circ H\circ f_x^{-(n+k)}.
$$
Also, it follows from (4) and (4b) along with (a) and (b) that
$$
H|_{U\cap f_*^{-n}(V)}=g^k\circ H\circ f_x^{-k}|_{U\cap f_*^{-n}(V)}.
$$
Therefore,
$$
\aligned
H_w|_{f^n(U\cap f_*^{-n}(V))}
&=g^{n+k}\circ H\circ f_x^{-(n+k)}|_{f^n(U\cap f_*^{-n}(V))} \\
&=g^n\circ\(g^k\circ H\circ f_x^{-k}|_{U\cap f_*^{-n}(V)}\)
         \circ f_*^{-n}|_{f^n(U\cap f_*^{-n}(V))} \\
&=g^n\circ H\circ f_*^{-n}|_{f^n(U\cap f_*^{-n}(V))} .
\endaligned
$$
Thus, 
\beq\lab{6_2017_12_07}
g^n\circ H|_{U\cap f_*^{-n}(V)}
=H_w\circ f^n|_{U\cap f_*^{-n}(V)}.
\eeq
Since $C\sbt U$ and since $f_*^{-n}(V)\spt f_*^{-n}(V\cap f^n(C))\ne\es$, we conclude from \eqref{5_2017_12_07}, that $U\cap f_*^{-n}(V)\spt f_*^{-n}(V\cap f^n(C))\ne\es$. Hence, we directly get from \eqref{6_2017_12_07}, that
\beq\lab{7_2017_12_07}
g^n\circ H|_{f_*^{-n}(V\cap f^n(C))}=H_w\circ f^n|_{f_*^{-n}(V\cap f^n(C))}
\eeq
Since $C\sbt U$ and because of \eqref{1dns153}, the respective maps $g^n\circ H|_C$  
and $H_w\circ f^n|_C$ are well defined. Therefore, since $f_*^{-n}(V\cap f^n(C))\sbt C$, we conlude from \eqref{7_2017_12_07} that 
$$
g^n\circ H|_C=H_w\circ f^n|_C.
$$
So, taking $\tilde H_{G_w}:=H_w$, we see that $G_w\in\cF$, contrary to maximality of $G$. We thus have that
$$
G\spt \C\sms\Om(f).
$$
Put
$$
H_f:=H_{\C\sms\Om(f)}:\C\sms\Om(f)\lra\C.
$$ 
By symmetry of the situation we also have now a holomorphic function $H_g:\C\sms \Om(g)\lra\C$ which extends $h^{-1}:J_g\lra J_f$ from neighborhood of $J(g)\sms \ov{{\rm PC}(g)}$. But then, invoking also the fact that $H^{-1}(\Om(g))=\Om(f)$ (as $h:J(f)\lra J(g)$ is a topological conjugacy and rationally indifferent periodic points are topologically determined), 
$$
H_g\circ H_f:\C\sms \Om(f)\lra\C
$$ 
is a well defined holomorphic function which is identity on the uncountable set $J(f)\sms \Om(f)$. Therefore,
\beq\label{1dns155}
H_g^{-1}\circ H_f
=\Id_{\C\sms \Om(f)}
\eeq
Hence, the holomorphic map $H_f:\C\sms\Om(f)\lra\C$ is 1--to--1. But since also, $\Om(f)$ and $\Om(g)$ are both finite sets contained in $\C$ and $H(\Om(f))=\Om(g)$, all points of $\Om(f)$ must be removable singularities of $H_f$. Thus, $H_f$ extends holomorphically to $\C$; whence the same also holds for $H_g$. Because of \eqref{1dns155}, we then also have
\beq\label{1dns155B}
H_g^{-1}\circ H_f
=\Id_{\C}.
\eeq
Thus, $H_f:\C\lra\C$ is a holomorphic isomorphism, so of the form
$$
\C\ni z\longmapsto az+b\in\C,
$$
and also
$$
H_f\circ f=g\circ H_f
$$ 
The proof of Lemma~\ref{l1j21} is complete. 
\endpf

\sp The implication {\bf{(6)}}$\Rightarrow${\bf{(1)}} of
Theorem~\ref{t1042706} now
readily follows from Lemma~\ref{l1j15}, Lemma~\ref{l1j19}, and
Lemma~\ref{l1j21}.  The proof of Theorem~\ref{t1042706} is complete.
\endpf

\backmatter

\appendix

\chapter[Appendix]{} 

\section[A Quick Review of Some Selected Facts from Complex Analysis] {A Quick Review of Some Selected Facts from Complex Analysis of One Complex Variable}\label{Review Complex Analysis}

In this appendix we collect for the convenience of the reader many basic and fundamental theorems of complex analysis. We provide no proofs but we give detailed references (arbitrarily chosen) where the proofs can be found. We use these theorems throughout the book without directly referring to them. 

\bthm\label{open} {\rm (Open Mapping Theorem)}\index{(N)}{Open
Mapping Theorem} The image of an open set under a non--constant
analytic map is an open set.
\ethm

For a proof see  e.g. Theorem~7.1  in  \cite{BK}.

\sp\bthm\label{Liouville} {\rm (Liouville
Theorem)}\index{(N)}{Liouville Theorem}  
If $f:\C\lra\C$ is a bounded analytic function, then $f$ is constant. 
\ethm

For a proof see  e.g. Theorem~5.10  in  \cite{BK}.

\sp\bthm\label{maximum} {\rm (Maximum  Modulus
Principle)}\index{(N)}{Maximum  Modulus Principle} Suppose that $D\sbt\C$ is an open connected set. Let $f:\C\lra\C$ be a holomorphic function with in a  continuous extension to the boundary $\partial{D}$. if
$$
M:=  \max\{|f(z)|: z \in \partial{D}\},
$$
then either 
$$
|f(z)| < M
$$ 
for every $z\in D$ or $f$ is a constant function with absolute value equal to $M$.
\ethm

For a proof see  \cite{Ga}.

\sp\bthm\label{Bloch}{\rm (Bloch's Theorem)}\index{(N)}{Bloch's
Theorem} There exists $\b>0$ such that if $f:\mathbb D\lra\C$ is a holomorphic function defined on the open unit disk $\mathbb D$ for which 
$$
f(0)=0
\  \  {\rm and} 
\  \   f'(0)=1,
$$
then there exists an open Euclidean disk $D \subset \mathbb D$ such that the restricted function $f|_D$ is one--to--one and $f(D)$
contains an open Euclidean disk of radius $ \beta$.
\ethm

For a proof see \cite{C}. Such a largest constant $\b>$ is called
Bloch's constant\index{(N)}{Bloch's  constant}.

\sp\bthm\label{Riemann} {\rm (Riemann  Mapping
Theorem)}\index{(N)}{Riemann  Mapping  Theorem} If $D$ is an open connected
simply connected  subset of $\C$ whose boundary (in $\C$) contains
at least two points, then for every $\xi\in D$ there exists a conformal  homeomorphism from the unit disk $\D\sbt\C$ onto $D$ sending $0$ to $\xi$. 
\ethm 

For a proof see \cite{C}. The conformal homeomorphism in the above theorem
is commonly called a Riemann mapping.

\sp\bdfn\label{locally-connected}\index{(N)}{locally  connected} A
topological space $X$ is called locally connected at a point $\xi\in X$  
if and only if for every open set $V$ containing $\xi$ there exists a connected set $\Ga$ such that $\xi\in \Int(\Ga)\sbt\Ga\sbt V$. The space $X$ is said to be locally connected if and only if it is locally connected at $\xi$ for all $\xi \in X$. Equivalently, as it is well known, the space $X$ is locally connected if and only if it has a topological basis consisting of open connected sets.
\edfn

\sp\bthm\label{Caratheodory}{\rm (Carath\'eodory Theorem)}\index{(N)}{
Carath\'eodory  Theorem} If $D$ is an open connected
simply connected subset of $\C$ whose boundary (in $\C$) contains
at least two points, then the boundary $\partial{D}$ is locally connected  if and only if some (equivalently, any) Riemann mapping from $\D$ onto $D$ extends  continuously to the boundary $\partial\D$. 
\ethm

For a proof see  e.g. Theorem~2.1  in  \cite{CG}.

\sp

\bthm\label{reflection}{\rm (Schwarz Reflection
Principle)}\index{(N)}{Schwarz Reflection Principle} Suppose that $G$ is an open connected subset of $\C$ such that 
$$
G^*:=\{z\in\C: \ov{z} \in G\}=G
$$
Let
$$
G_+:=\{z\in G: \im(z)>0\},
\  \
G_-:=\{z\in G: \im(z)<0\}
\  \  {\rm and,} \  \
G_0:=\{z\in G: \im(z)=0\}.
$$
If $f:G_+\cup G_0\lra\C$ is a continuous function such that $f|_{G_+}$ is holomorphic, and if
$$
f(G_0)\sbt \R,
$$
then there exists a holomorphic function $g:G\lra \C$ such that
$$
g|_{G_+\cup G_0}=f.
$$
\ethm

For proof see  e.g. Theorem~7.8 in \cite{BK} or \cite{C}.

\sp

\bthm\label{Picard}{\rm(Great Picard Theorem)}\index{(N)}{Great
Picard Theorem} 
Suppose that $G$ is an open subset of $\C$ and $\xi$ is an isolated point of the boundary of $G$. If $f:G\lra \C$ is a holomorphic function that has an essential singularity at $\xi$, then for every open neighborhood $V$ of $\xi$ in $\C$, the set $\C\sms f(V\cap G)$ contains at most one point.
\ethm

For a proof see  \cite{C}.

\sp\bdfn 
Suppose that $G$ is an open subset of $\C$. Let $\mathcal F$ be a family of meromorphic functions from $G$ to $\oc$. The family  $\mathfrak F$ is called
normal \index{(N)}{normal family} if and only if every sequence $\(f_n\)_{n=1}^\infty$ in $ \mathfrak F$ contains a subsequence that converges uniformly in the spherical metric on compact subsets  of $G$.
\edfn

\sp \bthm\label{Montel}{\rm (Montel's Criterion of
Normality)}\index{(N)}{Montel's Criterion of Normality} 
Suppose that $G$ is an open subset of $\C$. Let
$\mathfrak F$ be a family of meromorphic functions on $D$.
If there are three values in $\oc$ that are omitted by every every function$f\in \mathfrak F$, i.e. if 
$$
\oc\sms\bu_{f\in \mathfrak F}f(G)
$$
consists of at least three points, then $ \mathfrak F$ is a normal family. 
\ethm

For a proof see e.g. Theorem~3.2 in  \cite{CG}.

\sp

\bthm\label{Uniformization} {\rm (Koebe-Poincar\'e's  Uniformization
Theorem)}\index{(N)}{Koebe-Poincar\'e's  Uniformization Theorem} The following statements hold:

\ben
\item Every connected simply connected Riemann surface is conformally equivalent to either the open unit disc $\mathbb D$, the complex plane ${\mathbb C}$, or the Riemann sphere $\oc$. 

\,

\item Any Riemann surface has one of these three, $\D$, $\C$, or $\oc$, as its universal covering surface. 

\,

\item The only Riemann surface having the Riemann sphere $\oc$ as its
universal covering surface is the sphere $\oc$ itself. 

\,

\item 
The only surface having the complex plane $\mathbb C$ as universal covering surface are the complex plane itself $\mathbb C$, the punctured complex
plane $\mathbb C \sms \{0\}$, and tori $\mT$.

\,

\item  All others Riemann  surface  have the open unit disc $\mathbb  D$  as universal covering surface. 
\een
\ethm

For proof see \cite{Ga}

\sp

\bthm\label{3.8.6jones} If $(f_n)_{n=1}^\infty$ is a sequence
of holomorphic functions defined on a an open connected set $D\sbt \mathbb C$ such that the product
$$
f:=\prod_{n=1}^\infty f_n 
$$ 
converges normally on all compact subsets of $D$, then the function  
$f:D\lra\C$ is holomorphic.
\ethm

For a proof see e.g. Theorem 3.8.6 in \cite{JS}.

\sp

\bthm\label{3.8.7jones}
Let $D$, $f$, and $(f_n)_{n=1}^\infty$ be as in Theorem~\ref{3.8.6jones}, and let $z \in D$. Then $f(z)=0$ if and only if $f_n(z)=0$ for some $ n\in \mathbb{N}$. In this case the set
$$
Z:=\{n\in \N: f_n(z)=0\}
$$
is finite and the order of $z$ as zero of $f$ is the sum of the orders of $z$ as zeros of the functions $f_n$, where $n\in Z$.
\ethm

For a proof see e.g. Theorem 3.8.7 in \cite{JS}.

\sp

\bthm\label{3.8.8jones}
If $D$, $f$, and $(f_n)_{n=1}^\infty$ are as in Theorem~\ref{3.8.6jones},
then  
$$
\sum_{n=1}^\infty \frac{f_n'}{f_n}
$$  
converges uniformly to $f'/f$ on all compact subsets of $D$.
\ethm

For a proof see e.g. Theorem 3.8.8 in \cite{JS}.

\sp

\bthm\label{monodromy} {\rm (Monodromy Theorem)}\index{(N)}
{Monodromy Theorem} With the terminology of \cite{JS} let $D$ be an open connected simply connected subset of
$\oc$, and let $(G,f)$ be a function element with $G \sbt
D$. If $(G,f)$ can be continued meromorphically along all paths in $D$ starting at some point $a \in G$, then there is a direct  meromorphic
continuation $(G,f) \sim (D,f)$.
\ethm

For proof see e.g. Theorem 4.5.3 in \cite{JS}.

\sp

\bdfn\label{beltrami} A Beltrami coefficient on an open set $U\sbt
\oc$  is a measurable function $\mu: U \lra  \mathbb C$
such that $|\mu|$ has essential supremum  $k <1$. This means that
the set  $\{z \in U; \, \, |\mu(z)| >k\}$  has zero Lebesgue measure
and $k$ is the smallest number with this property. 

If $f: U \lra f(U)$ is a quasiconformal  homeomorphism, then the function $\mu_f: U \lra  \mathbb C$, given by the formula
$$
\mu_f(z):=\frac{\ov{\partial}{f(z)}}{\partial{f}(z)},
$$ 
is called the Beltrami coefficient of $f$. 
\edfn

\sp

\bthm\label{ABB}{\rm (Ahlfors--Bers--Bojarski
Theorem)}\index{(N)}{Ahlfors-Bers-Bojarski Theorem} Let $\La \sbt
\mathbb C^n$ be an open set and $ \mu: \oc\times \La\lra
\mathbb D$ be a measurable function  satisfying:
\begin{itemize}
\item [(a)] There exists $k<1$ such that $| \mu(z, \l)|\leq k$ for all $ \l  \in \La$ and for Lebesgue almost all $z \in \oc$,

\,

\item [(b)] The function $\La\ni \l \lmt \mu(z,\l)$ is a holomorphic for Lebesgue almost all $z \in \oc$.
\end{itemize}
Then there  exists a unique function $F: \oc \times \La
\lra \oc$ such that

\,

\begin{enumerate}
\item  $F(0,\l)=0$, $F(1, \l)=1$, $F(\infty, \l)=\infty$,

\,

\item For every $\l \in \La$ the map $ \oc\ni z \lmt  F(z,\l)\in\oc$ is a
quasiconformal homeomorphism whose Beltrami coefficient is $
\mu(\cdot, \l)$,

\,

\item The map $\La\ni \l \lmt F(z,\l)\in\oc$ is holomorphic for Lebesgue almost every $z\in\oc$.
\end{enumerate}
\ethm

For a proof see \cite{AB}, \cite{Boj}.

\chapter[Appendix]{}

\section{Proof of Sullivan's Non--Wandering Theorem for Speiser class $\cS$}

As the title of this section indicates we shall prove in it the Non--Wandering Theorem listed in the main body of the book as Theorem~\ref{Sullivan}. This theorem has been conjectured by P. Fatou in \cite{F0}
and \cite{F1} for rational functions and was proved in this form by D. Sullivan in \cite{Su0}. In fact it holds for all meromorphic functions in Speiser class $\cS$ and was proved by I.N. Baker, J. Kotus, and Y. L\"u in \cite{BKL4}. 
See this paper and references therein for more historical and bibliographical information. Examples of wandering components for analytic self--maps of $\mathbb C^*$ were provided by J. Kotus in \cite{K2} and by I.N.~Baker, J.~Kotus and Y.L\"u  in \cite{BKL2}. We provide  below the appropriate statement, i.e Theorem~\ref{Sullivan}, and its original proof copied from \cite{BKL4} without any substantial processing it. We start with the following auxiliary fact.

\blem\label{l1in{bkl4}} If a meromorphic function $f:\C\lra\oc$ is in  Speiser class $\cS$ and its Fatou set $F(f)$ has some
wandering component, then there exists $U$, a wandering connected component of $F(f)$, which is simply connected, for each integer $k\ge 0$ the set $f^k(U)$ is a connected component of $F(f)$, and the restriction
$$
f|_U^k:U\lra \C
$$
is a 1-to-1 map.
\elem

\fr{\sl Proof.} Let $U_0$ be a wandering component of $F(f)$ and for every integer $n\ge 0$ let $U_n$ be the connected component of $F(f)$ which contains $f^n(U_0)$. Since the set $\Sing(f^{-1})$ is finite, replacing $U_0$ by $U_k$ with some sufficiently large integer $k\ge 0$, we may assume that
\beq\label{620190905}
U_n\cap \Sing(f^{-1})=\es
\eeq
for all integers $n \geq 0$. Hence , the restriction
$$
f|_{U_n}:U_n \lra U_{n+1}
$$
is a (regular) covering map; in particular surjective. It thus suffices to show that $U_0$ (and, since it satisfied the same assumptions, each set $U_n$, for all $n\geq 0$) is simply connected.

Seeking contradiction suppose that $U_0$ is not simply connected. This means that $U_0$ contains the image of an analytic Jordan curve $\g:S^1\to U_0$ which is not null--homotopic (not homotopic to a constant) in
$U_0$. If for some integer $n\ge 0$ the closed path
$$
\g_n:=f^n\circ\g:S^1\lra U_n
$$
is null homotopic in $U_n$, then the corresponding homotopy lifts, by the covering map $f|_{U_0}^n:U_0\to U_n$, to a homotopy between $\g$ and a constant loop in $U_0$. So, $\g$ is null homotopic in
$U_0$, and this contradiction shows that $\g_n$ is not null homotopic in $U_n$ for any integer $n\ge 0$.

For each integer $n\ge 0$ the path $\g_n$ may have self--intersections, but its image is a finite union of analytic arcs, which are all disjoint apart from common end--points; each such arc may be though traversed more than once. Abusing slightly notation for the ease of exposition, we will in the sequel denote the image $\g_n(S^1)$ just by $\g_n$; this will not cause any confusion. Let
\beq\label{120190905}
\varepsilon:=\min\big\{\dist_s(a,b): a, b\in \Sing(f^{-1}) \  {\rm and} \ a\ne b\big\}\in (0,1)
\eeq
By Proposition~\ref{p1_2017_07-26} the sequence $\(f|_{U_0}^n\)_{n=0}^\infty$ has only  constant limit functions. Hence, the spherical
diameters $\diam_s(\g_n\)$ converge to $0$ as $n \to \infty$. Therefore, there exists an integer $N\ge 0$ such that
\beq\label{220190905}
\diam_s(\g_n\) < \varepsilon
\eeq
for all integers $n \geq N$. For every $n \geq 0$, the set
$\oc\sms \g_n$ has finitely many connected components exactly one of which, called unbounded, contains $\infty$. Denote it by $D_\infty(n)$. The remaining components of $\oc\sms \g_n$ will be  called bounded. The union of $\g_n$ and all bounded components of $\oc\sms \g_n$ is equal to
\beq\label{420190905}
\tilde{\g}_n:= \oc \sms D_\infty(n).
\eeq
Note that
\beq\label{320190905}
\diam_s(\^\g_n\)=\diam_s(\g_n\).
\eeq
We claim that
$$
f(\tilde{\g}_n)\sbt \tilde{\g}_{n+1}.
$$
Fix $n\geq N$ and let $D$ be a bounded component of $\oc\sms \g_n$. Then $\partial{D}$ is a piecewise analytic Jordan curve contained in $\g_n$, whence
\beq\label{520190905}
f(\bd D)\sbt f(\g_n)=\g_{n+1}.
\eeq
Because of \eqref{120190905}, \eqref{220190905}, and \eqref{320190905}, the intersection
$$
\tilde{\g}_{n+1}\cap \Sing(f^{-1})
$$
contains at most one element. Suppose first that it does contain one element and denote it by $a$. Since $\g_{n+1}\sbt U_{n+1}$, it follows from
\eqref{620190905} and \eqref{520190905} that
$$
a\notin f(\bd D)\sbt\g_{n+1}\sbt U_{n+1}.
$$
Now, the loop $f(\bd D)$, treated as part of the loop $\g_{n+1}$, is homotopic to $k\d$ in the open connected set $B(a,2
\varepsilon)\sms\{a\}$, with some $k \in \mathbb Z$, where $\d$ is the loop given by the formula
$$
[0,2\pi]\ni t\lmt \d(t):=a+\eta e^{it}
\in B(a,\varepsilon)\sms\{a\}
\sbt B(a,2
\varepsilon)\sms\{a\}
$$
with some/any $0<\eta<\varepsilon$. If $\xi \in
\bd D$, then the branch $g$ of $f^{-1}$ with $g(f(\xi))=\xi$ can be continued throughout $B(a,2 \varepsilon)\sms\{a\}$ and lifts the homotopy $f(\bd D)\sim k \d$ to a homotopy of $\bd D$ to a closed curve $\beta$ which projects to $k \d$ by $f$. Furthermore, since the branch $g$ of the meromorphic function $f$ has an isolated singularity at $a$, this point is either a logarithmic branch point or an algebraic branch point.

 If $k \neq 0$, $g$ has a an algebraic singularity at $a$, since
 $g(z)$ returns  to its initial value after $z$ makes $k$ circuits
 of $\d$. In this  case we may lift the homotopy
 $a+(1-s)(f(\bd D)-a)$, $ 0 \leq s \le 1 $. Hence, with 
$$
\Delta:=B(a, 2\varepsilon),
$$ 
we have proved the following. 

\sp Claim~1: $\bd D$ is homotopic in $A$ to a constant path $\a$, where $A$ is the component of $f^{-1}(\Delta)$ containing $\xi$.
Thus $f$ is holomorphic in $A$ and  
$$
\diam_s(f(A)) \leq  4 \varepsilon.
$$
If $k=0$, the same Claim~1 holds with $\Delta :=B(a,2 \varepsilon)\sms\{a\}$.  If $\tilde{\g}_{n+1}$ does not intersect $\Sing(f^{-1})$, Claim~1 holds with $\Delta$ being any sufficiently small simply--connected neighborhood of
   $\tilde{\g}_{n+1}$. Thus Claim~1 holds in all cases.

\sp We now will show that $D \sbt A$. Indeed, if $z \notin A$, then the winding number of $\bd D$ satisfies $n(\partial{D},z)= n
(\a,z)=0$. Hence $z \notin D$ Thus $D\sbt A$, whence
$$
\diam_s(f(D) \leq \diam_s(f(A)\leq 4 \varepsilon.
$$
Since $f:\C\lra\oc$ is meromorphic, we have
$$
\partial{f(D)}\sbt  f(\partial{D})\sbt \tilde{\g}_{n+1}.
$$
The set $f(D)$ cannot intersect $D_\infty(n+1)$, for otherwise $f(D)$ would include $D_\infty(n+1)$ and would had spherical diameter exceeding $4\varepsilon$. Thus,
$$
f(D)\sbt \tilde{\g}_{n+1}.
$$
Since this holds for every bounded component $D$, we conclude that
$$
f(\tilde{\g}_{n}) \sbt \tilde{\g}_{n+1}.
$$
Therefore by induction,
$$
f^j(\tilde{\g}_n)\sbt \tilde{\g}_{n+j}
$$
for all integers $n \geq N$ and $j\geq 0$. But this, along with \eqref{220190905} and \eqref{320190905}, implies that the family
$$
\lt(f|_{\Int_{\C}(\tilde{\g}_n)}^j\rt)_{j=0}^\infty
$$
is normal. Hence,
$$
\Int_{\C}(\tilde{\g}_n)\sbt F(f).
$$
Since, in addition, $\bd \tilde{\g}_n\sbt \g_n\sbt U_n\sbt F(f)$, we thus get that
$$
\tilde{\g}_n\sbt F(f).
$$
Since also, $U_n$ is a connected of $F(f)$, we thus conclude that
$$
\tilde{\g}_n\sbt U_n.
$$
But $\tilde{\g}_n$ intersects the Julia set $J(f)$ since $\g_n $ is not  contractible to a point in $U_n$. This contradiction shows that $U_0$ is indeed simply--connected and completes the proof of Lemma~\ref{l1in{bkl4}}. \endpf

\sp We will also need the following.

\blem\lab{Baker Lemma 7.9}
There is a constant $K>1$ such that if $H:\oc\lra\oc$ is a $K$--quasiconformal homeomorphism, which fixes the three points $0$, $1$, and $\infty$, then
$$
|H(z)-z|<\frac12
$$
for every $z=e^{i\th}$ with $\lt[\frac13\pi\le \th\le \frac53\pi\rt]$.
\elem

\bpf
For every $K\ge 1$ the family $\cF_K$ of all quasiconformal homeomorphisms fixing the three points $0$, $1$, and $\infty$ is normal and if $K\to 1$, then the maps of $\cF_K$ converge to the identity map uniformly on $\lt\{e^{i\th}:\lt[\frac13\pi\le \th\le \frac53\pi\rt]\rt\}$. So, the proof is complete.
\epf

\sp As an immediate consequence of this lemma, we get the following. 

\blem\lab{Baker Lemma 7.9B}
If $\e>0$ is small enough, then for any two quasiconformal homeomorphisms $f,g:\oc\lra\oc$, which both fix $0$, $1$, and $\infty$, with respective dilatations $\mu_f$ and $\mu_g$ satisfying $\|\mu_f\|_\infty, \|\mu_g\|_\infty<\e$, one has
$$
|g\circ f(z)-z|<\frac12
$$
for every $z=e^{i\th}$ with $\lt[\frac13\pi\le \th\le \frac53\pi\rt]$.
\elem

\sp The main result of this appendix is the following.

\bthm
No functions in Speiser class $\mathcal S$ have wandering domains. 
\ethm

\bpf
The proof is long and involved. We split it into five parts.

\sp {\bf I}. Suppose for a contrary that for some meromorphic function in $f\in
\mathcal S$ its Fatou set $F(f)$ has a wandering component $U$, which due to
Lemma~\ref{l1in{bkl4}} can be assumed to be the one produced in this lemma.  Let 
$$
p:=\#\Sing(f^{-1})\in\N.
$$ 
Take an integer $M> 2p+8$ and define
$$
T:=\big\{t \in \mathbb R^M: |t_i|<\d \ {\rm for \ all} \ 1\leq i \leq M \big\},
$$
where $0<\d<\pi/2$ is a sufficiently small constant whose required smallness will be determined later in the course of the proof. Take $ \a_i$ so that
$$
\d < \a_1<\ldots  <\a_M< 2\pi - \d
$$
and also all the sets
$$
I_j:=[\a_i-\d,\a_i+\d]
$$
are mutually disjoint. Define $\psi: [0,2 \pi]\lra\R$ the following formula:
$$
\psi_j(\theta)
:=\left\{\begin{array}{ll}
    \d \exp\lt(-\frac{\d^2}{(\theta-\a_j)^2-\d^2}\rt) \ & {\rm if} \ \theta \in I_j,\\
     0 \ & {\rm if} \  \theta \in  [0, 2 \pi ] \sms I_j.
\end{array}\right.
$$
Then $\psi_j $ is $C^\infty$ and
$$
\|\psi_j\|_\infty\leq \d/2
\  \  \  {\rm and} \  \  \
\|\psi_j'\|_\infty\leq M',
$$
with some constant $M'$ independent of $j$ and $\d$. Set
$$
\psi(t, \theta)
:=\theta+\sum_{j=1}^ Mt_j\psi_j(\theta),\,\,\, t \in T.
$$
Then with 
$$
\psi_\theta:=\frac{\bd \psi}{\bd \th}
$$
we have that
$$
1-MM'\d \leq \psi_\theta\leq 1+ MM'\d.
$$
Assume $\d>0$ to be so small that
$$
MM'\d<1.
$$
Then the function $\psi$ is increasing in $\theta$. For every $t\in T$ define the map 
$$
\phi_{t}:\ov{\D}\lra \ov{\D}
$$
by the following formula:
$$
\phi_{t}(z)= r\exp(i \psi( t, \theta)), 
$$
where $z=re^{i\th}$ with $\th\in[0,2\pi]$. Note that 
$$
\phi_{t}(\D)=\D
$$
and moreover, $\phi_{t}$ is a quasiconformal homeomorphism self of the unit disk $\D$ onto itself with dilatation
$$
\nu_{t}(z)=- \frac{e^{2i\theta}\sum_{j=1}^Mt_j\psi_j'(\theta)}{2+\sum_{j=1}^Mt_j\psi_j'(\theta)}.
$$
Thus,
$$
\|\nu_{t}\|_\infty\leq\frac{ M M'\d }{2-  MM'\d}\leq  M M'\d< \varepsilon,
$$
where $\varepsilon>0$ comes from Lemma~\ref{Baker Lemma 7.9B} and we assumed that
$$
\d\in\lt(0,\frac{\varepsilon}{3 MM'}\rt).
$$
Note that for each $t\in T$, 
$$
\phi_{t}(0)=0.
$$ 
Note that $\psi_j (\theta)=0$ for each $j=1,2,\ld,M$ and every $\theta \in [0,2\pi]\sms \bigcup_{i=1}^M I_i$. Hence, $\psi(t, \theta)=\theta$ for all $t\in T$ and all $\theta \in [0,2\pi]\sms \bigcup_{i=1}^M I_i$. Therefore 
$$
\phi_{t}(e^{i\theta})=
\exp(i\psi(t, \theta))=\e^{i\theta}
$$
for all such $t$ and $\th$. Since $[0,2\pi]\sms \bigcup_{i=1}^M I_i$ contains a non--empty open interval, denote it by $S$, we thus conclude that 
$$
\phi_{t}|_{S}=\Id_S
$$
for all $t\in T$. Furthermore, 

\centerline{$\phi_t|_{\D}=\phi_s|_{\D}$ implies that $t=s$. }

\sp {\bf II.} Since $U$ is simply connected, because of the Riemann Mapping Theorem,  there exists a conformal homeomorphism
$$
h:\D \lra U
$$
Then the homeomorphism
$$
\chi:=h\circ\phi_{t}\circ h^{-1}: U \lra U
$$
 is quasiconformal with dilatation
$$
\mu_t=\nu_t\circ h^{-1} \frac{\ov{(h^{-1})'}}{(h^{-1})'},
$$
which is measurable and satisfies
\beq\lab{120200410}
\|\mu_t\|_\infty\leq \varepsilon.
\eeq
With
$$
P:=\bigcup_{n=0}^\infty f^{-n}(\infty)
$$
we say that two points $x$ and $y$ in $\mathbb C \sms P$ are equivalent if and only if
$$
f^m(x)=f^n(y)
$$
for some $m,n \in\mathbb N$. By Lemma~\ref{l1in{bkl4}} each equivalence class $[x]$ of this relation intersects $U$ in at most one point. So, we can extend $\mu_{t}$ to $\mathbb C$ by the following formula:
$$
\mu_{t}(z):=
\begin{cases}
\frac{\mu_{t}(z)(f^m)'(x)\ov{(f^n)'}(z)}{(f^m)'(x)\ov{(f^n)'}(z)} &
 {\rm if}\,\,
  f^n(z)=f^m(x),\, {\rm with \ some} \,\, x\in U \, {\rm and} \  m,n \in \mathbb N,\\
 0 & {\rm if} \,\,[z]\cap U=\es\,\, {\rm or} \,\, z \in P. 
 \end{cases}
$$
Then $\mu_{t}$ is measurable with 
\beq\lab{120200414}
\|\mu_t\|_\infty< \varepsilon
\eeq
and
\begin{equation}\label{3in{bkl4}}
\mu_t(f(z))=\mu_t(z)\frac{f'(z)}{\ov{f'(z)}}
\,\,\, \  {\rm Lebesgue} \ a.e.
\end{equation}
Fix different finite $b_1, b_2\notin\Sing(f^{-1})$. By Ahlfors--Bers--Bojarski Theorem (see \cite{Boj}) there is a unique quasiconformal   homeomorphism
\beq\label{120190904}
\Phi_{t}:\oc\lra\oc
\eeq
which fixes $b_1$, $b_2$, and $\infty$, and has dilatation $\mu_{t}$. Furthermore (see e.g. \cite{Boj}), for each $z\in\C$ the function
$$
T\ni {t} \lmt \Phi_{t}(z)\in\C
$$
belongs to $C^1(T)$. It also follows from (\ref{3in{bkl4}}) that
\beq\label{220190904}
f_{t}=\Phi_{t}\circ f\circ \Phi_{t}^{-1}:\oc\lra\oc
\eeq
is meromorphic in $\mathbb C$ for each ${t}\in T$. Moreover, the restriction
$$
f_{t}: \mathbb C \sms\Phi_{t}(f^{-1}(\Sing(f^{-1})))\lra  \mathbb C \sms
\Phi_{t}(\Sing(f^{-1}))
$$
is an (unbranched) covering map.

\sp {\bf III}. Since $M> 2p+4$, there exists a homeomorphic embedding
$$
\om:I=[0,\sigma_0]\lra T  
$$
such that the following restriction of the $C^1$--map from the formula \eqref{120190904},
$$
\om(I)\ni t \lmt \Phi_{t}\Big|_{\Sing(f^{-1})\cup \{a_1,a_2\}}
$$
is constant. For  every $\sigma \in I$ abbreviate
$$
\Phi_\sigma = \Phi_{\om(\sigma)} \  \  {\rm and} \  \  f_\sigma=
f_{\om(\sigma)}.
$$
We shall prove the following.

\blem\label{l2in{bkl4}}
For all $\sg, \tau\in I=[0,\sg_0]$ we have that
$$
f_\sigma=f_\tau.
$$
\elem

\bpf Fix $a_1\in f^{-1}(b_1)$, $a_2\in f^{-1}(b_2)$. We already know that both maps
$$\begin{aligned}
&f_{\sg}\circ \Phi_\sigma^{-1}: \mathbb C \sms
\Phi_{\sigma}(f^{-1}(\Sing(f^{-1})))\lra \ov{ \mathbb C} \sms \Sing(f^{-1})\\
{\rm and}\\
&f_{\tau}\circ  \Phi_\tau^{-1}: \mathbb C \sms
\Phi_{\tau}(f^{-1}(\Sing(f^{-1})))\lra \ov{ \mathbb C} \sms \Sing(f^{-1})
\end{aligned}$$
are unbranched coverings, while for every $\lambda \in I$, the composition
$\Phi_\lambda^{-1}\circ \Phi_\sigma$ is a quasiconformal self--map of the set
$\oc\sms \Sing(f^{-1})$ which lifts to a homeomorphism
$$
\psi_{\sigma \lambda}: \mathbb C \sms \Phi_{\sigma}(f^{-1}(\Sing(f^{-1})))\lra
\mathbb C\sms \Phi_\tau(f^{-1}(\Sing(f^{-1})))
$$
such that $\psi_{\sigma
\lambda}(\Phi_\sigma(a_1))=\Phi_\sg(a_1)=\Phi_\tau(a_1)$ since
$\Phi_\lambda^{-1}\circ\Phi_\sigma(b_1)=b_1$. We have a commuting
diagram:
\[\begin{tikzcd}
 \mathbb C \sms  \Phi_{\sigma}(f^{-1}(\Sing(f^{-1})))\arrow{r}{f\circ\Phi^{-1}_\sigma } \arrow[swap]{d}{\psi_{\sigma \lambda}} & \oc\sms \Sing(f^{-1})\arrow{d}{\Phi_{\lambda}^{-1}\circ\phi_\sigma} \\
{ \mathbb C \sms  \Phi_{\sigma}(f^{-1}(\Sing(f^{-1})))} \arrow{r}{f\circ\Phi^{-1}_\tau} & {\oc\sms \Sing(f^{-1})}
\end{tikzcd}
\]
Since $\Phi_{\lambda}(a_2)= \Phi_\sigma(a_2)$ and
$\Phi_\lambda(b_2)=\Phi_\sigma(b_2)$, we  have $\psi_{\sigma\lambda}\circ
\Phi_{\sigma}(a_{p+2}) \in \Phi_{\tau} (f^{-1}(b_2))$, which is a
discrete set. Therefore, since the function
$$
I\ni\l\lmt\psi_{\tau \lambda}(\Phi_{\sigma}(a_{p+2}))\in\C
$$
is continuous, it is constant. Now, $\psi_{\sigma
\sigma}=\Phi_{\tau}\circ\Phi_{\sigma}^{-1}$, whence
$$
\psi_{\sigma \lambda}(
\Phi_{\sigma}(a_{p+2}))= \psi_{\sigma \sigma}(
\Phi_{\sigma}(a_{p+2}))=\Phi_{\tau}(a_{p+2})=\Phi_{\sigma}(a_{p+2}).
$$
Taking $\lambda=\tau$ the commutative diagram above gives
$$
f_\sigma=\Phi_\sigma \circ f\circ \Phi_{\sigma}^{-1}=\Phi_\tau\circ f\circ
\Phi_{\tau}^{-1}\circ\psi_{\sigma \tau}=f_\tau\circ\psi_{\sigma \tau}.
$$
Thus, $\psi$ is analytic except, possibly, at isolated points, hence in fact it is a conformal automorphism of $\oc$. It  must clearly fix $\infty$ and also $ \Phi_{\sigma}(a_1), \Phi_\sigma(a_2)$, so
it is the identity. The proof of Lemma~\ref{l2in{bkl4}} is complete.
\epf

\sp {\bf IV.} For every $\sigma \in I$ set
$$
\Om_\sigma:=\Phi_0^{-1}\circ \Phi_\sigma
$$
Lemma~\ref{l2in{bkl4}} is then equivalent to saying that
$$
\Om_\sigma\circ f = f\circ\Om_\sigma.
$$
Given $\sg\in I$, for each $q \in \mathbb N$ the map $\Om_{\sigma}$ permutes the periodic points of period $q$ of $f$. Since in addition
$\Om_0$ is an identity map and since for every $z\in\C$, the map $I\ni\sg\lmt \Om_\sigma(z)\in\C$ is continuous, we conclude that $\Om_\sg$
fixes all of them. Hence, $\Om_\sg$ fixes all periodic points of $f$. Since the set of all such points is dense in the Julia set $J(f)$ of $f$, we thus conclude that
$$
\Om_\sg|_{J(f)}=\Id_{J(f)}.
$$
In particular, 
\beq\lab{520200413}
\Om_\sg|_{\bd U}=\Id_{\bd U}.
\eeq
We also have 
$$
\Om_\sg(U)\sbt U.
$$
Hence, with the map $h:\D \lra U$ defined in Part II of our proof, we get that the following composition is well defined.
$$
G_\sg:=h^{-1}\circ\Om_{\sigma}\circ h:\D\lra\D
$$
is a quasiconformal homeomorhism. Now we are in position to formulate the following lemma due to I.N. Baker (see \cite{Ba4}). We also provide its proof. 

\sp\blem~\label{l3in{bkl4}}
For every $\sg\in I$, the quasiconformal homeomorhism $G_\sg:\D\lra\D$ has a continuous extension to $\ov{\D}$ which is identity on $\partial{\D}$.
\elem

\bpf Let $\rho_\D$ and $\rho_U$ be the hyperbolic (Poincar\'e) metrics respectively on the unit disk $\D=\{z\in\C:|z|<1\}$ and the set $U$. Since Poincar\'e metrics are conformally invariant, using the definition of $G_\sg$, we have for every $z\in \D$ that
\beq\lab{620200413}
\rho_\D(z, G_\sg(z))
=\rho_U(h(z),h(G_\sg(z)))
=\rho_U\(h(z),\Omega_\sigma(h(z))\).
\eeq
Denote by $\alpha $ a nearest point of $\partial{U}$ to $h(z)$ and set 
$$
r:=|\alpha-h(z)|.
$$
Assuming that $z$ lies sufficiently close to $\bd\D$, we see that there exists  a point $\beta\in\partial{U}$ whose distance from $\alpha$ is $r$. Since \eqref{120200414} holds, since $\Omega_\sigma$ fixes the three points $\alpha$, $\beta$, and $\infty$, (see \eqref{520200413}), and since the angle of the triangle $h(z)$, $\a$, $\b$, at the vertex $\a$ is at least $\pi/3$, Lemma~\ref{Baker Lemma 7.9B}, after an appropriate translation, rotation, and rescaling by $1/r$, with $h(z)$ corresponding to $z$, $\a$ corresponding to $0$, and  $\b$ corresponding to $1$, gives that
$$
|\Omega_\sigma(h(z))- h(z)| < \frac{1}{2}r.
$$
Since also
\beq\lab{720200413}
\rho_U\(h(z),\Omega_\sigma(h(z))\)
\leq  \rho_W\(h(z), \Omega_\sigma(h(z))\), 
\eeq
where $W$ is the disk $B_\C(h(z),r)$, we conclude that $\rho_U\(h(z),\Omega_\sigma(h(z))\)$ is bounded above by an absolute constant, say $M$, independent of $z\in U$ sufficiently close to $\bd U$. Hence, invoking \eqref{620200413} and \eqref{720200413}, we get that
$$
\rho_\D(z, G_\sg(z))\le M
$$
for every $z\in\D$ sufficiently close to $\bd\D$. Therefore, 
$$
\lim_{\D\ni z\to \xi}G_\sg(z)=\xi
$$
for every $\xi\in\bd\D$ and the convergence is uniform with respect to $\xi\in\D$. The proof is complete. 
\epf

\sp {\bf V.} We now can conclude the proof of Theorem~\ref{Sullivan}.
By Lemma~\ref{l3in{bkl4}} 
$$
\Omega_\sigma=h\circ G_\sigma \circ h^{-1}: U\lra U
$$ 
has continuous extension to $\ov{U}$ which is identity on $\partial{U}$ for $\sigma \in I$.
Since $\Omega_\sigma=\Phi^{-1}_0\circ \Phi_\sigma$, we have that 
$$ 
\Phi_\sigma= \Phi_0\circ \Omega_\sigma
$$
for all $\sigma \in I$.
Therefore, all the images $\Phi_\sigma(U)$, $\sigma \in I$, are
they same. Denote their common value by $V$.  
Let $H:\D\lra V$ be a conformal homeomorhism resulting from the Riemann Mapping Theorem. Then Lemma~\ref{l3in{bkl4}} shows that
$$
h^{-1}\circ\Phi^{-1}_{\om(0)}\circ H \circ H^{-1}\circ \Phi_{\om(\sigma)}\circ h,
$$ 
restricted $\partial{\D}$, is the identity map. Thus
$$
H^{-1}\circ\Phi_{\om(\sigma)}\circ h=H^{-1}\circ\Phi_{\sg(0)}\circ h
$$ 
holds on $\partial{\D}$ and both sides are maps of $\D$ to $
\D$ with respective dilatation $\mu_{\om(\sigma)}$ and $\mu_{\om(0)}$. These delatations are the same on $\D$ as the dilatations of $\phi_{\om(\sigma)}$ and $\phi_{\om(0)}$ are. Hence, there exists  M\"obius maps $L_\sigma,L_0 : \D \lra \D$ such that
$$
H^{-1}\circ\Phi_{\om(\sigma)}\circ h
=L_\sigma\circ \Phi_{\om(\sigma)} 
\quad \text{and}\quad 
H^{-1}\circ\Phi_{\om(0)}\circ h=
L_0\circ \Phi_{\om(0)}.
$$
By Lemma~\ref{l3in{bkl4}} 
$$
(L_0\circ\phi_{\sg(0)})^{-1}\circ (L_\sigma\circ\phi_{\om(\sigma)})=\Id
$$ 
on $\partial{D}$. Thus,
$$
L_0\phi_{\om(0)}=L_\sigma\phi_{\om(\sigma)}
$$ 
on $\partial{D}$. Since $\phi_{\om(\sg)}=\Id$ on a non--empty open subset of $\partial{D}$, we have that $L_0=L_\sigma$ which implies $\om(0)=\om(\sigma)$, against the construction. The existence of wandering $U$ thus led to a
contradiction. The proof of Theorem~\ref{Sullivan} is complete.
\epf

\backmatter


\bibliographystyle{amsalpha}


\Printindex{(S)}{Symbols} \Printindex{(N)}{Index}

\end{document}